\documentclass[12pt]{book}

\usepackage{amssymb,amsmath,amscd, amsthm,stmaryrd,verbatim, mathrsfs,latexsym, makeidx}
\makeindex

\begin{document}

\baselineskip=18pt

\newcommand{\la}{\langle}
\newcommand{\ra}{\rangle}
\newcommand{\psp}{\vspace{0.4cm}}
\newcommand{\pse}{\vspace{0.2cm}}
\newcommand{\ptl}{\partial}
\newcommand{\dlt}{\delta}
\newcommand{\sgm}{\sigma}
\newcommand{\al}{\alpha}
\newcommand{\be}{\beta}
\newcommand{\G}{\Gamma}
\newcommand{\gm}{\gamma}
\newcommand{\vs}{\varsigma}
\newcommand{\Lmd}{\Lambda}
\newcommand{\lmd}{\lambda}
\newcommand{\td}{\tilde}
\newcommand{\vf}{\varphi}
\newcommand{\yt}{Y^{\nu}}
\newcommand{\wt}{\mbox{wt}\:}
\newcommand{\der}{\mbox{Der}\:}
\newcommand{\ad}{\mbox{ad}\:}
\newcommand{\stl}{\stackrel}
\newcommand{\ol}{\overline}
\newcommand{\ul}{\underline}
\newcommand{\es}{\epsilon}
\newcommand{\dmd}{\diamond}
\newcommand{\clt}{\clubsuit}
\newcommand{\vt}{\vartheta}
\newcommand{\ves}{\varepsilon}
\newcommand{\dg}{\dagger}
\newcommand{\tr}{\mbox{Tr}\:}
\newcommand{\ga}{{\cal G}({\cal A})}
\newcommand{\hga}{\hat{\cal G}({\cal A})}
\newcommand{\Edo}{\mbox{End}\:}
\newcommand{\for}{\mbox{for}}
\newcommand{\kn}{\mbox{ker}\:}
\newcommand{\Dlt}{\Delta}
\newcommand{\rad}{\mbox{Rad}}
\newcommand{\rta}{\rightarrow}
\newcommand{\mbb}{\mathbb}
\newcommand{\rd}{\mbox{Res}\:}
\newcommand{\stc}{\stackrel{\circ}}
\newcommand{\stdt}{\stackrel{\bullet}}
\newcommand{\lra}{\Longrightarrow}
\newcommand{\f}{\varphi}
\newcommand{\rw}{\rightarrow}
\newcommand{\op}{\oplus}
\newcommand{\co}{\Omega}
\newcommand{\dar}{\Longleftrightarrow}
\newcommand{\sta}{\theta}
\newcommand{\wht}{\widehat}
\newcommand{\msr}{\mathscr}
\newcommand{\mfk}{\mathfrak}

\title{\Large \bf Representations of Lie Algebras and \\ Partial Differential Equations}
\vspace{6in}

\author{{\large \bf  Xiaoping Xu}\\ \\{\large Institute of Mathematics}
\\{\large Academy of Mathematics and System Sciences}\\{\large Chinese Academy of Sciences}
\\{\large 55 Zhongguancun Dong Lu}
\\{\large Beijing 100190, P. R. China}}

\vspace{2in}

\date{ 2016}

\thispagestyle{empty} \maketitle

\pagenumbering{roman}

\chapter*{Preface}
\addcontentsline{toc}{chapter}{\numberline{}Preface}
\markboth{PREFACE}{PREFACE}

Symmetry is an important phenomenon in the natural world. Lie
algebra is not purely an abstract mathematics but a fundamental tool
of studying symmetry. In fact, Norwegian mathematician Sorphus Lie
introduced Lie group and Lie algebra in 1874 in order to study the
symmetry of differential equations. Lie algebras are the
infinitesimal structures (bones) of Lie groups, which are symmetric
manifolds. Lie theory has extensive and important applications in
many other fields of mathematics, such as geometry, topology, number
theory, control theory, integrable systems, operator theory, and
stochastic process etc. Representations of finite-dimensional
semisimple Lie algebras play fundamental roles in quantum mechanics.
 The controllability property of the unitary propagator of an $N$-level
quantum mechanical system subject to a single control field can be
described in terms of the structure theory of semisimple Lie
algebras. Moreover, Lie algebras were used to explain the
degeneracies encountered in  genetic codes as the result of a
sequence of symmetry breakings that have occurred during its
evolution. The structures  and representations of simple Lie
algebras are connected with  solvable quantum many-body system in
one-dimension. Our research also showed that the representation
theory of Lie algebras is connected with algebraic coding theory.

The existing classical books on finite-dimensional Lie algebras,
such as the ones by Jacobson and by Humphreys, purely focus on the
algebraic structures of semisimple Lie algebras and their
finite-dimensional representations. Explicit irreducible
representations of simple Lie algebras had not been addressed
extensively. Moreover, the relations of Lie algebras with the other
subjects had not been narrated that much. It seems to us that a book
on Lie algebras with more extensive view is needed in coupling with
modern development of mathematics, sciences and technology.

 This book is mainly an
exposition of the author's works and his joint works with his former
students on explicit representations of finite-dimensional simple
Lie algebras, related partial differential equations, linear
orthogonal algebraic codes, combinatorics and algebraic varieties.
Various oscillator generalizations of the classical representation
theorem on harmonic polynomials are presented. New functors from the
representation category of a simple Lie algebra to that of another
simple Lie algebra are given. Partial differential equations play
key roles in solving certain representation problems. The weight
matrices of the minimal and adjoint representations over the simple
Lie algebras of types E and F are proved to generate ternary
orthogonal linear codes with large minimal distances. New
multi-variable hypergeometric functions related to the root systems
of simple Lie algebras are introduced in connection with  quantum
many-body system in one-dimension. Certain equivalent combinatorial
properties on representation formulas are found. Irreducibility of
representations are proved directly related to algebraic varieties.

This book is self-contained  with the minimal prerequisite of
calculus and linear algebra. It consists of three parts. The first
part is mainly the classical structure and finite-dimensional
representation theory of finite-dimensional semisimple Lie algebras,
 written with Humphreys' book ``Introduction to Lie Algebras and
Representation Theory" as the main reference, where we give more
examples and direct constructions of simple Lie algebras of
exceptional types, revise some arguments and prove some new
statements. This part serves as the preparation for the main
context. The second part is our explicit representation theory of
finite-dimensional simple Lie algebras. Many of the irreducible
representations in this part are infinite-dimensional and some of
them are even not of highest-weight type. Certain important natural
representation problems are solved by means of solving partial
differential equations. In particular, some of our irreducible
presentations are completely characterized by invariant partial
differential equations. New representation functors are constructed
from inhomogeneous oscillator representations of simple Lie
algebras, which give fractional representations of the corresponding
Lie groups, such as the projective representations of special linear
Lie groups and the conformal representations of complex orthogonal
groups. The third part is an extension of the second part. First we
give supersymmetric generalizations of the classical representation
theorem on harmonic polynomials. They can also be viewed as certain
supersymmetric Howe dualities. Then we present our theory of Lie
theoretic codes. Finally, we talk about root related integrable
systems and our new multi-variable hypergeometric functions, which
are natural multi-variable analogues of classical Gauss
hypergeometric functions. The corresponding hypergeometric partial
differential equations are found.

Part of this book has been taught for many times at the University
of Chinese Academy of Sciences. The book can serve as a research
reference book for mathematicians and scientists. It can also be
treated as a text book for students  after a proper selection of
materials.

\vspace{1cm}

\hfill Xiaoping Xu\\
Beijing, P. R. China\\
2016

\newpage

\chapter*{Introduction}
\addcontentsline{toc}{chapter}{\numberline{}Introduction}
\markboth{INTRODUCTION}{INTRODUCTION}

Algebraic study of partial differential equations traces back to
Norwegian mathematician Sophus Lie [Lie], who invented the powerful
tool of continuous groups (known  as Lie groups) in 1874 in order to
study symmetry of differential equations. Lie's idea has been
carried on mainly by the mathematicians in the former states of
Soviet Union, East Europe and some mathematicians in North America.
Now it has become an important mathematical field known as ``group
analysis of differential equations," whose main objective is to find
symmetry group of  differential equations, related conservation laws
and similarity solutions. One of the classical long-standing
problems in the area was to determine all invariant partial
differential equation under the natural representations of classical
groups (It was also mentioned by Olver [Op]). The problem was
eventually settled down theoretically by the author [X5], where we
found the complete set of functional generators for the differential
invariants of classical groups.

Gel'fand, Dikii and Dorfman [GDi1, GDi2, GDo1-GDo3] introduced in
1970s a theory of Hamiltonian operators in order to study the
integrability of nonlinear evolution partial differential equations.
Our first experience with partial differential equation was in the
works [X3, X4, X9, X10, X11] on the structure of Hamiltonian
operators and their supersymmetric generalizations. In particular,
we [X10] proved that Lie conformal algebras in the sense of Kac
[Kv2] are equivalent to linear Hamiltonian operators as mathematical
structures. Since Borcherds' vertex algebras are determined by their
positive parts (modes in physics terminology) of vertex operators
that form Lie conformal algebras (e.g., cf. [Kv2, X6]), our result
essentially established an equivalence between  vertex algebras and
Hamiltonian operators. The result was later generalized by Barakat,
De Sole and Kac [BDK] to vertex poisson algebras. Moreover, we used
the techniques and thoughts from Lie algebras to solve various
physical partial differential equations such as the equation of
transonic gas flows, the equation
 of geopotential  forecast, the nonlinear Schr\"{o}dinger equations,
the Navier-Stokes equations and the classical boundary layer
equations (cf. [X21] and the references therein). In this book, we
present an explicit representation theory of finite-dimensional
simple Lie algebras and show that partial differential equations can
be used to solve representation problems of Lie algebras.

 Abstractly, a Lie algebra ${\msr
 G}$ is a vector space with a bilinear map $[\cdot,\cdot]:{\msr G}\times{\msr G}\rta{\msr
 G}$ such that
$$[u,v]=-[v,u],\;\;[u,[v,w]]=[[u,v],w]+[v,[u,w]],\qquad
u,v,w\in{\msr G}.\eqno(0.1)$$ A Lie algebra ${\msr G}$ is called
{\it simple}\index{simple Lie algebra} if it does not contain any
nonzero ``invariant" proper subspace ${\cal I}$; that is£º
$$[u,{\cal I}]\subset {\cal I},\qquad u\in {\msr G}.\eqno(0.2)$$
Finite-dimensional complex simple Lie algebras were classified by
Killing and Cartan in later nineteen century. For a vector space
$V$, we denote by $V^\ast$ the space of linear functions on $V$. It
turns out that such a Lie algebra ${\msr G}$  must contain a
subspace $H$, called {\it toral Cartan subalgebra},\index{toral
Cartan subalgebra} such that
$${\msr G}=\bigoplus_{\al\in H^\ast}{\msr G}_\al,\qquad{\msr
G}_\al=\{u\in{\msr G}\mid [h,u]=\al(h)u\;\for\;h\in H\},\eqno(0.3)$$
and ${\msr G}_0=H$, where
$$\Phi=\{\al\in H^\ast\setminus\{0\}\mid {\msr G}_\al\neq\{0\}\}\eqno(0.4)$$
is called the {\it root system of ${\msr G}$}\index{root system}.
Moreover,
$$\dim {\msr G}_\al=1\qquad\for\;\;\al\in\Phi\eqno(0.5)$$ and
there exists a subset $\Phi_+$ of $\Phi$ containing a basis
$\Pi=\{\al_1,...,\al_n\}$ of $H^\ast$ such that
$$\Phi=-\Phi_+\bigcup\Phi_+,\;\;-\Phi_+\bigcap\Phi_+=\emptyset\eqno(0.6)$$
and for any $\be\in\Phi_+$,
$$\be=\sum_{i=1}^nk_i\al_i\;\;\mbox{with}\;\;0\leq
k_i\in\mbb{Z},\eqno(0.7)$$where $\mbb{Z}$ is the ring of integers.
Furthermore, there exists $h_i\in [{\msr G}_{\al_i},{\msr
G}_{-\al_i}]$ such that
$$\al_i(h_i)=2\qquad\for\;\;i=1,2,...,n.\eqno(0.8)$$
It turns out that the structure of ${\msr G}$ is completely
determined by the matrix $(\al_i(h_j))_{n\times n}$, which is called
the {\it Cartan matrix}\index{Cartan matrix} of ${\msr G}$.
Consequently, there are only nine classes of finite-dimensional
simple Lie algebras, which are called of types $A_n, B_n, C_n, D_n,
G_2,F_4,E_6,E_7,E_8$. The first four are infinite series, called
``classical Lie algebras." The last five are fixed simple Lie
algebras, called ``exceptional types." Using the root lattices
$\Lmd_r=\sum_{\al\in\Pi}\mbb{Z}\al$ of type $A$ and the Coxeter
automorphism of $\Lmd_r$, we constructed in [X1] two families of
self-dual complex lattices that are also real self-dual (unimodular)
lattices. These lattices can be used in studying geometry of numbers
and conformal field theory.

Since every finite-dimensional complex simple Lie algebra must
contain a toral Cartan subalgebra, it is natural to ask if there
exists a complex simple Lie algebra that does not contain any toral
Cartan subalgebra.  In [X7, X8], we constructed six families of
infinite-dimensional complex simple Lie algebras without any toral
Cartan subalgebras.

For a vector space $M$, we denote by $\Edo M$ the space of all
linear transformations on $M$. A {\it representation} $\nu$ of a Lie
algebra ${\msr G}$ on $M$ is a linear map from ${\msr G}$ to $\Edo
M$ such that
$$\nu([u,v])=\nu(u)\nu(v)-\nu(v)\nu(u)\qquad\for\;\;u,v\in {\msr
G}.\eqno(0.9)$$ The space $M$ is called a ${\msr G}$-{\it
module}.\index{module} In this book, we sometimes use the notions:
$$\nu(\xi)=\xi|_M,\;\;\nu(\xi)(w)=\xi(w)\qquad\for\;\;\xi\in{\msr G},\;w\in
M.\eqno(0.10)$$ A subspace $\msr N$ of a  ${\msr G}$-module $M$ is
called a {\it submodule}\index{submodule} of $M$ if
$$\xi(w)\in \msr N\qquad\for\;\;\xi\in{\msr G},\;w\in \msr N.\eqno(0.11)$$
If $M$ dose not contain any proper nonzero submodule, we say that
$M$ is an {\it irreducible}\index{irreducible} ${\msr G}$-module.
Let ${\msr G}$ be a Lie algebra with a decomposition (0.3)
satisfying (0.4), (0.6) and (0.7). A ${\msr G}$-module $M$ is called
a {\it weight module}\index{weight module} if
$$M=\bigoplus_{\mu\in H^\ast}M_\mu,\;\;M_\mu=\{w\in M\mid
h(w)=\mu(h)w\}.\eqno(0.12)$$ The set
$$\Lmd(M)=\{\mu\in H^\ast\mid M_\mu\neq \{0\}\}\eqno(0.13)$$
is called the {\it weight set}\index{weight set} of $M$ and the
elements in $\Lmd(M)$ are called the {\it weights of
$M$}.\index{weight} A vector in $M_\mu$ with $\mu\in\Lmd(M)$ is
called a {\it weight vector}.\index{weight vector} A nonzero weight
vector $w$ is called a {\it singular vector}\index{singular vector}
if
$$\xi(w)=0\qquad \for\;\;\xi\in \bigcup_{\be\in\Phi_+}{\msr
G}_\be.\eqno(0.14)$$ If $M$ is generated by a singular vector $v$
with weight $\lmd$, we call $M$ a {\it highest-weight ${\msr
G}$-module}\index{highest-weight module} and $\lmd$ the {\it highest
weight} of $M$. \index{highest weight}

Finite-dimensional representations of a complex finite-dimensional
simple Lie algebra were essentially determined by Cartan in early
twenty century. The approaches were simplified and further developed
by Weyl in 1920s. Let ${\msr G}$ be a finite-dimensional complex
simple Lie algebra and take the settings in (0.3)-(0.8). Set
$$\Lmd^+=\{\lmd\in H^\ast\mid 0\leq \lmd
(h_i)\in\mbb{Z}\;\for\;i=1,...,n\}.\eqno(0.15)$$ It turns out that
any finite-dimensional ${\msr G}$-module  is a
 direct sum of its irreducible submodules and any finite-dimensional irreducible
  ${\msr G}$-module is a highest-weight module with its
  highest weight in $\Lmd^+$. Conversely, for any element $\lmd\in\Lmd^+$,
  there exists a unique finite-dimensional irreducible ${\msr G}$-module with highest weight $\lmd$.
The above conclusion is now called {\it Weyl's theorem on completely
reducibility}. \index{Weyl's theorem on complete reducibility} For
any finite-dimensional ${\msr G}$-module $M$,
$$\mu(h_i)\in\mbb{Z}\qquad\for\;\;i\in\{1,...,n\},\;\mu\in\Lmd(M).\eqno(0.16)$$
Write $\Lmd(M)\setminus\{0\}=\{\mu_1,\mu_2,...,\mu_m\}$. We found in
[X20] that the set
$$\{(\mu_1(h_i),\mu_2(h_i),...,\mu_m(h_i))\mid
i=1,...,n\}\eqno(0.17)$$ spans a binary or ternary orthogonal code
for certain ${\msr G}$ and $M$. In particular, we showed that when
${\msr G}$ is of type $F_4,\;E_6,\; E_7$ and $E_8$, and $M$ is of
minimal dimension or $M={\msr G}$ with left multiplication as the
representation, the set (0.17) spans a ternary orthogonal code with
large minimal distance $d$, which can be used to correct $\llbracket
d/2 \rrbracket$ errors in information technology.

A module of a finite-dimensional simple Lie algebra is called {\it
cuspidal}\index{cuspital} if it is not induced from its proper
``parabolic subalgebras." Infinite-dimensional irreducible weight
modules of finite-dimensional simple Lie algebras with
finite-dimensional weight subspaces had been intensively studied by
the authors in [BBL, BFL, BHL, BL1, BL2, Fs, Fv, Mo]. In particular,
Fernando [Fs] proved that such modules must be cuspidal or
parabolically induced. Moreover, such cuspidal modules exist only
for special linear Lie algebras and symplectic Lie algebras. A
similar result was independently obtained by Futorny [Fv]. Mathieu
[Mo] proved that these cuspidal modules
 are irreducible components in the tensor
modules of their multiplicity-free modules with finite-dimensional
modules. Thus the structures of irreducible weight modules of
finite-dimensional simple Lie algebra with finite-dimensional weight
subspaces were essentially determined by Fernando's result in [Fs]
and Methieu's  result in [Mo]. However, explicit presentations of
such modules were almost unknown.

An important feature of this book is the connection between
representations of simple Lie algebras and partial differential
equations. Let $E_{r,s}$ be the square matrix with 1 as its
$(r,s)$-entry and 0 as the others. Denote by $\mbb{R}$ the field of
real numbers.  Let $n\geq 3$ be an integer. The compact orthogonal
Lie algebra $o(n,\mbb{R})=\sum_{1\leq r<s\leq
n}\mbb{R}(E_{r,s}-E_{s,r}),$ whose representation  on the polynomial
algebra ${\msr A}=\mbb{R}[x_1,...,x_n]$ is given by
$(E_{r,s}-E_{s,r})|_{\msr A}=x_r\ptl_{x_s}-x_s\ptl_{x_r}$. Denote by
${\msr A}_k$ the subspace of homogeneous polynomials in ${\msr A}$
with degree $k$. Recall that the {\it Laplace
operator}\index{Laplace operator}
$\Dlt=\ptl_{x_1}^2+\cdots+\ptl_{x_n}^2$ and its corresponding
invariant $\eta=x_1^2+x_2^2+\cdots+x_n^2$. It is well known that the
subspaces
$${\msr H}_k=\{f\in{\msr A}_k\mid
\Dlt(f)=0\}\eqno(0.18)$$of {\it harmonic polynomials}\index{harmonic
polynomial} form  irreducible $o(n,\mbb{R})$-submodules and
$${\msr A}=\bigoplus_{i,k=0}^\infty \eta^i{\msr H}_k\eqno(0.19)$$
is a direct sum of irreducible submodules. In other words, the
irreducible submodules are characterized by the Laplace operator
$\Dlt$ and its dual invariant $\eta$ gives the complete reducibility
of the polynomial algebra ${\msr A}$. The above conclusion is called
the {\it classical theorem on harmonic polynomials}.\index{classical
theorem on harmonic polynomials}

The {\it Navier equations}\index{Navier equations}
$$\iota_1\Delta (\vec u)+(\iota_1+\iota_2)(\nabla^t\cdot\nabla )(\vec
u)=0\eqno(0.20)$$ are used to describe the deformation of a
homogeneous, isotropic and linear elastic medium in the absence of
body forces, where $\vec u$ is an $n$-dimensional vector-valued
function, $\nabla=(\ptl_{x_1},\ptl_{x_2},...,\ptl_{x_n})$ is the
gradient operator,  $\iota_1$ and $\iota_2$ are Lam\'{e} constants
 with  $\iota_1>0$,  $2\iota_1+\iota_2>0$ and $\iota_1+\iota_2\neq0$.
In fact, $\nabla^t\cdot\nabla$ is the well-known Hessian operator.
 Mathematically, the
  above system is a natural vector  generalization of
  the Laplace equation. In [X16], we
found methods of solving linear flag partial differential equations
for polynomial solutions. Cao [Cb] used a method of us to prove that
the subspaces of homogeneous polynomial-vector solutions are exactly
direct sums of three explicitly given irreducible $o(n,\mbb
R)$-submodules if $n\neq 4$ and of four explicitly given irreducible
$o(n,\mbb R)$-submodules if $n=4$. This gave a vector generalization
of the classical theorem on harmonic polynomials.

Denote by ${\msr B}=\mbb{C}[x_1,...,x_n,y_1,...,y_n]$ the polynomial
algebra in $x_1,...,x_n$ and $y_1,...,y_n$ over the field $\mbb C$
of complex numbers. Fix
 $n_1,n_2\in\{1,...,n\}$ with $n_1\leq n_2$. We redefine
 $$\mbox{deg}\:x_i=-1,\;\mbox{deg}\:x_j=1,\;\mbox{deg}\:y_r=1,\;\mbox{deg}\:y_s=-1\eqno(0.21)$$
for $i\in\{1,...,n_1\},\;j\in\{n_1+1,..,n\},\;r\in\{1,...,n_2\}$ and
$s\in\{n_2+1,...,n\}$. Denote by ${\msr B}_{\la\ell_1,\ell_2\ra}$
the subspace of homogeneous polynomials with degree $\ell_1$ in
$\{x_1,...,x_n\}$ and degree $\ell_2$ in $\{y_1,...,y_n\}$. We
deform the Laplace operator $\sum_{i=1}^n\ptl_{x_i}\ptl_{y_i}$ to
the operator
$$\td\Dlt=-\sum_{i=1}^{n_1}x_i\ptl_{y_i}+\sum_{r=n_1+1}^{n_2}\ptl_{x_r}\ptl_{y_r}-\sum_{s=n_2+1}^n
y_s\ptl_{x_s}\eqno(0.22)$$ and its dual $\sum_{i=1}^nx_iy_i$ to the
operator
$$\eta=\sum_{i=1}^{n_1}y_i\ptl_{x_i}+\sum_{r=n_1+1}^{n_2}x_ry_r+\sum_{s=n_2+1}^n
x_s\ptl_{y_s}.\eqno(0.23)$$
 Define
$${\msr H}_{\la\ell_1,\ell_2\ra}=\{f\in {\msr B}_{\la
\ell_1,\ell_2\ra}\mid \td\Dlt(f)=0\}.\eqno(0.24)$$  Luo and the
author [LX1] constructed a new representation on ${\msr B}$ for the
simple Lie algebra ${\msr G}$ of type $A_{n-1}$ such that the
operators in (0.22) and (0.23) commute with the elements in ${\msr
G}|_{\msr B}$ and used a  method in [X16] to prove that ${\msr
H}_{\la\ell_1,\ell_2\ra}$ with
 $\ell_1,\ell_2\in\mbb{Z}$ such that
$\ell_1+\ell_2\leq n_1-n_2+1-\dlt_{n_1,n_2}$ are irreducible ${\msr
G}$-modules.  Moreover, ${\msr
B}_{\la\ell_1,\ell_2\ra}=\bigoplus_{m=0}^\infty\eta^m({\msr
H}_{\la\ell_1-m,\ell_2-m\ra})$ is a decomposition of irreducible
${\msr G}$-submodules. This establishes a $\mbb{Z}^2$-graded
analogue of the classical theorem on harmonic polynomials for the
simple Lie algebra of type $A_{n-1}$.

Set
$${\msr B}_{\la k\ra}=\sum_{\ell\in\mbb{Z}}{\msr B}_{\la
\ell,k-\ell\ra},\;\;{\msr H}_{\la k\ra}=\sum_{\ell\in\mbb{Z}}{\msr
H}_{\la \ell,k-\ell\ra}\qquad\for\;\;k\in\mbb{Z}.\eqno(0.25)$$ Luo
and the author [LX2] extended the above representation of ${\msr G}$
on ${\msr B}$ to a representation of the simple Lie algebra ${\msr
G}_1$ of type $D_n$ on ${\msr B}$ and proved that ${\msr H}_{\la
k\ra}$ with $n_1-n_2+1-\dlt_{n_1,n_2}\geq k\in\mbb{Z}$ are
irreducible ${\msr G}_1$-modules and ${\msr B}_{\la
k\ra}=\bigoplus_{i=0}^\infty\eta^i({\msr H}_{\la k-2i\ra})$ is a
decomposition of irreducible ${\msr G}_1$-submodules. This gives a
$\mbb{Z}$-graded analogue of the classical theorem on harmonic
polynomials for the simple Lie algebra of type $D_n$.

 Let ${\msr
B}'=\mbb{C}[x_0,x_1,...,x_n,y_1,...,y_n]$ be the polynomial algebra
in $x_0,x_1,...,x_n$ and $y_1,...,y_n$. We define
$\mbox{deg}\:x_0=1$ and take (0.21). Write the deformed Laplace
operator
$$\td\Dlt'=\ptl_{x_0}^2-2\sum_{i=1}^{n_1}x_i\ptl_{y_i}+2\sum_{r=n_1+1}^{n_2}\ptl_{x_r}\ptl_{y_r}-2\sum_{s=n_2+1}^n
y_s\ptl_{x_s}\eqno(0.26)$$ and its dual operator
$$\eta'=x_0^2+2\sum_{i=1}^{n_1}y_i\ptl_{x_i}+2\sum_{r=n_1+1}^{n_2}x_ry_r+2\sum_{s=n_2+1}^n
x_s\ptl_{y_s}.\eqno(0.27)$$ Set
$${\msr B}'_{\la k\ra}=\sum_{i=0}^\infty {\msr B}_{\la
k-i\ra}x_0^i,\qquad {\msr H}'_{\la k\ra}=\{f\in {\msr B}'_{\la
k\ra}\mid \td\Dlt'(f)=0\}.\eqno(0.28)$$ Luo and the author [LX2]
extend the above representation of ${\msr G}_1$ on ${\msr B}$ to a
representation of the simple Lie algebra ${\msr G}_2$ of type $B_n$
on ${\msr B}'$ and proved that ${\msr H}'_{\la k\ra}$ with
$k\in\mbb{Z}$ are irreducible ${\msr G}_2$-modules and ${\msr
B}'=\bigoplus_{k\in\mbb{Z}}\bigoplus_{i=0}^\infty(\eta')^i({\msr
H}'_{\la k\ra})$ is a decomposition of irreducible ${\msr
G}_2$-submodules. This is  a $\mbb{Z}$-graded analogue of the
classical theorem on harmonic polynomials for the simple Lie algebra
of type $B_n$.

When $n_1<n_2$, the bases of the subspaces ${\msr
H}_{\la\ell_1,\ell_2\ra},\;{\msr H}_{\la k\ra}$ and ${\msr H}'_{\la
k\ra}$ had been obtained. Using Fourier transformation, we can
identify the subspaces ${\msr H}_{\la\ell_1,\ell_2\ra},\;{\msr
H}_{\la k\ra}$ and ${\msr H}'_{\la k\ra}$ with the subspaces of
homogeneous solutions for the corresponding usual Laplace equations
in certain spaces of generalized functions. In [LX3], Luo and the
author extended the above representation of ${\msr G}$ on ${\msr B}$
to a representation of the simple Lie algebra ${\msr G}_3$ of type
$C_n$ on ${\msr B}$ and proved that if $n_1<n_2$ or $k\neq 0$, the
subspace ${\msr B}_{\la k\ra}$ with $k\in\mbb{Z}$ are irreducible
${\msr G}_3$-modules, and when $n_1=n_2$, the subspace ${\msr
B}_{\la 0\ra}$ is a direct sum of two
 explicitly given irreducible ${\msr G}_3$-submodules.

 The above irreducible modules except $n_1=n_2=n$ in the case of
 type $A_{n-1}$ are explicit infinite-dimensional weight modules of
 finite-dimensional weight subspaces. They are not unitary. The
 irreducible modules for ${\msr G}_1,{\msr G}_2$ and ${\msr G}_3$
 are generically neither of highest-weight type nor of unitary type.
 In [LX4],
 Luo and the author generalized the above results for ${\msr
 G},{\msr G}_1$ and ${\msr G}_2$ to those for the corresponding Lie
 superalgebras. Bai [B1, B2] proved that
some of irreducible representations in [LX1-LX3] have distinguished
Gelfand-Kirillov dimensions.

Dickson [D] (1901) first realized that there exists an
$E_6$-invariant trilinear form on its 27-dimensional basic
irreducible module, which corresponds to the unique fundamental
cubic polynomial invariant and constant-coefficient differential
operator. We proved in [X17] that the space of homogeneous
polynomial solutions with degree $m$ for the  cubic Dickson
invariant differential equation is exactly a direct sum of
$\llbracket m/2 \rrbracket+1$ explicitly determined irreducible
$E_6$-submodules and the whole polynomial algebra is a free module
over the polynomial algebra in the Dickson invariant generated by
these solutions. Thus we obtained a cubic $E_6$-generalization of
the classical theorem on harmonic polynomials.

  Next we want to describe our new multi-variable hypergeometric
functions of type A and their connections with the representations
of type-A simple Lie algebras and the Calogero-Sutherland model. For
$a\in\mbb{C}$ and $0<n\in\mbb{Z}$, we denote
$$(a)_n=a(a+1)\cdots (a+n-1),\qquad (a)_0=1.\eqno(0.29)$$
The well-known Gauss hypergeometric function is
$$_2F_1(\al,\be;\gm;z)=\sum_{n=0}^\infty
\frac{(\al)_n(\be)_n}{n!(\gm)_n}z^n.\eqno(0.30)$$ \index{Gauss
hypergeometric function} Many well known elementary functions and
orthogonal polynomials are special cases of $_2F_1(\al,\be;\gm;z)$
(e.g., cf. [WG, AAR, X21]). It satisfies the classical
hypergeometric equation
$$z(1-z){y'}'+[\gm-(\al+\be+1)z]y'-\al\be y=0.\eqno(0.31)$$
\index{classical hypergeometric equation} Denoting
$D=z\frac{d}{dz}$, we can rewrite the above equation as
$$(\gm+D)\frac{d}{dz}(y)=(\al+D)(\be+D)(y).\eqno(0.32)$$
Moreover, it has the differential property:
$$\frac{d}{dz}\:
_2F_1(\al,\be;\gm;z)=\frac{\al\be}{\gm}\:
_2F_1(\al+1,\be+1;\gm+1;z)\eqno(0.33)$$ Furthermore,  the Euler
integral representation
$$
_2F_1(\al,\be;\gm;z)=\frac{\G(\gm)}{\G(\be)\G(\gm-\be)}\int_0^1t^{\be
-1}(1-t)^{\gm-\be-1}(1-zt)^{-\al}dt\eqno(0.34)$$ holds in the $z$
plane cut along the real axis from 1 to $\infty$.

 Let $\mbb N$ be the addtive semigroup of nonnegative integers and
 let
$$\G_A=\sum_{1\leq p<q\leq n}\mbb{N}\es_{q,p}\eqno(0.35)$$
be the additive semigroup  of  rank $n(n-1)/2$ with $\es_{q,p}$ as
base elements.
 For $\al=\sum_{1\leq p<q\leq n}\al_{q,p}\es_{q,p}\in \G_A$, we denote
$$\al_{1\ast}=\al_n^\ast=0,\;\;\al_{k\ast}=\sum_{r=1}^{k-1}\al_{k,r},\;\;
\al_l^\ast=\sum_{s=l+1}^n\al_{s,l}\eqno(0.36)$$
 Given $\vt\in\mbb{C}\setminus\{-k\mid 0<k\in\mbb Z\}$ and $\tau_r\in\mbb{C}$ with $r\in\{1,...,n\}$,
 we defined $(n(n-1)/2)$-variable {\it hypergeometric function of type
 A} in [X14]
 by \index{hypergeometric function of type A}
$${\msr X}_A(\tau_1,..,\tau_n;\vt)\{z_{j,k}\}=\sum_{\be\in\G_A}\frac{\left[\prod_{s=1}^{n-1}
(\tau_s-\be_{s\ast})_{_{\be_s^\ast}}\right](\tau_n)_{_{\be_{n\ast}}}}
{\be!(\vt)_{_{\be_{n\ast}}}}z^\be,\eqno(0.37)$$ where
$$\be!=\prod_{1\leq k<j\leq n}\be_{j,k}!,\qquad z^\be=\prod_{1\leq k<j\leq n}
z_{j,k}^{\be_{j,k}}.\eqno(0.38)$$ The variables $\{z_{j,k}\}$
correspond to the negative root vectors of the simple Lie algebra of
type $A_{n-1}$ and the above definition has essentially used the
structure of its root system.

 Our functions ${\msr
X}_A(\tau_1,..,\tau_n;\vt)\{z_{j,k}\}$ are indeed natural
generalizations of the Gauss hypergeometric function
$_2F_1(\al,\be;\gm;z)$. Denote
$${\msr D}_{p\ast}=\sum_{r=1}^{p-1}z_{p,r}\ptl_{z_{p,r}},\;\; {\msr D}_{q}^\ast=
\sum_{s=q+1}^nz_{s,q}\ptl_{z_{s,q}}\;\;\for\;\;p\in\{2,...,n\},\;q\in\{1,...,n-1\}.\eqno(0.39)$$
We have the following analogous system of partial differential
equations of (0.32):
$$(\tau_{r_2}-1-{\msr D}_{r_2\ast}+{\msr D}_{r_2}^\ast)\ptl_{z_{r_2,r_1}}({\msr X}_A)
=(\tau_{r_2}-1-{\msr D}_{r_2\ast})(\tau_{r_1}-{\msr
D}_{r_1\ast}+{\msr D}_{r_1}^\ast) ({\msr X}_A)\eqno(0.49)$$ for
$1\leq r_1<r_2\leq n-1$, and
$$(\vt+{\msr D}_{n\ast})\ptl_{z_{n,r}}({\msr X}_A)=(\tau_n+{\msr D}_{n\ast})
(\tau_r-{\msr D}_{r\ast}+{\msr D}_r^\ast)({\msr X}_A)\eqno(0.41)$$
for $r\in\{1,...,n-1\}$. Define
$$\left(\begin{array}{ccccc}1&0&0&\cdots&0\\
P_{[1,2]}&
1&0&\cdots&0\\ P_{[1,3]}&P_{[2,3]}&1&\ddots&\vdots\\
\vdots&\vdots&\ddots&\ddots&0\\
P_{[1,n]}&P_{[2,n]}&\cdots&P_{[n-1,n]}&
1\end{array}\right)=\left(\begin{array}{ccccc}1&0&0&\cdots&0\\
z_{2,1}&
1&0&\cdots&0\\ z_{3,1}&z_{3,2}&1&\ddots&\vdots\\
\vdots&\vdots&\ddots&\ddots&0\\ z_{n,1}&z_{n,2}&\cdots&z_{n,n-1}&
1\end{array}\right)^{-1}\eqno(0.42)$$ and treat $P_{[i,i]}=1$. Then
$P_{[i,j]}$ are polynomials in $\{z_{r,s}\mid 1\leq s<r\leq n\}$.
 Moreover, for $1\leq
r_1<r_2\leq n-1$ and $r\in \{1,...,n-1\}$, we have the differential
property:
$$\ptl_{z_{r_2,r_1}}({\msr X}_A)=\sum_{s=1}^{r_1}\tau_sP_{[s,r_1]}{\msr X}_A[s,r_2],\eqno(0.43)$$
$$\ptl_{z_{n,r}}({\msr X}_A)=\frac{\tau_n}{\vt}\sum_{s=1}^r\tau_sP_{[s,r]}{\msr
X}_A[s,n], \eqno(0.44)$$ where ${\msr X}_A[i,j]$ is obtained from
${\msr X}_A$ by changing $\tau_i$ to $\tau_i+1$ and $\tau_j$ to
$\tau_j-1$ for $1\leq i<j\leq n-1$ and ${\msr X}_A[k,n]$ is obtained
from ${\msr X}_A$ by changing $\tau_i$ to $\tau_i+1$, $\tau_n$ to
$\tau_n+1$ and $\vt$ to $\vt+1$ for $k\in\ol{1,n-1}$. The system of
(0.43) and (0.44) is an analogue of (0.33). Furthermore, if
$\mbox{\it Re}\:\tau_n>0$ and $\mbox{\it Re}\:(\vt-\tau_n)>0$, then
we have the Euler integral representation
 $${\msr X}_A=\frac{\G(\vt)}{\G(\vt-\tau_n)\G(\tau_n)}\int_0^1\left[\prod_{r=1}^{n-1}
(\sum_{s=r}^{n-1}P_{[r,s]}+tP_{[r,n]})^{-\tau_r}\right]t^{\tau_n-1}(1-t)^{\vt-\tau_n-1}dt
\eqno(0.45)$$ on the region
$[P_{[r,n]}/(\sum_{s=r}^{n-1}P_{[r,s]})]\not\in (-\infty,-1)$ for
$r\in\{1,...,n-1\}$. Expression (0.45) is exactly an analogue of
(0.34).

The Calogero-Sutherland model\index{Calogero-Sutherland model} is an
exactly solvable quantum many-body system in one-dimension
 (cf. [Cf, Sb]), whose Hamiltonian is given by
$$H_{CS}=\sum_{\iota=1}^n\ptl_{x_\iota}^2+K\sum_{1\leq p<q\leq n}\frac{1}{\sinh^2(x_p-x_q)},\eqno(0.46)$$
where $K$ is a constant. Set
$$\xi_{r_2,r_1}^A=\prod_{s=r_1}^{r_2-1}\frac{e^{2x_{r_2}}}{e^{2x_{r_2}}-e^{2x_s}}\qquad
 \for\;\;1\leq r_1<r_2\leq n.\eqno(0.47)$$
Take $(\lmd_1,...,\lmd_n)\in\mbb{C}^n$  such that
$$\lmd_1-\lmd_2=\cdots
=\lmd_{n-2}-\lmd_{n-1}=\mu\;\;\mbox{and}\;\;\lmd_{n-1}-\lmd_n=\sgm
\not\in\mbb{N}, \eqno(0.48)$$ for some constants $\mu$ and $\sgm$.
Based on Etingof's work [Ep] related to representations of simple
Lie algebra of type $A_{n-1}$, we proved in [X14] that
$$\prod_{r=1}^ne^{2(\lmd_r+(n+1)/2-r)x_r}
{\msr
X}_A(\mu+1,..,\mu+1,-\mu;-\sgm)\{\xi_{r_2,r_1}^A\}\eqno(0.49)$$ is a
solution of the Calogero-Sutherland model.

We also introduced in [X14] new multi-variable hypergeometric
functions related to the simple Lie algebras of types $B_n,C_n,D_n$,
where the functions of type $C$ give rise to the solutions of the
corresponding Olshanesky-Perelomov model (cf. [OP]).

Partial differential equations are also used by us [X12, X17, X18,
X19, X22] to explicitly determine singular vectors in certain weight
modules of finite-dimensional simple Lie algebras and the structures
of the modules (Some of the results were abstractly determined by
Brion [B] before us). In particular,  we obtained combinatorial
identities on the dimensions of an infinite family of irreducible
modules of $F_4,\;E_6$ and $E_7$, respectively.
Differential-operator representations of Lie algebras are called
{\it oscillator representations} in physics terminology.\index{
oscillator representations} In the oscillator representation of the
simple Lie algebra of type $F_4$ lifted from its basic irreducible
module, we showed in [X18] that the number of irreducible submodules
contained in the space of homogeneous harmonic polynomials with
degree $k\geq 2$ is $\geq \llbracket
k/3\rrbracket+\llbracket(k-2)/3\rrbracket+2$. Moreover, in the
polynomial algebra over the 56-dimensional basic irreducible module
of the simple Lie algebra $E_7$,  we  found in [X22] two
three-parameter families of irreducible submodules in the solution
space of the Cartan's fourth-order $E_7$-invariant partial
differential equation.

For any $\lmd\in H^\ast$, there exists a unique highest-weight
${\msr G}$-module $M(\lmd)$ with highest weight $\lmd$ such that any
other highest-weight ${\msr G}$-module with highest weight $\lmd$
must be a quotient module of $M(\lmd)$. The module $M(\lmd)$ is
called a {\it Verma module} (cf. [V1, V2]). Suppose that ${\msr G}$
is of type $A_{n-1}$. The Weyl group of ${\msr G}$ in this case is
exactly the group $S_n$ of permutations on $\{1,2,...,n\}$. In
[X19], we derived a system of variable-coefficient second-order
linear partial differential equations that determines the singular
vectors in $M(\lmd)$, and a differential-operator representation of
$S_n$ on
 the related space of truncated power series. We proved that the solution
 space of the system of partial differential equations is exactly
 spanned by $\{\sgm(1)\mid \sgm\in S_n\}$. Moreover, the singular
 vectors in $M(\lmd)$ are given by those
 $\sgm(1)$ that are polynomials.
Xiao [XW1] found the explicit formula for $\sgm(1)$ when $\sgm$ is
the reflection associated with a positive root and determined the
condition for it to be a polynomial. A similar result as that in
[X19] for the simple Lie algebra of type $C_2$ was proved by us
[X12]. Xiao [XW2] generalized our result of $A_{n-1}$ in [X19] and
that of $C_2$ in [X12] to all finite-dimensional simple Lie
algebras.

 The second important feature of this
book is the construction of irreducible representations of simple
Lie algebras from their natural inhomogeneous oscillator
representations, which give fractional representations of the
corresponding Lie groups.

Suppose that $n\geq 2$ is an integer. Let ${\msr G}$ be the simple
Lie algebra of types $A_n,\;B_n,\;D_n,\\ E_6$ and $E_7$. It is known
that it has the following subalgebra decomposition
$${\msr G}={\msr
G}_-\oplus{\msr G}_0\oplus {\msr G}_+\eqno(0.50)$$ with
$$[{\msr G}_+,{\msr G}_-]\subset {\msr G}_0,\qquad[{\msr G}_0,{\msr
G}_\pm]\subset{\msr G}_\pm,\qquad[{\msr G}_\pm,{\msr
G}_\pm]=\{0\},\eqno(0.51)$$ where $\msr G_0$ is a direct sum of a
one-dimensional trivial Lie algebra with the simple Lie algebra
$\msr G_0'$ of types $A_{n-1},\;B_{n-1},\;D_{n-1},\;D_5$ and $E_6$,
respectively. Moreover, ${\msr G}_-$ forms the minimal natural
irreducible $\msr G_0'$-module with respect to the adjoint
representation of ${\msr G}$ if $\msr G$ is not of type $E_6$, and
${\msr G}_-$ gives the ``spin representation" of $\msr G_0'$ when
$\msr G$ is  of type $E_6$. Moreover, ${\msr G}_+$ is its dual $\msr
G_0'$-module of $\msr G_-$. We lift the representation of $\msr
G_0'$ on $\msr G_-$ to an oscillator representation of ${\msr G}_0$
on the algebra $\msr A=\mbb C[x_1,...,x_m]$ of polynomial functions
$\msr G_-$, where $m=\dim \msr G_-$.

Denote by $\msr A_k$ the subspace of homogeneous polynomials with
degree $k$. We extend the oscillator
 representation of ${\msr G}_0$ to an inhomogeneous oscillator
 representation of ${\msr G}$ by fixing $\msr G_+|_{\msr
 A}=\sum_{i=1}^m\mbb C\ptl_{x_i}$ and solving $\msr G_-|_{\msr
 A}$ in $\sum_{i=1}^m\msr A_2 \ptl_{x_i}$. When $\msr G$ is of type
 $A_n$, the extended representation is exactly the one induced from
 the corresponding projective transformations (cf. [ZX]). If $\msr G$ is of types
 $B_n$ and $D_n$, the extended representation is exactly the one induced from
 the corresponding conformal transformations (cf. [XZ]). In the cases of $E_6$
 and $E_7$, we found the corresponding fractional representations of
 the corresponding Lie groups (cf. [X23, X24]).

 We define a Lie algebra
 $$\msr G_{\msr A}={\msr G}|_{\msr A}\oplus \msr A\msr
 G_0\eqno(0.52)$$
 with the Lie bracket
 $$[d_1+f_1u,d_2+f_2v]=(d_1d_2-d_2d_1)+(d_1(f_2)v-d_2(f_1)u+f_1f_2[u,v])\eqno(0.53)$$
 for $d_1,d_2\in {\msr G}|_{\msr A}$, $f_1,f_2\in\msr A$ and
 $u,v\in \msr G_0$. Using Shen's idea of mixed product (cf. [Sg]),
 we found a Lie algebra monomorphism $\iota: {\msr G}\rta \msr G_{\msr
 A}$ such that $\iota(\msr G)\not\subset\msr G|_{\msr A}$. Let $M$ be any $\msr G_0'$-module. We
 extend it to a $\msr G_0$-module by letting the central element
 take a constant map. Define a $\msr G_{\msr A}$-module structure on $\wht {M}=\msr
 A\otimes_{\mbb C}M$ by
 $$(d+fu)(g\otimes w)=d(g)\otimes w+fg\otimes u(w)\eqno(0.54)$$
for $d\in
 {\msr G}|_{\msr A},\;f,g\in\msr A,\;u\in\msr G_0$ and $w\in
 M.$ Then $\wht {M}$ becomes a ${\msr G}$-module with the
 action
 $$\xi(\varpi)=\iota(\xi)(\varpi)\qquad\for\;\;\xi\in {\msr
 G},\;\varpi\in \wht {M}.\eqno(0.55)$$

The map $M\mapsto \widehat {M}$ defines a new functor from $\msr
G_0'$-{\bf Mod} to ${\msr
 G}$-{\bf Mod}. When $M$ is an irreducible $\msr G_0'$-module, $\mbb C[\msr G_-](1\otimes M)$ is  naturally an irreducible
$\msr G$-module due to the fact $\msr G_+|_{\msr
 A}=\sum_{i=1}^m\mbb C\ptl_{x_i}$, where $\mbb C[\msr G_-]$ is the polynomial algebra generated by $\msr G_-$. Thus the less explicit functor
 $M\mapsto \mbb C[\msr G_-](1\otimes M)$ gives a polynomial extension
from irreducible $\msr G_0'$-modules to irreducible ${\msr
 G}$-modules. In particular, it can be used to construct Gel'fand-Zetlin type bases for $E_6$ from those for
 $D_5$ and then construct  Gel'fand-Zetlin type bases for $E_7$ from those
 for $E_6$ (cf. [GZ1, GZ2, Ma]). When $M$ is a
 finite-dimensional irreducible $\msr G_0'$-module, Zhao and the author [ZX, XZ] determined the condition
 for $\wht {M}=\mbb C[\msr G_-](1\otimes M)$ if $\msr G$ is of types $A_n,\;B_n$ and $D_n$.  The condition
 for $E_6$ and $E_7$ was found by the author in [X23, X24]. In these
 works,  the idea of Kostant's characteristic identities in
[Kb] played a key role. Our approaches heavily depend on the
explicit decompositions of $\msr A$ and the tensor module of $\msr
G_-$ with any finite-dimensional $\msr G_0'$-module into direct sums
of irreducible $\msr G_0'$-submodules. This is a reason why we
adopted the case-by-case approach.

The above works also yield a one-parameter ($c$) family  of
inhomogeneous first-order differential operator representations of
$\msr G$. Letting these operators act on the space of
exponential-polynomial functions that depend on a parametric vector
$\vec a\in \mbb C^n$, we proved in [X26-X28] that the space forms an
irreducible $\msr G$-module  for any $c\in\mbb C$ if $\vec a$ is not
on an explicitly given algebraic variety when $\msr G$ is not of
type $A_n$. The equivalent combinatorial properties of the
representations played key roles. In the case of $A_n$, we [X25] got
new non-weight irreducible modules of the simple Lie algebras of
type $A_n$ and $C_m$ if $n=2m$ is an even integer. By partially
swapping differential operators and multiplication operators, we
obtain more general oscillator representations of $\msr G$. If
$c\not\in\mbb Z$, we obtained in [X25] explicit infinite-dimensional
weight modules with finite-dimensional weight subspaces for the
simple Lie algebras of types $A_n$ and $C_m$ when $n=2m$ is an even
integer. When $c\not\in\mbb Z/2$, we [X26] found  explicit
infinite-dimensional weight modules with finite-dimensional weight
subspaces for the simple Lie algebras of types $B_n$ and $D_n$.
These weight modules are not of highest type.

Below we give a chapter-by-chapter introduction. Throughout this
book, $\mbb{F}$ is always a field with characteristic 0 such as the
field $\mbb{Q}$ of rational numbers, the field $\mbb{R}$ of real
numbers and the field $\mbb{C}$ of complex numbers. All the vector
spaces are assumed over $\mbb{F}$ unless they are specified.

The first part of this book is the classical theory of
finite-dimensional Lie algebras and their representations, which
serves as the preparation of our later main contents.

 In Chapter 1, we give  basic concepts and examples of Lie
algebras. Moreover, Engel's theorem on nilpotent Lie algebras and
Lie's theorem on solvable Lie algebras are proved. Furthermore, we
derive the Jordan-Chevalley decomposition of a linear transformation
and use it to show Cartan's criterion on the solvability.

In Chapter 2, we first introduce the Killing form and prove that the
semisimplicity of a finite-dimensional Lie algebra over $\mbb{C}$ is
equivalent to the nondegeneracy of its Killing form. Then we use the
Killing form to derive the decomposition of a finite-dimensional
semisimple Lie algebra over $\mbb{C}$ into a direct sum of simple
ideals. Moreover, it is showed that a derivation  of ${\msr G}$ must
be a left multiplication operator of $\msr G$ for such a Lie
algebra. Furthermore, we study the completely reducible modules of a
Lie algebra and proved the Weyl's theorem of complete reducibility.
The equivalence of the complete reducibility of real and complex
modules is also given. Cartan's root-space decomposition of a
finite-dimensional semisimple Lie algebra over $\mbb{C}$ is derived.
In particular, we prove that such a Lie algebra is generated by two
elements. The complete reducibility of finite-dimensional modules of
$A_1$ plays an important role in proving the properties of the
corresponding root systems.

In Chapter 3, we start with the axiom of root system and
 give the root systems of special linear algebras, orthogonal Lie
 algebras and symplectic Lie algebras. Then we derive some basic
 properties of root systems; in particular, the existence of the bases of root systems.
 As finite symmetries of root systems, the Weyl
 groups are introduced and studied in detail.  The classification and explicit constructions of root
 systems are presented. The automorphism groups of roots systems are
 determined. As a preparation for later representation theory of Lie
 algebras, the corresponding weight lattices and their saturated
 subsets are investigated.

In Chapter 4, we show that  the structure of a finite-dimensional
semisimple Lie algebra over $\mbb{C}$ is completely determined by
its root system. Moreover, we prove that any two Cartan subalgebras
of such a Lie algebra ${\msr G}$ are conjugated under the group of
inner automorphisms of ${\msr G}$. In particular, the automorphism
group of ${\msr G}$ is determined when it is simple. Furthermore, we
give explicit constructions of the simple Lie algebras of
exceptional types.

 By
Weyl's Theorem, any finite-dimensional representation of a
finite-dimensional semisimple Lie algebra over $\mbb{C}$ is
completely reducible. The main goal in Chapter 5 is to study
finite-dimensional irreducible representations of a
finite-dimensional semisimple Lie algebra over $\mbb{C}$. First we
introduce the universal enveloping algebra of a Lie algebra and
prove the Poincar\'{e}-Birkhoff-Witt (PBW) Theorem on its basis.
Then we use the universal enveloping algebra of a finite-dimensional
semisimple Lie algebra ${\msr G}$ to construct the Verma modules of
${\msr G}$. Moreover, we prove that any finite-dimensional ${\msr
G}$-module is the quotient of a Verma module modulo its maximal
proper submodule, whose generators are explicitly given.
Furthermore, the Weyl's character formula of a finite-dimensional
irreducible ${\msr G}$-module is derived and the dimensional formula
of the module is determined. Finally, we decompose the tensor module
of two finite-dimensional irreducible ${\msr G}$-modules into a
direct sum of irreducible ${\msr G}$-submodules in terms of
characters.

Part II is about the explicit representations of finite-dimensional
simple Lie algebra over $\mbb{C}$, which is the main content of this
book.

In Chapter 6, we give various explicit representations of the simple
Lie algebras of type $A$. First we present a fundamental lemma of
solving flag partial differential equations for polynomial
solutions, which was due to our work [X16]. Then we present the
canonical bosonic and fermionic oscillator representations over
their minimal natural modules and minimal orthogonal modules.
Moreover, we determine the structure of the noncanonical oscillator
representations obtained from the canonical bosonic oscillator
representations over their minimal natural modules and minimal
orthogonal modules by partially swapping differential operators and
multiplication operators. The case over their minimal natural
modules was due to Howe [Hr4]. The results in the case over their
minimal orthogonal modules are generalizations of the classical
theorem on harmonic polynomials. The results  were due to Luo and
the author [LX1].  We construct a functor from the category of
$A_{n-1}$-modules to the category of $A_n$-modules, which is related
to $n$-dimensional projective transformations. This work was due to
Zhao and the author [ZX]. Finally, we present multi-parameter
families of irreducible projective oscillator representations of the
algebras given in [X25].

Chapter 7 is devoted to natural explicit representations of the
simple Lie algebra of type $D$. First we present the canonical
bosonic and fermionic oscillator representations over their minimal
natural modules. Then we determine the structure of the noncanonical
oscillator representations obtained from the above bosonic
representations by partially swapping differential operators and
multiplication operators, which are generalizations of the classical
theorem on harmonic polynomials (cf. [LX2]).  Furthermore, we speak
about a functor from the category of $D_n$-modules to the category
of $D_{n+1}$-modules, which is related to $2n$-dimensional conformal
transformations (cf. [XZ]). In addition, we present multi-parameter
families of irreducible conformal oscillator representations of the
algebras given in [X26].

Chapter 8 is about  natural explicit representations of the simple
Lie algebra of type $B$. We give the canonical  bosonic and
fermionic oscillator representations over their minimal natural
modules. Moreover, we determine the structure of the noncanonical
oscillator representations obtained from the above bosonic
representations by partially swapping differential operators and
multiplication operators, which are generalizations of the classical
theorem on harmonic polynomials (cf. [LX2]). Furthermore, we present
a  functor from the category of $B_n$-modules to the category of
$B_{n+1}$-modules, which is related to $(2n+1)$-dimensional
conformal transformations (cf. [XZ]). Besides, we present
multi-parameter families of irreducible conformal oscillator
representations of the algebras given in [X26].

Chapter 9 determines the structure of the canonical bosonic and
fermionic oscillator representations  of the simple Lie algebras of
type $C$ over their minimal natural modules. Moreover, we study the
noncanonical oscillator representations obtained from the above
bosonic representations by partially swapping differential operators
and multiplication operators, and obtain a two-parameter family of
new infinite-dimensional irreducible representations (cf. [LX3]).
Finally, we present multi-parameter families of irreducible
projective oscillator representations of the algebras given in
[X25].

In Chapter 10, we determine the structure of the canonical bosonic
and fermionic oscillator representations of the simple Lie algebra
of type $G_2$ over its 7-dimensional module. Moreover, we use
partial differential equations to  find the explicit irreducible
decomposition of the space of polynomial functions on 26-dimensional
basic irreducible module of the simple Lie algebra of type $F_4$
(cf. [X18]).

Chapter 11 studies explicit representations of the simple Lie
algebra of type $E_6$. First we prove the cubic $E_6$-generalization
of the classical theorem on harmonic polynomials given [X17]. Then
we study the functor from the module category of $D_5$ to the module
category of $E_6$ developed in [X23]. Finally, we give a family of
inhomogeneous oscillator representations of the simple Lie algebra
of type $E_6$ on a space of exponential-polynomial functions and
prove that their irreducibility is related to an explicit given
algebraic variety (cf. [X27]).

Explicit representations of the simple Lie algebra of type $E_7$ are
given in Chapter 12. By solving certain partial differential
equations, we find the explicit decomposition of the polynomial
algebra over the 56-dimensional basic irreducible module of the
simple Lie algebra $E_7$ into a sum of irreducible submodules (cf.
[X22]). Then we study the functor from the module category of $E_6$
to the module category of $E_7$ developed in [X24]. Moreover, we
construct a family of irreducible inhomogeneous oscillator
representations of the simple Lie algebra of type $E_7$ on a space
of exponential-polynomial functions, related to an explicitly given
algebraic variety (cf. [X28]).

Part III is an extension of Part II.

In Chapter 13, we first establish two-parameter $\mbb{Z}^2$-graded
supersymmetric oscillator generalizations of the classical theorem
on harmonic polynomials for the general linear Lie superalgebra.
Then we extend the result to two-parameter $\mbb{Z}$-graded
supersymmetric oscillator generalizations of the classical theorem
on harmonic polynomials for the ortho-symplectic Lie superalgebras.
This chapter is a reformulation of Luo and the author's work [LX4].

Linear codes with large minimal distances are important error
correcting codes in information theory. Orthogonal codes have more
applications in the other fields of mathematics. In Chapter 14, we
study the binary and ternary orthogonal codes generated by the
weight matrices of finite-dimensional  modules of
  simple Lie algebras. The Weyl groups of the Lie algebras act on these codes isometrically.
  It turns out
 that certain weight matrices of the simple Lie algebras of types $A$ and $D$ generate doubly-even  binary orthogonal codes
 and ternary orthogonal codes with large minimal distances. Moreover, we prove that the
 weight matrices of $F_4$, $E_6$, $E_7$
 and $E_8$ on their minimal irreducible modules and adjoint modules all generate ternary orthogonal codes
 with large minimal distances.
 In determining the minimal distances, we have used the Weyl groups and branch rules of
 the irreducible representations of the related simple Lie algebras.
 The above results are taken from [X20].

 In Chapter 15, we prove that certain variations of the classical Weyl
 functions are solutions of the Calogero-Sutherland model and its generalizations---the Olshanesky-Perelomov model in
 various cases. New multi-variable hypergeometric
functions related to the root systems of classical simple Lie
algebras are introduced. In particular, those of type $A$ give rise
to solutions of the Calogero-Sutherland model based on Etingof's
work [Ep] and those of type $C$ yield solutions of the
Olshanesky-Perelomov model of type C based on Etingof and Styrkas'
work [ES]. The differential properties and multi-variable
hypergeometric equations for these multi-variable hypergeometric
functions are given. The  Euler integral representations of the
type-A functions are found. These results come from the author's
work [X14].\psp

{\bf Acknowledgements}\hspace{0.3cm} The research in this book was
partly supported by the National Natural Science Foundation of China
(Grant No. 11171324).

\section*{Notational Conventions}
\addcontentsline{toc}{section}{\numberline{}Notational
Conventions}

$\mbb{N}$: $\{0,1,2,3,...\}$, the set of nonnegative integers.\\
$\mbb{Z}$: the ring of integers.\\
$\mbb{Q}$: the field of rational numbers.\\
$\mbb{R}$: the field of real numbers.\\
$\mbb{C}$: the field of complex numbers.\\
$\mbb{F}$: a field with characteristic 0, such as, $\mbb{Q},\mbb{R},\mbb{C}$.\\
$\ol{i,i+j}$: $\{i,i+1,i+2,...,i+j\}$, an index set.\\
$\dlt_{i,j}=1$ if $i=j$, $0$ if $i\neq j$.\\
${\msr A}$: an associative algebra.\\
${\msr G}$: a Lie algebra.\\
$\Pi$: the set of positive simple roots.\\
$\Phi$: the set of roots.\\
$\Omega$:  Casimier operator.\\
$\omega$: Casimir element.\\
${\msr W}$: the Weyl group.\\
$\mbb{F}[x_1,x_2,...,x_n]$: the algebra of polynomials in
$x_1,x_2,...,x_n$.\\
$M(\lmd)$: the Verma module of highest weight $\lmd$.\\
$V(\lmd)$: the  irreducible module with highest weight $\lmd$.\\
 $\der {\msr U}$: the derivation Lie algebra of the
algebra ${\msr
U}$. \\
$\mbb{W}_n$: the Witt algebra; that is,
$\der\mbb{F}[x_1,x_2,...,x_n]$.
\\
$\sgm_i$: the $i$th simple reflection.\\
$\Edo V$: the algebra of linear transformations on the vector space
$V$.\\

\pagebreak

\thispagestyle{empty}

\newpage

\tableofcontents

\newpage

\pagenumbering{arabic}

\part{Fundament of Lie Algebras}

\chapter{Preliminary of Lie Algebras}

In this chapter, we give  basic concepts and examples of Lie
algebras. Moreover, Engel's theorem on nilpotent Lie algebras and
Lie's theorem on solvable Lie algebras are proved. Furthermore, we
derive the Jordan-Chevalley decomposition of a linear transformation
and use it to show Cartan's criterion on the solvability.

\section{Basic Definitions and Examples}

In this section, we give definitions and examples of algebras.

 An {\it algebra}\index{algebra} ${\msr A}$ is a
vector space with a bilinear map $m(\cdot,\cdot):{\msr A}\times
{\msr A}\rta {\msr A}$; that is,
$$m(au+bv,w)=am(u,v)+bm(v,w),\;\;m(w,au+bv)=am(w,u)+bm(w,v)\eqno(1.1.1)$$
for $a,b\in\mbb{F}$ and $u,v,w\in{\msr A}$. The map $m(\cdot,\cdot)$
is called an {\it algebraic operation}\index{algebraic operation}.
The algebra ${\msr A}$ is called {\it associative}\index{
associative} if
$$m(u,m(v,w))=m(m(u,v),w)\qquad\for\;\;u,v,w\in{\msr
A}.\eqno(1.1.2)$$ In this case, we denote $m(u,v)=u\cdot v$ or
simply $uv$. The above equation becomes
$$u\cdot (v\cdot w)=(u\cdot v)\cdot w\qquad\for\;\;u,v,w\in{\msr
A}.\eqno(1.1.3)$$ \psp

{\bf Example 1.1.1}. Denote by $\Edo V$ the space of all linear
transformations (endomorphisms)  on a vector space $V$. Then $\Edo
V$ becomes an associative algebra with respect to the composition of
two linear transformations; that is,
$$ (T_1\cdot T_2)(u)=T_1(T_2(u))\qquad\for\;\;T_1,T_2\in\Edo
V,\;u\in V.\eqno(1.1.4)$$ The other familiar examples of associative
algebras are the algebra $M_{n\times n}(\mbb{F})$ of all $n\times n$
matrices with entries in $\mbb{F}$ and the algebra
$\mbb{F}[x_1,x_2,...,x_n]$ of all the polynomials in variables
$\{x_1,x_2,...,x_n\}$.

 \vspace{0.1cm}

{\bf Example 1.1.2}. A {\it group}\index{group} $G$ is a set with a
map $\cdot: G\times G\rta G$ and a spacial element $1_G$ such that
$$ (g_1\cdot g_2)\cdot g_3=(g_1\cdot g_2)\cdot
g_3,\;\;1_G\cdot g_1=g_1\cdot 1_G=g_1\qquad\for\;\;g_1,g_2,g_3\in
G,\eqno(1.1.5)$$ and for any $g\in G$, there exists an element
$g^{-1}\in G$ for which
$$g\cdot g^{-1}=g^{-1}\cdot g=1_G.\eqno(1.1.6)$$
If
$$g_1\cdot g_2=g_2\cdot g_1\;\;\for\;\;g_1,g_2\in G,\eqno(1.1.7)$$
we call ${\msr G}$ an {\it abelian group}.\index{abelian group}

Given a group $G$. Denote by $\mbb{F}[G]$ the vector space with a
basis $\{\Im_g\mid g\in G\}$. Define an algebraic operation ``$\cdot
$" on $\mbb{F}[G]$ by:
$$\Im_{g_1}\cdot\Im_{g_2}=\Im_{g_1\cdot
g_2}\qquad\for\;\;g_1,g_2\in G.\eqno(1.1.8)$$ Then $\mbb{F}[G]$
becomes an associative algebra with respect to (1.1.8), which is
called the {\it group algebra} of $G$.\psp

A {\it Lie algebra} ${\msr G}$\index{Lie algebra} is a vector space
with an algebraic operation $[\cdot,\cdot]$ satisfying the following
two axioms:
$$[u,v]=-[v,u]\qquad(\mbox{Skew Symmetry}),\eqno(1.1.9)$$\index{skew symmetry}
$$[[u,v],w]+[[v,w],u]+[[w,u],v]=0\qquad(\mbox{Jacobi
Identity})\eqno(1.1.10)$$\index{Jacobi identity} for $u,v,w\in {\msr
G}$. The algebraic operation $[\cdot,\cdot]$ is called the {\it Lie
bracket}\index{Lie bracket} of ${\msr G}$. Note
$$[u,u]=-[u,u]\lra [u,u]=0\qquad \for\;\;u\in{\msr
G}.\eqno(1.1.11)$$

 \psp

{\bf Example 1.1.3}. Let $({\msr A},\cdot)$ be an associative
algebra. Define another algebraic operation $[\cdot,\cdot]$ on
${\msr A}$:
$$[u,v]=u\cdot v-v\cdot u,\eqno(1.1.12)$$ which is called the {\it
commutator}\index{commutator} of $u$ and $v$ in ${\msr A}$. The skew
symmetry is obvious.  Moreover, \begin{eqnarray*}\hspace{2cm}&
&[[u,v],w]=[u\cdot v-v\cdot u,w]=(u\cdot v)\cdot w\\ & &-(v\cdot
u)\cdot w-w\cdot(u\cdot v)+w\cdot(v\cdot
u).\hspace{5.3cm}(1.1.13)\end{eqnarray*} Hence
\begin{eqnarray*}&&[[u,v],w]+[[v,w],u]+[[w,u],v]=
(u\cdot v)\cdot w-(v\cdot u)\cdot w\\ & &-w\cdot(u\cdot
v)+w\cdot(v\cdot u)+(v\cdot w)\cdot u-(w\cdot v)\cdot
u-u\cdot(v\cdot w)\\ & &+u\cdot(w\cdot v)+ (w\cdot u)\cdot v-(u\cdot
w)\cdot v-v\cdot(w\cdot u)+v\cdot(u\cdot
v)=0.\hspace{2.2cm}(1.1.14)\end{eqnarray*} Thus $({\msr
A},[\cdot,\cdot])$ forms a Lie algebra. In particular, we have Lie
algebras $(\Edo V,[\cdot,\cdot])$ and $(\mbb{F}[G],[\cdot,\cdot])$.

Indeed, the above Lie bracket can be generalized as follows. Take
any fixed nonzero element $\xi\in {\msr A}$ and define an algebraic
operation $[\cdot,\cdot]_\xi$ on ${\msr A}$ by
$$[u,v]_\xi=u\cdot\xi\cdot v-v\cdot\xi\cdot
u\qquad\for\;\;u,v\in{\msr A}.\eqno(1.1.15)$$ It can be verified
that $({\msr A},[\cdot,\cdot]_\xi)$ forms a Lie algebra (exercise).
This Lie algebra structure is useful in the study of
infinite-dimensional simple Lie algebras and has also appeared in
the other fields of mathematics. \psp

{\bf Example 1.1.4}. There are so-called ``Sin algebras"\index{Sin
algebra} appeared in geometry and dynamics. We assume
$\mbb{F}=\mbb{R}$, the field of real numbers. Let $(\G,+)$ be an
abelian group. Suppose that $\phi(\cdot,\cdot):\G\times
\G\rightarrow
 \mbb{R}$ is a skew-symmetric $\mbb{Z}$-bilinear map; that
 is,
 $$\phi(\al,\be)=-\phi(\be,\al),\qquad\phi(\al+\be,\gm)=\phi(\al,\gm)
 +\phi(\be,\gm).\eqno(1.1.16)$$
On the group algebra $\mbb{R}[\G]$, we define another algebraic
operation $[\cdot,\cdot]$
 defined by
$$[\Im_{\al},\Im_{\be}]=\sin \phi(\al,\be)\;\Im_{\al+\be}\qquad\for\;\;\al,\be\in\G.
\eqno(1.1.17)$$ Then $(\mbb{R}[\G],[\cdot,\cdot])$ forms a Lie
algebra (exercise), which is called a {\it Sin algebra}. \psp

\section{Lie Subalgebras}

In this section, we present derivation Lie algebras, Lie algebras of
Block type, special linear Lie algebras, orthogonal Lie algebras,
symplectic Lie algebras and co-adjoint Poisson structures.

A subspace $\msr L$ of a Lie algebra ${\msr G}$ is called a {\it
subalgebra}\index{Lie subalgebra} of ${\msr G}$ if
$$ [u,v]\in \msr L\qquad \mbox{for any}\;\;u,v\in \msr L.\eqno(1.2.1)$$
Let $({\msr U},\circ)$ be any algebra. A linear transformation $d$
on ${\msr U}$ is called a {\it derivation}\index{derivation} of
${\msr U}$ if
$$ d(u\circ v)=d(u)\circ v+u\circ d(v)\qquad\for\;\;u,v\in{\msr
U}.\eqno(1.2.2)$$
 Denote
$$\der {\msr U}=\mbox{the set of all derivations on}\;{\msr
U}.\eqno(1.2.3)$$ Then $\der {\msr U}$ forms a subspace of $\Edo
{\msr U}$. Moreover, for $d_1,d_2\in \der{\msr U}$ and $u,v\in{\msr
U}$, we have \begin{eqnarray*} \hspace{2cm}& & d_1d_2(u\circ
v)=d_1(d_2(u)\circ v+u\circ d_2(v))=d_1d_2(u)\circ v\\ &
&+d_2(u)\circ d_1(v)+d_1(u)\circ d_2(v)+u\circ
d_1d_2(v).\hspace{4.3cm}(1.2.4)\end{eqnarray*} Thus
\begin{eqnarray*} & &[d_1,d_2](u\circ v)\\&=& (d_1d_2-d_2d_1)(u\circ
v)=(d_1d_2-d_2d_1)(u\circ v)=d_1d_2(u\circ v)-d_2d_1(u\circ v)\\
&=&d_1d_2(u)\circ v+d_2(u)\circ d_1(v)+d_1(u)\circ d_2(v)+u\circ
d_1d_2(v)\\ & &-[d_2d_1(u)\circ v+d_1(u)\circ d_2(v)+d_2(u)\circ
d_1(v)+u\circ d_2d_1(v)]\\ &=&(d_1d_2(u)-d_2d_1(u))\circ
v+u\circ(d_1d_2(v)-d_2d_1(v))\\ &=&[d_1,d_2](u)\circ v+u\circ
[d_1,d_2](v),\hspace{8.4cm}(1.2.5)\end{eqnarray*}or equivalently,
$[d_1,d_2]\in\der{\msr U}$. Hence $\der {\msr U}$ forms a Lie
subalgebra of $\Edo {\msr U}$. \psp

{\bf Example 1.2.1}. Let ${\msr A}=\mbb{F}[x_1,x_2,...,x_n]$ be the
polynomial algebra in $n$ variables. It can be verified (exercise)
that
$$ \mbb W_n=\der {\msr A}=\{\sum_{i=1}^nf_i\ptl_{x_i}\mid f_i\in
\mbb{F}[x_1,x_2,...,x_n]\}.\eqno(1.2.6)$$ Note that
$$f_1\ptl_{x_i}(f_2\ptl_{x_j}(g))=f_1\ptl_{x_i}(f_2)\ptl_{x_j}(g)
+f_1f_2\ptl_{x_i}\ptl_{x_j}(g).\eqno(1.2.7)$$ Thus
\begin{eqnarray*}& &[f_1\ptl_{x_i},f_2\ptl_{x_j}](g)\\&=&
f_1\ptl_{x_i}(f_2\ptl_{x_j}(g))-f_2\ptl_{x_j}(f_1\ptl_{x_i}(g))\\
&=&f_1\ptl_{x_i}(f_2)\ptl_{x_j}(g)
+f_1\ptl_{x_i}\ptl_{x_j}(g)-[f_2\ptl_{x_j}(f_1)\ptl_{x_i}(g)
+f_2\ptl_{x_j}\ptl_{x_i}(g)]\\
&=&f_1\ptl_{x_i}(f_2)\ptl_{x_j}(g)-f_2\ptl_{x_j}(f_1)\ptl_{x_i}(g)
=(f_1\ptl_{x_i}(f_2)\ptl_{x_j}-f_2\ptl_{x_j}(f_1)\ptl_{x_i})(g).\hspace{1.7cm}(1.2.8)
\end{eqnarray*}
Therefore,
$$[f_1\ptl_{x_i},f_2\ptl_{x_j}]=f_1\ptl_{x_i}(f_2)\ptl_{x_j}-f_2\ptl_{x_j}(f_1)\ptl_{x_i}.
\eqno(1.2.9)$$ The Lie algebra $\mbb W_n$ is called a {\it Witt
algebra}.\index{Witt algebra}

\psp

{\bf Example 1.2.2}. Let ${\msr A}$ be an associative algebra
satisfying the commutativity:
$$ uv=vu\qquad\for\;\;u,v\in{\msr A}.\eqno(1.2.10)$$
Such an algebra is called a {\it commutative associative
algebra}.\index{commutative associative algebra}
 Suppose that we have two commuting derivations
$\ptl_1$ and $\ptl_2$ on ${\msr A}$. Define another algebraic
operation $[\cdot,\cdot]$ on ${\msr A}$ by:
$$[u,v]=\ptl_1(u)\ptl_2(v)-\ptl_2(u)\ptl_1(v)+u\ptl_1(v)-\ptl_1(u)v\qquad\for\;\;u,v\in{\msr
A}.\eqno(1.2.11)$$ Then $({\msr A},[\cdot,\cdot])$ forms a Lie
algebra (exercise), which is called a {\it Lie algebra of Block
type}\index{Block algebra}. When ${\msr A}$ is the group algebra of
an abelian group, the Lie algebra was discovered by Block [Br] in
1958 and the general case was obtained by the author [X7] in 1999.
Chen [Cl] generalized the algebras in [X7] to more general Lie
algebras of twisted Hamiltonian type.\psp

The Jacobi identity can be rewritten as
$$ [u,[v,w]]=[[u,v],w]+[v,[u,w]],\eqno(1.2.12)$$
which shows that the left multiplication operator is a derivation.
Let ${\msr A}$ be an associative algebra. A linear map $\vf:{\msr
A}\rta \mbb{F}$ is called a {\it trace map}\index{trace map} if
$$ \vf(uv)=\vf(vu)\qquad\for\;\;u,v\in{\msr A}.\eqno(1.2.13)$$
Recall that the trace of an $n\times n$ matrix $B=(b_{i,j})_{n\times
n}$ is defined by
$$\tr B=\sum_{i=1}b_{i,i},\eqno(1.2.14)$$
which is a trace map of the algebra $M_{n\times n}(\mbb{F})$ of all
$n\times n$ matrices. For a group algebra $\mbb{F}[G]$ defined in
Example 1.1.2, we have the following trace map:
$$\tr \Im_g=\dlt_{g,1_G}\qquad\for\;\;g\in G.\eqno(1.2.15)$$

Suppose that ${\msr A}$ is an associative algebra with a trace map
$\tr$. Recall that ${\msr A}$ forms a Lie algebra with respect to
the commutator (1.1.12). Moreover, the set
$$s({\msr A})=\{u\in{\msr A}\mid \tr u=0\}\eqno(1.2.16)$$
forms a Lie subalgebra, because
$$\tr ([u,v])=\tr (uv-vu)=\tr uv-\tr vu=0\eqno(1.2.17)$$
by (1.2.13) for $u,v\in{\msr A}$. When we view $\Edo V$ as a Lie
algebra with respect to the commutator, we call it a {\it general
linear Lie algebra}\index{general linear Lie algebra} and redenote
it as $gl(V)$. If $\dim V=n$, $\Edo V\cong M_{n\times n}(\mbb{F})$.
We use $gl(n,\mbb{F})$ to denote $M_{n\times n}(\mbb{F})$ as the Lie
algebra with respect to the commutator. The Lie subalgebra
$$sl(n,\mbb{F})=\{B\in gl(n,\mbb{F})\mid \tr B=0\}\eqno(1.2.18)$$
is called a {\it special linear Lie algebra}.\index{special linear
Lie algebra}

An {\it involutive anti-automorphism}\index{involutive
anti-automorphism} $\sgm$ of an associative algebra ${\msr A}$ is a
linear transformation on ${\msr A}$ such that
$$\sgm^2=\mbox{Id}_{\msr
A},\;\;\sgm(uv)=\sgm(v)\sgm(u).\eqno(1.2.19)$$ Set
$${\msr A}^\sgm_-=\{u\in{\msr A}\mid \sgm(u)=-u\}.\eqno(1.2.20)$$
Then ${\msr A}^\sgm_-$ forms a Lie subalgebra of ${\msr A}$ with
respect to the commutator. In fact,
$$\sgm([u,v])=\sgm(uv-vu)=\sgm(uv)-\sgm(vu)=
\sgm(v)\sgm(u)-\sgm(u)\sgm(v)=-[\sgm(u),\sgm(v)]\eqno(1.2.21)$$ for
$u,v\in{\msr A}$. The transpose of matrix is an involutive
anti-automorphism of $M_{n\times n}(\mbb{F})$. Thus we have the
following Lie subalgebra
$$o(n,\mbb{F})=\mbox{the set of all}\;n\times
n\;\mbox{skew-symmetric matrices with entries
in}\;\mbb{F},\eqno(1.2.22)$$ which is called an {\it orthogonal Lie
algebra}.\index{orthogonal Lie algebra} In fact, $o(n,\mbb{R})$ is
the Lie algebra of rotation group. To see this, let us show the case
$n=2$. Set
$$A=\left(\begin{array}{cc}0&1\\
-1&0\end{array}\right).\eqno(1.2.23)$$ Then $o(2,\mbb{R})=\mbb{R}A$
and
$$ A^2=\left(\begin{array}{cc}0&1\\
-1&0\end{array}\right)\left(\begin{array}{cc}0&1\\
-1&0\end{array}\right)=\left(\begin{array}{cc}-1&0\\
0&-1\end{array}\right)=-I_2.\eqno(1.2.24)$$ Thus
\begin{eqnarray*} \qquad e^{\theta A}&=&\sum_{n=0}^{\infty}\frac{(\theta
A)^n}{n!}=\sum_{k=0}^{\infty}\frac{\theta^{2k}A^{2k}}{(2k)!}+\sum_{k=0}^{\infty}\frac{\theta^{2k+1}A^{2k+1}}
{(2k+1)!}\\
&=&\sum_{k=0}^{\infty}\frac{(-1)^k\theta^{2k}}{(2k)!}I_2+\sum_{k=0}^{\infty}\frac{(-1)^k
\theta^{2k+1}} {(2k+1)!}A\\ &=&\cos \theta\:I_2+\sin\theta\:A=
\left(\begin{array}{cc}\cos \theta&\sin\theta\\
-\sin\theta&\cos
\theta\end{array}\right),\hspace{5cm}(1.2.25)\end{eqnarray*} which
is a rotation linear transformation on the plane.

Suppose that $n=2k$ is even integer. Denote by $I_k$ the $k\times k$
identity matrix. Use up index ``t" to denote the transpose of
matrix. Define
$$ A^\dg=\left(\begin{array}{cc}0&I_k\\
-I_k&0\end{array}\right)A^t\left(\begin{array}{cc}0&-I_k\\
I_k&0\end{array}\right)\qquad\for\;\;A\in M_{n\times
n}(\mbb{F}).\eqno(1.2.26)$$ It is easy to check that $\dg$ is an
involutive anti-automorphism of $M_{n\times n}(\mbb{F})$. The
corresponding Lie subalgebra is:
$$sp(2k,\mbb{F})=\left\{\left(\begin{array}{cc}A&B\\
C&-A^t\end{array}\right)\mid A,B,C\in M_{k\times
k}(\mbb{F}),\;B^t=B,\;C^t=C\right\},\eqno(1.2.27)$$ which is called
a {\it symplectic Lie algebra}.\index{symplectic  algebra}

On the group algebra $\mbb{F}[G]$, we have the following involutive
anti-automorphism $\sgm$:
$$\sgm(\Im_g)=\Im_{g^{-1}}\qquad\for\;\;g\in G.\eqno(1.2.28)$$
It seems there is not much known on the Lie algebra
$\mbb{F}[G]^\sgm_-$.

Suppose that $\{u_1,u_2,...,u_n\}$ is a basis of a Lie algebra
${\msr G}$. We write
$$ [u_i,u_j]=\sum_{k=1}^nc_{i,j}^ku_k\qquad
\for\;\;i,j\in\{1,2,...,n\}.\eqno(1.2.29)$$ The constants
$\{c_{i,j}^k\}$ are called the {\it structure
constants}\index{structure constants}. Then the skew symmetry
(1.1.9) becomes
$$ c_{i,j}^k=-c_{j,i}^k\eqno(1.2.30)$$
and the Jacobi identity (1.1.10) can be written to:
$$
\sum_{r=1}^n(c_{i,j}^rc_{r,k}^s+c_{j,k}^rc_{r,i}^s+c_{k,i}^rc_{r,j}^s)=0.\eqno(1.2.31)$$
We define an algebraic operation $\{\cdot,\cdot\}$ on
$\mbb{F}[x_1,x_2,...,x_n]$ by
$$\{f,g\}=\sum_{i,j,k=1}^nc_{i,j}^k\ptl_{x_i}(f)\ptl_{x_j}(g)x_k.\eqno(1.2.32)$$
Obviously, $\{\cdot,\cdot\}$ is skew-symmetric. Moreover, for
$f,g,h\in \mbb{F}[x_1,x_2,...,x_n]$, we have
\begin{eqnarray*}\{\{f,g\},h\}&=&\sum_{i,j,k=1}^nc_{i,j}^k\{\ptl_{x_i}(f)\ptl_{x_j}(g)x_k,h\}
\\
&=&\sum_{i,j,k,r,s,t=1}^nc_{i,j}^kc_{r,s}^t\ptl_{x_r}(\ptl_{x_i}(f)\ptl_{x_j}(g)x_k)
\ptl_{x_s}(h)x_t\\
&=&\sum_{i,j,r,s,t=1}^nc_{i,j}^rc_{r,s}^t\ptl_{x_i}(f)\ptl_{x_j}(g)\ptl_{x_s}(h)x_t
+\sum_{i,j,k,r,s,t=1}^n c_{i,j}^kc_{r,s}^t\\ & &\times
(\ptl_{x_r}\ptl_{x_i}(f)\ptl_{x_j}(g)\ptl_{x_s}(h)
+\ptl_{x_i}(f)\ptl_{x_r}\ptl_{x_j}(g)\ptl_{x_s}(h))x_kx_t.\hspace{2cm}(1.2.33)\end{eqnarray*}
Thus
\begin{eqnarray*}& &\{\{f,g\},h\}+\{\{g,h\},f\}+\{\{h,f\},g\}\\
&=&\sum_{i,j,r,s,t=1}^nc_{i,j}^rc_{r,s}^t\ptl_{x_i}(f)\ptl_{x_j}(g)\ptl_{x_s}(h)x_t
+\sum_{i,j,k,r,s,t=1}^nc_{i,j}^kc_{r,s}^t\\ && \times
(\ptl_{x_r}\ptl_{x_i}(f)\ptl_{x_j}(g)\ptl_{x_s}(h)
+\ptl_{x_i}(f)\ptl_{x_r}\ptl_{x_j}(g)\ptl_{x_s}(h))x_kx_t\\ &
&+\sum_{i,j,r,s,t=1}^nc_{i,j}^rc_{r,s}^t\ptl_{x_i}(g)\ptl_{x_j}(h)\ptl_{x_s}(f)x_t
+\sum_{i,j,k,r,s,t=1}^n c_{i,j}^kc_{r,s}^t\\ & &\times
(\ptl_{x_r}\ptl_{x_i}(g)\ptl_{x_j}(h)\ptl_{x_s}(f)
+\ptl_{x_i}(g)\ptl_{x_r}\ptl_{x_j}(h)\ptl_{x_s}(f))x_kx_t\\ & &
+\sum_{i,j,r,s,t}c_{i,j}^rc_{r,s}^t\ptl_{x_i}(h)\ptl_{x_j}(f)\ptl_{x_s}(g)x_t+\sum_{i,j,k,r,s,t}
c_{i,j}^kc_{r,s}^t\\ & &\times
(\ptl_{x_r}\ptl_{x_i}(h)\ptl_{x_j}(f)\ptl_{x_s}(g)
+\ptl_{x_i}(h)\ptl_{x_r}\ptl_{x_j}(f)\ptl_{x_s}(g))x_kx_t\hspace{6cm}\end{eqnarray*}\begin{eqnarray*}
 &=&
\sum_{i,j,r,s,t}(c_{i,j}^rc_{r,s}^t+c_{j,s}^rc_{r,i}^t+c_{s,i}^rc_{r,j}^t)\ptl_{x_i}(f)
\ptl_{x_j}(g)\ptl_{x_s}(h)x_t
\\ &&+
\sum_{i,j,k,r,s,t}c_{i,j}^kc_{r,s}^t
(\ptl_{x_r}\ptl_{x_i}(f)\ptl_{x_j}(g)\ptl_{x_s}(h)+\ptl_{x_i}(h)\ptl_{x_r}\ptl_{x_j}(f)
\ptl_{x_s}(g))x_kx_t
\\ &
&+\sum_{i,j,k,r,s,t}c_{i,j}^kc_{r,s}^t(\ptl_{x_i}(f)\ptl_{x_r}\ptl_{x_j}(g)\ptl_{x_s}(h)
+\ptl_{x_r}\ptl_{x_i}(g)\ptl_{x_j}(h)\ptl_{x_s}(f))x_kx_t\\
&&+\sum_{i,j,k,r,s,t}c_{i,j}^kc_{r,s}^t(\ptl_{x_i}(g)\ptl_{x_r}\ptl_{x_j}(h)\ptl_{x_s}(f)+
\ptl_{x_r}\ptl_{x_i}(h)\ptl_{x_j}(f)\ptl_{x_s}(g))x_kx_t\\ &=&
\sum_{i,j,k,r,s,t}c_{i,j}^k(c_{r,s}^t+c_{s,r}^t)\ptl_{x_r}\ptl_{x_i}(f)\ptl_{x_j}(g)\ptl_{x_s}(h)
x^kx_t\\ &
&+\sum_{i,j,k,r,s,t}c_{i,j}^k(c_{r,s}^t+c_{s,r}^t)\ptl_{x_r}\ptl_{x_i}(g)\ptl_{x_j}(h)
\ptl_{x_s}(f) x^kx_t\\
&&+\sum_{i,j,k,r,s,t}c_{i,j}^k(c_{r,s}^t+c_{s,r}^t)\ptl_{x_r}\ptl_{x_i}(h)\ptl_{x_j}(f)
\ptl_{x_s}(g) x^kx_t=0.\hspace{3.8cm}(1.2.34)\end{eqnarray*} Thus
$\{\cdot,\cdot\}$ is a Lie bracket, which is called {\it co-adjoint
Poisson structure}\index{co-adjoint Poisson structure}
 over ${\msr G}$. The structure has been used in
symplectic geometry and integrable systems.

\section{Ideals}

In this section, we define an ideal of a Lie algebra and simple Lie
algebra. Moreover, we present the basic technique of proving
simplicity.

If $U$ and $V$ are subspaces of a Lie algebra ${\msr G}$, we define
$$[U,V]=\mbox{Span}\{[u,v]\mid u\in U,\;v\in V\}\subset {\msr
G}.\eqno(1.3.1)$$ A subspace ${\cal I}$ of a Lie algebra ${\msr G}$
is called an {\it ideal}\index{ideal} if
$$[{\msr G},{\cal I}]\subset {\cal I}.\eqno(1.3.2)$$
In particular, an ideal is a subalgebra. Obvious examples of ideals
are $\{0\},\;{\msr G}$ and $[{\msr G},{\msr G}]$, which is called
{\it derived algebra}\index{derived algebra} of ${\msr G}$ . Define
the center of ${\msr G}$:
$$Z({\msr G})=\{u\in{\msr G}\mid [u,v]=0\;\mbox{for any}\;v\in
{\msr G}\}.\eqno(1.3.3)$$ Easily see that $Z({\msr G})$ is an ideal
of ${\msr G}$. \psp

{\bf Example 1.3.1}. Set
$${\msr K}_n=\mbb{F}\ptl_{x_1}+\sum_{i=0}^n\mbb{F}x_1^i\ptl_{x_2}\subset
\der \mbb{F}[x_1,x_2].\eqno(1.3.4)$$ Then
$$[\ptl_{x_1},\ptl_{x_1}]=0,\;[\ptl_{x_1},x_1^i\ptl_{x_2}]=ix^{i-1}\ptl_{x_2},\;\;
[x_1^i\ptl_{x_2},x_1^j\ptl_{x_2}]=0.\eqno(1.3.5)$$ Thus ${\msr K}_n$
forms a Lie algebra. Moreover,
$${\cal I}_k=\sum_{i=0}^k\mbb{F}x_1^i\ptl_{x_2},\eqno(1.3.6)$$
is an ideal of ${\msr K}_n$ if $k\leq n$. Furthermore,
$$[{\msr K}_n,{\msr K}_n]={\cal I}_{n-1},\qquad Z({\msr
K}_n)=\mbb{F}\ptl_{x_2}.\eqno(1.3.7)$$ \vspace{0.1cm}

A Lie algebra ${\msr G}$ is called {\it simple}\index{simple} if
$\dim {\msr G}>1$ and the only ideals of ${\msr G}$ are $\{0\}$ and
${\msr G}$. \psp

{\bf Lemma 1.3.1}. {\it Let $T$ be a linear transformation on a
vector space $V$. Suppose that $U$ is a subspace of $V$ such that
$T(U)\subset U$, and $\{u_1,u_2,...,u_k\}$ are eigenvectors of $T$
corresponding to distinct eigenvalues
$\{\lmd_1,\lmd_2,...,\lmd_k\}$. If $\sum_{i=1}^ku_i\in U$, then each
$u_i\in U$ for $i=1,2,...,k$.}

{\it Proof}. Note
$$v_j=T^j(\sum_{i=1}^ku_i)=\sum_{i=1}^k\lmd_i^ju_i\in
U.\eqno(1.3.8)$$ Consider the system
$$\sum_{i=1}^k\lmd_i^ju_i=v_j,\qquad j=0,1,...,k-1.\eqno(1.3.9)$$
Since the coefficient determinant is a Vandermonde determinant
$$\left|\begin{array}{cccc}1&1&\cdots& 1\\
\lmd_1&\lmd_2&\cdots&\lmd_k\\ \vdots&\vdots&\vdots&\vdots\\
\lmd_1^{k-1}&\lmd_2^{k-1}&\cdots
&\lmd_k^{k-1}\end{array}\right|=\prod_{1\leq i<j\leq
k}(\lmd_j-\lmd_i)\neq 0,\eqno(1.3.10)$$ we can solve (1.3.9) for
$\{u_1,u_2,...,u_k\}$ as linear combinations of
$\{v_0,v_1,...,v_{k-1}\}\subset U$. Thus $u_i\in U.\qquad\Box$\psp

{\bf Example 1.3.2}. Denote by $E_{i,j}$ the square matrix with 1 as
its $(i,j)$-entry and 0 as the others. Note
$$sl(2,\mbb{F})=\mbb{F}E_{1,2}+\mbb{F}(E_{1,1}-E_{2,2})+\mbb{F}E_{2,1}\eqno(1.3.11)$$
with the Lie bracket determined by
$$[E_{1,2},E_{2,1}]=E_{1,1}-E_{2,2},\;\;[E_{1,1}-E_{2,2},E_{1,2}]=2E_{1,2},
\;\;[E_{1,1}-E_{2,2},E_{2,1}]=-2E_{2,1}.\eqno(1.3.12)$$ Define a
linear transformation $T$ on $sl(2,\mbb{F})$ by
$$T(u)=[E_{1,1}-E_{2,2},u]\qquad \for\;\;u\in
sl(2,\mbb{F}).\eqno(1.3.13)$$ Then $E_{1,2}$ is an eigenvector of
$T$ with eigenvalue 2, $E_{2,1}$ is an eigenvector of $T$ with
eigenvalue $-2$ and $E_{1,1}-E_{2,2}$ is an eigenvector of $T$ with
eigenvalue 0. Suppose ${\cal I}$ is a nonzero ideal of
$sl(2,\mbb{F})$. Then $T({\cal I})\subset {\cal I}$. Moreover, there
exists a nonzero element
$$u=a_1E_{1,2}+a_2(E_{1,1}-E_{2,2})+a_3E_{2,1}\in {\cal
I}.\eqno(1.3.14)$$ By the above lemma,
$$a_1E_{1,2},a_2(E_{1,1}-E_{2,2}),a_3E_{2,1}\in {\cal
I}.\eqno(1.3.15)$$ Hence one of the elements in
$\{E_{1,2},E_{1,1}-E_{2,2},E_{2,1}\}$ in ${\cal I}$. If one of
$E_{1,2}$ and $E_{2,1}$ in ${\cal I}$, then the first equation in
(1.3.12) implies $E_{1,1}-E_{2,2}\in {\cal I}$. Thus we can always
assume $E_{1,1}-E_{2,2}\in {\cal I}$. The second and third equations
in (1.3.12) imply $E_{1,2},E_{2,1}\in{\cal I}$. Therefore, ${\cal
I}=sl(2,\mbb{F})$; that is, $sl(2,\mbb{F})$ is simple.\psp

In order to prove simplicity for more Lie algebras, we need to
generalize Lemma 1.3.1.\psp

{\bf Lemma 1.3.2}. {\it Let $H$ be a  subspace of $\Edo V$ of
commuting linear transformations on a vector space $V$. Denote
$$H^\ast=\mbox{the set of all linear maps (functions)
from}\;H\;\mbox{to}\;\mbb{F},\eqno(1.3.16)$$ which forms a vector
space with linear operation:
$$(af_1+bf_2)(h)=af_1(h)+bf_2(h)\qquad\for\;\;h\in
H.\eqno(1.3.17)$$ Suppose that $U$ is a subspace of $V$ such that
$$h(U)\subset U\qquad\for\;\;h\in H,\eqno(1.3.18)$$
and $\{u_1,u_2,...,u_k\}\subset V$ are vectors such that there exist
distinct elements $\{\al_1,\al_2,...,\al_k\}$ in $H^\ast$ for which
$$h(u_i)=\al_i(h)u_i\qquad \for\;\;h\in H.\eqno(1.3.19)$$
If $\sum_{i=1}^ku_i\in U$, then each $u_i\in U$ for $i=1,2,...,k$.}

{\it Proof}. Set
$$ \bar{H}=\{h\in H\mid
\al_1(h)=\al_2(h)=\cdots=\al_k(h)=0\}.\eqno(1.3.20)$$ Then
$$\dim H/\bar{H}\leq k.\eqno(1.3.21)$$
Take a subspace $H_1$ of $H$ such that
$$H=H_1\oplus \bar{H},\eqno(1.3.22)$$
which implies $\dim H_1\leq k$, and
$\{\al_1|_{H_1},\al_2|_{H_1},...,\al_k|_{H_1}\}$ are distinct.
Replacing $H$ by $H_1$, we can assume that $\dim H$ is finite.

We prove it by induction on $\dim H$. The lemma is equivalent to
Lemma 1.3.1 when $\dim H=1$. Suppose that the lemma holds for $\dim
H\leq r$. Assume $\dim H=r+1$. We take a nonzero vector $T\in H$ and
an $r$-dimensional subspace $\td{H}$ of $H$ such that
$$H=\mbb{F}T+\td{H}.\eqno(1.3.23)$$
Suppose that $\{\lmd_1,\lmd_2,...,\lmd_s\}$ is the maximal subset of
distinct elements in $\{\al_1(T),\al_2(T),...,\\ \al_k(T)\}$. Denote
$$J_i=\{j\in\{1,2,...,k\}\mid \al_j(T)=\lmd_i\}\eqno(1.3.24)$$
and set
$$ v_i=\sum_{j\in J_i}u_j.\eqno(1.3.25)$$
Then
$$ \sum_{r=1}^sv_r=\sum_{i=1}^ku_i\in
U,\;\;T(v_r)=\lmd_rv_r.\eqno(1.3.25)$$ By Lemma 1.3.1, we have
$$v_1,v_2,...,v_s\in U.\eqno(1.3.26)$$
Since $\{\al_1,\al_2,...,\al_k\}$ are distinct,
$\{\al_j|_{\td{H}}\mid j\in J_i\}$ are distinct. By induction,
$$u_j\in U\qquad \for\;\;j\in J_i.\eqno(1.3.27)$$
Thus the lemma holds.$\qquad\Box$\psp

For two integers $m_1,m_2$, we use the notation:
$$\ol{m_1,m_2}=\left\{\begin{array}{ll}\{m_1,m_1+1,...,m_2\}&\mbox{if}\;\;m_1\leq m_2\\
\emptyset & \mbox{if}\;\;m_1> m_2.\end{array}\right.\eqno(1.3.28)$$

 \pse

{\bf Example 1.3.3}. Let $n>2$.  The Lie algebra
$$sl(n,\mbb{F})=\sum_{1\leq i<
j\leq
n}(\mbb{F}E_{i,j}+\mbb{F}E_{j,i})+\sum_{i=1}^{n-1}\mbb{F}(E_{i,i}-E_{i+1,i+1}).\eqno(1.3.29)$$
Set
$$H=\sum_{i=1}^{n-1}\mbb{F}(E_{i,i}-E_{i+1,i+1})=\{\sum_{i=1}^na_iE_{i,i}\mid
a_i\in\mbb{F},\;\sum_{i=1}^na_i=0\}.\eqno(1.3.30)$$ For
$i\in\ol{1,n}$, we define $\es_i\in H^\ast$ by
$$\es_i(\sum_{j=1}^na_jE_{j,j})=a_i\qquad\for\;\;\sum_{j=1}^na_jE_{j,j}\in
H.\eqno(1.3.31)$$ Since
$$[\sum_{r=1}^na_rE_{r,r},E_{i,j}]=(a_i-a_j)E_{i,j}\qquad\for\;\;\sum_{r=1}^na_rE_{r,r}\in
H,\eqno(1.3.32)$$ we have
$$[h,E_{i,j}]=(\es_i-\es_j)(h)E_{i,j}\qquad\for\;\;h\in
H,\;i,j\in\ol{1,n}.\eqno(1.3.33)$$

Let ${\cal I}$ be a nonzero ideal of $sl(n,\mbb{F})$. Since
$\{0,\es_i-\es_j\mid i\neq j\}$ are distinct linear functions on
$H$, we have some $E_{i,j}\in {\cal I}$ with $i\neq j$ or $0\neq
h\in {\cal I}\bigcap H$ by Lemma 1.3.2, (1.3.30), (1.3.31) and
(1.3.33). Since the center
$$Z(gl(n,\mbb{F}))=\mbb{F}\sum_{i=1}^nE_{i,i}\not\subset
sl(n,\mbb{F}),\eqno(1.3.34)$$ the case $0\neq h\in {\cal I}\bigcap
H$ imply that there exists some $E_{r,s}$ such that
$$0\neq [h,E_{r,s}]=(\es_r-\es_s)(h)E_{r,s}\in {\cal
I}.\eqno(1.3.35)$$ Thus we can always assume some $E_{i,j}\in{\cal
I}$ with $i\neq j$.

For any $i,j\neq l\in\ol{1,n}$, we have
$$[E_{l,i},E_{i,j}]=E_{l,j},\;[E_{i,j},E_{j,l}]=E_{i,l}\in{\cal
I}.\eqno(1.3.36)$$ Furthermore, if $i,j,l\neq k\in\ol{1,n}$,
$$[E_{l,j},E_{j,k}]=E_{l,k}\in{\cal I}.\eqno(1.3.37)$$
Observe
$$[E_{i,j},E_{j,i}]=E_{i,i}-E_{j,j}\in{\cal I}.\eqno(1.3.38)$$
So
$$ E_{j,i}=\frac{1}{2}[E_{j,i},E_{i,i}-E_{j,j}]\in{\cal
I}.\eqno(1.3.39)$$ Similarly as (1.3.36), we have
$$E_{l,i},E_{j,l}\in{\cal I}.\eqno(1.3.40)$$
Thus
$$E_{r,s}\in {\cal I}\qquad\for\;\;r,s\in\ol{1,n},\;r\neq
s.\eqno(1.3.41)$$ Finally,
$$[E_{i,i+1},E_{i+1,i}]=E_{i,i}-E_{i+1,i+1}\in {\cal
I}.\eqno(1.3.42)$$ Hence ${\cal I}=sl(n,\mbb{F})$. Therefore,
$sl(n,\mbb{F})$ is simple.\psp

We leave reader to prove that the symplectic algebra
$sp(2k,\mbb{F})$ is simple when $k>1$.

\section{Homomorphisms}

In this section, we discuss homomorphisms between Lie algebras.

A {\it homomorphism}\index{homomorphism}
 $\nu$ from a Lie algebra ${\msr G}_1$
to a Lie algebra ${\msr G}_2$ is a linear map such that
$$\nu([u,v])=[\nu(u),\nu(v)]\qquad\for\;\;u,v\in{\msr
G}_1.\eqno(1.4.1)$$ Immediately, we have
$$\nu({\msr G}_1)=\{\nu(u)\mid u\in{\msr G}_1\}\;\mbox{forms a
subalgebra of}\;{\msr G}_2.\eqno(1.4.2)$$ Set
$$\kn\nu =\{u\in{\msr G}_1\mid \nu(u)=0\}.\eqno(1.4.3)$$
If $u\in\kn\nu$, then
$$\nu([u,v])=[\nu(u),\nu(v)]=[0,\nu(v)]=0.\eqno(1.4.4)$$
Thus $[u,v]\in\kn\nu$. So $\kn\nu$ is an ideal of ${\msr G}_1$. The
map $\nu$ is called {\it injective} if $\kn\nu=\{0\}$ and surjective
if $\nu({\msr G}_1)={\msr G}_2$. Moreover, $\nu$ is called an {\it
isomorphism}\index{isomorphism} if it is both injective and
surjective (i.e. {\it bijective}). In this case, we say that ${\msr
G}_1$ is {\it isomorphic} to ${\msr G}_2$ via $\nu$, denoted as
$${\msr G}_1\stl{\nu}{\cong}{\msr G}_2.\eqno(1.4.5)$$
Isomorphic Lie algebras are viewed as the same algebra.\psp

{\bf Example 1.4.1}. Note that
$$[\ptl_x,x^2\ptl_x]=2x\ptl_x,\;\;[2x\ptl_x,\ptl_x]=-2\ptl_x,\;\;[2x\ptl_x,x^2\ptl_x]
=2x^2\ptl_x\eqno(1.4.6)$$ Thus
$$\mbb{F}\ptl_x+\mbb{F}x\ptl_x+\mbb{F}x^2\ptl_x\eqno(1.4.7)$$
forms a Lie algebra isomorphic to $sl(2,\mbb{F})$ with the
isomorphism $\nu$ given by:
$$\nu(\ptl_x)=-E_{2,1},\;\;\nu(2x\ptl_x)=E_{1,1}-E_{2,2},\;\;\nu(x^2\ptl_x)
=E_{1,2}.\eqno(1.4.8)$$

In general, the space
$$\sum_{i,j=1}^n\mbb{F}x_i\ptl_{x_j}\eqno(1.4.9)$$
forms a Lie algebra isomorphic to $gl(n,\mbb{F})$ with the
isomorphism $\tau$ defined by:
$$\tau(x_i\ptl_j)=E_{i,j}\qquad\for\;\;i,j\in\ol{1,n}.\eqno(1.4.10)$$
\vspace{0.1cm}

{\bf Example 1.4.2}. The symplectic Lie
algebra\begin{eqnarray*}\hspace{1cm}sp(2n,\mbb{F})&=&
\sum_{i,j=1}^n\mbb{F}(E_{i,j}-E_{n+j,n+i})+\sum_{i=1}^n(\mbb{F}E_{i,n+i}+\mbb{F}E_{n+i,i})\\
& &+\sum_{1\leq i<j\leq n
}[\mbb{F}(E_{i,n+j}+E_{j,n+i})+\mbb{F}(E_{n+i,j}+E_{n+j,i})].\hspace{2cm}(1.4.11)\end{eqnarray*}
View the elements of $\mbb{F}[x_1,x_2,...,x_n]$ also as the
multiplication operators on $\mbb{F}[x_1,x_2,...,x_n]$. We have the
following Lie subalgebra
$$\sum_{i,j=1}^n(\mbb{F}(x_i\ptl_{x_j}+\dlt_{i,j}/2)+\mbb{F}x_ix_j
+\mbb{F}\ptl_{x_i}\ptl_{x_j})\eqno(1.4.12)$$ of
$\Edo\mbb{F}[x_1,x_2,...,x_n]$, which is isomorphic to
$sp(2n,\mbb{F})$ with the isomorphism $\nu$ defined by
$$\nu(x_i\ptl_{x_j}+\dlt_{i,j}/2)=E_{i,j}-E_{n+j,n+i},\eqno(1.4.13)$$
$$\nu(\ptl_{x_i}\ptl_{x_j})=E_{n+i,j}+E_{n+j,i},\;\;\nu(x_ix_j)=
-(E_{i,n+j}+E_{j,n+i})\eqno(1.4.14)$$ for $i,j\in\ol{1,n}$. \psp

A homomorphism from a Lie algebra ${\msr G}$ to $gl(V))$ is called a
{\it representation}\index{representation} on $V$. The isomorphisms
in the above examples are representations of the Lie algebras on
polynomial algebras. For any element $u$ in  Lie algebra ${\msr G}$,
we define $\ad u\in gl({\msr G})$ by
$$ \ad u(v)=[u,v]\qquad\for\;\;v\in {\msr G}.\eqno(1.4.15)$$
Then Jacobi identity implies
$$\ad u([v,w])=[\ad u(v),w]+[v,\ad u(w)]\qquad\for\;\;v,w\in{\msr
G}.\eqno(1.4.16)$$ Thus $\ad u$ is a derivation of ${\msr G}$, which
is called an {\it inner derivation}\index{inner derivation}
 of ${\msr G}$. Moreover,
\begin{eqnarray*}[\ad u,\ad v](w)&=& (\ad u\:\ad v-\ad v\:\ad
u)(w)=\ad(u)([v,w])-\ad v([u,w])\\
&=&[u,[v,w]]-[v,[u,w]] =[u,[v,w]]+[v,[w,u]]\\ &
=&-[w,[u,v]]=[[u,v],w]=\ad
[u,v](w).\hspace{4.2cm}(1.4.17)\end{eqnarray*} Hence
$$[\ad u,\ad v]=\ad [u,v].\eqno(1.4.18)$$
Thus $\ad$ is a representation of ${\msr G}$ on ${\msr G}$, which is
called the {\it adjoint representation} of ${\msr G}$. Furthermore,
$$ \kn\ad=Z({\msr G}).\eqno(1.4.19)$$
We leave the reader to show that $\ad {\msr G}$ is an ideal of
$\der{\msr G}$.

Suppose that ${\cal I}$ is an ideal of a Lie algebra ${\msr G}$.
Define the cosets
$$ u+{\cal I}=\{u+v\mid v\in {\cal I}\}\qquad\for\;\;u\in{\msr
G}.\eqno(1.4.20)$$ The quotient space
$${\msr G}/{\cal I}=\{u+{\cal I}\mid u\in{\msr
G}\}.\eqno(1.4.21)$$ Define an algebraic operation $[\cdot,\cdot]$
on ${\msr G}/{\cal I}$ by
$$ [u+{\cal I},v+{\cal I}]=[u,v]+{\cal I}\qquad\for\;\;u,v\in{\msr
G},\eqno(1.4.22)$$ which is independent of the choices of
representative elements in the cosets. Then $({\msr G}/{\cal
I},[\cdot,\cdot])$ forms a Lie algebra. If $\nu$ is a homomorphism
from a Lie algebra ${\msr G}_1$ to a Lie algebra ${\msr G}_2$, then
$${\msr G}_1/\kn \nu\cong \nu({\msr G}_1).\eqno(1.4.23)$$
Let $\msr L$ be a subalgebra of a Lie algebra ${\msr G}$ and let
${\cal I}$ be an ideal of ${\msr G}$. Then $\msr L+{\cal I}$ forms a
Lie subalgebra and ${\cal I}$ is an ideal of $\msr L+{\cal I}$.
Moreover,
$$(\msr L+{\cal I})/{\cal I}\cong \msr L/(L\bigcap {\cal I}).\eqno(1.4.24)$$
\vspace{0.1cm}

{\bf Example 1.4.3}. In Example 1.3.1,
$${\msr K}_n/{\cal I}_k\cong {\msr K}_{n-k-1}\eqno(1.4.25)$$
through the isomorphism:
$$ \ptl_{x_1}+{\cal I}_k\rta
\ptl_{x_1},\;\;\frac{1}{i!}x^i_1\ptl_{x_2}+{\cal I}_k\rta
\frac{1}{(i-k-1)!}x^{i-k-1}_1\ptl_{x_2},\;\;i\in\ol{k+1,n}.\eqno(1.4.26)$$

An isomorphism from a Lie algebra ${\msr G}$ to itself is called an
{\it automorphism}. Suppose that $\ptl$ is a derivation of a Lie
algebra ${\msr G}$ and $e^\ptl$ is well defined (e.g., $\ptl$ is
nilpotent). Note
\begin{eqnarray*}e^\ptl([u,v])&=&\sum_{n=0}^{\infty}\frac{1}{n!}\ptl^n([u,v])
=\sum_{n=0}^{\infty}\frac{1}{n!}\sum_{i=0}^n
\left(\!\!\!\begin{array}{c}n\\
i\end{array}\!\!\!\right)[\ptl^i(u),\ptl^{n-i}(v)]\\
&=&\sum_{n=0}^{\infty}\frac{1}{n!}\sum_{i=0}^n
\frac{n!}{i!(n-i)!}[\ptl^i(u),\ptl^{n-i}(v)]=\sum_{n=0}^{\infty}
\sum_{i=0}^n\frac{1}{i!(n-i)!}[\ptl^i(u),\ptl^{n-i}(v)]
\\
&=&[\sum_{i=0}^{\infty}\frac{1}{i!}\ptl^i(u),\sum_{j=0}^{\infty}\frac{1}{j!}\ptl^j(v)]
=[e^{\ptl}(u),e^{\ptl}(v)].\hspace{4.8cm}(1.4.27)\end{eqnarray*} So
$e^\ptl$ is an automorphism of ${\msr G}$.
 If $e^{\ad u}$ is defined for $u\in
{\msr G}$, we called $e^{\ad u}$ an {\it inner
automorphism}.\index{inner automorphism}

Let ${\msr A}$ be an associative algebra and let $u$ be an element
in ${\msr A}$ such that $e^u$ is well defined (e.g., $u$ is
nilpotent). Define the left multiplication operator $L_u$ and the
right multiplication operator $R_u$ by
$$L_u(v)=uv,\qquad R_u(v)=vu\qquad\for\;\;v\in{\msr
A}.\eqno(1.4.28)$$ Thus
$$\ad u(v)=uv-vu=(L_u-R_u)(v)\qquad\for\;\;v\in{\msr
A}.\eqno(1.4.29)$$ Moreover,
$$L_uR_u=R_uL_u\lra e^{\ad
u}=e^{L_u-R_u}=e^{L_u}e^{-R_u}.\eqno(1.4.30)$$ Thus
$$e^{\ad u}(v)=e^{L_u}e^{-R_u}(v)=e^uve^{-u}\eqno(1.4.31)$$
( A very useful formula!).

Let $n=2k$ be an even integer. We define a map $\ast: M_{n\times
n}(\mbb{F}) \rta M_{n\times n}(\mbb{F})$ by
$$ A^\ast=\left(\begin{array}{cc}0&I_k\\
I_k& 0\end{array}\right)A^t\left(\begin{array}{cc} 0&I_k\\
I_k&0\end{array}\right)\qquad\for\;\;A\in M_{n\times
n}(\mbb{F}).\eqno(1.4.32)$$ Then $\ast$ is an involutive
anti-isomorphism. Thus we have a Lie subalgebra
\begin{eqnarray*}o'(2k,\mbb{F})&=&\{A\in M_{n\times n}(\mbb{F})\mid
A^\ast=-A\}\\ &=&\left\{\left(\begin{array}{cc}A&B\\
 C&-A^t\end{array}\right)\mid A,B,C\in M_{k\times
k}(\mbb{F}),\;B^t=-B,C^t=-C\right\},\hspace{1.2cm}(1.4.33)\end{eqnarray*}
 which is simple
when $k>2$, and is called a Lie algebra of type $D_k$ for $k>3$. The
algebra  $o'(2k,\mbb{R})$ is the Lie algebra of the Lie group
preserving the metric defined by the symmetric matrix
$$\left(\begin{array}{cc}0&I_k\\
I_k& 0\end{array}\right).\eqno(1.4.34)$$

Suppose $n=2k+1$. We define a map $\ast: M_{n\times n}(\mbb{F}) \rta
M_{n\times n}(\mbb{F})$ by
$$ A^\ast=\left(\begin{array}{ccc}1&0&0\\ 0&0&I_k\\ 0
&I_k&0\end{array}\right)A^t\left(\begin{array}{ccc}1&0&0\\
0&0&I_k\\ 0 &I_k&0\end{array}\right)\qquad\for\;\;A\in M_{n\times
n}(\mbb{F}).\eqno(1.4.35)$$ Then $\ast$ is an involutive
anti-isomorphism. Thus we have a Lie subalgebra
\begin{eqnarray*} & &o'(2k+1,\mbb{F})=\{A\in M_{n\times n}(\mbb{F})\mid
A^\ast=-A\}=\{\left(\begin{array}{ccc}0&-\vec b^t& -\vec a^t\\
\vec a&A&B\\ \vec b& C&-A^t\end{array}\right)\\ & &\mid\vec a,\vec
b\in M_{k\times 1}(\mbb{F});\;A,B,C\in M_{k\times
k}(\mbb{F}),\;B^t=-B,C^t=-C\},\hspace{3.2cm}(1.4.36)\end{eqnarray*}
is simple when $k>1$, and is called a Lie algebra of type $B_k$ for
$k>1$. The algebra $o'(2k+1,\mbb{R})$ is the Lie algebra of the Lie
group preserving the metric defined by the symmetric matrix:
$$\left(\begin{array}{ccc}1&0&0\\ 0&0&I_k\\
0&I_k& 0\end{array}\right).\eqno(1.4.37)$$

Let ${\msr G}_1$ and ${\msr G}_2$ be Lie algebras. The Lie bracket
on
$${\msr G}_1\oplus {\msr G}_2=\{(u_1,u_2)\mid u_1\in {\msr G}_1,
\;u_2\in {\msr G}_2\}\eqno(1.4.38)$$ is defined by
$$[(u_1,u_2),(v_1,v_2)]=([u_1,v_1],[u_2,v_2])\qquad  u_1,v_1\in {\msr G}_1,
\;u_2,v_2\in {\msr G}_2.\eqno(1.4.39)$$
 We leave the reader to prove
$$o(2k,\mbb{C})\cong o'(2k,\mbb{C}),\qquad o(2k+1,\mbb{C})\cong
o'(2k+1,\mbb{C}).\eqno(1.4.40)$$ Thereby, prove
$$ o(3,\mbb{C})\cong sl(2,\mbb{C}),\;\;o(4,\mbb{C})\cong
sl(2,\mbb{C})\oplus sl(2,\mbb{C}),\;\;o(6,\mbb{C})\cong
sl(4,\mbb{C}).\eqno(1.4.41)$$

\section{Nilpotent and solvable Lie Algebras}

In this section, we define nilpotent and solvable Lie algebras.
Moreover, we prove the Engel's theorem, Lie's theorem,
Jordan-Chevalley decomposition theorem and Cartan's criterion of
solvability.

Let ${\msr G}$ be a Lie algebra. The {\it central
series}\index{central series} of ${\msr G}$ is
defined as: ${\msr G}^0={\msr G},\;{\msr G}^1=[{\msr G},{\msr G}],\\
{\msr G}^{i+1}=[{\msr G},{\msr G}^i]$. Then
$${\msr G}^0\supset {\msr G}^1\supset {\msr G}^2\supset
\cdots\supset {\msr G}^i\supset\cdots\eqno(1.5.1)$$ forms a
descending series of ideals of ${\msr G}$, invariant under any
automorphism $\tau$; that is,
$$\tau({\msr G}^i)={\msr G}^i.\eqno(1.5.2)$$
If ${\msr G}^n=\{0\}$ for some positive integer $n$, then we say
that ${\msr G}$ is {\it nilpotent}.\index{nilpotent} For any vector
space $V$, we define a trivial Lie bracket on $V$ by
$$[u,v]=0\qquad\for\;\;u,v\in V.\eqno(1.5.3)$$
We call $(V,[\cdot,\cdot])$ a {\it trivial Lie
algebra}\index{trivial Lie algebra}  ~or {\it abelian Lie
algebra}.\index{abelian Lie algebra}
 It is obviously a nilpotent Lie algebra.

\psp

{\bf Example 1.5.1}. The Heisenberg Lie algebra ${\msr H}_n$ is a
vector space with a basis $\{u_i,v_i,\kappa\mid i=1,2,...,n\}$ whose
Lie bracket is determined  by
$$
[u_i,v_j]=\dlt_{i,j}\kappa,\;\;[u_i,u_j]=[v_i,v_j]=[\kappa,u_i]=[\kappa,v_i]=0.\eqno(1.5.4)$$
This algebra is the symmetry of harmonic oscillators in physics.
Note
$${\msr H}^1=\mbb{F}\kappa,\qquad {\msr
H}_n^2=\{0\}.\eqno(1.5.5)$$ So it is a nilpotent Lie algebra.\psp

{\bf Example 1.5.2}. Recall the Lie algebra ${\msr K}_n$ defined in
Example 1.3.1. We have
$${\msr K}_n^i={\cal I}_{n-i}\eqno(1.5.6)$$
(cf. (1.3.7)). In particular,
$${\msr K}^n=\mbb{F}\ptl_{x_2},\qquad{\msr
K}^{n+1}=\{0\}.\eqno(1.5.7)$$ Thus ${\msr K}_n$ is a nilpotent Lie
algebra.\psp

{\bf Example 1.5.3}. The subspace of strict up-triangular matrices
$${\msr
N}_n=\left\{\left(\begin{array}{ccccc}0&a_{1,2}&a_{1,3}&\cdots
&a_{1,n}\\ 0& 0&a_{2,3}&\ddots&\vdots\\
\vdots&\ddots&\ddots&\ddots&a_{n-2,n}\\ 0&\cdots &0&0&a_{n-1,n}\\
0&\cdots &0&0&0\end{array}\right)\mid a_{i,j}\in \mbb{F}\right\}=
\sum_{1\leq i<j\leq n}\mbb{F}E_{i,j}\eqno(1.5.8)$$ forms a Lie
subalgebra of $gl(n,\mbb{F})$. Set
$${\msr E}_r=\sum_{j-i>r}\mbb{F}E_{i,j}.\eqno(1.5.9)$$
Then
$${\msr N}_n^r={\msr E}_r.\eqno(1.5.10)$$
In particular,
$${\msr N}_n^{n-1}=\mbb{F}E_{1,n},\;\;{\msr
N}^n=\{0\}.\eqno(1.5.11)$$ Thus ${\msr N}_n$ is a nilpotent Lie
algebra. \psp

{\bf Theorem 1.5.1}. {\it Let $V\neq \{0\}$ be a finite-dimensional
space and let ${\msr G}$ be a subalgebra of $gl(V)$, consisting of
nilpotent linear transformations. There exists nonzero $v\in V$ such
that ${\msr G}(v)=\{0\}$.}

{\it Proof}. We prove the conclusion by induction on $\dim {\msr
G}$. If $\dim {\msr G}=1$, we write ${\msr G}=\mbb{F}T$. Take any
eigenvector $v$ of $T$, we have ${\msr G}(v)=\{0\}$.

Assume that the theorem holds for $\dim {\msr G}\leq k$. Consider
the case $\dim {\msr G}=k+1$. Take a maximal proper Lie subalgebra
${\msr K}$ of ${\msr G}$. For $u\in {\msr G}$, we want to prove that
$\ad u$ is nilpotent on $gl(V)$. Note that
$$(\ad
u)^n=(L_u-R_u)^n=\sum_{i=0}^n{n\choose
i}(-1)^iL_u^{n-i}R_u^i=\sum_{i=0}^n{n\choose
i}(-1)^iL_{u^{n-i}}R_{u^i}\eqno(1.5.12)$$ by (1.4.28) and (1.4.29).
If $u^m=0$, then $(\ad u)^{2m-1}=0$ by the above equation.

Define the action of ${\msr K}$ on the quotient space ${\msr
G}/{\msr K}$ by
$$\ad\xi(u+{\msr K})=[\xi,u]+{\msr K}\qquad\for\;\;\xi\in{\msr K},\;u\in{\msr
G}.\eqno(1.5.13)$$ Then
$$\ad {\msr K}|_{{\msr G}/{\msr K}}\eqno(1.5.14)$$
forms a Lie subalgebra of $gl({\msr G}/{\msr K})$. If $\xi^m=0$,
then $(\ad\xi)^{2m-1}=0$ as we have shown in the above. Furthermore,
$$\dim \ad {\msr K}|_{{\msr G}/{\msr K}}\leq\dim {\msr K}\leq
k.\eqno(1.5.15)$$ Thus there exists $\xi_0\in {\msr G}\setminus
{\msr K}$ such that
$$\ad\xi (\xi_0+{\msr K})={\msr K}\qquad\for\;\;\xi\in{\msr
K},\eqno(1.5.16)$$ which implies
$$[\xi,\xi_0]\in {\msr K}\qquad\for\;\;\xi\in{\msr
K}.\eqno(1.5.17)$$ Hence $\mbb{F}\xi_0+{\msr K}$ forms a Lie
subalgebra of ${\msr G}$. Since ${\msr K}$ is maximal, we have
$${\msr G}=\mbb{F}\xi_0+{\msr K}.\eqno(1.5.18)$$

Set
$$V_1=\{u\in V\mid {\msr K}(u)=\{0\}\}.\eqno(1.5.19)$$
Then $V_1\neq \{0\}$ by induction. For any $u\in V_1$ and $\xi\in
{\msr K}$, we have
$$\xi \xi_0 (u)=[\xi,\xi_0](u)+\xi_0\xi(u)=0\eqno(1.5.20)$$
by (1.5.16); that is, $\xi_0 (u)\in V_1$. Hence $\xi_0(V_1)\subset
V_1$. Take any eigenvector $v$ of $\xi_0|_{V_1}$ and we have ${\msr
G}(v)=\{0\}.\qquad\Box$\psp

The above theorem and induction show that ${\msr G}$ is isomorphic
to a Lie algebra of strict up-triangular matrices (cf. (1.5.8)).
Hence ${\msr G}$ is a nilpotent Lie algebra.\psp

{\bf Theorem 1.5.2 (Engel)}.\index{Engel's Theorem} {\it If all the
elements in a finite-dimensional Lie algebra ${\msr G}$ are ad
nilpotent, then ${\msr G}$ is nilpotent}.

{\it Proof}. By Theorem 1.5.1, $\ad {\msr G}$ is nilpotent. Hence
there exists a positive integer $m$ such that $(\ad {\msr
G})^m=\ad{\msr G}^m=\{0\}$. So ${\msr G}^m\in Z({\msr G})$.
Therefore, ${\msr G}^{m+1}=[{\msr G},{\msr G}^m]\subset [{\msr
G},Z({\msr G})]=\{0\}.\qquad\Box$\psp

{\bf Corollary 1.5.3}. {\it Let ${\msr G}$ be a finite-dimensional
nilpotent Lie algebra and let ${\cal I}$ be a nonzero ideal of
${\msr G}$. Then ${\cal I}\bigcap Z({\msr G})\neq \{0\}$.}

{\it Proof}. Note that $\ad{\msr G}|_{\cal I}$ is a Lie algebra of
nilpotent elements. Thus there exists $0\neq v\in {\cal I}$ such
that $[{\msr G},v]=\{0\}$. So $0\neq v\in {\cal I}\bigcap Z({\msr
G}).\qquad\Box$\psp

 The {\it derived series}\index{derived series} of a Lie algebra ${\msr G}$ is defined as: ${\msr
G}^{(0)}={\msr G},\;{\msr G}^{(1)}=[{\msr G},{\msr G}],\;{\msr
G}^{(i+1)}=[{\msr G}^{(i)},{\msr G}^{(i)}]$.  The Lie algebra ${\msr
G}$ is called {\it solvable}\index{solvable} if ${\msr
G}^{(n)}=\{0\}$ for some positive integer $n$. Since ${\msr
G}^{(i)}\subset {\msr G}^i$, any nilpotent Lie algebra is
solvable.\psp

{\bf Properties:}

{\it (1) Let ${\cal I}$ be an ideal of ${\msr G}$. Then ${\msr G}$
is solvable if and only if ${\cal I}$ and ${\msr G}/{\cal I}$ are
solvable}.

{\it Proof}. Suppose that ${\cal I}$ and ${\msr G}/{\cal I}$ are
solvable. There exist positive integer $m$ and $n$ such that
$${\cal I}^{(m)}=\{0\},\qquad ({\msr G}/{\cal
I})^{(n)}=\{0\}.\eqno(1.5.21)$$ Thus ${\msr G}^{(n)}\subset {\cal
I}$ and ${\msr G}^{(m+n)}\subset {\cal I}^{(m)}=\{0\}$.

(2) {\it If ${\cal I}$ and ${\cal J}$ are solvable ideals of ${\msr
G}$, so is ${\cal I}+{\cal J}$}.

{\it Proof}. This is because $({\cal I}+{\cal J})/{\cal J}\cong
{\cal I}/({\cal I}\bigcap{\cal J})$ and  ${\cal I}\bigcap{\cal J}$
are solvable.\psp

We define the {\it radical}\index{radical!of Lie algebra} of ${\msr
G}$ to be the unique maximal solvable ideal, denoted by
$\mbox{Rad}\:{\msr G}$. The Lie algebra ${\msr G}$ is called {\it
semisimple}\index{semisimple} of $\mbox{Rad}\:{\msr G}=\{0\}$.\psp

{\bf Example 1.5.4}.  The {\it Diamond Lie algebra}\index{Diamond
Lie algebra} is a vector space ${\msr G}_\diamondsuit$ with a basis
$\{h,u_+,u_-,\kappa\}$ whose Lie bracket is given by
$$[h,u_{\pm}]=\pm u_{\pm},\;\;[u_+,u_-]=\kappa,\;\;[{\msr G}_\diamondsuit,\kappa]=\{0\}.\eqno(1.5.22)$$ Then
$${\msr G}_\diamondsuit^{(1)}=[{\msr G}_\diamondsuit,{\msr G}_\diamondsuit]=\mbb{F}u_++\mbb{F}u_-+\mbb{F}\kappa,
\eqno(1.5.23)$$
$${\msr G}_\diamondsuit^{(2)}=\mbb{F}\kappa,\qquad {\msr G}_\diamondsuit^{(3)}=\{0\}.\eqno(1.5.24)$$ So ${\msr G}_\diamondsuit$ is a solvable Lie
algebra. Moreover, we define a symmetric bilinear form on ${\msr
G}_\diamondsuit$ by
$$\la u_+,u_-\ra=\la h,\kappa\ra=1,\eqno(1.5.25)$$
$$\la u,v\ra=0\qquad \mbox{for the other
pairs}\;u,v\in\{h,u_+,v_-,\kappa\}.\eqno(1.5.26)$$ Then
$\la\cdot,\cdot\ra$ is invariant:
$$\la [u,v],w\ra=\la u,[v,w]\ra\qquad\for\;\;u,v,w\in{\msr
G}\eqno(1.5.27)$$ (exercise). Thus the diamond Lie algebra ${\msr
G}_\diamondsuit$ is a solvable Lie algebra with a nondegenerate
symmetric invariant bilinear form.\psp

{\bf Example 1.5.5}. The Lie algebra of up-triangular matrices:
$${\msr
B}_n=\left\{\left(\begin{array}{ccccc}a_{1,1}&a_{1,2}&a_{1,3}&\cdots&
a_{1,n}\\ 0& a_{2,2}&a_{2,3}&\ddots&\vdots\\
\vdots&\ddots&\ddots&\ddots&a_{n-2,n}\\ 0&\cdots &0&a_{n-1,n-1}&a_{n-1,n}\\
0&\cdots &0&0&a_{n,n}\end{array}\right)\mid a_{i,j}\in
\mbb{F}\right\}= \sum_{1\leq i\leq j\leq
n}\mbb{F}E_{i,j}\eqno(1.5.28)$$ forms a Lie subalgebra of
$gl(n,\mbb{F})$. Moreover, $[{\msr B}_n,{\msr B}_n]={\msr N}_n$ is a
nilpotent ideal (cf. Example 1.5.3). Thus ${\msr B}_n$ is a solvable
Lie algebra. \psp

{\bf Theorem 1.5.4 (Lie)}.\index{Lie's Theorem} {\it Suppose that
$\mbb{F}=\mbb{C}$ and $V$ is a nonzero finite-dimensional space. Let
${\msr G}$ be a solvable algebra of $gl(V)$. There exists a common
eigenvector for all the elements in ${\msr G}$. In particular,
${\msr G}$ is isomorphic to a Lie algebra of up-triangular
matrices.}

{\it Proof}. We prove it by induction on $\dim{\msr G}$. When
$\dim{\msr G}=1$, ${\msr G}=\mbb{F}T$. Take any eigenvector $v$ of
$T$. Then $v$ is a common eigenvector all the elements in ${\msr
G}$.

 Assume that the theorem hold for $\dim {\msr G}\leq k$. Consider the case
 $\dim{\msr G}=k+1$. Since ${\msr G}$ is solvable, $[{\msr
 G},{\msr G}]={\msr G}^{(1)}\neq {\msr G}$. Take any codimension one subspace
 ${\msr K}\supset {\msr G}^{(1)}$ of ${\msr G}$. Then ${\msr K}$
 is a solvable ideal of ${\msr G}$. By induction, there exists a
 common eigenvector $v$ for all the elements in ${\msr K}$. Hence
 there exists a linear function $\al$ on ${\msr K}$ such that
 $$ \xi(v)=\al(\xi)v\qquad\for\;\;\xi\in{\msr K}.\eqno(1.5.29)$$
 Therefore, the subspace
 $$V_\al=\{v\in V\mid \xi(v)=\al(\xi)v\;\for\;\xi\in{\msr
 K}\}\neq\{0\}.\eqno(1.5.30)$$

Write ${\msr G}=\mbb{F}\xi_0+{\msr K}$. For any $0\neq w\in V_\al$,
we have
$$
\xi\xi_0(w)=[\xi,\xi_0](w)+\xi_0\xi(w)=\al(\xi)\xi_0(w)+\al([\xi,\xi_0])w.\eqno(1.5.31)$$
We want to prove $\al([\xi,\xi_0])=0$. If $\xi_0(w)\in\mbb{F}w$, we
are done. Otherwise, since $\dim V$ is finite, there exists a
positive integer $n$ such that
$\{w,\xi_0(w),\xi_0^2(w),...,\xi_0^n(w)\}$ are linear independent
and
$$ \xi_0^{n+1}(w)\in W=\sum_{i=0}^n\mbb{F}\xi_0^i(w).\eqno(1.5.32)$$
Thus $\xi_0(W)\in W$. Set
$$W_i=\sum_{r=0}^i\mbb{F}\xi_0^r(w).\eqno(1.5.33)$$
For any $\xi\in{\msr K}$, we have $\xi(w)=\al(\xi)w$. Suppose
$$\xi(\xi_0^i(w))\equiv
\al(\xi)\xi_0^i(w)\;\;(\mbox{mod}\:W_{i-1})\eqno(1.5.34)$$ for some
$i<n$. Now
\begin{eqnarray*}\hspace{2cm}&&\xi(\xi_0^{i+1}(w))=\xi\xi_0(\xi_0^i(w))=[\xi,\xi_0](\xi_0^i(w))+\xi_0\xi(\xi_0^i(w))
\\ &\equiv&
\al([\xi,\xi_0])\xi_0^i(w)+\al(\xi)\xi_0^{i+1}(w)\;\;(\mbox{mod}\:\xi_0
(W_{i-1})+W_{i-1})\\
&\equiv&\al(\xi)\xi_0^{i+1}(w)\;\;(\mbox{mod}\:W_i).\hspace{7.4cm}(1.5.35)\end{eqnarray*}
By induction, (1.5.34) holds for any $i\leq n$. Thus $\xi (W)\subset
W$ and
$$\tr \xi|_W=(n+1)\al(\xi)\qquad \mbox{for any}\;\;\xi\in{\msr
K}.\eqno(1.5.36)$$ Since $[\xi,\xi_0]\in{\msr K}$, we have
$$\tr [\xi,\xi_0]|_W=(n+1)\al([\xi,\xi_0])=\tr
[\xi|_W,\xi_0|_W]=0;\eqno(1.5.37)$$ that is, $\al([\xi,\xi_0])=0$.
Therefore,
$$\xi(\xi_0 w)=\al(\xi)\xi_0 w\eqno(1.5.38)$$
for $\xi\in{\msr K}$ by (1.5.31). So $\xi_0 (V_\al)\subset V_\al$.
Any eigenvector of $\xi_0|_{V_\al}$ is a common eigenvector of the
elements in ${\msr G}.\qquad\Box$ \psp

The above theorem shows that ${\msr G}$ becomes a Lie subalgebra of
up-triangular matrices under a certain basis. Thus $[{\msr G},{\msr
G}]$ is a Lie algebra of strict up-triangular matrices. Hence
$[{\msr G},{\msr G}]$ is nilpotent.

Let $V$ be a vector space with a basis $\{v_i\mid i\in\mbb{Z}\}$.
Define $T\in gl(V)$ by
$$ T(v_i)=v_{i+1}\qquad\for\;\;i\in\mbb{Z}.\eqno(1.5.39)$$
Then $T$ does not have any eigenvector (exercise). So Lie's theorem
fail if we allow $\dim V=\infty$.

\psp

{\bf Corollary 1.5.5}. {\it If ${\msr G}$ is a finite-dimensional
solvable Lie algebra over $\mbb{C}$. Then $[{\msr G},{\msr G}]$ is
nilpotent}.

{\it Proof}. Note $\ad {\msr G}$ is a solvable Lie subalgebra of
$gl({\msr G})$. Thus
$$\ad [{\msr G},{\msr G}]=[\ad{\msr G},\ad{\msr G}]\eqno(1.5.40)$$
is a nilpotent Lie subalgebra. Since $\kn\ad=Z({\msr G})$, $[{\msr
G},{\msr G}]$ is nilpotent.$\qquad\Box$\psp

 An element $T\in \Edo V$ is called {\it semisimple}\index{semisimple transformation} if it can be
diagonalized under a certain basis of $V$.\psp

{\bf Theorem 1.5.6 (Jordan-Chevalley
Decomposition)}.\index{Jordan-Chevalley Decomposition} {\it Suppose
$\mbb{F}=\mbb{C}$. Let $V$ be a finite-dimensional vector space.
Given $T\in \Edo V$. There exist unique semisimple element $T_s$ and
nilpotent element $T_n$ in $\Edo V$ such that $T=T_s+T_n$ and
$[T_s,T_n]=0$. Moreover, $T_s$ and $T_n$ can be written as
polynomials of $T$ without constant.}

{\it Proof}. Suppose that $\{\lmd_1,\lmd_2,...,\lmd_k\}$ are
distinct eigenvalues of $T$. Then the characteristic polynomial
 is of the form $\prod_{i=1}^k(t-\lmd_i)^{n_i}$. Set
 $$V_i=\{u\in V\mid (T-\lmd_i)^{n_i}(u)=0\}.\eqno(1.5.41)$$
 We have
 $$V=\bigoplus_{i=1}^kV_i.\eqno(1.5.42)$$

 Suppose that $0$ is an eigenvalue of $T$. Set
 $$ f_i(t)=\sum_{i\neq j}(t-\lmd_j)^{n_j}.\eqno(1.5.43)$$
 Then the maximal common divisor of $\{f_1(t),f_2(t),...,f_k(t)\}$
 is 1. By long division, there exists polynomials
 $\{h_1(t),h_2(t),...,h_k(t)\}$ such that
 $$ \sum_{i=1}^th_i(t)f_i(t)=1.\eqno(1.5.44)$$
 Hence
 $$h_i(t)f_i(t)\equiv
 \dlt_{i,j}\;\;(\mbox{mod}\:(t-\lmd_j)^{n_j}).\eqno(1.5.45)$$
 Let
 $$P(t)=\sum_{i=1}^k\lmd_ih_i(t)f_i(t).\eqno(1.5.46)$$
 By (1.5.46), we get
 $$P(t)\equiv \lmd_i\;\;(\mbox{mod}\:(t-\lmd_i)^{n_i}).\eqno(1.5.47)$$
If $0$ is not an eigenvalue, then we set
$$f_0(t)=\prod_{j=1}^k(t-\lmd_j)^{n_j},\;\;f_i(t)=t\sum_{i\neq
j}(t-\lmd_j)^{n_j},\qquad 0<i\leq k.\eqno(1.5.48)$$ By long
division, there exists polynomials
 $\{h_0,h_1(t),h_2(t),...,h_k(t)\}$ such that
 $$ \sum_{i=0}^kh_i(t)f_i(t)=1.\eqno(1.5.49)$$
 Now (1.5.45) holds. We define $P(t)$ by (1.5.46) again. Equation
 (1.5.48) holds. In any case, we always have $P(0)=0$.

 Set
 $$ T_s=P(T),\qquad T_n=T-P(T).\eqno(1.5.50)$$
 Note $P(T)(V_i)\subset V_i$ and
 $$T_s|_{Vi}=\lmd_i\mbox{Id}_{V_i}.\eqno(1.5.51)$$
 Denote
 $$m=\mbox{max}\{n_1,n_2,...,n_k\}.\eqno(1.5.52)$$
 Then
 $$T_n^m=(T-T_s)^m=0\eqno(1.5.53)$$
 by (1.5.41), (1.5.42) and (1.5.51).

 Suppose that $T=A+B$ such that $A$ is a semisimple element in
 $gl(V)$, $B$ is a nilpotent element in $gl(V)$ and $AB=BA$. Note
 $AT=A(A+B)=(A+B)A=TA$, which implies $BT=TB$. Thus
 $T_sA=AT_s$ and $T_nB=BT_n$.  Since both $T_s$ and $A$ are
 diagonalizable and they commute, $T_s-A$ is semsimple. On the
 other hand, $B-T_n$ is nilpotent. Thus
 $$T_s+T_n=A+B\lra T_s-A=B-T_n\lra  T_s-A=B-T_n=0,\eqno(1.5.54)$$
 because it is both semisimple and nilpotent.$\qquad\Box$
\psp

 {\bf Corollary 1.5.7}. {\it If $T=T_s+T_n$ is the Jordan-Chevalley
decomposition, so is $\ad T=\ad T_s+\ad T_n$}.

 {\it Proof}. This is because $\ad T_s$ is semisimple, $\ad T_n$
 is nilpotent and $[\ad T_s,\ad T_n]=\ad [T_s,T_n]=0.\qquad\Box$
\psp

{\bf Lemma 1.5.8}. {\it If $({\msr A},\circ)$ is any
finite-dimensional algebra, the semisimple and nilpotent parts of an
element in $\der {\msr A}$ are also in $\der {\msr A}$.}

{\it Proof}. Let $T\in\der {\msr A}$. For any $\lmd\in\mbb{F}$, we
define
$$ {\msr A}_\lmd=\{u\in {\msr A}\mid (T-\lmd)^n(u)=0\;\mbox{for
some positive integer}\;n\}.\eqno(1.5.55)$$ Then
$${\msr A}=\bigoplus{\msr A}_\lmd,\qquad
T_s|_{{\msr A}_\lmd}=\lmd\mbox{Id}_{{\msr A}_\lmd}.\eqno(1.5.56)$$
For $u\in \msr A_\lmd$ and $v\in \msr A_\mu$, there exist positive
integers $m$ and $n$ such that
$$(T-\lmd)^m(u)=(T-\mu)^n(v)=0.\eqno(1.5.57)$$
Since $T$ is a derivation, we have
\begin{eqnarray*}& &(T-\lmd-\mu)(u\circ v)=T(u\circ v)-(\lmd+\mu)u\circ v
\\ &=&T(u)\circ v+u\circ T(v)-(\lmd+\mu)u\circ v=(T-\lmd)(u)\circ
v+u\circ (T-\mu)(v).\hspace{1.2cm}(1.5.58)\end{eqnarray*} The
Leibniz rule shows
$$(T-\lmd-\mu)^{m+n}(u\circ v)=\sum_{i=0}^{m+n}{m+n\choose
i}(T-\lmd)^i(u)\circ (T-\mu)^{m+n-i}(v)=0.\eqno(1.5.59)$$ Thus
$u\circ v\in{\msr A}_{\lmd+\mu}$. Hence
$$T_s(u\circ v)=(\lmd+\mu)u\circ v=(\lmd u)\circ v+u\circ(\mu
v)=T_s(u)\circ v+u\circ T_s(v)\eqno(1.5.60)$$ by (1.5.56). Hence
$T_s$ is a derivation and so is $T_n=T-T_s.\qquad\Box$\psp

In the rest of this section, we assume that $\mbb{F}=\mbb{C}$. The
trace of linear transformation is trace of its corresponding matrix
under a basis.\psp

 {\bf Lemma 1.5.9}. {\it Let $V$ be a finite-dimensional space and let
 ${\msr A}\subset {\msr B}$ be two subspaces of $\Edo V$. Set
 $M=\{T\in gl(V)\mid [T,{\msr B}]\subset {\msr A}\}$. Fix
 $\xi\in M$. If $\tr \xi T'=0$ for any $T'\in M$, then $\xi$ is
 nilpotent.}

{\it Proof}. Let $\xi=\xi_s+\xi_n$ be  the Jordan-Chevalley
decomposition. Under a certain basis $\{v_1,v_2,...,v_m\}$,
$\xi_s=\mbox{diag}(\lmd_1,\lmd_2,...,\lmd_m)$ and $gl(V)=M_{m\times
m}(\mbb{F})$. Set
$$H=\sum_{i=1}^m\mbb{Q}\lmd_i,\eqno(1.5.61)$$
a subspace of  $\mbb{F}$ over $\mbb{Q}$. Let $f:H\rta \mbb{Q}$ be
any $\mbb{Q}$-linear map and define
$$
T_0=\mbox{diag}(f(\lmd_1),f(\lmd_2),...,f(\lmd_m)).\eqno(1.5.62)$$
We want to prove $T_0\in M$. Note
$$\ad \xi_s(E_{i,j})=(\lmd_i-\lmd_j)E_{i,j},\qquad \ad
T_0(E_{i,j})=(f(\lmd_i)-f(\lmd_j))E_{i,j}.\eqno(1.5.63)$$ Suppose
that $\{a_1,a_2,...,a_s\}$ is the maximal set of distinct elements
in $\{\lmd_i-\lmd_j\mid 1\leq i,j\}$. Define
$$\zeta(t)=\sum_{i=1}^sf(a_i)\frac{\prod_{i\neq
j\in\{1,2,...,s\}}(t-a_j)}{\prod_{i\neq
j\in\{1,2,...,s\}}(a_i-a_j)}.\eqno(1.5.64)$$ Then
$$\zeta(a_i)=f(a_i).\eqno(1.5.65)$$
Since $0\in\{a_1,a_2,...,a_s\}$, $t|\zeta(t)$. Furthermore,
$$\zeta(\lmd_i-\lmd_j)=f(\lmd_i-\lmd_j)=f(\lmd_i)-f(\lmd_j).\eqno(1.5.66)$$
Thus
$$\ad T_0=\zeta (\ad \xi_s).\eqno(1.5.67)$$
By Theorem 1.5.6 and Corollary 1.5.7, $\ad \xi_s$ is a polynomial of
$\ad\xi$ without constant. So is $\ad T_0$. Therefore, $T_0\in M$.

Now
$$\tr T_0\xi=\sum_{i=1}^mf(\lmd_i)\lmd_i=0\eqno(1.5.68)$$
by our assumption. Hence
$$0=f(0)=f(\sum_{i=1}^mf(\lmd_i)\lmd_i)=\sum_{i=1}^mf(\lmd_i)f(\lmd_i)=
\sum_{i=1}^mf(\lmd_i)^2. \eqno(1.5.69)$$ Thus all $f(\lmd_i)=0$.
Since $f$ is any linear function, we have $H=\{0\}$ or equivalently,
$\xi_s=0.\qquad\Box$\psp

{\bf Theorem 1.5.10 (Cartan's Criterion)}.\index{Cartan's Criterion}
 {\it Let $V$ be a
finite-dimensional space and let ${\msr G}$ be a subalgebra of
$gl(V)$. If $\tr uv=0$ for any $u\in[{\msr G},{\msr G}]$ and $v\in
{\msr G}$, then ${\msr G}$ is solvable.}

{\it Proof}. In the above theorem, we let ${\msr A}=[{\msr G},{\msr
G}]$ and ${\msr B}={\msr G}$. Then
$$M=\{T\in gl(V)\mid [T,{\msr G}]\subset [{\msr G},{\msr
G}]\}.\eqno(1.5.70)$$ Take any $u,v\in{\msr G}$ and $T\in M$. We
have
\begin{eqnarray*} \hspace{2cm}& &\tr [u,v]T=\tr(uv-vu)T=\tr (uvT-vuT)=\tr
uvT-\tr vuT\\ &=&\tr uvT-\tr uTv=\tr (uvT-uTv)=\tr u(vT-Tv)\\ &
=&-\tr u[T,v]=-\tr [T,v]u=0.\hspace{6.2cm}(1.5.71)\end{eqnarray*}
 By
linearity, we have
$$\tr wT=0\qquad\mbox{for any}\;w\in[{\msr G},{\msr
G}].\eqno(1.5.72)$$
 According to the above theorem, all the elements in $[{\msr G},{\msr G}]$
 are nilpotent. So $[{\msr G},{\msr G}]$ is a nilpotent Lie
 algebra. Therefore, ${\msr G}$ is solvable.$\qquad\Box$\psp

{\bf Corollary 1.5.11}. {\it If $\tr\ad\xi\:\ad u=0$ for any $\xi\in
[{\msr G},{\msr G}]$ and $u\in {\msr G}$. Then ${\msr G}$ is
solvable.}

{\it Proof}. By the above theorem, $\ad {\msr G}$ is
solvable.$\qquad\Box$\psp

 In [X13],
we introduced a family of nilpotent Lie algebras associated with a
root tree in connection with evolution partial differential
equations. L. Luo [Ll,L2] generalized our algebras to those
associated with more general directed graphs.

\chapter{Semisimple Lie Algebras}

In this chapter, we first introduce the Killing form and prove that
the semisimplicity of a finite-dimensional Lie algebra over
$\mbb{C}$ is equivalent to the nondegeneracy of its Killing form.
Then we use the Killing form to derive the decomposition of a
finite-dimensional semisimple Lie algebra over $\mbb{C}$ into a
direct sum of simple ideals. Moreover, it is showed that $\der{\msr
G}=\ad {\msr G}$ for such a Lie algebra. Furthermore, we study the
completely reducible modules of a Lie algebra and proved the Weyl's
theorem of complete reducibility. The equivalence of the complete
reducibility of real and complex modules is also given. Cartan's
root-space decomposition of a finite-dimensional semisimple Lie
algebra over $\mbb{C}$ is derived. In particular, we prove that such
a Lie algebra is generated by two elements. The completely
reducibility of finite-dimensional modules of $sl(2,\mbb{C})$ plays
an important role in proving the properties of the corresponding
root systems.

\section{Killing Form}

In this section, we assume that Lie algebras are  finite
dimensional. Moreover, we define the Killing form of a Lie algebra
and use it to characterize semisimple Lie algebras.

Let ${\msr G}$ be  a Lie algebra. Define
$$\kappa(u,v)=\tr \ad u\:\ad v\qquad\for\;\;u,v\in{\msr
G}.\eqno(2.1.1)$$ Then $\kappa(\cdot,\cdot)$ is a symmetric bilinear
form that is called {\it Killing form}.\index{Killing form}
 Moreover,
\begin{eqnarray*}\kappa([u,v],w)&=&\tr
\ad[u,v]\:\ad w=\tr [\ad u,\ad v]\:\ad w\\ &=&\tr (\ad u\:\ad v-\ad
v\:\ad u)\ad w=\tr (\ad u\:\ad v\:\ad w-\ad v\:\ad
u\:\ad w)\\ &=&\tr (\ad u\:\ad v\:\ad w)-\tr(\ad v\:\ad u\:\ad w)\\
&=&\tr (\ad u\:\ad v\:\ad w)-\tr(\ad u\:\ad w\:\ad v)\\ &=&\tr (\ad
u\:\ad v\:\ad w-\ad u\:\ad w\:\ad v)\\ &=&\tr \ad u(\ad v\:\ad w-\ad
w\:\ad v)=\tr \ad u[\ad v,\ad w]\\ &=&\tr \ad
u\:\ad[v,w]=\kappa(u,[v,w])\hspace{6.5cm}(2.1.2)\end{eqnarray*} for
$u,v,w\in{\msr G}$. Thus $\kappa$ is {\it associative} ({\it
invariant}).\index{invariant form}\psp

{\bf Lemma 2.1.1}. {\it Let ${\cal I}$ be an ideal of a Lie algebra
${\msr G}$. If $\kappa(\cdot,\cdot)$ is the Killing form of ${\msr
G}$, then $\kappa(|_{\cal I},|_{\cal I})$ is the Killing form of
${\cal I}$.}

{\it Proof}. Assume that ${\cal I}$ is a proper nonzero ideal. Take
subspace ${\msr K}$ of ${\msr G}$ such that
$${\msr G}={\msr K}\oplus {\cal I}.\eqno(2.1.3)$$
Let ${\msr S}_1$ be a basis of ${\msr K}$ and let ${\msr S}_2$ be a
basis of ${\cal I}$. Then ${\msr S}={\msr S}_1\bigcup {\msr S}_2$
forms a basis of ${\msr G}$. For each $u\in {\cal I}$, the matrix
form
$$\mbox{ad}_{\msr G}u=\left(\begin{array}{cc}0&A_u\\
0&B_u\end{array}\right)\eqno(2.1.4)$$ with respect to the basis
$\msr S$ and the matrix form $\mbox{ad}_{\cal I}u=B_u$ with respect
to $\msr S_2$. Thus
$$\ad_{\msr G}u\:\ad_{\msr G}v=\left(\begin{array}{cc}0&A_u
\\ 0&B_u\end{array}\right)\left(\begin{array}{cc}0&A_v
\\ 0&B_v\end{array}\right)=\left(\begin{array}{cc}0&A_uB_v\\
0&B_uB_v\end{array}\right)\qquad\for\;\;u,v\in {\msr
I}.\eqno(2.1.5)$$ Hence
$$\kappa(u,v)=\tr\mbox{ad}_{\msr G}u\:\mbox{ad}_{\msr G}v=\tr
B_uB_v=\tr \mbox{ad}_{\cal I}u\:\mbox{ad}_{\cal I}v\eqno(2.1.6)$$
for $u,v\in {\cal I}.\qquad\Box$\psp

The {\it radical} of a symmetric bilinear form $\be(\cdot,\cdot)$ on
a vector space $V$ is defined by\index{radical!of symmetric form}
$$\rad\:\be=\{u\in V\mid \be(u,v)=0\;\mbox{for any}\;v\in
V\}.\eqno(2.1.7)$$ If $\rad\;\be=\{0\}$, we call $\be(\cdot,\cdot)$
{\it nondegenerate}.\index{nondegenerate}
 Recall the radical of a Lie algebra
${\msr G}$ is the unique maximal solvable ideal. The algebra ${\msr
G}$ is called semisimple if $\rad\:{\msr G}=\{0\}$. In the rest of
this section, we assume the base field $\mbb{F}=\mbb{C}$.\psp

{\bf Theorem 2.1.2}. {\it Let ${\msr G}$ be a Lie algebra. Then
${\msr G}$ is semisimple if and only if the Killing form is
nondegenerate}.

{\it Proof}. Assume that ${\msr G}$ is semisimple. For any
$u\in\rad\:\kappa$ and $v,w\in{\msr G}$, we have
$$\kappa([u,v],w)=\kappa(u,[v,w])=0.\eqno(2.1.8)$$
Hence $[u,v]\in\rad\:\kappa$. So $\rad\:\kappa$ is an ideal of
${\msr G}$. By Cartan's Criterion, $\ad \rad\:\kappa$ is solvable,
and so is $\rad\:\kappa$ because $\kn\ad=Z({\msr G})=\{0\}$. Hence
$\rad\:\kappa\subset \rad\:{\msr G}=\{0\}$.

Suppose that $\rad\:\kappa=\{0\}$ and $R=\rad\:{\msr G}\neq\{0\}$.
Let $n$ be the minimal positive integer such that $R^{(n)}=\{0\}$.
Then $R^{(n-1)}$ is an abelian ideal of ${\msr G}$. For any element
$u\in R^{(n-1)}$ and $v,w\in {\msr G}$, we have
$$(\ad u\:\ad v)^2(w)=[u,[v,[u,[v,w]]]]\in
[R^{(n-1)},R^{(n-1)}]=\{0\}.\eqno(2.1.9)$$ Thus
$$(\ad u\:\ad v)^2=0\lra \tr \ad u\:\ad v=0\sim
\kappa(u,v)=0.\eqno(2.1.10)$$ So $R^{(n-1)}\subset
\rad\:\kappa=\{0\}$, which leads to a contradiction. Therefore,
$\rad\:{\msr G}=\{0\}.\qquad\Box$\psp

{\bf Theorem 2.1.3}. {\it For any  semisimple Lie algebra ${\msr
G}$, there exist ideals ${\msr G}_1,{\msr G}_2,...,{\msr G}_k$ of
${\msr G}$, which are simple Lie algebras, such that
$${\msr G}=\bigoplus_{i=1}^k{\msr G}_i.\eqno(2.1.11)$$
The set $\{{\msr G}_1,{\msr G}_2,...,{\msr G}_k\}$ enumerates all
simple ideals of ${\msr G}$}.

{\it Proof}. Let ${\msr G}_1$ be a minimal nonzero ideal of ${\msr
G}$. Since ${\msr G}_1$ is not solvable by the semisimplicity of
${\msr G}$, $\dim {\msr G}_1>1$. Set
$${\msr G}'=\{u\in {\msr G}\mid \kappa(u,v)=0\;\mbox{for
any}\;v\in {\msr G}_1\}.\eqno(2.1.12)$$ For $u\in{\msr
G}',\;v\in{\msr G}_1$ and $w\in {\msr G}$, we have $[w,v]\in{\msr
G}_1$ and
$$\kappa([u,w],v)=\kappa(u,[w,v])=0.\eqno(2.1.13)$$
Hence $[u,w]\in {\msr G}'$. So ${\msr G}'$ is an ideal of ${\msr
G}$. Since $\kappa(|_{{\msr G}_1},|_{{\msr G}_1})$ is also the
Killing form of ${\msr G}_1$ by Lemma 2.1.1, the radical of the
Killing form of ${\msr G}_1$ is exactly ${\msr G}'\bigcap {\msr
G}_1$. If ${\msr G}'\bigcap {\msr G}_1\neq \{0\}$, then it is an
ideal of ${\msr G}$ included in ${\msr G}_1$. By minimality, ${\msr
G}'\bigcap {\msr G}_1={\msr G}_1$, which implies that the Killing
form of ${\msr G}_1$ is zero. By Cartan's Criterion, ${\msr G}_1$ is
a solvable ideal, which is absurd. Hence ${\msr G}'\bigcap {\msr
G}_1=\{0\}.$ Thus $\kappa(|_{{\msr G}_1},|_{{\msr G}_1})$ is
nondegenerate. Therefore, we have
$${\msr G}={\msr G}_1\oplus {\msr G}'.\eqno(2.1.14)$$
Observe that $[{\msr G}_1,{\msr G}']\subset {\msr G}_1\bigcap {\msr
G}'=\{0\}$. Thus any ideals of ${\msr G}_1$ and ${\msr G}'$ are
ideals of ${\msr G}$. In particular, ${\msr G}_1$ is a simple Lie
algebra and ${\msr G}'$ is semisimple if ${\msr G}'\neq \{0\}$. Now
we work on ${\msr G}'$. By induction, there exist simple ideals
${\msr G}_2,...,{\msr G}_k$ of ${\msr G}'$ such that ${\msr
G}'=\bigoplus_{i=2}^k{\msr G}_i$. So we have (2.1.11).

If ${\msr L}$ is any simple ideal of ${\msr G}$, then $[{\msr
L},{\msr L}]$ is a nonzero ideal of ${\msr L}$. Hence $[{\msr
L},{\msr L}]={\msr L}$, which implies ${\msr L}=[{\msr G},{\msr
L}]$. Now
$${\msr L}=[{\msr G},{\msr L}]=\big[\bigoplus_{i=1}^k{\msr G}_i,{\msr L}\big].\eqno(2.1.15)$$
So $[{\msr G}_i,{\msr L}]\neq \{0\}$ for some $i\in\{1,2,...,k\}$.
But $[{\msr G}_i,{\msr L}]$ is an ideal of ${\msr G}$ that is also
an ideal of ${\msr G}_i$ and ${\msr L}$. Thus ${\msr L}=[{\msr
G}_i,{\msr L}]={\msr G}_i.\qquad\Box$\psp

Let ${\cal I}$ be a nonzero ideal of ${\msr G}$. Suppose
$$u=\sum_{i=1}^k u_i\in{\cal I}.\eqno(2.1.16)$$
If $u_i\neq 0$, then $[u,{\msr G}_i]=[u_i,{\msr G}_i]\neq 0$, which
implies $\{0\}\neq {\msr G}_i\bigcap {\cal I}={\msr G}_i$.
Therefore, ${\cal I}$ is a direct sum of some simple ideals.\psp

{\bf Theorem 2.1.4}. {\it If ${\msr G}$ is a semisimple Lie algebra,
then $\der\: {\msr G}=\ad\: {\msr G}.$}

{\it Proof}. We define a symmetric associative bilinear form $\be$
on $\der{\msr G}\subset \Edo G$ by
$$\be(d_1,d_2)=\tr d_1d_2\qquad\for\;\;d_1,d_2\in\der{\msr
G}.\eqno(2.1.17)$$ Note that
$$(\ad{\msr
G})^{\bot}=\{d\in\der {\msr G}\mid \be(d,\ad u)=0\;\mbox{for
any}\;u\in{\msr G}\}\eqno(2.1.18)$$ is an ideal of $\der{\msr G}$
because  $\ad{\msr G}$ is an ideal of $\der{\msr G}$. Since
$\be(|_{\ad {\msr G}},|_{\ad {\msr G}})=\kappa(\cdot,\cdot)$, the
nondegenerate Killing form of ${\msr G}$, we have $\ad{\msr
G}\bigcap(\ad{\msr G})^{\bot}=\{0\}$. Thus
$$\der {\msr G}=\ad{\msr G}\oplus (\ad{\msr G})^{\bot}.\eqno(2.1.19)$$ For
any $d\in (\ad{\msr G})^{\bot}$ and $u\in {\msr G}$, we have $[d,\ad
u]=\ad d(u)\in\ad{\msr G}\bigcap(\ad{\msr G})^{\bot}=\{0\}$. Thus
$d(u)\in Z({\msr G})=\{0\}$. So $d=0$. Hence $(\ad{\msr
G})^{\bot}=\{0\}$. Therefore $\der{\msr G}=\ad{\msr G}$ by (2.1.19).
$\qquad\Box$\psp

For any element $u$ in a semisimple Lie algebra ${\msr G}$, there
exist elements $u_s,u_n\in{\msr G}$ such that $\ad u=\ad u_s+\ad
u_n$ is the Jordan-Chevalley decomposition of $\ad u$. Since
$Z({\msr G})=\{0\}$, we have $u=u_s+u_n$, which is called the {\it
abstract Jordan decomposition} of $u$.

As an exercise, prove that a finite-dimensional Lie algebra ${\msr
G}$ is solvable if and only if $[{\msr G},{\msr G}]\subset
\rad\:\kappa(\cdot,\cdot)$.

\section{Modules}

In this section, we define a module of a Lie algebra and discuss the
complete reducibility. Moreover, the dual module and tensor module
are introduced. The lifting from a given representation to an
oscillator representation is given.

Recall that a representation $\nu$ of ${\msr G}$ on a vector space
$V$ is a Lie algebra homomorphism from ${\msr G}$ to $gl(V)$. It is
more convenient to work on the elements of ${\msr G}$ directly
acting on $V$. We denote
$$\xi(u)=\nu(\xi)(u)\qquad\for\;\;\xi\in{\msr G},\;u\in
V.\eqno(2.2.1)$$ The conditions for $\nu$ to be a homomorphism
becomes
$$ \xi_1(au+bv)=a\xi_1(u)+b\xi_1(v),\qquad
(a\xi_1+b\xi_2)(u)=a\xi_1(u)+b\xi_2(u),\eqno(2.2.2)$$
$$\xi_1(\xi_2(u))-\xi_2(\xi_1(u))=[\xi_1,\xi_2](u)\eqno(2.2.3)$$
for $a,b\in\mbb{F},\;\xi_1,\xi_2\in{\msr G}$ and $u,v\in V$. We call
$V$ a ${\msr G}$-{\it module}.\index{module} Conversely, a map from
${\msr G}\times V\rta V$: $(\xi,u)\mapsto \xi(u)$ satisfying (2.2.2)
and (2.2.3) gives a representation $\nu$ of ${\msr G}$ on $V$
defined by (2.2.1).  A {\it submodule}\index{submodule} $U$ of a
${\msr G}$-module $V$ is a subspace such that ${\msr G}(U)\subset
U$. A {\it trivial} ${\msr G}$-module $V$ is a vector space with the
${\msr G}$-action: $\xi(u)=0$ for any $\xi\in{\msr G}$ and $u\in V$.

A {\it homomorphism}\index{module homomorphism} $\phi$ from a ${\msr
G}$-module $V$ to a ${\msr G}$-module $W$ is a linear map from $V$
to $W$ such that
$$\phi(\xi(u))=\xi(\phi(u))\qquad\for\;\;\xi\in{\msr G},\;u\in
V.\eqno(2.2.4)$$ In this case, $\phi(V)$ forms a submodule of $W$
and $\kn\phi$ is a submodule of $V$. Moreover, $\phi(V)\cong
V/\kn\phi$. The map $\phi$ is called an {\it
isomorphism}\index{isomorphism} if it is bijective. A ${\msr
G}$-module $V$ is called {\it irreducible}\index{irreducible} if it
has exactly two submodules $\{0\}$ and $V$. We do not view $\{0\}$
as an irreducible ${\msr G}$-module. A module is called {\it
completely reducible}\index{completely reducible} if it is direct
sum of irreducible submodules. A {\it complement}\index{module
complement} of a submodule $W$ of $V$ is a submodule $W'$ such that
$V=W\oplus W'$.\psp

{\bf Lemma 2.2.1}. {\it A module $V$ is completely reducible if and
only if any submodule has a complement.}

{\it Proof}. Suppose that $V\neq \{0\}$ and any submodule of $V$ has
a complement. Take any nonzero minimal submodule $W_1$ of $V$, which
must be irreducible. If $V=W_1$, it is done. Otherwise, it has a
nonzero complement $W_1'$. Suppose that $W_2$ is any nonzero proper
submodule of $W_1'$. Then $W_1+W_2$ is a module. By assumption,
there exists a complement $U$ of $W_1+W_2$. This $U$ may not be in
$W_1'$. Define the projection ${\msr P}: V\rta W_1'$ by
$${\msr P}(u+u')=u'\qquad\for\;\;u\in W_1,\;u'\in
W_1'.\eqno(2.2.5)$$ Then ${\msr P}$ is a Lie algebra module
homomorphism. In fact, ${\msr P}(U)$ is a submodule of $W_1'$.

If $u\in {\msr P}(U)\bigcap W_2$, then there exists $v\in W_1$ such
that $v+u\in U$. But $v+u\in W_1+W_2$. Thus $v+u=0$, or
equivalently, $u=v=0$. Hence ${\msr P}(U)\bigcap W_2=\{0\}$.
Moreover, any element $w\in W_1'$ can be written as $w=w_1+w_2+w'$
with $w_1\in W_1,\;w_2\in W_2$ and $w'\in U$. Note $w'-{\msr
P}(w')\in W_1$. We have
$$ w-w_2-{\msr P}(w')=w_1+w'-{\msr P}(w')\in W_1\bigcap
W_1'=\{0\};\eqno(2.2.6)$$ that is, $w=w_2+{\msr P}(w')$. Hence
$$ W_1'=W_2\oplus {\msr P}(U),\eqno(2.2.7)$$
or equivalently, any submodule of $W_1'$ has a complement in $W_1'$.
By induction, $W_1'=\bigoplus_{i=2}^kW_i$ is a direct sum of
irreducible submodules $W_2,W_3,...,W_k$ in $W_1'$. So is
$V=\bigoplus_{i=1}^kW_i$.

Assume that $V=\bigoplus_{i=1}^kW_i$ is a direct sum of irreducible
submodules $W_1,W_2,...,W_k$ in $V$ and $U$ is any nonzero proper
submodule of $V$. Then there exists some $W_i\not\subset U$. Now
$$\bar{V}=V/W_i=\bigoplus_{i\neq j}(W_j+W_i)/W_i\eqno(2.2.8)$$ is
completely reducible. By induction, $(U+W_i)/W_i$ has a complement
$U'/W_i$ in $\bar{V}$, where $U'$ is a submodule of $V$ containing
$W_i$. If $u\in U\bigcap  U'$, then $u+W_i\in [(U+W_i)/W_i]\bigcap
(U'/W_i)=W_i$. So $u\in W_i\bigcap U=\{0\}$. That is, $U\bigcap
U'=\{0\}$. For any element $v\in V$, there exists $v_1\in U$ and
$v_2\in U'$ such that
$$ v+W_i=(v_1+W_i)+(v_2+W_i)=v_1+v_2+W_i;\eqno(2.2.9)$$
equivalently, there exists $w\in W_i$ such that
$$v=v_1+v_2+w.\eqno(2.2.10)$$
But $W_i\subset U'$. We have $v_2+w\in U'$. Thus $v\in U+U'$.
Therefore, $V=U\oplus U'$; that is, $U'$ is a complement of
$U.\qquad\Box$\psp

Recall that $V^\ast$ denotes the space of linear functions on $V$.
Suppose that $V$ is a module of a Lie algebra ${\msr G}$. We define
an action of ${\msr G}$ on $V^\ast$  by
$$\xi(f)(u)=-f(\xi(u))\qquad\for\;\;\xi\in{\msr G},\;f\in
V^\ast,\;u\in V.\eqno(2.2.11)$$ Note
\begin{eqnarray*}& &[\xi_1,\xi_2](f)(u)\\ &=&f(-[\xi_1,\xi_2](u))
=f(-\xi_1(\xi_2(u))+\xi_2(\xi_1(u)))\\
&=&-f(\xi_1(\xi_2(u)))+f(\xi_2(\xi_1(u)))=-\xi_1(f)(\xi_2(u))+\xi_2(f)(\xi_1(u))\\
&=&-\xi_2(\xi_1(f))(u)+\xi_1(\xi_2(f))(u)=(\xi_1(\xi_2(f))-\xi_2(\xi_1(f)))(u)
\hspace{3.2cm}(2.2.12)\end{eqnarray*} for $\xi_1,\xi_2\in{\msr
G},\;f\in V^\ast$ and $u\in V$. Thus
$$[\xi_1,\xi_2](f)=\xi_1(\xi_2(f))-\xi_2(\xi_1(f)).\eqno(2.2.13)$$
Hence $V^\ast$ forms a ${\msr G}$-module, which is called the {\it
dual (contragredient) module}\index{dual module} of $V$. In general,
$V\not\cong V^\ast$ as ${\msr G}$-modules. It is in general a
difficult problem of determining if they are isomorphic.

A bilinear form $\be:V\times V\rta \mbb{F}$ on a ${\msr G}$-module
$V$ is called ${\msr G}$-{\it invariant} if
$$\be(\xi(u),v)=-\be(u,\xi(v))\qquad\for\;\;\xi\in{\msr
G},\;u,v\in V.\eqno(2.2.14)$$ In fact, we have:\psp

{\bf Lemma 2.2.2}. {\it A finite-dimensional module $V$ of a Lie
algebra ${\msr G}$ is isomorphic to its dual module $V^\ast$ if and
only if it has a nondegenerate invariant bilinear form}.

{\it Proof}. Suppose that $V$ is isomorphic to $V^\ast$ via $\sgm$.
We define a bilinear form $\be$ on $V$ by
$$\be(u,v)=\sgm(u)(v)\qquad\for\;\;u,v\in V.\eqno(2.2.15)$$
Then
$$\be(\xi(u),v)=\sgm(\xi(u))(v)=\xi(\sgm(u))(v)=-\sgm(u)(\xi(v))=-\be(u,\xi(v)).\eqno(2.2.16)$$
for $\xi\in{\msr G}$ and $u,v\in V$. So $\be$ is a nondegenerate
invariant bilinear form.

Assume that $V$ has a nondegenerate invariant bilinear form $\be$.
Equation (2.2.15) define a linear map $\sgm:V\rta V^\ast$ and
(2.2.16) implies that $\sgm$ is a module homomorphism. The
nondegeneracy implies that $\sgm$ is  bijective.$\qquad\Box$\psp

Next, we want to lift a given representation to an oscillator
(differential-operator) representation. Let $\{u_i\mid i\in I\}$ be
a basis of a module $V$ of a Lie algebra ${\msr G}$. For any $\xi\in
{\msr G}$, we write
$$\xi(u_i)=\sum_{j\in I}\nu_{i,j}(\xi)u_j.\eqno(2.2.17)$$
Then $\nu_{i,j}\in {\msr G}^\ast$. We define an action of ${\msr G}$
on the polynomial algebra $\mbb{F}[x_i\mid i\in I]$ by:
$$\xi (f)=\sum_{i,j\in
I}\nu_{i,j}(\xi)x_j\ptl_{x_i}(f).\eqno(2.2.18)$$ By (2.2.3) and
(2.2.17), we have
$$\sum_{j\in
I}(\nu_{i,j}(\xi_2)\nu_{j,k}(\xi_1)-\nu_{i,j}(\xi_1)\nu_{j,k}(\xi_2))
=\nu_{i,k}([\xi_1,\xi_2])\qquad\for\;\;\xi_1,\xi_2\in{\msr
G}.\eqno(2.2.19)$$ On the other hand,
\begin{eqnarray*}& &[\sum_{i,j\in
I}\nu_{i,j}(\xi_1)x_j\ptl_{x_i},\sum_{r,s\in
I}\nu_{r,s}(\xi_2)x_s\ptl_{x_r}]\\ &=&\sum_{i,j,r,s\in
I}\nu_{i,j}(\xi_1)\nu_{r,s}(\xi_2)x_j\ptl_{x_i}(x_s)\ptl_{x_r}-
\sum_{i,j,r,s\in
I}\nu_{i,j}(\xi_1)\nu_{r,s}(\xi_2)x_s\ptl_{x_r}(x_j)\ptl_{x_i}\\
&=&\sum_{i,j,r\in I}\nu_{i,j}(\xi_1)\nu_{r,i}(\xi_2)x_j\ptl_{x_r}-
\sum_{i,j,s\in I}\nu_{i,j}(\xi_1)\nu_{j,s}(\xi_2)x_s\ptl_{x_i}\\
&=& \sum_{r,j\in I}\sum_{i\in I} (\nu_{r,i}(\xi_2)\nu_{i,j}(\xi_1)-
\nu_{r,i}(\xi_1)\nu_{i,j}(\xi_2))x_j\ptl_{x_r}\\ &=&\sum_{r,j\in
I}\nu_{r,j}([\xi_1,\xi_2])x_j\ptl_{x_r}\hspace{10cm}(2.2.20)\end{eqnarray*}
for $\xi_1,\xi_2\in{\msr G}$. Thus $\mbb{F}[x_i\mid i\in I]$ becomes
a ${\msr G}$-module. It is important to understand this module
structure from the information of the module $V$. For instance, any
one-dimensional trivial submodule gives an invariant of the
corresponding Lie group. To find invariants is equivalent to solve
the system of partial differential equations:
$$\sum_{i,j\in
I}\nu_{i,j}(\xi)x_i\ptl_{x_j}(f)=0,\qquad\xi\in{\msr
G}.\eqno(2.2.21)$$ The above defined module is isomorphic to the
``symmetric tensor of $V$."

Let $V$ and $W$ be modules of a Lie algebra ${\msr G}$. We define an
action of ${\msr G}$ on $V\otimes W$ by
$$ \xi(u\otimes v)=\xi(u)\otimes v+u\otimes
\xi(v)\qquad\for\;\;\xi\in {\msr G},\;u\in V,\;v\in
W.\eqno(2.2.22)$$ Note
\begin{eqnarray*}\hspace{1cm}& &(\xi_1\xi_2-\xi_2\xi_1)(u\otimes v)\\&=&\xi_1\xi_2(u\otimes v)
-\xi_2\xi_1(u\otimes v)\\ &=&\xi_1(\xi_2(u)\otimes v+u\otimes
\xi_2(v))-\xi_2(\xi_1(u)\otimes v+u\otimes \xi_1(v))\\
&=&\xi_1(\xi_2(u))\otimes v+\xi_1(u)\otimes \xi_2(v)+\xi_2(u)\otimes
\xi_1(v)+u\otimes \xi_1(\xi_2(v))\\ & &-(\xi_2(\xi_1(u))\otimes
v+\xi_2(u)\otimes \xi_1(v)+\xi_1(u)\otimes \xi_2(v)+u\otimes
\xi_2(\xi_1(v))\\ &=&(\xi_1(\xi_2(u))-\xi_2(\xi_1(u)))\otimes
v+u\otimes(\xi_1(\xi_2(v))-\xi_2(\xi_1(v)))\\
&=&[\xi_1,\xi_2](u)\otimes v+u\otimes[\xi_1,\xi_2](v)
=[\xi_1,\xi_2](u\otimes v)\hspace{4.4cm}(2.2.23)\end{eqnarray*} for
$\xi_1,\xi_2\in{\msr G},\;u\in V$ and $v\in W$. Thus $V\otimes W$
becomes a ${\msr G}$-module.

Denote by $\mbox{Hom}_{\mbb{C}}(V,W)$ the space of all linear maps
from $V$ to $W$. Define $\tau:V^\ast\otimes W\rta
\mbox{Hom}_{\mbb{C}}(V,W)$ by
$$\tau(\sum_{i=1}^kf_i\otimes v_i)(u)=\sum_{i=1}^kf_i(u) v_i\qquad
\for\;\;u\in V.\eqno(2.2.24)$$ If $V$ and $W$ are modules of a Lie
algebra ${\msr G}$ and $\dim W<\infty$, then $\tau$ is a linear
isomorphism, $V^\ast\otimes W$ forms a ${\msr G}$-module and $\tau$
induces a ${\msr G}$-module structure on
$\mbox{Hom}_{\mbb{C}}(V,W)$. In fact,
$$\xi(g)(u)=\xi(g(u))-g(\xi(u))\qquad\for\;\;\xi\in{\msr G},\;g\in
\mbox{Hom}_{\mbb{C}}(V,W),\;u\in V.\eqno(2.2.25)$$

As an exercise, write out the module structure of
$\mbox{Hom}_{\mbb{C}}(V,W)$ for $sl(2,\mbb{F})$ if (1) $V$ is the
adjoint module and $W$ is the $2$-dimensional module of defining
$sl(2,\mbb{F})$; (2) exchange the positions of two modules. Are
these two modules irreducible?

\section{Real and Complex Modules}

In this section, we prove the equivalence of real and complex
complete reducibility. In particular, we find the connection between
real and complex finite-dimensional irreducible modules.

 Let $V$
be a module of a Lie algebra ${\msr G}$ and let $U$ be a nonzero
proper submodule of $V$. We say that $V$ {\it splits} at $U$ if $U$
has a complementary submodule in $V$.\psp

{\bf Lemma 2.3.1}. {\it A ${\msr G}$-module $V$ splits at a
submodule $U$ if and only if there exits a ${\msr G}$-module
homomorphism $\tau:V\rta U$ such that $\tau|_U=\mbox{\it Id}_U$.}

{\it Proof}. If $V$ splits at a submodule $U$, then $U$ has a
complement $U'$ such that $V=U\oplus U'$. The projection $\tau$
defined by
$$ \tau(u+u')=u\qquad\for\;\;u\in U,\;u'\in U',\eqno(2.3.1)$$
is the required homomorphism. Conversely, if there exits a ${\msr
G}$-module homomorphism $\tau:V\rta U$ such that $\tau|_U=\mbox{\it
Id}_U$, then
$$v=\tau(v)+(1-\tau)(v)\qquad\for\;\;v\in V,\eqno(2.3.2)$$
which implies that $(1-\tau)(V)$ is a complement of $U$ because
$1-\tau$ is a ${\msr G}$-module homomorphism from $V$ to
$V$.$\qquad\Box$\psp

 Let $V$ be a vector space over $\mbb{R}$ with a basis $\{v_i\mid
i\in I\}$. The complexification of $V$ is a vector space over
$\mbb{C}$ with $\{v_i\mid i\in I\}$ as a basis; that is,
$$V_{\mbb{C}}=\sum_{i\in I}\mbb{C}v_i.\eqno(2.3.3)$$
Moreover, for $u=\sum_{i\in I}a_iv_i\in V$ and $c\in\mbb{C}$, we
define
$$cu=\sum_{i\in I}a_icv_i.\eqno(2.3.4)$$

Suppose that ${\msr G}$ is a finite-dimensional real Lie algebra.
The Lie bracket on ${\msr G}_{\mbb{C}}$ defined by
$$[\sum_{i\in I}a_iu_i,\sum_{\in I}jb_jv_j]=\sum_{i,j\in I}a_ib_j[u_i,v_j],\qquad
u_i,v_j\in{\msr G},\;a_i,b_j\in\mbb{C}.\eqno(2.3.5)$$ In particular,
${\msr G}$ is a real subalgebra of ${\msr G}_{\mbb{C}}$. In fact,
${\msr G}_{\mbb{C}}={\msr G}+\sqrt{-1}{\msr G}$. For the Killing
form $\kappa_{\msr G}$ of ${\msr G}$ and the Killing form
$\kappa_{{\msr G}_\mbb{C}}$ of ${\msr G}_\mbb{C}$,
$\rad\:\kappa_{{\msr G}_\mbb{C}}=\mbb{C}(\rad\:\kappa_{\msr G})$ by
linear algebra. In particular, $\rad\:\kappa_{\msr G}$ is a solvable
ideal of ${\msr G}$ because $\rad\:\kappa_{{\msr G}_\mbb{C}}$ is.
Thus ${\msr G}$ is semisimple if and only $\kappa_{\msr G}$ is
nondegenerate. By the proof of Theorem 2.1.3, ${\msr G}$ is a direct
sum of simple ideals if it is semisimple. Moreover, ${\msr G}$ is
semisimple if and only if ${\msr G}_\mbb{C}$ is.

If $V$ is a finite-dimensional (real) ${\msr G}$-module, then
$V_{\mbb{C}}$ becomes a (complex) ${\msr G}_{\mbb{C}}$-module with
the action:
$$(\sum_ia_i\xi_i)(\sum_j
b_jv_j)=\sum_{i,j}a_ib_j\xi_i(v_j),\qquad\for\;\;a_i,b_j\in\mbb{C},\;\xi_i\in{\msr
G},\;v_j\in V.\eqno(2.3.6)$$ We can also view $V_{\mbb{C}}$ as a
real ${\msr G}$-module. In fact, $V_{\mbb{C}}=V+\sqrt{-1}V$ is a
decomposition of real ${\msr G}$-submodules. As real ${\msr
G}$-modules, $V\cong \sqrt{-1}V$.\pse

{\bf Theorem 2.3.2}. {\it The space $V$ is a completely reducible
real ${\msr G}$-module if and only if $V_{\mbb{C}}$ is a completely
reducible complex ${\msr G}_{\mbb{C}}$-module}.

{\it Proof}. Suppose that $V$ is a completely reducible real ${\msr
G}$-module. Then $V=\bigoplus_{i=1}^kV_i$ is a direct sum of real
irreducible ${\msr G}$-submodules. Moreover,
$$V_{\mbb{C}}=\bigoplus_{i=1}^k(V_i\oplus
\sqrt{-1}V_i)\eqno(2.3.7)$$ is a direct sum of $2k$ real irreducible
${\msr G}$-submodules. So $V_{\mbb{C}}$ is a completely reducible
real ${\msr G}$-module. Let $U$ be a complex ${\msr
G}_{\mbb{C}}$-submodule of $V_{\mbb{C}}$. Note that $U$ is also a
real ${\msr G}$-submodule of $V_{\mbb{C}}$. Thus there exists a real
${\msr G}$-module homomorphism $\tau:V_{\mbb{C}}\rta U$ such that
$\tau|_U=\mbox{Id}_U$. Define
$$\td{\tau}(v)=\frac{1}{2}(\tau(v)-\sqrt{-1}\tau(\sqrt{-1}v))\qquad\for\;\;v\in
V.\eqno(2.3.8)$$ Note
\begin{eqnarray*}\td{\tau}(\sqrt{-1}v)&=&\frac{1}{2}(\tau(\sqrt{-1}v)+\sqrt{-1}\tau(v))
\\
&=&\frac{1}{2}\sqrt{-1}(-\sqrt{-1}\tau(\sqrt{-1}v)+\tau(v))=\sqrt{-1}\td{\tau}(v).
\hspace{4.1cm}(2.3.9)\end{eqnarray*} for $v\in V$. So $\td{\tau}$ is
a $\mbb{C}$-linear map from $V_{\mbb{C}}$ to $U$. Moreover,
$$\xi(\td{\tau}(v))=\frac{1}{2}(\tau(\xi(v))-\sqrt{-1}\tau(\sqrt{-1}\xi(v)))=\td{\tau}(\xi(v))
\eqno(2.3.10)$$ for $v\in V$ and $\xi\in{\msr G}$ by (2.3.6). Thus
$\td{\tau}$ is a ${\msr G}_{\mbb{C}}$-module homomorphism due to
${\msr G}_{\mbb{C}}=\mbb{C}{\msr G}$ by (2.3.4). Furthermore, for
$u\in U$, we have $\sqrt{-1}u\in U$ because $U$ is a complex ${\msr
G}_{\mbb{C}}$-submodule of $V_{\mbb{C}}$. Hence $\tau(u)=u$ and
$\tau(\sqrt{-1}u)=\sqrt{-1}u$. Thus
$$\td{\tau}(u)=\frac{1}{2}(u-\sqrt{-1}(\sqrt{-1}u))=u;\eqno(2.3.11)$$
that is, $\td{\tau}|_U=\mbox{Id}_U$. By the above lemma,
$V_{\mbb{C}}$ is completely reducible.

Conversely, we assume that $V_{\mbb{C}}$ is a completely reducible
complex ${\msr G}_{\mbb{C}}$-module. Let $W$ be a real ${\msr
G}$-submodule of $V$. Then $W_{\mbb{C}}$ forms a ${\msr
G}_{\mbb{C}}$-submodule. There exists a ${\msr G}_{\mbb{C}}$-module
homomorphism $\bar{\tau}:V_{\mbb{C}}\rta W_{\mbb{C}}$ such that
$\bar{\tau}|_{W_{\mbb{C}}}=\mbox{Id}_{W_{\mbb{C}}}$. So
$\bar{\tau}|_V$ is real ${\msr G}$-module homomorphism from $V$ to
$W_{\mbb{C}}$ because ${\msr G}\subset {\msr G}_{\mbb{C}}$. But in
general, $\bar{\tau}(V)\not\subset W$. Define
$$\nu(w_1+\sqrt{-1}w_2)=w_1\qquad w_1,w_2\in W.\eqno(2.3.12)$$
Then $\nu:W_{\mbb{C}}\rta W$ is a real ${\msr G}$-module
homomorphism. Therefore, $\tau=\nu\bar{\tau}$ is a real ${\msr
G}$-module homomorphism from $V$ to $W$ and $\tau|_W=\mbox{Id}_W$.
By the above lemma, $V$ is  a completely reducible real ${\msr
G}$-module.$\qquad\Box$\pse

Let $M=\bigoplus_{i=1}^kM_i$ be a completely reducible module of a
Lie algebra ${\msr G}$ over an arbitrary field, where $M_i$'s are
irreducible ${\msr G}$-submodules. For $i\in\ol{1,k}$, we define the
projection ${\msr P}_i:M\rta M_i$ by
$${\msr P}_i(\sum_{r=1}^ku_r)=u_i\qquad \for\;\;u_r\in
M_r.\eqno(2.3.13)$$ Then ${\msr P}_i$'s are ${\msr G}$-module
homomorphisms. Let $U$ be any irreducible ${\msr G}$-submodule of
$M$. If ${\msr P}_i(U)\neq\{0\}$, then ${\msr P}_i|_U$ is a ${\msr
G}$-module isomorphism from $U$ to $M_i$ by the irreducibility of
$U$ and $M_i$; in particular, ${\msr P}_i(U)=M_i$. Set
$$I=\{r\in\ol{1,k}\mid {\msr P}_r(U)\neq\{0\}\}.\eqno(2.3.14)$$
For $r,s\in I$, the map ${\msr P}_s({\msr P}_r|_U)^{-1}$ is a ${\msr
G}$-module isomorphism from $M_r$ to $M_s$. Re-indexing
$\{M_1,M_2,...,M_k\}$ if necessary, we can assume $1\in I$. Now
$$u=\sum_{r\in I}{\msr P}_r(u)={\msr P}_1(u)+\sum_{1\neq r\in
I}[{\msr P}_r({\msr P}_1|_U)^{-1}]({\msr P}_1(u)) \qquad\for\;\;u\in
U.\eqno(2.3.15)$$ Thus
$$U=\{v+\sum_{1\neq r\in
I}[{\msr P}_r({\msr P}_1|_U)^{-1}](v)\mid v\in M_1\}.\eqno(2.3.16)$$

Suppose that $M=\bigoplus_{j=1}^{k'}V_j$ is another decomposition of
irreducible ${\msr G}$-submodules. Then ${\msr P}_i(V_1)\neq \{0\}$
for some $i\in\ol{1,k}$. Re-indexing $\{M_1,M_2,...,M_k\}$ if
necessary, we can assume ${\msr P}_1(V_1)\neq \{0\}$. So ${\msr
P}_1|_{V_1}$ is a ${\msr G}$-module isomorphism from $V_1$ to $M_1$.
Moreover,
$$V=V_1\oplus\bigoplus_{i=2}^kM_i.\eqno(2.3.17)$$
In particular,
$$\bigoplus_{i=2}^k(M_i+V_1)/V_1=V/V_1=\bigoplus_{r=2}^{k'}(V_i+V_1)/V_1\eqno(2.3.18)$$
as ${\msr G}$-modules. By induction on $k$, we have $k-1=k'-1$;
equivalently, $k=k'$. This shows that although a completely
reducible ${\msr G}$-module may have more than one decompositions
into irreducible submodules, the numbers of irreducible submodules
in the decompositions are equal.

Let $V$ be a finite-dimensional real irreducible module of a real
Lie algebra ${\msr G}$. Then $V_{\mbb{C}}=\bigoplus_{i=1}^kW_i$ is a
direct sum of irreducible ${\msr G}_{\mbb{C}}$-submodules. Since
$W_i$ are ${\msr G}$-modules and $V_{\mbb{C}}=V\oplus \sqrt{-1}V$ is
a direct sum of two irreducible real ${\msr G}$-submodules, $W_i$
contain a real ${\msr G}$-submodule isomorphic to $V$.  Let $M$ be a
finite-dimensional irreducible ${\msr G}_{\mbb{C}}$-module. Take a
minimal nonzero real ${\msr G}$-submodule $U$ of $M$. Then $U$ is an
irreducible real ${\msr G}$-module and so is $\sqrt{-1}U$. If
$U\bigcap \sqrt{-1}U\neq \{0\}$, then $U=\sqrt{-1}U$. So $U$ is also
a ${\msr G}_\mbb{C}$-submodule. Hence $M=U$ is an irreducible real
${\msr G}$-module. Assume $U\bigcap \sqrt{-1}U= \{0\}$. Then
$U+\sqrt{-1}U$ is a ${\msr G}_\mbb{C}$-submodule, and so
$M=U+\sqrt{-1}U$. Denote by ${\msr P}$ the corresponding projection
from $M$ to $U$. Let $U'$ be any nonzero proper ${\msr G}$-submodule
such that ${\msr P}(U')\neq\{0\}$. Then ${\msr P}(U')=U$ by the
irreducibility of $U$. If $\kn {\msr P}|_{U'}=\{0\}$, then
$U'\stl{{\msr P}|_{U'}}{\cong}U$ is an irreducible real ${\msr
G}$-module. Suppose $\kn {\msr P}|_{U'}\neq\{0\}$. Since $\kn{\msr
P}=\sqrt{-1}U$, we have $U'\bigcap \sqrt{-1}U\neq\{0\}$. By the
irreducibility of $\sqrt{-1}U$, $U'\supset \sqrt{-1}U$. The facts
${\msr P}(U')=U$ and $M=U+\sqrt{-1}U$ imply $U'=M$, which
contradicts the properness of $U'$. Note that any nonzero proper
${\msr G}$-submodule $U'$ such that ${\msr P}(U')=\{0\}$ must be a
submodule of $\sqrt{-1}U$. So $U'=\sqrt{-1}U$ by the irreducibility
of $\sqrt{-1}U$. In summary, we have: \psp

{\bf Theorem 2.3.3}. {\it Let ${\msr G}$ be a real Lie algebra. Any
finite-dimensional irreducible real ${\msr G}$-module must be a real
${\msr G}$-submodule of a finite-dimensional irreducible complex
${\msr G}_\mbb{C}$-module. A finite-dimensional irreducible complex
${\msr G}_\mbb{C}$-module is either an irreducible real ${\msr
G}$-module or any of its proper ${\msr G}$-submodules is an
irreducible real ${\msr G}$-module.}\psp

Thus we can find all finite-dimensional irreducible real ${\msr
G}$-modules from finite-dimensional irreducible complex ${\msr
G}_{\mbb{C}}$-modules.

\section{Weyl's Theorem}

In this section, we prove the Weyl's theorem of complete
reducibility.

Assume the base field $\mbb{F}=\mbb{C}$. We will prove that any
finite-dimensional module of a finite-dimensional semisimple Lie
algebra is completely reducible. A representation is called
irreducible if the associated module is irreducible. First, we
have:\psp

{\bf Lemma 2.4.1 (Schur)}.\index{Schur's Lemma} {\it Let ${\msr G}$
be a Lie algebra and let $\phi:{\msr G}\rta gl(V)$ be a
finite-dimensional irreducible representation. If $T\in gl(V)$
commutes with every element in $\phi({\msr G})$, then
$T=\lmd\mbox{\it Id}_V$ for some $\lmd\in\mbb{C}$.}

{\it Proof}. Since we assume $\mbb{F}=\mbb{C}$ and $V$ is
finite-dimensional, $T$ has an eigenvector $u$ whose eigenvalue
denoted as $\lmd$. Denote by $\mbb{Z}_+$ the set of positive
integers. Set
$$U=\mbox{Span}\{u,\xi_1(\cdots(\xi_i(u))\cdots)\mid
i\in\mbb{Z}_+,\;\xi_j\in{\msr G}\}.\eqno(2.4.1)$$ Then $U$ forms a
nonzero submodule of $V$. The irreducibility of $V$ implies $V=U$.
Since
\begin{eqnarray*}\hspace{2cm}T(\xi_1(\cdots(\xi_i(u))\cdots))&=&\xi_1(T(\xi_2(\cdots(\xi_i(u))\cdots)))
\\ &=&\xi_1(\cdots(\xi_i(T(u)))\cdots)\\
&=&\lmd\xi_1(\cdots(\xi_i(u))\cdots),\hspace{4.7cm}(2.4.2)\end{eqnarray*}
We have $T=\lmd\mbox{Id}_V.\qquad\Box$\psp

Suppose that $\phi:{\msr G}\rta gl(V)$ is a finite-dimensional
faithful ($\kn\phi=\{0\}$) representation of a semisimple Lie
algebra ${\msr G}$. Define
$$\be(\xi,\zeta)=\tr \phi(\xi)\phi(\zeta)\qquad \for\;\;\xi,\zeta\in{\msr G}.\eqno(2.4.3)$$ By Cartan's Criterion,
$\phi(\rad\:\be)$ is a solvable ideal of $\phi({\msr G})$, which is
semisimple. Thus $\rad\: \be=\{0\}$, and there exists an orthonormal
basis $\{\xi_1,\xi_2,...,\xi_n\}$ of ${\msr G}$ with respect to
$\be$. Set
$$\co_\phi=\sum_{i=1}^n\phi(\xi_i)^2\in\Edo V,\eqno(2.4.4)$$
 which is called a {\it Casimier operator}.\index{Casimier operator} Write
 $$[\xi_i,\xi_j]=\sum_{k=1}^nc_{i,j}^k\xi_k.\eqno(2.4.5)$$
 Then
 $$c_{i,j}^k=\be([\xi_i,\xi_j],\xi_k)=-\be(\xi_j,[\xi_i,\xi_k])=-c_{i,k}^j.\eqno(2.4.6)$$
 Thus
 \begin{eqnarray*}[\phi(\xi_i),\co_\phi]&=&[\phi(\xi_i),\sum_{j=1}^n\phi_j(u)^2]
 =\sum_{j=1}^n[\phi(\xi_i),\phi(\xi_j)^2]\\
 &=&\sum_{j=1}^n([\phi(\xi_i),\phi(\xi_j)]\phi(\xi_j)+\phi(\xi_j)[\phi(\xi_i),\phi(\xi_j)])\\
 &=&\sum_{j=1}^n(\phi([\xi_i,\xi_j])\phi(\xi_j)+\phi(\xi_j)\phi([\xi_i,\xi_j]))
 \\ &=&\sum_{j,k=1}^nc_{i,j}^k(\phi(\xi_k)\phi(\xi_j)+\phi(\xi_j)\phi(\xi_k))\\
 &=&\sum_{j,k=1}^nc_{i,j}^k\phi(\xi_k)\phi(\xi_j)+\sum_{j,k=1}^nc_{i,j}^k\phi(\xi_j)\phi(\xi_k)
\\ &=&\sum_{j,k=1}^nc_{i,k}^j\phi(\xi_j)\phi(\xi_k)+\sum_{j,k=1}^nc_{i,j}^k\phi(\xi_j)\phi(\xi_k)
\\
&=&\sum_{j,k=1}^n(c_{i,j}^k+c_{i,k}^j)\phi(\xi_j)\phi(\xi_k)=0.\hspace{6.5cm}(2.4.7)\end{eqnarray*}
So we have
$$\phi(u)\co_\phi=\co_\phi\phi(u)\qquad\mbox{for}\;\;u\in{\msr
G}.\eqno(2.4.8)$$ If $\phi$ is irreducible, then
$\co_\phi=\lmd_\phi\mbox{Id}_V$ for some $\lmd_\phi\in \mbb{C}$ by
Lemma 2.4.1. In fact,
$$\lmd_\phi\dim V=\tr \co_\phi=\sum_{i=1}^n\tr
\phi(\xi_i)^2=\sum_{i=1}^n\be(\xi_i,\xi_i)=n.\eqno(2.4.9)$$ Hence
$$\lmd_\phi=n/\dim V=\dim{\msr G}/\dim V.\eqno(2.4.10)$$
Furthermore, we have:\psp

{\bf Lemma 2.4.2}. {\it Let  $\phi:{\msr G}\rta gl(V)$ be a
finite-dimensional  representation of a semisimple Lie algebra
${\msr G}$. Then $\phi({\msr G})\in sl(V)$. In particular,
one-dimensional module of ${\msr G}$ must be trivial.}

{\it Proof}. This follows from $\phi({\msr G})=\phi([{\msr G},{\msr
G}])=[\phi({\msr G}),\phi({\msr G})]$ and the fact $\tr [A,B]=0$ for
any $A,B\in gl(V).\qquad\Box$\psp

{\bf Lemma 2.4.3}. {\it Let  $\phi:{\msr G}\rta gl(V)$ be a
finite-dimensional  representation of a semisimple Lie algebra
${\msr G}$. Suppose that $V$ has an irreducible submodule $U$ of
codimension one. Then $V$ splits at $U$.}

{\it Proof}. If $\phi({\msr G})=\{0\}$, the lemma holds trivially.
Suppose $\phi({\msr G})\neq\{0\}$. Since $\kn\phi$ is either $\{0\}$
or a direct sum of some of simple ideals of ${\msr G}$, $\phi({\msr
G})\cong {\msr G}/\kn\phi$ is semisimple. Replacing ${\msr G}$ by
$\phi({\msr G})$, we may assume $\phi$ is faithful. By Lemma 2.4.2,
$V/U$ is a trivial module. Thus
$$\xi(V)\subset U\qquad \for\;\;\xi\in {\msr G}.\eqno(2.4.11)$$
If $U$ is a trivial submodule, then $V$ is two-dimensional and
$\phi({\msr G})$ is isomorphic to the Lie algebra of strict
up-triangular matrices, which is absurd. Thus $\phi({\msr G})|_U\neq
\{0\}$ is semisimple.

Expression (2.4.11) implies
$$\be(\xi,\zeta)=\tr\phi(\xi)\phi(\zeta)=\tr
\phi(\xi)|_U\phi(\zeta)|_U\qquad\for\;\;\xi,\zeta\in{\msr
G}.\eqno(2.4.12)$$ This shows
$$\co_\phi|_U=(\dim\phi({\msr
G})|_U/\dim U)\mbox{Id}_U.\eqno(2.4.13)$$ Hence
$$\tau=\frac{\dim U}{\dim\phi({\msr
G})|_U}\co_\phi\eqno(2.4.14)$$ is a module homomorphism from $V$ to
$U$ with $\tau|_U=\mbox{Id}_U.$ So $V$ splits at $U$ by Lemma
2.3.1.$\qquad\Box$\psp

{\bf Lemma 2.4.4}. {\it Let  $\phi:{\msr G}\rta gl(V)$ be a
finite-dimensional  representation of a semisimple Lie algebra
${\msr G}$. Suppose that $V$ has a submodule $U$ of codimension one.
Then $V$ splits at $U$.}

{\it Proof}. We prove it by induction on $\dim V$. It holds for
$\dim V=1$. Suppose that it holds for $\dim V<k$. Consider the case
$\dim V=k$. We assume $U\neq\{0\}$. If $U$ is irreducible, this is
the above lemma. Suppose that $U$ has a nonzero proper submodule
$W$. Then $U/W$ is a ${\msr G}$-submodule of $V/W$ with codimension
one. By induction, there exists a submodule $W'\supset W$ of $V$
such that
$$V/W=W'/W\oplus U/W.\eqno(2.4.15)$$
Note $\dim W'/W=1$. By induction, there exists a one-dimensional
submodule $U'$ of $W'$ such that $W'=U'\oplus W$. According to
(2.4.15), $U'\not\subset U$. Hence $U'$ is a complement of $U$ in
$V.\qquad\Box$\psp

{\bf Theorem 2.4.5 (Weyl)}.\index{Weyl's Theorem} {\it Any nonzero
finite-dimensional module $V$ of a semisimple Lie algebra ${\msr G}$
is completely reducible}.

{\it Proof}. Suppose that $U$ is a nonzero proper submodule of $V$.
Then $\mbox{Hom}_{\mbb{C}}(V,U)$ forms a ${\msr G}$-module. If
$f\in\mbox{Hom}_{\mbb{C}}(V,U)$ such that $f|_U=a\mbox{Id}_U$ with
$a\in\mbb{C}$, then
$$\xi(f)(u)=\xi(f(u))-f(\xi(u))=\xi(au)-a\xi(u)=0\qquad\for\;\;u\in
U.\eqno(2.4.16)$$  Thus
$$W=\{f\in\mbox{Hom}_{\mbb{C}}(V,U)\mid
f|_U\in\mbb{C}\mbox{Id}_U\},\;\;W_1=\{f\in\mbox{Hom}_{\mbb{C}}(V,U)\mid
f|_U=0\}\eqno(2.4.17)$$ are submodules. Moreover, $\dim W/W_1=1$. By
Lemma 2.4.4, $W_1$ has a complement $W_1'$ in $W$, which is a
one-dimensional trivial module. Hence there exists $\tau\in W_1'$
such that $\tau|_U=\mbox{Id}_U$ and $\xi(\tau)=0$ for any
$\xi\in{\msr G}$; that is,
$$0=\xi(\tau)(v)=\xi(\tau(v))-\tau(\xi(v))\qquad\for\;\;v\in
V.\eqno(2.4.18)$$ So $\tau$ is a ${\msr G}$-homomorphism. By Lemma
2.3.1, $V$ splits at $U$. Therefore, $V$ is completely reducible by
Lemma 2.2.1.$\qquad\Box$\psp

{\bf Theorem 2.4.6}. {\it If ${\msr G}\subset gl(V)$ is a semisimple
Lie subalgebra and dim\:$V<\infty$, then ${\msr G}$ includes the
semisimple part and nilpotent part of its any element. In
particular, the abstract and usual Jordan decompositions coincide}.

{\it Proof}. By the above theorem, $V=\bigoplus_{i=1}^kV_i$ is a
direct sum of irreducible submodules. Set
$${\msr L}=\{T\in gl(V)\mid [T,{\msr G}]\subset{\msr G};\;T(V_i)\subset
V_i,\;\tr T|_{V_i}=0,\;i=1,2,...,k\}.\eqno(2.4.19)$$ Since ${\msr
G}|_{V_i}$ is a homomorphic image of the semisimple Lie algebra
${\msr G}$, it is a semisimple Lie algebra. Thus we have $\tr
\xi|_{V_i}=0$ for $\xi\in{\msr G}$ and $i=1,...,k$ by Lemma 2.4.2.
Hence ${\msr G}\subset {\msr L}$. For any $T\in {\msr L}$,
$\ad_{\msr G} T\in \der {\msr G}=\ad{\msr G}$ (cf. Theorem 2.1.4).
So there exists an element $\xi'\in{\msr G}$ such that $\ad_{\msr G}
T=\ad_{\msr G} \xi'$. Set
$$T_1=T-\xi'.\eqno(2.4.20)$$
Then
$$[T_1,\xi]=0\qquad\for\;\;\xi\in{\msr G}.\eqno(2.4.21)$$
In particular, $T_1|_{V_i}$ is a constant map $\lmd \mbox{Id}_{V_i}$
by Schur's Lemma (Lemma 2.4.1). Now $0=\tr T_1|_{V_i}=\lmd(\dim
V_i)$ by (2.4.19), or equivalently, $\lmd=0$. Hence $T_1=0$, that
is, $T=\xi'\in {\msr G}$. This shows ${\msr L}={\msr G}$.

Given $\xi\in {\msr G}$. Let $\xi=\xi_s+\xi_n$ be the usual
Jordan-Chevalley decomposition of $\xi$. Recall that $\xi_n$ is a
polynomial of $\xi$ without constant. Hence $\xi_n(V_i)\subset V_i$.
Since $\mbox{ad}_{gl(V)}
\xi=\mbox{ad}_{gl(V)}\xi_s+\mbox{ad}_{gl(V)}\xi_n$ is also a
Jordan-Chevalley decomposition, we have $[\xi_n,{\msr G}]\subset
{\msr G}$. Note that $\xi_n|_{V_i}$ is nilpotent because $\xi_n$ is.
So $\tr\xi_n|_{V_i}=0$; that is, $\xi_n\in {\msr L}={\msr G}$.
Moreover, $\xi_s=\xi-\xi_n\in{\msr G}$. Note that the usual
Jordan-Chevalley decomposition $\xi=\xi_s+\xi_n$ implies that
$\mbox{ad}_{\msr G}\xi_s$ is semisimple, $\mbox{ad}_{\msr G}\xi_n$
is nilpotent and $[\mbox{ad}_{\msr G}\xi_s,\mbox{ad}_{\msr
G}\xi_n]=0$. By the uniqueness of the Jordan-Chevalley
decomposition, $\xi=\xi_s+\xi_n$ is also the abstract
Jordan-Chevalley decomposition. $\qquad\Box$\psp

{\bf Corollary 2.4.7}. {\it Let  $\phi:{\msr G}\rta gl(V)$ be a
finite-dimensional  representation of a semisimple Lie algebra
${\msr G}$. Given $\xi\in {\msr G}$. If $\xi=\xi_s+\xi_n$ is the
abstract Jordan-Chevalley  decomposition, then
$\phi(\xi)=\phi(\xi_s)+\phi(x_n)$ is the usual Jordan-Chevalley
decomposition}.

{\it Proof}. If $\phi({\msr G})\neq \{0\}$, it is semisimple.
Moreover, the usual Jordan-Chevalley decomposition
$\phi(\xi)=\phi(\xi)_s+\phi(\xi)_n$ is also the abstract
Jordan-Chevalley decomposition of $\phi(\xi)$ by the above theorem.
If $\xi=\xi_s+\xi_n$ is the abstract Jordan-Chevalley decomposition,
then $\ad\xi=\ad\xi_s+\ad\xi_n$ is the usual  Jordan-Chevalley
decomposition. Since $\phi(\msr G)\cong {\msr G}/\ker\:\phi$, $\ad
\phi(\xi_s)$ is semisimple and $\ad\phi(\xi_n)$ is nilpotent.
Moreover,
$[\ad\phi(\xi_n),\ad\phi(\xi_n)]=\ad\phi([\xi_s,\xi_n])=0$. Hence
$\ad\phi(\xi)=\ad\phi(\xi_s)+\ad\phi(\xi_n)$ is the Jordan-Chevalley
decomposition of $\phi(\xi)$ by the uniqueness. Thus
$\phi(\xi)=\phi(\xi_s)+\phi(\xi_n)$ is the abstract Jordan-Chevalley
decomposition of $\phi(\xi)$. Therefore, $\phi(\xi_s)=\phi(\xi)_s$
and $\phi(\xi_n)=\phi(\xi)_n.\qquad\Box$ \psp

Let ${\msr G}$ be a simple Lie algebra and  let $\be,\gm$ be two
 symmetric associative bilinear forms on ${\msr G}$. If $\be\neq
 0$, then $\gm=a\be$ for some $a\in\mbb{C}$ (exercise).

 We can identify $({\msr G}\otimes {\msr G})^\ast$ with the space
 of bilinear forms on ${\msr G}$. An associative form $\be$ is an
 element in $({\msr G}\otimes {\msr G})^\ast$ such that ${\msr
 G}(\be)=\{0\}$ (exercise), where we view ${\msr G}$ as the adjoint
 module of ${\msr G}$.

\section{Root Space Decomposition}

In this section, we prove the Cartan's root space decomposition of a
semisimple Lie algebra.

We always assume that ${\msr G}$ is a finite-dimensional semisimple
Lie algebra over $\mbb{C}$. Recall the abstract Jordan decomposition
of any element $\xi=\xi_s+\xi_n$ in ${\msr G}$, where
$\xi_s,\xi_n\in{\msr G}$ and $\ad\xi=\ad\xi_s+\ad\xi_n$ is the usual
Jordan decomposition. By Engel's Theorem, there exist elements in
${\msr G}$ that possess nonzero semisimple parts. Thus ${\msr G}$
has nonzero semisimple elements. A Lie subalgebra of ${\msr G}$
consisting of semisimple elements is called {\it toral}.\index{toral
subalgebra}\psp

{\bf Lemma 2.5.1}. {\it A toral subalgebra is abelian}.

{\it Proof}. Let ${\cal T}$ be a toral subalgebra of ${\msr G}$.
Given $\xi\in {\cal T}$. The operator $\ad_{\cal T} \xi$ is
diagonalizable. Suppose that $\zeta$ is an eigenvector of $\ad_{\cal
T}\xi$ with eigenvalue $\lmd$; that is $[\xi,\zeta]=\lmd v$. On the
other hand, $\ad_{\cal T}v$ is diagonalizable implies that ${\cal
T}$ has a basis $\{\zeta_1,\zeta_2,...,\zeta_n\}$ consisting
eigenvectors of $\ad_{\cal T}\zeta$ with corresponding eigenvalues
$\{\lmd_1,\lmd_2,...,\lmd_n\}$. In particular, we can take
$\zeta_1=\zeta$ with $\lmd_1=0$. Write
$$\xi=\sum_{i=1}^na_i\zeta_i,\qquad a_i\in\mbb{C}.\eqno(2.5.1)$$
Then
$$\lmd
\zeta=[\xi,\zeta]=[\sum_{i=1}^na_i\zeta_i,\zeta]=-\sum_{i=2}^n\lmd_ia_i\zeta_i.\eqno(2.5.2)$$
By linear independence of $\{\zeta,\zeta_2,...,\zeta_n\}$, $\lmd=0$.
Thus all the eigenvalues of $\ad_{\cal T} u$ are 0. Hence $\ad_{\cal
T}u=0.\qquad\Box$\psp

Let $H$ be a maximal toral subalgebra of ${\msr G}$. For instance,
we take the algebra of all diagonal elements in $sl(n,\mbb{C})$.
Since $\{\ad h\mid h\in H\}$ are commuting diagonalizable operators,
by linear algebra,
$${\msr G}=\bigoplus_{\al\in H^\ast}{\msr G}_\al,\qquad{\msr
G}_\al=\{\xi\in{\msr G}\mid [h,\xi]=\al(h)\xi\;\for\;h\in
H\}.\eqno(2.5.3)$$ \vspace{0.1cm}

{\bf Lemma 2.5.2}. {\it For $\al,\be\in H^\ast$, $[{\msr
G}_\al,{\msr G}_\be]\subset{\msr G}_{\al+\be}$. In particular, for
any $u\in{\msr G}_\al$ with $\al\neq 0$, $\ad u$ is nilpotent.
Moreover,
$$\kappa({\msr G}_\al,{\msr G}_\be)=\{0\}\qquad\for\;\;\al,\be\in
H^\ast,\;\al\neq-\be,\eqno(2.5.4)$$ which implies that
$\kappa(|_{{\msr G}_0},|_{{\msr G}_0})$ is nondegenerate.}

{\it Proof}. For $h\in H,\;\xi\in{\msr G}_\al$ and $\zeta\in{\msr
G}_\be$, we have:
$$[h,[\xi,\zeta]]=[[h,\xi],\zeta]+[\xi,[h,\zeta]]=(\al(h)+\be(h))[\xi,\zeta]=(\al+\be)(h)[\xi,\zeta].\eqno(2.5.5)$$
Thus $[{\msr G}_\al,{\msr G}_\be]\subset {\msr G}_{\al+\be}$.
Suppose that $\al,\be\in H^\ast$ such that $\al\neq -\be$. Given
$h\in H,\;\xi\in{\msr G}_\al$ and $\zeta\in{\msr G}_\be$, we get
$$\al(h)\kappa(\xi,\zeta)=\kappa([h,\xi],\zeta)=-\kappa(\xi,[h,\zeta])=-\be(h)\kappa(\xi,\zeta).\eqno(2.5.6)$$
Thus
$$(\al+\be)(h)\kappa(\xi,\zeta)=0\qquad\mbox{for any}\;\;h\in
H.\eqno(2.5.7)$$ Thus $\kappa(\xi,\zeta)=0$; that is, (2.5.4) holds.
For $0\neq \xi\in {\msr G}_0$, we have
$$\kappa(\xi,{\msr G}_\al)=\{0\}\qquad\for\;\;0\neq\al\in
H^\ast.\eqno(2.5.8)$$ Since $\kappa$ is nondegenerate,
$\kappa(\xi,{\msr G}_0)\neq\{0\}$. So $\kappa(|_{{\msr
G}_0},|_{{\msr G}_0})$ is nondegenerate.$\qquad\Box$\psp

{\bf Lemma 2.5.3}.  ${\msr G}_0=H$.

{\it Proof}. For any $\xi\in {\msr G}_0$, we have the abstract
Jordan decomposition $\xi=\xi_s+\xi_n$. Note $\ad
\xi(H)=[\xi,H]=\{0\}$. Since $\ad\xi_s$ is a polynomial of $\ad\xi$
without constant, we have $[\xi_s,H]=\{0\}$. Thus $H+\mbb{C}\xi_s$
forms a toral subalgebra. But $H$ is a maximal toral subalgebra, we
have $\xi_s\in H$, which implies $\xi_n=\xi-\xi_s\in {\msr G}_0$
because $H\subset{\msr G}_0$. Since $\ad\xi_s|_{{\msr G}_0}=0$ and
$\ad\xi_n$ is nilpotent, $\ad\xi|_{{\msr G}_0}$ is nilpotent. By
Engel's Theorem, ${\msr G}_0$ is a nilpotent Lie algebra.

Suppose that $h\in H$ satisfies $\kappa(h,H)=\{0\}$. For any
$\xi\in{\msr G}_0$, let $\xi=\xi_s+\xi_n$ be the abstract Jordan
decomposition. Since $[\ad h,\ad\xi_n]=\ad[h,\xi_n]=0$, $\ad h\:\ad
\xi_n$ is a nilpotent element. Thus $\kappa(h,\xi_n)=0$. So
$\kappa(h,\xi)=\kappa(h,\xi_s)+\kappa(h,\xi_n)=0$; that is $h\in
\rad\;\kappa(|_{{\msr G}_0},|_{{\msr G}_0})=\{0\}$. Therefore,
$\kappa(|_H,|_H)$ is nondegenerate. For any $h\in H$ and
$\xi_1,\xi_2\in{\msr G}_0$,
$$\kappa(h,[\xi_1,\xi_2])=\kappa([h,\xi_1],\xi_2)=0.\eqno(2.5.9)$$
Hence
$$\kappa(H,[{\msr G}_0,{\msr G}_0])=\{0\}.\eqno(2.5.10)$$
If $[{\msr G}_0,{\msr G}_0]\neq\{0\}$, we take $0\neq \xi\in Z({\msr
G}_0)\bigcap [{\msr G}_0,{\msr G}_0]$. Let $\xi=\xi_s+\xi_n$ be the
abstract Jordan decomposition. Then $\xi_s\in H$ and
$\kappa(H,\xi_n)=\{0\}$. Hence
$$\kappa (H,\xi_s)=\kappa(H,\xi)=\{0\}\eqno(2.5.11)$$
by (2.5.10). The nondegeneracy of $\kappa(|_H,|_H)$ implies
$\xi_s=0$. So $\xi$ is ad-nilpotent. If $[{\msr G}_0,{\msr
G}_0]=\{0\}$ and $H\neq {\msr G}_0$, then $u_n$ is a nonzero
ad-nilpotent element for any $u\in{\msr G}_0\setminus H$. In
summary, if $H\neq{\msr G}_0$, we have a nonzero ad-nilpotent
element $\xi\in Z({\msr G}_0)$.  For any $\xi'\in{\msr G}_0$,
$\ad\xi'\:\ad\xi$ is nilpotent. So $\kappa(\xi',\xi)=0$; that is,
$\xi\in\rad\:\kappa(|_{{\msr G}_0},|_{{\msr G}_0})=\{0\}$, which is
absurd. Therefore $H={\msr G}_0.\qquad\Box$\psp

Set
$$\Phi=\{\al\in H^\ast\mid \al\neq 0,\;{\msr
G}_\al\neq\{0\}\}.\eqno(2.5.12)$$ The elements of $\Phi$ are called
{\it roots}\index{root} of ${\msr G}$ and $\Phi$ is called the {\it
root system} of ${\msr G}$. Since $\kappa(|_H,|_H)$ is
nondegenerate, for any $\al\in H^\ast$, there exists a unique
$t_\al\in H$ such that
$$\al(h)=\kappa(t_\al,h)\qquad\for\;\;h\in H.\eqno(2.5.13)$$
Moreover, $t_{-\al}=-t_\al$.\psp

{\bf Lemma 2.5.4}. {\it (a) $\Phi$ spans $H^\ast$.

(b) If $\al\in \Phi$, then $-\al\in\Phi$.

(c) Let $\al\in\Phi,\;\xi\in{\msr G}_\al,\;\zeta\in{\msr G}_{-\al}$.
Then $[\xi,\zeta]=\kappa(\xi,\zeta)t_\al$. In particular, $[{\msr
G}_\al,{\msr G}_{-\al}]=\mbb{C}t_\al$.

(d) If $\al\in\Phi$ and $0\neq \xi_\al\in{\msr G}_\al$, there exists
$\zeta_\al\in{\msr G}_{-\al}$ such that
$\{\xi_\al,t_\al,\zeta_\al\}$ spans a Lie subalgebra isomorphic to
$sl(2,\mbb{C})$}.

{\it Proof}. (a) If $\Phi$ does not span $H^\ast$, then there exists
$0\neq h\in H$ such that $\al(h)=0$ for any $\al\in \Phi$. By
(2.5.3), $h\in Z({\msr G})$, which is absurd.

(b) By (2.5.4), $\kappa({\msr G}_\al,{\msr G}_{-\al})\neq\{0\}$ if
$\al\in\Phi$.

(c) By Lemmas 2.5.2 and 2.5.3, $[\xi,\zeta]\in H$. For any $h\in H$,
we have
$$\kappa(h,[\xi,\zeta])=\kappa([h,\xi],\zeta)=\al(h)\kappa(\xi,\zeta)=
\kappa(h,t_\al)\kappa(\xi,\zeta),\eqno(2.5.14)$$ or equivalently,
$$\kappa(h,[\xi,\zeta]-\kappa(\xi,\zeta)t_\al)
=\kappa(h,[\xi,\zeta])-\kappa(h,t_\al)\kappa(\xi,\zeta)=0.\eqno(2.5.15)$$
Since $\kappa(|_H,|_H)$ is nondegenerate, we have
$[\xi,\zeta]-\kappa(\xi,\zeta)t_\al=0$. $[{\msr G}_\al,{\msr
G}_{-\al}]=\mbb{C}t_\al$ because $\kappa({\msr G}_\al,{\msr
G}_{-\al})\neq\{0\}$.

(d) By (2.5.3) and the nongeneracy of $\kappa$, there exist
$\zeta\in{\msr G}_{-\al}$ such that $\kappa(\xi_\al,\zeta)=1$. So
$$[\xi_\al,\zeta]=t_\al\eqno(2.5.16)$$
by (c). If $\al(t_\al)=0$, then
$$[t_\al,\xi_\al]=\al(t_\al)\xi_\al=0=-\al(t_\al)\zeta=[t_\al,\zeta].\eqno(2.5.17)$$
Thus ${\msr K}=\mbb{C}\xi_\al+\mbb{C}t_\al+\mbb{C}\zeta$ forms a
solvable Lie subalgebra of ${\msr G}$ and $t_\al\in [{\msr K},{\msr
K}]$. By Lie's theorem, $\ad t_\al$ is nilpotent. But $t_\al\in H$
is ad-semisimple. Hence $\ad t_\al=0$. Since $Z({\msr G})=\{0\}$, we
have $t_\al=0$. Thereby, $\al=0$, which is absurd. Hence $0\neq
\al(t_\al)=\kappa(t_\al,t_\al)$.

 Set
$$h_\al=\frac{2}{\al(t_\al)}t_\al,\qquad
\zeta_\al=\frac{2}{\al(t_\al)}\zeta.\eqno(2.5.18)$$ Then
$$[\xi_\al,\zeta_\al]=h_\al,\;\;[h_\al,\xi_\al]=2\xi_\al,\;\;[h_\al,\zeta_\al]=-2\zeta_\al.
\eqno(2.5.19)$$ Hence
$${\msr
S}_\al=\mbb{C}\xi_\al+\mbb{C}h_\al+\mbb{C}\zeta_\al\eqno(2.5.20)$$
forms a Lie subalgebra isomorphic to $sl(2,\mbb{C})$ via the linear
map determined by
$$\xi_\al\mapsto E_{1,2},\;\;h_\al\mapsto
E_{1,1}-E_{2,2},\;\;\zeta_\al\mapsto E_{2,1}.\eqno(2.5.21)$$ By
(2.5.18), $h_{-\al}=-h_\al.\qquad\Box$ \psp

{\bf Remark 2.5.5}. In [X7, X8], we constructed six families of
infinite-dimensional simple Lie algebras ${\msr G}$ without any
toral subalgebra $H$ such that  ${\msr G}_0=H$ (cf. (2.5.3)).

\section{Properties of Roots and Root Subspaces}

In this section, we  use finite-dimensional representations of
$sl(2,\mbb{C})$ to study the properties of roots and root subspaces
of finite-dimensional semisimple Lie algebra ${\msr G}$. In
particular, we prove that such a Lie algebra is generated by two
elements.

Again we assume the base field $\mbb{F}=\mbb{C}$. For convenience,
we also take $\ptl_x=d/dx$. For any nonnegative integer $n$, set
$$V(n)=\sum_{i=0}^n\mbb{C}x^i\subset \mbb{C}[x].\eqno(2.5.1)$$
Recall
$$sl(2,\mbb{C})=\mbb{C}E_{1,2}+\mbb{C}(E_{1,1}-E_{2,2})+\mbb{C}E_{2,1}.
\eqno(2.6.2)$$ We define an action of $sl(2,\mbb{C})$ on $V(n)$ by
$$E_{1,2}|_{V(n)}=x^2\ptl_x-nx,\;\;E_{2,1}|_{V(n)}=-\ptl_x,\;\;(E_{1,1}-E_{2,2})|_{V(n)}
=2x\ptl_x-n.\eqno(2.6.3)$$ It can be verified that $V(n)$ forms an
irreducible $sl(2,\mbb{C})$-module (exercise). For convenience, we
denote
$$e=E_{1,2},\;\;f=E_{2,1},\;\;h=E_{1,1}-E_{2,2}.\eqno(2.6.4)$$
For any $sl(2,\mbb{C})$-module $W$, we denote
$$W_a=\{w\in W\mid
h(w)=aw\}\qquad\for\;\;a\in\mbb{C}.\eqno(2.6.5)$$ For instance,
$$V(n)_{2i-n}=\mbb{C}x^i,\qquad
V(n)=\bigoplus_{i=0}^nV(n)_{2i-n}.\eqno(2.6.6)$$ \pse

{\bf Theorem 2.6.1}. {\it Any $(n+1)$-dimensional  irreducible
$sl(2,\mbb{C})$-module is isomorphic to $V(n)$. A finite-dimensional
$sl(2,\mbb{C})$-module $W$ is a direct sum of $\dim W_0+\dim W_1$
irreducible submodules}.

{\it Proof}. Let $V$ be any $(n+1)$-dimensional irreducible
irreducible $sl(2,\mbb{C})$-module. Note that ${\msr
B}=\mbb{C}h+\mbb{C}e$  forms a solvable Lie subalgebra of
$sl(2,\mbb{C})$. By Lie's Theorem, there exists a common eigenvector
$v$ of $h$ and $e$ in $V$; that is,
$$h(v)=\lmd v,\;\;e(v)=\mu v,\qquad\lmd,\mu\in\mbb{C}.\eqno(2.6.7)$$
But
$$2\mu v=2e(v)=[h,e](v)=h(e(v))-e(h(v))=\mu\lmd v-\lmd\mu
v=0.\eqno(2.6.8)$$ Thus $\mu=0$; that is, $e(v)=0$.

Recall that $\mbb{N}$ is the set of nonnegative integers.
 Note
$$h(f^i(v))=(\lmd-2i)f^i(v)\qquad\for\;\;i\in\mbb{N}.\eqno(2.6.9)$$
Moreover,
\begin{eqnarray*}\hspace{2cm}e(f^{i+1}(v))&=&e(f((f^i(v))))=[e,f](f^i(v))+f(e((f^i(v))))\\
&=&(\lmd-2i)f^i(v)+f(e((f^i(v))).
\hspace{4.7cm}(2.6.10)\end{eqnarray*} By induction, we have
$$e(f^{i+1}(v))=((i+1)\lmd-2\sum_{r=0}^ir)f^i(v)=(i+1)(\lmd-i)f^i(v).\eqno(2.6.11)$$
The above equation shows
$$ f^i(v)\neq 0\qquad \mbox{if}\;\;\lmd\not
\in\{0,1,...,i-1\}.\eqno(2.6.12)$$ On the other hand, (2.6.9)
implies that
$$\{v,f(v),...,f^i(v)\mid \lmd\not
\in\{0,1,...,i-1\}\}\eqno(2.6.13)$$ is a set of eigenvectors of $h$
with distinct eigenvalues. So it is linearly independent. Since
$\dim V=n+1$, $\lmd\in\{0,1,...,n\}$. By (2.6.11),
$e(f^{\lmd+1}(v))=0$. Hence
$$\sum_{i=1}^{\infty}\mbb{C}f^{\lmd+i}(v)\eqno(2.6.14)$$
forms a proper $sl(2,\mbb{C})$ submodule of $V$. The irreducibility
of $V$ implies
$$ f^{\lmd+1}(v)=0.\eqno(2.6.15)$$
Now $\sum_{i=0}^{\lmd}\mbb{C}f^i(v)$ forms a nonzero submodule of
$V$. Thus
$$V=\sum_{i=0}^{\lmd}\mbb{C}f^i(v)\eqno(2.6.16)$$
and $\{v,f(v),...,f^{\lmd}(v)\}$ is a basis by (2.6.9). Therefore,
$\lmd =n$.

By (2.6.3) and (2.6.4),
$$f^i(x^n)=(-1)^in(n-1)\cdots (n-i+1)x^{n-i}\qquad\for\;\;
i\in\{0,1,...,n\}.\eqno(2.6.17)$$ Define a linear map $\sgm:V(n)\rta
V$ by
$$\sgm(f^i(x^n))=f^i(v)\qquad\for\;\; i\in\{0,1,...,n\}.\eqno(2.6.18)$$
The equations (2.6.9), (2.6.11) and (2.6.15) show that  the linear
map $\sgm$ is $sl(2,\mbb{C})$-module isomorphism.

Expressions (2.6.5) and (2.6.6) imply that
$$ \dim V_0=0,\qquad\dim V_1=1\qquad\mbox{if}\;\;n\in
2\mbb{N}+1\eqno(2.6.19)$$ and
$$ \dim V_0=1,\qquad\dim V_1=0\qquad\mbox{if}\;\;n\in
2\mbb{N}.\eqno(2.6.20)$$  Now Weyl's Theorem tell us that any
finite-dimensional $sl(2,\mbb{C})$-module $W$ is a direct sum of
$\dim W_0+\dim W_1$ irreducible submodules.$\qquad\Box$\psp

Let ${\msr G}$ be a finite-dimensional semisimple Lie algebra and
let $H$ be a maximal toral subalgebra of ${\msr G}$. Then we have
the root space decomposition:
$${\msr G}=H+\sum_{\al\in\Phi}{\msr G}_\al,\eqno(2.6.21)$$
where $\Phi$ is the root system of ${\msr G}$. Recall Lemma 2.5.4.
\psp

{\bf Lemma 2.6.2}. {\it (a) $\dim {\msr G}_\al=1$  and
$\mbb{C}\al\bigcap \Phi=\{\al,-\al\}$ for any $\al\in\Phi$. In
particular, ${\msr S}_\al={\msr G}_\al+\mbb{C}h_\al+{\msr G}_{-\al}$
(cf. (2.5.20)).

(b) If $\al,\be\in\Phi$, then $\be(h_\al)\in\mbb{Z}$ (which is
called {\bf Cartan integer})\index{Cartan integer} and
$\be-\be(h_\al)\al\in\Phi$.

(c) If $\al,\be,\al+\be\in\Phi$, then $[{\msr G}_\al,{\msr
G}_\be]={\msr G}_{\al+\be}$.

(d) Let $\al,\be\in\Phi,\;\be\neq \pm\al$. Let $r,q$ be
(respectively) the largest integers such that $\be-r\al$ and
$\be+q\al$ are roots. Then all $\be+i\al\in\Phi$ for $-r\leq i\leq
q$.

(e) ${\msr G}$ is generated by $\{{\msr G}_\al\mid \al\in\Phi\}$.}

{\it Proof}. Denote $\dim H=n$. Let $\al\in \Phi$ be any element.
Recall the Lie subalgebra ${\msr S}_\al$ defined in (2.5.20). By the
above lemma, all the eigenvalues of $\ad h_\al$ are integers. If
$a\al\in\Phi$, then
$$a\al(h_\al)=2a\in\mbb{Z}\lra a\in\frac{\mbb{Z}}{2}.\eqno(2.6.22)$$
Note that
$${\msr L}=H+\sum_{m\in\mbb{Z}}{\msr G}_{m\al}\eqno(2.6.23)$$
forms an $\ad {\msr S}_\al$-submodule of ${\msr G}$. Moreover, the
eigenvalues of $\ad h_\al$ in ${\msr L}$ are even integers by
(2.6.22). Thus ${\msr L}$ is direct sum of $n$ irreducible ${\msr
S}_\al$-submodules. Set
$$H'=\{h\in H\mid \al(h)=0\}.\eqno(2.6.24)$$
For any $h\in H'$, we have
$$[h,\xi_\al]=\al(h)\xi_\al=0,\;\;[h,\zeta_\al]=-\al(h)\zeta_\al=0,\;\;[h,h_\al]=0\eqno(2.6.25)$$
(cf. Lemma 2.5.4). Since ${\msr
S}_\al=\mbb{C}\xi_\al+\mbb{C}h_\al+\mbb{C}\zeta_\al$, $H'$ is a
direct sum of $(n-1)$ one-dimensional irreducible $\ad{\msr
S}_\al$-submodules. Thus ${\msr S}_\al+H'$ is a direct sum of $n$
irreducible $\ad{\msr S}_\al$-submodules. Therefore, the complement
of ${\msr S}_\al+H'$ in ${\msr L}$ must be zero; that is, ${\msr
L}=H'+{\msr S}_\al$. So ${\msr S}_\al={\msr
G}_\al+\mbb{C}h_\al+{\msr G}_{-\al},\;\dim{\msr
G}_\al=\dim\mbb{C}\xi_\al=1$ and $\mbb{Z}\al\bigcap
\Phi=\{\al,-\al\}$. In particular, $2\al\not\in\Phi$ for any
$\al\in\Phi$. This shows $\al/2\not\in\Phi$; that is, ${\msr
G}_{\al/2}=\{0\}$. Since $\bar{\msr L}=\sum_{m\in\mbb{Z}}{\msr
G}_{m\al+\al/2}$ forms an $\ad{\msr S}_\al$-submodule, the
eigenvalues of $\ad h_\al$ in $\bar{\msr L}$ are odd integers and
the eigenspace of the eigenvalue 1 is ${\msr G}_{\al/2}=\{0\}$, we
must have $\bar{\msr L}=\{0\}$ by the above lemma. According to
(2.6.22), we get $\mbb{C}\al\bigcap \Phi=\{\al,-\al\}$ for any
$\al\in\Phi$. This proves (a).

Let $\al,\be\in\Phi,\;\be\neq\pm\al$. Then
$$ {\msr K}=\sum_{i\in\mbb{Z}}{\msr G}_{\be+i\al}\eqno(2.6.26)$$
forms an $\ad{\msr S}_\al$-submodule of ${\msr G}$. Since ${\msr
G}_\be$ is an eigenspace of $\ad h_\al$ with the eigenvalue
$\be(h_\al)$, we must have $\be(h_\al)\in\mbb{Z}$. By (a), $\dim
{\msr G}_{\be+i\al}\leq 1$. If $\be+i\al\in\Phi$, then ${\msr
G}_{\be+i\al}$ is an eigenspace of $\ad h_\al$ with the eigenvalue
$\be(h_\al)+2i$. Thus
$$\dim\{\xi\in {\msr K}\mid [h_\al,\xi]=0\}+\dim\{\zeta\in {\msr K}\mid
[h_\al,\zeta]=\zeta\}\leq 1.\eqno(2.6.27)$$ By the above lemma,
${\msr K}$ is an irreducible $\ad{\msr S}_\al$-module. Let $r,q$ be
(respectively) the largest integers such that $\be-r\al$ and
$\be+q\al$ are roots. Then all $\be+i\al\in\Phi$ for $-r\leq i\leq
q$ by (2.6.6) and the above lemma. This proves (d). If
$\al+\be\in\Phi$, then $[\xi_\al,{\msr G}_\be]={\msr G}_{\al+\be}$
by (2.6.11). So (c) holds. Expression (2.6.6) shows
$$-(\be(h_\al)-2r)=\dim{\msr K}-1=\be(h_\al)+2q,\eqno(2.6.28)$$
which implies $\be(h_\al)=r-q$. Furthermore,
$$\be-\be(h_\al)\al=\be-(r-q)\al=\be+q\al-r\al\in\Phi.\eqno(2.6.29)$$
Thus (b) holds. Since $\Phi$ spans $H^\ast$, $\{t_\al\mid
\al\in\Phi\}$ spans $H$. The equations in (2.5.18) and (2.5.19)
yield (e).$\qquad\Box$\psp

{\bf Proposition 2.6.3}. {\it The semisimple Lie algebra ${\msr G}$
is generated by two elements}.

{\it Proof}. Denote
$$\Psi=\{\al-\be\mid \al,\be\in\Phi,\;\al\neq \be\}.\eqno(2.6.30)$$
Take
$$h_0\in H\setminus \bigcup_{\gm\in\Psi}\{h\in
H\mid\gm(h)=0\}.\eqno(2.6.31)$$ For each $\al\in\Phi$, we pick
$0\neq\xi_\al\in{\msr G}_\al$. Then $\{\xi_\al\mid \al\in\Phi\}$ are
eigenvectors of $\ad h_0$ with distinct eigenvalues
$\{\al(h_0)\mid\al\in\Phi\}$. Set
$$\xi=\sum_{\al\in\Phi}\xi_\al.\eqno(2.6.32)$$
Let ${\msr L}$ be the Lie subalgebra of ${\msr G}$ generated by
$h_0$ and $\xi$. Then
$$\ad h_0 ({\msr L})\subset {\msr L},\;\; \ad h_0(\xi)\in{\msr L}.\eqno(2.6.33)$$ By Lemma 1.3.1,
$$\xi_\al\in {\msr L}\qquad\for\;\;\al\in \Phi.\eqno(2.6.34)$$
By Lemma 2.6.2 (e), ${\msr L}={\msr G}.\qquad\Box$\psp

Define
$$(\al,\be)=\kappa(t_\al,t_\be)\qquad \for\;\;\al,\be\in
H^\ast\eqno(2.6.35)$$ (cf. (2.5.13)). Then $(\cdot,\cdot)$ is a
nondegenerate symmetric bilinear form on $H^\ast$. By the arguments
above (2.5.18), $(\al,\al)=\al(t_\al)\neq 0$ for $\al\in\Phi$.
Moreover, we set
$$\la\be,\al\ra=\frac{2(\be,\al)}{(\al,\al)}\qquad\for\;\;\al\in
\Phi,\;\be\in H^\ast.\eqno(2.6.36)$$ By (2.5.18) and Lemma 2.6.2
(b),
$$\la\al,\be\ra\in\mbb{Z},\;\;\be-\la\be,\al\ra\al\in\Phi\qquad\for\;\;
\al,\be\in\Phi. \eqno(2.6.37)$$ Since $\Phi$ spans $H^\ast$, $\Phi$
contains a basis $\{\al_1,...,\al_n\}$ of $H^\ast$. For any
$\be\in\Phi$, $\be=\sum_{i=1}^nc_i\al_i$ with $c_i\in\mbb{C}$. Note
$$\sum_{i=1}^nc_i\la\al_i,\al_j\ra=\la\be,\al_j\ra,\qquad
j=1,2,...,n.\eqno(2.6.38)$$ View $\{c_1,c_2,...,c_n\}$ as unknowns.
The coefficients in the above linear systems are integers. The
coefficient determinant
$$|(\la\al_i,\al_j\ra)_{n\times
n}|=\frac{2^n}{\prod_{r=1}^n(\al_r,\al_r)}|((\al_i,\al_j))_{n\times
n}|\neq 0\eqno(2.6.39)$$ by the nondegeneracy of $(\cdot,\cdot)$ on
$H^\ast$. Solving (2.6.38) for $\{c_i\}$ in terms of
$\{\la\al_i,\al_j\ra,\la\be,\al_j\ra\}$, we obtain
$$c_i\in\mbb{Q},\qquad i=1,2,...,n.\eqno(2.6.40)$$
Thus the $\mbb{Q}$-subspace
$$\G=\sum_{\al\in\Phi}\mbb{Q}\al\eqno(2.6.41)$$
of $H^\ast$ is $n$-dimensional.

By (2.6.21),
$$(\lmd,\mu)=\tr\ad t_\lmd \ad t_\mu=\sum_{\al\in
\Phi}\al(t_\lmd)\al(t_\mu)=\sum_{\al\in
\Phi}(\al,\lmd)(\al,\mu)\qquad\for\;\;\lmd,\mu\in
H^\ast.\eqno(2.6.42)$$ In particular,
$$(\be,\be)=\sum_{\al\in\Phi}(\al,\be)^2\qquad\for\;\;\be\in\Phi,\eqno(2.6.43)$$
equivalently,
$$\frac{1}{(\be,\be)}=\frac{1}{4}\sum_{\al\in\Phi}\left(\frac{2(\al,\be)}{(\be,\be)}\right)^2
=\frac{1}{4}\sum_{\al\in\Phi}\la\al,\be\ra^2\in\mbb{Q}.\eqno(2.6.44)$$
So $(\be,\be)\in\mbb{Q}$. For $\al,\be\in\phi$, we have
$\al-\la\al,\be\ra\be\in\Phi$ by the above lemma. Hence
$$\la\al,\be\ra(\al,\be)=\frac{1}{2}[(\al,\al)+\la\al,\be\ra^2(\be,\be)
-(\al-\la\al,\be\ra\be,\al-\la\al,\be\ra\be)]\in\mbb{Q}.\eqno(2.6.45)$$
Therefore,
$$(\al,\be)\in\mbb{Q}\qquad\for\;\;\al,\be\in\Phi.\eqno(2.6.46)$$
Thus $(\cdot,\cdot)$ is a positive definite $\mbb{Q}$-valued
symmetric $\mbb{Q}$-bilinear form on $\G$ by (2.6.42) and (2.6.46).
Extend $(\cdot,\cdot)$ on
$${\msr E}=\G_{\mbb{R}}=\mbb{R}\otimes_{\mbb{Q}}\G\eqno(2.6.47)$$
$\mbb{R}$-bilinearly. Then ${\msr E}$ is isomorphic to the Euclidean
space $\mbb{R}^n$. Now we obtain  a main theorem.\psp

{\bf Theorem 2.6.4}. {\it Let ${\msr G}$ be a finite-dimensional
semisimple Lie algebra. There exists a maximal toral subalgebra $H$,
whose dimension denoted by $n$. With respect to $H$, ${\msr G}$ has
a root space decomposition ({\bf Cartan decomposition})\index{Cartan
root space decomposition} (2.6.21), where
 $\dim {\msr G}_\al=1$ and if $\al,\be,\al+\be\in\Phi$, then $[{\msr G}_\al,{\msr
G}_\be]={\msr G}_{\al+\be}$.  Moreover, $\Phi$ can be identified
with a finite subset of  the Euclidean space $\mbb{R}^n$ satisfying:
(a) $\Phi$ spans $\mbb{R}^n$ and $0\not\in\Phi$; (b)
$\mbb{R}\al\bigcap \Phi=\{\al,-\al\}$ for $\al\in\Phi$; (c) If
$\al,\be\in\Phi$, then $2(\be,\al)/(\al,\al)\in\mbb{Z}$ and
$\be-[2(\be,\al)/(\al,\al)]\al\in\Phi.$}\psp

As an exercise, find all finite-dimensional real irreducible modules
for $o(3,\mbb{R})$ and $o(4,\mbb{R})$ defined in (1.2.22).

\chapter{Root Systems}

 In this chapter, we start with the axiom of root system and
 give the root systems of special linear algebras, orthogonal Lie
 algebras and symplectic Lie algebras. Then we derive some basic
 properties of root systems; in particular, the existence of the bases of root systems.
 As finite symmetries of root systems, the Weyl
 groups are introduced and studied in detail. Classification and explicit construction of root
 systems are presented. The automorphism groups of roots systems are
 determined. As a preparation for later representation theory of Lie
 algebras, the corresponding weight lattices and their saturated
 subsets are investigated.

\section{Definitions, Examples and Properties}

In this section, we present the basic definitions on root systems,
examples and properties.

An $n$-dimensional Euclidean space ${\msr E}$ is the vector space
$\mbb{R}^n$ with a positive definite symmetric $\mbb{R}$-bilinear
form $(\cdot,\cdot)$ (called {\it inner product}).\index{inner
product} For $\al,\be\in {\msr E}$ with $\al\neq 0$, we define
$$\la\be,\al\ra=\frac{2(\be,\al)}{(\al,\al)}\eqno(3.1.1)$$
and
$$P_\al=\{\gm\in{\msr E}\mid (\gm,\al)=0\}.\eqno(3.1.2)$$
Note that the number $\la\be,\al\ra$  and $P_\al$ remain unchanged
if we replace $(\cdot,\cdot)$ by $c(\cdot,\cdot)$ for any $0\neq
c\in\mbb{R}$. Moreover, $(\al,\al)>0$ implies
$${\msr E}=\mbb{R}\al+P_\al.\eqno(3.1.3)$$
Furthermore, we define a linear transformation $\sgm_\al$ by
$$\sgm_\al(\be)=\be-\la\be,\al\ra\al\qquad\for\;\;\be\in{\msr
E}.\eqno(3.1.4)$$ Then $\sgm_\al|_{P_\al}=\mbox{Id}_{P_\al}$ and
$$\sgm_\al(\al)=\al-\la\al,\al\ra\al=\al-\frac{2(\al,\al)}{(\al,\al)}\al=-\al.\eqno(3.1.5)$$
So $\sgm_\al$ is the reflection with respect to the hyperplane
$P_\al$. Moreover, $\sgm_{a\al}=\sgm_\al$ for any $0\neq
a\in\mbb{R}$. For $\be,\gm\in{\msr E}$, we write $\be=\be_1+\be_2$
and $\gm=\gm_1+\gm_2$ with $\be_1,\gm_1\in\mbb{R}\al$ and
$\be_2,\gm_2\in P_\al$, and get
$$\la\sgm_\al(\be),\sgm_\al(\gm)\ra=\la
-\be_1+\be_2,-\gm_1+\gm_2\ra=\la\be_1,\gm_1\ra+\la\be_2,\gm_2\ra=\la\be,\gm\ra.
\eqno(3.1.6)$$ Thus $\sgm_\al$ is an isometry.\psp

A {\it root system}\index{root system} $\Phi$ in ${\msr E}$ is a
finite set spanning ${\msr E}$ such that: (1) $0\not\in\Phi$ and
$\mbb{R}\al\bigcap\Phi=\{\al,-\al\}$ for $\al\in\Phi$; (2)
$\la\be,\al\ra\in\mbb{Z}$ and $\sgm_\al(\Phi)=\Phi$ for
$\al,\be\in\Phi$. \psp

Recall that $h_\al$ in a semisimple Lie algebra ${\msr
G}=H+\sum_{\al\in\Phi}{\msr G}_\al$ is uniquely determined by
$$[{\msr G}_\al,{\msr G}_{-\al}]=\mbb{C}h_\al,\;\;\ad h_\al|_{{\msr
G}_\al}=2\mbox{Id}_{{\msr G}_\al}.\eqno(3.1.7)$$ Moreover,
$$\la\be,\al\ra=\be(h_\al)\eqno(3.1.8)$$
by (2.5.18), (2.5.19), (2.6.35) and (2.6.36).\psp

{\bf Example 3.1.1}. The special linear Lie algebra
$$sl(n,\mbb{C})=\sum_{1\leq i<j\leq
n}(\mbb{C}E_{i,j}+\mbb{C}E_{j,i})+\sum_{i=1}^{n-1}\mbb{C}(E_{i,i}-E_{i+1,i+1})\eqno(3.1.9)$$
and recall the notion of indices: $\ol{1,n}=\{1,2,...,n\}.$ Take
$$H=\sum_{i=1}^{n-1}\mbb{C}(E_{i,i}-E_{i+1,i+1})=\{\sum_{i=1}^na_iE_{i,i}\mid
a_i\in \mbb{C},\;\sum_{j=1}^na_j=0\}.\eqno(3.1.10)$$ For
$i\in\ol{1,n}$, we define $\ves_i\in H^\ast$ by
$$\ves_i(\sum_{j=1}^na_jE_{j,j})=a_i.\eqno(3.1.11)$$
Then the root system  of $sl(n,\mbb{C})$ is
$$\Phi_{A_{n-1}}=\{\ves_i-\ves_j\mid i,j\in\ol{1,n},\;i\neq j\}\eqno(3.1.12)$$
and $${\msr
G}_{\ves_i-\ves_j}=\mbb{C}E_{i,j},\;\;h_{\ves_i-\ves_j}=E_{i,i}-E_{j,j},\;\;\;i,j\in\ol{1,n},\;i\neq
j.\eqno(3.1.13)$$ \vspace{0.1cm}

{\bf Example 3.1.2}. The even orthogonal Lie algebra
\begin{eqnarray*}\hspace{1cm} & &o(2n,\mbb{C})=\sum_{i,j=1}^n(\mbb{C}(E_{i,j}-E_{n+j,n+i})\\ & &+
\sum_{1\leq i< j\leq
n}[\mbb{C}(E_{i,n+j}-E_{j,n+i})+\mbb{C}(E_{n+j,i}-E_{n+i,j})].
\hspace{3.9cm}(3.1.14)\end{eqnarray*} Take
$$H=\sum_{i=1}^n\mbb{C}(E_{i,i}-E_{n+i,n+i})\eqno(3.1.15)$$
and define
$$\ves_i(\sum_{j=1}^na_j(E_{j,j}-E_{n+j,n+j}))=a_i\qquad\for\;\;i\in\ol{1,n}.
\eqno(3.1.16)$$
 Then the
root system of $o(2n,\mbb{C})$ is
$$\Phi_{D_n}=\{\ves_i-\ves_j,\pm(\ves_i+\ves_j)\mid
i,j\in\ol{1,n},\;i\neq j\}.\eqno(3.1.17)$$ Moreover,
$${\msr G}_{\ves_i-\ves_j}=\mbb{C}(E_{i,j}-E_{n+j,n+i}),\;\;{\msr
G}_{\ves_i+\ves_j}=\mbb{C}(E_{i,n+j}-E_{j,n+i}),\eqno(3.1.18)$$
$${\msr G}_{-\ves_i-\ves_j}=\mbb{C}(E_{n+i,j}-E_{n+j,i})\eqno(3.1.19)$$ and
$$h_{\ves_i-\ves_j}=E_{i,i}-E_{j,j}-E_{n+i,n+i}+E_{n+j,n+j},\eqno(3.1.20)$$
$$h_{\ves_i+\ves_j}=-h_{-\ves_i-\ves_j}=E_{i,i}+E_{j,j}-E_{n+i,n+i}
-E_{n+j,n+j},\eqno(3.1.21)$$ for $i,j\in\ol{1,n},\;i\neq j$.\psp

{\bf Example 3.1.3}. The symplectic Lie algebra
\begin{eqnarray*} \hspace{1cm}& &sp(2n,\mbb{C})=\sum_{i,j=1}^n\mbb{C}(E_{i,j}-E_{n+j,n+i})\\ & &
+\sum_{1\leq i\leq j\leq
n}[\mbb{C}(E_{i,n+j}+E_{j,n+i})+\mbb{C}(E_{n+j,i}+E_{n+i,j})].\hspace{3.9cm}(3.1.22)
\end{eqnarray*}
Take $H$ in (3.1.15) and $\ves_i$ in (3.1.16).  Then the root system
of $sp(2n,\mbb{C})$ is
$$\Phi_{C_n}=\{\ves_i-\ves_j,\pm(\ves_i+\ves_j),\pm 2\ves_i\mid
i,j\in\ol{1,n},\;i\neq j\}.\eqno(3.1.23)$$ Moreover,
$${\msr G}_{\ves_i-\ves_j}=\mbb{C}(E_{i,j}-E_{n+j,n+i}),\;\;{\msr
G}_{\ves_i+\ves_j}=\mbb{C}(E_{i,n+j}+E_{j,n+i}),\eqno(3.1.24)$$
$${\msr G}_{-\ves_i-\ves_j}=\mbb{C}(E_{n+i,j}+E_{n+j,i}),\;\;{\msr G}_{2\ves_i}=\mbb{C}E_{i,n+i},\;\; {\msr
G}_{-2\ves_i}=\mbb{C}E_{n+i,i},\eqno(3.1.25)$$
$$h_{2\ves_i}=-h_{-2\ves_i}=E_{i,i}-E_{n+i,n+i},\eqno(3.1.26)$$
 (3.1.20) and (3.1.21) hold, for $i,j\in\ol{1,n},\;i\neq j$.\psp

{\bf Example 3.1.4}. The odd orthogonal Lie algebra
$$o(2n+1,\mbb{C})=o(2n,\mbb{C})+\sum_{i=1}^n[\mbb{C}(E_{0,i}-E_{n+i,0})+\mbb{C}(E_{0,n+i}-E_{i,0})].\eqno(3.1.27)$$
Take $H$ in (3.1.15) and $\ves_i$ in (3.1.16).  Then the root system
of $o(2n+1,\mbb{C})$ is
$$\Phi_{B_n}=\{\ves_i-\ves_j,\pm(\ves_i+\ves_j),\pm\ves_i\mid
i,j\in\ol{1,n},\;i\neq j\}.\eqno(3.1.28)$$
 Moreover, (3.1.18)-(3.1.21) hold and
$${\msr G}_{-\ves_i}=\mbb{C}(E_{0,i}-E_{n+i,0}),\;\;{\msr
G}_{\ves_i}=\mbb{C}(E_{0,n+i}-E_{i,0}),\eqno(3.1.29)$$
$$h_{\ves_i}=-h_{-\ves_i}=2(E_{i,i}
-E_{n+i,n+i}),\eqno(3.1.30)$$  for $i\in\ol{1,n}$.\psp

Recall that the symmetric bilinear form $(\cdot,\cdot)$ on $H^*$ for
a semismiple Lie algebra ${\msr G}$ is induced from its Killing form
$\kappa$. When ${\msr G}$ is simple, any invariant (associative)
symmetric bilinear form is a scalar multiple of $\kappa$.  In a root
system $\Phi$, only integers $\la\be,\al\ra$ for $\al,\be\in\Phi$
play roles in its structure. But the number $\la\be,\al\ra$ and
$P_\al$ remain unchanged if we replace $(\cdot,\cdot)$ by
$c(\cdot,\cdot)$ for any $0\neq c\in\mbb{R}$. In the above four
examples, the bilinear form $\zeta(A,B)=\tr AB$ is a scalar multiple
of $\kappa$. Replacing $\kappa$ by $\zeta$ in Example 3.1.1, and by
$\zeta/2$ in Examples 3.1.2, 3.1.3 and 3.1.4, we have
$$(\ves_i,\ves_j)=\dlt_{i,j}\qquad\for\;\;i,j\in\ol{1,n}.\eqno(3.1.31)$$\pse

For any $\al,\be$ in the Euclidean space ${\msr E}$, we have
$$|(\al,\be)|\leq
||\al||\:||\be||,\qquad\mbox{where}\;\;||\al||=\sqrt{(\al,\al)},\;||\be||=\sqrt{(\be,\be)}.
\eqno(3.1.32)$$ We define the angle $\theta$ between two nonzero
vectors $\al$ and $\be$ by
$$\cos\theta=\frac{(\al,\be)}{
||\al||\:||\be||}.\eqno(3.1.33)$$ Then
$$\la\al,\be\ra\la\be,\al\ra=4\left(\frac{(\al,\be)}{
||\al||\:||\be||}\right)^2=4\cos^2\theta.\eqno(3.1.34)$$ Moreover,
$$\frac{||\al||^2}{||\be||^2}=\frac{\la\al,\be\ra}{\la\be,\al\ra}\qquad\mbox{if}
\;\;(\al,\be)\neq 0. \eqno(3.1.35)$$ Thus we have:\psp

{\bf Lemma 3.1.1}. {\it Let $\al,\be$ be two elements in a root
system such that $\al\neq\pm\be$ and $||\al||\leq ||\be||$. If
$(\al,\be)\neq 0$, then one of the following statements holds}:
$$\begin{tabular}{|c|c|c|c|c|}\hline\mbox{\it
Cases}&$\la\al,\be\ra$&$\la\be,\al\ra$&$\theta$&$||\al||^2/||\be||^2$\\\hline\hline
(1)&1&1&$\pi/3$&1\\ \hline (2)& -1&-1&$2\pi/3$&1\\ \hline
(3)&1&2&$\pi/4$&1/2\\ \hline (4)&-1&-2&$3\pi/4$&1/2 \\
\hline(5)&1&3&$\pi/6$&1/3\\ \hline(6)&-1&-3&$5\pi/6$&1/3\\
\hline\end{tabular}.\eqno(3.1.36)$$

{\it Proof}. By (3.1.34), (3.1.35) and the fact
$0<\la\al,\be\ra\la\be,\al\ra\in\mbb{Z}.\qquad\Box$\psp

{\bf Lemma 3.1.2}. {\it Let $\al,\be$ be two roots in a root system
$\Phi$ such that $\al\neq\pm\be$. If $(\al,\be)>0$, then
$\al-\be\in\Phi$. When $(\al,\be)<0$, $\al+\be\in\Phi$}.

{\it Proof}. If $(\al,\be)>0$, then $\la\al,\be\ra=1$ or
$\la\be,\al\ra=1$ by (3.1.36). Moreover,
$$\la\al,\be\ra=1\lra
\sgm_\be(\al)=\al-\la\al,\be\ra\be=\al-\be\in\Phi,\eqno(3.1.37)$$
$$\la\be,\al\ra=1\lra
\sgm_\al(\be)=\be-\la\be,\al\ra\al=\be-\al\in\Phi\lra
\al-\be=-(\be-\al)\in\Phi.\eqno(3.1.38)$$

When $(\al,\be)<0$, then $\la\al,\be\ra=-1$ or $\la\be,\al\ra=-1$.
Moreover,
$$\la\al,\be\ra=-1\lra
\sgm_\be(\al)=\al-\la\al,\be\ra\be=\al+\be\in\Phi,\eqno(3.1.39)$$
$$\la\be,\al\ra=-1\lra
\sgm_\al(\be)=\be-\la\be,\al\ra\al=\be+\al\in\Phi.\qquad\Box\eqno(3.1.40)$$
\vspace{0.1cm}

{\bf Lemma 3.1.3}. {\it Let $\al,\be$ be two roots in a root system
$\Phi$ such that $\al\neq\pm\be$. Let $r$ and $q$ be the largest
integers such that $\be-r\al,\be+q\al\in\Phi$. Then
$\be+i\al\in\Phi$ for any $-r\leq i\leq q$, and
$\la\be,\al\ra=r-q$.}

{\it Proof}. Suppose that $\be+i\al\not\in\Phi$ for some $-r<i<q$.
There exists integers $s,p$ with $-r\leq s<i<p\leq q$ such that
$$\be+s\al,\be+p\al\in\Phi\;\;\mbox{and}\;\;\be+(s+1)\al,\be+(p-1)\al\not\in\Phi.
\eqno(3.1.41)$$ By the above lemma,
$$(\be+s\al,\al)\geq 0\geq (\be+p\al,\al)\lra s(\al,\al)\geq
p(\al,\al)\sim s>p,\eqno(3.1.42)$$ which is absurd. Thus
$\be+i\al\in\Phi$ for any $-r\leq i\leq q$.

Note that
$$\sgm_\al(\be+i\al)=\be-\la\be,\al\ra\al-i\al\in\Phi\qquad\for\;\;-r\leq
i\leq q.\eqno(3.1.43)$$ Thus $\sgm_\al$ flips the $\al$-string
through $\be$:
$$(\be-r\al,...,\be,...,\be+q\al)\mapsto
(\be-(\la\be,\al\ra+q)\al,....,\be,...,\be+(r-\la\be,\al\ra)\al).\eqno(3.1.44)$$
Hence
$$\la\be,\al\ra+q=r\lra
\la\be,\al\ra=r-q.\qquad\Box\eqno(3.1.45)$$ \vspace{0.1cm}

Suppose that $\Phi$ is a root system in a Euclidean space ${\msr
E}$. A subset $\Pi$ of $\Phi$ is called a {\it base}\index{base} if
$\Pi$ is a basis of ${\msr E}$ and $\be=\sum_{\al\in\Pi}k_\al\al$
for any $\be\in \Phi$ such that all $ k_\al\in\mbb{N}$ or all $
-k_\al\in\mbb{N}$, where $\mbb{N}$ is the set of all nonnegative
integers.

 A vector $\gm\in{\msr E}\setminus \bigcup_{\al\in\Phi}P_\al$ is called {\it
 regular}. Since $\dim P_\al=\dim{\msr E}-1$, ${\msr E}\neq
 \bigcup_{\al\in\Phi}P_\al$. So regular vector exists. For a
 regular vector $\gm\in{\msr E}$, we define
 $$\Phi^+(\gm)=\{\al\in\Phi\mid (\gm,\al)>0\},\;\;\Phi^-(\gm)=-\Phi^+(\gm).\eqno(3.1.46)$$
 Then
 $$\Phi=\Phi^+(\gm)\bigcup \Phi^-(\gm).\eqno(3.1.47)$$
 An element $\be\in \Phi^+(\gm)$ is called {\it decomposable} if
 there exist $\be_1,\be_2\in  \Phi^+(\gm)$ such that
 $\be=\be_1+\be_2$. In this case,
 $(\be,\gm)=(\be_1,\gm)+(\be_2,\gm)$. In particular,
 $(\be_1,\gm),(\be_2,\gm)<(\be,\gm)$. Thus the elements
 $\al\in\Phi^+(\gm)$ with minimal $(\al,\gm)$ are indecomposable
 elements. Set
 $$\Pi(\gm)=\mbox{the set of all indecomposble elements
 of}\;\;\Phi^+(\gm).\eqno(3.1.48)$$
Since $\{(\al,\gm)\mid\al\in\Phi\}$ is a finite set of positive
numbers, any element $\be\in\Phi^+(\gm)$ can be written as
$$\be=\sum_{\al\in\Pi(\gm)}k_\al\al,\qquad
k_\al\in\mbb{N},\eqno(3.1.49)$$ by induction on $(\be,\gm)$.
 \psp

{\bf Theorem 3.1.4}. {\it  A subset of $\Phi$ is base if and only if
it is of the form $\Pi(\gm)$}.

{\it Proof}. Let $\al,\be\in \Pi(\gm)$ such that $\al\neq \be$. If
$\al-\be\in\Phi$, then either $\al-\be\in\Phi^+(\gm)$ or
$\be-\al=-(\al-\be)\in\Phi^+(\gm)$, which is impossible because
$$\al=(\al-\be)+\be,\qquad\be=(\be-\al)+\al.\eqno(3.1.50)$$
By Lemma 3.1.2,
$$(\al,\be)\leq 0\qquad\for\;\;\al,\be\in\Pi(\gm),\;\al\neq\be.\eqno(3.1.51)$$ Suppose
$$\sum_{\al\in\Pi(\gm)}a_\al\al=0\qquad\mbox{with}\;\;a_\al\in\mbb{R}.\eqno(3.1.52)$$
Write
$$\ves=\sum_{\al\in\Pi(\gm),\;a_\al\geq
0}a_\al\al=-\sum_{\be\in\Pi(\gm),\;a_\be\leq0}a_\be\be.\eqno(3.1.53)$$
By (3.1.53), we have
\begin{eqnarray*}\hspace{1cm}(\ves,\ves)&=&(\sum_{\al\in\Pi(\gm),\;a_\al\geq
0}a_\al\al,-\sum_{\be\in\Pi(\gm),\;a_\be\leq 0}a_\be\be)\\ &
=&\sum_{\al\in\Pi(\gm),\;a_\al\geq
0}\;\sum_{\be\in\Pi(\gm),\;a_\be\leq 0}a_\al(-a_\be)(\al,\be)\leq
0.\hspace{4.2cm}(3.1.54)\end{eqnarray*} Thus $(\ves,\ves)=0$; that
is,
$$\sum_{\al\in\Pi(\gm),\;a_\al\geq
0}a_\al\al=-\sum_{\be\in\Pi(\gm),\;a_\be\leq
0}a_\be\be=\ves=0.\eqno(3.1.55)$$ Hence
$$\sum_{\al\in\Pi(\gm),\;a_\al\geq
0}a_\al(\gm,\al)=-\sum_{\be\in\Pi(\gm),\;a_\be\leq
0}a_\be(\gm,\be)=0,\eqno(3.1.56)$$ or equivalently,
$$a_\al=0\qquad\for\;\;\al\in\Pi(\gm).\eqno(3.1.57)$$
Therefore, $\Pi(\gm)$ is linearly independent.  Thus $\Pi(\gm)$ is a
base by (3.1.47) and (3.1.49) because $\Phi$ spans ${\msr E}$.

Conversely, we  assume that $\Pi$ is a base of $\Phi$. Denote
$$\Pi=\{\al_1,\al_2,...,\al_n\}.\eqno(3.1.58)$$
So $\dim{\msr E}=n$. Set
$${\msr E}_i=\{\be\in{\msr E}\mid (\be,\al_j)=0\;\for\;\;i\neq
j\in\ol{1,n}\}.\eqno(3.1.59)$$ Then $\dim {\msr E}_i=1$ and we can
take $\gm_i\in{\msr E}_i$ such that $(\gm_i,\al_i)>0$. Denote
$$\gm=\sum_{i=1}^n\gm_i.\eqno(3.1.60)$$
We have
$$(\gm,\al_i)>0\qquad\for\;\;i\in\ol{1,n};\eqno(3.1.61)$$ that is, $\Pi\subset\Phi^+(\gm)$.
If some $\al_i=\be_1+\be_2$ with $\be_1,\be_2\in\Phi^+(\gm)$, then
$$\be_1=\sum_{j=1}^nl_j\al_j,\;\;\be_2=\sum_{j=1}^nk_j\al_j,\qquad
l_j,k_j\in\mbb{N},\eqno(3.1.62)$$ by  (3.1.61) and the definition of
base. Note
$$\al_i=\sum_{j=1}^n(l_j+k_j)\al_j\lra \sum_{j=1}^n(l_i+k_i)=1\lra
\be_1=0\;\mbox{or}\;\be_2=0,\eqno(3.1.63)$$ which is absurd. Hence
all $\al_i$ are indecomposable; that is, $\Pi\subset \Pi(\gm)$. But
$|\Pi(\gm)|=\dim{\msr E}=|\Pi|$. Therefore,
$\Pi=\Pi(\gm).\qquad\Box$.\psp

If a base $\Pi=\Pi(\gm)$, we simply denote $\Phi^+=\Phi^+(\gm)$. We
call the elements in $\Phi^+$ {\it positive roots}. Correspondingly,
the elements in $\Phi^-=-\Phi^+$ are called {\it negative roots}.
The elements in $\Pi$ are called {\it simple roots}.\psp

{\bf Lemma 3.1.5}. {\it If $\al$ is positive but not simple, then
there exists $\be\in\Pi$ such that $\al-\be\in\Phi^+$.}

{\it Proof}. Otherwise, we have $(\al,\be)\leq 0$ for any
$\be\in\Pi$. By the arguments in the above theorem,
$\{\be\}\bigcup\Pi$ is linearly independent, which constricts to the
fact that $\Pi$ is base.$\qquad\Box$\psp

This lemma shows:\psp

{\bf Corollary 3.1.6}. {\it Any positive root
$\al=\sum_{i=1}^k\be_i$ with $\be_i\in\Pi$ such that
$\sum_{i=1}^s\be_i$ are positive roots for $s\in\ol{1,k}$.} \psp

For any $\al$ in a root system $\Phi$, we define
$$\check{\al}=\frac{2\al}{(\al,\al)}.\eqno(3.1.64)$$
Then
$$\check{\Phi}=\{\check{\al}\mid\al\in\Phi\}\eqno(3.1.65)$$
forms a root system (exercise).

If $\Phi$ satisfies all the conditions of root system except
$\mbb{R}\al\bigcap\Phi=\{\pm\al\}$, then the only possible multiples
of $\al\in\Phi$ in $\Phi$ are $\{\pm \al/2,\pm\al,\pm 2\al\}$ and
$\{\al\in\Phi\mid 2\al\not\in\Phi\}$ forms a root system (exercise).

Let $\al,\be$ be two elements in a root system $\Phi$ such that
$\al\neq\pm \be$. Suppose that the $\al$-string through $\be$ is
$(\be-r\al,...,\be+q\al)$ and  the $\be$-string through $\al$ is
$(\al-r'\be,...,\al+q\be)$. As an exercise, prove
$$\frac{q(r+1)}{(\be,\be)}=\frac{q'(r'+1)}{(\al,\al)}.\eqno(3.1.66)$$

\section{Weyl Groups}

In this section, we  introduce the Weyl group associated to each
root system and use it to study the structure of the root system.

Let $\Phi$ be a root system in a Euclidean space ${\msr E}$. Denote
by $GL({\msr E})$ the group of all invertible linear transformations
on ${\msr E}$. The {\it Weyl group}\index{Weyl group} associated to
$\Phi$ is
$${\msr W}=\mbox{the subgroup of}\;GL({\msr E})\;\mbox{generated
by}\:\{\sgm_\al\mid \al\in\Phi\}.\eqno(3.2.1)$$ Since
$\sgm_\al(\Phi)=\Phi$ and $\Phi$ spans ${\msr E}$, ${\msr W}$ is a
subgroup of the permutation group on $\Phi$. So ${\msr W}$ is a
finite group.\psp

{\bf Lemma 3.2.1}. {\it Let $S$ be a finite subset spanning ${\msr
E}$. Suppose $\sgm_\al(S)=S$ for some $0\neq\al\in{\msr E}$. If an
element $\sgm\in GL({\msr E})$ satisfies that
$\sgm(\al)=-\al,\;\sgm(S)=S$ and $\sgm|_P=\mbox{Id}_P$ for  some
hyperplane $P$ in ${\msr E}$, then $\sgm=\sgm_\al$ and $P=P_\al$.}

{\it Proof}. Note that
$${\msr E}=P+\mbb{R}\al=P_\al+\mbb{R}\al.\eqno(3.2.2)$$
So
$$\sgm|_{{\msr E}/\mbb{R}\al}=\mbox{Id}_{{\msr E}/\mbb{R}\al},\;\;
\sgm_\al|_{{\msr E}/\mbb{R}\al}=\mbox{Id}_{{\msr
E}/\mbb{R}\al}.\eqno(3.2.3)$$ Note
$$\sgm\sgm_\al(\al)=\sgm(-\al)=\al,\;\;\sgm\sgm_\al|_{{\msr E}/\mbb{R}\al}=\mbox{Id}_{{\msr
E}/\mbb{R}\al}.\eqno(3.2.4)$$ Take a basis ${\msr B}$ of $P$. Then
${\msr B}\bigcup\{\al\}$ forms a basis of ${\msr E}$, under which
$$\sgm\sgm_\al=\left(\begin{array}{cc}I&A\\
0&1\end{array}\right),\eqno(3.2.5)$$ where $I$ is an identity
matrix. Since $\sgm\sgm_\al(S)=S$, it is of finite order; say, $k$.
Then
$$\left(\begin{array}{cc}I&0\\
0&1\end{array}\right)=(\sgm\sgm_\al)^k=\left(\begin{array}{cc}I&A\\
0&1\end{array}\right)^k=\left(\begin{array}{cc}I&kA\\
0&1\end{array}\right).\eqno(3.2.6)$$ Thus $A=0$, or equivalently,
$\sgm\sgm_\al=1$. Hence $\sgm=\sgm_\al^{-1}=\sgm_\al$. So
$P=(1+\sgm)(V)=(1+\sgm_\al)(V)=P_\al.\qquad\Box$\psp

{\bf Lemma 3.2.2}. {\it Let $\Phi$ be a root system in a Euclidean
space ${\msr E}$. If an element $\sgm\in GL({\msr E})$ satisfies
$\sgm(\Phi)=\Phi$, then $\sgm\sgm_\al\sgm^{-1}=\sgm_{\sgm(\al)}$ for
any $\al\in \Phi$, and $\la\sgm(\al),\sgm(\be)\ra=\la \al,\be\ra$
for any $\al,\be\in \Phi$  (such $\sgm$ is called an {\bf
automorphism} of $\Phi$).}

{\it Proof}. First we have
$$\sgm\sgm_\al\sgm^{-1}(\sgm(\al))=\sgm\sgm_\al(\al)=\sgm(-\al)=-\sgm(\al).\eqno(3.2.7)$$
Moreover, for any $\gm\in P_\al$, we have
$$\sgm\sgm_\al\sgm^{-1}(\sgm(\gm))=\sgm\sgm_\al(\gm)=\sgm(\gm);\eqno(3.2.8)$$
that is,
$\sgm\sgm_\al\sgm^{-1}|_{\sgm(P_\al)}=\mbox{Id}_{\sgm(P_\al)}$.
Since $\sgm\sgm_\al\sgm^{-1}(\Phi)=\Phi$ and $\sgm(P_\al)$ is a
hyperplane, we have $\sgm\sgm_\al\sgm^{-1}=\sgm_{\sgm(\al)}$ by the
above lemma.

Now for any $\be\in\Phi$, we have
$$\sgm\sgm_\al\sgm^{-1}(\sgm(\be))=\sgm(\be-\la\be,\al\ra\al)=\sgm(\be)-\la\be,\al\ra\sgm(\al),
\eqno(3.2.9)$$ which is equal to
$$\sgm_{\sgm(\al)}(\sgm(\be))=\sgm(\be)-\la\sgm(\be),\sgm(\al)\ra\sgm(\al).\eqno(3.2.10)$$
Thus
$$\la\sgm(\be),\sgm(\al)\ra=\la\be,\al\ra.\qquad\Box\eqno(3.2.11)$$
\vspace{0.1cm}

The hyperplanes $\{P_\al\mid\al\in\Phi\}$ partition ${\msr E}$ into
finite many regions; the connected components of ${\msr
E}\setminus\bigcup_{\al\in \Phi}P_\al$ are called the (open) {\it
Weyl chambers} of ${\msr E}$.\index{Weyl chamber} Each regular
vector $\gm$ belongs a unique Weyl chamber, denoted by ${\msr
C}(\gm)$. If $\Pi=\Pi(\gm)$ is a base, we call ${\msr C}(\gm)={\msr
C}(\Pi)$ the {\it fundamental Weyl chamber relative to} $\Pi$. In
fact, we have
\begin{eqnarray*}\hspace{2cm}{\msr C}(\gm)&=&\{\be\in{\msr
E}\mid(\be,\al)>0\;\for\;\al\in\Phi^+(\gm)\}\\&=&\{\be\in{\msr
E}\mid(\be,\al)>0\;\for\;\al\in\Pi^+(\gm)\}.\hspace{4.4cm}(3.2.12)\end{eqnarray*}

\psp

{\bf Lemma 3.2.3}. {\it For $\sgm\in {\msr W}$, $\sgm({\msr
C}(\gm))={\msr C}(\sgm(\gm))$ and $\sgm(\Pi(\gm))=\Pi(\sgm(\gm))$.}

{\it Proof}.  Observe that
$$(\sgm(\gm),\sgm(\al))=(\gm,\al)>0\qquad \for\;\;\al\in
\Phi^+(\gm).\eqno(3.2.13)$$  Thus $\sgm(\gm)$ is regular by (3.1.46)
and (3.1.47), and $\sgm(\Phi^+(\gm))=\Phi^+(\sgm(\gm))$. Hence
$\sgm({\msr C}(\gm))={\msr C}(\sgm(\gm))$ by (3.2.12). Since $\sgm$
is linear, $\{\sgm(\al)\mid\al\in\Pi(\gm)\}$ are indecomposable
elements in $\Phi^+(\sgm(\gm))$. Thus
$\sgm(\Pi(\gm))=\Pi(\sgm(\gm)).\qquad\Box$\psp

 Now we fix a base $\Pi$ of the root system $\Phi$. Then
 $\Phi^+=\{\be\in\Phi\mid \be=\sum_{\al\in\Pi}k_\al\al,\;k_\al\geq
 0\}$.\psp

{\bf Lemma 3.2.4}. {\it For $\al\in\Pi$,
$\sgm_\al(\Phi^+\setminus\{\al\})=\Phi^+\setminus\{\al\}$}.

{\it Proof}. Let $\gm\in\Phi^+\setminus\{\al\}$. Then
$\gm=\sum_{\be\in\Pi}k_\be\be$ with $k_\be\in\mbb{N}$ and
$k_{\be'}>0$ for some $\al\neq\be'\in\Pi$. Note
$$\sgm_\al(\gm)=\sum_{\al\neq\be}k_\be\be+(k_\al-\la\gm,\al\ra)\al,\eqno(3.2.14)$$
whose coefficient of $\be'$ is still $k_{\be'}>0$. Thus
$\sgm_\al(\gm)\in\Phi^+\setminus\{\al\}.\qquad\Box$ \psp

{\bf Corollary 3.2.5}. {\it Set
$$\rho=\frac{1}{2}\sum_{\al\in\Phi^+}\al.\eqno(3.2.15)$$
Then
$$\sgm_\al(\rho)=\rho-\al\qquad\for\;\;\al\in\Pi,\eqno(3.2.16)$$
or equivalently,
$$\la\rho,\al\ra=1\qquad \;\;\al\in\Pi.\eqno(3.2.17)$$}
\vspace{0.1cm}

We define partial ordering $\prec$ on ${\msr E}$ by
$$\al\prec\be\;\mbox{if}\;\be-\al=\sum_{\al\in\Pi}a_\al\al\;\mbox{such
that}\;0\leq a_\al\in\mbb{R}, \eqno(3.2.18)$$ where $\al,\be\in{\msr
E}$.\psp

{\bf Proposition 3.2.6}. {\it Let $\al,\be\in \Phi$ be two positive
distinct roots such that $\al\prec\be$. Then there exist
$\gm_1,...,\gm_s\in\Pi$ such that $\be-\al=\sum_{i=1}^s\gm_i$ and
$$\al+\sum_{i=1}^r\gm_i\in\Phi\qquad
\for\;\;r=1,...,s.\eqno(3.2.19)$$}

{\it Proof}. Write $\be-\al=\sum_{\gm\in\Pi'}k_\gm\gm$ with
$k_\gm>0$ and $\Pi'\subset\Pi$. We prove the proposition by
induction on $\sum_{\gm\in\Pi'}k_\gm$. If $(\al,\gm_0)<0$ or
$(\be,\gm_0)>0$ for some $\gm_0\in\Pi'$, then $\al+\gm_0\in\Phi$ or
$\be-\gm_0\in\Phi$ by Lemma 3.1.2. The conclusion holds for the pair
$(\al+\gm_0,\be)$ or $(\al,\be-\gm_0)$ by induction, which implies
that the conclusion holds for the pair $(\al,\be)$. Suppose
$$(\al,\gm)\geq 0,\;\;(\be,\gm)\leq 0\qquad\for\;\;\gm\in\Pi'.\eqno(3.2.20)$$
Then
$$0\leq (\be-\al,\be-\al)=(\be-\al,\sum_{\gm\in\Pi'}k_\gm\gm)=
\sum_{\gm\in\Pi'}k_\gm[(\be,\gm)-(\al,\gm)]\leq 0.\eqno(3.2.21)$$
Thus $(\be-\al,\be-\al)=0$; that is, $\be-\al=0$, which contradicts
$\al\neq \be.\qquad\Box$\psp

The above proposition had been used in studying the structure of
finite-dimensional simple Lie superalgebras. It is also useful in
the representations of simple Lie algebras.\psp

{\bf Lemma 3.2.7}. {\it Let $\al_1,\al_2,...,\al_k\in\Pi$ (not
necessarily distinct). Write $\sgm_i=\sgm_{\al_i}$. If
$\sgm_1\cdots\sgm_{k-1}(\al_k)\prec 0$, then for some
$s\in\ol{1,k-1}$,
$\sgm_1\cdots\sgm_k=\sgm_1\cdots\sgm_{s-1}\sgm_{s+1}\cdots\sgm_{k-1}$.}

{\it Proof}. Set
$$\be_{k-1}=\al_k,\;\;\be_i=\sgm_{i+1}\cdots\sgm_{k-1}(\al_k)\qquad\for\;\;i\in\{0,1,...,k-2\}.
\eqno(3.2.22)$$ Then there exists $s\in\ol{1,k-1}$ such that
$\be_s\succ 0$ and $\be_{s-1}=\sgm_{\al_s}(\be_s)\prec 0$. By the
above lemma, $\be_s=\al_s$. Thus
$$\sgm_{s+1}\cdots\sgm_{k-1}\sgm_k(\sgm_{s+1}\cdots\sgm_{k-1})^{-1}=\sgm_s
\eqno(3.2.23)$$ by Lemma 3.2.2. Hence
$$\sgm_{s+1}\cdots\sgm_{k-1}\sgm_k=\sgm_s\sgm_{s+1}\cdots\sgm_{k-1},\eqno(3.2.24)$$
$$\sgm_s\sgm_{s+1}\cdots\sgm_{k-1}\sgm_k=\sgm_{s+1}\cdots\sgm_{k-1},\eqno(3.2.25)$$or equivalently,
$\sgm_1\cdots\sgm_k=\sgm_1\cdots\sgm_{s-1}\sgm_{s+1}\cdots\sgm_{k-1}.\qquad\Box$\psp

{\bf Corollary 3.2.8}. {\it If $\sgm=\sgm_1\sgm_2\cdots\sgm_k$ is a
product  of minimal number of the reflections of simple roots, then
$\sgm(\al_k)\prec 0$}.

{\it Proof}. Otherwise,
$$\sgm(\al_k)=\sgm_1\sgm_2\cdots\sgm_{k-1}(-\al_k)=-\sgm_1\sgm_2\cdots\sgm_{k-1}(\al_k)\succ
0.\eqno(3.2.26)$$ Thus $\sgm_1\sgm_2\cdots\sgm_{k-1}(\al_k)\prec 0$.
By the above lemma, the product $\sgm=\sgm_1\sgm_2\cdots\sgm_k$ is
not minimal.$\qquad\Box$ \psp

 {\bf Theorem 3.2.9}. {\it Let $\Pi$ be a base of a root
system $\Phi$.

(a) If $\gm\in{\msr E}$ is regular, there exists $\sgm\in{\msr W}$
such that $(\sgm(\gm),\al)>0$ for $\al\in\Pi$ (so ${\msr W}$ acts
transitively on Weyl chambers).

(b) ${\msr W}$ acts transitively on bases.

(c) If $\al\in\Phi$, there exists $\sgm\in{\msr W}$ such that
$\sgm(\al)\in\Pi$.

(d) ${\msr W}$ is generated by $\{\sgm_\al\mid\al\in\Pi\}$.

(e) If $\sgm(\Pi)=\Pi,\;\sgm\in{\msr W},$ then $\sgm=1$}.

{\it Proof}. Let ${\msr W}'$ be the subgroup of ${\msr W}$ generated
by $\{\sgm_\al\mid\al\in\Pi\}$. Since ${\msr W}$ is finite, ${\msr
W}'$ is finite.

(a) Take $\sgm\in{\msr W}'$ such that $(\sgm(\gm),\Pi)$ is maximal
in $\{(\tau(\gm),\Pi)\mid\tau\in{\msr W}'\}$. For any $\al\in\Pi$,
we have
$$(\sgm(\gm),\Pi)\geq
(\sgm_\al\sgm(\gm),\Pi)=(\sgm(\gm),\sgm_\al(\Pi))=(\sgm(\gm),\Pi-\al)=(\sgm(\gm),\Pi)
-(\sgm(\gm),\al),\eqno(3.2.27)$$ which implies $(\sgm(\gm),\al)\geq
0$. Since $\sgm(\gm)$ is regular, we have
$$(\sgm(\gm),\al)>0\qquad\for\;\;\al\in\Pi.\eqno(3.2.28)$$
Thus we have
$$\sgm(\gm)\in{\msr C}(\Pi)\eqno(3.2.29)$$
by (3.2.12). Hence $\sgm({\msr C}(\gm))={\msr C}(\sgm(\gm))={\msr
C}(\Pi)$; that is,  ${\msr W}$ acts transitively on Weyl chambers.

(b) By (a) and (3.2.12), $\sgm(\Phi^+(\gm))=\Phi^+$. Since $\sgm$ is
linear, it maps indecomposable elements to indecomposable elements;
that is, $\sgm(\Pi(\gm))=\Pi(\sgm(\gm))=\Pi$.

(c) Since $\sgm_\al(\al)=-\al$, we may assume $\al\succ 0$. Suppose
${\msr W}'(\al)\bigcap \Pi=\emptyset$. By Lemma 3.2.4, ${\msr
W}'(\al)$ is a finite set of positive roots. Choose a minimal
element $\gm\in {\msr W}'(\al)$ with respect the partial ordering
$\prec$ in (3.2.18). Then
$$\sgm_\be(\gm)=\gm-\la\gm,\be\ra\be\not\prec\gm\qquad\for\;\;\be\in\Pi.
\eqno(3.2.30)$$ Thus
$$\la\gm,\be\ra\be\prec 0\sim (\gm,\be)\leq 0 \qquad\for\;\;\be\in\Pi.
\eqno(3.2.31)$$ According to the proof of the arguments in
(3.1.56)-(3.1.62), $\{\gm,\Pi\}$ is an independent set (the
positivity of $\gm$ is needed (cf. (3.1.61) and (3.1.62))), which is
absurd.

(d) For any $\al\in\Phi$, we have $\sgm_\al=\sgm_{-\al}$. Thus we
assume that $\al$ is a positive root. By (c), there exists
$\sgm\in{\msr W}'$ such that $\sgm(\al)=\be\in\Pi$. Thus
$\al=\sgm^{-1}(\be)$ and
$$\sgm_\al=\sgm_{\sgm^{-1}(\be)}=\sgm^{-1}\sgm_\be\sgm\in{\msr
W}'\eqno(3.2.32)$$ by Lemma 3.2.2. Since ${\msr W}$ is generated by
$\{\sgm_\al\mid\al\in\Phi\}$, we have ${\msr W}={\msr W}'$.

(e) Suppose $\sgm\neq 1$. Then we can write
$\sgm=\sgm_1\cdots\sgm_k$ as a product of minimal number of the
reflections of simple roots and $\sgm_k=\sgm_{\al_k}$. By Corollary
3.2.8, then $\sgm(\al_k)\prec 0$, which contradicts
$\sgm(\Pi)=\Pi.\qquad\Box$\psp

For $\sgm\in{\msr W}$, if $\sgm=\sgm_1\sgm_2\cdots\sgm_k$ is a
product of  minimal number of the reflections of simple roots, we
called the product $\sgm_1\sgm_2\cdots\sgm_k$ {\it reduced} and call
$k$ the {\it length of $\sgm$}, denoted as $\ell(\sgm)$. Moreover,
we define
$$n(\sgm)=|\{\al\in\Phi^+\mid \sgm(\al)\prec 0\}|.\eqno(3.2.33)$$
\vspace{0.1cm}

{\bf Lemma 3.2.10}. $\ell(\sgm)=n(\sgm)$.

{\it Proof}. We prove it by induction on $\ell(\sgm)$. The case
$\ell(\sgm)=1$ is Lemma 3.2.4. Suppose the lemma holds for
$\ell(\sgm)<k$ for some positive integer $k$. Consider
$\ell(\sgm)=k$. Write $\sgm=\sgm_1\sgm_2\cdots\sgm_k$ with
$\sgm_i=\sgm_{\al_i}$ and $\al_i\in\Pi$.  Let $\sgm'=\sgm_1\cdots
\sgm_{k-1}$. Then $\ell(\sgm')=k-1$. Set
$$S=\{\al\in\Phi^+\mid \sgm'(\al)\prec 0\}.\eqno(3.2.34)$$
By assumption, $|S|=k-1$. According to Lemma 3.2.7,
$\sgm'(\al_k)\succ 0$. Thus
$$S=\{\be\in \Phi^+\setminus\{\al_k\}\mid\sgm'(\be)\prec0\}.\eqno(3.2.35)$$
Since $\sgm_k(\Phi^+\setminus\{\al_k\})=\Phi^+\setminus\{\al_k\}$ by
Lemma 3.2.4, we have
$$\sgm_k(S)=\{\sgm_k(\be)\in \Phi^+\setminus\{\al_k\}\mid\sgm'(\be)\prec 0\}=\{\be'\in \Phi^+\setminus\{\al_k\}\mid
\sgm(\be')\prec 0\}.\eqno(3.2.36)$$ Therefore,
$$\sgm_k(S)\bigcup\{\al\}=\{\be\in \Phi^+\mid\sgm(\be)\prec 0\},\eqno(3.2.37)$$
which implies $n(\sgm)=k.\qquad\Box$\psp

{\bf Remark 3.2.11}. According to (3.2.37), if
$\sgm=\sgm_1\sgm_2\cdots\sgm_k$ is a reduced product of the
reflections of simple roots with $\sgm_i=\sgm_{\al_i}$, then
$$\{\be\in\Phi^+\mid\sgm(\be)\prec
0\}=\{\al_k,\sgm_k\cdots\sgm_{i+1}(\al_i)\mid
i\in\ol{1,k-1}\}.\eqno(3.2.38)$$ \vspace{0.1cm}

{\bf Lemma 3.2.12}. {\it Let $\lmd,\mu\in\ol{{\msr C}(\Pi)}$. If
$\sgm(\lmd)=\mu$, then $\lmd=\mu$.}

{\it Proof}. If $\ell(\sgm)\neq 0$, let
$\sgm=\sgm_1\sgm_2\cdots\sgm_k$ be a reduced product of the
reflections of simple roots with $\sgm_i=\sgm_{\al_i}$. Then
$\sgm(\al_k)\prec 0$. Thus
$$0\geq (\mu,\sgm(\al_k))=(\sgm^{-1}(\mu),\al_k)=(\lmd,\al_k)\geq
0.\eqno(3.2.39)$$ Thus $(\lmd,\al_k)=0$, or equivalently,
$\sgm_k(\lmd)=\lmd$. Now $\sgm_1\cdots\sgm_{k-1}(\lmd)=\mu$. By
induction on $\ell(\sgm)$, we have $\lmd=\mu.\qquad\Box$\psp

{\bf Lemma 3.2.13}. {\it For any $0\neq\gm\in{\msr E}$, the orbit
$\{\sgm(\gm)\mid\sgm\in{\msr W}\}$ has a unique maximal element
$\be$ with respect to the partial ordering $\prec$ defined in
(3.2.18). Moreover, $\be\in\ol{{\msr C}(\Pi)}$. In particular,}
$$\sgm(\gm)\prec \gm\qquad\mbox{\it for any}\;\;\gm\in\ol{{\msr
C}(\Pi)}.\eqno(3.2.40)$$

{\it Proof}. Suppose that $\be=\sgm_1(\gm)$ is a maximal element in
$\{\sgm(\gm)\mid\sgm\in{\msr W}\}$. For any $\al\in\Pi$,
$$\sgm_\al(\be)=\be-\la\be,\al\ra\al\not\succ\be\lra
\la\be,\al\ra\geq 0\sim (\be,\al)\geq 0.\eqno(3.2.41)$$ Thus
$\be\in\ol{{\msr C}(\Pi)}$.

 Suppose that $\be'=\sgm_2(\gm)$ is also a  maximal element in
$\{\sgm(\gm)\mid\sgm\in{\msr W}\}$. Then we also have
$\be'\in\ol{{\msr C}(\Pi)}$. Since $\sgm_1\sgm_2^{-1}(\be')=\be$, we
have $\be=\be'$ by the above lemma. If $\gm\in\ol{{\msr C}(\Pi)}$,
we have $\be=\gm$ again by the above lemma.$\qquad\Box$\psp

The above two lemmas shows the ${\msr W}$-orbit space
$${\msr E}/{\msr W}\cong \ol{{\msr C}(\Pi)}.\eqno(3.2.42)$$
In other words, there exists a one-to-one correspondence between the
functions on $ \ol{{\msr C}(\Pi)}$ and ${\msr W}$-invariant
functions on ${\msr E}$.

 A root system $\Phi$ is called {\it reducible} if it is a
union of two nonempty subsets $\Phi_1$ and $\Phi_2$ such that
$$(\al,\be)=0\qquad\for\;\;\al\in\Phi_1,\;\be\in\Phi_2.\eqno(3.2.43)$$
\vspace{0.1cm}

{\bf Lemma 3.2.14}. {\it Let $\Pi$ be a base of a roots system
$\Phi$. Then $\Phi$ is reducible if and only if $\Pi$ is a union of
two nonempty subsets $\Pi_1$ and $\Pi_2$ such that}
$$(\al,\be)=0\qquad\mbox{\it for}\;\;\al\in\Pi_1,\;\be\in\Pi_2.\eqno(3.2.44)$$

{\it Proof}. Suppose that  $\Phi$ is a union of two nonempty subsets
$\Phi_1$ and $\Phi_2$ such that (3.2.43) holds. If
$\Pi\subset\Phi_1$, then the elements in $\Phi_2$ are orthogonal to
a basis of ${\msr E}$. So $\Phi_2=\{0\}$, which is absurd. Thus
$\Pi\not\subset \Phi_i$. Set $\Pi_i=\Pi_i\bigcap\Phi_i$. Now $\Pi$
is a union of two nonempty subsets $\Pi_1$ and $\Pi_2$ such that
(3.2.44) holds.

Assume that $\Pi$ is a union of two nonempty subsets $\Pi_1$ and
$\Pi_2$ such that (3.2.44) holds. Set
$$\Phi_1=\{\sgm(\al)\mid\al\in\Pi_1,\;\sgm\in{\msr W}\},\;\;\Phi_2=\{\sgm(\al)\mid\al\in\Pi_2,
\;\sgm\in{\msr W}\}.\eqno(3.2.45)$$ For $\al\in\Pi_1$ and
$\be\in\Pi_2$, we have $\sgm_\al(\be)=\be-\la\be,\al\ra\al=\be$ and
$$\sgm_\al\sgm_\be\sgm_\al^{-1}=\sgm_{\sgm_\al(\be)}=\sgm_\be\sim
\sgm_\al\sgm_\be=\sgm_\be\sgm_\al.\eqno(3.2.46)$$ Let
$${\msr W}_i=\mbox{the subgroup of}\;{\msr W}\;\mbox{generated
by}\;\{\sgm_\al\mid\al\in\Pi_i\}.\eqno(3.2.47)$$ Since ${\msr W}$ is
generated by $\{\sgm_\al\mid\al\in\Pi\}$ by Theorem 3.2.9 (d), we
have
$${\msr W}={\msr W}_1{\msr W}_2\eqno(3.2.48)$$
by (3.2.46). Thus
$$\Phi_i=\{\sgm(\al)\mid\al\in\Pi_i,\;\sgm\in{\msr
W}_i\}.\eqno(3.2.49)$$ For $\al_i\in\Pi_i$ and $\sgm_i\in{\msr
W}_i$, we have
$$\sgm_2(\al_1)=\al_1,\qquad \sgm_1(\al_2)=\al_2.\eqno(3.2.50)$$
Hence
$$(\sgm_1(\al_1),\sgm_2(\al_2))=(\sgm_2^{-1}(\al_1),\sgm_1^{-1}(\al_2))=(\al_1,\al_2)=0.
\eqno(3.2.51)$$ So (3.2.44) holds. By Theorem 3.2.9 (c),
$\Phi=\Phi_1\bigcup\Phi_2.\qquad\Box$\psp

{\bf Lemma 3.2.15}. {\it Let $\Phi$ be an irreducible root system.
Relative to the partial ordering $\prec$, there exists a unique
maximal root $\be$. Moreover, $\be=\sum_{\al\in\Pi}k_\al\al$ with
all $k_\al>0$.}

{\it Proof}. Let $\be=\sum_{\al\in\Pi}k_\al\al$ be a maximal root.
Then $\be\succ 0$. Set
$$\Pi_1=\{\al\in\Pi\mid k_\al>0\},\qquad\Pi_2=\{\al\in\Pi\mid
k_\al=0\}.\eqno(3.2.52)$$ Then $\Pi=\Pi_1\bigcup\Pi_2$ and
$\Pi_1\neq\emptyset$. If $\Pi_2\neq\emptyset$, there exist
$\al_1\in\Pi_1$ and  $\al_2\in\Pi_2$ such that $(\al_1,\al_2)<0$ by
(3.1.51), the above lemma and the irreducibility of $\Phi$. Again
(3.1.51) implies $(\be,\al_2)<0$. According to Lemma 3.1.2,
$\be\prec\be+\al_2\in\Phi$, which contradicts the maximality of
$\be$. Thus $\Pi_2=\emptyset$.

Suppose that $\be'$ is another maximal root. By Lemma 3.1.2,
$(\be',\al)\geq 0$ for $\al\in\Pi$. Since $\be'$ can not be
orthogonal to all the elements of $\Pi$, there exists $\al_1\in\Pi$
such that $(\be',\al_1)>0$. This implies
$$(\be,\be')=\sum_{\al\in\Pi}k_\al(\al,\be')>0.\eqno(3.2.53)$$
Lemma 3.1.2 tells us that $\be-\be'$ is a root. Without loss of
generality, we can assume $\be-\be'\succ 0$. Then
$$\be=(\be-\be')+\be'\succ \be',\eqno(3.2.54)$$
which contradicts the maximality of $\be'.\qquad\Box$\psp

{\bf Lemma 3.2.16}. {\it Let $\Phi$ be an irreducible root system.
Then ${\msr W}$ acts irreducibly on ${\msr E}$. In particular, the
${\msr W}$-orbit of a root $\al$ spans ${\msr E}$.}

{\it Proof}. Let $W$ be a nonzero ${\msr W}$-submodule of ${\msr E}$
and let $W^\bot$ be its orthogonal complement in ${\msr E}$. Suppose
$W^\bot\neq\{0\}$. For any $\al\in\Phi$, if $\al\not\in W^\bot$,
there exists $\gm\in W$ such that $(\gm,\al)\neq 0$. Thus
$$\sgm_\al(\gm)=\gm-\la\gm,\al\ra\al\in W\lra \la\gm,\al\ra\al\in W\lra
\al\in W.\eqno(3.2.55)$$ Symmetrically, if a root $\be\not\in W$,
then $\be\in W^\bot$. Since $\Phi$ spans ${\msr E}$, we have
$$\Phi_1=\Phi\bigcap W\neq\emptyset,\qquad \Phi_2=\Phi\bigcap
W^\bot\neq\emptyset\eqno(3.2.56)$$ and $\Phi=\Phi_1\bigcup\Phi_2$,
which leads a contradiction to the irreducibility of
$\Phi.\qquad\Box$\psp

{\bf Lemma 3.2.17}. {\it Let $\Phi$ be an irreducible root system.
There are only at most two root lengths  occurring in $\Phi$, and
all roots with the same length are conjugate under ${\msr W}$ (in
the same ${\msr W}$-orbit).}

{\it Proof}. For any $\al,\be\in\Phi$,  $(\al,\sgm(\be))\neq 0$ for
some $\sgm\in{\msr W}$ by the above lemma. Replacing $\be$ by
$\sgm(\be)$, we may assume $(\al,\be)\neq 0$ and $(\al,\al)\leq
(\be,\be)$. By Lemma 3.1.1,
$$\frac{(\be,\be)}{(\al,\al)}\in\{1,2,3\}.\eqno(3.2.57)$$
If we have $\al_1,\al_2,\al_3\in\Phi$ with
$(\al_1,\al_1)<(\al_2,\al_2)<(\al_3,\al_3)$. By (3.2.57), the pairs
$(\al_1,\al_2)$ and $(\al_1,\al_3)$ give
$$(\al_2,\al_2)=2(\al_1,\al_1),\qquad(\al_3,\al_3)=3(\al_1,\al_1).\eqno(3.2.58)$$
But
$$\frac{(\al_3,\al_3)}{(\al_2,\al_2)}=\frac{3}{2},\eqno(3.2.59)$$
which contradicts (3.2.57). Thus there are only at most two root
lengths occurring in $\Phi$.

Suppose we have two  elements $\al,\be\in\Phi$ with the same length.
Replacing $\be$ by $\sgm(\be)$ for some $\sgm\in{\msr W}$ if
necessary, we may assume $(\al,\be)\neq 0$ and $\al\neq\be$. By
Lemma 3.1.1, $\la\al,\be\ra=\la\be,\al\ra=\pm 1$. Replacing $\be$ by
$\sgm_\be(\be)=-\be$ if necessarily, we may assume
$\la\al,\be\ra=\la\be,\al\ra=1$. Now
$$\sgm_\al\sgm_\be\sgm_\al(\be)=\sgm_\al\sgm_\be(\be-\al)=\sgm_\al(-\al)=\al.
\qquad\Box\eqno(3.2.60)$$ \vspace{0.1cm}

 {\bf Lemma 3.2.18}. {\it Let $\Phi$ be an irreducible root system
 with two distinct lengths. Then the maximal root $\be$ is a long
 root.}

 {\it Proof}. By Lemma 3.1.2, $\be\in\ol{{\msr C}(\Pi)}$.
 Let $\al\in\Phi$ be any root such that $||\al||\neq||\be||$. According to Theorem 3.2.9 (a),
 there exists $\sgm\in{\msr W}$ such that $\sgm(\al)\in \ol{{\msr
 C}(\Pi)}$. Replacing $\al$ by $\sgm(\al)$, we can assume $\al\in\ol{{\msr C}(\Pi)}$.
Since $\be-\al\succ 0$, we have
$$(\be,\be-\al)\geq 0\lra (\be,\be)\geq (\be,\al)\eqno(3.2.61)$$
and
$$(\be-\al,\al)\geq 0\lra (\be,\al)\geq (\al,\al).\eqno(3.2.62)$$
Thus $(\be,\be)\geq (\al,\al).\qquad\Box$\psp

The following fact is useful in representation theory and Lie
superalgebras.\psp

 {\bf Lemma 3.2.19}. {\it Let $\al$ and $\be$ be two positive roots
 of a root system $\Phi$ such that $\al\prec\be$. Then there exist roots $\al=\gm_0,\gm_1,...,
 \gm_k=\be$ such that $\gm_i-\gm_{i-1}\in\Pi$ for $i=1,...,k$.}

{\it Proof}. We write
$$\be-\al=\sum_{r\in S}m_r\al_r,\qquad m_r>0,\eqno(3.2.63)$$
where $S$ is a subset of $\{1,...,n\}$. We prove by induction on
$$d(\be,\al)=\sum_{r\in S}m_r.\eqno(3.2.64)$$
If $d(\be,\al)=1$, then the lemma trivially holds. Suppose that the
lemma holds for $d(\be,\al)<m$. Assume $d(\be,\al)=m$. If
$(\al_s,\al)<0$ for some $r\in S$, then $\al+\al_r\in \Phi$ by Lemma
3.1.2. Note $\al+\al_r\prec\be$ and $d(\be,\al+\al_r)=m-1$. So there
exists $\gm_1=\al+\al_r,\gm_2,...,\gm_k=\be$ such that
$\gm_{i+1}-\gm_i\in\Pi$ for $i=1,...,k-1$. Then
$\al=\gm_0,\gm_1,..., \gm_k=\be$ satisfy the Lemma.

When $(\be,\al_r)>0$ for some $r\in S$,  then $\be-\al_r\in \Phi$ by
Lemma 3.1.2.  Note $\al\prec\be-\al_r$ and $d(\be-\al_r,\al)=m-1$.
So there exists $\gm_0=\al,\gm_1,...,\gm_{k-1}=\be-\al_r$ such that
$\gm_i-\gm_{i-1}\in\Pi$ for $i=1,...,k-1$. Then
$\al=\gm_0,\gm_1,..., \gm_k=\be$ satisfy the Lemma.

Finally we consider the case
$$(\al,\al_r)\geq 0,\qquad(\be,\al_r)\leq 0\qquad\for\;\;r\in
S.\eqno(3.2.65)$$ But
$$0<(\be-\al,\be-\al)=\sum_{i\in S}m_i(\be-\al,\al_i)=\sum_{i\in
S}m_i(\be,\al_i)-\sum_{i\in S}m_i(\al,\al_i)\leq 0,\eqno(3.2.66)$$
which leads a contradiction.$\qquad\Box$\psp

In the following examples, we let ${\msr E}$ be a vector space over
$\mbb{R}$ with a basis $\{\ves_1,\ves_2,...,\ves_n\}$ and an inner
product determined by
$$(\ves_i,\ves_j)=\dlt_{i,j}\qquad\for\;\;i,j\in\ol{1,n}.\eqno(3.2.67)$$
\vspace{0.1cm}

{\bf Example 3.2.1}. Recall the root system $\Phi_{A_{n-1}}$  in
(3.1.12) of $sl(n,\mbb{C})$. It is a root system in the subspace
$${\msr E}'=\{\sum_{i=1}^na_i\ves_i\mid
a_i\in\mbb{R},\;\sum_{r=1}^na_r=0\}.\eqno(3.2.68)$$ To understand
the Weyl group of $\Phi_{A_{n-1}}$, we consider its action on ${\msr
E}$. For $\ves_i-\ves_j\in\Phi_{A_{n-1}}$ with $i\neq j$,
$$\sgm_{\ves_i-\ves_j}(\ves_r)=\left\{\begin{array}{ll}\ves_r&\mbox{if}\;r\neq
i,j\\ \ves_j&\mbox{if}\;r=i\\
\ves_i&\mbox{if}\;\;r=j\end{array}\right.\qquad\for\;\;r\in\ol{1,n}.\eqno(3.2.69)$$
Let $S_n$ be the permutation group  on $\{1,2,...,n\}$. Define a
representation $\nu:S_n\rta GL({\msr E})$ by
$$\nu(\iota)(\ves_r)=\ves_{\iota(r)}\qquad\for\;\;\iota\in
S_n,\;r\in\ol{1,n}.\eqno(3.2.70)$$ Then
$$S_n\stl{\nu}{\cong}{\msr W}_{A_{n-1}},\eqno(3.2.71)$$
the Weyl group of $\Phi_{A_{n-1}}$.

 Since $\sgm_{\ves_i}=\sgm_{2\ves_i}$, the Weyl
groups of $\Phi_{B_n}$ in (3.1.28) and $\Phi_{C_n}$ in (3.1.27) are
the same. Note
$$\sgm_{\ves_i}(\ves_r)=(-1)^{\dlt_{i,r}}\ves_r,\;\;\sgm_{\ves_i+\ves_j}(\ves_r)=
\left\{\begin{array}{ll}\ves_r&\mbox{if}\;r\neq
i,j\\ -\ves_j&\mbox{if}\;r=i\\
-\ves_i&\mbox{if}\;\;r=j\end{array}\right.\qquad\for\;\;r\in\ol{1,n}\eqno(3.2.72)$$
Let $\mbb{Z}_2=\mbb{Z}/2\mbb{Z}$ be the field of two elements and
let $\Upsilon$ be a vector space over $\mbb{Z}_2$ with a basis
$\{\ves_1,\ves_2,...,\ves_n\}$. Define
$$\iota(\vs)=\sum_{i=1}^na_i\ves_{\iota(i)}\qquad\for\;\;\iota\in
S_n,\;\vs=\sum_{i=1}^na_i\ves_i\in\Upsilon.\eqno(3.2.73)$$ Viewing
$\Upsilon$ as an additive group, we form a semi-product of $S_n$
with $\Upsilon$:
$$G=S_n\ltimes\Upsilon ,\eqno(3.2.74)$$
whose group multiplication is defined by:
$$(\iota_1,\vs_1)\cdot(\iota_2,\vs_2)=(\iota_1\iota_2,\vs_1+\iota_1(\vs_2)).\eqno(3.2.75)$$
Moreover, we define $\nu:G\rta GL({\msr E})$ by
$$\nu[(\iota,\sum_{i=1}^na_i\ves_i,)](\ves_r)=(-1)^{a_r}\ves_{\iota(r)}\qquad\for\;\;r\in\ol{1,n}.
\eqno(3.2.76)$$ By (3.2.69) and (3.2.72), we have
$$G\stl{\nu}{\cong}{\msr W}_{B_n}={\msr W}_{C_n},\eqno(3.2.77)$$
the Weyl groups of $\Phi_{B_n}$ and $\Phi_{C_n}$, respectively.

Set
$$\Upsilon'=\{\sum_{i=1}^na_i\ves_i\in\Upsilon\mid
\sum_{r=1}^na_i=0\},\;\;G'= S_n\ltimes \Upsilon'.\eqno(3.2.78)$$
Then $G'$ is a subgroup of $G$. Moreover,
$$G'\stl{\nu}{\cong}{\msr W}_{D_n},\eqno(3.2.79)$$
the Weyl groups of $\Phi_{D_n}$ in (3.1.17).\psp

The followings are exercises:

(1) If $\sgm\in{\msr W}$ can be written as  a product of $k$
reflections of simple roots ({\it simple reflections}), prove
$k\equiv \ell(\sgm)\;(\mbox{mod}\;2)$.

(2) Prove that there exists a unique element $\sgm$ such that
$\sgm(\Phi^+)=-\Phi^+=\Phi^-$. Moreover, any reduced expression for
$\sgm$ must involve all $\{\sgm_\al\mid\al\in\Pi\}$. What is
$\ell(\sgm)$?

(3) Prove that the only reflections in ${\msr W}$ are
$\{\sgm_\al\mid\al\in\Phi\}$.

\section{Classification}

In this section, we  classify all the irreducible root systems.

Let $\Phi$ be an irreducible root system with a base
$$\Pi=\{\al_1,\al_2,...,\al_n\}.\eqno(3.3.1)$$
The matrix
$$CM(\Pi)=(\la\al_i,\al_j\ra)_{n\times n}\eqno(3.3.2)$$
is called the {\it Cartan matrix}\index{Cartan matrix} of the root
system $\Phi$.
 Since $\Pi$ is a basis of the
corresponding Euclidean space ${\msr E}$, each simple reflection
$\sgm_{\al_i}$ is determined by the integers
$$\{\la\al_1,\al_i\ra,\la\al_2,\al_i\ra,...,\la\al_n,\al_i\ra\};\eqno(3.3.3)$$
that is, the $i$th column of $CM(\Pi)$. Recall that the Weyl group
${\msr W}$ is generated by $\{\sgm_{\al_1},...,\sgm_{\al_n}\}$
(Theorem 3.2.8 (d)). Thus the root system $\Phi={\msr W}(\Pi)$ is
completely determined by $CM(\Pi)$. Observe that the diagonal
entries of $CM(\Pi)$ are 2. So $\Phi$ is determined by
$$\{\la\al_i,\al_j\ra\mid i,j\in\ol{1,n},\;i\neq
j\}.\eqno(3.3.4)$$

Recall
$$(\al_i,\al_j)\leq 0,\qquad\la\al_i,\al_j\ra
\la\al_j,\al_i\ra=0,1,2,3\qquad\for\;\;i,j\in\ol{1,n},\;i\neq j
\eqno(3.3.5)$$ (cf. Lemma 3.1.1 and (3.1.51)). Set
$$\mfk{e}_i=\frac{\al_i}{||\al_i||}\qquad\for\;\;i\in\ol{1,n},\eqno(3.3.5)$$
where $||\al_i||=\sqrt{(\al_i,\al_i)}$. Now
$\{\mfk{e}_1,...,\mfk{e}_n\}$ are unit vectors and
$$\la\mfk{e}_i,\mfk{e}_j\ra=2(\mfk{e}_i,\mfk{e}_j)=2\left(\frac{\al_i}{||\al_i||},
\frac{\al_j}{||\al_j||}\right)=\frac{2(\al_i,\al_j)}{||\al_i||\:||\al_j||}.\eqno(3.3.6)$$
Thus
$$\la\mfk{e}_i,\mfk{e}_j\ra\la\mfk{e}_j,\mfk{e}_i\ra=
\left(\frac{2(\al_i,\al_j)}{||\al_i||\:||\al_j||}\right)^2=\la\al_i,\al_j\ra
\la\al_j,\al_i\ra.\eqno(3.3.7)$$ Therefore, (3.3.5) is equivalent to
$$(\mfk{e}_i,\mfk{e}_j)\leq 0,\qquad\la\mfk{e}_i,\mfk{e}_j\ra
\la\mfk{e}_j,\mfk{e}_i\ra=4(\mfk{e}_i,\mfk{e}_j)^2=0,1,2,3\qquad\for\;\;i,j\in\ol{1,n},\;i\neq
j \eqno(3.3.8)$$

In general, a linearly independent subset of unit vectors
$\{\mfk{e}_1,\mfk{e}_2,...,\mfk{e}_n\}$ in a Euclidean space is
called an {\it admissible set} if $(3.3.8)$ holds. \index{admissible
set} Note that any subset of an admissible set is still an
admissible set. If an admissible set ${\msr S}$ can not be written
as a union of two orthogonal proper subsets, then we call ${\msr S}$
irreducible.

Let ${\msr S}=\{\mfk{e}_1,\mfk{e}_2,...,\mfk{e}_n\}$ be an
irreducible admissible set. Define the {\it Coxeter
graph}\index{Coxeter graph} of ${\msr S}$ to be a graph having $n$
vertices corresponding to the vectors in ${\msr S}$, in which the
$i$th vertex is jointed to $j$th vertex $(i\neq j)$ by
$\la\mfk{e}_i,\mfk{e}_j\ra \la\mfk{e}_j,\mfk{e}_i\ra$ edges. We
denote the Coxeter graph of ${\msr S}$ by $CG({\msr S})$. Note that
$$0<(\sum_{i=1}^n\mfk{e}_i,\sum_{i=1}^n\mfk{e}_i)=\sum_{i=1}^n(\mfk{e}_i,\mfk{e}_i)+\sum_{1\leq
i<j\leq n}2(\mfk{e}_i,\mfk{e}_j)=n+\sum_{1\leq i<j\leq
n}2(\mfk{e}_i,\mfk{e}_j).\eqno(3.3.9)$$ Since
$$2(\mfk{e}_i,\mfk{e}_j)=0,-1,-\sqrt{2},-\sqrt{3}\eqno(3.3.10)$$
and the graph is connected by the irreducibility, we have
$$\mbox{the number of connected pairs of vertices in}\;CG({\msr
S})\;\mbox{is}\; n-1=|{\msr S}|-1.\eqno(3.3.11)$$ Suppose there is a
cycle in $CG({\msr S})$. Denote by ${\msr S}'$ the subset of vectors
corresponding to the vertices in the cycle. Then ${\msr S}'$ is also
an irreducible admissible set. But the  number of connected pairs of
vertices in $CG({\msr S}')$ is greater than or equal to $|{\msr
S}'|$, which leads contradiction to (3.3.11) with ${\msr S}$
replaced by ${\msr S}'$. Thus
$$CG({\msr S})\;\mbox{contains no cycles}.\eqno(3.3.12)$$

\setlength{\unitlength}{3pt}

 If $(\mfk{e}_i,\mfk{e}_j)\neq 0$, we say that $\mfk{e}_i$ is {\it
connected} to $\mfk{e}_j$. For $\mfk{e}\in{\msr S}$, let
$\{\vs_1,...,\vs_k\}$ are all the vectors connected to $\mfk{e}$. By
(3.3.12), no pair of $\vs_i$ and $\vs_j$ is connected; that is,
$$(\vs_i,\vs_j)=0\qquad\for\;\;i,j\in\ol{1,k},\;i\neq
j.\eqno(3.3.13)$$ Set
$$\al=\mfk{e}-\sum_{i=1}^k(\mfk{e},\vs_i)\vs_i.\eqno(3.3.14)$$
Then $\al\neq 0$ and
$$(\al,\vs_i)=0\qquad\for\;\;i\in\ol{1,k}.\eqno(3.3.15)$$
Thus
$$0<(\al,\al)=(\al,\mfk{e}-\sum_{i=1}^k(\mfk{e},\vs_i)\vs_i)=(\al,\mfk{e})
=1-\sum_{i=1}^k(\mfk{e},\vs_i)^2;\eqno(3.3.16)$$ that is,
$$\sum_{i=1}^k(\mfk{e},\vs_i)^2<1\sim
\sum_{i=1}^k4(\mfk{e},\vs_i)^2<4.\eqno(3.3.17)$$ By (3.3.8),
$$\mbox{no more than three edges can oringinate at a vertex
in}\;CG({\msr S}),\eqno(3.3.18)$$ which implies that if $CG({\msr
G})$ contains a triple edges, then it must be
\begin{picture}(17,3.4)\put(1,0){\circle{2}}\put(1,1.2){\line(1,0){14}}
\put(2,0.285){\line(1,0){12}}\put(1,-0.8){\line(1,0){14}}\put(15,0){\circle{2}}\end{picture}.
Moreover, $CG({\msr S})$ contains no following subgraphs:
\begin{center}\begin{picture}(100,5)\put(1,0){\circle{2}}\put(1,1.2){\line(1,0){14}}
\put(1,-0.8){\line(1,0){14}}\put(15,0){\circle{2}}\put(15,1.2){\line(1,0){14}}
\put(15,-0.8){\line(1,0){14}}\put(29,0){\circle{2}}
\put(40,0){\circle{2}}\put(40.5,1.2)
{\line(1,0){14}}\put(40.5,-0.8){\line(1,0){14}}\put(54.5,0){\circle{2}}\put(55.1,0.5)
{\line(3,1){12}}\put(67.5,3.7){\circle{2}}\put(55.1,-0.5)
{\line(3,-1){12}}\put(67.5,-3.7){\circle{2}}
\put(80,3.7){\circle{2}}\put(80,-3.7){\circle{2}}
\put(80.8,4.6){\line(3,-1){12}}\put(80.8,-4.6){\line(3,1){12}}\put(94,0){\circle{2}}
\put(94.6,0.5){\line(3,1){12}}\put(107.2,3.7){\circle{2}}\put(94.6,-0.5)
{\line(3,-1){12}}\put(107.2,-3.7){\circle{2}}
\end{picture}\end{center}
\psp

 Suppose that $CG({\msr S})$ contains a simple chain:
\begin{picture}(45,5)\put(1,0){\circle{2}}\put(2,0){\line(1,0){12}}\put(1,-5){$\mfk{e}_i$}
\put(15,0){\circle{2}}\put(15,-5){$\mfk{e}_{i+1}$}\put(20,0){...}\put(28,0){\circle{2}}
\put(28,-5){$\mfk{e}_{j-1}$}
\put(29,0){\line(1,0){12}}\put(42,0){\circle{2}}\put(42,-5){$\mfk{e}_j$}\end{picture}.
Set\psp

$$\mfk{e}=\sum_{r=i}^j\mfk{e}_r.\eqno(3.3.19)$$
Then
$$(\mfk{e},\mfk{e})=j-i+1+\sum_{r=i}^{j-1}2(\mfk{e}_r,\mfk{e}_{r+1})=j-i+1-(j-i)=1,\eqno(3.3.20)$$
$$(\mfk{e}_r,\mfk{e})=(\mfk{e}_r,\mfk{e}_i)\;\mbox{or}\;(\mfk{e}_r,\mfk{e}_j)\qquad\for\;\;r<i\;\mbox{or}\;r>j.
\eqno(3.3.21)$$ Thus $\{\mfk{e}\}\bigcup{\msr
S}\setminus\{\mfk{e}_i,...,\mfk{e}_j\}$ is an irreducible admissible
set. By the arguments in the above paragraph, $CG({\msr S})$
contains no subgraphs of the forms:
\begin{center}\begin{picture}(67,10)\put(1,0){\circle{2}}\put(1,1.2){\line(1,0){14}}
\put(1,-0.8){\line(1,0){14}}\put(15,0){\circle{2}}\put(16,0){\line(1,0){12}}\put(29.5,0){\circle{2}}
\put(33,0){...}\put(38,0){\circle{2}}\put(39,0){\line(1,0){12}}\put(52,0){\circle{2}}\put(52,1.2){\line(1,0){14}}
\put(52,-0.8){\line(1,0){14}}\put(66,0){\circle{2}}\end{picture}

\begin{picture}(67,10)\put(1,0){\circle{2}}\put(1.5,1.2)
{\line(1,0){14}}\put(1.5,-0.8){\line(1,0){14}}\put(15,0){\circle{2}}\put(16,0){\line(1,0){12}}
\put(29.5,0){\circle{2}}
\put(33,0){...}\put(38,0){\circle{2}}\put(39,0){\line(1,0){12}}\put(52,0){\circle{2}}\put(52.6,0.5)
{\line(3,1){12}}\put(65,3.7){\circle{2}}\put(52.6,-0.5)
{\line(3,-1){12}}\put(65,-3.7){\circle{2}}\end{picture}

\begin{picture}(67,12)\put(1,3.7){\circle{2}}\put(1,-3.7){\circle{2}}
\put(1.8,4.6){\line(3,-1){12}}\put(1.8,-4.6){\line(3,1){12}}\put(15,0){\circle{2}}\put(16,0)
{\line(1,0){12}}\put(29.5,0){\circle{2}}
\put(33,0){...}\put(38,0){\circle{2}}\put(39,0){\line(1,0){12}}\put(52,0){\circle{2}}\put(52.6,0.5)
{\line(3,1){12}}\put(65,3.7){\circle{2}}\put(52.6,-0.5)
{\line(3,-1){12}}\put(65,-3.7){\circle{2}}\end{picture}
\end{center}\psp

The remaining possible cases of  $CG({\msr G})$ are:
\begin{center}
\begin{picture}(17,10)\put(1,0){\circle{2}}\put(1,1.2){\line(1,0){14}}
\put(2,0.285){\line(1,0){12}}\put(1,-0.8){\line(1,0){14}}\put(15,0){\circle{2}}
\put(15,-5){$\mfk{e}_2$}\put(1,-5){$\mfk{e}_1$}\end{picture}

\begin{picture}(50,12)\put(1,0){\circle{2}}\put(2,0){\line(1,0){12}}\put(1,-5){$\mfk{e}_1$}
\put(15,0){\circle{2}}\put(15,-5){$\mfk{e}_2$}\put(23,0){...}\put(33,0){\circle{2}}
\put(33,-5){$\mfk{e}_{n-1}$}
\put(34,0){\line(1,0){12}}\put(47,0){\circle{2}}\put(47,-5){$\mfk{e}_n$}
\end{picture}

\begin{picture}(102,12)\put(1,0){\circle{2}}\put(1,-5){$\mfk{e}_1$}\put(2,0){\line(1,0){12}}
\put(15,0){\circle{2}}\put(15,-5){$\mfk{e}_2$}\put(16,0){\line(1,0){12}}\put(30,0){\circle{2}}
\put(30,-5){$\mfk{e}_3$}
\put(37,0){...}\put(45,0){\circle{2}}\put(45,-5){$\mfk{e}_r$}\put(45,1.2){\line(1,0){13.6}}
\put(45,-0.8){\line(1,0){13.6}}\put(58.5,0){\circle{2}}\put(58.5,-5){$\vs_s$}\put(59.5,0)
{\line(1,0){12}}\put(72.5,0){\circle{2}}\put(72.5,-5){$\vs_{s-1}$}
\put(80,0){...}\put(88.5,0){\circle{2}}\put(88.5,-5){$\vs_2$}\put(89.5,0){\line(1,0){12}}
\put(102.5,0){\circle{2}} \put(102.5,-5){$\vs_1$}\end{picture}

\begin{picture}(88,25)\put(1,0){\circle{2}}\put(1,-5){$\mfk{e}_1$}\put(2,0){\line(1,0){12}}
\put(16,0){\circle{2}}\put(16,-5){$\mfk{e}_2$}
\put(24,0){...}\put(32,0){\circle{2}}\put(32,-5){$\mfk{e}_{p-1}$}\put(33,0){\line(1,0){12}}
\put(46,0){\circle{2}}\put(46,3){$\psi$}\put(46.6,0.5)
{\line(3,1){12}}\put(60,4.5){\circle{2}}\put(60,0){$\xi_{r-1}$}\put(64,5.8){.}
\put(66,6.46){.}\put(68,7.14){.}\put(72,8.5){\circle{2}}\put(72,4){$\xi_2$}
\put(72.6,8.8){\line(3,1){13.6}}
\put(87,13.7){\circle{2}}\put(87,8.3){$\xi_1$} \put(46.6,-0.5)
{\line(3,-1){12}}\put(60,-4.5){\circle{2}}\put(60,-9){$\vs_{q-1}$}\put(64,-5.8){.}
\put(66,-6.46){.}\put(68,-7.14){.}\put(72,-8.5){\circle{2}}\put(72,-13.5){$\vs_2$}
\put(72.6,-8.8){\line(3,-1){13.6}}
\put(87,-13.7){\circle{2}}\put(87,-18.7){$\vs_1$}
\end{picture}
\end{center}\vspace{1cm}

In the third case, we let
$$\mfk{e}=\sum_{i=1}^ri\mfk{e}_i,\qquad\vs=\sum_{i=1}^si\vs_s.\eqno(3.3.22)$$
Moreover,
$$2(\mfk{e}_i,\mfk{e}_{i+1})=-1=2(\vs_i,\vs_{i+1}),\qquad
4(\mfk{e}_r,\vs_s)^2=2.\eqno(3.3.23)$$ Thus
$$(\mfk{e},\mfk{e})=\sum_{i=1}^ri^2-\sum_{s=1}^{r-1}(i+1)i=r^2-\sum_{i=1}^{r-1}i=\frac{r(r+1)}{2},\eqno(3.3.24)$$
$$(\vs,\vs)=\frac{s(s+1)}{2},\;\;(\mfk{e},\vs)^2=r^2s^2(\mfk{e}_r,\vs_s)^2=\frac{r^2s^2}{2}.
\eqno(3.3.25)$$ By the Schwartz inequality, we have:
$$(\mfk{e},\vs)^2<(\mfk{e},\mfk{e})(\vs,\vs)\sim\frac{r^2s^2}{2}<\frac{rs(r+1)(s+1)}{4}\eqno(3.3.26)$$
$$\sim 2rs<(r+1)(s+1)\sim(r-1)(s-1)<2\sim
r=1\;\mbox{or}\;s=1\;\mbox{or}\;r=s=2.\eqno(3.3.27)$$

Suppose that we have the fourth case with $p\geq q\geq r>1$. Let
$$\mfk{e}=\sum_{i=1}^{p-1}i\mfk{e}_i,\;\;\vs=\sum_{i=1}^{q-1}i\vs_i,\;\;
\xi=\sum_{i=1}^{r-1}i\xi_i.\eqno(3.3.28)$$ By (3.3.24), we have
$$(\mfk{e},\mfk{e})=\frac{p(p-1)}{2},\;\;(\vs,\vs)=\frac{q(r-1)}{2},\;\;(\xi,\xi)
=\frac{r(r-1)}{2}.\eqno(3.3.29)$$ Moreover, $\mfk{e},\;\vs$ and
$\xi$ are mutually orthogonal. Set
$$\be=\psi-\frac{(\psi,\mfk{e})}{(\mfk{e},\mfk{e})}\mfk{e}-\frac{(\psi,\vs)}{(\vs,\vs)}\vs-\frac{(\psi,\xi)}{(\xi,\xi)}\xi.
\eqno(3.3.30)$$ Then $0\neq\be $ is orthogonal to $\mfk{e},\;\vs$
and $\xi$. Furthermore,
$$0<(\be,\be)=(\be,\psi)=1-\frac{(\psi,\mfk{e})^2}{(\mfk{e},\mfk{e})}
-\frac{(\psi,\vs)^2}{(\vs,\vs)}-\frac{(\psi,\xi)^2}{(\xi,\xi)}.\eqno(3.3.31)$$

On the other hand,
$$4(\psi,\mfk{e})^2=4(p-1)^2(\psi,\mfk{e}_{p-1})^2=(p-1)^2,\eqno(3.3.32)$$
$$4(\psi,\vs)^2=(q-1)^2,\;\;4(\psi,\xi)^2=(r-1)^2.\eqno(3.3.33)$$ Thus
$$\frac{(p-1)^2}{(\mfk{e},\mfk{e})}
+\frac{(q-1)^2}{(\vs,\vs)}+\frac{(r-1)^2}{(\xi,\xi)}<4\sim
\frac{p-1}{p}+\frac{q-1}{q}+\frac{r-1}{r}<2\eqno(3.3.34)$$
$$\lra \frac{1}{p}+\frac{1}{q}+\frac{1}{r}>1.\eqno(3.3.35)$$
In particular,
$$\frac{3}{r}>1\sim r<3\lra r=2.\eqno(3.3.36)$$
Now (3.3.35) becomes
$$\frac{1}{p}+\frac{1}{q}>\frac{1}{2}\lra
\frac{2}{q}>\frac{1}{2}\sim q<4\sim q=2,3.\eqno(3.3.37)$$ If $q=3$,
then (3.3.35) becomes
$$\frac{1}{p}>\frac{1}{6}\sim p<6\sim p=2,3,4,5.\eqno(3.3.38)$$
In summary, the fourth case yields:
$$(p,q,r)=(p,2,2),\;(3,3,2),\;(4,3,2),\;(5,3,2).\eqno(3.3.39)$$

This gives all possible $CG({\msr S})$. Whenever a double edge or a
triple edge occurs in $CG(\Pi)$, we add an arrow pointing to the
shorter of two roots. The resulting figure is called the {\it Dynkin
diagram}\index{Dynkin diagram} of $\Phi$.\psp

{\bf Theorem 3.3.1}. {\it Let $\Phi$ be an irreducible root system.
The Dynkin diagram of $\Phi$ must be one of the followings:}

\begin{picture}(70,8)\put(2,0){$A_n$:}
\put(21,0){\circle{2}}\put(22,0){\line(1,0){12}}\put(21,-5){1}
\put(35,0){\circle{2}}\put(35,-5){2}\put(43,0){...}\put(53,0){\circle{2}}
\put(53,-5){n-1}
\put(54,0){\line(1,0){12}}\put(67,0){\circle{2}}\put(67,-5){n}
\end{picture}

\begin{picture}(80,12)\put(2,0){$B_n$:}\put(21,0){\circle{2}}\put(21,-5){1}\put(22,0)
{\line(1,0){12}} \put(35,0){\circle{2}}\put(35,-5){2}
\put(42,0){...}\put(51,0){\circle{2}}\put(52,0){\line(1,0){12}}
\put(51,-5){n-2}\put(65,0){\circle{2}}\put(65,-5){n-1}\put(65,1.2){\line(1,0){13.6}}
\put(65,-0.8){\line(1,0){13.6}}\put(78.5,0){\circle{2}}\put(78.5,-5){n}\put(72,-1){$\ra$}
\end{picture}

\begin{picture}(80,12)\put(2,0){$C_n$:}\put(21,0){\circle{2}}\put(21,-5){1}\put(22,0)
{\line(1,0){12}} \put(35,0){\circle{2}}\put(35,-5){2}
\put(42,0){...}\put(51,0){\circle{2}}\put(52,0){\line(1,0){12}}
\put(51,-5){n-2}\put(65,0){\circle{2}}\put(65,-5){n-1}\put(65,1.2){\line(1,0){13.6}}
\put(65,-0.8){\line(1,0){13.6}}\put(78.5,0){\circle{2}}\put(78.5,-5){n}\put(72,-1){$\la$}
\end{picture}

\begin{picture}(81,15)\put(2,0){$D_n$}\put(21,0){\circle{2}}\put(21,-5){1}\put(22,0){\line(1,0)
{12}} \put(35,0){\circle{2}}\put(35,-5){2}
\put(43,0){...}\put(52,0){\circle{2}}\put(52,-5){n-3}\put(53,0){\line(1,0){12}}
\put(66,0){\circle{2}}\put(66,-5){n-2}\put(66.6,0.5)
{\line(3,1){12}}\put(80,4.5){\circle{2}}\put(80,0){n-1}\put(66.6,-0.5)
{\line(3,-1){12}}\put(80,-4.5){\circle{2}}\put(80,-9){n}
\end{picture}

\begin{picture}(80,23)
\put(2,0){$E_6$:}\put(21,0){\circle{2}}\put(21,
-5){1}\put(22,0){\line(1,0){12}}\put(35,0){\circle{2}}\put(35,
-5){3}\put(36,0){\line(1,0){12}}\put(49,0){\circle{2}}\put(49,
-5){4}\put(49,1){\line(0,1){10}}\put(49,12){\circle{2}}\put(52,10){2}\put(50,0){\line(1,0){12}}
\put(63,0){\circle{2}}\put(63,-5){5}\put(64,0){\line(1,0){12}}\put(77,0){\circle{2}}\put(77,
-5){6}
\end{picture}

\begin{picture}(93,23)
\put(2,0){$E_7$:}\put(21,0){\circle{2}}\put(21,
-5){1}\put(22,0){\line(1,0){12}}\put(35,0){\circle{2}}\put(35,
-5){3}\put(36,0){\line(1,0){12}}\put(49,0){\circle{2}}\put(49,
-5){4}\put(49,1){\line(0,1){10}}\put(49,12){\circle{2}}\put(52,10){2}\put(50,0){\line(1,0){12}}
\put(63,0){\circle{2}}\put(63,-5){5}\put(64,0){\line(1,0){12}}\put(77,0){\circle{2}}\put(77,
-5){6}\put(78,0){\line(1,0){12}}\put(91,0){\circle{2}}\put(91,
-5){7}
\end{picture}

\begin{picture}(110,23)
\put(2,0){$E_8$:}\put(21,0){\circle{2}}\put(21,
-5){1}\put(22,0){\line(1,0){12}}\put(35,0){\circle{2}}\put(35,
-5){3}\put(36,0){\line(1,0){12}}\put(49,0){\circle{2}}\put(49,
-5){4}\put(49,1){\line(0,1){10}}\put(49,12){\circle{2}}\put(52,10){2}\put(50,0){\line(1,0){12}}
\put(63,0){\circle{2}}\put(63,-5){5}\put(64,0){\line(1,0){12}}\put(77,0){\circle{2}}\put(77,
-5){6}\put(78,0){\line(1,0){12}}\put(91,0){\circle{2}}\put(91,
-5){7}\put(92,0){\line(1,0){12}}\put(105,0){\circle{2}}\put(105,
-5){8}
\end{picture}

\begin{picture}(60,12)\put(2,0){$F_4$:}
\put(21,0){\circle{2}}\put(21,-5){1}\put(22,0){\line(1,0){12}}
\put(35,0){\circle{2}}\put(35,-5){2}\put(35,1.2){\line(1,0){13.6}}
\put(35,-0.8){\line(1,0){13.6}}\put(41,-1){$\ra$}\put(48.5,0){\circle{2}}\put(48.5,-5){3}\put(49.5,0)
{\line(1,0){12}}\put(62.5,0){\circle{2}}\put(62.5,-5){4}
\end{picture}

\begin{picture}(27,10)\put(2,0){$G_2$:}\put(21,0){\circle{2}}\put(21,1.2){\line(1,0){14}}
\put(22,0.285){\line(1,0){12}}\put(21,-0.8){\line(1,0){14}}\put(35,0){\circle{2}}
\put(35,-5){2}\put(21,-5){1}\put(27,-1){$\ra$}\end{picture}

\vspace{1cm}

As an exercise, calculate the determinants of the Cartan matrices of
the above root systems, which are as follows:
$$A_n:\;n+1;\;B_n:\;2;\;C_n:\;2;\;D_n:\;4;\;E_6:\;3;\;E_7:\;2;\;E_8,\;F_4,\;G_2:
\;1.\eqno(3.3.40)$$

\section{Automorphisms, Constructions and Weights}

In this section, we construct all the irreducible root systems and
their associated weight lattices. Moreover, their automorphism
groups are determined and the saturated subsets of weight lattices
are investigated.

Let $\Phi$ be a root system in a Euclidean space ${\msr E}$ and let
$\Pi=\{\al_1,\al_2,...,\al_n\}$ be a base of $\Phi$. The additive
subgroup
$$\Lmd_r=\sum_{i=1}^n\mbb{Z}\al_i\subset{\msr E}\eqno(3.4.1)$$
is called the {\it root lattice} of $\Phi$.\index{root lattice}
Moreover,
$$\Lmd=\{\gm\in{\msr E}\mid
\la\gm,\al\ra\in\mbb{Z}\;\for\;\al\in\Phi\}=\{\gm\in{\msr E}\mid
\la\gm,\al\ra\in\mbb{Z}\;\for\;\al\in\Pi\}\eqno(3.4.2)$$ is called
the {\it weight lattice} of $\Phi$,\index{weight lattice} whose
elements are called {\it weights}.\index{weight} Obviously,
$\Lmd\supset \Lmd_r$. We call $\Lmd/\Lmd_r$ the {\it fundamental
group} of $\Phi$.\index{fundamental group}

Since the inner product of ${\msr E}$ is nondegenerate, there exist
$\lmd_1,\lmd_2,...,\lmd_n\in{\msr E}$ such that
$$\la\lmd_i,\al_j\ra=\dlt_{i,j}\qquad\for\;\;i,j\in\ol{1,n}.\eqno(3.4.3)$$
We call $\{\lmd_1,...,\lmd_n\}$ {\it fundamental dominant weights
relative to} $\Pi$. Moreover,
$$\Lmd =\sum_{i=1}^n\mbb{Z}\lmd_i.\eqno(3.4.4)$$
Write
$$\al_i=\sum_{r=1}^nm_{i,r}\lmd_r\qquad\for\;\;i\in\ol{1,n}.\eqno(3.4.5)$$
Then
$$\la\al_i,\al_j\ra=\la\sum_{r=1}^nm_{i,r}\lmd_r,\al_j\ra=m_{i,j}.\eqno(3.4.6)$$
Thus
$$\al_i=\sum_{j=1}^n\la\al_i,\al_j\ra\lmd_j\qquad\for\;\;i\in\ol{1,n}.\eqno(3.4.7)$$
By the fundamental theorem of finitely generated abelian groups,
$$|\Lmd/\Lmd_r|=|(\la\al_i,\al_j\ra)_{n\times n}|\eqno(3.4.8)$$
is finite. Moreover,
$$\lmd_i\in\frac{1}{|(\la\al_i,\al_j\ra)_{n\times
n}|}\Lmd_r\qquad\for\;\;i\in\ol{1,n}.\eqno(3.4.9)$$ By (3.2.15) and
(3.2.17),
$$\rho=\frac{1}{2}\sum_{\al\in\Phi^+}\al=\sum_{i=1}^n\lmd_i.\eqno(3.4.10)$$
\pse

{\bf Example 3.4.1}. Let ${\msr E}$ be the Euclidean space with a
basis $\{\ves_1,\ves_2,...,\ves_n\}$ such that
$(\ves_i,\ves_j)=\dlt_{i,j}$. The root system of $A_{n-1}$ and a
base are as follows:
$$\Phi_{A_{n-1}}=\{\ves_i-\ves_j\mid i,j\in\ol{1,n},\;i\neq
j\},\;\;\Pi_{A_{n-1}}=\{\ves_i-\ves_{i+1}\mid
i\in\ol{1,n-1}\}.\eqno(3.4.11)$$ Set
$$Q=\sum_{i=1}^n\mbb{Z}\ves_i.\eqno(3.4.12)$$
The root lattice of $\Phi_{A_{n-1}}$:
$$\Lmd_r(A_{n-1})=\{\sum_{i=1}^na_i\ves_i\in Q\mid
\sum_{i=1}^na_i=0\}\eqno(3.4.13)$$ and
$$\lmd_i=\frac{1}{n}(\sum_{r=1}^i(n-i)\ves_r-\sum_{s=i+1}^ni\ves_s)\qquad
\for\;\;i\in\ol{1,n-1}.\eqno(3.4.14)$$ The group
$\Lmd(\Phi_{A_{n-1}})/\Lmd_r(\Phi_{A_{n-1}})\cong
\mbb{Z}_n=\mbb{Z}/n\mbb{Z}$, the cyclic group of order
$n$(exercise).

The root system of $B_n$ is
$$\Phi_{B_n}=\{\ves_i-\ves_j,\pm\ves_i,\pm(\ves_i+\ves_j)\mid i,j\in\ol{1,n},\;i\neq
j\}\eqno(3.4.15)$$ with a base
$$\Pi_{B_n}=\{\ves_1-\ves_2,...,\ves_{n-1}-\ves_n,\ves_n\}.\eqno(3.4.16)$$
The root lattice
$$\Lmd_r(B_n)=Q\eqno(3.4.17)$$
and
$$\lmd_i=\sum_{r=1}^i\ves_r,\;\;\lmd_n=\frac{1}{2}\sum_{s=1}^n\ves_s,\qquad\for\;\;i\in\ol{1,n-1}.
\eqno(3.4.18)$$ Thus the weight lattice:
$$\Lmd(B_n)=Q\bigcup (\lmd_n+Q).\eqno(3.4.19)$$
The fundamental group: $\Lmd(B_n)/\Lmd_r(B_n)\cong \mbb{Z}_2$.

The root system of $C_n$ is
$$\Phi_{C_n}=\{\ves_i-\ves_j,\pm 2\ves_i,\pm(\ves_i+\ves_j)\mid i,j\in\ol{1,n},\;i\neq
j\}\eqno(3.4.20)$$ with a base
$$\Pi_{C_n}=\{\ves_1-\ves_2,...,\ves_{n-1}-\ves_n,2\ves_n\}.\eqno(3.4.21)$$
The root lattice
$$\Lmd_r(C_n)=\{\sum_{i=1}^na_i\ves_i\in Q\mid \sum_{i=1}^na_i\in 2\mbb{Z}\}\eqno(3.4.22)$$
and
$$\lmd_i=\sum_{r=1}^i\ves_i,\;\;\lmd_n=\sum_{s=1}^n\ves_s,\qquad\for\;\;i\in\ol{1,n-1}.
\eqno(3.4.23)$$ Thus the weight lattice:
$$\Lmd(C_n)=Q.\eqno(3.4.24)$$
The fundamental group: $\Lmd(C_n)/\Lmd_r(C_n)\cong \mbb{Z}_2$.

The root system of $D_n$ is
$$\Phi_{D_n}=\{\ves_i-\ves_j,\pm(\ves_i+\ves_j)\mid i,j\in\ol{1,n},\;i\neq
j\}\eqno(3.4.25)$$ with a base
$$\Pi_{D_n}=\{\ves_1-\ves_2,...,\ves_{n-1}-\ves_n,\ves_{n-1}+\ves_n\}.\eqno(3.4.26)$$
The root lattice and weight lattice:
$$\Lmd_r(D_n)=\Lmd_r(C_n).\eqno(3.4.27)$$
$$\lmd_i=\sum_{r=1}^i\ves_i\;\;\lmd_{n-1}=\frac{1}{2}(\sum_{r=1}^{n-1}\ves_i-\ves_n),\;\;
\lmd_n=\frac{1}{2}\sum_{s=1}^n\ves_s,\qquad\for\;\;i\in\ol{1,n-2}.
\eqno(3.4.28)$$ Thus the weight lattice
$$\Lmd(D_n)=\Lmd(B_n)=Q\bigcup (\lmd_n+Q).\eqno(3.4.29)$$
The fundamental group:
$$\Lmd(D_n)/\Lmd_r(D_n)\cong \left\{\begin{array}{ll}\mbb{Z}_4&\mbox{if}\;n\in 2\mbb{Z}+1\\
 \mbb{Z}_2\times\mbb{Z}_2&\mbox{if}\;n\in
 2\mbb{Z}\end{array}\right.\eqno(3.4.30)$$
\psp

{\bf Example 3.4.2}. The root and weight lattices of $E_8$:
$$\Lmd_r(E_8)=\Lmd(D_8)=\Lmd(E_8),\eqno(3.4.31)$$
which is a self-dual integral even lattice. The root system
$$\Phi_{E_8}=\left\{\frac{1}{2}\sum_{i=1}^8(-1)^{\iota_i}\ves_i\mid
\iota_i\in\mbb{Z}_2,\;\sum_{r=1}^8\iota_r=0\;\mbox{in}\;\mbb{Z}_2\right\}\bigcup
\Phi_{D_8}\eqno(3.4.32)$$ with a base
$$\Pi_{E_8}=\left\{\ves_1-\ves_2,\frac{1}{2}(\sum_{s=4}^8\ves_s-\sum_{r=1}^3\ves_r),\ves_{i+1}-\ves_{i+2}
\mid i\in\ol{1,6}\right\}.\eqno(3.4.33)$$ Moreover,
$$\lmd_1=\frac{1}{2}(3\ves_1+\sum_{i=2}^8\ves_i),\;\;\lmd_2=\sum_{i=1}^8\ves_i,\;\;\lmd_3
=2(\ves_1+\ves_2)+\sum_{i=3}^8\ves_i,\eqno(3.4.34)$$
$$\lmd_4=\frac{1}{2}(5(\ves_1+\ves_2+\ves_3)+3\sum_{i=4}^8\ves_i),\;\;
\lmd_5=\sum_{i=1}^4(2\ves_i+\ves_{4+i}),\eqno(3.4.35)$$
$$\lmd_6=\frac{1}{2}(3\sum_{i=1}^5\ves_i+\ves_6+\ves_7+\ves_8),\;\;\lmd_7=\sum_{i=1}^6\ves_i,
\;\;\lmd_8=\frac{1}{2}(\sum_{i=1}^7\ves_i-\ves_8).\eqno(3.4.36)$$

Observe
$$\Pi_{E_7}=\left\{\ves_1-\ves_2,\frac{1}{2}(\sum_{r=1}^3\ves_r-\sum_{s=4}^8\ves_s),\ves_{i+1}-\ves_{i+2}
\mid i\in\ol{1,5}\right\}\eqno(3.4.37)$$ and
$$\Pi_{E_6}=\left\{\ves_1-\ves_2,\frac{1}{2}(\sum_{r=1}^3\ves_r-\sum_{s=4}^8\ves_s),\ves_{i+1}-\ves_{i+2}
\mid i\in\ol{1,4}\right\}.\eqno(3.4.38)$$ Moreover, the fundamental
groups:
$$\Lmd(E_7)/\Lmd_r(E_7)\cong\mbb{Z}_2,\qquad
\Lmd(E_6)/\Lmd_r(E_6)\cong\mbb{Z}_3\eqno(3.4.39)$$ (exercise).\psp

{\bf Example 3.4.3}. The root and weight lattices of $F_4$:
$$\Lmd_r(F_4)=\Lmd(F_4)=\Lmd(D_4).\eqno(3.4.40)$$
Moreover,
$$\Phi_{F_4}=\Phi_{B_4}\bigcup\left\{\frac{1}{2}\sum_{i=1}^4(-1)^{\iota_i}\ves_i\mid
\iota_i\in\mbb{Z}_2\right\}\eqno(3.4.41)$$ with a base
$$\Pi_{F_4}=\left\{\ves_2-\ves_3,\ves_3-\ves_4,\ves_4,\frac{1}{2}(\ves_1-\ves_2-\ves_3-\ves_4)\right\}.
\eqno(3.4.42)$$ Moreover,
$$\lmd_1=\ves_1+\ves_2,\;\;\lmd_2=2\ves_1+\ves_2+\ves_3,\;\;
\lmd_3=\frac{1}{2}(3\ves_1+\ves_2+\ves_3+\ves_4),\;\;\lmd_4=\ves_1.\eqno(3.4.43)$$

{\bf Example 3.4.4}.  The root and weight lattices of $G_2$:
$$\Lmd_r(G_2)=\Lmd(G_2)=\Lmd_r(A_2).\eqno(3.4.44)$$
Moreover,
$$\Phi_{G_2}=\{\pm(2\ves_i-\ves_j-\ves_k),\ves_i-\ves_j\mid\{i,j,k\}=\{1,2,3\}\}\eqno(3.4.45)$$
with a base
$$\Pi_{G_2}=\{\ves_2+\ves_3-2\ves_1,\ves_1-\ves_2\}.\eqno(3.4.46)$$
Moreover,
$$\lmd_1=2\ves_3-\ves_1-\ves_2,\qquad\lmd_2=\ves_3-\ves_2.\eqno(3.4.47)$$
\pse

Let $\Phi$ be a root system in a Euclidean space ${\msr E}$ and let
$\Pi$ be a base of $\Phi$. Recall that $\Lmd$ is the set of weights
of $\Phi$. Set
$$\Lmd^+=\Lmd\bigcap \ol{{\msr C}(\Pi)}\eqno(3.4.48)$$
and call the elements in $\Lmd^+$ {\it dominant
weights}.\index{dominant weight} Moreover, we call the elements of
$\Lmd\bigcap{\msr C}(\Pi)$ {\it strong dominant
weights}.\index{strong dominant weight} We define another partial
ordering $\lhd$ on $\Lmd$ by
$$\lmd\lhd\mu\;\mbox{if}\;\mu-\lmd=\sum_{\al\in\Pi}a_\al\al\;\mbox{such
that}\;0\leq a_\al\in\mbb{Z}, \eqno(3.4.49)$$ where
$\lmd,\mu\in\Lmd$. By (3.1.4) and (3.4.2), two partial orderings
coincide on ${\msr W}(\lmd)$ for any $\lmd\in\Lmd$. According to
Lemma 3.2.12, any weight is conjugated to a dominant weight under
the Weyl group ${\msr W}$. If $\lmd$ is a dominant weight, then
$\sgm(\lmd)\lhd \lmd$ for $\sgm\in{\msr W}.$

If $\lmd$ is a strong dominant weight and $\sgm(\lmd)=\lmd$ with
$\sgm\in{\msr W}$, then $\sgm=1$ by Theorem 3.2.8 (e).\psp

{\bf Lemma 3.4.1}. {\it For $\lmd\in\Lmd^+$, $|\{\mu\in\Lmd^+\mid
\mu\lhd\lmd\}|<\infty$.}

{\it Proof}. Note
$$0\leq (\lmd,\lmd-\mu)=(\lmd,\lmd)-(\lmd,\mu)\eqno(3.4.50)$$
and
$$0\leq (\lmd-\mu,\mu)=(\lmd,\mu)-(\mu,\mu).\eqno(3.4.51)$$
Adding two inequalities, we obtain
$$0\leq (\lmd,\lmd)-(\mu,\mu)\sim (\mu,\mu)\leq
(\lmd,\lmd).\eqno(3.4.52)$$ Since the intersection of $\Lmd$ with
the sphere of radius $||\lmd||$ is a finite set, so is its subset
$\{\mu\in\Lmd^+\mid \mu\lhd\lmd\}.\qquad\Box$\psp

A subset $\mbb{S}$ of $\Lmd$ is called {\it
saturated}\index{saturated} if $\mu-i\al\in\mbb{S}$ for
$\mu\in\mbb{S},\;\al\in\Phi$ and $i$  between $0$ and
$\la\mu,\al\ra$. By definition ${\msr W}(\mbb{S})=\mbb{S}$. If
$\mbb{S}$ has a unique maximal weight $\lmd$, we call $\lmd$ the
{\it highest weight} of $\mbb{S}$.\index{highest weight} Lemma
3.2.12 says that $\lmd\in\Lmd^+$. By Lemma 3.2.12 and Lemma 3.4.1,
we have:\psp

{\bf Lemma 3.4.2}. {\it A saturated set of weights with a highest
weight must be finite}.\psp

{\bf Lemma 3.4.3}. {\it Let $\mbb{S}$ be a saturated subset with a
highest weight $\lmd$. If $\lmd\rhd\mu\in\Lmd^+$ , then
$\mu\in\mbb{S}$.}

{\it Proof}. Assume $\mu\neq\lmd$. Suppose
$\lmd'=\mu+\sum_{\al\in\Pi}k_\al\al\in\mbb{S}$ with $0\leq
k_\al\in\mbb{Z}$ and $\sum_{\al\in\Pi}k_\al>0$. For instance $\lmd$
is such a $\lmd'$. Since
$$0<(\sum_{\al\in\Pi}k_\al\al,\sum_{\al\in\Pi}k_\al\al)=
\sum_{\be\in\Pi}k_\be(\sum_{\al\in\Pi}k_\al\al,\be),\eqno(3.4.53)$$
there exists $\gm\in\Pi$ such that
$$k_\gm(\sum_{\al\in\Pi}k_\al\al,\gm)>0\sim
k_\gm,\;(\sum_{\al\in\Pi}k_\al\al,\gm)>0.\eqno(3.4.54)$$ So
$$\la\lmd',\gm\ra=\la\mu,\gm\ra+\la\sum_{\al\in\Pi}k_\al\al,\gm\ra>0.\eqno(3.4.55)$$
By definition, $\lmd'-\gm\in\mbb{S}$. An induction on
$\sum_{\al\in\Pi}k_\al$ shows $\mu\in\mbb{S}.\qquad\Box$\psp

{\bf Lemma 3.4.4}. {\it The subset $\mbb{S}$ is a saturated subset
of weights with a highest weight $\lmd$ if and only if
$$\lmd\in \Lmd^+,\;\;\mbb{S}={\msr W}(\{\mu\in\Lmd^+\mid \mu\lhd
\lmd\}).\eqno(3.4.56)$$}

{\it Proof}. Set $\mbb{S}_\lmd={\msr W}(\{\mu\in\Lmd^+\mid \mu\lhd
\lmd\})$. If $\mbb{S}$ is a saturated subset of weights with a
highest weight $\lmd$, then (3.4.56) holds by Lemmas 3.2.12 and
3.4.3.

Suppose $\lmd\in\Lmd^+$. Then ${\msr W}(\mbb{S}_\lmd)=\mbb{S}_\lmd$
and $\lmd$ is the highest weight of $\mbb{S}_\lmd$. For any
$\lmd\rhd\mu\in\Lmd^+,\;\sgm\in {\msr W},\; \al\in \Phi$ and $i$
between $0$ and $\la\sgm(\mu),\al\ra$, there exists $\tau\in{\msr
W}$ such that $\tau(\sgm(\mu)-i\al)\in\Lmd^+$. Note
$\tau\sgm(\mu)\lhd \mu\lhd \lmd$ by Lemma 3.2.12. First, we assume
$\tau(\al)\in\Phi^+$. If $\la\sgm(\mu),\al\ra\geq 0$, then $i\geq
0$. So
$$\tau(\sgm(\mu)-i\al)=\tau\sgm(\mu)-i\tau(\al)\lhd
\tau\sgm(\mu)\lhd \lmd.\eqno(3.4.57)$$ When $\la\sgm(\mu),\al\ra<0$,
we have $i<0$ and
$$\tau\sgm(\mu)-i\tau(\al)\lhd\tau\sgm(\mu)-\la\sgm(\mu),\al\ra\tau(\al)=
\tau\sgm_\al\sgm(\mu)\lhd \mu\lhd \lmd\eqno(3.4.58)$$ by Lemma
3.2.12. Hence $\tau(\sgm(\mu)-i\al)\in \mbb{S}_\lmd$. Thus
$\sgm(\mu)-i\al\in \mbb{S}_\lmd$.

Assume $\tau(\al)\in \Phi^-$. Then
$$\sgm(\mu)-i\al=\sgm(\mu)-(-i)(-\al)\eqno(3.4.59)$$
and $-i$ is between $0$ and $\la\sgm(\mu),-\al\ra$. Replacing $\al$
by $-\al$ in the above arguments, we have
$$\sgm(\mu)-i\al=\sgm(\mu)-(-i)(-\al)\in \mbb{S}_\lmd.\eqno(3.4.60)$$
Therefore, $\mbb{S}_\lmd$ is saturated.$\qquad\Box$\psp

The automorphism group of $\Phi$:
$$\mbox{Aut}\:\Phi=\{\sgm\in GL({\msr E})\mid
\sgm(\Phi)=\Phi\}.\eqno(3.4.61)$$ By Lemma 3.2.2, ${\msr W}$ is a
normal subgroup of $\mbox{Aut}\:\Phi$ and
$$\la\sgm(\al),\sgm(\be)\ra=\la\al,\be\ra\qquad \for\;\;\sgm\in
\mbox{Aut}\:\Phi,\;\al,\be\in\Phi.\eqno(3.4.62)$$ Since $\Phi$ spans
${\msr E}$, the above equality holds for any $\al\in{\msr E} $ and
$\be\in\Phi$. So $\sgm$ maps regular elements to regular elements.
Thus $\sgm(\Pi)$ is a base of $\Phi$ for any $\sgm\in
\mbox{Aut}\:\Phi$. By Theorem 3.2.8 (b), there exists
$\sgm_1\in{\msr W}$ such that $\sgm_1\sgm(\Pi)=\Pi$. Hence
$\mbox{Aut}\:\Phi/{\msr W}$ is completely determined by
$$\{\sgm\in GL({\msr E})\mid
\sgm(\Pi)=\Pi,\;\la\sgm(\al),\sgm(\be)\ra=\la\al,\be\ra\;\for\;\;\al,\be\in\Pi\}
\eqno(3.4.63)$$ because $\Phi={\msr W}(\Pi)$. By Theorem 3.3.1, we
have:\psp

{\bf Theorem 3.4.5}. {\it Let $\Phi$ be an irreducible root system.
Then $\mbox{Aut}\:\Phi/{\msr W}$ is completely determined as
follows:
$$\mbox{Aut}\:\Phi_{A_n}/{\msr W}\cong\mbox{Aut}\:\Phi_{E_6}/{\msr
W}\cong\mbb{Z}_2,\eqno(3.4.64)$$
$$\mbox{Aut}\:\Phi_{D_4}/{\msr W}\cong S_3,\;\;\mbox{Aut}\:\Phi_{D_n}/{\msr W}\cong
\mbb{Z}_2\qquad\for\;\;n>4,\eqno(3.4.65)$$ and
$\mbox{Aut}\:\Phi/{\msr W}=\{1\}$ for the other types of
$\Phi$.}\psp

Exercise: prove the Weyl group of $G_2$ is isomorphic the dihedral
group ${\msr D}_6$.

\chapter{Isomorphisms, Conjugacy and Exceptional Types}

In this chapter, we show that  the structure of a finite-dimensional
semisimple Lie algebra over $\mbb{C}$ is completely determined by
its root system. Moreover, we prove that any two Cartan subalgebras
of such a Lie algebra ${\msr G}$ are conjugated under the group of
inner automorphisms of ${\msr G}$. In particular, the automorphism
group of ${\msr G}$ is determined when it is simple. Furthermore, we
give explicit constructions of the simple Lie algebras of
exceptional types.

\section{Isomorphisms}

In this section, we  prove that two semisimple Lie algebras are
isomorphic if their root systems are.

Let ${\msr G}$ be a finite-dimensional semisimple Lie algebra over
$\mbb{C}$ and let $H$ be a maximal toral subalgebra of ${\msr G}$.
Then we have the Contan root decomposition:
$${\msr G}=H+\sum_{\al\in\Phi}{\msr G}_\al,\eqno(4.1.1)$$
where $\Phi\subset H^\ast$ is a root system of ${\msr G}$.\psp

{\bf Lemma 4.1.1}. {\it If ${\msr G}$ is simple, then $\Phi$ is
irreducible}.

{\it Proof}. Suppose
$\Phi=\Phi_1\bigcup\Phi_2,\;\Phi_1\neq\emptyset\neq\Phi_2$ and
$\Phi_1\bot\Phi_2$. For any $\al_1\in\Phi_1$ and $\al_2\in\Phi_2$,
$\al_1+\al_2$ is not orthogonal to $\Phi_1$ and $\Phi_2$. So
$\al_1+\al_2\not\in\Phi$. Hence $[{\msr G}_{\al_1},{\msr
G}_{\al_2}]=\{0\}$.  Let ${\msr G}_i$ be the Lie algebra generated
by $\{{\msr G}_\al\mid\al\in\Phi_i\}$. Then ${\msr G}={\msr
G}_1+{\msr G}_2$ and $[{\msr G}_1,{\msr G}_2]=\{0\}$, which
contradicts the simplicity of ${\msr G}.\qquad\Box$\psp

 {\bf Corollary 4.1.2}. {\it In general, if ${\msr G}={\msr
G}_1\oplus{\msr G}_2\oplus\cdots\oplus {\msr G}_k$ is the
decomposition of ${\msr G}$ into simple ideals, then $H_i=H\bigcap
{\msr G}_i$ is a maximal toral subalgebra of ${\msr G}_i$ whose
corresponding root system $\Phi_i$ are irreducible components of
$\Phi$; that is, $\Phi=\bigcup_{i=1}^k\Phi_i$.}\psp

{\bf Lemma 4.1.3}. {\it If ${\msr G}$ is simple, then it is
generated by $\{{\msr G}_\al,{\msr G}_{-\al}\mid\al\in\Pi\}$}.

{\it Proof}. According to Corollary 3.1.6, any positive root
$\al=\sum_{i=1}^k\be_i$ with $\be_i\in\Pi$ such that
$\sum_{i=1}^s\be_i$ are positive roots for $s\in\ol{1,k}$. By Lemma
2.6.2 (c),
$${\msr G}_\al=[{\msr G}_{\be_k},[{\msr G}_{\be_{k-1}},...,[{\msr
G}_{\be_2},{\msr G}_{\be_1}]...]].\eqno(4.1.2)$$ Symmetrically,
$${\msr G}_{-\al}=[{\msr G}_{-\be_k},[{\msr G}_{-\be_{k-1}},...,[{\msr
G}_{-\be_2},{\msr G}_{-\be_1}]...]].\eqno(4.1.3)$$ Since ${\msr G}$
is generated by $\{{\msr G}_\al\mid\al\in\Phi\}$ by Lemma 2.6.2 (e),
it is generated by $\{{\msr G}_\al,{\msr
G}_{-\al}\mid\al\in\Pi\}.\qquad\Box$\psp

{\bf Theorem 4.1.4}. {\it Suppose that ${\msr G}$ is simple. Let
${\msr G}'$ be another finite-dimensional simple Lie algebra and
$H'$ is a maximal toral subalgebra ${\msr G}'$ with corresponding
root system $\Phi'$. If $\Phi\cong\Phi'$, then ${\msr G}\cong {\msr
G}'$.}

{\it Proof}. Let ${\msr E}$ be the Euclidean space associated with
$\Phi$ and let ${\msr E}'$ be the Euclidean space associated with
$\Phi'$. Suppose that $\tau:\al\mapsto \al'$ be an isomorphism from
$\Phi$ to $\Phi'$; that is, $\tau$ is a linear map from ${\msr E}$
to ${\msr E}'$ such that $\tau(\Phi)=\Phi'$ and
$$\la\al',\be'\ra=\la\al,\be\ra\qquad\for\;\;\al,\be\in\Phi.\eqno(4.1.4)$$
Note that $\Pi'=\tau(\Pi)$ is a base of $\Phi'$.

Let ${\msr L}={\msr G}\oplus {\msr G}'$ be the direct sum of two Lie
algebras. So ${\msr L}$ is semisimple, and ${\msr G},\;{\msr G}'$
are the only simple ideals. Take
$$\{\xi_\al\in{\msr
G}_\al,\;\zeta_\al\in{\msr G}_{-\al},\xi'_{\al'}\in{\msr
G}'_{\al'},\zeta'_{\al'}\in{\msr G}'_{-\al'}
\mid\al\in\Pi\}\eqno(4.1.5)$$ such that
$$[\xi_\al,\zeta_\al]=h_{\al},\;\;[h_\al,\xi_\al]=2\xi_\al,\;\;[h_\al,\zeta_\al]
=-2\zeta_\al\eqno(4.1.6)$$ and
$$[\xi'_{\al'},\zeta'_{\al'}]=h_{\al'}',\;\;[h_{\al'}',\xi'_{\al'}]=2\xi_{\al'},\;\;[h_{\al'}',
\zeta'_{\al'}] =-2\zeta'_{\al'}\eqno(4.1.7)$$ for $\al\in\Pi$. Let
$${\msr K}=\mbox{the Lie subalgebra generated
by}\;\;\{\xi_\al+\xi'_{\al'},\zeta_\al+\zeta'_{\al'}\mid\al\in\Pi\}.\eqno(4.1.8)$$
Define ${\msr P},{\msr P}'\in\Edo {\msr L}$ by
$${\msr P}(u+u')=u,\;\;{\msr P}'(u+u')=u'\qquad\for\;\;u\in{\msr
G},\;u'\in{\msr G}'.\eqno(4.1.9)$$ Then
$${\msr P}({\msr K})={\msr G},\;\;{\msr P}'({\msr K})={\msr
G}'\eqno(4.1.10)$$ by Lemma 4.1.3.

Let $\be$ be the maximal root in $\Phi$ relative to $\Pi$. Then
$\be'$ is the maximal root in $\Phi'$ relative to $\Pi'$. Take any
$0\neq \zeta\in{\msr G}_{\be}$ and $0\neq \zeta'\in{\msr
G}'_{\be'}$. Set
$${\msr
M}=\mbox{span}\{\zeta+\zeta',\ad(\zeta_{\al_1}+\zeta'_{\al_1'})\cdots\ad(\zeta_{\al_{_k}}+\zeta'_{\al_{_k}'})(\zeta+\zeta')
\mid\al_i\in\Pi,\;k\in\mbb{Z}_+\}.\eqno(4.1.11)$$ Note
$$[\xi_\al+\xi'_{\al'},\zeta_\gm+\zeta'_{\gm'}]=\dlt_{\al,\gm}(h_\al+h'_{\al'})\qquad\for\;\;
\al,\gm\in\Pi\eqno(4.1.12)$$ by (4.1.6) and (4.1.7), because
$\al-\gm$ is not a root. Moreover,
$$[h_\al+h'_{\al'},\zeta_\gm+\zeta'_{\gm'}]=\gm(h_\al)\zeta_\gm+\gm'(h'_{\al'})\zeta'_{\gm'}
=\la\gm,\al\ra\zeta_\gm+\la\gm',\al'\ra\zeta'_{\gm'}=\la\gm,\al\ra(\zeta_\gm+\zeta'_{\gm'})
\eqno(4.1.13)$$ and
$$[h_\al+h'_{\al'},\zeta+\zeta']=\la\be,\al\ra(\zeta+\zeta')\eqno(4.1.14)$$
by (4.1.4). Furthermore,
$$[\xi_\al+\xi'_{\al'},\zeta+\zeta']=0\qquad\for\;\;\al\in\Pi,\eqno(4.1.15)$$
due to $\be+\al\not\in\Phi$ and $\be'+\al'\not\in\Phi'$ by the
maximality of $\be$. Equations (4.1.12)-(4.1.15) imply that
$$[{\msr K},{M}]\subset{M}.\eqno(4.1.16)$$

Since
$${M}\bigcap({\msr G}_{\be}+{\msr
G}'_{\be'})=\mbb{F}(\zeta+\zeta'),\eqno(4.1.17)$$ we have
$${M}\neq {\msr L},{\msr G},{\msr G}'.\eqno(4.1.18)$$
Thus ${M}$ is not an ideal of ${\msr L}$. This means ${\msr D}\neq
{\msr L}$. If ${\msr K}\bigcap {\msr G}\neq\{0\}$, then
$$[{\msr K},{\msr K}\bigcap {\msr G}]=[{\msr G},{\msr K}\bigcap {\msr
G}]\subset{\msr K}\bigcap {\msr G}\eqno(4.1.19)$$
 by (4.1.10). So ${\msr K}\bigcap {\msr
G}$ is a nonzero ideal of ${\msr G}$. By the simplicity of ${\msr
G}$, ${\msr G}={\msr K}\bigcap {\msr G}$, which implies ${\msr
L}={\msr K}$ by (4.1.10). This contradicts ${\msr K}\neq {\msr L}$.
Thus ${\msr K}\bigcap {\msr G}=\{0\}$. Symmetrically, ${\msr
K}\bigcap {\msr G}'=\{0\}$. These facts show that the map
$$u\mapsto \Im(u)\;\;\mbox{determined by}\;\;u+\Im(u)\in{\msr K},\qquad
u\in{\msr G},\eqno(4.1.20)$$ is a Lie isomorphism from ${\msr G}$ to
${\msr G}'.\qquad\Box$ \psp

Note
$$\Im(\xi_\al)=\xi'_{\al'},\;\;\Im(\zeta_\al)=\zeta'_{\al'},\;\Im(h_\al)=h'_{\al'}\eqno(4.1.21)$$
by (4.1.6)-(4.1.8). In particular, any automorphism $\sgm$ of $\Phi$
gives rise to an automorphism $\Im_\sgm$ (not unique) of ${\msr G}$
such that
$$\Im_\sgm(H)=H,\;\;\Im_\sgm({\msr G}_\al)={\msr
G}_{\sgm(\al)}\qquad\for\;\;\al\in \Phi.\eqno(4.1.22)$$ Note that
$-\mbox{Id}_{\msr E}\in\mbox{Aut}\:\Phi$. We have a unique
automorphism $\vt$ of ${\msr G}$ such that
$$\vt(\xi_\al)=-\zeta_\al,\;\;\vt(\zeta_\al)=-\xi_\al,\;\;\vt(h_\al)=-h_\al
\qquad\for\;\;\al\in\Pi.\eqno(4.1.23)$$ The automorphism $\vt$ is
called a Cartan {\it involution}. For any $\al\in\Pi$, $\ad\xi_\al$
and $\ad\zeta_\al$ are nilpotent by Lemma 2.5.2. The linear
transformation
$$\tau_\al=e^{\ad \xi_\al}e^{-\ad\zeta_\al}e^{\ad
\xi_\al}\eqno(4.1.24)$$ is an automorphism of ${\msr G}$ by (1.4.27)
satisfying
$$\tau_\al({\msr G}_\be)={\msr G}_{\sgm_\al(\be)}\qquad
\for\;\;\be\in\Phi\eqno(4.1.25)$$ and $\tau_\al=\vt$ when ${\msr
G}=sl(2,\mbb{C})$ (exercise).

\section{Cartan Subalgebras}

In this section, we give characterizations of Cartan subalgebras.
Moreover, we prove a correspondence between Cartan subalgebras under
homomorphisms.

 Let ${\msr G}$ be a finite-dimensional Lie algebra
over $\mbb{C}$. For any $\xi\in{\msr G}$, we have
$${\msr G}=\bigoplus_{a\in\mbb{F}}{\msr G}_\al(\ad\xi),\eqno(4.2.1)$$ where
$${\msr G}_a(\ad\xi)=\{u\in{\msr
G}\mid(\ad\xi-a)^m(u)=0\;\;\mbox{for
some}\;m\in\mbb{N}\}.\eqno(4.2.2)$$ Moreover,
$$[{\msr G}_a(\ad\xi),{\msr G}_b(\ad\xi)]={\msr G}_{a+b}(\ad\xi)\qquad\for\;\;a,b\in\mbb{F}\eqno(4.2.3)$$ by
(1.5.60). In particular, ${\msr G}_0(\ad\xi)$ forms a Lie subalgebra
of ${\msr G}$, which is called an {\it Engel
subalgebra}.\index{Engel subalgebra}

For a subalgebra ${\msr K}$ of ${\msr G}$, we define the {\it
normalizer}\index{normalizer} of ${\msr K}$ in ${\msr G}$ by
$$N_{\msr G}({\msr K})=\{\xi\in{\msr G}\mid [\xi,{\msr
K}]\subset{\msr K}\}.\eqno(4.2.4)$$ Then $N_{\msr G}({\msr K})$ is
the unique maximal subalgebra of ${\msr G}$ that contains ${\msr K}$
as an ideal.\psp

{\bf Lemma 4.2.1}. {\it If ${\msr K}$ is a subalgebra of ${\msr G}$
containing an Engel subalgebra. Then $N_{\msr G}({\msr K})={\msr
K}$.}

{\it Proof}. Suppose ${\msr K}\supset{\msr G}_0(\ad\xi)$ for some
$\xi\in{\msr G}$. Obviously, $\xi\in{\msr G}_0(\ad\xi)\subset{\msr
K}$. If $N_{\msr G}({\msr K})/{\msr K}\neq\{0\}$, then $\ad\xi$ acts
on it as an invertible linear transformation because ${\msr
G}_0(\ad\xi)\subset{\msr K}$. This contradicts $[\xi,N_{\msr
G}({\msr K})]\subset{\msr K}.\qquad\Box$\psp

We call  a subalgebra ${\msr K}$ of ${\msr G}$ {\it
self-normalizing}\index{self-normalizing} if $N_{\msr G}({\msr
K})={\msr K}$.\psp

{\bf Lemma 4.2.2}. {\it Let ${\msr K}$ be a subalgebra of ${\msr
G}$. Take $\upsilon\in {\msr K}$ such that ${\msr G}_0(\ad\upsilon)$
is a minimal element in $\{{\msr G}_0(\ad\xi)\mid \xi\in{\msr K}\}$
with respect to inclusion. If ${\msr K}\subset {\msr
G}_0(\ad\upsilon)$, then $ {\msr G}_0(\ad\upsilon)\subset {\msr
G}_0(\ad\xi)$ for any $\xi\in{\msr K}$.}

{\it Proof}. Let $0\neq \xi\in{\msr K}$ be a fixed element. Suppose
$$\dim{\msr G}=n,\qquad\dim  {\msr
G}_0(\ad\upsilon)=m.\eqno(4.2.5)$$ For any $c\in\mbb{F}$, we have
$$\ad(\upsilon+c\xi)({\msr
G}_0(\ad\upsilon))\subset {\msr G}_0(\ad\upsilon).\eqno(4.2.6)$$ Let
$$f(t,c)=t^m+f_1(c)t^{m-1}+\cdots +f_m(c)\eqno(4.2.7)$$
be the characteristic polynomial of $\ad(\upsilon+c\xi)|_{{\msr
G}_0(\ad\upsilon)}$ and let
$$g(t,c)=t^{n-m}+g_1(c)t^{n-m-1}+\cdots +g_{n-m}(c)\eqno(4.2.8)$$
be the characteristic polynomial of $\ad(\upsilon+c\xi)|_{{\msr G}/
{\msr G}_0(\ad\upsilon)}$. Then the characteristic polynomial of
$\ad(\upsilon+c\xi)$ on ${\msr G}$ is $f(t,c)g(t,c)$. Moreover, as
polynomials of $c$,
$$\mbox{deg}\:f_i(c)\leq i,\qquad \mbox{deg}\:g_j(c)\leq j.\eqno(4.2.9)$$

Since $\ad \upsilon|_{{\msr G}/ {\msr G}_0(\ad\upsilon)}$ is
invertible, we have $g_{n-m}(0)\neq 0$. So $g_{n-m}(c)$ is a nonzero
polynomial.  There exists distinct $c_1,c_2,...,c_{n+1}\in \mbb{C}$
such that
$$g_{n-m}(c_i)\neq 0\qquad\for\;\;i\in\ol{1,n+1}.\eqno(4.2.10)$$
Thus
$${\msr G}_0(\ad(\upsilon+c_i\xi))\subset {\msr G}_0(\ad\upsilon).\eqno(4.2.11)$$
By the minimality of ${\msr G}_0(\ad\upsilon)$, we have
$${\msr G}_0(\ad(\upsilon+c_i\xi))={\msr
G}_0(\ad\upsilon)\qquad\for\;\;i\in\ol{1,n+1}, \eqno(4.2.12)$$ or
equivalently,
$$f(t,c_i)=t^m\qquad\for\;\;i\in\ol{1,n+1}. \eqno(4.2.13)$$
By (4.2.7),
$$f_j(c_i)=0\qquad\for\;\;i\in\ol{1,n+1},\;j\in\ol{1,m}.\eqno(4.2.14)$$
Then (4.2.9) and (4.2.14) imply that
$$f_1(c)\equiv \cdots\equiv f_m(c)\equiv 0\;\;\mbox{as polynomials
of}\;c.\eqno(4.2.15)$$ So
$${\msr G}_0(\ad\upsilon)\subset{\msr
G}_0(\ad(\upsilon+c\xi))\qquad\mbox{for any}\;\xi\in{\msr
K},\;c\in\mbb{F}.\eqno(4.2.16)$$ But $\upsilon +{\msr K}={\msr K}$
because $\upsilon\in{\msr K}$. Therefore,  $ {\msr
G}_0(\ad\upsilon)\subset {\msr G}_0(\ad\eta)$ for any $\eta\in{\msr
K}.\qquad\Box$\psp

A {\it Cartan subalgebra}\index{Cartan subalgebra} (abbreviated CSA)
of a Lie algebra ${\msr G}$ is a self-normalizing nilpotent
subalgebra.\psp

{\bf Theorem 4.2.3}. {\it A subalgebra $H$ of ${\msr G}$ is a CSA if
and only if it is a minimal Engel subalgebra of ${\msr G}$.}

{\it Proof}. Suppose that $H={\msr G}_0(\ad\upsilon)$ is a minimal
Engel subalgebra. It is self-normalizing by Lemma 4.2.1. Taking
${\msr K}=H$ in the above lemma, we have
$$H\subset{\msr G}_0(\ad\xi)\qquad\mbox{for any}\;\;\xi\in
H,\eqno(4.2.17)$$ or equivalently, $\ad\xi|_H$ is nilpotent for any
$\xi\in H$. By Engel's Theorem, $H$ is nilpotent.

Conversely, we assume that $H$ is a CSA. Let ${\msr
G}_0(\ad\upsilon)$ be a minimal element in $\{{\msr
G}_0(\ad\xi)\mid\xi\in H\}$. Since $H$ is nilpotent, we have
$$H\subset {\msr G}_0(\ad\upsilon).\eqno(4.2.18)$$
By the above lemma with ${\msr K}=H$, $\ad\xi|_{{\msr
G}_0(\ad\upsilon)}$ is nilpotent and so is $\ad\xi|_{{\msr
G}_0(\ad\upsilon)/H}$ for any $\xi\in H\subset{\msr
G}_0(\ad\upsilon)$. If $H\neq{\msr G}_0(\ad\upsilon)$, then there
exists $u\in {\msr G}_0(\ad\upsilon)\setminus H$ such that
$$[H,u+H]\subset H\eqno(4.2.19)$$
by Theorem 1.5.1. So $u\in N_{\msr G}(H)\setminus H$, which
contradicts that $H$ is self-normalizing. Thus $H={\msr
G}_0(\ad\upsilon)$.

If ${\msr G}_0(\ad u)\subset H$ for some $u\in {\msr G}$, then $u\in
H$. Since $H$ is nilpotent, we have $H\subset{\msr G}_0(\ad u)$.
Thus $H={\msr G}_0(\ad\upsilon)$ is a minimal Engel
subalgebra.$\qquad\Box$\psp

{\bf Corollary 4.2.4}. {\it Let ${\msr G}$ be a finite-dimensional
semisimple Lie algebra over $\mbb{C}$. Then the CSA's of ${\msr G}$
are precisely the maximal toral subalgebras of ${\msr G}$.}

{\it Proof}. Suppose that $H $ is a maximal toral subalgebra of
${\msr G}$. Then we have the Cartan decomposition
$${\msr G}=H+\sum_{\al\in\Phi}{\msr G}_\al,\eqno(4.2.20)$$
which implies that $H$ is self-normalizing. So $H$ is a CSA.

Conversely, we assume that $H$ is a CSA of ${\msr G}$. By the above
theorem, $H={\msr G}_0(\ad\upsilon)$ is a minimal Engel subalgebra.
Let $\upsilon=\upsilon_s+\upsilon_n$ be the abstract Jordan
decomposition. Then
$$C_{\msr G}(\upsilon_s)=\{u\in{\msr G}\mid
[\upsilon_s,u]=0\}={\msr G}_0(\ad\upsilon)=H.\eqno(4.2.21)$$ On the
other hand, there exists a maximal toral subalgebra $H'$ containing
$\upsilon_s$. But $H'\subset C_{\msr G}(\upsilon_s)=H$ is a CSA and
it is also a minimal Engel suablgebra of ${\msr G}$ by the above
theorem. By the minimality of ${\msr G}_0(\ad\upsilon)=H$, we have
$H=H'.\qquad\Box$\psp

In the above, we have $\upsilon=\upsilon_s$, which is called a {\it
regular semisimple} element.\index{regular semisimple}\psp

{\bf Lemma 4.2.5}. {\it Let $\phi: {\msr G}\rta {\msr G}'$ be an
epimorphism of Lie algebras. If $H$ is a CSA of ${\msr G}$, then
$\phi(H)$ is a CSA of ${\msr G}'$. Conversely, if $H'$ is a CSA of
${\msr G}'$, then any CSA of $\phi^{-1}(H')$ is a CSA of ${\msr G}$.
}

{\it Proof}. Set $\kn\phi={\msr K}$. Then ${\msr G}'\cong {\msr
G}/{\msr K}$. Obviously, $\phi(H)$ is a nilpotent subalgebra of
${\msr G}'$ because $H$ is nilpotent. Let $\xi'\in N_{{\msr
G}'}(\phi(H))$. Take an element $\xi\in\phi^{-1}(\xi')$. Then
$$\phi([\xi,H+{\msr K}])=[\xi',\phi(H)]\subset
\phi(H).\eqno(4.2.22)$$ Thus
$$\xi\in N_{\msr G}(H+{\msr K}).\eqno(4.2.23)$$
Since $H$ is an Engel subalgebra by Theorem 4.2.3, we have
$$N_{\msr G}(H+{\msr K})=H+{\msr K}\eqno(4.2.24)$$
by Lemma 4.2.1. So $\xi'=\phi(\xi)\in \phi(H+{\msr K})=\phi(H)$;
that is, $\phi(H)$ is self-normalizing. Thus $\phi(H)$ is a CSA of
${\msr G}'$.

Suppose that $H'$ is a CSA of ${\msr G}'$ and $H$ is a CSA of
$\phi^{-1}(H')$. By the above argument, $\phi(H)$ is a CSA of $H'$.
So $\phi(H)$ is an Engel subalgebra of $H'$. Since $H'$ is
nilpotent, we have $\phi(H)=H'$. If $\xi\in {\msr G}$ such that
$[\xi,H]\subset H$, then $[\phi(\xi),H']\subset H'$. Thus
$\phi(\xi)\in H'$ because $H'$ is self-normalizing. Hence $\xi\in
N_{\phi^{-1}(H')}(H)=H$; that is, $H=N_{\msr G}(H).$ Therefore, $H$
is a CSA of ${\msr G}$ because it is nilpotent as a CSA of
$\phi^{-1}(H').\qquad\Box$\psp

Exercises:

1. Let ${\msr G}$ be a finite-dimensional semisimple Lie algebra
over $\mbb{C}$. Prove that a nonzero semsimple element
$\upsilon\in{\msr G}$ is regular if and only if $\upsilon$ lies in
exactly one CSA of ${\msr G}$.

2. Prove that a CSA is a maximal nilpotent subalgebra, but a maximal
nilpotent subalgebra may not be a CSA.

\section{Conjugacy Theorems}

In this section, we always assume that ${\msr G}$ is a
finite-dimensional Lie algebra over $\mbb{C}$. We want to prove that
all CSA's of ${\msr G}$ are conjugated under the group
$\mbox{Int}\:{\msr G}$ of inner automorphisms of ${\msr G}$.

An element $u\in{\msr G}$ is called {\it strongly nilpotent} if
$u\in{\msr G}_a(\ad\xi)$ for some $\xi\in{\msr G}$ and $0\neq
a\in\mbb{C}$. Set
$${\msr N}({\msr G})=\mbox{the set of all strongly nilpotent elements
in}\;{\msr G}.\eqno(4.3.1)$$ Then ${\msr N}({\msr G})$ is invariant
under $\mbox{Aut}\:{\msr G}$. Denote
$${\msr E}({\msr G})=\mbox{the subgroup of}\;\mbox{Aut}\:{\msr
G}\;\mbox{generated by}\;\{e^{\ad u}\mid u\in {\msr N}({\msr
G})\}.\eqno(4.3.2)$$ For any $\tau\in\mbox{Aut}\:{\msr G}$ and $u\in
{\msr N}({\msr G})$, we have $\tau(u)\in {\msr N}({\msr G})$ and
$$\tau e^{\ad u}\tau^{-1}=e^{\ad \tau(u)}\in {\msr E}({\msr
G}).\eqno(4.3.3)$$ Thus ${\msr E}({\msr G})$ is a normal subgroup of
$\mbox{Aut}\:{\msr G}$.

If ${\msr K}$ is a subalgebra of ${\msr G}$, then
$${\msr N}({\msr K})\subset {\msr N}({\msr G}).\eqno(4.3.4)$$
Set
$${\msr E}({\msr G};{\msr K})=\mbox{the subgroup of}\;\mbox{Aut}\:{\msr
G}\;\mbox{generated by}\;\{e^{\ad_{_{\msr G}} u}\mid u\in {\msr
N}({\msr K})\},\eqno(4.3.5)$$ which is a subgroup of ${\msr E}({\msr
G})$.

Let $\phi:{\msr G}\rta{\msr G}'$ be a Lie algebra epimorphism. For
any $\xi\in{\msr G}$, we have $\phi({\msr G}_a(\ad\xi))\subset {\msr
G}'_a(\ad\phi(\xi))$. Since ${\msr G}=\sum_{b\in\mbb{C}}{\msr
G}_b(\ad\xi)$ and $\phi$ is surjective, we have
$$\phi({\msr G}_a(\ad\xi))={\msr
G}'_a(\ad\phi(\xi))\qquad\for\;\;a\in\mbb{C}.\eqno(4.3.6)$$ Thus
$$\phi({\msr N}({\msr G}))={\msr N}({\msr G}').\eqno(4.3.7)$$
For any $u\in{\msr N}({\msr G})$,
$$e^{\ad\phi(u)}(\phi(v))=\phi(e^{\ad u}(v)).\eqno(4.3.8)$$
Hence for any $\tau\in{\msr E}({\msr G})$, there exists a unique
$\tau'\in {\msr E}({\msr G}')$ such that
$$\tau'(\phi(v))=\phi(\tau(v))\qquad\for\;\;v\in{\msr
G}.\eqno(4.3.9)$$ Thus we have:\psp

{\bf Lemma 4.3.1}. {\it  The map $\tau\mapsto \tau'$ is a group
epimorphism from ${\msr E}({\msr G})$ to ${\msr E}({\msr G}')$}.\psp

{\bf Theorem 4.3.2}. {\it If ${\msr G}$ is solvable, then any two
CSA's are conjugated under ${\msr E}({\msr G})$.}

{\it Proof}. Let $H$ and $H'$ be two CSA's of ${\msr G}$. We prove
by induction on $\dim{\msr G}$. If ${\msr G}$ is nilpotent, then
$H=H'={\msr G}$ (This includes $\dim{\msr G}=1$. So the theorem
holds). Assume that the theorem holds for $\dim{\msr G}<k$ with
$1<k\in\mbb{Z}$. Consider $\dim{\msr G}=k$ and ${\msr G}$ is not
nilpotent. By Lie's theorem, there exists a common eigenvector
$v\in{\msr G}$ for $\ad{\msr G}$. So $\mbb{C}v$ is an ideal of
${\msr G}$. By Lemma 4.2.5, $H+\mbb{C}v$ and $H'+\mbb{C}v$ are two
CSA's of ${\msr G}/\mbb{C}v$. According to the above lemma, there
exists $\tau\in{\msr E}({\msr G})$ such that
$$\tau(H')+\mbb{C}v=H+\mbb{C}v.\eqno(4.3.10)$$

Now $\tau(H')$ and $H$ are two CSA's of the Lie subalgebra ${\msr
K}=H+\mbb{C}v$. If ${\msr K}\neq {\msr G}$, then $\dim{\msr K}<k$.
By assumption, there exists $\sgm\in{\msr E}({\msr G};{\msr K})$
such that $\sgm\tau(H')=H$. We are done.

Suppose ${\msr K}={\msr G}$. So $H+\mbb{C}v={\msr G}=H'+\mbb{C}v$.
By Theorem 4.2.3, $H={\msr G}_0(\ad\xi)$, which implies $[\xi,v]=av$
with $0\neq a\in\mbb{C}$. Thus $v\in{\msr N}({\msr G})$. Write
$\xi=\zeta+bv$ with $\zeta\in H'$. Then
$$\zeta=\xi-bv=e^{\ad a^{-1}bv}(\xi).\eqno(4.3.11)$$
Hence
$$e^{\ad a^{-1}bv}(H)={\msr G}_0(\ad\zeta)\supset H'.\eqno(4.3.12)$$
Since $e^{\ad a^{-1}bv}(H)$ and $H'$ are both CSA's, they are
minimal Engel subalgebras. Thus $e^{\ad
a^{-1}bv}(H)=H'.\qquad\Box$\psp

A maximal solvable subalgebra of ${\msr G}$ is called a {\it Borel
subalgebra}.\index{Borel subalgebra}

 \psp

{\bf Lemma 4.3.3}. {\it A Borel subalgebra ${\msr B}$ is
self-normalizing}.

{\it Proof}. For any $u\in N_{\msr G}({\msr B})$, ${\msr
B}+\mbb{C}u$ is also a solvable subalgebra. By maximality,
$u\in{\msr B}.\qquad\Box$\psp

Recall that $\rad\: {\msr G}$ is the unique maximal solvable ideal
of ${\msr G}$. For any solvable subalgebra ${\msr K}$, ${\msr
K}+\rad\:{\msr G}$ is still a solvable subalgebra. Thus any Borel
subalgebra contains $\rad\:{\msr G}$. By Lemma 4.3.1, we have:\psp

{\bf Lemma 4.3.4}. {\it If  all the Borel subalgebras of ${\msr
G}/\rad\;{\msr G}$ are conjugated under ${\msr E}({\msr
G}/\rad\;{\msr G})$, then all the Borel subalgebras of ${\msr G}$
are conjugated under ${\msr E}({\msr G})$.}\psp

Note that ${\msr G}/\rad\:{\msr G}$ is semisimple.\psp

{\bf Theorem 4.3.5}. {\it All the Borel subalgebras of ${\msr G}$
are conjugated under ${\msr E}({\msr G})$.}

{\it Proof}.  Set
$$m_0({\msr G})=\max\{\dim{\msr K}\mid {\msr K}\;\mbox{is a Borel
subalgebra of}\; {\msr G}\}.\eqno(4.3.13)$$ Let ${\msr B}$ and
${\msr B}'$ be two Borel subalgebras of ${\msr G}$. Define
$$m({\msr B},{\msr B}')=m_0({\msr G})-\dim {\msr B}\bigcap{\msr
B}'.\eqno(4.3.14)$$ We prove the theorem by induction on $\dim{\msr
G}$ and on $m({\msr B},{\msr B}')$. The theorem holds if $\dim{\msr
G}=1$ or $m({\msr B},{\msr B}')=0$, which implies ${\msr B}={\msr
B}'$. Suppose that the theorem holds for $\dim{\msr G}<k_1$ or
$m({\msr B},{\msr B})<k_2$. Now we assume $\dim{\msr G}=k_1$ and
${\msr B}\neq{\msr B}'$ with $m({\msr B},{\msr B}')=k_2$. If $\rad\;
{\msr G}\neq \{0\}$, then the theorem holds by Lemma 4.3.4 because
any two Borel subalgebras of ${\msr G}/\rad\:{\msr G}$ are
conjugated under ${\msr E}({\msr G}/\rad\:{\msr G})$ due to
$\dim{\msr G}/\rad\:{\msr G}<k$.

Now we assume that ${\msr G}$ is semsimple and $H$ is a maximal
toral subalgebra of ${\msr G}$ with root system $\Phi$. Let $\Pi$ be
a base of $\Phi$. With respect to $\Pi$, we have positive roots and
negative roots. Set
$${\msr B}(\Pi)=H+\sum_{\al\in\Phi^+}{\msr G}_\al,\;\;
{\msr N}(\Pi)=\sum_{\al\in\Phi^+}{\msr G}_\al.\eqno(4.3.15)$$
Suppose $\Pi=\{\al_1,\al_2,...,\al_n\}$. For any
$\be=\sum_{i=1}^nb_i\al_i\in \Phi$, we have
$$\mbox{ht}\:\be=\sum_{i=1}^nb_i.\eqno(4.3.16)$$
Let
$$\ell=\max\{\mbox{ht}\:\be\mid\be\in\Phi\}.\eqno(4.3.17)$$Set
$${\msr N}_i(\Pi)=\sum_{\be\in\Phi^+;\mbox{\tiny ht}\:\be\geq
i}{\msr G}_\be.\eqno(4.3.18)$$ Then we have
$$[{\msr B}(\Pi),{\msr B}(\Pi)]={\msr N}(\Pi),\;\;({\msr
N}(\Pi))^i\subset {\msr N}_i(\Pi),\eqno(4.3.19)$$ where $({\msr
N}(\Pi))^i=[{\msr N}(\Pi),({\msr N}(\Pi))^{i-1}]$ is the central
series. In particular, $({\msr N}(\Pi))^{\ell+1}=\{0\}$. So ${\msr
N}(\Pi)$ is nilpotent and ${\msr B}(\Pi)$ is solvable. If ${\msr K}$
is a subalgebra strictly containing ${\msr B}(\Pi)$, then
$${\msr K}\supset {\msr G}_\al+\mbb{C}h_\al+{\msr G}_{-\al}\cong
sl(2,\mbb{C})\qquad\mbox{for some}\;\;\al\in\Phi.\eqno(4.3.20)$$
Thus ${\msr K}$ is not solvable. Hence ${\msr B}(\Pi)$ is a Borel
subalgebra, which is called {\it standard} relative to $H$. To prove
that any two Borel subalgebras are conjugated under ${\msr E}({\msr
G})$, it is enough to prove that any Borel subalgebra ${\msr B}'$ is
conjugated to ${\msr B}(\Pi)$ under ${\msr E}({\msr G})$. So we may
assume ${\msr B}={\msr B}(\Pi)$ and simply denote ${\msr N}={\msr
N}(\Pi)$. Note that
$${\msr N}=\mbox{the set of all}\;\mbox{ad}_{\msr
G}\mbox{-nilpotent elements in}\;{\msr B}.\eqno(4.3.21)$$

If $u$ is an element of a Borel subalgebra $\ol{\msr B}$ and
$u=u_s+u_n$ is the abstract Jordan decomposition, then
$[u_s,\ol{\msr B}]\subset\ol{\msr B}$ because $\ad u_s$ is a
polynomial of $\ad u$ without constant. Hence $\mbb{C}u_s+\ol{\msr
B}$ is a solvable subalgebra. By maximality of $\bar{\msr B}$,
$u_s\in\ol{\msr B}$ and also $u_n=u-u_s\in\ol{\msr B}$. Thus ${\msr
B}\bigcap{\msr B}'$ contains the semsimple and nilpotent parts of
its elements. Moreover,
$$N'={\msr N}\bigcap{\msr B}'=\mbox{the set of all}\;\mbox{ad}_{\msr
G}\mbox{-nilpotent elements in}\;{\msr B}\bigcap{\msr B}'
\eqno(4.3.22)$$ by (4.3.21). According to (4.3.19),
$$N'\supset [{\msr B}\bigcap{\msr B}',{\msr B}\bigcap{\msr
B}'].\eqno(4.3.23)$$ So $N'$ is an ideal of ${\msr B}\bigcap{\msr
B}'$.\pse

{\bf Case 1}. $N'\neq \{0\}$. \psp

Set
$${\msr K}=N_{\msr G}(N').\eqno(4.3.24)$$
Then ${\msr B}\bigcap{\msr B}'\subset {\msr K}$ because $N'$ is an
ideal of ${\msr B}\bigcap{\msr B}'$.
  Note that ${\msr B}\neq{\msr B}'$ implies ${\msr B}/{\msr
B}\bigcap{\msr B}'\neq\{0\}$ and ${\msr B}'/{\msr B}\bigcap{\msr
B}'\neq\{0\}$. Since $\ad_{\msr G}N'$ consists of nilpotent
elements,  there exist $u\in {\msr B}\setminus{\msr B}\bigcap{\msr
B}'$ and $u'\in{\msr B}'\setminus{\msr B}\bigcap{\msr B}'$ such that
$$[u,N'],[u',N']\subset {\msr B}\bigcap {\msr
B}'.\eqno(4.3.25)$$ By Lie's Theorem, the elements $\ad_{\msr
G}[{\msr B}',{\msr B}']$ are nilpotent. So
$$[u,N'],[u',N']\subset N'\eqno(4.3.26)$$
by (4.3.22). Thus
$${\msr K}\bigcap {\msr B}\setminus{\msr B}\bigcap {\msr B}'\neq
\emptyset,\;\;{\msr K}\bigcap {\msr B}'\setminus{\msr B}\bigcap
{\msr B}'\neq \emptyset.\eqno(4.3.27)$$ Let ${\msr C}\supset {\msr
K}\bigcap{\msr B}$ and ${\msr C}'\supset {\msr K}\bigcap{\msr B}'$
be  Borel subalgebras of ${\msr K}$. Since ${\msr G}$ is semisimple,
$N'$ can not be an ideal of ${\msr G}$ because it is nilpotent. So
$\dim{\msr K}<k$. By assumption, there exists $\nu\in{\msr E}({\msr
G};{\msr K})$ such that $\nu({\msr C})={\msr C}'$.

Take Borel subalgebras ${\msr C}_1\supset{\msr C}$ and ${\msr
C}_2\supset{\msr C}'$ of ${\msr G}$. Then
$${\msr C}_1\bigcap{\msr B}\supset{\msr K}\bigcap{\msr
B},\;\nu({\msr C}_1)\bigcap{\msr C}_2\supset {\msr C}'\supset {\msr
K}\bigcap{\msr B}',\;{\msr C}_2\bigcap{\msr B}'\supset{\msr
K}\bigcap{\msr B}'.\eqno(4.3.28)$$ Thus
$$m({\msr B},{\msr C}_1),\;m(\nu({\msr C}_1),{\msr C}_2),\;m({\msr
C}_2,{\msr B}')<m_1.\eqno(4.3.29)$$ By assumption, there exist
$\sgm_1,\sgm_2,\sgm_3\in{\msr E}({\msr G})$ such that
$$\sgm_1({\msr B})={\msr C}_1,\;\;\sgm_2\nu({\msr C}_1)={\msr
C}_2,\;\;\sgm_3({\msr C}_2)={\msr B}'.\eqno(4.3.30)$$ Therefore,
$$\sgm_3\sgm_2\nu\sgm_1({\msr B})={\msr B}'.\eqno(4.3.31)$$
\vspace{0.1cm}

{\bf Case 2}. $N'=\{0\}$ and ${\msr B}\bigcap{\msr B}'\neq\{0\}$.
\pse

In this case, $T={\msr B}\bigcap{\msr B}'$ is a subalgebra
consisting  of semismiple elements. So $T$ is abelian by Lemma
2.4.1. Let ${\msr C}$ be a CSA of
$$C_{\msr B}(T)=\{u\in{\msr B}\mid [u,T]=\{0\}\}.\eqno(4.3.32)$$
Since ${\msr C}$ is self-normalizing and $[T,{\msr C}]=\{0\}$, we
have $T\subset {\msr C}$. For any $v\in N_{\msr B}({\msr C})$ and
$h\in T$, we have $[h,v]\in{\msr C}$ and $(\ad h)^2(v)=0$. The
semisimplicity of $h$ implies $[h,v]=0$. Hence $v\in C_{\msr B}(T)$.
So $v\in N_{C_{\msr B}(T)}({\msr C})={\msr C}$. Thus ${\msr C}$ is a
self-normalizing nilpotent subalgebra of ${\msr B}$; that is, a CSA
of ${\msr B}$. On the other hand, $H$ is a CSA of ${\msr B}$. By
Theorem 4.3.2, there exists $\tau\in{\msr E}({\msr G};{\msr B})$
such that $\tau({\msr C})=H$. Note $\tau({\msr B})={\msr B}$. So
$$\tau({\msr B}')\bigcap {\msr B}=\tau({\msr B}')\bigcap\tau({\msr
B})=\tau({\msr B}'\bigcap{\msr B})=\tau(T)\subset\tau({\msr
C})=H.\eqno(4.3.33)$$
 Replacing
${\msr B}'$ by $\tau({\msr B}')$, we may assume $T\subset H$.

Suppose that $T\neq H$. Note $H\subset C_{\msr G}(T)=\{u\in{\msr
G}\mid [u,T]=\{0\}\}$. If ${\msr B}'\subset  C_{\msr G}(T)$, we take
a Borel subalgebra $\ol{\msr B}\supset H$ of $C_{\msr G}(T)$. Since
$Z({\msr G})=\{0\}$, we have $\dim C_{\msr G}(T)<k$. By assumption,
there exists $\sgm\in{\msr E}({\msr G};C_{\msr G}(T))$ such that
$\sgm({\msr B}')=\ol{\msr B}$. Replacing ${\msr B}'$ by $\sgm({\msr
B}')$, we may assume $T=H$. When ${\msr B}'\not\subset C_{\msr
G}(T)$, there exist $u\in {\msr B}'$ and $h\in T$ such that
$[h,u]=u$. Set
$${\msr C}=H+\sum_{\al\in\Phi;\;0<\al(h)\in\mbb{Z}}{\msr
G}_{\al}.\eqno(4.3.34)$$ Then ${\msr C}$ is a solvable subalgebra
and $u\in {\msr C}$. Take  a Borel subalgebra ${\msr C}_1\supset
{\msr C}$ of ${\msr G}$. We have $\dim{\msr C}_1\bigcap{\msr
B}'>\dim T=\dim{\msr B}\bigcap {\msr B}'$. So $m({\msr B}',{\msr
C}_1)<m_1$. By assumption, there exists $\varrho\in {\msr E}({\msr
G})$ such that $\varrho({\msr B}')={\msr C}_1$. Replacing ${\msr
B}'$ by $\varrho({\msr B}')$, we may again assume $T=H$.

Now we only need to deal with the case $T=H$. Since ${\msr
B}'\neq{\msr B}$, we have ${\msr G}_{-\al}\subset{\msr B}'$ for some
$\al\in\Phi^+$. Take $\xi\in{\msr G}_{\al}$ and $\zeta\in{\msr
G}_{-\al}$ such that
$$[\xi,\zeta]=h_\al,\;\;[h_\al,\xi]=2\xi.\eqno(4.3.35)$$
Denote
$$\tau_\al=e^{\ad\xi}e^{-\ad\zeta}e^{\ad\xi}\in{\msr E}({\msr
G}).\eqno(4.3.36)$$ Then
$$\tau_\al({\msr G}_{-\al})={\msr
G}_{\al},\;\;\tau_\al(H)=H.\eqno(4.3.37)$$ Now $m({\msr B},
\tau_\al({\msr B}'))<m_1$. By assumption, there exists $\sgm\in{\msr
E}({\msr G})$ such that $\sgm\tau_\al({\msr B}')={\msr B}$.\psp

{\bf Case 3}. ${\msr B}\bigcap {\msr B}'=\{0\}$.\psp

Let $T$ be a maximal toral subalgebra of ${\msr B}'$. If $T=\{0\}$,
then ${\msr B}'$ consists of ad-nilpotent elements because it
contains the semisimple parts of all its elements. By Engel's
Theorem and Lemma 4.3.3, ${\msr B}'$ is a CSA. But Corollary 4.2.4
says that ${\msr B}'$ is a toral subalgebra. Thus ${\msr B}'=\{0\}$,
which is absurd. Thus $T\neq\{0\}$. Take $\bar{H}\supset T$ to be a
maximal toral subalgebra of ${\msr G}$. Take any standard Borel
subalgebra $\ol{\msr B}$ of ${\msr G}$ relative to $\bar{H}$. Then
${\msr B}'\bigcap \ol{\bf B}\neq\{0\}$. By Case 2, ${\msr B}'$ is
conjugated to $\ol{\msr B}$. So
$$\dim{\msr B}'=\dim\ol{\msr B}=\frac{1}{2}(\dim{\msr G}+\dim
\bar{H})>\frac{1}{2}\dim{\msr G}.\eqno(4.3.38)$$ On the other hand,
$$\dim{\msr B}=\frac{1}{2}(\dim{\msr G}+\dim
H)>\frac{1}{2}\dim{\msr G}.\eqno(4.3.39)$$ Thus
$$\dim{\msr G}\geq \dim{\msr B}+\dim{\msr B}'>\frac{1}{2}\dim{\msr
G}+\frac{1}{2}\dim{\msr G}=\dim{\msr G},\eqno(4.3.40)$$ which is
absurd. Therefore, this case does not exist. $\qquad\Box$\psp

{\bf Corollary 4.3.6}. {\it Any Borel subalgebra ${\msr B}$ of a
semisimple Lie algebra ${\msr G}$ is standard}.

{\it Proof}. Take  a maximal toral subalgebra $H$ of ${\msr G}$ and
let $\Pi$ be a base of the corresponding root system $\Phi$. Then
there exists a $\sgm\in{\msr E}({\msr G})$ such that $\sgm({\msr
B}(\Pi))={\msr B}$ by the above theorem (cf. (4.3.15)). Note that
$\sgm(H)\subset {\msr B}$ is also a maximal toral subalgebra of
${\msr G}$. Replacing $H$ by $\sgm(H)$, we may assume ${\msr
B}\supset H$. Again there exists $\sgm_1\in{\msr E}({\msr G})$ such
that $\sgm_1({\msr B}(\Pi))={\msr B}$. Now $H$ and $\sgm_1(H)$ are
CSA's of ${\msr B}$. By Lemma 4.3.2, there exists $\tau\in {\msr
E}({\msr G};{\msr B})$ such that $\tau\sgm_1(H)=H$.

 Since $\tau({\msr
B})={\msr B}$, we have $\tau\sgm_1({\msr B}(\Pi))={\msr B}$. Denote
$\nu=\tau\sgm_1$. So we have $\nu(H)=H$ and $\nu({\msr
B}(\Pi))={\msr B}$. Define the action of $\nu$ on $H^\ast$ by
$$\nu(\al)(h)=\al(\nu^{-1}(h))\qquad\for\;\;\al\in
H^\ast.\eqno(4.3.41)$$ For any $\al\in\Phi$ and $u\in{\msr G}_\al$,
we have
$$[h,\nu(u)]=\nu([\nu^{-1}(h),u])=\al(\nu^{-1}(h))\nu(u)\qquad\for\;\;h\in
H.\eqno(4.3.42)$$ Thus
$$\nu({\msr G}_\al)={\msr
G}_{\nu(\al)}\qquad\for\;\;\al\in\Phi.\eqno(4.3.43)$$ Hence $\nu
(\Phi)=\Phi$. Since $\nu$ is linear, $\nu(\Pi)$ is also a base of
$\Phi$. Thus ${\msr B}={\msr B}(\nu(\Pi))$ is a standard Borel
subalgebra.$\qquad\Box$\psp

{\bf Corollary 4.3.7}. {\it Any two CSA's of ${\msr G}$ are
conjugated under ${\msr E}({\msr G})$}.

{\it Proof}. Let $H$ and $H'$ be two CSA's of ${\msr G}$. Take Borel
subalgebras ${\msr B}\supset H$ and ${\msr B}'\supset H'$. By
Theorem 4.3.5, there exists $\tau\in{\msr E}({\msr G})$ such that
$\tau({\msr B}')={\msr B}$. Now $H$ and $\tau(H')$ are two CSA's in
the solvable Lie algebra ${\msr B}$. According to Lemma 4.3.2, there
exists $\sgm\in{\msr E}({\msr G};{\msr B})$ such that
$\sgm\tau(H')=H.\qquad\Box$.\psp

Let $H$ and $H'$ be two maximal toral subalgebras (which are CSA's
by Corollary 4.2.4) of a finite-dimensional semisimple Lie algebra
${\msr G}$. Their associated root systems are denoted by $\Phi$ and
$\Phi'$, respectively. For $\al\in \Phi$ and $\al'\in \Phi'$, we
denote
$${\msr G}_\al=\{u\in {\msr G}\mid [h,u]=\al(h)u\;\for\;h\in
H\},\eqno(4.3.44)$$
$$ {\msr G}'_{\al'}=\{u\in {\msr G}\mid
[h',u]=\al'(h')u\;\for\;h'\in H'\}\eqno(4.3.45)$$ and take $h_\al\in
[{\msr G}_\al,{\msr G}_{-\al}]$ and $h'_{\al'}\in [{\msr
G}'_{\al'},{\msr G}'_{-\al'}]$ such that
$$\ad h_\al|_{{\msr G}_\al}=2\mbox{Id}_{{\msr G}_\al},\;\;
\ad h'_{\al'}|_{{\msr G}'_{\al'}}=2\mbox{Id}_{{\msr
G}'_{\al'}}.\eqno(4.3.46)$$ Then
$$\be(h_\al)=\la \be,\al\ra,\;\;\be'(h'_{\al'})=\la \be',\al'\ra\;\;\for\;\;\al,\be\in \Phi,\;\al',\be'\in \Phi'.
\eqno(4.3.47)$$

By the above corollary, there exists $\nu\in{\msr E}((\msr G))$ such
that $\nu(H)=H'$. We define $\nu: H^\ast\rta (H')^\ast$ by
$$\nu(\al)(h')=\al(\nu^{-1}(h'))\qquad\for\;\al\in H^\ast,\;h'\in
H'.\eqno(4.3.48)$$ For any $\al\in \Phi$, $v\in{\msr G}_\al$ and
$h'\in H'$, we have
\begin{eqnarray*}\hspace{2cm}
[h',\nu(v)]&=&[\nu(\nu^{-1}(h')),\nu(v)]=\nu([\nu^{-1}(h'),v])\\
&=&\al(\nu^{-1}(h'))\nu(v)=\nu(\al)(h')\nu(v).
\hspace{4.5cm}(4.3.49)\end{eqnarray*} Thus
$$\nu({\msr G}_\al)={\msr
G}'_{\nu(\al)}\qquad\for\;\;\al\in\Phi.\eqno(4.3.50)$$ In
particular, $\nu(\Phi)=\Phi'.$ Since $\nu$ is an automorphism,
$$\nu(h_\al)=h'_{\nu(\al)}\qquad\for\;\;\al\in\Phi.\eqno(4.3.51)$$
Now for $\al,\be\in\Phi$ and $0\neq u\in {\msr G}_\be$, we have
$$\la\be,\al\ra\nu(u)=\nu([h_\al,u])=[\nu(h_\al),\nu(u)]=\nu(\be)(h'_{\nu(\al)})\nu(u)=\la\nu(\be),\nu(\al)\ra\nu(u).
\eqno(4.3.52)$$ Therefore,
$$\la\nu(\be),\nu(\al)\ra=\la\be,\al\ra\qquad\for\;\;\al,\be\in\Phi;\eqno(4.3.53)$$
that is, $\nu:\Phi\rta\Phi'$ is an isomorphism of two root systems.
This shows that the root system of ${\msr G}$ is independent of the
choice of its maximal toral subalgebra.

The dimension of a CSA of a finite-dimensional semisimple Lie
algebra ${\msr G}$ is called the {\it rank} of ${\msr G}$. For
instance, $sl(n,\mbb{C})$ for $n>1$ is a finite-dimensional simple
Lie algebra of rank $n-1$. While the ranks of
$so(2n,\mbb{C}),\;so(2n+1,\mbb{C})$ and $sp(2n,\mbb{C})$ are $n$.
The Lie algebra $sl(2,\mbb{C})$ is the unique rank-1
finite-dimensional simple Lie algebra. There are three rank-2
finite-dimensional simple Lie algebras:
$sl(3,\mbb{C}),\;o(5,\mbb{C})\cong sp(4,\mbb{C})$ and the Lie
algebra of type $G_2$.

Finally, we want to study the structure of $\mbox{Aut}\:{\msr G}$
when ${\msr G}$ is a finite-dimensional simple Lie algebra over
$\mbb{C}$. Note that $sl(2,\mbb{C})$ has the following
representation on the polynomial algebra ${\msr
A}=\mbb{C}[x_1,x_2]$:
$$E_{1,2}|_{\msr A}=x_1\ptl_{x_2},\;\;E_{2,1}|_{\msr
A}=x_2\ptl_{x_1},\;\;(E_{1,1}-E_{2,2})|_{\msr
A}=x_1\ptl_{x_1}-x_2\ptl_{x_2}.\eqno(4.3.54)$$ For $0\neq
a\in\mbb{C}$, we define the operator
$$\sgm(a)=e^{aE_{1,2}}e^{-a^{-1}E_{2,1}}e^{aE_{1,2}}e^{E_{1,2}}e^{-E_{2,1}}e^{E_{1,2}}\eqno(4.3.55)$$
on ${\msr A}$. For any $f(x_1,x_2)\in{\msr A}$,
\begin{eqnarray*}& &e^{aE_{1,2}}e^{-a^{-1}E_{2,1}}e^{aE_{1,2}}(f(x_1,x_2))\\&=&e^{ax_1\ptl_{x_2}}e^{-a^{-1}x_2\ptl_{x_1}}
e^{ax_1\ptl_{x_2}}(f(x_1,x_2))=e^{ax_1\ptl_{x_2}}e^{-a^{-1}x_2\ptl_{x_1}}
(f(x_1,x_2+ax_1))\\&=&e^{ax_1\ptl_{x_2}}(f(x_1-a^{-1}x_2,ax_1))=f(-a^{-1}x_2,ax_1)\hspace{5.6cm}(4.3.56)\end{eqnarray*}
by Taylor's expansion. Thus
$$\sgm(a)(f(x_1,x_2))=f(-a^{-1}x_1,-ax_2).\eqno(4.3.57)$$
Hence
$$(E_{1,1}-E_{2,2})(x_1^{\ell_1}x_2^{\ell_2})=(\ell_1-\ell_2)x_1^{\ell_1}x_2^{\ell_2},\;\;
\sgm(a)(x_1^{\ell_1}x_2^{\ell_2})=(-a)^{\ell_2-\ell_1}x_1^{\ell_1}x_2^{\ell_2}\eqno(4.3.58)$$
for $\ell_1,\ell_2\in\mbb{N}$. It is easy to verify that the
subspace ${\msr A}_m$ of homogeneous polynomials with degree $m$ in
${\msr A}$ forms an $(m+1)$-dimensional irreducible
$sl(2,\mbb{C})$-submodule for $m\in\mbb{N}$. By (4.3.58) and Weyl's
theorem of complete reducibility (cf. Theorem 2.4.5), we have:\psp

{\bf Lemma 4.3.8}. {\it Let $u$ be an element of a
finite-dimensional $sl(2,\mbb{C})$-module such that
$(E_{1,1}-E_{2,2})(u)=ku$ with $k\in\mbb{Z}$. Then
$\sgm(a)(u)=(-a)^{-k}u$.}\psp

 Take a maximal toral subalgebra $H$ of ${\msr G}$ and let
$\Pi$ be a base of the corresponding root system $\Phi$. For any
$\tau\in\mbox{Aut}\:{\msr G}$, $\tau({\msr B}(\Pi))$ is a Borel
subalgebra of ${\msr G}$. By Theorem 4.3.5, there exists
$\sgm\in{\msr E}({\msr G})$ such that $\sgm\tau({\msr B}(\Pi))={\msr
B}(\Pi)$. Since $\sgm\tau(H)$ and $H$ are two CSA's in ${\msr B}$.
By Lemma 4.3.2, there exists $\tau_1\in{\msr E}({\msr G};{\msr
B}(\Pi))$ such that $\tau_1\sgm\tau(H)=H$. Denote
$\nu=\tau_1\sgm\tau$. By the proof of Corollary 4.3.6, (4.3.43)
holds. According to Theorem 3.2.8 (e), (4.1.24) and (4.1.25), there
exists $\sgm_1\in {\msr W}$ and $\Im_{\sgm_1} \in{\msr E}({\msr G})$
such that $\sgm_1\nu(\Pi)=\Pi$ and
$$\Im_{\sgm_1}\nu({\msr G}_\al)={\msr
G}_{\sgm_1\nu(\al)}\qquad\mbox{for}\;\;\al\in\Phi.\eqno(4.3.59)$$
This shows that
$$\tau\in \Im_\iota{\msr E}({\msr G}),\eqno(4.3.60)$$
where $\iota$ is an automorphism of $\Phi$ such that
$\iota(\Pi)=\Pi$ and $\Im_\iota\in \mbox{Aut}\;{\msr G}$ such that
$\Im_\iota(H)=H$ and $\Im_\iota({\msr G}_\al)={\msr G}_{\iota(\al)}$
for $\al\in\Phi$.

If $\vt$ is another automorphism of ${\msr G}$ such that $\vt(H)=H$
and $\vt({\msr G}_\al)={\msr G}_{\iota(\al)}$ for $\al\in\Phi$, then
$$\vt^{-1}\Im_\iota(h)=h,\;\;\vt^{-1}\Im_\iota({\msr G}_\al)={\msr
G}_{\al}\qquad\for\;\;h\in H,\;\al\in\Phi.\eqno(4.3.61)$$ Suppose
$\Pi=\{\al_1,\al_2,...,\al_n\}$. We get
$$\vt^{-1}\Im_\iota|_{{\msr G}_{\al_i}}=b_i\mbox{Id}_{{\msr
G}_{\al_i}},\;\;\vt^{-1}\Im_\iota|_{{\msr
G}_{\al_i}}=b_i^{-1}\mbox{Id}_{{\msr G}_{\al_i}},\qquad 0\neq
b_i\in\mbb{C},\eqno(4.3.62)$$ for $i\in\ol{1,n}$. By Lemma 4.1.3,
$\vt^{-1}\Im_\iota$ is completely determined by (4.3.62). Take
$\xi_1\in{\msr G}_{\al_i}$ and $\zeta_i\in{\msr G}_{-\al_i}$ such
that $[\xi_i,\zeta_i]=h_{\al_i}$ for $i\in\ol{1,n}$. Moreover, we
define
$$\sgm_i(a_i)=e^{a_i\ad\xi_i}e^{-a^{-1}\ad\zeta_i}e^{a_i\ad\xi_i}
e^{\ad\xi_i}e^{-\ad\zeta_i}e^{\ad\xi_i}\in{\msr E}({\msr G}),
\;\;a_i\in\mbb{C}.\eqno(4.3.63)$$ By Lemma 4.3.8,
$\sgm_i(a_i)|_H=\mbox{Id}_H$ and
$$\sgm_i(a_i)|_{{\msr G}_{\al_j}}=(-a_i)^{-\la\al_j,\al_i\ra}\mbox{Id}_{{\msr
G}_{\al_j}}\qquad\for\;\;i,j\in\ol{1,n}.\eqno(4.3.64)$$ Since the
Cartan matric $(\la\al_h,\al_i\ra)_{n\times n}$ is nondegenerate,
there exist $a_1,a_2,...,a_n\in\mbb{C}$ such that
$$\prod_{i=1}^n(-a_i)^{-\la\al_j,\al_i\ra}=e^{-\sum_{i=1}^n\la\al_j,\al_i\ra\ln (-a_i)}=e^{\ln b_j}=b_j\qquad\for\;\;j\in\ol{1,n}.\eqno(4.3.65)$$
In this case,
$$\vt^{-1}\Im_\iota=\sgm_1(a_1)\sgm_2(a_2)\cdots\sgm_n(a_n)\in{\msr
E}({\msr G}).\eqno(4.3.66)$$

For each automorphism $\iota$ of $\Phi$ induced by an automorphism
of the Dynkin diagram of ${\msr G}$, we fix an automorphism
$\Im_\iota$ of ${\msr G}$ such that such that $\Im_\iota(H)=H$ and
$\Im_\iota({\msr G}_\al)={\msr G}_{\iota(\al)}$ for $\al\in\Phi$.
Then
$$\tau\in\Im_\iota{\msr E}({\msr G})\eqno(4.3.67)$$
for some automorphism $\iota$ of $\Phi$ induced by an automorphism
of the Dynkin diagram of ${\msr G}$. Therefore,
$$\mbox{Aut}\:{\msr G}/{\msr E}({\msr G})\cong\mbox{the autmorphism
group of the Dynkin diagram of}\;{\msr G}.\eqno(4.3.68)$$
 In particular,
$$\mbox{Aut}\:{\msr G}={\msr E}({\msr G})\;\;\mbox{if}\;{\msr
G}\;\mbox{is of types:}\;\;B_n, C_n,E_8, E_7,F_4,
G_2.\eqno(4.3.69)$$ \vspace{0.1cm}

Let $\mbb{B}({\msr G})$ be the set of Borel subalgebras of a Lie
algebra ${\msr G}$.\psp

{\bf Lemma 4.3.9}. {\it We have $\bigcap_{{\msr B}\in\mbb{B}({\msr
G})}{\msr B}=\rad\; {\msr G}$.}

{\it Proof}. Since any Borel subalgebra contains $\rad\: {\msr G}$,
the conclusion is equivalent to $\bigcap_{\ol{\msr B}\in\mbb{
B}({\msr G}/\mbox{\tiny Rad}\:{\msr G})}\ol{\msr B}=\{0\}$. So we
only need to consider the case when ${\msr G}$ is semisimple. Let
$H$ be a maximal toral subalgebra ${\msr G}$ with root system
$\Phi$. Take a base $\Pi$ of $\Phi$. We have the following two Borel
subalgebras:
$${\msr B}_1=H+\sum_{\al\in\Phi^+}{\msr G}_\al,\;\;{\msr B}_2=H+\sum_{\al\in\Phi^-}{\msr
G}_\al.\eqno(4.3.70)$$ Note
$${\msr B}_1\bigcap{\msr B}_2=H.\eqno(4.3.71)$$
Let $\mbb{H}$ be the set of maximal toral subalgebras of ${\msr G}$.
By (4.3.71),
$$\bigcap_{{\msr B}\in\mbb{B}({\msr G})}{\msr B}\subset \bigcap_{H'\in\mbb{H}}H'.\eqno(4.3.72)$$

For any $\al\in\Phi$, take $0\neq \xi_\al\in{\msr G}_\al$. We have
$$e^{\ad\xi_\al}(h)=h-\al(h)\xi_\al\qquad\for\;\;h\in
H.\eqno(4.3.73)$$ Thus
$$e^{\ad\xi_\al}(H)\bigcap H=\{h\in H\mid
\al(h)=0\}.\eqno(4.3.74)$$ Since $\{e^{\ad\xi_\al}(H)\mid
\al\in\Phi\}$ are maximal toral subalgebras of ${\msr G}$, we have
$$\bigcap_{{\msr B}\in\mbb{B}({\msr G})}{\msr B}\subset
H\bigcap\bigcap_{\al\in\Phi}e^{\ad\xi_\al}(H)=\{h\in
H\mid\al(h)=0\;\mbox{for any}\;\al\in\Phi\}=\{0\}\eqno(4.3.75)$$ by
(4.3.72).$\qquad\Box$\psp

Exercises:

1. Prove that a regular semisimple element of a semisimple Lie
algebra lies in only finitely many Borel subalgebras.

2. Prove that the intersection of two Borel subalgebras in a
semisimple Lie algebra contains a CSA (Hint: if ${\msr B}$ and
${\msr B}'$ are Borel subalgebras and ${\msr N}=[{\msr B},{\msr
B}],\;{\msr N}'=[{\msr B}',{\msr B}']$, then ${\msr B}^\perp={\msr
N},\;{\msr N}^\perp={\msr B}$ and ${{\msr B}'}^\perp={\msr
N}',\;{{\msr N}'}^\perp={\msr B}'$. Moreover, ${\msr B}\bigcap{\msr
N}'\subset {\msr N}$ and  ${\msr B}'\bigcap{\msr N}\subset {\msr
N}'$).

\section{Simple Lie Algebra of Exceptional Types}

 In this section, we construct finite-dimensional complex simple Lie algebra of exceptional
 types.

 We  use ${\msr G}^X$ to denote the simple Lie algebra of type
$X$. It is well known that the simple Lie algebra of type $G_2$ is
exactly the derivation Lie
 algebra of a special 8-dimensional algebra with an identity element, which is called {\it octonion
 algebra}. Since any derivation annihilates the identity element
 (exercise), the octonion
 algebra gives a 7-dimensional representation of ${\msr G}^{G_2}$,
 which is presented as follows. Denote by
$\mbb{Z}_3=\mbb{Z}/3\mbb{Z}$ and identify $i+3\mbb{Z}$ with $i$ for
$i\in\{1,2,3\}$. Write
\begin{eqnarray*}\hspace{1.2cm}o(7,\mbb{C})&=&\sum_{i,j\in\mbb{Z}_3}[\mbb{C}(E_{i,j}-E_{j',i'})
+\mbb{C}(E_{i,j'}-E_{j,i'}) +\mbb{C}(E_{i',j}-E_{j',i})]\\ &
&+\sum_{i\in\mbb{Z}_3}[\mbb{C}(E_{0,i}-E_{i',0})+\mbb{C}(E_{0,i'}-E_{i,0})].\hspace{4.3cm}(4.4.1)\end{eqnarray*}
Then the Lie subalgebra
$${\msr K}=\sum_{i,j\in\mbb{Z}_3;i\neq
j}\mbb{C}(E_{i,j}-E_{j',i'})+\sum_{r=1,2}\mbb{C}(E_{r,r}-E_{r+1,r+1}-E_{r',r'}+E_{(r+1)',(r+1)'})\eqno(4.4.2)$$
 is isomorphic to $sl(3,\mbb{C})$. Denote
$$C_r=E_{r,(r+2)'}-E_{r+2,r'}+\sqrt{2}(E_{0,r+1}-E_{(r+1)',0}),\eqno(4.4.3)$$
$$C_r'=E_{r',r+2}-E_{(r+2)',r}+\sqrt{2}(E_{0,(r+1)'}-E_{r+1,0})\eqno(4.4.4)$$
for $r\in\mbb{Z}_3$. Then
$$[C_r,C_{r+1}]=2C_{r+2}',\qquad[C_r',C_{r+1}']=2C_{r+2},\eqno(4.4.5)$$
$$[C_r,C_r']=E_{r',r'}-E_{r,r}+E_{(r+2)',(r+2)'}-E_{r+2,r+2}+2(E_{r+1,r+1}-E_{(r+1)',(r+1)'})
\eqno(4.4.6)$$
$$[C_r,C_{r+1}']=3(E_{r+2,r+1}-E_{(r+1)',(r+2)'}) \eqno(4.4.7)$$ for
$r\in\mbb{Z}_3$. Moreover, $\sum_{r=1}^3\mbb{C}C_r$ and
$\sum_{r=1}^3\mbb{C}C_r'$  form $(\ad {\msr K})$-submodules, which
are isomorphic the three dimensional natural $sl(3,\mbb{C})$-module
and its dual, respectively. Thus
$${\msr G}^{G_2}={\msr
K}+\sum_{r=1}^3(\mbb{C}C_r+\mbb{C}C_r')\eqno(4.4.8)$$ forms a Lie
subalgebra of $o(7,\mbb{C})$, which is a simple Lie algebra of type
$G_2$.

 Note that  if $|\be(h_\al)|<|\al(h_\be)|$, then
 $|\la\be,\al\ra|<|\la\al,\be\ra|$, which implies
 $(\al,\al)>(\be,\be)$. Thus we set
$$h_1=E_{1,1}-E_{2,2}-E_{1',1'}
+E_{2',2'},\;h_2=-2E_{1,1}+E_{2,2}+E_{3,3}+2E_{1',1'}-E_{2',2'}-E_{3',3'},\eqno(4.4.9)$$
Take the Cartan subalgebra
$$H=\mbb{C}h_1+\mbb{C}h_2.\eqno(4.4.10)$$
Then
$${\msr
G}^{G_2}_{\al_1}=\mbb{C}(E_{1,2}-E_{2',1'}),\;\;{\msr
G}^{G_2}_{\al_2}=\mbb{C}C_3,\;\;{\msr
G}^{G_2}_{\al_1+\al_2}=\mbb{C}C_1,\eqno(4.4.11)$$
$${\msr G}^{G_2}_{\al_1+2\al_2}=\mbb{C}C_2',\;\;{\msr
G}^{G_2}_{\al_1+3\al_2}=\mbb{C}(E_{3,1}-E_{1',3'}),\;\; {\msr
G}^{G_2}_{2\al_1+3\al_2}=\mbb{C}(E_{3,2}-E_{2',3'}),\eqno(4.4.12)$$
$${\msr
G}^{G_2}_{-\al_1}=\mbb{C}(E_{2,1}-E_{1',2'}),\;\;{\msr
G}^{G_2}_{-\al_2}=\mbb{C}C_3',\;\;{\msr
G}^{G_2}_{-\al_1-\al_2}=\mbb{C}C_1',\eqno(4.4.13)$$
$${\msr G}^{G_2}_{-\al_1-2\al_1}=\mbb{C}C_2,\;\;{\msr
G}^{G_2}_{-\al_1-3\al_2}=\mbb{C}(E_{1,3}-E_{3',1'}),\;\; {\msr
G}^{G_2}_{-2\al_1-3\al_2}=\mbb{C}(E_{2,3}-E_{3',2'}).\eqno(4.4.14)$$
\pse

Let $\Lmd_r$ be the root lattices of types $X=E_6,\;E_7$ and $E_8$
(cf. (3.4.33), (3.4.37) and (3.4.38)). Recall that $(\cdot,\cdot)$
is a symmetric $\mbb Z$-bilinear form on $\Lmd_r$ such that
$(\al,\al)=2$ for $\al\in \Phi_X$. In fact,
$$0<(\sum_{i=1}^nk_i\al_i,\sum_{i=1}^nk_i\al_i)=2(\sum_{p=1}^nk_p^2+\sum_{i<
j}k_ik_j(\al_i,\al_j))\eqno(4.4.15)$$ for
$0\neq\sum_{i=1}^nk_i\al_i\in \Lmd_r$. So $(\al,\al)$ is a positive
even integer for any $0\neq\al\in \Lmd_r$. For $\al,\be\in\Phi$,
$$\al+\be\in\Phi\Rightarrow
2=(\al+\be,\al+\be)=(\al,\al)+(\be,\be)+2(\al ,\be)=2+2+2(\al
,\be),\eqno(4.4.16)$$ So $(\al,\be)=-1$. By Lemma 3.1.2,
 $$\al+\be\in\Phi\Leftrightarrow(\al,\be)=-1\qquad\for\;\;\al,\be\in\Phi.\eqno(4.4.17)$$ Denote
$\Pi=\{\al_1,\al_2,...,\al_n\}$ according to the Dynkin diagrams of
$E_6,\;E_7$ and $E_8$ in Theorem 3.3.1. Define
$F(\cdot,\cdot):\;\Lmd_r\times \Lmd_r\rta \{\pm 1\}$ by
$$F(\sum_{i=1}^nk_i\al_i,\sum_{j=1}^nl_j\al_j)=(-1)^{\sum_{i=1}^nk_il_i+\sum_{i>j}k_il_j
(\al_i,\al_j)},\qquad k_i,l_j\in\mbb{Z}.\eqno(4.4.18)$$ Then for
$\al,\be,\gm\in \Lmd_r$,
$$F(\al+\be,\gm)=F(\al,\gm)F(\be,\gm),\;\;F(\al,\be+\gm)=F(\al,\be)F(\al,\gm),
\eqno(4.4.19)$$
$$F(\al,\be)F(\be,\al)^{-1}=(-1)^{(\al,\be)},\;\;F(\al,\al)=(-1)^{(\al,\al)/2}.
\eqno(4.4.20)$$ In particular,
$$F(\al,\be)=-F(\be,\al)\qquad
\mbox{if}\;\;\al,\be,\al+\be\in\Phi\eqno(4.4.21)$$ by (4.4.17).

Denote
$$H=\sum_{i=1}^n\mbb{C}\al_i.\eqno(4.4.22)$$
 Define
$${\msr
G}^X=H\oplus\bigoplus_{\al\in\Phi}\mbb{C}E_{\al},\eqno(4.4.23)$$
 an algebraic operation $[\cdot,\cdot]$ and a symmetric bilinear
 form $(\cdot,\cdot)$ on ${\msr G}$ by:
 $$[H,H]=0,\;\;[h,E_{\al}]=-[E_{\al},h]=(h,\al)E_{\al},\;\;[E_{\al},E_{-\al}]=-\al,
 \eqno(4.4.24)$$
 $$[E_{\al},E_{\be}]=\left\{\begin{array}{ll}0&\mbox{if}\;\al+\be\not\in\Phi,\\
 F(\al,\be)E_{\al+\be}&\mbox{if}\;\al+\be\in\Phi.\end{array}\right.\eqno(4.4.25)$$
\pse

{\bf Theorem 4.4.1} {\it The pair $({\msr G}^X,[\cdot,\cdot])$ forms
a simple Lie algebra of type $X$.}

{\bf Proof}. The skew symmetry of $[\cdot,\cdot]$ follows from
(4.4.21). Recall the Jacobi identity:
$$[[x,y],z]+[[y,z],x]+[[z,x],y]=0.\eqno(4.4.26)$$
It is trivial if $x,y,z\in H$. For $x,y\in H$ and $z=E_{\al}$ with
$\al\in \Phi$, (4.4.26) becomes
$$0+(y,\al)(-(x,\al))E_{\al}+(-(x,\al))(-(y,\al))E_{\al}=0,\eqno(4.4.27)$$
which obviously holds. When $x\in H,\;y=E_{\al}$ and $z=E_{\be}$
with $\al,\be\in\Phi$, (4.4.26) gives
$$(x,\al)[E_{\al},E_{\be}]+(-(x,\al+\be))[E_{\al},E_{\be}]+(-(\be,x))[E_{\be},E_{\al}]=0,
\eqno(4.4,28)$$ which holds because of the skew symmetry. Observe
that (4.4.26) holds whenever two of $\{x,y,z\}$ are equal.

 Consider
$x=E_{\al},\;y=E_{-\al}$ and $z=E_{\be}$ with $\al,\be\in\Phi$ and
$\be\neq\al,-\al$. If $\al+\be\in\Phi$, then $(\al,\be)=-1$, and
$(-\al,\be)=1$. So $\be-\al\not\in\Phi$ by (4.4.17). Now (4.4.26)
yields
$$(-\al,\be)E_{\be}+0+F(\be,\al)F(\be+\al,-\al)E_{\be}=0,\eqno(4.4.29)$$
or equivalently, $ 1+F(\be,\al)F(\be+\al,-\al)=0.$ This holds
because
\begin{eqnarray*}\hspace{1cm}F(\be,\al)F(\be+\al,-\al)&=&
F(\be,\al)F(\be,-\al)F(\al,-\al)
\\ &=&F(\be,0)F(\al,\al)^{-1}=-1\hspace{4.6cm}(4.4.30)\end{eqnarray*} by (4.4.20) and
(4.4.21). Assume $\al+\be,\be-\al\not\in\Phi$, or equivalently,
$$(\be\pm\al,\be\pm\al)=4\pm
2(\be,\al)>2\Rightarrow\pm(\be,\al)>-1\sim
-1<(\be,\al)<1.\eqno(4.4.31)$$ Thus we have $(\be,\al)=0$ because
$(\be,\al)\in\mbb{Z}$. So all three terms in (4.4.26) are 0.

Suppose  that $x=E_{\al},\;y=E_{\be},\;z=E_{\gm}$ with
$\al,\be,\gm\in\Phi$ such that $\al\neq\pm\be,\pm\gm$ and
$\be\neq\pm \gm$, and also at least one of the terms in (4.4.26) is
nonzero. We may assume $\al+\be\in\Phi$ and $\al+\be+\gm=0$ or
$\al+\be+\gm\in\Phi$. First we have $(\al,\be)=-1$ by (4.4.17).
Observe that
$$F(\be,-\al)=F(\be,\al)^{-1}=F(-\be,\al)=-F(\al,\be)\eqno(4.4.32)$$
and
$$F(\be,-\be)=F(-\al,\al)=-1\eqno(4.4.33)$$
by (4.4.20). If $\al+\be+\gm=0$, then $\gm=-\al-\be$ and (4.4.26)
gives
$$F(\al,\be)(-(\al+\be))+F(\be,-(\al+\be))\al+F(-(\al+\be),\al)\be=0,\eqno(4.4.34)$$
or equivalently,
$$F(\al,\be)=F(\be,-(\al+\be))=F(-(\al+\be),\al),\eqno(4.4.35)$$
which holds by (4.4.18)-(4.4.20).  Assume $\al+\be+\gm\in\Phi$. So
$$(\al+\be,\gm)=(\al,\gm)+(\be,\gm)=-1.\eqno(4.4.36)$$
Since $[(\al,\gm)]^2<(\al,\al)(\gm,\gm)=4$, we have $(\al,\gm)=0,\pm
1$. Similarly, $(\be,\gm)=0,\pm 1$. Thus (4.4.36) forces
$(\al,\gm)=-1,\;(\be,\gm)=0$ or $(\al,\gm)=0,\;(\be,\gm)=-1.$ We may
assume $(\al,\gm)=-1$ and $ (\be,\gm)=0$. Then (4.4.26) becomes
$$F(\al,\be)F(\al+\be,\gm)E_{\al+\be+\gm}+0
+F(\gm,\al)F(\al+\gm,\be)E_{\al+\be+\gm}=0,\eqno(4.4.37)$$ or
equivalently, \begin{eqnarray*} \hspace{2cm}&
&F(\al,\be)F(\al+\be,\gm) +F(\gm,\al)F(\al+\gm,\be)=0\\
& \sim &F(\al+\be,\gm) +F(\gm,\al)F(\gm,\be)=0\\ &\sim&
F(\al+\be,\gm)
+F(\gm,\al+\be)=0,\hspace{5.9cm}(4.4.38)\end{eqnarray*} which is
implied by (4.4.21). This completes the proof of the Jacobi identity
(4.4.26). The simplicity follows from the fact that $\Phi$ is an
irreducible root system. $\qquad\Box$\psp

Now to understand the simple Lie algebras of types $E_6,\;E_7$ and
$E_8$ explicitly, we need to write all the positive roots as a
linear combinations of simple positive roots. In terms of the order
in the Dynkin diagram of $E_8$ in Theorem 3.3.1, we have:
$$\al_1=\ves_1-\ves_2,\qquad\al_2=\frac{1}{2}(\sum_{s=4}^8\ves_s-\sum_{r=1}^3\ves_r),\eqno(4.4.39)$$
$$\al_{i+2}=\ves_{i+1}-\ves_{i+2}\qquad\for\;\;
i\in\ol{1,6},\eqno(4.4.40)$$ where $\{\ves_i\mid i\in\ol{1,8}\}$ is
the standard basis in the Euclidean space $\mbb{R}^8$. Recall the
root system
$$\Phi_{E_8}=\left\{\frac{1}{2}\sum_{i=1}^8(-1)^{\iota_i}\ves_i\mid
\iota_i\in\mbb{Z}_2,\;\sum_{r=1}^8\iota_r=0\;\mbox{in}\;\mbb{Z}_2\right\}\bigcup
\Phi_{D_8}.\eqno(4.4.41)$$

Let $\Phi^+_{E_8}$ be the set of  positive roots. For convenience,
 we also denote
$$\al_{(k_1,...,k_r)}=\sum_{s=1}^rk_s\al_s,\;k_r\neq
0\eqno(4.4.42)$$ and
$$E_{(k_1,...,k_r)}=E_{\al_{(k_1,...,k_r)}}\qquad\mbox{if}\;\;\al_{(k_1,...,k_r)}\in\Phi_{E_8}^+.\eqno(4.4.43)$$
Moreover, we write
$$(aE_\be+bE_\gm)'=aE'_\be+bE'_\gm=aE_{-\be}+bE_{-\gm}\qquad\for\;\;a,b\in\mbb F,\;\be,\gm\in
\Phi_{E_8}^+.\eqno(4.4.44)$$ Take the convention
$\sum_{r=i}^j\al_r=0$ if $j<i$. We calculate
$$\ves_1-\ves_j=\al_1+\sum_{r=3}^j\al_r\qquad\for\;\;j\in\ol{2,8},\eqno(4.4.45)$$
$$\ves_i-\ves_j=\sum_{r=i+1}^j\al_r\qquad\for\;\;2\leq i<j\leq 8,\eqno(4.4.46)$$
$$\ves_1+\ves_j=\al_{(2,2,3,4,3,2,1)}+\sum_{s=j+1}^8\al_s
\qquad\for\;\;j\in\ol{2,8},\eqno(4.4.47)$$
$$\ves_7+\ves_8=\al_{(1,2,2,3,2,1)},\qquad
\frac{1}{2}(\sum_{r=2}^8\ves_r-\ves_1)=\al_{(1,3,3,5,4,3,2,1)}\eqno(4.4.48)$$
$$\ves_i+\ves_j=\al_{(1,2,2,3,2,1)}+\sum_{r=i+1}^7\al_r+\sum_{s=j+1}^8\al_s
\qquad\for\;\;2\leq i<j\leq 8,\eqno(4.4.49)$$
$$\frac{1}{2}(\sum_{r\neq 1,
j,k}\ves_r-\ves_1-\ves_j-\ves_k)=\sum_{s=2}^j\al_s+\sum_{t=4}^k\al_t\qquad\for\;\;2\leq
j<k\leq 8,\eqno(4.4.50)$$
$$\frac{1}{2}(\sum_{r\neq i}\ves_r-\ves_i)=\al_{(2,2,3,5,4,3,2,1)}+
\sum_{r=2}^i\al_r\qquad\for\;\;i\in\ol{2,8},\eqno(4.4.51)$$
$$\frac{1}{2}(\sum_{r\neq i,
j,k}\ves_r-\ves_i-\ves_j-\ves_k)=\sum_{\iota=1}^i\al_\iota+
\sum_{s=3}^j\al_s+\sum_{t=4}^k\al_t\eqno(4.4.52)$$ for $2\leq i<
j<k\leq 8$. \psp

{\bf Proposition 4.4.3}. {\it The positive roots of type $E_8$ in
terms of linear combinations of $\Pi_{E_8}$ are given in (4.4.39),
(4.4.40) and (4.4.45)-(4.4.52). Moreover, $|\Phi_{E_8}^+|=120$. In
particular},
$$\mbox{dim}\:{\msr G}^{E_8}=248.\eqno(4.4.53)$$
\pse

 According to Example 3.4.2 and the above proposition,
the positive roots of type $E_7$ in terms of linear combinations of
$\Pi_{E_7}$ are:
$$\al_{(2,2,3,4,3,2,1)},\;\;\{\al_{(1,2,2,3,2,1)}+\sum_{r=i+1}^7\al_r
\mid i\in\ol{2,7}\},\eqno(4.4.54)$$
$$\{\al_1+\sum_{r=3}^j\al_r\mid j\in\ol{2,7}\},\;\;
\{\sum_{r=i+1}^j\al_r\mid 2\leq i<j\leq 7\},\eqno(4.4.55)$$
$$\{\sum_{s=2}^j\al_s+\sum_{t=4}^k\al_t\mid 2\leq j<k\leq
7\},\eqno(4.4.56)$$
$$\{\sum_{\iota=1}^i\al_\iota+
\sum_{s=3}^j\al_s+\sum_{t=4}^k\al_t\mid 2\leq i< j<k\leq
7\}.\eqno(4.4.57)$$ Thus
$$|\Phi_{E_7}^+|=63,\qquad\mbox{dim}\:{\msr G}^{E_7}=133.\eqno(4.4.58)$$
 Moreover, the positive roots of type $E_6$ in
terms of linear combinations of $\Pi_{E_6}$ are:
$$\{\al_1+\sum_{r=3}^j\al_r\mid j\in\ol{2,6}\}\bigcup
\{\sum_{r=i+1}^j\al_r\mid 2\leq i<j\leq 6\},\eqno(4.4.59)$$
$$\al_{(1,2,2,3,2,1)},\;\;\{\sum_{s=2}^j\al_s+\sum_{t=4}^k\al_t\mid 2\leq j<k\leq
6\}\eqno(4.4.60)$$ and
$$\{\sum_{\iota=1}^i\al_\iota+
\sum_{s=3}^j\al_s+\sum_{t=4}^k\al_t\mid 2\leq i< j<k\leq
6\}.\eqno(4.4.61)$$ Hence
$$|\Phi_{E_6}^+|=36,\qquad\mbox{dim}\:{\msr G}^{E_6}=78.\eqno(4.4.62)$$
\pse

Finally we give a construction of the simple Lie algebra of type
$F_4$. From the Dynkin diagram of $E_6$ in Theorem 3.3.1, we have
the following automorphism of $\Phi_{E_6}$:
$$\sgm(\sum_{i=1}^6k_i\al_i)=k_6\al_1+k_2\al_2+k_5\al_3+k_4\al_4+k_3\al_5+k_1\al_6\eqno(4.4.63)$$
for $\sum_{i=1}^6k_i\al_i\in \Phi_{E_6}$. In the construction of the
Lie algebra ${\msr G}$ in (4.4.22)-(4.4.25), we take
$\Lmd_r=\Lmd_r(E_6)$ and replace the $F$ in (4.4.18) by
$$F(\sum_{i=1}^6k_i\al_i,\sum_{j=1}^6l_j\al_j)=(-1)^{\sum_{i=1}^6k_il_i+
k_1l_3+k_4l_2+k_3l_4+k_5l_4+k_6l_5},\qquad
k_i,l_j\in\mbb{Z}.\eqno(4.4.64)$$ Then (4.4.19)-(4.4.21) holds and
$$F(\sgm(\al),\sgm(\be))=F(\al,\be)\qquad\for\;\;\al,\be\in
\Lmd_r(E_6).\eqno(4.4.65)$$ Thus we have the following automorphism
$\hat\sgm$ of the Lie algebra ${\msr G}$:
$$\hat\sgm(\sum_{i=1}^6b_i\al_i)=\sum_{i=1}^6b_i\sgm(\al_i),\qquad
b_i\in\mbb{C},\eqno(4.4.66)$$
$$\hat\sgm(E_\al)=E_{\sgm(\al)}\qquad\for\;\;\al\in\Phi_{E_6}.\eqno(4.4.67)$$

 Using the notation in (4.4.43), we write
$$F_{(1)}=E_{\al_2},\;\;F_{(0,1)}=E_{\al_4},\;\;
 F_{(0,0,1)}=E_{\al_3}+E_{\al_5},\;\;F_{(0,0,0,1)}=E_{\al_1}+E_{\al_6},
\eqno(4.4.68)$$
$$F_{(1,1)}=E_{(0,1,0,1)},\;\;
F_{(0,1,1)}=E_{(0,0,1,1)}+E_{(0,0,0,1,1)},\;\;F_{(0,0,1,1)}=E_{(1,0,1)}+E_{(0,0,0,0,1,1)},\eqno(4.4.69)$$
$$F_{(0,0,1,1)}=E_{(1,0,1)}+E_{(0,0,0,0,1,1)},\;\;F_{(1,1,1)}=
E_{(0,1,1,1)}+E_{(0,1,0,1,1)}, \eqno(4.4.70)$$
$$F_{(0,1,1,1)}=E_{(1,0,1,1)}+E_{(0,0,0,1,1,1)},
\;\;F_{(0,1,2)}=E_{(0,0,1,1,1)},\;\;F_{(1,1,2)}=E_{(0,1,1,1,1)},
\eqno(4.4.71)$$
$$F_{(0,1,2,1)}=E_{(1,0,1,1,1)}+
E_{(0,0,1,1,1,1)},\;\;F_{(1,1,1,1)}=E_{(1,1,1,1)}+
E_{(0,1,0,1,1,1)},\eqno(4.4.72)$$
$$F_{(1,2,2)}=E_{(0,1,1,2,1)},\;\;
F_{(1,1,2,1)}=E_{(1,1,1,1,1)}+
E_{(0,1,1,1,1,1)},\;\;F_{(0,1,2,2)}=E_{(1,0,1,1,1,1)},
\eqno(4.4.73)$$
$$F_{(1,2,2,1)}=E_{(1,1,1,2,1)}
+E_{(0,1,1,2,1,1)},\; F_{(1,1,2,2)}=E_{(1,1,1,1,1,1)},\;
F_{(1,2,2,2)}=E_{(1,1,1,2,1,1)},\eqno(4.4.74)$$
$$F_{(1,2,3,1)}=E_{(1,1,2,2,1)}
+E_{(0,1,1,2,2,1)},\;\;F_{(1,2,3,2)}=E_{(1,1,2,2,1,1)}
+E_{(1,1,1,2,2,1)},\eqno(4.4.75)$$
$$F_{(1,2,4,2)}=E_{(1,1,2,2,2,1)},\;\;F_{(1,3,4,2)}=E_{(1,1,2,3,2,1)},\;\;
F_{(2,3,4,2)}=E_{(1,2,2,3,2,1)} .\eqno(4.4.76)$$ Moreover, we set
$$h_1=\al_2,\qquad h_2=\al_4,\qquad h_3=\al_3+\al_5,\qquad h_4=\al_1+\al_6.\eqno(4.4.77)$$

The Dynkin diagram of $F_4$ is

\begin{picture}(60,12)\put(2,0){$F_4$:}
\put(21,0){\circle{2}}\put(21,-5){1}\put(22,0){\line(1,0){12}}
\put(35,0){\circle{2}}\put(35,-5){2}\put(35,1.2){\line(1,0){13.6}}
\put(35,-0.8){\line(1,0){13.6}}\put(41,-1){$\ra$}\put(48.5,0){\circle{2}}\put(48.5,-5){3}\put(49.5,0)
{\line(1,0){12}}\put(62.5,0){\circle{2}}\put(62.5,-5){4}
\end{picture}
\vspace{0.6cm}

\noindent In order to make notation distinguishable, we add a bar on
the roots in the root system $\Phi_{F_4}$ of type $F_4$. In
particular, we let $\{\bar\al_1,\bar\al_2,\bar\al_3,\bar\al_4\}$ be
the simple positive roots corresponding to the above  Dynkin diagram
of $F_4$, where $\bar\al_1,\bar\al_2$ are long roots and
$\bar\al_3,\bar\al_4$ are short roots. Let $\Phi_{F_4}^+$ be the set
of positive roots of $F_4$. Denote
$$S_{F_4}=\{(k_1,...,k_r)\mid \bar\al_{(k_1,...,k_r)}=\sum_{i=1}^rk_i\bar\al_i\in\Phi_{F_4}^+,\;k_r\neq
0\}.\eqno(4.4.78)$$ Moreover, we denote
$$F'_{(k_1,...,k_r)}=(F_{(k_1,...,k_r)})'\qquad\for\;\;(k_1,...,k_r)\in
S_{F_4}\eqno(4.4.79)$$ (cf. (4.4.44)). Then
 the simple Lie algebra of
type $F_4$ is
$${\msr G}^{F_4}=\{u\in {\msr
G}^{E_6}\mid\hat\sgm(u)=u\}=\sum_{\varpi\in S_{F_4}}(\mbb
FF_\varpi+\mbb F F'_\varpi)+\sum_{i=1}^4\mbb Fh_i.\eqno(4.4.80)$$ In
fact, $F_\varpi$ is a root vector of root $\bar\al_\varpi$,
$F'_\varpi$ is a root vector of root $-\bar\al_\varpi$, and
$H_{F_4}=\sum_{i=1}^4\mbb Fh_i$ is the corresponding Cartan
subalgebra of $\msr G^{F_4}$.

\chapter{Hight-Weight Representation Theory}

In this chapter, we always assume that the base field
$\mbb{F}=\mbb{C}$. By Weyl's Theorem, any finite-dimensional
representation of a finite-dimensional semisimple Lie algebra is
completely reducible. The main goal in this chapter is to study
finite-dimensional irreducible representation of a
finite-dimensional semisimple Lie algebra. First we introduce the
universal enveloping algebra of a Lie algebra and prove the
Poincar\'{e}-Birkhoff-Witt (PBW) Theorem on its basis. Then we use
the universal enveloping algebra of a finite-dimensional semisimple
Lie algebra ${\msr G}$ to construct the Verma module of ${\msr G}$,
which is the maximal  module of a given highest weight. Moreover, we
prove that any finite-dimensional irreducible ${\msr G}$-module is
the quotient of a Verma module modulo its maximal proper submodule,
whose generators are explicitly given. Furthermore, the Weyl's
character formula of a finite-dimensional irreducible ${\msr
G}$-module is derived and the dimensional formula of the module is
determined. Finally, we decompose the tensor module of two
finite-dimensional irreducible ${\msr G}$-modules into a direct sum
of irreducible ${\msr G}$-submodules in terms of characters.

\section{Universal Enveloping Algebras}

In this section, we construct the universal enveloping algebra of a
Lie algebra and prove the Poincar\'{e}-Birkhoff-Witt Theorem on its
basis.

Let $V$ be vector space. The {\it free associative algebra generated
by} $V$ is an associative algebra ${\msr A}(V)$ and a linear
injective map $\iota:V\rta {\msr A}(V)$ such that if $\vf: V\rta
{\msr B}$ is linear map from $V$ to an associative algebra ${\msr
B}$, then there exists a unique algebra homomorphism $\phi: {\msr
A}(V)\rta {\msr B}$ for which $\phi(1_{{\msr A}(V)})=1_{\msr B}$ and
$\phi|_V=\vf$. If $S$ is a basis of $V$, we also say that ${\msr
A}(V)$ is a free associative algebra generated by $S$.

The algebra ${\msr A}(V)$ is indeed the {\it tensor algebra} over
$V$ as follows. Suppose that $S=\{u_i\mid i\in I\}$ is a basis of
$V$. We set ${\cal I}(V)$ to be the vector space with
$$\{1,u_{_{l_1}}\otimes u_{_{l_2}}\otimes\cdots\otimes
u_{_{l_s}}\mid 0<s\in\mbb{Z},\;l_1,...,l_s\in I\}\;\mbox{as a
basis}\eqno(5.1.1)$$ and define an associative algebraic operation
$\cdot$ on ${\cal I}(V)$ by $1\cdot 1=1$,
$$1\cdot (u_{_{l_1}}\otimes\cdots\otimes u_{_{l_r}})=(u_{_{l_1}}\otimes\cdots\otimes u_{_{l_r}})\cdot 1=u_{_{l_1}}\otimes\cdots\otimes
u_{_{l_r}}\eqno(5.1.2)$$ and
$$(u_{_{l_1}}\otimes\cdots\otimes u_{_{l_r}})\cdot (u_{_{k_1}}\otimes\cdots\otimes
u_{_{k_s}}) =u_{_{l_1}}\otimes\cdots\otimes u_{_{l_r}}\otimes
u_{_{k_1}}\otimes\cdots\otimes u_{_{k_s}}.\eqno(5.1.3)$$ Given
$$v_i=\sum_{j\in I}a_{i,j}u_j\in V\qquad
\for\;\;i\in\ol{1,r},\eqno(5.1.4)$$ we define
$$v_1\otimes v_2\otimes\cdots\otimes v_r=\sum_{j_1,...,j_r\in
I}a_{_{1,j_1}}a_{_{2,j_2}}\cdots a_{_{r,j_r}}u_{_{j_1}}\otimes
u_{_{j_1}}\otimes\cdots \otimes u_{_{j_r}}.\eqno(5.1.5)$$
 As
associative algebras, ${\msr A}(V)\cong {\cal I}(V)$. The map
$\iota: V\rta {\cal I}$ is the inclusion.

Denote by
$${\msr J}_0(V)=\mbox{the ideal of}\;{\cal I}(V)\;\mbox{generated
by}\;\{u\otimes v-v\otimes u\mid u,v\in V\}.\eqno(5.1.6)$$ The {\it
symmetric tensor algebra over} $V$ (or {\it free commutative
associative algebra generated by} $V$) is defined as the quotient
algebra
$${\msr S}(V)={\cal I}(V)/{\msr J}_0(V).\eqno(5.1.7)$$
In fact, we have the following algebra isomorphism: $1_{{\msr
S}(V)}\mapsto 1$ and
$$u_{_{j_1}}\otimes u_{_{j_2}}\cdots \otimes u_{_{j_r}}+{\msr
J}_0(V)\mapsto x_{_{j_1}} x_{_{j_2}}\cdots x_{_{j_r}}\eqno(5.1.8)$$
from ${\msr S}(V)$ to the polynomial algebra $\mbb{C}[x_i\mid i\in
I]$.

Suppose that ${\msr G}$ is a Lie algebra. We have the tensor algebra
${\cal I}({\msr G})$. Set
$${\msr J}_1({\msr G})=\mbox{the ideal of}\;{\cal I}({\msr G})\;\mbox{generated
by}\;\{u\otimes v-v\otimes u-[u,v]\mid u,v\in {\msr G}
\}.\eqno(5.1.9)$$ Define
$$U({\msr G})={\cal I}({\msr G})/{\msr J}_1({\msr
G}).\eqno(5.1.10)$$ For convenience, we identify $u+{\msr J}_1({\msr
G})$ in $U({\msr G})$ with $u\in{\msr G}$. Suppose that $\{u_i\mid
i\in I\}$ is a basis of ${\msr G}$ and $I$ has a total ordering
$\prec$. If $I=\{1,2,...,n\}$ is finite, we can take a total
ordering just usual $<$. Recall that $\mbb{Z}_+$ is the set of
positive integers.

\psp

{\bf Theorem 5.1.1 (Poincar\'{e}-Birkhoff-Witt (PBW) Theorem)}.
\index{PBW Theorem} {\it The set
$$\{1,u_{i_1}^{m_1}u_{i_2}^{m_2}\cdots
u_{i_r}^{m_r}\mid r,m_1,...,m_r\in\mbb{Z}_+,\;i_1,...,i_r\in
I;\;i_1\prec i_2\prec\cdots\prec i_r\}\eqno(5.1.11)$$ forms a basis
of} $U({\msr G})$.

 {\it
Proof}. We define ${\msr U}_0=\mbb{C}1,\;{\msr U}_1={\msr G}$ and
$${\msr U}_i=\mbox{Span}\:\{1,u_{i_1}u_{i_2}\cdots u_{i_r}\mid
i_s\in I,\;i\geq r\in\mbb{Z}_+\}.\eqno(5.1.12)$$ Then we have
$$U({\msr G})=\bigcup_{i=0}^{\infty}{\msr
U}_i.\eqno(5.1.13)$$ Let ${\msr U}'$ be the subspace of $ U({\msr
G})$ spanned by (5.1.11). Now ${\msr U}_0,{\msr U}_1\subset {\msr
U}'$. Suppose that ${\msr U}_k\subset {\msr U}'$.

Note
$${\msr U}_{k+1}=\mbox{Span}\:\{u_{i_1}\cdots u_{i_{_{k+1}}}\mid
i_s\in I\}+{\msr U}_k.\eqno(5.1.14)$$ Suppose $i_{s+1}\prec i_s$ for
some $s\in\ol{1,k}$. Then
\begin{eqnarray*}& &u_{i_1}\cdots
u_{i_{s-1}}u_{i_s}u_{i_{s+1}}u_{i_{s+2}}\cdots u_{_{k+1}}\\ &
=&u_{i_1}\cdots u_{i_{s-1}}u_{i_{s+1}}u_{i_s}u_{i_{s+2}}\cdots
u_{_{k+1}}+u_{i_1}\cdots
u_{i_{s-1}}(u_{i_s}u_{i_{s+1}}-u_{i_{s+1}}u_{i_s})u_{i_{s+2}}\cdots
u_{_{k+1}}\\ &=&u_{i_1}\cdots
u_{i_{s-1}}u_{i_{s+1}}u_{i_s}u_{i_{s+2}}\cdots
u_{_{k+1}}+u_{i_1}\cdots
u_{i_{s-1}}[u_{i_s},u_{i_{s+1}}]u_{i_{s+2}}\cdots u_{_{k+1}}
\\ &\equiv&u_{i_1}\cdots
u_{i_{s-1}}u_{i_{s+1}}u_{i_s}u_{i_{s+2}}\cdots
u_{_{k+1}}\;\;(\mbox{mod}\:{\msr
U}_k).\hspace{6cm}(5.1.15)\end{eqnarray*} This shows that we can
change the orders of adjacent factors in the product $u_{i_1}\cdots
u_{i_{_k}}$ modulo ${\msr U}_k$. After finite steps of exchanging
the positions of adjacent factors, we get
$$u_{i_1}\cdots u_{i_{_{k+1}}}=u_{_{j_1}}^{m_1}\cdots
u_{_{j_s}}^{m_s}\;\;(\mbox{mod}\:{\msr U}_k),\eqno(5.1.16)$$ where
$$j_1,...,j_s\in I,\;j_1\prec j_2\prec\cdots\prec
j_s,\;\;m_1,...,m_s\in\mbb{Z}_+,\;\;\sum_{r=1}^sm_r=k+1.\eqno(5.1.17)$$
By (5.1.11), $u_{_{j_1}}^{m_1}\cdots u_{_{j_s}}^{m_s}\in{\msr U}'$.
Since ${\msr U}_{\:k}\subset {\msr U}'$ by assumption, we have
$$u_{i_1}\cdots u_{i_{_{k+1}}}\in {\msr U}'\qquad\mbox{for any}\;i_1,...,i_{k+1}\in
I.\eqno(5.1.18)$$ Thus ${\msr U}_{\:k+1}\subset {\msr U}'$. By
induction, $U({\msr G})={\msr U}'$.

By (5.1.16), ${\msr U}_{\:k+1}/{\msr U}_{\;k}$ is isomorphic to the
subspace of homogeneous polynomials of degree $k+1$ in
$\mbb{C}[x_i\mid i\in{\cal I}]$. Therefore, (5.1.11) is linearly
independent. $\qquad\Box$ \psp

Let ${\msr A}$ be an associative algebra. Suppose that $\vf: {\msr
G}\rta {\msr A}$ is a Lie algebra homomorphism, where the Lie
bracket of ${\msr A}$ is the commutator. There exists a unique
associative algebra homomorphism $\vf_1:{\cal I}({\msr G})\rta {\msr
A}$ such that $\vf_1(1)=1$ and $\vf_1|_{\msr G}=\vf$. Moreover,
$$\vf_1(u\otimes v-v\otimes
u-[u,v])=\vf(u)\vf(v)-\vf(v)\vf(u)-\vf([u,v])=0\eqno(5.1.19)$$ for
$u,v\in{\msr G}$ because $\vf$ is a Lie algebra homomorphism. Hence
${\msr J}_1({\msr G})\subset \kn\vf_1$. So we have a unique
associative algebra homomorphism $\phi:U({\msr G})\rta {\msr A}$
defined by
$$\phi(w+{\msr J}_1({\msr G}))=\vf_1(w)\qquad\for\;\;w\in{\msr
I}({\msr G})\eqno(5.1.20)$$ with $\phi|_{\msr G}=\vf$.

Suppose that $\ol{\msr U}$ is an associative algebra with an
injective linear map $\iota:{\msr G}\rta \ol{\msr U}$ such that
given a Lie algebra homomorphism $\vf'$ from ${\msr G}$ to an
associative ${\msr A}$, there exists a unique associative algebra
homomorphism $\phi': \ol{\msr U}\rta {\msr A}$  for which
$\phi'\iota=\vf'$. Such an algebra $\ol{\msr U}$ is called an {\it
universal enveloping algebra} of ${\msr G}$. By the above paragraph,
there exists an associative algebra homomorphism $\tau:U({\msr
G})\rta \ol{\msr U}$ such that $\tau|_{\msr G}=\iota$. Since the
inclusion map from ${\msr G}$ to $U({\msr G})$ is an obvious Lie
algebra monomorphism, there exists an associative algebra
homomorphism $\tau':\ol{\msr U}\rta U({\msr G})$ such that
$\tau'\iota=\mbox{Id}_{\msr G}$.

Now $\tau\tau':\ol{\msr U}\rta \ol{\msr U}$ is an associative
algebra endomorphism such that $\tau\tau'\iota=\iota$. Since
$\mbox{Id}_{\ol{\msr U}}\iota=\iota$, we have
$\tau\tau'=\mbox{Id}_{\ol{\msr U}}$ by the uniqueness from the
universality of $\ol{\msr U}$. On the other hand, $\tau'\tau:
U({\msr G})\rta U({\msr G})$ an associative algebra endomorphism
such that $\tau'\tau|_{\msr G}=\tau'\iota=\mbox{Id}_{\msr G}$. Since
$U({\msr G})$ is generated by ${\msr G}$, we have
$\tau'\tau=\mbox{Id}_{U({\msr G})}$. Thus $\ol{\msr U}\cong U({\msr
G})$; that is, $U({\msr G})$ is the unique universal enveloping
algebra of ${\msr G}$.

\section{Highest-Weight Modules}

In this section, we use the universal enveloping algebra of a
finite-dimensional semisimple Lie algebra ${\msr G}$ to construct
the Verma module of ${\msr G}$. Moreover, we prove that any
finite-dimensional irreducible ${\msr G}$-module is the quotient of
a Verma module modulo its maximal proper submodule.

Let ${\msr G}$ be a Lie algebra (may not be finite-dimensional) with
a toral Cartan subalgebra $H$ and two subalgebras ${\msr G}_{\pm}$
such that
$${\msr G}={\msr G}_-\oplus H\oplus {\msr G}_+,\qquad [H,{\msr
G}_\pm]\subset {\msr G}_\pm.\eqno(5.2.1)$$ By PBW Theorem,
$$U({\msr G})=U({\msr G}_-)U(H)U({\msr
G}_+).\eqno(5.2.2)$$ The universality of $U({\msr G})$ implies that
a vector space forms a module of a Lie algebra ${\msr G}$ if and
only if it forms a module of the associative algebra $U({\msr G})$.
Suppose that $V$ is a ${\msr G}$-module. For any $\lmd\in H^\ast$,
we define
$$V_\lmd=\{v\in V\mid h(v)=\lmd(h)v\;\for\;h\in H\}.\eqno(5.2.3)$$
Any nonzero vector in the {\it weight subspace}\index{weight
subspace} $V_\lmd$ is called a {\it weight vector}\index{weight
vector} with {\it weight}\index{weight} $\lmd$. A weight vector $v$
is called {\it singular}\index{singular} if
$${\msr G}_+(v)=\{0\}.\eqno(5.2.4)$$
If $v$ is a singular vector, then
$$U({\msr G})(v)=U({\msr G}_-)(v)\eqno(5.2.5)$$
is a submodule by (5.2.2). A ${\msr G}$-module generated by a
singular vector $v$ of weight $\lmd$ is called a {\it highest-weight
module}.\index{highest-weight module} In this case, $v$ is called a
{\it highest weight vector}\index{highest weight vector} with
highest weight $\lmd$. In particular,  an irreducible module
containing a singular vector is a highest weight module. The module
$V$ is called a {\it weight module}\index{weight module} if
$$V=\sum_{\lmd\in H^\ast}V_\lmd.\eqno(5.2.6)$$
A highest-weight module is naturally a weight module. One of
fundamental problems in the representation theory of Lie algebra is
to investigate weight modules. For instance, an interesting but
difficult problem is to find all the singular vectors in a weight
module.

Observe that
$${\msr B}=H+{\msr G}_+\eqno(5.2.7)$$
is a Lie subalgebra of ${\msr G}$. For $\lmd\in H^\ast$, we define a
one-dimensional ${\msr B}$-module $\mbb{C}v_\lmd$ by
$${\msr
G}_+(v_\lmd)=\{0\},\;\;h(v_\lmd)=\lmd(h)v_\lmd\qquad\for\;\;h\in
H.\eqno(5.2.8)$$ Since ${\msr G}={\msr G}_-\oplus {\msr B}$, we have
$U({\msr G})=U({\msr G}_-)U({\msr B})$ by PBW Theorem. Now we form
an {\it induced $U({\msr G})$-module}
$$M(\lmd)=U({\msr G})\otimes_{U({\msr B})}\mbb{C}v_\lmd\cong
U({\msr G}_-)\otimes_{\mbb{C}}\mbb{C}v_\lmd\;\;(\mbox{as vector
spaces}),\eqno(5.2.9)$$ which is also a module of the Lie algebra
${\msr G}$. In fact,
$$\xi_1(u\otimes v_\lmd)=\xi_1u\otimes v_\lmd,\;\;\xi_2(u\otimes
v_\lmd)=[\xi_2,u]\otimes v_\lmd+u\otimes
\xi_2(v_\lmd)\eqno(5.2.10)$$ for $\xi_1\in{\msr G}_-,\;\xi_2\in{\msr
B}$ and $u\in U({\msr G}_-)$. Identify $1\otimes v_\lmd$ with
$v_\lmd$. The module $M(\lmd)$ is called a {\it Verma
module}.\index{Verma module} Obviously, $M(\lmd)$ is a highest
weight module with $v_\lmd$ as a highest-weight vector and $\lmd$ as
the highest weight.

Note that
$$U({\msr G}_-)=\bigoplus_{\mu\in H^\ast}U({\msr
G}_-)_\mu,\qquad U({\msr G}_-)_\mu=\{u\in U({\msr G}_-)\mid
[h,u]=\mu(h)u\;\for\;h\in H\}.\eqno(5.2.11)$$ Then
$$M(\lmd)_{\lmd+\mu}=\{w\in M(\lmd)\mid
h(w)=(\lmd+\mu)(h)w\;\for\;h\in H\}=U({\msr G}_-)_\mu
v_\lmd.\eqno(5.2.12)$$ Assume
$$U({\msr G}_-)_0=\mbb{C}1.\eqno(5.2.13)$$
Then
$$M(\lmd)_{\lmd}=\mbb{C}v_\lmd.\eqno(5.2.14)$$
Since any proper submodule $U$ of $M(\lmd)$ must be a weight module,
we have
$$U\subset \sum_{0\neq\mu\in
H^\ast}M(\lmd)_{\lmd+\mu}.\eqno(5.2.15)$$ Thus the sum  of all
proper submodules of $M(\lmd)$ is the unique maximal proper
submodule $N(\lmd)$ of $M(\lmd)$. The quotient module
$$V(\lmd)=M(\lmd)/N(\lmd)\eqno(5.2.16)$$
is an irreducible highest-weight ${\msr G}$-module.\psp

{\bf Fundamental Problem}: Determine $N(\lmd)$ and the formula of
finding the dimensions of weight subspaces of $V(\lmd)$ (so-called
``character of $V(\lmd)$").\psp

Let ${\msr G}$ be a finite-dimensional semisimple Lie algebra and
let $H$ be a CSA of ${\msr G}$ with root system $\Phi$. Fix a base
$\Pi=\{\al_1,\al_2,...,\al_n\}$ of $\Phi$, and so we have positive
and negative roots. Set
$${\msr G}_\pm=\sum_{\al\in\Phi^\pm}{\msr G}_\al.\eqno(5.2.17)$$
Then ${\msr G}_\pm$ are nilpotent subalgebras of ${\msr G}$ and
(5.2.1) holds. Moreover, ${\msr B}={\msr B}(\Pi)$ (cf. (5.2.7)) is
the standard Borel subalgebra.

Let $V$ be any finite-dimensional irreducible ${\msr G}$-module. By
Lie's Theorem, there exists a common eigenvector $v$ of ${\msr B}$.
Since ${\msr G}_+=[{\msr B},{\msr B}]$, we have ${\msr
G}_+(v)=\{0\}$; that is, $v$ is a singular vector. Hence
$$V=U({\msr G})(v)=U({\msr G}_-)(v).\eqno(5.2.18)$$
Let $\lmd$ be the weight of $v$. Define
$$\tau(w\otimes v_\lmd)=wv\qquad \for\;\;w\in U({\msr
G}_-).\eqno(5.2.19)$$ The map $\tau$ is a Lie algebra module
epimorphism from $M(\lmd)$ to $V$ by (5.2.9). By (5.2.17), (5.2.13)
holds. Now $\kn\tau\subset N(\lmd)$. Furthermore,
$$M(\lmd)/\kn\tau\cong V\eqno(5.2.20)$$
is irreducible, which implies that $\kn\tau$ is a maximal proper
submodule of $M(\lmd)$. Thus $N(\lmd)=\kn\tau$, or equivalently,
$$V\cong V(\lmd).\eqno(5.2.21)$$

Take $0\neq \xi_i\in{\msr G}_{\al_i},\;\zeta_i\in{\msr G}_{-\al_i}$
such that
$$[\xi_i,\zeta_i]=h_{\al_i}=h_i,\;\;[h_i,\xi_i]=2\xi_i\qquad\for\;\;i\in\ol{1,n}.\eqno(5.2.22)$$
Then
$${\msr S}_i=\mbb{C}\xi_i+\mbb{C}h_i+\mbb{C}\zeta_i\eqno(5.2.23)$$
forms a Lie subalgebra isomorphic to $sl(2,\mbb{C})$. By the proof
of Theorem 2.6.1, $v$ generates a finite-dimensional irreducible
${\msr S}_i$-module. Thus
$$0\leq
\lmd(h_i)=\la\lmd,\al_i\ra\in\mbb{Z}\qquad\for\;\;i\in\ol{1,n};\eqno(5.2.24)$$
that is, $\lmd$ is a dominant integral weight. Therefore, we
obtain:\psp

{\bf Lemma 5.2.1}. {\it Any finite-dimensional irreducible ${\msr
G}$-module is of the form $V(\lmd)$ for some  dominant integral
weight $\lmd$.}\psp

 Now given a dominant weight $\lmd$, we want to
determine $N(\lmd)$ and prove $V(\lmd)$ is a finite-dimensional
irreducible ${\msr G}$-module.\psp

{\bf Lemma 5.2.2}. {\it The elements
$$v_{i,\lmd}=\zeta_i^{\la\lmd,\al_i\ra+1}v_\lmd\qquad\for\;\;i\in\ol{1,n}\eqno(5.2.25)$$
are singular vectors in $M(\lmd)$.}

{\it Proof}. Since ${\msr G}_+$ is generated by
$\{\xi_1,...,\xi_n\}$, we only need to prove
$$\xi_j(v_{i,\lmd})=0\qquad\for\;\;j\in\ol{1,n},\eqno(5.2.26)$$
which is obvious if $i\neq j$ because $\xi_j\zeta_i=\zeta_i\xi_j$ in
$U({\msr G})$. Moreover,
\begin{eqnarray*}\qquad\xi_i(v_{i,\lmd})&=&\sum_{s=0}^{\la\lmd,\al_i\ra}\zeta_i^{\la\lmd,\al_i\ra-s}[\xi_i,\zeta_i]
\zeta_i^sv_\lmd=\sum_{s=0}^{\la\lmd,\al_i\ra}\zeta_i^{\la\lmd,\al_i\ra-s}h_i
\zeta_i^sv_\lmd\\
&=&\sum_{s=0}^{\la\lmd,\al_i\ra}(\la\lmd,\al_i\ra-2s)\zeta_i^{\la\lmd,\al_i\ra-s}
\zeta_i^sv_\lmd=0.\qquad\Box\hspace{4.3cm}(5.2.27)\end{eqnarray*}
\vspace{0.1cm}

Set
$$N'=\sum_{i=1}^nU({\msr G})v_{i,\lmd}=\sum_{i=1}^nU({\msr
G}_-)v_{i,\lmd}.\eqno(5.2.28)$$ Then $N'$ forms a ${\msr
G}$-submodule of $M(\lmd)$. Form a quotient ${\msr G}$-module:
$$V'=M(\lmd)/N'.\eqno(5.2.29)$$
We can identify $v_\lmd$ with its image in $V'$. Given
$\al,\be\in\Phi$. Assume that $q$ is the largest nonnegative integer
such that $\be+q\al\in\Phi$. By Lemmas 3.1.1 and 3.1.2,
$$0\leq \la\be+q\al,\al\ra\leq 3.\eqno(5.2.30)$$
Moreover, by Lemma 3.1.3,
$$\be+q\al-(\la\be+q\al,\al\ra+1)\al\not\in\Phi\lra\be+(q-4)\al\not\in\Phi\lra
q\leq 3.\eqno(5.2.31)$$ Thus
$$(\ad u_\al)^4(u_\be)=0\qquad\for\;\;u_\al\in{\msr
G}_\al,\;u_\be\in{\msr G}_\be.\eqno(5.2.32)$$ Hence
\begin{eqnarray*}&
&\xi_i^{4r}(\zeta_{\be_1}\zeta_{\be_2}\cdots\zeta_{\be_r}v_\lmd)
=\sum_{m_1+\cdots+m_{r+1}=4r}\frac{(4r)!}{m_1!m_2!\cdots m_{r+1}!}\\
& &(\ad\xi_i)^{m_1}
(\zeta_{\be_1})(\ad\xi_i)^{m_2}(\zeta_{\be_2})\cdots(\ad\xi_i)^{m_r}
(\zeta_{\be_r})\xi_i^{m_{r+1}}v_\lmd=0\hspace{3.9cm}(5.2.33)\end{eqnarray*}
and
\begin{eqnarray*}&
&\zeta_i^{4r+\la\lmd,\al_i\ra}(\zeta_{\be_1}\zeta_{\be_2}\cdots\zeta_{\be_r}v_\lmd)=
\sum_{m_1+\cdots+m_{r+1}=4r+\la\lmd,\al_i\ra}\frac{(4r+\la\lmd,\al_i\ra)!}{m_1!m_2!\cdots
m_{r+1}!}\\ & &
(\ad\zeta_i)^{m_1}(\zeta_{\be_1})(\ad\zeta_i)^{m_2}(\zeta_{\be_2})\cdots(\ad\zeta_i)^{m_r}
(\zeta_{\be_r})
\zeta_i^{m_{r+1}}v_\lmd=0\hspace{3.8cm}(5.2.34)\end{eqnarray*} for
$i\in\ol{1,n}$ and  $\zeta_{\be_j}\in{\msr G}_{\be_i}$ with
$\be_j\in\Phi^-$. Therefore, $\{\xi_i,\zeta_i\mid i\in\ol{1,n}\}$
are locally nilpotent on $V'$. So
$$\hat{\tau}_i=e^{\xi_i}e^{-\zeta_i}e^{\xi_i}\qquad\for\;\;i\in\ol{1,n}\eqno(5.2.35)$$
are linear automorphisms of $V'$.

Recall
$$V'_\mu=\{w\in V'\mid h(w)=\mu(h)w\;\for\;h\in
H\}\qquad\for\;\;\mu\in H^\ast.\eqno(5.2.36)$$ Set
$$\mbb{S}(\lmd)=\{\mu\in H^\ast\mid V'_\mu\neq 0\},\eqno(5.2.37)$$
the set of weights of $V'$. Then
$$\mu\lhd \lmd\qquad\for\;\;\mu\in\mbb{S}(\lmd)\eqno(5.2.38)$$
because $V'=U({\msr G}_-)v_\lmd$. Moreover, as operators on $V'$,
$$
\hat{\tau}_ih=\hat{\tau}_ih\hat{\tau}_i^{-1}\hat{\tau}_i=
e^{\xi_i}e^{-\zeta_i}e^{\xi_i}h
e^{-\xi_i}e^{\zeta_i}e^{-\xi_i}\hat{\tau}_i
=[e^{\ad\xi_i}e^{-\ad\zeta_i}e^{\ad\xi_i}(h)]\hat{\tau}_i=\tau_{\al_i}(h)\hat{\tau}_i\eqno(5.2.39)$$
by (1.4.31) (cf. (4.1.24) and (4.1.25)). So
$$\hat{\tau}_i(V'_{\mu})=V'_{\sgm_{\al_i}(\mu)}\qquad\for\;\;\mu\in
H^\ast\eqno(5.2.40)$$ by (5.2.39). For each $\sgm\in{\msr W}$, we
pick an expression
$\sgm=\sgm_{\al_{i_1}}\sgm_{\al_{i_2}}\cdots\sgm_{\al_{i_r}}$ and
define
$$\hat{\tau}_\sgm=\hat{\tau}_{i_1}\hat{\tau}_{i_2}\cdots\hat{\tau}_{i_r}.\eqno(5.2.41)$$
By (5.2.40), we get
$$\hat{\tau}_\sgm(V'_{\mu})=V'_{\sgm(\mu)}\qquad\for\;\;\mu\in
\mbb{S}(\lmd).\eqno(5.2.42)$$ So
$\sgm(\mbb{S}(\lmd))=\mbb{S}(\lmd)$.

Recall that $\Lmd^+$ denotes the set of dominant integral weights.
Set $$\mbb{S}^+(\lmd)=\mbb{S}(\lmd)\bigcap \Lmd^+.\eqno(5.2.43)$$
According to Lemma 3.4.1 and (5.2.38), $\mbb{S}^+(\lmd)$ is a finite
set. Moreover, any element in $\mbb{S}(\lmd)$ is conjugated to an
element in $\mbb{S}^+(\lmd)$ under the Weyl group ${\msr W}$ by
Lemma 3.2.13. Thus
$$|\mbb{S}(\lmd)|\leq |{\msr W}||\mbb{S}^+(\lmd)|<\infty.\eqno(5.2.44)$$

Given $0\lhd\gm=\sum_{i=1}^nk_i\al_i$, we define ${\msr P}(\gm)$ to
be the number of distinct sets of nonnegative integers
$\{l_\al\mid\al\in\Phi^+\}$ such that
$\sum_{\al\in\Phi^+}l_\al\al=\gm$. Then
$$\dim M(\lmd)_{\mu}={\msr P}(\lmd-\mu)\qquad\for\;\;\lmd\succ\mu\in
H^\ast\eqno(5.2.45)$$ (exercise). Hence
$$\dim V'_\mu\leq \dim
M(\lmd)_{\mu}={\msr
P}(\lmd-\mu)\qquad\for\;\;\mu\in\Pi(\lmd).\eqno(5.2.46)$$
Expressions (5.2.44) and (5.2.46) imply that $V'$ is
finite-dimensional. By Weyl's Theorem, $V'$ is a direct sum of
irreducible submodules, which are weight modules. Since
$$\dim V'_\lmd=\dim\mbb{C}v_\lmd=1,\eqno(5.2.47)$$
$v_\lmd$ must be in one of irreducible summands, say $U$. But
$$V'=U({\msr G})v_\lmd\subset U.\eqno(5.2.48)$$
Thus $V'=U$ is a finite-dimensional irreducible ${\msr G}$-module.
In summary, we obtain:\psp

{\bf Theorem 5.2.3}. {\it The set $\{V(\lmd)\mid \lmd\in\Lmd^+\}$
are all finite-dimensional irreducible ${\msr G}$-modules. Moreover,
$$N(\lmd)=\sum_{i=1}^nU({\msr
G}_-)v_{i,\lmd}.\eqno(5.2.49)$$ The set $\mbb{S}(\lmd)$ is a
saturated set of weights.}\psp

Suppose that
$${\msr G}={\msr G}_1\oplus{\msr G}_2\oplus\cdots\oplus {\msr
G}_s\eqno(5.2.50)$$ is a direct sum of simple ideals. Take a CSA
$H_i$ of ${\msr G}_i$ with root system $\Phi_i$ for $i\in\ol{1,s}$.
Then
$$H=H_1\oplus H_2\oplus\cdots\oplus H_s\eqno(5.2.51)$$
is a CSA of ${\msr G}$ with root system
$$\Phi=\bigcup_{i=1}^s\Phi_i.\eqno(5.2.52)$$
Pick a base $\Pi_i$ for $\Phi_i$ for $i\in\ol{1,s}$. The set
$$\Pi=\bigcup_{i=1}^s\Pi_i\eqno(5.2.53)$$
is a base for $\Phi$. Moreover,
$$\Phi^\pm=\bigcup_{i=1}^s\Phi_i^\pm.\eqno(5.2.54)$$

For $\lmd_i\in H^\ast_i$, denote
$$\lmd=\bigoplus_{i=1}^s\lmd_i\in H^\ast.\eqno(5.2.55)$$
The highest weight vector $v_\lmd$ for ${\msr G}$ is equivalent to
$$v_\lmd=v_{\lmd_1}\otimes v_{\lmd_2}\otimes\cdots\otimes
v_{\lmd_s},\eqno(5.2.56)$$ where $v_{\lmd_i}$ are the highest weight
vectors of ${\msr G}_i$. By PBW Theorem,
$$U({\msr G})\cong U({\msr G}_1)\otimes_{\mbb{C}}U({\msr G}_2)
\otimes_{\mbb{C}}\cdots \otimes_{\mbb{C}}U({\msr
G}_s).\eqno(5.2.57)$$ In particular, the Verma module:
$$M(\lmd)\cong M(\lmd_1)\otimes_{\mbb{C}}
M(\lmd_2)\otimes_{\mbb{C}}\cdots
\otimes_{\mbb{C}}M(\lmd_s)\eqno(5.2.58)$$ with the action:
$$(\sum_{i=1}^su_i)(w_1\otimes w_2\otimes\cdots\otimes
w_s)=\sum_{i=1}^sw_1\otimes\cdots\otimes w_{i-1}\otimes
u_i(w_i)\otimes w_{i+1}\otimes \cdots\otimes w_s\eqno(5.2.59)$$ for
$$u_i\in{\msr G}_i,\;\;w_i\in
M(\lmd_i)\qquad\for\;\;i\in\ol{1,s}.\eqno(5.2.60)$$ Moreover, it can
be proved that the maximal proper submodule:
$$N(\lmd)\cong \sum_{i=1}^kM(\lmd_1)\otimes_{\mbb{C}}\cdots
M(\lmd_{i-1})\otimes_{\mbb{C}}N(\lmd_i)\otimes_{\mbb{C}}
M(\lmd_{i+1})\otimes_{\mbb{C}}\cdots\otimes_{\mbb{C}}M(\lmd_s).\eqno(5.2.61)$$
Thus the irreducible highest weight module: $$ V(\lmd)\cong
V(\lmd_1)\otimes_{\mbb{C}}
V(\lmd_2)\otimes_{\mbb{C}}\cdots\otimes_{\mbb{C}}V(\lmd_s),\eqno(5.2.62)$$
where $V(\lmd_i)$ is the irreducible highest weight ${\msr
G}_i$-module with the highest weight $\lmd_i$.\psp

{\bf Remark 5.2.4}. Following Lepowsky [Lj], the Verma module can be
generalized as follows. Let
$$\msr G=\msr G_-\oplus G_0\oplus G_+\eqno(5.2.63)$$be a Lie algebra
that is a direct sum of its subalgebras $\msr G_\pm$ and $\msr G_0$
such that
$$[\msr G_-,\msr G_+]\subset \msr G_0,\qquad [\msr G_0,\msr G_\pm]\subset \msr
G_\pm.\eqno(5.2.64)$$ Denote
$$\msr B=\msr G_0+\msr G_+.\eqno(5.2.65)$$
For any $\msr G_0$-module $V$, we extend it to a $\msr B$-module by
letting
$$\msr G_+(V)=\{0\}.\eqno(5.2.66)$$Then we have the induced $\msr G$-module
$$M(V)=U(\msr G)\otimes_{U(\msr B)}V\cong U(\msr G_-)\otimes_{\mbb
C}V\;\;(\mbox{as vector spaces}),\eqno(5.2.67)$$ which is called a
{\it generalized Verma module}.\index{generalized Verma module}

\section{Formal Characters}

In this section, we introduce the notion of formal character of a
weight module of a finite-dimensional semisimple Lie algebra.

Let ${\msr G}$ be a finite-dimensional semisimple Lie algebra and
let $H$ be a CSA of ${\msr G}$ with root system $\Phi$. Denote by
${\msr F}(H^\ast)$  the set of $\mbb{C}$-valued functions on
$H^\ast$, which becomes a vector space with the linear operation:
$$(af+bg)(\lmd)=af(\lmd)+bg(\lmd)\qquad\for\;\;f,g\in{\msr
F}(H^\ast),\;a,b\in\mbb{C}.\eqno(5.3.1)$$ Set
$$\mbox{supp}\:f=\{\lmd\in H^\ast\mid f(\lmd)\neq
0\}\qquad\for\;\;f\in{\msr F}(H^\ast).\eqno(5.3.2)$$ Given
$\mu_1,...,\mu_r\in H^\ast$, we denote
$$S(\mu_1,...,\mu_r)=\bigcup_{i=1}^r(\mu_i-\sum_{s=1}^n\mbb{N}\al_s).\eqno(5.3.3)$$
Define
$${\msr X}(H^\ast)=\{f\in{\msr F}(H^\ast)\mid\mbox{supp}\:f\subset
S(\mu_1,...,\mu_r)\;\mbox{for some}\;\mu_1,...,\mu_r\in
H^\ast\}.\eqno(5.3.4)$$

It can be verified that ${\msr X}(H^\ast)$ is a subspace of ${\msr
F}(H^\ast)$ (exercise). Moreover, we define an algebraic operation
$\ast$ ({\it convolution}) on ${\msr X}(H^\ast)$ by:
$$(f\ast g)(\lmd)=\sum_{\mu_1,\mu_2\in
H^\ast;\;\mu_1+\mu_2=\lmd}f(\mu_1)g(\mu_2)\qquad\for\;\;f,g\in {\msr
X}(H^\ast),\;\lmd\in H^\ast.\eqno(5.3.5)$$ The above equation makes
sense due to (5.3.4). The convolution $\ast$ is commutative and
associative. Given $\lmd\in H^\ast$, we define $\ves_\lmd\in{\msr
X}(H^\ast)$ by
$$\ves_\lmd(\mu)=\left\{\begin{array}{ll}1&\mbox{if}\;\mu=\lmd\\
0&\mbox{if}\;\mu\neq\lmd.\end{array}\right.\eqno(5.3.6)$$ Then
$$\ves_\lmd\ast\ves_\mu=\ves_{\lmd+\mu}\qquad\for\;\;\lmd,\mu\in
H^\ast.\eqno(5.3.7)$$ In fact, $\ves_0$ is the identity element of
${\msr X}(H^\ast)$.

Define an action of the Weyl group ${\msr W}$ on ${\msr X}(H^\ast)$
by
$$\sgm(f)(\lmd)=f(\sgm^{-1}(\lmd))\qquad\for\;\;f\in{\msr
X}(H^\ast),\;\sgm\in{\msr W},\;\lmd\in H^\ast.\eqno(5.3.8)$$ In
particular,
\begin{eqnarray*}\hspace{1.2cm}\sgm(\ves_\lmd)(\mu)&=&\ves_\lmd(\sgm^{-1}(\mu))
=\left\{\begin{array}{ll}1&\mbox{if}\;\sgm^{-1}(\mu)=\lmd\\
0&\mbox{if}\;\sgm^{-1}(\mu)\neq\lmd\end{array}\right.\\ &=&
\left\{\begin{array}{ll}1&\mbox{if}\;\mu=\sgm(\lmd)\\
0&\mbox{if}\;\mu\neq\sgm(\lmd)\end{array}\right.=\ves_{\sgm(\lmd)}(\mu).
\hspace{5.8cm}(5.3.9)\end{eqnarray*} So
$\sgm(\ves_\lmd)=\ves_{\sgm(\lmd)}$.

 Let $V=\bigoplus_{\mu\in
H^\ast}V_\mu$ be a weight ${\msr G}$-module such that $\dim
V_\mu<\infty$ for $\mu\in H^\ast$ ($V$ may not be
finite-dimensional). Identify $e_\mu$ with $\ves_\mu$ for
$\mu\in\Lmd$ and define the {\it formal character} of $V$ by
$$\mbox{ch}_V=\sum_{\mu\in H^\ast}m_{_V}(\mu)\ves_\mu,\qquad
m_{_V}(\mu)=\dim V_\mu.\eqno(5.3.10)$$ For any $\lmd\in H^\ast$, we
have
$$\mbox{ch}_{M(\lmd)}\in{\msr X}(H^\ast),\eqno(5.3.11)$$
where $M(\lmd)$ is the Verma module (cf. (5.2.9)). Furthermore, we
define  the {\it Kostant function}:
$$p=\mbox{ch}_{M(0)}=\prod_{\al\in\Phi^+}(\sum_{i=0}^\infty\ves_{-i\al})\eqno(5.3.12)$$
and the {\it Weyl function}:
$$q=\prod_{\al\in\Phi^+}(\ves_{\al/2}-\ves_{-\al/2})=\ves_{\rho}\ast
 \prod_{\al\in\Phi^+}(\ves_0-\ves_{-\al}).\eqno(5.3.13)$$

Note that
$$\mbox{ch}_{M(\lmd)}=\ves_\lmd\ast p.\eqno(5.3.14)$$
 Moreover,
$$(\ves_0-\ves_{-\al})\ast
\sum_{i=0}^\infty\ves_{-i\al}=\sum_{i=0}^\infty\ves_{-i\al}-\sum_{i=0}^\infty\ves_{-(i+1)\al}=\ves_0.
\eqno(5.3.15)$$ Thus
$$q\ast p=\ves_{\rho}\ast
\prod_{\al\in\Phi^+}[(\ves_0-\ves_{-\al})\ast(\sum_{i=0}^\infty\ves_{-i\al})]=\ves_\rho.
\eqno(5.3.16)$$ By (5.3.14) and (5.3.16), we obtain:\psp

{\bf Lemma 5.3.1}. {\it For any $\lmd\in H^\ast$},
$$q\ast\mbox{ch}_{M(\lmd)}=\ves_{\lmd+\rho}.\eqno(5.3.17)$$

  Recall that the
weight lattice of ${\msr G}$ is defined by:
$$\Lmd=\{\lmd\in H^\ast\mid
\la\lmd,\al\ra\in\mbb{Z}\;\for\;\al\in\Phi\}.\eqno(5.3.18)$$ Take a
base $\Pi$ of $\Phi$. The set of dominant integral weight
$$\Lmd^+=\{\lmd\in \Lmd\mid
0\leq\la\lmd,\al\ra\;\for\;\al\in\Pi\}.\eqno(5.3.19)$$ If $U$ and
$V$ both are finite-dimensional ${\msr G}$-modules, we have
$$(U\otimes V)_\lmd=\sum_{\mu_1,\mu_2\in\Lmd;\;\mu_1+\mu_2=\lmd}U_{\mu_1}\otimes
V_{\mu_2}\qquad\for\;\;\lmd\in\Lmd\eqno(5.3.20)$$ by (2.2.22). Thus
we have
$$\mbox{ch}_{U\otimes
V}=\mbox{ch}_U\cdot\mbox{ch}_V.\eqno(5.3.21)$$
 For any
$\lmd\in\Lmd^+$, the space $V(\lmd)$ is the finite-dimensional
irreducible module with highest weight $\lmd$. By (5.2.41) and
(5.2.42), we have
$$m_{_{V(\lmd)}}(\sgm(\mu))=m_{_{V(\lmd)}}(\mu)\qquad\for\;\;\sgm\in{\msr
W},\;\mu\in\Lmd,\eqno(5.3.22)$$ which is equivalent to:
$$\sgm(\mbox{ch}_{V(\lmd)})=\mbox{ch}_{V(\lmd)}\qquad\for\;\;\sgm\in{\msr
W}.\eqno(5.3.23)$$ So $\mbox{ch}_{V(\lmd)}$ is ${\msr W}$-invariant.
\psp

{\bf Remark 5.3.2}. We can treat $\ves_\mu$ as the following
function on $H$:
$$\ves_\mu(h)=e^{\mu(h)}\qquad\for\;\;h\in H.\eqno(5.3.24)$$
If $V$ is a finite-dimensional ${\msr G}$-module, then
$$\mbox{ch}_V(h)=\mbox{tr}\:e^h\qquad\for\;\;h\in H.\eqno(5.3.25)$$

\section{ Weyl's Character Formula}

In this section, we continue to study the characters of
finite-dimensional irreducible modules of  a finite-dimensional
semisimple Lie algebra ${\msr G}$. In particular, we prove  Weyl's
character formula and use it to prove the tensor formula of two
finite-dimensional $\msr G$-modules.

Pick an orthonormal basis $\{h_i\in\ol{1,n}\}$ of the Cartan
subalgebra $H$ of ${\msr G}$ with respect to its Killing form
$\kappa$. For any $\al\in\Phi^+$, we take $\xi_\al\in{\msr G}_\al$
and $\zeta_\al\in{\msr G}_{-\al}$ such that
$\kappa(\xi_\al,\zeta_\al)=1$. Then
$$\omega=\sum_{i=1}^nh_i^2+\sum_{\al\in\Phi^+}(\xi_\al\zeta_\al+\zeta_\al\xi_\al)\eqno(5.4.1)$$
is a quadratic central element of the universal enveloping algebra
$U({\msr G})$; that is,
$$ u\omega=\omega u\qquad\for\;\;u\in U({\msr G})\eqno(5.4.2)$$
(exercise). We call $\omega$ a {\it Casimier element}\index{Casimier
element}. Indeed, for any representation $\phi$ of ${\msr G}$,
$\phi(\omega)=\omega_\phi$ the Casimier operator defined in (2.4.4).

Let $V$ be a highest-weight ${\msr G}$-module with a highest-weight
vector $v$ of weight $\lmd$. According to Lemma 2.5.4,
$[\xi_\al,\zeta_\al]=t_\al$ for $\al\in\Phi^+$. Thus
$$\omega(v)=[\sum_{i=1}^nh_i^2+\sum_{\al\in\Phi^+}\xi_\al\zeta_\al](v)\\
=[\sum_{i=1}^n\lmd(h_i)^2+\sum_{\al\in\Phi^+}\lmd(t_\al)]v=[(\lmd,\lmd)+2(\rho,\lmd)]v\eqno(5.4.3)$$
(cf. (3.2.15)).  Since $V=U({\msr G})(v)$, (5.4.2) yields
$$\omega|_V=[(\lmd,\lmd)+2(\rho,\lmd)]\mbox{Id}_V.\eqno(5.4.4)$$
Suppose that $u$ is a singular vector of $V$ with weight $\mu$. Then
we have
$$\omega(u)=[(\mu,\mu)+2(\rho,\mu)]u.\eqno(5.4.5)$$
Hence
$$(\mu,\mu)+2(\rho,\mu)=(\lmd,\lmd)+2(\rho,\lmd)\sim
(\rho+\mu,\rho+\mu)=(\rho+\lmd,\rho+\lmd).\eqno(5.4.6)$$

Recall the base $\Pi=\{\al_1,\al_2,...,\al_n\}$.
  Given $\lmd\in
H^\ast$, we denote
$$B_\lmd=\{\mu\in\lmd-\sum_{i=1}^n\mbb{N}\al_i\mid \mu\neq\lmd,\;
(\rho+\mu,\rho+\mu)=(\rho+\lmd,\rho+\lmd)\},\eqno(5.4.7)$$ where
$\mbb{N}$ is the additive semigroup of nonnegative integers. The
nondegeneracy of the bilinear form $(\cdot,\cdot)$ (cf. Theorem
2.1.2 and (2.6.35)) implies that $B_\lmd$ is a finite set. Moreover,
$$B_{\lmd'}\subset
B_\lmd\qquad\for\;\;\lmd'\in B_\lmd.\eqno(5.4.8)$$

Let $V$ be a highest weight ${\msr G}$-module with the highest
weight $\lmd$. We define $m_{_V}=0$ if $B_\lmd=\emptyset$ and
otherwise,
$$m_{_V}=\sum_{\mu\in B_\lmd}\dim V_\mu.\eqno(5.4.9)$$
\psp

{\bf Lemma 5.4.1}. {\it A highest weight ${\msr G}$-module $V$ with
highest weight $\lmd$ has a sequence of submodules:
$$V_0=\{0\}\subset V_1\subset V_2\subset \cdots\subset
V_r=V\eqno(5.4.10)$$ such that $V_r/V_{r-1}=V(\lmd)$ and}
$$V_i/V_{i-1}\cong V(\mu_i)\;\;\mbox{\it with
}\;\;\mu_i\in B_\lmd\qquad\mbox{\it
for}\;\;\;\;i\in\ol{1,r-1}.\eqno(5.4.11)$$

{\it Proof}. Suppose that $V$ has a proper nonzero submodule $U$,
which is also a weight module. Let $0\neq v\in U$ be a weight vector
with weight $\mu=\lmd-\sum_{i=1}^nm_i\al_i$, where
$\Pi=\{\al_1,\al_2,...,\al_n\}$. Denote
$$m=\sum_{i=1}^nm_i.\eqno(5.4.12)$$
Take $0\neq \xi_i\in{\msr G}_{\al_i}$ for $i\in\ol{1,n}$. Since
$V_\lmd\not\subset U$,
$$\xi_{i_1}\xi_{i_2}\cdots\xi_{i_m}(v)=0\qquad
\for\;\;i_r\in\ol{1,n}.\eqno(5.4.13)$$ If $v$ is not a singular
vector, then there exist $j_1,...,j_s\in\ol{1,n}$ such that
$\xi_{j_1}\xi_{j_2}\cdots\xi_{j_s}(v)\in U$ is a singular vector by
(5.4.14). So $U$ has a singular vector $u$ with weight $\mu\in
B_\lmd$.

If $m_{_V}=0$, then $V$ has a unique singular vector up to a scalar.
Thus $V=V(\lmd)$ is irreducible and the lemma holds. Suppose that it
holds for $V$ with $m_{_V}<k$. Assume that $m_{_V}=k$ and $V$ is not
irreducible. Take a singular vector $u$ of $V$ with weight $\mu\in
B_\lmd$. Set
$$W=U({\msr G}_-)u.\eqno(5.4.14)$$
Now both $V/W$ and $W$ are highest weight modules with
$m_{_{V/W}},m_{_W}<k$. Moreover, $\mu\in B_\lmd$. By assumption, $W$
has a sequence of submodules
$$W_0=\{0\}\subset W_1\subset \cdots \subset W_{r_1}=W\eqno(5.4.15)$$
such that $W_i/W_{i-1}\cong V(\mu_i)$ with $\mu_i\in\ B_\mu$ for
$i\in\ol{1,r_1-1}$ and $W_{r_1}/W_{r_1-1}\cong V(\mu)$. Furthermore,
$V$ contains a sequence of submodules
$$V_0'=W\subset V_1'\subset\cdots\subset V'_{r_2}=V\eqno(5.4.16)$$
such that
$$V_{r_2}'/V_{r_2-1}'\cong (V_{r_2}'/W)/(V_{r_2-1}'/W)\cong
V(\lmd)\eqno(5.4.17)$$ and
$$V_j'/V_{j-1}'\cong (V_j'/W)/(V_{j-1}'/W)\cong
V(\mu'_j)\eqno(5.4.18)$$ with $\mu'_j\in B_\lmd$ for
$j\in\ol{1,r_2}$. Now the sequence:
$$W_0=\{0\}\subset W_1\subset \cdots \subset W_{r_1}=W\subset V_1'\subset\cdots\subset
V'_{r_2}=V\eqno(5.4.19)$$ is the required.$\qquad\Box$\psp

By the above lemma,
$$\mbox{ch}_V=\sum_{i=1}^r\mbox{ch}_{V_i/V_{i-1}}=\mbox{ch}_{V(\lmd)}+\sum_{i=1}^{r-1}
\mbox{ch}_{V(\mu_i)}\eqno(5.4.20)$$ with $\mu_i\in B_\lmd$. In
particular,
$$\mbox{ch}_{M(\lmd)}=\mbox{ch}_{V(\lmd)}+\sum_{\mu\in B_\lmd}c(\lmd,\mu)
\mbox{ch}_{V(\mu)},\qquad c(\lmd,\mu)\in\mbb{N},\eqno(5.4.21)$$ for
any $\lmd\in H^\ast$.

Fix $\lmd\in H^\ast$. If $\mu\in B_\lmd$ is minimal with respect to
``$\lhd$" in (3.4.49), then $B_\mu=\emptyset$. So $M(\mu)=V(\mu)$,
which implies $\mbox{ch}_{M(\mu)}=\mbox{ch}_{V(\mu)}$. If $\mu$ is
not minimal, we have
$$\mbox{ch}_{V(\mu)}=\mbox{ch}_{M(\mu)}-\sum_{\mu'\in B_\mu}c(\mu,\mu')
\mbox{ch}_{V(\mu')}.\eqno(5.4.22)$$ By induction,
$$\mbox{ch}_{V(\mu)}=\mbox{ch}_{M(\mu)}+\sum_{\mu'\in B_\mu}d(\mu,\mu')
\mbox{ch}_{M(\mu')},\qquad d(\mu,\mu')\in\mbb{Z}.\eqno(5.4.23)$$ In
particular,
$$\mbox{ch}_{V(\lmd)}=\mbox{ch}_{M(\lmd)}+\sum_{\mu\in B_\lmd}a(\mu)
\mbox{ch}_{M(\mu)},\qquad a(\mu)\in\mbb{Z}.\eqno(5.4.24)$$ Thus
$$q\ast\mbox{ch}_{V(\lmd)}=\ves_{\lmd+\rho}+\sum_{\mu\in B_\lmd}a(\mu)
\ves_{\mu+\rho}\eqno(5.4.25)$$ by (5.3.17).

Note that Lemma 3.2.10 implies
$$\sgm(q)=(-1)^{\ell(\sgm)}q\qquad\for\;\;\sgm\in{\msr
W}.\eqno(5.4.26)$$  In fact,
$$\mbox{sn}(\sgm)=\mbox{the determinant
of}\;\sgm=(-1)^{\ell(\sgm)}.\eqno(5.4.27)$$ So
$$\sgm(q)=\mbox{sn}(\sgm)q\qquad\for\;\;\sgm\in{\msr
W}.\eqno(5.4.28)$$

Observe that ${\msr W}$ is an automorphism subgroup of the algebra
$({\msr X}(H^\ast),\ast)$. Suppose $\lmd\in\Lmd^+$. Then
$\mbox{ch}_{V(\lmd)}$ is ${\msr W}$-invariant. Applying
$\sgm\in{\msr W}$ to (5.4.25), we obtain
\begin{eqnarray*}\hspace{2cm}& &\mbox{sn}(\sgm)(\ves_{\lmd+\rho}+\sum_{\mu\in B_\lmd}a(\mu)
\ves_{\mu+\rho})\\ &
=&\sgm(q)\ast\sgm(\mbox{ch}_{V(\lmd)})=\sgm(q\ast\mbox{ch}_{V(\lmd)})\\
&=& \ves_{\sgm(\lmd+\rho)}+\sum_{\mu\in B_\lmd}a(\mu)
\ves_{\sgm(\mu+\rho)}.\hspace{7cm}(5.4.29)\end{eqnarray*} On the
other hand, if $\mu \in B_\lmd$ and $\rho+\mu\in \Lmd^+$, then
$$
(\rho+\lmd,\rho+\lmd)-(\rho+\mu,\rho+\mu)=(\rho,\lmd-\mu)+(\lmd,\lmd-\mu)+(\rho+\mu,\lmd-\mu)>0,\eqno(5.4.30)$$
which contradicts (5.4.7). Thus
$$\{\rho+\mu\mid \mu\in
B_\lmd\}\bigcap\Lmd^+=\emptyset.\eqno(5.4.31)$$ For any $\mu'\in
B_\lmd$, there exists $\sgm\in {\msr W}$ such that
$\sgm(\rho+\mu')\in\Lmd^+$ by Lemma 3.2.13. Then (5.4.29) and
(5.4.31) imply
$$\sgm(\rho+\mu')=\rho+\lmd\;\;\mbox{and}\;\;a(\mu')=\mbox{sn}(\sgm).\eqno(5.4.32)$$
Since $\lmd+\rho$ is strongly dominant (regular),
$\{\sgm(\lmd+\rho)\mid\sgm\in{\msr W}\}$ are distinct by Theorem
3.2.9 (e).  Observe that $\mbox{sn}(\sgm)=\mbox{sn}(\sgm^{-1})$.
Therefore, we have:\psp

{\bf Theorem 5.4.2 (Weyl)}.\index{Weyl's Character Formula} {\it For
$\lmd\in\Lmd^+$,}
$$q\ast\mbox{ch}_{V(\lmd)}=\sum_{\sgm\in{\msr
W}}\mbox{sn}(\sgm)\ves_{\sgm(\lmd+\rho)}.\eqno(5.4.33)$$
\vspace{0.1cm}

Taking $\lmd=0$ in the above equation, we get the {\it Weyl
denominator identity}:
$$\prod_{\al\in\Phi^+}(\ves_{\al/2}-\ves_{-\al/2})=
\sum_{\sgm\in{\msr
W}}\mbox{sn}(\sgm)\ves_{\sgm(\rho)}.\eqno(5.4.34)$$ \vspace{0.1cm}

 {\bf Example 5.4.1}. Let
$n>1$ be a integer and let ${\msr E}$ be the Euclidean space with a
basis $\{\ves_1,\ves_2,...,\ves_n\}$ such that
$(\ves_i,\ves_j)=\dlt_{i,j}$. The root system of $sl(n,\mbb{C})$ and
a base are as follows:
$$\Phi=\{\ves_i-\ves_j\mid i,j\in\ol{1,n},\;i\neq
j\},\;\;\Pi=\{\ves_i-\ves_{i+1}\mid i\in\ol{1,n-1}\}.\eqno(5.4.35)$$
Denote
$$\Lmd^+_1=\{\sum_{i=1}^{n-1}m_i\ves_i\mid
m_i\in\mbb{N},\;m_1\geq m_2\geq\cdots\geq m_{n-1}\}.\eqno(5.4.36)$$
Define
$$|\mu|=\sum_{i=1}^{n-1}m_i\qquad\for\;\;\mu=\sum_{i=1}^{n-1}m_i\ves_i\in\Lmd_1^+\eqno(5.4.37)$$
and
$${\bf 1}=\sum_{i=1}^n\ves_i.\eqno(5.4.38)$$
Then
$$\Lmd^+=\left\{\mu-\frac{|\mu|}{n}{\bf 1}\mid
\mu\in\Lmd_1^+\right\}\eqno(5.4.39)$$ (exercise, cf.
(3.4.11)-(3.4.14)). In particular,
$$\rho=\frac{1}{2}\sum_{1\leq i\leq j\leq
n}(\ves_i-\ves_j)=\frac{1}{2}\sum_{i=1}^n(n+1-2i)\ves_i=\sum_{i=1}^{n-1}(n-i)\ves_i
-\frac{n-1}{2}{\bf 1}.\eqno(5.4.40)$$ As we showed in Example 3.2.1,
the Weyl group of $sl(n,\mbb{C})$ isomorphic to the permutation
group $S_n$ on $\{1,2,...,n\}$. Note
$$\sgm(\ves_i)=\ves_{\sgm(i)}\qquad\for\;\;i\in\ol{1,n},\;\sgm\in
S_n.\eqno(5.4.41)$$

Denote
$$x_i=\ves_{\ves_i}\qquad\for\;\;i\in\ol{1,n}.\eqno(5.4.42)$$
Moreover, we write
$$x^\mu=x_1^{a_1}x_2^{a_2}\cdots
x_n^{a_n}\qquad\for\;\;\mu=\sum_{i=1}^na_i\ves_i\in H^\ast
.\eqno(5.4.43)$$ So we can treat
$$\ves_\mu= x^\mu.\eqno(5.4.44)$$
Thus
$$\sgm(x_1^{a_1}x_2^{a_2}\cdots
x_n^{a_n})=x_{\sgm(1)}^{a_1}x_{\sgm(2)}^{a_2}\cdots
x_{\sgm(n)}^{a_n}\qquad\for\;\;\sgm\in
S_n,\;a_i\in\mbb{C}.\eqno(5.4.45)$$ In particular,
$$\sgm(x^{a{\bf 1}})=\sgm(x_1^ax_2^a\cdots x_n^a)=x^{a{\bf 1}}\qquad\for\;\;\sgm\in
S_n.\eqno(5.4.46)$$

Now \begin{eqnarray*}q&=&\prod_{1\leq i< j\leq
n}(x_i^{1/2}x_j^{-1/2}-x_i^{-1/2}x_j^{1/2})=[\prod_{1\leq i< j\leq
n}(x_ix_j)^{-1/2}][\prod_{1\leq i< j\leq n}(x_i-x_j)]\\ &
=&x^{-(n-1){\bf 1}/2}\prod_{1\leq i< j\leq
n}(x_i-x_j)\hspace{8.8cm}(5.4.47)\end{eqnarray*} and
\begin{eqnarray*}\qquad  \sum_{\sgm\in{\msr
W}}\mbox{sn}(\sgm)\ves_{\sgm(\rho)}&=&\sum_{\sgm\in
S_n}\mbox{sn}(\sgm)\sgm(x^{-(n-1){\bf 1}/2}x_1^{n-1}x_2^{n-2}\cdots
x_{n-1})\\ &=& x^{-(n-1){\bf 1}/2}\sum_{\sgm\in
S_n}\mbox{sn}(\sgm)x_{\sgm(1)}^{n-1}x_{\sgm(2)}^{n-2}\cdots
x_{\sgm(n-1)}\\ &=&x^{-(n-1){\bf 1}/2}\left|\begin{array}{cccc}x_1^{n-1}&x_2^{n-1}&\cdots&x_n^{n-1}\\
x_1^{n-2}&x_2^{n-2}&\cdots&x_n^{n-2}\\ \vdots&\vdots&
\vdots&\vdots\\
1&1&\cdots&1\end{array}\right|\hspace{3.7cm}(5.4.48)\end{eqnarray*}
by (5.4.40). So the Weyl denominator identity (5.4.34) is equivalent
to the identity of  Vandermonde determinant:
$$\prod_{1\leq i< j\leq
n}(x_i-x_j)=\left|\begin{array}{cccc}x_1^{n-1}&x_2^{n-1}&\cdots&x_n^{n-1}\\
x_1^{n-2}&x_2^{n-2}&\cdots&x_n^{n-2}\\ \vdots&\vdots&
\vdots&\vdots\\
1&1&\cdots&1\end{array}\right|.\eqno(5.4.49)$$

For $\mu=\sum_{i=1}^{n-1}m_i\ves_i\in\Lmd_1^+$, the Schur's
symmetric polynomial
\begin{eqnarray*}s(\mu)&=&\frac{\sum_{\sgm\in
S_n}\mbox{sn}(\sgm)x_{\sgm(1)}^{m_1+n-1}x_{\sgm(2)}^{m_2+n-2}\cdots
x_{\sgm(n-1)}^{m_{n-1}}}{\prod_{1\leq i< j\leq n}(x_i-x_j)}\\
&=&\frac{1}{\prod_{1\leq i< j\leq
n}(x_i-x_j)}\left|\begin{array}{cccc}x_1^{m_1+n-1}&x_2^{m_1+
n-1}&\cdots&x_n^{m_1+n-1}\\
x_1^{m_2+n-2}&x_2^{m_2+n-2}&\cdots&x_n^{m_2+n-2}\\ \vdots&\vdots&
\vdots&\vdots\\
x_1^{m_{n-1}}&x_2^{m_{n-1}}&\cdots&x_n^{m_{n-1}}\end{array}\right|.\hspace{2.6cm}(5.4.50)\end{eqnarray*}

Let $\lmd=\mu-|\mu|{\bf 1}/n\in\Lmd^+$ with $\mu\in\Lmd_1^+$. By
(5.4.33), (5.4.40), (5.4.46), (5.4.47) and (5.4.50),
$$\mbox{ch}_{V(\lmd)}=x^{-|\mu|{\bf 1}/n}s(\mu).\eqno(5.4.51)$$
\pse

Let $\lmd',{\lmd'}'\in\Lmd^+$. Since $V(\lmd')\otimes V({\lmd'}')$
is completely reducible,
$$\mbox{ch}_{V(\lmd')}\ast
\mbox{ch}_{V({\lmd'}')}=\sum_{\lmd\in\Lmd^+}n(\lmd)\mbox{ch}_{V(\lmd)},\qquad
n(\lmd)\in\mbb{N}.\eqno(5.4.52)$$ Our goal is to find a formula for
$n(\lmd)$. Multiplying $q$ to both sides of the above equation, we
get
$$\mbox{ch}_{V(\lmd')}\ast(\sum_{\sgm\in{\msr
W}}\mbox{sn}(\sgm)\ves_{\sgm({\lmd'}'+\rho)})=\sum_{\lmd\in\Lmd^+}\sum_{\sgm\in{\msr
W}}n(\lmd)\mbox{sn}(\sgm)\ves_{\sgm(\lmd+\rho)}\eqno(5.4.53)$$ by
(5.4.33).

For any $\mu\in\Lmd$, we define
$${\msr W}(\mu)\bigcap \Lmd^+=\{\bar{\mu}\}\eqno(5.4.54)$$ by Lemma 3.2.12 and
$${\msr W}_\mu=\{\sgm\in{\msr W}\mid \sgm(\mu)=\mu\}.\eqno(5.4.55)$$
If ${\msr W}_\mu\neq\{1\}$, then $\mu$ is not regular by Theorem
3.2.9 (e). Thus there exists $\al\in\Phi$ such that $(\mu,\al)=0$.
Hence $\sgm_\al\in{\msr W}_\mu$. Note $\mbox{sn}(\sgm_\al)=-1$,
${\msr W}_\mu={\msr W}_\mu\sgm_\al$ and so
$$\sum_{\tau\in {\msr W}_\mu}\mbox{sn}(\tau)=\sum_{\tau\in {\msr W}_\mu}\mbox{sn}(\tau\sgm_\al)
=\sum_{\tau\in {\msr
W}_\mu}\mbox{sn}(\tau)\mbox{sn}(\sgm_\al)=-\sum_{\tau\in {\msr
W}_\mu}\mbox{sn}(\tau).\eqno(5.4.56)$$ Thus
$$\sum_{\tau\in {\msr W}_\mu}\mbox{sn}(\tau)=0.\eqno(5.4.57)$$
Write
$${\msr W}=\bigcup_{i=1}^r\tau_i{\msr W}_\mu,\qquad
\tau_i\tau_j^{-1}\not\in{\msr W}_\mu\;\;\mbox{if}\;i\neq
j.\eqno(5.4.58)$$ Then
$$\sum_{\sgm\in{\msr W}}\mbox{sn}(\sgm)\ves_{\sgm(\mu)}=\sum_{i=1}^r
\sum_{\tau\in{\msr W}_\mu}\mbox{sn}(\tau_i\tau)\ves_{\tau_i(\mu)}=
\sum_{i=1}^r(\sum_{\tau\in {\msr
W}_\mu}\mbox{sn}(\tau))\mbox{sn}(\tau_i)\ves_{\tau_i(\mu)}=0.\eqno(5.4.59)$$
Define
$$\psi(\mu)=\left\{\begin{array}{ll}0&\mbox{if}\;{\msr
W}_\mu\neq\{1\}\\ \mbox{sn}(\sgm)&\mbox{if}\;{\msr
W}_\mu=\{1\}\;\mbox{and}\;\sgm(\mu)=\bar{\mu}\end{array}\right.
\qquad\for\;\;\mu\in\Lmd.\eqno(5.4.60)$$

Since
$$m_{_{V(\lmd')}}(\sgm(\mu))=m_{_{V(\lmd')}}(\mu)\qquad\for\;\;\sgm\in{\msr
W},\;\mu\in\mbb{S}(\lmd'),\eqno(5.4.61)$$ we have
$$\mbox{ch}_{V(\lmd')}=\sum_{\mu\in\mbb{S}(\lmd')}m_{_{V(\lmd')}}(\mu)\ves_{\sgm(\mu)}\qquad\for\;\;\sgm\in{\msr
W}.\eqno(5.4.62)$$ So
$$\mbox{ch}_{V(\lmd')}\ast(\sum_{\sgm\in{\msr
W}}\mbox{sn}(\sgm)\ves_{\sgm({\lmd'}'+\rho)})=\sum_{\mu\in\mbb{S}(\lmd')}\sum_{\sgm\in{\msr
W}}m_{_{V(\lmd')}}(\mu)\mbox{sn}(\sgm)\ves_{\sgm({\lmd'}'+\rho+\mu)}.\eqno(5.4.63)$$
Note that ${\msr W}_{{\lmd'}'+\rho+\mu}=\{1\}$ if and only if
$\ol{{\lmd'}'+\rho+\mu}\in\Lmd^++\rho$.
 Comparing the coefficients of $\ves_{\lmd+\rho}$ in (5.4.53) and
 using (5.4.59),
we get:\psp

{\bf Theorem 5.4.3}. {\it For $\lmd',{\lmd'}'\in\Lmd^+$,

$$\mbox{ch}_{V(\lmd')}\ast
\mbox{ch}_{V({\lmd'}')}=\sum_{\mu\in\mbb{S}(\lmd'),\;\ol{{\lmd'}'+\rho+\mu}\in\Lmd^++\rho}m_{_{V(\lmd')}}(\mu)
\psi({\lmd'}'+\rho+\mu)
\mbox{ch}_{V(\ol{{\lmd'}'+\rho+\mu}-\rho)}.\eqno(5.4.64)$$} \pse

Exercise:

Calculate (5.4.64) for $sl(n,\mbb{C})$ with $\lmd'=\lmd_1$ and
${\lmd'}'=\lmd_2$.

\section{Dimensional Formula}

In this section, we  find the dimensional formula of a
finite-dimensional irreducible module of the finite-dimensional
semisimple Lie algebra ${\msr G}$.

For $\al\in\Phi^+$, we define $\ptl_\al\in\mbox{Der}\:{\msr
X}(H^\ast)$ by:
$$\ptl_\al(\ves_\lmd)=(\lmd,\al)\ves_\lmd\qquad\for\;\;\lmd\in
H^\ast,\eqno(5.5.1)$$ where $(\al,\lmd)=\kappa(t_\al,t_\lmd)$. Set
$$\ptl=\prod_{\al\in\Phi^+}\ptl_\al.\eqno(5.5.2)$$
Note that \begin{eqnarray*}\ptl(\sum_{\sgm\in{\msr
W}}\mbox{sn}(\sgm)\ves_{\sgm(\lmd+\rho)})&=&\sum_{\sgm\in{\msr
W}}\mbox{sn}(\sgm)\prod_{\al\in\Phi^+}(\sgm(\lmd+\rho),\al)\ves_{\sgm(\lmd+\rho)}
\\ &=&\sum_{\sgm\in{\msr
W}}\mbox{sn}(\sgm)\prod_{\al\in\Phi^+}(\lmd+\rho,\sgm^{-1}(\al))\ves_{\sgm(\lmd+\rho)}
\\ &=&\sum_{\sgm\in{\msr
W}}\prod_{\al\in\Phi^+}(\lmd+\rho,\al)\ves_{\sgm(\lmd+\rho)}\\
&=&\prod_{\al\in\Phi^+}(\lmd+\rho,\al)\sum_{\sgm\in{\msr
W}}\ves_{\sgm(\lmd+\rho)}\hspace{5.8cm}(5.5.3)\end{eqnarray*} by
Lemma 3.2.10.

If $f=\sum_{\mu\in H^\ast}a_\mu\ves_\mu\in{\msr X}(H^\ast)$ has only
finite number of $a_\mu\neq 0$, we define
$$\jmath(f)=\sum_{\mu\in H^\ast}a_\mu.\eqno(5.5.4)$$
If $g\in{\msr X}(H^\ast)$ has the same property as $f$, then
$$\jmath(f\ast g)=\jmath(f)\jmath(g).\eqno(5.5.5)$$
In particular,
$$\jmath[\ptl(\sum_{\sgm\in{\msr
W}}\mbox{sn}(\sgm)\ves_{\sgm(\lmd+\rho)})]=|{\msr
W}|\prod_{\al\in\Phi^+}(\lmd+\rho,\al)\eqno(5.5.6)$$ by (5.5.3).
Taking $\lmd=0$, we have
$$\jmath(\ptl(q))=|{\msr
W}|\prod_{\al\in\Phi^+}(\rho,\al)\eqno(5.5.7)$$ by (5.4.33).

On the other hand,
$q=\prod_{\be\in\Phi^+}(\ves_{\be/2}-\ves_{-\be/2})$. If $S$ is a
proper subset of $\Phi^+$, then $(\prod_{\al\in S}\ptl_\al)(q)$ is a
sum of the terms with a factor $(\ves_{\be/2}-\ves_{-\be/2})$ for
some $\be\in\Phi$. Since $\jmath(\ves_{\be/2}-\ves_{-\be/2})=0$, we
have
$$\jmath((\prod_{\al\in S}\ptl_\al)(q))=0\eqno(5.5.8)$$
by (5.5.5). This shows
$$\jmath(\ptl(q\ast\mbox{ch}_{V(\lmd)}))=\jmath(\ptl(q))
\jmath(\mbox{ch}_{V(\lmd)})=|{\msr
W}|[\prod_{\al\in\Phi^+}(\rho,\al)]\dim V(\lmd).\eqno(5.5.9)$$
Expressions (5.4.33), (5.5.6) and (5.5.9) yield:\psp

{\bf Theorem 5.5.1}. {\it For $\lmd\in\Lmd^+$,}
$$\mbox{\it
dim}\:V(\lmd)=\frac{\prod_{\al\in\Phi^+}(\lmd+\rho,\al)}{\prod_{\al\in\Phi^+}(\rho,\al)}=
\frac{\prod_{\al\in\Phi^+}\la\lmd+\rho,\al\ra}{\prod_{\al\in\Phi^+}\la\rho,\al\ra}.\eqno(5.5.10)$$
\vspace{0.1cm}

{\bf Example 5.5.1}. As in Example 3.4.1, the simple roots of
$sl(n,\mbb{C})$ are:
$$\al_i=\ves_i-\ves_{i+1}\qquad\for\;\;i\in\ol{1,n-1}.\eqno(5.5.11)$$
The positive roots:
$$\ves_i-\ves_j=\sum_{r=i}^{j-1}\al_r\qquad\for\;\;1\leq i<j\leq
n.\eqno(5.5.12)$$ For $\lmd=\sum_{i=1}^{n-1}m_i\lmd_i\in\Lmd^+$, we
have
$$\la\lmd+\rho,\ves_i-\ves_j\ra=\sum_{r=i}^{j-1}m_r+j-i.\eqno(5.5.13)$$
Thus
$$\dim V(\lmd)=\prod_{1\leq i<j\leq
n}\frac{\sum_{r=i}^{j-1}m_r+j-i}{j-i}.\eqno(5.5.14)$$

\part{Explicit Representations}

\chapter{Representations of Special Linear Algebras}

In this chapter, we give various explicit representations of special
linear Lie algebras. In Section 6.1, we present a fundamental lemma
of solving flag partial differential equations for polynomial
solutions, which was due to our work [X16]. In Section 6.2, we talk
about the canonical  bosonic oscillator representations  and
fermionic oscillator representations
 of special linear Lie algebras over their minimal natural modules and minimal orthogonal modules.
In Section 6.3, we give a general set-up of noncanonical oscillator
representations of special linear Lie algebras and Howe's result on
the representations obtained from the canonical bosonic oscillator
representations over minimal modules  by partially swapping
differential operators and multiplication operators (cf. [Hr4]).
Sections 6.4, 6.5 and 6.6 are devoted  to determining the structure
of the noncanonical oscillator representations obtained from the
canonical bosonic oscillator representations over minimal orthogonal
modules by partially swapping differential operators and
multiplication operators, which are generalizations of the classical
theorem on harmonic polynomials. The results in Sections 6.4, 6.5
and 6.6 were due to Luo and the author [LX1]. In Section 6.7, we
construct a functor from the category of $A_{n-1}$-modules to the
category of $A_n$-modules, which is related to $n$-dimensional
projective transformations. This work was due to Zhao and the author
[ZX]. As a consequence of Section 6.7,  we obtain a one-parameter
family of inhomogeneous first-order differential operator
representations of special linear Lie algebras. By partially
swapping differential operators and multiplication operators, we
obtain in Section 6.8 more general differential operator
representations. Letting these differential operators act on the
corresponding polynomial algebra and the space of
exponential-polynomial functions, we construct  multi-parameter
families of explicit infinite-dimensional irreducible
representations. These results are taken from our work [X25].

\section{Fundamental Lemma on Polynomial Solutions}

In this section, we present a fundamental lemma of solving for
polynomial solutions of linear partial differential equations. In
terms of representations, it will be useful in obtaining a basis for
irreducible oscillator representations.\psp

{\bf Lemma 6.1.1}. {\it Let ${\msr B}$ be a commutative associative
algebra and let ${\msr A}$ be a free ${\msr B}$-module
 generated  by a filtrated subspace $V=\bigcup_{r=0}^\infty V_r$
(i.e., $V_r\subset V_{r+1}$). Let $T_1$ be a linear operator on
$\msr B\oplus {\msr A}$ with a right inverse $T_1^-$ such that
$$T_1({\msr B}),\;T_1^-({\msr B})\subset{\msr B},\qquad
T_1(\zeta_1\zeta_2)=T_1(\zeta_1)\zeta_2,\qquad
T_1^-(\zeta_1\zeta_2)=T_1^-(\zeta_1)\zeta_2 \eqno(6.1.1)$$ for
$\zeta_1 \in {\msr B},\;\zeta_2\in V$, and let $T_2$ be a linear
operator on ${\msr A}$ such that $T_2(V_0)=\{0\}$,
$$ T_2(V_{r+1})\subset {\msr B}V_r,\;\;
T_2(f\zeta)=fT_2(\zeta) \qquad\for\;\; r\in\mbb{N},\;\;f\in{\msr
B},\;\zeta\in{\msr A}.\eqno(6.1.2)$$ Then we
have \begin{eqnarray*}\hspace{1cm}&&\{f\in{\msr A}\mid (T_1+T_2)(f)=0\}\\
& =&\mbox{Span}\{ \sum_{i=0}^\infty(-T_1^-T_2)^i(hg)\mid g\in
V,\;h\in {\msr B};\;T_1(h)=0\}, \hspace{3.4cm}(6.1.3)\end{eqnarray*}
where the summation is finite due to (6.1.2). Moreover, the operator
$\sum_{i=0}^\infty(-T_1^-T_2)^iT_1^-$ is a right inverse of
$T_1+T_2$}.

{\it Proof}. For $h\in {\msr B}$ such that $T_1(h)=0$ and $g\in V$,
we have
\begin{eqnarray*}& &(T_1+T_2)(\sum_{i=0}^\infty(-T_1^-T_2)^i(hg))\\
&=&T_1(hg)-\sum_{i=1}^\infty T_1[T_1^-T_2(-T_1^-T_2)^{i-1}(hg)]+
\sum_{i=0}^\infty T_2[(-T_1^-)^i(hg)]\\
&=&T_1(h)g-\sum_{i=1}^\infty (T_1T_1^-)T_2(-T_1^-T_2)^{i-1}(hg)+
\sum_{i=0}^\infty T_2(-T_1^-T_2)^i(hg)\\
&=&-\sum_{i=1}^\infty T_2(-T_1^-T_2)^{i-1}(hg)+ \sum_{i=0}^\infty
T_2(-T_1^-T_2)^i(hg)=0\hspace{4.5cm} (6.1.4)\end{eqnarray*} by
(6.1.1).  Set $V_{-1}=\{0\}$. For $k\in\mbb{N}$, we take
$\{\psi_i\mid i\in I_k\}\subset V_k$ such that
$$\{\psi_i+V_{k-1}\mid i\in I_k\}\;\;\mbox{forms a basis
of}\;\;V_k/V_{k-1},\eqno(6.1.5)$$ where $I_k$ is an index set.
 Let
$${\msr A}^{(m)}={\msr B}V_m=\sum_{s=0}^m\;\sum_{i\in I_s}{\msr
B}\psi_i.\eqno(6.1.6)$$ Obviously,
$$T_1({\msr
A}^{(m)}),\;T_1^-({\msr A}^{(m)}),\; T_2({\msr A}^{(m+1)})\subset
{\msr A}^{(m)}\qquad\for\;\;m\in\mbb{N}\eqno(6.1.7)$$ by (6.1.1) and
(6.1.2), and
$${\msr A}=\bigcup_{m=0}^\infty {\msr A}^{(m)}.\eqno(6.1.8)$$

  Suppose $\phi\in {\msr
A}^{(m)}$ such that $(T_1+T_2)(\phi)=0$. If $m=0$, then
$$\phi=\sum_{i\in I_0}h_i\psi_i,\qquad h_i\in {\msr
B}.\eqno(6.1.9)$$ Now
$$0=(T_1+T_2)(\phi)=\sum_{i\in I_0}T_1(h_i)\psi_i+\sum_{i\in
I_0}h_iT_2(\psi_i)=\sum_{i\in I_0}T_1(h_i)\psi_i,\eqno(6.1.10)$$
Since $T_1(h_i)\in {\msr B}$ by (6.1.1) and ${\msr A}$ is a free
${\msr B}$-module generated by $V$, we have $T_1(h_i)=0$ for $i\in
I_0$. Denote by ${\msr S}$ the right hand side of the equation
(6.1.3). Then
$$\phi=\sum_{i\in I_0}\;\sum_{m=0}^\infty(-T_1^-T_2)^m(h_i\psi_i)\in {\msr
S}.\eqno(6.1.11)$$

Suppose $m>0$. We write
$$\phi=\sum_{i\in I_m}h_i\psi_i+\phi',\qquad h_i\in{\msr
B},\;\phi'\in {\msr A}^{(m-1)}.\eqno(6.1.12)$$ Then
$$0=(T_1+T_2)(\phi)=\sum_{i\in
I_m}T_1(h_i)\psi_i+T_1(\phi')+T_2(\phi).\eqno(6.1.13)$$ Since
$T_1(\phi')+T_2(\phi)\in {\msr A}^{(m-1)}$, we have $T_1(h_i)=0$ for
$i\in I_m$. Now
$$\phi-\sum_{i\in
I_m}\sum_{k=0}^\infty(-T_1^-T_2)^k(h_i\psi_i)=\phi'-\sum_{i\in
I_m}\sum_{k=1}^\infty(-T_1^-T_2)^k(h_i\psi_i)\in {\msr
A}^{(m-1)}\eqno(6.1.14)$$ and (6.1.4) implies
$$(T_1+T_2)(\phi-\sum_{i\in
I_m}\sum_{k=0}^\infty(-T_1^-T_2)^k(h_i\psi_i))=0.\eqno(6.1.15)$$ By
induction on $m$, $$\phi-\sum_{i\in
I_m}\sum_{k=0}^\infty(-T_1^-T_2)^k(h_i\psi_i)\in {\msr
S}.\eqno(6.1.16)$$ Therefore, $\phi\in {\msr S}.$

For any $f\in {\msr A}$, we have:
\begin{eqnarray*} \hspace{1cm} & &(T_1+T_2)(\sum_{i=0}^\infty(-T_1^-T_2)^iT_1^-)(f)
\\ &=&f-\sum_{i=1}^\infty
T_2(-T_1^-T_2)^{i-1}T_1^-(f)+\sum_{i=0}^\infty
T_2(-T_1^-T_2)^iT_1^-(f)=f.\hspace{2.1cm}(6.1.17)\end{eqnarray*}
Thus
 the operator
$\sum_{i=0}^\infty(-T_1^-T_2)^iT_1^-$ is a right inverse of
$T_1+T_2.\qquad\Box$\psp

 We remark that the above operator $T_1$ and $T_2$ may not
 commute.

Define
$$x^\al=x_1^{\al_1}x_2^{\al_2}\cdots
x_n^{\al_n}\;\;\for\;\;\al=(\al_1,...,\al_n)\in\mbb{N}^{\:n}.\eqno(6.1.18)$$
Moreover, we denote
$$\es_i=(0,...,0,\stl{i}{1},0,...,0)\in \mbb{N}^{\:n}.\eqno(6.1.19)$$
 For each
$i\in\ol{1,n}$, we define the linear operator $\int_{(x_i)}$ on
${\msr A}$ by:
$$\int_{(x_i)}(x^\al)=\frac{x^{\al+\es_i}}{\al_i+1}\;\;\for\;\;\al\in
\mbb{N}^{\:n}.\eqno(6.1.20)$$ Furthermore, we let
$$\int_{(x_i)}^{(0)}=1,\qquad\int_{(x_i)}^{(m)}=\stl{m}{\overbrace{\int_{(x_i)}\cdots\int_{(x_i)}}}
\qquad\for\; \;0<m\in\mbb{Z}\eqno(6.1.21)$$ and denote
$$\ptl^{\al}=\ptl_{x_1}^{\al_1}\ptl_{x_2}^{\al_2}\cdots
\ptl_{x_n}^{\al_n},\;\;
\int^{(\al)}=\int_{(x_1)}^{(\al_1)}\int_{(x_2)}^{(\al_2)}\cdots
\int_{(x_n)}^{(\al_n)}\qquad\for\;\;\al\in
\mbb{N}^{\:n}.\eqno(6.1.22)$$ Obviously, $\int^{(\al)}$ is a right
inverse of $\ptl^\al$ for $\al\in \mbb{N}^{\:n}.$ We remark that
$\int^{(\al)}\ptl^\al\neq 1$ if $\al\neq 0$ due to $\ptl^\al(1)=0$.

Consider the wave equation in Riemannian space with a nontrivial
conformal group (cf. [I1]):
$$u_{tt}-u_{x_1x_1}-\sum_{i,j=2}^ng_{i,j}(x_1-t)u_{x_ix_j}=0,\eqno(6.1.23)$$
where we assume that $g_{i,j}(z)$ are one-variable polynomials.
Change variables:
$$z_0=x_1+t,\qquad z_1=x_1-t.\eqno(6.1.24)$$
Then
$$\ptl_t^2=(\ptl_{z_0}-\ptl_{z_1})^2,\qquad
\ptl_{x_1}^2=(\ptl_{z_0}+\ptl_{z_1})^2.\eqno(6.1.25)$$ So the
equation (6.1.23) changes to:
$$2\ptl_{z_0}\ptl_{z_1}+
\sum_{i,j=2}^ng_{i,j}(z_1)u_{x_ix_j}=0.\eqno(6.1.26)$$ Denote
$$T_1=2\ptl_{z_0}\ptl_{z_1},\qquad
T_2=\sum_{i,j=2}^ng_{i,j}(z_1)\ptl_{x_i}\ptl_{x_j}.\eqno(6.1.27)$$
Take $T_1^-=\frac{1}{2}\int_{(z_0)}\int_{(z_1)}$, and
$${\msr B}=\mbb{F}[z_0,z_1],\qquad V=\mbb{F}[x_2,...,x_n],\qquad
V_r=\{f\in V\mid\mbox{deg}\;f\leq r\}.\eqno(6.1.28)$$ Then the
conditions in Lemma 6.1.1 hold. Thus we have:
 \psp

{\bf Theorem 6.1.2}. {\it The space of all polynomial solutions for
the equation (6.1.23) is:
\begin{eqnarray*} \hspace{2cm}& &\mbox{Span}\:\{\sum_{m=0}^\infty(-2)^{-m}
(\sum_{i,j=2}^n\int_{(z_0)}\int_{(z_1)}g_{i,j}(z_1)\ptl_{x_i}\ptl_{x_j})^m(f_0g_0+f_1g_1)\\
& & \mid
f_0\in\mbb{F}[z_0],\;f_1\in\mbb{F}[z_1],\;g_0,g_1\in\mbb{F}[x_2,...,x_n]\}
\hspace{4cm}(6.1.29)\end{eqnarray*} with $z_0,z_1$ defined in
(6.1.24).}\psp

Let $m_1,m_2,...,m_n$ be positive integers. According to Lemma
6.1.1, the set
\begin{eqnarray*}\hspace{1.9cm}& &\{\sum_{k_2,...,k_n=0}^\infty(-1)^{k_2+\cdots+k_n}{k_2+\cdots+k_k\choose
k_2,...,k_n} \int_{(x_1)}^{((k_2+\cdots +k_n)m_1)}(x_1^{\ell_1})\\
& &\times\ptl_{x_2}^{k_2m_2}(x_2^{\ell_2})\cdots
\ptl_{x_n}^{k_nm_n}(x_n^{\ell_n})\mid
\ell_1\in\ol{0,m_1-1},\;\ell_2,...,\ell_n\in\mbb{N}\}\hspace{1.7cm}
(6.1.30)\end{eqnarray*} forms a basis of the space of polynomial
solutions for the equation
$$(\ptl_{x_1}^{m_1}+\ptl_{x_2}^{m_2}+\cdots+\ptl_{x_n}^{m_n})(u)=0\eqno(6.1.31)$$
in  ${\msr A}$.

Next we give an example of applying Lemma 6.1.1 iteratively . Let
$$f_i\in\mbb{F}[x_1,...,x_i]\qquad \for\;\;i\in\ol{1,n-1}.\eqno(6.1.32)$$
Consider the equation:
$$(\ptl_{x_1}^{m_1}+f_1\ptl_{x_2}^{m_2}+\cdots+f_{n-1}\ptl_{x_n}^{m_n})(u)=0
\eqno(6.1.33)$$
 Denote
$$d_1=\ptl_{x_1}^{m_1},\;\;
d_r=\ptl_{x_1}^{m_1}+f_1\ptl_{x_2}^{m_2}+\cdots+f_{r-1}\ptl_{x_r}^{m_r}\qquad\for\;\;r
\in\ol{2,n}.\eqno(6.1.34)$$ We will successively apply Lemma 6.1.1
with $T_1=d_r,\;T_2=f_r\ptl_{x_{r+1}}^{m_{r+1}}$ and ${\msr
B}=\mbb{F}[x_1,...,x_r],\\ V=\mbb{F}[x_{r+1}]$,
$$V_k=\sum_{s=0}^k\mbb Fx_{r+1}^s.\eqno(6.1.35)$$
 Take a right inverse
$d_1^-=\int_{(x_1)}^{(m_1)}$.  Suppose that we have found a right
inverse $d_s^-$ of $d_s$ for some $s\in\ol{1,n-1}$ such that
$$x_id_s^-=d_s^-x_i,\;\;
\ptl_{x_i}d_s^-=d_s^-\ptl_{x_i}\qquad\for\;\;i\in\ol{s+1,n}.\eqno(6.1.36)$$
Lemma 6.1.1 enables us to take
$$d_{s+1}^-=\sum_{i=0}^\infty(-d_s^-f_s)^id_s^-\ptl_{x_{s+1}}^{im_{s+1}}\eqno(6.1.37)$$
as a right inverse of $d_{s+1}$.  Obviously, $$
x_id_{s+1}^-=d_{s+1}^-x_i,\;\;
\ptl_{x_i}d_{s+1}^-=d_{s+1}^-\ptl_{x_i}\qquad\for\;\;i\in\ol{s+2,n}
\eqno(6.1.38)$$ according to (6.1.34). By induction, we have found a
right inverse $d_s^-$ of $d_s$ such that (6.1.36) holds for each
$s\in\ol{1,n}$.

We set
$${\msr S}_r=\{g\in \mbb{F}[x_1,...,x_r]\mid
d_r(g)=0\}\qquad\for\;\;r\in\ol{1,k}.\eqno(6.1.39)$$ By (6.1.34),
$${\msr S}_1=\sum_{i=0}^{m_1-1}\mbb{F}x_1^i.\eqno(6.1.40)$$
Suppose that we have found ${\msr S}_r$ for some $r\in \ol{1,n-1}$.
Given $h\in {\msr S}_r$ and $\ell\in \mbb{N}$, we define
$$\sgm_{r+1,\ell}(h)=\sum_{i=0}^\infty
(-d_r^-f_r)^i(h)\ptl_{x_{r+1}}^{im_{r+1}}(x_{r+1}^{\ell}),\eqno(6.1.41)$$
which is actually a finite summation. Lemma 6.1.1 says
$${\msr S}_{r+1}=\sum_{\ell=0}^\infty \sgm_{r+1,\ell}({\msr
S}_r).\eqno(6.1.42)$$ By  induction, we obtain:\psp

{\bf Theorem 6.1.3}. {\it The set
$$\{\sgm_{n,\ell_n}\sgm_{n-1,\ell_{n-1}}\cdots\sgm_{2,\ell_2}(x_1^{\ell_1})\mid
\ell_1\in\ol{0,m_1-1},\;\ell_2,...,\ell_n\in\mbb{N}\}\eqno(6.1.43)$$
forms a basis of the polynomial solution space ${\msr S}_n$ of the
partial differential equation (6.1.33).}\psp

\section{Canonical Oscillator Representations}

In this section, we present the canonical oscillator representations
of special linear Lie algebras.

 Recall that we denote by $E_{i,j}$ the square matrix
with 1 as its $(i,j)$-entry and 0 as the others.  The special linear
Lie algebra
$$sl(n,\mbb F)=\sum_{1\leq i<j\leq
n}(\mbb FE_{i,j}+\mbb FE_{j,i})+\sum_{r=1}^{n-1}\mbb
C(E_{r,r}-E_{r+1,r+1})\eqno(6.2.1)$$ with the Lie bracket:
$$[A,B]=AB-BA\qquad \for\;\;A,B\in sl(n,\mbb F).\eqno(6.2.2)$$
 Set
$$h_i=E_{i,i}-E_{i+1,i+1},\qquad i=1,2,...,n-1.\eqno(6.2.3)$$
The subspace
$$H=\sum_{i=1}^{n-1}\mbb Fh_i\eqno(6.2.4)$$
forms a Cartan subalgebra of $sl(n,\mbb F)$. We choose
$$\{E_{i,j}\mid 1\leq i<j\leq n\}\;\;\mbox{as positive root vectors}.\eqno(6.2.5)$$
In particular, we have
$$\{E_{i,i+1}\mid i=1,2,...,n-1\}\;\;\mbox{as positive simple root vectors}.\eqno(6.2.6)$$
Accordingly,
$$\{E_{i,j}\mid 1\leq j<i\leq n\}\;\;\mbox{are negative root vectors}\eqno(6.2.7)$$
and we have
$$\{E_{i+1,i}\mid i=1,2,...,n-1\}\;\;\mbox{as negative simple root vectors}.\eqno(6.2.8)$$
In particular,
$$sl(n,\mbb F)_+=\sum_{1\leq i<j\leq
n}\mbb FE_{i,j}\;\;\mbox{and}\;\;sl(n,\mbb F)_-=\sum_{1\leq i<j\leq
n}\mbb FE_{j,i}\eqno(6.2.9)$$ are the nilpotent subalgebra of
positive root vectors and the nilpotent subalgebra of negative root
vectors, respectively. Recall that an $sl(n,\mbb F)$-module $V$ is
called a {\it weight module} if for any $h\in H$, $h|_V$ is
diagonalizable. A {\it singular vector}\index{singular vector} of
$v$ of $V$ is a weight vector annihilated by the elements in
$sl(n,\mbb F)_+$. The fundamental weights $\lmd_i\in H^\ast$ are
$$\lmd_i(h_r)=\dlt_{i,r}.\eqno(6.2.10)$$

For $\al=(\al_1,\al_2,...,\al_n)\in\mbb{N}^{\:n}$, we define
$$|\al|=\sum_{i=1}^n\al_n.\eqno(6.2.11)$$
Denote ${\msr A}=\mbb F[x_1,x_2,...,x_n]$. Define a representation
of $sl(n,\mbb F)$ on ${\msr A}$ by:
$$E_{i,j}|_{\msr
A}=x_i\ptl_j\qquad\for\;\;i,j\in\ol{1,n}.\eqno(6.2.12)$$ For any
nonnegative integer $m$, we set
$${\msr A}_m=\mbox{Span}\:\{x^\al\mid
\al\in\mbb{N}^{\:n},\;|\al|=m\}.\eqno(6.2.13)$$ Note that a
polynomial $f\in {\msr A}$ is a solution of the system
$$0=E_{i,i+1}(f)=x_i\ptl_{x_{i+1}}(f)\qquad\for\;\;i\in\ol{1,n-1}\eqno(6.2.14)$$
if and only if $f$ is a one-variable polynomial in $x_1$. For
$m\in\mbb{N}$, ${\msr A}_m$ is a finite-dimensional $sl(n,\mbb
F)$-module and it has a unique singular vector $x_1^m$ up to scalar
multiple. Since ${\msr A}_m$ is weight module and any submodule must
contain a singular vector, we obtain:\psp

{\bf Proposition 6.2.1}. {\it The subspace ${\msr A}_m$ is an
irreducible $sl(n,\mbb F)$-module with highest weight $m\lmd_1$.}
\psp

Let $\Psi$ be the exterior algebra generated by
$\{\sta_1,\sta_2,...,\sta_n\}$; that is, an associative algebra with
the defining relations:
$$\sta_i\sta_j=-\sta_j\sta_i.\eqno(6.2.15)$$
Set
$$\Psi_r=\mbox{Span}\:\{\sta_{i_1}\sta_{i_2}\cdots\sta_{i_r}\mid
i_1,...,i_r\in\ol{1,n}\}\eqno(6.2.16)$$
 For each $i\in\ol{1,n}$, we
define a linear transformation $\ptl_{\sta_i}$ by:
$$\ptl_{\sta_i}(\sta_j)=\dlt_{i,j}\eqno(6.2.17)$$
and $$
\ptl_{\sta_i}(\sta_{j_1}\cdots\sta_{j_r})=\sum_{s=1}^r(-1)^{s-1}\sta_{j_1}\cdots
\sta_{j_{s-1}}\ptl_{\sta_i}(\sta_{j_s})\sta_{j_{s+1}}\cdots\sta_{j_r}.
\eqno(6.2.18)$$
 Now $sl(n,\mbb F)$ has a representation on $\Psi$ defined by:
$$E_{i,j}|_{\Psi}=\sta_i\ptl_{\sta_j}.\eqno(6.2.19)$$

Observe that $\sta_1\sta_2\cdots\sta_r$ is a singular vector of
$\Psi_r$ and
$$\sta_{i_1}\sta_{i_2}\cdots\sta_{i_r}=E_{i_1,1}E_{i_2,2}\cdots
E_{i_r,r}(\sta_1\sta_2\cdots\sta_r)\eqno(6.2.20)$$ for any $1\le
i_1<i_2<\cdots<i_r\leq n$. So $\Psi_r$ is an $sl(n,\mbb F)$-module
generated by $\sta_1\sta_2\cdots\sta_r$.\psp

{\bf Proposition 6.2.2}. {\it The subspace $\Psi_r$ forms an
irreducible $sl(n,\mbb F)$-module with highest  weight
$\lmd_r$.}\psp

Let
$${\cal Q}=\mbb F(x_1,...,x_n,y_1,...,y_n),\eqno(6.2.21)$$
the space of rational functions in $x_1,...,x_n,y_1,...,y_n$. Define
the {\it bosonic orthogonal oscillator representation}\index{
bosonic orthogonal oscillator representation}  of $sl(n,\mbb F)$ on
${\cal Q}$ via
$$E_{i,j}|_{\cal
Q}=x_i\ptl_{x_j}-y_j\ptl_{y_i}\qquad\for\;\;i,j\in\ol{1,n}.\eqno(6.2.22)$$
Set
$$\eta=\sum_{i=1}^n x_iy_i.\eqno(6.2.23)$$
Then
$$\xi(\eta)=0\qquad\for\;\;\xi\in sl(n,\mbb F).\eqno(6.2.24)$$
\pse

{\bf Lemma 6.2.3}. {\it Any singular function in ${\cal Q}$ is a
rational function in $x_1,y_n,\eta$}.

{\it Proof}. Let $f\in{\cal Q}$ be a singular function. We can write
$$f=g(x_1,...,x_{n-1},\eta,y_1,...,y_n)\eqno(6.2.25)$$ as a rational
functions in $x_1,...,x_{n-1},\eta,y_1,...,y_n$. By (6.2.23) and
(6.2.24), we have:
$$E_{i,n}(f)=(x_i\ptl_{x_n}-y_n\ptl_{y_i})(g)=y_ng_{y_i}=0\qquad\for\;\;i\in\ol{1,n-1},\eqno(6.2.26)$$
or equivalently,
$$g_{y_i}=0\qquad\for\;\;i\in\ol{1,n-1}.\eqno(6.2.27)$$
Thus (6.2.22), (6.2.24) and (6.2.27) imply
$$E_{1,i}(g)=(x_1\ptl_{x_i}-y_i\ptl_{y_1})(g)=x_1g_{x_i}=0\qquad\for\;\;i\in\ol{2,n-1};\eqno(6.2.28)$$
that is,
$$g_{x_i}=0\qquad\for\;\;i\in\ol{2,n-1},\eqno(6.2.29)$$
Therefore, $g$ is independent of $x_2,...,x_{n-1}$ and
$y_1,...,y_{n-1}.\qquad\Box$ \psp

Denote ${\msr B}=\mbb F[x_1,...,x_n,y_1,...,y_n]\subset {\cal Q}$.
Set
$${\msr
B}_{\ell_1,\ell_2}=\sum_{\al,\be\in\mbb{N}^{\:n};\;|\al|=\ell_1,\;|\be|=\ell_2}
\mbb Fx^\al
y^\be\qquad\for\;\;\ell_1,\ell_2\in\mbb{N}.\eqno(6.2.30)$$ Then
${\msr B}_{\ell_1,\ell_2}$ is a finite-dimensional $sl(n,\mbb
F)$-submodule by (6.2.22) and ${\msr
B}=\bigoplus_{\ell_1,\ell_2=0}^\infty{\msr B}_{\ell_1,\ell_2}$.
 The function
$x_1^{\ell_1}y_n^{\ell_2}$ is a singular function of weight
$\ell_1\lmd_1+\ell_2\lmd_{n-1}$. According to the above lemma, any
singular polynomial in ${\msr B}_{\ell_1,\ell_2}$ must be of the
form $ax_1^{\ell_1-i}y_n^{\ell_2-i}\eta^i$ for some $0\neq a\in\mbb
F$ and $i\in\mbb{N}$. Define
$$V_{\ell_1,\ell_2}=\mbox{the submodule
generated by}\;x_1^{\ell_1}y_n^{\ell_2}.\eqno(6.2.31)$$ Since ${\msr
B}_{\ell_1,\ell_2}$ is a weight module, Weyl's Theorem 2.3.6 of
complete reducibility implies that ${\msr B}_{\ell_1,\ell_2}$ is a
direct sum of its irreducible submodules, which are generated by its
singular polynomials. So $V_{\ell_1,\ell_2}$ is an irreducible
highest weight module with the highest weight
$\ell_1\lmd_1+\ell_2\lmd_{n-1}$ and
$${\msr B}_{\ell_1,\ell_2}=V_{\ell_1,\ell_2}\oplus \eta
{\msr B}_{\ell_1-1,\ell_2-1}\eqno(6.2.32)$$ as a direct sum of two
$sl(n,\mbb F)$-submodules, where we treat $V_{i,j}=\{0\}$ if
$\{i,j\}\not\subset\mbb{N}$.

Denote
$$\Dlt=\sum_{i=1}^n\ptl_{x_i}\ptl_{y_i}.\eqno(6.2.33)$$
It can be verified that
$$\xi\Dlt=\Dlt\xi\qquad\for\;\;\xi\in sl(n,\mbb F),\eqno(6.2.34)$$
as operators on ${\cal Q}$. Set
$${\msr H}_{\ell_1,\ell_2}=\{f\in {\msr B}_{\ell_1,\ell_2}\mid
\Dlt(f)=0\}\eqno(6.2.35)$$
 Since $\Dlt(x_1^{\ell_1}y_n^{\ell_2})=0$,
we have $$V_{\ell_1,\ell_2}\subset {\msr
H}_{\ell_1,\ell_2}\eqno(6.2.36)$$ by (6.2.31) and (6.2.34). On the
other hand,
$${\msr B}_{\ell_1,\ell_2}=\bigoplus_{i=0}^\infty
\eta^iV_{\ell_1-i,\ell_2-i}\eqno(6.2.37)$$ by (6.2.32) and
induction. Note
$$\Dlt\eta=n+\eta\Dlt+\sum_{i=1}^n(x_i\ptl_{x_i}+y_i\ptl_{y_i})\eqno(6.2.38)$$
as operators on ${\cal Q}$. Thus
$$\Dlt(\eta^i
g)=\sum_{r=1}^i(n+\ell_1+\ell_2-2r)(\eta^{i-1}g)=i(n+\ell_1+\ell_2-i+1)
\eta^{i-1}g \eqno(6.2.39)$$ for $i\in\mbb{N}+1,\;g\in
V_{\ell_1-i,\ell_2-i}.$ Hence
$${\msr H}_{\ell_1,\ell_2}\bigcap \eta {\msr
B}_{\ell_1-1,\ell_2-1}=\{0\}.\eqno(6.2.40)$$ Therefore,
$$V_{\ell_1,\ell_2}={\msr H}_{\ell_1,\ell_2}\eqno(6.2.41)$$
by (6.2.32) and (6.2.36). Now Lemma 6.1.1 gives:\psp

{\bf Theorem 6.2.4}. {\it The subspace ${\msr H}_{\ell_1,\ell_2}$ is
an irreducible  $sl(n,\mbb F)$-module with  highest weight
$\ell_1\lmd_1+\ell_2\lmd_{n-1}$. Moreover, it has a basis}
$$\left\{\sum_{i=0}^\infty (-1)^i\frac{(x_1y_1)^i(\sum_{i=2}^n\ptl_{x_i}\ptl_{y_i})^i(x^\al y^\be)}
{\prod_{r=1}^i(\al_1+i)(\be_1+i)}\mid\al,\be\in\mbb{N};\al_1\be_1=0,|\al|=\ell_1,|\be|=\ell_2\right\}.
\eqno(6.2.42)$$ \psp

The above result was due to the author [X16]. According to (6.2.38),
the space $\mbb F\Dlt+\mbb F[\Dlt,\eta]+\mbb F\eta$ forms a Lie
subalgebra of linear operators on $\msr B$, which is isomorphic to
$sl(2,\mbb F)$. The above result can be interpreted as an
$(sl(2),sl(n))$-Howe duality (cf. [Hr1-Hr4]).

 Consider the exterior
algebra $\check{\msr A}$ generated by
$\{\sta_1,...,\sta_n,\vt_1,...,\vt_n\}$ (cf. (6.2.15)). We define
the {\it fermionic orthogonal oscillator
representation}\index{fermionic orthogonal oscillator
representation} of $sl(n,\mbb F)$ on $\check{\msr A}$ by
$$E_{i,j}|_{\check{\msr
A}}=\sta_i\ptl_{\sta_j}-\vt_j\ptl_{\vt_i}\qquad\for\;\;i,j\in\ol{1,n}.\eqno(6.2.43)$$
Denote
$$\Theta_1=\sum_{i=1}^n\mbb F\sta_i,\qquad
\Theta_2=\sum_{i=1}^n\mbb F\vt_i.\eqno(6.2.44)$$ For
$\ell_1,\ell_2\in\ol{1,n}$, we define
$$\check{\msr
A}_{\ell_1,\ell_2}=\Theta_1^{\ell_1}\Theta_2^{\ell_2}.\eqno(6.2.45)$$
Then $\check{\msr A}_{\ell_1,\ell_2}$ is a finite-dimensional
$\check{\msr G}$-module and
$$\check{\msr
A}=\bigoplus_{\ell_1,\ell_2=0}^n\check{\msr
A}_{\ell_1,\ell_2}.\eqno(6.2.46)$$ Moreover, we define an ordering:
$$\sta_1\prec\sta_2\prec\cdots\prec\sta_n\prec\vt_n\prec\vt_{n-1}\prec\cdots\prec\vt_1.\eqno(6.2.47)$$
On the basis
\begin{eqnarray*}\qquad& &\{\sta_{i_1}\cdots\sta_{i_r}\vt_{j_1}\cdots\vt_{j_s}\mid
r,s\in\ol{0,n};\\ & &1\leq i_1<i_2<\cdots<i_r\leq n;n\geq
j_1>j_2>\cdots j_s\geq 1\}\hspace{3.6cm}(6.2.48)\end{eqnarray*} of
$\check{\msr A}$, we define the partial ordering ``$\prec$"
lexically.

Write
$$\check\Dlt=\sum_{r=1}^n\ptl_{\sta_r}\ptl_{\vt_r},\qquad
\check\eta=\sum_{r=1}^n\sta_r\vt_r.\eqno(6.2.49)$$ For
$r\in\ol{1,n}$, we define
$$\vec\sta_r=\sta_1\cdots\sta_r,\qquad
\vec\vt_r=\vt_n\cdots\vt_r.\eqno(6.2.50)$$ For convenience, we let
$$\vec\sta_0=1=\vec\vt_{n+1}.\eqno(6.2.51)$$
It can easily proved that a minimal term of any singular vector in
$\check{\msr A}_{\ell_1,\ell_2}$ is of the form
$\vec\sta_r\vec\vt_s$ for some $r\in\ol{0,n}$ and $s\in\ol{r+1,n+1}$
or
$$\vec\sta_r\sta_{r+1}\cdots\sta_{s_1}\vec\vt_{s_2}\vt_{s_1}\cdots\vt_{r+1}\eqno(6.2.52)$$
for some $0\leq r<s_1<s_2\leq n+1$.  By comparing minimal terms, we
can prove that
$$\{\check\eta^\ell\vec\sta_r\vec\vt_s\mid 0\leq r<s\leq
n+1;\ell\in\ol{0,s-r-1};r+\ell=\ell_1;\ell+n-s+1=\ell_2\}\eqno(6.2.53)$$
is the set of all $sl(n,\mbb F)$-singular vectors in $\check{\msr
A}_{\ell_1,\ell_2}$. Let $V_{r,s}$ be the finite-dimensional
irreducible $sl(n,\mbb F)$-submodule generated by
$\vec\sta_r\vec\vt_s\in\check{\msr A}_{r,n+1-s}$. Since $\check{\msr
A}$ is a weight module, Theorem 2.4.5 of Weyl's complete
reducibility gives
$$\check{\msr A}=\bigoplus_{0\leq r<s\leq
n+1}\bigoplus_{\ell=0}^{s-r-1}\check\eta^\ell V_{r,s}\eqno(6.2.54)$$
is a direct sum of irreducible $sl(n,\mbb F)$-submodules.

Define
$$\check{\msr H}=\{f\in\check{\msr
A}\mid\check\Dlt(f)=0\}.\eqno(6.2.55)$$
 Note
$$E_{r,s}\check\Dlt=\check\Dlt E_{r,s},\;\;
E_{r,s}\check\eta=\check\eta E_{r,s}\;\;\mbox{on}\;\;\check{\msr A}
\eqno(6.2.56)$$ for $r,s\in\ol{1,n}$. Moreover,
$$\check\Dlt\check\eta=\check\eta\check\Dlt-n+\sum_{r=1}^n(\sta_r\ptl_{\sta_r}
+\vt_r\ptl_{\vt_r})\eqno(6.2.57)$$ by (6.2.49). Furthermore,
$$\check\Dlt(\vec\sta_r\vec\vt_s)=0\qquad\mbox{if}\;\;r<s.\eqno(6.2.58)$$
Hence
$$V_{r,s}\subset \check{\msr H}\qquad\for\;\;0\leq r<s\leq
n+1\eqno(6.2.59)$$ by (6.2.56). Suppose $0\leq r+1<s\leq n+1$ and
$\ell\in\ol{1,s-r-1}$. For any $f\in V_{r,s}$, we have
$$\check\Dlt(\check\eta^\ell f)=(\sum_{p=0}^{\ell-1}(2p+r+1-s))\check\eta^{\ell-1} f
=\ell(\ell+r-s)\check\eta^{\ell-1} f.\eqno(6.2.60)$$ Therefore,
$$\check{\msr H}=\bigoplus_{0\leq r<s\leq
n+1} V_{r,s}.\eqno(6.2.61)$$ In particular,
$$\check{\msr H}_{r,n+1-s}=\{f\in\check{\msr
A}_{r,n+1-s}\mid\check\Dlt(f)=0\}=\check{\msr
A}_{r,n+1-s}\bigcap\check{\msr H}= V_{r,s}\eqno(6.2.62)$$  for
$0\leq r<s\leq n+1$ and
$$\check{\msr A}_{\ell_1,\ell_2}\bigcap\check{\msr
H}=\{0\}\qquad\mbox{if}\;\;\ell_1+\ell_2\geq n+1.\eqno(6.2.63)$$
Treat $\lmd_0=\lmd_n=0$.  We have:\psp

{\bf Theorem 6.2.5}. {\it For $0\leq r<s\leq n+1$, $\check{\msr
H}_{r,n+1-s}$ is a finite-dimensional irreducible $sl(n,\mbb
F)$-module with the highest-weight vector $\vec\sta_r\vec\vt_s$ of
weight $\lmd_r+\lmd_{s-1}$. Moreover,
$$\check{\msr A}=\bigoplus_{0\leq r<s\leq
n+1}\bigoplus_{\ell=0}^{s-r-1}\check\eta^\ell \check{\msr
H}_{r,n+1-s}.\eqno(6.2.64)$$}\pse

The above theorem was due to Luo and the author's work [LX4].
 According to (6.2.57),
the space $\mbb F\check{\Dlt}+\mbb F[\check{\Dlt},\check{\eta}]+\mbb
F\check{\eta}$ forms a Lie subalgebra of linear operators on
$\check{\msr A}$, which is isomorphic to $sl(2,\mbb F)$. The above
result can also be interpreted as an $(sl(2),sl(n))$-Howe duality.

\section{Noncanonical   Representations I: General }

In this section, we study the noncanonical oscillator
representation\index{oscillator representation} of deformed from the
canonical ones in (6.2.12) and (6.2.22).

Fix $1\leq n_1<n$. Note the symmetry
$$[\ptl_{x_r},x_r]=1=[-x_r,\ptl_{x_r}].\eqno(6.3.1)$$
Changing operators $\ptl_{x_r}\mapsto -x_r$ and $x_r\mapsto
\ptl_{x_r}$ in (6.2.12) for $r\in\ol{1,n_1}$, we obtain the
followings representation of $sl(n,\mbb{F})$ (also $gl(n,\mbb F)$)
on ${\msr A}=\mbb{F}[x_1,x_2,...,x_n]$:
$$E_{i,j}|_{\msr A}=\left\{\begin{array}{ll}-x_j\ptl_{x_i}-\delta_{i,j}&\mbox{if}\;
i,j\in\ol{1,n_1};\\ \ptl_{x_i}\ptl_{x_j}&\mbox{if}\;i\in\ol{1,n_1},\;j\in\ol{n_1+1,n};\\
-x_ix_j &\mbox{if}\;i\in\ol{n_1+1,n},\;j\in\ol{1,n_1};\\
x_i\partial_{x_j}&\mbox{if}\;i,j\in\ol{n_1+1,n}.
\end{array}\right.\eqno(6.3.2)$$
For any $k\in\mbb{Z}$, we denote
$${\msr A}_{\la
k\ra}=\mbox{Span}\:\{x^\al\mid\al\in\mbb{N}^n;\sum_{r=n_1+1}^n\al_r-\sum_{i=1}^{n_1}\al_i=k\}.
\eqno(6.3.3)$$ Then
$${\msr A}=\bigoplus_{k\in\mbb{Z}}{\msr A}_{\la k\ra}.\eqno(6.3.4)$$
Set
$$\td D=\sum_{r=n_1+1}^nx_r\ptl_{x_r}-\sum_{i=1}^{n_1}x_i\ptl_{x_i}.\eqno(6.3.5)$$
We have
$${\msr A}_{\la k\ra}=\{f\in{\msr A}\mid\td D(f)=kf\}.\eqno(6.3.6)$$
Moreover, we have
$$\td D E_{i,j}=E_{j,i}\td D\;\;\mbox{on}\;\;{\msr
A}\qquad\for\;i,j\in\ol{1,n}.\eqno(6.3.7)$$ Thus ${\msr A}_{\la
k\ra}$ forms a ${\msr G}$-module for any subalgebra ${\msr G}$ of
$gl(n,\mbb{F})$.

For $\al\in\mbb{N}^n$, we denote
$$\alpha!=\prod_{i=1}^n\al_i!.\eqno(6.3.8)$$
Define a symmetric bilinear form $(\cdot|\cdot)$ on ${\msr A}$ by
$$(x^\alpha|x^\beta)=\dlt_{\alpha,\beta}(-1)^{\sum_{i=1}^{n_1}\alpha_i}\alpha!\qquad\for\;\;
\al,\be\in \mbb{N}^n.\eqno(6.3.9)$$ Then we have: \psp

{\bf Lemma 6.3.1}. {\it For any $A\in gl(n,\mbb{F})$ and
$f,g\in\mbb{\msr A}$, we have
$$(A|g)=(f|A^t(g)),\eqno(6.3.10)$$
where $A^t$ denote the transpose of the matrix $A$.}

{\it Proof.} Let $\al,\be\in \mbb{N}^n$. For $i,j\in\ol{1,n_1}$,
$$(E_{i,j}(x^\al)|x^\be)=-\al_i(x^{\al+\es_j-\es_i}|x^\be)
-\dlt_{i,j}(x^\al|x^\be)\eqno(6.3.11)$$ and
$$(x^\al|E_{j,i}(x^\be))=-\be_j(x^\al|x^{\be+\es_i-\es_j})
-\dlt_{i,j}(x^\al|x^\be)\eqno(6.3.12)$$ by (6.3.2). Note
\begin{eqnarray*}\hspace{1cm}\al_i(x^{\al+\es_j-\es_i}|x^\be)&=&
\dlt_{\alpha+\es_j-\es_i,\beta}(-1)^{\sum_{r=1}^{n_1}\alpha_r}(\al_j+1)\alpha!
\\ &=&\be_j\dlt_{\alpha,\beta+\es_i-\es_j}(-1)^{\sum_{r=1}^{n_1}\alpha_r}\alpha!
=\be_j(x^\al|x^{\be+\es_i-\es_j})\hspace{2.3cm}(6.3.13)\end{eqnarray*}
by (6.3.9). Hence
$$(E_{i,j}(x^\al)|x^\be)=(x^\al|E_{j,i}(x^\be)).\eqno(6.3.14)$$
If $i,j\in\ol{n_1+1,n}$, then (6.3.13) holds and so does (6.3.14).

Consider $i\in\ol{1,n_1}$ and $j\in\ol{n_1+1,n}$.
$$(E_{i,j}(x^\al)|x^\be)=\al_i\al_j(x^{\al-\es_i-\es_j}|x^\be)=-
\dlt_{\alpha-\es_i-\es_j,\beta}(-1)^{\sum_{r=1}^{n_1}\alpha_r}\alpha!\eqno(6.3.15)$$
and
$$(x^\al|E_{j,i}(x^\be))=-(x^\al|x^{\be+\es_i+\es_j})=
-\dlt_{\alpha,\beta+\es_i+\es_j}(-1)^{\sum_{r=1}^{n_1}\alpha_r}\alpha!\eqno(6.3.16)$$
by (6.3.2) and (6.3.9). So (6.3.14) holds. Therefore, (6.3.10) holds
by the symmetry of the form.$\qquad\Box$\psp

For a subspace $V$ of ${\msr A}$, the {\it radical} of $V$ with
respect to the symmetric bilinear form $(\cdot|\cdot)$
is:\index{radical!of symmetric form}
$${\msr R}_V=\{v\in V\mid (v|u)=0\;\mbox{for any}\;u\in
V\}.\eqno(6.3.17)$$ If ${\msr R}_V=\{0\}$, we call $V$ {\it
nondegerate}.\index{nondegerate}

 Let ${\msr G}$ be simple Lie subalgebra of
$gl(n,\mbb{F})$ such that the transpose
$$A^t\in{\msr G}\qquad\mbox{if}\;A\in{\msr G}.\eqno(6.3.18)$$
Suppose that $H$ is a Cartan subalgebra of ${\msr G}$ and
$\Pi=\{\al_1,\al_2,....,\al_\ell\}$ is a base of the corresponding
root system $\Phi$. Take $0\neq \xi_i\in{\msr G}_{\al_i}$ and assume
$\zeta_i=\xi_i^t\in{\msr G}_{-\al_i}$ for $i\in\ol{1,\ell}$. Set
$${\msr G}_\pm=\sum_{\al\in\Phi_\pm}{\msr G}_\al.\eqno(6.3.19)$$
Then ${\msr G}_+$ is a subalgebra generated by $\{\xi_i\mid
i\in\ol{1,\ell}\}$, ${\msr G}_-$ is a subalgebra generated by
$\{\zeta_i\mid i\in\ol{1,\ell}\}$ and
$$[\xi_i,\zeta_j]=0,\;\;[\xi_i,\zeta_i]\in
H\qquad\for\;\;i,j\in\ol{1,\ell},\;i\neq j.\eqno(6.3.20)$$

Assume that ${\msr A}$ forms a weight ${\msr G}$-module with respect
to $H$.  A linear transformation (operator) $T$ on ${\msr A}$ is
called {\it locally nilpotent}\index{locally nilpotent} if for any
$f\in {\msr A}$, there exists a positive integer $k$ such that
$T^k(f)=0$. Moreover, an element $g\in{\msr A}$ is called {\it
nilpotent with respect to} ${\msr G}_+$ if there exist a positive
integer $m$ such that
$$u_1\cdots u_m(g)=0\qquad\mbox{for
any}\;u_1,...,u_m\in{\msr G}_+.\eqno(6.3.21)$$ A subspace $V$ of
${\msr A}$ is called {\it nilpotent with respect to} ${\msr G}_+$ if
all its elements are nilpotent with respect to ${\msr G}_+$.  If the
elements of ${\msr G}_+|_{\msr A}$ are locally nilpotent and
$${\msr G}_+({\msr A}_i)\subset \sum_{r=0}^i{\msr
A}_r\qquad\mbox{for any}\;\;i\in\mbb{N},\eqno(6.3.22)$$ then any
element of ${\msr A}$ is nilpotent with respect to ${\msr G}_+$ by
Theorem 1.5.1.\psp

{\bf Lemma 6.3.2}. {\it If a submodule $M$ of ${\msr A}$ is
nilpotent with respect to ${\msr G}_+$, $M$ contains a unique
singular vector $v$ (up to a scaler) and $(v|v)\neq0$, then $M$ is
an irreducible summand of ${\msr A}$}.

{\it Proof}. First, $V=U(\mathcal{G})(v)$ is an irreducible
submodule by the uniqueness of singular vector. For any $g\in M$,
there exists $m\in \mbb{N}+1$ such that (6.3.21) holds. Set
$$V_1=\mbox{Span}\:\{v,\zeta_{i_1}\cdots\zeta_{i_s}(v)\mid
m>s\in\mbb{N}+1;\;i_1,...,i_s\in\ol{1,\ell}\}\eqno(6.3.23)$$ and
$$V_2=\mbox{Span}\:\{\zeta_{i_1}\cdots\zeta_{i_s}(v)\mid
m\leq s\in\mbb{N}+1;\;i_1,...,i_s\in\ol{1,\ell}\}.\eqno(6.3.24)$$
Then
$$V=U({\msr G}_-)(v)=V_1\oplus V_2\eqno(6.3.25)$$
because ${\msr G}_-$ is  generated by $\{\zeta_i\mid
i\in\ol{1,\ell}\}$. Note that ${\msr R}_V$ is a submodule of $V$ by
Lemma 6.3.1 and the assumption (6.3.18). If ${\msr R}_V\neq \{0\}$,
then it contains a singular vector  under the nilpotent assumption,
which must be $v$ by the uniqueness. This contradicts  the
assumption $(v|v)\neq0$. Therefore ${\msr R}_V=\{0\}$; that is $V$
is nondegenerate.

According to Lemma 6.3.1
$$(\zeta_{i_1}\cdots\zeta_{i_s}(v)|\zeta_{j_1}\cdots\zeta_{j_r}(v))=
(\xi_{j_r}\cdots\xi_{j_1}\zeta_{i_1}\cdots\zeta_{i_s}(v)|v)=0\;\;\mbox{if}\;\;s<r,
\eqno(6.3.26)$$ by (6.3.20) and
$$(\zeta_{j_1}\cdots\zeta_{j_r}(v)|g)=(v|\xi_{j_r}\cdots\xi_{j_1}(g))=0\;\;\mbox{when}\;\;r\geq
m\eqno(6.3.27)$$ by the assumption of (6.3.21). Hence
$$(V_1|V_2)=\{0\},\;\;(V_2|g)=\{0\}.\eqno(6.3.28)$$
Thus $V_1$ and $V_2$ are nondegenerate. Since $\dim V_1$ is finite,
there exists $g'\in V_1$ such that $g-g'\in V_1^\perp$ by linear
algebra. Moreover, $g-g'\in V_2^\perp$ by (6.3.28). Thus $g-g'\in
V^\perp$. Thus we have
$$M=V\oplus M\bigcap V^\perp.\eqno(6.3.29)$$

 By  Lemma 6.3.1 and the assumption (6.3.18), $V^\perp$ is a
submodule of ${\msr A}$. So $M\bigcap V^\perp$ is a submodule of
$M$, which contains $v$ if it is nonzero. The assumption
$(v|v)\neq0$ leads $M\bigcap V^\perp=\{0\}$. Therefore, $M=V$ is
irreducible.$\qquad\Box$\psp

Let $f\in{\msr A}$ be a singular vector of $sl(n,\mbb{F})$. Then
$$E_{i,n_1}(f)=-x_{n_1}\ptl_{x_i}(f)=0\qquad\for\;\;i\in\ol{1,n_1-1}\eqno(6.3.30)$$
and
$$E_{n_1+1,j}(f)=x_{n_1+1}\ptl_{x_j}(f)=0\qquad\for\;\;j\in\ol{n_1+2,n}\eqno(6.3.31)$$
by (6.3.2). Thus $f=g(x_{n_1},x_{n_1+1})$ only depends on $x_{n_1}$
and $x_{n_1+1}$. Moreover,
$$0=E_{n_1,n_1+1}(f)=\ptl_{x_{n_1}}\ptl_{x_{n_1+1}}(g).\eqno(6.3.32)$$
This shows that $g$ is a one-variable function in $x_{n_1}$ or
$x_{n_1+1}$. Since $f$ is a weight vector, $f=x_{n_1}^{m_1}$ or
$x_{n_1+1}^{m_2}$ up to a scalar multiple for $m_1\in\mbb{N}$ and
$0\neq m_2\in\mbb{N}$. Note $x_{n_1}^{m_1}\in {\msr A}_{\la
-m_1\ra}$ and $x_{n_1+1}^{m_2}\in {\msr A}_{\la m_2\ra}$. So ${\msr
A}_{\la k\ra}$ has a unique singular vector for any $k\in\mbb{Z}$.
By the above lemma, ${\msr A}_{\la k\ra}$ is an irreducible
$sl(n,\mbb{F})$-submodule.

For any graded subspace $V$ of ${\msr A}$ (i.e.
$V=\bigoplus_{k=0}^\infty V\bigcap {\msr A}_k$), we define it
$q$-{\it dimension}\index{$q$-dimension}
 by
$$\dim_qV=\sum_{k=0}^\infty (\dim V\bigcap {\msr A}_k)q^k.\eqno(6.3.33)$$
 Recall that $\lmd_i$ is the $i$th fundamental
weight of $sl(n+1,\mbb F)$. For convenience, we treat
$\lmd_0=0$.\psp

{\bf Theorem 6.3.3}. {\it Let $m_1,m_2\in\mbb{N}$ with $m_1>0$,
${\msr A}_{\la -m_1\ra}$ is an irreducible highest-weight
$sl(n,\mbb{F})$-submodule with highest weight
$m_1\lmd_{n_1-1}-(m_1+1)\lmd_{n_1}$ and ${\msr A}_{\la m_2\ra}$ is
an irreducible highest-weight $sl(n,\mbb{F})$-submodule with highest
weight $-(m_2+1)\lmd_{n_1}+m_2\lmd_{n_1+1}$. Moreover,
\begin{eqnarray*}\qquad \sum_{m\in\mbb{Z}}z^m\mbox{\it ch}_{{\msr A}_{\la
m\ra}}
&=&\frac{\ves_{-\lmd_{n_1}}}{(1-z^{-1}\ves_{-\lmd_1})(1-z\ves_{-\lmd_{n-1}})}\\
& &\times\frac{1}{
[\prod_{i=2}^{n_1}(1-z^{-1}\ves_{\lmd_{i-1}-\lmd_i})]
[\prod_{r=n_1+1}^{n-1}(1-z\ves_{\lmd_r-\lmd_{r-1}})]},\hspace{1.2cm}(6.3.34)\end{eqnarray*}
$$\sum_{m\in\mbb{Z}}(\dim_q{\msr A}_{\la
m\ra})z^m=\frac{1}{(1-z^{-1}q)^{n_1}(1-zq)^{n-n_1}}.\eqno(6.3.35)$$}\pse

The above result was essentially due to Howe [Hr4], where he used
the symmetric tensor over several copies of the natural module and
its dual  to construct the representation. Since the structure
constants are integers, the result works for any field with
characteristic 0.

Fix $n_1,n_2\in\ol{1,n}$ such that $n_1\leq n_2$. Recall that ${\cal
Q}$ is the space of rational functions in $x_1,...,x_n,y_1,...,y_n$.
Changing operators $\ptl_{x_r}\mapsto -x_r,\;
 x_r\mapsto\ptl_{x_r}$  for $r\in\ol{1,n_1}$ and $\ptl_{y_s}\mapsto -y_s,\;
 y_s\mapsto\ptl_{y_s}$  for $s\in\ol{n_2+1,n}$ in the representation
 (6.2.22) of $sl(n,\mbb F)$, we get a new representation of $sl(n,\mbb{F})$ on ${\cal Q}$ determined
by
$$E_{i,j}|_{\cal
Q}=E_{i,j}^x-E_{j,i}^y\qquad\for\;\;i,j\in\ol{1,n}\eqno(6.3.36)$$
with
$$E_{i,j}^x=\left\{\begin{array}{ll}-x_j\ptl_{x_i}-\delta_{i,j}&\mbox{if}\;
i,j\in\ol{1,n_1};\\ \ptl_{x_i}\ptl_{x_j}&\mbox{if}\;i\in\ol{1,n_1},\;j\in\ol{n_1+1,n};\\
-x_ix_j &\mbox{if}\;i\in\ol{n_1+1,n},\;j\in\ol{1,n_1};\\
x_i\partial_{x_j}&\mbox{if}\;i,j\in\ol{n_1+1,n}
\end{array}\right.\eqno(6.3.37)$$
and
$$E_{i,j}^y=\left\{\begin{array}{ll}y_i\ptl_{y_j}&\mbox{if}\;
i,j\in\ol{1,n_2};\\ -y_iy_j&\mbox{if}\;i\in\ol{1,n_2},\;j\in\ol{n_2+1,n};\\
\ptl_{y_i}\ptl_{y_j} &\mbox{if}\;i\in\ol{n_2+1,n},\;j\in\ol{1,n_2};\\
-y_j\partial_{y_i}-\delta_{i,j}&\mbox{if}\;i,j\in\ol{n_2+1,n}.
\end{array}\right.\eqno(6.3.38)$$
Recall $\td D$ in (6.3.5) and define
$$\td D'=\sum_{i=1}^{n_2}y_i\ptl_{y_i}-\sum_{r=n_2+1}^ny_r\ptl_{y_r}.\eqno(6.3.39)$$
We set
$$\td\Dlt=\sum_{i=1}^{n_1}x_i\ptl_{y_i}-\sum_{r=n_1+1}^{n_2}\ptl_{x_r}\ptl_{y_r}+\sum_{s=n_2+1}^n
y_s\ptl_{x_s}\eqno(6.3.40)$$ and
$$\eta=\sum_{i=1}^{n_1}y_i\ptl_{x_i}+\sum_{r=n_1+1}^{n_2}x_ry_r+\sum_{s=n_2+1}^n
x_s\ptl_{y_s}.\eqno(6.3.41)$$ Then
$$TE_{i,j}|_{\cal Q}=E_{i,j}|_{\cal Q}T\qquad\for\;\;T=\td D,\td D',\td\Dlt,\eta;\;i,j\in\ol{1,n}.\eqno(6.3.42)$$ Moreover,
$$[\td D,\td\Dlt]=[\td D',\td\Dlt]=-\td\Dlt,\qquad
[\td D,\eta]=[\td D',\eta]=\eta.\eqno(6.3.43)$$ Furthermore, we have
$$E_{i,r}|_{\cal
Q}=-x_r\ptl_{x_i}-y_r\ptl_{y_i}\qquad\for\;\;1\leq i<r\leq
n_1,\eqno(6.3.44)$$
$$E_{i,n_1+s}|_{\cal
Q}=\ptl_{x_i}\ptl_{x_{n_1+s}}-y_{n_1+s}\ptl_{y_i}\qquad\for\;\;i\in\ol{1,n_1},\;s\in\ol{1,n_2-n_1},
\eqno(6.3.45)$$
$$E_{r,s}|_{\cal Q}=x_r\ptl_{x_s}-y_s\ptl_{y_r}\qquad\for\;\;n_1<
r<s\leq n_2,\eqno(6.3.46)$$
$$E_{n_2,n_2+1}|_{\cal Q}=x_{n_2}\ptl_{x_{n_2+1}}-\ptl_{y_{n_2}}\ptl_{y_{n_2+1}},\eqno(6.3.47)$$
$$E_{i,r}|_{\cal
Q}=x_i\ptl_{x_r}+y_i\ptl_{y_r}\qquad\for\;\;n_2+1\leq i<r\leq
n.\eqno(6.3.48)$$ Thus $\msr B=\mbb F[x_1,...,x_n,y_1,...,y_n]$ is
nilpotent with respect to $sl(n,\mbb F)_+$ in (6.2.9).

 Denote
$$\zeta_1=x_{n_1-1}y_{n_1}-x_{n_1}y_{n_1-1},\;\;
\zeta=\sum_{r=n_1+1}^{n_2}x_ry_r,\;\;\zeta_2=x_{n_2+1}y_{n_2+2}-x_{n_2+2}y_{n_2+1}.
\eqno(6.3.49)$$ Then
$$E_{i,j}(\zeta_1)=0\qquad\for\;\;1\leq i<j\leq n_1,\eqno(6.3.50)$$
$$E_{p,q}(\zeta)=0\qquad\for\;\;n_1+1\leq p<q\leq
n_2,\eqno(6.3.51)$$
$$E_{r,s}(\zeta_2)=0\qquad\for\;\;n_2+1\leq r<s\leq
n.\eqno(6.3.52)$$

 We will process according to three cases.\psp

\section{Noncanonical Representations II: $n_1+1<n_2$ }

In this section, we study the representation of $sl(n,\mbb{F})$
given in (6.3.36)-(6.3.38) with $n_1+1<n_2$.

 Assume $n_1+1<n_2<n$. Suppose that $f\in{\cal Q}$ is a singular vector. By Lemma 6.2.3,
$f$ can be written as a rational function in
$$\{x_i,y_r,\zeta_1,\zeta,\zeta_2\mid
n_2+2\neq i\in\ol{1,n_1+1}\bigcup\ol{n_2+1,n},\;n_1-1\neq
r\in\ol{1,n_1}\bigcup\ol{n_2,n}\}.\eqno(6.4.1)$$ Note that (6.3.44)
and (6.3.50) give
$$E_{n_1-1,n_1}(f)=-(x_{n_1}\ptl_{x_{n_1-1}}+y_{n_1}\ptl_{y_{n_1-1}})(f)=-x_{n_1}f_{x_{n_1-1}}=0.\eqno(6.4.2)$$
Moreover, (6.3.48) and (6.3.52) yield
$$E_{n_2+1,n_2+2}(f)=(x_{n_2+1}\ptl_{x_{n_2+2}}+y_{n_2+1}\ptl_{y_{n_2+2}})(f)=y_{n_2+1}f_{y_{n_2+2}}=0.
\eqno(6.4.3)$$ So $f$ is independent of $x_{n_1-1}$ and $y_{n_2+2}$
in terms of (6.4.1). For $i\in\ol{1,n_1-2}$, we have
\begin{eqnarray*}\qquad E_{i,n_1-1}(f)&=&-(x_{n_1-1}\ptl_{x_i}+y_{n_1-1}\ptl_{y_i})(f)\\
&=& -x_{n_1-1}f_{x_i}-y_{n_1-1}f_{y_i}\\ &=& -x_{n_1-1}(f_{x_i}
+x_{n_1}^{-1}y_{n_1}f_{y_i})
+x_{n_1}^{-1}\zeta_1f_{y_i}=0\hspace{4.2cm}(6.4.4)\end{eqnarray*}
 by (6.3.44), (6.4.49) and (6.4.50). Since both
 $f_{x_i}+x_{n_1}^{-1}y_{n_1}f_{y_i}$ and
 $x_{n_1}^{-1}\zeta_1f_{y_i}$ are independent of $x_{n_1-1}$ in terms of (6.4.1),
 we have $f_{y_i}=0$, and so $f_{y_i}=0$. Thus $f$ is
 independent of $\{x_i,y_i\mid i\in\ol{1,n_1-1}\}$ in terms of (6.4.1). Similarly, we can
 prove that $f$ is
 independent of $\{x_i,y_i\mid i\in\ol{n_2+1,n}\}$ in terms of (6.4.1). Therefore, $f$
 only depends on
 $$\{x_{n_1},x_{n_1+1},x_{n_2+1},y_{n_1},y_{n_2},y_{n_2+1},\zeta_1,\zeta,\zeta_2\}.
 \eqno(6.4.5)$$

According to (6.3.45) and (6.3.49), $E_{n_1,n_1+1}|_{\cal
Q}=\ptl_{x_{n_1}}\ptl_{x_{n_1+1}}-y_{n_1+1}\ptl_{y_{n_1}}$ and
\begin{eqnarray*}\qquad\qquad& &E_{n_1,n_1+1}(f)=f_{x_{n_1}x_{n_1+1}}-y_{n_1-1}f_{\zeta_1x_{n_1+1}}\\ & &
+y_{n_1+1}(f_{x_{n_1}\zeta}-y_{n_1-1}
f_{\zeta_1\zeta}-f_{y_{n_1}}-x_{n_1-1}f_{\zeta_1})=0.
\hspace{3.9cm}(6.4.6)\end{eqnarray*} Applying $E_{n_1+1,n_2}|_{\cal
Q}=x_{n_1+1}\ptl_{x_{n_2}}-y_{n_2}\ptl_{y_{n_1+1}}$  in (6.3.45) to
the above equation, we get
$$f_{x_{n_1}\zeta}-y_{n_1-1}f_{\zeta_1\zeta}-f_{y_{n_1}}-x_{n_1-1}f_{\zeta_1}=0\eqno(6.4.7)$$
by  (6.3.52). According to (6.3.49),
$$x_{n_1-1}=y_{n_1}^{-1}\zeta_1+x_{n_1}y_{n_1}^{-1}y_{n_1-1}.\eqno(6.4.8)$$
Substituting it into (6.4.7), we get
$$f_{x_{n_1}\zeta}-y_{n_1-1}(f_{\zeta_1\zeta}+y_{n_1}^{-1}x_{n_1}f_{\zeta_1})-f_{y_{n_1}}
-y_{n_1}^{-1}\zeta_1f_{\zeta_1}=0.\eqno(6.4.9)$$
 Since $f$ is independent of $y_{n_1-1}$ in terms of (6.4.5), we have
$$f_{\zeta_1\zeta}+y_{n_1}^{-1}x_{n_1}f_{\zeta_1}=0.\eqno(6.4.10)$$
Thus
$$f_{\zeta_1}=e^{-y_{n_1}^{-1}x_{n_1}\zeta}g\eqno(6.4.11)$$
for some rational function $g$ in the variables of (6.4.5) except
$\zeta$, i.e., $g_{\zeta}=0$. But $f$ is a rational function in the
variables of (6.4.5) and so is $f_{\zeta_1}$. Hence (6.4.11) forces
$f_{\zeta_1}=0$; that is, $f$ is independent of $\zeta_1$ in terms
of (6.4.5). Similarly, we can prove that $f$ is independent of
$\zeta_2$  in terms of (6.4.5). Now $f$ only depends on
$$\{x_{n_1},x_{n_1+1},x_{n_2+1},y_{n_1},y_{n_2},y_{n_2+1},\zeta\}.\eqno(6.4.12)$$

Next we only consider $f\in{\msr
B}=\mbb{F}[x_1,...,x_n,y_1,...,y_n]$. Since
$\zeta=\sum_{i=n_1+1}^{n_2}x_iy_i$, $f$ must be a polynomial in the
variables (6.4.12). Now (6.4.6) and (6.4.7) are equivalent to
$$f_{x_{n_1}x_{n_1+1}}=0,\qquad
f_{x_{n_1}\zeta}-f_{y_{n_1}}=0.\eqno(6.4.13)$$ Similarly, we can
prove
$$f_{y_{n_2}y_{n_2+1}}=0,\qquad
f_{y_{n_2+1}\zeta}-f_{x_{n_2+1}}=0.\eqno(6.4.14)$$ Set
$$\phi(m_1,m_2)=\sum_{i=0}^{\infty}
\frac{y_{n_1}^i(\ptl_{x_{n_1}}\ptl_{\zeta})^i(x_{n_1}^{m_1}\zeta^{m_2})}{i!}\qquad\for
\;\;m_1,m_2\in\mbb{N}.\eqno(6.4.15)$$ By Lemma 6.1.1 with
$T_1=\ptl_{y_{n_1}}\;T_1^-=\int_{(y_{n_1})}$ (cf. (6.1.20)) and
$T_2=-\ptl_{x_{n_1}}\ptl_{\zeta}$, the polynomial solution space of
the second equation in (6.4.13) is
$$[\sum_{m_1=0}^\infty\sum_{m_2=0}^\infty\mbb{F}\phi(m_1,m_2)]
[\mbb{F}[x_{n_1+1},x_{n_2+1},y_{n_2},y_{n_2+1}]].\eqno(6.4.16)$$ The
first equation in (6.4.13) says that $f$ can not contain the
monomials with $x_{n_1}x_{n_1+1}$ as a factor. Moreover,
$\phi(0,m_2)=\zeta^{m_2}$. Thus the polynomial solution space of
(6.4.13) is
$$[\mbb{F}[x_{n_1+1},\zeta]+\sum_{m_1=1}^\infty\sum_{m_2=0}^\infty\mbb{F}\phi(m_1,m_2)]
[\mbb{F}[x_{n_2+1},y_{n_2},y_{n_2+1}]].\eqno(6.4.17)$$

Denote
$$\psi(m_1,m_2)=\sum_{i=0}^{\infty}\frac{x_{n_2+1}^i(\ptl_{y_{n_2+1}}\ptl_{\zeta})^i(y_{n_2+1}^{m_1}\zeta^{m_2})}{i!}\qquad\for
\;\;m_1,m_2\in\mbb{N},\eqno(6.4.18)$$
\begin{eqnarray*}\qquad\phi(m_1,m_2,m_3)&=&\sum_{r=0}^\infty\frac{x_{n_2+1}^r(\ptl_{y_{n_2+1}}\ptl_{\zeta})^r(\phi(m_1,m_2)y_{n_2+1}^{m_3})}{r!}
\\ &=&\sum_{i,r=0}^\infty\frac{y_{n_i}^ix_{n_2+1}^r\ptl_{x_{n_1}}^i\ptl_{y_{n_2+1}}^r\ptl_{\zeta}^{i+r}
(x_{n_1}^{m_1}\zeta^{m_2}y_{n_2+1}^{m_3})}{i!r!}.\hspace{2.8cm}(6.4.19)\end{eqnarray*}
By Lemma 6.1.1 with
$T_1=\ptl_{x_{n_2+1}},\;T_1^-=\int_{(x_{n_2+1})}$ (cf. (6.1.20)) and
$T_2=-\ptl_{y_{n_2+1}}\ptl_{\zeta}$, the polynomial solution space
of the system (6.4.13) and the second equation in (6.4.14) is
$$[\sum_{m_1,m_2=0}^\infty\mbb{F}[x_{n_1+1}]\psi(m_1,m_2)+\sum_{m_1=1}^\infty\sum_{m_2,m_3=0}^\infty\mbb{F}\phi(m_1,m_2,m_3)]
[\mbb{F}[y_{n_2}]].\eqno(6.4.20)$$ The first equation in (6.4.14)
says that $f$ can not contain the monomials with $y_{n_2}y_{n_2+1}$
as a factor. Moreover, $\psi(0,m_2)=\zeta^{m_2}$ and
$\phi(m_1,m_2,0)=\phi(m_1,m_2)$. Therefore the polynomial solution
space of the system (6.4.13) and  (6.4.14) is
\begin{eqnarray*}\qquad& &\mbb{F}[x_{n_1+1},y_{n_2},\zeta]+\sum_{m_1,m_3=1}^\infty\sum_{m_2=0}^\infty \mbb{F}
\phi(m_1,m_2,m_3)
\\ &&+\sum_{m_1=1}^\infty\sum_{m_2=0}^\infty(\mbb{F}[y_{n_2}]\phi(m_1,m_2)+
\mbb{F}[x_{n_1+1}]\psi(m_1,m_2)).\hspace{4cm}(6.4.21)\end{eqnarray*}
According to (6.3.41),
$$x_{n_1+1}^{m_1}y_{n_2}^{m_2}\zeta^{m_3}=\eta^{m_3}(x_{n_1+1}^{m_1}y_{n_2}^{m_2}),
\eqno(6.4.22)$$
$$\eta^{m_2}(x_{n_1}^{m_1}y_{n_2}^{m_3})=(\zeta+y_{n_1}\ptl_{x_{n_1}})^{m_2}(x_{n_1}^{m_1}y_{n_2}^{m_3})
=\phi(m_1,m_2)y_{n_2}^{m_3},\eqno(6.4.23)$$
$$\eta^{m_2}(y_{n_2+1}^{m_1}x_{n_1+1}^{m_3})=(\zeta+x_{n_2+1}\ptl_{y_{n_2+1}})^{m_2}
(y_{n_2+1}^{m_1}x_{n_1+1}^{m_3})
=\psi(m_1,m_2)x_{n_1+1}^{m_3},\eqno(6.4.24)$$
$$\eta^{m_2}(x_{n_1}^{m_1}y_{n_2+1}^{m_3})=(\zeta+y_{n_1}\ptl_{x_{n_1}}
+x_{n_2+1}\ptl_{y_{n_2+1}})^{m_2}(x_{n_1}^{m_1} y_{n_2+1}^{m_3})
=\phi(m_1,m_2,m_3).\eqno(6.4.25)$$ It can be verified that
$\{\eta^{m_1}(x_i^{m_2}y_j^{m_3})\mid
m_1,m_2,m_3\in\mbb{N};i=n_1,n_1+1;j=n_2,n_2+1\}$ are singular
vectors. By (6.4.19)-(6.4.25), we have:\psp

{\bf Lemma 6.4.1}. {\it The  nonzero vectors in
$$\{\mbb{F}[\eta](x_i^{m_1}y_j^{m_2})\mid
m_1,m_2\in\mbb{N};i=n_1,n_1+1;j=n_2,n_2+1\}\eqno(6.4.26)$$ are all
the singular vectors of $sl(n,\mbb{F})$ in ${\msr
 B}=\mbb{F}[x_1,...,x_{n_1},y_1,...,y_{n_2}]$.

 Similarly, when $n_2=n$ and $n_1\leq n-2$,
 the  nonzero vectors in
$$\{\mbb{F}[\eta](x_i^{m_1}y_n^{m_2})\mid
m_1,m_2\in\mbb{N};i=n_1,n_1+1\}\eqno(6.4.27)$$ are all the singular
vectors of $sl(n,\mbb{F})$ in ${\msr B}$.}\psp

Denote
$${\msr H}=\{f\in{\msr B}\mid\td\Dlt(f)=0\}.\eqno(6.4.28)$$
By (6.3.42), ${\msr H}$ forms an $sl(n,\mbb{F})$-submodule.  Set
$${\msr B}_{\la \ell_1,\ell_2\ra}=\mbox{Span}\{x^\al
y^\be\mid\al,\be\in\mbb{N}^n;\sum_{r=n_1+1}^n\al_r-\sum_{i=1}^{n_1}\al_i=\ell_1;
\sum_{i=1}^{n_2}\be_i-\sum_{r=n_2+1}^n\be_r=\ell_2\}\eqno(6.4.29)$$
for $\ell_1,\ell_2\in\mbb{Z}$. Then
$${\msr B}_{\la \ell_1,\ell_2\ra}=\{f\in{\msr B}\mid\td D(f)=\ell_1f;\td D'(f)=\ell_2f\}.
\eqno(6.4.30)$$ By (6.3.42), ${\msr B}_{\la \ell_1,\ell_2\ra}$ forms
an $sl(n,\mbb{F})$-submodule, and so does
$${\msr H}_{\la \ell_1,\ell_2\ra}={\msr B}_{\la
\ell_1,\ell_2\ra}\bigcap {\msr H}.\eqno(6.4.31)$$

Next (6.3.40) and (6.3.41) imply
$$[\td\Dlt,\eta]=n_1-n_2-\td D-\td D',\;\;\td\Dlt(x_i^{m_1}y_j^{m_2})=0\eqno(6.4.32)$$ for
$m_1,m_2\in\mbb{N},\;i=n_1,n_1+1$ and $j=n_2,n_2+1$. Thus
$$x_{n_1+1}^{m_1}y_{n_2}^{m_2}\in {\msr H}_{\la m_1,m_2\ra},
\qquad x_{n_1+1}^{m_1}y_{n_2+1}^{m_2}\in {\msr H}_{\la
m_1,-m_2\ra},\eqno(6.4.33)$$
$$x_{n_1}^{m_1}y_{n_2}^{m_2}\in {\msr H}_{\la -m_1,m_2\ra},
\qquad x_{n_1}^{m_1}y_{n_2+1}^{m_2}\in {\msr H}_{\la
-m_1,-m_2\ra}.\eqno(6.4.34)$$ For any $g\in {\msr
H}_{\la\ell_1,\ell_2\ra}$and $0<m\in\mbb{N}$, we have
$$\eta^m(g)\in{\msr B}_{\ell_1+m,\ell_2+m}\eqno(6.4.35)$$
and
$$\td \Dlt(\eta^m(g))=m(-\ell_1-\ell_2+n_1-n_2-m+1)\eta^{m-1}(g).\eqno(6.4.36)$$
Thus
$$\td \Dlt(\eta^m(g))=0\;\;\mbox{if and only if}\;\;\ell_1+\ell_2\leq
n_1-n_2\;\;\mbox{and}\;\;m=n_1-n_2-\ell_1-\ell_2+1.\eqno(6.4.37)$$
If so,
$$\eta^m(g)\in{\msr
H}_{n_1-n_2-\ell_2+1,n_1-n_2-\ell_1+1}.\eqno(6.4.38)$$ Note
$$(n_1-n_2-\ell_2+1)+(n_1-n_2-\ell_1+1)=n_1-n_2+2+(n_1-n_2-\ell_1-\ell_2)\geq
n_1-n_2+2.\eqno(6.4.39)$$

Let $f_{\la\ell_1,\ell_2\ra}\in {\msr H}_{\la\ell_1,\ell_2\ra}$ be a
singular vector in (6.4.33) and (6.4.34). Then the singular vectors
in ${\msr H}$ are nonzero weight vectors in
$$\mbox{Span}\{f_{\la\ell_1,\ell_2\ra},\eta^{n_1-n_2+1-r_1-r_2}(f_{\la
r_1,r_2\ra})\mid\ell_1,\ell_2,r_1,r_2\in\mbb{Z};r_1+r_2\leq
n_1-n_2\}\eqno(6.4.40)$$ by Lemma 6.4.1, where
$$\eta^{n_1-n_2+1-r_1-r_2}(f_{\la r_1,r_2\ra})\in {\msr H}_{\la
n_1-n_2+1-r_2,n_1-n_2+1-r_1\ra}.\eqno(6.4.41)$$  Thus we get:\psp

{\bf Lemma 6.4.2}. {\it When $n_1+1<n_2<n,$
  we have
$${\msr H}_{\la \ell_1,\ell_2\ra}\;\mbox{has a unique singular
vector if}\;\;\ell_1+\ell_2\leq n_1-n_2+1\eqno(6.4.42)$$ and
$${\msr H}_{\la \ell_1,\ell_2\ra}\;\mbox{has exactly two singular
vectors if}\;\;\ell_1+\ell_2> n_1-n_2+1.\eqno(6.4.43)$$ In the case
$n_1+1<n_2=n,$ ${\msr B}_{\la\ell_1,\ell_2\ra}=0$ if $\ell_2<0$, and
for $\ell_1\in{Z}$ and $\ell_2\in\mbb{N}$, $${\msr H}_{\la
\ell_1,\ell_2\ra}\;\mbox{has a unique singular vector
if}\;\;\ell_1\geq n_1-n+2\;\mbox{or}\;\ell_1+\ell_2\leq
n_1-n+1,\eqno(6.4.44)$$
$${\msr H}_{\la
\ell_1,\ell_2+1\ra}\;\mbox{has exactly two singular vector
if}\;\mbox{and}\;\;n_1-n+1-\ell_2\leq\ell_1\leq
n_1-n+1.\eqno(6.4.45)$$}\psp

According to the bilinear form $(\cdot|\cdot)$ in Section 6.3 with
$n$ replaced by $2n$, $\ol{1,n_1}$ replaced by
$\ol{1,n_1}\bigcup\ol{n+n_2+1,2n}$ and $x_{n+i}$ replaced by $y_i$
for $i\in\ol{1,n}$, we have
$$(x^{\al}y^\be|x^{\al'}y^{\be'})=\dlt_{\al,\al'}
\dlt_{\be,\be'}(-1)^{\sum_{i=1}^{n_1}\al_i+\sum_{r=n_2+1}^n\be_r}\al!\be!\qquad\for
\;\;\al,\be,\al',\be'\in\mbb{N}^n.\eqno(6.4.46)$$  It can be
verified that
$$(\td\Dlt(x^{\al}y^\be)|x^{\al'}y^{\be'})=
-(x^{\al}y^\be|\eta(x^{\al'}y^{\be'})).\eqno(6.4.47)$$ We remind the
concerned Lie algebra ${\msr G}=sl(n,\mbb{F})$ in this section. The
Cartan subalgebra $H$ is given in (6.2.4). Suppose that
$f_{\la\ell_1,\ell_2\ra}\in {\msr H}_{\la\ell_1,\ell_2\ra}$ is a
singular vector in (6.4.33) and (6.4.34). Then
$$(f_{\la\ell_1,\ell_2\ra}|f_{\la\ell_1,\ell_2\ra})\neq
0\eqno(6.4.48)$$ and
$$(f_{\la\ell_1,\ell_2\ra}|f_{\la\ell_1',\ell_2'\ra})=0\qquad\mbox{if}\;\;(\ell_1,\ell_2)\neq
(\ell_1',\ell_2').\eqno(6.4.49)$$ Recall that ${\msr
G}_+=\sum_{1\leq i<j\leq n}\mbb{F}E_{i,j}$ is the subalgebra spanned
by the positive root vectors and ${\msr G}_-=\sum_{1\leq i<j\leq
n}\mbb{F}E_{j,i}$ is the subalgebra spanned by the negative root
vectors. Moreover,
$${\msr G}_-^t={\msr G}_+.\eqno(6.4.50)$$
According to (6.3.46)-(6.3.50), ${\msr B}$ is nilpotent with respect
to ${\msr G}_+$. Thus all ${\msr H}_{\la\ell_1,\ell_2\ra}$ with
$\ell_1+\ell_2\leq n_1-n_2+1$ are irreducible
$sl(n,\mbb{F})$-submodules by Lemma 6.3.2 and Lemma 6.4.2, and so
are $\eta^m({\msr H}_{\la\ell_1,\ell_2\ra})$ for any $m\in\mbb{N}$
by (6.3.44).

We extend the transpose to an algebraic anti-isomorphism on $U({\msr
G})$ by $1^t=1$ and
$$(A_1A_2\cdots A_r)^t=A_r^t\cdots A_2^tA_1^t\qquad \for\;\; A_i\in{\msr
G}.\eqno(6.4.51)$$  By the irreducibility,
$$ {\msr H}_{\la\ell_1,\ell_2\ra}=U({\msr G}_-)(f_{\la\ell_1,\ell_2\ra})\qquad\mbox{if}\;\;\ell_1+\ell_2\leq
n_1-n_2+1.\eqno(6.4.52)$$

Let $\ell_1,\ell_2,\ell_1',\ell_2'\in\mbb{Z}$ such that
$\ell_1+\ell_2,\ell_1'+\ell_2'\leq n_1-n_2+1$ and
$(\ell_1,\ell_2)\neq (\ell_1',\ell_2')$. Then
$$(w(f_{\la\ell_1,\ell_2\ra})|f_{\la\ell_1',\ell_2'\ra})=
(f_{\la\ell_1,\ell_2\ra}|w^t(f_{\la\ell_1',\ell_2'\ra}))=0\qquad\for\;\;w\in
U({\msr G}_-){\msr G}-\eqno(6.4.53)$$ by Lemma 6.3.1. Since
$f_{\la\ell_1,\ell_2\ra}$ is a weight vector, we have
$$U(H)(f_{\la\ell_1,\ell_2\ra})\subset\mbb{F}f_{\la\ell_1,\ell_2\ra}.\eqno(6.4.54)$$
Thus for any $w_1,w_2\in U({\msr G}_-)$,
$$(w_1(f_{\la\ell_1,\ell_2\ra})|w_2(f_{\la\ell_1',\ell_2'\ra}))=
(w_2^tw_1(f_{\la\ell_1,\ell_2\ra})|f_{\la\ell_1',\ell_2'\ra})=c(f_{\la\ell_1,\ell_2\ra})|f_{\la\ell_1',\ell_2'\ra})
\eqno(6.4.55)$$ for some $c\in\mbb{F}$ by (6.4.53) and (6.4.54).
Hence (6.4.52) implies
$$({\msr H}_{\la\ell_1,\ell_2\ra}|{\msr
H}_{\la\ell_1',\ell_2'\ra})=\{0\}.\eqno(6.4.56)$$

For $f\in {\msr H}_{\la\ell_1,\ell_2\ra},\;g\in {\msr B}$ and
$m,m'\in\mbb{N}$ such that $m\leq m'$, we find
\begin{eqnarray*} \qquad(\eta^m(f)|\eta^{m'}(g))&=& (-1)^{m'}(\td
\Dlt^{m_1}\eta^m(f)|g)\\
&=&(-1)^{m'}\dlt_{m',m}m![\prod_{r=0}^{m-1}
(-\ell_1-\ell_2+n_1-n_2-r)](f|g)\hspace{1.2cm}(6.4.57)\end{eqnarray*}
by (6.4.36) and (6.4.47). In particular, the singular vectors
$\eta^{n_1-n_2+1-r_1-r_2}(f_{\la r_1,r_2\ra})$ for
$r_1,r_2\in\mbb{Z}$ with $r_1+r_2\leq n_1-n_2$ are isotropic
polynomials. Moreover, for $m,m'\in\mbb{N}$ and
$\ell_1,\ell_2,\ell_1',\ell_2'\in\mbb{Z}$ such that
$\ell_1+\ell_2,\ell_1'+\ell_2'\leq n_1-n_2+1$,
$$(\eta^m({\msr H}_{\la\ell_1,\ell_2\ra})|\eta^{m'}({\msr
H}_{\la\ell_1',\ell_2'\ra}))=\{0\}\qquad\mbox{if}\;\;(m,\ell_1,\ell_1)\neq
(m',\ell_1',\ell_2').\eqno(6.4.58)$$ On the other hand, for
$0<m\in\mbb N$,
\begin{eqnarray*}\qquad\qquad & &(\eta^m(f_{\la\ell_1-m,\ell_2-m\ra})|\eta^m
(f_{\la\ell_1-m,\ell_2-m\ra}))\\ &=&(-1)^mm![\prod_{r=0}^{m-1}
(2m-\ell_1-\ell_2+n_1-n_2-r)]\\&
&\times(f_{\la\ell_1-m,\ell_2-m\ra}|
f_{\la\ell_1-m,\ell_2-m\ra})\neq
0\hspace{6.2cm}(6.4.59)\end{eqnarray*} by (6.4.57). Since the
radical of $(\cdot|\cdot)$ on $\eta^m({\msr
H}_{\la\ell_1-m,\ell_2-m\ra})$ is a proper submodule by Lemma 6.3.1,
the irreducibility of $\eta^m({\msr H}_{\la\ell_1-m,\ell_2-m\ra})$
implies that
$$(\cdot|\cdot)\;\;\mbox{is nondegenerate rewtricted to}\;\;\eta^m({\msr
H}_{\la\ell_1-m,\ell_2-m\ra}).\eqno(6.4.60)$$

Fix $\ell_1,\ell_2\in\mbb{Z}$ with $\ell_1+\ell_2\leq n_1-n_2+1$.
Set
$$\hat{\msr B}_{\la\ell_1,\ell_2\ra}=\sum_{m=0}^\infty\eta^m({\msr
H}_{\la\ell_1-m,\ell_2-m\ra}).\eqno(6.4.61)$$ By (6.4.58) and
(6.4.60), the above sum is a direct sum and $(\cdot|\cdot)$ is
nondegenerate restricted to $\hat{\msr B}_{\la\ell_1,\ell_2\ra}$.
 For any
$g\in \msr B_{\la\ell_1,\ell_2\ra}$, there exists an $m_0\in\mbb{N}$
such that
$$\td\Dlt^{m_0+1}(g)=0\eqno(6.4.62)$$
by (6.3.40). So
$$(g|\eta^m({\msr
H}_{\la\ell_1-m,\ell_2-m\ra}))=(-1)^m(\td\Dlt^m(g)|{\msr
H}_{\la\ell_1-m,\ell_2-m\ra})=\{0\}\;\;\for\;\;m_0< m\in\mbb
N.\eqno(6.4.63)$$ By the arguments in (6.3.23)-(6.3.29), there
exists
$$g'\in\sum_{m=0}^{m_0}\eta^m({\msr
H}_{\la\ell_1-m,\ell_2-m\ra})\;\;\mbox{such
that}\;\;(g-g'|\sum_{m=0}^{m_0}\eta^m({\msr
H}_{\la\ell_1-m,\ell_2-m\ra}))=\{0\}.\eqno(6.4.64)$$ On the other
hand,
$$(g-g'|\sum_{k=m_0+1}^\infty\eta^k({\msr
H}_{\la\ell_1-k,\ell_2-k\ra})=\{0\}\eqno(6.4.65)$$ by (6.4.58) and
(6.4.63). Thus $g-g'\in \hat{\msr B}_{\la\ell_1,\ell_2\ra}^\perp$.
Therefore
$${\msr B}_{\la\ell_1,\ell_2\ra}=\hat{\msr
B}_{\la\ell_1,\ell_2\ra}\oplus (\hat{\msr
B}_{\la\ell_1,\ell_2\ra}^\perp\bigcap{\msr
B}_{\la\ell_1,\ell_2\ra}).\eqno(6.4.66)$$ If $\hat{\msr
B}_{\la\ell_1,\ell_2\ra}^\perp\bigcap{\msr
B}_{\la\ell_1,\ell_2\ra}\neq \{0\}$, then it contains a singular
vector, which must be of the form
$\eta^{m_1}(f_{\la\ell_1-m_1,\ell_2-m_1\ra})$ for some
$m_1\in\mbb{N}$ Lemma 6.4.1. This contradicts (6.4.59). Thus
$\hat{\msr B}_{\la\ell_1,\ell_2\ra}^\perp\bigcap{\msr
B}_{\la\ell_1,\ell_2\ra}= \{0\}$, or equivalently
$${\msr B}_{\la\ell_1,\ell_2\ra}=\bigoplus_{m=0}^\infty\eta^m({\msr
H}_{\la\ell_1-m,\ell_2-m\ra})\eqno(6.4.67)$$ is completely
reducible. Applying (6.4.67) to ${\msr
B}_{\la\ell_1-1,\ell_2-1\ra}$, we have
$${\msr B}_{\ell_1,\ell_2}={\msr H}_{\la\ell_1,\ell_2\ra}\oplus
\eta({\msr B}_{\la\ell_1-1,\ell_2-1\ra})\qquad\mbox{if}\;\;
\ell_1+\ell_2\leq n_1-n_2+1.\eqno(6.4.68)$$

Assume $n_1+1<n_2=n$. For $\ell_1\in\mbb{Z}$ and $\ell_2\in\mbb{N}$
such that $\ell_1\geq n_1-n+2\;\mbox{or}\;\ell_1+\ell_2\leq
n_1-n+1$, all  ${\msr H}_{\ell_1,\ell_2}$  are irreducible
submodules of ${\msr B}_{\ell_1,\ell_2}$ by Lemma 6.3.2, (6.4.44)
and (6.4.45).  In summary, we have:\psp

{\bf Theorem 6.4.3}. {\it Suppose $n_1+1<n_2$. For
$\ell_1,\ell_2\in\mbb{Z}$ with $\ell_1+\ell_2\leq n_1-n_2+1$ ,
${\msr H}_{\la\ell_1,\ell_2\ra}$ is an infinite-dimensional
irreducible highest-weight $sl(n,\mbb{F})$-module and
$${\msr B}_{\la\ell_1,\ell_2\ra}=\bigoplus_{m=0}^\infty\eta^m({\msr
H}_{\la\ell_1-m,\ell_2-m\ra})\eqno(6.4.69)$$ is an orthogonal
decomposition of irreducible submodules. In particular, ${\msr
B}_{\la\ell_1,\ell_2\ra}={\msr H}_{\la\ell_1,\ell_2\ra}\oplus
\eta({\msr B}_{\la\ell_1-1,\ell_2-1\ra})$. The symmetric bilinear
form restricted to $\eta^m({\msr H}_{\la\ell_1-m,\ell_2-m\ra})$ is
nondegenerate. If $n_2<n$, all ${\msr H}_{\la\ell_1,\ell_2\ra}$ for
$\ell_1,\ell_2\in\mbb{Z}$ with $\ell_1+\ell_2> n_1-n_2+1$ have
exactly two singular vectors.

 Assume $n_2=n$. Then ${\msr
B}_{\la\ell,0\ra}={\msr H}_{\la\ell,0\ra}$ with $\ell\in\mbb{Z}$ are
 infinite-dimensional irreducible highest-weight $sl(n,\mbb{F})$-modules.  All ${\msr
H}_{\la\ell_1,\ell_2\ra}$ for $\ell_1\in\mbb{Z}$ and
$\ell_2\in\mbb{N}$ such that $\ell_1\geq n_1-n+2$ are also
 infinite-dimensional irreducible highest-weight $sl(n,\mbb{F})$-modules. Moreover, for
$\ell_2\in 1+\mbb{N}$ and $n_1-n_2+1+\ell_2\leq \ell_1\in\mbb{Z}$,
the orthogonal decompositions in (6.4.69) hold. Furthermore, ${\msr
H}_{\la\ell_1,\ell_2+1\ra}$ for $\ell_1\in\mbb{Z}$ and
$\ell_2\in\mbb{N}$ such that $n_1-n+1-\ell_2\leq \ell_1\leq n_1-n+1$
have exactly two singular vectors.}\psp

Indeed, we have more detailed information.  Suppose $n_1+1<n_2<n$.
For $m_1,m_2\in\mbb{N}$ with $m_1+m_2\geq n_2-n_1-1$, ${\msr
H}_{\la-m_1,-m_2\ra}$ has  a highest-weight vector
$x_{n_1}^{m_1}y_{n_2+1}^{m_2}$ of weight
$m_1\lmd_{n_1-1}-(m_1+1)\lmd_{n_1}-(m_2+1)\lmd_{n_2}+m_2(1-\dlt_{n_2,n-1})\lmd_{n_2+1}$.
When $m_1,m_2\in\mbb{N}$ with $m_2-m_1\geq n_2-n_1-1$, ${\msr
H}_{\la m_1,-m_2\ra}$ has a highest-weight vector
$x_{n_1+1}^{m_1}y_{n_2+1}^{m_2}$ of weight
$-(m_1+1)\lmd_{n_1}+m_1\lmd_{n_1+1}-(m_2+1)\lmd_{n_2}+m_2(1-\dlt_{n_2,n-1})\lmd_{n_2+1}$.
If $m_1,m_2\in\mbb{N}$ with $m_1-m_2\geq n_2-n_1-1$, ${\msr H}_{\la
-m_1,m_2\ra}$ is has a highest-weight vector
$x_{n_1}^{m_1}y_{n_2}^{m_2}$ of weight
$m_1\lmd_{n_1-1}-(m_1+1)\lmd_{n_1}+m_2\lmd_{n_2-1}-(m_2+1)\lmd_{n_2}$.

Assume $n_1+1<n_2=n$. When $m_1,m_2\in\mbb{N}$,  ${\msr H}_{\la
m_1,m_2\ra}$  has a highest-weight vector $x_{n_1+1}^{m_1}y_n^{m_2}$
of weight $-(m_1+1)\lmd_{n_1}+m_1\lmd_{n_1+1}+m_2\lmd_{n-1}$. If
$m_1,m_2\in\mbb{N}$ with $m_1\leq n-n_1-2$ or $m_2-m_1\leq n_1-n+1$,
${\msr H}_{\la -m_1,m_2\ra}$ has a highest-weight vector
$x_{n_1}^{m_1}y_n^{m_2}$ of weight
$m_1\lmd_{n_1-1}-(m_1+1)\lmd_{n_1}+m_2\lmd_{n-1}$.

By Lemma 6.1.1, ${\msr H}_{\la\ell_1,\ell_2\ra}$ has a basis
\begin{eqnarray*}\qquad& &\big\{\sum_{i=0}^\infty\frac{(x_{n_1+1}y_{n_1+1})^i(\td\Dlt+\ptl_{x_{n_1}+1}\ptl_{y_{n_1}+1})^i(x^\al
y^\be)}{\prod_{r=1}^i(\al_{n_1+1}+r)(\be_{n_1+1}+r)}
\mid\al,\be\in\mbb{N}^n;\\ &
&\al_{n_1+1}\be_{n_1+1}=0;\sum_{i=1}^{n_1}\al_i+\ell_1=\sum_{r=n_1+1}^n\al_r;
\sum_{i=1}^{n_2}\be_i=\sum_{r=n_2+1}^n\be_r+\ell_2
\big\}.\hspace{1.7cm}(6.4.70)\end{eqnarray*} \psp

\section{Noncanonical   Representations III: $n_1+1=n_2$ }

In this section, we study the representation of $sl(n,\mbb{F})$
given in (6.3.36)-(6.3.38) with $n_1+1=n_2$.

First we consider the subcase $n_2<n$. Suppose that $f\in{\cal Q}$
is a singular vector. Note $\zeta=x_{n_1+1}y_{n_1+1}$ and
$y_{n_2}=y_{n_1+1}$ in this case. According to the arguments in
(6.4.1)-(6.4.5), $f$ is a rational function in
$$\{x_{n_1},x_{n_1+1},x_{n_1+2},y_{n_1},y_{n_1+2},\zeta,\zeta_1,\zeta_2\}.\eqno(6.5.1)$$
By (6.4.6),
$$f_{x_{n_1}x_{n_1+1}}-y_{n_1-1}f_{\zeta_1x_{n_1+1}}\\
+y_{n_1+1}(f_{x_{n_1}\zeta}-y_{n_1-1}
f_{\zeta_1\zeta}-f_{y_{n_1}}-x_{n_1-1}f_{\zeta_1})=0.\eqno(6.5.2)$$
Substituting
$$x_{n_1-1}=y_{n_1}^{-1}\zeta_1+x_{n_1}y_{n_1}^{-1}y_{n_1-1},\;\;y_{n_1+1}=x_{n_1+1}^{-1}\zeta\eqno(6.5.3)$$
into (6.5.2), we get
\begin{eqnarray*}\qquad\qquad&
&f_{x_{n_1}x_{n_1+1}}+x_{n_1+1}^{-1}\zeta(f_{x_{n_1}\zeta}-f_{y_{n_1}}-y_{n_1}^{-1}\zeta_1f_{\zeta_1})\\
&&-y_{n_1-1}[f_{\zeta_1x_{n_1+1}}+x_{n_1+1}^{-1}\zeta(f_{\zeta_1\zeta}+x_{n_1}y_{n_1}^{-1}f_{\zeta_1})]=0.\hspace{4cm}
(6.5.4)\end{eqnarray*}

Since $f$ is independent of $y_{n_1-1}$ in terms of (6.5.1), we have
$$f_{\zeta_1x_{n_1+1}}+x_{n_1+1}^{-1}\zeta(f_{\zeta_1\zeta}+x_{n_1}y_{n_1}^{-1}f_{\zeta_1})]=0;\eqno(6.5.5)$$
equivalently
$$x_{n_1+1}f_{\zeta_1x_{n_1+1}}+\zeta(f_{\zeta_1\zeta}+x_{n_1}y_{n_1}^{-1}f_{\zeta_1})=0.\eqno(6.5.6)$$
Change variable
$$f_{\zeta_1}=e^{-x_{n_1}y_{n_1}^{-1}\zeta}g.\eqno(6.5.7)$$
Then
$$x_{n_1+1}f_{\zeta_1x_{n_1+1}}=x_{n_1+1}e^{-x_{n_1}y_{n_1}^{-1}\zeta}g_{x_{n_1+1}},\;\;
f_{\zeta_1\zeta}+x_{n_1}y_{n_1}^{-1}f_{\zeta_1}=e^{-x_{n_1}y_{n_1}^{-1}\zeta}g_\zeta.\eqno(6.5.8)$$
Thus (6.5.6) is equivalent to
$$x_{n_1+1}g_{x_{n_1+1}}+\zeta g_\zeta=0,\eqno(6.5.9)$$
which says that the total degree of $g$ with respect to $x_{n_1+1}$
and $\zeta$ is zero. In particular, $x_{n_1+1}/\zeta$ is a solution.
Changing variables $(x_{n_1+1},\zeta)\mapsto
(x_{n_1+1}/\zeta,\zeta)$, we can prove
$$g=h(x_{n_1+1}/\zeta)\eqno(6.5.10)$$
for some one-variable function $h(z)$ (it can also be derived by the
method of characteristic line in partial differential equation).
Since
$$f_{\zeta_1}=e^{-x_{n_1}y_{n_1}^{-1}\zeta}h(x_{n_1+1}/\zeta)\eqno(6.5.11)$$
is rational in $x_{n_1+1}$, $h(z)$ must be a rational function in
$z$. Now $f_{\zeta_1}$ and $h(x_{n_1+1}/\zeta)$ are rational in
$\zeta$. Equation (6.5.11) forces $h(z)=0$ and so $f_{\zeta_1}=0$.
Thus $f$ is independent of $\zeta_1$ in terms of (6.5.1).
Symmetrically, $f$ is independent of $\zeta_2$ in terms of (6.5.1).
 Hence
$f$ a rational function in
$$\{x_{n_1},x_{n_1+1},x_{n_1+2},y_{n_1},y_{n_1+2},\zeta\}.
 \eqno(6.5.12)$$
The fact $\zeta=x_{n_1+1}y_{n_1+1}$ allows us to rewrite $f$ as
 a rational function in
$$\{x_{n_1},x_{n_1+1},x_{n_1+2},y_{n_1},y_{n_1+1},y_{n_1+2}\}.
 \eqno(6.5.13)$$

Now $f$ is a singular vector if and only if it is a weight vector
satisfying the following system of differential equations
$$E_{n_1,n_1+1}(f)=(\ptl_{x_{n_1}}\ptl_{x_{n_1+1}}-y_{n_1+1}\ptl_{y_{n_1}})(f)=0,\eqno(6.5.14)$$
$$E_{n_1+1,n_1+2}(f)=(x_{n_1+1}\ptl_{x_{n_1+2}}-\ptl_{y_{n_1+1}}\ptl_{y_{n_1+2}})(f)=0.
\eqno(6.5.15)$$  Since the differential operators in the above two
equations do not commute, we have to solve them according to the
proof of Theorem 1.5.1. Note
$$E_{n_1,n_1+2}|_{\cal Q}=[E_{n_1,n_1+1}|_{\cal
Q},E_{n_1+1,n_1+2}|_{\cal Q}]=
\ptl_{x_{n_1}}\ptl_{x_{n_1+2}}-\ptl_{y_{n_1}}\ptl_{y_{n_1+2}}.\eqno(6.5.16)$$
 So
$$(\ptl_{x_{n_1}}\ptl_{x_{n_1+2}}-\ptl_{y_{n_1}}\ptl_{y_{n_1+2}})(f)=0.\eqno(6.5.17)$$

For our purpose of representation, we only consider $f$ is a
polynomial in $\{x_i,y_i\mid i=n_1,n_1+1,n_1+2\}$. Now
$$\td\Dlt=\sum_{i=1}^{n_1}x_i\ptl_{y_i}-\ptl_{x_{n_1+1}}\ptl_{y_{n_1+1}}+\sum_{s=n_1+2}^n
y_s\ptl_{x_s}\eqno(6.5.18)$$ and
$$\eta=\sum_{i=1}^{n_1}y_i\ptl_{x_i}+x_{n_1+1}y_{n_1+1}+\sum_{s=n_1+2}^n
x_s\ptl_{y_s}.\eqno(6.5.19)$$ Set
\begin{eqnarray*}\qquad\phi_{m_1,m_2,m_3}&=&[\prod_{s=1}^{m_2}(m_1+s)]\sum_{i=0}^\infty
\frac{x_{n_1}^{m_1+i}x_{n_1+2}^i(\ptl_{y_{n_1}}\ptl_{y_{n_1+2}})^i
(y_{n_1}^{m_2}y_{n_1+2}^{m_3})} {i!\prod_{r=1}^i(m_1+r)}\\
&=&\sum_{i=0}^\infty{m_2\choose i}
\ptl_{x_{n_1}}^{m_2-i}(x_{n_1}^{m_1+m_2})x_{n_1+2}^iy_{n_1}^{m_2-i}\ptl_{y_{n_1+2}}^i
(y_{n_1+2}^{m_3})\\
&=&(y_{n_1}\ptl_{x_{n_1}}+x_{n_1+2}\ptl_{y_{n_1+2}})^{m_2}
(x_{n_1}^{m_1+m_2}y_{n_1+2}^{m_3})\hspace{4.1cm}(6.5.20)\end{eqnarray*}
and
\begin{eqnarray*}\qquad\psi_{m_1,m_2,m_3}&=&\frac{(m_1+m_2)!\prod_{s=1}^{m_1}(m_3+s)}{m_1!}
\sum_{i=0}^\infty
\frac{x_{n_1}^ix_{n_1+2}^{m_1+i}(\ptl_{y_{n_1}}\ptl_{y_{n_1+2}})^i
(y_{n_1}^{m_2}y_{n_1+2}^{m_3})} {i!\prod_{r=1}^i(m_1+r)}
\\&=&\sum_{i=0}^{m_2}{m_2\choose i}
\frac{(m_1+m_2)!x_{n_1}^ix_{n_1+2}^{m_1+i}
y_{n_1}^{m_2-i}\ptl_{y_{n_1+2}}^{m_1+i}(y_{n_1+2}^{m_1+m_3})}
{(m_1+i)!}
\\&=&\sum_{i=0}^{m_2}{m_2\choose m_2-i}
\frac{(m_1+m_2)!x_{n_1}^i
y_{n_1}^{m_2-i}(x_{n_1+2}\ptl_{y_{n_1+2}})^{m_1+i}(y_{n_1+2}^{m_1+m_3})}
{(m_1+i)!}
\\&=&\sum_{i=0}^{m_2}
\frac{(m_1+m_2)!(y_{n_1}\ptl_{x_{n_1}})^{m_2-i}(x_{n_1}^{m_2})
(x_{n_1+2}\ptl_{y_{n_1+2}})^{m_1+i}(y_{n_1+2}^{m_1+m_3})}
{(m_2-i)!(m_1+i)!}
\\&=&\sum_{r=0}^\infty
\frac{(m_1+m_2)!(y_{n_1}\ptl_{x_{n_1}})^r(x_{n_1}^{m_2})
(x_{n_1+2}\ptl_{y_{n_1+2}})^{m_1+m_2-r}(y_{n_1+2}^{m_1+m_3})}
{r!(m_1+m_2-r)!}\\
&=&(y_{n_1}\ptl_{x_{n_1}}+x_{n_1+2}\ptl_{y_{n_1+2}})^{m_1+m_2}
(x_{n_1}^{m_2}y_{n_1+2}^{m_1+m_3}).\hspace{3.4cm}(6.5.21)\end{eqnarray*}
By Lemma 6.1.1 with
$T_1=\ptl_{x_{n_1}}\ptl_{x_{n_1+2}},\;T_1^-=\int_{(x_{n_1+1})}\int_{(x_{n_1+2})}$
(cf. (6.1.20)) and $T_2=\ptl_{y_{n_1}}\ptl_{y_{n_1+2}}$, the
polynomial solution space of (6.5.17) is
$$\msr S=\mbox{Span}\{\phi_{m_1,m_2,m_3}x_{n_1+1}^{m_4}y_{n_1+1}^{m_5},\psi_{m_1+1,m_2,m_3}x_{n_1+1}^{m_4}y_{n_1+1}^{m_5}\mid
m_i\in\mbb {N}\}.\eqno(6.5.22)$$

On the other hand,
$$\ptl_{x_{n_1}}(\phi_{0,m_2,m_3})=m_2\psi_{1,m_2-1,m_3-1},\eqno(6.5.23)$$
$$\ptl_{x_{n_1}}(\phi_{m_1,m_2,m_3})=(m_1+m_2)\phi_{m_1-1,m_2,m_3}\qquad\mbox{if}\;m_1>0,\eqno(6.5.24)$$
$$
\ptl_{x_{n_1}}(\psi_{m_1,m_2,m_3})=m_2\psi_{m_1+1,m_2-1,m_3-1},\eqno(6.5.25)$$
$$\ptl_{y_{n_1}}(\phi_{m_1,m_2,m_3})=m_2(m_1+m_2)\phi_{m_1,m_2-1,m_3},\eqno(6.5.26)$$
$$\ptl_{y_{n_1}}(\psi_{m_1,m_2,m_3})=m_2(m_1+m_2)\psi_{m_1,m_2-1,m_3}.\eqno(6.5.27)$$

In Lemma 6.1.1, $\msr A$ is a free $\msr B$-module, and so the
condition (6.1.1) yields
$$\{hg\mid g\in V,\;h\in\msr B;\;T_1(h)=0\}=\{f\in\msr A\mid
T_1(f)=0\}.\eqno(6.5.28)$$ Thus
$$\{f\in\msr A\mid (T_1+T_2)(f)=0\}=(\sum_{i=0}^\infty(-T_1^-T_2)^i)[\{g\in\msr A\mid
T_1(g)=0\}].\eqno(6.5.29)$$ Next we want to solve (6.5.14) in the
space $\msr S$ in (6.5.22). Let $T_1=\ptl_{x_{n_1}}\ptl_{x_{n_1+1}}$
and $T_2=-y_{n_1+1}\ptl_{y_{n_1}}$. First we want to find the
solution space
$$\{g\in\msr S\mid T_1(g)=0\}=\{g\in\msr S\mid
\ptl_{x_{n_1}}\ptl_{x_{n_1+1}}(g)=0\}.\eqno(6.5.30)$$ Note that
$$\ptl_{x_{n_1}}(\phi_{m_1,m_2,m_3})=0\dar
m_1=m_2m_3=0\eqno(6.5.31)$$ and
$$\ptl_{x_{n_1}}(\psi_{m_1,m_2,m_3})=0\dar m_2m_3=0.\eqno(6.5.32)$$
Moreover,
$$\phi_{0,m_2,0}=m_2!y_{n_1}^{m_2},\;\;\phi_{0,0,m_3}=y_{n_1+2}^{m_3},\eqno(6.5.33)$$
$$\psi_{m_1,m_2,0}=x_{n_1+2}^{m_1}y_{n_1}^{m_2},\;\;\psi_{m_1,0,m_3}=x_{n_1+2}^{m_1}y_{n_+2}^{m_3}.\eqno(6.5.34)$$
Since
$$\ptl_{x_{n_1+1}}(\phi_{m_1,m_2,m_3})=\ptl_{x_{n_1+1}}(\phi_{m_1,m_2,m_3})=0\eqno(6.5.35)$$
by (6.5.20) and (6.5.21), we have
\begin{eqnarray*}
\{g\in\msr S\mid
T_1(g)=0\}&=&\mbox{span}\{\phi_{m_1,m_2,m_3}y_{n_1+1}^{m_4},\psi_{m_1,m_2,m_3}y_{n_1+1}^{m_4},
y_{n_1}^{m_1}x_{n_1+1}^{m_2}y_{n_1+1}^{m_3}x_{n_1+2}^{m_4},\\
& &x_{n_1+1}^{m_1}y_{n_1+1}^{m_2}x_{n_1+2}^{m_3}y_{n_1}^{m_4}\mid
m_1,m_2,m_3,m_4\in\mbb N\}.\hspace{2.3cm}(6.5.36)\end{eqnarray*}

By  (6.5.23)-(6.5.25), we take $T_1^-$ as follows:
$$T_1^-(\phi_{m_1,m_2,m_3}x_{n_1+1}^{m_4}y_{n_1+1}^{m_5})=\frac{1}{(m_1+m_2+1)(m_4+1)}\phi_{m_1+1,m_2,m_3}x_{n_1+1}^{m_4+1}y_{n_1+1}^{m_5},
\eqno(6.5.37)$$
$$T_1^-(\psi_{1,m_2,m_3}x_{n_1+1}^{m_4}y_{n_1+1}^{m_5})
=\frac{1}{(m_2+1)(m_4+1)}\phi_{0,m_2+1,m_3+1}x_{n_1+1}^{m_4+1}y_{n_1+1}^{m_5},\eqno(6.5.38)$$
$$T_1^-(\psi_{m_1+2,m_2,m_3}x_{n_1+1}^{m_4}y_{n_1+1}^{m_5})=\frac{1}{(m_2+1)(m_4+1)}\psi_{m_1+1,m_2+1,m_3+1}
x_{n_1+1}^{m_4+1}y_{n_1+1}^{m_5}\eqno(6.5.39)$$ for $m_i\in\mbb N$
with $i\in\ol{1,5}$. In particular,
$$T_1^-\ptl_{y_{n_1}}(\phi_{m_1,m_2,m_3}x_{n_1+1}^{m_4}y_{n_1+1}^{m_5})=
\frac{m_2}{m_4+1}\phi_{m_1+1,m_2-1,m_3}x_{n_1+1}^{m_4+1}y_{n_1+1}^{m_5},\eqno(6.5.40)$$
$$T_1^-\ptl_{y_{n_1}}(\psi_{1,m_2,m_3}x_{n_1+1}^{m_4}y_{n_1+1}^{m_5})=
\frac{m_2+1}{m_4+1}\phi_{0,m_2,m_3+1}x_{n_1+1}^{m_4+1}y_{n_1+1}^{m_5},\eqno(6.5.41)$$
$$T_1^-\ptl_{y_{n_1}}(\psi_{m_1+2,m_2,m_3}x_{n_1+1}^{m_4}y_{n_1+1}^{m_5})=
\frac{m_1+m_2+2}{m_4+1}\psi_{m_1+1,m_2,m_3+1}x_{n_1+1}^{m_4+1}y_{n_1+1}^{m_5}\eqno(6.5.42)$$
for $m_i\in\mbb N$ with $i\in\ol{1,5}$.

By Lemma 6.1.1, applying the operator
$\sum_{i=0}^\infty(T_1^-y_{n_1+1}\ptl_{y_{n_1}})^i$ to (6.5.36) will
lead to the solution space of (6.5.14) in $\msr S$. We calculate
\begin{eqnarray*}& &\sum_{i=0}^\infty
(T_1^-y_{n_1+1}\ptl_{y_{n_1}})^i(\phi_{m_1,m_2,m_3})\\ &=&
\sum_{i=0}^{m_2}{m_2\choose
i}\phi_{m_1+i,m_2-i,m_3}(x_{n_1+1}y_{n_1+1})^i\\
&=&\sum_{i=0}^{m_2}{m_2\choose
i}(x_{n_1+1}y_{n_1+1})^i(y_{n_1}\ptl_{x_{n_1}}+x_{n_1+2}\ptl_{y_{n_1+2}})^{m_2-i}
(x_{n_1}^{m_1+m_2}y_{n_1+2}^{m_3})
\\&=&\eta^{m_2}(x_{n_1}^{m_1+m_2}y_{n_1+2}^{m_3}),
\hspace{10.2cm}(6.5.43)\end{eqnarray*}
\begin{eqnarray*}& &\sum_{i=0}^\infty
(T_1^-y_{n_1+1}\ptl_{y_{n_1}})^i(\psi_{m_1,m_2,m_3})\\ &=&
\sum_{i=0}^{m_1-1}{m_1+m_2\choose
i}\psi_{m_1-i,m_2,m_3+i}(x_{n_1+1}y_{n_1+1})^i\\
& &+\sum_{r=0}^{m_2}{m_1+m_2\choose
m_1+r}\phi_{r,m_2-r,m_1+m_3}(x_{n_1+1}y_{n_1+1})^{m_1+r}
\\ &=& \sum_{i=0}^{m_1-1}{m_1+m_2\choose
i}(x_{n_1+1}y_{n_1+1})^i(y_{n_1}\ptl_{x_{n_1}}+x_{n_1+2}\ptl_{y_{n_1+2}})^{m_1+m_2-i}
(x_{n_1}^{m_2}y_{n_1+2}^{m_1+m_3})\\
& &+\sum_{r=0}^{m_2}{m_1+m_2\choose
m_1+r}(x_{n_1+1}y_{n_1+1})^{m_1+r}(y_{n_1}\ptl_{x_{n_1}}+x_{n_1+2}\ptl_{y_{n_1+2}})^{m_2-r}
(x_{n_1}^{m_2}y_{n_1+2}^{m_1+m_3})
\\&=&\eta^{m_1+m_2}(x_{n_1}^{m_2}y_{n_1+2}^{m_1+m_3})
\hspace{10cm}(6.5.44)\end{eqnarray*}
 by (6.5.19). Moreover,
\begin{eqnarray*}\;\;\qquad&&[\prod_{j=1}^{m_1}(m_2+j)]\sum_{r=0}^\infty
(T_1^-y_{n_1+1}\ptl_{y_{n_1}})^r(y_{n_1}^{m_1}x_{n_1+1}^{m_2}y_{n_1+1}^{m_3})\\
&=&[\prod_{j=1}^{m_1}(m_2+j)]\sum_{r=0}^{m_1}{m_1\choose r}
\frac{y_{n_1}^{m_1-r}x_{n_1}^rx_{n_1+1}^{m_2+r}y_{n_1+1}^{m_3+r}}{\prod_{s=1}^r(m_2+s)}
\\&=&\sum_{r=0}^{m_1}{m_1+m_2\choose
m_1-r}(y_{n_1}\ptl_{x_{n_1}})^{m_1-r}(x_{n_1}^{m_1})x_{n_1+1}^{m_2+r}y_{n_1+1}^{m_3+r}
\\&=&\sum_{s=0}^{m_1+m_2}{m_1+m_2\choose
s}(y_{n_1}\ptl_{x_{n_1}})^s(x_{n_1}^{m_1})x_{n_1+1}^{m_1+m_2-r}y_{n_1+1}^{m_1+m_3-r}
\\
&=&\eta^{m_1+m_2}(x_{n_1}^{m_1}y_{n_1+1}^{m_3-m_2}),\hspace{8.7cm}(6.5.45)\end{eqnarray*}
Note that $\eta^{m_2}(x_{n_1}^{m_1}y_{n_1+1}^{m_3}y_{n_1+2}^{m_4})$
is a solution of the system of (6.5.14) and (6.5.17) for any
$m_1,m_2,m_3,m_4\in\mbb N$. Therefore, the polynomial solution space
of the system of (6.5.14) and (6.5.17) is
\begin{eqnarray*}
\qquad\;\msr
S_1&=&\mbox{span}\{\eta^{m_2}(x_{n_1}^{m_1}y_{n_1+1}^{m_3}y_{n_1+2}^{m_4}),
\eta^{m_1+m_2}(x_{n_1}^{m_1}y_{n_1+1}^{m_3-m_2}x_{n_1+2}^{m_4}),\\
& &x_{n_1+1}^{m_1}y_{n_1+1}^{m_2}x_{n_1+2}^{m_3}y_{n_1+2}^{m_4}\mid
m_1,m_2,m_3,m_4\in\mbb
N\}.\hspace{4.2cm}(6.5.46)\end{eqnarray*}Recall
$E_{n_1+1,n_1+2}|_{\msr
B}=x_{n_1+1}\ptl_{x_{n_1+2}}-\ptl_{y_{n_1+1}}\ptl_{y_{n_1+2}}$. So
the solution space of the equation (6.5.15) in $\msr S_1$ is the
span of
$$\mbox{the solutions of (6.5.15) in}\;\;
\mbox{span}\{\eta^{m_2}(x_{n_1}^{m_1}y_{n_1+1}^{m_3}y_{n_1+2}^{m_4})\mid
m_i\in\mbb N;m_1>0\},\eqno(6.5.47)$$
\begin{eqnarray*}\qquad& &\mbox{the solutions of (6.5.15) in}\;\;
\mbox{span}\{\eta^{m_1+m_2}(x_{n_1}^{m_1}y_{n_1+1}^{m_3-m_2}x_{n_1+2}^{m_4})\\
& &\hspace{2cm}\mid m_i\in\mbb
N;m_1>0;m_3-m_2<0\;\mbox{or}\;m_4>0\}\hspace{3.5cm}(6.5.48)\end{eqnarray*}
and
$$\mbox{the solutions of (6.5.15) in}\;\;
\mbox{span}\{x_{n_1+1}^{m_1}y_{n_1+1}^{m_2}x_{n_1+2}^{m_3}y_{n_1+2}^{m_4}\mid
m_i\in\mbb N\}.\eqno(6.5.49)$$ Moreover,
$$E_{n_1+1,n_1+2}[\eta^{m_2}(x_{n_1}^{m_1}y_{n_1+1}^{m_3}y_{n_1+2}^{m_4})]
=\eta^{m_2}(E_{n_1+1,n_1+2}(x_{n_1}^{m_1}y_{n_1+1}^{m_3}y_{n_1+2}^{m_4}))=0\eqno(6.5.50)$$
if and only if $m_3m_4=0$, and
\begin{eqnarray*}\qquad\qquad&
&E_{n_1+1,n_1+2}[\eta^{m_1+m_2}(x_{n_1}^{m_1}y_{n_1+1}^{m_3-m_2}x_{n_1+2}^{m_4})]\\&=&
\eta^{m_1+m_2}(E_{n_1+1,n_1+2}(x_{n_1}^{m_1}y_{n_1+1}^{m_3-m_2}x_{n_1+2}^{m_4}))=0\hspace{4.6cm}(6.5.51)\end{eqnarray*}
if and only if $m_4=0$.

Let
$T_1=\ptl_{y_{n_1+1}}\ptl_{y_{n_1+2}},\;T_1^-=\int_{(y_{n_1+1})}\int_{(y_{n_1+2})}$
and $T_2=-x_{n_1+1}\ptl_{x_{n_1+2}}$ in Lemma 6.1.1. Then (6.5.49)
is determined by
\begin{eqnarray*}\qquad&&[\prod_{j=1}^{m_3}(m_2+j)]\sum_{r=0}^\infty
(T_1^-x_{n_1+1}\ptl_{x_{n_1+2}})^r(x_{n_1+1}^{m_1}y_{n_1+1}^{m_2}x_{n_1+2}^{m_3})\\
&=&[\prod_{j=1}^{m_3}(m_2+j)]\sum_{r=0}^{m_3}{m_3\choose r}
\frac{(x_{n_1+1}^{m_1+r}y_{n_1+1}^{m_2+r}y_{n_1+2}^rx_{n_1+2}^{m_3-r})}{\prod_{s=1}^r(m_2+s)}
\\&=&\sum_{r=0}^{m_3}{m_2+m_3\choose m_3-r}
x_{n_1+1}^{m_1+r}y_{n_1+1}^{m_2+r}x_{n_1+2}^{m_3-r}\ptl_{y_{n_2+2}}^{m_3-r}(y_{n_1+2}^{m_3})
\\&=&\sum_{k=0}^{m_2+m_3}{m_2+m_3\choose k}
x_{n_1+1}^{m_1+m_2-k}y_{n_1+1}^{m_2+m_3-k}x_{n_1+2}^k\ptl_{y_{n_2+2}}^k(y_{n_1+2}^{m_3})
\\
&=&\eta^{m_2+m_3}(x_{n_1+1}^{m_1-m_3}y_{n_1+2}^{m_3}),\hspace{8.6cm}(6.5.52)\end{eqnarray*}
\begin{eqnarray*}\;\;\qquad&&[\prod_{j=1}^{m_3}(m_2+j)]\sum_{r=0}^\infty
(T_1^-x_{n_1+1}\ptl_{x_{n_1+2}})^r(x_{n_1+1}^{m_1}y_{n_2+1}^{m_2}x_{n_1+2}^{m_3})\\
&=&[\prod_{j=1}^{m_3}(m_2+j)]\sum_{r=0}^{m_3}{m_3\choose r}
\frac{(x_{n_1+1}^{m_1+r}y_{n_1+1}^ry_{n_1+2}^{m_2+r}x_{n_1+2}^{m_3-r})}{\prod_{s=1}^r(m_2+s)}
\\&=&\sum_{r=0}^{m_3}{m_3\choose m_r}
x_{n_1+1}^{m_1+r}y_{n_1+1}^{m_2+r}x_{n_1+2}^{m_3-r}\ptl_{y_{n_2+2}}^{m_3-r}(y_{n_1+2}^{m_2+m_3})
\\
&=&\eta^{m_3}(x_{n_1+1}^{m_1}y_{n_1+2}^{m_2+m_3}).\hspace{8.9cm}(6.5.53)\end{eqnarray*}
Therefore, the singular vectors in ${\msr B}$ are nonzero vectors in
\begin{eqnarray*} & &\mbox{Span}\{\eta^{m_2}(x_i^{m_1}y_j^{m_3}),
\eta^{m_1+m_2}(x_{n_1}^{m_2}y_{n_1+1}^{m_3-m_1}),\eta^{m_1+m_2}(y_{n_1+2}^{m_2}x_{n_1+1}^{m_3-m_1})
\\ & &\qquad\;\;\mid m_r\in\mbb{N};
(i,j)=(n_1,n_1+1),(n_1,n_1+2),(n_1+1,n_1+2)\}.\hspace{2.1cm}(6.5.54)\end{eqnarray*}

According to (6.4.36),
$$\td\Dlt[\eta^{m_1+m_2}(x_{n_1}^{m_2}y_{n_1+1}^{m_3-m_1})]=-(m_1+m_2)m_3
\eta^{m_1+m_2-1}(x_{n_1}^{m_2}y_{n_1+1}^{m_3-m_1}).\eqno(6.5.55)$$
 Thus we find a singular
$$\eta^{m_1+m_2}(x_{n_1}^{m_2}y_{n_1+1}^{-m_1})\in {\msr H}_{\la
m_1,m_2\ra}\eqno(6.5.56)$$ of new type if $m_1,m_2\geq 1$.
Symmetrically, $\eta^{m_1+m_2}(y_{n_1+2}^{m_2}x_{n_1+1}^{-m_1})\in
{\msr H}_{\la m_2,m_1\ra}$ is a singular vector. Recall the singular
vectors
$$f_{\la-m_1,-m_2\ra}=x_{n_1}^{m_1}y_{n_1+2}^{m_2}\in{\msr
H}_{\la -m_1,-m_2\ra},\qquad
f_{\la-m_1,m_2\ra}=x_{n_1}^{m_1}y_{n_1+1}^{m_2}\in{\msr H}_{\la
-m_1,m_2\ra},\eqno(6.5.57)$$
$$f_{\la m_1,-m_2\ra}=x_{n_1+1}^{m_1}y_{n_1+2}^{m_2}\in{\msr
H}_{\la m_1,-m_2\ra}.\eqno(6.5.58)$$ Moreover, we have the singular
vectors
$$\eta^{-\ell_1-\ell_2}(f_{\la\ell_1,\ell_2\ra})\in {\msr
H}_{\la -\ell_2,-\ell_1\ra}\qquad
\for\;\;\ell_1,\ell_2\in\mbb{Z}\;\mbox{with}\;\ell_1+\ell_2\leq
-1\eqno(6.5.59)$$ by (6.4.35)-(6.4.39). Therefore, any singular
vector in ${\msr H}$ (cf. (6.4.28)) is a nonzero weight vector in
\begin{eqnarray*}& & \mbox{Span}\{f_{\la\ell_1,\ell_2\ra},\eta^{-\ell_1'-\ell_2'}
(f_{\la\ell_1',\ell_2'\ra}),\eta^{m_1+m_2}(x_{n_1}^{m_2}y_{n_1+1}^{-m_1}),\eta^{m_1+m_2}(y_{n_1+2}^{m_2}x_{n_1+1}^{-m_1})\\
&
&\qquad\;\;\mid\ell_1,\ell_2,\ell_1',\ell_2'\in\mbb{Z},\;m_1,m_2\in\mbb{N}+1;\ell_1\leq
0\;\mbox{or}\;\ell_2\leq 0;\ell_1'+\ell_2'\leq
-1\}.\hspace{1.1cm}(6.5.60)\end{eqnarray*}

Assume $n_2=n$. We similarly find that the singular vectors in
${\msr B}$ are nonzero vectors in
$$\mbox{Span}\{\eta^{m_2}(x_{n-1}^{m_1}y_n^{m_3}),
x_n^{m_1}y_n^{m_2},\eta^{m_1+m_2}(x_{n-1}^{m_2}y_n^{m_3-m_1})\mid
m_i\in\mbb{N}\}.\eqno(6.5.61)$$ In particular, any singular vector
in ${\msr H}$ (cf. (6.4.28)) is a nonzero weight vector in
\begin{eqnarray*}\qquad& &
\mbox{Span}\{x_{n-1}^{m_1}y_n^{m_2},x_n^{m_1},\eta^{m_1+1}
(x_{n-1}^{m_1+m_2+1}y_n^{m_2}),\\
&&\qquad\;\;\eta^{m_1+m_2+2}(x_{n-1}^{m_2+1}y_n^{-m_1-1})\mid
m_1,m_2\in\mbb{N}\}.\hspace{5cm}(6.5.62)\end{eqnarray*} By the
arguments of (6.4.46)-(6.4.68), we have:\psp

{\bf Theorem 6.5.1}. {\it Suppose $n_1+1=n_2$. For
$\ell_1,\ell_2\in\mbb{Z}$ with $\ell_1+\ell_2\leq 0$ or $n_2=n$ and
$0\leq\ell_2\leq\ell_1$,  ${\msr H}_{\la\ell_1,\ell_2\ra}$ is an
infinite-dimensional irreducible highest-weight
$sl(n,\mbb{F})$-module and
$${\msr B}_{\la\ell_1,\ell_2\ra}=\bigoplus_{m=0}^\infty\eta^m({\msr
H}_{\la\ell_1-m,\ell_2-m\ra})\eqno(6.5.63)$$ is an orthogonal
decomposition of irreducible submodules. In particular, ${\msr
B}_{\la\ell_1,\ell_2\ra}={\msr H}_{\la\ell_1,\ell_2\ra}\oplus
\eta({\msr B}_{\la\ell_1-1,\ell_2-1\ra})$. The symmetric bilinear
form restricted to $\eta^m({\msr H}_{\la\ell_1-m,\ell_2-m\ra})$ is
nondegerate.

Assume $n_2<n$. For $m_1,m_2\in\mbb{N}+1$, ${\msr H}_{\la
m_1,m_2\ra}$ has exactly three singular vectors. All the submodules
${\msr H}_{\la\ell_1,\ell_2\ra}$ for $\ell_1,\ell_2\in\mbb{Z}$ such
$\ell_1+\ell_2>0$ and $\ell_1\ell_2\leq 0$ have exactly two singular
vectors. Consider $n_2=n$. For $m_1,m_2\in\mbb{N}$ with $m_1<m_2$,
${\msr H}_{\la m_1,m_2\ra}$ and ${\msr H}_{-m_2-m_1-1,m_1+1}$ are
also infinite-dimensional irreducible highest-weight
$sl(n,\mbb{F})$-modules. All submodules ${\msr H}_{\la
-m_1,m_1+m_2+1\ra}$ with $m_1,m_2\in\mbb{N}$ have  exactly two
singular vectors}.\psp

Indeed, we have more detailed information.  Suppose $n_2<n$. For
$m_1,m_2\in\mbb{N}$, ${\msr H}_{\la-m_1,-m_2\ra}$ has  a
highest-weight vector $x_{n_1}^{m_1}y_{n_1+2}^{m_2}$ of weight
$m_1\lmd_{n_1-1}-(m_1+1)\lmd_{n_1}-(m_2+1)\lmd_{n_1+1}+m_2\lmd_{n_1+2}$.
When $m_1,m_2\in\mbb{N}$ with $m_2-m_1\geq 0$, ${\msr H}_{\la
m_1,-m_2\ra}$ has a highest-weight vector
$x_{n_1+1}^{m_1}y_{n_1+2}^{m_2}$ of weight
$-(m_1+1)\lmd_{n_1}+(m_1-m_2-1)\lmd_{n_1+1}+m_2\lmd_{n_1+2}$. If
$m_1,m_2\in\mbb{N}$ with $m_1-m_2\geq 0$, ${\msr H}_{\la
-m_1,m_2\ra}$ is has a highest-weight vector
$x_{n_1}^{m_1}y_{n_1+1}^{m_2}$ of weight
$m_1\lmd_{n_1-1}+(m_2-m_1-1)\lmd_{n_1}-(m_2+1)\lmd_{n_1+1}$.

Assume $n_2=n$. For $m_1,m_2\in\mbb{N}$ with $m_2\leq m_1$, ${\msr
H}_{\la -m_1,m_2\ra}$ has a highest-weight vector
$x_{n-1}^{m_1}y_n^{m_2}$ of weight
$m_1\lmd_{n-2}+(m_2-m_1-1)\lmd_{n-1}$. Moreover,
 ${\msr H}_{\la
m,0\ra}$ has a highest-weight vector $x_{n-1}^m$ of weight
$m\lmd_{n-2}-(m+1)\lmd_{n-1}$ for $m\in\mbb{Z}$. For
$m_1,m_2\in\mbb{N}+1$, ${\msr H}_{\la m_1,m_2\ra}$ has a
highest-weight vector $\eta^{m_1+m_2}(x_{n-1}^{m_2}y_n^{-m_1})$ of
weight $m_2\lmd_{n-2}+(m_1-m_2-1)\lmd_{n-1}$. Again ${\msr
H}_{\la\ell_1,\ell_2\ra}$ has a basis of the form (6.4.70). \psp

\section{Noncanonical  Representations VI: $n_1=n_2$}

In this section, we continue the discussion from last section.

Recall $n\geq 2$. In this case,
$$\td\Dlt=\sum_{i=1}^{n_1}x_i\ptl_{y_i}+\sum_{s=n_1+1}^n
y_s\ptl_{x_s}\eqno(6.6.1)$$ and
$$\eta=\sum_{i=1}^{n_1}y_i\ptl_{x_i}+\sum_{s=n_1+1}^n
x_s\ptl_{y_s}.\eqno(6.6.2)$$

First we consider the subcase $1<n_1<n-1$. Suppose that $f\in{\cal
Q}$ is a singular vector. According to the arguments in
(6.4.1)-(6.4.4), $f$ is a rational function in
$$\{x_{n_1},x_{n_1+1},y_{n_1},y_{n_1+1},\zeta_1,\zeta_2\}
 \eqno(6.6.3)$$
(cf. (6.4.5)). Note
$$E_{n_1,n_1+1}|_{\cal
Q}=\ptl_{x_{n_1}}\ptl_{x_{n_1+1}}-\ptl_{y_{n_1}}\ptl_{y_{n_1+1}}.\eqno(6.6.4)$$
Now $E_{n_1,n_1+1}(f)=0$ implies
$$(\ptl_{x_{n_1}}\ptl_{x_{n_1+1}}-\ptl_{y_{n_1}}\ptl_{y_{n_1+1}})(f)=0,\eqno(6.6.5)$$
or equivalently,
\begin{eqnarray*}\qquad& &
(x_{n_1-1}x_{n_1+2}-y_{n_1-1}y_{n_1+2})f_{\zeta_1\zeta_2}-y_{n_1-1}f_{\zeta_1x_{n_1+1}}
-x_{n_1-1}f_{\zeta_1y_{n_1+1}}\\ &
&+y_{n_1+2}f_{\zeta_2x_{n_1}}+x_{n_1+2}f_{\zeta_2y_{n_1}}+f_{x_{n_1}x_{n_1+1}}-
f_{y_{n_1}y_{n_1+1}}=0.\hspace{3.7cm}(6.6.6)\end{eqnarray*}
According (6.3.49),
$$y_{n_1-1}=x_{n_1}^{-1}y_{n_1}x_{n_1-1}-x_{n_1}^{-1}\zeta_1,\qquad
y_{n_1+2}=x_{n_1+1}^{-1}\zeta_2+x_{n_1+1}^{-1}y_{n_1+1}x_{n_1+2}.\eqno(6.6.7)$$
Substituting (6.6.7) into (6.6.6), the coefficient of
$x_{n_1-1}x_{n_1+2}$ implies $f_{\zeta_1\zeta_2}=0$. Thus
$$f=g+h\qquad\mbox{with}\;\;g_{\zeta_2}=h_{\zeta_1}=0.\eqno(6.6.8)$$
Now (6.6.6) becomes
\begin{eqnarray*}& &
x_{n_1}^{-1}\zeta_1g_{\zeta_1x_{n_1+1}}-(x_{n_1}^{-1}y_{n_1}g_{\zeta_1x_{n_1+1}}+g_{\zeta_1y_{n_1+1}})x_{n_1-1}
+(x_{n_1+1}^{-1}y_{n_1+1}h_{\zeta_2x_{n_1}}+h_{\zeta_2y_{n_1}})x_{n_1+2}
\\ & &+x_{n_1+1}^{-1}\zeta_2h_{\zeta_2x_{n_1}}+g_{x_{n_1}x_{n_1+1}}-
g_{y_{n_1}y_{n_1+1}}+h_{x_{n_1}x_{n_1+1}}-
h_{y_{n_1}y_{n_1+1}}=0,\hspace{2.4cm}(6.6.9)\end{eqnarray*} which
implies
$$x_{n_1}^{-1}y_{n_1}g_{\zeta_1x_{n_1+1}}+g_{\zeta_1y_{n_1+1}}=0,\qquad
x_{n_1+1}^{-1}y_{n_1+1}h_{\zeta_2x_{n_1}}+h_{\zeta_2y_{n_1}}=0.\eqno(6.6.10)$$

For the representation purpose, we assume that $g$ is polynomial in
$\zeta_1$ with $g|_{\zeta_1=0}=0$ and $h$ is polynomial in
$\zeta_2$. The first assumption and the first equation in (6.6.10)
yield
$$y_{n_1}g_{x_{n_1+1}}+x_{n_1}^{-1}g_{y_{n_1+1}}=0.\eqno(6.6.11)$$
Note
$$\zeta_3=x_{n_1}y_{n_1+1}-x_{n_1+1}y_{n_1}\eqno(6.6.12)$$
is a solution of (6.6.11). Write $g$ as a function in
$\{x_{n_1},x_{n_1+1},\;y_{n_1},\zeta_1,\zeta_3\}$. Then (6.6.11)
gives $g_{x_{n_1+1}}$=0. So
$$g\;\;\mbox{is a function
in}\;\;x_{n_1},y_{n_1},\zeta_1,\zeta_3.\eqno(6.6.13)$$ Moreover,
(6.6.9) says
$$-x_{n_1}^{-1}y_{n_1}\zeta_1g_{\zeta_1\zeta_3}-y_{n_1}g_{x_{n_1}\zeta_3}-x_{n_1}g_{y_{n_1}\zeta_3}=0.\eqno(6.6.14)$$
Again we can assume that $g=\hat g+\td g$ is polynomial in
$x_{n_1},y_{n_1},\zeta_1,\zeta_3$ with $\hat g|_{\zeta_1=0}=\hat
g|_{\zeta_3=0}=0$ and $\td g_{\zeta_3}=0$. Then (6.6.14) is
equivalent to
$$y_{n_1}\zeta_1\hat g_{\zeta_1}+x_{n_1}y_{n_1}\hat g_{x_{n_1}}+x_{n_1}^2\hat g_{y_{n_1}}=0.\eqno(6.6.15)$$
Observe that $\zeta_1/x_{n_1}$ and $x_{n_1}^2-y_{n_1}^2$ are
solutions of (6.6.15). Write $\hat g$ as a function in
$\{x_{n_1},\zeta_1/x_{n_1},x_{n_1}^2-y_{n_1}^2,\zeta_3\}$. Then
(6.6.15) yields $\hat g_{x_{n_1}}=0$. This shows that $\hat g$ is a
function in $\{\zeta_1/x_{n_1},x_{n_1}^2-y_{n_1}^2,\zeta_3\}$.
 If
$\hat g\in\msr B$, then $\hat g$ is independent of $\zeta_1$. The
assumption  $\hat g|_{\zeta_1=0}=0$ implies $\hat g=0$. So the
polynomial solution of $g$ must be a polynomial in
$x_{n_1},y_{n_2},\zeta_1$ with $g_{\zeta_1}\neq 0$. Similarly, if
$h_{\zeta_2}\neq 0$ and $h|_{\zeta_2=0}=0$, the polynomial solution
of $h$ must be a polynomial in $x_{n+1},y_{n+1},\zeta_2$. Assume
$h_{\zeta_2}=0$. Then
$$h_{x_{n_1}x_{n_1+1}}-h_{y_{n_1}y_{n_1+1}}=0.\eqno(6.6.16)$$

By (6.6.2) and (6.5.17), (6.5.20), (6.5.21) with $n_1+2$ replaced by
$n_1+1$, the polynomial solution of $h$ must be in
$$\mbox{Span}\{\eta^{m_3}(x_{n_1}^{m_1}y_{n_1+1}^{m_2})\mid m_1,m_2,m_3\in\mbb{N}\}.\eqno(6.6.17)$$
Therefore, we have:\psp

{\bf Lemma 6.6.1}, {\it  Any singular vector in ${\msr B}$ must be a
nonzero weight vector in}
$$\mbox{\it Span}\{x_{n_1}^{m_1}y_{n_1}^{m_2}\zeta_1^{m_3+1},
x_{n_1+1}^{m_1}y_{n_1+1}^{m_2}\zeta_2^{m_3+1},
\eta^{m_3}(x_{n_1}^{m_1}y_{n_1+1}^{m_2})\mid
m_i\in\mbb{N}\}.\eqno(6.6.18)$$

 Note
$$x_{n_1}^{m_1}y_{n_1}^{m_2}\zeta_1^{m_3+1}\in{\msr
B}_{\la-m_1-m_3-1,m_2+m_3+1\ra},\eqno(6.6.19)$$
$$x_{n_1+1}^{m_1}y_{n_1+1}^{m_2}\zeta_2^{m_3+1}\in {\msr B}_{\la
m_1+m_3+1,-m_2-m_3-1\ra}.\eqno(6.6.20)$$ Moreover,
$$\td\Dlt(x_{n_1}^{m_1}y_{n_1}^{m_2}\zeta_1^{m_3+1})=
m_2x_{n_1}^{m_1+1}y_{n_1}^{m_2-1}\zeta_1^{m_3+1}=0\dar
m_2=0\eqno(6.6.21)$$ and
$$\td\Dlt(x_{n_1+1}^{m_1}y_{n_1+1}^{m_2}\zeta_2^{m_3+1})
=m_1x_{n_1+1}^{m_1-1}y_{n_1+1}^{m_2+1}\zeta_2^{m_3+1}=0\dar
m_1=0\eqno(6.6.22)$$ by (6.3.49) and (6.6.1). Furthermore,
$$x_{n_1}^{m_1}y_{n_1}^{m_2}\zeta_1^{m_3+1}=
\frac{\eta^{m_2}(x_{n_1}^{m_1+m_2}\zeta_1^{m_3+1})}{\prod_{r=1}^{m_2}(m_1+r)},\;\;
x_{n_1+}^{m_1}y_{n_1+1}^{m_2}\zeta_2^{m_3+1}=
\frac{\eta^{m_1}(y_{n_1+1}^{m_1+m_2}\zeta_2^{m_3+1})}{\prod_{r=1}^{m_1}(m_2+r)}
\eqno(6.6.23)$$ by (6.6.2). Indeed,
$$\eta^{m_1+1}(x_{n_1}^{m_1}\zeta_1^{m_2})=\eta^{m_1+1}(y_{n_1+1}^{m_1}\zeta_2^{m_2})=0
\qquad\for\;\;m_1,m_2\in\mbb{N}.\eqno(6.6.24)$$

Since $x_{n_1}^{m_1}y_{n_1+1}^{m_2}\in{\msr H}_{\la-m_1,-m_2\ra}$,
 (6.4.36) says that $\eta^m(x_{n_1}^{m_1}y_{n_1+1}^{m_2})$ with
$m>0$ is a singular vector only if $m=m_1+m_2+1$. But
$\eta^{m_1+m_2+1}(x_{n_1}^{m_1}y_{n_1+1}^{m_2})=0$ by (6.6.2). Thus
we have:\psp

{\bf Lemma 6.6.2}. {\it Any singular vector in ${\msr H}$ (cf.
(6.6.28)) is a nonzero weight vector in
$$\mbox{Span}\{x_{n_1}^{m_1}\zeta_1^{m_2+1},y_{n_1+1}^{m_1}\zeta_2^{m_2+1},
x_{n_1}^{m_1}y_{n_1+1}^{m_2}\mid
m_1,m_2\in\mbb{N}\}.\eqno(6.6.25)$$}\pse

 Since ${\msr B}$ is nilpotent
with respect to ${\msr G}_+$, any nonzero submodule of ${\msr B}$
has a singular vector. The above fact implies
$${\msr
H}_{\la\ell_1,\ell_2\ra}=\{0\}\qquad\for\;\;\ell_1,\ell_2\in\mbb{Z}\;\mbox{such
that}\;\ell_1+\ell_2>0.\eqno(6.6.26)$$

Observe that
$$\sum_{i=0}^{m_1}(-1)^{m_1-i}i!(m_2+i)!=\sum_{i+0}^{m_1-1}[(m_1-1)!(m_1+m_2)!+(-1)^ii!(m_2+i)!]>0\eqno(6.6.27)$$
for $m_1,m_2\in\mbb{N}.$ Thus
\begin{eqnarray*}& &(x_{n_1}^{m_1}\zeta_1^{m_2}|x_{n_1}^{m_1}\zeta_1^{m_2})\\&=&(\sum_{i=0}^{m_2}{m_2\choose
i}(-1)^ix_{n_1-1}^{m_2-i}x_{n_1}^{m_1+i}y_{n_1-1}^iy_{n_1}^{m_2-i}|\sum_{i=0}^{m_2}{m_2\choose
i}(-1)^ix_{n_1-1}^{m_2-i}x_{n_1}^{m_1+i}y_{n_1-1}^iy_{n_1}^{m_2-i})\\
&=&(-1)^{m_1}m_2!\sum_{i=0}^{m_2}(-1)^{m_2-i}(m_1+i)!(m_2-i))!\neq
0\hspace{5.1cm}(6.6.28)\end{eqnarray*} by (6.6.27). Similarly,
$(y_{n_1+1}^{m_1}\zeta_2^{m_2}|y_{n_1}^{m_1}\zeta_2^{m_2})\neq 0$.

Next we assume $n_1=n_2=1$ and $n\geq 3$. By the arguments in the
above, a singular vector in ${\msr B}$ must be a nonzero weight
vector in
$$\mbox{Span}\{x_{n_1+1}^{m_1}y_{n_1+1}^{m_2}\zeta_2^{m_3+1},
\eta^{m_3}(x_{n_1}^{m_1}y_{n_1+1}^{m_2})\mid
m_i\in\mbb{N}\}.\eqno(6.6.29)$$ Thus any singular vector in ${\msr
H}$ (cf. (6.4.28)) is a nonzero weight vector in
$$\mbox{Span}\{y_{n_1+1}^{m_1}\zeta_2^{m_2+1},
x_{n_1}^{m_1}y_{n_1+1}^{m_2}\mid m_1,m_2\in\mbb{N}\}.\eqno(6.6.30)$$
The above fact implies
$${\msr
H}_{\la\ell_1,\ell_2\ra}=\{0\}\qquad\for\;\;\ell_1,\ell_2\in\mbb{Z}\;\mbox{such
that}\;\ell_1+\ell_2>0\;\mbox{or}\;\;\ell_2>0.\eqno(6.6.31)$$

Consider the subcase $n_1=n_2=n-1$ and $n\geq 3$. A singular vector
in ${\msr B}$ must be a nonzero weight vector in
$$\mbox{Span}\{x_{n_1}^{m_1}y_{n_1}^{m_2}\zeta_1^{m_3+1},
\eta^{m_3}(x_{n_1}^{m_1}y_{n_1+1}^{m_2})\mid
m_i\in\mbb{N}\}.\eqno(6.6.32)$$ Thus any singular vector in ${\msr
H}$ (cf. (6.4.28)) is a nonzero weight vector in
$$\mbox{Span}\{x_{n_1}^{m_1}\zeta_1^{m_2+1},
x_{n_1}^{m_1}y_{n_1+1}^{m_2}\mid m_1,m_2\in\mbb{N}\}.\eqno(6.6.33)$$
The above fact implies
$${\msr
H}_{\la\ell_1,\ell_2\ra}=\{0\}\qquad\for\;\;\ell_1,\ell_2\in\mbb{Z}\;\mbox{such
that}\;\ell_1+\ell_2>0\;\mbox{or}\;\;\ell_1>0.\eqno(6.6.34)$$

Suppose $n_1=n_2=1$ and $n=2$. A singular vector in ${\msr B}$ must
be a nonzero weight vector in
$$\mbox{Span}\{\eta^{m_3}(x_1^{m_1}y_2^{m_2})\mid
m_i\in\mbb{N}\}.\eqno(6.6.35)$$ Thus any singular vector in ${\msr
H}$ (cf. (6.4.28)) is a nonzero weight vector in
$$\mbox{Span}\{x_1^{m_1}y_2^{m_2}\mid m_1,m_2\in\mbb{N}\}.\eqno(6.6.36)$$
The above fact implies
$${\msr
H}_{\la\ell_1,\ell_2\ra}=\{0\}\qquad\for\;\;\ell_1,\ell_2\in\mbb{Z}\;\mbox{such
that}\;\ell_1>0\;\mbox{or}\;\ell_2>0.\eqno(6.6.37)$$

Finally, we assume $n_1=n_2=n$. A singular vector in ${\msr B}$ must
be a nonzero weight vector in
$$\mbox{Span}\{x_{n_1}^{m_1}y_{n_1}^{m_2}\zeta_1^{m_3}\mid m_i\in\mbb{N}\}.\eqno(6.6.38)$$
 Thus any singular vector in
${\msr H}$ (cf. (6.4.29)) is a nonzero weight vector in
$$\mbox{Span}\{x_{n_1}^{m_1}\zeta_1^{m_2}\mid m_1,m_2\in\mbb{N}\}.\eqno(6.6.39)$$
The above fact implies
$${\msr
H}_{\la\ell_1,\ell_2\ra}=\{0\}\qquad\for\;\;\ell_1,\ell_2\in\mbb{Z}\;\mbox{such
that}\;\ell_1+\ell_2>0.\eqno(6.6.40)$$

By the arguments of (6.4.46)-(6.4.68),  we obtain:\psp

{\bf Theorem 6.6.3}. {\it Suppose $n_1=n_2$. Let
 $\ell_1,\ell_2\in\mbb{Z}$ such that $\ell_2\geq 0$ when $n_1=n$, or $\ell_1+\ell_2\leq 0$ and: (a)
 $\ell_2\leq 0 $  if $n_1=1$ and $n\geq 3$;
 (b) $\ell_1\leq 0$ if $n_1=n-1$ and $n\geq 3$; (c) $\ell_1,\ell_2\leq
 0$ when $n_1=1$ and $n=2$. Then
   ${\msr H}_{\la\ell_1,\ell_2\ra}$ is an
irreducible highest-weight $sl(n,\mbb{F})$-module and
$${\msr B}_{\la\ell_1,\ell_2\ra}=\bigoplus_{m=0}^\infty\eta^m({\msr
H}_{\la\ell_1-m,\ell_2-m\ra})\eqno(6.6.41)$$ is an orthogonal
decomposition of irreducible submodules. The symmetric bilinear form
restricted to $\eta^m({\msr H}_{\la\ell_1-m,\ell_2-m\ra})$ is
nondenerate.
 In particular, ${\msr
B}_{\la\ell_1,\ell_2\ra}={\msr H}_{\la\ell_1,\ell_2\ra}\oplus
\eta({\msr B}_{\la\ell_1-1,\ell_2-1\ra})$. If the conditions fails,
${\msr H}_{\la\ell_1,\ell_2\ra}=\{0\}$}.\psp

Suppose $n_1<n-1$. Let $m_1,m_2\in\mbb{N}$. The subspace ${\msr
H}_{\la-m_1,-m_2\ra}$ has  a highest-weight vector
$x_{n_1}^{m_1}y_{n_1+1}^{m_2}$ of weight
$m_1(1-\dlt_{1,n_1})\lmd_{n_1-1}-(m_1+m_2+2)\lmd_{n_1}+m_2\lmd_{n_1+1}$.
If $n_1\geq 2$, the subspace ${\msr H}_{\la-m_1-m_2-1,m_2+1\ra}$ has
a highest-weight vector $x_{n_1}^{m_1}\zeta_1^{m_2+1}$ of weight
$(m_2+1)\lmd_{n_1-2}-m_1\lmd_{n_1-1}-(m_1+m_2+3)\lmd_{n_1}$. The
subspace ${\msr H}_{\la m_1+1,-m_2-m_1-1\ra}$ has  a highest-weight
vector $y_{n_1+1}^{m_2}\zeta_2^{m_1+1}$ of weight
$-(m_1+m_2+3)\lmd_{n_1}+m_2\lmd_{n_1+1}-(m_1+1)(1-\dlt_{n_1,n-2})\lmd_{n_1+2}$.

Consider $n_1=n-1$. The subspace ${\msr H}_{\la-m_1,-m_2\ra}$ has a
highest-weight vector $x_{n_1}^{m_1}y_{n_1+1}^{m_2}$ of weight
$m_1(1-\dlt_{n,2})\lmd_{n-2}-(m_1+m_2+2)\lmd_{n-1}$. If $n\geq 3$,
the subspace ${\msr H}_{\la-m_1-m_2-1,m_2+1\ra}$ has  a
highest-weight vector $x_{n_1}^{m_1}\zeta_1^{m_2+1}$ of weight
$(m_2+1)(1-\dlt_{n,3})\lmd_{n-3}-m_1\lmd_{n-2}-(m_1+m_2+3)\lmd_{n-1}$.

Assume $n_1=n$. The subspace ${\msr H}_{\la-m_1-m_2,m_2\ra}$ has a
highest-weight vector $x_n^{m_1}\zeta_1^{m_2}$ of weight
$m_2(1-\dlt_{n,2})\lmd_{n-2}+m_1\lmd_{n-1}$.\psp

Now we want to find an explicit expression for ${\msr
H}_{\la\ell_1,\ell_2\ra}$ when it is irreducible. First we assume
$n_1<n$. Note
$$-E_{n_1+r,i}|_{\msr
B}=x_ix_{n_1+r}-y_iy_{n_r}\qquad\for\;\;i\in\ol{1,n_1},\;r\in\ol{1,n-n_1}\eqno(6.6.42)$$
by (6.3.36)-(6.3.38). Thus
$$(x_ix_{n_1+r}-y_iy_{n_r})({\msr
H}_{\la\ell_1,\ell_2\ra})\subset{\msr
H}_{\la\ell_1,\ell_2\ra}.\eqno(6.6.43)$$ For $m_1,m_2\in\mbb{N}$, we
set
\begin{eqnarray*}\qquad{\msr H}_{\la
-m_1,-m_2\ra}'&=&\mbox{Span}\{[\prod_{r=1}^{n_1}x_r^{l_r}]
[\prod_{s=1}^{n-n_1}y_{n_1+s}^{k_s}][\prod_{r=1}^{n_1}\prod_{s=1}^{n-n_1}(x_rx_{n_1+s}-
y_ry_{n_1+s})^{l_{r,s}}]\\& &\qquad\;\;\mid
l_r,k_s,l_{r,s}\in\mbb{N};\sum_{r=1}^{n_1}l_r=m_1;\sum_{s=1}^{n-n_1}k_s=m_2\}.\hspace{2.2cm}(6.6.44)\end{eqnarray*}
By (6.6.1), ${\msr H}_{\la m_1,m_2\ra}'\subset{\msr H}_{\la
-m_1,-m_2\ra}$. Moreover,
$$E_{i,r}|_{\msr
B}=-x_r\ptl_{x_i}-y_r\ptl_{y_i}-\dlt_{i,r}\qquad\for\;\;i,r\in\ol{1,n_1},\eqno(6.6.45)$$
$$E_{n_1+j,n_1+s}|_{\msr
B}=x_{n_1+j}\ptl_{x_{n_1+s}}+y_{n_1+j}\ptl_{y_{n_1+s}}+\dlt_{j,s}\qquad\for\;\;j,s\in\ol{1,n-n_1}\eqno(6.6.46)$$
and
$$E_{i,n_1+s}|_{\msr
B}=\ptl_{x_i}\ptl_{x_{n_1+s}}-\ptl_{y_i}\ptl_{y_{n_1+s}}\qquad\for\;\;i\in\ol{1,n_1},\;s\in\ol{1,n-n_1}
\eqno(6.6.47)$$ by (6.3.36)-(6.3.38). It is obvious that
$$E_{p,p'}({\msr H}_{\la
-m_1,-m_2\ra}'),E_{n_1+q,n_1+q'}({\msr H}_{\la
-m_1,-m_2\ra}'),E_{n_1+q,p'}({\msr H}_{\la-m_1,-m_2\ra}')\subset
{\msr H}_{\la -m_1,-m_2\ra}'\eqno(6.6.48)$$ for $p,p'\in\ol{1,n_1}$
and $q,q'\in\ol{1,n-n_1}$ by (6.6.42), (6.6.45) and (6.6.46).
Furthermore,
\begin{eqnarray*}& &E_{p,n_1+q}\{[\prod_{r=1}^{n_1}x_r^{l_r}]
[\prod_{s=1}^{n-n_1}y_{n_1+s}^{k_s}][\prod_{r=1}^{n_1}\prod_{s=1}^{n-n_1}(x_rx_{n_1+s}-
y_ry_{n_1+s})^{l_{r,s}}]\}\\
&=&(\ptl_{x_p}\ptl_{x_{n_1+q}}-\ptl_{y_p}\ptl_{y_{n_1+q}})\{[\prod_{r=1}^{n_1}x_r^{l_r}]
[\prod_{s=1}^{n-n_1}y_{n_1+s}^{k_s}][\prod_{r=1}^{n_1}\prod_{s=1}^{n-n_1}(x_rx_{n_1+s}-
y_ry_{n_1+s})^{l_{r,s}}]\}\\
&=&\big(l_p\sum_{\iota=1}^{n_1}\frac{l_{\iota,q}x_\iota}{x_p(x_\iota
x_{n_1+q}-y_\iota
y_{n_1+q})}+k_q\sum_{\iota_1=1}^{n-n_1}\frac{l_{p,\iota_1}y_{n_1+\iota_1}}{y_{n_1+q}(x_p
x_{n_1+\iota_1}-y_p y_{n_1+\iota_1})}\\ & &+ \frac{l_{p,q}}{x_p
x_{n_1+q}-y_p y_{n_1+q}}
+\sum_{\iota_1=1}^{n_1}\sum_{\iota_2=1}^{n-n_1}\frac{l_{\iota_1,q}l_{p,\iota_2}(x_{\iota_1}x_{n_1+\iota_2}
-y_{\iota_1}y_{n_1+\iota_2})}{(x_{\iota_1}x_{n_1+q}
-y_{\iota_1}y_{n_1+q})(x_px_{n_1+\iota_2}
-y_py_{n_1+\iota_2})}\big)\\& &\times
 [\prod_{r=1}^{n_1}x_r^{l_r}]
[\prod_{s=1}^{n-n_1}y_{n_1+s}^{k_s}][\prod_{r=1}^{n_1}\prod_{s=1}^{n-n_1}(x_rx_{n_1+s}-
y_ry_{n_1+s})^{l_{r,s}}]\in {\msr H}_{\la
-m_1,-m_2\ra}'.\hspace{1.9cm}(6.6.49)\end{eqnarray*} Thus ${\msr
H}_{\la -m_1,-m_2\ra}'$ is a nonzero submodule of ${\msr H}_{\la
-m_1,-m_2\ra}$, which is irreducible. So we have ${\msr H}_{\la
-m_1,-m_2\ra}'={\msr H}_{\la -m_1,-m_2\ra}$.

Next we assume $n_1>1$. We let
\begin{eqnarray*}& &{\msr H}_{\la
-m_1-m_2,m_2\ra}'\\&=&\mbox{Span}\{[\prod_{r=1}^{n_1}x_r^{l_r}]
[\prod_{1\leq p<q\leq
n_1}(x_py_q-x_qy_p)^{k_{p,q}}][\prod_{r=1}^{n_1}\prod_{s=1}^{n-n_1}(x_rx_{n_1+s}-
y_ry_{n_1+s})^{l_{r,s}}]\\& &\qquad\;\;\mid
l_r,k_{p,q},l_{r,s}\in\mbb{N};\sum_{r=1}^{n_1}l_r=m_1;\sum_{1\leq
p<q\leq n_1}k_{p,q}=m_2\}.\hspace{3.6cm}(6.6.50)\end{eqnarray*} For
convenience, we treat $k_{p,q}=-k_{q,p}$. To prove ${\msr H}_{\la
-m_1-m_2,m_2\ra}'={\msr H}_{\la -m_1-m_2,m_2\ra}$, it is enough to
calculate
\begin{eqnarray*}& &
\ptl_{x_i}[\prod_{1\leq p<q\leq
n_1}(x_py_q-x_qy_p)^{k_{p,q}}]\ptl_{x_{n_1+j}}[\prod_{r=1}^{n_1}\prod_{s=1}^{n-n_1}(x_rx_{n_1+s}-
y_ry_{n_1+s})^{l_{r,s}}]\\ & &-\ptl_{y_i}[\prod_{1\leq p<q\leq
n_1}(x_py_q-x_qy_p)^{k_{p,q}}]\ptl_{y_{n_1+j}}[\prod_{r=1}^{n_1}\prod_{s=1}^{n-n_1}(x_rx_{n_1+s}-
y_ry_{n_1+s})^{l_{r,s}}]\\ &=&\big(\sum_{i\neq
\iota_1\in\ol{1,n_1}}\sum_{\iota_2=1}^{n_1}\frac{k_{i,i_1}l_{\iota_2,j}(x_{\iota_1}y_{\iota_2}
-x_{\iota_2}y_{\iota_1})}{(x_iy_{\iota_1}-x_{\iota_1}y_i)(x_{\iota_2}x_{n_1+j}-y_{\iota_2}y_{n_1+j})}\big)
\\ & &\times
[\prod_{1\leq p<q\leq
n_1}(x_py_q-x_qy_p)^{k_{p,q}}][\prod_{r=1}^{n_1}\prod_{s=1}^{n-n_1}(x_rx_{n_1+s}-
y_ry_{n_1+s})^{l_{r,s}}]\hspace{3cm}(6.6.51)\end{eqnarray*} for
$i\in\ol{1,n_1}$ and $j\in\ol{1,n_1}$ by (6.6.45)-(6.6.47) and
(6.6.49).

Symmetrically, similar expression can obtained for ${\msr H}_{\la
m_2, -m_1-m_2\ra}$. In summary, we have:\psp

{\bf Theorem 6.6.4}. {\it Let $m_1,m_2\in\mbb{N}$. If $n_1=n_2<n$,
we have
\begin{eqnarray*}\qquad{\msr H}_{\la
-m_1,-m_2\ra}&=&\mbox{Span}\{[\prod_{r=1}^{n_1}x_r^{l_r}]
[\prod_{s=1}^{n-n_1}y_{n_1+s}^{k_s}][\prod_{r=1}^{n_1}\prod_{s=n_1+1}^n(x_rx_s-
y_ry_s)^{l_{r,s}}]\\& &\qquad\;\;\mid
l_r,k_s,l_{r,s}\in\mbb{N};\sum_{r=1}^{n_1}l_r=m_1;\sum_{s=1}^{n-n_1}k_s=m_2\}.\hspace{2.1cm}(6.6.52)\end{eqnarray*}
When $2\leq n_1=n_2<n$,
\begin{eqnarray*}&&{\msr H}_{\la
-m_1-m_2,m_2\ra}\\&=&\mbox{Span}\{[\prod_{r=1}^{n_1}x_r^{l_r}]
[\prod_{1\leq p<q\leq
n_1}(x_py_q-x_qy_p)^{k_{p,q}}][\prod_{r=1}^{n_1}\prod_{s=n_1+1}^n(x_rx_s-
y_ry_s)^{l_{r,s}}]\\& &\qquad\;\;\mid
l_r,k_{p,q},l_{r,s}\in\mbb{N};\sum_{r=1}^{n_1}l_r=m_1;\sum_{1\leq
p<q\leq n_1}k_{p,q}=m_2\}.\hspace{3.7cm}(6.6.53)\end{eqnarray*} In
the case $n_1=n_2<n-1$,
\begin{eqnarray*}& &{\msr H}_{\la
m_2,-m_1-m_2\ra}=\mbox{Span}\{ [\prod_{n_1+1\leq p<q\leq
n}(x_py_q-x_qy_p)^{k_{p,q}}][\prod_{r=1}^{n_1}\prod_{s=n_1+1}^n(x_rx_s-
y_ry_s)^{l_{r,s}}]\\& &\times
[\prod_{r=1}^{n-n_1}y_{n_1+r}^{l_r}]\mid
l_r,k_{p,q},l_{r,s}\in\mbb{N};\sum_{r=1}^{n_1}l_r=m_1;\sum_{n_1+1\leq
p<q\leq n}k_{p,q}=m_2\}.\hspace{2.35cm}(6.6.54)\end{eqnarray*}
 If $n_1=n_2=n$,
\begin{eqnarray*}{\msr H}_{\la
-m_1-m_2,m_2\ra}&=&\mbox{Span}\{[\prod_{r=1}^nx_r^{l_r}]
[\prod_{1\leq p<q\leq n}(x_py_q-x_qy_p)^{k_{p,q}}]\\&
&\qquad\;\;\mid
l_r,k_{p,q}\in\mbb{N};\sum_{r=1}^nl_r=m_1;\sum_{1\leq p<q\leq
n}k_{p,q}=m_2\}.\hspace{2.2cm}(6.6.55)\end{eqnarray*}}

\section{Extensions of the Projective Representations}

In this section, we construct a functor from the category of
$A_{n-1}$-modules to the category of $A_n$-modules, which is related
to $n$-dimensional projective transformations.

Denote by $GL(n+1,\mbb{F})$ the group of $(n+1)\times (n+1)$
invertible matrices.  A projective transformation on $\mathbb{F}^n$
for a field $\mathbb{F}$ is given by
$$u\mapsto \frac{Au+\vec b}{\vec c\:^t u+d}\qquad\mbox{for}\;\;u\in
\mathbb{F}^n,\eqno(6.7.1)$$ where all the vector in $\mathbb{F}^n$
are in column form and
$$\left(\begin{array}{cc}A&\vec b\\ \vec c\:^t&
d\end{array}\right)\in GL(n+1,\mbb{F}).\eqno(6.7.2)$$ It is
well-known that a transformation of mapping straight lines to lines
must be a projective transformation. The group of projective
transformations is the fundamental group of $n$-dimensional
projective geometry.  Physically, the group  with $n=4$ and
$\mathbb{F}=\mathbb{R}$ consists of all the transformations of
keeping free particles including light signals moving with constant
velocities along straight lines (e.g., cf. [GWZ1-2]). Based on the
embeddings of the poincar\'{e} group and De Sitter group into the
projective group with $n=4$ and $\mathbb{F}=\mathbb{R}$, Guo, Wu and
Zhou [GWZ1-2] proposed three kinds of special relativity.

 Now we consider the Lie
algebra
$$sl(n+1,\mbb{F})=\sum_{1\leq i<j\leq
n+1}(\mbb{F}E_{i,j}+\mbb{F}E_{j,i})+\sum_{r=1}^n\mbb{F}(E_{r,r}-E_{r+1,r+1}).\eqno(6.7.3)$$
Recall that
$$D=\sum_{r=1}^nx_r\ptl_{x_r}\eqno(6.7.4)$$
is the degree operator on the polynomial algebra ${\msr
A}=\mbb{F}[x_1,...,x_n]$.  Note
$$[\partial_{x_{j}},x_iD]=\left\{\begin{array}{llll}
x_{i}\partial_{x_{j}},  & i \neq j,\\
D+x_{i}\partial_{x_{i}}, & i =j,\end{array}\right.\eqno(6.7.5)$$
$$[x_{i}\partial_{x_{j}},x_kD]=\delta_{j,k}x_iD,
\;\;
[x_{i}\partial_{x_{j}},\partial_{x_{k}}]=-\delta_{i,k}\partial_{x_{j}}\eqno(6.7.6)$$
for $i,j,k\in\ol{1,n}$. Thus we  have the following inhomogeneous
representation of $sl(n+1,\mbb{F})$ on ${\msr A}$ determined by
$$E_{i,j}|_{\msr A}=x_i\ptl_{x_j},\qquad E_{i,n+1}|_{\msr A}=x_iD,\eqno(6.7.7)$$
$$E_{n+1,i}|_{\msr A}=-\ptl_{x_i},\;\;(E_{i,i}-E_{n+1,n+1})|_{\msr A}=D+x_i\ptl_{x_i}\eqno(6.7.8)$$
for $i,j\in\ol{1,n}$ with $i\neq j$.
 The above
representation is exactly the one of $sl(n+1,\mbb{F})$ induced by
the projective transformations in (6.7.1).

Recall the classical Witt algebra $\mbb W_n=\sum_{i=1}^n{\msr
A}\ptl_{x_i}$ with the following Lie bracket:
$$[\sum\limits_{i=1}^{n}f_{i}\ptl_{x_i},\sum\limits_{j=1}^{n}g_j\ptl_{x_j}]=
\sum\limits_{i,j=1}^{n}(f_{j}\ptl_{x_j}
(g_{i})-g_{j}\ptl_{x_j}(f_{i}))\ptl_{x_i}\eqno(6.7.9)$$ (cf. Example
1.2.1). Acting on the entries of the elements of the Lie algebra
$gl(n,{\msr A})$, $\mbb W_n$ becomes a Lie subalgebra of the Lie
algebras of derivations of $gl(n,{\msr A})=\sum_{i,j=1}^n{\msr
A}E_{i,j}$. In particular,
$$\widehat{\mbb W}_n=\mbb W_n\oplus gl(n,{\msr A})\eqno(6.7.10)$$
becomes a Lie algebra with the Lie bracket
$$[d_1\oplus A_1,d_2\oplus A_2]=[d_1,d_2]\oplus[d_1(A_2)-d_2(A_1)+[A_1,A_2]]\eqno(6.7.11)$$
for $d_1,d_2\in\mbb W_n$ and $A_1,A_2\in gl(n,{\msr A})$.\psp

{\bf Lemma 6.7.1 (Shen)}. {\it The  map $\Im$ given by
$$\Im(\sum\limits_{i=1}^{n}f_{i}\ptl_{x_i})=
\sum\limits_{i=1}^{n}f_{i}\ptl_{x_i}\oplus\sum_{i,j=1}^n\ptl_{x_i}(f_j)
E_{i,j}\eqno(6.7.12)$$ is
 a Lie algebra monomorphism from $\mbb W_n$ to $\widehat{\mbb
 W}_n.$}

{\it Proof}. Note that
\begin{eqnarray*} & &\Im([\sum_{i=1}^nf_i\ptl_{x_i},\sum_{j=1}^ng_j\ptl_{x_j}])
\\&=&\Im(\sum\limits_{i,k=1}^{n}(f_{k}\ptl_{x_k}
(g_{i})-g_{k}\ptl_{x_k}(f_{i}))\ptl_{x_i})\\
&=&\sum_{i,k=1}^{n}(f_{k}\ptl_{x_k}
(g_{i})-g_{k}\ptl_{x_k}(f_{i}))\ptl_{x_i})\oplus\sum_{i,j,k=1}^n\ptl_{x_i}[f_{k}\ptl_{x_k}(g_{j})-g_{k}\ptl_{x_k}(f_{j})]E_{i,j}\hspace{1.5cm}(6.7.13)
\end{eqnarray*}
and
\begin{eqnarray*} &
&[\Im(\sum_{i=1}^nf_i\ptl_{x_i}),\Im(\sum_{j=1}^ng_j\ptl_{x_j})]\\
&=&[\sum_{i=1}^nf_i\ptl_{x_i}\oplus\sum_{i,r=1}^n\ptl_{x_i}(f_r)E_{i,r},\sum_{j=1}^ng_j\ptl_{x_j}\oplus\sum_{j,s=1}^n\ptl_{x_j}(g_s)E_{j,s}]
\\ &=&\sum_{i,k=1}^{n}(f_{k}\ptl_{x_k}
(g_{i})-g_{k}\ptl_{x_k}(f_{i}))\ptl_{x_i})\oplus[\sum_{i,j,s=1}^nf_i\ptl_{x_i}\ptl_{x_j}(g_s)E_{j,s}-\sum_{i,j,r}\sum_{i,r=1}^ng_j
\ptl_{x_j}\ptl_{x_i}(f_r)E_{i,r}\\ &
&+\sum_{i,r,s=1}^n\ptl_{x_i}(f_r)\ptl_{x_r}(g_s)E_{i,s}-\sum_{j,r,s=1}^n\ptl_{x_j}(g_s)\ptl_{x_s}(f_r)E_{j,r}].\hspace{4.1cm}(6.7.14)
\end{eqnarray*}
Moreover,
\begin{eqnarray*} &
&\sum_{i,j,s=1}^nf_i\ptl_{x_i}\ptl_{x_j}(g_s)E_{j,s}-\sum_{i,j,r}\sum_{i,r=1}^ng_j
\ptl_{x_j}\ptl_{x_i}(f_r)E_{i,r}\\ &
&+\sum_{i,r,s=1}^n\ptl_{x_i}(f_r)\ptl_{x_r}(g_s)E_{i,s}-\sum_{j,r,s=1}^n\ptl_{x_j}(g_s)\ptl_{x_s}(f_r)E_{j,r}\\
&=&\sum_{i,j,k=1}^n[f_k\ptl_{x_k}\ptl_{x_i}(g_j)-g_k
\ptl_{x_k}\ptl_{x_i}(f_j)+\ptl_{x_i}(f_k)\ptl_{x_k}(g_j)-\ptl_{x_i}(g_k)\ptl_{x_k}(f_j)]E_{i,j}
\\ &=&\sum_{i,j,k=1}^n\ptl_{x_i}[f_{k}\ptl_{x_k}(g_{j})-g_{k}\ptl_{x_k}(f_{j})]E_{i,j}\hspace{7.3cm}(6.7.15)
\end{eqnarray*}
Thus
$$\Im([\sum_{i=1}^nf_i\ptl_{x_i},\sum_{j=1}^ng_j\ptl_{x_j}])=[\Im(\sum_{i=1}^nf_i\ptl_{x_i}),\Im(\sum_{j=1}^ng_j\ptl_{x_j})].\qquad\Box\eqno(6.7.16)$$
\psp

The above lemma was an important discovery of Shen [Sg].  Now we
have the Lie algebra monomorphism $\nu: sl(n+1,\mbb{F})\rta
\widehat{\mbb W}_n$ given by
$$\nu(\xi)=\Im(\xi|_{\msr A})\qquad\for\;\;\xi\in
sl(n+1,\mbb{F}).\eqno(6.7.17)$$ In particular,
$$\nu(E_{i,n+1})=\Im(x_iD)=x_iD\oplus
(\sum_{r=1}^n(x_iE_{r,r}+x_rE_{i,r}))\qquad\for\;\;i\in\ol{1,n}\eqno(6.7.18)$$
and
$$\nu(E_{i,j})=x_i\ptl_{x_j}\oplus
E_{i,j},\;\nu(E_{i,i}-E_{j,j})=(x_i\ptl_{x_i}-x_j\ptl_{x_j})\oplus
(E_{i,i}-E_{j,j})\eqno(6.7.19)$$ for $i,j\in\ol{1,n}$ with $i\neq
j$.

Let $M$ be an  $gl(n,\mbb{F})$-module. Set
$$\wht M={\msr A}\otimes_{\mbb{F}}M.\eqno(6.7.20)$$
According to (6.7.11), $\wht M$ becomes a $ \widehat{\mbb
W}_n$-module with the action
$$(d\oplus fA)(g\otimes v)=d(g)\otimes v+fg\otimes
A(v)\eqno(6.7.21)$$for $d\in\mbb W_n,\;f,g\in{\msr A},\;A\in
gl(n,\mbb{F})$ and $v\in M.$ Thus we have the following
$sl(n+1,\mbb{F})$-module structure on $\wht M$:
$$\xi(g\otimes v)=\nu(\xi)(g\otimes v)\qquad\for\;\;\xi\in
sl(n+1,\mbb{F}),\;g\in{\msr A},\;v\in M.\eqno(6.7.22)$$

For convenience, we denote
$${\msr G}_0=sl(n,\mbb{F})+\mbb{F}(E_{n,n}-E_{n+1,n+1}),\;{\msr
G}_+=\sum_{i=1}^n\mbb{F}E_{n+1,i},\;{\msr
G}_-=\sum_{i=1}^n\mbb{F}E_{i,n+1}.\eqno(6.7.23)$$ Then  ${\msr
G}_\pm$ are abelian Lie subalgebras of $sl(n+1,\mbb{F})$. Moreover,
$$\nu(E_{n+1,i})=-\ptl_{x_i}\qquad\for\;\;i\in\ol{1,n}.\eqno(6.7.24)$$
Thus
$${\msr G}_+(1\otimes M)=\{0\}.\eqno(6.7.25)$$
Furthermore,
$${\msr
G}_0=sl(n,\mbb{F})+\mbb{F}(\sum_{r=1}^nE_{r,r}-nE_{n+1,n+1})\eqno(6.7.26)$$
and
$$\nu(\sum_{r=1}^nE_{r,r}-nE_{n+1,n+1})=(n+1)(D\oplus\sum_{r=1}^nE_{r,r})\in\wht{\mbb
W}_n.\eqno(6.7.27)$$ So
$$U({\msr G}_0)(1\otimes M)=1\otimes U(sl(n,\mbb{F}))(M)=1\otimes
M.\eqno(6.7.28)$$ Therefore, the $sl(n+1,\mbb{F})$-submodule
$$U({\msr G})(1\otimes M)=U({\msr G}_-)(1\otimes M).\eqno(6.7.29)$$
By (6.7.24), we have: \psp

{\bf Proposition 6.7.2}. {\it The map $M\mapsto U({\msr
G}_-)(1\otimes M)$ gives rise to a functor from the category of
$gl(n,\mbb{F})$-modules to the category of
$sl(n+1,\mbb{F})$-modules. In particular, it maps irreducible
$gl(n,\mbb{F})$-modules to irreducible
$sl(n+1,\mbb{F})$-modules.}\psp

From another point view, $U({\msr G}_-)(1\otimes M)$ is a polynomial
extension from $gl(n,\mbb{F})$-module $M$ to an
$sl(n+1,\mbb{F})$-module. Next we want to study when $\wht M=U({\msr
G}_-)(1\otimes M)$. Set
$$E_i=E_{i,i}-\frac{1}{1+n}\sum_{r=1}^{n+1}E_{r,r}\qquad\for\;\;i\in\ol{1,n}.\eqno(6.7.30)$$
Then we have a Lie algebra isomorphism $\vs: gl(n,\mbb{F})\rta \msr
G_0$ determined by
$$\vs(E_{r,s})=E_{r,s},\;\;\vs(E_{i,i})=E_i\qquad\for\;\;i,r,s\in\ol{1,n},\;r\neq
s.\eqno(6.7.31)$$ Moreover,
$$E_i=E_{i,i}-E_{n+1,n+1}-\frac{1}{1+n}(\sum_{r=1}^nE_{r,r}-nE_{n+1,n+1})\qquad\for
\;i\in\ol{1,n}.\eqno(6.7.32)$$ According to (6.7.8),
$$E_i|_{\msr A}=x_i\ptl_{x_i}\qquad\for
\;i\in\ol{1,n}.\eqno(6.7.33)$$ Thus
$$\nu(E_i)=x_i\ptl_{x_i}\oplus E_{i,i}\qquad\for
\;i\in\ol{1,n}.\eqno(6.7.34)$$
 We calculate
$$[E_i,E_{j,n+1}]=\dlt_{i,j}E_{j,n+1}\qquad\for
\;i,j\in\ol{1,n}.\eqno(6.7.35)$$ In the rest of this section, we use
$\vs$ as the identification of $gl(n,\mbb{F})$ with $\msr G_0$. In
particular,
$$E_{i,j}|_{\msr
A}=x_i\ptl_{x_j}\qquad\for\;\;i,j\in\ol{1,n}\eqno(6.7.36)$$ by
(6.7.7) and (6.7.33). So $\wht M$ is a usual tensor module of
$gl(n,\mbb{F})$.

 Denote
$$E^\al=E_{1,n+1}^{\al_1}E_{2,n+1}^{\al_2}\cdots E_{n,n+1}^{\al_n}\qquad\for\;\al=(\al_1,...,\al_n)\in\mbb{N}^n.
\eqno(6.7.37)$$ We write
$$\wht M_k={\msr A}_k\otimes_{\mbb{F}}M,\;(U({\msr G}_-)(1\otimes
M))_k=\wht M _k\bigcap (U({\msr G}_-)(1\otimes
M))\qquad\for\;\;k\in\mbb{N}.\eqno(6.7.38)$$ According to (6.7.18),
$$(U({\msr G}_-)(1\otimes
M))_k=\sum_{\al\in\mbb{N}^n,\;|\al|=k}E^\al(1\otimes
M),\eqno(6.7.39)$$ where $|\al|=\sum_{r=1}^n\al_r$. Define a linear
map $\vf: \wht M\rta U({\msr G}_-)(1\otimes M)$ by
$$\vf(x^\al\otimes v)=E^\al(1\otimes
v)\qquad\for\;\;\al\in\mbb{N}^n,\;v\in M.\eqno(6.7.40)$$ Then
(6.7.19) and (6.7.35) implies that $\vf$ is $gl(n,\mbb{F})$-module
homomorphism.

Note that
$$\omega=\sum_{i,
j=1}^nE_{i,j}E_{j,i}\in U(gl(n,\mbb F))\eqno(6.7.41)$$ is a Casimier
element of $gl(n,\mbb{F})$. Set
$$\td\omega=\frac{1}{2}(\omega|_{\wht M}-\omega|_{\msr A}\otimes
\mbox{Id}_M-\mbox{Id}_{\msr A}\otimes\omega|_M)=\sum_{i,j=1}^n
E_{i,j}|_{\msr A}\otimes E_{j,i}|_M.\eqno(6.7.42)$$ Note that
$\sum_{r=1}^nE_{r,r}$ is the central element of $gl(n,\mbb{F})$.
\psp

{\bf Lemma 6.7.3}. {\it If $(\sum_{r=1}^nE_{r,r})|_M=c\;\mbox{\it
Id}_M$ for $c\in\mbb{F}$, then  we have $\vf|_{\wht
M_1}=(\td\omega+c) |_{\wht M_1}$.}

{\it Proof}. For $s\in\ol{1,n}$ and $v\in M$,
\begin{eqnarray*}
\hspace{1cm}\vf(x_s\otimes v)&=&\nu(E_{s,n+1})(1\otimes v)
=[P_s\oplus (\sum_{r=1}^n(x_sE_{r,r}+x_rE_{s,r})](1\otimes
v)\\&=&cx_s\otimes v+\sum_{r=1}^nx_r\otimes
E_{s,r}(v).\hspace{6cm}(6.7.43)\end{eqnarray*} Moreover,
$$\td\omega(x_s\otimes v)=\sum_{i,j=1}^nx_i\ptl_{x_j}(x_s)\otimes
E_{j,i}(v) =\sum_{i=1}^nx_i\otimes
E_{s,i}(v).\qquad\Box\eqno(6.7.44)$$ \pse

Fix the Cartan subalgebra $H=\sum_{i=1}^n\mbb{F}E_{i,i}$ of
$gl(n,\mbb{F})$ and define $\ves_i\in H^\ast$ by
$$\ves_i(E_{j,j})=\dlt_{i,j}\qquad\for\;\;i,j\in\ol{1,n}.\eqno(6.7.45)$$
For any $\mu=\sum_{i=1}^n\mu_i\ves_i,\mu'=\sum_{i=1}^n\mu'_i\ves_i$,
we define
$$(\mu,\mu')=\sum_{i=1}^n\mu_i\mu'_i,\;\;|\mu|=\sum_{i=1}^n\mu_i.\eqno(6.7.46)$$

Set
$$\Lmd^+=\{\mu=\sum_{i=1}^n\mu_i\ves_i\in
H^\ast\mid\mu_i-\mu_{i+1}\in\mbb{N}\;\for\;i\in\ol{1,n-1}\}.\eqno(6.7.47)$$
Let $V(\mu)$ be a highest-weight irreducible $gl(n,\mbb F)$-module
with $\mu$ as the highest weight.  According to Section 5.2, any
finite-dimensional $gl(n,\mbb{F})$-module is a highest-weight
irreducible module $V(\mu)$ for some $\mu\in \Lmd^+$ when $\mbb F$
is algebraically closed. Moreover,
$$(\sum_{r=1}^nE_{r,r})|_{V(\mu)}=|\mu|\:\mbox{Id}_{V(\mu)}\eqno(6.7.48)$$
by its action on the highest-weight vector. By (3.2.15),
$$\rho=\frac{1}{2}\sum_{r=1}^n(n-2r+1)\ves_r.\eqno(6.7.49)$$
Consider the action of $\omega$ on the highest-weight vector, we
have:
$$\omega|_{V(\mu)}=(\mu,\mu+2\rho)\:\mbox{Id}_{V(\mu)}\eqno(6.7.50)$$

As a $gl(n,\mbb{F})$-module,
$${\msr A}_1\cong V(\ves_1),\eqno(6.7.51)$$
whose weights are
$$\{\ves_1,\ves_2,...,\ves_n\}.\eqno(6.7.52)$$
Fix $\mu\in\Lmd$ and take
$$M=V(\mu).\eqno(6.7.53)$$
Moreover, we define
$${\cal I}(\mu)=\{i_1,i_2,...,i_{s+1}\}\subset\ol{1,n}\;\;\mbox{such
that}\;\;1=i_1<i_2<\cdots <i_{s+1}=n+1,\eqno(6.7.54)$$
$$\mu_{i_r}=\mu_{i_{r+1}-1}\;\;\for\;\;r\in\ol{1,s}\;\;\mbox{and}\;\;
\mu_{i_r}>\mu_{i_{r+1}}\;\;\mbox{if}\;\;r<s.\eqno(6.7.55)$$
According to Theorem 5.4.3, we have:\psp

{\bf Lemma 6.7.4 (Pieri's formula)}. \index{Pieri's formula! for
$sl(n)$}{\it As $gl(n,\mbb{F})$-modules,
$$\wht {V(\mu)}_1={\msr
A}_1\otimes_{\mbb{F}}V(\mu)\cong\bigoplus_{r=1}^sV(\mu+\ves_{i_r}),\eqno(6.7.56)$$
where ${\cal I}(\mu)=\{i_1,i_2,...,i_{s+1}\}$.} \psp

By  (6.7.50), the eigenvalues of $\omega$ on $\wht {V(\mu)}_1$ are:
$$\{(\mu+\ves_{i_1},\mu+\ves_{i_1}+2\rho),....,(\mu+\ves_{i_s},\mu+\ves_{i_s}+2\rho)\}.\eqno(6.7.57)$$
On the other hand, $\omega|_{\msr A_1}\otimes \mbox{Id}_{V(\mu)}$
has the only eigenvalue $(\ves_1,\ves_1+2\rho)=n$ and
$\mbox{Id}_{\msr A_1}\otimes\omega|_{V(\mu)}$ has the only
eigenvalue $(\mu,\mu+2\rho)$. Calculate
$$(\mu+\ves_{i_r},\mu+\ves_{i_r}+2\rho)-(\mu,\mu+2\rho)-n=2(\mu_{i_r}+1-i_r).\eqno(6.7.58)$$
Thus we get: \psp

{\bf Lemma 6.7.5}. {\it The eigenvalues of $\td\omega$ on $\wht
{V(\mu)}_1$ are}
$$\{\mu_{i_r}+1-i_r\mid r\in\ol{1,s}\}.\eqno(6.7.59)$$

For $f\in{\msr A}$, we define the action
$$f(g\otimes v)=fg\otimes v\qquad\for\;\;g\in{\msr A},\;v\in
M.\eqno(7.3.60)$$ Now we have:\psp

{\bf Theorem 6.7.6}. {\it The $sl(n+1,\mbb{F})$-module
$\wht{V(\mu)}$ is irreducible if and only if
$|\mu|+\mu_{i_s}-i_s\not\in-\mbb{N}-1$.}

{\it Proof}. Note
$$\vf[\wht{V(\mu)}_k]= (U({\msr G}_-)(1\otimes
V(\mu)))_k\qquad\for\;\;k\in\mbb{N}.\eqno(6.7.61)$$ By Proposition
6.7.2, $\wht{V(\mu)}$ is irreducible if and only if $\vf$ is
injective. Observe that
$$\wht{V(\mu)}_0=(U({\msr G}_-)(1\otimes
V(\mu))_0=1\otimes V(\mu)\eqno(6.7.62)$$ and
$$\vf(1\otimes v)=1\otimes v\qquad\for\;\;v\in
V(\mu).\eqno(6.7.63)$$ So $\vf|_{\wht{V(\mu)}_0}$ is injective.

Suppose that $\vf|_{\wht{V(\mu)}_k}$ is injective for some
$k\in\mbb{N}$. So
$$\wht{V(\mu)}_k= (U({\msr G}_-)(1\otimes
V(\mu)))_k.\eqno(6.7.64)$$ Thus $\vf|_{\wht{V(\mu)}_{k+1}}$ is
injective if and only if
$$\sum_{i=1}^nE_{i,n+1}[\wht{V(\mu)}_k]=\wht{V(\mu)}_{k+1}.\eqno(6.7.65)$$
For $f\in{\msr A}_k$ and $v\in V(\mu)$, we have
\begin{eqnarray*}\qquad\qquad E_{i,n+1}(f\otimes v)&=&(x_iD\oplus
(\sum_{r=1}^n(x_iE_{r,r}+x_rE_{i,r})))(f\otimes v)\\&=& kx_if\otimes
v+f\sum_{r=1}^n(x_iE_{r,r}+x_rE_{i,r})(1\otimes v)
\\ &=&f[(k+\vf)(x_i\otimes v)]\hspace{6.3cm}(6.7.66)\end{eqnarray*} by (6.7.18). Thus (6.7.64)
holds if and only if
$$(k+\vf)|_{\wht{V(\mu)}_1}\;\;\mbox{is injective}.\eqno(6.7.67)$$

According to (6.7.48) and Lemmas 6.7.3, 6.7.5, the eigenvalues of
$(k+\vf)|_{\wht{V(\mu)}_1} $ are
$$\{k+|\mu|+\mu_{i_r}+1-i_r\mid r\in\ol{1,s}\}.\eqno(6.7.68)$$
Therefore, $\vf$ is injective if and only if
$$0\not\in \{k+|\mu|+\mu_{i_r}+1-i_r\mid
r\in\ol{1,s},\;k\in\mbb{N}\}.\eqno(6.7.69)$$ Observe that
$$\mu_{i_r}-\mu_{i_s}+i_s-i_r\in\mbb{N}\eqno(6.7.70)$$
by (6.7.47) and (6.7.54). Suppose
$|\mu|+\mu_{i_s}-i_s\not\in-\mbb{N}-1$. Then
$$|\mu|+\mu_{i_r}-i_r=|\mu|+\mu_{i_s}-i_s+\mu_{i_r}-\mu_s+i_s-i_r\not\in-\mbb{N}-1\eqno(6.7.71)$$
for $r\in\ol{1,s}$ by (6.7.70). This shows that
$$(6.7.69)\;\;\mbox{holds for any $k\in\mbb N$ if and only
if}\;\;|\mu|+\mu_{i_s}-i_s\not\in-\mbb{N}-1.\qquad\Box\eqno(6.7.72)$$
\psp

Let $\lmd_i$ be the $i$th fundamental weight of $sl(n+1,\mbb F)$
with respect to the order $n+1,1,2,...,n$.
 Since
$\nu(E_{1,1}-E_{n+1,n+1})=(D+x_1\ptl_{1,1})\oplus
(\sum_{r=1}^nE_{r,r}+E_{1,1})$ by (6.7.8), (6.7.12) and (6.7.17),
the irreducible $sl(n+1,\mbb F)$-module $U(\msr G_-)(V(\mu))$ is a
highest-weight module with highest weight
$-(|\mu|+\mu_1)\lmd_1+\sum_{r=1}^{s-1}(\mu_{i_r}-\mu_{i_{r+1}})\lmd_{i_r+1}$,
and so is $\wht{V(\mu)}$ when
$|\mu|+\mu_{i_s}-i_s\not\in-\mbb{N}-1$. In particular, $U(\msr
G_-)(V(\mu))$ is a finite-dimensional irreducible $sl(n+1,\mbb
F)$-module if and only if $|\mu|+\mu_1\in-\mbb N$.

 Let ${\msr B}=\mbb{F}[y_1,...,y_n,z_1,...,z_n]$. Given $c\in \mbb{F}$, we define an action
 of $gl(n,\mbb{F})$ on ${\msr B}$ by
 $$E_{i,j}|_{\msr
 B}=c\dlt_{i,j}+y_i\ptl_{y_j}-z_j\ptl_{z_i}\qquad\for\;\;i,j\in\ol{1,n}.\eqno(6.7.73)$$
For $\ell_1,\ell_2\in\mbb{N}$, we set
$${\msr
B}_{\ell_1,\ell_2}=\sum_{\al,\be\in\mbb{N}^{\:n};\;|\al|=\ell_1,\;|\be|=\ell_2}
\mbb{F}y^\al z^\be.\eqno(6.7.74)$$ Moreover, we define
$${\msr H}_{\ell_1,\ell_2}=\{f\in {\msr
B}_{\ell_1,\ell_2}\mid
(\sum_{i=1}^n\ptl_{y_i}\ptl_{z_i})(f)=0\}.\eqno(6.7.75)$$ According
to Theorem 6.2.4, ${\msr H}_{\ell_1,\ell_2}$ forms a
finite-dimensional irreducible $gl(n,\mbb{F})$-module and
$y_1^{\ell_1}z_n^{\ell_2}$ is a highest-weight vector with the
weight $\ell_1\ves_1-\ell_2\ves_n+c\sum_{i=1}^n\ves_i$. By the above
Theorem, we get:\psp

{\bf Corollary 6.7.7}. {\it Suppose $\ell_1>0$.  The
$sl(n+1,\mbb{F})$-module $\wht{{\msr H}_{\ell_1,\ell_2}}$ is
irreducible if and only if $c\not\in -(\mbb{N}+\ell_1-1)/(n+1)$ when
$\ell_2=0$ and $c\not\in (2\ell_2+n-\ell_1-1-\mbb{N})/(n+1)$ when
$\ell_2\neq 0$.} \psp

Let $\Psi$ be the exterior algebra generated by
$\{\sta_1,\sta_2,...,\sta_n\}$ (cf. (6.2.15)) and take the settings
(6.2.16)-(6.2.18). Given $c\in\mbb{F}$, we define an action of
$gl(n,\mbb{F})$ on $\Psi$ by
$$E_{i,j}|_{\Psi}=c\dlt_{i,j}+\sta_i\ptl_{\sta_j}.\eqno(6.7.76)$$
For any $k\in\ol{1,n-1}$, $\Psi_k$ is a finite-dimensional
irreducible $gl(n,\mbb{F})$-module and $\sta_1\sta_2\cdots\sta_r$ is
a highest-weight vector with the weight
$\sum_{r=1}^k\ves_r+c\sum_{i=1}^n\ves_i$. By Theorem 6.7.6, we
have:\psp

{\bf Corollary 6.7.8}. {\it The $sl(n+1,\mbb{F})$-module
$\wht{\Psi_k}$ is irreducible if and only if
$c\not\in-\mbb{N}/(n+1)$.} \psp

Suppose that $M=\mbb{F}v$ is a one-dimensional
$gl(n,\mbb{F})$-module. Then
$$E_{i,j}|_M=c\dlt_{i,j}\:\mbox{Id}_M\qquad\for\;\;i,j\in\ol{1,n}.\eqno(6.7.77)$$
We identify ${\msr A}$ with $\wht M$ via
$$f\leftrightarrow f\otimes v\qquad\for\;\;f\in{\msr A}.\eqno(6.7.78)$$
Denote
$$\kappa=(n+1)c.\eqno(6.7.79)$$
Recall that $D=\sum_{s=1}^nx_s\ptl_{x_s}$ is the degree operator on
${\msr A}=\mbb{F}[x_1,...,x_n]$. Then we have the following
one-parameter generalization $\pi_\kappa$ of the projective
representation of $sl(n+1,\mbb{F})$ in (6.7.7) and (6.7.8):
$$\pi_\kappa(E_{i,j})=x_i\ptl_{x_j},\qquad\pi_\kappa(E_{i,n+1})=x_i(D+\kappa),\eqno(6.7.80)$$
$$\pi_\kappa(E_{n+1,i})=-\ptl_{x_i},\;\;\pi_\kappa(E_{i,i}-E_{n+1,n+1})=D+\kappa+x_i\ptl_{x_i}\eqno(6.7.81)$$
for $i,j\in\ol{1,n}$ with $i\neq j$.
 According to Theorem 6.7.6, we have: \psp

{\bf Corollary 6.7.9}. {\it The representation $\pi_\kappa$ of
$sl(n+1,\mbb{F})$ on $\msr A$ is irreducible if and only if
$\kappa\not\in -\mbb{N}$.}\psp

For $\ell\in\mbb N$, we set
$$\msr A_{(\ell)}=\sum_{i=0}^\ell\msr A_i.\eqno(6.7.82)$$
Reordering $\ol{1,n+1}$ as $\{n+1,1,2,...,n\}$, we get directly by
(6.7.80) and (6.7.81):\psp

{\bf Corollary 6.7.10}. {\it If $\kappa=-\ell$ for some $\ell\in\mbb
N$, then $\msr A_{(\ell)}$ is a finite-dimensional irreducible
$sl(n+1,\mbb{F})$-module with highest weight $\ell\lmd_1$ and $\msr
A/\msr A_{(\ell)}$ is an infinite-dimensional highest-weight
irreducible $sl(n+1,\mbb{F})$-module with highest weight
$-(\ell+2)\lmd_1+(\ell+1)\lmd_2.$}

\section{Projective Oscillator Representations}

In this section, we study oscillator generalizations of the
representation $\pi_\kappa$ of $sl(n+1,\mbb{F})$ in (6.7.80) and
(6.7.81).

  Note the symmetry
$$[\ptl_{x_r},x_r]=1=[-x_r,\ptl_{x_r}].\eqno(6.8.1)$$
Fix $n_1\in\ol{0,n-1}$. We define $\pi_{\kappa,0}=\pi_\kappa$.
Changing operators $\ptl_{x_r}\mapsto -x_r$ and $x_r\mapsto
\ptl_{x_r}$ for $r\in \ol{1,n_1}$ in (6.7.80) and (6.7.81) when
$n_1>0$, we get another differential-operator representation
$\pi_{\kappa,n_1}$ of $sl(n+1,\mbb{F})$. We call $\pi_{\kappa,n_1}$
{\it projective oscillator representations} in terms of physics
terminology.\index{projective oscillator representation!for $sl(n)$}
 For
$\vec a=(a_1,a_2,...,a_n)^t\in\mbb{F}^n$, we denote $\vec a\cdot\vec
x=\sum_{i=1}^na_ix_i$. Recall ${\msr A}=\mbb{F}[x_1,x_2,...,x_n]$.
Set
$${\msr A}_{\vec a}=\{fe^{\vec a\cdot\vec
x}\mid f\in{\msr A}\}.\eqno(6.8.2)$$ Denote by
$\pi_{\kappa,n_1}^{\vec a}$ the representation $\pi_{\kappa,n_1}$ of
$sl(n+1,\mbb{F})$ on ${\msr A}_{\vec a}$. Corollary 6.7.9 says that
the representation $\pi_{\kappa,0}^{\vec 0}$ of $sl(n+1,\mbb{F})$ is
irreducible if and only if $\kappa\not\in -\mbb{N}$. The following
is the main theorem in this section.\psp

{\bf Theorem 6.8.1}. {\it The representation $\pi_{\kappa,n_1}^{\vec
0}$ is irreducible for any $\kappa\in\mbb F\setminus \mbb Z$, and
the underlying module ${\msr A}$ is an infinite-dimensional weight
$sl(n+1,\mbb{F})$-module with finite-dimensional weight subspaces.
If $a_i\neq 0$ for some $i\in \ol{n_1+1,n}$, then the representation
$\pi_{\kappa,n_1}^{\vec a}$ of $sl(n+1,\mbb{F})$ is always
irreducible for any $\kappa\in\mbb F$. When  $n_1>1,\;\vec a\neq 0$
and $a_i=0$  for any $i\in \ol{n_1+1,n}$, the representation
$\pi_{\kappa,n_1}^{\vec a}$ of $sl(n+1,\mbb{F})$ is irreducible for
$\kappa\in\mbb F\setminus \mbb Z$.}\psp

 We will prove Theorem 6.8.1 case by case.\psp

{\it Case 1. The representation $\pi_{\kappa,n_1}^{\vec 0}$ with
$0<n_1<n$.} \psp

As in Section 6.3, we set
$$\td D=\sum_{r=n_1+1}^nx_r\ptl_{x_r}-\sum_{i=1}^{n_1}x_i\ptl_{x_i}.\eqno(6.8.3)$$
Then the representation $\pi_{\kappa,n_1}^{\vec 0}$ of
$sl(n+1,\mbb{F})$ is the representation $\pi_{\kappa,n_1}$ on ${\msr
A}$ with
$$\pi_{\kappa,n_1}(E_{i,j})=\left\{\begin{array}{ll}-x_j\ptl_{x_i}-\delta_{i,j}&\mbox{if}\;
i,j\in\ol{1,n_1};\\ \ptl_{x_i}\ptl_{x_j}&\mbox{if}\;i\in\ol{1,n_1},\;j\in\ol{n_1+1,n};\\
-x_ix_j &\mbox{if}\;i\in\ol{n_1+1,n},\;j\in\ol{1,n_1};\\
x_i\partial_{x_j}&\mbox{if}\;i,j\in\ol{n_1+1,n},
\end{array}\right.\eqno(6.8.4)$$
$$\pi_{\kappa,n_1}(E_{i,n+1})=\left\{\begin{array}{ll}(\td D+\kappa-n_1-1)\ptl_{x_i}&\mbox{if}\;\;i\leq
n_1,\\ x_i(\td D+\kappa-n_1)&\mbox{if}\;\;i>
n_1,\end{array}\right.\eqno(6.8.5)$$
$$\pi_{\kappa,n_1}(E_{n+1,i})=\left\{\begin{array}{ll}x_i&\mbox{if}\;\;i\leq
n_1,\\-\ptl_{x_i}&\mbox{if}\;\;i>
n_1,\end{array}\right.\eqno(6.8.6)$$
$$\pi_{\kappa,n_1}(E_{n,n}-E_{n+1,n+1})=\td D-n_1+\kappa+x_n\ptl_{x_n}\eqno(6.8.7)$$
Recall
$${\msr A}_{\la
k\ra}=\mbox{Span}\:\{x^\al=\prod_{i=1}x_i^{\al_i}\mid\al=(\al_1,...,\al_n)\in\mbb{N}^n;\sum_{i=1}^{n_1}\al_i-\sum_{r=n_1+1}^n\al_r=k\}.
\eqno(6.8.8)$$ Then ${\msr A}=\bigoplus_{k\in\mbb{Z}}{\msr A}_{\la
k\ra}$ and
$${\msr A}_{\la k\ra}=\{f\in{\msr A}\mid
\td D(f)=kf\}.\eqno(6.8.9)$$ Note that
$${\msr G}_0=\sum_{1\leq i<j\leq
n}(\mbb{F}E_{i,j}+\mbb{F}E_{j,i})+\sum_{r=1}^{n-1}\mbb{F}(E_{r,r}-E_{r+1,r+1})\eqno(6.8.10)$$
is a Lie subalgebra of $sl(n+1,\mbb{F})$ isomorphic to
$sl(n,\mbb{F})$.
 \psp

{\bf Lemma 6.8.2}. {\it The representation $\pi_{\kappa,n_1}^{\vec
0}$ of $sl(n+1,\mbb{F})$ is irreducible for any $\kappa\in\mbb
F\setminus \mbb{Z}$.}

{\it Proof}. Let $k$ be any integer.  For any $0\neq f\in {\msr
A}_{\la k\ra}$, we have
$$0\neq E_{n+1,1}(f)=x_1f\in {\msr A}_{\la k-1\ra}\eqno(6.8.11)$$
by (6.8.6), and
$$0\neq E_{n,n+1}(f)=(k+\kappa-n_1)x_nf\in {\msr A}_{\la
k+1\ra}\eqno(6.8.12)$$ by (6.8.5). Let ${M}$ be a nonzero
$sl(n+1,\mbb{F})$-submodule of ${\msr A}$. If $k_1,k_2\in\mbb{Z}$
with $k_1\neq k_2$, then the highest weights of ${\msr A}_{\la
k_1\ra}$ and ${\msr A}_{\la k_2\ra}$ are different as ${\msr
G}_0$-modules by Theorem 6.3.3. So ${\msr A}_{\la k_0\ra}\subset
{M}$ for some $k_0\in\mbb{Z}$. Moreover, (6.8.11) and (6.8.12) imply
${\msr A}_{\la k\ra}\subset {M}$ for any $k\in\mbb{Z}$. Hence ${\msr
M}={\msr A}$. $\qquad\Box$ \psp

Expressions (6.8.4)-(6.8.7) imply the above representation is not of
highest-weight type. Moreover, ${\msr A}$ is a weight
$sl(n+1,\mbb{F})$-module with finite-dimensional weight subspaces.

\psp

{\it Case 2. The representation $\pi_{\kappa,0}^{\vec a}$ with $\vec
0\neq \vec a\in\mbb{F}^n$.} \psp

In this case,
$$E_{n+1,i}(fe^{\vec a\cdot\vec
x})=-(\ptl_{x_i}+a_i)(f)e^{\vec a\cdot\vec
x}\qquad\for\;\;i\in\ol{1,n},\;f\in{\msr A}.\eqno(6.8.13)$$ Thus
$$(E_{n+1,i}+a_i)(fe^{\vec a\cdot\vec
x})=-\ptl_{x_i}(f)e^{\vec a\cdot\vec
x}\qquad\for\;\;i\in\ol{1,n},\;f\in{\msr A}.\eqno(6.8.14)$$ The
second result in this section.\psp

{\bf Lemma 6.8.3}. {\it The representation $\pi_{\kappa,0}^{\vec a}$
with $\vec 0\neq \vec a\in\mbb{F}^n$ is an irreducible
representation of $sl(n+1,\mbb{F})$ for any $\kappa\in\mbb{F}$.}

{\it Proof}. Let $\msr{A}_k$ be the subspace of homogeneous
polynomials with degree $k$. Set
$$\msr{A}_{\vec a,k}=\msr{A}_ke^{\vec a\cdot\vec
x}\qquad\for\;k\in\mbb{N}.\eqno(6.8.15)$$ Without loss of
generality, we assume $a_1\neq 0$. Let  ${M}$ be a nonzero
$sl(n+1,\mbb{F})$-submodule of ${\msr A}_{\vec a}$. Take any $0\neq
fe^{\vec a\cdot\vec x}\in  M$ with $f\in \msr{A}$. By (6.8.14),
$$\ptl_{x_i}(f)e^{\vec a\cdot\vec
x}\in M\qquad\for\;\;i\in\ol{1,n}.\eqno(6.8.16)$$ By induction, we
have $e^{\vec a\cdot\vec x}\in M$; that is, $\msr{A}_{\vec
a,0}\subset M$.

Suppose $\msr{A}_{\vec a,\ell}\subset M$ for some $\ell\in\mbb{N}$.
For any $ge^{\vec a\cdot\vec x}\in \msr{A}_{\vec a,\ell}$,
$$E_{i,1}(ge^{\vec a\cdot\vec x})=x_i(\ptl_{x_1}+a_1)(g)e^{\vec a\cdot\vec x}
=a_1x_ige^{\vec a\cdot\vec x}+x_i\ptl_{x_1}(g)e^{\vec a\cdot\vec
x}\in M\qquad\for\;\;i\in\ol{2,n}\eqno(6.8.17)$$ by (6.7.80). Since
$x_i\ptl_{x_1}(g)e^{\vec a\cdot\vec x}\in \msr{A}_{\vec
a,\ell}\subset M$, we have
$$x_ige^{\vec a\cdot\vec x}\in M\qquad\for\;\;i\in\ol{2,n}.\eqno(6.8.18)$$

On the other hand,
$$(E_{1,1}-E_{2,2})(ge^{\vec a\cdot\vec x})=a_1x_1ge^{\vec a\cdot\vec
x}+(x_1\ptl_{x_1}-x_2\ptl_{x_2}-a_2x_2)(g)e^{\vec a\cdot\vec x}\in
M\eqno(6.8.19)$$ by (6.7.81).  Our assumption says that
$(x_1\ptl_{x_1}-x_2\ptl_{x_2})(g)e^{\vec a\cdot\vec x}\in
\msr{A}_{\vec a,\ell}\subset M$. According to (6.8.18),
$-a_2x_2(g)e^{\vec a\cdot\vec x}\in M$. Therefore,
$$x_1ge^{\vec a\cdot\vec
x}\in M.\eqno(6.8.20)$$ Expressions (6.8.19) and (6.8.20) imply
$\msr{A}_{\vec a,\ell+1}\subset M$. By induction, $\msr{A}_{\vec
a,\ell}\subset M$ for any $\ell\in\mbb{N}$. So $\msr{A}_{\vec a}=
M$. Hence $\msr{A}_{\vec a}$ is an irreducible
$sl(n+1,\mbb{F})$-module. $\qquad\Box$\psp

{\it Case 3. The representation $\pi_{\kappa,n_1}^{\vec a}$ with
$n_1>0$, and $a_i\neq 0$ for some $i\in \ol{n_1+1,n}$.}\psp

The following is the third result in this section.\psp

{\bf Lemma 6.8.4}. {\it Under the above assumption, the
representation $\pi_{\kappa,n_1}^{\vec a}$ is an irreducible
representation of $sl(n+1,\mbb{F})$ for any $\kappa\in\mbb F$.}

{\it Proof}. Let ${M}$ be a nonzero $sl(n+1,\mbb{F})$-submodule of
${\msr A}_{\vec a}$. By (6.8.6) and (6.8.13)-(6.8.16), there exists
$0\neq fe^{\vec a\cdot\vec x}\in M$ with
$f\in\mbb{F}[x_1,...,x_{n_1}]$. By symmetry, we can assume $a_n\neq
0$. According to (6.8.4),
$$E_{i,n}(fe^{\vec a\cdot\vec
x})=(\ptl_{x_i}+a_i)(\ptl_{x_n}+a_n)(f)e^{\vec a\cdot\vec
x}=a_ia_nfe^{\vec a\cdot\vec x}+a_n\ptl_{x_i}(f)e^{\vec a\cdot\vec
x}\;\;\for\;i\in\ol{1,n_1}.\eqno(6.8.21)$$ Thus
$$(a_n^{-1}E_{i,n}-a_i)(fe^{\vec a\cdot\vec
x})=\ptl_{x_i}(f)e^{\vec a\cdot\vec x}\in
M\;\;\for\;i\in\ol{1,n_1}.\eqno(6.8.22)$$ By induction on the degree
of $f$, we get $e^{\vec a\cdot\vec x}\in M$; that is, $\msr{A}_{\vec
a,0}\subset M$.

 The arguments in (6.8.17)-(6.8.20) yield
$$\mbb{F}[x_{n_1+1},...,x_n]e^{\vec a\cdot\vec
x}\subset M.\eqno(6.8.23)$$ According to (6.8.6),
$$E_{n+1,1}^{\ell_1}\cdots E_{n+1,n_1}^{\ell_{n_1}}(\mbb{F}[x_{n_1+1},...,x_n]e^{\vec a\cdot\vec
x})=x_1^{\ell_1}\cdots
x_{n_1}^{\ell_{n_1}}(\mbb{F}[x_{n_1+1},...,x_n]e^{\vec a\cdot\vec
x})\subset M\eqno(6.8.24)$$ for $\ell_i\in\mbb{N}$ with
$i\in\ol{1,n_1}.$ Thus ${\msr A}_{\vec a}= M$. So  $\msr{A}_{\vec
a}$ is an irreducible $sl(n+1,\mbb{F})$-module. $\qquad\Box$\psp

{\it Case 4. $n_1>1,\;\vec a\neq\vec 0$ and $a_i= 0$ for any
$i\in\ol{n_1+1,n}$.}\psp

The following is the fourth result in this section.\psp

{\bf Lemma 6.8.5}. {\it Under the above assumption, the
representation $\pi_{\kappa,n_1}^{\vec a}$ is an irreducible
representation of $sl(n+1,\mbb{F})$ for any $\kappa\in\mbb
F\setminus\mbb Z$.}

{\it Proof}.  Applying the transformation
$$\vec x\mapsto T\vec x ,\;\;A\mapsto TAT^{-1},\;\;E_{n,n}-E_{n+1,n+1}\mapsto
E_{n,n}-E_{n+1,n+1}.\eqno(6.8.25)$$
$$(E_{1,n+1},...,E_{n,n+1})\mapsto
(E_{1,n+1},...,E_{n,n+1})T^{-1},\eqno(6.8.26)$$
$$(E_{n+1,1},...,E_{n+1,n})\mapsto
(E_{n+1,n},...,E_{n+1,n})T^{-1}\eqno(6.8.27)$$ with $A\in
sl(n,\mbb{F})$ for some $n\times n$ orthogonal  matrix $T$, we can
assume $a_1\neq 0$ and $a_i=0$ for $i\in\ol{2,n}$.

 Let ${M}$ be a nonzero $sl(n+1,\mbb{F})$-submodule
of ${\msr A}_{\vec a}$. Take any $0\neq fe^{\vec a\cdot\vec x}\in
M$.
 Note that
$$E_{1,2}(fe^{\vec a\cdot\vec x})=-x_2(\ptl_{x_1}+a_1)(f)e^{\vec a\cdot\vec
x}\in M\eqno(6.8.28)$$ by (6.8.4) and
$$E_{n+1,2}(fe^{\vec a\cdot\vec x})=x_2fe^{\vec a\cdot\vec x}\in M\eqno(6.8.29)$$ by
(6.8.6). Thus
$$x_2\ptl_{x_1}(f)e^{\vec a\cdot\vec
x}\in M.\eqno(6.8.30)$$ Repeatedly applying (6.8.29) if necessary,
we can assume $f\in\mbb{F}[x_2,....,x_n]$. We apply the arguments in
the proof of Lemma 6.8.2 to the Lie subalgebra
$$\msr{L}=\sum_{2\leq i<j\leq
n+1}(\mbb{F}E_{i,j}+\mbb{F}E_{j,i})+\sum_{r=2}^n\mbb{F}(E_{r,r}-E_{r+1,r+1})\eqno(6.8.31)$$
and obtain
$$\mbb{F}[x_2,....,x_n]e^{\vec a\cdot\vec x}\subset  M\eqno(6.8.32)$$
when $\kappa\not\in\mbb Z$. According to (6.8.6),
$$E_{n+1,1}^\ell(\mbb{F}[x_2,....,x_n]e^{\vec a\cdot\vec x})=x_1^\ell(\mbb{F}[x_2,....,x_n]e^{\vec a\cdot\vec x})
\subset  M\eqno(6.8.33)$$ for any $\ell\in\mbb{N}$. Therefore $
M=\msr{A}_{\vec a}$. So $\msr{A}_{\vec a}$ is an irreducible
$sl(n+1,\mbb{F})$-module. $\qquad\Box$\psp

With the representation $\pi_{\kappa,n_1}^{\vec a}$ with $\vec
a\neq\vec 0$, ${\msr A}_{\vec a}$ is not a weight
$sl(n+1,\mbb{F})$-module. Now Theorem 6.8.1 follows from Lemmas
6.8.2-6.8.5.

\chapter{Representations of Even Orthogonal Lie Algebras}

In this chapter, we focus on  the natural explicit representations
of even orthogonal Lie algebras. In Section 7.1, we study the
canonical  bosonic  and fermionic oscillator representations over
their minimal natural modules.  The spin representations are also
presented.  In Section 7.2, we determine the structure of the
noncanonical oscillator representations obtained from the above
bosonic representations by partially swapping differential operators
and multiplication operators, which are generalizations of the
classical theorem on harmonic polynomials. The results were due to
Luo and the author [LX2].  In Section 7.3, we construct a new
functor from the category $o(2n,\mbb C)$-{\bf mod} to the category
$o(2n+2,\mbb C)$-{\bf mod}, which is an extension of the conformal
representation of $o(2n+2,\mbb C)$. Moreover, we find the condition
for the functor to map a finite-dimensional irreducible $o(2n,\mbb
C)$-module to an infinite-dimensional irreducible $o(2n+2,\mbb
C)$-module. This work was due to Zhao and the author [XZ]. The work
in Section 7.3 gives rise to a one-parameter ($c$) family  of
inhomogeneous first-order differential operator representations of
$o(2n+2,\mbb C)$. Letting these operators act on the space of
exponential-polynomial functions that depend on  a parametric vector
$\vec a\in \mbb C^{2n}$, we prove in Section 7.4 that the space
forms an irreducible $o(2n+2,\mbb C)$-module for any $c\in\mbb C$ if
$\vec a$ is not on a certain hypersurface. By partially swapping
differential operators and multiplication operators, we obtain more
general differential operator representations of $o(2n+2,\mbb C)$ on
the polynomial algebra $\msr B$ in $2n$ variables. We prove that
$\msr B$ forms an infinite-dimensional irreducible weight
$o(2n+2,\mbb C)$-module with finite-dimensional weight subspaces if
$c\not\in\mbb Z/2$. These results are taken from [X26].

\section{Canonical Oscillator Representations}

In this section , we present the canonical  bosonic  and fermionic
oscillator representations of even orthogonal Lie algebras over
their minimal natural modules.

Let $n>1$ be an integer.  Recall that the split even orthogonal Lie
algebra
\begin{eqnarray*}\qquad o(2n,\mbb{F})&=&\sum_{1\leq
p<q\leq n}[\mbb{F}(E_{p,n+q}-E_{q,n+p})+
\mbb{F}(E_{n+p,q}-E_{n+q,p})]\\ &
&+\sum_{i,j=1}^n\mbb{F}(E_{i,j}-E_{n+j,n+i}).\hspace{7cm}(7.1.1)\end{eqnarray*}
We take the subspace
$$H=\sum_{i=1}^n\mbb{F}(E_{i,i}-E_{n+i,n+i})\eqno(7.1.2)$$
as a Cartan subalgebra and define $\{\ves_i\mid
i\in\ol{1,n}\}\subset H^\ast$ by
$$\ves_i(E_{j,j}-E_{n+j,n+j})=\dlt_{i,j}.\eqno(7.1.3)$$
The inner product $(\cdot,\cdot)$ on the $\mbb{Q}$-subspace
$$L_\mbb{Q}=\sum_{i=1}^n\mbb{Q}\ves_i\eqno(7.1.4)$$
is given by
$$(\ves_i,\ves_j)=\dlt_{i,j}\qquad\for\;\;i,j\in\ol{1,n}.\eqno(7.1.5)$$
Then the root system of $o(2n,\mbb{F})$ is
$$\Phi_{D_n}=\{\pm \ves_i\pm\ves_j\mid1\leq i<j\leq
n\}.\eqno(7.1.6)$$ We take the set of positive roots
$$\Phi_{D_n}^+=\{\ves_i\pm\ves_j\mid1\leq i<j\leq
n\}.\eqno(7.1.7)$$ In particular,
$$\Pi_{D_n}=\{\ves_1-\ves_2,...,\ves_{n-1}-\ves_n,\ves_{n-1}+\ves_n\}\;\;
\mbox{is the set of positive simple roots}.\eqno(7.1.8)$$

Recall the set of dominate integral weights
$$\Lmd^+=\{\mu\in L_\mbb{Q}\mid
(\ves_{n-1}+\ves_n,\mu),(\ves_i-\ves_{i+1},\mu)\in\mbb{N}\;\for\;i\in\ol{1,n-1}\}.\eqno(7.1.9)$$
According to (7.1.5),
$$\Lmd^+=\{\mu=\sum_{i=1}^n\mu_i\ves_i\mid
\mu_i\in\mbb{Z}/2;\mu_i-\mu_{i+1},\mu_{n-1}+\mu_n\in\mbb{N}\}.\eqno(7.1.10)$$
Note that if $\mu\in\Lmd^+$, then $\mu_{n-1}\geq |\mu_n|$.  For
$\lmd\in\Lmd^+$, we denote by $V(\lmd)$ the finite-dimensional
irreducible $o(2n,\mbb{F})$-module with highest weight $\lmd$.

Let ${\msr B}=\mbb{F}[x_1,...,x_n,y_1,...,y_n]$. Recall
(6.2.22)-(6.2.42). The canonical bosonic oscillator
representation\index{canonical bosonic oscillator representation!of
$o(2n,\mbb F)$} of $o(2n,\mbb{F})$ is given by
$$(E_{i,j}-E_{n+j,n+i})|_{\msr B}=x_i\ptl_{x_j}-y_j\ptl_{x_i},\;(E_{n+i,j}-E_{n+j,i})|_{\msr B}
=y_i\ptl_{x_j}-y_j\ptl_{x_i},\eqno(7.1.11)$$
$$(E_{i,n+j}-E_{j,n+i})|_{\msr B}
=x_i\ptl_{y_j}-x_j\ptl_{y_i}\eqno(7.1.12)$$ for $i,j\in\ol{1,n}$.
The positive root vectors of $o(2n,\mbb{F})$ are
$$\{E_{i,j}-E_{n+j,n+i},E_{i,n+j}-E_{j,n+i}\mid 1\leq i<j\leq
n\}.\eqno(7.1.13)$$ It can be verified that $\eta$ in (6.2.23) is an
$o(2n,\mbb{F})$-invariant and $\Dlt$ in (6.2.33)  is an
$o(2n,\mbb{F})$-invariant differential operator, or equivalently,
$$\xi\eta=\eta\xi,\;\xi\Dlt=\Dlt\xi\;\;\;\mbox{on}\;\;{\msr B}\;\;\for\;\;\xi\in
o(2n,\mbb{F}).\eqno(7.1.14)$$

Recall the notion ${\msr B}_{\ell_1,\ell_2}$ in (6.2.30) and set
$${\msr B}_k=\sum_{\ell=0}^k{\msr B}_{\ell,k-\ell},\;\;{\msr
H}_k=\{f\in{\msr
B}_k\mid\Dlt(f)=0\}\qquad\for\;\;k\in\mbb{N}.\eqno(7.1.15)$$ \pse

{\bf Theorem 7.1.1}. {\it For $k\in\mbb{N}$, the subspace ${\msr
H}_k$ forms a finite-dimensional irreducible
$o(2n,\mbb{F})$-submodule and $x_1^k$ is a highest-weight vector
with the weight $k\ves_1$.  Moreover, it has a basis
$$\left\{\sum_{i=0}^\infty (-1)^i\frac{(x_1y_1)^i(\sum_{i=2}^n\ptl_{x_i}\ptl_{y_i})^i(x^\al y^\be)}
{\prod_{r=1}^i(\al_1+i)(\be_1+i)}\mid\al,\be\in\mbb{N};\al_1\be_1=0,|\al|+|\be|=k\right\}.
\eqno(7.1.16)$$ Furthermore,
$${\msr B}=\bigoplus_{k_1,k_2}\eta^{k_1}{\msr
H}_{k_2}\eqno(7.1.17)$$ is a direct sum of irreducible
$o(2n,\mbb{F})$-submodules.}

{\it Proof}. According to (6.2.35),
$${\msr H}_k=\sum_{\ell=0}^k{\msr H}_{\ell,k-\ell}.\eqno(7.1.18)$$
Thus (7.1.16) follows from (6.2.42) and (7.1.17) follows from
(6.2.37) and (6.2.41). Moreover, the second equation in (7.1.14)
implies that ${\msr H}_k$ forms an $o(2n,\mbb{F})$-submodule. Note
that $o(2n,\mbb{F})$ contains the Lie subalgebra
$$\sum_{i,j=1}^n\mbb{F}(E_{i,j}-E_{n+j,n+i})\cong
gl(n,\mbb{F})\eqno(7.1.19)$$ and the representation of
$o(2n,\mbb{F})$ in (7.1.11) and (7.1.12) is essentially an extension
of that for $sl(n,\mbb{F})$ given in (6.2.22). Thus the only
possible highest-weight vectors for $o(2n,\mbb{F})$ in ${\msr H}_k$
are:
$$\{x_1^\ell y_n^{k-\ell}\mid \ell\in\ol{0,k}\}\eqno(7.1.20)$$
by Lemma 6.2.3, (6.2.39) and (7.1.13). Note
$$(E_{1,2n}-E_{n,n+1})(x_1^\ell
y_n^{k-\ell})=(x_1\ptl_{y_n}-x_n\ptl_{y_1})(x_1^\ell
y_n^{k-\ell})=(k-\ell)x_1^\ell y_n^{k-\ell-1}.\eqno(7.1.21)$$ Thus
${\msr H}_k$ has a unique (up to a scalar multiple) singular vector
$x_1^k$ of weight $k\ves_1$. By Weyl's Theorem 2.3.6 of complete
reducibility if $\mbb F=\mbb C$ or Lemma 6.3.2 with $n_1=0$, it is
irreducible. The first equation in (7.1.14) shows that all
$\eta^{k_1}{\msr H}_k$ are irreducible $o(2n,\mbb{F})$-submodules.
$\qquad\Box$ \psp

Consider the exterior algebra $\check{\msr A}$ generated by
$\{\sta_1,...,\sta_n,\vt_1,...,\vt_n\}$ (cf. (6.2.15)) and take the
settings in (6.2.44)-(6.2.51). We have the {\it fermionic
osicillator representation}\index{fermionic osicillator
representation!of of $o(2n,\mbb{F})$} $o(2n,\mbb{F})$ on
$\check{\msr A}$ by
$$(E_{i,j}-E_{n+j,n+i})|_{\check{\msr
A}}=\sta_i\ptl_{\sta_j}-\vt_j\ptl_{\vt_i},\;\;(E_{n+i,j}-E_{n+j,i})|_{\check{\msr
A}} =\vt_i\ptl_{\sta_j}-\vt_j\ptl_{\sta_i},\eqno(7.1.22)$$
$$(E_{i,n+j}-E_{j,n+i})|_{\check{\msr
A}} =\sta_i\ptl_{\vt_j}-\sta_j\ptl_{\vt_i}\eqno(7.1.23)$$ for
$i,j\in\ol{1,n}$. By (7.1.19), (7.1.22) and (7.1.23), the above
representation of $o(2n,\mbb{F})$  is essentially an extension of
that for $sl(n,\mbb{F})$ given in (6.2.43). Set
$$\check{\msr A}_k=\sum_{\ell=0}^{\min\{k,n\}}\check{\msr
A}_{\ell,k-\ell}\qquad\for\;\;k\in\ol{0,2n}\eqno(7.1.24)$$ (cf.
(6.2.45)). According to (6.2.53), the $o(2n,\mbb{F})$-singular
vectors of $\check{\msr A}_k$ are in the set
$$\{\check\eta^\ell\vec\sta_r\vec\vt_s\mid 0\leq r<s\leq
n+1;\ell\in\ol{0,s-r-1};r+2\ell+n-s+1=k\}\eqno(7.1.25)$$ (cf.
(6.2.50) and (6.2.51)).

For $0<i<j\leq n$,
$$[E_{i,n+j}-E_{j,n+i},\check\eta]=\sta_j\sta_i-\sta_i\sta_j=2\sta_j\sta_i\eqno(7.1.26)$$
by (6.2.49) and (7.1.23).
 If $0< r+1<s\leq n$ and $\ell\in\ol{0,s-r-1}$, then
$$(E_{r+1,n+s}-E_{s,n+r+1})(\check\eta^\ell\vec\sta_r\vec\vt_s)=
2\ell\sta_s\sta_{r+1}\check\eta^{\ell-1}\vec\sta_r\vec\vt_s
 +(-1)^{n-s} \check\eta^\ell\vec\sta_{r+1}\vec\vt_{s+1}\neq 0\eqno(7.1.27)$$
 When $1\leq r+1=s\leq n-1$,
$$(E_{n,n+s}-E_{s,2n})(\vec\sta_{s-1}\vec\vt_s)=\vec\sta_{s-1}\vec\vt_{s+1}\sta_n
-\vec\sta_s\vt_{n-1}\cdots \vt_s\neq 0.\eqno(7.1.28)$$

Note that $\eta^\ell\vec\sta _r=0$ if $\ell >n-r$ by (6.2.49). For
$r\in\ol{1,n}$ and $\ell\in\ol{0,n-r}$, (7.1.23) yields
$$(E_{n-1,2n}-E_{n,2n-1})(\check\eta^\ell\vec\sta_r)=2\ell\sta_n\sta_{n-1}
\eta^{\ell-1}\vec\sta_r=0\Leftrightarrow\ell=0\;\;\mbox{or}\;\;\ell\geq
n-r.\eqno(7.1.29)$$ Furthermore,
$$(E_{n-1,2n}-E_{n,2n-1})(\vec\sta_{n-1}\vt_n)=\vec\sta_{n-1}\sta_{n-1}=0.\eqno(7.1.30)$$
Since $\{E_{i,i+1}-E_{n+i+1,n+i},E_{n-1,2n}-E_{n,2n-1}\mid
i\in\ol{1,n-1}\}$ is the set of simple positive root vectors,  the
singular vectors of $o(2n,\mbb{F})$ in $\check{\msr A}$ are
$$\{\vec\sta_n,\vec\sta_{n-1}\vt_n,\vec\sta_r,\check\eta^{n-r}\vec\sta_r\mid
r\in\ol{0,n-1}\}.\eqno(7.1.31)$$ By Weyl's Theorem 2.3.6 of complete
reducibility if $\mbb F=\mbb C$ or an analogue of Lemma 6.3.2, we
have: \psp

 {\bf Theorem 7.1.2}. {\it For $r\in\ol{1,n-1}$, $\check{\msr A}_r$
 and $\check{\msr A}_{2n-r}$ are finite-dimensional irreducible
 $o(2n,\mbb{F})$-submodules with highest weight
 $\sum_{i=1}^r\ves_i$. Moreover, $\check{\msr A}_n$ is the direct
 sum of a finite-dimensional irreducible
 $o(2n,\mbb{F})$-submodule with  highest weight
 $\sum_{i=1}^n\ves_i$ and a finite-dimensional irreducible
 $o(2n,\mbb{F})$-submodule with  highest weight
 $\sum_{i=1}^{n-1}\ves_i-\ves_n$. Furthermore, $\check{\msr A}_{2n}$ and $\check{\msr A}_0$ are isomorphic  to the
one-dimensional trivial $o(2n,\mbb{F})$-module.}
 \psp

Let $\Psi$ be the exterior algebra generated by
$\{\sta_1,\sta_2,...,\sta_n\}$ (cf. (6.2.15)). Take the settings as
those in (6.2.16)-(6.2.18). Observe that
$$[\sta_i\sta_j,\ptl_{\sta_j}\ptl_{\sta_i}]=\sta_i\ptl_{\sta_i}+\sta_j\ptl_{\sta_j}-1\eqno(7.1.32)$$
and
$$[\sta_i\sta_j,\ptl_{\sta_j}\ptl_{\sta_l}]=\sta_i\ptl_{\sta_l}\eqno(7.1.33)$$
for distinct $i,j,l\in\ol{1,n}$. Thus we have the following
representation of $o(2n,\mbb{F})$:
$$(E_{i,j}-E_{n+j,n+i})|_\Psi=\sta_i\ptl_{\sta_j}-\frac{\dlt_{i,j}}{2}\qquad\for\;\;i,j\in\ol{1,n},\eqno(7.1.34)$$
$$(E_{r,n+s}-E_{s,n+r})|_\Psi=\sta_r\sta_s,\qquad
(E_{n+r,s}-E_{n+s,r})|_\Psi=\ptl_{\sta_r}\ptl_{\sta_s}\eqno(7.1.35)$$
for $1\leq r<s\leq n$. The above representation is called the {\it
spin representation} of $o(2n,\mbb{F})$.\index{spin
representation!of $o(2n,\mbb{F})$}

 Recall the notion $\Psi_r$ defined in
(6.2.16) and we set
$$\Psi_{(0)}=\sum_{i=0}^{\llbracket n/2
\rrbracket}\Psi_{2i},\qquad\Psi_{(1)}=\sum_{s=0}^{\llbracket (n-1)/2
\rrbracket}\Psi_{2s+1}.\eqno(7.1.36)$$ Then
$$\Psi=\Psi_{(0)}\oplus\Psi_{(1)}.\eqno(7.1.37)$$
Moreover, we have:\psp

{\bf Theorem 7.1.3}. {\it The subspaces $\Psi_{(0)}$ and
$\Psi_{(1)}$ form  $2^{n-1}$-dimensional irreducible
$o(2n,\mbb{F})$-submodules. If $n$ is even, then
$\sta_1\cdots\sta_n$ is a highest-weight vector of $\Psi_{(0)}$ with
the weight $(\sum_{i=1}^n\ves_i)/2$ and $\sta_1\cdots\sta_{n-1}$ is
a highest-weight vector of $\Psi_{(1)}$
 with the weight $(\sum_{i=1}^{n-1}\ves_i-\ves_n)/2$.
When $n$ is odd, then $\sta_1\cdots\sta_n$ is a highest-weight
vector of $\Psi_{(1)}$ with  weight $(\sum_{i=1}^n\ves_i)/2$ and
$\sta_1\cdots\sta_{n-1}$ is a highest-weight vector of $\Psi_{(0)}$
 with weight $(\sum_{i=1}^{n-1}\ves_i-\ves_n)/2$.}

\section{Noncanonical Oscillator Representations}

In this section, we determine the structure of the noncanonical
oscillator representations of even orthogonal Lie algebras obtained
from the above bosonic representations by partially swapping
differential operators and multiplication operators, which are
generalizations of the classical theorem on harmonic polynomials.

Let ${\msr B}=\mbb{F}[x_1,...,x_n,y_1,...,y_n]$. Recall the
canonical oscillator representation of $o(2n,\mbb{F})$ is given in
(7.1.11) and (7.1.12). Fix
 $n_1,n_2\in\ol{1,n}$ with $n_1\leq n_2$.
  Swapping operators $\ptl_{x_r}\mapsto -x_r,\;
 x_r\mapsto
\ptl_{x_r}$  for $r\in\ol{1,n_1}$ and $\ptl_{y_s}\mapsto -y_s,\;
 y_s\mapsto\ptl_{y_s}$  for $s\in\ol{n_2+1,n}$ in the canonical
oscillator representation (7.1.11) and (7.1.12), we get another
noncanonical oscillator representation of
 $o(2n,\mbb{F})$ on
  ${\msr B}$ determined by
$$(E_{i,j}-E_{n+j,n+i})|_{\msr
B}=E_{i,j}^x-E_{j,i}^y\eqno(7.2.1)$$ with
$$E_{i,j}^x=\left\{\begin{array}{ll}-x_j\ptl_{x_i}-\delta_{i,j}&\mbox{if}\;
i,j\in\ol{1,n_1};\\ \ptl_{x_i}\ptl_{x_j}&\mbox{if}\;i\in\ol{1,n_1},\;j\in\ol{n_1+1,n};\\
-x_ix_j &\mbox{if}\;i\in\ol{n_1+1,n},\;j\in\ol{1,n_1};\\
x_i\partial_{x_j}&\mbox{if}\;i,j\in\ol{n_1+1,n};
\end{array}\right.\eqno(7.2.2)$$
and
$$E_{i,j}^y=\left\{\begin{array}{ll}y_i\ptl_{y_j}&\mbox{if}\;
i,j\in\ol{1,n_2};\\ -y_iy_j&\mbox{if}\;i\in\ol{1,n_2},\;j\in\ol{n_2+1,n};\\
\ptl_{y_i}\ptl_{y_j} &\mbox{if}\;i\in\ol{n_2+1,n},\;j\in\ol{1,n_2};\\
-y_j\partial_{y_i}-\delta_{i,j}&\mbox{if}\;i,j\in\ol{n_2+1,n};
\end{array}\right.\eqno(7.2.3)$$
and
$$E_{i,n+j}|_{\msr B}=\left\{\begin{array}{ll}
\ptl_{x_i}\ptl_{y_j}&\mbox{if}\;i\in\ol{1,n_1},\;j\in\ol{1,n_2};\\
-y_j\ptl_{x_i}&\mbox{if}\;i\in\ol{1,n_1},\;j\in\ol{n_2+1,n};\\
x_i\ptl_{y_j}&\mbox{if}\;i\in\ol{n_1+1,n},\;j\in\ol{1,n_2};\\
-x_iy_j&\mbox{if}\;i\in\ol{n_1+1,n},\;j\in\ol{n_2+1,n};\end{array}\right.\eqno(7.2.4)$$
and
$$E_{n+i,j}|_{\msr B}=\left\{\begin{array}{ll}
-x_jy_i&\mbox{if}\;j\in\ol{1,n_1},\;i\in\ol{1,n_2};\\
-x_j\ptl_{y_i}&\mbox{if}\;j\in\ol{1,n_1},\;i\in\ol{n_2+1,n};\\
y_i\ptl_{x_j}&\mbox{if}\;j\in\ol{n_1+1,n},\;i\in\ol{1,n_2};\\
\ptl_{x_j}\ptl_{y_i}&\mbox{if}\;j\in\ol{n_1+1,n},\;i\in\ol{n_2+1,n}.\end{array}\right.\eqno(7.2.5)$$
Recall
$$\td\Dlt=\sum_{i=1}^{n_1}x_i\ptl_{y_i}-\sum_{r=n_1+1}^{n_2}\ptl_{x_r}\ptl_{y_r}+\sum_{s=n_2+1}^n
y_s\ptl_{x_s}\eqno(7.2.6)$$ and its dual
$$\eta=\sum_{i=1}^{n_1}y_i\ptl_{x_i}+\sum_{r=n_1+1}^{n_2}x_ry_r+\sum_{s=n_2+1}^n
x_s\ptl_{y_s}.\eqno(7.2.7)$$ Moreover,
$$[\td\Dlt,\eta]=n_1-n_2-\wht D,\;\;\;\wht D=
-\sum_{i=1}^{n_1}x_i\ptl_{x_i}
+\sum_{r=n_1+1}^nx_r\ptl_{x_r}+\sum_{j=1}^{n_2}y_j\ptl_{y_j}-\sum_{s=n_2+1}^ny_s\ptl_{y_s}\eqno(7.2.8)$$
and as operators on ${\msr B}$,
$$\xi\td\Dlt=\td\Dlt\xi,\;\;\xi\eta=\eta\xi,\;\;\xi \wht D=\wht D\xi\qquad\for\;\;\xi\in
o(2n,\mbb{F}).\eqno(7.2.9)$$

 Set
$${\msr B}_{\la k\ra}=\mbox{Span}\{x^\al
y^\be\mid\al,\be\in\mbb{N}^n;\sum_{r=n_1+1}^n\al_r-\sum_{i=1}^{n_1}\al_i+
\sum_{i=1}^{n_2}\be_i-\sum_{r=n_2+1}^n\be_r=k\}\eqno(7.2.10)$$ for
$k\in\mbb{Z}$. Then
$${\msr B}_{\la k\ra}=\{u\in{\msr B}\mid \wht D(u)=k
u\}\eqno(7.2.11)$$ forms an $o(2n,\mbb{F})$-submodule by the third
equation (7.2.9). Define
$${\msr H}_{\la k\ra}=\{f\in {\msr B}_{\la
k\ra}\mid \td\Dlt(f)=0\}.\eqno(7.2.12)$$ According to the first
equation (7.2.9), ${\msr H}_{\la k\ra}$ forms an
$o(2n,\mbb{F})$-submodule.
 The
following is our main theorem:\psp

{\bf Theorem 7.2.1}. {\it For any $n_1-n_2+1-\dlt_{n_1,n_2}\geq
k\in\mbb{Z}$, ${\msr H}_{\la k\ra}$ is an infinite-dimensional
irreducible $o(2n,\mbb{F})$-submodule and ${\msr B}_{\la
k\ra}=\bigoplus_{i=0}^\infty\eta^i({\msr H}_{\la k-2i\ra})$ is a
decomposition of irreducible submodules. In particular, ${\msr
B}_{\la k\ra}={\msr H}_{\la k\ra}\oplus \eta({\msr B}_{\la
k-2\ra})$. The module ${\msr H}_{\la k\ra}$ under the assumption is
of highest-weight type only if $n_2=n$, and its highest weight is
$-k\lmd_{n_1-1}+(k-1)\lmd_{n_1}$ if $n_1>1$ and $(k-1)\lmd_1$ if
$n_1=1$.}

{\it Proof}. We will prove it case by case. According to (7.1.19),
the representation of $o(2n,\mbb{F})$ given in (7.2.1)-(7.2.5) is an
extension of that for $sl(n,\mbb{F})$ given in (6.3.36)-(6.3.38) on
${\msr B}$. Moreover,
$${\msr B}_{\la k\ra}=\bigoplus_{\ell\in\mbb{Z}}{\msr
B}_{\la\ell,k-\ell\ra},\qquad{\msr H}_{\la
k\ra}=\bigoplus_{\ell\in\mbb{Z}}{\msr
H}_{\la\ell,k-\ell\ra}\eqno(7.2.13)$$ by (6.4.29), (6.4.31),
(7.2.10) and (7.2.12). \psp

 {\it Case 1}.  $n_1+1<n_2$ {\it and}
$n_1-n_2+1\geq k\in\mbb{Z}$. \psp

According to (6.4.40) and (6.4.41), the $sl(n,\mbb{F})$-singular
vectors in ${\msr H}_{\la k\ra}$ are: for $m_1,m_2\in\mbb{N}$,
$$x_{n_1}^{m_1}y_{n_2+1}^{m_2}\qquad\mbox{with}\;-(m_1+m_2)=k,\eqno(7.2.14)$$
$$x_{n_1+1}^{m_1}y_{n_2+1}^{m_2}\qquad\mbox{with}\;m_1-m_2=k,\eqno(7.2.15)$$
$$x_{n_1}^{m_1}y_{n_2}^{m_2}\qquad\mbox{with}\;-m_1+m_2=k.\eqno(7.2.16)$$

Note
$$(E_{n+n_2+1,n_1}-E_{n+n_1,n_2+1})|_{\msr
B}=-x_{n_1}\ptl_{y_{n_2+1}}-y_{n_1}\ptl_{x_{n_2+1}}\eqno(7.2.17)$$
by (7.2.5). So
$$(E_{n+n_2+1,n_1}-E_{n+n_1,n_2+1})^{m_2}(x_{n_1}^{m_1}y_{n_2+1}^{m_2})=(-1)^{m_2}m_2!x_{n_1}^{-k}\eqno(7.2.18)$$
for the vectors in (7.2.14). Moreover,
$$(E_{n+n_2+1,n_1+1}-E_{n+n_1+1,n_2+1})|_{\msr
B}=\ptl_{x_{n_1+1}}\ptl_{y_{n_2+1}}-y_{n_1+1}\ptl_{x_{n_2+1}}\eqno(7.2.19)$$
again by (7.2.5), which implies
$$(E_{n+n_2+1,n_1+1}-E_{n+n_1+1,n_2+1})^{m_1}(x_{n_1+1}^{m_1}y_{n_2+1}^{m_2})=m_1![\prod_{r=0}^{m_1-1}(m_2-r)]
y_{n_2+1}^{-k}\eqno(7.2.20)$$ for the vectors in (7.2.15).
Furthermore,
$$(E_{n_1,n+n_2}-E_{n_2,n+n_1})|_{\msr
B}=\ptl_{x_{n_1}}\ptl_{y_{n_2}}-x_{n_2}\ptl_{y_{n_1}}\eqno(7.2.21)$$
by (7.2.4), which implies
$$(E_{n_1,n+n_2}-E_{n_2,n+n_1})^{m_2}(x_{n_1}^{m_1}y_{n_2}^{m_2})=
m_2![\prod_{r=0}^{m_2-1}(m_1-r)] x_{n_1}^{-k}\eqno(7.2.22)$$ for the
vectors in (7.2.16).

On the other hand,
$$(E_{n_1,n+n_2+1}-E_{n_2+1,n+n_1})|_{\msr
B}=-y_{n_2+1}\ptl_{x_{n_1}}-x_{n_2+1}\ptl_{y_{n_1}}\eqno(7.2.23)$$
by (7.2.4), which implies
$$(E_{n_1,n+n_2+1}-E_{n_2+1,n+n_1})^{m_2}(x_{n_1}^{-k})=(-1)^{m_2}[\prod_{r=0}^{m_2-1}(-k-r)]
x_{n_1}^{m_1}y_{n_2+1}^{m_2}\eqno(7.2.24)$$ for the vectors in
(7.2.14). Moreover,
$$(E_{n_1+1,n+n_2+1}-E_{n_2+1,n+n_1+1})|_{\msr
B}=-x_{n_1+1}y_{n_2+1}-x_{n_2+1}\ptl_{y_{n_1+1}}\eqno(7.2.25)$$ by
(7.2.4), which implies
$$(E_{n_1+1,n+n_2+1}-E_{n_2+1,n+n_1+1})^{m_1}(y_{n_2+1}^{-k})=(-1)^{m_1}
x_{n_1+1}^{m_1}y_{n_2+1}^{m_2} \eqno(7.2.26)$$ for the vectors in
(7.2.15). Furthermore,
$$(E_{n+n_2,n_1}-E_{n+n_1,n_2})|_{\msr
B}=-x_{n_1}y_{n_2}-y_{n_1}\ptl_{x_{n_2}}\eqno(7.2.27)$$ by (7.2.5),
which implies
$$(E_{n+n_2,n_1}-E_{n+n_1,n_2})^{m_2}(x_{n_1}^{-k})=
(-1)^{m_2}x_{n_1}^{m_1}y_{n_2}^{m_2}\eqno(7.2.28)$$ for the vectors
in (7.2.16). Thus for any two vectors in (7.2.14)-(7.2.16), there
exists an element in the universal enveloping algebra
$U(o(2n,\mbb{F}))$ which carries one to another. On the other hand,
the vectors in (7.2.14)-(7.2.16) have distinct weights. Thus any
nonzero $o(2n,\mbb F)$-submodule of ${\msr H}_{\la k\ra}$ must
contain one of the vectors in (7.2.14)-(7.2.16). Hence all the
vectors in (7.2.14)-(7.2.16) are in the submodule by
(7.2.17)-(7.2.28). Therefore, the submodule must contains all ${\msr
H}_{\la\ell,k-\ell\ra}$ for $\ell\in\mbb{Z}$ by Theorem 6.4.3, or
equivalently, it is  equal to ${\msr H}_{\la k\ra}$ due to (7.2.13).
So ${\msr H}_{\la k\ra}$ is irreducible. By (7.2.25) and (7.2.27),
${\msr H}_{\la k\ra}$ is not of highest-weight type. For any
$m\in\mbb{N}$, $n_1-n_2+1\geq k-m$ and so ${\msr H}_{\la k-m\ra}$ is
an irreducible $o(2n,\mbb{F})$-submodule.

Since $\td\Dlt$ is locally nilpotent, for any $0\neq u\in {\msr
B}_{\la k\ra}$, there exists an element $\kappa(u)\in\mbb{N}$ such
that
$$\td\Dlt^{\kappa(u)}(u)\neq
0\;\;\mbox{and}\;\;\td\Dlt^{\kappa(u)+1}(u)=0.\eqno(7.2.29)$$ Set
$$V=\sum_{i=0}^\infty\eta^i({\msr
H}_{\la k-2i\ra}).\eqno(7.2.30)$$ Given $0\neq u\in {\msr B}_{\la
k\ra}$, $\kappa(u)=0$ implies $u\in {\msr H}_{\la k\ra}\subset V$.
Suppose that $u\in V$ whenever $\kappa(u)<r$ for some positive
integer $r$. Assume $\kappa(u)=r$. First
$$v=\td\Dlt^r(u)\in {\msr
H}_{\la k-2r\ra}\subset V.\eqno(7.2.31)$$ Note
$$\td\Dlt^r[\eta^r(v)]=r![\prod_{i=1}^r(n_1-n_2-k+r+i)]v\eqno(7.2.32)$$
by (7.2.8) and (7.2.11). Thus we  have either
$$u=\frac{1}{r![\prod_{i=1}^r(n_1-n_2-k+r+i)]}\eta^r(v)\in V\eqno(7.2.33)$$
or
$$\kappa\left(u-\frac{1}{r![\prod_{i=1}^r(n_1-n_2-k+r+i)]}\eta^r(v)\right)<r.\eqno(7.2.34)$$
By induction,
$$u-\frac{1}{r![\prod_{i=1}^r(n_1-n_2-k+r+i)]}\eta^r(v)\in V,\eqno(7.2.35)$$
which implies $u\in V$. Therefore, we have $V={\msr B}_{\la k\ra}$.
Since the weight of any $sl(n,\mbb{F})$-singular vector in
$\eta^i({\msr H}_{\la k-2i\ra})$ is different from  the weight of
any $sl(n,\mbb{F})$-singular vector in $\eta^j({\msr H}_{\la
k-2j\ra})$ when $i\neq j$,  the sums in Theorem 7.2.1 are direct
sums. \psp

{\it Case 2}. $n_1+1=n_2$ {\it and} $0\geq k\in\mbb{Z}$.\psp

Assume $n_1+1=n_2<n$. By (6.5.60), (7.2.8) and (7.2.11), the
$sl(n,\mbb{F})$-singular vectors in ${\msr H}_{\la k\ra}$ are those
in (7.2.14)-(7.2.16). So Theorem 7.2.1 holds by the arguments in
Case 1. Suppose $n_1<n_2=n$. According to  (6.5.62), (7.2.8) and
(7.2.11), the $sl(n,\mbb{F})$-singular vectors in ${\msr H}_{\la
k\ra}$ are those in (7.2.16). Expressions (7.2.22)
 and (7.2.28)-(7.2.35) imply the conclusions in the Theorem 7.2.1.
\psp

{\it Case 3}. $n_1=n_2$ {\it and} $0\geq k\in\mbb{Z}$.\psp

 Recall
$$\zeta_1=x_{n_1-1}y_{n_1}-x_{n_1}y_{n_1-1},\;\;
\;\zeta_2=x_{n_1+1}y_{n_1+2}-x_{n_1+2}y_{n_1+1}. \eqno(7.2.36)$$
Suppose $n_1=n_2<n-1$. Lemma 6.6.2 tells us that the
$sl(n,\mbb{F})$-singular vectors in ${\msr H}_{\la k\ra}$ are those
in (7.2.14) and
$$x_{n_1}^{-k}\zeta_1^{m+1}\qquad\for\;\;m\in\mbb{N},\eqno(7.2.37)$$
$$y_{n_1+1}^{-k}\zeta_2^{m+1}\qquad\for\;\;m\in\mbb{N}.\eqno(7.2.38)$$
Again all the singular vectors have distinct weights. If ${\msr N}$
is a nonzero submodule of ${\msr H}_{\la k\ra}$, then ${\msr N}$
must contain one of the above $sl(n,\mbb{F})$-singular vectors. If
${\msr N}$ contains a singular vector in (7.2.14), then
$x_{n_1}^{-k}\in N$ by (7.2.18). Suppose
$x_{n_1}^{-k}\zeta_1^{m+1}\in N$ for some $m\in\mbb{N}$.
 Note
$$(E_{n_1-1,n+n_1}-E_{n_1,n+n_1-1})|_{\msr
B}=\ptl_{x_{n_1-1}}\ptl_{y_{n_1}} -\ptl_{x_{n_1}}\ptl_{y_{n_1-1}}
\eqno(7.2.39)$$ by (7.2.4). Thus
\begin{eqnarray*}\qquad&
&(E_{n_1-1,n+n_1}-E_{n_1,n+n_1-1})^{m+1}(x_{n_1}^{-k}\zeta_1^{m+1})
\\ &=&\left[\sum_{r=0}^{m+1}(-1)^r{m+1\choose
r}(\ptl_{x_{n_1-1}}\ptl_{y_{n_1}})^{m+1-r}(\ptl_{x_{n_1}}\ptl_{y_{n_1-1}})^r\right]
\\ &&\times\left[\sum_{s=0}^{m+1}(-1)^s{m+1\choose
s}(x_{n_1-1}y_{n_1})^{m+1-s}x_{n_1}^{-k+s}y_{n_1-1}^s\right]
\\&=&\left(\sum_{r=0}^{m+1}{m+1\choose
r}^2[(m+1-r!)]^2r![\prod_{i=1}^r(-k+i)]\right)x_{n_1}^{-k}
\\ &=&[(m+1)!]^2\left(\sum_{r=0}^{m+1}{-k+r\choose
r}\right)x_{n_1}^{-k}\in N.\hspace{5.8cm}(7.2.40)\end{eqnarray*} So
we have $x_{n_1}^{-k}\in {\msr N}$ again. Symmetrically,
$y_{n_1+1}^{-k}\in{\msr N}$ if
 $y_{n_1+1}^{-k}\zeta_2^{m+1}\in {\msr N}$ for some $m\in\mbb{N}$.
Observe
$$(E_{n+n_1,n_1+1}-E_{n+n_1+1,n_1})|_{\msr
B}=y_{n_1}\ptl_{x_{n_1+1}}+x_{n_1}\ptl_{y_{n_1+1}}\eqno(7.2.41)$$ by
(7.2.5). Thus
$$(E_{n+n_1,n_1+1}-E_{n+n_1+1,n_1})^{-k}(y_{n_1+1}^{-k})=(-k)!x_{n_1}^{-k}\in{\msr
N}.\eqno(7.2.42)$$ Thus we always have  $x_{n_1}^{-k}\in {\msr N}$.

According to (7.2.24), ${\msr N}$ contains all the singular vectors
in (7.2.14). Observe that
$$(E_{n+n_1-1,n_1}-E_{n+n_1,n_1-1})|_{\msr B}=\zeta_1,\;\;
(E_{n_1+2,n+n_1+1}-E_{n_1+1,n+n_1+2})|_{\msr
B}=\zeta_2\eqno(7.2.43)$$ as multiplication operators on ${\msr B}$
by (7.2.4) and (7.2.5). Thus
$$(E_{n+n_1-1,n_1}-E_{n+n_1,n_1-1})^{m+1}(x_{n_1}^{-k})=x_{n_1}^{-k}\zeta_1^{m+1},\eqno(7.2.44)$$
$$(E_{n_1+2,n+n_1+1}-E_{n_1+1,n+n_1+2})^{m+1}(x_{n_1}^{-k})=x_{n_1}^{-k}\zeta_2^{m+1}\in
{\msr N}.\eqno(7.2.45)$$ Note that
$$(E_{n_1+1,n+n_1}-E_{n_1,n+n_1+1})|_{\msr
B}=x_{n_1+1}\ptl_{y_{n_1}}+y_{n_1+1}\ptl_{x_{n_1}}\eqno(7.2.46)$$ by
(7.2.4). So
$$\frac{1}{(-k)!}(E_{n_1+1,n+n_1}-E_{n_1,n+n_1+1})^{-k}(x_{n_1}^{-k}\zeta_2^{m+1})
=y_{n_1+1}^{-k}\zeta_2^{m+1}\in{\msr N}.\eqno(7.2.47)$$ Thus ${\msr
N}$ contains all the $sl(n,\mbb{F})$-singular vectors in ${\msr
H}_{\la k\ra}$, which implies that it contains all ${\msr
H}_{\la\ell,k-\ell\ra}\subset {\msr H}_{\la k\ra}$ by Theorem 6.6.3.
So ${\msr N}={\msr H}_{\la k\ra}$ by (7.2.13); that is, ${\msr
H}_{\la k\ra}$ is an irreducible $o(2n,\mbb{F})$-module. By
(7.2.29)-(7.2.35), the direct sums in Theorem 7.2.1 holds. Observe
that Theorem 7.2.1 under the subcase $n_1=n-1,n_2=n$ is implied by
(6.6.30), (6.6.33), (6.6.36), (6.6.39) and the related partial
arguments in the above. This completes the proof of Theorem 7.2.1
$\qquad\Box$\psp

Suppose $n_1<n_2$. Taking $T_1=-\ptl_{x_{n_1+1}}\ptl_{y_{n_1+1}}$,
$T_1^-=-\int_{(x_{n_1+1})}\int_{(y_{n_1+1})}$ and
$T_2=\td\dlt+\ptl_{x_{n_1}+1}\ptl_{y_{n_1}+1}$ in Lemma 6.1.1,
${\msr H}_{\la k\ra}$ has a basis
\begin{eqnarray*}\qquad& &\big\{\sum_{i=0}^\infty\frac{(x_{n_1+1}y_{n_1+1})^i(\td\Dlt+\ptl_{x_{n_1}+1}\ptl_{y_{n_1}+1})^i(x^\al
y^\be)}{\prod_{r=1}^i(\al_{n_1+1}+r)(\be_{n_1+1}+r)}
\mid\al,\be\in\mbb{N}^n;\\ &
&\al_{n_1+1}\be_{n_1+1}=0;-\sum_{i=1}^{n_1}\al_i+\sum_{r=n_1+1}^n\al_r+
\sum_{i=1}^{n_2}\be_i-\sum_{r=n_2+1}^n\be_r=k
\big\}\hspace{2cm}(7.2.48)\end{eqnarray*}

Finally, we want to find an expression for ${\msr H}_{\la k\ra}$ for
$0\geq k\in\mbb{Z}$ when $n_1=n_2$. First we assume $n_1=n_2=1$ and
$n=2$. According to (6.6.36), (6.6.52) and (7.2.13),
$${\msr H}_{\la
-k\ra}=\mbox{Span}\{[x_1^r y_2^s(x_1x_2- y_1y_2)^l\mid
r,s,l\in\mbb{N};r+s=k\}.\eqno(7.2.49)$$ Next we consider the subcase
$n_1=n_2=1$ and $n\geq 3$. According to (6.6.30), (6.6.52), (6.6.54)
and (7.2.13)
\begin{eqnarray*}& &{\msr H}_{\la
-k\ra}\\ &=&\mbox{Span}\{[\prod_{r=2}^ny_r^{\hat l_r}][\prod_{2\leq
p<q\leq n}(x_py_q-x_qy_p)^{\hat k_{p,q}}][\prod_{s=2}^n(x_1x_s-
y_1y_s)^{\hat l_s}],x_1^l [\prod_{s=2}^ny_s^{k_s}]\\
&&\times[\prod_{s=2}^n(x_1x_s- y_1y_s)^{l_s}]\mid l,k_s,l_s,\hat
l,\hat k_{p,q},\hat l_s\in\mbb{N};l+\sum_{s=2}^nk_s=\sum_{r=2}^n\hat
l_r=k\}.\hspace{1.5cm}(7.2.50)\end{eqnarray*}

 Assume $1<n_1=n_2<n-1$.  By (6.6.25), (6.6.52)-(6.6.54) and (7.2.13), we have
\begin{eqnarray*}& &{\msr H}_{\la
-k\ra}\\ &=&\mbox{Span}\{[\prod_{r=1}^{n_1}x_r^{l_r'}] [\prod_{1\leq
p<q\leq
n_1}(x_py_q-x_qy_p)^{k_{p,q}'}][\prod_{r=1}^{n_1}\prod_{s=n_1+1}^n(x_rx_s-
y_ry_s)^{l_{r,s}'}],\\
& &[\prod_{r=1}^{n-n_1}y_{n_1+r}^{\hat l_r}][\prod_{n_1+1\leq
p<q\leq n}(x_py_q-x_qy_p)^{\hat
k_{p,q}}][\prod_{r=1}^{n_1}\prod_{s=n_1+1}^n(x_rx_s- y_ry_s)^{\hat
l_{r,s}}],\\ & &[\prod_{r=1}^{n_1}x_r^{l_r}]
[\prod_{s=1}^{n-n_1}y_{n_1+s}^{k_s}][\prod_{r=1}^{n_1}\prod_{s=1}^{n-n_1}(x_rx_{n_1+s}-
y_ry_{n_1+s})^{l_{r,s}}]\mid l_r,k_s,l_{r,s},l_r',k_{p,q}',\\ &
&l_{r,s}',\hat l_r,\hat k_{p,q},\hat
l_{r,s}\in\mbb{N};\sum_{r=1}^{n_1}l_r+\sum_{s=1}^{n-n_1}k_s=\sum_{r=1}^{n_1}l_r'=\sum_{r=1}^{n-n_1}\hat
l_r=k\}.\hspace{3.2cm}(7.2.51)\end{eqnarray*}

 Consider the subcase
$n_1=n_2=n-1$ and $n\geq 3$. By (6.6.33), (6.6.52), (6.6.53) and
(7.2.13), we obtain
\begin{eqnarray*}& &{\msr H}_{\la
-k\ra}\\ &=&\mbox{Span}\{[\prod_{r=1}^{n-1}x_r^{l_r'}] [\prod_{1\leq
p<q\leq n-1}(x_py_q-x_qy_p)^{k_{p,q}'}][\prod_{r=1}^{n-1}(x_rx_n-
y_ry_n)^{\bar l_r'}],[\prod_{r=1}^{n-1}x_r^{l_r}] y_n^{\hat k}\\
&&\times[\prod_{r=1}^{n-1}(x_rx_n- y_ry_n)^{\bar l_r}]\mid l_r,\hat
k,\bar l_r,l_r',k_{p,q}',\bar
l_r'\in\mbb{N};\sum_{r=1}^{n-1}l_r+\hat
k=\sum_{r=1}^{n-1}l_r'=k\}.\hspace{1.2cm}(7.2.52)\end{eqnarray*} At
last, we assume $n_1=n_2=n$. By (6.6.39), (6.6.53) and (7.2.13),
$${\msr H}_{\la -k\ra}=\mbox{Span}\{\prod_{r=1}^nx_r^{l_r}] [\prod_{1\leq
p<q\leq n}(x_py_q-x_qy_p)^{k_{p,q}}]\mid
l_r,k_{p,q}\in\mbb{N};\sum_{r=1}^nl_r=k\}.\eqno(7.2.53)$$

\section{Extensions of the Conformal Representation}

In  this section, we construct a new functor from the category
$o(2n,\mbb C)$-{\bf mod} to the category $o(2n+2,\mbb C)$-{\bf mod},
which is an extension of the conformal representation of
$o(2n+2,\mbb C)$.

The $n$-dimensional conformal group with respect to Euclidean metric
$(\cdot,\cdot)$ is generated by the translations, rotations,
dilations and special conformal transformations
$$\vec x\mapsto\frac{\vec x-(\vec x,\vec x)\vec b}{(\vec b,\vec b)
(\vec x,\vec x)-2(\vec b,\vec x)+1}.\eqno(7.3.1)$$ The above
transformations yield a nonhomogeneous representation of the real
Lie algebra $o(2,n)$. We assume the base field $\mbb F=\mbb C$.

For $\lmd\in\Lmd^+$, we denote by $V(\lmd)$ the finite-dimensional
irreducible $o(2n,\mbb{C})$-module with highest weight $\lmd$. The
$2n$-dimensional natural module of $o(2n,\mbb{C})$ is $V(\ves_1)$
with the weights $\{\pm\ves_i\mid i\in\ol{1,n}\}.$ Recall
$$\rho=\frac{1}{2}\sum_{\nu\in\Phi_{D_n}^+}\nu.\eqno(7.3.2)$$
Then
$$(\rho,\nu)=1\qquad\for\;\;\nu\in\Pi_{D_n}\eqno(7.3.3)$$
by (3.2.17). Expression (7.1.7) yields
$$\rho=\sum_{i=1}^{n-1}(n-i)\ves_i.\eqno(7.3.4)$$
 The following
result can be derived from  Theorem 5.4.3:\psp

{\bf Lemma 7.3.1 (Pieri's formula)}. \index{Pieri's formula!for
$o(2n)$}{\it Given $\mu\in\Lmd^+$ with ${\msr
S}(\mu)=\{i_1,i_2,...,i_{s+1}\}$ (cf. (6.7.54) and (6.7.55)),
$$V(\ves_1)\otimes_\mbb{C}V(\mu)\cong \bigoplus_{r=1}^s
(V(\mu+\ves_{i_r})\oplus V(\mu-\ves_{i_{r+1}-1}))\eqno(7.3.5)$$ if
$\mu_{n-1}+\mu_n>0$ and
$$V(\ves_1)\otimes_\mbb{C}V(\mu)\cong \bigoplus_{r=1}^{s-2+\dlt_{\mu_{_n},0}}
V(\mu-\ves_{i_{r+1}-1})\oplus \bigoplus_{r=1}^s
V(\mu+\ves_{i_r})\eqno(7.3.6)$$ when $\mu_{n-1}+\mu_n=0$.}\psp

Denote by $U({\msr G})$ the universal enveloping algebra of a Lie
algebra ${\msr G}$. The algebra $U({\msr G})$  can be imbedded into
the tensor algebra $U({\msr G})\otimes U({\msr G})$ by the
associative algebra homomorphism $\mfk{d}: U({\msr G}) \rightarrow
U({\msr G})\otimes_\mbb{C} U({\msr G})$ determined  by
$$\mfk{d}(u)=u\otimes 1 +1 \otimes u \qquad \mbox{ for} \ u\in
\msr{G}.\eqno(7.3.7)$$

Note that the Casimir element of $o(2n,\mbb{C})$ is
\begin{eqnarray*}\omega&=&\sum_{1\leq
i<j\leq
n}[(E_{i,n+j}-E_{j,n+i})(E_{n+j,i}-E_{n+i,j})+(E_{n+j,i}-E_{n+i,j})(E_{i,n+j}-E_{j,n+i})]
\\ & &+\sum_{i,j=1}^n(E_{i,j}-E_{n+j,n+i})(E_{j,i}-E_{n+i,n+j})\in
U(o(2n,\mbb C)).\hspace{3.8cm}(7.3.8)\end{eqnarray*} Set
$$\td\omega=\frac{1}{2}(\mfk{d}(\omega)-\omega\otimes 1-1\otimes
\omega)\in U(o(2n,\mbb C))\otimes_\mbb{C} U(o(2n,\mbb
C)).\eqno(7.3.9)$$ By (7.3.8),
\begin{eqnarray*}\qquad\qquad\;\;\;\td\omega\!&=&\!\sum_{1\leq
i<j\leq n}[(E_{i,n+j}-E_{j,n+i})\otimes(E_{n+j,i}-E_{n+i,j})\\
& &+(E_{n+j,i}-E_{n+i,j})\otimes(E_{i,n+j}-E_{j,n+i})]
\\ & &+\sum_{i,j=1}^n(E_{i,j}-E_{n+j,n+i})\otimes(E_{j,i}-E_{n+i,n+j}).\hspace{3.8cm}(7.3.10)\end{eqnarray*}

For any $\mu\in \Lmd^+$, we have
$$\omega|_{V(\mu)}=(\mu+2\rho,\mu)\mbox{Id}_{V(\mu)}.\eqno(7.3.11)$$
Denote
$$\ell^+_i(\mu)=\dim
V(\mu+\ves_i)\qquad\mbox{if}\;\;\mu+\ves_i\in\Lmd^+\eqno(7.3.12)$$
and
$$\ell^-_i(\mu)=\dim
V(\mu-\ves_i)\qquad\mbox{if}\;\;\mu-\ves_i\in\Lmd^+.\eqno(7.3.13)$$
Observe that
$$(\mu+\ves_i+2\rho,\mu+\ves_i)-(\mu+2\rho,\mu)-(\ves_1+2\rho,\ves_1)=2(\mu_i+1-i)\eqno(7.3.14)$$
and
$$(\mu-\ves_i+2\rho,\mu-\ves_i)-(\mu+2\rho,\mu)-(\ves_1+2\rho,\ves_1)=2(1+i-2n-\mu_i)\eqno(7.3.15)$$
for $\mu=\sum_{r=1}^n\mu_r\ves_r$ by (7.3.10). Moreover, the algebra
$U(o(2n,\mbb C))\otimes_{\mbb C}U(o(2n,\mbb C))$ acts on
$V(\ves_1)\otimes_\mbb{C}V(\mu)$ by
$$(\xi_1\otimes \xi_2)(v\otimes
u)=\xi_1(v)\otimes\xi_2(u)\;\;\for\;\;\xi_1,\xi_2\in U(o(2n,\mbb
F)),\;v\in V(\ves_1),\;u\in V(\mu).\eqno(7.3.16)$$ By Lemma 7.3.1,
(7.3.9) and (7.3.12)-(7.3.15), we obtain:\psp

{\bf Lemma 7.3.2}. {\it Let $\mu=\sum_{i=1}^n\mu_i\ves_i\in\Lmd^+$
with ${\msr S}(\mu)=\{i_1,i_2,...,i_{s+1}\}$ (cf. (6.7.54) and
(6.7.55)). If $\mu_{n-1}+\mu_n>0$, the characteristic polynomial of
$\td\omega|_{V(\ves_1)\otimes_\mbb{C}V(\mu)}$ is
$$\prod_{r=1}^s[(t-\mu_{i_r}+i_r-1)^{\ell^+_{i_r}(\mu)}
(t+\mu_{i_r}+2n-i_{r+1})^{\ell^-_{i_{r+1}-1}(\mu)}].\eqno(7.3.17)$$
When $\mu_{n-1}+\mu_n=0$, the characteristic polynomial of
$\td\omega|_{V(\ves_1)\otimes_\mbb{C}V(\mu)}$ is
$$[\prod_{r=1}^s[(t-\mu_{i_r}+i_r-1)^{\ell^+_{i_r}(\mu)}]
[\prod_{\iota=1}^{s-2+\dlt_{\mu_n,0}}(t+\mu_{i_\iota}+2n-i_{\iota+1})^{\ell^-_{i_{\iota+1}-1}(\mu)}].\eqno(7.3.18)$$}
\pse

We remark that the above lemma is equivalent to an explicit version
of Kostant's characteristic identity (cf. [Kb]).

Denote
$$A_{i,j}=E_{i,j}-E_{n+1+j,n+1+i},\;B_{i,j}=E_{i,n+1+j}-E_{j,n+1+i},\;C_{i,j}=E_{n+1+i,j}-E_{n+1+j,i}\eqno(7.3.19)$$
for $i,j\in\ol{1,n+1}$. Then
 the split orthogonal Lie algebra
$$o(2n+2,\mbb C)=\sum_{i,j=1}^{n+1}\mbb CA_{i,j}+\sum_{1\leq i<j\leq
n+1}(\mbb CB_{i,j}+\mbb CC_{j,i})\eqno(7.3.20)$$ Recall
$${\msr A}=\mbb C[x_1,x_2,...,x_{2n}].\eqno(7.3.21)$$
Set
$$\eta=\sum_{i=1}^nx_ix_{n+i},\qquad
D=\sum_{r=1}^{2n}x_r\partial_{x_r}.\eqno(7.3.22)$$ Replace $n$ by
$2n$ and take the product $(\vec x,\vec
y)=(1/2)\sum_{i=1}^n(x_iy_{n+i}+x_{n+i}y_i)$ in (7.3.1).
 Using the derivatives of
one-parameter conformal transformations just as one derives the Lie
algebra of a Lie group, we get the following conformal
representation of $o(2n+2,\mbb C)$ determined by
$$A_{i,j}|_{\msr A}=x_i\ptl_{x_j}-x_{n+j}\ptl_{x_{n+i}},\qquad
B_{i,j}|_{\msr
A}=x_i\ptl_{x_{n+j}}-x_j\ptl_{x_{n+i}},\eqno(7.3.23)$$
$$C_{i,j}|_{\msr
A}=x_{n+i}\ptl_{x_j}-x_j\ptl_{x_{n+i}},\qquad A_{n+1,n+1}|_{\msr
A}=-D,\eqno(7.3.24)$$
$$A_{n+1,i}|_{\msr A}=\ptl_{x_i},\qquad B_{i,n+1}|_{\msr
A}=-\ptl_{x_{n+i}},\eqno(7.3.25)$$
$$A_{i,n+1}|_{\msr A}=-x_iD+\eta\ptl_{x_{n+i}},\;\;C_{n+1,i}|_{\msr
A}=x_{n+i}D-\eta\ptl_{x_i}\eqno(7.3.26)$$ for $i,j\in\ol{1,n}$.

Note that
$$\msr L=\sum_{i,j=1}^n\mbb CA_{i,j}+\sum_{1\leq i<j\leq
n}(\mbb CB_{i,j}+\mbb CC_{j,i})\eqno(7.3.27)$$ forms a Lie
subalgebra of $o(2n+2,\mbb C)$ that is isomorphic to $o(2n,\mbb C)$.
For convenience, we make the identification of $o(2n,\mbb C)$ with
$\msr L$ as follows:
$$E_{i,j}-E_{n+j,n+i}\leftrightarrow A_{i,j},\;\;
E_{i,n+j}-E_{j,n+i}\leftrightarrow
B_{i,j},\;\;E_{n+i,j}-E_{n+j,i}\leftrightarrow
C_{i,j}\eqno(7.3.28)$$ for $i,j\in\ol{1,n}$. Recall the Witt algebra
${\mbb W}_{2n}=\sum_{i=1}^{2n}{\msr A}\ptl_{x_i}$, and Shen's
monomorphism $\Im: \mbb W_{2n}\rta\wht{\mbb W}_{2n}=\mbb
W_{2n}\oplus gl(2n,{\msr A})$ (cf. (6.7.12)) given by
$$\Im(\sum_{i=1}^{2n}f_i\ptl_{x_i})=\sum_{i=1}^{2n}f_i\ptl_{x_i}\oplus\sum_{i,j=1}^n\ptl_{x_i}(f_j)E_{i,j}.
\eqno(7.3.29)$$ Now we have the Lie algebra monomorphism $\nu:
o(2n+2,\mbb C)\rta \widehat{\mbb W}_{2n}$ given by
$$\nu(\xi)=\Im(\xi|_{\msr A})\qquad\for\;\;\xi\in
o(2n+2,\mbb C).\eqno(7.3.30)$$ Denote
$$I_{2n}=\sum_{r=1}^{2n}E_{r,r}.\eqno(7.3.31)$$
According to the identification (7.3.28), we have
$$\nu(A_{i,j})=(x_i\ptl_{x_j}-x_{n+j}\ptl_{x_{n+i}})\oplus
A_{i,j},\;\;\nu(B_{i,j})=(x_i\ptl_{x_{n+j}}-x_j\ptl_{x_{n+i}})\oplus
B_{i,j},\eqno(7.3.32)$$
$$\nu(C_{i,j})=(x_{n+i}\ptl_{x_j}-x_{n+j}\ptl_{x_i})\oplus
C_{i,j},\;\;\nu(A_{n+1,n+1})=-(D\oplus I_{2n}),\eqno(7.3.33)$$
$$\nu(A_{n+1,i})=\ptl_{x_i},\qquad\nu(B_{i,n+1})=-\ptl_{x_{n+i}},\eqno(7.3.34)$$
$$\nu(A_{i,n+1})=(-x_iD+\eta\ptl_{x_{n+i}})\oplus[\sum_{p=1}^n(x_{n+p}B_{p,i}-
x_pA_{i,p})-x_iI_{2n}],\eqno(7.3.35)$$
$$\nu(C_{n+1,i})=(x_{n+i}D-\eta\ptl_{x_i})\oplus[\sum_{p=1}^n(x_pC_{i,p}-x_{n+p}A_{p,i})+x_{n+i}I_{2n}]
\eqno(7.3.36)$$ for $i,j\in\ol{1,n}$.

Observe that
$$\wht{\mbb W}^o_{2n}=\wht{\mbb W}_{2n}\oplus[o(2n,{\msr A})+{\msr
A}I_{2n}]\eqno(7.3.37)$$ form a Lie subalgebra of $\wht{\mbb
W}_{2n}$ and
$$\nu(o(2n+2,\mbb C))\subset \wht{\mbb W}^o_{2n}.\eqno(7.3.38)$$
Let  $M$ be an $o(2n,\mbb C)$-module and let $c\in\mbb C$ be a fixed
constant. Then
$$\widehat M={\msr A}\otimes_{\mbb C}M\eqno(7.3.39)$$
becomes a $\wht{\mbb W}^o_{2n}$-module with the action:
$$(d+f_1A+f_2I_{2n})(g\otimes
v)=(d(g)+cf_2g)\otimes v+f_1g\otimes A(v)\eqno(7.3.40)$$ for
$f_1,f_2,g\in{\msr A},\;A\in o(2n,\mbb C)$ and $v\in M$. Moreover,
we make $\widehat M$ an $o(2n+2,\mbb C)$-module with the action:
$$\xi(\varpi)=\nu(\xi)(\varpi)\qquad\for\;\;\xi\in o(2n+2,\mbb C),\;\varpi\in
\widehat M.\eqno(7.3.41)$$

Denote
$$\msr G=o(2n+2,\mbb C),\qquad \msr G_0=o(2n,\mbb C)+\mbb
CA_{n+1,n+1}\eqno(7.3.42)$$
$$\msr G_+=\sum_{i=1}^n(\mbb C A_{n+1,i}+\mbb CB_{i,n+1}),\qquad
\msr G_-=\sum_{i=1}^n(\mbb C A_{i,n+1}+\mbb
CC_{n+1,i})\eqno(7.3.43)$$ (cf. (7.3.28)). Then $\msr G_\pm$ are
abelian Lie subalgebras of $\msr G$ and $\msr G_0$ is a Lie
subalgebra of $\msr G$. Moreover,
$$\msr G=\msr G_-+\msr G_0+\msr G_+,\qquad[\msr G_0,\msr
G_\pm]\subset \msr G_\pm.\eqno(7.3.44)$$ By (7.3.32)-(7.3.34),
$$\msr G_+(1\otimes M)=\{0\},\;\;U(\msr G_0)(1\otimes M)=1\otimes
M.\eqno(7.3.45)$$ Thus
$$U(\msr G)(1\otimes M)=U(\msr G_-)(1\otimes M).\eqno(7.3.46)$$
Using (7.3.34), we can prove:\psp

{\bf Proposition 7.3.3}. {\it The map $M\mapsto U({\msr
G}_-)(1\otimes M)$ gives rise to a functor from the category of
$o(2n,\mbb C)$-modules to the category of $o(2n+2,\mbb C)$-modules.
In particular, it maps irreducible $o(2n,\mbb C)$-modules to
irreducible $o(2n+2,\mbb C)$-modules.}\psp

From another point view, $U({\msr G}_-)(1\otimes M)$ is a polynomial
extension from $o(2n,\mbb C)$-module $M$ to an $o(2n+2,\mbb
C)$-module. Next we want to study when $\wht M=U({\msr
G}_-)(1\otimes M)$. Write
$$x^{\alpha}=\prod_{i=1}^{2n}x_{i}^{\alpha_i}
,\;\;E^\al=\prod_{r=1}^n[(-A_{r,n+1})^{\al_r}C_{n+1,r}^{\al_{n+r}}]
\qquad\for\;\;\al=(\al_1,\al_2,...,\al_{2n})\in\mbb{N}^{2n}.\eqno(7.3.47)$$
 For $k\in\mbb{N}$, we set
$${\msr A}_k=\mbox{Span}\{x^\al\mid
\al \in\mbb{N}^{2n}, \ |\al|=k\},\;\; \widehat M_k={\msr
A}_k\otimes_\mbb{C}M\eqno(7.3.48)$$ (recall
$|\al|=\sum_{i=1}^{2n}\al_i$) and
 \begin{eqnarray*}\hspace{1cm}(U(\msr G_-)(1\otimes
M))_k&=&U(\msr G_-)(1\otimes M)\bigcap \widehat M_k\\ &=&
\mbox{Span}\{ E^\al(1\otimes M)\mid  \al \in\mbb{N}^{2n}, \
|\al|=k\}\hspace{2.7cm}(7.3.49)\end{eqnarray*} by (7.3.35) and
(7.3.36). Moreover,
$$(U(\msr G_-)(1\otimes
M))_0=\widehat M_0=1\otimes M.\eqno(7.3.50)$$ Furthermore,
 $$\widehat M=\bigoplus\limits_{k=0}^\infty\widehat M_k,\qquad
 U(\msr G_-)(1\otimes M)=\bigoplus\limits_{k=0}^\infty(U(\msr G_-)(1\otimes
M))_k.\eqno(7.3.51)$$

According to (7.3.26), we define a $\msr G_0$-module homomorphism
$\vf:\wht M\rta  U(\msr G_-)(1\otimes M)$ by
$$\vf(x^\al\otimes v)=E^\al(1\otimes
v)\qquad\for\;\;\al\in\mbb{N}^{2n},\;v\in M.\eqno(7.3.52)$$ Then we
have
$$\vf(\wht M_k)=(U(\msr G_-)(1\otimes
M))_k\qquad\for\;\;k\in\mbb{N}.\eqno(7.3.53)$$ Under the
identification (7.3.28), we have:\psp

{\bf Lemma 7.3.4}. {\it We have $\vf|_{\wht
M_1}=(c+\td\omega)|_{\wht M_1}$ (cf. (7.3.8)-(7.3.10)).}\psp

{\it Proof}.  Let $i\in \ol{1,n}$ and $v\in M$. Expressions (7.3.35)
and (7.3.36) give
$$\vf(x_i\otimes v)=-A_{i,n+1}(1\otimes v)
=\sum_{p=1}^n(x_p\otimes A_{i,p}(v)+x_{n+p}\otimes
B_{i,p}(v))+cx_i\otimes v,\eqno(7.3.54)$$ $$\vf(x_{n+i}\otimes
v)=C_{n+1,i}(1\otimes v)=\sum_{p=1}^n(x_p\otimes
C_{i,p}(v)-x_{n+p}\otimes A_{p,i}(v))+cx_{n+i}\otimes
v.\eqno(7.3.55)$$

On the other hand, (7.3.10), (7.3.28), (7.3.32) and (7.3.33) yield
\begin{eqnarray*}\qquad\td\omega(x_i\otimes v)&=&
\sum_{1\leq p<q\leq n}[(B_{p,q}\otimes C_{q,p})(x_i\otimes v)
+(C_{q,p}\otimes B_{p,q})(x_i\otimes v)]
\\ & &+\sum_{r,s=1}^n(A_{r,s}\otimes A_{s,r})(x_i\otimes v)
\\ &=&\sum_{q=1}^nx_{n+q}\otimes B_{i,q}(v)+\sum_{r=1}^nx_r\otimes
A_{i,r}(v), \hspace{4.4cm}(7.3.56)\end{eqnarray*}
\begin{eqnarray*}\qquad\td\omega(x_{n+i}\otimes v)&=&
\sum_{1\leq p<q\leq n}[(B_{p,q}\otimes C_{q,p})(x_{n+i}\otimes v)
+(C_{q,p}\otimes B_{p,q})(x_{n+i}\otimes v)]
\\ & &+\sum_{r,s=1}^n(A_{r,s}\otimes A_{s,r})(x_{n+i}\otimes v)
\\ &=&\sum_{p=1}^nx_p\otimes C_{i,p}(v)-\sum_{s=1}^nx_{n+s}\otimes
A_{s,i}(v).\hspace{4cm}(7.3.57)\end{eqnarray*} Comparing the above
four expressions, we get the conclusion in the lemma.
$\qquad\Box$\psp

For $f\in{\msr A}$, we define the action
$$f(g\otimes v)=fg\otimes v\qquad\for\;\;g\in{\msr A},\;v\in
M.\eqno(7.3.58)$$ Then we have the $o(2n,\mbb{C})$-invariant
operator
$$T=\sum_{i=1}^n[\nu(C_{n+1,i})x_i-\nu(A_{i,n+1})x_{n+i}]\eqno(7.3.59)$$
on $\wht M$.

{\bf Lemma 7.3.5}. {\it We have $T|_{\wht
M_k}=(2c+2-2n+k)\eta|_{\wht M_k}$}.

{\it Proof}. Let $f$ be a homogeneous polynomial with degree $k$ and
let $v\in M$. According to (7.3.22), we find
\begin{eqnarray*}\qquad& &\sum_{i=1}^n[(x_iD-\eta\partial_{x_{n+i}})(x_{n+i}f)+
(x_{n+i}D-\eta\partial_{x_i})(x_if)]
\\ &=&2(k+1)\eta f-2n\eta f-\eta D(f)=(k+2-2n)\eta f.\hspace{4.4cm}(7.3.60)
\end{eqnarray*}
Moreover, the skew-symmetry of $C_{i,p}$ and $B_{p,i}$ implies
$$\sum_{i=1}^n[\sum_{p=1}^n(x_pC_{i,p}-x_{n+p}A_{p,i})+x_{n+i}I_{2n}]x_i
=-\sum_{i,p=1}^nx_ix_{n+p}A_{p,i}+\eta I_{2n},\eqno(7.3.61)$$
$$-\sum_{i=1}^n[\sum_{p=1}^n(x_{n+p}B_{p,i}-
x_pA_{i,p})-x_iI_{2n}]x_{n+i}=\sum_{i,p=1}^nx_{n+i}x_pA_{i,p}+\eta
I_{2n}.\eqno(7.3.62)$$ Now the lemma follows from (7.3.35) and
(7.3.36). $\qquad\Box$ \psp

 For $ 0\neq \mu=\sum_{i=1}^n\mu_i\ves_i\in\Lmd^+$ with ${\msr
S}(\mu)=\{i_1,i_2,...,i_{s+1}\}$ (cf. (6.7.54) and (6.7.55)), we
define
$$\Theta(\mu)=\left\{\begin{array}{ll}
\mu_1+n-1-\mbb{N}&\mbox{if}\;\mu_{n-1}=-\mu_n>0\;\mbox{and}\;s=2,
\\\mu_1+2n-i_2-\mbb{N}&\mbox{otherwise}.\end{array}\right.\eqno(7.3.63)$$

{\bf Theorem 7.3.6}.  {\it For $0\neq\mu\in\Lmd^+$, the
$o(2n+2,\mbb{C})$-module $\widehat{V(\mu)}$ defined by
(7.3.32)-(7.3.41) is irreducible if $c\in
\mbb{C}\setminus\{n-1-\mbb{N}/2,\Theta(\mu)\}.$}

{\it Proof}. By Proposition 7.3.3, it is enough to prove that the
homomorphism $\vf$ defined in (7.3.47) and (7.3.52) satisfies
$\vf(\wht{V(\mu)})=\wht{V(\mu)}$.  According to (7.3.53), we only
need to prove
$$\vf(\wht{V(\mu)}_k)=\wht{V(\mu)}_k\eqno(7.3.64)$$
for any $k\in\mbb{N}$. We will prove it by induction on $k$.

When $k=0$, (7.3.64) holds by (7.3.50) and (7.3.52). Consider $k=1$.
Write $\mu=\sum_{i=1}^n\mu_i\ves_i\in\Lmd^+$ with ${\msr
S}(\mu)=\{i_1,i_2,...,i_{s+1}\}$. According to Lemma 7.3.2 and Lemma
7.3.4 with $M=V(\mu)$, the eigenvalues of $\vf|_{\wht{V(\mu)}_1}$
are
$$c+\mu_{i_r}-i_r+1,\;c-\mu_{i_r}-2n+i_{r+1}\;\;\for\;\;r\in\ol{1,s}\;\;\mbox{if}\;\;\mu_{n-1}+\mu_n>0\eqno(7.3.65)$$
and
$$c+\mu_{i_r}-i_r+1,\;c-\mu_{i_\iota}-2n+i_{\iota+1}\;\;\for\;\;r\in\ol{1,s},\;\iota\in\ol{1,s-2+\dlt_{\mu_n,0}}\eqno(7.3.66)$$
when $\mu_{n-1}+\mu_n=0.$ Recall that $\mu_i\in\mbb{Z}/2$ for
$i\in\ol{1,n}$,
$$\mu_\iota-\mu_{\iota+1}\in\mbb{N}\;\;\for\;\;\iota\in\ol{1,n-1},\qquad
\;\;\mu_{n-1}\geq |\mu_n|. \eqno(7.3.67)$$
 So $c\not\in\Theta(\mu)$ implies that the eigenvalues of
 $\vf|_{\wht{V(\mu)}_1}$ are nonzero.
 Thus (7.3.64) holds for
$k=1$.

Suppose that (7.3.64) holds for $k\leq \ell$ with $\ell\geq 1$.
Consider $k=\ell+1$. Note that
\begin{eqnarray*}\qquad\vf(\widehat{V(\mu)}_{
\ell+1})&=&\sum_{i=1}^{2n}\vf(x_i\widehat{V(\mu)}_\ell)
=\sum_{i=1}^n[A_{i,n+1}[\vf(\widehat{V(\mu)}_
\ell)]+C_{n+1,i}[\vf(\widehat{V(\mu)}_ \ell)]]\\ &=&
\sum_{i=1}^n[A_{i,n+1}(\widehat{V(\mu)}_
\ell)+C_{n+1,i}(\widehat{V(\mu)}_\ell)]\hspace{4.4cm}(7.3.68)\end{eqnarray*}
by the inductional assumption.  To prove (7.3.64) with $k=\ell+1$ is
equivalent to prove
$$\sum_{i=1}^n[A_{i,n+1}(\widehat{V(\mu)}_\ell)+C_{n+1,i}(\widehat{V(\mu)}_\ell)]=\widehat{V(\mu)}_{
\ell+1}.\eqno(7.3.69)$$

For any $u\in \widehat{V(\mu)}_{ \ell-1}$, Lemma 7.3.5 says that
$$\sum_{i=1}^n[C_{n+1,i}(x_iu)-A_{i,n+1}(x_{n+i}u)]=(2c+2-2n+\ell)\eta
u.\eqno(7.3.70)$$ Since $c\not\in n-1-\mbb{N}/2$, $2c+2-2n+\ell\neq
0$ and  (7.3.70) gives
$$\eta u\in \sum_{i=1}^n[A_{i,n+1}(\widehat{V(\mu)}_\ell)+C_{n+1,i}(\widehat{V(\mu)}_\ell)]\qquad\for\;\;u\in  \widehat{V(\mu)}_{
\ell-1}.\eqno(7.3.71)$$

Let $g\otimes v\in \widehat{V(\mu)}_{\la \ell\ra}$. According to
(7.3.35), (7.3.36) and Lemma 7.3.4,
\begin{eqnarray*} -A_{i,n+1}(g\otimes v)&=&x_iD(g)\otimes
v-\eta\ptl_{x_{n+i}}(g)\otimes v+x_ig\otimes I_{2n}(v)\\
& &+\sum_{p=1}^n(x_{n+p}g\otimes
B_{i,n+p}(v)+x_pg\otimes A_{i,p}(v))\\
&=&\ell x_ig\otimes
v-\eta\ptl_{x_{n+i}}(g)\otimes v+g[x_i\otimes I_{2n}(v)\\
& &+\sum_{p=1}^n(x_{n+p}\otimes B_{i,n+p}(v)+x_p\otimes A_{i,p}(v))]
\\&=&\ell x_ig\otimes
v-\eta\ptl_{x_{n+i}}(g)\otimes v-g[A_{i,n+1}(1\otimes v)]
\\&=&\ell x_ig\otimes
v-\eta\ptl_{x_{n+i}}(g)\otimes v+g\vf(x_i\otimes v)
\\&=&-\eta\ptl_{x_{n+i}}(g)\otimes v +g[(\ell+c+\td\omega)(x_i\otimes
v)]\hspace{4cm}(7.3.72)\end{eqnarray*} and
\begin{eqnarray*} C_{n+i,i}(g\otimes v)&=&x_{n+i}D(g)\otimes
v-\eta\ptl_{x_i}(g)\otimes v+x_{n+i}g\otimes I_{2n}(v)\\
& &+\sum_{p=1}^n(x_pg\otimes C_{i,p}(v)-x_{n+p}g\otimes
A_{p,i}(v))\\&=&\ell x_{n+i}g\otimes
v-\eta\ptl_{x_i}(g)\otimes v+g[x_{n+i}\otimes I_{2n}(v)\\
& &+\sum_{p=1}^n(x_p\otimes C_{i,p}(v)-x_{n+p}\otimes A_{p,i}(v))]
\\&=&\ell x_{n+i}g\otimes v-\eta\ptl_{x_i}(g)\otimes
v+g[C_{n+1,i}(1\otimes v)] \\&=&\ell x_{n+i}g\otimes
v-\eta\ptl_{x_i}(g)\otimes v+g\vf(x_{n+i}\otimes v)
\\&=&-\eta\ptl_{x_i}(g)\otimes v +g[(\ell+c+\td\omega)(x_{n+i}\otimes
v)]\hspace{4.2cm}(7.3.73)\end{eqnarray*} for $i\in\ol{1,n}$. Since
(7.3.71) says that
$$\eta\ptl_{x_r}(g)\otimes v
\in\sum_{i=1}^n[A_{i,n+1}(\widehat{V(\mu)}_\ell)+C_{n+1,i}(\widehat{V(\mu)}_\ell)])\eqno(7.3.74)$$
for $r\in\ol{1,2n}$, Expressions (7.3.72) and (7.3.73) show
$$g[(\ell+c+\td\omega)(x_r\otimes
v)]\in\sum_{i=1}^n[A_{i,n+1}(\widehat{V(\mu)}_\ell)+C_{n+1,i}(\widehat{V(\mu)}_\ell)]\eqno(7.3.75)$$
for $r\in\ol{1,2n}$ and $g\in{\msr A}_\ell$.

According to Lemmas 7.3.2, the eigenvalues of
$(\ell+c+\td\omega)|_{\widehat{V(\mu)}_{\la 1\ra}}$ are among
$$\ell+c+\mu_{i_r}-i_r+1,\;\ell+c-\mu_{i_r}-2n+i_{r+1}\;\;\for\;\;r\in\ol{1,s}.\eqno(7.3.76)$$
Again $$-\ell-\mu_{i_r}+i_r-1,\;-\ell+\mu_{i_r}+2n-i_{r+1}\in
\mu_1+2n-i_2-\mbb{N}\;\;\for\;\;i\in\ol{1,s}.\eqno(7.3.77)$$ If
$c\not\in \mu_1+2n-i_2-\mbb{N}$, then all the eigenvalues of
$(\ell+c+\td\omega)|_{\widehat{V(\mu)}_1}$ are nonzero.  In the case
$\mu_n=-\mu_{n-1}>0$ and $s=2$, $\mu_1=\mu_{n-1}=-\mu_n$ and the
eigenvalues $(\ell+c+\td\omega)|_{\widehat{V(\mu)}_1}$ are
$c+\mu_1+\ell$ and $c+\mu_n-n+1+\ell=c-\mu_1-n+1+\ell$, which are
not equal to 0 because of $c\not\in \Theta(\mu)=\mu_1+n-1-\mbb{N}$.
 Hence
$$(\ell+c+\td\omega)(\widehat{V(\mu)}_1)=\widehat{V(\mu)}_1.\eqno(7.3.78)$$ By (7.3.75) and (7.3.78),
$$g(\widehat{V(\mu)}_{\la 1\ra})\subset \sum_{i=1}^n[A_{i,n+1}(\widehat{V(\mu)}_\ell)+C_{n+1,i}(\widehat{V(\mu)}_\ell)]
\qquad\for\;\;g\in{\msr A}_\ell,\eqno(7.3.79)$$ or equivalently,
(7.3.64) holds for $k=\ell+1$ by (7.3.68). By induction, (7.3.64)
holds for any $k\in\mbb{N}.\qquad\Box$\psp

Let $\lmd_i$ be the $i$th fundamental weight of $o(2n+2,\mbb C)$
with respect to the following simple positive roots
$$\{\ves_{n+1}-\ves_n,\ves_n-\ves_{n-1},...,\ves_2-\ves_1,\ves_2+\ves_1
\}.\eqno(7.3.80)$$ By the Weyl group's actions on $o(2n,\mbb C)$ and
$V(\mu)$, $\wht{V(\mu)}$ is a highest-weight $o(2n+2,\mbb C)$-module
with highest weight
$-(c+\mu_1)\lmd_1+\sum_{r=1}^{s-1}(\mu_{i_r}-\mu_{i_{r+1}})\lmd_{i_{r+1}}+(\mu_{n-1}+\mu_n)\lmd_{n+1}$.

 Up to this stage, we do not known if the
condition in Theorem 7.3.6 is necessary for the
$o(2n+2,\mbb{C})$-module $\widehat{V(\mu)}$ to be irreducible if
$\mu\neq 0$. We will deal with  the case $\mu=0$ in next section.

\section{ Conformal
Oscillator Representations}

In this section, we  study  the $o(2n+2,\mbb{C})$-module
$\widehat{V(0)}$ and its  oscillator generalizations.

In order to use the results in Sections 7.1 and 7.2, we redenote
$$y_i=x_{n+i}\qquad\for\;\;i\in\ol{1,n}\eqno(7.4.1)$$
and use ${\msr B}=\mbb{C}[x_1,...,x_n,y_1,...,y_n]$ to replace
${\msr A}$ in last section. Fix $c\in\mbb{C}$ and identify
$\widehat{V(0)}={\msr B}\otimes v_0$ with ${\msr B}$ by
$$f\otimes v_0\leftrightarrow f\qquad\for\;\;f\in{\msr
B},\eqno(7.4.2)$$ where $V(0)=\mbb{C}v_0$. Then we have the
following one-parameter generalization $\pi_c$ of the conformal
representation of $o(2n+2,\mbb{C})$:
$$\pi_c(A_{i,j})=x_i\ptl_{x_j}-y_j\ptl_{x_i},\;\pi_c(B_{i,j})=x_i\ptl_{y_j}-x_j\ptl_{y_i},\;
\pi_c(C_{i,j})=y_i\ptl_{x_j}-y_j\ptl_{x_i},\eqno(7.4.3)$$
$$\pi_c(A_{n+1,i})=\ptl_{x_i},\;\;\pi_c(A_{n+1,n+1})=-D-c,\;\;\pi_c(B_{i,n+1})=-\ptl_{y_i},\eqno(7.4.4)$$
$$\pi_c(A_{i,n+1})=\eta\ptl_{y_i}-x_i(D+c),\;\;\pi_c(C_{n+1,i})=y_i(D+c)-\eta\ptl_{x_i}\eqno(7.4.5)$$
The following result is taken from Zhao and the author's work
[XZ].\psp

 {\bf Theorem 7.4.1}. {\it The representation $\pi_c$ of
$o(2n+2,\mbb{C})$ on $\msr B$ is irreducible if and only if
$c\not\in-\mbb{N}$. When $c=0$, the  representation $\pi_c$ of
$o(2n+2,\mbb{C})$ on $\msr B$ is  the natural conformal
representation of $o(2n+2,\mbb{C})$ given in (7.3.23)-(7.3.26) in
terms (7.4.1) with $\msr A$ replaced by $\msr B$. The subspace
$\mbb{C}1_{\msr B}$ forms a trivial $o(2n+2,\mbb{C})$-submodule of
the conformal module ${\msr B}$ and the quotient space ${\msr
B}/\mbb{C}1_{\msr B}$ forms an irreducible
$o(2n+2,\mbb{C})$-module.}

{\it Proof}. Recall
$$\Dlt=\sum_{i=1}^n\ptl_{x_i}\ptl_{y_i},\;\;D=\sum_{r=1}^nx_r\partial_{x_r}+\sum_{s=1}^ny_s\ptl_{y_s},\;\;\eta=\sum_{i=1}^nx_iy_i\eqno(7.4.6)$$
in terms of (7.4.1) and the identification (7.3.28) .  Moreover,
$\msr B_k$ denotes the subspace of homogeneous polynomials in $\msr
B$ with degree $k$ and
$${\msr H}_k=\{f\in{\msr
A}_k\mid\Dlt(f)=0\}\qquad\for\;\;k\in\mbb{N}.\eqno(7.4.7)$$
According to Theorem 7.1.1,  ${\msr H}_k$ is an irreducible
$o(2n,\mbb{C})$-submodule with the highest-weight vector $x_1^k$ of
weight $k\ves_1$ for any $k\in\mbb N$, and
$$\msr B=\bigoplus_{m,k=0}^\infty \eta^m{\msr H}_k\eqno(7.4.8)$$ is a direct sum of irreducible
$o(2n,\mbb{C})$-submodules. On the other hand, $U(\msr G_-)(1_{\msr
B})$ forms an irreducible $o(2n+2,\mbb{C})$-submodule of $\msr B$
due to the first and third equations in (7.4.4). If $c=-\ell$ for
some $\ell\in\mbb N$, then
$$U(\msr G_-)(1_{\msr B})\bigcap \msr B_{\ell+1}\subset \eta\msr
B_{\ell-1}\eqno(7.4.9)$$ by (7.4.5). So $c\not\in-\mbb{N}$ is a
necessary condition for the representation $\pi_c$ of
$o(2n+2,\mbb{C})$ on $\msr B$ to be irreducible.

Next we assume $c\not\in -\mbb{N}-1$.  Let $M$ be a nonzero
$o(2n+2,\mbb{C})$-submodule of $\msr B$ such that
$$M\not\subset
\mbb{C}1_{\msr B}\;\;\mbox{if}\;\; c= 0.\eqno(7.4.10)$$ Repeatedly
applying
 the first and third equations in (7.4.4) if necessary, we have $1_{\msr B}\in M$. Note
$$-A_{i,n+1}(1_{\msr B})=cx_i,\;\;C_{n+1,i}(1_{\msr B})=cy_i\qquad\for\;\;i\in\ol{1,n}.\eqno(7.4.11)$$ Thus
$$\msr B_1\subset M\eqno(7.4.12)$$ if $c\neq 0$.
 When $c=0$, (7.4.12) also holds by the first and third equations in
 (7.4.4) and the fact that $\msr B_1$ is an irreducible $o(2n,\mbb{C})$-submodule. Suppose that
$$\msr B_k\subset M\qquad\for\;\;k<\ell\eqno(7.4.13)$$
with $2\leq \ell\in\mbb{N}$. According to (7.4.8),
$$\msr B_\ell=\bigoplus_{m=0}^{\llbracket \ell/2
\rrbracket}\eta^m{\msr H}_{\ell-2m}.\eqno(7.4.14)$$
 Moreover,
$$[\Dlt,\eta]=n+D.\eqno(7.4.15)$$Set
$$\msr B_ {\ell,r}=\bigoplus_{m=0}^r\eta^m{\msr H}_{\ell-2m}\eqno(7.4.16)$$
for $r\in\ol{0,\llbracket \ell/2 \rrbracket}$.  Then
$$\msr B_ {\ell,r}=\{w\in \msr B_ \ell\mid \Dlt^{r+1}(w)=0\}\eqno(7.4.17)$$ and
$$\Dlt^r(\msr B_ {\ell,r})={\msr H}_{\ell-2r}\eqno(7.4.18)$$ by (7.4.16).

Since
$$-A_{1,n+1}(x_1^{\ell-1})=(c+\ell-1)x_1^\ell\in
M\eqno(7.4.19)$$ and $c\not\in-\mbb{N}-1$, we have
$$x_1^\ell\in M\eqno(7.4.20)$$ by (7.4.5).
Hence
$${\msr H}_\ell\subset
M\eqno(7.4.21)$$ because ${\msr H}_\ell$ is an irreducible
$o(2n,\mbb{C})$-submodule generated by $x_1^\ell$. Recall $n\geq 2$
by our assumption. For $r\in\ol{1,\llbracket \ell/2 \rrbracket}$,
$$-A_{2,n+1}(x_1^{\ell-r-1}y_1^r)=(c+\ell-1)x_1^{\ell-r-1}x_2y_1^r
\in M\eqno(7.4.22)$$  by (7.4.5).  So $x_1^{\ell-r-1}x_2y_1^r\in M$.
Moreover,
$$\Dlt^r(x_1^{\ell-r-1}x_2y_1^r)=r![\prod_{s=1}^r(\ell-r-s)]x_1^{\ell-2r-1}x_2 \in
{\msr H}_{\ell-2r}.\eqno(7.4.23)$$ Observe that (7.4.14) is a direct
sum of irreducible $o(2n,\mbb{C})$-submodules with distinct highest
weights. So
$$\msr B_\ell\bigcap M=\bigoplus_{m=0}^{\llbracket \ell/2
\rrbracket}(\eta^m{\msr H}_{\ell-2m})\bigcap M.\eqno(7.4.24)$$ By
(7.4.17)-(7.4.23), $(\eta^r{\msr H}_{\ell-2r})\bigcap M$ is a
nonzero $o(2n,\mbb{C})$-submodule. Since $\eta^r{\msr H}_{\ell-2r}$
is an irreducible $o(2n,\mbb{C})$-module, we have $\eta^r{\msr
H}_{\ell-2r}=(\eta^r{\msr H}_{\ell-2r})\bigcap M$. Therefore, $\msr
B_\ell\subset M$. By induction, $\msr B_k\subset M$ for any
$k\in\mbb{N}$, or equivalently, $M=\msr B$. This proves the theorem.
$\qquad\Box$\psp

With respect to (7.3.80), $\msr B$ is a highest-weight irreducible
$o(2n+2,\mbb C)$-module with highest weight $-c\lmd_1$ when
$c\not\in -\mbb N$. If $c=0$, then  $\msr B/\mbb C1_{\msr B}$ is a
highest-weight irreducible $o(2n+2,\mbb C)$-module with highest
weight  $-2\lmd_1+\lmd_2$. The rest of this section is taken from
our work [X26].

For $\vec a=(a_1,a_2,...,a_n)^t,\;\vec
b=(b_1,b_2,...,b_n)^t\in\mbb{C}^n$, we put
$$\vec a\cdot\vec
x=\sum_{i=1}^na_ix_i,\qquad\vec b\cdot\vec
y=\sum_{i=1}^nb_iy_i.\eqno(7.4.25)$$ Set
$${\msr B}_{\vec a,\vec b}=\{fe^{\vec a\cdot\vec
x+\vec b\cdot\vec y}\mid f\in{\msr B}\}.\eqno(7.4.26)$$
 Denote by $\pi_{c,\vec a,\vec b}$ the
representation $\pi_c$ of $o(2n+2,\mbb C)$ on $\msr B_{\vec a,\vec
b}$.\psp

{\bf Theorem 7.4.2}. {\it The representation $\pi_{c,\vec a,\vec b}$
of $o(2n+2,\mbb{C})$ is irreducible for any $c\in\mbb{C}$ if
$\sum_{i=1}^na_ib_i\neq 0$.}

{\it Proof}. By symmetry, we may assume $a_1\neq 0$. Let ${M}$ be a
nonzero $o(2n+2,\mbb{C})$-submodule of ${\msr B}_{\vec a,\vec b}$.
Take any $0\neq fe^{\vec a\cdot\vec x+\vec b\cdot\vec y}\in  M$ with
$f\in \msr{B}$.   Set
$$\msr B_{\vec a,\vec b,k}=\msr{B}_ke^{\vec a\cdot\vec
x+\vec b\cdot\vec y}\qquad\for\;k\in\mbb{N}.\eqno(7.4.27)$$
 According to (7.4.4),
$$(A_{n+1,i}-a_i)(fe^{\vec a\cdot\vec x+\vec b\cdot\vec
y})=\ptl_{x_i}(f)e^{\vec a\cdot\vec x+\vec b\cdot\vec
y},\;\;-(B_{i,n+1}+b_i)(fe^{\vec a\cdot\vec x+\vec b\cdot\vec
y})=\ptl_{y_i}(f)e^{\vec a\cdot\vec x+\vec b\cdot\vec
y}\eqno(7.4.28)$$ for $i\in\ol{1,n}$. Repeatedly applying (7.4.28),
we obtain $e^{\vec a\cdot\vec x+\vec b\cdot\vec y}\in  M$.
Equivalently, $\msr{B}_{\vec a,\vec b,0}\subset M$.

Suppose $\msr{B}_{\vec a,\vec b,\ell}\subset M$ for some
$\ell\in\mbb{N}$.  Take any $ge^{\vec a\cdot\vec x+\vec b\cdot\vec
y}\in \msr{B}_{\vec a,\vec b,\ell}$. Since
$$(x_i\ptl_{x_1}-y_1\ptl_{y_i})(g)e^{\vec a\cdot\vec x+\vec b\cdot\vec
y},(y_i\ptl_{x_1}-y_1\ptl_{x_i})(g)e^{\vec a\cdot\vec x+\vec
b\cdot\vec y}\in \msr{B}_{\vec a,\vec b,\ell}\subset
M,\eqno(7.4.29)$$ we have
$$A_{i,1}(ge^{\vec a\cdot\vec x+\vec b\cdot\vec y})\equiv (a_1x_i-b_iy_1)ge^{\vec a\cdot\vec x+\vec b\cdot\vec y}
\equiv 0\;\;(\mbox{mod}\;M)\eqno(7.4.30)$$ and
$$C_{i,1}(ge^{\vec a\cdot\vec x+\vec b\cdot\vec y})\equiv
(a_1y_i-a_iy_1)ge^{\vec a\cdot\vec x+\vec b\cdot\vec y} \equiv
0\;\;(\mbox{mod}\;M)\eqno(7.4.31)$$ for $i\in\ol{1,n}$ by (7.4.3).
On the other hand,
$$(D+c)(g)e^{\vec a\cdot\vec x+\vec
b\cdot\vec y}\in \msr{B}_{\vec a,\vec b,\ell}\subset
M,\eqno(7.4.32)$$ and so the second equation in (7.4.4) implies
gives
$$-A_{n+1,n+1}(ge^{\vec a\cdot\vec x+\vec b\cdot\vec y})\equiv
[\sum_{i=1}^n(a_ix_i+b_iy_i)]ge^{\vec a\cdot\vec x+\vec b\cdot\vec
y} \equiv 0\;\;(\mbox{mod}\;M)\eqno(7.4.33)$$

Substituting (7.4.30) and (7.4.31) into (7.4.33), we
get$$(\sum_{i=1}^na_ib_i)y_1ge^{\vec a\cdot\vec x+\vec b\cdot\vec y}
\equiv 0\;\;(\mbox{mod}\;M).\eqno(7.4.34)$$ Equivalently,
$y_1ge^{\vec a\cdot\vec x+\vec b\cdot\vec y}\in M$. Substituting it
to (7.4.30) and (7.4.31), we obtain
$$x_ige^{\vec a\cdot\vec x+\vec
b\cdot\vec y},y_ige^{\vec a\cdot\vec x+\vec b\cdot\vec y}\in
M\eqno(7.4.35)$$ for $i\in\ol{1,n}$. Therefore, $\msr{B}_{\vec
a,\vec b,\ell+1}\subset M$. By induction, $\msr{B}_{\vec a,\vec b,
\ell}\subset M$ for any $\ell\in\mbb{N}$. So $\msr{B}_{\vec a,\vec
b}= M$. Hence $\msr{B}_{\vec a,\vec b}$ is an irreducible
$o(2n+2,\mbb{C})$-module. $\qquad\Box$\psp

In the above theorem, $\msr B_{\vec a,\vec b}$ is not a weight
$o(2n+2,\mbb{C})$-module.

Fix $n_1,n_2\in\ol{1,n}$ with $n_1\leq n_2$. We take (7.2.6),
(7.2.8) and write
$$\td\eta=\sum_{i=1}^{n_1}y_i\ptl_{x_i}+\sum_{r=n_1+1}^{n_2}x_ry_r+\sum_{s=n_2+1}^n
x_s\ptl_{y_s}.\eqno(7.4.36)$$ Moreover, we denote
$$\td c=c+n_2-n_1-n.\eqno(7.4.37)$$
 Changing operators $\ptl_{x_r}\mapsto -x_r,\;
 x_r\mapsto
\ptl_{x_r}$  for $r\in\ol{1,n_1}$ and $\ptl_{y_s}\mapsto -y_s,\;
 y_s\mapsto\ptl_{y_s}$  for $s\in\ol{n_2+1,n}$ in the above
 representation of $o(2n+2,\mbb{C})$ (7.4.3)-(7.4.5),
we have the following representation $\pi_c^{n_1,n_2}$ of the Lie
algebra $o(2n+2,\mbb{C})$  determined by
$$\pi_c^{n_1,n_2}(A_{i,j})=E_{i,j}^x-E_{j,i}^y\eqno(7.4.38)$$ with
$$E_{i,j}^x=\left\{\begin{array}{ll}-x_j\ptl_{x_i}-\delta_{i,j}&\mbox{if}\;
i,j\in\ol{1,n_1};\\ \ptl_{x_i}\ptl_{x_j}&\mbox{if}\;i\in\ol{1,n_1},\;j\in\ol{n_1+1,n};\\
-x_ix_j &\mbox{if}\;i\in\ol{n_1+1,n},\;j\in\ol{1,n_1};\\
x_i\partial_{x_j}&\mbox{if}\;i,j\in\ol{n_1+1,n};
\end{array}\right.\eqno(7.4.39)$$
and
$$E_{i,j}^y=\left\{\begin{array}{ll}y_i\ptl_{y_j}&\mbox{if}\;
i,j\in\ol{1,n_2};\\ -y_iy_j&\mbox{if}\;i\in\ol{1,n_2},\;j\in\ol{n_2+1,n};\\
\ptl_{y_i}\ptl_{y_j} &\mbox{if}\;i\in\ol{n_2+1,n},\;j\in\ol{1,n_2};\\
-y_j\partial_{y_i}-\delta_{i,j}&\mbox{if}\;i,j\in\ol{n_2+1,n};
\end{array}\right.\eqno(7.4.40)$$
and
$$\pi_c^{n_1,n_2}(E_{i,n+1+j})=\left\{\begin{array}{ll}
\ptl_{x_i}\ptl_{y_j}&\mbox{if}\;i\in\ol{1,n_1},\;j\in\ol{1,n_2};\\
-y_j\ptl_{x_i}&\mbox{if}\;i\in\ol{1,n_1},\;j\in\ol{n_2+1,n};\\
x_i\ptl_{y_j}&\mbox{if}\;i\in\ol{n_1+1,n},\;j\in\ol{1,n_2};\\
-x_iy_j&\mbox{if}\;i\in\ol{n_1+1,n},\;j\in\ol{n_2+1,n};\end{array}\right.\eqno(7.4.41)$$
$$\pi_c^{n_1,n_2}(E_{n+1+i,j})=\left\{\begin{array}{ll}
-x_jy_i&\mbox{if}\;j\in\ol{1,n_1},\;i\in\ol{1,n_2};\\
-x_j\ptl_{y_i}&\mbox{if}\;j\in\ol{1,n_1},\;i\in\ol{n_2+1,n};\\
y_i\ptl_{x_j}&\mbox{if}\;j\in\ol{n_1+1,n},\;i\in\ol{1,n_2};\\
\ptl_{x_j}\ptl_{y_i}&\mbox{if}\;j\in\ol{n_1+1,n},\;i\in\ol{n_2+1,n};\end{array}\right.\eqno(7.4.42)$$
$$\pi_c^{n_1,n_2}(A_{n+1,n+1})=-\wht D-\td c;\eqno(7.4.43)$$
$$\pi_c^{n_1,n_2}(A_{n+1,i})=\left\{\begin{array}{ll}-x_i&\mbox{if}\;\;i\in\ol{1, n_1};\\
\ptl_{x_i}&\mbox{if}\;\;i\in\ol{n_1+1,n};\end{array}\right.\eqno(7.4.44)$$
$$\pi_c^{n_1,n_2}(B_{i,n+1})=\left\{\begin{array}{ll}-\ptl_{y_i}&\mbox{if}\;\;\in\ol{1, n_2};\\
y_i&\mbox{if}\;\;i\in\ol{n_2+1,n};\end{array}\right. \eqno(7.4.45)$$
$$\pi_c^{n_1,n_2}(A_{i,n+1})=\left\{\begin{array}{ll} \td\eta\ptl_{y_i}-(\wht D+\td c-1)\ptl_{x_i}&\mbox{if}\;\;i\in\ol{1, n_1};\\
 \td\eta\ptl_{y_i}-x_i(\wht D+\td c)&\mbox{if}\;\;i\in\ol{n_1+1, n_2};\\
 - \td\eta y_i-x_i(\wht D+\td
 c)&\mbox{if}\;\;i\in\ol{n_2+1,n};\end{array}\right.\eqno(7.4.46)$$
$$\pi_c^{n_1,n_2}(C_{n+1,i})
=\left\{\begin{array}{ll} \td\eta x_i+y_i(\wht D+\td c)&\mbox{if}\;\;i\in\ol{1, n_1};\\
- \td\eta\ptl_{x_i}+y_i(\wht D+\td c)&\mbox{if}\;\;i\in\ol{n_1+1, n_2};\\
 -\td\eta\ptl_{x_i}+(\wht D+\td
 c-1)\ptl_{y_i}&\mbox{if}\;\;i\in\ol{n_2+1,n}\end{array}\right.\eqno(7.4.47)$$
for $i,j\in\ol{1,n}$.

 Recall that
 as operators on ${\msr B}$, (7.2.9) holds with $\eta$ replaced by
 $\td\eta$. We take (7.2.10) and (7.2.12). Moreover, we have
 (7.2.11).\psp

{\bf Theorem 7.4.3}. {\it The representation $\pi_c^{n_1,n_2}$ of
$o(2n+2,\mbb{C})$ on $\msr B$ is irreducible if $c\not\in
\mbb{Z}/2$.}

{\it Proof}. Let $M$ be a nonzero $o(2n+2,\mbb{C})$-submodule of
$\msr B$. By (7.2.11) and (7.4.43),
$$M=\bigoplus_{k\in\mbb{Z}}\msr B_{\la k\ra}\bigcap
M.\eqno(7.4.48)$$ Thus $\msr B_{\la k\ra}\bigcap M\neq\{0\}$ for
some $k\in \mbb{Z}$. If $k>n_1-n_2+1-\dlt_{n_1,n_2}$, then
$$
\{0\}\neq(-x_1)^{k-(n_1-n_2+1-\dlt_{n_1,n_2})}(\msr B_{\la
k\ra}\bigcap M) =A_{n+1,1}^{k-(n_1-n_2+1-\dlt_{n_1,n_2})}(\msr
A_{\la k\ra}\bigcap M)\eqno(7.4.49)$$ by (7.4.44),  which implies
$\msr B _{\la n_1-n_2+1-\dlt_{n_1,n_2} \ra}\bigcap M\neq \{0\}$.
Thus we can assume $k\leq n_1-n_2+1-\dlt_{n_1,n_2}$. Observe that
the Lie subalgebra
$$\msr L=\sum_{i,j=1}^n\mbb CA_{i,j}\cong
gl(n,\mbb{C}).\eqno(7.4.50)$$ By Theorem 7.2.1, $\msr B_{\la
k\ra}=\bigoplus_{i=0}^\infty\td\eta^i(\msr H_{\la k-2i\ra})$ is a
decomposition of irreducible $o(2n,\mbb C)$-submodules. Moreover,
the weight sets of $\msr L$-singular vectors in $\td\eta^i(\msr
H_{\la k-2i\ra})$ are distinct  by Sections 6.4-6.6. Hence
$$\td\eta^i(\msr H_{\la k-2i\ra})\subset M\;\;\mbox{for some}\;\;i\in\mbb{N}.\eqno(7.4.51)$$

Observe that
$$x_1^{-k+2i}\in \msr H_{\la k-2i\ra}.\eqno(7.4.52)$$
By (7.4.35) and (7.4.44),
$$i!(-1)^i(\prod_{r=1}^i(-k+i+r))x_1^{-k+i}=B_{1,n+1}^i(\eta^i_{n_1,n_2}(x_1^{-k+2i}))\in
M.\eqno(7.4.53)$$ Thus
$$\msr H_{\la k-i\ra}\subset M.\eqno(7.4.54)$$
So we can just assume
$$\msr H_{\la k\ra}\subset M.\eqno(7.4.55)$$
According to (7.4.43),
$$x_1^{-k+s}=(-1)^sA_{n+1,1}^s(x_1^{-k})\in
M\qquad\for\;\;s\in\mbb{N}.\eqno(7.4.56)$$ So Theorem 7.2.1 gives
$$\msr H_{\la k-s\ra}\subset M\qquad\for\;\;s\in\mbb{N}.\eqno(7.4.57)$$
For any $r\in k-\mbb{N}$, we suppose
$\td\eta^s(x_1^{-r+s}),\td\eta^s(x_1^{-r+s+1})\in M$ for some
$s\in\mbb N$. Applying (7.4.47) to it, we get
$$C_{n+1,1}[\td\eta^s(x_1^{-r+s})]=\td\eta^{s+1}(x_1^{-r+s+1})
+(r+\td c)\td\eta^s(y_1x_1^{-r+s}) \in M.\eqno(7.4.58)$$ By (7.4.36)
and (7.4.47),
$$C_{n+1,i}[\td\eta^s(x_1^{-r+s+1})]=(r-1+\td
c)\td\eta^s(y_ix_1^{-r+s+1})\in M\eqno(7.4.59)$$ for
$i\in\ol{n_1+1,n_2}$. According to  (7.4.36) and (7.4.46),
$$A_{i,n+1}[\td\eta^s(y_ix_1^{-r+s+1})]=\td\eta^{s+1}(x_1^{-r+s+1})
-(r+\td c)\td\eta^s(x_iy_ix_1^{-r+s+1})\in M\eqno(7.4.60)$$ for
$i\in\ol{n_1+1,n_2}$. Again (7.4.36), and
$(-r+s+1)\times(7.4.58)-\sum_{i=n_1+1}^{n_2}(7.4.60)$ lead to
$$ (1+s+\td
c-n_2+n_1)\td\eta^{s+1}(x_1^{-r+s+1})\in M\Rightarrow
\td\eta^{s+1}(x_1^{-r+s+1})\in M\;\;\for\;r\in k-\mbb
N.\eqno(7.4.61)$$ Replacing $r$ by $r-1$ in the above expression, we
get $\td\eta^{s+1}(x_1^{-r+s+2})\in M$. By induction,
$$\td\eta^\ell(x_1^{-r+\ell})\in M\qquad\for\;\;r\in k-\mbb N,\;\ell\in\mbb
N.\eqno(7.4.62)$$ Since $\td\eta^\ell(\msr H_{\la r-\ell\ra})\ni
\td\eta^\ell(x_1^{-r+\ell})$ is an irreducible $o(2n,\mbb C)$-module
by Theorem 7.2.1, we have
$$\td\eta^\ell(\msr H_{\la
r-\ell\ra})\subset M \qquad\for\;\;\ell\in\mbb N.\eqno(7.4.63)$$
Taking $r=m-\ell$ with $m\in k-\mbb N$, we get
$$\td\eta^\ell(\msr H_{\la
m-2\ell\ra})\subset M \qquad\for\;\;\ell\in\mbb N.\eqno(7.4.64)$$
According to Theorem 7.2.1,
$$\msr B_{\la
m\ra}=\bigoplus_{\ell=0}^\infty\td\eta^\ell(\msr H_{\la
m-2\ell\ra})\subset M\qquad\for\;\;m\in k-\mbb N.\eqno(7.4.65)$$

Expression (7.4.46) gives
$$\pi_c^{n_1,n_2}(A_{i,n+1})y_i=\left\{\begin{array}{ll} \td\eta(y_i\ptl_{y_i}+1)-y_i\ptl_{x_i}(\wht D+\td c+1)&\mbox{if}\;\;i\in\ol{1, n_1},\\
 \td\eta(y_i\ptl_{y_i}+1)-x_iy_i(\wht D+\td c+1)&\mbox{if}\;\;i\in\ol{n_1+1,n_2},\end{array}\right.
 \eqno(7.4.66)$$
$$\pi_c^{n_1,n_2}(A_{j,n+1})\ptl_{y_j}=- \td\eta y_j\ptl_{y_j}-x_j\ptl_{y_j}(\wht D+\td
 c+1)\qquad\for\;\;j\in\ol{n_2+1,n}.\eqno(7.4.67)$$
Moreover, (7.4.47) yields
$$\pi_c^{n_1,n_2}(C_{n+1,r})\ptl_{x_r}=\td\eta x_r\ptl_{x_i} +y_r\ptl_{x_r}(\wht D+\td
c+1)\qquad\for\;\;r\in\ol{1,n_1},\eqno(7.4.68)$$
\begin{eqnarray*}\qquad& &\pi_c^{n_1,n_2}(C_{n+1,s})x_s \\&=&\left\{\begin{array}{ll}
- \td\eta(x_s\ptl_{x_s}+1)+x_sy_s(\wht D+\td c+1)&\mbox{if}\;\;s\in\ol{n_1+1,n_2},\\
 -\td\eta(x_s\ptl_{x_s}+1)+x_s\ptl_{y_s}(\wht D+\td
 c+1)&\mbox{if}\;\;s\in\ol{n_2+1,n}.\end{array}\right.\hspace{2.6cm}(7.4.69)\end{eqnarray*}
Thus
\begin{eqnarray*}\hspace{3cm}& &\sum_{i=1}^{n_2}\pi_c^{n_1,n_2}(A_{i,n+1})y_i+\sum_{j=n_2+1}^n\pi_c^{n_1,n_2}(A_{j,n+1})\ptl_{y_j}
\\
&&-\sum_{r=1}^{n_1}\pi_c^{n_1,n_2}(C_{n+1,r})\ptl_{x_r}-\sum_{s=n_1+1}^n\pi_c^{n_1,n_2}(C_{n+1,s})x_s\\
&=&\td\eta(-\wht D+n_2+n-n_1-2(\td
c+1))\hspace{4.8cm}(7.4.70)\end{eqnarray*} as operators on $\msr B$.
Suppose that $\msr B_{\la \ell-s\ra}\subset M$ for some
$k\leq\ell\in\mbb{Z}$ and any $s\in\mbb{N}$. For any $f\in \msr
B_{\la \ell-1\ra}$, we apply the above equation to it and get
$$(1-\ell+n_2+n-n_1-2(\td
c+1))\td\eta(f)\in M.\eqno(7.4.71)$$ Since $c\not\in \mbb Z/2$, we
have
$$\td\eta(f)\in M.\eqno(7.4.72)$$

Now for any $g\in\msr B_{\la \ell\ra}$, we have $\ptl_{y_1}(g)\in
\msr B_{\la \ell-1\ra}$. By (7.4.46),
$$A_{1,n+1}(g)=\td\eta(\ptl_{y_1}(g))-(\ell+\td
c)\ptl_{x_1}(g)\in M.\eqno(7.4.73)$$ Moreover, (7.4.72) and (7.4.73)
yield
$$\ptl_{x_1}(g)\in M\qquad\for\;\;g\in \msr B_{\la
\ell\ra}.\eqno(7.4.74)$$ Since
$$\ptl_{x_1}(\msr B_{\la
\ell\ra})=\msr B_{\la \ell+1\ra},\eqno(7.4.75)$$ we obtain
$$\msr B_{\la \ell+1\ra}\subset M.\eqno(7.4.76)$$
By induction on $\ell$, we find
$$\msr B_{\la \ell\ra}\subset M\qquad\for\;\;\ell\in\mbb{Z},\eqno(7.4.77)$$
or equivalently, $\msr B=\bigoplus_{\ell\in\mbb{Z}}\msr B_{\la
\ell\ra}=M$. Thus $\msr B$ is an irreducible $o(2n+2,\mbb
C)$-module.$\qquad\Box$\psp

{\bf Remark 7.4.4}. The above irreducible representation depends on
the three parameters $c\in \mbb{F}$ and $n_1,n_2\in\ol{1,n}$. It is
not highest-weight type because of the mixture of multiplication
operators and differential operators in
 (7.4.39))-(7.4.42) and
(7.4.44)-(7.4.47). Since $\msr B$ is not completely reducible as an
$\msr L$-module by Sections 6.4-6.6 when $n\geq 2$ and $n_1<n$,
$\msr B$ is not a unitary $o(2n+2,\mbb{C})$-module. Expression
(7.4.43) shows that $\msr B$ is a weight $o(2n+2,\mbb{C})$-module
with finite-dimensional weight subspaces.

\chapter{Representations of Odd Orthogonal Lie Algebras}

In this chapter, we focus on  the natural explicit representations
of odd orthogonal Lie algebras. In Section 8.1, we study the
canonical bosonic and fermionic oscillator representations over
their minimal natural modules. The spin representations are also
presented. In Section 8.2, we determine the structure of the
noncanonical oscillator representations obtained from the above
bosonic representations by partially swapping differential operators
and multiplication operators, which are generalizations of the
classical theorem on harmonic polynomials. The results were due to
Luo and the author [LX2].  In Section 8.3, we construct a new
functor from the category $o(2n+1,\mbb C)$-{\bf mod} to the category
$o(2n+3,\mbb C)$-{\bf mod}, which is an extension of the conformal
representation of $o(2n+3,\mbb C)$. Moreover, we find the condition
for the functor to map a finite-dimensional irreducible $o(2n+1,\mbb
C)$-module to an infinite-dimensional irreducible $o(2n+3,\mbb
C)$-module. This work was due to Zhao and the author [XZ]. The work
in Section 8.3 gives rise to a one-parameter ($c$) family  of
inhomogeneous first-order differential operator representations of
$o(2n+3,\mbb C)$. Letting these operators act on the space of
exponential-polynomial functions that depend on  a parametric vector
$\vec a\in \mbb C^{2n+1}$, we prove in Section 8.4 that the space
forms an irreducible $o(2n+3,\mbb C)$-module for any $c\in\mbb C$ if
$\vec a$ is not on a certain hypersurface. By partially swapping
differential operators and multiplication operators, we obtain more
general differential operator representations of $o(2n+3,\mbb C)$ on
the polynomial algebra $\msr B'$ in $2n+1$ variables. We prove that
$\msr B'$ forms an infinite-dimensional irreducible weight
$o(2n+3,\mbb C)$-module with finite-dimensional weight subspaces if
$c\not\in\mbb Z/2$. These results are taken from our work [X26].

\section{Canonical Oscillator Representations}
In this section, we study the canonical bosonic and fermionic
oscillator representations of odd orthogonal Lie algebras over their
minimal natural modules and the spin representations.

Let $n\geq 1$ be an integer.  The odd orthogonal Lie algebra
$$o(2n+1,\mbb F)=o(2n,\mbb F)+\sum_{i=1}^n[\mbb F(E_{0,i}-E_{n+i,0})
+\mbb F(E_{0,n+i}-E_{i,0})]\eqno(8.1.1)$$ (cf. (7.1.1)).  We take
(7.1.2) as a Cartan subalgebra and use the settings in
(7.1.3)-(7.1.5). Then the root system of $o(2n+1,\mbb F)$ is
$$\Phi_{B_n}=\{\pm \ves_i\pm\ves_j,\pm\ves_r\mid1\leq i<j\leq
n;r\in\ol{1,n}\}.\eqno(8.1.2)$$ We take the set of positive roots
$$\Phi_{B_n}^+=\{\ves_i\pm\ves_j,\ves_r\mid1\leq i<j\leq
n,r\in\ol{1,n}\}.\eqno(8.1.3)$$ In particular,
$$\Pi_{B_n}=\{\ves_1-\ves_2,...,\ves_{n-1}-\ves_n,\ves_n\}\;\;\mbox{is the set of positive simple roots}.\eqno(8.1.4)$$

Recall the set of dominant integral weights
$$\Lmd^+=\{\mu\in L_\mbb{Q}\mid
2(\ves_n,\mu),(\ves_i-\ves_{i+1},\mu)\in\mbb{N}\;\for\;i\in\ol{1,n-1}\}.\eqno(8.1.5)$$
According to (7.1.5),
$$\Lmd^+=\{\mu=\sum_{i=1}^n\mu_i\ves_i\mid
\mu_i\in\mbb{N}/2;\mu_i-\mu_{i+1}\in\mbb{N}\;\for\;i\in\ol{1,n-1}\}.\eqno(8.1.6)$$

Let ${\msr B}'=\mbb{F}[x_0,x_1,...,x_n,y_1,...,y_n]$. Recall The
canonical oscillator representation of $o(2n+1,\mbb{F})$ is given by
$$(E_{i,j}-E_{n+j,n+i})|_{{\msr B}'}=x_i\ptl_{x_j}-y_j\ptl_{x_i},\;(E_{n+i,j}-E_{n+j,i})|_{{\msr B}'}
=y_i\ptl_{x_j}-y_j\ptl_{x_i},\eqno(8.1.7)$$
$$(E_{i,n+j}-E_{j,n+i})|_{{\msr B}'}
=x_i\ptl_{y_j}-x_j\ptl_{y_i},\;(E_{0,i}-E_{n+i,0})|_{{\msr
B}'}=x_0\ptl_{x_i}-y_i\ptl_{x_0}, \eqno(8.1.8)$$
$$(E_{0,n+i}-E_{i,0})|_{{\msr
B}'}=x_0\ptl_{y_i}-x_i\ptl_{x_0}\eqno(8.1.9)$$ for $i,j\in\ol{1,n}$.
The positive root vectors of $o(2n+1,\mbb{F})$ are
$$\{E_{i,j}-E_{n+j,n+i},E_{i,n+j}-E_{j,n+i},E_{0,n+r}-E_{r,0}\mid 1\leq i<j\leq
n;r\in\ol{1,n}\}.\eqno(8.1.10)$$

Denote
$$\Dlt'=\ptl_{x_0}^2+2\sum_{i=1}^n\ptl_{x_i}\ptl_{y_i},\;\;\eta'=x^2_0+2\sum_{i=1}^nx_iy_i.
\eqno(8.1.11)$$ It can be verified that $\eta'$  is an
$o(2n+1,\mbb{F})$-invariant and $\Dlt'$  is an
$o(2n+1,\mbb{F})$-invariant differential operator, or equivalently,
$$\xi\eta'=\eta'\xi,\;\xi\Dlt'=\Dlt'\xi'\;\;\;\mbox{on}\;\;{\msr B}'\;\;\for\;\;\xi'\in
o(2n+1,\mbb{F}).\eqno(8.1.12)$$ Let ${\msr B}'_k$ be the subspace of
homogeneous polynomials in ${\msr B}'$ with degree $k$. Set
$${\msr
H}_k'=\{f\in{\msr
B}_k'\mid\Dlt'(f)=0\}\qquad\for\;\;k\in\mbb{N}.\eqno(8.1.13)$$ \pse

{\bf Theorem 8.1.1}. {\it For $k\in\mbb{N}$, the subspace ${\msr
H}_k'$ forms a finite-dimensional irreducible
$o(2n+1,\mbb{F})$-submodule and $x_1^k$ is a highest-weight vector
with the weight $k\ves_1$.  Moreover, it has a basis
$$\left\{\sum_{i=0}^\infty (-2)^i\frac{x_0^{2i+\iota}(\sum_{i=1}^n\ptl_{x_i}\ptl_{y_i})^i(x^\al y^\be)}
{(2i+\iota)!}\mid\iota=0,1;\al,\be\in\mbb{N}^n;|\al|+|\be|+\iota=k\right\}.
\eqno(8.1.14)$$ Furthermore,
$${\msr B}'=\bigoplus_{k_1,k_2}{\eta'}^{k_1}{\msr
H}_{k_2}'\eqno(8.1.15)$$ is a direct sum of irreducible
$o(2n+1,\mbb{F})$-submodules.}

{\it Proof}. Take $\eta=\sum_{i=1}^nx_iy_i$. According to Theorem
7.1.1, the $o(2n+1,\mbb{F})$-singular vectors in ${\msr B}'$ are in
$$\sum_{i,j,k=0}^\infty \mbb{F}x_0^ix_1^j\eta^k.\eqno(8.1.16)$$
Since
$$\eta=\frac{1}{2}(\eta'-x_0^2),\eqno(8.1.17)$$
the $o(2n+1,\mbb{F})$-singular vectors in ${\msr B}'$ are in
$$\sum_{i,j,k=0}^\infty \mbb{F}x_0^ix_1^j{\eta'}^k.\eqno(8.1.18)$$

Suppose that
$$g=\sum_{i,j,k}a_{i,j,k}x_0^ix_1^j{\eta'}^k\;\;\mbox{is a singular
vector}.\eqno(8.1.19)$$ Then
$$(E_{0,n+1}-E_{1,0})(g)=-x_1\sum_{i,j,k}ia_{i,j,k}x_0^{i-1}x_1^j{\eta'}^k=0\eqno(8.1.20)$$
by the first equation in (8.1.12). So $g$ is independent of $x_0$.
Thus the $o(2n+1,\mbb{F})$-singular vectors in ${\msr B}_k'$ are
$$\{{\eta'}^{\ell}x_1^{k-2\ell}\mid\ell\in\mbb{N}\}.\eqno(8.1.21)$$
Note
$$[\Dlt',\eta']=2(1+2n)+4D',\;\;D'=x_0\ptl_{x_0}+\sum_{i=1}^n(x_i\ptl_{x_i}+y_i\ptl_{y_i}).\eqno(8.1.22)$$
Hence
$$\Dlt'[{\eta'}^{\ell}x_1^{k-2\ell}]=2\ell[1+2(n+k-1)]{\eta'}^{\ell-1}x_1^{k-2\ell}.\eqno(8.1.23)$$
Thus ${\msr H}_k'$ has a unique $o(2n+1,\mbb{F})$-singular vector
$x_1^k$. According to the second equation in (8.1.12), ${\msr H}_k'$
is a finite-dimensional $o(2n+1,\mbb{F})$-module. By Weyl's Theorem
2.3.6 of complete reducibility if $\mbb F=\mbb C$ or more generally
by Lemma 6.3.2 with $n_1=0$, it is irreducible. The equation
(8.1.15) follows from the first equation in (8.1.12) and the fact
that $o(2n+1,\mbb{F})$-singular vector are polynomials in $x_1$ and
$\eta'$.  Expression (8.1.14) follows from Lemma 6.1.1 with
$T_1=\ptl_{x_0}^2,\;T_1^-=\int_{(x_0)}^2$ and
$T_2=2\sum_{i=1}^n\ptl_{x_i}\ptl_{y_i}$. $\qquad\Box$\psp

Consider the exterior algebra $\check{\msr A}'$ generated by
$\{\sta_0,\sta_1,...,\sta_n,\vt_1,...,\vt_n\}$ (cf. (6.2.15)) and
take the settings in (6.2.44)-(6.2.51). Define a representation of
$o(2n+1,\mbb{F})$ on $\check{\msr A}'$ by
$$(E_{i,j}-E_{n+j,n+i})|_{\check{\msr
A}'}=\sta_i\ptl_{\sta_j}-\vt_j\ptl_{\vt_i},\;\;(E_{n+i,j}-E_{n+j,i})|_{\check{\msr
A}'} =\vt_i\ptl_{\sta_j}-\vt_j\ptl_{\sta_i},\eqno(8.1.24)$$
$$(E_{i,n+j}-E_{j,n+i})|_{\check{\msr
A}'} =\sta_i\ptl_{\vt_j}-\sta_j\ptl_{\vt_i},\;\;
(E_{0,i}-E_{n+i,0})|_{\check{\msr
A}'}=\sta_0\ptl_{\sta_i}-\vt_i\ptl_{\sta_0}, \eqno(8.1.25)$$
$$(E_{0,n+i}-E_{i,0})|_{\check{\msr
A}'}=\sta_0\ptl_{\vt_i}-\sta_i\ptl_{\sta_0}\eqno(8.1.26)$$ for
$i,j\in\ol{1,n}$. By (8.1.24) and (8.1.25), the above representation
of $o(2n+1,\mbb{F})$ is essentially an extension of that for
$o(2n,\mbb{F})$ given in (7.1.22) and (7.1.23). Recall the notion
$\check{\msr A}_k$ in (7.1.24). We set
$$\check{\msr A}_0'=\{1\},\;\;\check{\msr A}_{k+1}'=\sta_0{\msr
A}_k+\check{\msr A}_{k+1}\qquad\for\;\;k\in\ol{0,2k}.\eqno(8.1.27)$$

  According to (7.1.29), the $o(2n+1,\mbb{F})$-singular
vectors of $\check{\msr A}_k'$ are in the set
$$\{\sta_0\vec\sta_{k-1},\vec\sta_k\}\eqno(8.1.28)$$ if $k<n$.
Applying (8.1.26), we get that $\vec\sta_k$ is the unique
$o(2n+1,\mbb{F})$-singular vector in $\check{\msr A}_k'$. By
(7.1.29) again, the $o(2n+1,\mbb{F})$-singular vectors of
$\check{\msr A}_n'$ are in the set
$$\{\vec\sta_n,\vec\sta_{n-1}\vt_n,\sta_0\vec\sta_{n-1}\}.\eqno(8.1.29)$$
Since
$$(E_{0,2n}-E_{n,0})(\vec\sta_{n-1}\vt_n)=\vec\sta_{n-1}\sta_0,\;\;
(E_{0,2n}-E_{n,0})(\sta_0\vec\sta_{n-1})=(-1)^n\vec\sta_n\eqno(8.1.30)$$
by (8.1.26), $\vec\sta_n$ is the unique $o(2n+1,\mbb{F})$-singular
vector in $\check{\msr A}_n'$.

Next (7.1.29) implies that the $o(2n+1,\mbb{F})$-singular vectors of
$\check{\msr A}_{n+1}'$ are in the set
$$\{\sta_0\vec\sta_n,\sta_0\vec\sta_{n-1}\vt_n,\check\eta\vec\sta_{n-1}\}.\eqno(8.1.31)$$
Since
$$(E_{0,2n}-E_{n,0})(\sta_0\vec\sta_{n-1}\vt_n)=(-1)^n\vec\sta_n\vt_n,\;\;
(E_{0,2n}-E_{n,0})(\check\eta\vec\sta_{n-1})=(-1)^n\sta_0\vec\sta_n\eqno(8.1.32)$$
by (6.2.49) and (8.1.26), $\sta_0\vec\sta_n$ is the unique
$o(2n+1,\mbb{F})$-singular vector in $\check{\msr A}_{n+1}'$. For
$k\in\ol{1,n-1}$, the $o(2n+1,\mbb{F})$-singular vectors of
$\check{\msr A}_{2n+1-k}'$ are in the set
$$\{\sta_0\check\eta^{n-k}\vec\sta_k,
\check\eta^{n+1-k}\vec\sta_{k-1}\}\eqno(8.1.33)$$ by (7.1.29). Note
$$\check\eta^{n-r}\vec\sta_r=(n-r)!\sta_{r+1}\vt_{r+1}\cdots
\sta_n\vt_n\vec\sta_r.\eqno(8.1.34)$$ So
$\sta_0\check\eta^{n-k}\vec\sta_k$ is an $o(2n+1,\mbb{F})$-singular
vector by (8.1.26). Moreover,
$$(E_{0,n+k}-E_{k,0})(\check\eta^{n+1-k}\vec\sta_{k-1})=(n+1-k)\check\eta^{n-k}\vec\sta_k\sta_0\eqno(8.1.35)$$
by (6.2.49) and (8.1.26). Thus $\sta_0\check\eta^{n-k}\vec\sta_k$ is
the unique $o(2n+1,\mbb{F})$-singular vector in $\check{\msr
A}_{2n+1-k}'$. Furthermore, $\check{\msr A}_{2n+1}'$ is the
one-dimensional trivial module.
 By Weyl's Theorem 2.3.6 of
complete reducibility or more generally by an analogue of Lemma
6.3.2, we have: \psp

 {\bf Theorem 8.1.2}. {\it For $r\in\ol{1,n}$, $\check{\msr A}_r'$
 and $\check{\msr A}_{2n+1-r}'$ are finite-dimensional irreducible
 $o(2n+1,\mbb{F})$-submodules with the highest weight
 $\sum_{i=1}^r\ves_i$. Moreover, $\check{\msr A}_{2n+1}'$ and $\check{\msr A}_0'$ are isomorphic to the
one-dimensional trivial $o(2n+1,\mbb{F})$-module.}\psp

Let $\Psi'$ be the associative algebra generated by
$\{\sta_0,\sta_1,...,\sta_n\}$ with the defining relations:
$$\sta_0^2=1,\;\;\sta_0\sta_i=-\sta_i\sta_0,\;\;\sta_i\sta_j=-\sta_j\sta_i\qquad\for\;\;i,j\in\ol{1,n}.\eqno(8.1.36)$$
The algebra $\Psi'$ is an extension  of the exterior algebra $\Psi$
(cf. (6.2.15)). Take the settings as those in (6.2.16)-(6.2.18),
where $i\in\ol{1,n}$ and $j,j_s\in\ol{0,n}$. Note
$$[\sta_0\ptl_{\sta_i},\sta_0\sta_i]=2\sta_i\ptl_{\sta_i}-1\qquad\for\;\;i\in\ol{1,n}.\eqno(8.1.37)$$
Thus we have the following representation of $o(2n+1,\mbb{F})$ on
$\Psi'$:
$$(E_{i,j}-E_{n+j,n+i})|_{\Psi'}=\sta_i\ptl_{\sta_j}-\frac{\dlt_{i,j}}{2}\qquad\for\;\;i,j\in\ol{1,n}.\eqno(8.1.38)$$
$$(E_{r,n+s},E_{s,n+r})|_{\Psi'}=\sta_r\sta_s,\qquad
(E_{n+r,s},E_{n+s,r})|_{\Psi'}=\ptl_{\sta_r}\ptl_{\sta_s}\eqno(8.1.39)$$
for $1\leq r<s\leq n$, and
$$(E_{0,i}-E_{n+i,0})|_{\Psi'}=\frac{1}{\sqrt{2}}\sta_0\ptl_{\sta_i},\;\;(E_{0,n+i}-E_{i,0})|_{\Psi'}=\frac{1}{\sqrt{2}}
\sta_0\sta_i\;\;\for\;\;i\in\ol{1,n}.\eqno(8.1.40)$$
 The above representation is called the {\it spin representation}
of $o(2n+1,\mbb{F})$,\index{spin representation!of
$o(2n+1,\mbb{F})$} which is an extension of the spin representation
of$o(2n,\mbb{F})$ given (7.1.32) and (7.1.33).

Set
$$\Psi'_\pm=\Psi(1\pm\sta_0).\eqno(8.1.41)$$
Then $\Psi'_\pm$ are $o(2n+1,\mbb{F})$-submodules. Indeed, we
 have:\psp

{\bf Theorem 8.1.3}. {\it The subspaces $\Psi_\pm'$ are
$2^n$-dimensional irreducible $o(2n+1,\mbb{F})$-submodules with the
highest-weight vector $\sta_1\sta_2\cdots\sta_n(1\pm\sta_0)$ of
weight
 $(\sum_{i=1}^n\ves_i)/2$.}

\section{Noncanonical Oscillator Representations}

In this section, we study the noncanonical oscillator
representations of $o(2n+1,\mbb{F})$ obtained from the bosonic
oscillator representation
 in (8.1.7)-(8.1.9) by partially swapping differential operators and multiplication
operators.

Recall that ${\msr B}'=\mbb{F}[x_0,x_1,...,x_n,y_1,...,y_n]$ and the
canonical oscillator representations of $o(2n+1,\mbb{F})$ given in
(8.1.7)-(8.1.9). Fix
 $n_1,n_2\in\ol{1,n}$ with $n_1\leq n_2$.
  Changing operators $\ptl_{x_r}\mapsto -x_r,\;
 x_r\mapsto
\ptl_{x_r}$  for $r\in\ol{1,n_1}$ and $\ptl_{y_s}\mapsto -y_s,\;
 y_s\mapsto\ptl_{y_s}$  for $s\in\ol{n_2+1,n}$ in the canonical
oscillator representation, we get a {\it noncanonical oscillator
representation}\index{noncanonical oscillator representation!of
 $o(2n+1,\mbb{F})$} of
 $o(2n+1,\mbb{F})$ on
  ${\msr B}'$ determined by
$$(E_{i,j}-E_{n+j,n+i})|_{{\msr
B}'}=E_{i,j}^x-E_{j,i}^y\eqno(8.2.1)$$ with
$$E_{i,j}^x=\left\{\begin{array}{ll}-x_j\ptl_{x_i}-\delta_{i,j}&\mbox{if}\;
i,j\in\ol{1,n_1};\\ \ptl_{x_i}\ptl_{x_j}&\mbox{if}\;i\in\ol{1,n_1},\;j\in\ol{n_1+1,n};\\
-x_ix_j &\mbox{if}\;i\in\ol{n_1+1,n},\;j\in\ol{1,n_1};\\
x_i\partial_{x_j}&\mbox{if}\;i,j\in\ol{n_1+1,n};
\end{array}\right.\eqno(8.2.2)$$
and
$$E_{i,j}^y=\left\{\begin{array}{ll}y_i\ptl_{y_j}&\mbox{if}\;
i,j\in\ol{1,n_2};\\ -y_iy_j&\mbox{if}\;i\in\ol{1,n_2},\;j\in\ol{n_2+1,n};\\
\ptl_{y_i}\ptl_{y_j} &\mbox{if}\;i\in\ol{n_2+1,n},\;j\in\ol{1,n_2};\\
-y_j\partial_{y_i}-\delta_{i,j}&\mbox{if}\;i,j\in\ol{n_2+1,n};
\end{array}\right.\eqno(8.2.3)$$ and
$$E_{i,n+j}|_{{\msr B}'}=\left\{\begin{array}{ll}
\ptl_{x_i}\ptl_{y_j}&\mbox{if}\;i\in\ol{1,n_1},\;j\in\ol{1,n_2};\\
-y_j\ptl_{x_i}&\mbox{if}\;i\in\ol{1,n_1},\;j\in\ol{n_2+1,n};\\
x_i\ptl_{y_j}&\mbox{if}\;i\in\ol{n_1+1,n},\;j\in\ol{1,n_2};\\
-x_iy_j&\mbox{if}\;i\in\ol{n_1+1,n},\;j\in\ol{n_2+1,n};\end{array}\right.\eqno(8.2.4)$$
$$E_{n+i,j}|_{{\msr B}'}=\left\{\begin{array}{ll}
-x_jy_i&\mbox{if}\;j\in\ol{1,n_1},\;i\in\ol{1,n_2};\\
-x_j\ptl_{y_i}&\mbox{if}\;j\in\ol{1,n_1},\;i\in\ol{n_2+1,n};\\
y_i\ptl_{x_j}&\mbox{if}\;j\in\ol{n_1+1,n},\;i\in\ol{1,n_2};\\
\ptl_{x_j}\ptl_{y_i}&\mbox{if}\;j\in\ol{n_1+1,n},\;i\in\ol{n_2+1,n}.\end{array}\right.\eqno(8.2.5)$$
$$ E_{0,i}|_{{\msr
B}'}=\left\{\begin{array}{ll}-x_0x_i&\mbox{if}\;i\in\ol{1,n_1};\\
x_0\ptl_{x_i}&\mbox{if}\;i\in\ol{n_1+1,n};\\
x_0\ptl_{y_{i-n}}&\mbox{if}\;i\in\ol{n+1,n+n_2};\\
-x_0y_{i-n}&\mbox{if}\;i\in\ol{n+n_2+1,2n};\end{array}\right.\eqno(8.2.6)$$
and
$$ E_{i,0}|_{{\msr
B}'}=\left\{\begin{array}{ll}\ptl_{x_0}\ptl_{x_i}&\mbox{if}\;i\in\ol{1,n_1};\\
x_i\ptl_{x_0}&\mbox{if}\;i\in\ol{n_1+1,n};\\
y_{i-n}\ptl_{x_0}&\mbox{if}\;i\in\ol{n+1,n+n_2};\\
\ptl_{x_0}\ptl_{y_{i-n}}&\mbox{if}\;i\in\ol{n+n_2+1,2n}.\end{array}\right.\eqno(8.2.7)$$
Correspondingly we obtain a  variation of the  Laplace operator:
$$\td\Dlt'=\ptl_{x_0}^2-2\sum_{i=1}^{n_1}x_i\ptl_{y_i}+2\sum_{r=n_1+1}^{n_2}\ptl_{x_r}\ptl_{y_r}-2\sum_{s=n_2+1}^n
y_s\ptl_{x_s}\eqno(8.2.8)$$ and its dual operator
$$\eta'=x_0^2+2\sum_{i=1}^{n_1}y_i\ptl_{x_i}+2\sum_{r=n_1+1}^{n_2}x_ry_r+2\sum_{s=n_2+1}^n
x_s\ptl_{y_s}.\eqno(8.2.9)$$ Set
$${\msr B}'_{\la k\ra}=\sum_{i=0}^\infty {\msr B}_{\la
k-i\ra}x_0^i,\qquad {\msr H}'_{\la k\ra}=\{f\in {\msr B}'_{\la
k\ra}\mid \td\Dlt'(f)=0\}\eqno(8.2.10)$$ \psp (cf. (7.2.10)).
 Denote
$$\td D'=x_0\ptl_{x_0}-\sum_{i=1}^{n_1}x_i\ptl_{x_i}
+\sum_{r=n_1+1}^nx_r\ptl_{x_r}+\sum_{j=1}^{n_2}y_j\ptl_{y_j}-\sum_{s=n_2+1}^ny_s\ptl_{y_s}.\eqno(8.2.11)$$
Then
$${\msr B}'_{\la k\ra}=\{f\in {\msr B}'\mid \td
D'(f)=kf\}.\eqno(8.2.12)$$
 A straightforward verification shows
$$\td\Dlt'\xi=\xi \td\Dlt',\;\xi\eta'=\eta'\xi,\;\xi\td D'=\td D'\xi\;\mbox{on}\;{\msr
B}'\;\;\for\;\;\xi\in o(2n+1,\mbb F).\eqno(8.2.13)$$ Thus all ${\msr
B}'_{\la k\ra}$ and ${\msr H}'_{\la k\ra}$ with $k\in\mbb{Z}$ are
$o(2n+1,\mbb F)$-submodules and ${\msr
B}'=\bigoplus_{k\in\mbb{Z}}{\msr B}'_{\la k\ra}$ forms a
$\mbb{Z}$-graded algebra.

We take $\td\Dlt$ in (7.2.6) and $\eta$ in (7.2.7). Then
$$\td\Dlt'=\ptl_{x_0}^2+2\td\Dlt,\qquad\eta'=x_0^2+2\eta. \eqno(8.2.14)$$
  For $\iota=0,1$, we define
$$\wht T_\iota=\sum_{i=0}^\infty\frac{(-2)^ix_0^{2i+\iota}}{(2i+\iota)!}\td\Dlt^i.\eqno(8.2.15)$$
 By Lemma 6.1.1 with $T_1=\ptl_{x_0}^2,\;T_1^-=\int_{(x_0)}^{(2)}$
 and $T_2=2\td\Dlt$, we obtain
$${\msr H}'_{\la k\ra}=\wht T_0({\msr B}_{\la k\ra})\oplus
\wht T_1({\msr B}_{\la k-1\ra}).\eqno(8.2.16)$$\pse

{\bf Lemma 8.2.1}. {\it Suppose $n_1<n_2$. For $k\in\mbb{N}$, ${\msr
H}_{\la k\ra}'$ is an $o(2n+1,\mbb{F})$-submodule generated by
$x_{n_1+1}^k$ and ${\msr H}_{\la -k\ra}'$ is an
$o(2n+1,\mbb{F})$-submodule generated by $x_{n_1}^k$.}

{\it Proof}.  Let $V$ be an $o(2n+1,\mbb{F})$-submodule generated by
$x_{n_1+1}^k$. Since $x_{n_1+1}^k\in{\msr H}'_{\la k\ra}$ that is an
$o(2n+1,\mbb{F})$-submodule, we have $V\subset {\msr H}'_{\la
k\ra}$.

Observe
$$ (E_{n_1+1,n+j}-E_{j,n+n_1+1})|_{{\msr
B}'}=-x_{n_1+1}y_j-x_j\ptl_{y_{n_1+1}}\qquad\for\;\;j\in\ol{n_2+1,n}
\eqno(8.2.17)$$ by (8.2.4). Repeatedly applying (8.2.17) to
$x_{n_1+1}^k$ with various $j\in\ol{n_2+1,n}$, we obtain
$$x_{n_1+1}^{k+\ell}\prod_{s=n_2+1}^ny_s^{\be_s}\in V\qquad
\for\;\;\ell,\be_{n_2+1},...,\be_n\in\mbb{N},\;\sum_{s=n_2+1}^n\be_s=\ell.\eqno(8.2.18)$$
 Note
$$(E_{n+i,n+n_1+1}-E_{n_1+1,i})|_{{\msr
B}'}=y_i\ptl_{y_{n_1+1}}+x_ix_{n_1+1}\qquad\for\;\;i\in\ol{1,n_1}\eqno(8.2.19)$$
by (8.2.1)-(8.2.3). Repeatedly applying (8.2.19) to (8.2.18) with
various $i\in\ol{1,n_1}$, we have
$$[\prod_{i=1}^{n_1+1}x_i^{\al_i}][\prod_{s=n_2+1}^ny_s^{\be_s}]\in
V\eqno(8.2.20)$$ for
$\al_1,...,\al_{n_1+1},\be_{n_2+1},...,\be_n\in\mbb{N}$ such that
$k+\sum_{r=1}^{n_1}\al_r+\sum_{s=n_2+1}^n\be_s=\al_{n_1+1}$.

Denote
$$I=\{0,\ol{n_1+1,n_2},\ol{n+n_1+1,n+n_2}\}.\eqno(8.2.21)$$
Observe the Lie subalgebra
$${\msr
G}=o(2n+1,\mbb{F})\bigcap(\sum_{i,j\in I}\mbb{F}E_{i,j})\cong
o(2(n_2-n_1)+1,\mbb{F}).\eqno(8.2.22)$$ Set
$$\bar{\msr H}_{\la\ell\ra}={\msr H}'_{\la\ell\ra}\bigcap
\mbb{F}[x_0,x_{n_1+1},...,x_{n_2},y_{n_1+1},...,y_{n_2}]\qquad\for\;\;\ell\in\mbb{N}.\eqno(8.2.23)$$
Then $\bar{\msr H}_{\la\ell\ra}$ forms an irreducible ${\msr
G}$-module by Theorem 7.1.1. Repeatedly applying ${\msr G}|_{{\msr
B}'}$ to (8.2.20), we get
$$\wht T_\iota(x^\al y^\be)\in V\eqno(8.2.24)$$
for $\al,\be\in\mbb{N}^n$ such that $\al_i=0$ if $i>n_2$, $\be_j=0$
if $j\leq n_1$,  and
$$\iota+\sum_{r=n_1+1}^{n_2}\al_r-\sum_{i=1}^{n_1}\al_i+\sum_{s=n_1+1}^{n_2}\be_s-\sum_{j=n_2+1}^n\be_j=k.\eqno(8.2.25)$$
Repeatedly applying (8.2.19) to (8.2.24) satisfying the above
conditions, we get that (8.2.24) holds for $\al,\be\in\mbb{N}^n$
such that $\al_i=0$ if $i>n_2$,  and
$$\iota+\sum_{r=n_1+1}^{n_2}\al_r-\sum_{i=1}^{n_1}\al_i+\sum_{s=1}^{n_2}\be_s-\sum_{j=n_2+1}^n\be_j=k.\eqno(8.2.26)$$

According (8.2.1)-(8.2.3),
$$(E_{i,n_1+1}-E_{n+n_1+1,n+i})|_{{\msr
B}'}=x_i\ptl_{x_{n_1+1}}+y_{n_1+1}y_i\qquad\for\;\;i\in\ol{n_2+1,n}.\eqno(8.2.27)$$
Repeatedly applying (8.2.27) to (8.2.24) satisfying (8.2.26), we get
(8.2.25) for any $\al,\be\in\mbb{N}^n$ such that
$$\iota+\sum_{r=n_1+1}^n\al_r-\sum_{i=1}^{n_1}\al_i+\sum_{s=1}^{n_2}\be_s-\sum_{j=n_2+1}^n\be_j=k.\eqno(8.2.28)$$
Thus $V={\msr H}_{\la k\ra}'$.

 Let $U$ be an
$o(2n+1,\mbb{F})$-submodule generated by $x_{n_1}^k$. Since
$x_{n_1}^k\in{\msr H}'_{\la -k\ra}$ that is an
$o(2n+1,\mbb{F})$-submodule, we have $U\subset {\msr H}'_{\la
-k\ra}$. Repeatedly applying (8.2.17) and (8.2.19) to $x_{n_1}^k$
with various $i\in\ol{1,n_1}$ and $j\in\ol{n_2+1,n}$, we have
$$x_{n_1}^k[\prod_{i=1}^{n_1+1}x_i^{\al_i}][\prod_{s=n_2+1}^ny_s^{\be_s}]\in
U\;
\for\;\;\al_1,...,\al_{n_1+1},\be_{n_2+1},...,\be_n\in\mbb{N}\eqno(8.2.29)$$
such that $\sum_{r=1}^{n_1}\al_r+\sum_{s=n_2+1}^n\be_s=\al_{n_1+1}$.
Note
$$(E_{n_1,i}-E_{n+i,n+n_1})|_{{\msr
B}'}=-x_i\ptl_{x_{n_1}}-y_i\ptl_{y_{n_1}}\qquad\for\;\;i\in\ol{n_1-1}\eqno(8.2.30)$$
by (8.2.1)-(8.2.3). Repeatedly applying (8.2.31) to (8.2.30), we
have
$$[\prod_{i=1}^{n_1+1}x_i^{\al_i}][\prod_{s=n_2+1}^ny_s^{\be_s}]\in
U\;
\for\;\;\al_1,...,\al_{n_1+1},\be_{n_2+1},...,\be_n\in\mbb{N}\eqno(8.2.31)$$
such that
$\sum_{r=1}^{n_1}\al_r+\sum_{s=n_2+1}^n\be_s=\al_{n_1+1}+k$.

Repeatedly applying ${\msr G}|_{{\msr B}'}$ to (8.2.31), we get
$$\wht T_\iota(x^\al y^\be)\in U\eqno(8.2.32)$$
for $\al,\be\in\mbb{N}^n$ such that $\al_i=0$ if $i>n_2$, $\be_j=0$
if $j\leq n_1$,  and
$$\iota+\sum_{r=n_1+1}^{n_2}\al_r-\sum_{i=1}^{n_1}\al_i+\sum_{s=n_1+1}^{n_2}\be_s-\sum_{j=n_2+1}^n\be_j=-k.\eqno(8.2.33)$$
Repeatedly applying (8.2.20) to (8.2.33) satisfying the above
conditions, we get that (8.2.32) holds for $\al,\be\in\mbb{N}^n$
such that $\al_i=0$ if $i>n_2$,  and
$$\iota+\sum_{r=n_1+1}^{n_2}\al_r-\sum_{i=1}^{n_1}\al_i+\sum_{s=1}^{n_2}\be_s-\sum_{j=n_2+1}^n\be_j=-k.\eqno(8.2.34)$$
Repeatedly applying (8.2.27) to (8.2.32) satisfying (8.2.34), we get
(8.2.32) for any $\al,\be\in\mbb{N}^n$ such that
$$\iota+\sum_{r=n_1+1}^n\al_r-\sum_{i=1}^{n_1}\al_i+\sum_{s=1}^{n_2}\be_s-\sum_{j=n_2+1}^n\be_j=-k.\eqno(8.2.35)$$
Thus $U={\msr H}_{\la -k\ra}'$. $\qquad\Box$\psp

{\bf Lemma 8.2.2}. {\it Suppose $n_1=n_2$. For $k\in\mbb{N}$, ${\msr
H}_{\la -k\ra}'$ is an $o(2n+1,\mbb{F})$-submodule generated by
$x_{n_1}^k$ and ${\msr H}_{\la k\ra}'$ is an
$o(2n+1,\mbb{F})$-submodule generated by $\wht T_1(y_{n_1}^{k-1})$
when $k>0$.}

{\it Proof}. In this case,
$$\td\Dlt=-\sum_{i=1}^{n_1}x_i\ptl_{y_i}-\sum_{s=n_1+1}^n
y_s\ptl_{x_s}\eqno(8.2.36)$$ Note
$$ (E_{n+n_1,0}-E_{0,n_1})|_{{\msr
B}'}=y_{n_1}\ptl_{x_0}+x_0x_{n_1} \eqno(8.2.37)$$ by (8.2.6) and
(8.2.7) . Thus for $u\in\mbb{\msr B}$,
\begin{eqnarray*}& &(E_{n+n_1,0}-E_{0,n_1})[\wht T_0(u)]\\
&=&(y_{n_1}\ptl_{x_0}+x_0x_{n_1})(\sum_{i=0}^\infty\frac{(-2)^ix_0^{2i}}{(2i)!}\td\Dlt^i(u))\\
&=&\sum_{i=1}^\infty\frac{(-2)^ix_0^{2i-1}}{(2i-1)!}y_{n_1}\td\Dlt^i(u)+\sum_{i=0}^\infty\frac{(-2)^ix_0^{2i+1}}{(2i)!}\td\Dlt^i(x_{n_1}u)\\
&=&\sum_{i=1}^\infty\frac{(-2)^iix_0^{2i-1}}{(2i-1)!}\td\Dlt^{i-1}(x_{n_1}u)+\sum_{i=1}^\infty\frac{(-2)^ix_0^{2i-1}}{(2i-1)!}
\td\Dlt^{i-1}[\td\Dlt(y_{n_1}u)]\\ & &+
\sum_{i=0}^\infty\frac{(-2)^ix_0^{2i+1}}{(2i)!}\td\Dlt^i(x_{n_1}u)
\\
&=&-\sum_{i=0}^\infty\frac{(-2)^i2(i+1)x_0^{2i+1}}{(2i+1)!}\td\Dlt^i(x_{n_1}u)-2\sum_{i=0}^\infty\frac{(-2)^ix_0^{2i+1}}{(2i+1)!}\td\Dlt^i[\td\Dlt(y_{n_1}u)]\\
& &+
\sum_{i=0}^\infty\frac{(-2)^ix_0^{2i+1}}{(2i)!}\td\Dlt^i(x_{n_1}u)
=-\wht T_1(x_{n_1}u)-2\wht
T_1[\td\Dlt(y_{n_1}u)]\hspace{3.7cm}(8.2.38)\end{eqnarray*} and
\begin{eqnarray*}\qquad& &(E_{n+n_1,0}-E_{0,n_1})[\wht T_1(u)]\\
&=&(y_{n_1}\ptl_{x_0}+x_0x_{n_1})(\sum_{i=0}^\infty\frac{(-2)^ix_0^{2i+1}}{(2i+1)!}\td\Dlt^i(u))\\
&=&\sum_{i=0}^\infty\frac{(-2)^ix_0^{2i}}{(2i)!}y_{n_1}\td\Dlt^i(u)+\sum_{i=0}^\infty\frac{(-2)^ix_0^{2i+2}}{(2i+1)!}\td\Dlt^i(x_{n_1}u)\\
&=&-\sum_{i=1}^\infty\frac{(-2)^{i-1}x_0^{2i}}{(2i-1)!}\td\Dlt^{i-1}(x_{n_1}u)+\sum_{i=0}^\infty\frac{(-2)^ix_0^{2i}}{(2i)!}\td\Dlt^i(y_{n_1}u)\\
& &+
\sum_{i=0}^\infty\frac{(-2)^ix_0^{2i+2}}{(2i+1)!}\td\Dlt^i(x_{n_1}u)
\\
&=&\sum_{i=0}^\infty\frac{(-2)^ix_0^{2i}}{(2i)!}\td\Dlt^i(y_{n_1}u)=\wht
T_0(y_{n_1}u).\hspace{6.6cm}(8.2.39)\end{eqnarray*} Symmetrically,
we have,
$$ (E_{n_1+1,0}-E_{0,n+n_1+1})|_{{\msr
B}'}=x_{n_1+1}\ptl_{x_0}+x_0y_{n_1+1} \eqno(8.2.40)$$ by (8.2.6) and
(8.2.7). So
$$(E_{n_1+1,0}-E_{0,n+n_1+1})[\wht T_0(u)]=-\wht T_1(y_{n_1+1}u)-2\wht T_1[\td\Dlt(x_{n_1+1}u)]\eqno(8.2.41)$$ and
$$(E_{n_1+1,0}-E_{0,n+n_1+1})[\wht T_1(u)]=\wht T_0(x_{n_1+1}u).\eqno(8.2.42)$$

Let $V$ be an $o(2n+1,\mbb{F})$-submodule generated by $x_{n_1}^k$.
Since $x_{n_1}^k\in{\msr H}'_{\la -k\ra}$ that is an
$o(2n+1,\mbb{F})$-submodule, we have $V\subset {\msr H}'_{\la
-k\ra}$.  Note
$$(E_{n_1+1,n+n_1}-E_{n_1,n+n_1+1})|_{{\msr
B}'}=x_{n_1+1}\ptl_{y_{n_1}}+y_{n_1+1}\ptl_{x_{n_1}}\eqno(8.2.43)$$
by (8.2.4). Repeatedly applying (8.2.43) to $x_{n_1}^k$, we have
$$x_{n_1}^{k-r}y_{n_1+1}^r=\wht T_0(x_{n_1}^{k-r}y_{n_1+1}^r)\in
V\qquad\for\;\;r\in\ol{0,k}.\eqno(8.2.44)$$ According to (8.2.38),
\begin{eqnarray*}& &(E_{n+n_1,0}-E_{0,n_1})[\wht T_0(x_{n_1}^{k-r}y_{n_1+1}^r)]
\\ &=&-\wht T_1(x_{n_1}^{k-r+1}y_{n_1+1}^r))-2\wht T_1[\td\Dlt(x_{n_1}^{k-r}y_{n_1}y_{n_1+1}^{r})]
=\wht T_1(x_{n_1}^{k-r+1}y_{n_1+1}^r) \in V.\hspace{1.7cm}(8.2.45)
\end{eqnarray*}
Similarly, (8.2.41) shows
$$(E_{n_1+1,0}-E_{0,n+n_1+1})[\wht T_0(x_{n_1}^{k-r}y_{n_1+1}^r)]=
\wht T_1(x_{n_1}^{k-r}y_{n_1+1}^{r+1}) \in V.\eqno(8.2.47)$$

Suppose that
$$\wht T_1(x_{n_1}^{m_1}y_{n_1}^{m_2}x_{n_1+1}^{m_3}y_{n_1+1}^{m_4})\in
V\qquad\for\;\;m_i\in\mbb{N},\;m_1+m_4-m_2-m_3=k+1\eqno(8.2.48)$$
and $ \sum_{i=1}^4m_i=\ell$. Then (8.2.39) implies
$$(E_{n+n_1,0}-E_{0,n_1})
[\wht T_1(x_{n_1}^{m_1}y_{n_1}^{m_2}x_{n_1+1}^{m_3}y_{n_1+1}^{m_4})]
=\wht
T_0(x_{n_1}^{m_1}y_{n_1}^{m_2+1}x_{n_1+1}^{m_3}y_{n_1+1}^{m_4})\in
V\eqno(8.2.49)$$ and (8.2.42) yields
$$(E_{n_1+1,0}-E_{0,n+n_1+1})[\wht T_1(x_{n_1}^{m_1}y_{n_1}^{m_2}x_{n_1+1}^{m_3}y_{n_1+1}^{m_4})]
=\wht
T_0(x_{n_1}^{m_1}y_{n_1}^{m_2}x_{n_1+1}^{m_3+1}y_{n_1+1}^{m_4})\in
V.\eqno(8.2.50)$$ For $i_1,i_2\in\mbb{N}$ such that $i_1+i_2=1$,
(8.2.36), (8.2.37), (8.2.49) and (8.2.50) give
\begin{eqnarray*}& &(E_{n+n_1,0}-E_{0,n_1})[\wht T_0(x_{n_1}^{m_1}y_{n_1}^{m_2+i_1}x_{n_1+1}^{m_3+i_2}y_{n_1+1}^{m_4})]
\\ &=&-\wht T_1(x_{n_1}^{m_1+1}y_{n_1}^{m_2+i_1}x_{n_1+1}^{m_3+i_2}y_{n_1+1}^{m_4})-2\wht T_1[\td\Dlt(x_{n_1}^{m_1}y_{n_1}^{m_2+i_1+1}x_{n_1+1}^{m_3+i_2}y_{n_1+1}^{m_4})]
\\
&=&(2(m_2+i_1)+1)\wht T_1(x_{n_1}^{m_1+1}y_{n_1}^{m_2+i_1}x_{n_1+1}^{m_3+i_2}y_{n_1+1}^{m_4})\\
&&+2(m_3+i_2)\wht T_1
(x_{n_1}^{m_1}y_{n_1}^{m_2+i_1+1}x_{n_1+1}^{m_3+i_2-1}y_{n_1+1}^{m_4+1})\in
V.\hspace{5.05cm}(8.2.51)
\end{eqnarray*}
On the other hand,
$$(E_{n_1+1,n_1}-E_{n+n_1,n+n_1+1})|_{{\msr
B}'}=-x_{n_1}x_{n_1+1}+y_{n_1}y_{n_1+1}\eqno(8.2.52)$$ by
(8.2.1)-(8.2.3). Thus
\begin{eqnarray*}& &(E_{n_1+1,n_1}-E_{n+n_1,n+n_1+1})
[\wht
T_1(x_{n_1}^{m_1}y_{n_1}^{m_2+i_1}x_{n_1+1}^{m_3+i_2-1}y_{n_1+1}^{m_4})]
\\ &=&-\wht T_1(x_{n_1}^{m_1+1}y_{n_1}^{m_2+i_1}x_{n_1+1}^{m_3+i_2}y_{n_1+1}^{m_4})
+\wht T_1
(x_{n_1}^{m_1}y_{n_1}^{m_2+i_1+1}x_{n_1+1}^{m_3+i_1-1}y_{n_1+1}^{m_4+1})\in
V.\hspace{1.5cm}(8.2.53)
\end{eqnarray*}
Solving the system (8.2.51) and (8.2.53) for $\wht
T_1(x_{n_1}^{m_1+1}y_{n_1}^{m_2+i_1}x_{n_1+1}^{m_3+i_2}y_{n_1+1}^{m_4})$,
we get
$$\wht T_1(x_{n_1}^{m_1+1}y_{n_1}^{m_2+i_1}x_{n_1+1}^{m_3+i_2}y_{n_1+1}^{m_4})\in
V.\eqno(8.2.54)$$

Symmetrically, (8.2.36), (8.2.41), (8.2.49) and (8.2.50) give
\begin{eqnarray*}& &(E_{n_1+1,0}-E_{0,n+n_1+1})[\wht T_0(x_{n_1}^{m_1}y_{n_1}^{m_2+i_1}x_{n_1+1}^{m_3+i_2}y_{n_1+1}^{m_4})]
\\ &=&-\wht T_1(x_{n_1}^{m_1}y_{n_1}^{m_2+i_1}x_{n_1+1}^{m_3+i_2}y_{n_1+1}^{m_4+1})-2\wht T_1[{\msr
D}(x_{n_1}^{m_1}y_{n_1}^{m_2+i_1}x_{n_1+1}^{m_3+i_2+1}y_{n_1+1}^{m_4})]
\\
&=&(2(m_3+i_2)+1)\wht T_1(x_{n_1}^{m_1}y_{n_1}^{m_2+i_1}x_{n_1+1}^{m_3+i_2}y_{n_1+1}^{m_4+1})\\
&&+2(m_2+i_1)\wht T_1
(x_{n_1}^{m_1+1}y_{n_1}^{m_2+i_1-1}x_{n_1+1}^{m_3+i_2+1}y_{n_1+1}^{m_4})\in
V\hspace{4.8cm}(8.2.55)
\end{eqnarray*}
and (8.2.52) yields
\begin{eqnarray*}& &(E_{n_1+1,n_1}-E_{n+n_1,n+n_1+1})
[\wht
T_1(x_{n_1}^{m_1}y_{n_1}^{m_2+i_1-1}x_{n_1+1}^{m_3+i_2}y_{n_1+1}^{m_4})]
\\ &=&-\wht T_1(x_{n_1}^{m_1+1}y_{n_1}^{m_2+i_1-1}x_{n_1+1}^{m_3+i_2+1}y_{n_1+1}^{m_4})
+\wht T_1
(x_{n_1}^{m_1}y_{n_1}^{m_2+i_1}x_{n_1+1}^{m_3+i_1}y_{n_1+1}^{m_4+1})\in
V.\hspace{1.5cm}(8.2.56)
\end{eqnarray*}
Solving the system (8.2.55) and (8.2.56) for $\wht
T_1(x_{n_1}^{m_1}y_{n_1}^{m_2+i_1}x_{n_1+1}^{m_3+i_2}y_{n_1+1}^{m_4+1})$,
we find
$$\wht T_1(x_{n_1}^{m_1}y_{n_1}^{m_2+i_1}x_{n_1+1}^{m_3+i_2}y_{n_1+1}^{m_4+1})\in
V.\eqno(8.2.57)$$ By induction on $\sum_{i=1}^4m_i$, we conclude
$$\wht T_\iota(x_{n_1}^{m_1}y_{n_1}^{m_2}x_{n_1+1}^{m_3}y_{n_1+1}^{m_4})\in
V\qquad\for\;\;m_i\in\mbb{N},\;m_1+m_4-m_2-m_3=k+\dlt_{\iota,1}.\eqno(8.2.58)$$

Set
$${\msr A}'=\mbb{F}[x_1,...,x_{n_1},y_1,...,y_{n_1}],\;\;{\msr A}^\ast=\mbb{F}[x_{n_1+1},...,x_n,
y_{n_1+1},...,y_n],\eqno(8.2.59)$$
$${\msr K}_1=\sum_{i,j=1}^{n_1}[\mbb{F}(E_{i,j}-E_{n+j,n+i})+
\mbb{F}(E_{i,n+j}-E_{j,n+i})+\mbb{F}(E_{n+i,j}-E_{n+j,i})]\eqno(8.2.60)$$
and
$${\msr K}_2=\sum_{i,j=n_1+1}^{n}[\mbb{F}(E_{i,j}-E_{n+j,n+i})+
\mbb{F}(E_{i,n+j}-E_{j,n+i})+\mbb{F}(E_{n+i,j}-E_{n+j,i})].\eqno(8.2.61)$$
Then ${\msr K}_1\cong gl(n_1,\mbb{F})$ and ${\msr K}_2\cong
gl(n-n_1,\mbb{F})$ are  Lie subalgebras of $o(2n+1,\mbb F)$ and
${\msr B}={\msr A}'{\msr A}^\ast$. Since
$${\msr A}'_{\la\ell_1,\ell_2\ra}={\msr A}'\bigcap
{\msr B}_{\la\ell_1,\ell_2\ra},\;\;{\msr
A}^\ast_{\la\ell_1,\ell_2\ra}={\msr A}^\ast\bigcap {\msr
B}_{\la\ell_1,\ell_2\ra}\eqno(8.2.62)$$ are of finite-dimensional
for any $\ell_1,\ell_2\in\mbb{Z}$ (cf. (6.4.29)) and
$${\msr A}'=\bigoplus_{\ell_1,\ell_2\in\mbb{Z}}{\msr
A}'_{\la\ell_1,\ell_2\ra},\;\;{\msr
A}^\ast=\bigoplus_{\ell_1,\ell_2\in\mbb{Z}}{\msr
A}^\ast_{\la\ell_1,\ell_2\ra},\eqno(8.2.63)$$ ${\msr A}'$ is a
${\msr K}_1$-module generated by its ${\msr K}_1$-singular vectors
and ${\msr A}^\ast$ is a ${\msr K}_2$-module generated by its ${\msr
K}_2$-singular vectors. According to (6.6.18), (6.6.29) and (6.6.32)
$$\{x_{n_1}^{m_1}y_{n_1}^{m_2}\zeta_1^{\dlt_{1,n_1}m_3}\mid
m_i\in\mbb{N}\}\eqno(8.2.64)$$ is the set of homogeneous ${\msr
K}_1$-singular vectors in ${\msr A}'$ and
$$\{x_{n_1+1}^{m_1}y_{n_1+1}^{m_2}\zeta_2^{\dlt_{1,n-n_1}m_3}\mid
m_i\in\mbb{N}\}\eqno(8.2.65)$$ is the set of homogeneous ${\msr
K}_2$ singular vectors in ${\msr A}^\ast$. Therefore, ${\msr A}$ is
a $({\msr K}_1+{\msr K}_2)$-module generated by
$$\{x_{n_1}^{m_1}y_{n_1}^{m_2}\zeta_1^{\dlt_{1,n_1}m_3}x_{n_1+1}^{m_4}y_{n_1+1}^{m_5}
\zeta_2^{\dlt_{1,n-n_1}m_6}\mid m_i\in\mbb{N}\}.\eqno(8.2.66)$$

 Observe
$$(E_{n+n_1-1,n_1}-E_{n+n_1,n_1-1})|_{{\msr B}'}=\zeta_1,\;\;
(E_{n_1+2,n+n_1+1}-E_{n_1+1,n+n_1+2})|_{{\msr
B}'}=\zeta_2\eqno(8.2.67)$$ as multiplication operators on ${\msr
B}'$ by (8.2.4) and (8.2.5). Repeatedly applying (8.2.67) to
(8.2.58), we obtain
$$\wht T_\iota(x_{n_1}^{m_1}y_{n_1}^{m_2}\zeta_1^{\dlt_{1,n_1}m_5}x_{n_1+1}^{m_3}y_{n_1+1}^{m_4}\zeta_2^{\dlt_{1,n-n_1}m_6})\in
V\eqno(8.2.68)$$ for $m_i\in\mbb{N}$ such that
$m_1+m_4-m_2-m_3=k+\dlt_{\iota,1}$. Applying $U({\msr K}_1+{\msr
K}_2)$ to (8.2.68), we have
$$\wht T_0({\msr A}_{\la-k\ra}),\wht T_1({\msr A}_{\la-k-1\ra})\subset
V.\eqno(8.2.69)$$ According to (8.2.17), ${\msr
H}'_{\la-k\ra}\subset V$. Thus ${\msr H}'_{\la-k\ra}= V$.

Let $U$ be an $o(2n+1,\mbb{F})$-submodule generated by $\wht
T_0(y_{n_1}^k)$ with $k>0$. Since $T_1(y_{n_1}^{k-1})\in{\msr
H}'_{\la k\ra}$ that is an $o(2n+1,\mbb{F})$-submodule, we have
$U\subset {\msr H}'_{\la k\ra}$. Repeatedly applying (8.2.43) to
$\wht T_1(y_{n_1}^{k-1})$, we have
$$\wht T_1(x_{n_1+1}^ry_{n_1}^{k-r-1})\in
U\qquad\for\;\;r\in\ol{0,k-1}.\eqno(8.2.70)$$ By the same arguments
as (8.2.47)-(8.2.69) with $V$ replaced by $U$ and $k$ replaced by
$-k$, we prove ${\msr H}'_{\la k\ra}=U.\qquad\Box$\psp

{\bf Theorem 8.2.3}. {\it  For any $k\in\mbb{Z}$, ${\msr H}'_{\la
k\ra}$ is an irreducible $o(2n+1,\mbb{F})$-submodule and ${\msr
B}'_{\la k\ra}=\bigoplus_{i=0}^\infty{\eta'}^i({\msr H}'_{\la
k-2i\ra})$ is a decomposition of irreducible submodules. } \psp

{\it Proof}. Note that
$${\msr K}=\sum_{i,j=1}^n\mbb{F}(E_{i,j}-E_{n+j,n+i})\eqno(8.2.71)$$
forms a Lie subalgebra of $o(2n+1,\mbb{F})$, which is isomorphic to
$gl(n,\mbb{F})$. Our representation (8.2.1)-(8.2.7) of
$o(2n+1,\mbb{F})$ is an extension the representation
(6.3.36)-(6.3.38) of $sl(n,\mbb F)$ with ${\cal Q}$ replaced by
$\msr B$. First we prove that ${\msr H}'_{\la m\ra}$ is an
irreducible $o(2n+1,\mbb{F})$-submodule for any $m\in\mbb{Z}$. We
divide it into two cases. \psp

{\it Case 1}. $n_1+1\leq n_2$.\psp

 Consider ${\msr H}'_{\la k\ra}$ with
$k\in\mbb{N}$. Let $U$ be any nonzero $o(2n+1,\mbb{F})$-submodule of
${\msr H}'_{\la k\ra}$.  Since $U$ is a ${\msr K}$-module, it must
contain a ${\msr K}$-singular vector. According to Lemma 6.4.1,
(6.5.54) and (6.5.61), $U$ must contain some vectors in
$$\{\mbb{F}[\eta](x_i^{m_1}y_j^{m_2})\mid
m_1,m_2\in\mbb{N};i=n_1,n_1+1;j=n_2,n_2+1\}\eqno(8.2.72)$$ if
$n_2<n$, or
$$\{\mbb{F}[\eta](x_i^{m_1}y_n^{m_2})\mid
m_1,m_2\in\mbb{N};i=n_1,n_1+1\}\eqno(8.2.73)$$when $n_2=n$.

 Note
$$(E_{n+1,n+n_2}-E_{n_2,1})|_{{\msr
B}'}=y_1\ptl_{y_{n_2}}+x_1x_{n_2}\eqno(8.2.74)$$ by (8.2.1)-(8.2.3).
Moreover,
$$(E_{n_1+1,n+1}-E_{1,n+n_1+1})|_{{\msr
B}'}=x_{n_1+1}\ptl_{y_1}-x_1\ptl_{y_{n_1+1}}\eqno(8.2.75)$$ by
(8.2.4).
 If
$\wht T_\iota(\eta^\ell(x_j^{k_1}y_{n_2}^{k_2}))\in U$ with
$j=n_1,n_1+1$, we have
\begin{eqnarray*}& &(k_2!)^{-2}(E_{n_1+1,n+1}-E_{1,n+n_1+1})^{k_2}(E_{n+1,n+n_2}-E_{n_2,1})^{k_2}
[\wht T_\iota(\eta^\ell(x_j^{k_1}y_{n_2}^{k_2}))]\\
&=&\wht T_\iota(\eta^\ell(x_j^{k_1}x_{n_1+1}^{k_2}))\in
U.\hspace{9.8cm}(8.2.76)\end{eqnarray*} We use (8.2.73) because
$n_2$ may be equal to $n_1+1$.  Consider the case $\wht
T_\iota(\eta^\ell(x_j^{k_1}y_{n_2+1}^{k_2}))\in U$ with
$j=n_1,n_1+1$. Observe
$$(E_{n+n_1,n_2+1}-E_{n+n_2+1,n_1})_{{\msr
B}'}=y_{n_1}\ptl_{x_{n_2+1}}+x_{n_1}\ptl_{y_{n_2+1}}\eqno(8.2.77)$$
by (8.2.5). Hence
$$\frac{1}{k_2!}(E_{n+n_1,n_2+1}-E_{n+n_2+1,n_1})^{k_2}
[\wht T_\iota(\eta^\ell(x_j^{k_1}y_{n_2+1}^{k_2}))]=\wht
T_\iota(\eta^\ell(x_j^{k_1}x_{n_1}^{k_2}))\in U.\eqno(8.2.78)$$
Therefore we always have
$$\wht T_\iota(\eta^\ell(x_{n_1}^{k_1}x_{n_1+1}^{k_2}))\in U.\eqno(8.2.79)$$

Note
$$(E_{n_1+1,0}-E_{0,n+n_1+1})|_{{\msr B}'}=x_{n_1+1}\ptl_{x_0}-x_0\ptl_{y_{n_1+1}}\eqno(8.2.80)$$
by (8.2.6) and (8.2.7). Thus (7.2.7) and (7.2.8) give
\begin{eqnarray*}&
&(E_{n_1+1,0}-E_{0,n+n_1+1})[\wht
T_0(\eta^\ell(x_{n_1}^{k_1}x_{n_1+1}^{k_2}))]
\\ &=&(x_{n_1+1}\ptl_{x_0}-x_0\ptl_{y_{n_1+1}})(\sum_{i=0}^\infty\frac{(-2)^ix_0^{2i}}{(2i)!}\td\Dlt^i(\eta^\ell(x_{n_1}^{k_1}x_{n_1+1}^{k_2})))\\
&=&\sum_{i=1}^\infty\frac{(-2)^ix_0^{2i-1}}{(2i-1)!}x_{n_1+1}\td\Dlt^i(\eta^\ell(x_{n_1}^{k_1}x_{n_1+1}^{k_2}))-\sum_{i=0}^\infty\frac{(-2)^ix_0^{2i+1}}{(2i)!}\td\Dlt^i(\ptl_{y_{n_1+1}}(\eta^\ell(x_{n_1}^{k_1}x_{n_1+1}^{k_2})))\\
&=&\sum_{i=1}^\infty\frac{(-2)^ix_0^{2i-1}}{(2i-1)!}\td\Dlt^i(\eta^\ell(x_{n_1}^{k_1}x_{n_1+1}^{k_2+1}))
-\sum_{i=1}^\infty\frac{(-2)^iix_0^{2i-1}}{(2i-1)!}\td\Dlt^{i-1}(\ptl_{y_{n_1+1}}\eta^\ell(x_{n_1}^{k_1}x_{n_1+1}^{k_2}))
\\& &-\ell\sum_{i=0}^\infty\frac{(-2)^ix_0^{2i+1}}{(2i)!}\td\Dlt^i(\eta^{\ell-1}(x_{n_1}^{k_1}x_{n_1+1}^{k_2+1})))
\\&=&\ell(n_2-n_1+\ell+k_2-k_1)\sum_{i=1}^\infty\frac{(-2)^ix_0^{2i-1}}{(2i-1)!}\td\Dlt^{i-1}(\eta^{\ell-1}(x_{n_1}^{k_1}x_{n_1+1}^{k_2+1}))\\
& &+\ell \wht T_1((\eta^{\ell-1}(x_{n_1}^{k_1}x_{n_1+1}^{k_2+1})))
\\ &=&\ell[1-2(n_2-n_1+\ell+k_2-k_1)]\wht T_1((\eta^{\ell-1}(x_{n_1}^{k_1}x_{n_1+1}^{k_2+1})))
\hspace{4.5cm}(8.2.81)
\end{eqnarray*}
and
\begin{eqnarray*}&
&(E_{n_1+1,0}-E_{0,n+n_1+1})[\wht
T_1(\eta^\ell(x_{n_1}^{k_1}x_{n_1+1}^{k_2}))]
\\ &=&(x_{n_1+1}\ptl_{x_0}-x_0\ptl_{y_{n_1+1}})(\sum_{i=0}^\infty\frac{(-2)^ix_0^{2i+1}}{(2i+1)!}\td\Dlt^i(\eta^\ell(x_{n_1}^{k_1}x_{n_1+1}^{k_2})))
\\&=&\sum_{i=0}^\infty\frac{(-2)^ix_0^{2i}}{(2i)!}x_{n_1+1}\td\Dlt^i(\eta^\ell(x_{n_1}^{k_1}x_{n_1+1}^{k_2}))
-\sum_{i=0}^\infty\frac{(-2)^ix_0^{2i+2}}{(2i+1)!}\td\Dlt^i(\ptl_{y_{n_1+1}}(\eta^\ell(x_{n_1}^{k_1}x_{n_1+1}^{k_2})))\\
&=&\sum_{i=0}^\infty\frac{(-2)^ix_0^{2i}}{(2i)!}\td\Dlt^i(\eta^\ell(x_{n_1}^{k_1}x_{n_1+1}^{k_2+1}))
+\sum_{i=1}^\infty\frac{(-2)^{i-1}x_0^{2i}}{(2i-1)!}\td\Dlt^{i-1}(\ptl_{y_{n_1+1}}(\eta^\ell(x_{n_1}^{k_1}x_{n_1+1}^{k_2})))
\\ & &-\ell\sum_{i=0}^\infty\frac{(-2)^ix_0^{2i+2}}{(2i+1)!}\td\Dlt^i(\eta^{\ell-1}(x_{n_1}^{k_1}x_{n_1+1}^{k_2+1})))=\wht
T_0(\eta^\ell(x_{n_1}^{k_1}x_{n_1+1}^{k_2+1})).
\hspace{3.1cm}(8.2.82)\end{eqnarray*} By (8.2.79), (8.2.81) and
(8.2.82), we get
$$\wht T_0(x_{n_1}^{k_1}x_{n_1+1}^{k_2})\in U\qquad \mbox{for
some}\;k_1,k_2\in\mbb{N}\;\mbox{such
that}\;k_2-k_1=k.\eqno(8.2.83)$$

Observe
$$(E_{n_1,n_1+1}-E_{n+n_1+1,n+n_1})_{{\msr B}'}=\ptl_{x_{n_1}}\ptl_{x_{n_1+1}}
-y_{n_1+1}\ptl_{y_{n_1}}\eqno(8.2.84)$$ by (8.2.1)-(8.2.3).
Repeatedly applying (8.2.84) to (8.2.83), we obtain
$x_{n_1+1}^k=T_0(x_{n_1+1}^k)\in U$. Thanks to Lemma 8.2.1, $U={\msr
H}'_{\la k\ra}$. So ${\msr H}'_{\la k\ra}$ is an irreducible
$o(2n+1,\mbb{F})$-submodule.

Let $V$ be any nonzero $o(2n+1,\mbb{F})$-submodule of ${\msr
H}'_{\la -k\ra}$ with $k\in\mbb{N}+1$. By the above arguments,
$x_{n_1}^k\in V$.  Thanks to Lemma 8.2.1, $V={\msr H}'_{\la -k\ra}$.
So ${\msr H}'_{\la -k\ra}$ is an irreducible
$o(2n+1,\mbb{F})$-submodule.\psp

{\it Case 2}. $n_1=n_2$.\psp

In this case
$$\eta=\sum_{i=1}^{n_1}y_i\ptl_{x_i}+\sum_{s=n_1+1}^n
x_s\ptl_{y_s}.\eqno(8.2.85)$$ Thus
$$x_{n_1}^{m_1}y_{n_1}^{m_2}\zeta_1^{m_3+1}=\frac{1}{\prod_{i=1}^{m_2}(m_1+i)}\eta^{m_2}(
x_{n_1}^{m_1+m_2}\zeta_1^{m_3+1})\eqno(8.2.86)$$ and
$$x_{n_1+1}^{m_1}y_{n_1+1}^{m_2}\zeta_2^{m_3+1}=
\frac{1}{\prod_{i=1}^{m_1}(m_2+i)}\eta^{m_1}(
y_{n_1+1}^{m_1+m_2}\zeta_2^{m_3+1}).\eqno(8.2.87)$$
 First we consider ${\msr H}'_{\la k\ra}$ with
$k\in\mbb{N}+1$. Let $U$ be any nonzero $o(2n+1,\mbb{F})$-submodule
of ${\msr H}'_{\la k\ra}$. Then  $U$ must contain a ${\msr
K}$-singular vector. Moreover, (8.2.77) becomes
$$(E_{n+n_1,n_1+1}-E_{n+n_1+1,n_1})_{{\msr
B}'}=y_{n_1}\ptl_{x_{n_1+1}}+x_{n_1}\ptl_{y_{n_1+1}}.\eqno(8.2.88)$$
If $\wht T_\iota(\eta^\ell(x_{n_1}^{k_1}y_{n_1+1}^{k_2})\in U$, then
$$\frac{1}{k_2!}(E_{n+n_1,n_1+1}-E_{n+n_1+1,n_1})^{k_2}
[\wht T_\iota(\eta^\ell(x_{n_1}^{k_1}y_{n_1+1}^{k_2}))]= \wht
T_\iota(\eta^\ell(x_{n_1}^{k_1+k_2})) \in U.\eqno(8.2.89)$$ When
$\wht T_\iota(\eta^\ell(x_{n_1}^{k_1}\zeta_1^{k_2}))\in U$, we have
$\wht T_\iota(\eta^\ell(x_{n_1}^{k_1}))\in U$ by (7.2.40).
Symmetrically, we have $\wht T_\iota(\eta^\ell(y_{n_1+1}^{k_1}))\in
U$ if $\wht T_\iota(\eta^\ell(y_{n_1+1}^{k_1}\zeta_2^{k_2})\in U$,
and (8.2.90) implies $0\neq \wht
T_\iota(\eta^\ell(x_{n_1}^{k_1}))\in U$. By (6.6.18), (6.6.29) and
(6.6.32), we always have some $\wht
T_\iota(\eta^\ell(x_{n_1}^{k_1}))\in U$. According to (8.2.85),
$$\eta^m(x_{n_1}^{m_1})=0\qquad\for\;m,m_1\in\mbb{N}\;\mbox{such
that}\;m>m_1.\eqno(8.2.90)$$ Thus $\ell\leq k_1$.

 Note that
$$ (E_{n_1,0}-E_{0,n+n_1})|_{{\msr
B}'}=\ptl_{x_0}\ptl_{x_{n_1}}-x_0\ptl_{y_{n_1}} \eqno(8.2.91)$$ by
(8.2.6) and (8.2.7). Moreover, (8.2.37), (8.2.85) and (8.2.91) imply
that
\begin{eqnarray*}& &(E_{n_1,0}-E_{0,n+n_1})[\wht T_0(\eta^\ell(x_{n_1}^{k_1}))]\\
&=&(\ptl_{x_0}\ptl_{x_{n_1}}-x_0\ptl_{y_{n_1}})(\sum_{i=0}^\infty\frac{(-2)^ix_0^{2i}}{(2i)!}\td\Dlt^i(\eta^\ell(x_{n_1}^{k_1})))\\
&=&-\sum_{i=1}^\infty\frac{(-2)^iix_0^{2i-1}}{(2i-1)!}\td\Dlt^{i-1}(\ptl_{y_{n_1}}(\eta^\ell(x_{n_1}^{k_1})))
+k_1\sum_{i=1}^\infty\frac{(-2)^ix_0^{2i-1}}{(2i-1)!}\td\Dlt^i(\eta^\ell(x_{n_1}^{k_1-1})))\\&
&-\sum_{i=0}^\infty\frac{(-2)^ix_0^{2i+1}}{(2i)!}\td\Dlt^i(\ptl_{y_{n_1}}(\eta^\ell(x_{n_1}^{k_1})))
\\ &=&\ell k_1\sum_{i=0}^\infty\frac{(-2)^ix_0^{2i+1}}{(2i+1)!}\td\Dlt^i(\eta^{\ell-1}(x_{n_1}^{k_1-1}))-2\ell
k_1(\ell-k_1)\sum_{i=0}^\infty\frac{(-2)^ix_0^{2i+1}}{(2i+1)!}\td\Dlt^i(\eta^{\ell-1}(x_{n_1}^{k_1-1})))\\
&=&\ell k_1[1-2(\ell-k_1)] \wht T_1(\eta^{\ell-1}(x_{n_1}^{k_1-1})),
 \hspace{7.6cm}(8.2.92)\end{eqnarray*}
\begin{eqnarray*}& &(E_{n_1,0}-E_{0,n+n_1})[\wht T_1(\eta^\ell(x_{n_1}^{k_1}))]\\
&=&(\ptl_{x_0}\ptl_{x_{n_1}}-x_0\ptl_{y_{n_1}})(\sum_{i=0}^\infty\frac{(-2)^ix_0^{2i+1}}{(2i+1)!}\td\Dlt^i(\eta^\ell(x_{n_1}^{k_1})))\\
&=&\sum_{i=1}^\infty\frac{(-2)^{i-1}x_0^{2i}}{(2i-1)!}\td\Dlt^{i-1}(\ptl_{y_{n_1}}(\eta^\ell(x_{n_1}^{k_1})))
+\sum_{i=0}^\infty\frac{(-2)^ix_0^{2i}}{(2i)!}\td\Dlt^i(\ptl_{x_{n_1}}(\eta^\ell(x_{n_1}^{k_1})))\\
&
&-\sum_{i=0}^\infty\frac{(-2)^ix_0^{2i+2}}{(2i+1)!}\td\Dlt^i(\ptl_{y_{n_1}}(\eta^\ell(x_{n_1}^{k_1})))=k_1\wht
T_0(\eta^\ell(x_{n_1}^{k_1-1}))
.\hspace{4.1cm}(8.2.93)\end{eqnarray*} Therefore, we have
$$\frac{1}{(k-1)!}\wht T_1(\eta^{k-1}(x_{n_1}^{k-1}))=\wht T_1(y_{n_1}^{k-1})\in U.\eqno(8.2.94)$$
By Lemma 8.2.2, $U={\msr H}'_{\la k\ra}$. Thus ${\msr H}'_{\la
k\ra}$
 is an irreducible $o(2n+1,\mbb{F})$-submodule.

Let $V$ be any nonzero $o(2n+1,\mbb{F})$-submodule of ${\msr
H}'_{\la -k\ra}$ with $k\in\mbb{N}$. By the above arguments
(8.2.85)-(8.2.93), $x_{n_1}^k\in V$.  Thanks to Lemma 8.2.2,
$V={\msr H}'_{\la -k\ra}$. So ${\msr H}'_{\la -k\ra}$ is an
irreducible $o(2n+1,\mbb{F})$-submodule.\psp

Fix $k\in\mbb{Z}$. Since $\td\Dlt'$ in (8.2.8) is locally nilpotent,
for any $0\neq u\in {\msr B}'_{\la k\ra}$, there exists an element
$\kappa(u)\in\mbb{N}$ such that
$$(\td\Dlt')^{\kappa(u)}(u)\neq
0\;\;\mbox{and}\;\;(\td\Dlt')^{\kappa(u)+1}(u)=0.\eqno(8.2.95)$$ Set
$$\Psi'=\sum_{i=0}^\infty{\eta'}^i({\msr
H}'_{\la k-2i\ra}).\eqno(8.2.96)$$ Given $0\neq u\in {\msr B}'_{\la
k\ra}$, $\kappa(u)=1$ implies $u\in {\msr H}'_{\la
k\ra}\subset\Psi'$. Suppose that $u\in \Psi'$ whenever $\kappa(u)<r$
for some positive integer $r$. Assume $\kappa(u)=r$. First
$$v=(\td\Dlt')^r(u)\in {\msr
H}'_{\la k-2r\ra}\subset\Psi'.\eqno(8.2.97)$$ Note that
$[\ptl_{x_0}^2,x_0^2]=2+4x_0\ptl_{x_0}$ and so
$$(\td\Dlt')^r[{\eta'}^r(v)]=2^rr![\prod_{i=1}^r[1+2(n_1-n_2-k+r+i)]]v\eqno(8.2.98)$$
by (7.2.8) and (8.2.14).

 Thus we  have either
$$u=\frac{1}{2^rr!\prod_{i=1}^r[1+2(n_1-n_2-k+r+i)]}{\eta'}^r(v)\in\Psi'\eqno(8.2.99)$$
or
$$\kappa\left(u-\frac{1}{2^rr!\prod_{i=1}^r[1+2(n_1-n_2-k+r+i)]}{\eta'}^r(v)\right)<r.\eqno(8.2.100)$$
By induction,
$$u-\frac{1}{2^rr!\prod_{i=1}^r[1+2(n_1-n_2-k+r+i)]}{\eta'}^r(v)\in\Psi',\eqno(8.2.101)$$
which implies $u\in\Psi'$. Therefore, we have $\Psi'={\msr B}'_{\la
k\ra}$. Since the weight of any ${\msr K}$-singular vector in
${\eta'}^i({\msr H}'_{\la k-2i\ra})$ is different from  the weight
of any ${\msr K}$-singular vector in ${\eta'}^j({\msr H}'_{\la
k-2j\ra})$ when $i\neq j$, the sums in Theorem 8.2.3 are direct
sums. This completes the proof of Theorem 8.2.3. $\qquad\Box$

\section{Extensions of the Conformal Representation}

This section is to study extensions of the conformal representation
of the odd orthogonal Lie algebra $o(2n+3,\mbb{C})$ via
$o(2n+1,\mbb{C})$-modules.

 For $\lmd\in\Lmd^+$, we denote by $V(\lmd)$ the
finite-dimensional irreducible $o(2n+1,\mbb C)$-module with highest
weight $\lmd$. The $(2n+1)$-dimensional natural module of
$o(2n+1,\mbb C)$ is $V(\ves_1)$  with weights $\{0,\pm\ves_i\mid
i\in\ol{1,n}\}$. The following result can be derived from  Theorem
5.4.3:\psp

{\bf Lemma 8.3.1 (Pieri's formula)}. {\it Given $\mu\in\Lmd^+$ with
${\msr S}(\mu)=\{i_1,i_2,...,i_{s+1}\}$ (cf. (6.7.55) and (6.7.56)),
$$V(\ves_1)\otimes_{\mbb C}V(\mu)\cong
(1-\dlt_{\mu_n,0})V(\mu)\oplus\bigoplus_{\iota=1}^{s-\dlt_{\mu_n,0}-\dlt_{\mu_n,1/2}}
V(\mu-\ves_{i_{\iota+1}-1}))\oplus\bigoplus_{r=1}^s
V(\mu+\ves_{i_r}).\eqno(8.3.1)$$}\pse

Note that the Casimir element of $o(2n+1,\mbb C)$ is
\begin{eqnarray*}\omega&=&\sum_{1\leq
i<j\leq
n}[(E_{i,n+j}-E_{j,n+i})(E_{n+j,i}-E_{n+i,j})+(E_{n+j,i}-E_{n+i,j})(E_{i,n+j}-E_{j,n+i})]
\\
&&+\sum_{i=1}^n[(E_{0,i}-E_{n+i,0})(E_{i,0}-E_{0,n+i})+(E_{i,0}-E_{0,n+i})(E_{0,i}-E_{n+i,0})]
\\ & &+\sum_{i,j=1}^n(E_{i,j}-E_{n+j,n+i})(E_{j,i}-E_{n+i,n+j})\in
U(o(2n+1,\mbb C)).\hspace{3cm}(8.3.2)\end{eqnarray*} Recall the
operator $\mfk{d}$ defined in (7.3.7). Set
$$\td\omega=\frac{1}{2}(\mfk{d}(\omega)-\omega\otimes 1-1\otimes
\omega)\in U(o(2n+1,\mbb C))\otimes_{\mbb C} U(o(2n+1,\mbb
C)).\eqno(8.3.3)$$ By (8.3.2),
\begin{eqnarray*}\td\omega\!&=&\!\sum_{1\leq
i<j\leq
n}[(E_{i,n+j}-E_{j,n+i})\otimes(E_{n+j,i}-E_{n+i,j})+(E_{n+j,i}-E_{n+i,j})\otimes(E_{i,n+j}-E_{j,n+i})]
\\
&&+\sum_{i=1}^n[(E_{0,i}-E_{n+i,0})\otimes(E_{i,0}-E_{0,n+i})+(E_{i,0}-E_{0,n+i})\otimes(E_{0,i}-E_{n+i,0})]
\\ & &+\sum_{i,j=1}^n(E_{i,j}-E_{n+j,n+i})\otimes(E_{j,i}-E_{n+i,n+j}).\hspace{5.9cm}(8.3.4)\end{eqnarray*}
Moreover, (7.3.11) also holds. Take the settings in (7.3.12) and
(7.3.13). Denote
$$\ell(\mu)=\dim V(\mu).\eqno(8.3.5)$$

Recall
$$\rho=\frac{1}{2}\sum_{\nu\in\Phi_{B_n}^+}\nu.\eqno(8.3.6)$$
Then
$$\frac{2(\rho,\nu)}{(\nu,\nu)}=1\qquad\for\;\;\nu\in\Pi_{B_n}\eqno(8.3.7)$$
by (3.2.17). Expression (8.1.3) yields
$$\rho=\sum_{i=1}^{n-1}(n-i+1/2)\ves_i.\eqno(8.3.8)$$
 Observe that
$$(\mu+2\rho,\mu)-(\mu+2\rho,\mu)-(\ves_1+2\rho,\ves_1)=-2n,\eqno(8.3.9)$$
$$(\mu+\ves_i+2\rho,\mu+\ves_i)-(\mu+2\rho,\mu)-(\ves_1+2\rho,\ves_1)=2(\mu_i+1-i)\eqno(8.3.10)$$
and
$$(\mu-\ves_i+2\rho,\mu-\ves_i)-(\mu+2\rho,\mu)-(\ves_1+2\rho,\ves_1)=2(i-2n-\mu_i)\eqno(8.3.11)$$
for $\mu=\sum_{r=1}^n\mu_r\ves_r$ by (8.3.8). Moreover, the algebra
$U(o(2n+1,\mbb C))\otimes_{\mbb C}U(o(2n+1,\mbb C))$ acts on
$V(\ves_1)\otimes_{\mbb C}V(\mu)$ by
$$(\xi_1\otimes \xi_2)(v\otimes
u)=\xi_1(v)\otimes\xi_2(u)\;\for\;\xi_1,\xi_2\in U(o(2n+1,\mbb
C)),\;v\in V(\ves_1),\;u\in V(\mu).\eqno(8.3.12)$$ Note that
(7.3.11) also holds for $o(2n+1,\mbb C)$.
  By Lemma 8.3.1, (7.3.11), (8.3.3) and
(8.3.9)-(8.3.11), we get:\psp

{\bf Lemma 8.3.2}. {\it Let $\mu=\sum_{i=1}^n\mu_i\ves_i\in\Lmd^+$
with ${\msr S}(\mu)=\{i_1,i_2,...,i_{s+1}\}$ (cf. (6.7.55) and
(6.7.56)). The characteristic polynomial of
$\td\omega|_{V(\ves_1)\otimes_{\mbb C}V(\mu)}$ is
\begin{eqnarray*}\hspace{2cm} & &(t+n)^{\ell(\mu)(1-\dlt_{\mu_n,0})}[\prod_{\iota=1}^{s-\dlt_{\mu_n,0}-\dlt_{\mu_n,1/2}}(t+\mu_{i_\iota}+2n
-i_{\iota+1})^{\ell^-_{i_{\iota+1}-1}(\mu)}]
\\ & &\times[\prod_{r=1}^{s-1}(t-\mu_{i_r}+i_r-1)^{\ell^+_{i_r}(\mu)}].\hspace{6.8cm}(8.3.13)\end{eqnarray*}}
\pse

We remark that the above lemma is also equivalent to a special
detailed version of Kostant's characteristic identity.

Set
$${\msr A}=\mbb C[x_0,x_1,x_2,...,x_{2n}].\eqno(8.3.14)$$
Then ${\msr A}$ forms an $o(2n+1,\mbb C)$-module with the action
determined via
$$E_{i,j}|_{\msr
A}=x_i\ptl_{x_j}\qquad\for\;\;i,j\in\ol{0,2n}.\eqno(8.3.15)$$ The
corresponding Laplace operator, its dual invariant and the degree
operator are
$$\Dlt=\ptl_{x_0}^2+2\sum_{r=1}^n\partial_{x_r}\partial_{x_{n+r}},\;\;\eta=\frac{1}{2}x_0^2+\sum_{i=1}^nx_ix_{n+i},\;\;D=\sum_{r=0}^{2n}x_r\partial_{x_r}.
\eqno(8.3.16)$$ Denote
$$A_{i,j}=E_{i,j}-E_{n+1+j,n+1+i},\;B_{i,j}=E_{i,n+1+j}-E_{j,n+1+i},\;C_{i,j}=E_{n+1+i,j}-E_{n+1+j,i},\eqno(8.3.17)$$
$$K_i=E_{0,i}-E_{n+1+i,0},\qquad
K_{n+1+i}=E_{0,n+1+i}-E_{i,0}\eqno(8.3.18)$$ for $i,j\in\ol{1,n+1}$.
Then the split odd orthogonal Lie algebra
$$o(2n+3,\mbb C)=\sum_{i,j=1}^{n+1}\mbb CA_{i,j}+\sum_{1\leq i<j\leq
n+1}(\mbb CB_{i,j}+\mbb CC_{j,i})+\sum_{i=1}^{n+1}(\mbb CK_i+\mbb
CK_{n+1+i}).\eqno(8.3.19)$$

Replace $n$ by $2n+1$ and take the product $(\vec x,\vec
y)=(1/2)(x_0y_0+\sum_{i=1}^n(x_iy_{n+i}+x_{n+i}y_i)$ in (7.3.1).
 Using the derivatives of
one-parameter conformal transformations just as one derives the Lie
algebra of a Lie group, we get the following conformal
representation of $o(2n+3,\mbb C)$ determined by
$$A_{i,j}|_{\msr A}=x_i\ptl_{x_j}-x_{n+j}\ptl_{x_{n+i}},\qquad
B_{i,j}|_{\msr
A}=x_i\ptl_{x_{n+j}}-x_j\ptl_{x_{n+i}},\eqno(8.3.20)$$
$$C_{i,j}|_{\msr
A}=x_{n+i}\ptl_{x_j}-x_j\ptl_{x_{n+i}},\qquad A_{n+1,n+1}|_{\msr
A}=-D,\eqno(8.3.21)$$
$$A_{n+1,i}|_{\msr A}=\ptl_{x_i},\qquad B_{i,n+1}|_{\msr
A}=-\ptl_{x_{n+i}},\eqno(8.3.22)$$
$$A_{i,n+1}|_{\msr A}=-x_iD+\eta\ptl_{x_{n+i}},\;\;C_{n+1,i}|_{\msr
A}=x_{n+i}D-\eta\ptl_{x_i},\eqno(8.3.23)$$
$$K_i|_{\msr A}=x_0\ptl_{x_i}-x_{n+i}\ptl_{x_0},\qquad
K_{n+1+i}|_{\msr A}=x_0\ptl_{x_{n+i}}-x_i\ptl_{x_0},\eqno(8.3.24)$$
$$K_{2n+2}|_{\msr A}=-\ptl_{x_0},\qquad K_{n+1}|_{\msr A}=x_0D-\eta\ptl_{x_0}\eqno(8.3.25)$$
 for $i,j\in\ol{1,n}$.

Note that
$$\msr L=\sum_{i,j=1}^n\mbb CA_{i,j}+\sum_{1\leq i<j\leq
n}(\mbb CB_{i,j}+\mbb CC_{j,i}+\sum_{i=1}^n(\mbb CK_i+\mbb
CK_{n+1+i}))\eqno(8.3.26)$$ forms a Lie subalgebra of $o(2n+3,\mbb
C)$ that is isomorphic to $o(2n+1,\mbb C)$. For convenience, we make
the identification of $o(2n+1,\mbb C)$ with $\msr L$ as follows:
$$E_{i,j}-E_{n+j,n+i}\leftrightarrow A_{i,j},\;\;
E_{i,n+j}-E_{j,n+i}\leftrightarrow
B_{i,j},\;\;E_{n+i,j}-E_{n+j,i}\leftrightarrow
C_{i,j},\eqno(8.3.27)$$
$$E_{0,i}-E_{n+i,0}\leftrightarrow K_i,\qquad E_{0,n+i}-E_{i,0}\leftrightarrow
K_{n+1+i}\eqno(8.3.28)$$ for $i,j\in\ol{1,n}$. Recall the Witt
algebra ${\mbb W}_{2n+1}=\sum_{i=0}^{2n}{\msr A}\ptl_{x_i}$, and
Shen's monomorphism $\Im: {\mbb W}_{2n+1}\rta\wht{\mbb
W}_{2n+1}={\mbb W}_{2n+1}\oplus gl(2n+1,{\msr A})$ (cf. (6.7.13))
given by
$$\Im(\sum_{i=0}^{2n}f_i\ptl_{x_i})=\sum_{i=0}^{2n}f_i\ptl_{x_i}\oplus\sum_{i,j=0}^n\ptl_{x_i}(f_j)E_{i,j}.
\eqno(8.3.29)$$ Now we have the Lie algebra monomorphism $\nu:
o(2n+3,\mbb C)\rta \widehat{\mbb W}_{2n+1}$ given by
$$\nu(\xi)=\Im(\xi|_{\msr A})\qquad\for\;\;\xi\in
o(2n+3,\mbb C).\eqno(8.3.30)$$ Denote
$$I_{2n+1}=\sum_{r=0}^{2n}E_{r,r}.\eqno(8.3.31)$$
According to the identification (8.3.27) and (8.3.28), we have
$$\nu(A_{i,j})=(x_i\ptl_{x_j}-x_{n+j}\ptl_{x_{n+i}})\oplus
A_{i,j},\;\;\nu(B_{i,j})=(x_i\ptl_{x_{n+j}}-x_j\ptl_{x_{n+i}})\oplus
B_{i,j},\eqno(8.3.32)$$
$$\nu(C_{i,j})=(x_{n+i}\ptl_{x_j}-x_{n+j}\ptl_{x_i})\oplus
C_{i,j},\;\;\nu(A_{n+1,n+1})=-(D\oplus I_{2n+1}),\eqno(8.3.33)$$
$$\nu(A_{n+1,i})=\ptl_{x_i},\qquad\nu(B_{i,n+1})=-\ptl_{x_{n+i}},\qquad\nu(K_{2n+2})=-\ptl_{x_0},\eqno(8.3.34)$$
$$\nu(A_{i,n+1})=(-x_iD+\eta\ptl_{x_{n+i}})\oplus[\sum_{p=1}^n(x_{n+p}B_{p,i}-
x_pA_{i,p})-x_iI_{2n+1}+x_0K_{n+1+i}],\eqno(8.3.35)$$
$$\nu(C_{n+1,i})=(x_{n+i}D-\eta\ptl_{x_i})\oplus[\sum_{p=1}^n(x_pC_{i,p}-x_{n+p}A_{p,i})+x_{n+i}I_{2n+1}-x_0K_i],
\eqno(8.3.36)$$
$$\nu(K_i)=(x_0\ptl_{x_i}-x_{n+i}\ptl_{x_0})\oplus K_i,\qquad
\nu(K_{n+1+i})=(x_0\ptl_{x_{n+i}}-x_i\ptl_{x_0})\oplus
K_{n+i+i},\eqno(8.3.37)$$
$$\nu(K_{n+1})=(x_0D-\eta\ptl_{x_0})\oplus [\sum_{s=1}^n(x_sK_s+x_{n+s}K_{n+1+s})+x_0I_{2n+1}] \eqno(8.3.38)$$
for $i,j\in\ol{1,n}$.

Observe that
$$\wht{\mbb W}^o_{2n+1}=\wht{\mbb W}_{2n+1}\oplus[o(2n+1,{\msr A})+{\msr
A}I_{2n+1}]\eqno(8.3.39)$$ form a Lie subalgebra of $\wht{\mbb
W}_{2n+1}$ and
$$\nu(o(2n+3,\mbb C))\subset \wht{\mbb W}^o_{2n+1}.\eqno(8.3.40)$$
Let  $M$ be an $o(2n+1,\mbb C)$-module and let $c\in\mbb C$ be a
fixed constant. Then
$$\widehat M={\msr A}\otimes_{\mbb C}M\eqno(8.3.41)$$
becomes a $\wht{\mbb W}^o_{2n+1}$-module with the action:
$$(d+f_1A+f_2I_{2n})(g\otimes
v)=(d(g)+cf_2g)\otimes v+f_1g\otimes A(v)\eqno(8.3.42)$$ for
$f_1,f_2,g\in{\msr A},\;A\in o(2n+1,\mbb C)$ and $v\in M$. Moreover,
we make $\widehat M$ an $o(2n+3,\mbb C)$-module with the action:
$$\xi(\varpi)=\nu(\xi)(\varpi)\qquad\for\;\;\xi\in o(2n+3,\mbb C),\;\varpi\in
\widehat M.\eqno(8.3.43)$$

Denote
$$\msr G=o(2n+3,\mbb C),\qquad \msr G_0=o(2n+1,\mbb C)+\mbb
CA_{n+1,n+1}\eqno(8.3.44)$$
$$\msr G_+=\mbb CK_{2n+2}+\sum_{i=1}^n(\mbb C A_{n+1,i}+\mbb
CB_{i,n+1}),\;\; \msr G_-=\mbb C K_{n+1}\sum_{i=1}^n(\mbb C
A_{i,n+1}+\mbb CC_{n+1,i})\eqno(8.3.45)$$ (cf. (8.3.27) and
(8.3.28)). Then $\msr G_\pm$ are abelian Lie subalgebras of $\msr G$
and $\msr G_0$ is a Lie subalgebra of $\msr G$. Moreover,
$$\msr G=\msr G_-+\msr G_0+\msr G_+,\qquad[\msr G_0,\msr
G_\pm]\subset \msr G_\pm.\eqno(8.3.46)$$ By (8.3.32)-(8.3.38),
$$\msr G_+(1\otimes M)=\{0\},\;\;U(\msr G_0)(1\otimes M)=1\otimes
M.\eqno(8.3.47)$$ Thus
$$U(\msr G)(1\otimes M)=U(\msr G_-)(1\otimes M).\eqno(8.3.48)$$
Using (8.3.34), we can prove:\psp

{\bf Proposition 8.3.3}. {\it The map $M\mapsto U({\msr
G}_-)(1\otimes M)$ gives rise to a functor from the category of
$o(2n+1,\mbb C)$-modules to the category of $o(2n+3,\mbb
C)$-modules. In particular, it maps irreducible $o(2n+1,\mbb
C)$-modules to irreducible $o(2n+3,\mbb C)$-modules.}\psp

From another point view, $U({\msr G}_-)(1\otimes M)$ is a polynomial
extension from $o(2n+1,\mbb C)$-module $M$ to an $o(2n+3,\mbb
C)$-module. Next we want to study when $\wht M=U({\msr
G}_-)(1\otimes M)$. Write
$$x^{\alpha}=\prod_{i=0}^{2n}x_{i}^{\alpha_i}
,\;\;E^\al=K_{n+1}^{\al_0}\prod_{r=1}^n[(-A_{r,n+1})^{\al_r}C_{n+1,r}^{\al_{n+r}}]
\eqno(8.3.49)$$for
$\al=(\al_0,\al_1,...,\al_{2n})\in\mbb{N}^{2n+1}$.
 For $k\in\mbb{N}$, we set
$${\msr A}_k=\mbox{Span}\{x^\al\mid
\al \in\mbb{N}^{2n+1}, \ |\al|=k\},\;\; \widehat M_k={\msr
A}_k\otimes_\mbb{C}M\eqno(8.3.50)$$ (recall
$|\al|=\sum_{i=0}^{2n}\al_i$) and
 \begin{eqnarray*}\hspace{1cm}(U(\msr G_-)(1\otimes
M))_k&=&U(\msr G_-)(1\otimes M)\bigcap \widehat M_k\\ &=&
\mbox{Span}\{ E^\al(1\otimes M)\mid  \al \in\mbb{N}^{2n+1}, \
|\al|=k\}\hspace{2.4cm}(8.3.51)\end{eqnarray*} by (8.3.36)-(8.3.38).
Moreover,
$$(U(\msr G_-)(1\otimes
M))_0=\widehat M_0=1\otimes M.\eqno(8.3.52)$$ Furthermore,
 $$\widehat M=\bigoplus_{k=0}^\infty\widehat M_k,\qquad
 U(\msr G_-)(1\otimes M)=\bigoplus\limits_{k=0}^\infty(U(\msr G_-)(1\otimes
M))_k.\eqno(8.3.53)$$

According to (8.3.23) and (8.3.25), we define a $\msr G_0$-module
homomorphism $\vf:\wht M\rta  U(\msr G_-)(1\otimes M)$ by
$$\vf(x^\al\otimes v)=E^\al(1\otimes
v)\qquad\for\;\;\al\in\mbb{N}^{2n+1},\;v\in M.\eqno(8.3.54)$$ Then
we have
$$\vf(\wht M_k)=(U(\msr G_-)(1\otimes
M))_k\qquad\for\;\;k\in\mbb{N}.\eqno(8.3.55)$$ Under the
identification (8.3.27) and (8.3.28), we have:\psp

{\bf Lemma 8.3.4}. {\it We have $\vf|_{\wht
M_1}=(c+\td\omega)|_{\wht M_1}$ (cf. (8.3.2)-(8.3.4)).}

 {\it Proof}.
Recall the identification (8.3.27) and (8.3.28). Moreover, $\widehat
M_1={\msr A}_1\otimes_\mbb CM$. Let $i\in \ol{1,n}$ and $v\in M$.
Expressions (8.3.38), (8.3.49) and (8.3.54) give
$$\vf(x_0\otimes v)=\sum_{s=1}^n[x_s\otimes
K_s(v)+x_{n+s}\otimes K_{n+1+s}(v)]+cx_0\otimes v.\eqno(8.3.56)$$
Moreover, (8.3.35), (8.3.49) and (8.3.54) imply
$$\vf(x_i\otimes v)=\sum_{p=1}^nx_{n+p}\otimes
B_{i,p}(v) -x_0\otimes K_{n+1+i}(v)+\sum_{q=1}^nx_q\otimes
A_{i,q}(v)+cx_i\otimes v \eqno(8.3.57)$$for $i\in\ol{1,n}$.
Furthermore, (8.3.36), (8.3.49) and (8.3.54) yield
\begin{eqnarray*}\qquad \vf(x_{n+i}\otimes v)&=&
-\sum_{p=1}^nx_{n+p}\otimes A_{p,i}(v)-x_0\otimes K_i(v)\\
& &+\sum_{q=1}^nx_q\otimes C_{i,q}(v)+cx_{n+i}\otimes
v\hspace{5.35cm}(8.3.58)\end{eqnarray*} for $i\in\ol{1,n}$.

On the other hand, (8.3.4), (8.3.15), (8.3.27) and (8.3.28) yield
$$\td\omega(x_0\otimes
v)=\sum_{i=1}^n[x_{n+i}\otimes K_{n+1+i}(v) +x_i\otimes
K_i(v)],\eqno(8.3.59)$$ $$\td\omega(x_i\otimes v)=
\sum_{p=1}^nx_{n+p}\otimes B_{i,p}(v)-x_0\otimes K_{n+1+i}(v)
+\sum_{r=1}^nx_r\otimes A_{i,r}(v),\eqno(8.3.60)$$
$$\td\omega(x_{n+i}\otimes v)=
\sum_{p=1}^nx_p\otimes C_{i,p}(v)-x_0\otimes K_i(v)
-\sum_{s=1}^nx_{n+s}\otimes A_{s,i}(v).\eqno(8.3.61)$$

 Comparing the above six
expressions, we get the conclusion in the lemma. $\qquad\Box$\psp

For $f\in{\msr A}$, we define the action
$$f(g\otimes v)=fg\otimes v\qquad\for\;\;g\in{\msr A},\;v\in
M.\eqno(8.3.62)$$ Then we have the $o(2n+1,\mbb{C})$-invariant
operator
$$T=\nu(K_{n+1})x_0+\sum_{i=1}^n[\nu(C_{n+1,i})x_i-\nu(A_{i,n+1})x_{n+i}]\eqno(8.3.63)$$
on $\wht M$. \psp

{\bf Lemma 8.3.5}. {\it We have
$T|_{\widehat{M}_k}=(2c-2n+k+2)\eta$}.

{\it Proof}. Let $f\in {\msr A}_k$ and  $v\in M$. According to
(8.3.38),
\begin{eqnarray*}\nu(K_{n+1})x_0(f\otimes v)&=&x_0\sum_{s=1}^n[x_sf\otimes K_s(v)+x_{n+s}f\otimes
K_{n+1+s}(v)]\\ & &+ [(k+1+c)x_0^2-\eta]f\otimes v-\eta
x_0\ptl_{x_0}(f)\otimes v.\hspace{2.7cm}(8.3.64)\end{eqnarray*}
Moreover, (8.3.35) gives
\begin{eqnarray*}\qquad& &-\nu(A_{i,n+1})x_{n+i}(f\otimes
v)\\&=&x_{n+i}\sum_{p=1}^n[ x_p\otimes A_{i,p}(v)-x_{n+p}\otimes
B_{p,i}(v)]-x_0x_{n+i}f\otimes K_{n+1+i}(v)\\ & &+
[(k+1+c)x_ix_{n+i}-\eta]f\otimes v-\eta
x_{n+i}\ptl_{x_{n+i}}(f)\otimes
v\hspace{3.7cm}(8.3.65)\end{eqnarray*} and (8.3.36) yields
\begin{eqnarray*}\qquad& &\nu(C_{n+1,i})x_i(f\otimes
v)\\&=&x_i\sum_{p=1}^n[ x_p\otimes C_{i,p}(v)-x_{n+p}\otimes
A_{p,i}(v)]-x_0x_if\otimes K_i(v)\\ & &+
[(k+1+c)x_ix_{n+i}-\eta]f\otimes v-\eta x_i\ptl_{x_i}(f)\otimes
v\hspace{4.5cm}(8.3.66)\end{eqnarray*} for $i\in\ol{1,n}$. Then the
lemma follows from (8.3.63)-(8.3.66) and  the skew-symmetry
$B_{i,j}=-B_{j,i}$ and $C_{i,j}=-C_{j,i}$.$\qquad\Box$\psp

 For $ 0\neq \mu=\sum_{i=1}^n\mu_i\ves_i\in\Lmd^+$ with ${\msr
S}(\mu)=\{i_1,i_2,...,i_{s+1}\}$ (cf. (6.7.54) and (6.7.55)), we
define
$$\Theta(\mu)=\left\{\begin{array}{ll}\emptyset
&\mbox{if}\;\mu=(\sum_{i=1}^n\ves_i)/2,
\\\mu_1+2n-i_2-\mbb{N}&\mbox{otherwise}.\end{array}\right.\eqno(8.3.67)$$
\pse

 {\bf Theorem 8.3.6}.  {\it For $0\neq\mu\in\Lmd^+$, the
 $o(2n+3,\mbb C)$-module $\widehat{V(\mu)}$
defined by (8.3.32)-(8.3.43) is irreducible if $c\in \mbb
C\setminus\{n-\mbb{N}/2,\Theta(\mu)\}.$}

{\it Proof}. By Lemma 8.3.3, it is enough to prove that the
homomorphism $\vf$ defined in (8.3.49) and (8.3.54) satisfies
$\vf(\widehat{V(\mu)})=\widehat{V(\mu)}$.  According to (8.3.55), we
only need to prove
$$\vf(\widehat{V(\mu)}_k)=\widehat{V(\mu)}_k\eqno(8.3.68)$$
for any $k\in\mbb{N}$. We will prove it by induction on $k$.

 When $k=0$, (8.3.68) holds by the definition (8.3.47). Consider
$k=1$. Write $\mu=\sum_{i=1}^n\mu_i\ves_i\in\Lmd^+$ with ${\msr
S}(\mu)=\{i_1,i_2,...,i_{s+1}\}$.  According to Lemma 8.3.2 and
Lemma 8.3.4 with $M=V(\mu)$, the eigenvalues of
$\vf|_{\widehat{V(\mu)}_{\la 1\ra}}$ are among
$$\{c-n,\;c+\mu_{i_r}-i_r+1,c-\mu_{i_r}-2n
+i_{r+1}\mid\for\;\;i\in\ol{1,s}\}.\eqno(8.3.69)$$
 Recall that $\mu_r\in\mbb{N}/2$ for
$r\in\ol{1,n}$,
$$\mu_{\iota+1}-\mu_\iota\in\mbb{N}\;\;\for\;\;\iota\in\ol{1,n-1}\eqno(8.3.70)$$
such that  (6.7.55) holds. So
$$-\mu_{i_r}+i_r-1,\;\mu_{i_r}+2n
-i_{r+1}\in
\mu_1+2n-i_2-\mbb{N}\;\;\for\;\;i\in\ol{1,s}.\eqno(8.3.71)$$ If
$c\not\in \mu_1+2n-n_1-\mbb{N}$ and $c\neq n$, then all the
eigenvalues of $\vf|_{\widehat{V(\mu)}_{\la 1\ra}}$ are nonzero.  In
the case $\mu=(\sum_{i=1}^n\ves_i)/2$, the eigenvalues
$\vf|_{\widehat{V(\mu)}_{\la 1\ra}}$ are $c-n$ and $c+1/2$, which
are not equal to 0 because of $c\not\in n-\mbb{N}/2$. Thus (8.3.68)
holds for $k=1$.

Suppose that (8.3.68) holds for $k\leq \ell$ with $\ell\geq 1$.
Consider $k=\ell+1$.  Note that
\begin{eqnarray*}\qquad\vf(\widehat{V(\mu)}_{
\ell+1})&=&\sum_{i=0}^{2n}\vf(x_i\widehat{V(\mu)}_\ell)
\\ &=&K_{n+1}[\vf(\widehat{V(\mu)}_
\ell)]+\sum_{i=1}^n[A_{i,n+1}[\vf(\widehat{V(\mu)}_
\ell)]+C_{n+1,i}[\vf(\widehat{V(\mu)}_ \ell)]]\\ &=&
K_{n+1}(\widehat{V(\mu)}_
\ell)+\sum_{i=1}^n[A_{i,n+1}(\widehat{V(\mu)}_
\ell)+C_{n+1,i}(\widehat{V(\mu)}_\ell)]\hspace{1.7cm}(8.3.72)\end{eqnarray*}
by the inductional assumption.  To prove (8.3.68) with $k=\ell+1$ is
equivalent to prove
$$K_{n+1}(\widehat{V(\mu)}_
\ell)+\sum_{i=1}^n[A_{i,n+1}(\widehat{V(\mu)}_\ell)+C_{n+1,i}(\widehat{V(\mu)}_\ell)]=\widehat{V(\mu)}_{
\ell+1}.\eqno(8.3.73)$$

For any $u\in \widehat{V(\mu)}_{ \ell-1}$, Lemma 8.3.5 says that
$$K_{n+1}(x_0u)+\sum_{i=1}^n[C_{n+1,i}(x_iu)-A_{i,n+1}(x_{n+i}u)]=(2c+2-2n+\ell)\eta
u.\eqno(8.3.74)$$ Since $c\not\in n-\mbb{N}/2$, $2c+2-2n+\ell\neq 0$
and  (8.3.73) gives
$$\eta u\in K_{n+1}(\widehat{V(\mu)}_
\ell)+\sum_{i=1}^n[A_{i,n+1}(\widehat{V(\mu)}_\ell)+C_{n+1,i}(\widehat{V(\mu)}_\ell)]\qquad\for\;\;u\in
\widehat{V(\mu)}_{ \ell-1}.\eqno(8.3.75)$$

Let $g\otimes v\in \widehat{V(\mu)}_{\la \ell\ra}$.
 According to
(8.3.35), (8.3.36), (8.3.38) and Lemma 8.3.4,
\begin{eqnarray*} -A_{i,n+1}(g\otimes v)&=&x_iD(g)\otimes
v-\eta\ptl_{x_{n+i}}(g)\otimes v+x_ig\otimes I_{2n+1}(v)-x_0g\otimes K_{n+1+i}(v)\\
& &+\sum_{p=1}^n(x_{n+p}g\otimes
B_{i,n+p}(v)+x_pg\otimes A_{i,p}(v))\\
&=&\ell x_ig\otimes
v-\eta\ptl_{x_{n+i}}(g)\otimes v+g[x_i\otimes I_{2n+1}(v)-x_0g\otimes K_{n+1+i}(v)\\
& &+\sum_{p=1}^n(x_{n+p}\otimes B_{i,n+p}(v)+x_p\otimes A_{i,p}(v))]
\\&=&\ell x_ig\otimes
v-\eta\ptl_{x_{n+i}}(g)\otimes v-g[A_{i,n+1}(1\otimes v)]
\\&=&\ell x_ig\otimes
v-\eta\ptl_{x_{n+i}}(g)\otimes v+g\vf(x_i\otimes v)
\\&=&-\eta\ptl_{x_{n+i}}(g)\otimes v +g[(\ell+c+\td\omega)(x_i\otimes
v)],\hspace{3.9cm}(8.3.76)\end{eqnarray*}
\begin{eqnarray*} C_{n+i,i}(g\otimes v)&=&x_{n+i}D(g)\otimes
v-\eta\ptl_{x_i}(g)\otimes v+x_{n+i}g\otimes I_{2n+1}(v)-x_0g\otimes K_i(v)\\
& &+\sum_{p=1}^n(x_pg\otimes C_{i,p}(v)-x_{n+p}g\otimes
A_{p,i}(v))\\&=&\ell x_{n+i}g\otimes
v-\eta\ptl_{x_i}(g)\otimes v+g[x_{n+i}\otimes I_{2n+1}(v)-x_0g\otimes K_i(v)\\
& &+\sum_{p=1}^n(x_p\otimes C_{i,p}(v)-x_{n+p}\otimes A_{p,i}(v))]
\\&=&\ell x_{n+i}g\otimes v-\eta\ptl_{x_i}(g)\otimes
v+g[C_{n+1,i}(1\otimes v)] \\&=&\ell x_{n+i}g\otimes
v-\eta\ptl_{x_i}(g)\otimes v+g\vf(x_{n+i}\otimes v)
\\&=&-\eta\ptl_{x_i}(g)\otimes v +g[(\ell+c+\td\omega)(x_{n+i}\otimes
v)]\hspace{4.2cm}(8.3.77)\end{eqnarray*} for $i\in\ol{1,n}$, and
\begin{eqnarray*} K_{n+1}(g\otimes v)&=&
x_0D(g)\otimes v-\eta\ptl_{x_0}(g)\otimes v+
x_0g\otimes I_{2n+1}(v)
\\
& &+\sum_{s=1}^n(x_sg\otimes K_s(v)+x_{n+s}g\otimes
K_{n+1+s}(v))\\&=&\ell x_0g\otimes v-\eta\ptl_{x_0}(g)\otimes v+
g[x_0\otimes I_{2n+1}(v)
\\
& &+\sum_{s=1}^n(x_s\otimes K_s(v)+x_{n+s}\otimes K_{n+1+s}(v))]
\\&=&\ell x_0g\otimes v-\eta\ptl_{x_0}(g)\otimes
v+g[K_{n+1}(1\otimes v)] \\&=&\ell x_0g\otimes
v-\eta\ptl_{x_0}(g)\otimes v+g\vf(x_0\otimes v)
\\&=&-\eta\ptl_{x_0}(g)\otimes v +g[(\ell+c+\td\omega)(x_0\otimes
v)].\hspace{4.5cm}(8.3.78)\end{eqnarray*}

 Since (8.3.75) says that
$$\eta\ptl_{x_r}(g)\otimes v
\in
K_{n+1}(\widehat{V(\mu)}_\ell)+\sum_{i=1}^n[A_{i,n+1}(\widehat{V(\mu)}_\ell)+C_{n+1,i}(\widehat{V(\mu)}_\ell)])\eqno(8.3.79)$$
for $r\in\ol{0,2n}$, Expressions (8.3.76)-(8.3.78) show
$$g[(\ell+c+\td\omega)(x_i\otimes
v)]\in
K_{n+1}(\widehat{V(\mu)}_\ell)+\sum_{i=1}^n[A_{i,n+1}(\widehat{V(\mu)}_\ell)+C_{n+1,i}(\widehat{V(\mu)}_\ell)]\eqno(8.3.80)$$
for $i\in\ol{0,2n}$ and $g\in{\msr A}_\ell$.

According to Lemmas 8.3.2, the eigenvalues of
$(\ell+c+\td\omega)|_{\widehat{V(\mu)}_{\la 1\ra}}$ are among
$$\ell+c+\mu_{i_r}-i_r+1,\;\ell+c-\mu_{i_r}-2n+i_{r+1}\;\;\for\;\;r\in\ol{1,s}.\eqno(8.3.81)$$
Again $$-\ell-\mu_{i_r}+i_r-1,\;-\ell+\mu_{i_r}+2n-i_{r+1}\in
\mu_1+2n-i_2-\mbb{N}\;\;\for\;\;i\in\ol{1,s}.\eqno(8.3.82)$$ If
$c\not\in \mu_1+2n-i_2-\mbb{N}$, then all the eigenvalues of
$(\ell+c+\td\omega)|_{\widehat{V(\mu)}_1}$ are nonzero.
 Hence
$$(\ell+c+\td\omega)(\widehat{V(\mu)}_1)=\widehat{V(\mu)}_1.\eqno(8.3.83)$$ By (8.3.80) and (8.3.83),
$$g(\widehat{V(\mu)}_{\la 1\ra})\subset \sum_{i=1}^n[A_{i,n+1}(\widehat{V(\mu)}_\ell)+C_{n+1,i}(\widehat{V(\mu)}_\ell)]
\qquad\for\;\;g\in{\msr A}_\ell,\eqno(8.3.84)$$ or equivalently,
(8.3.68) holds for $k=\ell+1$ by (8.3.84). By induction, (8.3.68)
holds for any $k\in\mbb{N}.\qquad\Box$\psp

Let $\lmd_i$ be the $i$th fundamental weight of $o(2n+3,\mbb C)$
with respect to the following simple positive roots
$$\{\ves_{n+1}-\ves_n,\ves_n-\ves_{n-1},...,\ves_2-\ves_1,\ves_1
\}.\eqno(8.3.85)$$ By the Weyl group's actions on $o(2n+1,\mbb C)$
and $V(\mu)$, $\wht{V(\mu)}$ is a highest-weight $o(2n+3,\mbb
C)$-module with highest weight
$-(c+\mu_1)\lmd_1+\sum_{r=1}^{s-1}(\mu_{i_r}-\mu_{i_{r+1}})\lmd_{i_{r+1}}+2\mu_n\lmd_{n+1}$.

 Up to this stage, we do not known if the
condition in Theorem 8.3.6 is necessary for the
$o(2n+3,\mbb{C})$-module $\widehat{V(\mu)}$ to be irreducible if
$\mu\neq 0$. We will deal with  the case $\mu=0$ in next section.

\section{ Conformal
Oscillator Representations}

In this section, we  study  the $o(2n+3,\mbb{C})$-module
$\widehat{V(0)}$ and its  oscillator generalizations.

In order to use the results in Sections 8.1 and 8.2, we redenote
$$y_i=x_{n+i}\qquad\for\;\;i\in\ol{1,n}\eqno(8.4.1)$$
and use ${\msr B}'=\mbb{C}[x_0,x_1,...,x_n,y_1,...,y_n]$ to replace
${\msr A}$ in last section. Recall
$$\Dlt'=\ptl_{x_0}^2+2\sum_{i=1}^n\ptl_{x_i}\ptl_{y_i},\;\;\eta'=x^2_0/2+\sum_{i=1}^nx_iy_i,
\eqno(8.4.2)$$
$$D'=\sum_{r=0}^nx_r\partial_{x_r}+\sum_{s=1}^ny_s\ptl_{y_s}.\eqno(8.4.3)$$
Fix $c\in\mbb{C}$ and identify $\widehat{V(0)}={\msr B}'\otimes v_0$
with $\msr B'$ by
$$f\otimes v_0\leftrightarrow f\qquad\for\;\;f\in{\msr
B}',\eqno(8.4.4)$$ where $V(0)=\mbb{C}v_0$. Then we have the
following one-parameter generalization $\pi_c$ of the conformal
representation of $o(2n+3,\mbb{C})$:
$$\pi_c(A_{i,j})=x_i\ptl_{x_j}-y_j\ptl_{x_i},\;\pi_c(B_{i,j})=x_i\ptl_{y_j}-x_j\ptl_{y_i},\;
\pi_c(C_{i,j})=y_i\ptl_{x_j}-y_j\ptl_{x_i},\eqno(8.4.5)$$
$$\pi_c(A_{n+1,i})=\ptl_{x_i},\;\;\pi_c(B_{i,n+1})=-\ptl_{y_i},\;\;\pi_c(K_{2n+2})=-\ptl_{x_0},\eqno(8.4.6)$$
$$\pi_c(A_{i,n+1})=\eta'\ptl_{y_i}-x_i(D'+c),\;\;\pi_c(C_{n+1,i})=y_i(D'+c)-\eta'\ptl_{x_i}\eqno(8.4.7)$$
$$\pi_c(K_i)=x_0\ptl_{x_i}-x_{n+i}\ptl_{x_0},\qquad
\pi_c(K_{n+1+i})=x_0\ptl_{x_{n+i}}-x_i\ptl_{x_0},\eqno(8.4.8)$$
$$\pi_c(A_{n+1,n+1})=-D'-c,,\qquad \pi_c(K_{n+1})=x_0(D'+c)-\eta'\ptl_{x_0}\eqno(8.4.9)$$
 for $i,j\in\ol{1,n}$.\psp

 {\bf Theorem 8.4.1}. {\it The representation $\pi_c$ of
$o(2n+2,\mbb{C})$ on $\msr B$ is irreducible if and only if
$c\not\in-\mbb{N}$. When $c=0$, the  representation $\pi_c$ of
$o(2n+3,\mbb{C})$ on $\msr B$ is  the natural conformal
representation of $o(2n+3,\mbb{C})$ given in (8.3.20)-(8.3.25) in
terms (8.4.1) with $\msr A$ replaced by $\msr B'$. The subspace
$\mbb{C}1_{\msr B'}$ forms a trivial $o(2n+3,\mbb{C})$-submodule of
the conformal module $\msr B'$ and the quotient space $\msr
B'/\mbb{C}1_{\msr B'}$ forms an irreducible
$o(2n+3,\mbb{C})$-module.}

{\it Proof}. Recall the identification (8.3.27) and (8.3.28) .
Moreover, $\msr B_k'$ denotes the subspace of homogeneous
polynomials in $\msr B'$ with degree $k$ and
$${\msr H}_k'=\{f\in\msr
B'_k\mid\Dlt'(f)=0\}\qquad\for\;\;k\in\mbb{N}.\eqno(8.4.10)$$
According to Theorem 8.1.1,  ${\msr H}_k'$ is an irreducible
$o(2n+1,\mbb{C})$-submodule with the highest-weight vector $x_1^k$
of weight $k\ves_1$ for any $k\in\mbb N$, and
$$\msr B'=\bigoplus_{m,k=0}^\infty {\eta'}^m{\msr H}_k'\eqno(8.4.11)$$ is a direct sum of irreducible
$o(2n+1,\mbb{C})$-submodules. On the other hand, $U(\msr
G_-)(1_{\msr B'})$ forms an irreducible $o(2n+3,\mbb{C})$-submodule
of $\msr B'$ due to (8.4.6). If $c=-\ell$ for some $\ell\in\mbb N$,
then
$$U(\msr G_-)(1_{\msr B'})\bigcap \msr B'_{\ell+1}\subset \eta'\msr
B_{\ell-1}'\eqno(8.4.12)$$ by (8.4.7) and (8.4.9). So
$c\not\in-\mbb{N}$ is a necessary condition for the representation
$\pi_c$ of $o(2n+3,\mbb{C})$ on $\msr B'$ to be irreducible.

Next we assume $c\not\in -\mbb{N}-1$.  Let $M$ be a nonzero
$o(2n+3,\mbb{C})$-submodule of $\msr B'$ such that
$$M\not\subset
\mbb{C}1_{\msr B'}\;\;\mbox{if}\;\; c= 0.\eqno(8.4.13)$$ Repeatedly
applying (8.4.6) if necessary, we have $1_{\msr B'}\in M$. Note
$$K_{n+1}(1_{\msr B'})=cx_0,\;\;-A_{i,n+1}(1_{\msr B'})=cx_i,\;\;C_{n+1,i}(1_{\msr B'})=cy_i\qquad\for\;\;i\in\ol{1,n}.\eqno(8.4.14)$$ Thus
$$\msr B_1'\subset M\eqno(8.4.15)$$ if $c\neq 0$.
 When $c=0$, (8.4.15) also holds by
 (8.4.6) and the fact that $\msr B_1'$ is an irreducible $o(2n+1,\mbb{C})$-submodule. Suppose that
$$\msr B_k'\subset M\qquad\for\;\;k<\ell\eqno(8.4.16)$$
with $2\leq \ell\in\mbb{N}$. According to (8.4.11),
$$\msr B_\ell'=\bigoplus_{m=0}^{\llbracket \ell/2
\rrbracket}{\eta'}^m{\msr H}_{\ell-2m}'.\eqno(8.4.17)$$
 Moreover,
$$[\Dlt',\eta']=1+2n+2D'.\eqno(8.4.18)$$Set
$$\msr B_ {\ell,r}'=\bigoplus_{m=0}^r{\eta'}^m{\msr H}_{\ell-2m}'\eqno(8.4.19)$$
for $r\in\ol{0,\llbracket \ell/2 \rrbracket}$.  Then
$$\msr B_ {\ell,r}'=\{w\in \msr B_ \ell'\mid {\Dlt'}^{r+1}(w)=0\}\eqno(8.4.20)$$ and
$${\Dlt'}^r(\msr B'_ {\ell,r})={\msr H}_{\ell-2r}'\eqno(8.4.21)$$ by (8.4.19).

Since
$$-A_{1,n+1}(x_1^{\ell-1})=(c+\ell-1)x_1^\ell\in
M\eqno(8.4.22)$$ and $c\not\in-\mbb{N}-1$, we have
$$x_1^\ell\in M\eqno(8.4.23)$$ by (8.4.7).
Hence
$${\msr H}_\ell'\subset
M\eqno(8.4.24)$$ because ${\msr H}_\ell'$ is an irreducible
$o(2n+1,\mbb{C})$-submodule generated by $x_1^\ell$. Furthermore,
$$K_{n+1}(x_1^{\ell-r-1}y_1^r)=(c+\ell-1)x_1^{\ell-r-1}x_0y_1^r\in
M. \eqno(8.4.25)$$
  So
$x_1^{\ell-r-1}x_0y_1^r\in M$. Moreover,
$${\Dlt'}^r(x_1^{\ell-r-1}x_0y_1^r)=2^rr![\prod_{s=1}^r(\ell-r-s)]x_1^{\ell-2r-1}x_0 \in
{\msr H}'_{\ell-2r}.\eqno(8.4.26)$$ Observe that (8.4.17) is a
direct sum of irreducible $o(2n+3,\mbb{C})$-submodules with distinct
highest weights. So
$$\msr B'_\ell\bigcap M=\bigoplus_{m=0}^{\llbracket \ell/2
\rrbracket}({\eta'}^m{\msr H}'_{\ell-2m})\bigcap M.\eqno(8.4.27)$$
By (8.4.20)-(8.4.26), $({\eta'}^r{\msr H}'_{\ell-2r})\bigcap M$ is a
nonzero $o(2n+1,\mbb{C})$-submodule. Since ${\eta'}^r{
H}_{\ell-2r}'$ is an irreducible $o(2n+1,\mbb{C})$-module, we have
${\eta'}^r{\msr H}'_{\ell-2r}=({\eta'}^r{\msr H}'_{\ell-2r})\bigcap
M$. Therefore, $\msr B'_\ell\subset M$. By induction, $\msr
B'_k\subset M$ for any $k\in\mbb{N}$, or equivalently, $M=\msr B'$.
This proves the theorem. $\qquad\Box$\psp

With respect to (8.3.85), $\msr B'$ is a highest-weight irreducible
$o(2n+3,\mbb C)$-module with highest weight $-c\lmd_1$ when
$c\not\in -\mbb N$. If $c=0$, then  $\msr B'/\mbb C1_{\msr B'}$ is a
highest-weight irreducible $o(2n+3,\mbb C)$-module with highest
weight $-2\lmd_1+\lmd_2$.

For $\vec a=(a_0,a_1,...,a_n)^t\in\mbb C^{n+1}$ and $\vec
b=(b_1,b_2,...,b_n)^t\in\mbb{C}^n$, we put
$$\vec a\cdot\vec
x=\sum_{i=0}^na_ix_i,\qquad\vec b\cdot\vec
y=\sum_{i=1}^nb_iy_i.\eqno(8.4.28)$$ Set
$${\msr B}'_{\vec a,\vec b}=\{fe^{\vec a\cdot\vec
x+\vec b\cdot\vec y}\mid f\in{\msr B'}\}.\eqno(8.4.29)$$
 Denote by $\pi_{c,\vec a,\vec b}$ the
representation $\pi_c$ of $o(2n+3,\mbb C)$ on $\msr B'_{\vec a,\vec
b}$.\psp

{\bf Theorem 8.4.2}. {\it If $a_0^2+2\sum_{i=1}^na_ib_i\neq 0$, then
representation $\pi_{c,\vec a,\vec b}$ of $o(2n+3,\mbb{C})$ is
irreducible  for any $c\in\mbb{C}$.}

{\it Proof}.  Set
$$\msr B'_{\vec a,\vec b,k}=\msr B'_ke^{\vec a\cdot\vec
x+\vec b\cdot\vec y}\qquad\for\;k\in\mbb{N}.\eqno(8.4.30)$$  Let
${M}$ be a nonzero $o(2n+3,\mbb{C})$-submodule of $\msr B'_{\vec
a,\vec b}$. Take any $0\neq fe^{\vec a\cdot\vec x+\vec b\cdot\vec
y}\in  M$ with $f\in \msr B$.  According to (8.4.6),
$$(A_{n+1,i}-a_i)(fe^{\vec a\cdot\vec x+\vec b\cdot\vec
y})=\ptl_{x_i}(f)e^{\vec a\cdot\vec x+\vec b\cdot\vec
y},\;\;-(B_{i,n+1}+b_i)(fe^{\vec a\cdot\vec x+\vec b\cdot\vec
y})=\ptl_{y_i}(f)e^{\vec a\cdot\vec x+\vec b\cdot\vec
y}\eqno(8.4.31)$$ for $i\in\ol{1,n}$, and
$$-(K_{2n+2}+a_0)(fe^{\vec a\cdot\vec x+\vec b\cdot\vec
y})=\ptl_{x_0}(f)e^{\vec a\cdot\vec x+\vec b\cdot\vec
y}.\eqno(8.4.32)$$
 Repeatedly applying (8.4.31) and (8.4.32), we obtain $e^{\vec
a\cdot\vec x+\vec b\cdot\vec y}\in  M$. Equivalently, $\msr B'_{\vec
a,\vec b,0}\subset M$. Suppose $\msr B'_{\vec a,\vec b,\ell}\subset
M$ for some $\ell\in\mbb{N}$. Let $ ge^{\vec a\cdot\vec x+\vec
b\cdot\vec y}$ be any element in $\msr B'_{\vec a,\vec b,\ell}$.\pse

{\it Case 1. $a_i\neq 0$ or $b_i\neq 0$ for some
$i\in\ol{1,n}$.}\pse

 By symmetry, we may assume $a_1\neq 0$. Since
$$(x_i\ptl_{x_1}-y_1\ptl_{y_i})(g)e^{\vec a\cdot\vec x+\vec b\cdot\vec
y},\;(y_i\ptl_{x_1}-y_1\ptl_{x_i})(g)e^{\vec a\cdot\vec x+\vec
b\cdot\vec y}\in  \msr B'_{\vec a,\vec b,\ell}\subset
M\eqno(8.4.33)$$ by inductional assumption,
 Expression (8.4.5) yields
$$A_{i,1}(ge^{\vec a\cdot\vec x+\vec b\cdot\vec y})\equiv (a_1x_i-b_iy_1)ge^{\vec a\cdot\vec x+\vec b\cdot\vec y}
\equiv 0\;\;(\mbox{mod}\;M)\eqno(8.4.34)$$ and
$$C_{i,1}(ge^{\vec a\cdot\vec x+\vec b\cdot\vec y})\equiv
(a_1y_i-a_iy_1)ge^{\vec a\cdot\vec x+\vec b\cdot\vec y} \equiv
0\;\;(\mbox{mod}\;M).\eqno(8.4.35)$$Moreover, the first equation in
(8.4.8) gives
$$K_1(ge^{\vec a\cdot\vec x+\vec b\cdot\vec y})\equiv
(a_1x_0-a_0y_1)ge^{\vec a\cdot\vec x+\vec b\cdot\vec y} \equiv
0\;\;(\mbox{mod}\;M)\eqno(8.4.36)$$ because
$$(x_0\ptl_{x_1}-y_1\ptl_{x_0})(g)e^{\vec a\cdot\vec x+\vec b\cdot\vec
y}\in\msr B'_{\vec a,\vec b,\ell}\subset M.\eqno(8.4.37)$$

 On the other hand,
$$(D'+c)(g)e^{\vec a\cdot\vec x+\vec b\cdot\vec
y}=(\ell+c)ge^{\vec a\cdot\vec x+\vec b\cdot\vec y}\in\msr B'_{\vec
a,\vec b,\ell}\subset M\eqno(8.4.38)$$ and the first equation in
(8.4.9) implies
$$-A_{n+1,n+1}(ge^{\vec a\cdot\vec x+\vec b\cdot\vec y})\equiv
[a_0x_0+\sum_{i=1}^n(a_ix_i+b_iy_i)]ge^{\vec a\cdot\vec x+\vec
b\cdot\vec y} \equiv 0\;\;(\mbox{mod}\;M).\eqno(8.4.39)$$
 Substituting (8.4.34)-(8.4.36) into (8.4.39), we get
$$(a_0^2+2\sum_{i=1}^na_ib_i)y_1ge^{\vec
a\cdot\vec x+\vec b\cdot\vec y} \equiv 0\;\;(\mbox{mod}\;
M).\eqno(8.4.40)$$ Equivalently, $y_1ge^{\vec a\cdot\vec x+\vec
b\cdot\vec y}\in M$. Substituting it to (8.4.34)-(8.4.36), we obtain
$$x_0ge^{\vec a\cdot\vec x+\vec
b\cdot\vec y},x_ige^{\vec a\cdot\vec x+\vec b\cdot\vec
y},y_ige^{\vec a\cdot\vec x+\vec b\cdot\vec y}\in M\eqno(8.4.41)$$
for $i\in\ol{1,n}$. Therefore, $\msr B'_{\vec a,\vec
b,\ell+1}\subset M$. By induction, $\msr B'_{\vec a,\vec b,
\ell}\subset M$ for any $\ell\in\mbb{N}$. So $\msr B'_{\vec a,\vec
b}= M$. Hence $\msr B'_{\vec a,\vec b}$ is an irreducible
$o(2n+3,\mbb{C})$-module.\pse

{\it Case 2. $a_0\neq 0$ and $a_i=b_i=0$ for $i\in\ol{1,n}$.}\pse

Under the above assumption,
$$K_i(ge^{\vec a\cdot\vec x+\vec b\cdot\vec
y})=(x_0\ptl_{x_i}-y_i\ptl_{x_0}-a_0y_i)(g)e^{\vec a\cdot\vec x+\vec
b\cdot\vec y}\in M\eqno(8.4.42)$$ and
$$K_{n+1+i}(ge^{\vec a\cdot\vec x+\vec b\cdot\vec
y})=(x_0\ptl_{y_i}-x_i\ptl_{x_0}-a_0x_i)(g)e^{\vec a\cdot\vec x+\vec
b\cdot\vec y}\in M\eqno(8.4.43)$$ for $i\in\ol{1,n}$. Note
$$(x_0\ptl_{x_i}-y_i\ptl_{x_0})(g)e^{\vec a\cdot\vec x+\vec
b\cdot\vec y}, (x_0\ptl_{y_i}-x_i\ptl_{x_0})(g)e^{\vec a\cdot\vec
x+\vec b\cdot\vec y}\in\msr B'_{\vec a,\vec b,\ell}\subset M
\eqno(8.3.44)$$ by the inductional assumption. Thus (8.4.42) and
(8.4.43) imply
$$y_ige^{\vec a\cdot\vec x+\vec
b\cdot\vec y},x_ige^{\vec a\cdot\vec x+\vec b\cdot\vec y}\in
M\qquad\for\;\;i\in\ol{1,n}.\eqno(8.4.45)$$ Now (8.4.39) yields
$x_0ge^{\vec a\cdot\vec x+\vec b\cdot\vec y}\in M$. So $\msr
B'_{\vec a,\vec b,\ell+1}\subset M$. By induction, $\msr B'=M$; that
is, $\msr B'$ is irreducible. $\qquad\Box$\psp

Fix $n_1,n_2\in\ol{1,n}$ with $n_1\leq n_2$. Set
$$\td D'=x_0\ptl_{x_0}-\sum_{i=1}^{n_1}x_i\ptl_{x_i}
+\sum_{r=n_1+1}^nx_r\ptl_{x_r}+\sum_{j=1}^{n_2}y_j\ptl_{y_j}-\sum_{s=n_2+1}^ny_s\ptl_{y_s},\eqno(8.4.46)$$
$$\td\Dlt'=\ptl_{x_0}^2-2\sum_{i=1}^{n_1}x_i\ptl_{y_i}+2\sum_{r=n_1+1}^{n_2}\ptl_{x_r}\ptl_{y_r}-2\sum_{s=n_2+1}^n
y_s\ptl_{x_s}\eqno(8.4.47)$$ and
$$\td\eta'=\frac{x_0^2}{2}+\sum_{i=1}^{n_1}y_i\ptl_{x_i}+\sum_{r=n_1+1}^{n_2}x_ry_r+\sum_{s=n_2+1}^n
x_s\ptl_{y_s}.\eqno(8.4.48)$$ Then the representation
$\pi_c^{n_1,n_2}$ of $o(2n+3,\mbb C)$ is determined as follows:
$\pi_c|_{o(2n+2,\mbb C)}$ is given by (7.4.38)-(7.4.47) with $\wht
D$ replaced by $\td D'$ and $\td\eta$ by $\td\eta'$, and
$$\pi_c^{n_1,n_2}(K_i)=\left\{\begin{array}{ll}-x_0x_i-y_i\ptl_{x_0}&\mbox{if}\;i\in\ol{1,n_1},\\
x_0\ptl_{x_i}-y_i\ptl_{x_0}&\mbox{if}\;i\in\ol{n_1+1,n_2},\\
x_0\ptl_{x_i}-\ptl_{x_0}\ptl_{y_i}&\mbox{if}\;i\in\ol{n_2+1,n},\end{array}\right.\eqno(8.4.49)$$
$$\pi_c^{n_1,n_2}(K_{n+1+i})=\left\{\begin{array}{ll}x_0\ptl_{y_i}-\ptl_{x_0}\ptl_{x_i}&\mbox{if}\;i\in\ol{1,n_1},\\
x_0\ptl_{y_i}-x_i\ptl_{x_0}&\mbox{if}\;i\in\ol{n_1+1,n_2},\\
-x_0y_i-x_i\ptl_{x_0}&\mbox{if}\;i\in\ol{n_2+1,n},\end{array}\right.\eqno(8.4.50)$$
$$\pi_c^{n_1,n_2}(K_{n+1})=x_0(\td D'+\td c)-\td\eta'\ptl_{x_0},\qquad
\pi_c^{n_1,n_2}(K_{2n+2})=-\ptl_{x_0}, \eqno(8.4.51)$$ where
$$\td c=c+n_2-n_1-n.\eqno(8.4.52)$$

Recall (8.2.10) and (8.2.12). We have the following third result in
this section:\psp

{\bf Theorem 8.4.3}. {\it The representation $\pi_c^{n_1,n_2}$ of
$o(2n+3,\mbb{C})$ on $\msr A$ is irreducible if $c\not\in
\mbb{Z}/2$.}

{\it Proof}. Let $M$ be a nonzero $o(2n+3,\mbb{C})$-submodule of
$\msr B$. Note
$$\pi_c^{n_1,n_2}(A_{n+1,n+1})=-\td D'-\td c.\eqno(8.4.53)$$
According to (8.2.12),
$$M=\bigoplus_{k\in\mbb{Z}}\msr B_{\la k\ra}\bigcap
M.\eqno(8.4.54)$$ Thus $\msr B_{\la k\ra}\bigcap M\neq\{0\}$ for
some $k\in \mbb{Z}$. On the other hand,
  $\msr B_{\la k\ra}'=\bigoplus_{i=0}^\infty(\td\eta')^i(\msr H'_{\la
k-2i\ra})$ is a decomposition of irreducible $o(2n+1,\mbb
C)$-submodules by Theorem 8.2.3 and (8.3.26). Moreover, the weights
of the $sl(n,\mbb C)$-singular vectors in $(\td\eta')^i(\msr H'_{\la
k-2i\ra})$ are distinct by Sections 6.4, 6.5 and 6.6. Hence
$$(\td\eta')^i(\msr H'_{\la k-2i\ra})\subset M\;\;\mbox{for some}\;\;i\in\mbb{N}.\eqno(8.4.55)$$
Theorem 8.2.1 and the arguments in (7.4.52)-(7.4.57) show that we
can assume
$$\msr H'_{\la k-s\ra}\subset M\qquad\for\;\;s\in\mbb{N}.\eqno(8.4.56)$$

Suppose $(\td\eta')^s(x_1^{-r+s})\in M$ for any $r\in k-\mbb N$ and
some $s\in\mbb N$. Then,
$$K_{n+1}((\td\eta')^s(x_1^{-r+s+1}))=(r-1+\td c)(\td\eta')^s(x_0x_1^{-r+s+1})\in
M\eqno(8.4.57)$$ by (8.4.48) and the first equation in (8.4.51),
which implies $\eta_{n_1,n_2}^s(x_0x_1^{-r+s+1})\in M.$ Moreover,
$$K_{n+1}((\td\eta')^s(x_0x_1^{-r+s+1}))=-(\td\eta')^{s+1}(x_1^{-r+s+1})+(r+\td c)
(\td\eta')^s(x_0^2x_1^{-r+s+1}).\eqno(8.4.58)$$ Now (7.4.59) and
(7.4.61) with $\td\eta$ replaced by $\td\eta'$,  and (8.4.58) lead
to
$$ (1/2+s+\td
c-n_2+n_1)(\td\eta')^{s+1}(x_1^{-r+s+1})\in M\Rightarrow
\eta_{n_1,n_2}^{s+1}(x_1^{-r+s+1})\in M.\eqno(8.4.59)$$ By
induction,
$$(\td\eta')^\ell(x_1^{-r+\ell})\in M\qquad\for\;\;\ell\in\mbb
N,\;r\in k-\mbb N.\eqno(8.4.60)$$  According to Theorem 8.2.1,
$$\msr B'_{\la
m\ra}=\bigoplus_{\ell=0}^\infty(\td\eta')^\ell(\msr H'_{\la
m-2\ell\ra})\subset M\qquad\for\;\;m\in k-\mbb N.\eqno(8.4.61)$$

Observe that
$$\pi_c^{n_1,n_2}(K_{n+1})x_0=x_0^2(D_{n_1,n_2}+\td c+1)-\eta_{n_1,n_2}(x_0\ptl_{x_0}+1)\eqno(8.4.62)$$
by (8.4.51). Then (8.4.62) and (7.4.66)-(7.4.69) with $\td\eta$
replaced by $\td\eta'$ and $\wht D$ replaced by $\td D'$ yield
\begin{eqnarray*}\hspace{2cm}& &-\pi_c^{n_1,n_2}(K_{n+1})x_0+\sum_{i=1}^{n_2}
\pi_c^{n_1,n_2}(A_{i,n+1})y_i+\sum_{j=n_2+1}^n\pi_c^{n_1,n_2}(A_{j,n+1})\ptl_{y_j}
\\
&&-\sum_{r=1}^{n_1}\pi_c^{n_1,n_2}(C_{n+1,r})\ptl_{x_r}-\sum_{s=n_1+1}^n\pi_c^{n_1,n_2}(C_{n+1,s})x_s\\
&=&\td\eta'(1-\td D'+n_2+n-n_1-2(\td
c+1))\hspace{5.2cm}(8.4.63)\end{eqnarray*}
 as operators on $\msr
B'$. The arguments in (7.4.71)-(7.4.77) show $M=\msr B'$; that is,
$\msr B'$ is an irreducible $o(2n+3,\mbb C)$-module. $\qquad\Box$
\psp

{\bf Remark 8.4.4}. The above irreducible representation depends on
the three parameters $c\in \mbb{C}$ and $n_1,n_2\in\ol{1,n}$. It is
not highest-weight type because of the mixture of multiplication
operators and differential operators in
 (7.4.39))-(7.4.42),
(7.4.44)-(7.4.47), (8.4.49) and (8.4.50) . Since $\msr B'$ is not
completely reducible as a $sl(n,\mbb C)$-module by Sections 6.4-6.6
when $n\geq 2$ and $n_1<n$, $\msr B$ is not a unitary
$o(2n+3,\mbb{C})$-module. Expression (8.4.53) shows that $\msr B'$
is a weight $o(2n+3,\mbb{C})$-module with finite-dimensional weight
subspaces.

\chapter{Representations of Symplectic Lie Algebras}

In this chapter, we study  the natural explicit representations of
symplectic Lie algebras. First in Section 9.1, we determine the
structure of the canonical bosonic and fermionic oscillator
representations  of the  Lie algebras over their minimal natural
modules. Moreover, the full quadratic operator oscillator
representations are presented. In Section 9.2, we study the
noncanonical oscillator representations obtained from the above
bosonic representations by partially swapping differential operators
and multiplication operators, and obtain a two-parameter family of
new infinite-dimensional irreducible representations. The results
were due to Luo and  the author [LX3].  In Section 9.3, we determine
the structures of the projective oscillator representations in
Section 6.8 restricted on symplectic Lie algebras, which were due to
our work [X25].

\section{Canonical Oscillator Representations}

In this section, we present the canonical bosonic and fermionic
oscillator representations of symplectic Lie algebras over their
minimal natural modules.

 Let $n>1$ be an integer. The symplectic Lie algebra
\begin{eqnarray*}\hspace{1cm}sp(2n,\mbb F)&=&
\sum_{i,j=1}^n\mbb{F}(E_{i,j}-E_{n+j,n+i})+\sum_{i=1}^n(\mbb FE_{i,n+i}+\mbb FE_{n+i,i})\\
& &+\sum_{1\leq i<j\leq n }[\mbb F(E_{i,n+j}+E_{n+j,i})+\mbb
F(E_{n+i,j}+E_{n+j,i})].\hspace{2.1cm}(9.1.1)\end{eqnarray*} We
again take
$$H=\sum_{i=1}^n\mbb{F}(E_{i,i}-E_{n+i,n+i})\eqno(9.1.2)$$ as a
Cartan subalgebra of $sp(2n,\mbb F)$ and the settings in
(7.1.3)-(7.1.5). Then the root system of $sp(2n,\mbb{F})$ is
$$\Phi_{C_n}=\{\pm \ves_i\pm\ves_j,\pm 2\ves_r\mid1\leq i<j\leq
n;r\in\ol{1,n}\}.\eqno(9.1.3)$$ We choose the set of positive roots
$$\Phi_{C_n}^+=\{\ves_i\pm\ves_j,2\ves_r\mid1\leq i<j\leq
n;r\in\ol{1,n}\}.\eqno(9.1.4)$$ In particular,
$$\Pi_{C_n}=\{\ves_1-\ves_2,...,\ves_{n-1}-\ves_n,2\ves_n\}\;\;
\mbox{is the set of positive simple roots}.\eqno(9.1.5)$$

Recall the set of dominant integral weights
$$\Lmd^+=\{\mu\in L_\mbb{Q}\mid
(\ves_n,\mu),(\ves_i-\ves_{i+1},\mu)\in\mbb{N}\;\for\;i\in\ol{1,n-1}.\eqno(9.1.6)$$
According to (7.1.5),
$$\Lmd^+=\{\mu=\sum_{i=1}^n\mu_i\ves_i\mid
\mu_i\in\mbb{N};\mu_i-\mu_{i+1}\in\mbb{N}\;\for\;i\in\ol{1,n-1}\}.\eqno(9.1.7)$$
For $\lmd\in\Lmd^+$, we denote by $V(\lmd)$ the finite-dimensional
irreducible $sp(2n,\mbb{F})$-module with highest weight $\lmd$.

Let ${\msr B}=\mbb{F}[x_1,...,x_n,y_1,...,y_n]$. Recall
(6.2.22)-(6.2.42). The canonical oscillator representation of
$sp(2n,\mbb{F})$ is given by
$$(E_{i,j}-E_{n+j,n+i})|_{\msr B}=x_i\ptl_{x_j}-y_j\ptl_{x_i},\;(E_{n+i,j}+E_{n+j,i})|_{\msr B}
=y_i\ptl_{x_j}+y_j\ptl_{x_i},\eqno(9.1.8)$$
$$(E_{i,n+j}+E_{j,n+i})|_{\msr B}
=x_i\ptl_{y_j}+x_j\ptl_{y_i}\eqno(9.1.9)$$ for $i,j\in\ol{1,n}$.

The positive root vectors of $sp(2n,\mbb{F})$ are
$$\{E_{i,j}-E_{n+j,n+i},E_{i,n+j}+E_{j,n+i},E_{r,n+r}\mid 1\leq i<j\leq
n;r\in\ol{1,n}\}.\eqno(9.1.10)$$ Recall the notion ${\msr
B}_{\ell_1,\ell_2}$ in (6.2.30) and ${\msr B}_k$ in (7.1.15). \psp

{\bf Theorem 9.1.1}. {\it For $k\in\mbb{N}$, the subspace ${\msr
B}_k$ forms a finite-dimensional irreducible
$sp(2n,\mbb{F})$-submodule and $x_1^k$ is a highest-weight vector
with the weight $k\ves_1$.}

{\it Proof}.  Note that $sp(2n,\mbb{F})$ contains the Lie subalgebra
$$\sum_{i,j=1}^n\mbb{F}(E_{i,j}-E_{n+j,n+i})\cong
gl(n,\mbb{F})\eqno(9.1.11)$$ and the representation of
$sp(2n,\mbb{F})$ in (9.1.8) and (9.1.9) is essentially an extension
of that for $sl(n,\mbb{F})$ given in (6.2.22). Thus the only
possible highest-weight vectors for $sp(2n,\mbb{F})$ in ${\msr B}_k$
are:
$$\{\eta^ix_1^\ell y_n^{k-\ell-2i}\mid i\in \ol{1,\llbracket k/2
\rrbracket}; \ell\in\ol{0,k-2i}\}\eqno(9.1.12)$$ by Lemma 6.2.3 and
(9.1.10), where $\eta=\sum_{i=1}^nx_iy_i$. Note
\begin{eqnarray*}\hspace{2cm}0&=&E_{n,2n}(\eta^ix_1^\ell
y_n^{k-\ell-2i})=x_n\ptl_{y_n}(\eta^ix_1^\ell y_n^{k-\ell-2i})
\\ &=&ix_n^2\eta^{i-1}x_1^\ell y_n^{k-\ell-2i}+(k-\ell-2i)x_n\eta^ix_1^\ell
y_n^{k-\ell-2i-1},\hspace{3cm}(9.1.13)\end{eqnarray*}
 which forces $i=0$ and $\ell=k$. Thus ${\msr B}_k$
has a unique (up to a scalar multiple) singular vector $x_1^k$ of
weight $k\ves_1$. By Weyl's Theorem 2.3.6 of complete reducibility
if $\mbb F=\mbb C$ or more generally by Lemma 6.3.2 with $n_1=0$, it
is irreducible. $\qquad\Box$ \psp

Consider the exterior algebra $\check{\msr A}$ generated by
$\{\sta_1,...,\sta_n,\vt_1,...,\vt_n\}$ (cf. (6.2.15)) and take the
settings in (6.2.44)-(6.2.51). Define a representation of
$sp(2n,\mbb{F})$ on $\check{\msr A}$ by
$$(E_{i,j}-E_{n+j,n+i})|_{\check{\msr
A}}=\sta_i\ptl_{\sta_j}-\vt_j\ptl_{\vt_i},\;\;(E_{n+i,j}+E_{n+j,i})|_{\check{\msr
A}} =\vt_i\ptl_{\sta_j}+\vt_j\ptl_{\sta_i},\eqno(9.1.14)$$
$$(E_{i,n+j}+E_{j,n+i})|_{\check{\msr
A}} =\sta_i\ptl_{\vt_j}+\sta_j\ptl_{\vt_i}\eqno(9.1.15)$$ for
$i,j\in\ol{1,n}$. By (9.1.11), (9.1.14) and (9.1.15), the above
representation of $sp(2n,\mbb{F})$  is essentially an extension of
that for $sl(n,\mbb{F})$ given in (6.2.43). Recall
$$\check\Dlt=\sum_{r=1}^n\ptl_{\sta_r}\ptl_{\vt_r},\qquad
\check\eta=\sum_{r=1}^n\sta_r\vt_r.\eqno(9.1.16)$$ It can be
verified that
$$\xi\check\Dlt=\check\Dlt\xi,\qquad
\xi\check\eta=\check\eta\xi\qquad\for\;\;\xi\in
sp(2n,\mbb{F})\eqno(9.1.17)$$ as operators on $\check{\msr A}$.

 Recall the notation
$\check{\msr A}_k$ in (7.1.24). We denote
$$\check{\msr H}_k=\{u\in \check{\msr A}_k\mid
\check\Dlt(u)=0\}.\eqno(9.1.18)$$ Then $\check{\msr H}_k$ forms an
$sp(2n,\mbb{F})$-submodule. According to (6.2.62),
$$\check{\msr H}_k=\sum_{\ell=0}^{\min\{k,n\}}\check{\msr
H}_{\ell,k-\ell}\qquad\for\;\;k\in\ol{0,2n}.\eqno(9.1.19)$$
Expressions (6.2.53) and (6.2.60) tell us that
 the $sp(2n,\mbb{F})$-singular
vectors of $\check{\msr H}_k$ are in the set
$$\{\vec\sta_r\vec\vt_s\mid 0\leq r<s\leq
n+1;r+n-s+1=k\}\eqno(9.1.20)$$ (cf. (6.2.50) and (6..2.51)). If
$s<n+1$, then
$$E_{s,n+s}(\vec\sta_r\vec\vt_s)=\sta_s\ptl_{\vt_s}(\vec\sta_r\vec\vt_s)=\vec\sta_r\vec\vt_{s+1}\sta_s\neq0.\eqno(9.1.21)$$
Thus $\check{\msr H}_k$ has an unique (up to a scalar multiple)
singular vector $\vec\sta_k$ of weight $\sum_{i=1}^k\ves_i$. By
(6.2.64), (9.1.21) and Weyl's Theorem 2.3.6 of completely
reducibility if $\mbb F=\mbb C$ or more generally by Lemma 6.3.2
with $n_1=0$, we have: \psp

 {\bf Theorem 9.1.2}. {\it For $k\in\ol{0,n}$, $\check{\msr
H}_k$ is a finite-dimensional irreducible $sp(2n,\mbb{F})$-module
with the highest-weight vector $\vec\sta_k$ of weight
$\sum_{i=1}^k\ves_i$. Moreover,
$$\check{\msr A}=\bigoplus_{k=0}^n\bigoplus_{\ell=0}^{\llbracket n-k/2\rrbracket}\check\eta^\ell \check{\msr
H}_k.\eqno(9.1.22)$$}\psp

We let ${\msr A}=\mbb{F}[x_1,...,x_n]$. We define a representation
of $sp(2n,\mbb{F})$ on ${\msr A}$ by
$$(E_{i,j}-E_{n+j,n+i})|_{\msr A}=x_i\ptl_{x_j}+\frac{1}{2}\dlt_{i,j}\eqno(9.1.23)$$
$$(E_{i,n+j}+E_{j,i+n})|_{\msr
A}=x_ix_j,\;\;(E_{n+i,j}+E_{n+j,i})|_{\msr
A}=-\ptl_{x_i}\ptl_{x_j}\eqno(9.1.24)$$ for $i,j\in\ol{1,n}$.
 Denote
$${\msr A}_{(0)}=\mbox{Span}\:\{x^\al\mid
\al\in\mbb{N}^{\:n},\;|\al|\;\mbox{is even}\},\;\;{\msr
A}_{(1)}=\mbox{Span}\:\{x^\al\mid
\al\in\mbb{N}^{\:n},\;|\al|\;\mbox{is odd}\}.\eqno(9.1.25)$$ \pse

{\bf Theorem 9.1.3}. {\it The subspace ${\msr A}_{(0)}$ is an
irreducible $sp(2n,\mbb{F})$-module and the identity  element
$1_{\msr A}$ is an highest weight vector with weight
$(1/2)\sum_{i=1}^n\ves_i$. Moreover,  ${\msr A}_{(1)}$ is also an
irreducible $sp(2n,\mbb{F})$-module and the identity $x_1$ element
is an highest weight vector with weight
$\ves_1+(1/2)\sum_{i=1}^n\ves_i$.}

\section{Noncanonical Oscillator Representations}

In this section, we study the noncanonical oscillator
representations of symplectic Lie algebras obtained from the
 the canonical bosonic oscillator representations  of sympletic  Lie algebras over
their minimal natural modules by partially swapping differential
operators and multiplication operators.

\subsection{Main Theorem}

Denote ${\msr B}=\mbb{F}[x_1,...,x_n,y_1,...,y_n]$.  Fix
 $n_1,n_2\in\ol{1,n}$ with $n_1\leq n_2$. Swapping operators $\ptl_{x_r}\mapsto -x_r,\;
 x_r\mapsto
\ptl_{x_r}$  for $r\in\ol{1,n_1}$ and $\ptl_{y_s}\mapsto -y_s,\;
 y_s\mapsto\ptl_{y_s}$  for $s\in\ol{n_2+1,n}$ in the canonical
oscillator representation (9.1.8) and (9.1.9), we have the
nonocanoical oscillator representation of the Lie algebra
$sp(2n,\mbb{F})$ on ${\msr B}$ determined by
$$(E_{i,j}-E_{n+j,n+i})|_{\msr
B}=E_{i,j}^x-E_{j,i}^y\eqno(9.2.1)$$ with
$$E_{i,j}^x=\left\{\begin{array}{ll}-x_j\ptl_{x_i}-\delta_{i,j}&\mbox{if}\;
i,j\in\ol{1,n_1};\\ \ptl_{x_i}\ptl_{x_j}&\mbox{if}\;i\in\ol{1,n_1},\;j\in\ol{n_1+1,n};\\
-x_ix_j &\mbox{if}\;i\in\ol{n_1+1,n},\;j\in\ol{1,n_1};\\
x_i\partial_{x_j}&\mbox{if}\;i,j\in\ol{n_1+1,n};
\end{array}\right.\eqno(9.2.2)$$
and
$$E_{i,j}^y=\left\{\begin{array}{ll}y_i\ptl_{y_j}&\mbox{if}\;
i,j\in\ol{1,n_2};\\ -y_iy_j&\mbox{if}\;i\in\ol{1,n_2},\;j\in\ol{n_2+1,n};\\
\ptl_{y_i}\ptl_{y_j} &\mbox{if}\;i\in\ol{n_2+1,n},\;j\in\ol{1,n_2};\\
-y_j\partial_{y_i}-\delta_{i,j}&\mbox{if}\;i,j\in\ol{n_2+1,n};
\end{array}\right.\eqno(9.2.3)$$
$$E_{i,n+j}|_{\msr B}=\left\{\begin{array}{ll}
\ptl_{x_i}\ptl_{y_j}&\mbox{if}\;i\in\ol{1,n_1},\;j\in\ol{1,n_2},\\
-y_j\ptl_{x_i}&\mbox{if}\;i\in\ol{1,n_1},\;j\in\ol{n_2+1,n},\\
x_i\ptl_{y_j}&\mbox{if}\;i\in\ol{n_1+1,n},\;j\in\ol{1,n_2},\\
-x_iy_j&\mbox{if}\;i\in\ol{n_1+1,n},\;j\in\ol{n_2+1,n};\end{array}\right.\eqno(9.2.4)$$
$$E_{n+i,j}|_{\msr B}=\left\{\begin{array}{ll}
-x_jy_i&\mbox{if}\;j\in\ol{1,n_1},\;i\in\ol{1,n_2},\\
-x_j\ptl_{y_i}&\mbox{if}\;j\in\ol{1,n_1},\;i\in\ol{n_2+1,n},\\
y_i\ptl_{x_j}&\mbox{if}\;j\in\ol{n_1+1,n},\;i\in\ol{1,n_2},\\
\ptl_{x_j}\ptl_{y_i}&\mbox{if}\;j\in\ol{n_1+1,n},\;i\in\ol{n_2+1,n}.\end{array}\right.\eqno(9.2.5)$$

As Section 7.2, we set
$${\msr B}_{\la k\ra}=\mbox{Span}\{x^\al
y^\be\mid\al,\be\in\mbb{N}^n;\sum_{r=n_1+1}^n\al_r-\sum_{i=1}^{n_1}\al_i+
\sum_{i=1}^{n_2}\be_i-\sum_{r=n_2+1}^n\be_r=k\}\eqno(9.2.6)$$ for
$k\in\mbb{Z}$. Then ${\msr B}=\bigoplus_{k\in\mbb{Z}}{\msr B}_{\la
k\ra}$ is a $\mbb{Z}$-graded space. Let
$$\wht D=\sum_{r=n_1+1}^nx_r\ptl_{x_r}-\sum_{i=1}^{n_1}x_i\ptl_{x_i}+
\sum_{i=1}^{n_2}y_i\ptl_{y_i}-\sum_{r=n_2+1}^ny_r\ptl_{y_r}.\eqno(9.2.7)$$
It is straightforward to verify
$${\msr B}_{\la k\ra}=\{f\in{\msr B}\mid \wht D(f)=kf\}\eqno(9.2.8)$$
and
$$[E_{i,j}^x,\wht D]=[E_{i,j}^y,\wht D]=[E_{i,n+j}|_{\msr
B},\wht D]=[E_{n+i,j}|_{\msr B},\wht D]=0\eqno(9.2.9)$$ for
$i,j\in\ol{1,n}$. Thus ${\msr B}_{\la k\ra}$ forms an
$sp(2n,\mbb{F})$-submodule for any $k\in\mbb{Z}$.

The bilinear form $(\cdot|\cdot)$ on ${\msr B}$ is defined by
(6.4.46).  The main theorem in this section is:\psp

{\bf Theorem 9.2.1}. {\it Let $k\in\mbb{Z}$. If $n_1<n_2$ or $k\neq
0$, the subspace ${\msr B}_{\la k\ra}$  is an irreducible
$sp(2n,\mbb{F})$-module. Moreover, it is a highest-weight module
only if $n_2=n$, in which case for $m\in\mbb{N}$, $x_{n_1}^{-m}$ is
a highest-weight vector of ${\msr B}_{\la -m\ra}$ with weight
$m\ves_{n_1}-\sum_{i=1}^{n_1}\ves_i$, $x_{n_1+1}^{m+1}$ is a
highest-weight vector of ${\msr B}_{\la m+1\ra}$ with weight
$(m+1)\ves_{n_1+1}-\sum_{i=1}^{n_1}\ves_i$
 if
$n_1<n$ and $y_n^{m+1}$ is a highest-weight vector of ${\msr B}_{\la
m+1\ra}$ with weight $-(m+1)\ves_n-\sum_{i=1}^{n_1}\ves_i$ when
$n_1=n$.

When $n_1=n_2$, the subspace ${\msr B}_{\la 0\ra}$ is a direct sum
of two irreducible $sp(2n,\mbb{F})$-submodules. If $n_1=n_2=n$, they
are highest-weight modules with a highest-weight vector $1$ of
weight $-\sum_{i=1}^n\ves_i$ and with a highest-weight vector
$x_{n-1}y_n-x_ny_{n-1}$ of weight
$-\ves_{n-1}-\ves_n-\sum_{i=1}^n\ves_i$, respectively.} \psp

When $n_1-n_2+1-\dlt_{n_1,n_2}< k\in\mbb{Z}$, ${\msr B}_{\la k\ra}$
is not completely reducible with respect to its Lie subalgebra in
(9.1.11) by Sections 6.4-6.6. Thus the above representation
(9.2.1)-(9.2.5) is not unitary. Observe that
$$(E_{i,i}-E_{n+i,n+i})|_{\msr
B}=\left\{\begin{array}{ll}-x_i\ptl_{x_i}-y_i\ptl_{y_i}-1&\mbox{if}\;\;1\leq i\leq n_1,\\
x_i\ptl_{x_i}-y_i\ptl_{y_i}&\mbox{if}\;\;n_1<i\leq n_2,\\
x_i\ptl_{x_i}+y_i\ptl_{y_i}+1&\mbox{if}\;\;n_2<i\leq
n\end{array}\right.\eqno(9.2.10)$$ by (9.2.1)-(9.2.3). So only
finite number of monomials $x_i^{\ell_1}y_i^{\ell_2}$ appear in a
weight subspace when $1\leq i\leq n_1$ or $n_2<i\leq n$. If
$n_1<i\leq n_2$, the degree of $x_i^{\ell_1}y_i^{\ell_2}$ is
$\ell_1+\ell_2$. This shows that all ${\msr B}_{\la k\ra}$ with
$k\in\mbb{Z}$ are infinite-dimensional weight
$sp(2n,\mbb{F})$-modules with finite-dimensional weight subspaces.

Due to the mixture of differential operators with multiplications
operators, the modules ${\msr B}_{\la k\ra}$ with $k\in\mbb{Z}$ are
cuspidal modules  without highest-weight vectors if $n_2<n$.

Note
$$(E_{i,n+j}+E_{j,n+i})|_{\msr
B}=\ptl_{x_i}\ptl_{y_j}+\ptl_{x_j}\ptl_{y_i},\;\;(E_{i,n+r}+E_{r,n+i})|_{\msr
B}=\ptl_{x_i}\ptl_{y_r}+x_r\ptl_{y_i},\eqno(9.2.11)$$
$$(E_{r,n+s}+E_{s,n+r})|_{\msr
B}=x_r\ptl_{y_s}+x_s\ptl_{y_r}\eqno(9.2.12)$$ for $i,j\in\ol{1,n_1}$
and $r,s\in\ol{n_1+1,n_2}$ by (9.2.4).
 Moreover,
$$(E_{i,j}-E_{n+j,n+i})|_{\msr
B}=-x_j\ptl_{x_i}-y_j\ptl_{y_i}-\dlt_{i,j}\eqno(9.2.13)$$ and
$$(E_{i,r}-E_{n+r,n+i})|_{\msr
B}=\ptl_{x_i}\ptl_{x_r}-y_r\ptl_{y_i}\eqno(9.2.14)$$ for
$i,j\in\ol{1,n_1}$ and $r\in\ol{n_1+1,n_2}$ by (9.2.1)-(9.2.3).
Recall
$$\eta=\sum_{i=1}^{n_1}y_i\ptl_{x_i}+\sum_{r=n_1+1}^{n_2}x_ry_r+\sum_{s=n_2+1}^n
x_s\ptl_{y_s}.\eqno(9.2.15)$$

\subsection{Proof of the Theorem When $n_2=n$}

In this subsection,  we assume $n_2=n$.\psp

Under the assumption, any element of ${\msr B}$ is nilpotent with
respect to $sp(2n,\mbb{F})_+$ by (6.3.21) and (9.2.11)-(9.2.14).
Thus any $sp(2n,\mbb{F})$-submodule of ${\msr B}$ must contain an
$sp(2n,\mbb{F})$-singular vector.

First we assume $n_1+1<n$. According to Lemma 6.4.1, the nonzero
weight vectors in
$$\mbox{Span}\{\eta^{m_3}(x_i^{m_1}y_n^{m_2})\mid
m_r\in\mbb{N};i=n_1,n_1+1\}\eqno(9.2.16)$$ are all the $sl(n,\mbb
F)$-singular vectors
 in ${\msr B}$. The singular vectors of $sp(2n,\mbb{F})$ in ${\msr
B}$ must be among them by (9.1.11). Moreover, the subalgebra
$sp(2n,\mbb{F})_+$ is generated by $sl(n,\mbb F)_+$ and $E_{n,2n}$.
According to (9.2.11), $E_{n,2n}|_{\msr B}=x_n\ptl_{y_n}$. Hence
$$E_{n,2n}(\eta^{m_3}(x_i^{m_1}y_n^{m_2}))=x_n[
m_3x_n\eta^{m_3-1}(x_i^{m_1}y_n^{m_2})+m_2\eta^{m_3}(x_i^{m_1}y_n^{m_2-1})]\eqno(9.2.17)$$
for $i=n_1,n_1+1$ by (9.2.15). Considering weights, we conclude that
the vectors $\{x_{n_1}^{m},x_{n_1+1}^{m+1}\mid m\in\mbb{N}\}$ are
all the singular vectors of $sp(2n,\mbb{F})$ in ${\msr B}$.
Furthermore,
$$x_{n_1}^{m}\in {\msr B}_{\la
-m\ra}\;\;\mbox{and}\;\;x_{n_1+1}^{m+1}\in {\msr B}_{\la
m+1\ra}\qquad\for\;\;m\in\mbb{N}.\eqno(9.2.18)$$ Thus each ${\msr
B}_{\la k\ra}$ has a unique non-isotropic singular vector for
$k\in\mbb{Z}$. By Lemma 6.3.2, all ${\msr B}_{\la k\ra}$ with
$k\in\mbb{Z}$ are irreducible highest-weight
$sp(2n,\mbb{F})$-submodules.

Consider the case $n_1+1=n$.  According to (6.5.61), the nonzero
weight vectors in
$$\mbox{Span}\{\eta^{m_2}(x_{n-1}^{m_1}y_n^{m_3}),
x_n^{m_1}y_n^{m_2},\eta^{m_1+m_2}(x_{n-1}^{m_2}y_n^{m_3-m_1})\mid
m_i\in\mbb{N}\}\eqno(9.2.19)$$ are all the $sl(n,\mbb F)$-singular
vectors
 in ${\msr B}$.
By (9.2.17) and considering weights, we again conclude that the
vectors $\{x_{n-1}^{m},x_n^{m+1}\mid m\in\mbb{N}\}$ are all the
singular vectors of $sp(2n,\mbb{F})$ in ${\msr B}$. Again all ${\msr
B}_{\la k\ra}$ with $k\in\mbb{Z}$ are irreducible highest-weight
$sp(2n,\mbb{F})$-submodules.

Suppose $n_1=n$. By (9.2.11), we have $E_{n,2n}|_{\msr
B}=\ptl_{x_n}\ptl_{y_n}$ in this case. According to (6.6.38), the
nonzero weight vectors in
$$\mbox{Span}\{x_n^{m_1}y_n^{m_2}\zeta_1^{m_3}\mid
m_i\in\mbb{N}\}\eqno(9.2.20)$$
 are all the $sl(n,\mbb F)$-singular vectors  in ${\msr B}$,
 where $\zeta_1=x_{n-1}y_n-x_ny_{n-1}$
in this case. Moreover,
\begin{eqnarray*}
& &E_{n,2n}(x_n^{m_1}y_n^{m_2}\zeta_1^{m_3})\\ &=&
m_1m_2x_n^{m_1-1}y_n^{m_2-1}\zeta_1^{m_3}+m_1m_3x_{n-1}x_n^{m_1-1}y_n^{m_2}\zeta_1^{m_3-1}
\\
&&-m_2m_3y_{n-1}x_n^{m_1}y_n^{m_2-1}\zeta_1^{m_3-1}
-m_3(m_3-1)x_{n-1}y_{n-1}x_n^{m_1}y_n^{m_2}\zeta_1^{m_3-2}
\\
&=&m_1(m_2+m_3)x_n^{m_1-1}y_n^{m_2-1}\zeta_1^{m_3}
+m_3(m_1-m_2-m_3+1)y_{n-1}x_n^{m_1}y_n^{m_2-1}\zeta_1^{m_3-1}
\\& &-m_3(m_3-1)y_{n-1}^2x_n^{m_1+1}y_n^{m_2-1}\zeta_1^{m_3-2}
.\hspace{7.2cm}(9.2.21)\end{eqnarray*}
 Considering weights, we again conclude that the
vectors $\{x_n^m,y_n^{m+1},\zeta_1\mid m\in\mbb{N}\}$ are all the
singular vectors of $sp(2n,\mbb{F})$ in ${\msr B}$. Moreover,
$$x_n^{m}\in {\msr B}_{\la -m\ra},\;\;\zeta_1\in{\msr B}_{\la 0\ra}\;\;\mbox{and}\;\;y_n^{m+1}
\in{\msr B}_{\la m+1\ra}\qquad\for\;\;m\in\mbb{N}.\eqno(9.2.22)$$
Thus each ${\msr B}_{\la k\ra}$ with $k\neq 0$ has a unique
non-isotropic singular vector for $k\in\mbb{Z}$. By Lemma 6.3.2, all
${\msr B}_{\la k\ra}$ with $0\neq k\in\mbb{Z}$ are irreducible
highest-weight $sp(2n,\mbb{F})$-submodules.

Set
$${\msr B}_{\la
0,1\ra}=\mbox{Span}\{[\prod_{1\leq r\leq s\leq
n}(x_ry_s+x_sy_r)^{l_{r,s}}]\mid
l_{r,s}\in\mbb{N}\}\eqno(9.2.23)$$and
$${\msr B}_{\la 0,2\ra}=\mbox{Span}\{[\prod_{1\leq r\leq s\leq
n}(x_ry_s+x_sy_r)^{l_{r,s}}](x_py_q-x_qy_p)\mid
l_{r,s}\in\mbb{N};1\leq p<q\leq n\}.\eqno(9.2.24)$$ Let $${\msr
G}'=\sum_{1\leq r\leq s\leq
n}\mbb{F}(E_{n+s,r}+E_{n+r,s})\eqno(9.2.25)$$
 and
$$\wht{\msr
G}=\sum_{i,j=1}^n\mbb{F}(E_{i,j}-E_{n+j,n+i})+\sum_{1\leq r\leq
s\leq n}\mbb{F}(E_{r,n+s}+E_{s,n+r}).\eqno(9.2.26)$$
 Then ${\msr G}'$ and
$\wht{\msr G}$ are Lie subalgebras of $sp(2n,\mbb{F})$ and
$sp(2n,\mbb{F})={\msr G}'\oplus \wht{\msr G}.$ By PBW Theorem
$$U(sp(2n,\mbb{F}))=U({\msr G}')U(\wht{\msr G}).\eqno(9.2.27)$$
Note
$$(E_{n+s,r}+E_{n+r,s})|_{\msr B}=-(x_ry_s+x_sy_r)\qquad\for\;\;r,s\in\ol{1,n}\eqno(9.2.28)$$
by (9.2.5). According to (9.2.11), (9.2.13) and (9.2.28),
$${\msr B}_{\la
0,1\ra}=U({\msr G}')(1)=U(sp(2n,\mbb{F}))(1)\eqno(9.2.29)$$ and
$${\msr B}_{\la
0,2\ra}=\sum_{1\leq p<q\leq n}U({\msr
G}')(x_py_q-x_qy_p)=U(sp(2n,\mbb{F}))(\zeta_1)\eqno(9.2.30)$$ are
$sp(2n,\mbb{F})$-submodules.

It is obvious, $1\not\in {\msr B}_{\la 0,2\ra}$. On the other hand,
$({\msr B}_{\la 0,1\ra}|x_{n-1}y_n-x_ny_{n-1})=\{0\}$. Hence
$x_{n-1}y_n-x_ny_{n-1}\not \in {\msr B}_{\la 0,1\ra}$. Thus ${\msr
B}_{\la 0,1\ra}$ and ${\msr B}_{\la 0,0\ra}$ have a unique
non-isotropic singular vector. By Lemma 6.3.2, they are irreducible.
Note that
$$\wht{\msr B}_{\la 0\ra}={\msr B}_{\la 0,1\ra}+ {\msr B}_{\la
0,2\ra}\eqno(9.2.31)$$ is a direct sum of
$sp(2n,\mbb{F})$-submodules and the restriction of the symmetric
bilinear form $(\cdot|\cdot)$ on $\wht{\msr B}_{\la 0\ra}$ is
nondegenerate because ${\msr B}_{\la 0,1\ra}$ and ${\msr B}_{\la
0,2\ra}$ are irreducible $sp(2n,\mbb{F})$-submodules whose singular
vectors are non-isotropic. Thus
$${\msr B}_{\la 0\ra}=\wht{\msr B}_{\la 0\ra}\oplus ((\wht{\msr B}_{\la
0\ra})^\bot\bigcap {\msr B}_{\la 0\ra})\eqno(9.2.32)$$ and $({\msr
B}'_{\la 0\ra})^\bot\bigcap {\msr B}_{\la 0\ra}$ is an
$sp(2n,\mbb{F})$-submodule. Thus $(\wht{\msr B}_{\la
0\ra})^\bot\bigcap {\msr B}_{\la 0\ra}$ contains an
$sp(2n,\mbb{F})$-singular vector. Since $1$ and
$x_{n-1}y_n-x_ny_{n-1}$ are the only singular vectors in ${\msr
B}_{\la 0\ra}$, it contains $1$ or $x_{n-1}y_n-x_ny_{n-1}$. This
contradicts the fact that $1$ and $x_{n-1}y_n-x_ny_{n-1}$ are
non-isotropic. Thus $(\wht{\msr B}_{\la 0\ra})^\bot\bigcap {\msr
B}_{\la 0\ra}=\{0\}$. This proves the theorem when $n_2=n$.

\subsection{Proof of the Theorem When $n_1<n_2<n$}

In this subsection, we assume $n_1<n_2<n$. \psp

{\it Case 1}. $n_1+1<n_2$\psp

Set
\begin{eqnarray*}\hspace{1cm}{\msr G}_1&=&
\sum_{i,j=1}^{n_1+1}\mbb{F}(E_{i,j}-E_{n+j,n+i})+\sum_{i=1}^{n_1+1}(\mbb{F}E_{i,n+i}+\mbb{F}E_{n+i,i})\\
& &+\sum_{1\leq i<j\leq n_1+1
}[\mbb{F}(E_{i,n+j}+E_{n+j,i})+\mbb{F}(E_{n+i,j}+E_{n+j,i})]\hspace{2.8cm}(9.2.33)\end{eqnarray*}
and
\begin{eqnarray*}\hspace{1cm}{\msr G}_2&=&
\sum_{i,j=n_1+2}^n\mbb{F}(E_{i,j}-E_{n+j,n+i})+\sum_{i=n_1+2}^n(\mbb{F}E_{i,n+i}+\mbb{F}E_{n+i,i})\\
& &+\sum_{n_1+2\leq i<j\leq n
}[\mbb{F}(E_{i,n+j}+E_{n+j,i})+\mbb{F}(E_{n+i,j}+E_{n+j,i})].\hspace{2.7cm}(9.2.34)\end{eqnarray*}
Then ${\msr G}_1=sp(2n_1+2,\mbb{F})$ and ${\msr G}_2\cong
sp(2(n-n_1)-2,\mbb{F})$ are Lie subalgebras of $sp(2n,\mbb{F})$.
Denote
$${M}^1=\mbb{F}[x_1,...,x_{n_1+1},y_1,...,y_{n_1+1}],\qquad
 M^2=\mbb{F}[x_{n_1+2},...,x_n,y_{n_1+2},...,y_n].\eqno(9.2.35)$$
Observe that ${M}^1$ is exactly the ${\msr G}_1$-module as ${\msr
B}$ in Subsection 9.2.2 with $n\rightarrow n_1+1$ and ${M}^2$ is
exactly the ${\msr G}_1$-module as ${\msr B}$ in Subsection 9.2.2
with  $n\rightarrow n-n_1-1$. Write
$$M^i_{\la k\ra}=M^i\bigcap\msr B_{\la
k\ra}.\eqno(9.2.36)$$ By Subsection 9.2.2, $\{M^i_{\la k\ra}\mid
k\in\mbb Z\}$ are irreducible $\msr G_i$-submodules.

 According to Lemma 6.4.1, the
nonzero weight vectors in
$$\mbox{Span}\{\eta^{m_3}(x_i^{m_1}y_j^{m_2})\mid
m_r\in\mbb{N};i=n_1,n_1+1;j=n_2,n_2+1\}\eqno(9.2.37)$$ are all the
$sl(n,\mbb F)$-singular vectors  in ${\msr B}$. Fix $k\in\mbb{N}$.
Then the $sl(n,\mbb F)$-singular vectors in ${\msr B}_{\la-k\ra}$
are
\begin{eqnarray*}\qquad\qquad& &\{\eta^{m_3}(x_{n_1}^{k+m_2+2m_3}y_{n_2}^{m_2}),
\eta^{m_3}(x_{n_1+1}^{m_1}y_{n_2+1}^{k+m_1+2m_3}),\\
& &\eta^{m_3}(x_{n_1}^{m_4}y_{n_2+1}^{m_5}) \mid
m_i\in\mbb{N};m_4+m_5-2m_3=k\}.\hspace{3.9cm}(9.2.38)\end{eqnarray*}
Let $M$ be a nonzero $sp(2n,\mbb{F})$-submodule of ${\msr
B}_{\la-k\ra}$. Then $M$ contains an $sl(n,\mbb F)$-singular vector
by (9.1.11). Suppose some
$\eta^{m_3}(x_{n_1}^{k+m_2+2m_3}y_{n_2}^{m_2})\in M$. We have
$E_{n_1,n+n_1}|_{\msr B}=\ptl_{x_{n_1}}\ptl_{y_{n_1}}$ and
$$E_{n_1,n+n_1}^{m_3}[\eta^{m_3}(x_{n_1}^{k+m_2+2m_3}y_{n_2}^{m_2})]=m_3![\prod_{r=1}^{2m_3}(k+m_2+r)]
x_{n_1}^{k+m_2}y_{n_2}^{m_2}\in M\eqno(9.2.39)$$ by (9.2.11) and
(9.2.15). Moreover, $(E_{n_1,n+n_2}+E_{n_2,n+n_1})|_{\msr
B}=\ptl_{x_{n_1}}\ptl_{y_{n_2}}+x_{n_2}\ptl_{y_{n_1}}$ and
$$(E_{n_1,n+n_2}+E_{n_2,n+n_1})^{m_2}(x_{n_1}^{k+m_2}y_{n_2}^{m_2})=m_2!
[\prod_{r=1}^{m_2}(k+r)] x_{n_1}^k\in M\eqno(9.2.40)$$ by (9.2.11).
Thus
$$x_{n_1}^k\in M.\eqno(9.2.41)$$

Assume some $\eta^{m_3}(x_{n_1+1}^{m_1}y_{n_2+1}^{k+m_1+2m_3})\in
M$. According to (9.2.5),
 $$(E_{n+i,j}+E_{n+j,i})|_{\msr B}=\ptl_{x_i}\ptl_{y_j}+\ptl_{x_j}\ptl_{y_i}\qquad\for\;\;i\in\ol{n_2+1,n}.\eqno(9.2.42)$$
So
$$E_{n+n_2+1,n_2+1}^{m_3}[\eta^{m_3}(x_{n_1+1}^{m_1}y_{n_2+1}^{k+m_1+2m_3})]=
m_3![\prod_{r=1}^{2m_3}(k+m_1+r)]x_{n_1+1}^{m_1}y_{n_2+1}^{k+m_1}\in
M.\eqno(9.2.43)$$ Moreover,
$$(E_{n+n_2+1,n_1+1}+E_{n+n_1+1,n_2+1})|_{\msr B}=\ptl_{x_{n_1+1}}
\ptl_{y_{n_2+1}}+y_{n_1+1}\ptl_{x_{n_2+1}}\eqno(9.2.44)$$ by
(9.2.5). Hence
$$(E_{n+n_2+1,n_1+1}+E_{n+n_1+1,n_2+1})^{m_1}(x_{n_1+1}^{m_1}y_{n_2+1}^{k+m_1})=m_1!
[\prod_{r=1}^{m_1}(k+r)] y_{n_2+1}^k\in M.\eqno(9.2.45)$$
Furthermore,
$$(E_{n+n_2+1,n_1}+E_{n+n_1,n_2+1})|_{\msr
B}=-x_{n_1}\ptl_{y_{n_2+1}}+y_{n_1}\ptl_{x_{n_2+1}}\eqno(9.2.46)$$
by (9.2.5). Thus
$$(E_{n+n_2+1,n_1}+E_{n+n_1,n_2+1})^k(y_{n_2+1}^k)=(-1)^kk!x_{n_1}^k\in
M.\eqno(9.2.47)$$ Thus (9.2.41) holds again.

Consider $\eta^{m_3}(x_{n_1}^{m_4}y_{n_2+1}^{m_5})$ for some
$m_3,m_3,m_4\in\mbb{N}$ such that $m_4+m_5-2m_3=k$.  Note that
$E_{n_1+1,n+n_1+1}|_{\msr B}=x_{n_1+1}\ptl_{y_{n_1+1}}$ by (9.2.12)
and
$$E_{n_1+1,n+n_1+1}^{m_3}[\eta^{m_3}(x_{n_1}^{m_4}y_{n_2+1}^{m_5})]
=m_3!x_{n_1+1}^{2m_3}x_{n_1}^{m_4}y_{n_2+1}^{m_5}\in
M.\eqno(9.2.48)$$ There exist $r_1,r_2\in\mbb{N}$ such that
$r_1+r_2=2m_3$ and $r_1\leq m_4,\;r_2\leq m_5$. Moreover,
$$(E_{n_1,n_1+1}-E_{n+n_1+1,n+n_1})|_{\msr
B}=\ptl_{x_{n_1}}\ptl_{x_{n_1+1}}-y_{n_1+1}\ptl_{y_{n_1}}\eqno(9.2.49)$$
by (9.2.1)-(9.2.3). So (9.2.44) and (9.2.49) yield
\begin{eqnarray*}\qquad &
&(E_{n_1,n_1+1}-E_{n+n_1+1,n+n_1})^{r_1}(E_{n+n_2+1,n_1+1}+E_{n+n_1+1,n_2+1})^{r_2}
(x_{n_1+1}^{2m_3}x_{n_1}^{m_4}y_{n_2+1}^{m_5})\\
&=&(2m_3)![\prod_{s_1=0}^{r_1-1}(m_4-s_1)]
[\prod_{s_2=0}^{r_2-1}(m_5-s_2)]x_{n_1}^{m_4-r_1}y_{n_2+1}^{m_5-r_2}\in
M.\hspace{3cm}(9.2.50)\end{eqnarray*} Furthermore, (9.2.5) yields
\begin{eqnarray*}\hspace{2cm}& &(E_{n+n_2+1,n_1}+E_{n+n_1,n_2+1})^{m_5-r_2}
(x_{n_1}^{m_4-r_1}y_{n_2+1}^{m_5-r_2})\\
&=&(-1)^{m_5-r_2}(m_5-r_2)!x_{n_1}^k\in
M.\hspace{6.3cm}(9.2.51)\end{eqnarray*} Thus we always have
$x_{n_1}^k\in M$.

Note that ${M}^1_{\la-k\ra}\ni x_{n_1}^k$ is an irreducible ${\msr
G}_1$-module (cf. (9.2.33) and (9.2.35)) by Subsection 9.2.2. So
$${M}^1_{\la-k\ra}\subset M.\eqno(9.2.52)$$
Observe that
$$(E_{n_1,n+n_2+1}+E_{n_2+1,n+n_1})|_{\msr
B}=-y_{n_2+1}\ptl_{x_{n_1}}+x_{n_2+1}\ptl_{y_{n_1}}\eqno(9.2.53)$$
by (9.2.4). Moreover,
$$(E_{n+n_2,n_1+1}+E_{n+n_1+1,n_2})|_{\msr
B}=y_{n_2}\ptl_{x_{n_1+1}}+y_{n_1+1}\ptl_{x_{n_2}}\eqno(9.2.54)$$ by
(9.2.5).

For any $r\in\mbb N$, we have $x_{n_1}^{k+r}x_{n_1+1}^{r}\in
M^1_{\la-k\ra}\subset M$. Take any $s_1\in\ol{0,k+r}$ and
$r_1\in\ol{0,r}$. Using (9.2.53) and (9.2.54), we get
\begin{eqnarray*}\qquad&&(E_{n+n_2,n_1+1}+E_{n+n_1+1,n_2})^{r_1}(E_{n_1,n+n_2+1}+E_{n_2+1,n+n_1})^s(x_{n_1}^{k+r}x_{n_1+1}^{r})
\\
&=&(-1)^s[\prod_{i_1=0}^{r_1}(r-i_1)][\prod_{i_2=0}^{s-1}(k+r-i_2)]x_{n_1}^{k+r-s}x_{n_1+1}^{r-r_1}y_{n_2}^{r_1}y_{n_2+1}^s\in
M.\hspace{1.9cm}(9.2.55)\end{eqnarray*} Thus
$$x_{n_1}^{r_1}x_{n_1+1}^{r_2}y_{n_2}^{s_1}y_{n_2+1}^{s_2}\in
M\;\;\mbox{whenever}\;\;r_2-r_1+s_1-s_2=-k.\eqno(9.2.56)$$ Applying
$U(\msr G_1)U(\msr G_2)$ to (9.2.56), we get
$$M^1_{\la k_1\ra}M^2_{\la k_2\ra}\subset M \;\;\mbox{whenever}\;\;k_1+k_2=-k.\eqno(9.2.57)$$
 Therefore,
${\msr B}_{\la -k\ra}=\sum_{k_1\in\mbb Z}M^1_{\la k_1\ra} M^2_{\la
-k-k_1\ra}=M$; that is, ${\msr B}_{\la -k\ra}$ is an irreducible
$sp(2n,\mbb{F})$-submodule.

Fix $0<k\in\mbb{N}$.  Then the singular vectors of ${\msr K}$ in
${\msr B}_{\la k\ra}$ are
\begin{eqnarray*}\qquad& &\{\eta^{m_2}(x_{n_1+1}^{k+m_1-2m_2}y_{n_2+1}^{m_1}),
\eta^{m_2}(x_{n_1}^{m_1}y_{n_2}^{k+m_1-2m_2}),\eta^{m_3}(x_{n_1+1}^{m_4}y_{n_2}^{m_5})
\\
& &\mid m_i\in\mbb{N};2m_2\leq
k+m_1;m_4+m_5+2m_3=k\}\hspace{4.9cm}(9.2.58)\end{eqnarray*} by
(9.2.39). Let $M$ be a nonzero $sp(2n,\mbb{F})$-submodule of ${\msr
B}_{\la k\ra}$. Then $M$ contains an $sl(n,\mbb F)$-singular vector
by (9.1.11). Suppose some
$\eta^{m_2}(x_{n_1+1}^{k+m_1-2m_2}y_{n_2+1}^{m_1})\in M$ with
$2m_2\leq k+m_1$. We have $E_{n_1+1,n+n_1+1}|_{\msr
B}=x_{n_1+1}\ptl_{y_{n_1+1}}$ and
$$E_{n_1+1,n+n_1+1}^{m_2}[\eta^{m_2}(x_{n_1+1}^{k+m_1-2m_2}y_{n_2}^{m_1})]=m_2!
x_{n_1+1}^{k+m_1}y_{n_2+1}^{m_1}\in M\eqno(9.2.59)$$ by (9.2.12) and
(9.2.15). Moreover, (9.2.5) gives
$$(E_{n+n_2+1,n_1+1}+E_{n+n_1+1,n_2+1})^{m_1}(x_{n_1+1}^{k+m_1}y_{n_2+1}^{m_1})=m_1!
[\prod_{r=1}^{m_1}(k+r)] x_{n_1+1}^k\in M.\eqno(9.2.60)$$ Thus
$$x_{n_1+1}^k\in M.\eqno(9.2.61)$$

Assume some $\eta^{m_2}(x_{n_1}^{m_1}y_{n_2}^{k+m_1-2m_2})\in M$
with $2m_2\leq k+m_1$. Observe $E_{n+n_2,n_2}=y_{n_2}\ptl_{x_{n_2}}$
by (9.2.5). So
$$E_{n+n_2,n_2}^{m_2}[\eta^{m_2}(x_{n_1}^{m_1}y_{n_2}^{k+m_1-2m_2})]=
m_2!x_{n_1}^{m_1}y_{n_2}^{k+m_1}\in M.\eqno(9.2.62)$$  Moreover,
(9.2.11) gives that $(E_{n_1,n+n_2}+E_{n_2,n+n_1})|_{\msr
B}=\ptl_{x_{n_1}}\ptl_{y_{n_2}}+x_{n_2}\ptl_{y_{n_1}}$ and
$$(E_{n_1,n+n_2}+E_{n_2,n+n_1})^{m_1}(x_{n_1}^{m_1}y_{n_2}^{k+m_1})=m_1!
[\prod_{r=1}^{m_1}(k+r)] y_{n_2}^k\in M.\eqno(9.2.63)$$ Furthermore,
(9.2.12) yields that $(E_{n_1+1,n+n_2}+E_{n_2,n+n_1+1})|_{\msr
B}=x_{n_1+1}\ptl_{y_{n_2}}+x_{n_2}\ptl_{y_{n_1+1}}$ and
$$(E_{n_1+1,n+n_2}+E_{n_2,n+n_1+1})^k(y_{n_2}^k)=k!x_{n_1+1}^k\in
M.\eqno(9.2.64)$$ Thus (9.2.61) holds again.

Consider $\eta^{m_3}(x_{n_1+1}^{m_4}y_{n_2}^{m_5})$ for some
$m_3,m_3,m_4\in\mbb{N}$ such that $m_4+m_5+2m_3=k$.  Note
$E_{n_1+1,n+n_1+1}=x_{n_1+1}\ptl_{y_{n_1+1}}$ by (9.2.12). So
$$E_{n_1+1,n+n_1+1}^{m_3}
[\eta^{m_3}(x_{n_1+1}^{m_4}y_{n_2}^{m_5})]=m_3!x_{n_1+1}^{m_4+2m_3}y_{n_2}^{m_5}\in
M.\eqno(9.2.65)$$ According to (9.2.12),
$$(E_{n_1+1,n+n_2}+E_{n_2,n+n_1+1})^{m_5}(x_{n_1+1}^{m_4+2m_3}y_{n_2}^{m_5})=m_5!x_{n_1+1}^k\in
M.\eqno(9.2.66)$$ Therefore, we always have $x_{n_1+1}^k\in M$.

 Observe that ${M}^1_{\la k\ra}\ni x_{n_1+1}^k$ is an irreducible ${\msr
G}_1$-module (cf. (9.2.34) and (9.2.35)) by Subsection 9.2.2. So
$${M}^1_{\la k \ra}\subset M.\eqno(9.2.67)$$
For any $r\in\mbb N$, we have $x_{n_1}^rx_{n_1+1}^{k+r}\in M^1_{\la
k\ra}\subset M$. Take any $s_1\in\ol{0,r}$ and $r_1\in\ol{0,k+r}$.
Using (9.2.53) and (9.2.54), we get
\begin{eqnarray*}\qquad&&(E_{n+n_2,n_1+1}+E_{n+n_1+1,n_2})^{r_1}(E_{n_1,n+n_2+1}+E_{n_2+1,n+n_1})^s(x_{n_1}^{r}x_{n_1+1}^{k+r})
\\
&=&(-1)^s[\prod_{i_1=0}^{r_1}(k+r-i_1)][\prod_{i_2=0}^{s-1}(r-i_2)]x_{n_1}^{r-s}x_{n_1+1}^{k+r-r_1}y_{n_2}^{r_1}y_{n_2+1}^s\in
M.\hspace{1.9cm}(9.2.68)\end{eqnarray*} Thus
$$x_{n_1}^{r_1}x_{n_1+1}^{r_2}y_{n_2}^{s_1}y_{n_2+1}^{s_2}\in
M\;\;\mbox{whenever}\;\;r_2-r_1+s_1-s_2=k.\eqno(9.2.69)$$ Applying
$U(\msr G_1)U(\msr G_2)$ to (9.2.69), we get
$$M^1_{\la k_1\ra}M^2_{\la k_2\ra}\subset M \;\;\mbox{whenever}\;\;k_1+k_2=k.\eqno(9.2.70)$$
 Therefore,
${\msr B}_{\la k\ra}=\sum_{k_1\in\mbb Z}M^1_{\la k_1\ra} M^2_{\la
k-k_1\ra}=M$; that is, ${\msr B}_{\la k\ra}$ is an irreducible
$sp(2n,\mbb{F})$-submodule.
 \psp

{\it Case 2}. $n_2=n_1+1$. \psp

 According to (6.5.54), the
nonzero weight vectors in
\begin{eqnarray*} & &\mbox{Span}\{\eta^{m_2}(x_i^{m_1}y_j^{m_3}),
x_{n_1+1}^{m_1}y_{n_1+1}^{m_2},\eta^{m_1+m_2}(x_{n_1}^{m_2}y_{n_1+1}^{m_3-m_1}),\eta^{m_1+m_2}(y_{n_1+2}^{m_2}x_{n_1+1}^{m_3-m_1})
\\ & &\qquad\;\;\mid m_r\in\mbb{N};
(i,j)=(n_1,n_1+1),(n_1,n_1+2),(n_1+1,n_1+2)\}.\hspace{2.2cm}(9.2.71)\end{eqnarray*}
are all the singular vectors of $sl(n,\mbb{F})$ in ${\msr B}$. Fix
$k\in\mbb{N}$. Then the $sl(n,\mbb{F})$-singular vectors  in ${\msr
B}_{\la-k\ra}$ are those in (9.2.38). Set
\begin{eqnarray*}\hspace{1cm}{\msr G}_3&=&
\sum_{i,j=n_1+1}^n\mbb{F}(E_{i,j}-E_{n+j,n+i})+\sum_{i=n_1+1}^n(\mbb{F}E_{i,n+i}+\mbb{F}E_{n+i,i})\\
& &+\sum_{n_1+1\leq i<j\leq n
}[\mbb{F}(E_{i,n+j}+E_{n+j,i})+\mbb{F}(E_{n+i,j}+E_{n+j,i})]\hspace{2.8cm}(9.2.72)\end{eqnarray*}
and
$$M^3=\mbb
F[x_{n_1+1},...,x_n,y_{n_1+1},...,y_n],\;\;M^3_{\la \ell\ra}=
M^3\bigcap \msr B_{\la \ell\ra}\;\;\for\;\;\ell\in\mbb
Z.\eqno(9.2.73)$$
 According to the arguments in
Case 1, (9.2.56) holds. Applying $U(\msr G_3)$ to (9.2.56), we have
$$x_{n_1}^{r_1}M^3_{\la r_2+s_1+s_2\ra}\subset M\eqno(9.2.73)$$
by Subsection 9.2.2. In particular,
$$x_{n_1}^{r_1}x_{n_1+1}^{r_2}M^2_{\la s_1+s_2\ra}\subset M.\eqno(9.2.74)$$
Applying $U(\msr G_1)$ to (9.2.74), we get (9.2.57). Thus $\msr
B_{\la-k\ra}$ is an irreducible $sp(2n,\mbb{F})$-submodule.

Let $0<k\in\mbb{N}$. Then the ${\msr K}$-singular vectors in ${\msr
B}_{\la k\ra}$ are
\begin{eqnarray*}\!\!\!\!\!\!\!\! &\{\eta^{m_2}(x_{n_1+1}^{k+m_1-2m_2}y_{n_1+2}^{m_1}),
\eta^{m_2}(x_{n_1}^{m_1}y_{n_1+1}^{k+m_1-2m_2}),
\eta^{m_5+m_6}(x_{n_1}^{m_6}y_{n_1+1}^{m_7-m_5}),\eta^{m_5+m_6}(y_{n_1+2}^{m_6}x_{n_1+1}^{m_7-m_5}),\hspace{0.6cm}
\\
 &\!\!\!\!\!x_{n_1+1}^{m_3}y_{n_1+1}^{m_4} \mid m_i\in\mbb{N};2m_2\leq
k+m_1;m_3+m_4=k=m_5+m_6+m_7\}\hspace{1.9cm}(9.2.75)\end{eqnarray*}
by (9.2.71).  Let $M$ be a nonzero $sp(2n,\mbb{F})$-submodule of
${\msr B}_{\la k\ra}$. As an $sl(n,\mbb F)$-module, $M$ contains an
$sl(n,\mbb F)$-singular vector. If
$x_{n_1+1}^{m_3}y_{n_1+1}^{m_4}\in M$ with $m_3+m_4=k$, then
$E_{n_1+1,n+n_1+1}|_{\msr B}=x_{n_1+1}\ptl_{y_{n_1+1}}$ and
$$E_{n_1+1,n+n_1+1}^{m_4}(x_{n_1+1}^{m_3}y_{n_1+1}^{m_4})=m_4!x_{n_1+1}^k\in M\lra x_{n_1+1}^k\in M\eqno(9.2.76)$$
by (9.2.12). Suppose some
$\eta^{m_5+m_6}(x_{n_1}^{m_6}y_{n_1+1}^{m_7-m_5})\in M$ with
$m_5+m_6+m_7=k$. According to (9.2.5),
$E_{n+n_1+1,n_1+1}=y_{n_1+1}\ptl_{x_{n_1+1}}$. So
$$E_{n+n_1+1,n_1+1}^{m_5+m_6}[\eta^{m_5+m_6}(x_{n_1}^{m_6}y_{n_1+1}^{m_7-m_5})]=(m_5+m_6)!
x_{n_1}^{m_6}y_{n_1+1}^{k+m_6}\in M.\eqno(9.2.77)$$ Moreover,
(9.2.4) yields that $(E_{n_1,n+n_1+1}+E_{n_1+1,n+n_1})|_{\msr
B}=\ptl_{x_{n_1}}\ptl_{y_{n_1+1}}+x_{n_1+1}\ptl_{y_{n_1}}$ and
$$(E_{n_1,n+n_1+1}+E_{n_1+1,n+n_1})^{m_6}(x_{n_1}^{m_6}y_{n_1+1}^{k+m_6})=m_6![\prod_{r=1}^{m_6}(k+r)]
y_{n_1+1}^k\in M. \eqno(9.2.78)$$ Assume some
$\eta^{m_5+m_6}(y_{n_1+2}^{m_6}x_{n_1+1}^{m_7-m_5})\in M$ with
$m_5+m_6+m_7=k$. By (9.2.12) and (9.2.15),
$$E_{n_1+1,n+n_1+1}^{m_5+m_6}[\eta^{m_5+m_6}(y_{n_1+2}^{m_6}x_{n_1+1}^{m_7-m_5})]=(m_5+m_6)!
y_{n_1+2}^{m_6}x_{n_1+1}^{k+m_6}\in M.\eqno(9.2.79)$$ Observe
$$(E_{n+n_1+2,n_1+1}+E_{n+n_1+1,n_1+2})|_{\msr B}=
\ptl_{x_{n_1+1}}\ptl_{y_{n_1+2}}+y_{n_1+1}\ptl_{x_{n_1+2}}\eqno(9.2.80)$$
by (9.2.5). Hence
$$(E_{n+n_1+2,n_1+1}+E_{n+n_1+1,n_1+2})^{m_6}(y_{n_1+2}^{m_6}x_{n_1+1}^{k+m_6})
=m_6![\prod_{r=1}^{m_6}(k+r)] x_{n_1+1}^k\in M.\eqno(9.2.81)$$

Expressions (9.2.59)-(9.2.66), (9.2.76), (9.2.78) and (9.2.81) show
that we always have $x_{n_1+1}^k\in M$. According to the arguments
in Case 1, (7.2.69) holds. Then (9.2.74) holds for
$r_2-r_1+s_1-s_2=k$. Applying $U(\msr G_1)$ to it, we get (9.2.70).
There ${\msr B}_{\la k\ra}$ is an irreducible
$sp(2n,\mbb{F})$-module.\psp

This completes the proof of the theorem when $n_1<n_2<n$.

\subsection{Proof of the Theorem When $n_1=n_2<n$}

In this subsection, we assume $n_1=n_2<n$.

Under the assumption,
$$\eta=\sum_{i=1}^{n_1}y_i\ptl_{x_i}+\sum_{s=n_2+1}^n
x_s\ptl_{y_s}.\eqno(9.2.82)$$ First we consider the subcase
$1<n_1<n-1$. Lemma 6.6.1 says that the nonzero weight vectors in
$$\mbox{Span}\{x_{n_1}^{m_1}y_{n_1}^{m_2}\zeta_1^{m_3+1},
x_{n_1+1}^{m_1}y_{n_1+1}^{m_2}\zeta_2^{m_3+1},
\eta^{m_3}(x_{n_1}^{m_1}y_{n_1+1}^{m_2})\mid
m_i\in\mbb{N}\}\eqno(9.2.83)$$ are all the $sl(n,\mbb F)$-singular
vectors
 in ${\msr B}$, where
$$\zeta_1=x_{n_1-1}y_{n_1}-x_{n_1}y_{n_1-1},\;\;\zeta_2=x_{n_1+1}y_{n_1+2}-x_{n_1+2}y_{n_1+1}.
\eqno(9.2.84)$$ Fix $k\in\mbb{N}+1$. Then the $sl(n,\mbb
F)$-singular vectors  in ${\msr B}_{\la-k\ra}$ are
\begin{eqnarray*}\qquad&
&\{x_{n_1}^{k+m_1}y_{n_1}^{m_1}\zeta_1^{m_2+1},
x_{n_1+1}^{m_1}y_{n_1+1}^{k+m_1}\zeta_2^{m_2+1},
\eta^{m_3}(x_{n_1}^{m_4}y_{n_1+1}^{m_5}) \\
& &\mid
m_i\in\mbb{N};m_4+m_5-2m_3=k\}.\hspace{7.5cm}(9.2.85)\end{eqnarray*}
Let $M$ be a nonzero $sp(2n,\mbb{F})$-submodule of ${\msr
B}_{\la-k\ra}$. As a ${\msr K}$-module, $M$ contains an $sl(n,\mbb
F)$-singular vector. Suppose some
$x_{n_1}^{k+m_1}y_{n_1}^{m_1}\zeta_1^{m_2+1}\in M$. Note
$E_{n_1,n+n_1}|_{\msr B}=\ptl_{x_{n_1}}\ptl_{y_{n_1}}$ by (9.2.11),
and so
\begin{eqnarray*} & &
E_{n_1,n+n_1}(x_{n_1}^{k+m_1}y_{n_1}^{m_1}\zeta_1^{m_2}) \\ &=&
(k+m_1)m_1x_{n_1}^{k+m_1-1}y_{n_1}^{m_1-1}\zeta_1^{m_2}-m_2(m_2-1)
x_{n_1}^{k+m_1}y_{n_1}^{m_1}x_{n_1-1}y_{n_1-1}\zeta_1^{m_2-2}
\\ & &+(k+m_1)m_2x_{n_1}^{k+m_1-1}y_{n_1}^{m_1}x_{n_1-1}\zeta_1^{m_2-1}
-m_1m_2x_{n_1}^{k+m_1}y_{n_1}^{m_1-1}y_{n_1-1}\zeta_1^{m_2-1}.
\hspace{1.2cm}(9.2.86)\end{eqnarray*} Moreover,
$$(E_{n_1-1,n_1}-E_{n+n_1,n+n_1-1})|_{\msr
B}=-(x_{n_1}\ptl_{x_{n_1-1}}+y_{n_1}\ptl_{y_{n_1-1}})\eqno(9.2.87)$$
by (9.2.1)-(9.2.3). Thus
\begin{eqnarray*} \qquad\qquad& &
(E_{n_1-1,n_1}-E_{n+n_1,n+n_1-1})^2E_{n_1,n+n_1}(x_{n_1}^{k+m_1}y_{n_1}^{m_1}\zeta_1^{m_2})
\\ &=&
-2m_2(m_2-1) x_{n_1}^{k+m_1+1}y_{n_1}^{m_1+1}\zeta_1^{m_2-2}\in
M.\hspace{4.7cm}(9.2.88)\end{eqnarray*} Hence
$$x_{n_1}^{k+m_1+1}y_{n_1}^{m_1+1}\zeta_1^{m_2-2}\in
M\qquad\mbox{if}\;\;m_2>1.\eqno(9.2.89)$$ Furthermore,
$$(E_{n_1-1,n_1}-E_{n+n_1,n+n_1-1})E_{n_1,n+n_1}(x_{n_1}^{k+m_1}y_{n_1}^{m_1}\zeta_1)
= -kx_{n_1}^{k+m_1}y_{n_1}^{m_1}\in M.\eqno(9.2.90)$$  So we always
have $x_{n_1}^{k+m}y_{n_1}^m\in M$ for some $m\in\mbb{N}$ by
induction on $m_2$.

 Observe
 $$E_{n_1,n+n_1}(x_{n_1}^{k+m}y_{n_1}^m)=\ptl_{x_{n_1}}\ptl_{y_{n_1}}(x_{n_1}^{k+m}y_{n_1}^m)=m![\prod_{r=1}^m
 (k+r)]x_{n_1}^k\eqno(9.2.91)$$
 by (9.2.11). Thus
 $$x_{n_1}^k\in M.\eqno(9.2.92)$$
Symmetrically, if some
$x_{n_1+1}^{m_1}y_{n_1+1}^{k+m_1}\zeta_2^{m_2+1}\in M$, we have
$y_{n_1+1}^k\in M$. But
$$(E_{n+n_1+1,n_1}+E_{n+n_1,n_1+1})|_{\msr
B}=-x_{n_1}\ptl_{y_{n_1+1}}+y_{n_1}\ptl_{x_{n_1+1}}\eqno(9.2.93)$$
by (9.2.5), which gives
$$(E_{n+n_1+1,n_1}+E_{n+n_1,n_1+1})^k(y_{n_1+1}^k)=(-1)^kk!x_{n_1}^k\in
M.\eqno(9.2.94)$$ Thus (9.2.92) holds again.

Assume that some $\eta^{m_3}(x_{n_1}^{m_4}y_{n_1+1}^{m_5})\in M$
with $m_4+m_5-2m_3=k$. Note there exists $r_1,r_2\in\mbb{N}$ such
that $r_1+r_2=m_3$ and $2r_1\leq m_4,\;2r_2\leq m_5$. Moreover, $
E_{n_1,n+n_1}|_{\msr B}=\ptl_{x_{n_1}}\ptl_{y_{n_1}}$ by (9.2.11)
and $ E_{n+n_1+1,n_1+1}|_{\msr B}=\ptl_{x_{n_1+1}}\ptl_{y_{n_1+1}}$
by (9.2.5). Thus
\begin{eqnarray*}\qquad & &E_{n_1,n+n_1}^{r_1}E_{n+n_1+1,n_1+1}^{r_2}
[\eta^{m_3}(x_{n_1}^{m_4}y_{n_1+1}^{m_5})]\\
&=&m_3![\prod_{s_1=0}^{2r_1-1}(m_4-s_1)]
[\prod_{s_2=0}^{2r_2-1}(m_5-s_2)]x_{n_1}^{m_4-2r_1}y_{n_1+1}^{m_5-2r_2}\in
M. \hspace{2.9cm}(9.2.95)\end{eqnarray*} Furthermore, (9.2.11) gives
$(E_{n+n_1+1,n_1}+E_{n+n_1,n_1+1})|_{\msr
B}=-x_{n_1}\ptl_{y_{n_1+1}}+y_{n_1}\ptl_{x_{n_1+1}}$ , and so
\begin{eqnarray*}\hspace{2cm}& &(E_{n+n_1+1,n_1}+E_{n+n_1,n_1+1})^{m_5-2r_2}
(x_{n_1}^{m_4-2r_1}y_{n_1+1}^{m_5-2r_2})\\
&=&(-1)^{m_5-2r_2}(m_5-2r_2)!x_{n_1}^k\in
M.\hspace{5.9cm}(9.2.96)\end{eqnarray*} Thus we always have
$x_{n_1}^k\in M$.

Now
$$(E_{n_1,n+n_1+1}+E_{n_1+1,n+n_1})|_{\msr
B}=-y_{n_1+1}\ptl_{x_{n_1}}+x_{n_1+1}\ptl_{y_{n_1}}\eqno(9.2.97)$$
by (9.2.11). For any $r\in\mbb{N}+1$,
$$\frac{(-1)^r}{\prod_{s=0}^{r-1}(k-r)}(E_{n_1,n+n_1+1}+E_{n_1+1,n+n_1})^r(x_{n_1}^k)
=x_{n_1}^{k-r}y_{n_1+1}^r\in M.\eqno(9.2.98)$$ Set
\begin{eqnarray*}\hspace{1cm}{\msr G}_4&=&
\sum_{i,j=1}^{n_1}\mbb{F}(E_{i,j}-E_{n+j,n+i})+\sum_{i=1}^{n_1}(\mbb{F}E_{i,n+i}+\mbb{F}E_{n+i,i})\\
& &+\sum_{1\leq i<j\leq n_1
}[\mbb{F}(E_{i,n+j}+E_{n+j,i})+\mbb{F}(E_{n+i,j}+E_{n+j,i})]\hspace{3.1cm}(9.2.99)\end{eqnarray*}
and
$$M^4=\mbb
F[x_1,...,x_{n_1},y_1,...,y_{n_1}],\;\;M^4_{\la \ell\ra}= M^4\bigcap
\msr B_{\la \ell\ra}\;\;\for\;\;\ell\in\mbb Z.\eqno(9.2.100)$$

If $k\geq 2$ and $r\in\ol{1,k-1}$, then
$${M}^4_{\la -k+r\ra}{M}^3_{\la -r\ra}=U({\msr G}_3)U({\msr
G}_4)(x_{n_1}^{k-r}y_{n_1+1}^r)\subset M\eqno(9.2.101)$$ (cf.
(9.2.72) and (9.2.73)) because ${M}^4_{\la -k+r\ra}$ is an
irreducible ${\msr G}_4$-module and ${M}^3_{\la -r\ra}$ is an
irreducible ${\msr G}_3$-module by Subsection 9.2.2. Moreover,
$${M}^4_{\la -k\ra}=U({\msr G}_4)(x_{n_1}^k),\;{M}^3_{\la -k\ra}=U({\msr G}_3)
(y_{n_1+1}^k)\subset M.\eqno(9.2.102)$$ Furthermore,
$${M}^4_{\la -k\ra}{M}^3_{\la 0\ra}=U({\msr G}_3)U({\msr
G}_4)(x_{n_1}^k)\subset M \;\;\mbox{if}\;\;n_1=n-1\eqno(9.2.103)$$
and
$${M}^4_{\la 0\ra}{M}^3_{\la -k\ra}=U({\msr G}_3)U({\msr G}_4)
(y_{n_1+1}^k)\subset M \;\;\mbox{if}\;\;n_1=1.\eqno(9.2.104)$$

Note
$$(E_{r,i}-E_{n+i,n+r})|_{\msr
B}=y_iy_r-x_ix_r\qquad\for\;\;i\in\ol{1,n_1},\;r\in\ol{n_1+1,n}\eqno(9.2.105)$$
by (9.2.1)-(9.2.3). In particular, if $k>1$ or $n_1=1$, we have
$$(E_{n_1+1,n_1}-E_{n+n_1,n+n_1+1})(x_{n_1}^k)=y_{n_1}x_{n_1}^ky_{n_1+1}-x_{n_1}^{k+1}x_{n_1+1}\in
M.\eqno(9.2.106)$$ Since
$$y_{n_1}x_{n_1}^ky_{n_1+1}\in {M}^4_{\la -k+1\ra}{M}^3_{\la
-1\ra}\subset M,\eqno(9.2.107)$$ we get $$ x_{n_1}^{k+1}x_{n_1+1}\in
M.\eqno(9.2.108)$$ Suppose $k=1$ and $n_1>1$. By (9.2.102),
$$\zeta_1x_{n_1}=(x_{n_1-1}y_{n_1}-x_{n_1}y_{n_1-1})x_{n_1}\in
M.\eqno(9.2.109)$$ Observe
$$(E_{n_1+1,n+n_1-1}+E_{n_1-1,n+n_1+1})|_{\msr
B}=x_{n_1+1}\ptl_{y_{n_1-1}}-y_{n_1+1}\ptl_{x_{n_1-1}}\eqno(9.2.110)$$
by (9.2.4). So
$$-(E_{n_1+1,n+n_1-1}+E_{n_1-1,n+n_1+1})(\zeta_1x_{n_1})=x_{n_1}^2x_{n_1+1}
-x_{n_1}y_{n_1}y_{n_1+1}\in M.\eqno(9.2.111)$$ On the other hand,
(9.2.5) gives
$$(E_{n+i,j}+E_{n+j,i})|_{\msr
B}=-(x_iy_j+x_jy_i)\qquad\for\;\;i,j\in\ol{1,n_1},\eqno(9.2.112)$$
which implies
$$-E_{n+n_1,n_1}(y_{n_1+1})=x_{n_1}y_{n_1}y_{n_1+1}\in
M.\eqno(9.2.113)$$ By (9.2.111), we have $x_{n_1}^2x_{n_1+1}\in M.$
So (9.2.108) always holds.

By Subsection 9.2.2,
$${M}^4_{\la -k-1\ra}{M}^3_{\la 1\ra}=U({\msr G}_3)U({\msr
G}_4)(x_{n_1}^{k+1}x_{n_1+1})\subset M.\eqno(9.2.114)$$ Suppose
$${M}^4_{\la -k-i\ra}{M}^3_{\la i\ra}\subset
M\eqno(9.2.115)$$ for $1\leq i\leq m$. Then
\begin{eqnarray*}\qquad\qquad& &(E_{n_1+1,n_1}-E_{n+n_1,n+n_1+1})(x_{n_1}^{k+m}x_{n_1+1}^m)
\\ &=&y_{n_1}x_{n_1}^{k+m}x_{n_1+1}^my_{n_1+1}-x_{n_1}^{k+m+1}x_{n_1+1}^{m+1}\in
M\hspace{4.7cm}(9.2.116)\end{eqnarray*} by (9.2.105). If $m>1$, we
have
$$y_{n_1}x_{n_1}^{k+m}x_{n_1+1}^my_{n_1+1}\in
{M}^4_{\la -k-(m-1)\ra}{M}^3_{\la m-1\ra}\subset M.\eqno(9.2.117)$$
Note
$$(E_{r,n+s}+E_{s,n+r})|_{\msr
B}=-(x_ry_s+x_sy_r)\qquad\for\;\;r,s\in\ol{n_1+1,n}\eqno(9.2.118)$$
by (9.2.4). If $m=1$, we have
$$y_{n_1}x_{n_1}^{k+1}x_{n_1+1}y_{n_1+1}
=-E_{n_1+1,n+n_1+1}(y_{n_1}x_{n_1}^{k+1})\subset
E_{n_1+1,n+n_1+1}({M}^1_{\la -k\ra})\subset M.\eqno(9.2.119)$$ Then
(9.2.116), (9.2.117) and (9.2.119) give
$$x_{n_1}^{k+m+1}x_{n_1+1}^{m+1}\in
M.\eqno(9.2.120)$$ Furthermore,
$${M}^4_{\la -k-m-1\ra}{M}^3_{\la m+1\ra}=U({\msr G}_3)U({\msr
G}_4)(x_{n_1}^{k+m+1}x_{n_1+m+1})\subset M.\eqno(9.2.121)$$ Thus
(9.2.115) holds for any $i\in\mbb{N}+1$. Symmetrically, we have
$${M}^4_{\la i\ra}{M}^3_{\la -k-i\ra}\subset M\qquad\for\;\;i\in\mbb{N}+1.\eqno(9.2.122)$$

Suppose $n_1<n-1$. Then $x_{n_1}^{k+1}x_{n_1+1}\zeta_2\in M$ by
(9.2.115). Moreover,
$$(k+1)y_{n_1}x_{n_1}^ky_{n_1+1}\zeta_2=-(k+1)E_{n+n_1,n_1}(x_{n_1}^{k-1}y_{n_1+1}\zeta_2)\in
M\eqno(9.2.123)$$ by (9.2.5) and (9.2.101). According to
(9.2.1)-(9.2.3),
$$(E_{n_1,n_1+1}-E_{n+n_1+1,n+n_1})|_{\msr
B}=\ptl_{x_{n_1}}\ptl_{x_{n_1+1}}-\ptl_{y_{n_1}}\ptl_{y_{n_1+1}}.\eqno(9.2.124)$$
Thus
\begin{eqnarray*} \qquad& &(E_{n_1,n_1+1}-E_{n+n_1+1,n+n_1})
[(x_{n_1}^{k+1}x_{n_1+1}-(k+1)y_{n_1}x_{n_1}^ky_{n_1+1})\zeta_2]
\\ &=& 3(k+1)x_{n_1}^k\zeta_2\in
M\hspace{9.2cm}(9.2.125)\end{eqnarray*} by (9.2.84). Hence
$${M}^4_{\la-k\ra}{M}^3_{\la 0\ra}=U({\msr G}_3)U({\msr
G}_4)(x_{n_1}^k)+U({\msr G}_3)U({\msr G}_4)(x_{n_1}^k\zeta_2)\subset
M\eqno(9.2.126)$$ by (9.2.92) and (9.2.125). Symmetrically,
$${M}^4_{\la 0 \ra}{M}^3_{\la -k\ra}\subset
M.\eqno(9.2.127)$$ By (9.2.101)-(9.2.104), (9.2.115), (9.2.121),
(9.2.122), (9.2.126) and (9.2.127),
$${M}^4_{\la -k-r \ra}{M}^3_{\la r\ra}\subset
M\qquad\for\;\;r\in\mbb{Z}.\eqno(9.2.128)$$ Therefore,
$${\msr B}_{\la -k\ra}=\bigoplus_{r\in\mbb{Z}}
{M}^1_{\la -k-r \ra}{M}^2_{\la r\ra}\subset M.\eqno(9.2.129)$$ We
get $M={\msr B}_{\la -k\ra}$; that is, ${\msr B}_{\la -k\ra}$ is an
irreducible $sp(2n,\mbb{F})$-module. We can similarly prove that
${\msr B}_{\la k\ra}$ is an irreducible $sp(2n,\mbb{F})$-module.

Finally, we study ${\msr B}_{\la 0\ra}$. We first consider the
generic case $1<n_1<n-1$. Set
\begin{eqnarray*} \qquad
{\msr B}_{\la 0,1\ra}&=&\mbox{Span}\{[\prod_{1\leq r\leq s\leq
n_1\;\mbox{or}\;n_1+1\leq r\leq s\leq n}(x_ry_s+x_sy_r)^{l_{r,s}}]\\
& &\qquad\;\;\times
[\prod_{p=1}^{n_1}\prod_{q=n_1+1}^n(x_px_q-y_py_q)^{k_{p,q}}] \mid
l_{r,s},k_{p,q}\in\mbb{N}\}\hspace{2.6cm}(9.2.130)\end{eqnarray*}
and \begin{eqnarray*}\qquad\qquad{\msr B}_{\la 0,2\ra}&=&\sum_{1\leq
r< s\leq n_1\;\mbox{or}\;n_1+1\leq r< s\leq n}{\msr B}_{\la
0,1\ra}(x_ry_s-x_sy_r)\\ & &+\sum_{p=1}^{n_1}\sum_{q=n_1+1}^n{\msr
B}_{\la 0,1\ra}(x_px_q+y_py_q).\hspace{5cm}(9.2.131)\end{eqnarray*}
 We
want to prove that ${\msr B}_{\la 0,1\ra}$ and ${\msr B}_{\la
0,2\ra}$ form $sp(2n,\mbb{F})$-submodules.

Let \begin{eqnarray*}\qquad {\msr G}^\ast&=&\sum_{1\leq r\leq s\leq
n_1}\mbb{F}(E_{n+s,r}+E_{n+r,s})+\sum_{n_1+1\leq p\leq q\leq
n}\mbb{F}(E_{p,n+q}+E_{q,n+p})\\ &
&+\sum_{r=1}^{n_1}\sum_{p=n_1+1}^n
\mbb{F}(E_{p,r}-E_{n+r,n+p})\hspace{6.5cm}(9.2.132)\end{eqnarray*}
and
\begin{eqnarray*}\td{\msr
G}&=&\sum_{i,j=1}^{n_1}\mbb{F}(E_{i,j}-E_{n+j,n+i})+\sum_{r,s=n_1+1}^n\mbb{F}
(E_{r,s}-E_{n+s,n+r})+\sum_{1\leq r\leq s\leq
n_1}\mbb{F}(E_{r,n+s}+E_{s,n+r})\\ & &+\sum_{n_1+1\leq p\leq q\leq
n}\mbb{F}(E_{n+q,p}+E_{n+p,q})+\sum_{r=1}^{n_1}\sum_{p=n_1+1}^n
[\mbb{F}(E_{r,p}-E_{n+p,n+r})\\
&
&+\mbb{F}(E_{r,n+p}+E_{p,n+r})+\mbb{F}(E_{n+r,p}+E_{n+p,r})].\hspace{5.4cm}(9.2.133)\end{eqnarray*}
Then ${\msr G}^\ast$ and $\td{\msr G}$ are Lie subalgebras of
$sp(2n,\mbb{F})$ and $sp(2n,\mbb{F})={\msr G}^\ast\oplus \td{\msr
G}.$ By PBW Theorem (Theorem 5.1.1), $U(sp(2n,\mbb{F}))=U({\msr
G}^\ast)U(\td{\msr G}).$

Expressions (9.2.105), (9.2.112) and (9.2.118) show
$$U({\msr G}^\ast)|_{\msr B}={\msr B}_{\la
0,1\ra}\;\;\mbox{as multiplication operators on}\;\;{\msr
B}.\eqno(9.2.134)$$ Moreover,
$$(E_{r,s}-E_{n+s,n+r})|_{\msr
B}=x_r\ptl_{x_s}+y_r\ptl_{y_s}+\dlt_{r,s},\eqno(9.2.135)$$
$$(E_{n+r,s}+E_{n+s,r})|_{\msr
B}=\ptl_{x_r}\ptl_{y_s}+\ptl_{x_s}\ptl_{y_r},\eqno(9.2.136)$$
$$(E_{n+r,i}+E_{n+i,r})|_{\msr
B}=-x_i\ptl_{y_r}+y_i\ptl_{x_r},\eqno(9.2.137)$$
$$(E_{i,n+r}+E_{r,n+i})|_{\msr
B}=-y_r\ptl_{x_i}+x_r\ptl_{y_i},\eqno(9.2.138)$$
$$(E_{i,r}-E_{n+r,n+i})|_{\msr
B}=\ptl_{x_i}\ptl_{x_r}-\ptl_{y_i}\ptl_{y_r}\eqno(9.2.139)$$ for
$i\in\ol{1,n_1}$ and $r,s\in\ol{n_1+1,n}$. According to (9.2.11),
(9.2.14), (9.2.133) and (9.2.135)-(9.2.139), $U(\td{\msr
G})(1)=\mbb{F}.$ Thus
$${\msr B}_{\la 0,1\ra}=U({\msr
G}^\ast)(1)=U(sp(2n,\mbb{F}))(1)\eqno(9.2.140)$$ forms an
$sp(2n,\mbb{F})$-submodule.

Let
$$W=\sum_{1\leq r< s\leq n_1\;\mbox{or}\;n_1+1\leq r< s\leq
n}\mbb{F}(x_ry_s-x_sy_r)+
\sum_{p=1}^{n_1}\sum_{q=n_1+1}^n\mbb{F}(x_px_q+y_py_q).\eqno(9.2.141)$$
By (9.2.11), (9.2.14), (9.2.133) and (9.2.135)-(9.2.139), we can
verify that $W$ forms an irreducible $\td{\msr G}$-submodule. Hence
$${\msr B}_{\la 0,2\ra}=U({\msr
G}^\ast)(W)= U(sp(2n,\mbb{F}))(W)\eqno(9.2.142)$$ forms an
$sp(2n,\mbb{F})$-submodule. Moreover,
$${\msr B}_{\la 0,1\ra}\bigcap W=\{0\}.\eqno(9.2.143)$$

Next we want to prove that ${\msr B}_{\la 0,1\ra}$ and ${\msr
B}_{\la 0,2\ra}$ are irreducible $sp(2n,\mbb{F})$-submodules.
According to (9.2.83), the $sl(n,\mbb F)$-singular vectors in ${\msr
B}_{\la 0 \ra}$ are
\begin{eqnarray*}\qquad&
&\{x_{n_1}^{m_1}y_{n_1}^{m_1}\zeta_1^{m_2+1},
x_{n_1+1}^{m_1}y_{n_1+1}^{m_1}\zeta_2^{m_2+1},
\eta^{m_3}(x_{n_1}^{m_4}y_{n_1+1}^{m_5}) \\
& &\mid
m_i\in\mbb{N};m_4+m_5=2m_3\}.\hspace{7.8cm}(9.2.144)\end{eqnarray*}
Let $M$ be a nonzero submodule of ${\msr B}_{\la 0,1\ra}$. Then $M$
contains an $sl(n,\mbb F)$-singular vector. Suppose some
$x_{n_1}^{m_1}y_{n_1}^{m_1}\zeta_1^{m_2}\in M$. By
(9.2.86)-(9.2.89), we can assume $m_2=0,1$. If $m_2=0$, (9.2.91)
yields $1\in M$. Then $M={\msr B}_{\la 0,1\ra}$ by (9.2.140).
Suppose $m_2=1$. We have $E_{n_1,n+n_1}|_{\msr
B}=\ptl_{x_{n_1}}\ptl_{y_{n_1}}$ by (9.2.11), and
$$E_{n_1,n+n_1}[x_{n_1}^{m_1}y_{n_1}^{m_1}\zeta_1]
=m_1(m_1+m_2)x_{n_1}^{m_1-1}y_{n_1}^{m_1-1}\zeta_1\eqno(9.2.145)$$
by (9.2.86). By induction on $m_1$, we have $\zeta_1\in M\subset
{\msr B}_{\la 0,1\ra}$, which contradicts (9.2.143). Similarly, if
some $x_{n_1+1}^{m_1}y_{n_1+1}^{m_1}\zeta_2^{m_2+1}\in M$, we have
$M= {\msr B}_{\la 0,1\ra}$. Assume some
$\eta^{m_3}(x_{n_1}^{m_4}y_{n_1+1}^{m_5})$ in $M$ with
$m_4+m_5=2m_3$. Note that $m_4$ and $m_5$ are both even or odd. If
$m_4=2r_1$ and $m_5=2r_2$ are even, then (9.2.95) gives $1\in M$, or
equivalently $M={\msr B}_{\la 0,1\ra}$. Suppose that $m_4=2r_1+1$
and $m_5=2r_2+1$ are odd. As (9.2.95), we find
$$\eta(x_{n_1}y_{n_1+1})=x_{n_1}x_{n_1+1}
+y_{n_1}y_{n_1+1}\in M\subset {\msr B}_{\la 0,1\ra}\eqno(9.2.146)$$
by (9.2.82), which contradicts (9.2.143) again. Thus we always have
$M={\msr B}_{\la 0,1\ra}$; that is, ${\msr B}_{\la 0,1\ra}$ is
irreducible. Similarly, we can prove that ${\msr B}_{\la 0,2\ra}$ is
irreducible.

If $n_1=1$ and $n=2$, we let $${\msr B}_{\la
0,1\ra}=\mbox{Span}\{[\prod_{i=1}^n(x_iy_i)^{m_i}]
(x_1x_2-y_1y_2)^{m_3}] \mid m_i\in\mbb{N}\}\eqno(9.2.147)$$
 and
${\msr B}_{\la 0,2\ra}={\msr B}_{\la 0,1\ra}(x_1x_2+y_1y_2).$ When
$n_1=1$ and $n>2$, we set
\begin{eqnarray*}\qquad\qquad {\msr B}_{\la
0,1\ra}&=&\mbox{Span}\{[(x_1y_1)^l\prod_{2\leq r\leq s\leq
n}(x_ry_s+x_sy_r)^{l_{r,s}}]\\ & &
[\prod_{q=2}^n(x_1x_q-y_1y_q)^{k_q}] \mid l,
l_{r,s},k_{q}\in\mbb{N}\}\hspace{4.5cm}(9.2.148)\end{eqnarray*} and
$${\msr B}_{\la
0,2\ra}=\sum_{2\leq r< s\leq n}{\msr B}_{\la
0,1\ra}(x_ry_s-x_sy_r)+\sum_{q=2}^n{\msr B}_{\la
0,1\ra}(x_1x_q+y_1y_q).\eqno(9.2.149)$$ In the case $1<n_1=n-1$, we
put
\begin{eqnarray*} \qquad
{\msr B}_{\la 0,1\ra}&=&\mbox{Span}\{(x_ny_n)^l[\prod_{1\leq r\leq
s\leq n-1}(x_ry_s+x_sy_r)^{l_{r,s}}]\\ & &\qquad\;\;\times
[\prod_{p=1}^{n-1}(x_px_n-y_py_n)^{k_{p}}] \mid
l,l_{r,s},k_{p}\in\mbb{N}\}\hspace{3.8cm}(9.2.150)\end{eqnarray*}
and
$${\msr B}_{\la
0,2\ra}=\sum_{1\leq r< s\leq n_1}{\msr B}_{\la
0,1\ra}(x_ry_s-x_sy_r)+\sum_{p=1}^{n_1}{\msr B}_{\la
0,1\ra}(x_px_n+y_py_n).\eqno(9.2.151)$$ By the arguments similar to
part of those in the above proof of the irreducibility of ${\msr
B}_{\la 0,1\ra}$ when $1<n_1<n-1$, we can show that ${\msr B}_{\la
0,1\ra}$ and ${\msr B}_{\la 0,2\ra}$ are irreducible.

Recall the symmetric bilinear form $(\cdot|\cdot)$ in (6.4.46). Now
$1$ is a non-isotropic element in
 ${\msr B}_{\la 0,1\ra}$ and $x_{n_1}x_{n_1+1}+y_{n_1}y_{n_1+1}$ a non-isotropic element in ${\msr B}_{\la 0,2\ra}$.
By Lemma 6.3.2, ${\msr B}_{\la 0,r\ra}\bigcap ({\msr B}_{\la
0,r\ra})^\bot$ is a proper $sp(2n,\mbb F)$-submodule of ${\msr
B}_{\la 0,r\ra}$, the irreducibility of ${\msr B}_{\la 0,r\ra}$
implies  ${\msr B}_{\la 0,r\ra}\bigcap ({\msr B}_{\la
0,r\ra})^\bot=\{0\}$; that is, the bilinear form $(\cdot|\cdot)$ is
nondegenerate on ${\msr B}_{\la 0,1\ra}$ and ${\msr B}_{\la
0,2\ra}$.

 Since $(1|{\msr B}_{\la
0,2\ra})=\{0\}$, ${\msr B}_{\la 0,1\ra}$ is orthogonal to ${\msr
B}_{\la 0,2\ra}$. Thus the symmetric bilinear form $(\cdot|\cdot)$
on ${\msr B}_{\la 0,1\ra}+{\msr B}_{\la 0,2\ra}$ is nondegenerate.
Recall that $\msr B_k$ is the subspace of homogeneous polynomials in
$\msr B$ with degree $k$ and $\msr B=\bigoplus_{k=0}^\infty\msr
B_k$. According to (6.4.46),
$$(\msr B_{k_1}|\msr B_{k_2})=\{0\}\qquad\mbox{if}\;\;k_1\neq
k_2.\eqno(9.2.152)$$ Moreover,
$${\msr B}_{\la 0\ra}=\bigoplus_{k=0}^\infty \msr B_k\bigcap{\msr B}_{\la
0\ra},\eqno(9.2.153)$$
$${\msr B}_{\la 0,r\ra}=\bigoplus_{k=0}^\infty \msr B_k\bigcap{\msr B}_{\la
0,r\ra}.\eqno(9.2.154)$$ The above three expressions show
$${\msr B}_{\la 0\ra}=({\msr B}_{\la 0,1\ra}+{\msr B}_{\la
0,2\ra})\oplus ({\msr B}_{\la 0,1\ra}+{\msr B}_{\la
0,2\ra})^\perp\bigcap {\msr B}_{\la 0\ra}.\eqno(9.2.155)$$ If
$({\msr B}_{\la 0,1\ra}+{\msr B}_{\la 0,2\ra})^\perp\bigcap {\msr
B}_{\la 0\ra}\neq\{0\}$, then it contains an $sl(n,\mbb F)$-singular
vector. Our above arguments in proving the irreducibility of ${\msr
B}_{\la 0,1\ra}$ show that it contains either ${\msr B}_{\la
0,1\ra}$ or ${\msr B}_{\la 0,2\ra}$, which is absurd. Therefore,
${\msr B}_{\la 0\ra}={\msr B}_{\la 0,1\ra}\oplus{\msr B}_{\la
0,2\ra}$ is an orthogonal decomposition of irreducible
$sp(2n,\mbb{F})$-submodules.\psp

This completes the proof of the theorem when $n_1=n_2<n$. Therefore,
the theorem holds by Subsections 9.2.2-9.2.4.\psp

\section{Projective Oscillator Representations}

Let us go back to Section 6.8 and  assume $n=2m+1>1$ is an odd
integer. In this section, we study the restrictions of the
projective representations
 on $sp(2m+2,\mbb F)$.

Set
$$A_{i,j}=E_{i,j}-E_{m+1+j,m+1+i},\;B^s_{i,j}=E_{i,m+1+j}+E_{j,m+1+i},\;C^s_{i,j}=
E_{m+1+i,j}+E_{m+1+j,i}\eqno(9.3.1)$$ for $i,j\in\ol{1,m+1}$. Then
 the symplectic
Lie algebra
$$sp(2m+2,\mbb{F})=\sum_{i,j=1}^{m+1}\mbb FA_{i,j}+\sum_{1\leq p\leq
q\leq m+1}(\mbb FB^s_{p,q}+\mbb FC^s_{q,p}).\eqno(9.3.2)$$ For
convenience, we rednote
$$x_0=x_{m+1},\;\;y_i=x_{m+1+i}\qquad\for\;i\in\ol{1,m}.\eqno(9.3.3)$$
In particular,
$$D=\sum_{i=0}^mx_i\ptl_{x_i}+\sum_{r=1}^my_r\ptl_{y_r}.\eqno(9.3.4)$$
According to (6.7.80) and (6.7.81), we have the representation
$\pi_\kappa$ of $sp(2m+2,\mbb{F})$:
$$\pi_\kappa(A_{i,j})=x_i\ptl_{x_j}-y_j\ptl_{y_i},\;\;\pi_\kappa(B^s_{i,j})
=x_i\ptl_{y_j}+x_j\ptl_{y_i},\;\;\pi_\kappa(C^s_{i,j})=y_i\ptl_{x_j}+y_j\ptl_{x_i},\eqno(9.3.5)$$
$$\pi_\kappa
(A_{i,m+1})=x_i\ptl_{x_0}+\ptl_{y_i},\;\;\pi_\kappa(C^s_{m+1,i})
=y_i\ptl_{x_0}-\ptl_{x_i},\;\;\pi_\kappa(C^s_{m+1,m+1})=-2\ptl_{x_0},\eqno(9.3.6)$$
$$\pi_\kappa(A_{m+1,i})=x_0\ptl_{x_i}-y_i(D+\kappa),\;\;\pi_\kappa(B^s_{i,m+1})=x_0\ptl_{y_i}+x_i(D+\kappa),
\eqno(9.3.7)$$
$$\pi_\kappa(B^s_{m+1,m+1})=2x_0(D+\kappa),\;\;\pi_\kappa(A_{m+1,m+1})=D+x_0\ptl_{x_0}+\kappa\eqno(9.3.8)$$
 for
$i,j\in\ol{1,m}$.

Denote $${\msr K}=\sum_{i,j=1}^m\mbb FA_{i,j}+\sum_{1\leq p\leq
q\leq m}(\mbb FB^s_{p,q}+\mbb FC^s_{q,p}),\eqno(9.3.9)$$ which is a
Lie subalgebra of $sp(2m+2,\mbb{F})$ isomorphic to $sp(2m,\mbb{F})$.
 Again we take
 $$\msr B=\mbb{F}[x_1,...,x_m,y_1,...,y_m],\;\;\msr
 B'=\mbb{F}[x_0,x_1,...,x_m,y_1,...,y_m].\eqno(9.3.10)$$
Let $\msr B_k$ be the subspace of homogeneous polynomials in $\msr
B$ with degree $k$ and let $\msr B_k'$ be the subspace of
homogeneous polynomials in $\msr B'$ with degree $k$. Set
$${\msr B}'_{(\ell)}=\sum_{i=0}^\ell {\msr
B}'_i\qquad\for\;\;\ell\in\mbb{N}.\eqno(9.3.11)$$ Take the Cartan
subalgebra $H=\sum_{i=1}^{m+1}\mbb{F}A_{i,i}$ of $sp(2m+2,\mbb{F})$.
Define $\{\ves_1,...,\ves_{m+1}\}\subset H^\ast$ by:
$$\ves_j(A_{i,i})=\dlt_{i,j}\qquad\for\;\;i,j\in\ol{1,m+1}.\eqno(9.3.12)$$
First we have:

 \psp

{\bf Theorem 9.3.1}. {\it If $\kappa\not\in -\mbb{N}$, the
representation $\pi_\kappa$ of $sp(2m+2,\mbb{F})$ on $\msr B'$ is a
highest-weight irreducible representation with highest-weight
 $-\kappa\lmd_1$. When $-\kappa=\ell\in\mbb{N}$, $\msr B'_{(\ell)}$ is
 a finite-dimensional irreducible
$sp(2m+2,\mbb{F})$-module with highest weight $\ell\lmd_n$ and
${\msr B}'/{\msr B}'_{(\ell)}$ is an irreducible highest weight
$sp(2m+2,\mbb F)$-module with highest weight
$-(\ell+2)\lmd_1+(\ell+1)\lmd_2$, where $\lmd_i$ is the $i$th
fundamental weight of $sp(2m+2,\mbb{F})$.}

{\it Proof}.  Observe that
$$\msr B'_k=\sum_{s=0}^kx_0^s{\cal
B}_{k-s}\qquad\for\;\;k\in\mbb{N}.\eqno(9.3.13)$$ Let $ M$ be a
nonzero $sp(2m+2,\mbb{F})$-submodule of $\msr B'$. Take any $0\neq
f\in M$. Repeatedly applying  (9.3.6) to $f$, we obtain $1\in
 M$. Note
$$(B^s_{m+1,m+1})^k(1)=2^k[\prod_{r=0}^{k-1}(r+\kappa)]x_0^k\in
 M\eqno(9.3.14)$$ by (9.3.8) and
$$(A_{1,m+1})^r(x_0^k)=[\prod_{i=0}^{r-1}(k-i)x_0^{k-r}x_1^r\qquad\for\;\;r\in\ol{1,k}\eqno(9.3.15)$$
by the first equation in (9.3.5). Suppose $c\not\in-\mbb{N}$. Then
(9.3.14) yields
$$x_0^k\in  M\qquad\for\;\;k\in\mbb{N}.\eqno(9.3.16)$$
Moreover, (9.3.15) with $k=r+p$ gives
$$x_0^rx_1^p\in V\qquad\for\;\;r,p\in\mbb{N}.\eqno(9.3.17)$$
Furthermore,
$$U({\cal K})(x_0^rx_1^p)=x_0^r{\msr B}_p\subset  M\eqno(9.3.18)$$
by Lemma 9.1.1. Thus
$$\msr B'=\sum_{r,p=0}^\infty x_0^r{\msr B}_p\subset
 M;\eqno(9.3.19)$$ that is, $ M=\msr B'$. So $\msr B'$ is
an irreducible $sp(2m+2,\mbb F)$-module and $1$ is its
highest-weight vector with weight $-c\lmd_1$ with respect to the
following simple positive roots
$$\{\ves_{n+1}-\ves_n,\ves_n-\ves_{n-1},...,\ves_2-\ves_1,2\ves_1
\}.\eqno(9.3.20)$$

Next we assume $c=-\ell$ with $\ell\in\mbb{N}$. Since
$$B^s_{m+1,m+1}|_{\msr B'_\ell}=0,\;\;
A_{m+1,i}|_{\msr B'_\ell}=x_0\ptl_{x_i},\;\;B_{m+1,i}|_{\msr
B'_\ell}=x_0\ptl_{y_i}\eqno(9.3.21)$$ by (9.3.7) and (9.3.8), $\msr
B'_{(\ell)}$ is
 a finite-dimensional $sp(2m+2,\mbb{F})$-module.
Let $ M$ be a nonzero $sp(2m+2,\mbb{F})$-submodule of $\msr
B'_{(\ell)}$. By (9.3.14),
$$x_0^k\in  M\qquad\for\;\;k\in\ol{0,\ell}.\eqno(9.3.22)$$
Moreover, (9.3.15) with $k=r+s$ gives
$$x_0^rx_1^s\in  M \qquad\for\;\;r,s\in\ol{0,\ell}\;\;\mbox{such that}\;\;r+s\leq \ell.\eqno(9.3.23)$$
 Thus
$${\msr B}'_{(\ell)}=\sum_{r=0}^\ell\sum_{s=0}^{\ell-r} x_0^r{\cal B}_s\subset
 M\eqno(9.3.24)$$ by Lemma 3.1; that is,
$ M={\msr B}'_{(\ell)}$. So ${\msr B}'_{(\ell)}$ is an irreducible
$sp(2n+2,\mbb{F})$-module and $1$ is again its highest-weight
vector.

Consider the quotient $sp(2n+2,\mbb{F})$-module ${\msr B}'/{\msr
B}'_{(\ell)}$. Let $W\supset {\msr B}'_{(\ell)}$ be an
$sp(2n+2,\mbb{F})$-submodule of ${\msr B}'$ such that $W\neq {\msr
B}'_{(\ell)}$. Take any $f\in W\setminus {\msr B}'_{(\ell)}$.
Repeatedly applying (9.3.6) to $f$ if necessary, we can assume $f\in
{\cal B}_{\ell+1}$. Since ${\cal B}_{\ell+1}$ is an irreducible
$\msr{K}$-module by Theorem 9.1.1, we have
$${\cal B}_{\ell+1}\subset W.\eqno(9.3.25)$$ In particular,
$x_1^{\ell+1}\in W$. According to (9.3.7),
$$(B^s_{1,m+1})^r(x_1^{\ell+1})=r!x_1^{\ell+1+r}\in
W\qquad \for\;\;0<r\in\mbb{Z}.\eqno(9.3.26)$$ Since ${\cal
B}_{\ell+1+r}\ni x_1^{\ell+1+r}$ is an irreducible $\msr{K}$-module
by Theorem 9.1.1, we have
$${\cal B}_{\ell+1+r}\subset W.\eqno(9.3.27)$$

Suppose that
$$x_0^r{\cal B}_p \subset
W\qquad\for\;\;r\in\ol{0,k}\;\mbox{and}\;p\in\mbb{N}\;\mbox{such
that}\;r+p\geq\ell+1.\eqno(9.3.28)$$ Fix such $r$ and $p$. Observe
$x_0^rx_1^{p-1}y_1\in x_0^r{\cal B}_p\subset W$. Using the second
equation in (9.3.7), we get
$$B^s_{1,m+1}(x_0^rx_1^{p-1}y_1)=(r+p-\ell)x_0^rx_1^py_1
+x_0^{r+1}x_1^{p-1}\in W.\eqno(9.3.29)$$ By the assumption (9.3.28),
$(r+p-\ell)x_0^rx_1^py_1\in x_0^r{\cal B}_{p+1}\subset W$. So
$$x_0^{r+1}x_1^{p-1}\in W\bigcap x_0^{r+1}{\cal
B}_{p-1}.\eqno(9.3.30)$$ Since $ x_0^{r+1}{\cal B}_{p-1}$ is an
irreducible $\msr{K}$-module by Theorem 9.1.1, we get
$$ x_0^{r+1}{\cal B}_{p-1}\subset W.\eqno(9.3.31)$$
By induction on $r$, we prove
$$x_0^r{\cal B}_p \subset
W\qquad\for\;\;r,p\in\mbb{N}\;\mbox{such that}\;r+p\geq
\ell+1.\eqno(9.3.32)$$ According to (9.3.13),
$$\sum_{k=\ell+1}^\infty\msr B'_k\subset W.\eqno(9.3.33)$$
Since $W\supset \msr B'_{(\ell)}$, we have $W={\msr B}'$. So $\msr
B'/\msr B'_{(\ell)}$ is an irreducible $sp(2m+2,\mbb{F})$-module.
Moreover, $x_n^{\ell+1}$ is a highest weight vector of weight
$-(\ell+2)\lmd_1+(\ell+1)\lmd_2$ with respect to (9.3.20).
$\qquad\Box$\psp

For $\vec a=(a_1,...,a_m)^t,\;\vec b=(b_1,b_2,...,b_m)^t\in\mbb
F^m$, we put
$$\vec a\cdot\vec
x=\sum_{i=1}^ma_ix_i,\qquad\vec b\cdot\vec
y=\sum_{i=1}^mb_iy_i.\eqno(9.3.34)$$ Set
$$\msr B'_{\vec a,\vec b}=\{fe^{\vec a\cdot\vec
x+\vec b\cdot\vec y}\mid f\in{\msr B'}\}.\eqno(9.3.35)$$
 Denote by $\pi_{\kappa,\vec a,\vec b}$ the
representation $\pi_\kappa$ of $sp(2m+m,\mbb F)$ on $\msr B'_{\vec
a,\vec b}$.
 Our second result in this section is:\psp

{\bf Theorem 9.3.2}. {\it The representation $\pi_{\kappa,\vec
a,\vec b}$ with $(\vec a,\vec b)\neq (\vec 0,\vec 0)$ is an
irreducible representation of $sp(2m+2,\mbb{F})$ for any
$\kappa\in\mbb{F}$.}

{\it Proof}. By symmetry, we may assume $a_1\neq 0$. Set
$$\msr B'_{\vec a,\vec b,k}=\msr B'_ke^{\vec a\cdot\vec
x+\vec b\cdot\vec y}\qquad\for\;k\in\mbb{N}.\eqno(9.3.36)$$ Let
${M}$ be a nonzero $sp(2m+2,\mbb{F})$-submodule of $\msr B'_{\vec
a,\vec b}$. Take any $0\neq fe^{\vec a\cdot\vec x}\in M$ with $f\in
\msr B'$. By the assumption $a_0=0$ and (9.3.6),
$$C^s_{m+1,m+1}(fe^{\vec a\cdot\vec x})=-2\ptl_{x_0}(f)e^{\vec a\cdot\vec
x} \in M,\eqno(9.3.37)$$
$$(A_{i,m+1}-b_i)(fe^{\vec a\cdot\vec
x})=[\ptl_{y_i}(f)+x_i\ptl_{x_0}(f)]e^{\vec a\cdot\vec x} \in
M,\eqno(9.3.38)$$
$$(C^s_{m+1,i}+a_i)(fe^{\vec a\cdot\vec
x})=[-\ptl_{x_i}(f)+y_i\ptl_{x_0}(f)]e^{\vec a\cdot\vec x} \in
M\eqno(9.3.39)$$ for $i\in\ol{1,m}$. Repeatedly applying
(9.3.37)-(9.3.39), we obtain $e^{\vec a\cdot\vec x}\in M$.
Equivalently, $\msr B'_{\vec a,\vec b,0}\subset M$.

Suppose $\msr B'_{\vec a,\vec b,\ell}\subset M$ for some
$\ell\in\mbb{N}$.  For any $ge^{\vec a\cdot\vec x}\in \msr B'_{\vec
a,\vec b,\ell}$,
$$A_{i,1}(ge^{\vec a\cdot\vec x})=[a_1x_i-b_iy_1+x_i\ptl_{x_1}-y_1\ptl_{x_i}](g)e^{\vec a\cdot\vec x}
\in M\eqno(9.3.40)$$  and
$$C^s_{i,1}(ge^{\vec a\cdot\vec x})=[a_1y_i+a_iy_1+y_i\ptl_{x_1}+y_1\ptl_{x_i}](g)e^{\vec a\cdot\vec x}
\in M\eqno(9.3.41)$$ by (9.3.5), where $i\in\ol{1,m}$. Since
$$(x_i\ptl_{x_1}-y_1\ptl_{x_i})(g)e^{\vec a\cdot\vec x},\;(y_i\ptl_{x_1}+y_1\ptl_{x_i})(g)e^{\vec a\cdot\vec
x}\in \msr B'_{\vec a,\vec b,\ell}\subset M,\eqno(9.3.42)$$ we have
$$(a_1x_i-b_iy_1)ge^{\vec a\cdot\vec x},\;(a_1y_i+a_iy_1)ge^{\vec a\cdot\vec x}
\in M\eqno(9.3.43)$$ for $i\in\ol{1,m}$. The above second equation
with $i=1$ gives
$$2a_1y_1ge^{\vec a\cdot\vec x}
\in M\Rightarrow y_1ge^{\vec a\cdot\vec x} \in M.\eqno(9.3.44)$$
Thus (9.3.43) yields
$$x_ige^{\vec a\cdot\vec x},y_ige^{\vec a\cdot\vec x}
\in M\qquad\for\;i\in\ol{1,m}.\eqno(9.3.45)$$

According to the first equation in (9.3.7),
\begin{eqnarray*}& &A_{m+1,1}(ge^{\vec a\cdot\vec
x})\\
&=&[a_1x_0-\sum_{i=1}^m(a_ix_i+b_iy_i)y_1+x_0\ptl_{x_1}-y_1(D+\kappa)](g)e^{\vec
a\cdot\vec x}\in M.\hspace{2.9cm}(9.3.46)\end{eqnarray*} Replacing
$ge^{\vec a\cdot\vec x}\in \msr B'_{\vec a,\vec b,\ell}$ by
$ge^{\vec a\cdot\vec x}\in \sum_{i=1}^m(x_i\msr B'_{\vec a,\vec
b,\ell}+y_i\msr B'_{\vec a,\vec b,\ell})$ in (9.3.40)-(9.3.45), we
obtain
$$x_iy_1ge^{\vec a\cdot\vec x},y_iy_1ge^{\vec a\cdot\vec x}
\in M\qquad\for\;i\in\ol{1,m}.\eqno(9.3.47)$$ Since $D(g)=\ell g$
and $x_0\ptl_{x_1}(g)e^{\vec a\cdot\vec x}\in \msr B'_{\vec a,\vec
b,\ell}\subset M$, we have
$$[-\sum_{i=1}^m(a_ix_i+b_iy_i)y_1+x_0\ptl_{x_1}-y_1(D+\kappa)](g)e^{\vec
a\cdot\vec x}\in M.\eqno(9.3.48)$$ Hence (9.3.46) yields
$x_0ge^{\vec a\cdot\vec x}\in M$. Therefore, $\msr B'_{\vec a,\vec
b,\ell+1}\subset M$. By induction, $\msr B'_{\vec a,\vec
b,\ell}\subset M$ for any $\ell\in\mbb{N}$. So $\msr B'_{\vec a,\vec
b}= M$. Hence $\msr B'_{\vec a,\vec b}$ is an irreducible
$sp(2m+2,\mbb{F})$-module. $\qquad\Box$\psp

Fix $m_1,m_2\in\ol{1,m}$ with $m_1\leq m_2$. Set
 $$\td D'=x_0\ptl_{x_0}+\sum_{r=m_1+1}^mx_r\ptl_{x_r}-\sum_{i=1}^{m_1}x_i\ptl_{x_i}+
\sum_{i=1}^{m_2}y_i\ptl_{y_i}-\sum_{r=m_2+1}^my_r\ptl_{y_r}\eqno(9.3.49)$$
and
$$\td \kappa=\kappa+m_2-m_1-m.\eqno(9.3.50)$$
Swapping operators $\ptl_{x_r}\mapsto -x_r,\;
 x_r\mapsto
\ptl_{x_r}$  for $r\in\ol{1,m_1}$ and $\ptl_{y_p}\mapsto -y_p,\;
 y_p\mapsto\ptl_{y_p}$  for $s\in\ol{n_2+1,n}$ in the
oscillator representation (9.3.5)-(9.3.8),
 we get the following representation $\pi_\kappa^{m_1,m_2}$ of the
Lie algebra $sp(2m+2,\mbb{F})$  determined by
$$\pi_\kappa^{m_1,m_2}(A_{i,j})=E_{i,j}^x-E_{j,i}^y\eqno(9.3.51)$$ with
$$E_{i,j}^x=\left\{\begin{array}{ll}-x_j\ptl_{x_i}-\delta_{i,j}&\mbox{if}\;
i,j\in\ol{1,m_1};\\ \ptl_{x_i}\ptl_{x_j}&\mbox{if}\;i\in\ol{1,m_1},\;j\in\ol{m_1+1,m};\\
-x_ix_j &\mbox{if}\;i\in\ol{m_1+1,m},\;j\in\ol{1,m_1};\\
x_i\partial_{x_j}&\mbox{if}\;i,j\in\ol{m_1+1,m}
\end{array}\right.\eqno(9.3.52)$$
and
$$E_{i,j}^y=\left\{\begin{array}{ll}y_i\ptl_{y_j}&\mbox{if}\;
i,j\in\ol{1,m_2};\\ -y_iy_j&\mbox{if}\;i\in\ol{1,m_2},\;j\in\ol{m_2+1,m};\\
\ptl_{y_i}\ptl_{y_j} &\mbox{if}\;i\in\ol{m_2+1,m},\;j\in\ol{1,m_2};\\
-y_j\partial_{y_i}-\delta_{i,j}&\mbox{if}\;i,j\in\ol{m_2+1,m};
\end{array}\right.\eqno(9.3.53)$$
$$\pi_\kappa^{m_1,m_2}(E_{i,m+1+j})=\left\{\begin{array}{ll}
\ptl_{x_i}\ptl_{y_j}&\mbox{if}\;i\in\ol{1,m_1},\;j\in\ol{1,m_2};\\
-y_j\ptl_{x_i}&\mbox{if}\;i\in\ol{1,m_1},\;j\in\ol{m_2+1,m};\\
x_i\ptl_{y_j}&\mbox{if}\;i\in\ol{m_1+1,m},\;j\in\ol{1,m_2};\\
-x_iy_j&\mbox{if}\;i\in\ol{m_1+1,m},\;j\in\ol{m_2+1,m};\end{array}\right.\eqno(9.3.54)$$
$$\pi_\kappa^{m_1,m_2}(E_{m+1+i,j})=\left\{\begin{array}{ll}
-x_jy_i&\mbox{if}\;j\in\ol{1,m_1},\;i\in\ol{1,m_2};\\
-x_j\ptl_{y_i}&\mbox{if}\;j\in\ol{1,m_1},\;i\in\ol{m_2+1,m};\\
y_i\ptl_{x_j}&\mbox{if}\;j\in\ol{m_1+1,m},\;i\in\ol{1,m_2};\\
\ptl_{x_j}\ptl_{y_i}&\mbox{if}\;j\in\ol{m_1+1,m},\;i\in\ol{m_2+1,m};\end{array}\right.\eqno(9.3.55)$$
 $$\pi_\kappa^{m_1,m_2}(C^s_{m+1,m+1})=-2\ptl_{x_0},\;\;\pi_\kappa^{m_1,m_2}(B^s_{m+1,m+1})=2x_0(\td
D+\td\kappa);\eqno(9.3.56)$$
$$\pi_\kappa^{m_1,m_2}(A_{i,m+1})=\left\{\begin{array}{ll}\ptl_{x_0}\ptl_{x_i}+\ptl_{y_i}&\mbox{if}
\;\;i\in\ol{1,m_1};
\\ x_i\ptl_{x_0}+\ptl_{y_i}&\mbox{if}\;\;i\in\ol{m_1+1,m_2};\\ x_i\ptl_{x_0}-y_i&
\mbox{if}\;\;i\in\ol{m_2+1,m},\end{array}\right.\eqno(9.3.57)$$
$$\pi_\kappa^{m_1,m_2}(C^s_{m+1,i})=\left\{\begin{array}{ll}y_i\ptl_{x_0}+x_i&\mbox{if}\;\;i\in\ol{1,m_1};
\\ y_i\ptl_{x_0}-\ptl_{x_i}&\mbox{if}\;\;i\in\ol{m_1+1,m_2},\\ \ptl_{x_0}\ptl_{y_i}-\ptl_{x_i}&
\mbox{if}\;\;i\in\ol{m_2+1,m};\end{array}\right.\eqno(9.3.58)$$
$$\pi_\kappa^{m_1,m_2}(A_{m+1,i})=\left\{\begin{array}{ll}-x_0x_i-y_i(\td D'+\td\kappa)
&\mbox{if}\;\;i\in\ol{1,m_1};
\\ x_0\ptl_{x_i}-y_i(\td D'+\td\kappa)&\mbox{if}\;\;i\in\ol{m_1+1,m_2};\\ x_0\ptl_{x_i}-(\td D'+\td
c-1)\ptl_{y_i} &
\mbox{if}\;\;i\in\ol{m_2+1,m};\end{array}\right.\eqno(9.3.59)$$
$$\pi_\kappa^{m_1,m_2}(B^s_{i, m+1})=\left\{\begin{array}{ll}x_0\ptl_{y_i}+(\td D'+\td
c-1)\ptl_{x_i} &\mbox{if}\;\;i\in\ol{1,m_1};
\\ x_0\ptl_{y_i}+x_i(\td D'+\td\kappa)&\mbox{if}\;\;i\in\ol{m_1+1,m_2};\\ -x_0y_i+x_i(\td D'+\td\kappa)
& \mbox{if}\;\;i\in\ol{m_2+1,m};\end{array}\right.\eqno(9.3.60)$$
$$\pi_\kappa^{m_1,m_2}(A_{m+1,m+1})=\td D'+x_0\ptl_{x_0}+\td\kappa
\eqno(9.3.61)$$ for $i,j\in\ol{1,m}$.

For $k\in\mbb Z$, we write
$${\msr B}_{\la k\ra}=\mbox{Span}\{x^\al
y^\be\mid\al,\be\in\mbb{N}^m;\sum_{r=m_1+1}^m\al_r-\sum_{i=1}^{m_1}\al_i+
\sum_{i=1}^{m_2}\be_i-\sum_{r=m_2+1}^n\be_r=k\}\eqno(9.3.62)$$ and
$${\msr B}'_{\la k\ra}=\sum_{i=0}^\infty {\msr B}_{\la
k-i\ra}x_0^i.\eqno(9.3.63)$$ Then
$${\msr B}'_{\la k\ra}=\{u\in\msr B\mid \td
D'(u)=ku\}.\eqno(9.3.64)$$ Recall the Lie subalgebra $\msr K$ in
(9.3.9). Suppose $m>1$. Theorem 9.2.1 says that if $m_1<m_2$ or
$k\neq 0$, the subspace ${\msr B}_{\la k\ra}$ is an irreducible
$\msr K$-module; when $m_1=m_2$, the subspace ${\msr B}_{\la 0\ra}$
is a direct sum of two irreducible $\msr K$-submodules.

Assume $m=m_1=m_2=1$, $\msr{B}=\mbb{F}[x_1,y_1]$ and
$$\pi_\kappa^{m_1,m_2}(\msr{K})=\mbb{F}(x_1\ptl_{x_1}+y_1\ptl_{y_1}+1)+\mbb{F}x_1y_1+\mbb{F}\ptl_{x_1}\ptl_{y_1}.
\eqno(9.3.65)$$ So all ${\msr B}_{\la k\ra}$ with $k\in\mbb{Z}$ are
irreducible ${\msr K}$-submodules.  The following is the third
result in this section.\psp

{\bf Theorem 9.3.3}. {\it The representation $\pi_\kappa^{m_1,m_2}$
of $sp(2m+2,\mbb{F})$ on $\msr B'$ is irreducible if
$\kappa\not\in\mbb{Z}$.}

{\it Proof}. Let $M$ be any nonzero $sp(2m+2,\mbb{F})$-submodule of
$\msr B'$. Repeatedly applying the first equation in (9.3.56) to
$M$, we get
$$M\bigcap {\msr B}\neq\{0\}.\eqno(9.3.66)$$
According to (9.3.61),
$${\msr B}_{\la k\ra}=\{f\in\msr B\mid
A_{m+1,m+1}(f)=(k+\td\kappa)f\}.\eqno(9.3.67)$$ Hence
$$M=\bigoplus_{k\in\mbb{Z}}M\bigcap {\msr B}_{\la
k\ra}.\eqno(9.3.68)$$ If $M\bigcap\msr B_{\la 0\ra}\neq \{0\}$, then
(9.3.58) gives
$$C^s_{m+1,1}(M\bigcap\msr B_{\la 0\ra})=
x_1(M\bigcap\msr B_{\la 0\ra})\subset M\bigcap\msr B_{\la
-1\ra}.\eqno(9.3.69)$$ Thus we always have $M\bigcap\msr B_{\la
k\ra}\neq\{0\}$ for some $0\neq k\in\mbb{Z}$. Note that $\msr B_{\la
k\ra}$ is an irreducible $\msr K$-module. So
$$\msr B_{\la k\ra}\subset M.\eqno(9.3.70)$$

Next (9.3.57) yields
$$\msr B_{\la k-r\ra}=(\ptl_{y_1})^r(\msr  B_{\la
k\ra})=A_{1,m+1}^r(\msr B_{\la k\ra})\subset
M\qquad\for\;\;r\in\mbb{N}.\eqno(9.3.71)$$ On the other hand,  if
 $\msr B_{\la \ell\ra}\subset \msr
M$, then the assumption $\kappa\not\in\mbb Z$ and the second
equation in (9.3.56) give
$$x_0^r\msr B_{\la \ell\ra}=(C^s_{m+1,m+1})^r(\msr B_{\la
\ell\ra})\subset M \qquad\for\;\;r\in\mbb{N}.\eqno(9.3.72)$$ Suppose
that for some $p\in\mbb{Z}$,
$$x_0^r\msr B_{\la p\ra},x_0^r{\msr B}_{\la p-1\ra}
\subset M \qquad\for\;\;r\in\mbb{N}.\eqno(9.3.73)$$ For any
$\ell\in\mbb{N}$,
$$x_0^\ell{\msr B}_{\la p+1\ra}=(\td D'+\td
c-1)\ptl_{x_1}(x_0^\ell{\msr B}_{\la
p\ra})=[B^s_{1,m+1}-x_0\ptl_{y_1}](x_0^\ell{\msr B}_{\la
p\ra})\eqno(9.3.74)$$ by (9.3.60). Note
$$x_0\ptl_{y_1}(x_0^\ell{\msr
B}_{\la p\ra})=x_0^{\ell+1}{\msr B}_{\la p-1\ra}\subset
M.\eqno(9.3.75)$$ Thus (9.3.74) leads to
$$x_0^\ell{\msr B}_{\la p+1\ra}\subset M.\eqno(9.3.76)$$

By (9.3.71)-(9.3.76) and induction on $p$, we prove
$$x_0^r{\msr B}_{\la k\ra}\subset
M\qquad\for\;\;x_0\in\mbb{N},\;k\in\mbb{Z}.\eqno(9.3.77)$$ So $
M=\msr B'$. Therefore, $\msr B'$ is an irreducible
$sp(2m+2,\mbb{F})$-module. $\qquad\Box$\psp

{\bf Remark 9.3.4}. The above irreducible representation depends on
the three parameters $\kappa\in \mbb{F}$ and $m_1,m_2\in\ol{1,n}$.
It is not highest-weight type because of the mixture of
multiplication operators and differential operators in (9.3.54),
(9.3.55) and (9.3.57)-(9.3.60). Since ${\msr B}$ is not completely
reducible as a module of the Lie subalgebra $\sum_{i,j=1}^m\mbb
FA_{i,j}$ by Sections 6.4-6.6 when $m\geq 2$ and $m_1<m$, $\msr B'$
is not a unitary $sp(2m+2,\mbb{F})$-module. The constraints on the
degrees of the monomials with a fixed weight via the operators
$\pi^{m_1,m_2}_\kappa(A_{i,i})$ in (3.51)-(3.53) with $i\in\ol{1,m}$
and $\pi^{m_1,m_2}_\kappa(A_{m+1,m+1})$ in (3.61) show that the
weight subspaces are finite-dimensional. Thus $\msr B'$ is a weight
$sp(2m+2,\mbb{F})$-module with finite-dimensional weight
subspaces.\psp

Next we assume that
$$b_{j_0}\neq 0\;\;\mbox{for
some}\;\;j_0\in\ol{1,m_1}\eqno(9.3.78)$$ and
$$a_{i_0}\neq 0\;\;\mbox{for
some}\;\;i_0\in\ol{m_2+1,m}\;\;\mbox{if}\;\;m_2<m.\eqno(9.3.79)$$
Under the assumption, we have the following fourth result in this
 section:\psp

{\bf Theorem 9.3.5}. {\it The representation
$\pi_{\kappa}^{m_1,m_2}$ of $sp(2m+2,\mbb{F})$ on $\msr B'_{\vec
a,\vec b}$ is irreducible for any $\kappa\in\mbb{F}$.}

{\it Proof}.  Let ${M}$ be a nonzero $sp(2m+2,\mbb{F})$-submodule of
$\msr B'_{\vec a,\vec b}$. Take any $0\neq fe^{\vec a\cdot\vec x}\in
M$ with $f\in \msr B'$.
  Repeatedly applying the first equation in
(9.3.56) to $fe^{\vec a\cdot\vec x}$ if necessary, we may assume
$f\in\msr B=\mbb{F}[x_1,...,x_m,y_1,...,y_m]$. Then (9.3.58) yields
$$(C^s_{m+1,i}+a_i)(fe^{\vec a\cdot\vec x})=-\ptl_{x_i}(f)e^{\vec a\cdot\vec
x}\in M\;\;\for\;\;i\in\ol{m_1+1,m}.\eqno(9.3.80)$$ Moreover,
(9.3.57) yields
$$(A_{i,m+1}-b_j)(fe^{\vec a\cdot\vec x})=\ptl_{y_j}(f)e^{\vec a\cdot\vec
x}\in M\;\;\for\;\;j\in\ol{1,m_2}.\eqno(9.3.81)$$ Repeatedly
applying (9.3.80) and (9.3.81) if necessary, we can assume
$$f\in\mbb{F}[x_1,...,x_{m_1},y_{m_2+1},...,y_m].\eqno(9.3.82)$$

According to (9.3.54),
$$(B^s_{i,j_0}-a_{j_0}b_i-a_ib_{j_0})(fe^{\vec a\cdot\vec x})=(b_{j_0}\ptl_{x_i}+b_i\ptl_{x_{j_0}})(f)e^{\vec a\cdot\vec
x}\in M\eqno(9.3.83)$$ for $i\in\ol{1,m_1}$. Taking $i=j_0$ in
(9.3.83), we get $\ptl_{x_{j_0}}(f)e^{\vec a\cdot\vec x}\in M$.
Substituting it to (9.3.83) for general $i$, we obtain
$$\ptl_{x_i}(f)e^{\vec a\cdot\vec
x}\in M\qquad\for\;\;i\in\ol{1,m_1}.\eqno(9.3.84)$$
 Moreover, (9.3.55) yields
$$(C^s_{j,i_0}-a_jb_{i_0}-a_{i_0}b_j)(fe^{\vec a\cdot\vec x})
=(a_{i_0}\ptl_{y_j}+a_j\ptl_{y_{i_0}})(f)e^{\vec a\cdot\vec x}\in
M\eqno(9.3.85)$$ for $j\in\ol{m_2+1,m}.$ Letting $j=i_0$ in
(9.3.85), we find $\ptl_{y_{i_0}}(f)e^{\vec a\cdot\vec x}\in M$.
Substituting it to (9.3.85) for general $j$, we get
$$\ptl_{y_j}(f)e^{\vec a\cdot\vec
x}\in M\qquad\for\;\;j\in\ol{m_2+1,m}.\eqno(9.3.86)$$
 Repeatedly applying
(9.3.84) and (9.3.86) if necessary, we obtain $e^{\vec a\cdot\vec
x}\in M$. Equivalently, $\msr B'_{\vec a,\vec a, 0}\subset M$ (cf.
(9.3.36)).

Suppose that for some $\ell\in\mbb{N}$, $\msr B'_{\vec a,\vec a,
k}\subset M$ whenever $\ell\geq k\in\mbb{N}$. For any $ge^{\vec
a\cdot\vec x}\in \msr B'_{\vec a,\vec b,\ell}$, (9.3.58) implies
$$(C^s_{m+1,i}-y_i\ptl_{x_0})(ge^{\vec a\cdot\vec
x})=x_ige^{\vec a\cdot\vec x}\in
M\qquad\for\;\;i\in\ol{1,m_1}\eqno(9.3.87)$$ and (9.3.57) leads to
$$(-A_{i,m+1}+x_i\ptl_{x_0})(ge^{\vec a\cdot\vec x})=y_jge^{\vec a\cdot\vec
x}\in M\qquad\for\;\;j\in \ol{m_2+1,m}.\eqno(9.3.88)$$ Moreover,
(9.3.54) gives
$$B^s_{i,m+1+j_0}(ge^{\vec a\cdot\vec
x})=[b_{j_0}x_i+x_i\ptl_{y_{j_0}}+(\ptl_{x_{j_0}}+a_{j_0})(\ptl_{y_i}+b_i)](g)e^{\vec
a\cdot\vec x}\in M\eqno(9.3.89)$$ if $i\in\ol{m_1+1,m_2}$, and
$$B^s_{i,j_0}(ge^{\vec a\cdot\vec
x})=[b_{j_0}x_i-a_{j_0}y_i+x_i\ptl_{y_{j_0}}-y_i\ptl_{x_{j_0}}](g)e^{\vec
a\cdot\vec x}\in M\eqno(9.3.90)$$ if $i\in\ol{m_2+1,m}$. Note that
the inductional assumption imply
$$[x_i\ptl_{y_{j_0}}+(\ptl_{x_{j_0}}+a_{j_0})(\ptl_{y_i}+b_i)](g)e^{\vec
a\cdot\vec x}\in M\eqno(9.3.91)$$  if $i\in\ol{m_1+1,m_2}$, and
$$[-a_{j_0}y_i+x_i\ptl_{y_{j_0}}-y_i\ptl_{x_{j_0}}](g)e^{\vec
a\cdot\vec x}\in M\eqno(9.3.92)$$ by (9.3.90) if $i\in\ol{m_2+1,m}$.
Thus
$$x_ige^{\vec a\cdot\vec
x}\in M\qquad\for\;\;i\in\ol{m_1+1,m}.\eqno(9.3.93)$$

On the other hand, (9.3.55) yields
$$C^s_{j,i_0}(ge^{\vec a\cdot\vec
x})=(a_{i_0}y_j-b_{i_0}x_j+y_j\ptl_{x_{i_0}}-x_j\ptl_{y_{i_0}})(g)e^{\vec
a\cdot\vec x}\in M\eqno(9.3.94)$$ if $j\in\ol{1,m_1}$, and
$$C^s_{j,i_0}(ge^{\vec a\cdot\vec
x})=[a_{i_0}y_j+y_j\ptl_{x_{i_0}}+(\ptl_{x_j}+a_j)(\ptl_{y_{i_0}}+b_{i_0})](g)e^{\vec
a\cdot\vec x}\in M\eqno(9.3.95)$$ if $j\in\ol{m_1+1,m_2}$. Observe
that the inductional assumption imply
$$(-b_{i_0}x_j+y_j\ptl_{x_{i_0}}-x_j\ptl_{y_{i_0}})(g)e^{\vec
a\cdot\vec x}\in M\eqno(9.3.96)$$ by (9.3.87) if $j\in\ol{1,m_1}$,
and
$$[y_j\ptl_{x_{i_0}}+(\ptl_{x_j}+a_j)(\ptl_{y_{i_0}}+b_{i_0})](g)e^{\vec
a\cdot\vec x}\in M\eqno(9.3.97)$$ if $j\in\ol{m_1+1,m_2}$. Hence
$$y_jge^{\vec a\cdot\vec x}\in M\qquad\for\;\;j\in
\ol{1,m_2}.\eqno(9.3.98)$$  Moreover, (9.3.60) yields
\begin{eqnarray*}\qquad& &B^s_{j_0,m+1}(ge^{\vec a\cdot\vec
x})\\&=&[b_{j_0}x_0+x_0\ptl_{y_{j_0}}+(\td
D'-\sum_{i=1}^{m_1}a_ix_i+\sum_{j=m_1+1}^ma_jx_j+\sum_{r=1}^{m_2}b_ry_r\\
& &-\sum_{s=m_2+1}^mb_sy_s +\td
\kappa+1)(a_{j_0}+\ptl_{x_{j_0}})](g)e^{\vec a\cdot\vec x}\in
M.\hspace{4.3cm}(9.3.99)\end{eqnarray*}  Note that
$$x_0\ptl_{y_{j_0}}(g)e^{\vec
a\cdot\vec x}\in\msr B'_{\vec a,\vec b,\ell}\subset M;\;\; (\td D
+\td\kappa+1)(\ptl_{x_{j_0}}(g))e^{\vec a\cdot\vec x}\in{\msr
B}'_{\vec a,\ell-1}\subset M.\eqno(9.3.100)$$ Now (9.3.87),
(9.3.88), (9.3.93) and (9.3.98)-(9.3.100) imply $x_0ge^{\vec
a\cdot\vec x}\in M$. Therefore, $\msr B'_{\vec a,\vec
b,\ell+1}\subset M$. By induction, $\msr B'_{\vec a,\vec
b,\ell}\subset M$ for any $\ell\in\mbb{N}$. So $\msr B'_{\vec a,\vec
b}=M$. Hence $\msr B'_{\vec a,\vec b}$ is an irreducible
$sp(2m+2,\mbb{F})$-module. $\qquad\Box$\psp

\chapter{Representations of $G_2$ and $F_4$}

In this chapter, we determine the structure of the canonical bosonic
and fermionic oscillator representations of the simple Lie algebra
of type $G_2$ over its 7-dimensional module. Moreover, we use
partial differential equations to  find the explicit irreducible
decomposition of the space of polynomial functions on 26-dimensional
basic irreducible module of the simple Lie algebra of type $F_4$
(cf. [X18]).

\section{Representations of $G_2$}

In this section, we study the structure of the canonical bosonic and
fermionic oscillator representations of the simple Lie algebra of
type $G_2$ over its 7-dimensional module.

Denote by $\mbb{Z}_3=\mbb{Z}/3\mbb{Z}$ and identify $i+3\mbb{Z}$
with $i$ for $i\in\{1,2,3\}$. Recall the simple Lie algebra
$${\msr
G}^{G_2}=
 \mbb{F}h_1+\mbb{F}h_2+\sum_{i,j\in\mbb Z_3,\;i\neq j}\mbb
 F(E_{i,j}-E_{j',i'})+\sum_{r\in\mbb Z_3}(\mbb FC_r+\mbb FC'_r)\eqno(10.1.1)$$
constructed in (4.4.1)-(4.4.10).  Let ${\cal Q}$ be the space of
rational function
 in $x_0,x_1,x_2,x_3,x_{1'},\\ x_{2'},x_{3'}$.  The corresponding oscillator representation is
 given by
$$h_1|_{\cal
Q}=x_1\ptl_{x_1}-x_2\ptl_{x_2}-x_{1'}\ptl_{x_{1'}}+x_{2'}\ptl_{x_{2'}},\eqno(10.1.2)$$
$$h_2|_{\cal
Q}=-2x_1\ptl_{x_1}+x_2\ptl_{x_2}+x_3\ptl_{x_3}+2x_{1'}\ptl_{x_{1'}}-x_{2'}\ptl_{x_{2'}}-x_{3'}\ptl_{x_{3'}},\eqno(10.1.3)$$
$$(E_{i,j}-E_{j',i'})|_{\cal
Q}=x_i\ptl_{x_j}-x_{j'}\ptl_{x_{i'}}\qquad\for\;\;i,j\in\mbb
Z_3,\;i\neq j,\eqno(10.1.4)$$
$$C_r|_{\cal
Q}=x_r\ptl_{x_{(r+2)'}}-x_{r+2}\ptl_{x_{r'}}+\sqrt{2}(x_0\ptl_{x_{r+1}}-x_{(r+1)'}\ptl_{x_0}),\eqno(10.1.5)$$
$$C_r'|_{\cal
Q}=x_{r'}\ptl_{x_{r+2}}-x_{(r+2)'}\ptl_{x_r}+\sqrt{2}(x_0\ptl_{x_{(r+1)'}}-x_{r+1}\ptl_{x_0})\eqno(10.1.6)$$
for $r\in\mbb{Z}_3$. Moreover, the vectors
$$C_1,\;C_2',\;C_3,\;E_{1,2}-E_{2',1'},\;E_{3,1}-E_{1',3'},\;E_{3,2}-E_{2',3'}\eqno(10.1.7)$$
are positive root vectors.

Define
$$\eta=x_0^2+2x_1x_{1'}+2x_2x_{2'}+2x_3x_{3'}.\eqno(10.1.8)$$
It can be verified that
$$\xi(\eta)=0\qquad\for\;\;\xi\in{\msr G}^{G_2}.\eqno(10.1.9)$$
\psp

{\bf Lemma 10.1.1}. {\it Any singular function in $\cal Q$ with
respect to ${\msr G}^{G_2}$ must be a rational function in $x_3$ and
$\eta$.}

{\it Proof}.
 Let $f$ be a ${\msr G}^{G_2}$-singular function in
$Q$. By Lemma 6.2.3 and
$$(E_{3,1}-E_{1',3'})(f)=0,\qquad
(E_{1,2}-E_{2',1'})(f)=0,\eqno(10.1.10)$$ we have
$$f=\vf(x_0,x_3,x_{2'},\eta)\eqno(10.1.11)$$
as a rational function in $x_0,x_3,x_{2'},\eta$.

Note
$$C_1(f)=[x_1\ptl_{x_{3'}}-x_3\ptl_{x_{1'}}+\sqrt{2}(x_0\ptl_{x_2}-x_{2'}\ptl_{x_0})]
(f) =-\sqrt{2}x_{2'}\ptl_{x_0}(\vf)=0.\eqno(10.1.12)$$ So $\vf$ is
independent of $x_0$. Moreover,
$$C_3(f)=[x_3\ptl_{x_{2'}}-x_2\ptl_{x_{3'}}+\sqrt{2}(x_0\ptl_{x_1}-x_{1'}\ptl_{x_0})]
=x_3\ptl_{x_{2'}}(\vf)=0. \eqno(10.1.13)$$ Thus $\vf$ is independent
of $x_{2'}$. Equivalently, $f$ is a function in $x_3$ and
$\eta.\hfill\Box$\psp

Recall that
$$\Dlt'=\ptl_{x_0}^2+2\ptl_{x_1}\ptl_{x_{1'}}+2\ptl_{x_2}\ptl_{x_{2'}}+2\ptl_{x_3}\ptl_{x_{3'}}
\eqno(10.1.14)$$ is an $o(7,\mbb{F})$-invariant differential
operator (cf. (8.1.11)) and ${\msr G}^{G_2}\subset o(7,\mbb{F})$. So
$$\xi\msr D'=\msr D'\xi\qquad\for\;\;\xi\in{\msr G}^{G_2}.\eqno(10.1.15)$$
Let ${\msr B}'=\mbb{F}[x_0,x_1,x_2,x_3,x_{1'}, x_{2'},x_{3'}]\subset
{\cal Q}$ be the polynomial algebra in $x_0,x_1,x_2,x_3,x_{1'},
\\ x_{2'},x_{3'}$. The above lemma implies that any ${\msr
G}^{G_2}$-singular vector in ${\msr B}'$ is a polynomial in $x_3$
and $\eta$. Denote by ${\msr B}'_k$ the subspace of homogeneous
polynomials in ${\msr B}'$ with degree $k$. Set
$${\msr H}'_k=\{f\in{\msr B}'_k\mid {\msr
D}'(f)=0\}.\eqno(10.1.16)$$ By Theorem 8.1.1, we have the following
theorem.\psp

{\bf Theorem 10.1.2}. {\it For any $k\in\mbb{N}$, the subspace
${\msr H}'_k$ forms an irreducible ${\msr G}^{G_2}$-module with
highest weight vector $x_3^k$ of weight $k\lmd_2$. Moreover,
$${\msr B}'=\bigoplus_{i,k=0}^\infty \eta^i{\msr
H}'_k\eqno(10.1.17)$$ is a direct sum of  irreducible ${\msr
G}_2$-modules.}\psp

Consider the exterior algebra $\check{\msr A}'$ generated by
$\{\sta_0,\sta_1,\sta_2,\sta_3,\sta_{1'}, \sta_{2'},\sta_{3'}\}$
(cf. (6.2.15)-(6.2.18)). Define a representation of ${\msr G}^{G_2}$
on $\check{\msr A}'$ by
$$h_1|_{\check{\msr A}'}=\sta_1\ptl_{\sta_1}-\sta_2\ptl_{\sta_2}-\sta_{1'}\ptl_{\sta_{1'}}+\sta_{2'}\ptl_{\sta_{2'}},
\eqno(10.1.18)$$
$$h_2|_{\check{\msr A}'}=-2\sta_1\ptl_{\sta_1}+\sta_2\ptl_{\sta_2}+\sta_3\ptl_{\sta_3}+2\sta_{1'}
\ptl_{\sta_{1'}}-\sta_{2'}\ptl_{\sta_{2'}}-\sta_{3'}\ptl_{\sta_{3'}},\eqno(10.1.19)$$
$$(E_{i,j}-E_{j',i'})|_{\check{\msr A}'}=\sta_i\ptl_{\sta_j}-\sta_{j'}\ptl_{\sta_{i'}}\qquad\for\;\;i,j\in\mbb
Z_3,\;i\neq j,\eqno(10.1.20)$$
$$C_r|_{\check{\msr A}'}=\sta_r\ptl_{\sta_{(r+2)'}}-\sta_{r+2}\ptl_{\sta_{r'}}+\sqrt{2}(\sta_0\ptl_{\sta_{r+1}}
-\sta_{(r+1)'}\ptl_{\sta_0}),\eqno(10.1.21)$$
$$C_r'|_{\check{\msr A}'}=\sta_{r'}\ptl_{\sta_{r+2}}-\sta_{(r+2)'}\ptl_{\sta_r}+\sqrt{2}(\sta_0\ptl_{\sta_{(r+1)'}}
-\sta_{r+1}\ptl_{\sta_0})\eqno(10.1.22)$$ for $r\in\mbb{Z}_3$.

Note that
$$\Theta=\sum_{i=0}^6\mbb{F}\sta_i\eqno(10.1.23)$$
is the 7-dimensional basic irreducible ${\msr G}^{G_2}$-module.
Define
$$\check{\msr
A}'_k=\Theta^k\qquad\for\;\;k\in\ol{0,7}.\eqno(10.1.24)$$ Then
$\check{\msr A}'_k$ forms a ${\msr G}^{G_2}$-module. By the Dynkin
diagram of $G_2$, any finite-dimensional ${\msr G}^{G_2}$-module is
isomorphic to its dual module which consists of its linear
functions. In particular,
$$\check{\msr A}'_k\cong \check{\msr A}'_{7-k}\;\;\mbox{as}\; {\msr
G}_2\mbox{-modules}.\eqno(10.1.25)$$ We only need to study the
structure ${\msr G}^{G_2}$-modules $\check{\msr A}'_2$ and
$\check{\msr A}'_3$.

Recall that
$$\check{\eta}=\sta_1\sta_{1'}+\sta_2\sta_{2'}+\sta_3\sta_{3'}\eqno(10.1.26)$$
is a ${\msr G}_0$-invariant. According to (6.2.53) with $n=3$ and
$\vt_i=\sta_{i'}$ for $i\in\ol{1,3}$, the $\msr G^{G_2}$-singular
vectors  in $\check{\msr A}'_2$ are in the set
$$\{a\sta_0\sta_3+b\sta_{1'}\sta_{2'},a\sta_0\sta_{2'}+b\sta_1\sta_3,\sta_3\sta_{2'},\check{\eta}\mid (0,0)\neq (a,b)\in\mbb{F}\}.
\eqno(10.1.27)$$Moreover,
$$C_3(a\sta_0\sta_3+b\sta_{1'}\sta_{2'})=-\sqrt{2}a\sta_{1'}\sta_3+b\sta_3\sta_{1'},\eqno(10.1.28)$$
$$C_3(a\sta_0\sta_{2'}+b\sta_1\sta_3)=a\sta_0\sta_3+\sqrt{2}b\sta_0\sta_3,
\eqno(10.1.29)$$ and
$$C_3(\check{\eta})=\sqrt{2}\sta_0\sta_{1'}+2\sta_1\sta_3\eqno(10.1.30)$$
by (10.1.21).  Thus $\sta_3\sta_{2'}$ is the only ${\msr
G}^{G_2}$-singular vector in $\check{\msr A}'_2$, which has weight
$\lmd_1$.

Again by (6.2.53),  the $\msr G^{G_2}$-singular vectors  in
$\check{\msr A}'_3$ are in the set
\begin{eqnarray*}
\hspace{2cm}& &\{a\sta_0\check\eta+b\sta_1\sta_2\sta_3
+c\sta_{1'}\sta_{2'}\sta_{3'},
a\sta_0\sta_1\sta_3+b\sta_{2'}\check\eta,a\sta_0\sta_{1'}\sta_{2'}+b\sta_1\check\eta,
\\& &\sta_0\sta_3\sta_{2'},\sta_1\sta_3\sta_{2'},\sta_3\sta_{1'}\sta_{2'}\mid
a,b,c\in\mbb{F}\}. \hspace{5.4cm}(10.1.31)\end{eqnarray*} Calculate
\begin{eqnarray*}
\hspace{2cm}& &C_3(a\sta_0\check\eta+b\sta_1\sta_2\sta_3
+c\sta_{1'}\sta_{2'}\sta_{3'})\\
&=&-\sqrt{2}a\sta_{1'}\check\eta+(2a+\sqrt{2}b)\sta_0\sta_2\sta_3+c\sta_{1'}(\sta_2\sta_{2'}+\sta_3\sta_{3'}),
\hspace{2.3cm}(10.1.32)\end{eqnarray*}
$$C_3(a\sta_0\sta_1\sta_3+b\sta_{2'}\check\eta)=-a\sqrt{2}\sta_1\sta_2\sta_{1'}
+b(\sta_3\check\eta+\sqrt{2}\sta_0\sta_{1'}\sta_{2'}+2\sta_2\sta_3\sta_{2'}),\eqno(10.1.33)$$
$$C_3(a\sta_0\sta_{1'}\sta_{2'}+b\sta_3\check\eta)=-a\sta_0\sta_3\sta_{1'}-b\sqrt{2}\sta_0\sta_3\sta_{1'},\eqno(10.1.34)$$
$$C_3(\sta_0\sta_3\sta_{2'})=\sqrt{2}\sta_3\sta_{1'}\sta_{2'},\;C_3(\sta_1\sta_3\sta_{2'})=
\sqrt{2}\sta_0\sta_3\sta_{2'},\;C_3(\sta_3\sta_{1'}\sta_{2'})=0\eqno(10.1.35)$$
by (10.1.21). Thus we have three $\msr G^{G_2}$-singular vectors
$\sta_0\check\eta+\sqrt{2}(-\sta_1\sta_2\sta_3+\sta_{1'}\sta_{2'}\sta_{3'})$
of weight 0, $\sqrt{2}\sta_0\sta_{1'}\sta_{2'}-\sta_3\check\eta$ of
weight $\lmd_2$ and $\sta_3\sta_{1'}\sta_{2'}$ of weight $2\lmd_2$
by (10.1.18) and (10.1.19). By Weyl's Theorem 2.3.6 of completely
reducibility, we have the following theorem. \psp

{\bf Theorem 10.1.2}. {\it The subspace $\check{\msr A}'_2$ is an
irreducible ${\msr G}^{G_2}$-module of highest weight $\lmd_1$. The
subspace $\check{\msr A}'_3$ is a direct sum of a one-dimensional
trivial ${\msr G}^{G_2}$-module, an irreducible ${\msr
G}^{G_2}$-module of highest weight $\lmd_2$ and an irreducible
${\msr G}^{G_2}$-module of highest weight $2\lmd_2$.}

 \section{Basic Oscillator Representation of $F_4$}

In this section, we present the bosonic oscillator tepresentation of
$F_4$ over its 26-dimensional irreducible module.

Let us go back to the construction of the simple Lie algebra $\msr
G^{F_4}$ of type $F_4$ via the root lattice construction of the
simple Lie algebra $\msr G^{E_6}$ of type $E_6$ in
(4.4.63)-(4.4.80).

Recall the Dynkin diagram of $E_6$:

\begin{picture}(80,23)
\put(2,0){$E_6$:}\put(21,0){\circle{2}}\put(21,
-5){1}\put(22,0){\line(1,0){12}}\put(35,0){\circle{2}}\put(35,
-5){3}\put(36,0){\line(1,0){12}}\put(49,0){\circle{2}}\put(49,
-5){4}\put(49,1){\line(0,1){10}}\put(49,12){\circle{2}}\put(52,10){2}\put(50,0){\line(1,0){12}}
\put(63,0){\circle{2}}\put(63,-5){5}\put(64,0){\line(1,0){12}}\put(77,0){\circle{2}}\put(77,
-5){6}
\end{picture}
\vspace{0.7cm}

Let $\Pi_{E_6}=\{\al_1,...,\al_6\}$ be the set of positive simple
roots corresponding to the above diagram and let $\Phi_{E_6}$ be the
root system. The simple Lie algebra $\msr G^{E_6}$ of type $E_6$ is
$${\msr
G}^{E_6}=H\oplus\bigoplus_{\al\in\Phi_{E_6}}\mbb{C}E_{\al},\;\;H=\sum_{i=1}^6\mbb
F\al,\eqno(10.2.1)$$ with the Lie bracket given in (4.4.24) and
(4.4.25) and $F(\cdot,\cdot)$ given in (4.4.64). In terms of
(4.4.68)-(4.4.80), the simple Lie algebra  of type $F_4$ is a Lie
subalgebra of ${\msr G}^{E_6}$:
$${\msr G}^{F_4}=\sum_{\varpi\in S_{F_4}}(\mbb FF_\varpi+\mbb F
F'_\varpi)+\sum_{i=1}^4\mbb Fh_i\eqno(10.2.2)$$(cf. (4.4.78),
(4.4.79)), whose Dynkin diagram is

\begin{picture}(60,12)\put(2,0){$F_4$:}
\put(21,0){\circle{2}}\put(21,-5){1}\put(22,0){\line(1,0){12}}
\put(35,0){\circle{2}}\put(35,-5){2}\put(35,1.2){\line(1,0){13.6}}
\put(35,-0.8){\line(1,0){13.6}}\put(41,-1){$\ra$}\put(48.5,0){\circle{2}}\put(48.5,-5){3}\put(49.5,0)
{\line(1,0){12}}\put(62.5,0){\circle{2}}\put(62.5,-5){4}
\end{picture}
\vspace{0.6cm}

Recall the notions in (4.4.42)-(4.4.44). Set
$$\xi_1=E_{(1,1,2,2,1,1)}-E_{(1,1,1,2,2,1)},\qquad
\xi_2=E_{(1,1,2,2,1)}-E_{(0,1,1,2,2,1)}, \eqno(10.2.3)$$
$$\xi_3=E_{(1,1,1,2,1)}-E_{(0,1,1,2,1,1)},\qquad
\xi_4=E_{(1,1,1,1,1)}-E_{(0,1,1,1,1,1)},\eqno(10.2.4)$$
$$\xi_5=E_{(1,1,1,1)}-E_{(0,1,0,1,1,1)},\qquad
\xi_6=E_{(1,0,1,1,1)}-E_{(0,0,1,1,1,1)}, \eqno(10.2.5)$$
$$\xi_7=E_{(0,1,1,1)}-E_{(0,1,0,1,1)},\;
\xi_8=E_{(1,0,1,1)}-E_{(0,0,0,1,1,1)},\;\xi_9=E_{(0,0,1,1)}-E_{(0,0,0,1,1)},
\eqno(10.2.6)$$
$$\xi_{10}=E_{(1,0,1)}-E_{(0,0,0,0,1,1)},\;\;\xi_{11}=E_{\al_3}-E_{\al_5},\;\;
\xi_{12}=E_{\al_1}-E_{\al_6},\eqno(10.2.7)$$
$$\xi_{13}=\al_1-\al_6,\qquad \xi_{14}=\al_3-\al_5,\qquad
\xi_{15}=E_{-\al_1}-E_{-\al_6}, \eqno(10.2.8)$$
$$\xi_{16}=E_{-\al_3}-E_{-\al_5},\;\;\xi_{17}=E'_{(1,0,1)}-E'_{(0,0,0,0,1,1)},\;\;
 \xi_{18}=E'_{(0,0,1,1)}-E'_{(0,0,0,1,1)}, \eqno(10.2.9)$$
$$\xi_{19}=E'_{(1,0,1,1)}-E'_{(0,0,0,1,1,1)},\;\ \xi_{20}=E'_{(0,1,1,1)}-E'_{(0,1,0,1,1)},
\eqno(10.2.10)$$
$$\xi_{21}=E'_{(1,0,1,1,1)}-E'_{(0,0,1,1,1,1)},\;\;
\xi_{22}=E_{(1,1,1,1)}-E_{(0,1,0,1,1,1)},\eqno(10.2.11)$$
$$
\xi_{23}=E'_{(1,1,1,1,1)}-E'_{(0,1,1,1,1,1)},\;\;\xi_{24}=E'_{(1,1,1,2,1)}-E'_{(0,1,1,2,1,1)},
\eqno(10.2.12)$$
$$\xi_{25}=E'_{(1,1,2,2,1)}-E'_{(0,1,1,2,2,1)},\;\;
\xi_{26}=E'_{(1,1,2,2,1,1)}-E'_{(1,1,1,2,2,1)}. \eqno(10.2.13)$$
Recall the automorphism $\hat\sgm$ of ${\msr G}^{E_6}$ defined in
(4.4.63), (4.4.66) and (4.4.67). Then the subspace
$$V=\sum_{i=1}^{26}\mbb F\xi_i=\{v\in{\msr G}^{E_6}\mid \hat\sgm
(v)=-v\}\eqno(10.2.14)$$ forms an irreducible $\msr G^{F_4}$-module
with respect to the adjoint representation of $\msr G^{E_6}$,
$\xi_1$ is a highest-weight vector of weight $\lmd_4$ and
$$ {\msr G}^{E_6}=\msr G^{F_4}\oplus V.\eqno(10.2.15)$$
 Write
$$[u,\xi_i]=\sum_{j=1}^{26}\vf_{i,j}(u)\xi_j\qquad\for\;\;u\in\msr
G^{F_4}.\eqno(10.2.16)$$ Set
$$\msr A=\mbb{F}[x_1,...,x_{26}]\eqno(10.2.17)$$
and define the {\it basic oscillator representation of} $\msr
G^{F_4}$ on $\msr A$ by
$$u(g)=\sum_{i,j=1}^{26}\vf_{i,j}(u)x_j\ptl_{x_i}(g)\qquad\for\;\;u\in\msr
G^{F_4},\;g\in \msr A\eqno(10.2.18)$$ (cf. (2.2.17)-(2.2.20)). Then
$\msr A$ forms a $\msr G^{F_4}$-module isomorphic to the symmetric
tensor $S(V)$ over $V$. More explicitly, we have
$$F_{(1)}|_{\msr A}=x_4\ptl_{x_6}+x_5\ptl_{x_8}+x_7\ptl_{x_9}-
x_{18}\ptl_{x_{20}}-x_{19}\ptl_{x_{22}}-x_{21}\ptl_{x_{23}},\eqno(10.2.19)$$
$$F_{(0,1)}|_{\msr A}=x_3\ptl_{x_4}+x_8\ptl_{x_{10}}+x_9\ptl_{x_{11}}
-x_{16}\ptl_{x_{18}}-x_{17}\ptl_{x_{19}}-x_{23}\ptl_{x_{24}},\eqno(10.2.20)$$
\begin{eqnarray*}\qquad F_{(0,0,1)}|_{\msr A}&=&-x_2\ptl_{x_3}-x_4\ptl_{x_5}-x_6\ptl_{x_8}
+x_{10}\ptl_{x_{12}}+x_{11}(\ptl_{x_{13}}-2\ptl_{x_{14}})\\
&
&-x_{14}\ptl_{x_{16}}-x_{15}\ptl_{x_{17}}+x_{19}\ptl_{x_{21}}+x_{22}\ptl_{x_{23}}+x_{24}\ptl_{x_{25}},
\hspace{2.35cm}(10.2.21)\end{eqnarray*}
\begin{eqnarray*}\hspace{1cm}F_{(0,0,0,1)}|_{\msr A}&=&-x_1\ptl_{x_2}-x_5\ptl_{x_7}-x_8\ptl_{x_9}-
x_{10}\ptl_{x_{11}}+x_{12}(\ptl_{x_{14}}-2\ptl_{x_{13}})
\\ &
&-x_{13}\ptl_{x_{15}}+x_{16}\ptl_{x_{17}}+x_{18}\ptl_{x_{19}}+x_{20}\ptl_{x_{22}}
+x_{25}\ptl_{x_{26}},\hspace{1.9cm}(10.2.22)\end{eqnarray*}
$$F_{(1,1)}|_{\msr A}=-x_3\ptl_{x_6}+x_5\ptl_{x_{10}}+x_7\ptl_{x_{11}}-x_{16}\ptl_{x_{20}}
-x_{17}\ptl_{x_{22}}+x_{21}\ptl_{x_{24}},\eqno(10.2.23)$$
\begin{eqnarray*}\hspace{1cm}F_{(0,1,1)}|_{\msr A}&=&x_2\ptl_{x_4}+x_3\ptl_{x_5}
+x_6\ptl_{10}+x_8\ptl_{x_{12}}+x_9(\ptl_{x_{13}}-2\ptl_{x_{14}})\\
&
&-x_{14}\ptl_{x_{18}}-x_{15}\ptl_{x_{19}}-x_{17}\ptl_{x_{21}}-x_{22}\ptl_{x_{24}}-x_{23}\ptl_{x_{25}},\hspace{2.1cm}
(10.2.24)\end{eqnarray*}
\begin{eqnarray*}\hspace{0.4cm}F_{(0,0,1,1)}|_{\msr A}&=&-x_1
\ptl_{x_3}+x_4\ptl_{x_7}+x_6\ptl_{x_{9}}-x_{10}(\ptl_{x_{13}}+\ptl_{x_{14}})-x_{11}\ptl_{x_{15}}
\\
&&+x_{12}\ptl_{x_{16}}-(x_{13}+x_{14})\ptl_{x_{17}}-x_{18}\ptl_{x_{21}}
-x_{20}\ptl_{x_{23}}+
x_{24}\ptl_{x_{26}},\hspace{1.15cm}(10.2.25)\end{eqnarray*}
\begin{eqnarray*}F_{(1,1,1)}|_{\msr A}&=&-x_2\ptl_{x_6}+x_3\ptl_{x_8}
+x_4\ptl_{x_{10}}+x_5\ptl_{x_{12}}+x_7(\ptl_{x_{13}}-2\ptl_{x_{14}})
-x_{14}\ptl_{x_{20}}\\ & &-x_{15}\ptl_{x_{22}}-x_{17}\ptl_{x_{23}}
-x_{19}\ptl_{x_{24}}+x_{21}\ptl_{x_{25}}, \hspace{4.9cm}(10.2.26)
\end{eqnarray*}
\begin{eqnarray*}F_{(0,1,1,1)}|_{\msr A}&=&x_1\ptl_{x_4}+
x_3\ptl_{x_7}-x_6\ptl_{x_{11}}-x_8(\ptl_{x_{13}}+\ptl_{x_{14}})-x_9\ptl_{x_{15}}+x_{12}\ptl_{x_{18}}
\\ & &-(x_{13}+x_{14})\ptl_{x_{19}}+x_{16}\ptl_{x_{21}}
-x_{20}\ptl_{x_{24}}-x_{23}\ptl_{x_{26}}, \hspace{3.25cm}(10.2.27)
\end{eqnarray*}
$$F_{(0,1,2)}|_{\msr A}=-x_2\ptl_{x_5}+x_6\ptl_{x_{12}}+x_9\ptl_{x_{16}}-x_{11}\ptl_{x_{18}}-x_{15}\ptl_{x_{21}}
+x_{22}\ptl_{x_{25}}, \eqno(10.2.28)$$
$$F_{(1,1,2)}|_{\msr A}=x_2\ptl_{x_8}
+x_4\ptl_{x_{12}}+x_7\ptl_{x_{16}}-x_{11}\ptl_{x_{20}}
-x_{15}\ptl_{x_{23}}-x_{19}\ptl_{x_{25}},\eqno(10.2.29)$$
\begin{eqnarray*}F_{(0,1,2,1)}|_{\msr A}&=&-x_1\ptl_{x_5}
+x_2\ptl_{x_{7}}+x_6(\ptl_{x_{14}}-2\ptl_{x_{13}})+x_8\ptl_{x_{16}}
+x_9\ptl_{x_{17}}\\
& &-x_{10}\ptl_{x_{18}}-x_{11}\ptl_{x_{19}}-x_{13}\ptl_{x_{21}}
-x_{20}\ptl_{x_{25}}+x_{22}\ptl_{x_{26}},\hspace{2.9cm}(10.2.30)
\end{eqnarray*}
\begin{eqnarray*}F_{(1,1,1,1)}|_{\msr A}&=&
-x_1\ptl_{x_6}-x_3\ptl_{x_9}-x_4\ptl_{x_{11}}-x_5(\ptl_{x_{13}}+\ptl_{x_{14}})-x_7\ptl_{x_{15}}
+x_{12}\ptl_{x_{20}}\\ &
&-(x_{13}+x_{14})\ptl_{x_{22}}+x_{16}\ptl_{x_{23}}+x_{18}\ptl_{x_{24}}
+x_{21}\ptl_{x_{26}},\hspace{3.2cm}(10.2.31)
\end{eqnarray*}
$$F_{(1,2,2)}|_{\msr A}=-x_2\ptl_{x_{10}}
+x_3\ptl_{x_{12}}+x_7\ptl_{x_{18}}-x_9\ptl_{x_{20}}-x_{15}\ptl_{x_{24}}+x_{17}\ptl_{x_{25}},\eqno(10.2.32)$$
\begin{eqnarray*}F_{(1,1,2,1)}|_{\msr A}&=&x_1\ptl_{x_8}-x_2\ptl_{x_9}
+x_4(\ptl_{x_{14}}-2\ptl_{x_{13}})+x_5\ptl_{x_{16}}+x_7\ptl_{x_{17}}-x_{10}\ptl_{x_{20}}\\
&&-x_{11}\ptl_{x_{22}}-x_{13}\ptl_{x_{23}}+x_{18}\ptl_{x_{25}}-x_{19}\ptl_{x_{26}},\hspace{4.6cm}(10.2.33)
\end{eqnarray*}
$$F_{(0,1,2,2)}|_{\msr A}=x_1\ptl_{x_7}
+x_6\ptl_{x_{15}}+x_8\ptl_{x_{17}}-x_{10}\ptl_{x_{19}}-x_{12}\ptl_{x_{21}}
-x_{20}\ptl_{x_{26}},\eqno(10.2.34)$$
\begin{eqnarray*}F_{(1,2,2,1)}|_{\msr A}&=&-x_1\ptl_{x_{10}}
+x_2\ptl_{x_{11}}+x_3(\ptl_{x_{14}}-2\ptl_{x_{13}})
+x_5\ptl_{x_{18}}+x_7\ptl_{x_{19}}-x_8\ptl_{x_{20}}\\
& &-x_9\ptl_{x_{22}}-x_{13}\ptl_{x_{24}} -x_{16}\ptl_{x_{25}}
+x_{17}\ptl_{x_{26}},\hspace{4.8cm}(10.2.35)
\end{eqnarray*}
$$F_{(1,1,2,2)}|_{\msr A}=
-x_1\ptl_{x_9}+x_4\ptl_{x_{15}}+x_5\ptl_{x_{17}}
-x_{10}\ptl_{x_{22}} -x_{12}\ptl_{x_{23}}
+x_{18}\ptl_{x_{26}},\eqno(10.2.36)$$
$$F_{(1,2,2,2)}|_{\msr A}=
x_1\ptl_{x_{11}}+x_3\ptl_{x_{15}}+x_5\ptl_{x_{19}}-x_8\ptl_{x_{22}}
-x_{12}\ptl_{x_{24}}-x_{16}\ptl_{x_{26}},\eqno(10.2.37)$$
\begin{eqnarray*}E_{(1,2,3,1)}|_{\msr A}&=&-x_1\ptl_{x_{12}}
-x_2(\ptl_{x_{14}}+\ptl_{x_{13}})-x_3\ptl_{x_{16}}
+x_4\ptl_{x_{18}}-x_6\ptl_{x_{20}}+x_7\ptl_{x_{21}}\\
& &-x_9\ptl_{x_{23}} +x_{11}\ptl_{x_{24}}
-(x_{13}+x_{14})\ptl_{x_{25}}
+x_{15}\ptl_{x_{26}},\hspace{3.35cm}(10.2.38)
\end{eqnarray*}
\begin{eqnarray*}F_{(1,2,3,2)}|_{\msr A}&=&x_1(\ptl_{x_{13}}-2\ptl_{x_{14}})
+x_2\ptl_{x_{15}} -x_3\ptl_{x_{17}}+x_4\ptl_{x_{19}}
+x_5\ptl_{x_{21}}-x_6\ptl_{x_{22}}\\
& &-x_8\ptl_{x_{23}}+x_{10}\ptl_{x_{24}}-x_{12}\ptl_{x_{25}}
-x_{14}\ptl_{x_{26}},\hspace{4.7cm}(10.2.39)
\end{eqnarray*}
$$F_{(1,2,4,2)}|_{\msr A}=
x_1\ptl_{x_{16}}-x_2\ptl_{x_{17}} +x_4\ptl_{x_{21}}-x_6\ptl_{x_{23}}
+x_{10}\ptl_{x_{25}}-x_{11}\ptl_{x_{26}},\eqno(10.2.40)$$
$$F_{(1,3,4,2)}|_{\msr A}=
-x_1\ptl_{x_{18}}+x_2\ptl_{x_{19}}-x_3\ptl_{x_{21}}
+x_6\ptl_{x_{24}}-x_8\ptl_{x_{25}}+x_9\ptl_{x_{26}},\eqno(10.2.41)$$
$$E_{(2,3,4,2)}|_{\msr A}=x_1\ptl_{x_{20}}
-x_2\ptl_{x_{22}} +x_3\ptl_{x_{23}} -x_4\ptl_{x_{24}}
+x_5\ptl_{x_{25}}-x_7\ptl_{x_{26}},\eqno(10.2.42)$$
\begin{eqnarray*}\hspace{1cm}h_1|_{\msr A}&=&x_4\ptl_{x_4}+x_5\ptl_{x_5}-x_6\ptl_{x_6}+x_7\ptl_{x_7}-x_8\ptl_{x_8}-x_9\ptl_{x_9}
+x_{18}\ptl_{x_{18}}\\
&&+x_{19}\ptl_{x_{19}}-x_{20}\ptl_{x_{20}}+x_{21}\ptl_{x_{21}}-x_{22}\ptl_{x_{22}}
-x_{23}\ptl_{x_{23}},\hspace{2.9cm}(10.2.43)
\end{eqnarray*}
\begin{eqnarray*}\hspace{1cm}h_2|_{\msr A}&=&x_3\ptl_{x_3}-x_4\ptl_{x_4}+x_8\ptl_{x_8}
+x_9\ptl_{x_9}-x_{10}\ptl_{x_{10}}-x_{11}\ptl_{x_{11}}
+x_{16}\ptl_{x_{16}}\\&
&+x_{17}\ptl_{x_{17}}-x_{18}\ptl_{x_{18}}-x_{19}\ptl_{x_{19}}+x_{23}\ptl_{x_{23}}
-x_{24}\ptl_{x_{24}},\hspace{2.9cm}(10.2.44)
\end{eqnarray*}
\begin{eqnarray*}\hspace{1cm}h_3|_{\msr A}&=&x_2\ptl_{x_2}-x_3\ptl_{x_3}+x_4\ptl_{x_4}-x_5\ptl_{x_5}+x_6\ptl_{x_6}-x_8\ptl_{x_8}
+x_{10}\ptl_{x_{10}}\\&
&+2x_{11}\ptl_{x_{11}}-x_{12}\ptl_{x_{12}}+x_{15}\ptl_{x_{15}}-2x_{16}\ptl_{x_{16}}-x_{17}\ptl_{x_{17}}
+x_{19}\ptl_{x_{19}}
\\
&&-x_{21}\ptl_{x_{21}}+x_{22}\ptl_{x_{22}}-x_{23}\ptl_{x_{23}}+x_{24}\ptl_{x_{24}}
-x_{25}\ptl_{x_{25}},\hspace{2.9cm}(10.2.45)
\end{eqnarray*}
\begin{eqnarray*}\hspace{1cm}h_4|_{\msr A}&=&x_1\ptl_{x_1}-x_2\ptl_{x_2}+x_5\ptl_{x_5}-x_7\ptl_{x_7}+x_8\ptl_{x_8}-x_9\ptl_{x_9}
+x_{10}\ptl_{x_{10}}\\&
&-x_{11}\ptl_{x_{11}}+2x_{12}\ptl_{x_{12}}-2x_{15}\ptl_{x_{15}}+x_{16}\ptl_{x_{16}}-x_{17}\ptl_{x_{17}}
+x_{18}\ptl_{x_{18}}
\\
&&-x_{19}\ptl_{x_{19}}+x_{20}\ptl_{x_{20}}-x_{22}\ptl_{x_{22}}+x_{25}\ptl_{x_{25}}
-x_{26}\ptl_{x_{26}}.\hspace{2.9cm}(10.2.46)
\end{eqnarray*}

Denote
$$\bar r=27-r\qquad\for\;\;r\in\ol{1,26}.\eqno(10.2.47)$$
Define an algebraic isomorphism $\tau$ on $\msr
A[\ptl_{x_1},...,\ptl_{x_{26}}]$ by
$$\tau(x_r)=x_{\bar r},\;\tau(\ptl_{x_r})=\ptl_{x_{\bar
r}}\qquad\for\;\; r\in \ol{1,26}\setminus\{13,14\}\eqno(10.2.48)$$
and
$$\tau(x_{13})=-x_{13},\;\tau(x_{14})=-x_{14},\;\tau(\ptl_{x_{13}})=-\ptl_{x_{13}},\;\tau(\ptl_{x_{14}})=
-\ptl_{x_{14}}.\eqno(10.2.49)$$ Then
$$F'_\varpi |_{\msr A}=\tau(F_\varpi |_{\msr A})\qquad\for\;\;\varpi\in S_{F_4}.\eqno(10.2.50)$$
 Thus we have given the full explicit formulas for the basic
oscillator representation of $F_4$.

\section{Decomposition of the  Representation of $F_4$}

In this section, we decompose the polynomial algebra $\msr A$ into a
direct sum of irreducible ${\msr G}^{F_4}$-submodules.

Suppose that
$$\eta_1=3\sum_{r=1}^{12}x_rx_{\bar
r}+ax_{13}^2+bx_{13}x_{14}+cx_{14}^2\eqno(10.3.1)$$ is an
$F_4$-invariant (cf . (10.2.47)), where $a,b,c$ are constants to be
determined. Then $F_{(1)}(\eta_1)=F_{(0,1)}(\eta_1)=0$ naturally
hold. Moreover,
$$0=F_{(0,0,1)}(\eta_1)=2(a-b)x_{11}x_{13}+[b-4c-3)]x_{11}x_{14},\eqno(10.3.2)$$
which gives
$$a=b=4c+3.\eqno(10.3.3)$$
Similarly, the constraint $E_{(0,0,0,1)}(\eta_1)=0$ yield
$b=c=4a+3$. So we have the quadratic invariant
$$\eta_1=3\sum_{r=1}^{12}x_rx_{\bar
r}-x_{13}^2-x_{13}x_{14}-x_{14}^2.\eqno(10.3.4)$$ This invariant
also gives a symmetric ${\msr G}^{F_4}$-invariant bilinear form on
$V$.

By (10.2.43)-(10.2.46), we try to find a quadratic singular vector
of the form
$$\zeta_1=x_1(a_1x_{13}+a_2x_{14})+a_3x_2x_{12}+a_4x_3x_{10}+a_5x_4x_8
+a_6x_5x_6,\eqno(10.3.5)$$ where $a_r$ are constants. Observe
$$0=F_{(1)}(\zeta_1)=(a_5+a_6)x_4x_5\lra
a_6=-a_5.\eqno(10.3.6)$$ Moreover,
$$0=F_{(0,1)}(\zeta_1)=(a_4+a_5)x_3x_8\lra
a_5=-a_4.\eqno(10.3.7)$$ Note
$$0=F_{(0,0,1)}(\zeta_1)=(a_1-2a_2)x_1x_{11}+(a_3-a_4)x_2x_{10}\lra
a_1=2a_2,\;\;a_3=a_4.\eqno(10.3.8)$$ Furthermore,
$$0=F_{(0,0,0,1)}(\zeta_1)=(a_2-2a_1-a_3)x_1x_{12}\lra
2a_1=a_2-a_3.\eqno(10.3.9)$$ Hence we have the singular vector
$$\zeta_1=x_1(2x_{13}+x_{14})-3x_2x_{12}-3x_3x_{10}+3x_4x_8
-3x_5x_6\eqno(10.3.10)$$ of weight $\lmd_4$. So it generates an
irreducible module that is isomorphic to the basic module $V$. Note
$$F'_{(1)}|_{\msr A}=-x_6\ptl_{x_4}-x_8\ptl_{x_5}-x_9\ptl_{x_7}+
x_{20}\ptl_{x_{18}}+x_{22}\ptl_{x_{19}}+x_{23}\ptl_{x_{21}},\eqno(10.3.11)$$
$$F'_{(0,1)}|_{\msr A}=-x_4\ptl_{x_3}-x_{10}\ptl_{x_8}-x_{11}\ptl_{x_{9}}
+x_{18}\ptl_{x_{16}}+x_{19}\ptl_{x_{17}}+x_{24}\ptl_{x_{23}},\eqno(10.3.12)$$
\begin{eqnarray*}\hspace{1cm}F'_{(0,0,1)}|_{\msr A}&=&x_3\ptl_{x_2}+x_5\ptl_{x_4}+x_8\ptl_{x_6}
-x_{12}\ptl_{x_{10}}+x_{16}(2\ptl_{x_{14}}-\ptl_{x_{13}})\\
&&+x_{14}\ptl_{x_{11}}+x_{17}\ptl_{x_{15}}-x_{21}\ptl_{x_{19}}-x_{23}\ptl_{x_{22}}-x_{25}\ptl_{x_{24}},
\hspace{2.1cm}(10.3.13)\end{eqnarray*}
\begin{eqnarray*}F'_{(0,0,0,1)}|_{\msr A}&=&x_2\ptl_{x_1}+x_7\ptl_{x_5}+x_9\ptl_{x_8}+
x_{11}\ptl_{x_{10}}+x_{15}(2\ptl_{x_{13}}-\ptl_{x_{14}})
\\ &
&+x_{13}\ptl_{x_{12}}-x_{17}\ptl_{x_{16}}-x_{19}\ptl_{x_{18}}-x_{22}\ptl_{x_{20}}
-x_{26}\ptl_{x_{25}}\hspace{3cm}(10.3.14)\end{eqnarray*} by
(10.2.48)-(10.2.50). To get a  basis of the module generated by
$\zeta_1$ compatible to $\{x_i\mid i\in\ol{1,26}\}$, we set
$$\zeta_2=F'_{(0,0,0,1)}(\zeta_1)=x_2(-x_{13}+x_{14})+3x_1x_{15}-3x_3x_{11}+3x_4x_9 -3x_6x_7,\eqno(10.3.15)$$
$$\zeta_3=F'_{(0,0,1)}(\zeta_2)=-x_3(x_{13}+2x_{14})+3x_1x_{17}+3x_2x_{16}
+3x_5x_9 -3x_7x_8,\eqno(10.3.16)$$
$$\zeta_4=-F'_{(0,1)}(\zeta_3)=-x_4(x_{13}+2x_{14})-3x_1x_{19}-3x_2x_{18} +3x_5x_{11}
-3x_7x_{10},\eqno(10.3.17)$$
$$\zeta_5=F'_{(0,0,1)}(\zeta_4)=x_5(-x_{13}+x_{14})+3x_1x_{21}
-3x_3x_{18}-3x_4x_{16} +3x_7x_{12},\eqno(10.3.18)$$
$$\zeta_6=-F'_{(1)}(\zeta_4)=-x_6(x_{13}+2x_{14})+3x_1x_{22}+3x_2x_{20}+3x_8x_{11}
-3x_9x_{10},\eqno(10.3.19)$$
$$\zeta_7=F'_{(0,0,0,1)}(\zeta_5)=x_7(2x_{13}+x_{14})+3x_2x_{21}
+3x_3x_{19}+3x_4x_{17}-3x_5x_{15}
 ,\eqno(10.3.20)$$
$$\zeta_8=-F'_{(1)}(\zeta_5)=x_8(-x_{13}+x_{14})
-3x_1x_{23}+3x_3x_{20}-3x_6x_{16} +3x_9x_{12},\eqno(10.3.21)$$
$$\zeta_9=-F'_{(1)}(\zeta_7)=x_9(2x_{13}+x_{14})-3x_2x_{23}
-3x_3x_{22} +3x_6x_{17}-3x_8x_{15},\eqno(10.3.22)$$
$$\zeta_{10}=-F'_{(0,1)}(\zeta_8)=x_{10}(-x_{13}+x_{14})
+3x_1x_{24}+3x_4x_{20}+3x_6x_{18} +3x_{11}x_{12},\eqno(10.3.23)$$
$$\zeta_{11}=-F'_{(0,1)}(\zeta_9)=x_{11}(2x_{13}+x_{14})+3x_2x_{24}
-3x_4x_{22} -3x_6x_{19}-3x_{10}x_{15},\eqno(10.3.24)$$
$$\zeta_{12}=-F'_{(0,0,1)}(\zeta_{10})=-x_{12}(x_{13}+2x_{14})+3x_1x_{25}
-3x_5x_{20}-3x_8x_{18} -3x_{10}x_{16},\eqno(10.3.25)$$
\begin{eqnarray*}\zeta_{13}=F'_{(0,0,0,1)}(\zeta_{12})&=&-x_{13}(x_{13}+2x_{14})
-3x_1x_{26}+3x_2x_{25}+3x_5x_{22}-3x_7x_{20}\\ &
&+3x_8x_{19}-3x_9x_{18}
+3x_{10}x_{17}-3x_{11}x_{16},\hspace{3.1cm}(10.3.26)
\end{eqnarray*}
\begin{eqnarray*}\zeta_{14}=F'_{(0,0,1)}(\zeta_{11})&=&x_{14}(2x_{13}+x_{14})-3x_2x_{25}+3x_3x_{24}+3x_4x_{23}
-3x_5x_{22}\\ & & +3x_6x_{21}-3x_8x_{19}-3x_{10}x_{17}+3x_{12}x_{15}
,\hspace{3.4cm}(10.3.27)
\end{eqnarray*}
$$\zeta_r=\tau(\zeta_{\bar r})\qquad\for\;\;r\in\ol{15,27},\eqno(10.3.28)$$
where $\tau$ is an algebra automorphism determined by (10.2.48) and
(10.2.49). The above construction shows that the map $x_r\mapsto
\zeta_r$ determine a module isomorphism from $V$ to the module
generated by $\zeta_1$. In particular,
$$u(x_i)=ax_j\Leftrightarrow u(\zeta_i)=a\zeta_j,\qquad
a\in\mbb{F},\;u\in\msr G^{F_4}.\eqno(10.3.29)$$

First
\begin{eqnarray*}\vt=(x_1\zeta_2-x_2\zeta_1)/3&=&x_1(-x_2x_{13}+x_1x_{15}-x_3x_{11}+x_4x_9
-x_6x_7)\\ & &+x_2(x_2x_{12}+x_3x_{10}-x_4x_8 +x_5x_6
)\hspace{3.7cm}(10.3.30)\end{eqnarray*} is a singular vector of
weight $\lmd_3$. Recall that the invariant $\eta_1$ in (10.3.4)
define an invariant bilinear form on $V$. Thus we have the following
cubic invariant
$$\eta_2=3\sum_{r=1}^{12}(\zeta_rx_{\bar r}+x_r\zeta_{\bar
r})-2x_{13}\zeta_{13}-x_{13}\zeta_{14}-x_{14}\zeta_{13}-2x_{14}\zeta_{14}.\eqno(10.3.31)$$
According to (10.3.10) and (10.3.15)-(10.3.28), we find
\begin{eqnarray*}\eta_2&=&9(1+\tau)[(x_2x_{12}+x_3x_{10}-x_4x_8+x_5x_6)x_{26}
+(x_3x_{11} -x_4x_9+x_6x_7)x_{25}\\ & &
+(x_7x_8-x_5x_9)x_{24}+x_{10}(x_4x_{23}+x_9x_{21})-x_{11}(x_5x_{23}+x_8x_{21}+x_{12}x_{17})
 \\ & &-x_{12}(x_7x_{22}+x_9x_{19})]
+2x_{13}^3+3x_{13}^2x_{14}-3x_{13}x_{14}^2-2x_{14}^3
+3x_1(2x_{13}+x_{14})x_{26}\\& &+3x_2(x_{13}+2x_{14})x_{25}
-3x_{13}[x_3x_{24}+x_4x_{23}+x_5x_{22}+x_6x_{21}-2x_7x_{20}+x_8x_{19}\\
&&-2x_9x_{18}+x_{10}x_{17}-2x_{11}x_{16}+x_{12}x_{15}]-3x_{14}[2x_3x_{24}+2x_4x_{23}
-x_5x_{22}\\& &+2x_6x_{21}-x_7x_{20}
-x_8x_{19}-x_9x_{18}-x_{10}x_{17} -x_{11}x_{16}+2x_{12}x_{15}]
,\hspace{1.8cm}(10.3.32)\end{eqnarray*} where $\tau$ is an algebra
automorphism defined in (10.2.47)-(10.2.49).\psp

{\bf Theorem 10.3.1}. {\it Any polynomial $f$ in ${\msr A}$
satisfying the system of partial differential equations
$$F_\varpi(f)=0\qquad\for\;\;\varpi\in S_{F_4}^+\eqno(10.3.33)$$ must
be a polynomial in $x_1,\zeta_1,\vt,\eta_1,\eta_2$. In particular,
the elements
$$\{x_1^{m_1}\zeta_1^{m_2}\vt^{m_3}\eta_1^{m_4}\eta_2^{m_5}\mid
m_1,m_2,m_3,m_4,m_5\in\mbb{N}\}\eqno(10.3.34)$$ are linearly
independent singular vectors and any singular vector is a linear
combination of those in (10.3.34) with the same weight. The weight
of $x_1^{m_1}\zeta_1^{m_2}\vt^{m_3}\eta_1^{m_4}\eta_2^{m_5}$ is
$m_3\lmd_3+(m_1+m_2)\lmd_4$.}

{\it Proof.} First we note
$$x_1x_{14}=\zeta_1-2x_1x_{13}+3x_2x_{12}+3x_3x_{10}-3x_4x_8
+3x_5x_6,\eqno(10.3.35)$$
$$3x_1x_{15}=\zeta_2+x_2(x_{13}-x_{14})+3x_3x_{11}-3x_4x_9 +3x_6x_7,\eqno(10.3.36)$$
$$3x_1x_{17}=\zeta_3+x_3(x_{13}+2x_{14})-3x_2x_{16}
-3x_5x_9 +3x_7x_8,\eqno(10.3.37)$$
$$3x_1x_{19}=3x_5x_{11}-\zeta_4-x_4(x_{13}+2x_{14})-3x_2x_{18}
-3x_7x_{10},\eqno(10.3.38)$$
$$3x_1x_{21}=\zeta_5+x_5(x_{13}-x_{14})
+3x_3x_{18}+3x_4x_{16} -3x_7x_{12},\eqno(10.3.39)$$
$$3x_1x_{22}=\zeta_6+x_6(x_{13}+2x_{14})-3x_2x_{20}-3x_8x_{11}
+3x_9x_{10},\eqno(10.3.40)$$
$$3x_1x_{23}=\zeta_8+x_8(x_{13}-x_{14})
-3x_3x_{20}+3x_6x_{16} -3x_9x_{12},\eqno(10.3.41)$$
$$3x_1x_{24}=\zeta_{10}+x_{10}(x_{13}-x_{14})
-3x_4x_{20}-3x_6x_{18}-3x_{11}x_{12},\eqno(10.3.42)$$
$$3x_2x_{25}+3x_1x_{26}=\eta_1-3\sum_{r=3}^{12}x_rx_{\bar
r}+x_{13}^2+x_{13}x_{14}+x_{14}^2,\eqno(10.3.43)$$
\begin{eqnarray*}\quad &&3[3(x_1x_{15}+x_3x_{11} -x_4x_9+x_6x_7)+x_2(x_{13}+2x_{14})]x_{25}
\\ & &+3[3(x_2x_{12}+x_3x_{10}-x_4x_8+x_5x_6)+x_1(2x_{13}+x_{14})]x_{26}
\\ &=&\eta_2-9x_1(x_{17}x_{24}-x_{19}x_{23}+x_{21}x_{22})-9x_2(x_{16}x_{24} -x_{18}x_{23}+x_{20}x_{21})\\ & &-3x_{13}^2x_{14}
-9(1+\tau)[(x_7x_8-x_5x_9)x_{24}+x_{10}(x_4x_{23}+x_9x_{21})-2x_{13}^3\\
& &-x_{11}(x_5x_{23}+x_8x_{21}+x_{12}x_{17})
-x_{12}(x_7x_{22}+x_9x_{19})]+3x_{13}x_{14}^2+2x_{14}^3
\\ & &+3x_{13}[x_3x_{24}+x_4x_{23}+x_5x_{22}+x_6x_{21}-2x_7x_{20}+x_8x_{19}-2x_9x_{18}
\\ & &+x_{10}x_{17}-2x_{11}x_{16}+x_{12}x_{15}]-3x_{14}[2x_3x_{24}+2x_4x_{23}
-x_5x_{22}\\ & &+2x_6x_{21}-x_7x_{20}
-x_8x_{19}-x_9x_{18}-x_{10}x_{17} -x_{11}x_{16}+2x_{12}x_{15}]
\hspace{2cm}(10.3.44)\end{eqnarray*} by (10.3.4), (10.3.10),
(10.3.15)-(10.3.19), (10.3.21), (10.3.23) and (10.3.32). Thus
$\{x_r\mid 16,18,20\neq r\in\ol{14,26}\}$ are rational functions in
$$\{x_r,\zeta_s,\eta_1,\eta_2\mid r\in\{\ol{1,13},16,18,20\};7,9\neq
s\in\ol{1,10}\}.\eqno(10.3.45)$$

Suppose that $f\in {\msr A}$ is a solution of (10.3.33). Write $f$
as a rational function $f_1$ in the variables of (10.3.45). In the
following calculations, we will always use (10.3.29). By (10.2.42),
$$0=F_{(2,3,4,2)}(f_1)=x_1\ptl_{x_{20}}(f_1).
\eqno(10.3.46)$$ So $f_1$ is independent of $x_{20}$. Moreover,
(10.2.41) gives
$$0=F_{(1,3,4,2)}(f_1)=-x_1\ptl_{x_{18}}(f_1).
\eqno(10.3.47)$$ Hence $f_1$ is independent of $x_{18}$.
Furthermore, (10.2.40) yields
$$0=F_{(1,2,4,2)}(f_1)=x_1\ptl_{x_{16}}(f_1).
\eqno(10.3.48)$$ Thus $f_1$ is independent of $x_{16}$. Successively
applying (10.2.39), (10.2.38), (10.2.37), (10.2.36) and (10.2.34) to
$f_1$, we obtain that $f_1$ is independent of $x_{13},x_{12},
x_{11}, x_9$ and $x_7$. Therefore, $f_1$ is a rational function in
$$\{x_r,\zeta_s,\eta_1,\eta_2\mid 7,9\neq r\in\ol{1,10};7,9\neq
s\in\ol{1,10}\}.\eqno(10.3.49)$$

By (10.2.25), (10.2.27), (10.2.30), (10.2.31), (10.2.33) and
(10.2.35),
$$0=F_{(0,0,1,1)}(f_1)=-x_1\ptl_{x_3}(f_1)-\zeta_1\ptl_{\zeta_3}(f_1),
\eqno(10.3.50)$$
$$0=F_{(0,1,1,1)}(f_1)=x_1\ptl_{x_4}(f_1)+\zeta_1\ptl_{\zeta_4}(f_1),
\eqno(10.3.51)$$
$$0=F_{(0,1,2,1)}(f_1)=-x_1\ptl_{x_5}(f_1)-\zeta_1\ptl_{\zeta_5}(f_1),
\eqno(10.3.52)$$
$$0=F_{(1,1,1,1)}(f_1)=-x_1\ptl_{x_6}(f_1)-\zeta_1\ptl_{\zeta_6}(f_1),
\eqno(10.3.53)$$
$$0=F_{(1,1,2,1)}(f_1)=x_1\ptl_{x_8}(f_1)+\zeta_1\ptl_{\zeta_8}(f_1),
\eqno(10.3.54)$$
$$0=F_{(1,2,2,1)}(f_1)=-x_1\ptl_{x_{10}}(f_1)-\zeta_1\ptl_{\zeta_{10}}(f_1).
\eqno(10.3.55)$$ Set
$$\eta_r=x_1\zeta_r-x_r\zeta_1,\;\;r\in\ol{3,6};\;\;\eta_7=x_1\zeta_8-x_8\zeta_1,
\;\;\eta_8=x_1\zeta_{10}-x_{10}\zeta_1.\eqno(10.3.56)$$ By the
characteristic method of solving linear partial differential
equations (e.g., cf. [X21]), we get that $f_1$ can be written as a
rational function $f_2$ in
$$\{x_r,\zeta_s,\eta_q\mid r,s=1,2;q\in\ol{1,8}\}.\eqno(10.3.57)$$

Next applying (10.2.28), (10.2.29) and (10.2.32) to $f_2$, we get
$$0=F_{(0,1,2)}(f_2)=-(x_1\zeta_2-\zeta_1x_2)\ptl_{\eta_5}(f_2)
=-3\vt \ptl_{\eta_5}(f_2),\eqno(10.3.58)$$
$$0=F_{(1,1,2)}(f_2)=3\vt\ptl_{\eta_7}(f_2),\qquad
0=F_{(1,2,2)}(f_2)=-3\vt\ptl_{\eta_8}(f_2)=0\eqno(10.3.59)$$ (cf.
(10.3.30)). Thus $f_2$ is independent of $\eta_5$, $\eta_7$ and
$\eta_8$. Furthermore, we apply (10.2.21), (10.2.24) and (10.2.26)
to $f_2$ and obtain
$$0=F_{(0,0,1)}(f_2)=-3\vt \ptl_{\eta_3}(f_2),\qquad 0=F_{(0,1,1)}(f_2)=3\vt
\ptl_{\eta_4}(f_2),\eqno(10.3.60)$$
$$0=F_{(1,1,1)}(f_2)=-3\vt
\ptl_{\eta_6}(f_2).\eqno(10.3.61)$$ Therefore, $f_2$ is a rational
function in $x_1,x_2,\zeta_1,\zeta_2,\eta_1,\eta_2$.  By (10.2.22),
$$0=F_{(0,0,0,1)}(f_2)=-x_1\ptl_{x_2}(f_2)-\zeta_1\ptl_{\zeta_2}(f_2).\eqno(10.3.62)$$
Again the characteristic method tell us that $f_2$ can be written as
a rational function $f_3$ in $x_1,\zeta_1,\vt,\eta_1,\eta_2$. Since
$f_2=f$ is a polynomial in $\{x_r\mid r\in\ol{1,26}\}$, Expressions
(10.3.30), (10.3.35), (10.3.43) and (10.3.44) imply that $f_2$ must
be a polynomial in $x_1,\zeta_1,\vt,\eta_1,\eta_2$. The other
statements follow directly.$\qquad\Box$\psp

Calculating the weights of the singular vectors in the above
theorem, we have:\psp

 {\bf Corollary 10.3.2}. {\it The space of polynomial ${\msr
G}^{F_4}$-invariants over its basic module is an subalgebra of
${\msr A}$ generated by $\eta_1$ and $\eta_2$}.\psp

Let $L(m_1,m_2,m_3,m_4,m_5)$ be the ${\msr G}^{F_4}$-submodule
generated by
$x_1^{m_1}\zeta_1^{m_2}\vt^{m_3}\eta_1^{m_4}\eta_2^{m_5}$. Note that
(10.2.2) is a Cartan root space decomposition $\msr G^{F_4}$ over
$\mbb F$.  Let ${\msr A}_k$ be the subspace of polynomials in ${\msr
A}$ with degree $k$. Then ${\msr A}_k$ is a finite-dimensional
weight ${\msr G}^{F_4}$-module by (10.2.43)-(10.2.46).
 Thus $L(m_1,m_2,m_3,m_4,m_5)$ is a
finite-dimensional irreducible ${\msr G}^{F_4}$-submodule with the
highest weight $m_3\lmd_3+(m_1+m_2)\lmd_4$. By Weyl's Theorem 2.3.6
of complete reducibility if $\mbb F=\mbb C$ or more generally by
Lemma 6.3.2 with $n_1=0$,
$${\msr A}=\bigoplus_{k=0}^\infty {\msr A}_k=\bigoplus_{m_1,m_2,m_3,m_4,m_5=0}^\infty
L(m_1,m_2,m_3,m_4,m_5).\eqno(10.3.63)$$ Denote by $d(k,l)$ the
dimension of the highest weight irreducible module with the weight
$k\lmd_3+l\lmd_4$. The above equation implies the following
combinatorial identity:
$$\frac{1}{(1-t)^{26}}=\frac{1}{(1-t^2)(1-t^3)}\sum_{k_1,k_2,k_3=0}^\infty
d(k_1,k_2+k_3)t^{3k_1+2k_2+k_3}.\eqno(10.3.64)$$ Multiplying
$(1-t)^2$ to the above equation, we obtain a new combinatorial
identity about twenty-four:
$$\frac{1}{(1-t)^{24}}=\frac{1}{(1+t)(1+t+t^2)}\sum_{k_1,k_2,k_3=0}^\infty
d(k_1,k_2+k_3)t^{3k_1+2k_2+k_3}.\eqno(10.3.65)$$ Equivalently, we
have:\psp

{\bf Corollary 10.3.3}. {\it The dimensions $d(k,l)$ of the
irreducible module with the weights $k\lmd_3+l\lmd_4$ are linearly
correlated by the following identity:}
$$(1+t)(1+t+t^2)=(1-t)^{24}\sum_{k_1,k_2,k_3=0}^\infty
d(k_1,k_2+k_3)t^{3k_1+2k_2+k_3}.\eqno(10.3.66)$$ \pse

Recall the quadratic invariants $\eta_1$ in (10.3.4). Dually we have
${\cal G}^{F_4}$-invariant Laplace operator
$$\Dlt_{F_4}=3\sum_{r=1}^{12}\ptl_{x_r}\ptl{x_{\bar
r}}-\ptl_{x_{13}}^2-\ptl_{x_{13}}\ptl_{x_{14}}-\ptl_{x_{14}}^2.\eqno(10.3.67)$$
Now the subspace of complex homogeneous harmonic polynomials with
degree $k$ is
$${\cal H}^{F_4}_k=\{f\in{\cal A}_k\mid
\Dlt_{F_4}(f)=0\}.\eqno(10.3.68)$$ Then ${\cal H}^{F_4}_1=V$. Assume
$k\geq 2$. Suppose that $k_1,k_2,m_1,m_2$ are nonnegative integers
such that
$$k_1+3k_2=m_1+3m_2+2=k.\eqno(10.3.69)$$
If $\Dlt_{F_4}(x_1^{k_1}\vt^{k_2})\neq 0$, then it is a singular
vector of degree $k-2$ with weight $k_2\lmd_3+k_1\lmd_4$. By Theorem
10.3.1, ${\cal A}_{k-2}$ does not contain a singular vector of such
weight. A contradiction. Thus $\Dlt_{F_4}(x_1^{k_1}\vt^{k_2})=0$. By
the same reason, $\Dlt_{F_4}(x_1^{m_1}\zeta_1\vt^{m_2})=0.$ Hence
the irreducible submodules
$$L(k_1,0,k_2,0,0),\;L(m_1,1,m_2,0,0,0)\subset {\cal
H}^{F_4}_k.\eqno(10.3.70)$$ This gives the following corollary:\psp

{\bf Corollary 10.3.4}. {\it The number of irreducible submodules
contained in the subspace ${\cal H}^{F_4}_k$ of complex homogeneous
harmonic polynomials with degree $k\geq 2$ is $\geq
 [\!|k/3|\!]+[\!|(k-2)/3|\!]+2$.}

\chapter{Representations of $E_6$}

 Firstly in this chapter, we prove that the space of homogeneous polynomial solutions with
degree $m$ in 27 variables for the  Dickson invariant differential
equation is exactly a direct sum of $\llbracket m/2 \rrbracket+1$
explicitly determined irreducible $E_6$-submodules and the whole
polynomial algebra is a free module over the polynomial algebra in
the Dickson invariant generated by these solutions. Thus we obtain a
cubic $E_6$-generalization of the classical theorem on harmonic
polynomials. The result was due to us [X17].

Secondly we construct a representation of the simple Lie algebra of
type $E_6$ on the polynomial algebra in 16 variables, which gives a
fractional representation of the corresponding Lie group on
16-dimensional space.  Using this representation and Shen's idea of
mixed product (cf. [Sg]), we construct a new functor from the
category of $D_5$-modules to the category of $E_6$-modules. A
condition for the functor to map a finite-dimensional irreducible
$D_5$-module to an infinite-dimensional irreducible $E_6$-module is
obtained. Our results  yield explicit constructions of certain
infinite-dimensional irreducible weight $E_6$-modules with
finite-dimensional weight subspaces.  In our approach, the idea of
Kostant's characteristic identities plays a key role. This part is
taken from [X23].

In the above work, we found a one-parameter ($\mfk{c}$) family of
inhomogeneous first-order differential operator representations of
the simple Lie algebra of type $E_6$ in $16$ variables.  Letting
these operators act on the space of exponential-polynomial functions
that depend on a parametric vector $\vec a\in \mbb
F^{16}\setminus\{\vec 0\}$, we prove that the space forms an
irreducible $E_6$-module for any constant $\mfk{c}$ if $\vec a$ is
not on an explicitly given projective algebraic variety. Certain
equivalent combinatorial properties of the spin oscillator
representation of $D_5$ play key roles in our proof. The result is
taken from our work [X27].

\section{ Basic Oscillator Representation}

In this section, we present the bosonic oscillator tepresentation of
$E_6$ over its 27-dimensional irreducible module.

 First we go back to the construction of the simple simple Lie algebra $\msr G^X$ in (4.4.15)-(4.4.25)
 with $X=E_6,\;E_7$ and $E_8$. Recall the
Dynkin diagram of $E_7$:

\begin{picture}(93,20)
\put(2,0){$E_7$:}\put(21,0){\circle{2}}\put(21,
-5){1}\put(22,0){\line(1,0){12}}\put(35,0){\circle{2}}\put(35,
-5){3}\put(36,0){\line(1,0){12}}\put(49,0){\circle{2}}\put(49,
-5){4}\put(49,1){\line(0,1){10}}\put(49,12){\circle{2}}\put(52,10){2}\put(50,0){\line(1,0){12}}
\put(63,0){\circle{2}}\put(63,-5){5}\put(64,0){\line(1,0){12}}\put(77,0){\circle{2}}\put(77,
-5){6}\put(78,0){\line(1,0){12}}\put(91,0){\circle{2}}\put(91,
-5){7}
\end{picture}
\vspace{0.7cm}

 \noindent Let $\{\al_i\mid i\in\ol{1,7}\}$ be the
simple positive roots corresponding to the vertices in the diagram,
and let $\Phi_{E_7}$ be the root system of $E_7$. The simple Lie
algebra of type $E_7$ is
 $$\msr G^{E_7}=H\oplus\bigoplus_{\al\in\Phi_{E_7}}\mbb FE_{\al},\qquad H=H_{E_7}=\sum_{i=1}^7\mbb F\al_i,\eqno(11.1.1)$$
with the Lie bracket given in (4.4.24) and (4.4.25). Note that the
Dynkin diagram of $E_6$ is a sub-diagram of that of $E_7$. Set
$$H_{E_6}=\sum_{i=1}^6\mbb F\al_i,\qquad
\Phi_{E_6}=\Phi_{E_7}\bigcap H_{E_6}.\eqno(11.1.2)$$ We take the
simple Lie algebra $\msr G^{E_6}$ of type $E_6$ as the Lie
subalgebra
$$\msr G^{E_6}=H_{E_6}\oplus\bigoplus_{\al\in
\Phi_{E_6}}\mbb FE_{\al}.\eqno(11.1.3)$$

Recall the notion in (4.4.42)-(4.4.44). Denote
$$\mfk b_1=E'_{(0,0,0,0,0,0,1)},\;\;\mfk b_2=E'_{(0,0,0,0,0,1,1)},
\;\;\mfk b_3=E'_{(0,0,0,0,1,1,1)},\eqno(11.1.4)$$
$$\mfk
b_4=E'_{(0,0,0,1,1,1,1)},,\;\;\mfk b_5=E'_{(0,0,1,1,1,1,1)},\;\;\mfk
b_6=E'_{(0,1,0,1,1,1,1)},\eqno(11.1.5)$$ $$\mfk
b_7=E'_{(0,1,1,1,1,1,1)}, \;\;\mfk b_8=E'_{(1,0,1,1,1,1,1)},
\;\;\mfk b_9=E'_{(0,1,1,2,1,1,1)},\eqno(11.1.6)$$ $$\mfk
b_{10}=E'_{(1,1,1,1,1,1,1)},\;\;\mfk
b_{11}=E'_{(0,1,1,2,2,1,1)},\;\; \mfk
b_{12}=E'_{(1,1,1,2,1,1,1)},\eqno(11.1.7)$$
$$\mfk b_{13}=E'_{(1,1,1,2,2,1,1)},\;\;\mfk b_{14}=E'_{(0,1,1,2,2,2,1)},\;\;
\mfk b_{15}=E'_{(1,1,2,2,1,1,1)},\eqno(11.1.8)$$
$$\mfk
b_{16}=E'_{(1,1,2,2,2,1,1)},\;\;\mfk
b_{17}=E'_{(1,1,1,2,2,2,1)},\;\;\mfk
b_{18}=E'_{(1,1,2,3,2,1,1)},\eqno(11.1.9)$$ $$\mfk
b_{19}=E'_{(1,1,2,2,2,2,1)}, \;\mfk
b_{20}=E'_{(1,2,2,3,2,1,1)},\;\;\mfk
b_{21}=E'_{(1,1,2,3,2,2,1)},\eqno(11.1.10)$$ $$\mfk
b_{22}=E'_{(1,1,2,3,3,2,1)},\;\;\mfk
b_{23}=E'_{(1,2,2,3,2,2,1)},\;\; \mfk
b_{24}=E'_{(1,2,2,3,3,2,1)},\eqno(11.1.11)$$
$$\mfk b_{25}=E'_{(1,2,2,4,3,2,1)},\;\;\mfk b_{26}=E'_{(1,2,3,4,3,2,1)},\;\;\mfk b_{27}=E'_{(2,2,3,4,3,2,1)},\eqno(11.1.12)$$
Then the subspace
$$V=\sum_{i=1}^{27}\mbb F\mfk b_i\eqno(11.1.13)$$ forms an irreducible $\msr G^{E_6}$-module
with respect to the adjoint representation of $\msr G^{E_7}$, $\mfk
b_1$ is a highest-weight vector of weight $\lmd_6$ and
$$\sum_{\al\in\Phi_{E_7}^-}\mbb FE_\al=\sum_{\be\in\Phi_{E_6}^-}\mbb FE_\be\oplus V.\eqno(11.1.14)$$

 Write
$$[u,\mfk b_i]=\sum_{j=1}^{27}\vf_{i,j}(u)\mfk b_j\qquad\for\;\;u\in\msr
G^{F_4}.\eqno(11.1.15)$$ Set
$$\msr A=\mbb{F}[x_1,...,x_{27}]\eqno(11.1.16)$$
and define the {\it basic oscillator representation $\mfk r$ of}
$\msr G^{E_6}$ \index{ basic oscillator representation of! $\msr
G^{E_6}$} on $\msr A$ by
$$\mfk r(u)=\sum_{i,j=1}^{27}\vf_{i,j}(u)x_j\ptl_{x_i}\eqno(11.1.17)$$ (cf. (2.2.17)-(2.2.20)). Then
$\msr A$ forms a $\msr G^{E_6}$-module isomorphic to the symmetric
tensor $S(V)$ over $V$. More explicitly, we have the following
representation formulas for the positive root vectors in $\msr
G^{E_6}$ by (4.4.59)-(4.4.61):
$$\mfk r(E_{\al_1})=x_5\ptl_{x_8}+x_7\ptl_{x_{10}}+x_9\ptl_{x_{12}}+x_{11}
\ptl_{x_{13}}+x_{14}\ptl_{x_{17}}+x_{26}\ptl_{x_{27}},\eqno(11.1.18)$$
$$\mfk r(E_{\al_2})=x_4\ptl_{x_6}+x_5\ptl_{x_7}+x_8\ptl_{x_{10}}
-x_{18}\ptl_{x_{20}}-x_{21}\ptl_{x_{23}}-x_{22}\ptl_{x_{24}},\eqno(11.1.19)$$
$$\mfk r(E_{\al_3})=x_4\ptl_{x_5}+x_6\ptl_{x_7}+x_{12}\ptl_{x_{15}}+x_{13}\ptl_{x_{16}}
+x_{17}\ptl_{x_{19}}-x_{25}\ptl_{x_{26}},\eqno(11.1.20)$$
$$\mfk r(E_{\al_4})=x_3\ptl_{x_4}-x_7\ptl_{x_9}-x_{10}\ptl_{x_{12}}
-x_{16}\ptl_{x_{18}}-x_{19}\ptl_{x_{21}}-x_{24}\ptl_{x_{25}},\eqno(11.1.21)$$
$$\mfk r(E_{\al_5})=x_2\ptl_{x_3}-x_9\ptl_{x_{11}}
-x_{12}\ptl_{x_{13}}-x_{15}\ptl_{x_{16}}-x_{21}\ptl_{x_{22}}-x_{23}\ptl_{x_{24}},\eqno(11.1.22)$$
$$\mfk r(E_{\al_6})=x_1\ptl_{x_2}-x_{11}\ptl_{x_{14}}-x_{13}\ptl_{x_{17}}
-x_{16}\ptl_{x_{19}}-x_{18}\ptl_{x_{21}}-x_{20}\ptl_{x_{23}},\eqno(11.1.23)$$
$$\mfk r(E_{(1,0,1)})=x_4\ptl_{x_8}+x_6\ptl_{x_{10}}-x_9\ptl_{x_{15}}-x_{11}\ptl_{x_{16}}-x_{14}\ptl_{x_{19}}
-x_{25}\ptl_{x_{27}},\eqno(11.1.24)$$
$$\mfk r(E_{(0,1,0,1)})=x_3\ptl_{x_6}+x_5\ptl_{x_9}+x_8\ptl_{x_{12}}+x_{16}\ptl_{x_{20}}
+x_{19}\ptl_{x_{23}}-x_{22}\ptl_{x_{25}},\eqno(11.1.25)$$
$$\mfk r(E_{(0,0,1,1)})=x_3\ptl_{x_5}+x_6\ptl_{x_9}-x_{10}\ptl_{x_{15}}+x_{13}\ptl_{x_{18}}
+x_{17}\ptl_{x_{21}}+x_{24}\ptl_{x_{26}},\eqno(11.1.26)$$
$$\mfk r(E_{(0,0,0,1,1)})=x_2\ptl_{x_4}-x_7\ptl_{x_{11}}-x_{10}\ptl_{x_{13}}+x_{15}\ptl_{x_{18}}
-x_{19}\ptl_{x_{22}}+x_{23}\ptl_{x_{25}},\eqno(11.1.27)$$
$$\mfk r(E_{(0,0,0,0,1,1)})=x_1\ptl_{x_3}-x_9\ptl_{x_{14}}-x_{12}\ptl_{x_{17}}
-x_{15}\ptl_{x_{19}}+x_{18}\ptl_{x_{22}}+x_{20}\ptl_{x_{24}},\eqno(11.1.28)$$
$$\mfk r(E_{(1,0,1,1)})=x_3\ptl_{x_8}+x_6\ptl_{x_{12}}+x_7\ptl_{x_{15}}-x_{11}\ptl_{x_{18}}
-x_{14}\ptl_{x_{21}} +x_{24}\ptl_{x_{27}},\eqno(11.1.29)$$
$$\mfk r(E_{(0,1,1,1)})=x_3\ptl_{x_7}-
x_4\ptl_{x_9}+x_8\ptl_{x_{15}}-x_{13}\ptl_{x_{20}}-x_{17}\ptl_{x_{23}}+x_{22}\ptl_{x_{26}},\eqno(11.1.30)$$
$$\mfk r(E_{(0,1,0,1,1)})=x_2\ptl_{x_6}+x_5\ptl_{x_{11}}+x_8\ptl_{x_{13}}-x_{15}\ptl_{x_{20}}
+x_{19}\ptl_{x_{x_{24}}}+x_{21}\ptl_{x_{25}},\eqno(11.1.31)$$
$$\mfk r(E_{(0,0,1,1,1)})=x_2\ptl_{x_5}+x_6\ptl_{x_{11}}-x_{10}\ptl_{x_{16}}-x_{12}\ptl_{x_{18}}
+x_{17}\ptl_{x_{22}}-x_{23}\ptl_{x_{26}},\eqno(11.1.32)$$
$$\mfk r(E_{(0,0,0,1,1,1)})=x_1\ptl_{x_4}-x_7\ptl_{x_{14}}-x_{10}\ptl_{x_{17}}+x_{15}\ptl_{x_{21}}+x_{16}\ptl_{x_{22}}-x_{20}\ptl_{x_{25}},
\eqno(11.1.33)$$
$$\mfk r(E_{(1,1,1,1)})=x_3\ptl_{x_{10}}-x_4\ptl_{x_{12}}-x_5\ptl_{x_{15}}+x_{11}\ptl_{x_{20}}+x_{14}\ptl_{x_{23}}
+x_{22}\ptl_{x_{27}},\eqno(11.1.34)$$
$$\mfk r(E_{(1,0,1,1,1)})=x_2\ptl_{x_8}+x_6\ptl_{x_{13}}+x_7\ptl_{x_{16}}+x_9\ptl_{x_{18}}-x_{14}\ptl_{x_{22}}-x_{23}\ptl_{x_{27}},
\eqno(11.1.35)$$
$$\mfk r(E_{(0,1,1,1,1)})=x_2\ptl_{x_7}-x_4\ptl_{x_{11}}+x_8\ptl_{x_{16}}+x_{12}\ptl_{x_{20}}
-x_{17}\ptl_{x_{24}}-x_{21}\ptl_{x_{26}},\eqno(11.1.36)$$
$$\mfk r(E_{(0,1,0,1,1,1)})=x_1\ptl_{x_6}+x_5\ptl_{x_{14}}
+x_8\ptl_{x_{17}}-x_{15}\ptl_{x_{23}}-x_{16}\ptl_{x_{24}}-x_{18}\ptl_{x_{25}},\eqno(11.1.37)$$
$$\mfk r(E_{(0,0,1,1,1,1)})=x_1\ptl_{x_5}+x_6\ptl_{x_{14}}-x_{10}\ptl_{x_{19}}-x_{12}\ptl_{x_{21}}-x_{13}\ptl_{x_{22}}
+x_{20}\ptl_{x_{26}},\eqno(11.1.38)$$
$$\mfk r(E_{(1,1,1,1,1})=x_2\ptl_{x_{10}}-x_4\ptl_{x_{13}}-x_5\ptl_{x_{16}}-x_9\ptl_{x_{20}}+x_{14}\ptl_{x_{24}}-x_{21}\ptl_{x_{27}},
\eqno(11.1.39)$$
$$\mfk r(E_{(1,0,1,1,1,1)})=x_1\ptl_{x_8}+x_6\ptl_{x_{17}}+x_7\ptl_{x_{19}}+x_9\ptl_{x_{21}}+x_{11}\ptl_{x_{22}}
+x_{20}\ptl_{x_{27}},\eqno(11.1.40)$$
$$\mfk r(E_{(0,1,1,2,1)})=x_2\ptl_{x_9}-x_3\ptl_{x_{11}}+x_8\ptl_{x_{18}}-x_{10}\ptl_{x_{20}}
-x_{17}\ptl_{x_{25}}+x_{19}\ptl_{x_{26}},\eqno(11.1.41)$$
$$\mfk r(E_{(0,1,1,1,1,1)}|_{\msr
A}=x_1\ptl_{x_7}-x_4\ptl_{x_{14}}+x_8\ptl_{x_{19}}+x_{12}\ptl_{x_{23}}+x_{13}\ptl_{x_{24}}
+x_{18}\ptl_{x_{26}},\eqno(11.1.42)$$
$$\mfk r(E_{(1,1,1,2,1)})=x_2\ptl_{x_{12}}-x_3\ptl_{x_{13}}
-x_5\ptl_{x_{18}}+x_7\ptl_{x_{20}}+x_{14}\ptl_{x_{25}}+x_{19}\ptl_{x_{27}},
\eqno(11.1.43)$$
$$\mfk r(E_{(1,1,1,1,1,1)})=x_1\ptl_{x_{10}}
-x_4\ptl_{x_{17}}-x_5\ptl_{x_{19}}-x_9\ptl_{x_{23}}-x_{11}\ptl_{x_{24}}+x_{18}\ptl_{x_{27}},
\eqno(11.1.44)$$
$$\mfk r(E_{(0,1,1,2,1,1)})=x_1\ptl_{x_9}-x_3\ptl_{x_{14}}+x_8\ptl_{x_{21}}-x_{10}\ptl_{x_{23}}
+x_{13}\ptl_{x_{25}}-x_{16}\ptl_{x_{26}},\eqno(11.1.45)$$
$$\mfk r(E_{(1,1,2,2,1)})=x_2\ptl_{x_{15}}-x_3\ptl_{x_{16}}
+x_4\ptl_{x_{18}}-x_6\ptl_{x_{20}}-x_{14}\ptl_{x_{26}}-x_{17}\ptl_{x_{27}},
\eqno(11.1.46)$$
$$\mfk r(E_{(1,1,1,2,1,1)})=x_1\ptl_{x_{12}}
-x_3\ptl_{x_{17}}-x_5\ptl_{x_{21}}+x_7\ptl_{x_{23}}-x_{11}\ptl_{x_{25}}-x_{16}\ptl_{x_{27}},
\eqno(11.1.47)$$
$$\mfk r(E_{(0,1,1,2,2,1)})=x_1\ptl_{x_{11}}
-x_2\ptl_{x_{14}}+x_8\ptl_{x_{22}}-x_{10}\ptl_{x_{24}}
-x_{12}\ptl_{x_{25}}+x_{15}\ptl_{x_{26}},\eqno(11.1.48)$$
$$\mfk r(E_{(1,1,2,2,1,1)})=x_1\ptl_{x_{15}}
-x_3\ptl_{x_{19}}+x_4\ptl_{x_{21}}-x_6\ptl_{x_{23}}+x_{11}\ptl_{x_{26}}+x_{13}\ptl_{x_{27}},
\eqno(11.1.49)$$
$$\mfk r(E_{(1,1,1,2,2,1)})=x_1\ptl_{x_{13}}
-x_2\ptl_{x_{17}}-x_5\ptl_{x_{22}}+x_7\ptl_{x_{24}}
+x_9\ptl_{x_{25}}+x_{15}\ptl_{x_{27}},\eqno(11.1.50)$$
$$\mfk r(E_{(1,1,2,2,2,1)})=x_1\ptl_{x_{16}}
-x_2\ptl_{x_{19}}+x_4\ptl_{x_{22}}-x_6\ptl_{x_{24}}
-x_9\ptl_{x_{26}}-x_{12}\ptl_{x_{27}},\eqno(11.1.51)$$
$$\mfk r(E_{(1,1,2,3,2,1)})=x_1\ptl_{x_{18}}
-x_2\ptl_{x_{21}}+x_3\ptl_{x_{22}}-x_6\ptl_{x_{25}}+x_7\ptl_{x_{26}}+x_{10}\ptl_{x_{27}},
\eqno(11.1.52)$$
$$\mfk r(E_{(1,2,2,3,2,1)})=x_1\ptl_{x_{20}}
-x_2\ptl_{x_{23}}+x_3\ptl_{x_{24}}-x_4\ptl_{x_{25}}+x_5\ptl_{x_{26}}+x_8\ptl_{x_{27}}.
\eqno(11.1.53)$$

Recall that we also view $\al_i$ as the elements of ${\msr G}^{E_6}$
(cf. (11.1.1)). Then
$$\mfk r(\al_r)=\sum_{i=1}^{27}a_{i,r}x_i\ptl_{x_i}\qquad\for\;\;r\in\ol{1,6}\eqno(11.1.54)$$
with $a_{i,r}$ given in the following table:
 \begin{center}{\bf \large Table 11.1.1}\end{center}
{\begin{center}\begin{tabular}{|r||r|r|r|r|r|r||r||r|r|r|r|r|r|}\hline
$i$&$a_{i,1}$&$a_{i,2}$&$a_{i,3}$&$a_{i,4}$&$a_{i,5}$&$a_{i,6}$&
$i$&$a_{i,1}$&$a_{i,2}$&$a_{i,3}$&$a_{i,4}$&$a_{i,5}$& $a_{i,6}$
\\\hline\hline 1&0&0&0&0&0&1&2&0&0&0&0&1&$-1$\\\hline
3&0&0&0&1&$-1$&0&4&$0$&1&1&$-1$&0&0
\\\hline 5&1&1&$-1$&0&0&0&6&0&$-1$&1&0&0&0
\\\hline 7&1&$-1$&$-1$&1&0&0&8&$-1$&1&0&$0$&0&0
\\\hline 9&1&0&0&$-1$&1&0&10&$-1$&$-1$&0&1&0&0
\\\hline 11&1&0&0&0&$-1$&1&12&$-1$&0&1&$-1$&1&0
\\\hline 13&$-1$&0&1&0&$-1$&1&14&1&0&0&0&0&$-1$
\\\hline 15&0&0&$-1$&0&1&0&16&0&0&$-1$&1&$-1$&1
\\\hline 17&$-1$&0&1&0&0&$-1$&18&0&1&0&$-1$&0&1
\\\hline 19&0&0&$-1$&1&0&$-1$&20&0&$-1$&0&0&0&1
\\\hline 21&0&1&0&$-1$&1&$-1$&22&0&1&0&0&$-1$&0
\\\hline 23&0&$-1$&0&0&1&$-1$&24&0&$-1$&0&1&$-1$&0
\\\hline 25&0&0&1&$-1$&0&0&26&1&0&$-1$&0&0&0
\\\hline 27&$-1$&0&0&0&0&0
&&&&&&&\\\hline\end{tabular}\end{center}}

We define a symmetric linear operation $\tau$ on the space
$\sum_{i,j=1}^{27}\mbb Fx_i\ptl_{x_j}$ by
$$\tau(x_i\ptl_{x_j})=x_j\ptl_{x_i}.\eqno(11.1.55)$$
Then
$$\mfk r(E_{-\al})=-\tau(\mfk r(E_{\al}))\qquad\for\;\;\al\in\Phi_{E_6}^+\eqno(11.1.56)$$
by the second equations in (4.4.19), (4.4.20) and (4.4.25).

\section{Decomposition of the  Oscillator Representation}

In this section, we decompose $\msr A$ into a direct sum of
irreducible $\msr G^{E_6}$-submodules and prove that $\msr A$ is a
free module over the algebra of $\msr G^{E_6}$-invariants generated
by the solutions of Dickson's cubic invariant partial differential
equations.

Recall that  a singular vector of ${\msr G}^{E_6}$ is a nonzero
weight vector annihilated by positive root vectors. According to
Table 11.1.1 and (11.1.18)-(11.1.23), we  find a singular vector
$$\zeta_1=x_1x_{14}+x_2x_{11}+x_3x_9+x_4x_7-x_5x_6
\eqno(11.2.1)$$ of weight $\lmd_1$, which generates a 27-dimensional
irreducible ${\msr G}^{E_6}$-module
$$U=\sum_{i=1}^{27}\mbb F\zeta_i\eqno(11.2.2)$$ with
$$\zeta_2=x_1x_{17}+x_2x_{13}+x_3x_{12}+x_4x_{10}-x_6x_8,
\eqno(11.2.3)$$
$$\zeta_3=x_1x_{19}+x_2x_{16}+x_3x_{15}+x_5x_{10}-x_7x_8,
\eqno(11.2.4)$$
$$\zeta_4=-x_1x_{21}-x_2x_{18}+x_4x_{15}-x_5x_{12}
+x_8x_9, \eqno(11.2.5)$$
$$\zeta_5=x_1x_{22}-x_3x_{18}-x_4x_{16}+x_5x_{13}
-x_8x_{11}, \eqno(11.2.6)$$
$$\zeta_6=x_1x_{23}+x_2x_{20}+x_6x_{15}-x_7x_{12}
+x_9x_{10}, \eqno(11.2.7)$$
$$\zeta_7=-x_1x_{24}+x_3x_{20}-x_6x_{16}+x_7x_{13}
-x_{10}x_{11}, \eqno(11.2.8)$$
$$\zeta_8=x_2x_{22}+x_3x_{21}+x_4x_{19}-x_5x_{17}
+x_8x_{14}, \eqno(11.2.9)$$
$$\zeta_9=x_1x_{25}+x_4x_{20}+x_6x_{18}-x_9x_{13}
+x_{11}x_{12}, \eqno(11.2.10)$$
$$\zeta_{10}=-x_2x_{24}-x_3x_{23}+x_6x_{19}-x_7x_{17}
+x_{10}x_{14}, \eqno(11.2.11)$$
$$\zeta_{11}=x_1x_{26}-x_5x_{20}-x_7x_{18}+x_9x_{16}
-x_{11}x_{15}, \eqno(11.2.12)$$
$$\zeta_{12}=x_2x_{25}-x_4x_{23}-x_6x_{21}+x_9x_{17}
-x_{12}x_{14}, \eqno(11.2.13)$$
$$\zeta_{13}=x_3x_{25}+x_4x_{24}+x_6x_{22}-x_{11}x_{17}
+x_{13}x_{14}, \eqno(11.2.14)$$
$$\zeta_{14}=x_1x_{27}-x_8x_{20}-x_{10}x_{18}+x_{12}x_{16}
-x_{13}x_{15}, \eqno(11.2.15)$$
$$\zeta_{15}=
x_2x_{26}+x_5x_{23}+x_7x_{21}-x_9x_{19}+x_{14}x_{15},
\eqno(11.2.16)$$
$$\zeta_{16}=x_3x_{26}-x_5x_{24}-x_7x_{22}+x_{11}x_{19}-x_{14}x_{16},
\eqno(11.2.17)$$
$$\zeta_{17}=
-x_2x_{27}-x_8x_{23}-x_{10}x_{21}+x_{12}x_{19}-x_{15}x_{17},
\eqno(11.2.18)$$
$$\zeta_{18}=x_4x_{26}+x_5x_{25}+x_9x_{22}-x_{11}x_{21}+x_{14}x_{18},
\eqno(11.2.19)$$
$$\zeta_{19}=-x_3x_{27}+x_8x_{24}+x_{10}x_{22}-x_{13}x_{19}+x_{16}x_{17},
\eqno(11.2.20)$$
$$\zeta_{20}=-x_6x_{26}-x_7x_{25}+x_9x_{24}-x_{11}x_{23}+x_{14}x_{20},
\eqno(11.2.21)$$
$$\zeta_{21}=-x_4x_{27}-x_8x_{25}-x_{12}x_{22}+x_{13}x_{21}-x_{17}x_{18},
\eqno(11.2.22)$$
$$\zeta_{22}=x_5x_{27}-x_8x_{26}+x_{15}x_{22}-x_{16}x_{21}+x_{18}x_{19},
\eqno(11.2.23)$$
$$\zeta_{23}=x_6x_{27}+x_{10}x_{25}-x_{12}x_{24}+x_{13}x_{23}-x_{17}x_{20},
\eqno(11.2.24)$$
$$\zeta_{24}=-x_7x_{27}+x_{10}x_{26}+x_{15}x_{24}-x_{16}x_{23}+x_{19}x_{20},
\eqno(11.2.25)$$
$$\zeta_{25}=-x_9x_{27}+x_{12}x_{26}+x_{15}x_{25}-x_{18}x_{23}+x_{20}x_{21},
\eqno(11.2.26)$$
$$\zeta_{26}=-x_{11}x_{27}+x_{13}x_{26}+x_{16}x_{25}-x_{18}x_{24}+x_{20}x_{22},
\eqno(11.2.27)$$
$$\zeta_{27}=-x_{14}x_{27}+x_{17}x_{26}+x_{19}x_{25}-x_{21}x_{24}+x_{22}x_{23}.
\eqno(11.2.28)$$

By (11.1.18)-(11.1.23) and (11.2.1)-(11.2.28), we calculate
$$E_{\al_1}|_U=\zeta_1\ptl_{\zeta_2}+\zeta_{11}\ptl_{\zeta_{14}}-\zeta_{15}\ptl_{\zeta_{17}}
-\zeta_{16}\ptl_{\zeta_{19}}-\zeta_{18}\ptl_{\zeta_{21}}-\zeta_{20}\ptl_{\zeta_{23}},\eqno(11.2.29)$$
$$E_{\al_2}|_U=\zeta_4\ptl_{\zeta_6}+\zeta_5\ptl_{\zeta_7}+\zeta_8\ptl_{\zeta_{10}}
-\zeta_{18}\ptl_{\zeta_{20}}-\zeta_{21}\ptl_{\zeta_{23}}-\zeta_{22}\ptl_{\zeta_{24}},\eqno(11.2.30)$$
$$E_{\al_3}|_U=\zeta_2\ptl_{\zeta_3}-\zeta_9\ptl_{\zeta_{11}}
-\zeta_{12}\ptl_{\zeta_{15}}-\zeta_{13}\ptl_{\zeta_{16}}-\zeta_{21}\ptl_{\zeta_{22}}-\zeta_{23}\ptl_{\zeta_{24}},
\eqno(11.2.31)$$
$$E_{\al_4}|_U=\zeta_3\ptl_{\zeta_4}+\zeta_7\ptl_{\zeta_9}+\zeta_{10}\ptl_{\zeta_{12}}
+\zeta_{16}\ptl_{\zeta_{18}}+\zeta_{19}\ptl_{\zeta_{21}}-\zeta_{24}\ptl_{\zeta_{25}},\eqno(11.2.32)$$
$$E_{\al_5}|_U=\zeta_4\ptl_{\zeta_5}+\zeta_6\ptl_{\zeta_7}+\zeta_{12}\ptl_{\zeta_{13}}+\zeta_{15}\ptl_{\zeta_{16}}
+\zeta_{17}\ptl_{\zeta_{19}}-\zeta_{25}\ptl_{\zeta_{26}},\eqno(11.2.33)$$
$$E_{\al_6}|_U=\zeta_5\ptl_{\zeta_8}+\zeta_7\ptl_{\zeta_{10}}+\zeta_9\ptl_{\zeta_{12}}+\zeta_{11}
\ptl_{\zeta_{15}}
-\zeta_{14}\ptl_{\zeta_{17}}-\zeta_{26}\ptl_{\zeta_{27}}.\eqno(11.2.34)$$

Define an anti-homomorphism $\dg$ from the associative algebra
$$\mbb{A}=\sum_{i_1,...,i_{27}=0}^\infty {\msr
A}\ptl_{x_1}^{i_1}\cdots\ptl_{x_{27}}^{i_{27}}\eqno(11.2.35)$$ of
differential operators to the associative algebra
$$\mbb{A}'=\sum_{\ell_1,...,\ell_{27},i_1,...,i_{27}=0}^\infty \mbb F\zeta_1^{\ell_1}\cdots\zeta_{27}^{\ell_{27}}
\ptl_{\zeta_1}^{i_1}\cdots\ptl_{\zeta_{27}}^{i_{27}}\eqno(11.2.36)$$
determined by
$$(x_1)^\dg=\ptl_{\zeta_{27}},\;\;(\ptl_{x_1})^\dg=\zeta_{27},\;\;(x_{27})^\dg=-\ptl_{\zeta_1},\;\;(\ptl_{x_{27}})^\dg=-\zeta_1,
\eqno(11.2.37)$$
$$(x_{14})^\dg=-\ptl_{\zeta_{14}},\;\;(\ptl_{x_{14}})^\dg=-\zeta_{14},\;\;(x_i)^\dg=\ptl_{\zeta_{28-i}},\;\;
(\ptl_{x_i})^\dg=\zeta_{28-i}\eqno(11.2.38)$$ for $14\neq
i\in\ol{2,26}$. Then the map $d\mapsto -d^\dg$ gives a Lie algebra
monomorphism from $\mfk r(\msr G^{E_6})\rta\mbb A'$. By
(11.1.18)-(11.1.23), (11.1.54)-(11.1.56), Table 11.1.1 and
(11.2.29)-(11.2.34), we have
$$u|_U=-(\mfk r(u))^\dg\qquad\for\;\;u\in\msr
G^{E_6}.\eqno(11.2.39)$$ In particular,
$$\al_r|_U=\sum_{i=1}^{27}b_{i,r}\zeta_i\ptl_{\zeta_i}\qquad\for\;\;r\in\ol{1,6}\eqno(11.2.40)$$
with $b_{i,r}$ given in the following table:
\begin{center}{\bf \large Table 11.2.1}\end{center}
\begin{center}\begin{tabular}{|r||r|r|r|r|r|r||r||r|r|r|r|r|r|}\hline
$i$&$b_{i,1}$&$b_{i,2}$&$b_{i,3}$&$b_{i,4}$&$b_{i,5}$&$b_{i,6}$&$i$&$b_{i,1}$&$b_{i,2}$&$b_{i,3}$&$b_{i,4}$&$b_{i,5}$&
$b_{i,6}$
\\\hline\hline 1&1&0&0&0&0&0&2&$-1$&0&1&0&0&0\\\hline 3&0&0&$-1$&1&0&0&4&$0$&1&0&$-1$&1&0
\\\hline 5&0&1&$0$&0&$-1$&1&6&$0$&$-1$&0&$0$&1&0 \\\hline
7&0&$-1$&$0$&1&$-1$&1&8&$0$&1&0&$0$&0&$-1$\\\hline
9&0&$0$&$1$&$-1$&0&1&10&$0$&$-1$&0&$1$&0&$-1$\\\hline
11&1&$0$&$-1$&$0$&0&1&12&$0$&0&1&$-1$&1&$-1$\\\hline
13&0&0&1&0&$-1$&0&14&$-1$&0&0&0&0&1\\\hline
15&1&0&$-1$&0&1&$-1$&16&1&0&$-1$&1&$-1$&0\\\hline
17&$-1$&0&0&0&1&$-1$&18&1&1&$0$&$-1$&0&0\\\hline
19&$-1$&0&0&1&$-1$&0&20&1&$-1$&$0$&0&0&0\\\hline
21&$-1$&1&1&$-1$&0&0&22&0&1&$-1$&0&0&0\\\hline
23&$-1$&$-1$&1&0&0&0&24&0&$-1$&$-1$&1&0&0\\\hline
25&0&0&0&$-1$&1&0&26&0&0&0&$0$&$-1$&1\\\hline 27&0&0&0&0&0&$-1$
&&&&&&&\\\hline\end{tabular}\end{center}\psp

According to Table 11.1.1, Table 11.2.1 and (11.1.18)-(11.1.23), we
find the following invariant
$$\sum_{14\neq
i\in\ol{1,26}}x_i\zeta_{28-i}-x_{14}\zeta_{14}-x_{27}\zeta_1=-3\eta\eqno(11.2.41)$$
with
\begin{eqnarray*}\eta &=&x_1(x_{17}x_{26}+x_{19}x_{25}-x_{21}x_{24}+x_{22}x_{23})
-x_{14}(x_1x_{27}-x_8x_{20}-x_{10}x_{18}+x_{12}x_{16} \\ &
&-x_{13}x_{15})
+x_2(x_{13}x_{26}+x_{16}x_{25}-x_{18}x_{24}+x_{20}x_{22})
-x_{11}(x_2x_{27}+x_8x_{23}+x_{10}x_{21}\\ &
&-x_{12}x_{19}+x_{15}x_{17})
+x_3(x_{12}x_{26}+x_{15}x_{25}-x_{18}x_{23}+x_{20}x_{21})
-x_9(x_3x_{27}-x_8x_{24}\\ &
&-x_{10}x_{22}+x_{13}x_{19}-x_{16}x_{17})
+x_4(x_{10}x_{26}+x_{15}x_{24}-x_{16}x_{23}+x_{19}x_{20})
-x_7(x_4x_{27}\\ &
&+x_8x_{25}+x_{12}x_{22}-x_{13}x_{21}+x_{17}x_{18})+x_5(x_{10}x_{25}-x_{12}x_{24}+x_{13}x_{23}-x_{17}x_{20})
\\ &
&+x_6(x_5x_{27}-x_8x_{26}+x_{15}x_{22}-x_{16}x_{21}+x_{18}x_{19}),
\hspace{4.7cm}(11.2.42)\end{eqnarray*} which is the Dickson's cubic
invariant (cf. [D]). \psp

{\bf Lemma 11.2.1}. {\it Any homogeneous singular vector in ${\msr
A}$ is a monomial in $x_1,\;\zeta_1$ and $\eta$}.

{\it Proof}. Note that
$$x_1x_{14}=\zeta_1-x_2x_{11}-x_3x_9-x_4x_7+x_5x_6
\eqno(11.2.43)$$
$$x_1x_{17}=\zeta_2-x_2x_{13}-x_3x_{12}-x_4x_{10}+x_6x_8,
\eqno(11.2.44)$$
$$x_1x_{19}=\zeta_3-x_2x_{16}-x_3x_{15}-x_5x_{10}
+x_7x_8, \eqno(11.2.45)$$
$$x_1x_{21}=-\zeta_4-x_2x_{18}+x_4x_{15}-x_5x_{12}
+x_8x_9, \eqno(11.2.46)$$
$$x_1x_{22}=\zeta_5+x_3x_{18}+x_4x_{16}-x_5x_{13}
+x_8x_{11}, \eqno(11.2.47)$$
$$x_1x_{23}=\zeta_6-x_2x_{20}-x_6x_{13}+x_7x_{12}
-x_9x_{10}, \eqno(11.2.48)$$
$$x_1x_{24}=-\zeta_7+x_3x_{20}-x_6x_{16}+x_7x_{13}
-x_{10}x_{11}, \eqno(11.2.49)$$
$$x_1x_{25}=\zeta_9-x_4x_{20}-x_6x_{18}+x_9x_{13}
-x_{11}x_{12}, \eqno(11.2.50)$$
$$x_1x_{26}=\zeta_{11}+x_5x_{20}+x_7x_{18}-x_9x_{16}
+x_{11}x_{15} \eqno(11.2.51)$$
 by (11.2.1), (11.2.3)-(11.2.7), (11.2.9), (11.2.11) and
(11.2.15). Moreover, (11.2.42) can be written as
\begin{eqnarray*}& &(x_5x_6-x_1x_{14}-x_2x_{11}-x_3x_9-x_4x_7)x_{27}
 \\&=&\eta-x_1(x_{17}x_{26}+x_{19}x_{25}-x_{21}x_{24}+x_{22}x_{23})
-x_{14}(x_8x_{20}+x_{10}x_{18}-x_{12}x_{16} \\ & &-x_{13}x_{15})
-x_2(x_{13}x_{26}+x_{16}x_{25}-x_{18}x_{24}+x_{20}x_{22})
+x_{11}(x_8x_{23}+x_{10}x_{21}\\ & &-x_{12}x_{19}+x_{15}x_{17})
-x_3(x_{12}x_{26}+x_{15}x_{25}-x_{18}x_{23}+x_{20}x_{21})
-x_9(x_8x_{24}\\ & &-x_{10}x_{22}-x_{13}x_{19}+x_{16}x_{17})
-x_4(x_{10}x_{26}+x_{15}x_{24}-x_{16}x_{23}+x_{19}x_{20})
\\ & &+x_7(x_8x_{25}+x_{12}x_{22}-x_{13}x_{21}+x_{17}x_{18})-x_5(x_{10}x_{25}-x_{12}x_{24}+x_{13}x_{23}\\ &
&-x_{17}x_{20})
+x_6(x_8x_{26}-x_{15}x_{22}+x_{16}x_{21}-x_{18}x_{19}).
\hspace{4.6cm}(11.2.52)\end{eqnarray*}

Let $f$ be any homogenous singular vector in ${\msr A}$. According
to the above equations, $f$ can be written as a rational function
$f_1$ in
$$\{x_i,\zeta_r,\eta\mid
i\in\{\ol{1,13},15,16,18,20\};\;r\in\{\ol{1,7},9,11\}\}.\eqno(11.2.53)$$
By (11.1.47)-(11.1.53) and (11.2.37)-(11.2.39), we have
$$0=E_{(1,1,1,2,1,1)}(f_1)=x_1\ptl_{x_{12}}(f_1), \eqno(11.2.54)$$
$$0=E_{(0,1,1,2,2,1)} (f_1)=x_1\ptl_{x_{11}}(f_1),\eqno(11.2.55)$$
$$0=E_{(1,1,2,2,1,1)}(f_1)=x_1\ptl_{x_{15}}(f_1),
\eqno(11.2.56)$$
$$0=E_{(1,1,1,2,2,1)}(f_1)=x_1\ptl_{x_{13}}(f_1),\eqno(11.2.57)$$
$$0=E_{(1,1,2,2,2,1)}(f_1)=x_1\ptl_{x_{16}}(f_1),\eqno(11.2.58)$$
$$0=E_{(1,1,2,3,2,1)}(f_1)=x_1\ptl_{x_{18}}(f_1),
\eqno(11.2.59)$$
$$0=E_{(1,2,2,3,2,1)}(f)=x_1\ptl_{x_{20}}(f).
\eqno(11.2.60)$$ So $f_1$ is independent of
$x_{11},x_{12},x_{13},x_{15},x_{16}, x_{18},x_{20}$; that is, $f_1$
is a rational function in
$$\{x_i,\zeta_r,\eta\mid
i\in\ol{1,10};\;r\in\{\ol{1,7},9,11\}\}.\eqno(11.2.61)$$

Next (11.1.40), (11.1.42)-(11.1.46) and (11.2.37)-(11.2.39) imply
that
$$0=E_{(1,0,1,1,1,1)}(f_1)=x_1\ptl_{x_8}(f_1),\eqno(11.2.62)$$
$$0=E_{(0,1,1,1,1,1)}(f_1)=x_1\ptl_{x_7}(f_1),\eqno(11.2.63)$$
$$0=E_{(1,1,1,2,1)}(f_1)=\zeta_1\ptl_{\zeta_9},
\eqno(11.2.64)$$
$$0=E_{(1,1,1,1,1,1)}(f_1)=x_1\ptl_{x_{10}}(f_1), \eqno(11.2.65)$$
$$0=E_{(0,1,1,2,1,1)}(f_1)=x_1\ptl_{x_9}(f_1),\eqno(11.2.66)$$
$$0=E_{(1,1,2,2,1)}(f_1)=-\zeta_1\ptl_{\zeta_{11}}(f_1).
\eqno(11.2.67)$$
 Hence $f_1$ is independent of
$x_7,x_8,x_9,x_{10},\zeta_9$ and $\zeta_{11}$; that is, $f_1$ is a
rational function in
$$\{x_i,\zeta_r,\eta\mid
i\in\ol{1,6},\;r\in\ol{1,7}\}.\eqno(11.2.68)$$

Expressions (11.1.18), (11.1.23), (11.1.24), (11.1.28), (11.1.29),
(11.1.33)-(11.1.35), (11.1.37)-(11.1.39) and (11.2.37)-(11.2.39)
yield that
$$0=E_{\al_1}(f_1)=\zeta_1\ptl_{\zeta_2}(f_1)\eqno(11.2.69)$$
$$0=E_{\al_6}(f_1)=x_1\ptl_{x_2}(f_1),\eqno(11.2.70)$$
$$0=E_{(1,0,1)}(f_1)=-\zeta_1\ptl_{\zeta_3}(f_1),\eqno(11.2.71)$$
$$0=E_{(0,0,0,0,1,1)}(f_1)=x_1\ptl_{x_3}(f_1),\eqno(11.2.72)$$
$$0=E_{(1,0,1,1)}(f_1)=\zeta_1\ptl_{\zeta_4}(f_1),\eqno(11.2.73)$$
$$0=E_{(0,0,0,1,1,1)}(f_1)=x_1\ptl_{x_4}(f_1), \eqno(11.2.74)$$
$$0=E_{(1,1,1,1)}(f_1)=\zeta_1\ptl_{\zeta_6}(f_1),\eqno(11.2.75)$$
$$0=E_{(1,0,1,1,1)}(f_1)=-\zeta_1\ptl_{\zeta_5}(f_1),
\eqno(11.2.76)$$
$$0=E_{(0,1,0,1,1,1)}(f_1)=x_1\ptl_{x_6}(f_1),\eqno(11.2.77)$$
$$0=E_{(0,0,1,1,1,1)}(f_1)=x_1\ptl_{x_5}(f_1),\eqno(11.2.78)$$
$$0=E_{(1,1,1,1,1}(f_1)=-\zeta_1\ptl_{\zeta_7}(f_1). \eqno(11.2.79)$$
Thus $f_1$ is independent of $\{x_i,\zeta_r\mid
i\in\ol{2,6},\;r\in\ol{2,7}\}$; that is, $f=f_1$ is a rational
function in $x_1,\;\zeta_1$ and $\eta$. By (11.2.1) and (11.2.42),
it must be a polynomial in $x_1,\;\zeta_1$ and $\eta$. Recall that
the weights of $x_1,\;\zeta_1$ and $\eta$ are $\lmd_6,\;\lmd_1$ and
$0$, respectively. The homogeneity of $f$ implies that it must be a
monomial in $x_1,\;\zeta_1$ and $\eta.\qquad\Box$\psp

Let $L(m_1,m_2,m_3)$ be the ${\msr G}^{E_6}$-submodule generated by
$x_1^{m_1}\zeta_1^{m_2}\eta^{m_3}$. Note that (11.1.3) is a Cartan
root space decomposition of $\msr G^{E_6}$ over $\mbb F$. Moreover,
(11.1.54) implies that ${\msr A}$ is a direct sum of weigh subspaces
of $\msr G^{E_6}$ and the subspaces of homogeneous polynomials are
finite-dimensional ${\msr G}^{E_6}$-submodules. Thus
$L(m_1,m_2,m_3)$ is a finite-dimensional irreducible submodule of
highest weight $m_1\lmd_1+m_2\lmd_6$. By Weyl's Theorem 2.3.6 of
complete reducibility if $\mbb F=\mbb C$ or more generally by Lemma
6.3.2 with $n_1=0$, we have
$${\msr A}=\sum_{m_1,m_2,m_2=0}^\infty L(m_1,m_2,m_3).\eqno(11.2.80)$$
Recall we denote by $V(\lmd)$ the finite-dimensional irreducible
module of highest weight $\lmd$. The above equation implies
$$\frac{1}{(1-q)^{27}}=\frac{1}{1-q^3}\sum_{m_1,m_2=0}^\infty (\dim
V(m_1\lmd_6+m_2\lmd_1))q^{m_1+2m_2}.\eqno(11.2.81)$$ Equivalently,
we have:\psp

{\bf Corollary 11.2.2}. {\it The following dimensional property of
irreducible ${\msr G}^{E_6}$-modules holds:}
$$(1-q)^{26}\sum_{m_1,m_2=0}^\infty (\mbox{\it dim}\:
V(m_1\lmd_1+m_2\lmd_6))q^{m_1+2m_2}=1+q+q^2.\eqno(11.2.82)$$ \pse

Set
$$W=\sum_{i=1}^{27}\mbb F\ptl_{x_i}.\eqno(11.2.83)$$
Then $W$ is isomorphic to the module of linear functions on $V$ via
$\ptl_{x_i}(x_j)=\dlt_{i,j}$. Indeed, the linear map determined by
$\ptl_{x_i}\mapsto (\ptl_{x_i})^\dg$ (cf. (11.2.37), (11.2.38)) is a
${\msr G}^{E_6}$-module isomorphism from $W$ to $U$ by (11.2.41). We
define a linear map $\Im:{\msr A}\rta
\mbb{R}[\ptl_{x_1},...,\ptl_{x_{27}}]$ by
$$\Im(x_1^{\al_1}x_2^{\al_2}\cdots x_{27}^{\al_{27}})=\ptl_{x_1}^{\al_1}\ptl_{x_2}^{\al_2}\cdots
\ptl_{x_{27}}^{\al_{27}}.\eqno(11.2.84)$$ Set
$$\msr D=\Im(\eta),\qquad{\msr
D}_1=\sum_{i=1}^{27}x_i\ptl_{x_i},\qquad{\msr
D}_2=\sum_{i=1}^{27}\zeta_i\Im(\zeta_i).\eqno(11.2.85)$$ Then ${\msr
D},\;\msr D_1$ and $\msr D_2$ are invariant differential operators;
that is, as operators on $\msr A$,
$$d\xi=\xi d\qquad\for\;\;d=\msr D,\msr D_1,\msr D_2;\;\xi\in{\msr G}^{E_6}.\eqno(11.2.86)$$

Note that Lemma 11.2.1 implies
$$V^2=L(2,0,0)+L(0,1,0).\eqno(11.2.87)$$
Symmetrically,
$$W^2=L'(0,2,0)+L'(1,0,0),\eqno(11.2.88)$$
where $L'(0,2,0)$ is a module generated by the highest weight vector
$\ptl_{27}^2$ with weight $2\lmd_1$ and $L'(1,0,0)$ is a module
generated by the highest weight vector $\Im(\zeta_{27})$ with weight
$\lmd_6$. Thus the subspace of invariants (the trivial submodule) in
$V^2W^2$ is two-dimensional. The trivial submodule of
$L(0,1,0)L'(1,0,0)$ is $\mbb F\msr D_2$. In $L(2,0,0)L'(0,2,0)$,
there exists an invariant $\msr D_3$ with a term
$x_1^2\ptl_{x_1}^2$. So any invariant in $V^2W^2$ must be in $\mbb
F\msr D_2+\mbb F\msr D_3$. In particular, the invariant differential
operator
$$[\msr D,\eta]=\msr D\eta-\eta\msr D=b_0+b_1\msr
D_1+b_2\msr D_2+b_3\msr D_3\eqno(11.2.89)$$ for some $b_s\in\mbb F$.
According to (11.2.42), $\eta$ does not contain $x_1^2$. So $b_3=0$.
Moreover, (11.2.42) also implies $b_0=45$.

According to (11.2.52), the coefficient of $x_{27}\ptl_{x_{27}}$ in
$[\msr D,\eta]$ must be $5$, which implies $b_1=5$. Observe that
there exists a unique monomial in $\eta$ containing $x_1x_{14}$,
which is $x_1x_{14}x_{27}$. Thus the coefficient of
$x_1x_{14}\ptl_{x_1}\ptl_{x_{14}}$ in $[\msr D,\eta]$ must be 1;
that is, $b_2=1$. So we have: \psp

{\bf Lemma 11.2.3}. {\it As operators on ${\msr A}$},
$$[\msr D,\eta]=45+5{\msr
D}_1+\msr D_2.\eqno(11.2.90)$$ \pse

Let $m_1$ and $m_2$ be nonnegative integers.
 If $\msr D(x_1^{m_1}\zeta_1^{m_2})\neq 0$, then it is also a
singular of degree $m_1+2m_2-3$ with the same weight
$m_1\lmd_6+m_2\lmd_1$. But Lemma 11.2.1 implies that any singular
vector with weight $m_1\lmd_6+m_2\lmd_1$ must has degree $\geq
m_1+2m_2$. This leads a contradiction. Thus
$${\msr
D}(x_1^{m_1}\zeta_1^{m_2})=0\qquad\for\;\;m_1,m_2\in\mbb{N}.\eqno(11.2.91)$$
Moreover, (11.2.76) implies
$$\msr D(L(m_1,m_2,0))=\{0\}\qquad\for\;\;m_1,m_2\in\mbb{N}.\eqno(11.2.92)$$
Since $\msr D_2(x_1^{m_1}\zeta_1^{m_2})$ is also a singular vector
of degree $m_1+2m_2$ with the same weight $m_1\lmd_6+m_2\lmd_1$, we
have
$${\msr
D}_2(x_1^{m_1}\zeta_1^{m_2})=cx_1^{m_1}\zeta_1^{m_2}\eqno(11.2.93)$$
for some $c\in\mbb F$. Let
$$x_i=0\qquad \for\;\;1,14\neq i\in\ol{1,27}\eqno(11.2.94)$$
in (11.2.83) and we get
\begin{eqnarray*}cx_1^{m_1+m_2}x_{14}^{m_2}&=&
\lim_{x_i\rta 0;\;8,10\neq
i\in\ol{2,11}}x_1x_{14}(\ptl_{x_1}\ptl_{x_{14}}+\ptl_{x_2}\ptl_{x_{11}}+\ptl_{{x_3}}\ptl_{x_9}
+\ptl_{x_4}\ptl_{x_7}-\ptl_{x_5}\ptl_{x_6})[x_1^{m_1}\\& &\times
(x_1x_{14}+x_2x_{11}+x_3x_9+x_4x_7-x_5x_6)^{m_2}]\\
&=&m_2(m_1+m_2+4)x_1^{m_1+m_2}x_{14}^{m_2}\hspace{6.1cm}(11.2.95)
\end{eqnarray*}
by (11.2.1)-(11.2.29); that is, $c=m_2(m_1+m_2+4)$. We get:\psp

{\bf Lemma 11.2.4}. {\it For $m_1,m_2\in\mbb{N}$,}
$$\msr D_2(x_1^{m_1}\zeta_1^{m_2})=m_2(m_1+m_2+4)x_1^{m_1}\zeta_1^{m_2}.\eqno(11.2.96)$$
\pse

According to Lemma 11.2.1,
$$V^4=L(4,0,0)+L(2,1,0)+L(1,0,1).\eqno(11.2.97)$$
Moreover, $L(1,0,1)=\eta V$. Thus the invariants in $V^4W$ are $\mbb
F\eta\msr D_1$. Hence
$$[\msr D_2,\eta]=c_1\eta+c_2\eta\msr D_1\qquad\mbox{for
some}\;\;c_1,c_2\in\mbb F.\eqno(11.2.98)$$ Letting the above
equation act on 1, we have
$$\msr D_2(\eta)=c_1\eta.\eqno(11.2.99)$$
By (11.2.1)-(11.2.28) and (11.2.42),
\begin{eqnarray*}& &c_1x_1x_{14}x_{27}=-\lim_{x_i\rta 0;\;14\neq
i\in\ol{2,16}}\eta=\lim_{x_i\rta 0;\;14\neq i\in\ol{2,16}}{\msr
D}_2(\eta)\\
&=&(x_1x_{14}\ptl_{x_1}\ptl_{x_{14}}+x_{14}x_{27}\ptl_{x_{14}}\ptl_{x_{27}}+x_1x_{27}\ptl_{x_1}\ptl_{x_{27}})(x_1x_{14}x_{27})
=3x_1x_{14}x_{27},\hspace{0.7cm}(11.2.100)\end{eqnarray*} So
$c_1=3$. Letting (11.2.98) act on $x_1$, we have:
$$\msr D_2(\eta x_1)=(3+c_2)\eta x_1.\eqno(11.2.101)$$
As (11.2.100),
\begin{eqnarray*}& &(3+c_2)x_1^2x_{14}x_{27}=-\lim_{x_i\rta 0;\;14\neq
i\in\ol{2,16}}(3+c_2)\eta x_1=-\lim_{x_i\rta 0;\;14\neq
i\in\ol{2,16}}{\msr
D}_2(\eta x_1)\\
&=&(x_1x_{14}\ptl_{x_1}\ptl_{x_{14}}+x_{14}x_{27}\ptl_{x_{14}}\ptl_{x_{27}}+x_1x_{27}\ptl_{x_1}\ptl_{x_{27}})(x_1^2x_{14}x_{27})
=4x_1^2x_{14}x_{27},\hspace{0.8cm}(11.2.102)\end{eqnarray*} Hence
$c_2=1$. We get: \psp

{\bf Lemma 11.2.5}. {\it As operators on ${\msr A}$,}
$$[\msr D_2,\eta]=\eta(3+\msr D_1).\eqno(11.2.103)$$
\pse

For $m,m_1,m_2\in\mbb{N}$ with $m>0$, we have
\begin{eqnarray*}\hspace{1cm}& &\msr D(\eta^mx_1^{m_1}\zeta_1^{m_2})=
[m(45+5m_1+m_2(m_1+m_2+4))\\ & &+ \sum_{s=1}^{m-1}
s(3+3(s+1)/2+m_1+2m_2)]
\eta^{m-1}x_1^{m_1}\zeta_1^{m_2}\neq0\hspace{3cm}(11.2.104)\end{eqnarray*}
by Lemmas 11.2.3-11.2.5. According to (11.2.80) and (11.2.104), we
have: \psp

{\bf Lemma 11.2.6}. {\it For any $0\neq f\in{\msr A}$},
$$\msr D(\eta f)\neq 0.\eqno(11.2.105)$$
\pse

The above lemma implies that
$$\{f\in{\msr A}\mid\msr D(f)\}=\sum_{m_1,m_2=0}^\infty
L(m_1,m_2,0).\eqno(11.2.106)$$ Recall that ${\msr A}_m$ is the
subspace of homogeneous polynomials of degree $m$ in ${\msr A}$.
Denote
$$\Phi_m=\{f\in{\msr A}_m\mid \msr(f)=0\}.\eqno(11.2.107)$$
In summary, we have the following version of the  main theorem. \psp

{\bf Theorem 11.2.7}. {\it The set
$\{x^{m_1}_1\zeta_1^{m_2}\eta^{m_3}\mid n_1,m_2,m_3\in\mbb{N}\}$ is
the set of all singular vectors in ${\msr A}$ up to a scalar
multiple. In particular, $\eta$ is the unique fundamental invariant
(up to constant) and  the identity
$$(1-q)^{26}\sum_{m_1,m_2=0}^\infty (\mbox{\it dim}\:
V(m_1\lmd_1+m_2\lmd_6))q^{m_1+2m_2}=1+q+q^2\eqno(11.2.108)$$ holds.
Furthermore,
$${\msr A}_k=\Phi_k\oplus \eta{\msr A}_{k-3}\qquad\mbox{\it
for}\;\; k\in\mbb{N}\eqno(11.2.109)$$ and
$$\Phi_m=\sum_{i=0}^{\llbracket m/2
\rrbracket}L(m-2i,i,0)\qquad{\it
for}\;\;m\in\mbb{N},\eqno(11.2.110)$$ where we treat ${\msr
A}_r=\{0\}$ if $r<0$.} \psp

As a consequence, the space of homogeneous polynomial solutions with
degree $m$ for the Dickson's invariant cubic partial differential
equation is exactly a direct sum of $\llbracket m/2 \rrbracket+1$
 irreducible $E_6$-submodules.

\section{Spin Oscillator Representation of $D_5$}

In this section, we give another construction of the spin oscillator
representation of $D_5$ which will be used later.

Let we go back to the construction of the simple simple Lie algebra
$\msr G^X$ in (4.4.15)-(4.4.25) with $X=E_6$ and $F(\cdot,\cdot)$
given in (4.4.64).

Recall the Dynkin diagram of $E_6$:

\begin{picture}(80,23)
\put(2,0){$E_6$:}\put(21,0){\circle{2}}\put(21,
-5){1}\put(22,0){\line(1,0){12}}\put(35,0){\circle{2}}\put(35,
-5){3}\put(36,0){\line(1,0){12}}\put(49,0){\circle{2}}\put(49,
-5){4}\put(49,1){\line(0,1){10}}\put(49,12){\circle{2}}\put(52,10){2}\put(50,0){\line(1,0){12}}
\put(63,0){\circle{2}}\put(63,-5){5}\put(64,0){\line(1,0){12}}\put(77,0){\circle{2}}\put(77,
-5){6}
\end{picture}
\vspace{0.7cm}

\noindent Let $\{\al_i\mid i\in\ol{1,6}\}$ be the simple positive
roots corresponding to the vertices in the diagram, and let
$\Phi_{E_6}$ be the root system of $E_6$. The simple Lie algebra of
type $E_6$ is
 $$\msr G^{E_6}=H\oplus\bigoplus_{\al\in\Phi_{E_6}}\mbb FE_{\al},\qquad H=H_{E_6}=\sum_{i=1}^6\mbb F\al_i,\eqno(11.3.1)$$
with the Lie bracket given in (4.4.24) and (4.4.25). Moreover, we
define a bilinear form $(\cdot|\cdot)$ on ${\msr G}^{E_6}$ by
$$(h_1|h_2)=(h_1,h_2),\;\; (h|E_{\al})=0,\;\;
 (E_{\al}|E_{\be})=-\dlt_{\al+\be,0}\eqno(11.3.2)$$
for $h_1,h_2\in H$ and $\al,\be\in \Phi_{E_6}$. It can be verified
that $(\cdot|\cdot)$ is a ${\msr G}^{E_6}$-invariant form; that is,
$$([u,v]|w)=(u|[v,w])\qquad\for\;\;u,v\in{\msr
G}^{E_6}.\eqno(11.3.3)$$

\subsection{Construction}

Let $E_{r,s}$ be the $10\times 10$ matrix with 1 as its
$(r,s)$-entry and $0$ as the others. Denote
$$A_{i,j}=E_{i,j}-E_{5+j,5+i},\;\;B_{i,j}=E_{i,5+j}-E_{j,5+i},\;\;C_{i,j}=E_{5+i,j}-E_{5+j,i}\eqno(11.3.4)$$
for $i,j\in\ol{1,5}$. Then the orthogonal Lie algebra
$$ o(10,\mbb F)=\sum_{i,j=1}^5
(\mbb FA_{i,j}+\mbb FB_{i,j}+\mbb FC_{i,j}).\eqno(11.3.5)$$ Write
$$Q^{D_5}=\sum_{i=1}^5\mbb{Z}\al_i,\qquad \Phi_{D_5}=\Phi_{E_6}\bigcap Q^{D_5}. \eqno(11.3.6)$$
Then
$${\msr G}^{D_5}=\sum_{i=1}^5\mbb F\al_i+\sum_{\be\in
\Phi_{D_5}}\mbb FE_\be\eqno(11.3.7)$$ forms a Lie subalgebra of
${\msr G}^{E_6}$, which is isomorphic to the orthogonal Lie algebra
$o(10,\mbb F)$.

 Denote by
$\Phi_{E_6}^+$ the set of positive roots of $E_6$ and by
$\Phi_{D_5}^+$ the set of positive roots of $D_5$. Recall the notion
$$\al_{(r_1,...,r_s)}=\sum_{i=1}^sr_i\al_i\in
\Phi_{E_6}^+\;\;\mbox{with}\;\;r_s\neq 0.\eqno(11.3.8)$$
 We find the
elements of $\Phi_{D_5}^+$: $\al_r\;(r\in\ol{1,5})$ and
$$\al_{(1,0,1)},\;\al_{(0,1,0,1)},\;\al_{(0,0,1,1)},\;\al_{(0,0,0,1,1)},\;
\al_{(1,0,1,1)},\;\al_{(0,1,1,1)},\;\al_{(0,1,0,1,1)},\;\al_{(0,0,1,1,1)},\eqno(11.3.9)$$
$$\al_{(1,1,1,1)},\,\al_{(1,0,1,1,1)},\al_{(0,1,1,1,1)},\al_{(1,1,1,1,1)},\al_{(0,1,1,2,1)},
\al_{(1,1,1,2,1)},\al_{(1,1,2,2,1)}.\eqno(11.3.10)$$
 Moreover, the elements
in $\Phi_{E_6}^+\setminus\Phi_{D_5}^+$ are:
$$\al_6,\;\al_{(0,0,0,0,1,1)},\;\al_{(0,0,0,1,1,1)},\;\al_{(0,0,1,1,1,1)},\;\al_{(0,1,0,1,1,1)},\;\al_{(0,1,1,1,1,1)},
\eqno(11.3.11)$$
$$
\al_{(1,0,1,1,1,1)},\;\al_{(1,1,1,1,1,1)},\;\al_{(0,1,1,2,1,1)},\;\al_{(1,1,1,2,1,1)},\;\al_{(0,1,1,2,2,1)},
\eqno(11.3.12)$$
$$\al_{(1,1,2,2,1,1)},\;\al_{(1,1,1,2,2,1)},\;\al_{(1,1,2,2,2,1)},\;\al_{(1,1,2,3,2,1)},\;\al_{(1,2,2,3,2,1)}.
\eqno(11.3.13)$$

 Recall the notions
$$E_{(r_1,...,r_s)}=E_{\al_{(r_1,...,r_s)}},\;\;(E_\al)'=E_{-\al}.\eqno(11.3.14)$$
Set
$$\xi_1=E_{\al_6},\;\xi_2=E_{(0,0,0,0,1,1)},\;\xi_3=E_{(0,0,0,1,1,1)},\;\xi_4=E_{(0,0,1,1,1,1)},\eqno(11.3.15)$$
$$\xi_5=E_{(0,1,0,1,1,1)},\;\xi_6=E_{(0,1,1,1,1,1)},\;
\xi_7=E_{(1,0,1,1,1,1)},\;\xi_8=E_{(1,1,1,1,1,1)},\eqno(11.3.16)$$
$$\xi_9=E_{(0,1,1,2,1,1)},\;\xi_{10}E_{(1,1,1,2,1,1)},\;\xi_{11}=E_{(0,1,1,2,2,1)},\;\xi_{12}=E_{(1,1,2,2,1,1)},
\eqno(11.3.17)$$
$$\xi_{13}=E_{(1,1,1,2,2,1)},\;\xi_{14}=E_{(1,1,2,2,2,1)},\;\xi_{15}=E_{(1,1,2,3,2,1)},\;\xi_{16}=E_{(1,2,2,3,2,1)}.
\eqno(11.3.18)$$
$$\eta_i=(\xi_i)'\qquad\for\;\;i\in\ol{1,16}.\eqno(11.3.19)$$
Write
$${\msr G}_+=\sum_{i=1}^{16}\mbb F\xi_i,\qquad{\msr
G}_-=\sum_{i=1}^{16}\mbb F\eta_i,\qquad{\msr G}_0={\msr
G}^{D_5}+\mbb{F}\al_6.\eqno(11.3.20)$$ It is straightforward to
verify that ${\msr G}_\pm$ are ableian Lie subalgebras of ${\msr
G}^{E_6}$, ${\msr G}_0$ is a reductive Lie subalgebra of ${\msr
G}^{E_6}$ and
$${\msr G}^{E_6}={\msr G}_-\oplus{\msr G}_0\oplus{\msr
G}_+.\eqno(11.3.21)$$ Moreover, ${\msr G}_\pm$ form irreducible
${\msr G}_0$-submodules with respect to the adjoint representation
of ${\msr G}^{E_6}$. Furthermore,
$$(\xi_i|\eta_j)=-\dlt_{i,j}\qquad\for\;\;i,j\in\ol{1,16}\eqno(11.3.22)$$
by (11.3.2). Expression (11.3.3) shows that ${\msr G}_+$ is
isomorphic to the dual ${\msr G}_0$-module of ${\msr G}_-$.

Set
$${\msr B}=\mbb F[x_1,x_2,...,x_{16}],\eqno(11.3.23)$$
the polynomial algebra in $x_1,x_2,...,x_{16}$. Write
$$[u,\eta_i]=\sum_{j=1}^{16}\vf_{i,j}(u)\eta_j\qquad\for\;\;i\in\ol{1,16},\;u\in{\msr
G}_0,\eqno(11.3.24)$$ where $\vf_{i,j}(u)\in\mbb F$. Define an
action of ${\msr G}_0$ on ${\msr B}$ by
$$u(f)=\sum_{i,j=1}^{16}\vf_{i,j}(u)x_j\ptl_{x_i}(f)\qquad\for\;\;u\in{\msr
G}_0,\;f\in{\msr B}.\eqno(11.3.25)$$ Then ${\msr B}$ forms a ${\msr
G}_0$-module and the subspace
$$V=\sum_{i=1}^{16}\mbb Fx_i\eqno(11.3.26)$$
forms a ${\msr G}_0$-submodule isomorphic to ${\msr G}_-$, where the
isomorphism is determined by $x_i\mapsto \eta_i$  for
$i\in\ol{1,16}$.

Denote by $\mbb{N}$ the set of nonnegative integers. Write
$$x^\al=\prod_{i=1}^{16}x_i^{\al_i},\;\;\ptl^\al=\prod_{i=1}^{16}\ptl_{x_i}^{\al_i}\qquad\for\;\;
\al=(\al_1,\al_2,...,\al_{16})\in\mbb{N}^{16}.\eqno(11.3.27)$$ Let
$$\mbb{B}=\sum_{\al\in\mbb{N}^{16}}{\msr B}\ptl^\al\eqno(11.3.28)$$
be the algebra of differential operators on ${\msr B}$. Then the
linear transformation $\tau$ determined by
$$\tau(x^\be\ptl^\gm)=x^\gm\ptl^\be\qquad\for\;\;\be,\gm\in\mbb{N}^{16}\eqno(11.3.29)$$
is an involutive anti-automorphism of $\mbb{B}$.

 According
to (4.4.24) and (4.4.25), we have the Lie algebra isomorphism $\nu:
o(10,\mbb F)\rta {\msr G}^{D_5}$ determined by the generators:
$$\nu(A_{1,2})=E_{\al_1},\;\nu(A_{2,3})=E_{\al_3},\;
\nu(A_{3,4})=E_{\al_4},\;\nu(A_{4,5})=E_{\al_5},\;\nu(B_{4,5})=E_{\al_2},\eqno(11.3.30)$$
$$
\nu(A_{2,1})=-E_{-\al_1},\;\nu(A_{3,2})=-E_{-\al_3},\;\;\nu(A_{4,3})=-E_{-\al_4},\;\;
\nu(A_{5,4})=-E_{-\al_5},\eqno(11.3.31)$$
$$\nu(C_{5,4})=-E_{-\al_2},\;\;\nu(A_{1,1})=\al_1+\al_3+\al_4+(
\al_2+\al_5)/2,\eqno(11.3.32)$$
$$\nu(A_{2,2})=\al_3+\al_4+(
\al_2+\al_5)/2,\;\;\nu(A_{3,3})=\al_4+\frac{1}{2}(
\al_2+\al_5),\eqno(11.3.33)$$
$$\nu(A_{4,4})=(
\al_2+\al_5)/2,\qquad\nu(A_{5,5})=( \al_2-\al_5)/2.\eqno(11.3.34)$$
Then ${\msr B}$ becomes an $o(10,\mbb F)$-module with respect to the
action
$$A(f)=\nu(A)(f)\qquad\for\;\;A\in o(10,\mbb F),\;f\in\mbb{\msr B}
.\eqno(11.3.35)$$

Thanks to (4.4.24), (4.4.25), (11.3.14)-(11.3.19), (11.3.24) and
(11.3.25), we have the following differential-operator
representation $\pi$ of $o(10,\mbb F)$ on $\msr B$:
$$\pi(A_{1,2})=x_4\ptl_{x_7}+x_6\ptl_{x_8}+x_9\ptl_{x_{10}}+x_{11}\ptl_{x_{13}},\eqno(11.3.36)$$
$$\pi(A_{2,3})=x_3\ptl_{x_4}+x_5\ptl_{x_6}+x_{10}\ptl_{x_{12}}+x_{13}\ptl_{x_{14}},\eqno(11.3.37)$$
$$\pi(A_{3,4})=-x_2\ptl_{x_3}-x_6\ptl_{x_9}-x_8\ptl_{x_{10}}+x_{14}\ptl_{x_{15}},\eqno(11.3.38)$$
$$\pi(A_{4,5})=-x_1\ptl_{x_2}+x_9\ptl_{x_{11}}+x_{10}\ptl_{x_{13}}+x_{12}\ptl_{x_{14}},\eqno(11.3.39)$$
$$\pi(B_{4,5})=-x_3\ptl_{x_5}-x_4\ptl_{x_6}-x_7\ptl_{x_8}+x_{15}\ptl_{x_{16}},\eqno(11.3.40)$$
$$\pi(A_{1,3})=-x_3\ptl_{x_7}-x_5\ptl_{x_8}+x_9\ptl_{x_{12}}+x_{11}\ptl_{x_{14}},\eqno(11.3.41)$$
$$\pi(A_{2,4})=x_2\ptl_{x_4}-x_5\ptl_{x_9}+x_8\ptl_{x_{12}}+x_{13}\ptl_{x_{15}},\eqno(11.3.42)$$
$$\pi(A_{3,5})=-x_1\ptl_{x_3}-x_6\ptl_{x_{11}}-x_8\ptl_{x_{13}}-x_{12}\ptl_{x_{15}},\eqno(11.3.43)$$
$$\pi(B_{3,5})=x_2\ptl_{x_5}-x_4\ptl_{x_9}-x_7\ptl_{x_{10}}+x_{14}\ptl_{x_{16}},\eqno(11.3.44)$$
$$\pi(A_{1,4})=-x_2\ptl_{x_7}+x_5\ptl_{x_{10}}+x_6\ptl_{x_{12}}+x_{11}\ptl_{x_{15}},\eqno(11.3.45)$$
$$\pi(A_{2,5})=x_1\ptl_{x_4}-x_5\ptl_{x_{11}}+x_8\ptl_{x_{14}}-x_{10}\ptl_{x_{15}},\eqno(11.3.46)$$
$$\pi(B_{2,5})=-x_2\ptl_{x_6}-x_3\ptl_{x_9}+x_7\ptl_{x_{12}}+x_{13}\ptl_{x_{16}},\eqno(11.3.47)$$
$$\pi(B_{3,4})=-x_1\ptl_{x_5}+x_4\ptl_{x_{11}}+x_7\ptl_{x_{13}}+x_{12}\ptl_{x_{16}},\eqno(11.3.48)$$
$$\pi(A_{1,5})=-x_1\ptl_{x_7}+x_5\ptl_{x_{13}}+x_6\ptl_{x_{14}}-x_9\ptl_{x_{15}},\eqno(11.3.49)$$
$$\pi(B_{1,5})=x_2\ptl_{x_8}+x_3\ptl_{x_{10}}+x_4\ptl_{x_{12}}+x_{11}\ptl_{x_{16}},\eqno(11.3.50)$$
$$\pi(B_{2,4})=x_1\ptl_{x_6}+x_3\ptl_{x_{11}}-x_7\ptl_{x_{14}}+x_{10}\ptl_{x_{16}},\eqno(11.3.51)$$
$$\pi(B_{1,4})=-x_1\ptl_{x_8}-x_4\ptl_{x_{14}}-x_3\ptl_{x_{13}}+x_9\ptl_{x_{16}},\eqno(11.3.52)$$
$$\pi(B_{2,3})=x_1\ptl_{x_9}-x_2\ptl_{x_{11}}+x_7\ptl_{x_{15}}-x_8\ptl_{x_{16}},\eqno(11.3.53)$$
$$\pi(B_{1,3})=-x_1\ptl_{x_{10}}+x_2\ptl_{x_{13}}+x_4\ptl_{x_{15}}-x_6\ptl_{x_{16}},\eqno(11.3.54)$$
$$\pi(B_{1,2})=x_1\ptl_{x_{12}}-x_2\ptl_{x_{14}}+x_3\ptl_{x_{15}}-x_5\ptl_{x_{16}},\eqno(11.3.55)$$
$$\pi(A_{j,i})=\tau(\pi(A_{i,j})),\;\;
\pi(C_{j,i})=\tau(\pi(B_{i,j}))\qquad\for\;\;1\leq i<j\leq
5,\eqno(11.3.56)$$
$$\pi(A_{r,r})=\sum_{i=1}^{16}(1/2+a_{r,i})x_i\ptl_{x_i},\qquad
r\in\ol{1,5},\eqno(11.3.57)$$ where $a_{r,i}$ are given in the
following table
\begin{center}{\bf \large Table 11.3.1}\end{center}
\begin{center}\begin{tabular}{|r||r|r|r|r|r|r|r|r|r|r|r|r|r|r|r|r|}\hline
$i$&1&2&3&4&5&6&7&8&9&10&11&12&13&14&15&16
\\\hline\hline
$a_{1,i}$&0&0&0&0&0&0&$-1$&$-1$&0&$-1$&0&$-1$&$-1$&$-1$&$-1$&$-1$
\\\hline
$a_{2,i}$&0&0&0&$-1$&0&$-1$&0&0&$-1$&0&$-1$&$-1$&0&$-1$&$-1$&$-1$
\\\hline$a_{3,i}$&0&0&$-1$&0&$-1$&0&0&0&$-1$&$-1$&$-1$&0&$-1$&0&$-1$&$-1$
\\\hline$a_{4,i}$&0&$-1$&0&0&$-1$&$-1$&0&$-1$&0&0&$-1$&0&$-1$&$-1$&0&$-1$
\\\hline$a_{5,i}$&$-1$&0&0&0&$-1$&$-1$&0&$-1$&$-1$&$-1$&0&$-1$&0&0&0&$-1$\\\hline
\end{tabular}\end{center}
Note that (11.3.36)-(11.3.55) are the representation formulas of all
the positive root vectors. In particular, $x_1$ is a highest-weight
vector of $V$ with weight $\lmd_4$, the forth fundamental weight of
$o(10,\mbb F)$, and $V$ gives a spin representation of $o(10,\mbb
C)$.
\subsection{Decomposition}

Recall that the representation of ${\msr G}^{D_5}$ on ${\msr A}$ is
given by (11.3.25) and the representation of $o(10,\mbb{C})$ is
given in (11.3.35). We calculate
$$\al_r|_{\msr
A}=\sum_{i=1}^{16}b_{r,i}x_i\ptl_{x_i}\qquad\for\;\;r\in\ol{1,5},\eqno(11.3.58)$$
where $b_{r,i}$ are given in the following table:
\begin{center}{\bf \large Table 11.3.2}\end{center}
\begin{center}\begin{tabular}{|r||r|r|r|r|r|r|r|r|r|r|r|r|r|r|r|r|}\hline
$i$&1&2&3&4&5&6&7&8&9&10&11&12&13&14&15&16
\\\hline\hline
$b_{1,i}$&0&0&0&1&0&1&$-1$&$-1$&1&$-1$&1&0&$-1$&0&0&0
\\\hline
$b_{2,i}$&0&0&1&1&$-1$&$-1$&1&$-1$&0&0&0&0&0&0&1&$-1$
\\\hline$b_{3,i}$&0&0&1&$-1$&1&$-1$&0&0&0&1&0&$-1$&1&$-1$&0&0
\\\hline$b_{4,i}$&0&1&$-1$&0&0&1&0&1&$-1$&$-1$&0&0&0&1&$-1$&0
\\\hline$b_{5,i}$&1&$-1$&0&0&0&0&0&0&1&1&$-1$&1&$-1$&$-1$&0&0\\\hline
\end{tabular}\end{center}
Recall that  a singular vector of $o(10,\mbb F)$ is a nonzero weight
vector annihilated by positive root vectors $\{A_{i,j},B_{i,j}\mid
1\leq i<j\leq 5\}$. Note that the weight of a singular vector in
${\msr B}$ must be dominate integral. Using the above table, we find
the singular vector
$$\zeta_1=x_1x_{11}+x_2x_9-x_3x_6+x_4x_5\eqno(11.3.59)$$
of weight $\lmd_1$, the first fundamental weight of $o(10,\mbb F)$.
Thus $\zeta_1$ generates the 10-dimensional natural $o(10,\mbb
F)$-module $U$. According to (11.3.36)-(11.3.40) and (11.3.56),
$$\pi(A_{2,1})=x_7\ptl_{x_4}+x_8\ptl_{x_6}+x_{10}\ptl_{x_9}+x_{13}\ptl_{x_{11}},\eqno(11.3.60)$$
$$\pi(A_{3,2})=x_4\ptl_{x_3}+x_6\ptl_{x_5}+x_{12}\ptl_{x_{10}}+x_{14}\ptl_{x_{13}},\eqno(11.3.61)$$
$$\pi(A_{4,3})=-x_3\ptl_{x_2}-x_9\ptl_{x_6}-x_{10}\ptl_{x_8}+x_{15}\ptl_{x_{14}},\eqno(11.3.62)$$
$$\pi(A_{5,4})=-x_2\ptl_{x_1}+x_{11}\ptl_{x_9}+x_{13}\ptl_{x_{10}}+x_{14}\ptl_{x_{12}},\eqno(11.3.63)$$
$$\pi(C_{5,4})=-x_5\ptl_{x_3}-x_6\ptl_{x_4}-x_8\ptl_{x_7}+x_{16}\ptl_{x_{15}}.\eqno(11.3.64)$$
We take
$$\zeta_2=A_{2,1}(\zeta_1)=x_1x_{13}+x_2x_{10}-x_3x_8+x_5x_7,\eqno(11.3.65)$$
$$\zeta_3=A_{3,2}(\zeta_2)=x_1x_{14}+x_2x_{12}-x_4x_8+x_6x_7,\eqno(11.3.66)$$
$$\zeta_4=A_{4,3}(\zeta_3)=x_1x_{15}-x_3x_{12}+x_4x_{10}-x_7x_9,\eqno(11.3.67)$$
$$\zeta_5=A_{5,4}(\zeta_4)=-x_2x_{15}-x_3x_{14}+x_4x_{13}-x_7x_{11},\eqno(11.3.68)$$
$$\zeta_{10}=C_{5,4}(\zeta_4)=x_1x_{16}+x_5x_{12}-x_6x_{10}+x_8x_9,\eqno(11.3.69)$$
$$\zeta_9=-C_{5,4}(\zeta_5)=x_2x_{16}-x_5x_{14}+x_6x_{13}-x_8x_{11},\eqno(11.3.70)$$
$$\zeta_8=-A_{4,3}(\zeta_9)=x_3x_{16}+x_5x_{15}+x_9x_{13}-x_{10}x_{11},\eqno(11.3.71)$$
$$\zeta_7=-A_{3,2}(\zeta_8)=-x_4x_{16}-x_6x_{15}-x_9x_{14}+x_{11}x_{12},\eqno(11.3.72)$$
$$\zeta_6=-A_{2,1}(\zeta_7)=x_7x_{16}+x_8x_{15}+x_{10}x_{14}-x_{12}x_{13}.\eqno(11.3.73)$$
Then $U=\sum_{i=1}^{10}\mbb F\zeta_i$ forms an $o(10,\mbb F)$-module
isomorphic to the 10-dimensional natural $o(10,\mbb F)$-module with
$\{\zeta_1,...,\zeta_{10}\}$ as the standard basis.\psp

{\bf Theorem 11.3.1}. {\it Any $o(10,\mbb F)$-singular vector in
$\msr B$ is a polynomial in $x_1$ and $\zeta_1$.}

{\it Proof}. Note
$$x_{11}=x_1^{-1}(\zeta_1-x_2x_9+x_3x_6-x_4x_5),\eqno(11.3.74)$$
$$x_{13}=x_1^{-1}(\zeta_2-x_2x_{10}+x_3x_8-x_5x_7),\eqno(11.3.75)$$
$$x_{14}=x_1^{-1}(\zeta_3-x_2x_{12}+x_4x_8-x_6x_7),\eqno(11.3.76)$$
$$x_{15}=x_1^{-1}(\zeta_4+x_3x_{12}-x_4x_{10}+x_7x_9),\eqno(11.3.77)$$
$$x_{16}=x_1^{-1}(\zeta_{10}-x_5x_{12}+x_6x_{10}-x_8x_9).\eqno(11.3.78)$$
Let $f$ be a singular vector in ${\msr B}$. Substituting
(11.3.74)-(11.3.78) into it, we can write
$$f=g(x_i,\zeta_1,\zeta_2,\zeta_3,\zeta_4,\zeta_{10}\mid 11,13,14,15,16\neq
i\in\ol{1,16}).\eqno(11.3.79)$$ By (11.3.39), (11.3.43), (11.3.46),
(11.3.48), (11.3.49),  (11.3.51)-(11.3.55), and the weights of
(11.3.59), (11.3.65)-(11.3.67) and (11.3.69),
$$A_{4,5}(f)=-x_1\ptl_{x_2}(g)=0\lra g_{x_2}=0,\eqno(11.3.80)$$
$$A_{3,5}(f)=-x_1\ptl_{x_3}(g)=0\lra g_{x_3}=0,\eqno(11.3.81)$$
$$A_{2,5}(f)=x_1\ptl_{x_4}(g)=0\lra g_{x_4}=0,\eqno(11.3.82)$$
$$B_{3,4}(f)=-x_1\ptl_{x_5}(g)=0\lra g_{x_5}=0,\eqno(11.3.83)$$
$$A_{1,5}(f)=-x_1\ptl_{x_7}(g)=0\lra g_{x_7}=0,\eqno(11.3.84)$$
$$B_{2,4}(f)=x_1\ptl_{x_6}(g)=0\lra g_{x_6}=0,\eqno(11.3.85)$$
$$B_{1,4}(f)=-x_1\ptl_{x_8}(g)=0\lra g_{x_8}=0,\eqno(11.3.86)$$
$$B_{2,3}(f)=x_1\ptl_{x_9}(g)=0\lra g_{x_9}=0,\eqno(11.3.87)$$
$$B_{1,3}(f)=-x_1\ptl_{x_{10}}(g)=0\lra g_{x_{10}}=0,\eqno(11.3.88)$$
$$B_{1,2}(f)=x_1\ptl_{x_{12}}(g)=0\lra g_{x_{12}}=0.\eqno(11.3.89)$$
Thus $f$ is a function in $x_1,\zeta_1,\zeta_2,\zeta_3,\zeta_4$ and
$\zeta_{10}$.

According to (11.3.36)-(11.3.38), (11.3.40), (11.3.59),
(11.3.65)-(11.3.67) and (11.3.69),
$$A_{1,2}(f)=\zeta_1\ptl_{\zeta_2}(g)=0\lra g_{\zeta_2}=0,\eqno(11.3.90)$$
$$A_{2,3}(f)=\zeta_2\ptl_{\zeta_3}(g)=0\lra g_{\zeta_3}=0,\eqno(11.3.91)$$
$$A_{3,4}(f)=\zeta_3\ptl_{\zeta_4}(g)=0\lra g_{\zeta_4}=0,\eqno(11.3.92)$$
$$B_{4,5}(f)=\zeta_4\ptl_{\zeta_{10}}(g)=0\lra g_{\zeta_{10}}=0.\eqno(11.3.93)$$
Hence $f$ is a function in $x_1$ and $\zeta_1$. Thanks to (11.3.79),
it must be a polynomial in $x_1$ and $\zeta_1.\qquad\Box$\psp

Let $V_{m_1,m_2}$ be the  $o(10,\mbb F)$-submodule generated by
$x_1^{m_1}\zeta^{m_2}$. Note that (11.3.7) is a Cartan root space
decomposition of $o(10,\mbb F)$. Since $\msr B$ is a weight
$o(10,\mbb F)$-module, Weyl's Theorem 2.3.6 of completely
reducibility yields that
$${\msr B}=\bigoplus_{m_1,m_2=0}^\infty V_{m_1,m_2}\eqno(11.3.94)$$
is a decomposition of irreducible $o(10,\mbb F)$-submodules. Denote
by $V(\lmd)$ the highest-weight irreducible $o(10,\mbb F)$-module
with the highest weight $\lmd$. The above equation leads to the
following combinatorial identity:
$$\sum_{m_1,m_2=0}^\infty (\dim
V(m_2\lmd_1+m_1\lmd_4))q^{m_1+2m_2}=\frac{1}{(1-q)^{16}},\eqno(11.3.95)$$
which was proved in (11.3.74)-(11.3.93)  by partial differential
equations.

\subsection{Symmetry and Equivalent Combinatorics}

Next we want to study the symmetry of $\ol{1,16}$ with respect to
the representation  in (11.3.36)-(11.3.57). Set
$$\mfk W=\{\pi(A_{i,j}),\pi(A_{j,i}),\pi(B_{i,j}),\pi(C_{j,i})\mid
1\leq i<j\leq 5\},\eqno(11.3.96)$$ which is the set of the
representations of root vectors in $o(10,\mbb F)$. Write
$$\mfk Z=\{\zeta_i\mid i\in\ol{1,10}\}\eqno(11.3.97)$$
(cf. (11.3.59) and (11.3.65)-(11.3.73)). First we have:\psp

{\bf Lemma 11.3.2}. {\it The set $\ol{1,16}$ is symmetric with
respect to the sets $\mfk W$ and $\mfk Z$.} \psp

{\it Proof}. Take $n=5$ in the spin representation of $o(2n,\mbb F)$
in (7.1.32)-(7.1.37) and Theorem 7.1.3. By (11.3.57) and Table
11.3.1, we have an $o(10,\mbb F)$-module isomorphism $\sgm: U\rta
\Psi_{(0)}$ such that $\sgm(x_1)=\sta_1\sta_2\sta_3\sta_4$.
According to (7.1.34), (7.1.35), (11.3.36)-(11.3.40) and (11.3.56),
the following sets are symmetric with respect to the representation
$\pi$ of $o(10,\mbb F)$:
$$S_1=\{x_1,x_2,x_3,x_4,x_7\},\;S_2=\{x_5,x_6,x_8,x_9,x_{10},x_{11},x_{12},x_{13},x_{14},x_{15}\}\eqno(11.3.98)$$
and $S_3=\{x_{16}\}$.

 We define an anti-automorphism
$\dg$ of the associative algebra
$$\mbb B=\sum_{i_1,...,i_{16}=0}^\infty\msr B\ptl_{x_1}^{i_1}\cdots
\ptl_{x_{16}}^{16}\eqno(11.3.99)$$ of differential operators by
$$(x_6)^\dg=-\ptl_{x_{10}},\;\;(\ptl_{x_6})^\dg=-x_{10},\;(x_{10})^\dg=-\ptl_{x_6},\;\;(\ptl_{x_{10}})^\dg=-x_6,
\eqno(11.3.100)$$
$$(x_i)^\dg=-\ptl_{x_{17-i}},\;(\ptl_{x_i})^\dg=-x_{17-i}\;\;\for\;\;i=4,13,
\eqno(11.3.101)$$
$$(x_7)^\dg=\ptl_{x_{11}},\;\;(\ptl_{x_7})^\dg=x_{11},\;(x_{11})^\dg=\ptl_{x_7},\;\;(\ptl_{x_{11}})^\dg=x_7,
\eqno(11.3.102)$$
$$(x_r)^\dg=\ptl_{x_{17-r}},\;(\ptl_{x_r})^\dg=x_{17-r}\;\;\for\;\;
r=1,2,3,5,8,9,12,14,15,16. \eqno(11.3.103)$$ It turns out that the
map $d\mapsto -d^\dg$ forms a Lie algebra automorphism of
$\pi(o(10,\mbb F))$ by (11.3.36)-(11.3.40), (11.3.56), (11.3.57) and
Table 11.3.1, which corresponds to an automorphism of $o(10,\mbb F)$
inducing the automorphism of its Dynkin diagram. Moreover,
$$-d^\dg\in\mfk W\qquad\mbox{for any}\;\;d\in\mfk
W.\eqno(11.3.104)$$

Define
$$\mfk s(d)=-\tau(d^\dg)\qquad\for\;\;d\in\mbb B.\eqno(11.3.105)$$
Then $\mfk s$ is an associative algebra automorphism of $\mbb B$
with order 2. In fact,
$$\mfk s(x_6)=-x_{10},\;\;\mfk s(x_4)=-x_{13},\;\;\mfk s(x_7)=x_{11},
\;\;\mfk s(x_r)^\dg=x_{17-r} \eqno(11.3.106)$$ for $r=1,2,3,5,8.$
Moreover,
$$\mfk s(\zeta_i)=\zeta_{5+i},\;\mfk
s(\zeta_r)=\zeta_r\qquad\for\;\;i\in\ol{1,4},\;r=5,10.\eqno(11.3.107)$$
This shows that $\{1,16\}$ and $\{2,15\}$ are symmetric with respect
to the sets $\mfk W$ and $\mfk Z$. Thus the lemma follows from the
conclusion of the first paragraph. $\qquad\Box$\psp

For any $i\in\ol{1,16}$, we define
$$\Upsilon_i=\{(r,s)\mid
\zeta_r\;\mbox{contains}\;x_ix_s\}.\eqno(11.3.108)$$ We have the
following equivalent combinatorial property:\psp

{\bf Lemma 11.3.3}. {\it Each polynomial $\zeta_r$ contains exactly
eight $x_i$'s. Moreover,
$$|\Upsilon_i|=5\qquad\for\;\;i\in\ol{16}.\eqno(11.3.109)$$}

{\it Proof}. According to (11.3.59) and  (11.3.65)-(11.3.73), we
calculate
$$\Upsilon_1=\{(1,11), (2,13), (3,14), (4,15), (10,16)\},\eqno(11.3.110)$$
Then (11.3.108) follows from Lemma 11.3.2. The first statement is
obtained by checking (11.3.59) and  (11.3.65)-(11.3.73)
one-by-one.$\qquad\Box$\psp

If we represent $\{x_i\mid i\in\ol{1,16}\}$ by 16 vertices and
represent $\{\zeta_r\mid r\in\ol{1,10}\}$ by 10 lines, then we
obtain a graph of $16$ vertices and 10 lines such that each line
contains 8 vertices and each vertex is exactly on 5 lines.

 For
$i\in\ol{1,16}$, we denote
$$I_i=\{r\in\ol{1,16}\mid\mbox{some}\;d\in\mfk W\;\mbox{contains}\;x_r\ptl_{x_i}\}\eqno(11.3.111)$$ and
$$\mfk W_i=\{d\in\mfk W\mid d\;\mbox{does not
contain}\;x_i\;\mbox{and}\;\ptl_{x_i}\}.\eqno(11.3.112)$$ Write
$$J_i=\{s\in\ol{1,16}\mid (r,s)\in\Upsilon_i\;\mbox{for
some}\;r\in\ol{1,16}\}.\eqno(11.3.113)$$ The following equivalent
combinatorial properties will be crucial to our main result in
Section 11.7.\psp

{\bf Lemma 11.3.4}. {\it For any $i\in\ol{1,16}$, we have
$$|I_i|=10,\;\;|\mfk W_i|=20,\;\;|J_i|=5.\eqno(11.3.114)$$
In fact,
$$I_i\bigcup J_i=\ol{1,16}\setminus\{i\}.\eqno(11.3.115)$$
Moreover, every element in $\mfk W_i$ contains exactly one $x_s$
with $s\in J_i$, and for any $r\in J_i$, $x_r$ is contained in
exactly four elements in $\mfk W_i$.}

{\it Proof}. By Lemma 11.3.2, we only need to prove it for $i=1$.
Note that the elements in $\mfk W$ containing $\ptl_{x_1}$ are
$$\pi(A_{5,4}),\;\;\pi(A_{5,3}),\;\;\pi(A_{5,2}),\;\;\pi(A_{5,1}),\;\;\pi(C_{4,3}),\eqno(11.3.116)$$
$$\pi(C_{4,2}),\;\;\pi(C_{4,1}),\;\;
\pi(C_{3,2}),\;\;\pi(C_{3,1}),\;\;\pi(C_{2,1})\eqno(11.3.117)$$ by
(11.3.36)-(11.3.56). Thus $I_1=\ol{2,10}\bigcup\{12\}$. On the other
hand,
$$J_1=\{11,13,14,15,16\}\eqno(11.3.118)$$ by (11.3.110).
Since $\tau(\mfk W)=\mfk W$ (cf. (11.3.29)), $|\mfk W_1|=20$. So
(11.3.114) and (11.3.115) holds for $i=1$.

The elements in $\mfk W_1$ containing $x_{11}$ are
$$\pi(A_{1,2}),\pi(A_{1,3}),\pi(A_{1,4}),\pi(B_{1,5}).\eqno(11.3.119)$$
 The elements in $\mfk W_1$
containing $x_{13}$ are
$$\pi(A_{2,1}),\pi(A_{2,3}),\pi(A_{2,4}),\pi(B_{2,5}).\eqno(11.3.120)$$
The elements in $\mfk W_1$ containing $x_{14}$ are
$$\pi(A_{3,1}),\pi(A_{3,2}),\pi(A_{3,4}),\pi(B_{3,5}).\eqno(11.3.121)$$
 The elements in $\mfk W_1$
containing $x_{15}$ are
$$\pi(A_{4,1}),\pi(A_{4,2}),\pi(A_{4,3}),\pi(B_{4,5}).\eqno(11.3.122)$$
The elements in $\mfk W_1$ containing $x_{16}$ are
$$\pi(B_{5,1}),\pi(B_{5,2}),\pi(B_{5,3}),\pi(B_{5,4}).\eqno(11.3.123)$$
Therefore every element in $\mfk W_1$ contains exactly one $x_s$
with $s\in J_1$, and for any $r\in J_1$, $x_r$ is contained in
exactly four elements in $\mfk W_1.\qquad\Box$

\section{Realization of $E_6$ in 16-Dimensional Space}

 In this section, we
want to find a differential-operator representation of ${\msr
G}^{E_6}$, or equivalently, a fractional representation on
16-dimensional space of the Lie group of type $E_6$.

Recall the Lie subalgebra $\msr G_0=\msr G^{D_5}\oplus\mbb F\al_6$
of $\msr G^{E_6}$. We can view the linear map $\nu$ determined in
(11.3.30)-(11.3.34) as an identification of $o(10,\mbb F)$ with
$\msr G^{D_5}$. We want to find the extension of the representation
$\pi$ of $o(10,\mbb F)$ given in (11.3.36)-(11.3.57) to $\msr
G^{E_6}$. According to (11.3.24) and (11.3.25), we calculate
$$\pi(\al_6)=-2x_1\ptl_{x_1}-\sum_{i=2}^{10}x_i\ptl_{x_i}-x_{12}\ptl_{x_{12}}.\eqno(11.4.1)$$
 Write
$$\widehat\al=2\al_1+4\al_3+6\al_4+3\al_2+5\al_5+4\al_6.\eqno(11.4.2)$$
Then
$$(\widehat\al,\al_r)=0\qquad\for\;\;r\in\ol{1,5}\eqno(11.4.3)$$
by the Dynkin diagram of $E_6$, where $(\cdot,\cdot)$ is a symmetric
$\mbb Z$-bilinear form on the root lattice of $E_6$ such that
$(\al,\al)=2$ for $\al\in \Phi_{E_6}$.

Thanks to (4.4.24),
$$[\wht\al,{\msr G}^{D_5}]=0.\eqno(11.4.4)$$
So $\wht\al$ is a central element of $\msr G_0$. By Schur's Lemma,
$\widehat\al|_V=c\sum_{i=1}^{16}x_i\ptl_{x_i}$. According to the
coefficients of $x_1\ptl_{x_1}$ in (11.3.58) with the data in Table
11.3.2 and (11.4.2), we have
$$\pi(\widehat\al)=-3D,\;\;\mbox{where}\;\;D=\sum_{i=1}^{16}x_i\ptl_{x_i}\eqno(11.4.5)$$
is the degree operator on ${\msr B}$.

Recall that the Lie bracket in the algebra $\mbb B$ (cf. (11.3.28))
is given by the commutator
$$[d_1,d_2]=d_1d_2-d_2d_1.\eqno(11.4.6)$$
Set
$$\msr T=\sum_{i=1}^{16}\mbb F\ptl_{x_i}.\eqno(11.4.7)$$
Then $\msr T$ forms an $o(10,\mbb F)$-module with respect to the
action
$$B(\ptl)=[\pi(B),\ptl]\qquad\for\;\;B\in
o(10,\mbb F),\;\ptl\in{\msr T}.\eqno(11.4.8)$$ On the other hand,
${\msr G}_\pm$ (cf. (11.3.14)-(11.3.20)) form $o(10,\mbb F)$-modules
with respect to the action
$$B(u)=[\nu(B),u]\qquad \for\;\;B\in
o(10,\mbb F),\;u\in{\msr G}_\pm\eqno(11.4.9)$$(cf.
(11.3.30)-(11.3.34)). According to (11.3.24) and (11.3.25), the
linear map determined by $\eta_i\mapsto x_i$ for $i\in\ol{1,16}$
gives an $o(10,\mbb F)$-module isomorphism from ${\msr G}_-$ to $V$.
Moreover, (11.3.3) and (11.3.22) imply that the linear map
determined by $\xi_i\mapsto \ptl_{x_i}$ for $i\in\ol{1,16}$ gives an
$o(10,\mbb F)$-module isomorphism from ${\msr G}_+$ to ${\msr T}$.
Hence we define the representation $\pi$ of ${\msr G}_+$ on $\msr B$
by
$$\pi(\xi_i)=\ptl_{x_i}\qquad\for\;\;i\in\ol{1,16}.\eqno(11.4.10)$$

Recall the Witt Lie subalgebra of $\mbb B$:
$${\mbb W}_{16}=\sum_{i=1}^{16}{\msr B}\ptl_{x_i}.\eqno(11.4.11)$$
Now we want to find the differential operators
$P_1,P_2,...,P_{16}\in {\mbb W}_{16}$ such that the following
representation $\pi$ matches the structure of ${\msr G}^{E_6}$:
$$\pi(\eta_i)=P_i\qquad\for\;\;i\in\ol{1,16}.\eqno(11.4.12)$$
Imposing
$$[\ptl_{x_1},P_1]=\pi([E_{\al_6},E_{-\al_6}])=-\pi(\al_6)
=2x_1\ptl_{x_1}+\sum_{i=2}^{10}x_i\ptl_{x_i}+x_{12}\ptl_{x_{12}},\eqno(11.4.13)$$
we take
$$P_1=x_1(\sum_{i=1}^{10}x_i\ptl_{x_i}+x_{12}\ptl_{x_{12}})+P_1',\eqno(11.4.14)$$
where $P_1'$ is a differential operator such that
$[\ptl_{x_1},P_1']=0$. Moreover,
$$[\ptl_{x_r},x_1(\sum_{i=1}^{10}x_i\ptl_{x_i}+x_{12}\ptl_{x_{12}})]=x_1\ptl_{x_r}\qquad\for\;\;r\in\{\ol{2,10},12\}.
\eqno(11.4.15)$$ Wanting $[\ptl_{x_r},P_1]\in \pi(o(10,\mbb F))$
(cf. (11.3.36)-(11.3.57)), we take
\begin{eqnarray*} P_1&=&
x_1(\sum_{i=1}^{10}x_i\ptl_{x_i}+x_{12}\ptl_{x_{12}})-(x_2x_9-x_3x_6+x_4x_5)\ptl_{x_{11}}
\\ &&-(x_2x_{10}-x_3x_8+x_5x_7)\ptl_{x_{13}}-(x_2x_{12}-x_4x_8+x_6x_7)\ptl_{x_{14}}
\\ &&+(x_3x_{12}-x_4x_{10}+x_7x_9)\ptl_{x_{15}}-(x_5x_{12}-x_6x_{10}+x_8x_9)\ptl_{x_{16}}
\\&=&
x_1D-\zeta_1\ptl_{x_{11}}-\zeta_2\ptl_{x_{13}}-\zeta_3\ptl_{x_{14}}-\zeta_4\ptl_{x_{15}}-\zeta_{10}\ptl_{x_{16}}
\hspace{4.3cm}(11.4.16)\end{eqnarray*} by (11.3.39), (11.3.43),
(11.3.46), (11.3.48), (11.3.49), (11.3.51)-(11.3.55), (11.3.59),
(11.3.65)-(11.3.67), (11.3.69) and (11.4.5). Then
$$[\ptl_{x_s},P_1]=\pi(\nu^{-1}([\xi_s,\eta_1]))\qquad\for\;\;s\in\ol{1,16}\eqno(11.4.17)$$ due to
(11.3.30)-(11.3.35).

Since $[E_{-\al_5},\eta_1]=\eta_2$ by (4.4.25), we take
\begin{eqnarray*}\qquad\qquad P_2&=&[\pi(\nu^{-1}(E_{-\al_5})),\pi(\eta_1)]=-[\pi(A_{5,4}),P_1]\\&=&
x_2D-\zeta_1\ptl_{x_9}-\zeta_2\ptl_{x_{10}}-\zeta_3\ptl_{x_{12}}+\zeta_5\ptl_{x_{15}}-\zeta_9\ptl_{x_{16}}
\hspace{2.9cm}(11.4.18)\end{eqnarray*} by (11.3.31) and (11.3.63).
Similarly, we take
$$
P_3=x_3D+\zeta_1\ptl_{x_6}+\zeta_2\ptl_{x_8}+\zeta_4\ptl_{x_{12}}+\zeta_5\ptl_{x_{14}}-\zeta_8\ptl_{x_{16}},\eqno
(11.4.19)$$
$$P_4=x_4D-\zeta_1\ptl_{x_5}+\zeta_3\ptl_{x_8}-\zeta_4\ptl_{x_{10}}-\zeta_5\ptl_{x_{13}}+\zeta_7\ptl_{x_{16}},
\eqno(11.4.20)$$
$$P_5=
x_5D-\zeta_1\ptl_{x_4}-\zeta_2\ptl_{x_7}-\zeta_{10}\ptl_{x_{12}}+\zeta_9\ptl_{x_{14}}-\zeta_8\ptl_{x_{15}},\eqno(11.4.21)$$
$$P_6=x_6D+\zeta_1\ptl_{x_3}-\zeta_3\ptl_{x_7}+\zeta_{10}\ptl_{x_{10}}-\zeta_9\ptl_{x_{13}}+\zeta_7\ptl_{x_{15}},
\eqno(11.4.22)$$
$$P_7=x_7D-\zeta_2\ptl_{x_5}-\zeta_3\ptl_{x_6}+\zeta_4\ptl_{x_9}+\zeta_5\ptl_{x_{11}}-\zeta_6\ptl_{x_{16}},
\eqno(11.4.23)$$ $$
P_8=x_8D+\zeta_2\ptl_{x_3}+\zeta_3\ptl_{x_4}-\zeta_{10}\ptl_{x_9}+\zeta_9\ptl_{x_{11}}-\zeta_6\ptl_{x_{15}},
\eqno(11.4.24)$$
$$P_9=x_9D-\zeta_1\ptl_{x_2}+\zeta_4\ptl_{x_7}-\zeta_{10}\ptl_{x_8}-\zeta_8\ptl_{x_{13}}+\zeta_7\ptl_{x_{14}},
\eqno(11.4.25)$$
$$P_{10}=x_{10}D-\zeta_2\ptl_{x_2}-\zeta_4\ptl_{x_4}+\zeta_{10}\ptl_{x_6}+\zeta_8\ptl_{x_{11}}-\zeta_6\ptl_{x_{14}},
\eqno(11.4.26)$$
$$
P_{11}=x_{11}D-\zeta_1\ptl_{x_1}+\zeta_5\ptl_{x_7}+\zeta_9\ptl_{x_8}+\zeta_8\ptl_{x_{10}}-\zeta_7\ptl_{x_{12}},
\eqno(11.4.27)$$
$$P_{12}=
x_{12}D-\zeta_3\ptl_{x_2}+\zeta_4\ptl_{x_3}-\zeta_{10}\ptl_{x_5}-\zeta_7\ptl_{x_{11}}+\zeta_6\ptl_{x_{13}},
\eqno(11.4.28)$$
$$P_{13}=
x_{13}D-\zeta_2\ptl_{x_1}-\zeta_5\ptl_{x_4}-\zeta_9\ptl_{x_6}-\zeta_8\ptl_{x_9}+\zeta_6\ptl_{x_{12}},
\eqno(11.4.29)$$ $$P_{14}=
x_{14}D-\zeta_3\ptl_{x_1}+\zeta_5\ptl_{x_3}+\zeta_9\ptl_{x_5}+\zeta_7\ptl_{x_9}-\zeta_6\ptl_{x_{10}},
\eqno(11.4.30)$$
$$P_{15}=x_{15}D-\zeta_4\ptl_{x_1}+\zeta_5\ptl_{x_2}-\zeta_8\ptl_{x_5}+\zeta_7\ptl_{x_6}-\zeta_6\ptl_{x_8},\eqno(11.4.31)$$
$$P_{16}=x_{16}D-\zeta_{10}\ptl_{x_1}-\zeta_9\ptl_{x_2}-\zeta_8\ptl_{x_3}+\zeta_7\ptl_{x_4}-\zeta_6\ptl_{x_7}.
\eqno(11.4.32)$$

Set
$${\msr P}=\sum_{i=1}^{16}\mbb FP_i,\qquad\mfk{C}_0=\pi(o(10,\mbb F))+\mbb FD\eqno(11.4.33)$$ (cf. (11.3.36)-(11.3.57)
 and (11.4.5)) and
$$\mfk{C}={\msr P}+\mfk{C}_0+{\msr T}\eqno(11.4.34)$$ (cf. (11.4.7)).
Then we have:\psp

{\bf Theorem 11.4.1}. {\it The space $\mfk{C}$ forms a Lie
subalgebra of the Witt algebra ${\mbb W}_{16}$ (cf. (11.4.11)).
Moreover, the linear map $\vt$ determined by
$$\vt(\xi_i)=\ptl_{x_i},\;\;\vt(\eta_i)=P_i,\;\;\vt(u)=\pi(\nu^{-1}(u))\qquad\for\;\;i\in\ol{1,16},\;u\in{\msr G}^{D_5}\eqno(11.4.35)$$
(cf. (11.3.30)-(11.3.34)) and
$$\vt(\al_6)=-2x_1\ptl_{x_1}-\sum_{i=2}^{10}x_i\ptl_{x_i}-x_{12}\ptl_{x_{12}}\eqno(11.4.36)$$
(cf. (11.4.1)) gives a Lie algebra isomorphism from ${\msr G}^{E_6}$
to $\mfk{C}$.}

{\it Proof}. Since $\msr T\cong {\msr G}_+$ as $\msr
G^{D_5}$-modules, we have
$${\msr G}_0+{\msr G}_+\stl{\vt}{\cong}\mfk{C}_0+\msr T\eqno(11.4.37)$$ as Lie algebras. Recall
that $U({\msr G})$ stands for the universal enveloping algebra of a
Lie algebra ${\msr G}$. Note that
$${\msr B}_-={\msr G}_0+{\msr G}_-,\qquad{\msr B}_+={\msr G}_0+{\msr G}_+\eqno(11.4.38)$$
are also  Lie subalgebras of ${\msr G}^{E_6}$ and
$${\msr G}^{E_6}={\msr B}_-\oplus {\msr G}_+={\msr G}_-\oplus {\msr B}_+.\eqno(11.4.39)$$
We define a one-dimensional ${\msr B}_-$-module $\mbb Fu_0$ by
$$w(u_0)=0\qquad\for\;\;w\in{\msr G}_-+\msr G^{D_5},\;\;\widehat\al(u_0)=48u_0\eqno(11.4.40)$$
(cf. (11.4.2)). Let
$$M=U({\msr G}^{E_6})\otimes_{{\msr B}_-}\mbb Fu_0\cong U({\msr G}_+)\otimes_{\mbb F} \mbb Fu_0\eqno(11.4.41)$$
be the induced ${\msr G}^{E_6}$-module.

Recall that $\mbb{N}$ is the set of nonnegative integers. Let
$$\wht B=\mbb F[\ptl_{x_1},\ptl_{x_2},...,\ptl_{x_{16}}].\eqno(11.4.42)$$
We define an action of the associative algebra $\mbb B$ (cf.
(11.3.28)) on $\wht B$ by
$$\ptl_{x_i}(\prod_{j=1}^{16}\ptl_{x_j}^{\be_j})=\ptl_{x_i}^{\be_i+1}\prod_{i\neq
j\in\ol{1,16}}\prod_{j=1}^{16}\ptl_{x_j}^{\be_j}\eqno(11.4.43)$$ and
$$x_i(\prod_{j=1}^{16}\ptl_{x_j}^{\be_j})=-\be_i\ptl_{x_i}^{\be_i-1}\prod_{i\neq
j\in\ol{1,16}}\prod_{j=1}^{16}\ptl_{x_j}^{\be_j}\eqno(11.4.44)$$ for
$i\in\ol{1,16}$. Since
$$[-x_i,\ptl_{x_j}]=[\ptl_{x_i},x_j]=\dlt_{i,j}\qquad\for\;\;i,j\in\ol{1,16},\eqno(11.4.45)$$
the above action gives an associative algebra representation of
$\mbb B$. Thus it also gives a Lie algebra representation of $\mbb
B$ (cf. (11.4.6)).
 It is straightforward to verify that
$$[d|_{\wht B},\ptl|_{\wht B}]=[d,\ptl]|_{\wht B}\qquad\for\;\;d\in\mfk{C}_0,\;\ptl\in{\msr T}.\eqno(11.4.46)$$

Define linear map $\mfk{F}: M\rta\wht B$ by
$$\mfk{F}(\prod_{i=1}^{16}\xi_i^{\ell_i}\otimes
u_0)=\prod_{i=1}^{16}\ptl_{x_i}^{\ell_i}\qquad(\ell_1,...,\ell_{16})\in\mbb{N}^{16}.\eqno(11.4.47)$$
According to (11.3.36)-(11.3.57), (11.4.43) and (11.4.44),
$$D(1)=-16,\;\;d(1)=0\qquad\for\;\;d\in \pi(o(10,\mbb F)).\eqno(11.4.48)$$
Moreover, (11.4.40), (11.4.41),  (11.4.43), (11.4.44) and (11.4.48)
imply
$$\mfk{F}(\xi(v))=\vt(\xi)[\mfk{F}(v)]\qquad\for\;\;\xi\in{\msr
G}_0,\;v\in M.\eqno(11.4.49)$$

Now (11.4.41) and (11.4.43) yield
$$\mfk{F}(w(u))=\vt(w)[\mfk{F}(u)]\qquad\for\;\;w\in{\msr
B}_+,\;\;u\in M.\eqno(11.4.50)$$ Thus  we have
$$\vs w|_{\Psi}\vs^{-1}=\vt(w)|_{{\msr A}'}\qquad\for\;\;w\in{\msr
B}_+.\eqno(11.4.51)$$ On the other hand, the linear map
$$\psi(v)=\mfk{F}(v|_{M})\mfk{F}^{-1}\qquad\for\;\; v\in{\msr
G}^{E_6}\eqno(11.4.52)$$ is a Lie algebra monomorphism from ${\msr
G}^{E_6}$ to $\mbb B|_{\wht B}$. According to (11.4.17) and
(11.4.45),
$$\psi(\eta_1)=P_1|_{\wht B}.\eqno(11.4.53)$$
By the constructions of $P_2,...,P_{16}$ in (11.4.18)-(11.4.32), we
have
$$\psi(\eta_i)=P_i|_{\wht B}\qquad\for\;\;i\in\ol{2,16}.\eqno(11.4.54)$$ Therefore, we have
$$\psi(v)=\vt(v)|_{\wht B}\qquad\for\;\;v\in{\msr
G}^{E_6}.\eqno(11.4.55)$$ In particular, $\mfk{C}|_{\wht
B}=\vt({\msr G}^{E_6})|_{\wht B}=\psi({\msr G}^{E_6})$ forms a Lie
algebra. Since the linear map $d\mapsto d|_{\wht B}$ for $d\in
\mfk{C}$ is injective, we have that $\mfk{C}$ forms a Lie subalgebra
of $\mbb B$ and $\vt$ is a Lie algebra isomorphism.$\qquad\Box$ \psp

By the above theorem, a Lie group of type $E_6$ is generated by the
linear transformations $\{e^{b\pi(u)}\mid b\in\mbb{R},\; u\in
o(10,\mbb{R})\}$ associated with (11.3.36)-(11.3.57), the real
translations and  dilations in $\mbb R^{16}$ with $x_r$ as the $r$th
coordinate function, and the fractional transformations
$\{e^{bP_s}:(x_1,x_2,...,x_{16})\mapsto
(e^{bP_s}(x_1),e^{bP_s}(x_2),...,e^{bP_s}(x_{16}))\mid
b\in\mbb{R},\; s\in\ol{1,16}\}$, where $e^{bP_s}(x_i)$ are of the
 forms as the following forms of the case $s=1$:
$$e^{bP_1}(x_i)=\frac{x_i}{1-bx_1},\;\;i\in\{\ol{1,10},12\},\eqno(11.4.56)$$
$$e^{bP_1}(x_{11})
=x_{11}-\frac{b(x_2x_9-x_3x_6+x_4x_5)}{1-bx_1},\eqno(11.4.57)$$
$$e^{bP_1}(x_{13})=x_{13}-\frac{b(x_2x_{10}-x_3x_8+x_5x_7)}{1-bx_1},\eqno(11.4.58)$$
$$e^{bP_1}(x_{14})=x_{14}-\frac{
b(x_2x_{12}-x_4x_8+x_6x_7)}{1-bx_1},\eqno(11.4.59)$$
$$e^{bP_1}(x_{15})=x_{15}+\frac{b(x_3x_{12}-x_4x_{10}+x_7x_9)}{1-bx_1},\eqno(11.4.60)$$
$$e^{bP_1}(x_{16})=x_{16}-\frac{b(x_5x_{12}-x_6x_{10}+x_8x_9)}{1-bx_1}.\eqno(11.4.61)$$

\section{Functor from $D_5$-Mod to $E_6$-Mod}

 In this section, we construct a functor from  the category of
$D_5$-modules to the category of $E_6$-modules.

Note that
$$o(10,{\msr B})=\sum_{i,j=1}^n({\msr B}A_{i,j}+{\msr B}B_{i,j}+{\msr
A}C_{i,j})\eqno(11.5.1)$$ (cf. (11.3.4)) forms a Lie subalgebra of
the matrix algebra $gl(10,{\msr B})$ over ${\msr B}$ with respect to
the commutator, i.e.
$$[fB_1,gB_2]=fg[B_1,B_2]\qquad\for\;\;f,g\in{\msr B},\;B_1,B_2\in
gl(10,\mbb F).\eqno(11.5.2)$$ Moreover, we define the Lie algebra
$${\msr K}=o(10,{\msr B})\oplus {\msr B}\kappa\eqno(11.5.3)$$
with the Lie bracket:
$$[\xi_1+f\kappa,\xi_2+g\kappa]=[\xi_1,\xi_2]\qquad\for\;\;\xi_1,\xi_2\in
o(10,{\msr B}),\;f,g\in{\msr B}.\eqno(11.5.4)$$ Similarly,
$gl(16,{\msr B})$ becomes a Lie algebra with the Lie bracket as that
in (11.5.2). Recall the Witt algebra ${\mbb
W}_{16}=\sum_{i=1}^{16}{\msr B}\ptl_{x_i}$, and Lemma 6.7.1 says
that there exists a monomorphism $\Im$ from the Lie algebra ${\mbb
W}_{16}$ to the Lie algebra of semi-product ${\mbb
W}_{16}+gl(16,{\msr B})$ defined by
$$\Im(\sum_{i=1}^{16}f_i\ptl_{x_i})=\sum_{i=1}^{16}f_i\ptl_{x_i}+\Im_1(\sum_{i=1}^{16}f_i\ptl_{x_i}),\;\;
\Im_1(\sum_{i=1}^{16}f_i\ptl_{x_i})=\sum_{i,j=1}^{16}\ptl_{x_i}(f_j)E_{i,j}.
\eqno(11.5.5)$$ According to our construction of $P_1$-$P_{16}$ in
(11.4.12)-(11.4.32),
$$\Im_1(P_i)=\sum_{r=1}^{16}x_r\Im_1(\pi([\xi_r,\eta_i]))\qquad\for\;\;i\in\ol{1,16}\eqno(11.5.6)$$
(cf. (11..3.36)-(11.3.57) and (11.4.4)). On the other hand,
$$\widehat{\msr K}={\mbb W}_{16}\oplus {\msr K}\eqno(11.5.7)$$
becomes a Lie algebra with the Lie bracket
\begin{eqnarray*}&&[d_1+f_1B_1+f_2\kappa,d_2+g_1B_2+g_2\kappa]\\
&=&[d_1,d_2]+f_1g_1[B_1,B_2]+d_1(g_2)B_2
-d_2(g_1)B_1+(d_1(g_2)-d_2(g_1))\kappa\hspace{1.9cm}(11.5.8)\end{eqnarray*}
for $f_1,f_2,g_1,g_2\in{\msr B},\;B_1,B_2\in o(10,\mbb F)$ and
$d_1,d_2\in{\mbb W}_{16}$. Note
$${\msr G}_0={\msr
G}^{D_5}\oplus\mbb F\widehat\al\eqno(11.5.9)$$ by (11.3.20) and
(11.4.2). So there exists a Lie algebra monomorphism $\mfk{k}:{\msr
G}_0\rta {\msr K}$ determined by
$$\mfk{k}(\widehat\al)=2\kappa,\;\;\mfk{k}(u)=\nu^{-1}(u)\qquad\for\;\;u\in{\msr
G}^{D_5}\eqno(11.5.10)$$ (cf. (11.3.30)-(11.3.34)). Since $\Im$ is a
Lie algebra monomorphism, our construction of $P_1$-$P_{16}$ in
(11.4.12)-(11.4.32) and (11.5.6) show that we have a Lie algebra
monomorphism $\iota: {\msr G}^{E_6}\rta \widehat{\msr K}$ given by
$$\iota(u)=\pi(u)+\mfk{k}(u)\qquad\for\;\;u\in{\msr
G}_0,\eqno(11.5.11)$$
$$\iota(\xi_i)=\ptl_{x_i},\;\;\iota(\eta_i)=P_i+\sum_{r=1}^{16}x_r\mfk{k}([\xi_r,\eta_i])\qquad\for
\;\;i\in\ol{1,16}.\eqno(11.5.12)$$

In order to calculate $\{\iota(\eta_1),...,\iota(\eta_{16})\}$
explicitly, we need the more formulas of $\nu$ on the positive root
vectors of $o(10,\mbb F)$ extended from (11.3.30)-(11.3.34). We
calculate
$$\nu(A_{3,5})=E_{(0,0,0,1,1)},\;\;\nu(A_{2,5})=-E_{(0,0,1,1,1)},\;\;\nu(B_{1,4})=E_{(1,1,1,1,1)},\eqno(11.5.13)$$
$$\nu(B_{3,4})=E_{(0,1,0,1,1)},\;\;\nu(B_{2,4}=-E_{(0,1,1,1,1)},\;\;\nu(A_{1,5})=E_{(1,0,1,1,1)},
\eqno(11.5.14)$$
$$\nu(B_{2,3})=-E_{(0,1,1,2,1)},\;\;\nu(B_{1,3})=E_{(1,1,1,2,1)},\;\;\nu(B_{1,2})
=-E_{(1,1,2,2,1)}.\eqno(11.5.15)$$ Using these formulas, we find
\begin{eqnarray*} \iota(\eta_1)&=&
P_1+x_1(\sum_{i=1}^4A_{i,i}-A_{5,5}-\kappa)/2
-x_2A_{4,5}-x_3A_{3,5}+x_4A_{2,5}-x_5B_{3,4}\\ & &+x_6B_{2,4}
-x_7A_{1,5}-x_8B_{1,4}+x_9B_{2,3}
-x_{10}B_{1,3}+x_{12}B_{1,2},\hspace{2.5cm}(11.5.16)\end{eqnarray*}
\begin{eqnarray*}
\iota(\eta_2)&=& P_2+x_2(\sum_{i\neq
4}(A_{i,i}-A_{4,4}-\kappa)/2-x_1A_{5,4} -x_3A_{3,4} +x_4A_{2,4}
+x_5B_{3,5}\\ & &-x_6B_{2,5}-x_7A_{1,4}+x_8B_{1,5}-x_{11}B_{2,3}
+x_{13}B_{1,3}-x_{14}B_{1,2},\hspace{2.4cm}(11.5.17)\end{eqnarray*}
\begin{eqnarray*}
\iota(\eta_3)&=& P_3+x_3(\sum_{i\neq
3}A_{i,i}-A_{3,3}-\kappa)/2-x_1A_{5,3}-x_2A_{4,3}
 +x_4A_{2,3}-x_5B_{4,5}\\ & &-x_7A_{1,3}-x_9B_{2,5}
+x_{10}B_{1,5}+x_{11}B_{2,4}
-x_{13}B_{1,4}+x_{15}B_{1,2},\hspace{2.35cm}(11.5.18)\end{eqnarray*}
\begin{eqnarray*}
\iota(\eta_4)&=& P_4+x_4(\sum_{i\neq
2}A_{i,i}-A_{2,2}-\kappa)/2+x_1A_{5,2}+x_2A_{4,2}
 +x_3A_{3,2} -x_6B_{4,5}\\ & &+x_7A_{1,2}-x_9B_{3,5}
+x_{11}B_{3,4}+x_{12}B_{1,5}
-x_{14}B_{1,4}+x_{15}B_{1,3},\hspace{2.3cm}(11.5.19)\end{eqnarray*}
\begin{eqnarray*} \iota(\eta_5)
 &=& P_5+x_5(A_{1,1}+A_{2,2}-\sum_{i=
3}^5A_{i,i}-\kappa)/2-x_1C_{4,3} +x_2C_{5,3}-x_3C_{5,4}
 +x_6A_{2,3}\\& &-x_8A_{1,3}-x_9A_{2,4}
+x_{10}A_{1,4}-x_{11}A_{2,5}
+x_{13}A_{1,5}-x_{16}B_{1,2},\hspace{2.3cm}(11.5.20)\end{eqnarray*}
\begin{eqnarray*} \iota(\eta_6)&=&
P_6+x_6(A_{1,1}+A_{3,3}-\sum_{i=2,4,5 }A_{i,i}-\kappa)/2
+x_1C_{4,2}-x_2C_{5,2}
 -x_4C_{5,4}+x_5A_{3,2}\\ & &+x_8A_{1,2}-x_9A_{3,4}
-x_{11}A_{3,5}+x_{12}A_{1,4}
+x_{14}A_{1,5}-x_{16}B_{1,3},\hspace{2.3cm}(11.5.21)\end{eqnarray*}
\begin{eqnarray*}
\iota(\eta_7)&=& P_7+x_7(\sum_{i
2}^5A_{i,i}-A_{1,1}-\kappa)/2-x_1A_{5,1}-x_2A_{4,1}
 -x_3A_{3,1}+x_4A_{2,1}\\ & &-x_8B_{4,5}-x_{10}B_{3,5}
+x_{12}B_{2,5}+x_{13}B_{3,4}-x_{14}B_{2,4}+x_{15}B_{2,3},\hspace{2.1cm}(11.5.22)\end{eqnarray*}
\begin{eqnarray*}
\iota(\eta_8)&=& P_8+x_8(A_{2,3}+A_{3,3}-\sum_{i=1,4,5
}A_{i,i}-\kappa)/2-x_1C_{4,1}+x_2C_{5,1}
 -x_5A_{3,1}+x_6A_{2,1}\\ & & -x_7C_{5,4}-x_{10}A_{3,4}
+x_{12}A_{2,4}-x_{13}A_{3,5}
+x_{14}A_{2,5}-x_{16}B_{2,3},\hspace{2.1cm}(11.5.23)\end{eqnarray*}
\begin{eqnarray*}
\iota(\eta_9) &=& P_9+x_9(A_{1,1}+A_{4,4} -\sum_{i=2,3,5
}A_{i,i}-\kappa)/2 +x_1C_{3,2}-x_3C_{5,2}
 -x_4C_{5,3}-x_5A_{4,2}\\& &-x_6A_{4,3}
 +x_{10}A_{1,2}+x_{11}A_{4,5}+x_{12}A_{1,3}
-x_{15}A_{1,5}+x_{16}B_{1,4},\hspace{2.1cm}(11.5.24)\end{eqnarray*}
\begin{eqnarray*}
\iota(\eta_{10}) &=& P_{10}+x_{10}(A_{2,2}+A_{4,4} -\sum_{i=1,3,5
}A_{i,i}-\kappa)/2-x_1C_{3,1}+x_3C_{5,1}
 +x_5A_{4,1}
 -x_7C_{5,3}\\ & &-x_8A_{4,3}+x_9A_{2,1}
+x_{12}A_{2,3}+x_{13}A_{4,5}
-x_{15}A_{2,5}+x_{16}B_{2,4},\hspace{2.1cm}(11.5.25)\end{eqnarray*}
\begin{eqnarray*}
\iota(\eta_{11}) &=&P_{11}+x_{11}(A_{1,1}+A_{5,5} -\sum_{i=2,3,4
}A_{i,i}-\kappa)/2 -x_2C_{3,2}+x_3C_{4,2}
 +x_4C_{4,3}-x_5A_{5,2}\\ & &-x_6A_{5,3}
 +x_9A_{5,4}+x_{13}A_{1,2}
+x_{14}A_{1,3}+x_{15}A_{1,4}+x_{16}B_{1,5},\hspace{2.1cm}(11.5.26)\end{eqnarray*}
\begin{eqnarray*}
\iota(\eta_{12})&=& P_{12}+x_{12}(A_{3,3}+A_{4,4} -\sum_{i=1,2,5
}A_{i,i}-\kappa)/2+x_1C_{2,1} +x_4C_{5,1}
 +x_6A_{4,1}
 +x_7C_{5,2}\\ & &+x_8A_{4,2}+x_9A_{3,1}
+x_{10}A_{3,2}+x_{14}A_{4,5}
-x_{15}A_{3,5}+x_{16}B_{3,4},\hspace{2.1cm}(11.5.27)\end{eqnarray*}
\begin{eqnarray*}
\iota(\eta_{13})&=& P_{13}+x_{13}(A_{2,2}+A_{5,5} -\sum_{i=1,3,4
}A_{i,i}-\kappa)/2 +x_2C_{3,1} -x_3C_{4,1}
 +x_5A_{5,1}+x_7C_{4,3}\\ & &-x_8A_{5,3}
 +x_{10}A_{5,4}+x_{11}A_{2,1}+x_{14}A_{2,3}
+x_{15}A_{2,4}+x_{16}B_{2,5},\hspace{2cm}(11.5.28)\end{eqnarray*}
\begin{eqnarray*}
\iota(\eta_{14})&=& P_{14}+x_{14}(A_{2,2}+A_{5,5} -\sum_{i=1,2,4
}A_{i,i}-\kappa)/2 -x_2C_{2,1} -x_4C_{4,1}
 +x_6A_{5,1}-x_7C_{4,2}\\& &+x_8A_{5,2}
 +x_{11}A_{3,1}+x_{12}A_{5,4}
+x_{13}A_{3,2}+x_{15}A_{3,4}+x_{16}B_{3,5},\hspace{2cm}(11.5.29)\end{eqnarray*}
\begin{eqnarray*}
\iota(\eta_{15})&=& P_{15}+x_{15}(A_{2,2}+A_{5,5} -\sum_{i=1,2,3
}A_{i,i}-\kappa)/2+x_3C_{2,1}+x_4C_{3,1}
 +x_7C_{3,2}-x_9A_{5,1}\\ & &-x_{10}A_{5,2}
 +x_{11}A_{4,1}-x_{12}A_{5,3}
+x_{13}A_{4,2}
+x_{14}A_{4,3}+x_{16}B_{4,5},\hspace{1.8cm}(11.5.30)\end{eqnarray*}
\begin{eqnarray*}
\iota(\eta_{16})&=& P_{16}-x_{16}(\sum_{i=1
}^5A_{i,i}+\kappa)/2-x_5C_{2,1} -x_6C_{3,1}
 -x_8C_{3,2}+x_9C_{4,1}\\ & &+x_{10}C_{4,2}
 +x_{11}C_{5,1}+x_{12}C_{4,3}
+x_{13}C_{5,2}+x_{14}C_{5,3}+x_{15}C_{5,4}.\hspace{2cm}(11.5.31)\end{eqnarray*}

Recall ${\msr B}=\mbb F[x_1,...,x_{16}]$. Let $M$ be an $o(10,\mbb
F)$-module and set
$$\widehat{M}={\msr B}\otimes_{\mbb F}M.\eqno(11.5.32)$$
We identify
$$f\otimes v=fv\qquad\for\;\;f\in{\msr B},\;v\in M.\eqno(11.5.33)$$
Recall the Lie algebra $\widehat{\msr K}$ defined via
(11.5.1)-(11.5.8). Fix $c\in\mbb F$. Then $\widehat M$ becomes a
$\widehat{\msr K}$-module with the action defined by
$$d(fv)=d(f)v,\;\;\kappa(fv)=cfv,\;\;(gB)(fv)=fg
B(v)\eqno(11.5.34)$$ for $d\in{\msr W}_{16},\;f,g\in{\msr B},\;v\in
M$ and $B\in o(10,\mbb F)$.

Since the linear map $\iota: {\msr G}^{E_6}\rta \widehat{\msr K}$
defined in (11.5.10)-(11.5.12) is a  Lie algebra monomorphism,
$\widehat{M}$ becomes a ${\msr G}^{E_6}$-module with the action
defined by
$$\xi(w)=\iota(\xi)(w)\qquad\for\;\;\xi\in{\msr
G}^{E_6},\;w\in\widehat{M}.\eqno(11.5.35)$$ In fact, we have:\psp

{\bf Theorem 11.5.1}. {\it The map $M\mapsto \widehat M$ gives a
functor from $o(10,\mbb F)$-{\bf Mod} to ${\msr G}^{E_6}$-{\bf Mod},
where the morphisms are $o(10,\mbb F)$-module homomorphisms and
${\msr G}^{E_6}$-module homomorphisms, respectively.}\psp

We remark that  the module $\widehat M$ is not a generalized Verma
module in general because it may not be equal to $U({\msr
G})(M)=U({\msr G}_-)(M)$. \psp

{\bf Proposition 11.5.2}. {\it If $M$ is an irreducible $o(10,\mbb
F)$-module, then $U({\msr G}_-)(M)$ is an irreducible ${\msr
G}^{E_6}$-module.}

{\it Proof.} Note that for any $i\in\ol{1,16}$, $f\in{\msr B}$ and
$v\in M$, (11.5.12), (11.5.34) and (11.5.35) imply
$$\xi_i(fv)=\ptl_{x_i}(f)v.\eqno(11.5.36)$$
Let $W$ be any nonzero ${\msr G}^{E_6}$-submodule. The above
expression shows that $W\bigcap M\neq\{0\}$. According to (11.5.34),
$W\bigcap M$ is an $o(10,\mbb F)$-submodule. By the irreducibility
of $M$, $M\subset W$. Thus $U({\msr G}_-)(M)\subset W$. So $U({\msr
G}_-)(M)=W$ is irreducible. $\qquad\Box$\psp

By the above proposition, the map $M\mapsto U({\msr G}_-)(M)$ is a
polynomial extension from  irreducible $o(10,\mbb F)$-modules to
irreducible ${\msr G}^{E_6}$-modules.  We can use it to derive
Gel'fand-Zetlin bases for $E_6$ from those for $o(10,\mbb F)$.

\section{Irreducibility of the Functor}

 In this section, we want to determine the irreducibility of the ${\msr
G}^{E_6}$-module $\widehat{M}$.

Note that $\widehat{M}$ can be viewed as an $o(10,\mbb F)$-module
with the representation $\iota(\nu(B))|_{\widehat{M}}$ (cf.
(11.3.30)-(11.3.34)). Indeed, (11.5.11) and (11.5.35) show
$$\nu(B(fv)=B(f)v+fB(v)\qquad\for\;\;B\in o(10,\mbb F),\;f\in{\msr
A},\;v\in M\eqno(11.6.1)$$ (cf. (11.3.36)-(11.3.57)). So
$\widehat{M}={\msr B}\otimes_{\mbb F}M$ is a tensor module of
$o(10,\mbb F)$. Write
$$\eta^\al=\prod_{i=1}^{16}\eta_{i}^{\alpha_i},\;\;
 |\al|=\sum_{i=1}^{16}\al_i\;\; \for\;\;
\al=(\al_1,\al_2,...,\al_{16})\in\mbb{N}^{16}\eqno(11.6.2)$$ (cf.
(11.3.14)-(11.3.19)).
 Recall the Lie subalgebras ${\msr G}_\pm$ and ${\msr G}_0$ of ${\msr G}^{E_6}$ defined
 in (11.3.20).  For $k\in\mbb{N}$, we set
$${\msr B}_k=\mbox{Span}\{x^\al\mid
\al \in\mbb{N}^{16};|\al|=k\},\;\; \widehat M_k={\msr
B}_kM\eqno(11.6.3)$$ (cf. (11.3.27), (11.5.34)) and
 $$(U({\msr G}_-)(
M))_k=\mbox{Span}\{\eta^\al( M)\mid \al \in\mbb{N}^{16}, \
|\al|=k\}. \eqno(11.6.4)$$
 Moreover,
$$(U({\msr G}_-)(
M))_0=\widehat M_0= M.\eqno(11.6.5)$$ Furthermore,
 $$\widehat M=\bigoplus\limits_{k=0}^\infty\widehat M_{
 k},\qquad
 U({\msr G}_-)(M)=\bigoplus\limits_{k=0}^\infty(U({\msr G}_-)(M))_k.\eqno(11.6.6)$$

Next we define a linear transformation $\vf$ on  $\widehat M$
determined by
$$\vf(x^\al v)=\eta^\al(
v)\qquad\for\;\;\al\in\mbb{N}^{16},\;v\in M.\eqno(11.6.7)$$ Note
that ${\msr B}_1=\sum_{i=1}^{16}\mbb Fx_i$ forms the
$16$-dimensional ${\msr G}_0$-module (equivalently, the $o(10,\mbb
F)$ spin module). According to (4.4.24)) and (4.4.25), ${\msr G}_-$
forms a ${\msr G}_0$-module with respect to the adjoint
representation, and the linear map from ${\msr B}_1$ to ${\msr G}_0$
determined by $x_i\mapsto \eta_i$ for $i\in\ol{1,16}$ gives a ${\msr
G}_0$-module isomorphism. Thus $\vf$ can also be viewed as a ${\msr
G}_0$-module homomorphism from $\widehat M$ to $U({\msr G}_0)(M)$.
Moreover,
$$\vf(\widehat M_k)=(U({\msr G}_-)(
M))_k\qquad\for\;\;k\in\mbb{N}.\eqno(11.6.8)$$

Note that the Casimir element of $o(10,\mbb F)$ is
$$\omega=\sum_{1\leq i<j\leq 5}(B_{i,j}C_{j,i}+C_{j,i}B_{i,j})
+\sum_{r,s=1}^5A_{r,s}A_{s,r}\in U(o(10,\mbb F))\eqno(11.6.9)$$ (cf.
(11.3.4)). The algebra $U(o(10,\mbb F))$ can be imbedded into the
tensor algebra $U(o(10,\mbb F))\otimes U(o(10,\mbb F))$ by the
associative algebra homomorphism $\mfk{d}: U(o(10,\mbb F))
\rightarrow U(o(10,\mbb F))\otimes_\mbb{F} U(o(10,\mbb F))$
determined  by
$$\mfk{d}(u)=u\otimes 1 +1 \otimes u \qquad \mbox{ for} \ u\in
o(10,\mbb F).\eqno(11.6.10)$$ Set
$$\td\omega=\frac{1}{2}(\mfk{d}(\omega)-\omega\otimes 1-1\otimes
\omega)\in U(o(10,\mbb F))\otimes_\mbb{F} U(o(2n,\mbb
F)).\eqno(11.6.11)$$ By (11.6.9),
$$\td\omega=\sum_{1\leq i<j\leq
5}(B_{i,j}\otimes C_{j,i}+C_{j,i}\otimes
B_{i,j})+\sum_{r,s=1}^5A_{r,s}\otimes A_{s,r}.\eqno(11.6.12)$$
Furthermore, $\td{\omega}$ acts on $\widehat M$ as an $o(10,\mbb
F)$-module homomorphism via
$$(B_1\otimes B_2)(fv)=B_1(f)B_2(v)\qquad\for\;\;B_1,B_2\in
o(10,\mbb F),\;f\in{\msr B},\;v\in M.\eqno(11.6.13)$$

{\bf Lemma 11.6.1}. {\it We have $\vf|_{\widehat M_
1}=(\td\omega-c/2)|_{\widehat M_1}$.}

{\it Proof.} For any $v\in M$,  (11.3.36)-(11.3.57), (11.5.16),
(11.6.12) and (11.6.13) give
\begin{eqnarray*}\td\omega(x_1v)&=&[-x_2A_{4,5}-x_3A_{3,5}
+x_4A_{2,5}-x_5B_{3,4}-x_7A_{1,5} +x_6B_{2,4}\\
& &-x_8B_{1,4}+x_9B_{2,3}-x_{10}B_{1,3}
+x_{12}B_{1,7}+x_1(\sum_{i=1}^4A_{i,i}-A_{5,5})]v
\\ &=& \eta_1(v)+(c/2)x_1v=(\vf+c/2)(x_1v),\hspace{5.9cm}(11.6.14)
\end{eqnarray*}
or equivalently, $\vf(x_1v)=(\td\omega-c/2)(x_1v)$. According to
(11.3.63), $-A_{5,4}(x_1)=x_2$. So
\begin{eqnarray*}& &\vf(x_2v)-\vf(x_1A_{5,4}(v))\\&=&
-A_{5,4}(\vf(x_1v)
=-A_{5,4}[(\td\omega-c/2)(x_1v)]=(c/2-\td\omega)A_{5,4}(x_1v)
\\
&=&(\td\omega-c/2)(x_2v-x_1A_{5,4}(v))
=(\td\omega-c/2)(x_2v)-(\td\omega-c/2)(x_1A_{5,4}(v))
\\ &=&(\td\omega-c/2)(x_2v)-\vf(x_1A_{5,4}(v)),\hspace{7.8cm}(11.6.15)\end{eqnarray*}
or equivalently, $\vf(x_2v)=(\td\omega-c/2)(x_2v)$.

Observe that
$$-A_{4,3}(x_2)=x_3,\;\;A_{3,2}(x_3)=x_4,\;\;-C_{5,4}(x_3)=x_5,\;\;-A_{3,1}(x_3)=x_7,\eqno(11.6.16)$$
$$C_{5,1}(x_2)=x_8,\;\;-C_{5,2}(x_2)=x_6,\;\;-C_{3,1}(x_1)=x_{10},\;\;-C_{3,2}(x_2)=x_{11},\eqno(11.6.17)$$
$$A_{3,2}(x_{10})=x_{12},\;\;C_{3,2}(x_1)=x_9,\;\;C_{3,1}(x_2)=x_{13},\eqno(11.6.18)$$
$$A_{3,2}(x_{13})=x_{14},\;\;
A_{4,3}(x_{14})=x_{15},\;\;C_{10,4}(x_{15})=x_{16}\eqno(11.6.19)$$
by (11.3.60)-(11.3.64). Using the argument similar to (11.6.14) and
induction, we can prove
$$\vf(x_iv)=(\td\omega-c/2)(x_iv)\qquad\for\;\;i\in\ol{1,16},\eqno(11.6.20)$$
or equivalently, the lemma holds. $\qquad\Box$\psp

We take (7.1.2)-(7.1.10) and (7.3.2)-(7.2.4) with $n=5$. For any
$\mu\in \Lmd^+$, we denote by $V(\mu)$ the finite-dimensional
irreducible $o(10,\mbb F)$-module with the highest weight $\mu$ and
have
$$\omega|_{V(\mu)}=(\mu+2\rho,\mu)\mbox{Id}_{V(\mu)}\eqno(11.6.21)$$
by (11.6.9).

Let $\mbb{Z}_2=\mbb{Z}/2\mbb{Z}=\{0,1\}$. According to (11.3.57) and
Table 11.3.1, the weight set of the $o(10,\mbb F)$-module ${\msr
B}_1$ is
$$\Lmd({\msr B}_1)=\{(1/2)\sum_{i=1}^5(-1)^{k_i}\ves_i\mid
k_i\in\mbb{Z}_2,\;\sum_{i=1}^5k_i=1\}.\eqno(11.6.22)$$ Fixing
$\lmd\in\Lmd^+$, we define
$$\mfk{H}(\lmd)=\{\lmd+\mu\mid\mu\in\Lmd({\msr
A}_1),\;\lmd+\mu\in\Lmd^+\}.\eqno(11.6.23)$$ \pse

{\bf Lemma 11.6.2}. {\it We have}:
$${\msr B}_1\otimes V(\lmd)\cong \bigoplus_{\lmd'\in
\mfk{H}(\lmd)}V(\lmd').\eqno(11.6.24)$$

{\it Proof}. Note that all the weight subspaces of ${\msr B}_1$ are
one-dimensional. Thus all the irreducible components of ${\msr
A}_1\otimes V(\lmd)$ are of multiplicity one. Since
$$\rho+\lmd+\mu\in\Lmd^+\qquad\for\;\;\mu\in \Lmd({\msr
A}_1),\eqno(11.6.25)$$ Theorem 5.4.3  says that $V(\lmd')$ is a
component of ${\msr B}_1\otimes V(\lmd)$ if and only if $\lmd'\in
\mfk{H}(\lmd).\qquad\Box$\psp

Recall
$$\mbox{the highest weight of}\;{\msr
B}_1=\frac{1}{2}(\ves_1+\ves_2+\ves_3+\ves_4-\ves_5)=\lmd_4,\eqno(11.6.26)$$
the forth fundamental weight of $o(10,\mbb F)$, by (11.3.57) and
Table 11.3.1. Thus the eigenvalues of
$\td{\omega}|_{\widehat{V(\lmd)}_1}$ are
$$\{[(\lmd'+2\rho,\lmd')-(\lmd+2\rho,\lmd)-(\lmd_4+2\rho,\lmd_4)]/2\mid\lmd'\in
\mfk{H}(\lmd)\}\eqno(11.6.27)$$ by (11.6.11) and (11.6.13). We
remark that the above fact is equivalent to special detailed version
of Kostant's characteristic identity (cf. [Kb]). Define
$$\ell_\omega(\lmd)=\min\{[(\lmd'+2\rho,\lmd')-(\lmd+2\rho,\lmd)-(\lmd_4+2\rho,\lmd_4)]/2\mid\lmd'\in
\mfk{H}(\lmd)\},\eqno(11.6.28)$$ which will be used to determine the
irreducibility of $\widehat{V(\lmd)}$. If
$\lmd'=\lmd+\lmd_4-\al\in\mfk{H}(\lmd)$ with $\al\in\Phi_{D_5}^+$,
then
$$(\lmd'+2\rho,\lmd')-(\lmd+2\rho,\lmd)-(\lmd_4+2\rho,\lmd_4)
=2[(\lmd,\lmd_4)+1-(\rho+\lmd+\lmd_4,\al)].\eqno(11.6.29)$$

Recall the differential operators $P_1,...,P_{16}$ given in
(11.4.16)-(11.4.32). We also view the elements of ${\msr B}$ as the
multiplication operators on ${\msr B}$. Recall $\zeta_1$ in
(11.3.59). It turns out that we need the following lemma in order to
determine the irreducibility of $\widehat{V(\lmd)}$.\psp

{\bf Lemma 11.6.3}. {\it As operators on} ${\msr B}$:
$$P_{11}x_1+P_1x_{11}+P_9x_2+P_2x_9-P_6x_3-P_3x_6+P_5x_4+P_4x_5=\zeta_1(D-6).\eqno(11.6.30)$$

{\it Proof}. According to  (11.4.16), (11.4.18)-(11.4.22), (11.4.25)
and (11.4.27), we find that
\begin{eqnarray*}\qquad&
&P_{11}x_1+P_1x_{11}+P_9x_2+P_2x_9-P_6x_3-P_3x_6+P_5x_4+P_4x_5\\
&=&-6\zeta_1+x_1P_{11}+x_{11}P_1+x_2P_9+x_9P_2-x_3P_6-x_6P_3+x_4P_5+x_5P_4\hspace{0.7cm}(11.6.31)\end{eqnarray*}
and
\begin{eqnarray*}&
&x_1P_{11}+x_{11}P_1+x_2P_9+x_9P_2-x_3P_6-x_6P_3+x_4P_5+x_5P_4 \\&=&
x_1(x_{11}D-\zeta_1\ptl_{x_1}+\zeta_5\ptl_{x_7}+\zeta_9\ptl_{x_8}+\zeta_8\ptl_{x_{10}}-\zeta_7\ptl_{x_{12}})
\\&
&+x_{11}(x_1D-\zeta_1\ptl_{x_{11}}-\zeta_2\ptl_{x_{13}}-\zeta_3\ptl_{x_{14}}-\zeta_4\ptl_{x_{15}}-\zeta_{10}\ptl_{x_{16}})
\\&&+x_2(x_9D-\zeta_1\ptl_{x_2}+\zeta_4\ptl_{x_7}-\zeta_{10}\ptl_{x_8}-\zeta_8\ptl_{x_{13}}+\zeta_7\ptl_{x_{14}})
\\&&+x_9(x_2D-\zeta_1\ptl_{x_9}-\zeta_2\ptl_{x_{10}}-\zeta_3\ptl_{x_{12}}+\zeta_5\ptl_{x_{15}}-\zeta_9\ptl_{x_{16}})
\\&&-x_3(x_6D+\zeta_1\ptl_{x_3}-\zeta_3\ptl_{x_7}+\zeta_{10}\ptl_{x_{10}}-\zeta_9\ptl_{x_{13}}+\zeta_7\ptl_{x_{15}})
\\&&-x_6(x_3D+\zeta_1\ptl_{x_6}+\zeta_2\ptl_{x_8}+\zeta_4\ptl_{x_{12}}+\zeta_5\ptl_{x_{14}}-\zeta_8\ptl_{x_{16}})\\
&&+x_4(x_5D-\zeta_1\ptl_{x_4}-\zeta_2\ptl_{x_7}-\zeta_{10}\ptl_{x_{12}}-\zeta_9\ptl_{x_{14}}-\zeta_8\ptl_{x_{15}})
\\&&+x_5(x_4D-\zeta_1\ptl_{x_5}+\zeta_3\ptl_{x_8}-\zeta_4\ptl_{x_{10}}-\zeta_5\ptl_{x_{13}}+\zeta_7\ptl_{x_{16}})
\\&=&2\zeta_1D-\zeta_1\sum_{i=1,2,3,4,5,6,9,11}x_i\ptl_{x_i}+(x_1\zeta_5+x_2\zeta_4+x_3\zeta_3-x_4\zeta_2)\ptl_{x_7}
\\&&+(x_1\zeta_9-x_2\zeta_{10}-x_6\zeta_2+x_5\zeta_3)\ptl_{x_8}
+(x_1\zeta_8-x_9\zeta_2-x_3\zeta_{10}-x_5\zeta_4)\ptl_{x_{10}}
\\&&-(x_1\zeta_7+x_9\zeta_3+x_6\zeta_4+x_4\zeta_{10})\ptl_{x_{12}}
-(x_{11}\zeta_2+x_2\zeta_8-x_3\zeta_9+x_5\zeta_5)\ptl_{x_{13}}
\\&&-(x_{11}\zeta_3-x_2\zeta_7+x_6\zeta_5-x_4\zeta_9)\ptl_{x_{14}}
-(x_{11}\zeta_4-x_9\zeta_5+x_3\zeta_7+x_4\zeta_8)\ptl_{x_{15}}
\\&&-(x_{11}\zeta_{10}+x_9\zeta_9-x_6\zeta_8-x_5\zeta_7)\ptl_{x_{16}}
=\zeta_1D. \qquad\Box\hspace{4.6cm}(11.6.32)\end{eqnarray*}\pse

We define the multiplication
$$f(gv)=(fg)v\qquad\for\;\;f,g\in{\msr B},\;v\in M.\eqno(11.6.33)$$
Then (11.5.16)-(11.5.21), (11.5.24) and (11.5.26) gives
\begin{eqnarray*}\qquad&&
\sum_{r=1}^{16}x_r\varrho([\xi_r,x_1\eta_{11}+x_{11}\eta_1+x_2\eta_9+x_9\eta_2-x_3\eta_6-x_6\eta_3+x_4\eta_5+x_5\eta_4])
\\&=&\sum_{i=1}^5\zeta_iA_{1,i}+\sum_{r=2}^5\zeta_{5+r}B_{1,r}-c\zeta_1
\hspace{7.3cm}(11.6.34)\end{eqnarray*} as operators on $\widehat M$
(cf. (11.5.33)), where $\zeta_i$ are defined in (11.3.59) and
(11.3.65)-(11.3.83).
 By Lemma 11.6.3, (11.5.5), (11.5.6) and (11.6.32),
\begin{eqnarray*}T_1&=&\iota(\eta_{11})x_1+\iota(\eta_1)x_{11}+\iota(\eta_9)x_2+\iota(\eta_2)x_9-\iota(\eta_6)x_3
-\iota(\eta_3)x_6+\iota(\eta_5)x_4+\iota(\eta_4)x_5
\\
&=&\zeta_1(D-c-6)+\sum_{i=1}^5\zeta_iA_{1,i}+\sum_{r=2}^5\zeta_{5+r}B_{1,r}
\hspace{5.8cm}(11.6.35)\end{eqnarray*}
 as operators on $\widehat M$. We define an $o(10,\mbb F)$-module
 structure on the space $\mbox{End}\:\widehat M$ of linear
 transformations on $\widehat M$ by
 $$B(T)=[\iota(\nu(B)),T]=\iota(\nu(B))T-T\iota(\nu(B))\qquad\for\;\;B\in
 o(10,\mbb F),\;T\in \mbox{End}\:\widehat M\eqno(11.6.36)$$
(cf. (11.6.1)). It can be verified that $T_1$ is an $o(10,\mbb
F)$-singular vector with weight $\ves_1$ in $\mbox{End}\:\widehat
M$. So it generates the 10-dimensional natural module.  Set
$$T_2=\zeta_2(D-c-6)+\sum_{i=1}^5\zeta_iA_{2,i}+
\sum_{r=1,3,4,5}\zeta_{5+r}B_{2,r}, \eqno(11.6.37)$$
$$T_3=\zeta_3(D-c-6)+\sum_{i=1}^5\zeta_iA_{3,i} +
\sum_{r=1,2,4,5}\zeta_{5+r}B_{3,r}, \eqno(11.6.38)$$
$$T_4=\zeta_4(D-c-6)+\sum_{i=1}^5\zeta_iA_{4,i} +
\sum_{r=1,2,3,5}\zeta_{5+r}B_{4,r},\eqno(11.6.39)$$
$$T_5=\zeta_5(D-c-6)+\sum_{i=1}^5\zeta_iA_{5,i}
-\sum_{r=1}^4\zeta_{5+r}B_{r,5},\eqno(11.6.40)$$
$$T_6=\zeta_6(D-c-6)-\sum_{i=2}^5\zeta_iC_{i,1}
-\sum_{r=1}^5\zeta_{5+r}A_{r,1},\eqno(11.6.41)$$
$$T_7=\zeta_7(D-c-6)-\sum_{i=1,3,4,5}\zeta_iC_{i,2}
-\sum_{r=1}^5\zeta_{5+r}A_{r,2},\eqno(11.6.42)$$
$$T_8=\zeta_8(D-c-6)+\sum_{i=1,2,4,5}\zeta_iC_{3,i}
-\sum_{r=1}^5\zeta_{5+r}A_{r,3},\eqno(11.6.43)$$
$$T_9=\zeta_9(D-c-6)+\sum_{i=1,2,3,5}\zeta_iC_{4,i}
-\sum_{r=1}^5\zeta_{5+r}A_{r,4},\eqno(11.6.44)$$
$$T_{10}=\zeta_{10}(D-c-6)+\sum_{i=1}^4\zeta_iC_{5,i}
-\sum_{r=1}^5\zeta_{5+r}A_{r,5}.\eqno(11.6.45)$$ Then ${\msr
T}=\sum_{i=1}^{10}\mbb FT_i$ forms the 10-dimensional natural module
of $o(10,\mbb F)$ with the standard basis $\{T_1,...,T_{10}\}$.

Denote
$$T'_i=T_i-\zeta_i(D-c-6)\qquad\for\;\;i\in\ol{1,10}.\eqno(11.6.46)$$
Easily see that  ${\msr T}'=\sum_{i=1}^{10}\mbb FT_i'$ forms the
10-dimensional natural module of $o(10,\mbb F)$ with the standard
basis $\{T_1',...,T_{10}'\}$. So we have the $o(10,\mbb F)$-module
isomorphism from $U=\sum_{i=1}^{10}\mbb F\zeta_i$ to ${\msr T}'$
determined by $\zeta_i\mapsto T_i'$ for $i\in\ol{1,10}$. The weight
set of $U$ is
$$\Lmd(U)=\{\pm\ves_1,...,\pm\ves_5\}.\eqno(11.6.47)$$
Let $\lmd\in\Lmd^+$. Denote
$$\mfk{H}'(\lmd)=\{\lmd+\mu\mid
\mu\in\Lmd(U),\;\;\lmd+\mu\in\Lmd^+\}.\eqno(11.6.48)$$ Take
$M=V(\lmd)$. It is known that
$$UV(\lmd)=U\otimes_{\mbb F}V(\lmd)\cong
\bigoplus_{\lmd'\in \mfk{H}'(\lmd)}V(\lmd').\eqno(11.6.49)$$

Given $\lmd'\in\mfk{H}'(\lmd)$, we pick a singular vector
$$u=\sum_{i=1}^{10}\zeta_iu_i\eqno(11.6.50)$$
of weight $\lmd'$ in $UV(\lmd)$, where $u_i\in V(\lmd)$. Moreover,
any singular vector of  weight $\lmd'$ in $UV(\lmd)$ is a scalar
multiple of $u$. Note that the vector
$$w=\sum_{i=1}^{10}T'_i(u_i)\eqno(11.6.51)$$
is also a  singular vector of  weight $\lmd'$ if it is not zero.
Thus
$$w=\flat_{\lmd'}u,\qquad \flat_{\lmd'}\in\mbb F.\eqno(11.6.52)$$
Set
$$\flat(\lmd)=\min\{\flat_{\lmd'}\mid\lmd'\in\mfk{H}'(\lmd)\}.\eqno(11.6.53)$$
\psp

{\bf Theorem 11.6.4}. {\it The ${\msr G}^{E_6}$-module
$\widehat{V(\lmd)}$ is irreducible if}
$$c\in\mbb F\setminus\{\flat(\lmd)-6+\mbb{N},2\ell_\omega(\lmd)+2\mbb{N}\}.\eqno(11.6.54)$$

{\it Proof}. Recall that the ${\msr G}^{E_6}$-submodule $U({\msr
G}_-)(V(\lmd))$ is irreducible by Proposition 11.5.2. It is enough
to prove $\widehat{V(\lmd)}=U({\msr G}_-)(V(\lmd))$. It is obvious
that
$$\widehat{V(\lmd)}_0=V(\lmd)=(U({\msr G}_-)(V(\lmd)))_0\eqno(11.6.55)$$ (cf. (11.6.3) and (11.6.4) with $M=V(\lmd)$).
Moreover, Lemma 11.6.1 with $M=V(\lmd)$,  (11.6.36) and (11.6.61)
imply that $\vf|_{\widehat{V(\lmd)}_1}$ is invertible, or
equivalently,
$$\widehat{V(\lmd)}_1=(U({\msr G}_-)(V(\lmd)))_1.\eqno(11.6.56)$$ Suppose that
$$\widehat{V(\lmd)}_i=(U({\msr G}_-)(V(\lmd)))_i\eqno(11.6.57)$$ for $i\in\ol{0,k}$ with $1\leq k\in\mbb{N}$.

For any $v\in V(\lmd)$ and $\al\in\mbb{N}^{16}$ such that
$|\al|=k-1$, we have
$$T_r(x^\al v)=x^\al[(|\al|-c-6)\zeta_r+T'_r](v)\in (U({\msr G}_-)(V(\lmd)))_{
k+1},\qquad r\in\ol{1,10}\eqno(11.6.58)$$ by (11.6.57) with
$i=k-1,k$. But
$$V'=\mbox{Span}\{[(|\al|-c-6)\zeta_r+T'_r](v)\mid r\in\ol{1,10},\;v\in
V(\lmd)\}\eqno(11.6.59)$$ forms an $o(10,\mbb F)$-submodule of
$UV(\lmd)$ with respect to the action in (11.6.1). Let $u$ be a
$o(10,\mbb F)$-singular vector in (11.6.50). Then
$$V'\ni
\sum_{r=1}^{10}[(|\al|-c-6)\zeta_r+T'_r](u_r)=(|\al|-c-6)u+w=(|\al|-c-6+\flat_{\lmd'})u\eqno(11.6.60)$$
by (11.6.51) and (11.6.52). Moreover, (11.6.53) and (11.6.54) yield
$u\in V'$. Since $UV(\lmd)$ is an $o(10,\mbb F)$-module generated by
all the singular vectors, we have $V'=UV(\lmd)$. So
$$x^\al UV(\lmd)\subset (U({\msr G}_-)(V(\lmd)))_{
k+1}.\eqno(11.6.61)$$ The arbitrariness of $\al$ implies
$$\zeta_r\widehat{V(\lmd)}_{ k-1}\subset (U({\msr G}_-)(V(\lmd)))_{
k+1}\qquad\for\;\;r\in\ol{1,10}.\eqno(11.6.62)$$

Given any $f\in{\msr B}_k$ and $v\in V(\lmd)$, we have
$$\zeta_r\ptl_{x_i}(f)v\in \zeta_r\widehat{V(\lmd)}_{ k-1}\subset (U({\msr G}_-)(V(\lmd)))_{
k+1}\qquad\for\;\;r\in\ol{1,10},\;i\in\ol{1,16}.\eqno(11.6.63)$$
Moreover,
\begin{eqnarray*}\qquad\eta_s(fv)&=&\iota(\eta_s)(fv)=P_s(fv)+f(\td\omega-c/2)(x_sv)\\
&\equiv&
f(k+\td\omega-c/2)(x_sv)\;\;(\mbox{mod}\;\sum_{r=1}^{10}\zeta_r\widehat{V(\lmd)}_{
k-1})\hspace{3.8cm}(11.6.64)\end{eqnarray*} for $s\in\ol{1,16}$ by
(11.4.16)-(11.4.32), (11.5.16)-(11.5.31) and Lemma 11.6.1. According
to (11.6.28), (11.6.54), (11.6.62) and (11.6.64), we get
$$x_sfv\in (U({\msr G}_-)(V(\lmd)))_{
k+1}\qquad\for\;\;s\in\ol{1,16}.\eqno(11.6.65)$$ Thus (11.6.57)
holds for $i=k+1$. By induction on $k$, (11.6.57) holds for any
$i\in\mbb{N}$; that is, $\widehat{V(\lmd)}=U({\msr
G}_-)(V(\lmd)).\qquad\Box$\psp

When $\lmd=0$, $V(0)$ is the one-dimensional trivial module and
$\ell_\omega(0)=\flat(0)=0$. So we have:\psp

{\bf Corollary 11.6.5}. {\it The ${\msr G}^{E_6}$-module
$\widehat{V(0)}$ is irreducible if} $c\in\mbb
F\setminus\{\mbb{N}-6\}.$ \psp

Next we consider the case $\lmd=k\ves_1=k\lmd_1$ for some positive
integer $k$, where $\lmd_1$ is the first fundamental weight. Note
$$\mfk{H}(k\ves_1)=\{\lmd_4+k\ves_1,\lmd_4+(k-1)\ves_1+\ves_5\}\eqno(11.6.66)$$ by
(11.6.22) and (11.6.23). Thus (11.6.28) and (11.6.29) give
$$\ell_{\omega}(k\ves_1)=-4-k/2.\eqno(11.6.67)$$
In order to calculate $\flat(k\ves_1)$, we give a realization of
$V(k\ves_1)$. Observe that we have a representation of $o(10,\mbb
F)$ on ${\msr B}=\mbb F[y_1,...,y_{10}]$ determined via
$$E_{i,j}|_{\msr
B}=y_i\ptl_{y_j}\qquad\for\;\;i,j\in\ol{1,10}.\eqno(11.6.68)$$
Denote by ${\msr B}_k$ the subspace of homogenous polynomials in
${\msr B}$ with degree $k$. Set
$${\msr H}_k=\{h\in{\msr B}_k\mid
(\sum_{i=1}^5\ptl_{y_i}\ptl_{y_{5+i}})(h)=0\}.\eqno(11.6.69)$$ Then
${\msr H}_k\cong V(k\ves_1)$ and $y_1^k$ is a highest-weight vector.

According to (11.6.47) and (11.6.48),
$$\mfk{H}'(k\ves_1)=\{(k+1)\ves_1,(k-1)\ves_1, k\ves_1+\ves_2\}.\eqno(11.6.70)$$
The vector $\zeta_1y_1^k$ is a singular vector in $U{\msr H}_k$ with
weight $(k+1)\ves_1$, where we take $M={\msr H}_k$ in the earlier
settings. By (11.6.35) and (11.6.46),
$$T_1'(y_1^k)=k\zeta_1y_1^k\lra \flat_{(k+1)\ves_1}=k.\eqno(11.6.71)$$
Moreover, $\zeta_1y_1^{k-1}y_2-\zeta_2y_1^k$  is a singular vector
in $U{\msr H}_k$ with weight $k\ves_1+\ves_2$. Moreover, (11.6.35),
(11.6.37) and (11.6.46) imply
\begin{eqnarray*}T'_1(y_1^{k-1}y_2)-T'_2(y_1^k)&=&(k-1)\zeta_1y_1^{k-1}y_2+\zeta_2y_1^k-k\zeta_1y_1^{k-1}y_2
\\ &=&\zeta_2y_1^k-\zeta_1y_1^{k-1}y_2=-(\zeta_1y_1^{k-1}y_2-\zeta_2y_1^k).\hspace{3.1cm}(11.6.72)\end{eqnarray*}
Thus $\flat_{k\ves_1+\ves_2}=-1$. Furthermore,
$$\varpi=(k+3)\sum_{i=1}^5[\zeta_iy_1^{k-1}y_{5+i}+\zeta_{5+i}y_1^{k-1}y_i]-(k-1)\zeta_1y_1^{k-2}\sum_{s=1}^5y_sy_{5+s}
\eqno(11.6.73)$$ is a singular vector in $U{\msr H}_k$ with weight
$(k-1)\ves_1$. Expressions (11.6.35) and (11.6.37)-(11.6.46) yield
\begin{eqnarray*}& &(k+3)\sum_{i=1}^5[T'_i(y_1^{k-1}y_{5+i})+T'_{5+i}(
y_1^{k-1}y_i)]-(k-1)T_1'(x^\al y_1^{k-2}\sum_{s=1}^5y_sy_{5+s})
\\ &=&(-8-k)x^\al\varpi\lra
\flat_{(k-1)\ves_1}=-8-k.\hspace{6.9cm}(11.6.74)
\end{eqnarray*}
Therefore, $\flat(k\ves_1)=-8-k.$ By Theorem 11.6.4 and (11.6.67),
we obtain:\psp

 {\bf Corollary 11.6.6}. {\it The ${\msr G}^{E_6}$-module $\widehat{V(k\lmd_1)}$ is
irreducible if} $c\in\mbb F\setminus\{\mbb{N}-14-k\}.$ \psp

By similar calculations, we obtain:\psp

 {\bf Corollary 11.6.7}. {\it (a) The ${\msr G}^{E_6}$-module $\widehat{V(\lmd_2)}$ is
irreducible if $c\in\mbb F\setminus\{\mbb{N}-16\}.$

(b) The ${\msr G}^{E_6}$-module $\widehat{V(\lmd_3)}$ is irreducible
if $c\in\mbb F\setminus\{\mbb{N}-15,-17,-19,-21\}.$

(c) The ${\msr G}^{E_6}$-module $\widehat{V(k\lmd_4)}$ is
irreducible if $c\in\mbb
F\setminus\{\mbb{N}-10-k/2,2\mbb{N}+k-12\}.$

(d) The ${\msr G}^{E_6}$-module $\widehat{V(k\lmd_5)}$ is
irreducible if $c\in\mbb
F\setminus\{\mbb{N}-10-k/2,2\mbb{N}-k-20\}$.}

\section{Representations on Exponential-Polynomial Functions}

In this section, we want to study representations of $\msr G^{E_6}$
on exponential-polynomial functions.

Recall
$$D=\sum_{i=1}^{16}x_i\ptl_{x_i}.\eqno(11.7.1)$$
Fix $\mfk{c}\in\mbb F$. Let $c=2\mfk{c}$ and identify
$\widehat{V(0)}={\msr B}\otimes v_0$ with ${\msr B}$ by
$$f\otimes v_0\leftrightarrow f\qquad\for\;\;f\in{\msr
B},\eqno(11.7.2)$$ where $V(0)=\mbb Fv_0$. Then we have the
following one-parameter  inhomogeneous first-order differential
operator representation $\pi_{\mfk{c}}$ of $\msr G^{E_6}$:
$$\pi_{\mfk{c}}(u)=\pi(\nu^{-1}(u))\qquad\for\;\;u\in \msr
G^{D_5}\eqno(11.7.3)$$(cf. (11.3.30)-(11.3.34) and
(11.3.36)-(11.3.57)),
$$\pi_{\mfk{c}}(\hat\al)=D+4\mfk{c},\qquad\pi_{\mfk{c}}(\xi_i)=\ptl_{x_i}\qquad\for\;\;i\in\ol{1,16},\eqno(11.7.4)$$
$$\pi_{\mfk{c}}(\eta_1)=
x_1(D-{\mfk{c}})-\zeta_1\ptl_{x_{11}}-\zeta_2\ptl_{x_{13}}-\zeta_3\ptl_{x_{14}}-\zeta_4\ptl_{x_{15}}-\zeta_{10}\ptl_{x_{16}},\eqno(11.7.5)$$
$$\pi_{\mfk{c}}(\eta_2)=x_2(D-\mfk{c})-\zeta_1\ptl_{x_9}-\zeta_2\ptl_{x_{10}}-\zeta_3\ptl_{x_{12}}+\zeta_5\ptl_{x_{15}}-\zeta_9\ptl_{x_{16}},
\eqno(11.7.6)$$
$$
\pi_{\mfk{c}}(\eta_3)=x_3(D-\mfk{c})+\zeta_1\ptl_{x_6}+\zeta_2\ptl_{x_8}+\zeta_4\ptl_{x_{12}}+\zeta_5\ptl_{x_{14}}-\zeta_8\ptl_{x_{16}},\eqno
(11.7.7)$$
$$\pi_{\mfk{c}}(\eta_4)=x_4(D-\mfk{c})-\zeta_1\ptl_{x_5}+\zeta_3\ptl_{x_8}-\zeta_4\ptl_{x_{10}}-\zeta_5\ptl_{x_{13}}+\zeta_7\ptl_{x_{16}},
\eqno(11.7.8)$$
$$\pi_{\mfk{c}}(\eta_5)=
x_5(D-\mfk{c})-\zeta_1\ptl_{x_4}-\zeta_2\ptl_{x_7}-\zeta_{10}\ptl_{x_{12}}+\zeta_9\ptl_{x_{14}}-\zeta_8\ptl_{x_{15}},\eqno(11.7.9)$$
$$\pi_{\mfk{c}}(\eta_6)=x_6(D-\mfk{c})+\zeta_1\ptl_{x_3}-\zeta_3\ptl_{x_7}+\zeta_{10}\ptl_{x_{10}}-\zeta_9\ptl_{x_{13}}+\zeta_7\ptl_{x_{15}},
\eqno(11.7.10)$$
$$\pi_{\mfk{c}}(\eta_7)=x_7(D-\mfk{c})-\zeta_2\ptl_{x_5}-\zeta_3\ptl_{x_6}+\zeta_4\ptl_{x_9}+\zeta_5\ptl_{x_{11}}-\zeta_6\ptl_{x_{16}},
\eqno(11.7.11)$$
$$\pi_{\mfk{c}}(\eta_8)=x_8(D-\mfk{c})+\zeta_2\ptl_{x_3}+\zeta_3\ptl_{x_4}-\zeta_{10}\ptl_{x_9}+\zeta_9\ptl_{x_{11}}-\zeta_6\ptl_{x_{15}},
\eqno(11.7.12)$$
$$\pi_{\mfk{c}}(\eta_9)=x_9(D-\mfk{c})-\zeta_1\ptl_{x_2}+\zeta_4\ptl_{x_7}-\zeta_{10}\ptl_{x_8}-\zeta_8\ptl_{x_{13}}+\zeta_7\ptl_{x_{14}},
\eqno(11.7.13)$$
$$\pi_{\mfk{c}}(\eta_{10})=x_{10}(D-\mfk{c})-\zeta_2\ptl_{x_2}-\zeta_4\ptl_{x_4}+\zeta_{10}\ptl_{x_6}+\zeta_8\ptl_{x_{11}}-\zeta_6\ptl_{x_{14}},
\eqno(11.7.14)$$
$$\pi_{\mfk{c}}(\eta_{11})=x_{11}(D-\mfk{c})-\zeta_1\ptl_{x_1}+\zeta_5\ptl_{x_7}+\zeta_9\ptl_{x_8}+\zeta_8\ptl_{x_{10}}-\zeta_7\ptl_{x_{12}},
\eqno(11.7.15)$$
$$\pi_{\mfk{c}}(\eta_{12})=
x_{12}(D-\mfk{c})-\zeta_3\ptl_{x_2}+\zeta_4\ptl_{x_3}-\zeta_{10}\ptl_{x_5}-\zeta_7\ptl_{x_{11}}+\zeta_6\ptl_{x_{13}},
\eqno(11.7.16)$$
$$\pi_{\mfk{c}}(\eta_{13})=
x_{13}(D-\mfk{c})-\zeta_2\ptl_{x_1}-\zeta_5\ptl_{x_4}-\zeta_9\ptl_{x_6}-\zeta_8\ptl_{x_9}+\zeta_6\ptl_{x_{12}},
\eqno(11.7.17)$$ $$\pi_{\mfk{c}}(\eta_{14})=
x_{14}(D-\mfk{c})-\zeta_3\ptl_{x_1}+\zeta_5\ptl_{x_3}+\zeta_9\ptl_{x_5}+\zeta_7\ptl_{x_9}-\zeta_6\ptl_{x_{10}},
\eqno(11.7.18)$$
$$\pi_{\mfk{c}}(\eta_{15})=x_{15}(D-\mfk{c})-\zeta_4\ptl_{x_1}+\zeta_5\ptl_{x_2}-\zeta_8\ptl_{x_5}+\zeta_7\ptl_{x_6}-\zeta_6\ptl_{x_8},
\eqno(11.7.19)$$
$$\pi_{\mfk{c}}(\eta_{16})=x_{16}(D-\mfk{c})-\zeta_{10}\ptl_{x_1}-\zeta_9\ptl_{x_2}-\zeta_8\ptl_{x_3}+\zeta_7\ptl_{x_4}-\zeta_6\ptl_{x_7}.
\eqno(11.7.20)$$

Let $\vec a=(a_1,...,a_{16})\in\mbb F^{16}\setminus\{\vec 0\}$.
Define
$$\vec a\cdot \vec x=\sum_{i=1}^{16}a_ix_i.\eqno(11.7.21)$$
Recall $\msr B=\mbb F[x_1,...,x_{16}]$ and set
$${\msr B}_{\vec a}=\{fe^{\vec a\cdot\vec
x}\mid f\in{\msr B}\}.\eqno(11.7.22)$$ Moreover, we write
$$\msr F=\{x_i\zeta_{r_1}\zeta_{r_2}\mid i\in\ol{1,16};\;
(r_1,s_1),(r_2,s_2)\in \Upsilon_i,\;r_1\neq r_2\}.\eqno(11.7.23)$$
By Lemma 11.3.3, $|\msr F|=160$. For any function
$f(x_1,...,x_{16})$, we define
$$f(\vec b)=f(b_1,...,b_{16})\qquad\for\;\;\vec b=(b_1,...,b_{16})\in\mbb
F^{16}.\eqno(11.7.24)$$ Now we define
$$\msr V=\{\vec b\in\mbb F^{16}\setminus\{\vec 0\}\mid f(\vec b)=0\;\for\;f\in\msr
F\},\eqno(11.7.25)$$which gives rise to a projective algebraic
variety. The following is our main theorem in this section.\psp

{\bf Theorem 11.7.1}. {\it With respect to the representation
$\pi_{\mfk{c}}$, ${\msr B}_{\vec a}$ forms an irreducible $\msr
G^{E_6}$-module for any $\mfk c\in\mbb F$ if $\vec a\not\in\msr V$.}

{\it Proof}. Let $\msr B_k$ be the subspace of homogeneous
polynomials with degree $k$. Set
$$\msr B_{\vec a,k}=\msr B_ke^{\vec a\cdot\vec
x}\qquad\for\;k\in\mbb{N}.\eqno(11.7.26)$$ Let  ${M}$ be a nonzero
$\msr G^{E_6}$-submodule of ${\msr B}_{\vec a}$. Take any $0\neq
fe^{\vec a\cdot\vec x}\in M$ with $f\in \msr B$. By the second
equation in (11.7.4),
$$(\xi_i-a_i)(fe^{\vec a\cdot\vec x})=\ptl_{x_i}(f)e^{\vec a\cdot\vec
x}\in M\qquad\for\;\;i\in\ol{1,16}.\eqno(11.7.27)$$ Repeatedly
applying (11.7.27) if necessarily, we obtain $e^{\vec a\cdot\vec
x}\in M$; that is, $\msr B_{\vec a,0}\subset M$.

Suppose $\msr B_{\vec a,\ell}\subset M$ for some $\ell\in\mbb{N}$.
Let $ge^{\vec a\cdot\vec x}$ be any element in $\msr B_{\vec
a,\ell}$. Note that
$$\pi_{\mfk{c}}(\msr
G^{E_6})=\pi(o(10,\mbb F))+\sum_{i=1}^{16}(\mbb
F\pi_{\mfk{c}}(\xi_i)+\mbb F\pi_{\mfk{c}}(\eta_i))+\mbb
F\pi_{\mfk{c}}(\hat\al)\eqno(11.7.28)$$ by (11.3.20), (11.3.21),
(11.4.2) and (11.7.3). Thus $M$ is also an $o(10,\mbb F)$-module.

By Lemma 11.3.2, we can assume that
$(x_1\zeta_{r_1}\zeta_{r_2})(\vec a)\neq 0$ for some
$r_1,r_2\in\{1,2,3,4,10\}$ with $r_1\neq r_2$ by (11.3.111). Under
the assumption, $a_1\neq 0$.

Applying (11.3.116) and (11.3.117) to $ge^{\vec a\cdot\vec x}$, we
get by (11.3.36)-(11.3.56) that
$$(-a_1x_2+a_9x_{11}+a_{10}a_{13}+a_{12}x_{14})ge^{\vec a\cdot\vec x}\equiv 0\;\;(\mbox{mod}\;M),\eqno(11.7.29)$$
$$(-a_1x_3-a_6x_{11}-a_8x_{13}-a_{12}x_{15})ge^{\vec a\cdot\vec x}\equiv 0\;\;(\mbox{mod}\;M),\eqno(11.7.30)$$
$$(a_1x_4-a_5x_{11}+a_8x_{14}-a_{10}x_{15})ge^{\vec a\cdot\vec x}\equiv 0\;\;(\mbox{mod}\;M),\eqno(11.7.31)$$
$$(-a_1x_5+a_4x_{11}+a_7x_{13}+a_{12}x_{16})ge^{\vec a\cdot\vec x}\equiv 0\;\;(\mbox{mod}\;M),\eqno(11.7.32)$$
$$(-a_1x_7+a_5x_{13}+a_6x_{14}-a_9x_{15})ge^{\vec a\cdot\vec x}\equiv 0\;\;(\mbox{mod}\;M),\eqno(11.7.33)$$
$$(a_1x_6+a_3x_{11}-a_7x_{14}+a_{10}x_{16})ge^{\vec a\cdot\vec x}\equiv 0\;\;(\mbox{mod}\;M),\eqno(11.7.34)$$
$$(-a_1x_8-a_4x_{14}-a_3x_{13}+a_9x_{16})ge^{\vec a\cdot\vec x}\equiv 0\;\;(\mbox{mod}\;M),\eqno(11.7.35)$$
$$(a_1x_9-a_2x_{11}+a_7x_{15}-a_8x_{16})ge^{\vec a\cdot\vec x}\equiv 0\;\;(\mbox{mod}\;M),\eqno(11.7.36)$$
$$(-a_1x_{10}+a_2x_{13}+a_4x_{15}-a_6x_{16})ge^{\vec a\cdot\vec x}\equiv 0\;\;(\mbox{mod}\;M),\eqno(11.7.37)$$
$$(a_1x_{12}-a_2x_{14}+a_3x_{15}-a_5x_{16})ge^{\vec a\cdot\vec x}\equiv 0\;\;(\mbox{mod}\;M).\eqno(11.7.38)$$

Multiplying $a_1$ to (11.3.119)-(11.3.125) and applying them to
$ge^{\vec a\cdot\vec x}$ by (11.3.36)-(11.3.56), we obtain
$$a_1(a_7x_4+a_8x_6+a_{10}x_9+a_{13}x_{11})ge^{\vec a\cdot\vec x}\equiv 0\;\;(\mbox{mod}\;M),\eqno(11.7.39)$$
$$a_1(-a_7x_3-a_8x_5+a_{12}x_9+a_{14}x_{11})ge^{\vec a\cdot\vec x}\equiv 0\;\;(\mbox{mod}\;M),\eqno(11.7.40)$$
$$a_1(-a_7x_2+a_{10}x_5+a_{12}x_6+a_{15}x_{11})ge^{\vec a\cdot\vec x}\equiv 0\;\;(\mbox{mod}\;M),\eqno(11.7.41)$$
$$a_1(a_8x_2+a_{10}x_3+a_{12}x_4+a_{16}x_{11})ge^{\vec a\cdot\vec x}\equiv 0\;\;(\mbox{mod}\;M),\eqno(11.7.42)$$
$$a_1(a_4x_7+a_6x_8+a_9x_{10}+a_{11}x_{13})ge^{\vec a\cdot\vec x}\equiv 0\;\;(\mbox{mod}\;M),\eqno(11.7.43)$$
$$a_1(a_4x_3+a_6x_5+a_{12}x_{10}+a_{14}x_{13})ge^{\vec a\cdot\vec x}\equiv 0\;\;(\mbox{mod}\;M),\eqno(11.7.44)$$
$$a_1(a_4x_2-a_9x_5+a_{12}x_8+a_{15}x_{13})ge^{\vec a\cdot\vec x}\equiv 0\;\;(\mbox{mod}\;M),\eqno(11.7.45)$$
$$a_1(-a_6x_2-a_9x_3+a_{12}x_7+a_{16}x_{13})ge^{\vec a\cdot\vec x}\equiv 0\;\;(\mbox{mod}\;M),\eqno(11.7.46)$$
$$a_1(-a_3x_7-a_5x_8+a_9x_{12}+a_{11}x_{14})ge^{\vec a\cdot\vec x}\equiv 0\;\;(\mbox{mod}\;M),\eqno(11.7.47)$$
$$a_1(a_3x_4+a_5x_6+a_{10}x_{12}+a_{13}x_{14})ge^{\vec a\cdot\vec x}\equiv 0\;\;(\mbox{mod}\;M),\eqno(11.7.48)$$
$$a_1(-a_3x_2-a_9x_6-a_{10}x_8+a_{15}x_{14})ge^{\vec a\cdot\vec x}\equiv 0\;\;(\mbox{mod}\;M),\eqno(11.7.49)$$
$$a_1(a_5x_2-a_9x_4-a_{10}x_7+a_{16}x_{14})ge^{\vec a\cdot\vec x}\equiv 0\;\;(\mbox{mod}\;M),\eqno(11.7.50)$$
$$a_1(-a_2x_7+a_5x_{10}+a_6x_{12}+a_{11}x_{15})ge^{\vec a\cdot\vec x}\equiv 0\;\;(\mbox{mod}\;M),\eqno(11.7.51)$$
$$a_1(a_2x_4-a_5x_9+a_8x_{12}+a_{13}x_{15})ge^{\vec a\cdot\vec x}\equiv 0\;\;(\mbox{mod}\;M),\eqno(11.7.52)$$
$$a_1(-a_2x_3-a_6x_9-a_8x_{10}+a_{14}x_{15})ge^{\vec a\cdot\vec x}\equiv 0\;\;(\mbox{mod}\;M),\eqno(11.7.53)$$
$$a_1(-a_5x_3-a_6x_4-a_8x_7+a_{16}x_{15})ge^{\vec a\cdot\vec x}\equiv 0\;\;(\mbox{mod}\;M),\eqno(11.7.54)$$
$$a_1(a_2x_8+a_3x_{10}+a_4x_{12}+a_{11}x_{16})ge^{\vec a\cdot\vec x}\equiv 0\;\;(\mbox{mod}\;M),\eqno(11.7.55)$$
$$a_1(-a_2x_6-a_3x_9+a_7x_{12}+a_{13}x_{16})ge^{\vec a\cdot\vec x}\equiv 0\;\;(\mbox{mod}\;M),\eqno(11.7.56)$$
$$a_1(a_2x_5-a_4x_9-a_7x_{10}+a_{14}x_{16})ge^{\vec a\cdot\vec x}\equiv 0\;\;(\mbox{mod}\;M),\eqno(11.7.57)$$
$$a_1(-a_3x_5-a_4x_6-a_7x_8+a_{15}x_{16})ge^{\vec a\cdot\vec x}\equiv 0\;\;(\mbox{mod}\;M).\eqno(11.7.58)$$

According to (11.7.31), (11.7.34) and (11.7.36),
$$a_1x_4ge^{\vec a\cdot\vec x}\equiv(a_5x_{11}-a_8x_{14}+a_{10}x_{15})ge^{\vec a\cdot\vec
x}\;\;(\mbox{mod}\;M),\eqno(11.7.59)$$
$$a_1x_6ge^{\vec a\cdot\vec x}\equiv (-a_3x_{11}+a_7x_{14}-a_{10}x_{16})ge^{\vec a\cdot\vec x}\;\;(\mbox{mod}\;M),
\eqno(11.7.60)$$
$$a_1x_9ge^{\vec a\cdot\vec x}\equiv (a_2x_{11}-a_7x_{15}+a_8x_{16})ge^{\vec a\cdot\vec x}\;\;(\mbox{mod}\;M).\eqno(11.7.61)$$
Substituting them into (11.7.39), we get
\begin{eqnarray*}\hspace{2cm}&
&[a_7(a_5x_{11}-a_8x_{14}+a_{10}x_{15})+a_8(-a_3x_{11}+a_7x_{14}-a_{10}x_{16})\\&
&+a_{10}(a_2x_{11}-a_7x_{15}+a_8x_{16})+a_1a_{13}x_{11}]ge^{\vec
a\cdot\vec x}\\
&=&[a_1a_{13}+a_2a_{10}-a_3a_8+a_5a_7]x_{11}ge^{\vec a\cdot\vec x}
\\ &=&\zeta_2(\vec a)x_{11}ge^{\vec a\cdot\vec x}\equiv 0\;\;(\mbox{mod}\;\msr
M)\hspace{6.4cm}(11.7.62)\end{eqnarray*} by (11.3.65) and (11.7.24).

Similarly we substitute (11.7.29)-(11.7.38) into
(11.7.40)-(11.7.58), and get by (11.3.59), (11.3.65)-(11.3.67) and
(11.3.69) that
$$\zeta_i(\vec a)x_{11}ge^{\vec a\cdot\vec x}\equiv 0\;\;(\mbox{mod}\;\msr
M)\qquad\for\;\;i=3,4,10;\eqno(11.7.63)$$
$$\zeta_i(\vec a)x_{13}ge^{\vec a\cdot\vec x}\equiv 0\;\;(\mbox{mod}\;\msr
M)\qquad\for\;\;i=1,3,4,10;\eqno(11.7.64)$$
$$\zeta_i(\vec a)x_{14}ge^{\vec a\cdot\vec x}\equiv 0\;\;(\mbox{mod}\;\msr
M)\qquad\for\;\;i=1,2,4,10;\eqno(11.7.65)$$
$$\zeta_i(\vec a)x_{15}ge^{\vec a\cdot\vec x}\equiv 0\;\;(\mbox{mod}\;\msr
M)\qquad\for\;\;i=1,2,3,10;\eqno(11.7.66)$$
$$\zeta_i(\vec a)x_{16}ge^{\vec a\cdot\vec x}\equiv 0\;\;(\mbox{mod}\;\msr
M)\qquad\for\;\;i=1,2,3,4.\eqno(11.7.67)$$

Since $\zeta_{r_1}(\vec a)\zeta_{r_2}(\vec a)\neq 0$ for some
$r_1,r_2\in\{1,3,4,5,10\}$ with $r_1\neq r_2$, (11.7.62)-(11.7.67)
imply
$$x_ige^{\vec a\cdot\vec x}\in\msr
M\qquad\for\;\;i=11,13,14,15,16.\eqno(11.7.68)$$ Substituting
(11.7.68) into (11.7.29)-(11.7.38), we obtain
$$x_ige^{\vec a\cdot\vec x}\in\msr
M\qquad\for\;\;i\in\ol{2,16}.\eqno(11.7.69)$$ By (11.3.57) and Table
11.3.1,
$$A_{1,1}(ge^{\vec a\cdot\vec x})\equiv a_1x_1ge^{\vec a\cdot\vec
x}/2\equiv 0 \;\;(\mbox{mod}\;M).\eqno(11.7.70)$$ So $x_1ge^{\vec
a\cdot\vec x}\in M$. Hence $\msr B_{\vec a,\ell+1}\subset M$. By
induction,
$$\msr B_{\vec
a,k}\subset M\qquad\for\;\;k\in\mbb N.\eqno(11.7.71)$$ Therefore,
$M=\sum_{k=0}^\infty\msr B_{\vec a,k}=\msr B_{\vec a}$. So $\msr
B_{\vec a}$ forms an irreducible $\msr G^{E_6}$-module. $\qquad\Box$

\chapter{Representations of $E_7$}

By solving certain partial differential equations, we find the
explicit decomposition of the polynomial algebra over the
56-dimensional basic irreducible module of the simple Lie algebra
$E_7$ into a direct sum of irreducible submodules, which was due to
our work [X22]. Moreover, we find a new representation of the simple
Lie algebra of type $E_7$ on the polynomial algebra in 27 variables,
which gives a fractional representation of the corresponding Lie
group on 27-dimensional space.  Using this representation and Shen's
idea of mixed product (cf. [Sg]), we construct a new functor from
the category of $E_6$-modules to the category of $E_7$-modules. A
condition for the functor to map a finite-dimensional irreducible
$E_6$-module to an infinite-dimensional irreducible $E_7$-module is
obtained. Our general frame also gives a direct polynomial extension
from irreducible $E_6$-modules to irreducible $E_7$-modules, which
can be used to derive Gel'fand-Zetlin bases for $E_7$ from those for
$E_6$ that can be obtained from those for $D_5$ in last chapter. Our
results also yield explicit constructions of certain
infinite-dimensional irreducible weight $E_7$-modules of
finite-dimensional weight subspaces.  In our approach, the idea of
Kostant's characteristic identities and the well-known Dickson's
$E_6$-invariant trilinear form play key roles. These results were
due to our work [X26].

In the above work, we fond a one-parameter ($\mfk{c}$) family of
inhomogeneous first-order differential operator representations of
the simple Lie algebra of type $E_7$ in $27$ variables.  Letting
these operators act on the space of exponential-polynomial functions
that depend on a parametric vector $\vec a\in \mbb
F^{27}\setminus\{\vec 0\}$, we prove that the space forms an
irreducible $E_7$-module for any constant $\mfk{c}$ if $\vec a$ is
not on an explicitly given projective algebraic variety. Certain
equivalent combinatorial properties of the basic oscillator
representation of $E_6$ play key roles in our proof. This part is
taken from our work [X28].

 \section{Basic Oscillator Representation of $E_7$}

In this section, we construct the oscillator representation of $E_7$
lifted from its 27-dimensional basic irreducible representation.

  First we go back to the construction of the simple
simple Lie algebra $\msr G^X$ in (4.4.15)-(4.4.25)
 with $X=E_6,\;E_7$ and $E_8$. Reacll the Dynkin diagram of $E_8$:

\begin{picture}(110,23)
\put(2,0){$E_8$:}\put(21,0){\circle{2}}\put(21,
-5){1}\put(22,0){\line(1,0){12}}\put(35,0){\circle{2}}\put(35,
-5){3}\put(36,0){\line(1,0){12}}\put(49,0){\circle{2}}\put(49,
-5){4}\put(49,1){\line(0,1){10}}\put(49,12){\circle{2}}\put(52,10){2}\put(50,0){\line(1,0){12}}
\put(63,0){\circle{2}}\put(63,-5){5}\put(64,0){\line(1,0){12}}\put(77,0){\circle{2}}\put(77,
-5){6}\put(78,0){\line(1,0){12}}\put(91,0){\circle{2}}\put(91,
-5){7}\put(92,0){\line(1,0){12}}\put(105,0){\circle{2}}\put(105,
-5){8}
\end{picture}
\vspace{0.7cm}

\noindent Let $\{\al_i\mid i\in\ol{1,8}\}$ be the simple positive
roots corresponding to the vertices in the diagram, and let
$\Phi_{E_8}$ be the root system of $E_8$. The simple Lie algebra of
type $E_8$ is
 $$\msr G^{E_8}=H\oplus\bigoplus_{\al\in\Phi_{E_8}}\mbb FE_{\al},\qquad H=H_{E_8}=\sum_{i=1}^8\mbb F\al_i,\eqno(12.1.1)$$
with the Lie bracket given in (4.4.24) and (4.4.25). Note that the
Dynkin diagram of $E_7$ is a sub-diagram of that of $E_8$. Set
$$H_{E_7}=\sum_{i=1}^7\mbb F\al_i,\qquad
\Phi_{E_7}=\Phi_{E_8}\bigcap H_{E_7}.\eqno(12.1.2)$$ We take the
simple Lie algebra $\msr G^{E_7}$ of type $E_7$ as the Lie
subalgebra
$$\msr G^{E_7}=H_{E_7}\oplus\bigoplus_{\al\in
\Phi_{E_7}}\mbb FE_{\al}.\eqno(12.1.3)$$

Recall the notion in (4.4.43). Denote
$$\varrho_1=E_{(2,3,4,6,5,4,3,1)},\quad \varrho_2=E_{(2,3,4,6,5,4,2,1)},\quad
\varrho_3=E_{(2,3,4,6,5,3,2,1)},\eqno(12.1.4)$$
$$\varrho_4=E_{(2,3,4,6,4,3,2,1)},\quad \varrho_5=E_{(2,3,4,5,4,3,2,1)},\quad
\varrho_6=E_{(2,3,3,5,4,3,2,1)},\eqno(12.1.5)$$
$$\varrho_7=E_{(2,2,4,5,4,3,2,1)},\quad \varrho_8=E_{(1,3,3,5,4,3,2,1)},\quad
\varrho_9=E_{(2,2,3,5,4,3,2,1)},\eqno(12.1.6)$$
$$\varrho_{10}=E_{(2,2,3,4,4,3,2,1)},\quad \varrho_{11}=E_{(1,2,3,5,4,3,2,1)},\quad
\varrho_{12}=E_{(2,2,3,4,3,3,2,1)}, \eqno(12.1.7)$$
$$ \varrho_{13}=E_{(1,2,3,4,4,3,2,1)},\quad \varrho_{14}=E_{(2,2,3,4,3,2,2,1)},\quad
\varrho_{15}=E_{(1,2,2,4,4,3,2,1)},\eqno(12.1.8)$$
$$\varrho_{16}=E_{(1,2,3,4,3,3,2,1)},\quad
\varrho_{17}=E_{(2,2,3,4,3,2,1,1)},\quad
\varrho_{18}=E_{(1,2,2,4,3,3,2,1)},\eqno(12.1.9)$$
$$\varrho_{19}=E_{(1,2,3,4,3,2,2,1)},\quad
\varrho_{20}=E_{(1,2,2,3,3,3,2,1)},\quad
\varrho_{21}=E_{(1,2,2,4,3,2,2,1)},\eqno(12.1.10)$$
$$\varrho_{22}=E_{(1,2,3,4,3,2,1,1)},
\quad \varrho_{23}=E_{(1,1,2,3,3,3,2,1)},\quad
\varrho_{24}=E_{(1,2,2,3,3,2,2,1)}, \eqno(12.1.11)$$
$$\varrho_{25}=E_{(1,2,2,4,3,2,1,1)},\quad
\varrho_{26}=E_{(1,1,2,3,3,2,2,1)},\quad
\varrho_{27}=E_{(1,2,2,3,2,2,2,1)}, \eqno(12.1.12)$$
$$\varrho_{28}=E_{(1,2,2,3,3,2,1,1)},\quad
\varrho_{29}=E_{(1,1,2,3,2,2,2,1)},\quad
\varrho_{30}=E_{(1,1,2,3,3,2,1,1)},
 \eqno(12.1.13)$$
$$\varrho_{31}=E_{(1,2,2,3,2,2,1,1)},\quad \varrho_{32}=E_{(1,1,2,2,2,2,2,1)},
\quad \varrho_{33}=E_{(1,1,2,3,2,2,1,1)}, \eqno(12.1.14)$$
$$ \varrho_{34}=E_{(1,2,2,3,2,1,1,1)},\quad
\varrho_{35}=E_{(1,1,1,2,2,2,2,1)},\quad
\varrho_{36}=E_{(1,1,2,2,2,2,1,1)},\eqno(12.1.15)$$
$$\varrho_{37}=E_{(1,1,2,3,2,1,1,1)},\quad \varrho_{38}=E_{(1,1,1,2,2,2,1,1)},\quad \varrho_{39}=E_{(1,1,2,2,2,1,1,1)},
\eqno(12.1.16)$$
$$\varrho_{40}=E_{(0,1,1,2,2,2,2,1)}, \quad \varrho_{41}=E_{(1,1,1,2,2,1,1,1)},
\quad \varrho_{42}=E_{(1,1,2,2,1,1,1,1)},\eqno(12.1.17)$$
$$\varrho_{43}=E_{(0,1,1,2,2,2,1,1)},\quad
\varrho_{44}=E_{(1,1,1,2,1,1,1,1)},\quad
\varrho_{45}=E_{(0,1,1,2,2,1,1,1)},\eqno(12.1.18)$$
$$\varrho_{46}=E_{(1,1,1,1,1,1,1,1)},\quad \varrho_{47}=E_{(0,1,1,2,1,1,1,1)},\quad
\varrho_{48}=E_{(0,1,1,1,1,1,1,1)},\eqno(12.1.19)$$
$$\varrho_{49}=E_{(1,0,1,1,1,1,1,1)},\quad
\varrho_{50}=E_{(0,1,0,1,1,1,1,1)},\quad
\varrho_{51}=E_{(0,0,1,1,1,1,1,1)},\eqno(12.1.20)$$
$$\varrho_{52}=E_{(0,0,0,1,1,1,1,1)},\quad
\varrho_{53}=E_{(0,0,0,0,1,1,1,1)},\quad
\varrho_{54}=E_{(0,0,0,0,0,1,1,1)},\eqno(12.1.21)$$
$$\varrho_{55}=E_{(0,0,0,0,0,0,1,1)},\quad
\varrho_{56}=E_{(0,0,0,0,0,0,0,1)}.\eqno(12.1.22)$$

Then the subspace
$$V=\sum_{i=1}^{56}\mbb F\varrho_i\eqno(12.1.23)$$ forms an irreducible $\msr G^{E_7}$-module
with respect to the adjoint representation of $\msr G^{E_8}$,
$\varrho_1$ is a highest-weight vector of weight $\lmd_7$. Set
$$\theta=\al(2,3,4,6,5,4,3,2)\eqno(12.1.24)$$ (cf. (4.4.42)). Then
$$[E_\theta,\msr G^{E_7}]=\{0\}\eqno(12.1.25)$$
and
$$\sum_{\al\in\Phi^+_{E_8}}\mbb{C}E_\al=\sum_{\be\in\Phi^+_{E_7}}\mbb{C}E_\be+V
+\mbb{C}E_\theta.\eqno(12.1.26)$$

Write
$$[u,\varrho_i]=\sum_{j=1}^{56}\vf_{i,j}(u)\varrho_j\qquad\for\;\;u\in\msr
G^{F_4}.\eqno(12.1.27)$$ Set
$$\msr B=\mbb{F}[x_1,...,x_{56}]\eqno(12.1.28)$$
and define the {\it basic oscillator representation of} $\msr
G^{E_7}$ on $\msr A$ by
$$u(g)=\sum_{i,j=1}^{56}\vf_{i,j}(u)x_j\ptl_{x_i}(g)\qquad\for\;\;u\in\msr
G^{E_7},\;g\in \msr B\eqno(12.1.29)$$ (cf. (2.2.17)-(2.2.20)). Then
$\msr B$ forms a $\msr G^{E_7}$-module isomorphic to the symmetric
tensor $S(V)$ over $V$. More explicitly, we have the following
representation formulas for the positive root vectors in $\msr
G^{E_7}$ by (4.4.54)-(4.4.58):
\begin{eqnarray*}E_{\al_1}|_{\msr B}&=&-x_6\ptl_{x_8}-x_9\ptl_{x_{11}}
-x_{10}\ptl_{x_{13}}-x_{12}\ptl_{x_{16}}-x_{14}\ptl_{x_{19}}-x_{17}\ptl_{x_{22}}
\\ &&+x_{35}\ptl_{x_{40}}+x_{38}\ptl_{x_{43}}+x_{41}\ptl_{x_{45}}
+x_{44}\ptl_{x_{47}}+x_{46}\ptl_{x_{48}}+x_{49}\ptl_{x_{51}},\hspace{1.9cm}(12.1.30)
\end{eqnarray*}
\begin{eqnarray*}E_{\al_2}|_{\msr B}&=&x_5\ptl_{x_7}+x_6\ptl_{x_9}
+x_8\ptl_{x_{11}}-x_{20}\ptl_{x_{23}}-x_{24}\ptl_{x_{26}}-x_{27}\ptl_{x_{29}}
\\ &&-x_{28}\ptl_{x_{30}}-x_{31}\ptl_{x_{33}}-x_{34}\ptl_{x_{37}}
+x_{46}\ptl_{x_{49}}+x_{48}\ptl_{x_{51}}+x_{50}\ptl_{x_{52}},\hspace{1.9cm}(12.1.31)
\end{eqnarray*}
\begin{eqnarray*}E_{\al_3}|_{\msr B}&=&-x_5\ptl_{x_6}-x_7\ptl_{x_9}
-x_{13}\ptl_{x_{15}}-x_{16}\ptl_{x_{18}}-x_{19}\ptl_{x_{21}}-x_{22}\ptl_{x_{25}}
\\ &&+x_{32}\ptl_{x_{35}}+x_{36}\ptl_{x_{38}}+x_{39}\ptl_{x_{41}}
+x_{42}\ptl_{x_{44}}+x_{48}\ptl_{x_{50}}+x_{51}\ptl_{x_{52}},
\hspace{1.9cm}(12.1.32)
\end{eqnarray*}
\begin{eqnarray*}E_{\al_4}|_{\msr B}&=&x_4\ptl_{x_5}-x_9\ptl_{x_{10}}
-x_{11}\ptl_{x_{13}}-x_{18}\ptl_{x_{20}}-x_{21}\ptl_{x_{24}}-x_{25}\ptl_{x_{28}}
\\ &&-x_{29}\ptl_{x_{32}}-x_{33}\ptl_{x_{36}}-x_{37}\ptl_{x_{39}}
-x_{44}\ptl_{x_{46}}-x_{47}\ptl_{x_{48}}+x_{52}\ptl_{x_{53}},
\hspace{1.9cm}(12.1.33)
\end{eqnarray*}
\begin{eqnarray*}E_{\al_5}|_{\msr B}&=&x_3\ptl_{x_4}-x_{10}\ptl_{x_{12}}
-x_{13}\ptl_{x_{16}}-x_{15}\ptl_{x_{18}}-x_{24}\ptl_{x_{27}}-x_{26}\ptl_{x_{29}}
\\ &&-x_{28}\ptl_{x_{31}}-x_{30}\ptl_{x_{33}}-x_{39}\ptl_{x_{42}}
-x_{41}\ptl_{x_{44}}-x_{45}\ptl_{x_{47}}+x_{53}\ptl_{x_{54}},
\hspace{1.9cm}(12.1.34)
\end{eqnarray*}
\begin{eqnarray*}E_{\al_6}|_{\msr B}&=&x_2\ptl_{x_3}-x_{12}\ptl_{x_{14}}
-x_{16}\ptl_{x_{19}}-x_{18}\ptl_{x_{21}}-x_{20}\ptl_{x_{24}}-x_{23}\ptl_{x_{26}}
\\ &&-x_{31}\ptl_{x_{34}}-x_{33}\ptl_{x_{37}}-x_{36}\ptl_{x_{39}}
-x_{38}\ptl_{x_{41}}-x_{43}\ptl_{x_{45}}+x_{54}\ptl_{x_{55}},
\hspace{1.9cm}(12.1.35)
\end{eqnarray*}
\begin{eqnarray*}E_{\al_7}|_{\msr B}&=&x_1\ptl_{x_2}-x_{14}\ptl_{x_{17}}
-x_{19}\ptl_{x_{22}}-x_{21}\ptl_{x_{25}}-x_{24}\ptl_{x_{28}}-x_{26}\ptl_{x_{30}}
\\ &&-x_{27}\ptl_{x_{31}}-x_{29}\ptl_{x_{33}}-x_{32}\ptl_{x_{36}}
-x_{35}\ptl_{x_{38}}-x_{40}\ptl_{x_{43}}+x_{55}\ptl_{x_{56}},
\hspace{1.9cm}(12.1.36)
\end{eqnarray*}
\begin{eqnarray*}E_{(1,0,1)}|_{\msr B}&=&-x_5\ptl_{x_8}-x_7\ptl_{x_{11}}
+x_{10}\ptl_{x_{15}}+x_{12}\ptl_{x_{18}}+x_{14}\ptl_{x_{21}}+x_{17}\ptl_{x_{25}}
\\ &&-x_{32}\ptl_{x_{40}}-x_{36}\ptl_{x_{43}}-x_{39}\ptl_{x_{45}}
-x_{42}\ptl_{x_{47}}+x_{46}\ptl_{x_{50}}+x_{49}\ptl_{x_{52}},
\hspace{1.4cm}(12.1.37)
\end{eqnarray*}
\begin{eqnarray*}E_{(0,1,0,1)}|_{\msr B}&=&-x_4\ptl_{x_7}-x_6\ptl_{x_{10}}
-x_8\ptl_{x_{13}}-x_{18}\ptl_{x_{23}}-x_{21}\ptl_{x_{26}}-x_{25}\ptl_{x_{30}}
\\ &&+x_{27}\ptl_{x_{32}}+x_{31}\ptl_{x_{36}}+x_{34}\ptl_{x_{39}}
+x_{44}\ptl_{x_{49}}+x_{47}\ptl_{x_{51}}+x_{50}\ptl_{x_{53}},
\hspace{1.1cm}(12.1.38)
\end{eqnarray*}
\begin{eqnarray*}E_{(0,0,1,1)}|_{\msr B}&=&x_4\ptl_{x_6}+x_7\ptl_{x_{10}}
-x_{11}\ptl_{x_{15}}+x_{16}\ptl_{x_{20}}+x_{19}\ptl_{x_{24}}+x_{22}\ptl_{x_{28}}
\\ &&+x_{29}\ptl_{x_{35}}+x_{33}\ptl_{x_{38}}+x_{37}\ptl_{x_{41}}
-x_{42}\ptl_{x_{46}}+x_{47}\ptl_{x_{50}}+x_{51}\ptl_{x_{53}},
\hspace{1.1cm}(12.1.39)
\end{eqnarray*}
\begin{eqnarray*}& &E_{(0,0,0,1,1)}|_{\msr B}=-x_3\ptl_{x_5}+x_9\ptl_{x_{12}}
+x_{11}\ptl_{x_{16}}-x_{15}\ptl_{x_{20}}+x_{21}\ptl_{x_{27}}+x_{25}\ptl_{x_{31}}
\\ &&-x_{26}\ptl_{x_{32}}-x_{30}\ptl_{x_{36}}+x_{37}\ptl_{x_{42}}
-x_{41}\ptl_{x_{46}}-x_{45}\ptl_{x_{48}}+x_{52}\ptl_{x_{54}},
\hspace{3.2cm}(12.1.40)
\end{eqnarray*}
\begin{eqnarray*}& &E_{(0,0,0,0,1,1)}|_{\msr B}=-x_2\ptl_{x_4}+x_{10}\ptl_{x_{14}}
+x_{13}\ptl_{x_{19}}+x_{15}\ptl_{x_{21}}-x_{20}\ptl_{x_{27}}-x_{23}\ptl_{x_{29}}
\\ &&+x_{28}\ptl_{x_{34}}+x_{30}\ptl_{x_{37}}-x_{36}\ptl_{x_{42}}
-x_{38}\ptl_{x_{44}}-x_{43}\ptl_{x_{47}}+x_{53}\ptl_{x_{55}},
\hspace{3.2cm}(12.1.41)
\end{eqnarray*}
\begin{eqnarray*}& &E_{(0,0,0,0,0,1,1)}|_{\msr B}=-x_1\ptl_{x_3}+x_{12}\ptl_{x_{17}}
+x_{16}\ptl_{x_{22}}+x_{18}\ptl_{x_{25}}+x_{20}\ptl_{x_{28}}+x_{23}\ptl_{x_{30}}
\\ &&-x_{27}\ptl_{x_{34}}-x_{29}\ptl_{x_{37}}-x_{32}\ptl_{x_{39}}
-x_{35}\ptl_{x_{41}}-x_{40}\ptl_{x_{45}}+x_{54}\ptl_{x_{56}},
\hspace{3.2cm}(12.1.42)
\end{eqnarray*}
\begin{eqnarray*}E_{(1,0,1,1)}|_{\msr B}&=&x_4\ptl_{x_8}+x_7\ptl_{x_{13}}
+x_9\ptl_{x_{15}}-x_{12}\ptl_{x_{20}}-x_{14}\ptl_{x_{24}}-x_{17}\ptl_{x_{28}}
\\ &&-x_{29}\ptl_{x_{40}}-x_{33}\ptl_{x_{43}}-x_{37}\ptl_{x_{45}}
+x_{42}\ptl_{x_{48}}+x_{44}\ptl_{x_{50}}+x_{49}\ptl_{x_{53}},
\hspace{1.1cm}(12.1.43)
\end{eqnarray*}
\begin{eqnarray*}E_{(0,1,1,1)}|_{\msr B}&=&-x_4\ptl_{x_9}+x_5\ptl_{x_{10}}
-x_8\ptl_{x_{15}}+x_{16}\ptl_{x_{23}}+x_{19}\ptl_{x_{26}}+x_{22}\ptl_{x_{30}}
\\ &&-x_{27}\ptl_{x_{35}}-x_{31}\ptl_{x_{38}}-x_{34}\ptl_{x_{41}}
+x_{42}\ptl_{x_{49}}-x_{47}\ptl_{x_{52}}+x_{48}\ptl_{x_{53}},
\hspace{1.1cm}(12.1.44)
\end{eqnarray*}
\begin{eqnarray*}& &E_{(0,1,0,1,1)}|_{\msr B}=x_3\ptl_{x_7}+x_6\ptl_{x_{12}}
+x_8\ptl_{x_{16}}-x_{15}\ptl_{x_{23}}+x_{21}\ptl_{x_{29}}+x_{24}\ptl_{x_{32}}
\\ &&+x_{25}\ptl_{x_{33}}+x_{28}\ptl_{x_{36}}-x_{34}\ptl_{x_{42}}
+x_{41}\ptl_{x_{49}}+x_{45}\ptl_{x_{51}}+x_{50}\ptl_{x_{54}},
\hspace{3.2cm}(12.1.45)
\end{eqnarray*}
\begin{eqnarray*}& &E_{(0,0,1,1,1)}|_{\msr B}=-x_3\ptl_{x_6}-x_7\ptl_{x_{12}}
+x_{11}\ptl_{x_{18}}+x_{13}\ptl_{x_{20}}-x_{19}\ptl_{x_{27}}-x_{22}\ptl_{x_{31}}
\\ &&+x_{26}\ptl_{x_{35}}+x_{30}\ptl_{x_{38}}-x_{37}\ptl_{x_{44}}
-x_{39}\ptl_{x_{46}}+x_{45}\ptl_{x_{50}}+x_{51}\ptl_{x_{54}},
\hspace{3.1cm}(12.1.46)
\end{eqnarray*}
\begin{eqnarray*}& &E_{(0,0,0,1,1,1)}|_{\msr B}=x_2\ptl_{x_5}-x_9\ptl_{x_{14}}
-x_{11}\ptl_{x_{19}}+x_{15}\ptl_{x_{24}}+x_{18}\ptl_{x_{28}}-x_{23}\ptl_{x_{32}}
\\ &&-x_{25}\ptl_{x_{34}}+x_{30}\ptl_{x_{39}}+x_{33}\ptl_{x_{42}}
-x_{38}\ptl_{x_{46}}-x_{43}\ptl_{x_{48}}+x_{52}\ptl_{x_{55}},
\hspace{3.2cm}(12.1.47)
\end{eqnarray*}
\begin{eqnarray*}& &E_{(0,0,0,0,1,1,1)}|_{\msr B}=x_1\ptl_{x_4}-x_{10}\ptl_{x_{17}}
-x_{13}\ptl_{x_{22}}-x_{15}\ptl_{x_{25}}+x_{20}\ptl_{x_{31}}+x_{23}\ptl_{x_{33}}
\\ &&+x_{24}\ptl_{x_{34}}+x_{26}\ptl_{x_{37}}-x_{32}\ptl_{x_{42}}
-x_{35}\ptl_{x_{44}}-x_{40}\ptl_{x_{47}}+x_{53}\ptl_{x_{56}},
\hspace{3.2cm}(12.1.48)
\end{eqnarray*}
\begin{eqnarray*}E_{(1,1,1,1)}|_{\msr B}&=&-x_4\ptl_{x_{11}}+x_5\ptl_{x_{13}}
+x_6\ptl_{x_{15}}-x_{12}\ptl_{x_{23}}-x_{14}\ptl_{x_{26}}-x_{17}\ptl_{x_{30}}
\\ &&+x_{27}\ptl_{x_{40}}+x_{31}\ptl_{x_{43}}+x_{34}\ptl_{x_{45}}
-x_{42}\ptl_{x_{51}}-x_{44}\ptl_{x_{52}}+x_{46}\ptl_{x_{53}},
\hspace{1.1cm}(12.1.49)
\end{eqnarray*}
\begin{eqnarray*}& &E_{(1,0,1,1,1)}|_{\msr B}=-x_3\ptl_{x_8}-x_7\ptl_{x_{16}}
-x_9\ptl_{x_{18}}-x_{10}\ptl_{x_{20}}+x_{14}\ptl_{x_{27}}+x_{17}\ptl_{x_{31}}
\\ &&-x_{26}\ptl_{x_{40}}-x_{30}\ptl_{x_{43}}+x_{37}\ptl_{x_{47}}+x_{39}\ptl_{x_{48}}
+x_{41}\ptl_{x_{50}}+x_{49}\ptl_{x_{54}}, \hspace{3.2cm}(12.1.50)
\end{eqnarray*}
\begin{eqnarray*}& &E_{(0,1,1,1,1)}|_{\msr B}=x_3\ptl_{x_9}-x_5\ptl_{x_{12}}
+x_8\ptl_{x_{18}}+x_{13}\ptl_{x_{23}}-x_{19}\ptl_{x_{29}}-x_{22}\ptl_{x_{33}}
\\ &&-x_{24}\ptl_{x_{35}}-x_{28}\ptl_{x_{38}}+x_{34}\ptl_{x_{44}}
+x_{39}\ptl_{x_{49}}-x_{45}\ptl_{x_{52}}+x_{48}\ptl_{x_{54}},
\hspace{3.2cm}(12.1.51)
\end{eqnarray*}
\begin{eqnarray*}& &E_{(0,1,0,1,1,1)}|_{\msr B}=-x_2\ptl_{x_7}-x_6\ptl_{x_{14}}
-x_8\ptl_{x_{19}}+x_{15}\ptl_{x_{26}}+x_{18}\ptl_{x_{30}}+x_{20}\ptl_{x_{32}}
\\ &&-x_{25}\ptl_{x_{37}}-x_{28}\ptl_{x_{39}}-x_{31}\ptl_{x_{42}}
+x_{38}\ptl_{x_{49}}+x_{43}\ptl_{x_{51}}+x_{50}\ptl_{x_{55}},
\hspace{3.2cm}(12.1.52)
\end{eqnarray*}
\begin{eqnarray*}& &E_{(0,0,1,1,1,1)}|_{\msr B}=x_2\ptl_{x_6}+x_7\ptl_{x_{14}}
-x_{11}\ptl_{x_{21}}-x_{13}\ptl_{x_{24}}-x_{16}\ptl_{x_{28}}+x_{23}\ptl_{x_{35}}
\\ &&+x_{22}\ptl_{x_{34}}-x_{30}\ptl_{x_{41}}-x_{33}\ptl_{x_{44}}
-x_{36}\ptl_{x_{46}}+x_{43}\ptl_{x_{50}}+x_{51}\ptl_{x_{55}},
\hspace{3.2cm}(12.1.53)
\end{eqnarray*}
\begin{eqnarray*}& &E_{(0,0,0,1,1,1,1)}|_{\msr B}=-x_1\ptl_{x_5}+x_9\ptl_{x_{17}}
+x_{11}\ptl_{x_{22}}-x_{15}\ptl_{x_{28}}-x_{18}\ptl_{x_{31}}-x_{21}\ptl_{x_{34}}
\\ &&+x_{23}\ptl_{x_{36}}+x_{26}\ptl_{x_{39}}+x_{29}\ptl_{x_{42}}
-x_{35}\ptl_{x_{46}}-x_{40}\ptl_{x_{48}}+x_{52}\ptl_{x_{56}},
\hspace{3.2cm}(12.1.54)
\end{eqnarray*}
\begin{eqnarray*}& &E_{(1,1,1,1,1)}|_{\msr B}=x_3\ptl_{x_{11}}-x_5\ptl_{x_{16}}
-x_6\ptl_{x_{18}}-x_{10}\ptl_{x_{23}}+x_{14}\ptl_{x_{29}}+x_{17}\ptl_{x_{33}}
\\ &&+x_{24}\ptl_{x_{40}}+x_{28}\ptl_{x_{43}}-x_{34}\ptl_{x_{47}}
-x_{39}\ptl_{x_{51}}-x_{41}\ptl_{x_{52}}+x_{46}\ptl_{x_{54}},
\hspace{3.2cm}(12.1.55)
\end{eqnarray*}
\begin{eqnarray*}& &E_{(1,0,1,1,1,1)}|_{\msr B}=x_2\ptl_{x_8}+x_7\ptl_{x_{19}}
+x_9\ptl_{x_{21}}+x_{10}\ptl_{x_{24}}+x_{12}\ptl_{x_{28}}-x_{17}\ptl_{x_{34}}
\\ &&-x_{23}\ptl_{x_{40}}+x_{30}\ptl_{x_{45}}+x_{33}\ptl_{x_{47}}
+x_{36}\ptl_{x_{48}}+x_{38}\ptl_{x_{50}}+x_{49}\ptl_{x_{55}},
\hspace{3.2cm}(12.1.56)
\end{eqnarray*}
\begin{eqnarray*}& &E_{(0,1,1,2,1)}|_{\msr B}=-x_3\ptl_{x_{10}}+x_4\ptl_{x_{12}}
-x_8\ptl_{x_{20}}+x_{11}\ptl_{x_{23}}+x_{19}\ptl_{x_{32}}-x_{21}\ptl_{x_{35}}
\\ &&+x_{22}\ptl_{x_{36}}-x_{25}\ptl_{x_{38}}-x_{34}\ptl_{x_{46}}
+x_{37}\ptl_{x_{49}}-x_{45}\ptl_{x_{53}}+x_{47}\ptl_{x_{54}},
\hspace{3.2cm}(12.1.57)
\end{eqnarray*}
\begin{eqnarray*}& &E_{(0,1,1,1,1,1)}|_{\msr B}=-x_2\ptl_{x_9}+x_5\ptl_{x_{14}}
-x_8\ptl_{x_{21}}-x_{13}\ptl_{x_{26}}-x_{16}\ptl_{x_{30}}-x_{20}\ptl_{x_{35}}
\\ &&+x_{22}\ptl_{x_{37}}+x_{28}\ptl_{x_{41}}+x_{31}\ptl_{x_{44}}
+x_{36}\ptl_{x_{49}}-x_{43}\ptl_{x_{52}}+x_{48}\ptl_{x_{55}},
\hspace{3.2cm}(12.1.58)
\end{eqnarray*}
\begin{eqnarray*}& &E_{(0,1,0,1,1,1,1)}|_{\msr B}=x_1\ptl_{x_7}+x_6\ptl_{x_{17}}
+x_8\ptl_{x_{22}}-x_{15}\ptl_{x_{30}}-x_{18}\ptl_{x_{33}}-x_{21}\ptl_{x_{37}}
\\ &&-x_{20}\ptl_{x_{36}}-x_{24}\ptl_{x_{39}}-x_{27}\ptl_{x_{42}}
+x_{35}\ptl_{x_{49}}+x_{40}\ptl_{x_{51}}+x_{50}\ptl_{x_{56}},
\hspace{3.2cm}(12.1.59)
\end{eqnarray*}
\begin{eqnarray*}& &E_{(0,0,1,1,1,1,1)}|_{\msr B}=-x_1\ptl_{x_6}-x_7\ptl_{x_{17}}
+x_{11}\ptl_{x_{25}}+x_{13}\ptl_{x_{28}}+x_{16}\ptl_{x_{31}}+x_{19}\ptl_{x_{34}}
\\ &&-x_{23}\ptl_{x_{38}}-x_{26}\ptl_{x_{41}}-x_{29}\ptl_{x_{44}}
-x_{32}\ptl_{x_{46}}+x_{40}\ptl_{x_{50}}+x_{51}\ptl_{x_{56}},
\hspace{3.2cm}(12.1.60)
\end{eqnarray*}
\begin{eqnarray*}& &E_{(1,1,1,2,1)}|_{\msr B}=-x_3\ptl_{x_{13}}+x_4\ptl_{x_{16}}
+x_6\ptl_{x_{20}}-x_9\ptl_{x_{23}}-x_{14}\ptl_{x_{32}}+x_{21}\ptl_{x_{40}}
\\ &&-x_{17}\ptl_{x_{36}}+x_{25}\ptl_{x_{43}}+x_{34}\ptl_{x_{48}}
-x_{37}\ptl_{x_{51}}-x_{41}\ptl_{x_{53}}+x_{44}\ptl_{x_{54}},
\hspace{3.2cm}(12.1.61)
\end{eqnarray*}
\begin{eqnarray*}& &E_{(1,1,1,1,1,1)}|_{\msr B}=-x_2\ptl_{x_{11}}+x_5\ptl_{x_{19}}
+x_6\ptl_{x_{21}}+x_{10}\ptl_{x_{26}}+x_{12}\ptl_{x_{29}}-x_{17}\ptl_{x_{37}}
\\ &&+x_{20}\ptl_{x_{40}}-x_{28}\ptl_{x_{45}}-x_{31}\ptl_{x_{47}}
-x_{36}\ptl_{x_{51}}-x_{38}\ptl_{x_{52}}+x_{46}\ptl_{x_{55}},
\hspace{3.2cm}(12.1.62)
\end{eqnarray*}
\begin{eqnarray*}& &E_{(1,0,1,1,1,1,1)}|_{\msr B}=-x_1\ptl_{x_8}-x_7\ptl_{x_{22}}
-x_9\ptl_{x_{25}}-x_{10}\ptl_{x_{28}}-x_{12}\ptl_{x_{31}}-x_{14}\ptl_{x_{34}}
\\ &&+x_{23}\ptl_{x_{43}}+x_{26}\ptl_{x_{45}}+x_{29}\ptl_{x_{47}}
+x_{32}\ptl_{x_{48}}+x_{35}\ptl_{x_{50}}+x_{49}\ptl_{x_{56}},
\hspace{3.2cm}(12.1.63)
\end{eqnarray*}
\begin{eqnarray*}& &E_{(0,1,1,2,1,1)}|_{\msr B}=x_2\ptl_{x_{10}}-x_4\ptl_{x_{14}}
+x_8\ptl_{x_{24}}-x_{11}\ptl_{x_{26}}+x_{16}\ptl_{x_{32}}-x_{18}\ptl_{x_{35}}
\\ &&-x_{22}\ptl_{x_{39}}+x_{25}\ptl_{x_{41}}-x_{31}\ptl_{x_{46}}
+x_{33}\ptl_{x_{49}}-x_{43}\ptl_{x_{53}}+x_{47}\ptl_{x_{55}},
\hspace{3.2cm}(12.1.64)
\end{eqnarray*}
\begin{eqnarray*}& &E_{(0,1,1,1,1,1,1)}|_{\msr B}=x_1\ptl_{x_9}-x_5\ptl_{x_{17}}
+x_8\ptl_{x_{25}}+x_{13}\ptl_{x_{30}}+x_{16}\ptl_{x_{33}}+x_{19}\ptl_{x_{37}}
\\ &&+x_{20}\ptl_{x_{38}}+x_{24}\ptl_{x_{41}}+x_{27}\ptl_{x_{44}}
+x_{32}\ptl_{x_{49}}-x_{40}\ptl_{x_{52}}+x_{48}\ptl_{x_{56}},
\hspace{3.2cm}(12.1.65)
\end{eqnarray*}
\begin{eqnarray*}& &E_{(1,1,2,2,1)}|_{\msr B}=-x_3\ptl_{x_{15}}+x_4\ptl_{x_{18}}
-x_5\ptl_{x_{20}}+x_7\ptl_{x_{23}}+x_{14}\ptl_{x_{35}}-x_{19}\ptl_{x_{40}}
\\ &&+x_{17}\ptl_{x_{38}}-x_{22}\ptl_{x_{43}}-x_{34}\ptl_{x_{50}}
+x_{37}\ptl_{x_{52}}-x_{39}\ptl_{x_{53}}+x_{42}\ptl_{x_{54}},
\hspace{3.2cm}(12.1.66)
\end{eqnarray*}
\begin{eqnarray*}& &E_{(1,1,1,2,1,1)}|_{\msr B}=x_2\ptl_{x_{13}}-x_4\ptl_{x_{19}}
-x_6\ptl_{x_{24}}+x_9\ptl_{x_{26}}-x_{12}\ptl_{x_{32}}+x_{18}\ptl_{x_{40}}
\\ &&+x_{17}\ptl_{x_{39}}-x_{25}\ptl_{x_{45}}+x_{31}\ptl_{x_{48}}
-x_{33}\ptl_{x_{51}}-x_{38}\ptl_{x_{53}}+x_{44}\ptl_{x_{55}},
\hspace{3.2cm}(12.1.67)
\end{eqnarray*}
\begin{eqnarray*}& &E_{(1,1,1,1,1,1,1)}|_{\msr B}=x_1\ptl_{x_{11}}-x_5\ptl_{x_{22}}
-x_6\ptl_{x_{25}}-x_{10}\ptl_{x_{30}}-x_{12}\ptl_{x_{33}}-x_{14}\ptl_{x_{37}}
\\ &&-x_{20}\ptl_{x_{43}}-x_{24}\ptl_{x_{45}}-x_{27}\ptl_{x_{47}}
-x_{32}\ptl_{x_{51}}-x_{35}\ptl_{x_{52}}+x_{46}\ptl_{x_{56}},
\hspace{3.2cm}(12.1.68)
\end{eqnarray*}
\begin{eqnarray*}& &E_{(0,1,1,2,2,1)}|_{\msr B}=-x_2\ptl_{x_{12}}+x_3\ptl_{x_{14}}
-x_8\ptl_{x_{27}}+x_{11}\ptl_{x_{29}}+x_{13}\ptl_{x_{32}}-x_{15}\ptl_{x_{35}}
\\ &&+x_{22}\ptl_{x_{42}}-x_{25}\ptl_{x_{44}}-x_{28}\ptl_{x_{46}}
+x_{30}\ptl_{x_{49}}-x_{43}\ptl_{x_{54}}+x_{45}\ptl_{x_{55}},
\hspace{3.2cm}(12.1.69)
\end{eqnarray*}
\begin{eqnarray*}& &E_{(0,1,1,2,1,1,1)}|_{\msr B}=-x_1\ptl_{x_{10}}+x_4\ptl_{x_{17}}
-x_8\ptl_{x_{28}}+x_{11}\ptl_{x_{30}}-x_{16}\ptl_{x_{36}}+x_{18}\ptl_{x_{38}}
\\ &&-x_{19}\ptl_{x_{39}}+x_{21}\ptl_{x_{41}}-x_{27}\ptl_{x_{46}}
+x_{29}\ptl_{x_{49}}-x_{40}\ptl_{x_{53}}+x_{47}\ptl_{x_{56}},
\hspace{3.2cm}(12.1.70)
\end{eqnarray*}
\begin{eqnarray*}& &E_{(1,1,2,2,1,1)}|_{\msr B}=x_2\ptl_{x_{15}}-x_4\ptl_{x_{21}}
+x_5\ptl_{x_{24}}-x_7\ptl_{x_{26}}+x_{12}\ptl_{x_{35}}-x_{16}\ptl_{x_{40}}
\\ &&-x_{17}\ptl_{x_{41}}+x_{22}\ptl_{x_{45}}-x_{31}\ptl_{x_{50}}
+x_{33}\ptl_{x_{52}}-x_{36}\ptl_{x_{53}}+x_{42}\ptl_{x_{55}},
\hspace{3.2cm}(12.1.71)
\end{eqnarray*}
\begin{eqnarray*}& &E_{(1,1,1,2,2,1)}|_{\msr B}=-x_2\ptl_{x_{16}}+x_3\ptl_{x_{19}}
+x_6\ptl_{x_{27}}-x_9\ptl_{x_{29}}-x_{10}\ptl_{x_{32}}+x_{15}\ptl_{x_{40}}
\\ &&-x_{17}\ptl_{x_{42}}+x_{25}\ptl_{x_{47}}+x_{28}\ptl_{x_{48}}
-x_{30}\ptl_{x_{51}}-x_{38}\ptl_{x_{54}}+x_{41}\ptl_{x_{55}},
\hspace{3.2cm}(12.1.72)
\end{eqnarray*}
\begin{eqnarray*}& &E_{(1,1,1,2,1,1,1)}|_{\msr B}=-x_1\ptl_{x_{13}}+x_4\ptl_{x_{22}}
+x_6\ptl_{x_{28}}-x_9\ptl_{x_{30}}+x_{12}\ptl_{x_{36}}+x_{14}\ptl_{x_{39}}
\\ &&-x_{18}\ptl_{x_{43}}-x_{21}\ptl_{x_{45}}+x_{27}\ptl_{x_{48}}
-x_{29}\ptl_{x_{51}}-x_{35}\ptl_{x_{53}}+x_{44}\ptl_{x_{56}},
\hspace{3.2cm}(12.1.73)
\end{eqnarray*}
\begin{eqnarray*}& &E_{(0,1,1,2,2,1,1)}|_{\msr B}=x_1\ptl_{x_{12}}-x_3\ptl_{x_{17}}
+x_8\ptl_{x_{31}}-x_{11}\ptl_{x_{33}}-x_{13}\ptl_{x_{36}}+x_{15}\ptl_{x_{38}}
\\ &&+x_{19}\ptl_{x_{42}}-x_{21}\ptl_{x_{44}}-x_{24}\ptl_{x_{46}}
+x_{26}\ptl_{x_{49}}-x_{40}\ptl_{x_{54}}+x_{45}\ptl_{x_{56}},
\hspace{3.2cm}(12.1.74)
\end{eqnarray*}
\begin{eqnarray*}& &E_{(1,1,2,2,2,1)}|_{\msr B}=-x_2\ptl_{x_{18}}+x_3\ptl_{x_{21}}
-x_5\ptl_{x_{27}}+x_7\ptl_{x_{29}}+x_{10}\ptl_{x_{35}}-x_{13}\ptl_{x_{40}}
\\ &&+x_{17}\ptl_{x_{44}}-x_{22}\ptl_{x_{47}}-x_{28}\ptl_{x_{50}}
+x_{30}\ptl_{x_{52}}-x_{36}\ptl_{x_{54}}+x_{39}\ptl_{x_{55}},
\hspace{3.2cm}(12.1.75)
\end{eqnarray*}
\begin{eqnarray*}& &E_{(1,1,2,2,1,1,1)}|_{\msr B}=-x_1\ptl_{x_{15}}+x_4\ptl_{x_{25}}
-x_5\ptl_{x_{28}}+x_7\ptl_{x_{30}}-x_{12}\ptl_{x_{38}}+x_{16}\ptl_{x_{43}}
\\ &&-x_{14}\ptl_{x_{41}}+x_{19}\ptl_{x_{45}}-x_{27}\ptl_{x_{50}}
+x_{29}\ptl_{x_{52}}-x_{32}\ptl_{x_{53}}+x_{42}\ptl_{x_{56}},
\hspace{3.2cm}(12.1.76)
\end{eqnarray*}
\begin{eqnarray*}& &E_{(1,1,1,2,2,1,1)}|_{\msr B}=x_1\ptl_{x_{16}}-x_3\ptl_{x_{22}}
-x_6\ptl_{x_{31}}+x_9\ptl_{x_{33}}+x_{10}\ptl_{x_{36}}-x_{15}\ptl_{x_{43}}
\\ &&-x_{14}\ptl_{x_{42}}+x_{21}\ptl_{x_{47}}+x_{24}\ptl_{x_{48}}
-x_{26}\ptl_{x_{51}}-x_{35}\ptl_{x_{54}}+x_{41}\ptl_{x_{56}},
\hspace{3.2cm}(12.1.77)
\end{eqnarray*}
\begin{eqnarray*}& &E_{(0,1,1,2,2,2,1)}|_{\msr B}=-x_1\ptl_{x_{14}}+x_2\ptl_{x_{17}}
-x_8\ptl_{x_{34}}+x_{11}\ptl_{x_{37}}+x_{13}\ptl_{x_{39}}-x_{15}\ptl_{x_{41}}
\\ &&+x_{16}\ptl_{x_{42}}-x_{18}\ptl_{x_{44}}-x_{20}\ptl_{x_{46}}
+x_{23}\ptl_{x_{49}}-x_{40}\ptl_{x_{55}}+x_{43}\ptl_{x_{56}},
\hspace{3.2cm}(12.1.78)
\end{eqnarray*}
\begin{eqnarray*}& &E_{(1,1,2,3,2,1)}|_{\msr B}=x_2\ptl_{x_{20}}-x_3\ptl_{x_{24}}
+x_4\ptl_{x_{27}}-x_7\ptl_{x_{32}}+x_9\ptl_{x_{35}}-x_{11}\ptl_{x_{40}}
\\ &&-x_{17}\ptl_{x_{46}}+x_{22}\ptl_{x_{48}}-x_{25}\ptl_{x_{50}}
+x_{30}\ptl_{x_{53}}-x_{33}\ptl_{x_{54}}+x_{37}\ptl_{x_{55}},
\hspace{3.1cm}(12.1.79)
\end{eqnarray*}
\begin{eqnarray*}& &E_{(1,1,2,2,2,1,1)}|_{\msr B}=x_1\ptl_{x_{18}}-x_3\ptl_{x_{25}}
+x_5\ptl_{x_{31}}-x_7\ptl_{x_{33}}-x_{10}\ptl_{x_{38}}+x_{13}\ptl_{x_{43}}
\\ &&+x_{14}\ptl_{x_{44}}-x_{19}\ptl_{x_{47}}-x_{24}\ptl_{x_{50}}
+x_{26}\ptl_{x_{52}}-x_{32}\ptl_{x_{54}}+x_{39}\ptl_{x_{56}},
\hspace{3.2cm}(12.1.80)
\end{eqnarray*}
\begin{eqnarray*}&
&E_{(1,1,1,2,2,2,1)}|_{\msr B}=-x_1\ptl_{x_{19}}+x_2\ptl_{x_{22}}
+x_6\ptl_{x_{34}}-x_9\ptl_{x_{37}}-x_{10}\ptl_{x_{39}}+x_{15}\ptl_{x_{45}}
\\ &&-x_{12}\ptl_{x_{42}}+x_{18}\ptl_{x_{47}}+x_{20}\ptl_{x_{48}}
-x_{23}\ptl_{x_{51}}-x_{35}\ptl_{x_{55}}+x_{38}\ptl_{x_{56}},
\hspace{3.2cm}(12.1.81)
\end{eqnarray*}
\begin{eqnarray*}& &E_{(1,2,2,3,2,1)}|_{\msr B}=-x_2\ptl_{x_{23}}+x_3\ptl_{x_{26}}
-x_4\ptl_{x_{29}}+x_5\ptl_{x_{32}}-x_6\ptl_{x_{35}}+x_8\ptl_{x_{40}}
\\ &&-x_{17}\ptl_{x_{49}}+x_{22}\ptl_{x_{51}}-x_{25}\ptl_{x_{52}}
+x_{28}\ptl_{x_{53}}-x_{31}\ptl_{x_{54}}+x_{34}\ptl_{x_{55}},
\hspace{3.3cm}(12.1.82)
\end{eqnarray*}
\begin{eqnarray*}& &E_{(1,1,2,3,2,1,1)}|_{\msr B}=-x_1\ptl_{x_{20}}+x_3\ptl_{x_{28}}
-x_4\ptl_{x_{31}}+x_7\ptl_{x_{36}}-x_9\ptl_{x_{38}}+x_{11}\ptl_{x_{43}}
\\ &&-x_{14}\ptl_{x_{46}}+x_{19}\ptl_{x_{48}}-x_{21}\ptl_{x_{50}}
+x_{26}\ptl_{x_{53}}-x_{29}\ptl_{x_{54}}+x_{37}\ptl_{x_{56}},
\hspace{3.3cm}(12.1.83)
\end{eqnarray*}
\begin{eqnarray*}& &E_{(1,1,2,2,2,2,1)}|_{\msr B}=-x_1\ptl_{x_{21}}+x_2\ptl_{x_{25}}
-x_5\ptl_{x_{34}}+x_7\ptl_{x_{37}}+x_{10}\ptl_{x_{41}}-x_{13}\ptl_{x_{45}}
\\ &&+x_{12}\ptl_{x_{44}}-x_{16}\ptl_{x_{47}}-x_{20}\ptl_{x_{50}}
+x_{23}\ptl_{x_{52}}-x_{32}\ptl_{x_{55}}+x_{36}\ptl_{x_{56}},
\hspace{3.3cm}(12.1.84)
\end{eqnarray*}
\begin{eqnarray*}& &E_{(1,2,2,3,2,1,1)}|_{\msr B}=x_1\ptl_{x_{23}}-x_3\ptl_{x_{30}}
+x_4\ptl_{x_{33}}-x_5\ptl_{x_{36}}+x_6\ptl_{x_{38}}-x_8\ptl_{x_{43}}
\\ &&-x_{14}\ptl_{x_{49}}+x_{19}\ptl_{x_{51}}-x_{21}\ptl_{x_{52}}
+x_{24}\ptl_{x_{53}}-x_{27}\ptl_{x_{54}}+x_{34}\ptl_{x_{56}},
\hspace{3.3cm}(12.1.85)
\end{eqnarray*}
\begin{eqnarray*}& &E_{(1,1,2,3,2,2,1)}|_{\msr B}=x_1\ptl_{x_{24}}-x_2\ptl_{x_{28}}
+x_4\ptl_{x_{34}}-x_7\ptl_{x_{39}}+x_9\ptl_{x_{41}}-x_{11}\ptl_{x_{45}}
\\ &&-x_{12}\ptl_{x_{46}}+x_{16}\ptl_{x_{48}}-x_{18}\ptl_{x_{50}}
+x_{23}\ptl_{x_{53}}-x_{29}\ptl_{x_{55}}+x_{33}\ptl_{x_{56}},
\hspace{3.3cm}(12.1.86)
\end{eqnarray*}
\begin{eqnarray*}& &E_{(1,2,2,3,2,2,1)}|_{\msr B}=-x_1\ptl_{x_{26}}+x_2\ptl_{x_{30}}
-x_4\ptl_{x_{37}}+x_5\ptl_{x_{39}}-x_6\ptl_{x_{41}}+x_8\ptl_{x_{45}}
\\ &&-x_{12}\ptl_{x_{49}}+x_{16}\ptl_{x_{51}}-x_{18}\ptl_{x_{52}}
+x_{20}\ptl_{x_{53}}-x_{27}\ptl_{x_{55}}+x_{31}\ptl_{x_{56}},
\hspace{3.3cm}(12.1.87)
\end{eqnarray*}
\begin{eqnarray*}& &E_{(1,1,2,3,3,2,1)}|_{\msr B}=-x_1\ptl_{x_{27}}+x_2\ptl_{x_{31}}
-x_3\ptl_{x_{34}}+x_7\ptl_{x_{42}}-x_9\ptl_{x_{44}}+x_{11}\ptl_{x_{47}}
\\ &&-x_{10}\ptl_{x_{46}}+x_{13}\ptl_{x_{48}}-x_{15}\ptl_{x_{50}}
+x_{23}\ptl_{x_{54}}-x_{26}\ptl_{x_{55}}+x_{30}\ptl_{x_{56}},
\hspace{3.2cm}(12.1.88)
\end{eqnarray*}
\begin{eqnarray*}&
&E_{(1,2,2,3,3,2,1)}|_{\msr B}=x_1\ptl_{x_{29}}-x_2\ptl_{x_{33}}
+x_3\ptl_{x_{37}}-x_5\ptl_{x_{42}}+x_6\ptl_{x_{44}}-x_8\ptl_{x_{47}}
\\ &&-x_{10}\ptl_{x_{49}}+x_{13}\ptl_{x_{51}}-x_{15}\ptl_{x_{52}}
+x_{20}\ptl_{x_{54}}-x_{24}\ptl_{x_{55}}+x_{28}\ptl_{x_{56}},
\hspace{3.1cm}(12.1.89)
\end{eqnarray*}
\begin{eqnarray*}&
&E_{(1,2,2,4,3,2,1)}|_{\msr B}=-x_1\ptl_{x_{32}}+x_2\ptl_{x_{36}}
-x_3\ptl_{x_{39}}+x_4\ptl_{x_{42}}-x_6\ptl_{x_{46}}+x_8\ptl_{x_{48}}
\\ &&-x_9\ptl_{x_{49}}+x_{11}\ptl_{x_{51}}-x_{15}\ptl_{x_{53}}
+x_{18}\ptl_{x_{54}}-x_{21}\ptl_{x_{55}}+x_{25}\ptl_{x_{56}},
\hspace{3.2cm}(12.1.90)
\end{eqnarray*}
\begin{eqnarray*}&
&E_{(1,2,3,4,3,2,1)}|_{\msr B}=-x_1\ptl_{x_{35}}+x_2\ptl_{x_{38}}
-x_3\ptl_{x_{41}}+x_4\ptl_{x_{44}}-x_5\ptl_{x_{46}}+x_8\ptl_{x_{50}}
\\ &&-x_7\ptl_{x_{49}}+x_{11}\ptl_{x_{52}}-x_{13}\ptl_{x_{53}}
+x_{16}\ptl_{x_{54}}-x_{19}\ptl_{x_{55}}+x_{22}\ptl_{x_{56}},
\hspace{3.2cm}(12.1.91)
\end{eqnarray*}
\begin{eqnarray*}&
&E_{(2,2,3,4,3,2,1)}|_{\msr B}=-x_1\ptl_{x_{40}}+x_2\ptl_{x_{43}}
-x_3\ptl_{x_{45}}+x_4\ptl_{x_{47}}-x_5\ptl_{x_{48}}+x_6\ptl_{x_{50}}
\\ &&-x_7\ptl_{x_{51}}+x_9\ptl_{x_{52}}-x_{10}\ptl_{x_{53}}
+x_{12}\ptl_{x_{54}}-x_{14}\ptl_{x_{55}}+x_{17}\ptl_{x_{56}}.
\hspace{3.5cm}(12.1.92)
\end{eqnarray*}

We define a symmetric linear operation $\tau$ on the space
$\sum_{i,j=1}^{56}\mbb{F}x_i\ptl_{x_j}$ by
$$\tau(x_i\ptl_{x_j})=x_j\ptl_{x_i}.\eqno(12.1.93)$$
Denote
$$\bar r=57-r\qquad\for\;\;r\in\ol{1,56}.\eqno(11.1.94)$$
 Then
$$E_{-\al}|_{\msr B}=-\tau(E_{\al}|_{\msr B})\qquad\for\;\;\al\in\Phi_{E_7}^+\eqno(12.1.95)$$
by (4.4.25). Moreover, (4.4.24) gives
$$\al_r|_V=\sum_{i=1}^{28}a_{i,r}(x_i\ptl_{x_i}-x_{\bar i}\ptl_{x_{\bar
i}})\qquad\for\;\;r\in\ol{1,7},\eqno(12.1.96)$$ where $a_{i,r}$ are
constants given by the following table:
\begin{center}{\bf \large Table 12.1.1}\end{center}
\begin{center}\begin{tabular}{|r||r|r|r|r|r|r|r||r||r|r|r|r|r|r|r|}\hline
$i$&$a_{i,1}$&$a_{i,2}$&$a_{i,3}$&$a_{i,4}$&$a_{i,5}$&$a_{i,6}$&$a_{i,7}$&$i$&$a_{i,1}$&$a_{i,2}$&$a_{i,3}$&$a_{i,4}$&$a_{i,5}$&
$a_{i,6}$&$a_{i,7}$
\\\hline\hline 1&0&0&0&0&0&0&1&2&0&0&0&0&0&1&$-1$\\\hline
3&0&0&0&0&1&$-1$&0&4&$0$&0&0&1&$-1$&0&0
\\\hline 5&0&1&1&$-1$&0&0&0&6&1&1&$-1$&0&0&0&0 \\\hline
7&0&$-1$&$1$&0&0&0&0&8&$-1$&1&0&$0$&0&0&0\\\hline
9&1&$-1$&$-1$&1&0&0&0&10&1&0&0&$-1$&1&0&0\\\hline
11&$-1$&$-1$&0&1&$0$&0&0&12&1&0&0&0&$-1$&1&0\\\hline
13&$-1$&0&1&$-1$&1&$0$&0&14&1&0&0&0&0&$-1$&1\\\hline
15&0&0&$-1$&0&1&0&0&16&$-1$&0&1&0&$-1$&1&0\\\hline
17&1&0&0&0&0&0&$-1$&18&0&0&$-1$&1&$-1$&1&0\\\hline
19&$-1$&0&1&0&0&$-1$&1&20&0&1&0&$-1$&0&1&0\\\hline
21&0&0&$-1$&1&0&$-1$&1&22&$-1$&0&1&0&0&0&$-1$\\\hline
23&0&$-1$&0&0&0&1&0&24&0&1&0&$-1$&1&$-1$&1\\\hline
25&0&0&$-1$&1&0&0&$-1$&26&0&$-1$&0&0&1&$-1$&1\\\hline
27&0&1&0&0&$-1$&0&1&28&0&1
&0&$-1$&1&0&$-1$\\\hline\end{tabular}\end{center}

\section{Constructions of Singular Vectors}

In this section, we construct special singular vectors in $\msr B$.

According to Table 12.1.1, we find the singular vector
$$\mfk{l}_1=x_1x_{17}+x_2x_{14}+x_3x_{12}+x_4x_{10}+x_5x_9-x_6x_7,
\eqno(12.2.1)$$ which generates an irreducible ${\msr
G}^{E_7}$-module $W$ isomorphic to the adjoint module. Set
$$\mbb{B}={\msr
B}[\ptl_{x_1},...,\ptl_{x_{56}}].\eqno(12.2.2)$$ In terms of
(11.1.94), we define an associative algebra isomorphism $\nu$ on
$\mbb{B}$ by
$$\nu(x_r)=x_{\bar r},\;\nu(\ptl_{x_r})=\ptl_{x_{\bar
r}}\qquad\for\;\; r\in \ol{1,56}.\eqno(12.2.3)$$
 Set
$$\mfk{l}_2=x_1x_{22}+x_2x_{19}+x_3x_{16}+x_4x_{13}+x_5x_{11}-x_7x_8,
\eqno(12.2.4)$$
$$\mfk{l}_3=x_1x_{25}+x_2x_{21}+x_3x_{18}+x_4x_{15}+x_6x_{11}-x_8x_9,
\eqno(12.2.5)$$
$$\mfk{l}_4=-x_1x_{28}-x_2x_{24}-x_3x_{20}+x_5x_{15}-x_6x_{13}+x_8x_{10},
\eqno(12.2.6)$$
$$\mfk{l}_5=x_1x_{30}+x_2x_{26}+x_3x_{23}+x_7x_{15}-x_9x_{13}+x_{10}x_{11},
\eqno(12.2.7)$$
$$\mfk{l}_6=x_1x_{31}+x_2x_{27}-x_4x_{20}-x_5x_{18}+x_6x_{16}-x_8x_{12},
\eqno(12.2.8)$$
$$\mfk{l}_7=-x_1x_{33}-x_2x_{29}+x_4x_{23}-x_7x_{18}+x_9x_{16}-x_{11}x_{12},
\eqno(12.2.9)$$
$$\mfk{l}_8=-x_1x_{34}+x_3x_{27}+x_4x_{24}+x_5x_{21}-x_6x_{19}+x_8x_{14},
\eqno(12.2.10)$$
$$\mfk{l}_9=-x_2x_{34}-x_3x_{31}-x_4x_{28}-x_5x_{25}+x_6x_{22}-x_8x_{17},
\eqno(12.2.11)$$
$$\mfk{l}_{10}=x_1x_{36}+x_2x_{32}+x_5x_{23}+x_7x_{20}-x_{10}x_{16}+x_{12}x_{13},
\eqno(12.2.12)$$
$$\mfk{l}_{11}=x_1x_{37}-x_3x_{29}-x_4x_{26}+x_7x_{21}-x_9x_{19}+x_{11}x_{14},
\eqno(12.2.13)$$
$$\mfk{l}_{12}=-x_1x_{38}-x_2x_{35}+x_6x_{23}+x_9x_{20}-x_{10}x_{18}+x_{12}x_{15},
\eqno(12.2.14)$$
$$\mfk{l}_{13}=x_2x_{37}+x_3x_{33}+x_4x_{30}-x_7x_{25}+x_9x_{22}-x_{11}x_{17},
\eqno(12.2.15)$$
$$\mfk{l}_{14}=-x_1x_{39}+x_3x_{32}-x_5x_{26}-x_7x_{24}+x_{10}x_{19}-x_{13}x_{14},
\eqno(12.2.16)$$
$$\mfk{l}_{15}=-x_2x_{39}-x_3x_{36}+x_5x_{30}+x_7x_{28}-x_{10}x_{22}+x_{13}x_{17},
\eqno(12.2.17)$$
$$\mfk{l}_{16}=x_1x_{41}-x_3x_{35}-x_6x_{26}-x_9x_{24}+x_{10}x_{21}-x_{14}x_{15},
\eqno(12.2.18)$$
$$\mfk{l}_{17}=x_1x_{42}+x_4x_{32}+x_5x_{29}+x_7x_{27}-x_{12}x_{19}+x_{14}x_{16},
\eqno(12.2.19)$$
$$\mfk{l}_{18}=x_2x_{41}+x_3x_{38}+x_6x_{30}+x_9x_{28}-x_{10}x_{25}+x_{15}x_{17},
\eqno(12.2.20)$$
$$\mfk{l}_{19}=x_1x_{43}+x_2x_{40}+x_8x_{23}+x_{11}x_{20}-x_{13}x_{18}+x_{15}x_{16},
\eqno(12.2.21)$$
$$\mfk{l}_{20}=x_2x_{42}-x_4x_{36}-x_5x_{33}-x_7x_{31}+x_{12}x_{22}-x_{16}x_{17},
\eqno(12.2.22)$$
$$\mfk{l}_{21}=-x_1x_{44}-x_4x_{35}+x_6x_{29}+x_9x_{27}-x_{12}x_{21}+x_{14}x_{18},
\eqno(12.2.23)$$
$$\mfk{l}_{22}=x_3x_{42}+x_4x_{39}+x_5x_{37}+x_7x_{34}-x_{14}x_{22}+x_{17}x_{19},
\eqno(12.2.24)$$
$$\mfk{l}_{23}=-x_1x_{45}+x_3x_{40}-x_8x_{26}-x_{11}x_{24}+x_{13}x_{21}-x_{15}x_{19},
\eqno(12.2.25)$$
$$\mfk{l}_{24}=-x_2x_{44}+x_4x_{38}-x_6x_{33}-x_9x_{31}+x_{12}x_{25}-x_{17}x_{18},
\eqno(12.2.26)$$
$$\mfk{l}_{25}=x_1x_{46}-x_5x_{35}-x_6x_{32}-x_{10}x_{27}+x_{12}x_{24}-x_{14}x_{20},
\eqno(12.2.27)$$
$$\mfk{l}_{26}=-x_2x_{45}-x_3x_{43}+x_8x_{30}+x_{11}x_{28}-x_{13}x_{25}+x_{15}x_{22},
\eqno(12.2.28)$$
$$\mfk{l}_{27}=-x_3x_{44}-x_4x_{41}+x_6x_{37}+x_9x_{34}-x_{14}x_{25}+x_{17}x_{21},
\eqno(12.2.29)$$
$$\mfk{l}_{28}=x_1x_{47}+x_4x_{40}+x_8x_{29}+x_{11}x_{27}-x_{16}x_{21}+x_{18}x_{19},
\eqno(12.2.30)$$
$$\mfk{l}_{29}=x_2x_{46}+x_5x_{38}+x_6x_{36}+x_{10}x_{31}-x_{12}x_{28}+x_{17}x_{20},
\eqno(12.2.31)$$
$$\mfk{l}_{30}=-x_1x_{48}+x_5x_{40}-x_8x_{32}-x_{13}x_{27}+x_{16}x_{24}-x_{19}x_{20},
\eqno(12.2.32)$$
$$\mfk{l}_{31}=x_2x_{47}-x_4x_{43}-x_8x_{33}-x_{11}x_{31}+x_{16}x_{25}-x_{18}x_{22},
\eqno(12.2.33)$$
$$\mfk{l}_{32}=x_3x_{46}-x_5x_{41}-x_6x_{39}-x_{10}x_{34}+x_{14}x_{28}-x_{17}x_{24},
\eqno(12.2.34)$$
$$\mfk{l}_{33}=x_1x_{49}-x_7x_{35}-x_9x_{32}+x_{10}x_{29}-x_{12}x_{26}+x_{14}x_{23},
\eqno(12.2.35)$$
$$\mfk{l}_{34}=-x_2x_{48}-x_5x_{43}+x_8x_{36}+x_{13}x_{31}-x_{16}x_{28}+x_{20}x_{22},
\eqno(12.2.36)$$
$$\mfk{l}_{35}=x_3x_{47}+x_4x_{45}+x_8x_{37}+x_{11}x_{34}-x_{19}x_{25}+x_{21}x_{22},
\eqno(12.2.37)$$
$$\mfk{l}_{36}=x_4x_{46}+x_5x_{44}+x_6x_{42}+x_{12}x_{34}-x_{14}x_{31}+x_{17}x_{27},
\eqno(12.2.38)$$
$$\mfk{l}_{37}=x_1x_{50}+x_6x_{40}+x_8x_{35}-x_{15}x_{27}+x_{18}x_{24}-x_{20}x_{21},
\eqno(12.2.39)$$
$$\mfk{l}_{38}=x_2x_{49}+x_7x_{38}+x_9x_{36}-x_{10}x_{33}+x_{12}x_{30}-x_{17}x_{23},
\eqno(12.2.40)$$
$$\mfk{l}_{39}=-x_3x_{48}+x_5x_{45}-x_8x_{39}-x_{13}x_{34}+x_{19}x_{28}-x_{22}x_{24},
\eqno(12.2.41)$$
$$\mfk{l}_{40}=-x_1x_{51}+x_7x_{40}-x_{11}x_{32}+x_{13}x_{29}-x_{16}x_{26}+x_{19}x_{23},
\eqno(12.2.42)$$
$$\mfk{l}_{41}=x_2x_{50}-x_6x_{43}-x_8x_{38}+x_{15}x_{31}-x_{18}x_{28}+x_{20}x_{25},
\eqno(12.2.43)$$
$$\mfk{l}_{42}=x_3x_{49}-x_7x_{41}-x_9x_{39}+x_{10}x_{37}-x_{14}x_{30}+x_{17}x_{26},
\eqno(12.2.44)$$
$$\mfk{l}_{43}=-x_4x_{48}-x_5x_{47}+x_8x_{42}+x_{16}x_{34}-x_{19}x_{31}+x_{22}x_{27},
\eqno(12.2.45)$$
$$\mfk{l}_{44}=x_1x_{52}+x_9x_{40}+x_{11}x_{35}+x_{15}x_{29}-x_{18}x_{26}+x_{21}x_{23},
\eqno(12.2.46)$$
$$\mfk{l}_{45}=-x_2x_{51}-x_7x_{43}+x_{11}x_{36}-x_{13}x_{33}+x_{16}x_{30}-x_{22}x_{23},
\eqno(12.2.47)$$
$$\mfk{l}_{46}=x_3x_{50}+x_6x_{45}+x_8x_{41}-x_{15}x_{34}+x_{21}x_{28}-x_{24}x_{25},
\eqno(12.2.48)$$
$$\mfk{l}_{47}=x_4x_{49}+x_7x_{44}+x_9x_{42}-x_{12}x_{37}+x_{14}x_{33}-x_{17}x_{29},
\eqno(12.2.49)$$
$$\mfk{l}_{48}=x_1x_{53}-x_{10}x_{40}-x_{13}x_{35}-x_{15}x_{32}+x_{20}x_{26}-x_{23}x_{24},
\eqno(12.2.50)$$
$$\mfk{l}_{49}=x_2x_{52}-x_9x_{43}-x_{11}x_{38}-x_{15}x_{33}+x_{18}x_{30}-x_{23}x_{25},
\eqno(12.2.51)$$
$$\mfk{l}_{50}=-x_3x_{51}+x_7x_{45}-x_{11}x_{39}+x_{13}x_{37}-x_{19}x_{30}+x_{22}x_{26},
\eqno(12.2.52)$$
$$\mfk{l}_{51}=x_4x_{50}-x_6x_{47}-x_8x_{44}+x_{18}x_{34}-x_{21}x_{31}+x_{25}x_{27},
\eqno(12.2.53)$$
$$\mfk{l}_{52}=x_5x_{49}-x_7x_{46}-x_{10}x_{42}+x_{12}x_{39}-x_{14}x_{36}+x_{17}x_{32},
\eqno(12.2.54)$$
$$\mfk{l}_{53}=x_1x_{54}+x_{12}x_{40}+x_{16}x_{35}+x_{18}x_{32}-x_{20}x_{29}+x_{23}x_{27},
\eqno(12.2.55)$$
$$\mfk{l}_{54}=x_2x_{53}+x_{10}x_{43}+x_{13}x_{38}+x_{15}x_{36}-x_{20}x_{30}+x_{23}x_{28},
\eqno(12.2.56)$$
$$\mfk{l}_{55}=x_3x_{52}+x_9x_{45}+x_{11}x_{41}+x_{15}x_{37}-x_{21}x_{30}+x_{25}x_{26},
\eqno(12.2.57)$$
$$\mfk{l}_{56}=-x_4x_{51}-x_7x_{47}+x_{11}x_{42}-x_{16}x_{37}+x_{19}x_{33}-x_{22}x_{29},
\eqno(12.2.58)$$
$$\mfk{l}_{57}=x_5x_{50}+x_6x_{48}+x_8x_{46}-x_{20}x_{34}+x_{24}x_{31}-x_{27}x_{28},
\eqno(12.2.59)$$
$$\mfk{l}_{58}=x_6x_{49}-x_9x_{46}+x_{10}x_{44}-x_{12}x_{41}+x_{14}x_{38}-x_{17}x_{35},
\eqno(12.2.60)$$
$$\mfk{l}_{59}=x_1x_{55}-x_{14}x_{40}-x_{19}x_{35}-x_{21}x_{32}+x_{24}x_{29}-x_{26}x_{27},
\eqno(12.2.61)$$
$$\mfk{l}_{60}=x_2x_{54}-x_{12}x_{43}-x_{16}x_{38}-x_{18}x_{36}+x_{20}x_{33}-x_{23}x_{31},
\eqno(12.2.62)$$
$$\mfk{l}_{61}=x_3x_{53}-x_{10}x_{45}-x_{13}x_{41}-x_{15}x_{39}+x_{24}x_{30}-x_{26}x_{28},
\eqno(12.2.63)$$
$$\mfk{l}_{62}=x_4x_{52}-x_9x_{47}-x_{11}x_{44}-x_{18}x_{37}+x_{21}x_{33}-x_{25}x_{29},
\eqno(12.2.64)$$
$$\mfk{l}_{63}=-x_5x_{51}+x_7x_{48}-x_{13}x_{42}+x_{16}x_{39}-x_{19}x_{36}+x_{22}x_{32},
\eqno(12.2.65)$$
\begin{eqnarray*}\hspace{1.5cm}\mfk{l}_{64}&=&x_1x_{56}+x_2x_{55}+x_{14}x_{43}+x_{17}x_{40}+x_{19}x_{38}
+x_{22}x_{35}\\ &
&+x_{21}x_{36}+x_{25}x_{32}-x_{24}x_{33}-x_{28}x_{29}+x_{26}x_{31}+x_{27}x_{30},
\hspace{1.8cm}(12.2.66)
\end{eqnarray*}
\begin{eqnarray*}\hspace{1.5cm}\mfk{l}_{65}&=&x_2x_{55}+x_3x_{54}+x_{12}x_{45}+x_{14}x_{43}+x_{16}x_{41}+x_{19}x_{38}
\\ & &+x_{18}x_{39}+x_{21}x_{36}-x_{20}x_{37}-x_{24}x_{33}+x_{23}x_{34}+x_{26}x_{31},\hspace{1.8cm}(12.2.67)\end{eqnarray*}
\begin{eqnarray*}\hspace{1.5cm}\mfk{l}_{66}&=&x_3x_{54}+x_4x_{53}+x_{10}x_{47}+x_{12}x_{45}+x_{13}x_{44}
+x_{15}x_{42}\\&&+x_{16}x_{41}+x_{18}x_{39}-x_{24}x_{33}-x_{27}x_{30}+x_{26}x_{31}+x_{28}x_{29},
\hspace{1.8cm}(12.2.68)\end{eqnarray*}
\begin{eqnarray*}\hspace{1.5cm}\mfk{l}_{67}&=&x_4x_{53}+x_5x_{52}+x_9x_{48}+x_{10}x_{47}+x_{11}x_{46}+x_{13}x_{44}\\
&&+x_{18}x_{39}+x_{20}x_{37}-x_{21}x_{36}-x_{24}x_{33}+x_{25}x_{32}+x_{28}x_{29},
\hspace{1.8cm}(12.2.69)\end{eqnarray*}
\begin{eqnarray*}\hspace{1.5cm}\mfk{l}_{68}&=&x_5x_{52}+x_6x_{51}+x_7x_{50}+x_8x_{49}+x_9x_{48}+x_{11}x_{46}\\
&&+x_{20}x_{37}+x_{23}x_{34}-x_{24}x_{33}-x_{26}x_{31}+x_{27}x_{30}+x_{28}x_{29},\hspace{1.8cm}(12.2.70)\end{eqnarray*}
\begin{eqnarray*}\hspace{1.5cm}\mfk{l}_{69}&=&x_5x_{52}-x_6x_{51}-x_7x_{50}+x_9x_{48}+x_{13}x_{44}-x_{15}x_{42}\\
& &-x_{16}x_{41}+x_{18}x_{39}+x_{19}x_{38}-x_{21}x_{36}-x_{22}x_{35}
 +x_{25}x_{32},\hspace{1.8cm}(12.2.71)\end{eqnarray*}
\begin{eqnarray*}\hspace{1.5cm}\mfk{l}_{70}&=&-x_6x_{51}+x_8x_{49}+x_9x_{48}-x_{10}x_{47}-x_{11}x_{46}+x_{13}x_{44}
\\ & &+x_{12}x_{45}-x_{14}x_{43}-x_{16}x_{41}+x_{17}x_{40}+x_{19}x_{38}-x_{22}x_{35},\hspace{1.8cm}(12.2.72)\end{eqnarray*}
\newpage

$$\mfk{l}_i=\nu(\mfk{l}_{134-i})\qquad\for\;\;i\in\ol{71,133}.\eqno(12.2.73)$$
Then
$$W=\sum_{i=1}^{133}\mbb F \mfk{l}_i.\eqno(12.2.74)$$
 Moreover,
$$\al_r|_W=\sum_{i=1}^{63}b_{i,r}(\mfk{l}_i\ptl_{\mfk{l}_i}-\mfk{l}_{134-i}\ptl_{\mfk{l}_{134-i}})
\qquad\for\;\;r\in\ol{1,7},\eqno(12.2.75)$$ where $b_{i,r}$ are
given by the following table.
\begin{center}{\bf \large Table 12.2.1}\end{center}
\begin{center}\begin{tabular}{|r||r|r|r|r|r|r|r||r||r|r|r|r|r|r|r|}\hline
$i$&$b_{i,1}$&$b_{i,2}$&$b_{i,3}$&$b_{i,4}$&$b_{i,5}$&$b_{i,6}$&$b_{i,7}$&$i$&$b_{i,1}$&$b_{i,2}$&$b_{i,3}$&$b_{i,4}$&$b_{i,5}$&
$b_{i,6}$&$b_{i,7}$
\\\hline\hline 1&1&0&0&0&0&0&0&2&$-1$&0&1&0&0&0&0\\\hline
3&0&0&$-1$&1&0&0&0&4&0&1 &0&$-1$&1&0&0
\\\hline 5&0&$-1$&0&0&1&0&0&6&0&1&0&0&$-1$&1&0 \\\hline
7&0&$-1$&0&1&$-1$&1&0&8&0&1&0&0&0&$-1$&1\\\hline
9&0&1&0&0&0&0&$-1$&10&0&0&1&$-1$&0&1&0\\\hline
11&0&$-1$&0&1&0&$-1$&1&12&1&0&$-1$&0&0&1&0\\\hline
13&0&$-1$&0&1&0&0&$-1$&14&0&0&1&$-1$&1&$-1$&1\\\hline
15&0&0&1&$-1$&1&0&$-1$&16&1&0&$-1$&0&1&$-1$&1\\\hline
17&0&0&1&0&$-1$&0&1&18&1&0&$-1$&0&1&0&$-1$\\\hline
19&$-1$&0&0&0&0&1&0&20&0&0&1&0&$-1$&1&$-1$\\\hline
21&1&0&$-1$&1&$-1$&0&1&22&0&0&1&0&0&$-1$&0\\\hline
23&$-1$&0&0&0&1&$-1$&1&24&1&0&$-1$&1&$-1$&1&$-1$\\\hline
25&1&1&0&$-1$&0&0&1&26&$-1$&0&0&0&1&0&$-1$\\\hline
27&1&0&$-1$&1&0&$-1$&0&28&$-1$&0&0&1&$-1$&0&1
\\\hline 29&1&1&0&$-1$&0&1&$-1$&30&$-1$&1&1&$-1$&0&0&1\\\hline
31&$-1$&0&0&1&$-1$&1&$-1$&32&1&1&0&$-1$&1&$-1$&0
\\\hline 33&1&$-1$&0&$0$&0&0&1&34&$-1$&1&1&$-1$&0&1&$-1$ \\\hline
35&$-1$&0&0&1&0&$-1$&0&36&1&1&0&0&$-1$&0&0\\\hline
37&0&1&$-1$&0&0&0&1&38&1&$-1$&0&$0$&0&1&$-1$\\\hline
39&$-1$&1&1&$-1$&1&$-1$&0&40&$-1$&$-1$&1&0&0&0&1\\\hline
41&0&1&$-1$&0&0&1&$-1$&42&1&$-1$&0&$0$&1&$-1$&0\\\hline
43&$-1$&1&1&0&$-1$&0&0&44&0&$-1$&$-1$&1&0&0&1\\\hline
45&$-1$&$-1$&1&0&0&1&$-1$&46&0&1&$-1$&0&1&$-1$&0\\\hline
47&1&$-1$&0&1&$-1$&0&0&48&0&0&0&$-1$&1&0&1\\\hline
49&0&$-1$&$-1$&1&0&1&$-1$&50&$-1$&$-1$&1&0&1&$-1$&0\\\hline
51&0&1&$-1$&1&$-1$&0&0&52&1&0&1&$-1$&0&0&0\\\hline
53&0&0&0&0&$-1$&1&1&54&$0$&0&0&$-1$&1&1&$-1$\\\hline
55&0&$-1$&$-1$&1&1&$-1$&0&56&$-1$&$-1$&1&1&$-1$&0&0
\\\hline 57&0&2&0&$-1$&0&0&0&58&2&0&$-1$&0&0&0&0\\\hline
59&0&0&0&0&0&$-1$&2&60&0&0&0&0&$-1$&2&$-1$\\\hline
61&0&0&0&$-1$&2&$-1$&0&62&0&$-1$&$-1$&2&$-1$&0&0\\\hline
63&$-1$&0&2&$-1$&0&0&0&64&0&0&0&0&0&0&0\\\hline
\end{tabular}\end{center}

 Note for $\al,\be\in \Phi_{E_7}$,
$$\al+\be\in \Phi_{E_7}\Leftrightarrow (\al,\be)=-1.\eqno(12.2.76)$$
Thus
$$E_{\al_r}(\mfk{l}_i)\neq 0 \Leftrightarrow
b_{i,r}=-1\qquad\for\;\;r\in\ol{1,7},\;i\in\ol{1,63}.\eqno(12.2.77)$$
Write
$$W_1=\sum_{i=1}^{63}\mbb{C}\mfk{l}_i.\eqno(12.2.78)$$
Based on the above fact and Table 12.2.1, we find
\begin{eqnarray*}E_{\al_1}|_{W_1}&=&-\mfk{l}_1\ptl_{\mfk{l}_2}-\mfk{l}_{12}\ptl_{\mfk{l}_{19}}
-\mfk{l}_{16}\ptl_{\mfk{l}_{23}}-\mfk{l}_{18}\ptl_{\mfk{l}_{26}}-\mfk{l}_{21}\ptl_{\mfk{l}_{28}}-\mfk{l}_{25}\ptl_{\mfk{l}_{30}}
\\&&-\mfk{l}_{24}\ptl_{\mfk{l}_{31}}-\mfk{l}_{29}\ptl_{\mfk{l}_{34}}-\mfk{l}_{27}\ptl_{\mfk{l}_{35}}
-\mfk{l}_{32}\ptl_{\mfk{l}_{39}}-\mfk{l}_{33}\ptl_{\mfk{l}_{40}}
\\ &&-\mfk{l}_{36}\ptl_{\mfk{l}_{43}}-\mfk{l}_{38}\ptl_{\mfk{l}_{45}}-\mfk{l}_{42}\ptl_{\mfk{l}_{50}}
-\mfk{l}_{47}\ptl_{\mfk{l}_{56}}-\mfk{l}_{52}\ptl_{\mfk{l}_{63}},\hspace{4.35cm}(12.2.79)
\end{eqnarray*}
\begin{eqnarray*}E_{\al_2}|_{W_1}&=&\mfk{l}_4\ptl_{\mfk{l}_5}+\mfk{l}_6\ptl_{\mfk{l}_7}
+\mfk{l}_8\ptl_{\mfk{l}_{11}}+\mfk{l}_9\ptl_{\mfk{l}_{13}}+\mfk{l}_{25}\ptl_{\mfk{l}_{33}}+\mfk{l}_{29}\ptl_{\mfk{l}_{38}}
\\&&+\mfk{l}_{30}\ptl_{\mfk{l}_{40}}+\mfk{l}_{32}\ptl_{\mfk{l}_{42}}+\mfk{l}_{37}\ptl_{\mfk{l}_{44}}
+\mfk{l}_{34}\ptl_{\mfk{l}_{45}}+\mfk{l}_{36}\ptl_{\mfk{l}_{47}}
\\ &&+\mfk{l}_{41}\ptl_{\mfk{l}_{49}}+\mfk{l}_{39}\ptl_{\mfk{l}_{50}}+\mfk{l}_{46}\ptl_{\mfk{l}_{55}}
+\mfk{l}_{43}\ptl_{\mfk{l}_{56}}+\mfk{l}_{51}\ptl_{\mfk{l}_{62}},\hspace{4.35cm}(12.2.80)
\end{eqnarray*}
\begin{eqnarray*}E_{\al_3}|_{W_1}&=&-\mfk{l}_2\ptl_{\mfk{l}_3}-\mfk{l}_{10}\ptl_{\mfk{l}_{12}}
-\mfk{l}_{14}\ptl_{\mfk{l}_{16}}-\mfk{l}_{15}\ptl_{\mfk{l}_{18}}-\mfk{l}_{17}\ptl_{\mfk{l}_{21}}-\mfk{l}_{20}\ptl_{\mfk{l}_{24}}
\\&&-\mfk{l}_{22}\ptl_{\mfk{l}_{27}}-\mfk{l}_{30}\ptl_{\mfk{l}_{37}}-\mfk{l}_{34}\ptl_{\mfk{l}_{41}}
-\mfk{l}_{40}\ptl_{\mfk{l}_{44}}-\mfk{l}_{39}\ptl_{\mfk{l}_{46}}
\\ &&-\mfk{l}_{45}\ptl_{\mfk{l}_{49}}-\mfk{l}_{43}\ptl_{\mfk{l}_{51}}-\mfk{l}_{50}\ptl_{\mfk{l}_{55}}
-\mfk{l}_{52}\ptl_{\mfk{l}_{58}}-\mfk{l}_{56}\ptl_{\mfk{l}_{62}},\hspace{4.35cm}(12.2.81)
\end{eqnarray*}
\begin{eqnarray*}E_{\al_4}|_{W_1}&=&\mfk{l}_3\ptl_{\mfk{l}_4}+\mfk{l}_7\ptl_{\mfk{l}_{10}}
+\mfk{l}_{11}\ptl_{\mfk{l}_{14}}+\mfk{l}_{13}\ptl_{\mfk{l}_{15}}+\mfk{l}_{21}\ptl_{\mfk{l}_{25}}+\mfk{l}_{24}\ptl_{\mfk{l}_{29}}
\\&&+\mfk{l}_{28}\ptl_{\mfk{l}_{30}}+\mfk{l}_{27}\ptl_{\mfk{l}_{32}}+\mfk{l}_{31}\ptl_{\mfk{l}_{34}}
+\mfk{l}_{35}\ptl_{\mfk{l}_{39}}+\mfk{l}_{44}\ptl_{\mfk{l}_{48}}
\\ &&+\mfk{l}_{47}\ptl_{\mfk{l}_{52}}+\mfk{l}_{49}\ptl_{\mfk{l}_{54}}+\mfk{l}_{51}\ptl_{\mfk{l}_{57}}
+\mfk{l}_{55}\ptl_{\mfk{l}_{61}}+\mfk{l}_{56}\ptl_{\mfk{l}_{63}},\hspace{4.35cm}(12.2.82)
\end{eqnarray*}
\begin{eqnarray*}E_{\al_5}|_{W_1}&=&\mfk{l}_4\ptl_{\mfk{l}_6}+\mfk{l}_5\ptl_{\mfk{l}_7}
+\mfk{l}_{14}\ptl_{\mfk{l}_{17}}+\mfk{l}_{15}\ptl_{\mfk{l}_{20}}+\mfk{l}_{16}\ptl_{\mfk{l}_{21}}+\mfk{l}_{18}\ptl_{\mfk{l}_{24}}
\\&&+\mfk{l}_{23}\ptl_{\mfk{l}_{28}}+\mfk{l}_{26}\ptl_{\mfk{l}_{31}}+\mfk{l}_{32}\ptl_{\mfk{l}_{36}}
+\mfk{l}_{39}\ptl_{\mfk{l}_{43}}+\mfk{l}_{42}\ptl_{\mfk{l}_{47}}
\\ &&+\mfk{l}_{46}\ptl_{\mfk{l}_{51}}+\mfk{l}_{48}\ptl_{\mfk{l}_{53}}+\mfk{l}_{50}\ptl_{\mfk{l}_{56}}
+\mfk{l}_{54}\ptl_{\mfk{l}_{60}}+\mfk{l}_{55}\ptl_{\mfk{l}_{62}},\hspace{4.35cm}(12.2.83)
\end{eqnarray*}
\begin{eqnarray*}E_{\al_6}|_{W_1}&=&\mfk{l}_6\ptl_{\mfk{l}_8}+\mfk{l}_7\ptl_{\mfk{l}_{11}}
+\mfk{l}_{10}\ptl_{\mfk{l}_{14}}+\mfk{l}_{12}\ptl_{\mfk{l}_{16}}+\mfk{l}_{20}\ptl_{\mfk{l}_{22}}+\mfk{l}_{19}\ptl_{\mfk{l}_{23}}
\\&&+\mfk{l}_{24}\ptl_{\mfk{l}_{27}}+\mfk{l}_{29}\ptl_{\mfk{l}_{32}}+\mfk{l}_{31}\ptl_{\mfk{l}_{35}}
+\mfk{l}_{34}\ptl_{\mfk{l}_{39}}+\mfk{l}_{38}\ptl_{\mfk{l}_{42}}
\\ &&+\mfk{l}_{41}\ptl_{\mfk{l}_{46}}+\mfk{l}_{45}\ptl_{\mfk{l}_{50}}+\mfk{l}_{49}\ptl_{\mfk{l}_{55}}
+\mfk{l}_{53}\ptl_{\mfk{l}_{59}}+\mfk{l}_{54}\ptl_{\mfk{l}_{61}},\hspace{4.35cm}(12.2.84)
\end{eqnarray*}
\begin{eqnarray*}E_{\al_7}|_{W_1}&=&\mfk{l}_8\ptl_{\mfk{l}_9}+\mfk{l}_{11}\ptl_{\mfk{l}_{13}}
+\mfk{l}_{14}\ptl_{\mfk{l}_{15}}+\mfk{l}_{16}\ptl_{\mfk{l}_{18}}+\mfk{l}_{17}\ptl_{\mfk{l}_{20}}+\mfk{l}_{21}\ptl_{\mfk{l}_{24}}
\\&&+\mfk{l}_{23}\ptl_{\mfk{l}_{26}}+\mfk{l}_{25}\ptl_{\mfk{l}_{29}}+\mfk{l}_{28}\ptl_{\mfk{l}_{31}}
+\mfk{l}_{30}\ptl_{\mfk{l}_{34}}+\mfk{l}_{33}\ptl_{\mfk{l}_{38}}
\\ &&+\mfk{l}_{37}\ptl_{\mfk{l}_{41}}+\mfk{l}_{40}\ptl_{\mfk{l}_{45}}+\mfk{l}_{44}\ptl_{\mfk{l}_{49}}
+\mfk{l}_{48}\ptl_{\mfk{l}_{54}}+\mfk{l}_{53}\ptl_{\mfk{l}_{60}}.\hspace{4.35cm}(12.2.85)
\end{eqnarray*}

We define a linear transformation $\mfk{s}$ on the space
$$\sum_{i,j=1}^{133}\mbb{C}\mfk{l}_i\ptl_{\mfk{l}_j}\eqno(12.2.86)$$
by
$$\mfk{s}(\mfk{l}_i\ptl_{\mfk{l}_j})=\mfk{l}_{134-j}\ptl_{\mfk{l}_{134-i}}.\eqno(12.2.87)$$
By symmetry, we have
$$E_{\al_1}|_W=(1-\mfk{s})(E_{\al_1}|_{W_1})+\mfk{l}_{70}\ptl_{\mfk{l}_{76}}
-\mfk{l}_{58}(2\ptl_{\mfk{l}_{70}}+\ptl_{\mfk{l}_{69}}),
\eqno(12.2.88)$$
$$E_{\al_2}|_W=(1+\mfk{s})(E_{\al_2}|_{W_1})+\mfk{l}_{68}\ptl_{\mfk{l}_{77}}+\mfk{l}_{57}(2\ptl_{\mfk{l}_{68}}+\ptl_{\mfk{l}_{67}}),
\eqno(12.2.89)$$
$$E_{\al_3}|_W=(1-\mfk{s})(E_{\al_3}|_{W_1})+\mfk{l}_{69}\ptl_{\mfk{l}_{71}}-
\mfk{l}_{63}(\ptl_{\mfk{l}_{70}}+2\ptl_{\mfk{l}_{69}}+\ptl_{\mfk{l}_{67}}),\eqno(12.2.90)$$
$$E_{\al_4}|_W=(1+\mfk{s})(E_{\al_4}|_{W_1})+\mfk{l}_{67}\ptl_{\mfk{l}_{72}}+\mfk{l}_{62}(\ptl_{\mfk{l}_{69}}
+\ptl_{\mfk{l}_{68}}
+2\ptl_{\mfk{l}_{67}}+\ptl_{\mfk{l}_{66}}),\eqno(12.2.91)$$
$$E_{\al_5}|_W=(1+\mfk{s})(E_{\al_5}|_{W_1})+\mfk{l}_{66}\ptl_{\mfk{l}_{73}}+\mfk{l}_{61}(\ptl_{\mfk{l}_{67}}
+2\ptl_{\mfk{l}_{66}}+\ptl_{\mfk{l}_{65}}),\eqno(12.2.92)$$
$$E_{\al_6}|_W=(1+\mfk{s})(E_{\al_6}|_{W_1})+\mfk{l}_{65}\ptl_{\mfk{l}_{74}}+
\mfk{l}_{60}(\ptl_{\mfk{l}_{66}}+2\ptl_{\mfk{l}_{65}}+\ptl_{\mfk{l}_{64}}),
\eqno(12.2.93)$$
$$E_{\al_7}|_W=(1+\mfk{s})(E_{\al_7}|_{W_1})+\mfk{l}_{64}\ptl_{\mfk{l}_{75}}+
\mfk{l}_{59}(\ptl_{\mfk{l}_{65}}+2\ptl_{\mfk{l}_{64}}).
\eqno(12.2.94)$$

 According to (12.1.30)-(12.1.36), (12.2.88)-(12.2.94)  Table 12.1.1 and Table 12.2.1, we find the following  singular
 vectors
\begin{eqnarray*}\vt&=&\frac{x_1}{2}\mfk{l}_{64}-x_2\mfk{l}_{59}+x_3\mfk{l}_{53}-
x_4\mfk{l}_{48}+x_5\mfk{l}_{44}-x_6\mfk{l}_{40}-x_7\mfk{l}_{37}\\&
&+x_8\mfk{l}_{33}+x_9\mfk{l}_{30}+x_{10}\mfk{l}_{28}-x_{11}\mfk{l}_{25}+x_{12}\mfk{l}_{23}
-x_{13}\mfk{l}_{21}+x_{14}\mfk{l}_{19}\\&
&+x_{15}\mfk{l}_{17}-x_{16}\mfk{l}_{16}+x_{18}\mfk{l}_{14}-x_{19}\mfk{l}_{12}
+x_{20}\mfk{l}_{11}+x_{21}\mfk{l}_{10} +x_{23}\mfk{l}_8\\& &
+x_{24}\mfk{l}_7+x_{26}\mfk{l}_6+x_{27}\mfk{l}_5
+x_{29}\mfk{l}_4+x_{32}\mfk{l}_3+x_{35}\mfk{l}_2+x_{40}\mfk{l}_1\hspace{3.8cm}(12.2.95)\end{eqnarray*}
of weight $\lmd_7$, and
$$\vs=\mfk{l}_1\mfk{l}_{19}-\mfk{l}_2\mfk{l}_{12}+\mfk{l}_3\mfk{l}_{10}-\mfk{l}_4
\mfk{l}_7+\mfk{l}_5\mfk{l}_6\eqno(12.2.96)$$ of weight $\lmd_6$.
Denote
$$I=\{\ol{1,6},\ol{9,12},15,\ol{17,19},21,22,\ol{26.29},\ol{33,36},\ol{40,44},\ol{50,55},58,63\}.\eqno(12.2.97)$$
We can similarly obtain the Cartan's  quartic
invariant\begin{eqnarray*}\eta &=&4(\sum_{i\in
I}\mfk{l}_i\mfk{l}_{134-i}-\sum_{r\in\ol{7,63}\setminus
I}\mfk{l}_r\mfk{l}_{134-r})+
4\mfk{l}_{70}(\mfk{l}_{70}-2\mfk{l}_{68}-3\mfk{l}_{69}+4\mfk{l}_{67}-3\mfk{l}_{66}\\
& &+2\mfk{l}_{65}-\mfk{l}_{64})+
\mfk{l}_{68}(7\mfk{l}_{68}+16\mfk{l}_{69}-24\mfk{l}_{67}+18\mfk{l}_{66}-12\mfk{l}_{65}+6\mfk{l}_{64})\\
&&+4\mfk{l}_{69}(3\mfk{l}_{69}-8\mfk{l}_{67}+6\mfk{l}_{66}
-4\mfk{l}_{65}+2\mfk{l}_{64})+4\mfk{l}_{67}(6\mfk{l}_{67}-9\mfk{l}_{66}+6\mfk{l}_{65}\\
&
&-3\mfk{l}_{64})+\mfk{l}_{66}(15\mfk{l}_{66}-20\mfk{l}_{65}+10\mfk{l}_{64})
+8\mfk{l}_{65}(\mfk{l}_{65}-\mfk{l}_{64})+3\mfk{l}_{64}^2.\hspace{3.1cm}(12.2.98)\end{eqnarray*}

\section{Decomposition of the Oscillator  Representation}

In this section, we want to prove the following theorem: \psp

{\bf Theorem 12.3.1}. {\it Any singular vector in ${\msr B}$ is a
polynomial in $x_1,\;\mfk{l}_1,\;\vt,\;\vs$ and $\eta$. Let
$L(n_1,n_2,n_3,n_4,n_5)$ be the irreducible submodule generated by
$\mfk{l}_1^{n_1}\vs^{n_2}x_1^{n_3}\vt^{n_4}\eta^{n_5}$ with highest
weight $n_1\lmd_1+n_2\lmd_6+(n_3+n_4)\lmd_7$. Then
$${\msr B}=\bigoplus_{n_1,n_2,n_3,n_4,n_5=0}^\infty
L(n_1,n_2,n_3,n_4,n_5).\eqno(12.3.1)$$ In particular,
\begin{eqnarray*}& &(1-q)^{55}\sum_{n_1,n_2,n_3,n_4=0}^\infty(\mbox{dim}\:
V(n_1\lmd_1+n_2\lmd_6+(n_3+n_4)\lmd_7))q^{2n_1+4n_2+n_3+3n_4}\\
&=&1+q+q^2+q^3.\hspace{10.8cm}(12.3.2)\end{eqnarray*}
  Let $\msr D$ be the invariant differential operator obtained
  from $\eta$ by changing $x_i$ to $\ptl_{x_i}$. Then
  $$\sum_{n_1,n_2,n_3=0}^\infty(L(n_1,n_2,n_3,0,0)+L(n_1,n_2,n_3,1,0))\subset
  \{f\in{\msr B}\mid \msr D(f)=0\}.\eqno(12.3.3)$$}

{\it Proof}. Let $f$ be a singular vector in ${\msr B}$. To
eliminate the extra variables from $f$, we need certain change of
variables. According to (12.2.1), (12.2.4)-(12.2.10),
(12.2.12)-(12.2.14), (12.2.16), (12.2.18), (12.2.19), (12.2.23),
(12.2.25), (12.2.27), (12.2.30), (12.2.32), (12.2.35), (12.2.39),
(12.2.42), (12.2.46), (12.2.50), (12.2.55), (12.2.61), (12.2.66) and
(12.2.96), we have
$$x_1x_{17}=\mfk{l}_1-x_2x_{14}-x_3x_{12}-x_4x_{10}-x_5x_9+x_6x_7.
\eqno(12.3.4)$$
$$x_1x_{22}=\mfk{l}_2-x_2x_{19}-x_3x_{16}-x_4x_{13}-x_5x_{11}+x_7x_8,
\eqno(12.3.5)$$
$$x_1x_{25}=\mfk{l}_3-x_2x_{21}-x_3x_{18}-x_4x_{15}-x_6x_{11}+x_8x_9,
\eqno(12.3.6)$$
$$x_1x_{28}=-\mfk{l}_4-x_2x_{24}-x_3x_{20}+x_5x_{15}-x_6x_{13}+x_8x_{10},
\eqno(12.3.7)$$
$$x_1x_{30}=\mfk{l}_5-x_2x_{26}-x_3x_{23}-x_7x_{15}+x_9x_{13}-x_{10}x_{11},
\eqno(12.3.8)$$
$$x_1x_{31}=\mfk{l}_6-x_2x_{27}+x_4x_{20}+x_5x_{18}-x_6x_{16}+x_8x_{12},
\eqno(12.3.9)$$
$$x_1x_{33}=-\mfk{l}_7-x_2x_{29}+x_4x_{23}-x_7x_{18}+x_9x_{16}-x_{11}x_{12},
\eqno(12.3.10)$$
$$x_1x_{34}=-\mfk{l}_8+x_3x_{27}+x_4x_{24}+x_5x_{21}-x_6x_{19}+x_8x_{14},
\eqno(12.3.11)$$
$$x_1x_{36}=\mfk{l}_{10}-x_2x_{32}-x_5x_{23}-x_7x_{20}+x_{10}x_{16}-x_{12}x_{13},
\eqno(12.3.12)$$
$$x_1x_{37}=\mfk{l}_{11}+x_3x_{29}+x_4x_{26}-x_7x_{21}+x_9x_{19}-x_{11}x_{14},
\eqno(12.3.13)$$
$$x_1x_{38}=-\mfk{l}_{12}-x_2x_{35}+x_6x_{23}+x_9x_{20}-x_{10}x_{18}+x_{12}x_{15},
\eqno(12.3.14)$$
$$x_1x_{39}=\mfk{l}_{14}+x_3x_{32}-x_5x_{26}-x_7x_{24}+x_{10}x_{19}-x_{13}x_{14},
\eqno(12.3.15)$$
$$x_1x_{41}=\mfk{l}_{16}+x_3x_{35}+x_6x_{26}+x_9x_{24}-x_{10}x_{21}+x_{14}x_{15},
\eqno(12.3.16)$$
$$x_1x_{42}=\mfk{l}_{17}-x_4x_{32}-x_5x_{29}-x_7x_{27}+x_{12}x_{19}-x_{14}x_{16},
\eqno(12.3.17)$$
\begin{eqnarray*}\hspace{1.8cm}x_1\mfk{l}_1x_{43}&=&
-\mfk{l}_1(x_2x_{40}+x_8x_{23}+x_{11}x_{20}-x_{13}x_{18}+x_{15}x_{16})
\\ & &+\vs+\mfk{l}_2\mfk{l}_{12}-\mfk{l}_3\mfk{l}_{10}+\mfk{l}_4\mfk{l}_7
-\mfk{l}_5\mfk{l}_6,\hspace{5cm}(12.3.18)\end{eqnarray*}
$$x_1x_{44}=-\mfk{l}_{21}-x_4x_{35}+x_6x_{29}+x_9x_{27}-x_{12}x_{21}+x_{14}x_{18},
\eqno(12.3.19)$$
$$x_1x_{45}=-\mfk{l}_{23}+x_3x_{40}-x_8x_{26}-x_{11}x_{24}+x_{13}x_{21}-x_{15}x_{19},
\eqno(12.3.20)$$
$$x_1x_{46}=\mfk{l}_{25}+x_5x_{35}+x_6x_{32}+x_{10}x_{27}-x_{12}x_{24}+x_{14}x_{20},
\eqno(12.3.21)$$
$$x_1x_{47}=\mfk{l}_{28}-x_4x_{40}-x_8x_{29}-x_{11}x_{27}+x_{16}x_{21}-x_{18}x_{19},
\eqno(12.3.22)$$
$$x_1x_{48}=-\mfk{l}_{30}+x_5x_{40}-x_8x_{32}-x_{13}x_{27}+x_{16}x_{24}-x_{19}x_{20},
\eqno(12.3.23)$$
$$x_1x_{49}=-\mfk{l}_{33}+x_7x_{35}+x_9x_{32}-x_{10}x_{29}+x_{12}x_{26}-x_{14}x_{23},
\eqno(12.3.24)$$
$$x_1x_{50}=\mfk{l}_{37}-x_6x_{40}-x_8x_{35}+x_{15}x_{27}-x_{18}x_{24}+x_{20}x_{21},
\eqno(12.3.25)$$
$$x_1x_{51}=-\mfk{l}_{40}+x_7x_{40}-x_{11}x_{32}+x_{13}x_{29}-x_{16}x_{26}+x_{19}x_{23},
\eqno(12.3.26)$$
$$x_1x_{52}=\mfk{l}_{44}-x_9x_{40}-x_{11}x_{35}-x_{15}x_{29}+x_{18}x_{26}-x_{21}x_{23},
\eqno(12.3.27)$$
$$x_1x_{53}=\mfk{l}_{48}+x_{10}x_{40}+x_{13}x_{35}+x_{15}x_{32}-x_{20}x_{26}+x_{23}x_{24},
\eqno(12.3.28)$$
$$x_1x_{54}=\mfk{l}_{53}-x_{12}x_{40}-x_{16}x_{35}-x_{18}x_{32}+x_{20}x_{29}-x_{23}x_{27},
\eqno(12.3.29)$$
$$x_1x_{55}=\mfk{l}_{59}+x_{14}x_{40}+x_{19}x_{35}+x_{21}x_{32}-x_{24}x_{29}+x_{26}x_{27}.
\eqno(12.3.30)$$

Moreover, (12.2.1) and (12.2.4)-(12.2.73) imply that when $x_i=0$
for $i\in\ol{3,54}$,
$$\vt=\frac{x_1}{2}(x_1x_{56}-x_2x_{55}),\eqno(12.3.31)$$
which is the homogeneous part of $\vt$ with degree 1 in $x_{55}$ and
$x_{56}$ (cf. (12.2.95)), and
$$\eta=3x_1^2x_{56}^2-6x_1x_2x_{55}x_{56}-5x_2^2x_{55}^2,\eqno(12.3.32)$$
which is  the homogeneous part of $\eta$ with degree 2 in $x_{55}$
and $x_{56}$ (cf. (12.2.98)). Under the assumption, we substitute
$x_1x_{56}=2x_1^{-1}\vt+x_2x_{55}$ into (12.3.32), we obtain
$$\eta=-8x_2^2x_{55}^2+12x_1^{-2}\vt^2.\eqno(12.3.33)$$
This shows that $f$ can be written as a function in
$\{x_i,\vt,\eta\mid i\in\ol{1,54}$. Moreover, using
(12.3.4)-(12.3.30), we obtain that $f$ is a function $f_1$ in
\begin{eqnarray*}\hspace{1.2cm}& &\{x_r,\mfk{l}_s,\vs,\vt,\eta\mid r\in \{\ol{1,35},40\}\setminus
\{17,22,25,28,30,31,33,34\},\\ &
&s\in\{\ol{1,17},21,23,25,28,30,33,37,40,44,48,53\}\setminus\{9,13,15\}\},\hspace{1.5cm}(12.3.34)\end{eqnarray*}
which is rational in the above variable except $\vt$ and $\eta$.

Denote
\begin{eqnarray*}\hspace{1.4cm}U&=&\sum_{i=1}^8\mbb{F}\mfk{l}_i+\sum_{13,15\neq
r\in\ol{10,17}}\mbb{F}\mfk{l}_r+\mbb{F}\mfk{l}_{19}+\mbb{F}\mfk{l}_{21}+\mbb{F}\mfk{l}_{23}+\mbb{F}\mfk{l}_{25}+\mbb{F}\mfk{l}_{28}
\\ &&+\mbb{F}\mfk{l}_{30}+\mbb{F}\mfk{l}_{33}+\mbb{F}\mfk{l}_{37}+\mbb{F}\mfk{l}_{40}
+\mbb{F}\mfk{l}_{44}+\mbb{F}\mfk{l}_{48}
+\mbb{F}\mfk{l}_{53}.\hspace{3.2cm}(12.3.35)\end{eqnarray*} By
(12.2.79)-(12.2.85), we get
$$E_{\al_1}|_U=-\mfk{l}_1\ptl_{\mfk{l}_2}-\mfk{l}_{12}\ptl_{\mfk{l}_{19}}
-\mfk{l}_{16}\ptl_{\mfk{l}_{23}}-\mfk{l}_{21}\ptl_{\mfk{l}_{28}}-\mfk{l}_{25}\ptl_{\mfk{l}_{30}}
-\mfk{l}_{33}\ptl_{\mfk{l}_{40}},\eqno(12.3.36)$$
$$E_{\al_2}|_U=\mfk{l}_4\ptl_{\mfk{l}_5}+\mfk{l}_6\ptl_{\mfk{l}_7}
+\mfk{l}_8\ptl_{\mfk{l}_{11}}+\mfk{l}_{25}\ptl_{\mfk{l}_{33}}
+\mfk{l}_{30}\ptl_{\mfk{l}_{40}}+\mfk{l}_{37}\ptl_{\mfk{l}_{44}},\eqno(12.3.37)$$
$$E_{al_3}|_U=-\mfk{l}_2\ptl_{\mfk{l}_3}-\mfk{l}_{10}\ptl_{\mfk{l}_{12}}
-\mfk{l}_{14}\ptl_{\mfk{l}_{16}}-\mfk{l}_{17}\ptl_{\mfk{l}_{21}}
-\mfk{l}_{30}\ptl_{\mfk{l}_{37}}-\mfk{l}_{40}\ptl_{\mfk{l}_{44}},\eqno(12.3.38)$$
$$E_{\al_4}|_U=\mfk{l}_3\ptl_{\mfk{l}_4}+\mfk{l}_7\ptl_{\mfk{l}_{10}}
+\mfk{l}_{11}\ptl_{\mfk{l}_{14}}+\mfk{l}_{21}\ptl_{\mfk{l}_{25}}
+\mfk{l}_{28}\ptl_{\mfk{l}_{30}}+\mfk{l}_{44}\ptl_{\mfk{l}_{48}},\eqno(12.3.39)$$
$$E_{\al_5}|_U=\mfk{l}_4\ptl_{\mfk{l}_6}+\mfk{l}_5\ptl_{\mfk{l}_7}
+\mfk{l}_{14}\ptl_{\mfk{l}_{17}}+\mfk{l}_{16}\ptl_{\mfk{l}_{21}}
+\mfk{l}_{23}\ptl_{\mfk{l}_{28}}+\mfk{l}_{48}\ptl_{\mfk{l}_{53}},\eqno(12.3.40)$$
$$E_{\al_6}|_U=\mfk{l}_6\ptl_{\mfk{l}_8}+\mfk{l}_7\ptl_{\mfk{l}_{11}}
+\mfk{l}_{10}\ptl_{\mfk{l}_{14}}+\mfk{l}_{12}\ptl_{\mfk{l}_{16}}+\mfk{l}_{19}\ptl_{\mfk{l}_{23}}
,\eqno(12.3.41)$$
$$E_{\al_7}|_U=0.\eqno(12.3.42)$$
It turns out that
$$E_{\al_i}(U)\subset U\qquad\for\;\;i\in\ol{1,7}.\eqno(12.3.43)$$
By induction, we can prove that
$$[E_\al,E_\be]|_U=[E_\al|_U,E_\be|_U]\qquad\for\;\;\al,\be\in\Phi_{E_7}^+.\eqno(12.3.44)$$
In particular,
$$E_\gm|_U=0\qquad\for\;\;\gm\in \Phi_{E_7}^+\setminus \Phi_{E_6}^+.\eqno(12.3.45)$$

According to (12.3.36)-(12.3.41) and (12.3.44), we have:
$$E_{(1,0,1)}|_U=\mfk{l}_1\ptl_{\mfk{l}_3}-\mfk{l}_{10}\ptl_{\mfk{l}_{19}}-\mfk{l}_{14}\ptl_{\mfk{l}_{23}}-\mfk{l}_{17}\ptl_{\mfk{l}_{28}}
+\mfk{l}_{25}\ptl_{\mfk{l}_{37}}+\mfk{l}_{33}\ptl_{\mfk{l}_{44}},
\eqno(12.3.46)$$
$$E_{(0,1,0,1)}|_U=-\mfk{l}_3\ptl_{\mfk{l}_5}+\mfk{l}_6\ptl_{\mfk{l}_{10}}
+\mfk{l}_8\ptl_{\mfk{l}_{14}}-\mfk{l}_{21}\ptl_{\mfk{l}_{33}}-\mfk{l}_{28}\ptl_{\mfk{l}_{40}}+\mfk{l}_{37}\ptl_{\mfk{l}_{48}},
\eqno(12.3.47)$$
$$E_{(0,0,1,1)}|_U=-\mfk{l}_2\ptl_{\mfk{l}_4}+\mfk{l}_7\ptl_{\mfk{l}_{12}}
+\mfk{l}_{11}\ptl_{\mfk{l}_{16}}
-\mfk{l}_{17}\ptl_{\mfk{l}_{25}}+\mfk{l}_{28}\ptl_{\mfk{l}_{37}}-\mfk{l}_{40}\ptl_{\mfk{l}_{48}},\eqno(12.3.48)$$
$$E_{(0,0,0,1,1)}|_U=\mfk{l}_3\ptl_{\mfk{l}_6}-\mfk{l}_5\ptl_{\mfk{l}_{10}}
+\mfk{l}_{11}\ptl_{\mfk{l}_{17}}-\mfk{l}_{16}\ptl_{\mfk{l}_{25}}-\mfk{l}_{23}\ptl_{\mfk{l}_{30}}
+\mfk{l}_{44}\ptl_{\mfk{l}_{53}},\eqno(12.3.49)$$
$$E_{(0,0,0,0,1,1)}|_U=\mfk{l}_4\ptl_{\mfk{l}_8}+\mfk{l}_5\ptl_{\mfk{l}_{11}}-\mfk{l}_{10}\ptl_{\mfk{l}_{17}}-\mfk{l}_{12}\ptl_{\mfk{l}_{21}}
-\mfk{l}_{19}\ptl_{\mfk{l}_{28}},\eqno(12.3.50)$$
$$E_{(1,0,1,1)}|_U=\mfk{l}_1\ptl_{\mfk{l}_4}+\mfk{l}_7\ptl_{\mfk{l}_{19}}
+\mfk{l}_{11}\ptl_{\mfk{l}_{23}}
-\mfk{l}_{17}\ptl_{\mfk{l}_{30}}-\mfk{l}_{21}\ptl_{\mfk{l}_{37}}+\mfk{l}_{33}\ptl_{\mfk{l}_{48}},\eqno(12.3.51)$$
$$E_{(0,1,1,1)}|_U=\mfk{l}_2\ptl_{\mfk{l}_5}+\mfk{l}_6\ptl_{\mfk{l}_{12}}
+\mfk{l}_8\ptl_{\mfk{l}_{16}}
+\mfk{l}_{17}\ptl_{\mfk{l}_{33}}-\mfk{l}_{28}\ptl_{\mfk{l}_{44}}-\mfk{l}_{30}\ptl_{\mfk{l}_{48}},\eqno(12.3.52)$$
$$E_{(0,1,0,1,1)}|_U=-\mfk{l}_3\ptl_{\mfk{l}_7}-\mfk{l}_4\ptl_{\mfk{l}_{10}}
+\mfk{l}_8\ptl_{\mfk{l}_{17}}+\mfk{l}_{16}\ptl_{\mfk{l}_{33}}+\mfk{l}_{23}\ptl_{\mfk{l}_{40}}
+\mfk{l}_{37}\ptl_{\mfk{l}_{53}},\eqno(12.3.53)$$
$$E_{(0,0,1,1,1)}|_U=-\mfk{l}_2\ptl_{\mfk{l}_6}-\mfk{l}_5\ptl_{\mfk{l}_{12}}
+\mfk{l}_{11}\ptl_{\mfk{l}_{21}}+\mfk{l}_{14}\ptl_{\mfk{l}_{25}}-\mfk{l}_{23}\ptl_{\mfk{l}_{37}}
-\mfk{l}_{40}\ptl_{\mfk{l}_{53}},\eqno(12.3.54)$$
$$E_{(0,0,0,1,1,1)}|_U=\mfk{l}_3\ptl_{\mfk{l}_8}-\mfk{l}_5\ptl_{\mfk{l}_{14}}-\mfk{l}_7\ptl_{\mfk{l}_{17}}+\mfk{l}_{12}\ptl_{\mfk{l}_{25}}
+\mfk{l}_{19}\ptl_{\mfk{l}_{30}},\eqno(12.3.55)$$
$$E_{(1,1,1,1)}|_U=-\mfk{l}_1\ptl_{\mfk{l}_5}+\mfk{l}_6\ptl_{\mfk{l}_{19}}
+\mfk{l}_8\ptl_{\mfk{l}_{23}}
+\mfk{l}_{17}\ptl_{\mfk{l}_{40}}+\mfk{l}_{21}\ptl_{\mfk{l}_{44}}+\mfk{l}_{25}\ptl_{\mfk{l}_{48}},\eqno(12.3.56)$$
$$E_{(1,0,1,1,1)}|_U=\mfk{l}_1\ptl_{\mfk{l}_6}-\mfk{l}_5\ptl_{\mfk{l}_{19}}
+\mfk{l}_{11}\ptl_{\mfk{l}_{28}}+\mfk{l}_{14}\ptl_{\mfk{l}_{30}}+\mfk{l}_{16}\ptl_{\mfk{l}_{37}}
+\mfk{l}_{33}\ptl_{\mfk{l}_{53}},\eqno(12.3.57)$$
$$E_{(0,1,1,1,1)}|_U=\mfk{l}_2\ptl_{\mfk{l}_7}-\mfk{l}_4\ptl_{\mfk{l}_{12}}
+\mfk{l}_8\ptl_{\mfk{l}_{21}}-\mfk{l}_{14}\ptl_{\mfk{l}_{33}}+\mfk{l}_{23}\ptl_{\mfk{l}_{44}}
-\mfk{l}_{30}\ptl_{\mfk{l}_{53}},\eqno(12.3.58)$$
$$E_{(0,1,0,1,1,1)}|_U=-\mfk{l}_3\ptl_{\mfk{l}_{11}}-\mfk{l}_4\ptl_{\mfk{l}_{14}}-\mfk{l}_6\ptl_{\mfk{l}_{17}}
-\mfk{l}_{12}\ptl_{\mfk{l}_{33}} -\mfk{l}_{19}\ptl_{\mfk{l}_{40}}
,\eqno(12.3.59)$$
$$E_{(0,0,1,1,1,1)}|_U=-\mfk{l}_2\ptl_{\mfk{l}_8}-\mfk{l}_5\ptl_{\mfk{l}_{16}}-\mfk{l}_7\ptl_{\mfk{l}_{21}}-\mfk{l}_{10}\ptl_{\mfk{l}_{25}}
+\mfk{l}_{19}\ptl_{\mfk{l}_{37}},\eqno(12.3.60)$$
$$E_{(1,1,1,1,1)}|_U=-\mfk{l}_1\ptl_{\mfk{l}_7}-\mfk{l}_4\ptl_{\mfk{l}_{19}}
+\mfk{l}_8\ptl_{\mfk{l}_{28}}-\mfk{l}_{14}\ptl_{\mfk{l}_{40}}-\mfk{l}_{16}\ptl_{\mfk{l}_{44}}
+\mfk{l}_{25}\ptl_{\mfk{l}_{53}},\eqno(12.3.61)$$
$$E_{(1,0,1,1,1,1)}|_U=\mfk{l}_1\ptl_{\mfk{l}_8}-\mfk{l}_5\ptl_{\mfk{l}_{23}}-\mfk{l}_7\ptl_{\mfk{l}_{28}}-\mfk{l}_{10}\ptl_{\mfk{l}_{30}}
-\mfk{l}_{12}\ptl_{\mfk{l}_{37}},\eqno(12.3.62)$$
$$E_{(0,1,1,1,1,1)}|_U=\mfk{l}_2\ptl_{\mfk{l}_{11}}-\mfk{l}_4\ptl_{\mfk{l}_{16}}-\mfk{l}_6\ptl_{\mfk{l}_{21}}+
\mfk{l}_{10}\ptl_{\mfk{l}_{33}}
-\mfk{l}_{19}\ptl_{\mfk{l}_{44}},\eqno(12.3.63)$$
$$E_{(0,1,1,2,1)}|_U=\mfk{l}_2\ptl_{\mfk{l}_{10}}+\mfk{l}_3\ptl_{\mfk{l}_{12}}
+\mfk{l}_8\ptl_{\mfk{l}_{25}}+\mfk{l}_{11}\ptl_{\mfk{l}_{33}}+\mfk{l}_{23}\ptl_{\mfk{l}_{48}}
+\mfk{l}_{28}\ptl_{\mfk{l}_{53}},\eqno(12.3.64)$$
$$E_{(1,1,1,2,1)}|_U=-\mfk{l}_1\ptl_{\mfk{l}_{10}}+\mfk{l}_3\ptl_{\mfk{l}_{19}}
+\mfk{l}_8\ptl_{\mfk{l}_{30}}+\mfk{l}_{11}\ptl_{\mfk{l}_{40}}-\mfk{l}_{16}\ptl_{\mfk{l}_{48}}
-\mfk{l}_{21}\ptl_{\mfk{l}_{53}},\eqno(12.3.65)$$
$$E_{(1,1,1,1,1,1)}|_U=-\mfk{l}_1\ptl_{\mfk{l}_{11}}-\mfk{l}_4\ptl_{\mfk{l}_{23}}-\mfk{l}_6\ptl_{\mfk{l}_{28}}
+\mfk{l}_{10}\ptl_{\mfk{l}_{40}} +\mfk{l}_{12}\ptl_{\mfk{l}_{44}}
,\eqno(12.3.66)$$
$$E_{(0,1,1,2,1,1)}|_U=\mfk{l}_2\ptl_{\mfk{l}_{14}}+\mfk{l}_3\ptl_{\mfk{l}_{16}}-\mfk{l}_6\ptl_{\mfk{l}_{25}}-
\mfk{l}_7\ptl_{\mfk{l}_{33}}
-\mfk{l}_{19}\ptl_{\mfk{l}_{48}},\eqno(12.3.67)$$
$$E_{(1,1,2,2,1)}|_U=-\mfk{l}_1\ptl_{\mfk{l}_{12}}-\mfk{l}_2\ptl_{\mfk{l}_{19}}
+\mfk{l}_8\ptl_{\mfk{l}_{37}}+\mfk{l}_{11}\ptl_{\mfk{l}_{44}}+\mfk{l}_{14}\ptl_{\mfk{l}_{48}}
+\mfk{l}_{17}\ptl_{\mfk{l}_{53}},\eqno(12.3.68)$$
$$E_{(1,1,1,2,1,1)}|_U=-\mfk{l}_1\ptl_{\mfk{l}_{14}}+\mfk{l}_3\ptl_{\mfk{l}_{23}}-\mfk{l}_6\ptl_{\mfk{l}_{30}}
-\mfk{l}_7\ptl_{\mfk{l}_{40}} +\mfk{l}_{12}\ptl_{\mfk{l}_{48}}
,\eqno(12.3.69)$$
$$E_{(0,1,1,2,2,1)}|_U=\mfk{l}_2\ptl_{\mfk{l}_{17}}+\mfk{l}_3\ptl_{\mfk{l}_{21}}+\mfk{l}_4\ptl_{\mfk{l}_{25}}+
\mfk{l}_5\ptl_{\mfk{l}_{33}}
-\mfk{l}_{19}\ptl_{\mfk{l}_{53}},\eqno(12.3.70)$$
$$E_{(1,1,2,2,1,1)}|_U=-\mfk{l}_1\ptl_{\mfk{l}_{16}}-\mfk{l}_2\ptl_{\mfk{l}_{23}}-\mfk{l}_6\ptl_{\mfk{l}_{37}}
-\mfk{l}_7\ptl_{\mfk{l}_{44}}
-\mfk{l}_{10}\ptl_{\mfk{l}_{48}},\eqno(12.3.71)$$
$$E_{(1,1,1,2,2,1)}|_U=-\mfk{l}_1\ptl_{\mfk{l}_{17}}+\mfk{l}_3\ptl_{\mfk{l}_{28}}+\mfk{l}_4\ptl_{\mfk{l}_{30}}+
\mfk{l}_5\ptl_{\mfk{l}_{40}}
+\mfk{l}_{12}\ptl_{\mfk{l}_{53}},\eqno(12.3.72)$$
$$E_{(1,1,2,2,2,1)}|_U=-\mfk{l}_1\ptl_{\mfk{l}_{21}}-\mfk{l}_2\ptl_{\mfk{l}_{28}}+\mfk{l}_4\ptl_{\mfk{l}_{37}}+
\mfk{l}_5\ptl_{\mfk{l}_{44}}
-\mfk{l}_{10}\ptl_{\mfk{l}_{53}},\eqno(12.3.73)$$
$$E_{(1,1,2,3,2,1)}|_U=-\mfk{l}_1\ptl_{\mfk{l}_{25}}-\mfk{l}_2\ptl_{\mfk{l}_{30}}-\mfk{l}_3\ptl_{\mfk{l}_{37}}+
\mfk{l}_5\ptl_{\mfk{l}_{48}}
+\mfk{l}_7\ptl_{\mfk{l}_{53}},\eqno(12.3.74)$$
$$E_{(1,2,2,3,2,1)}|_U=-\mfk{l}_1\ptl_{\mfk{l}_{33}}-\mfk{l}_2\ptl_{\mfk{l}_{40}}-\mfk{l}_3\ptl_{\mfk{l}_{44}}-
\mfk{l}_4\ptl_{\mfk{l}_{48}}
-\mfk{l}_6\ptl_{\mfk{l}_{53}}.\eqno(12.3.75)$$

By (12.1.83)-(12.1.92) and (12.3.45), we have
$$0=E_{(1,1,2,3,2,1,1)}(f_1)=-x_1\ptl_{x_{20}}(f_1),\eqno(12.3.76)$$
$$0=E_{(1,1,2,2,2,2,1)}(f_1)=-x_1\ptl_{x_{21}}(f_1),\eqno(12.3.77)$$
$$0=E_{(1,2,2,3,2,1,1)}(f_1)=x_1\ptl_{x_{23}}(f_1),\eqno(12.3.78)$$
$$0=E_{(1,1,2,3,2,2,1)}(f_1)=x_1\ptl_{x_{24}}(f_1),\eqno(12.3.79)$$
$$0=E_{(1,2,2,3,2,2,1)}(f_1)=-x_1\ptl_{x_{26}}(f_1),\eqno(12.3.80)$$
$$0=E_{(1,1,2,3,3,2,1)}(f_1)=-x_1\ptl_{x_{27}}(f_1),\eqno(12.3.81)$$
$$0=E_{(1,2,2,3,3,2,1)}(f_1)=x_1\ptl_{x_{29}}(f_1),\eqno(12.3.82)$$
$$0=E_{(1,2,2,4,3,2,1)}(f_1)=-x_1\ptl_{x_{32}}(f_1),\eqno(12.3.83)$$
$$0=E_{(1,2,3,4,3,2,1)}(f_1)=-x_1\ptl_{x_{35}}(f_1),\eqno(12.3.84)$$
$$0=E_{(2,2,3,4,3,2,1)}(f_1)=-x_1\ptl_{x_{40}}(f_1).\eqno(12.3.85)$$
So $f_1$ is independent of
$\{x_{20},x_{21},x_{23},x_{24},x_{26},x_{27},x_{29},x_{32},x_{35},x_{40}\}$;
that is, $f_1$ is  a  function  in
\begin{eqnarray*}\hspace{2cm}& &\{x_r,\mfk{l}_s,\vt,\vs\mid 17\neq r\in \ol{1,19};\;s\in\{\ol{1,17},21,23,\\
&&25,28,30,33,37,40,44,48,53\}\setminus\{9,13,15\}\}.\hspace{3.7cm}(12.3.86)\end{eqnarray*}

Expressions (12.1.36), (12.1.42), (12.1.48), (12.1.54), (12.1.59),
(12.1.60), (12.1.63), (12.1.64), (12.1.68), (12.1.70), (12.1.73),
(12.1.74), (12.1.76)-(12.1.78), (12.1.80), (12.1.81) and (12.3.45)
imply
$$0=E_{\al_7}(f_1)=x_1\ptl_{x_2}(f_1),\eqno(12.3.87)$$
$$0=E_{(0,0,0,0,0,1,1)}(f_1)=-x_1\ptl_{x_3}(f_1),\eqno(12.3.88)$$
$$0=E_{(0,0,0,0,1,1,1)}(f_1)=x_1\ptl_{x_4}(f_1),\eqno(12.3.89)$$
$$0=E_{(0,0,0,1,1,1,1)}(f_1)=-x_1\ptl_{x_5}(f_1),\eqno(12.3.90)$$
$$0=E_{(0,1,0,1,1,1,1)}(f_1)=x_1\ptl_{x_7}(f_1),\eqno(12.3.91)$$
$$0=E_{(0,0,1,1,1,1,1)}(f_1)=-x_1\ptl_{x_6}(f_1),\eqno(12.3.92)$$
$$0=E_{(1,0,1,1,1,1,1)}(f_1)=-x_1\ptl_{x_8}(f_1),\eqno(12.3.93)$$
$$0=E_{(0,1,1,1,1,1,1)}(f_1)=x_1\ptl_{x_9}(f_1),\eqno(12.3.94)$$
$$0=E_{(1,1,1,1,1,1,1)}(f_1)=x_1\ptl_{x_{11}}(f_1),\eqno(12.3.95)$$
$$0=E_{(0,1,1,2,1,1,1)}(f_1)=-x_1\ptl_{x_{10}}(f_1),\eqno(12.3.96)$$
$$0=E_{(1,1,1,2,1,1,1)}(f_1)=-x_1\ptl_{x_{13}}(f_1),\eqno(12.3.97)$$
$$0=E_{(0,1,1,2,2,1,1)}(f_1)=x_1\ptl_{x_{12}}(f_1),\eqno(12.3.98)$$
$$0=E_{(1,1,2,2,1,1,1)}(f_1)=-x_1\ptl_{x_{15}}(f_1),\eqno(12.3.99)$$
$$0=E_{(1,1,1,2,2,1,1)}(f_1)=x_1\ptl_{x_{16}}(f_1),\eqno(12.3.100)$$
$$0=E_{(0,1,1,2,2,2,1)}(f_1)=-x_1\ptl_{x_{14}}(f_1),\eqno(12.3.101)$$
$$0=E_{(1,1,2,2,2,1,1)}(f_1)=x_1\ptl_{x_{18}}(f_1),\eqno(12.3.102)$$
$$0=E_{(1,1,1,2,2,2,1)}(f_1)=-x_1\ptl_{x_{19}}(f_1).\eqno(12.3.103)$$
Hence $f_1$ is independent of $\{x_r,x_{18},x_{19}\mid\ol{2,16}\}$;
that is, $f_1$ is  a  function in
$$\{x_1,\mfk{l}_s,\vt,\vs\mid
s\in\{\ol{1,17},21,23,25,28,30,33,37,40,44,48,53\}\setminus\{9,13,15\}\}.\eqno(12.3.104)$$

Now by (12.3.36)-(12.3.41), (12.3.53), (12.3.59), (12.3.64),
(12.3.66), (12.3.68), (12.3.70), (12.3.71) and (12.3.73)-(12.3.75),
we obtain
$$\mfk{l}_1\ptl_{\mfk{l}_2}(f_1)+\mfk{l}_{16}\ptl_{\mfk{l}_{23}}(f_1)+\mfk{l}_{21}\ptl_{\mfk{l}_{28}}(f_1)
+\mfk{l}_{25}\ptl_{\mfk{l}_{30}}(f_1)
+\mfk{l}_{33}\ptl_{\mfk{l}_{40}}(f_1)=0,\eqno(12.3.105)$$
$$\mfk{l}_4\ptl_{\mfk{l}_5}(f_1)+\mfk{l}_6\ptl_{\mfk{l}_7}(f_1)
+\mfk{l}_8\ptl_{\mfk{l}_{11}}(f_1)+\mfk{l}_{25}\ptl_{\mfk{l}_{33}}(f_1)
+\mfk{l}_{30}\ptl_{\mfk{l}_{40}}(f_1)+\mfk{l}_{37}\ptl_{\mfk{l}_{44}}(f_1)=0,\eqno(12.3.106)$$
$$\mfk{l}_2\ptl_{\mfk{l}_3}(f_1)+\mfk{l}_{10}\ptl_{\mfk{l}_{12}}(f_1)
+\mfk{l}_{14}\ptl_{\mfk{l}_{16}}(f_1)+\mfk{l}_{17}\ptl_{\mfk{l}_{21}}(f_1)
+\mfk{l}_{30}\ptl_{\mfk{l}_{37}}(f_1)+\mfk{l}_{40}\ptl_{\mfk{l}_{44}}(f_1)=0,\eqno(12.3.107)$$
$$\mfk{l}_3\ptl_{\mfk{l}_4}(f_1)+\mfk{l}_7\ptl_{\mfk{l}_{10}}(f_1)
+\mfk{l}_{11}\ptl_{\mfk{l}_{14}}(f_1)+\mfk{l}_{21}\ptl_{\mfk{l}_{25}}(f_1)
+\mfk{l}_{28}\ptl_{\mfk{l}_{30}}(f_1)+\mfk{l}_{44}\ptl_{\mfk{l}_{48}}(f_1)=0,\eqno(12.3.108)$$
$$\mfk{l}_4\ptl_{\mfk{l}_6}(f_1)+\mfk{l}_5\ptl_{\mfk{l}_7}(f_1)
+\mfk{l}_{14}\ptl_{\mfk{l}_{17}}(f_1)+\mfk{l}_{16}\ptl_{\mfk{l}_{21}}(f_1)
+\mfk{l}_{23}\ptl_{\mfk{l}_{28}}(f_1)+\mfk{l}_{48}\ptl_{\mfk{l}_{53}}(f_1)=0,\eqno(12.3.109)$$
$$\mfk{l}_6\ptl_{\mfk{l}_8}(f_1)+\mfk{l}_7\ptl_{\mfk{l}_{11}}(f_1)
+\mfk{l}_{10}\ptl_{\mfk{l}_{14}}(f_1)+\mfk{l}_{12}\ptl_{\mfk{l}_{16}}(f_1)+\mfk{l}_{19}\ptl_{\mfk{l}_{23}}(f_1)
=0,\eqno(12.3.110)$$
$$-\mfk{l}_3\ptl_{\mfk{l}_7}(f_1)-\mfk{l}_4\ptl_{\mfk{l}_{10}}(f_1)
+\mfk{l}_8\ptl_{\mfk{l}_{17}}(f_1)+\mfk{l}_{16}\ptl_{\mfk{l}_{33}}(f_1)+\mfk{l}_{23}\ptl_{\mfk{l}_{40}}(f_1)
+\mfk{l}_{37}\ptl_{\mfk{l}_{53}}(f_1)=0,\eqno(12.3.111)$$
$$\mfk{l}_3\ptl_{\mfk{l}_{11}}(f_1)+\mfk{l}_4\ptl_{\mfk{l}_{14}}(f_1)+\mfk{l}_6\ptl_{\mfk{l}_{17}}(f_1)
+\mfk{l}_{12}\ptl_{\mfk{l}_{33}}(f_1)
+\mfk{l}_{19}\ptl_{\mfk{l}_{40}}(f_1) =0,\eqno(12.3.112)$$
$$\mfk{l}_2\ptl_{\mfk{l}_{10}}(f_1)+\mfk{l}_3\ptl_{\mfk{l}_{12}}(f_1)
+\mfk{l}_8\ptl_{\mfk{l}_{25}}(f_1)+\mfk{l}_{11}\ptl_{\mfk{l}_{33}}(f_1)+\mfk{l}_{23}\ptl_{\mfk{l}_{48}}(f_1)
+\mfk{l}_{28}\ptl_{\mfk{l}_{53}}(f_1)=0,\eqno(12.3.113)$$
$$\mfk{l}_2\ptl_{\mfk{l}_{14}}(f_1)+\mfk{l}_3\ptl_{\mfk{l}_{16}}(f_1)-\mfk{l}_6\ptl_{\mfk{l}_{25}}(f_1)-
\mfk{l}_7\ptl_{\mfk{l}_{33}}(f_1)
-\mfk{l}_{19}\ptl_{\mfk{l}_{48}}(f_1) =0,\eqno(12.3.114)$$
$$-\mfk{l}_1\ptl_{\mfk{l}_{12}}(f_1)
+\mfk{l}_8\ptl_{\mfk{l}_{37}}(f_1)+\mfk{l}_{11}\ptl_{\mfk{l}_{44}}(f_1)+\mfk{l}_{14}\ptl_{\mfk{l}_{48}}(f_1)
+\mfk{l}_{17}\ptl_{\mfk{l}_{53}}(f_1)=0,\eqno(12.3.115)$$
$$\mfk{l}_2\ptl_{\mfk{l}_{17}}(f_1)+\mfk{l}_3\ptl_{\mfk{l}_{21}}(f_1)+\mfk{l}_4\ptl_{\mfk{l}_{25}}(f_1)+
\mfk{l}_5\ptl_{\mfk{l}_{33}}(f_1)
-\mfk{l}_{19}\ptl_{\mfk{l}_{53}}(f_1) =0,\eqno(12.3.116)$$
$$\mfk{l}_1\ptl_{\mfk{l}_{16}}(f_1)+\mfk{l}_2\ptl_{\mfk{l}_{23}}(f_1)+\mfk{l}_6\ptl_{\mfk{l}_{37}}(f_1)
+\mfk{l}_7\ptl_{\mfk{l}_{44}}(f_1)
+\mfk{l}_{10}\ptl_{\mfk{l}_{48}}(f_1) =0,\eqno(12.3.117)$$
$$-\mfk{l}_1\ptl_{\mfk{l}_{21}}(f_1)-\mfk{l}_2\ptl_{\mfk{l}_{28}}(f_1)+\mfk{l}_4\ptl_{\mfk{l}_{37}}(f_1)+
\mfk{l}_5\ptl_{\mfk{l}_{44}}(f_1)
-\mfk{l}_{10}\ptl_{\mfk{l}_{53}}(f_1) =0,\eqno(12.3.118)$$
$$\mfk{l}_1\ptl_{\mfk{l}_{25}}(f_1)+\mfk{l}_2\ptl_{\mfk{l}_{30}}(f_1)+\mfk{l}_3\ptl_{\mfk{l}_{37}}(f_1)-
\mfk{l}_5\ptl_{\mfk{l}_{48}}(f_1) -\mfk{l}_7\ptl_{\mfk{l}_{53}}(f_1)
=0,\eqno(12.3.119)$$
$$\mfk{l}_1\ptl_{\mfk{l}_{33}}(f_1)+\mfk{l}_2\ptl_{\mfk{l}_{40}}(f_1)+\mfk{l}_3\ptl_{\mfk{l}_{44}}(f_1)+
\mfk{l}_4\ptl_{\mfk{l}_{48}}(f_1)+\mfk{l}_6\ptl_{\mfk{l}_{53}}(f_1)
=0.\eqno(12.3.120)$$

Thanks to (12.3.72), we have
$$-\mfk{l}_1\ptl_{\mfk{l}_{17}}(f_1)+\mfk{l}_3\ptl_{\mfk{l}_{28}}(f_1)+\mfk{l}_4\ptl_{\mfk{l}_{30}}(f_1)+
\mfk{l}_5\ptl_{\mfk{l}_{40}}(f_1)
+\mfk{l}_{12}\ptl_{\mfk{l}_{53}}(f_1) =0.\eqno(12.3.121)$$ Moreover,
$\mfk{l}_1\times (12.3.116)+\mfk{l}_2\times (12.3.121)$ gives
\begin{eqnarray*}\hspace{1cm}&&\mfk{l}_3(\mfk{l}_1\ptl_{\mfk{l}_{21}}(f_1)+\mfk{l}_2\ptl_{\mfk{l}_{28}}(f_1))
+\mfk{l}_4(\mfk{l}_1\ptl_{\mfk{l}_{25}}(f_1)+\mfk{l}_2\ptl_{\mfk{l}_{30}}(f_1))\\
& & +\mfk{l}_1\mfk{l}_5\ptl_{\mfk{l}_{33}}(f_1)
+\mfk{l}_2\mfk{l}_5\ptl_{\mfk{l}_{40}}(f_1)
+(\mfk{l}_2\mfk{l}_{12}-\mfk{l}_1\mfk{l}_{19})\ptl_{\mfk{l}_{53}}(f_1)
=0.\hspace{3cm}(12.3.122)\end{eqnarray*} Furthermore,
$(12.3.122)+\mfk{l}_3\times (12.3.118)-\mfk{l}_4\times (12.3.119)$
yields
\begin{eqnarray*}\hspace{1cm}&&\mfk{l}_5(\mfk{l}_1\ptl_{\mfk{l}_{33}}(f_1)+\mfk{l}_2\ptl_{\mfk{l}_{40}}(f_1)+
\mfk{l}_3\ptl_{\mfk{l}_{44}}(f_1)+\mfk{l}_4\ptl_{\mfk{l}_{48}}(f_1))
\\& &+(\mfk{l}_4\mfk{l}_7-\mfk{l}_3\mfk{l}_{10}
+\mfk{l}_2\mfk{l}_{12}-\mfk{l}_1\mfk{l}_{19})\ptl_{\mfk{l}_{53}}(f_1)
=0.\hspace{5.5cm}(12.3.123)\end{eqnarray*} Note that
$\mfk{l}_5\times (12.3.120)-(12.3.123)$ shows
$$(\mfk{l}_5\mfk{l}_6-\mfk{l}_4\mfk{l}_7+\mfk{l}_3\mfk{l}_{10}
-\mfk{l}_2\mfk{l}_{12}+\mfk{l}_1\mfk{l}_{19})\ptl_{\mfk{l}_{53}}(f_1)=0\lra
\ptl_{\mfk{l}_{53}}(f_1)=0.\eqno(12.3.124)$$

According to (12.3.69),
$$-\mfk{l}_1\ptl_{\mfk{l}_{14}}(f_1)+\mfk{l}_3\ptl_{\mfk{l}_{23}}(f_1)-\mfk{l}_6\ptl_{\mfk{l}_{30}}(f_1)
-\mfk{l}_7\ptl_{\mfk{l}_{40}}(f_1)
+\mfk{l}_{12}\ptl_{\mfk{l}_{48}}(f_1) =0.\eqno(12.3.125)$$ Moreover,
$\mfk{l}_1\times (12.3.114)+\mfk{l}_2\times (12.3.125)$ gives
\begin{eqnarray*}\hspace{0.7cm}&&\mfk{l}_3(\mfk{l}_1\ptl_{\mfk{l}_{16}}(f_1)+\mfk{l}_2\ptl_{\mfk{l}_{23}}(f_1))-
\mfk{l}_6(\mfk{l}_1\ptl_{\mfk{l}_{25}}(f_1)+\mfk{l}_2\ptl_{\mfk{l}_{30}}(f_1))\\
&&-\mfk{l}_7(\mfk{l}_1\ptl_{\mfk{l}_{33}}(f_1)+\mfk{l}_2\ptl_{\mfk{l}_{40}}(f_1))
+(\mfk{l}_2\mfk{l}_{12}-\mfk{l}_1\mfk{l}_{19})\ptl_{\mfk{l}_{48}}(f_1)
=0.\hspace{3.3cm}(12.3.126)\end{eqnarray*} Furthermore,
$(12.3.126)-\mfk{l}_3\times (12.3.117)+\mfk{l}_6\times
(12.3.119)+\mfk{l}_7\times (12.3.120)$ yields
$$(\mfk{l}_2\mfk{l}_{12}-\mfk{l}_3\mfk{l}_{10}+\mfk{l}_4\mfk{l}_7-\mfk{l}_5\mfk{l}_6-\mfk{l}_1\mfk{l}_{19})\ptl_{\mfk{l}_{48}}(f_1)
=0\lra \ptl_{\mfk{l}_{48}}(f_1)=0.\eqno(12.3.127)$$

Next (12.3.65) gives
$$-\mfk{l}_1\ptl_{\mfk{l}_{11}}(f_1)-\mfk{l}_4\ptl_{\mfk{l}_{23}}(f_1)-\mfk{l}_6\ptl_{\mfk{l}_{28}}
+\mfk{l}_{10}\ptl_{\mfk{l}_{40}}(f_1)
+\mfk{l}_{12}\ptl_{\mfk{l}_{44}}(f_1) =0.\eqno(12.3.128)$$ Moreover,
$\mfk{l}_1\times (12.3.112)+\mfk{l}_3\times (12.3.128)$ gives
\begin{eqnarray*}\qquad&&\mfk{l}_4(\mfk{l}_1\ptl_{\mfk{l}_{14}}(f_1)-\mfk{l}_3\ptl_{\mfk{l}_{23}}(f_1))+
\mfk{l}_6(\mfk{l}_1\ptl_{\mfk{l}_{17}}(f_1)-\mfk{l}_3\ptl_{\mfk{l}_{28}}(f_1))\\
& &+
\mfk{l}_{12}(\mfk{l}_1\ptl_{\mfk{l}_{33}}(f_1)+\mfk{l}_3\ptl_{\mfk{l}_{44}}(f_1))
+(\mfk{l}_3\mfk{l}_{10}+\mfk{l}_1\mfk{l}_{19})\ptl_{\mfk{l}_{40}}(f_1)
=0.\hspace{3cm}(12.3.129)\end{eqnarray*} Furthermore,
$(12.3.129)-\mfk{l}_{12}\times (12.3.120)+\mfk{l}_6\times
(12.3.121)+\mfk{l}_4\times (12.3.125)$ yields
$$(\mfk{l}_3\mfk{l}_{10}-\mfk{l}_2\mfk{l}_{12}-\mfk{l}_4\mfk{l}_7+\mfk{l}_5\mfk{l}_6+\mfk{l}_1\mfk{l}_{19})\ptl_{\mfk{l}_{40}}(f_1)
=0\lra \ptl_{\mfk{l}_{40}}(f_1)=0.\eqno(12.3.130)$$

Using (12.3.63), we get
$$\mfk{l}_2\ptl_{\mfk{l}_{11}}(f_1)-\mfk{l}_4\ptl_{\mfk{l}_{16}}(f_1)-\mfk{l}_6\ptl_{\mfk{l}_{21}}(f_1)+
\mfk{l}_{10}\ptl_{\mfk{l}_{33}}(f_1)
-\mfk{l}_{19}\ptl_{\mfk{l}_{44}}(f_1) =0.\eqno(12.3.131)$$
 Moreover,
$\mfk{l}_2\times (12.3.112)-\mfk{l}_3\times (12.3.131)$ gives
\begin{eqnarray*}\hspace{1cm}&&
\mfk{l}_4(\mfk{l}_2\ptl_{\mfk{l}_{14}}(f_1)+\mfk{l}_3\ptl_{\mfk{l}_{16}}(f_1))+
\mfk{l}_6(\mfk{l}_2\ptl_{\mfk{l}_{17}}(f_1)+\mfk{l}_3\ptl_{\mfk{l}_{21}}(f_1))
\\ & &+\mfk{l}_3\mfk{l}_{19}\ptl_{\mfk{l}_{44}}(f_1)+(\mfk{l}_2\mfk{l}_{12}-\mfk{l}_3\mfk{l}_{10})\ptl_{\mfk{l}_{33}}(f_1)
=0.\hspace{5.2cm}(12.3.132)\end{eqnarray*} Furthermore,
$(12.3.132)-\mfk{l}_4\times (12.3.114)-\mfk{l}_6\times
(12.3.116)-\mfk{l}_{19}\times (12.3.120)$ yields
$$(\mfk{l}_2\mfk{l}_{12}-\mfk{l}_1\mfk{l}_{19}-\mfk{l}_3\mfk{l}_{10}+\mfk{l}_4
\mfk{l}_7-\mfk{l}_5\mfk{l}_6)\ptl_{\mfk{l}_{33}}(f_1)=0\lra
\ptl_{\mfk{l}_{33}}(f_1)=0.\eqno(12.3.133)$$ Substituting the second
equations in (12.3.124), (12.3.127), (12.3.130) and (12.3.133) into
(12.3.120), we get $\ptl_{\mfk{l}_{44}}(f_1)=0.$

Note that (12.3.62) gives
$$\mfk{l}_1\ptl_{\mfk{l}_8}(f_1)-\mfk{l}_5\ptl_{\mfk{l}_{23}}(f_1)-\mfk{l}_7\ptl_{\mfk{l}_{28}}(f_1)-\mfk{l}_{10}\ptl_{\mfk{l}_{30}}(f_1)
-\mfk{l}_{12}\ptl_{\mfk{l}_{37}}(f_1)=0.\eqno(12.3.134)$$ Moreover,
$\mfk{l}_1\times(12.3.110)-\mfk{l}_6\times(12.3.134)$ yields
\begin{eqnarray*}\hspace{1cm}& &\mfk{l}_1(\mfk{l}_7\ptl_{\mfk{l}_{11}}(f_1)
+\mfk{l}_{10}\ptl_{\mfk{l}_{14}}(f_1)+\mfk{l}_{12}\ptl_{\mfk{l}_{16}}(f_1))+\mfk{l}_6\mfk{l}_7\ptl_{\mfk{l}_{28}}(f_1)\\
&&+\mfk{l}_6\mfk{l}_{10}\ptl_{\mfk{l}_{30}}(f_1)+(\mfk{l}_1\mfk{l}_{19}+\mfk{l}_5\mfk{l}_6)\ptl_{\mfk{l}_{23}}(f_1)
+\mfk{l}_6\mfk{l}_{12}\ptl_{\mfk{l}_{37}}(f_1)=0.\hspace{2.8cm}(12.3.135)\end{eqnarray*}
Furthermore, $(12.3.135)+\mfk{l}_7\times (12.3.128)$ says
\begin{eqnarray*}\hspace{1cm}& &\mfk{l}_1(
\mfk{l}_{10}\ptl_{\mfk{l}_{14}}(f_1)+\mfk{l}_{12}\ptl_{\mfk{l}_{16}}(f_1))+\mfk{l}_6\mfk{l}_{10}\ptl_{\mfk{l}_{30}}(f_1)\\
&&+(\mfk{l}_1\mfk{l}_{19}-\mfk{l}_4\mfk{l}_7+\mfk{l}_5\mfk{l}_6)\ptl_{\mfk{l}_{23}}(f_1)
+\mfk{l}_6\mfk{l}_{12}\ptl_{\mfk{l}_{37}}(f_1)=0\hspace{4.4cm}(12.3.136)\end{eqnarray*}
by the fact $\ptl_{\mfk{l}_{40}}(f_1)=\ptl_{\mfk{l}_{44}}(f_1)=0.$
Now $(12.3.136)+\mfk{l}_{10}\times (12.3.125)$ implies
$$\mfk{l}_1\mfk{l}_{12}\ptl_{\mfk{l}_{16}}(f_1)+(\mfk{l}_1\mfk{l}_{19}+\mfk{l}_3\mfk{l}_{10}-\mfk{l}_4\mfk{l}_7+\mfk{l}_5\mfk{l}_6)\ptl_{\mfk{l}_{23}}(f_1)
+\mfk{l}_6\mfk{l}_{12}\ptl_{\mfk{l}_{37}}(f_1)=0.\eqno(12.3.137)$$
Observe that now (12.3.117) becomes
$$\mfk{l}_1\ptl_{\mfk{l}_{16}}(f_1)+\mfk{l}_2\ptl_{\mfk{l}_{23}}(f_1)+\mfk{l}_6\ptl_{\mfk{l}_{37}}(f_1)=0.\eqno(12.3.138)$$
Besides, $(12.3.137)-\mfk{l}_{12}\times (12.3.137)$ gives
$$(\mfk{l}_1\mfk{l}_{19}-\mfk{l}_2\mfk{l}_{12}+\mfk{l}_3\mfk{l}_{10}-\mfk{l}_4\mfk{l}_7+\mfk{l}_5\mfk{l}_6)\ptl_{\mfk{l}_{23}}(f_1)
=0\lra \ptl_{\mfk{l}_{23}}(f_1)=0.\eqno(12.3.139)$$

According to (12.3.61),
$$-\mfk{l}_1\ptl_{\mfk{l}_7}(f_1)
+\mfk{l}_8\ptl_{\mfk{l}_{28}}(f_1)=0.\eqno(12.3.140)$$ Moreover,
(12.3.111) now becomes
$$-\mfk{l}_3\ptl_{\mfk{l}_7}(f_1)-\mfk{l}_4\ptl_{\mfk{l}_{10}}(f_1)
+\mfk{l}_8\ptl_{\mfk{l}_{17}}(f_1)=0,\eqno(12.3.141)$$ Subtracting
$\mfk{l}_1\times (12.3.141)$ from $\mfk{l}_3\times(12.3.140)$, we
get
$$\mfk{l}_1\mfk{l}_4\ptl_{\mfk{l}_{10}}(f_1)
+\mfk{l}_8(-\mfk{l}_1\ptl_{\mfk{l}_{17}}(f_1)
+\mfk{l}_3\ptl_{\mfk{l}_{28}}(f_1))=0.\eqno(12.3.142)$$ Moreover,
$(12.3.142)-\mfk{l}_8\times (12.3.121)$ gives
$$\mfk{l}_1\ptl_{\mfk{l}_{10}}(f_1)=\mfk{l}_8\ptl_{\mfk{l}_{30}}(f_1).\eqno(12.3.143)$$
Besides, now (12.3.115), (12.3.117), (12.3.125) and (12.3.128)
become
$$\mfk{l}_1\ptl_{\mfk{l}_{12}}(f_1)
=\mfk{l}_8\ptl_{\mfk{l}_{37}}(f_1),\qquad
\mfk{l}_1\ptl_{\mfk{l}_{16}}(f_1)=-\mfk{l}_6\ptl_{\mfk{l}_{37}}(f_1),\eqno(12.3.144)$$
$$\mfk{l}_1\ptl_{\mfk{l}_{11}}(f_1)=-\mfk{l}_6\ptl_{\mfk{l}_{28}},\qquad\mfk{l}_1\ptl_{\mfk{l}_{14}}(f_1)=
-\mfk{l}_6\ptl_{\mfk{l}_{30}}(f_1). \eqno(12.3.145)$$

Playing with (12.3.60), we find
$$-\mfk{l}_2\ptl_{\mfk{l}_8}(f_1)-\mfk{l}_5\ptl_{\mfk{l}_{16}}(f_1)-\mfk{l}_7\ptl_{\mfk{l}_{21}}(f_1)
-\mfk{l}_{10}\ptl_{\mfk{l}_{25}}(f_1)
+\mfk{l}_{19}\ptl_{\mfk{l}_{37}}(f_1)=0.\eqno(12.3.146)$$ On the
other hand, substituting (12.3.144) and (12.3.145) into
$\mfk{l}_1\times (12.3.110)$, we obtain
$$\mfk{l}_1\ptl_{\mfk{l}_8}(f_1)-\mfk{l}_7\ptl_{\mfk{l}_{28}}(f_1)
-\mfk{l}_{10}\ptl_{\mfk{l}_{30}}(f_1)-\mfk{l}_{12}\ptl_{\mfk{l}_{37}}(f_1)
=0.\eqno(12.3.147)$$ Calculating $\mfk{l}_1\times (12.3.146)$ by
(12.3.118), (12.3.119) and the second equation in (12.3.144), we
find \begin{eqnarray*}\hspace{2cm} &
&-\mfk{l}_1\mfk{l}_2\ptl_{\mfk{l}_8}(f_1)
+\mfk{l}_2\mfk{l}_{10}\ptl_{\mfk{l}_{30}}(f_1)+\mfk{l}_2\mfk{l}_7\ptl_{\mfk{l}_{28}}(f_1)
\\&&+(\mfk{l}_1\mfk{l}_{19}+\mfk{l}_3\mfk{l}_{10}-\mfk{l}_4\mfk{l}_7
+\mfk{l}_5\mfk{l}_6)\ptl_{\mfk{l}_{37}}(f_1)=0.\hspace{4.6cm}(12.3.148)\end{eqnarray*}
Moreover, $\mfk{l}_2\times(12.3.147)+(12.3.148)$ yields
$$(\mfk{l}_1\mfk{l}_{19}-\mfk{l}_2\mfk{l}_{12}+\mfk{l}_3\mfk{l}_{10}-\mfk{l}_4\mfk{l}_7+\mfk{l}_5\mfk{l}_6)\ptl_{\mfk{l}_{37}}(f_1)=0
\lra \ptl_{\mfk{l}_{37}}(f_1)=0.\eqno(12.3.149)$$ By (12.3.144),
$$\ptl_{\mfk{l}_{12}}(f_1)=\ptl_{\mfk{l}_{16}}(f_1)=0.\eqno(12.3.150)$$
Furthermore, (12.3.118) and (12.3.119) become
$$\mfk{l}_1\ptl_{\mfk{l}_{21}}(f_1)=-\mfk{l}_2\ptl_{\mfk{l}_{28}}(f_1),\qquad
\mfk{l}_1\ptl_{\mfk{l}_{25}}(f_1)=-\mfk{l}_2\ptl_{\mfk{l}_{30}}(f_1).\eqno(12.3.151)$$

Based on (12.3.57), we have
$$\mfk{l}_1\ptl_{\mfk{l}_6}(f_1)
+\mfk{l}_{11}\ptl_{\mfk{l}_{28}}(f_1)+\mfk{l}_{14}\ptl_{\mfk{l}_{30}}(f_1)=0.\eqno(12.3.152)$$
 Substituting (12.3.140) and (12.3.151) into $\mfk{l}_1\times (12.3.109)$,
we obtain
$$\mfk{l}_1\mfk{l}_4\ptl_{\mfk{l}_6}(f_1)
+\mfk{l}_1\mfk{l}_{14}\ptl_{\mfk{l}_{17}}(f_1)
+(\mfk{l}_1\mfk{l}_{23}-\mfk{l}_2\mfk{l}_{16}+\mfk{l}_5\mfk{l}_8)\ptl_{\mfk{l}_{28}}(f_1)=0.\eqno(12.3.153)$$
Moreover, we get
$$\mfk{l}_1\ptl_{\mfk{l}_{17}}(f_1)=\mfk{l}_3\ptl_{\mfk{l}_{28}}(f_1)+\mfk{l}_4\ptl_{\mfk{l}_{30}}(f_1)
\eqno(12.3.154)$$ by substituting (12.3.151) into (12.3.116). Again
we substitute (12.3.154) into (12.3.153):
$$\mfk{l}_1\mfk{l}_4\ptl_{\mfk{l}_6}(f_1)
+\mfk{l}_4\mfk{l}_{14}\ptl_{\mfk{l}_{30}}(f_1)
+(\mfk{l}_1\mfk{l}_{23}-\mfk{l}_2\mfk{l}_{16}+\mfk{l}_3\mfk{l}_{14}+\mfk{l}_5\mfk{l}_8)\ptl_{\mfk{l}_{28}}(f_1)=0.\eqno(12.3.155)$$
So $(12.3.155)-\mfk{l}_4\times (12.3.152)$ gives
$$(\mfk{l}_1\mfk{l}_{23}-\mfk{l}_2\mfk{l}_{16}+\mfk{l}_3\mfk{l}_{14}-
\mfk{l}_4\mfk{l}_{11}+\mfk{l}_5\mfk{l}_8)\ptl_{\mfk{l}_{28}}(f_1)=0\lra
\ptl_{\mfk{l}_{28}}(f_1)=0.\eqno(12.3.156)$$ By (12.3.140),
(12.3.145) and (12.3.151),
$$\ptl_{\mfk{l}_7}(f_1)=\ptl_{\mfk{l}_{11}}(f_1)=\ptl_{\mfk{l}_{21}}(f_1)=0.\eqno(12.3.157)$$
Now (12.3.56), (12.3.147) and (12.3.154) imply
$$\ptl_{\mfk{l}_5}(f_1)=0,\qquad \mfk{l}_1\ptl_{\mfk{l}_8}(f_1)=
\mfk{l}_{10}\ptl_{\mfk{l}_{30}}(f_1),\qquad
\mfk{l}_1\ptl_{\mfk{l}_{17}}(f_1)=\mfk{l}_4\ptl_{\mfk{l}_{30}}(f_1).\eqno(12.3.158)$$

Applying (12.3.55) to $f_1$, we have
$$\mfk{l}_3\ptl_{\mfk{l}_8}(f_1)-\mfk{l}_5\ptl_{\mfk{l}_{14}}(f_1)-\mfk{l}_7\ptl_{\mfk{l}_{17}}(f_1)+\mfk{l}_{12}\ptl_{\mfk{l}_{25}}(f_1)
+\mfk{l}_{19}\ptl_{\mfk{l}_{30}}(f_1)=0.\eqno(12.3.159)$$
Substituting (12.3.145), (12.3.151) and (12.3.158) into
$\mfk{l}_1\times(12.3.159$), we obtain
$$(\mfk{l}_3\mfk{l}_{10}+\mfk{l}_5\mfk{l}_6-\mfk{l}_4\mfk{l}_7-\mfk{l}_2\mfk{l}_{12}
+\mfk{l}_{19})\ptl_{\mfk{l}_{30}}(f_1)=0\lra
\ptl_{\mfk{l}_{30}}(f_1)=0.\eqno(12.3.161)$$ By (12.3.143),
(12.3.145), (12.3.151), (12.3.152) and (12.3.158), we get
$$\ptl_{\mfk{l}_6}(f_1)=\ptl_{\mfk{l}_8}(f_1)=\ptl_{\mfk{l}_{10}}(f_1)=\ptl_{\mfk{l}_{14}}(f_1)=\ptl_{\mfk{l}_{17}}(f_1)
=\ptl_{\mfk{l}_{25}}(f_1)=0.\eqno(12.3.162)$$ Up to this stage, we
have proved that $f_1$ is a function in
$\{x_1,\mfk{l}_1,\mfk{l}_2,\mfk{l}_3,\mfk{l}_4,\vt,\vs,\eta\}$.
Finally, (12.3.105), (12.3.107) and (12.3.108) give
$$\ptl_{\mfk{l}_2}(f_1)=\ptl_{\mfk{l}_3}(f_1)=\ptl_{\mfk{l}_4}(f_1)=0.\eqno(12.3.163)$$
This proves that $f_1$ is a function in
$\{x_1,\mfk{l}_1,\vt,\vs,\eta\}$, which is rational in
$x_1,\mfk{l}_1$ and $\vs$.

Note that $\vt$ and $\eta$ are functionally independent over the
rational functions in $\{x_i\mid i\in\ol{3,54}\}$ by (12.3.31) and
(12.3.32). Moreover,  the right hand sides of (12.3.31) and
(12.3.32) are coprime polynomials. Since $f_1=f$ is a polynomial in
$\{x_i\mid i\in\ol{1,56}\}$, there exists a positive integer $n$
such that $\ptl_{x_i}^n(f)=0$ for any $i\in\ol{1,56}$. Applying
$\ptl_{x_{56}}^n,\;\ptl_{x_{55}}^n,\;\ptl_{x_{43}}^n,\;\ptl_{x_{17}}^n$
and $\ptl_{x_1}^n$ to $f_1$ in the above order, we find that $f=f_1$
must be a polynomial in $\{x_1,\mfk{l}_1,\vt,\vs,\eta\}$ by
(12.2.1), (12.2.4)-(12.2.9), (12.2.12), (12.2.14), (12.3.17),
(12.3.30) and (12.3.31). Since each subspace of homogeneous
polynomials in ${\msr B}$ is a finite-dimensional ${\msr
G}^{E_7}$-module,  we get (12.3.1) by the Weyl's theorem of complete
reducibility and the fact that all finite-dimensional irreducible
${\msr G}^{E_7}$-modules are of highest-weight type. Moreover,
(12.3.1) gives
\begin{eqnarray*}& &(\sum_{r=0}^\infty q^4)\sum_{n_1,n_2,n_3,n_4=0}^\infty(\mbox{dim}\:
L(n_1,n_2,n_3,n_4,0))q^{2n_1+4n_2+n_3+3n_4}\\
&=&\frac{1}{(1-q)^{56}};\hspace{11.3cm}(12.3.164)\end{eqnarray*}
that is,
\begin{eqnarray*}& &\frac{1}{1-q^4}\sum_{n_1,n_2,n_3,n_4=0}^\infty(\mbox{dim}\:
L(n_1,n_2,n_3,n_4,0))q^{2n_1+4n_2+n_3+3n_4}\\
&=&\frac{1}{(1-q)^{56}};\hspace{11.2cm}(12.3.165)\end{eqnarray*}
which is equivalent to (12.3.2).

Observe that $\msr D(\mfk{l}_1^{m_1}\vs^{m_2}x_1^{m_3}\vt^{m_4})$ is
also a singular vector of weight
$m_1\lmd_1+m_2\lmd_6+(m_3+m_4)\lmd_7$ if it is nonzero. So is a
linear combination of the elements
$\mfk{l}_1^{m_1}\vs^{m_2}x_1^{n_3}\vt^{n_4}\eta^{n_5}$ such that
$$n_3+n_4=m_3+m_4,\qquad n_3+3n_4+4n_5=m_3+3m_4-4.\eqno(12.3.166)$$
Thus
$$m_4=2+n_4+2n_5\geq 2.\eqno(12.3.167)$$
This shows
$$\msr D(\mfk{l}_1^{m_1}\vs^{m_2}x_1^{m_3})={\msr
D}(\mfk{l}_1^{m_1}\vs^{m_2}x_1^{m_3}\vt)=0\eqno(12.3.168)$$ for any
nonnegative integers $m_1,m_2$ and $m_3$, or equivalently, (12.3.3)
holds. The proof of our  theorem is completed. $\qquad\Box$

\section{Realization of $E_7$ in 27-Dimensional Space}

In this section, we find a new representation of the simple Lie
algebra of type $E_7$ on the polynomial algebra in 27 variables,
which gives a fractional representation of the corresponding Lie
group on 27-dimensional space.

Let us go back to the construction of the basic oscillator
representation $\mfk r$ of $\msr G^{E_6}$ via the lattice
construction of $\msr G^{E_7}$ in Section 11.1.

Take the notations in (4.4.42)-(4.4.44) and (11.1.4)-(11.1.12). For
any $i\in\ol{1,27}$,  we write
$$\mfk a_i=E_\al\qquad\mbox{if}\;\;\mfk b_i=E'_\al.\eqno(12.4.1)$$
Set
$${\msr G}_-=\sum_{i=1}^{27}\mbb F\mfk b_i,\qquad {\msr G}_0={\msr
G}^{E_6}+\mbb F\al_7,\qquad {\msr G}_+=\sum_{i=1}^{27}\mbb F\mfk
a_i.\eqno(12.4.2)$$ Then ${\msr G}_\pm$ are abelian subalgebras of
${\msr G}^{E_7}$ and ${\msr G}_0$ is a maximal reductive Lie
subalgebra of ${\msr G}^{E_7}$. Moreover,
$$[{\msr G}_+,{\msr G}_-]\subset {\msr G}_0,\qquad[{\msr G}_0,{\msr
G}_\pm]\subset{\msr G}_\pm,\qquad{\msr G}^{E_7}={\msr
G}_-\oplus{\msr G}_0\oplus {\msr G}_+.\eqno(12.4.3)$$ Denote by
$\lmd_i$ the $i$th fundamental weight of ${\msr G}^{E_6}$. With
respect to the adjoint representation of ${\msr G}^{E_7}$, ${\msr
G}_+$ forms an irreducible ${\msr G}^{E_6}$-module with highest
weight $\lmd_1$ and ${\msr G}_-$ forms an irreducible ${\msr
G}^{E_6}$-module with highest weight $\lmd_6$.

Take the oscillator representation of $\msr G^{E_6}$ on
$${\msr A}=\mbb F[x_1,x_2,...,x_{27}]\eqno(12.4.4)$$
given in (11.1.18)-(11.1.56). Let
$$\widehat\al=2\al_1+3\al_2+4\al_3+6\al_4+5\al_5+4\al_6+3\al_7.\eqno(12.4.5)$$
Then
$$(\widehat\al,\al_r)=0\qquad\for\;\;r\in\ol{1,6}\eqno(12.4.6)$$
by the Dynkin diagram of $E_7$, where $(\cdot,\cdot)$ is a symmetric
$\mbb Z$-bilinear form on $\Lmd_r$ such that $(\al,\al)=2$ for
$\al\in \Phi_{E_7}$. Thanks to (4.4.24),
$$[\widehat\al,{\msr G}^{D_6}]=0.\eqno(12.4.7)$$
By Schur's Lemma, $\wht\al|_V=c\mbox{Id}_V$ for some constant $c$,
where $V=\msr G_-$ (cf. (11.1.13)). According to (11.1.15) and
(11.1.17), $\mfk r(\wht\al)=\wht\al|_{\msr
A}=c\sum_{i=1}^{16}x_i\ptl_{x_i}$. Expressions (4.4.24),
(11.1.4)-(11.1.12), (11.1.15), (11.1.17) and the Dynkin diagram of
$E_7$ yields
$$\mfk r(\al_7)=\al_7|_{\msr
A}=-2x_1\ptl_{x_1}-\sum_{i=2}^{13}x_i\ptl_{x_i}-x_{15}\ptl_{x_{15}}
-x_{16}\ptl_{x_{16}}-x_{18}\ptl_{x_{18}}
-x_{20}\ptl_{x_{20}}.\eqno(12.4.8)$$ According to (12.4.6),
(12.4.8), the coefficients of $x_1\ptl_{x_1}$ in (11.1.54) and Table
11.1.1, we have that
$$\mfk r(\wht\al)=-2D,\qquad\mbox{where}\;\;D=\sum_{i=1}^{16}x_i\ptl_{x_i}\eqno(12.4.9)$$
is the degree operator on ${\msr A}$.

Write
$$\msr T=\sum_{i=1}^{27}\mbb F\ptl_{x_i}.\eqno(12.4.10)$$
Then $\msr T$ forms a ${\msr G}^{E_6}$-module with respect to the
action
$$u(\ptl)=[\mfk r(u),\ptl]\qquad\for\;\;u\in{\msr G}^{E_6},\;\ptl\in\td\Dlt.\eqno(12.4.11)$$ On the other hand,
${\msr G}_\pm$ (cf. (11.1.4)-(11.1.12), (12.4.1) and (12.4.2)) form
${\msr G}^{E_6}$-modules with respect to the adjoint representation.
 According to (11.1.15)
and (11.1.17), the linear map determined by $\mfk b_i\mapsto x_i$
for $i\in\ol{1,27}$ gives a ${\msr G}^{E_6}$-module monomorphism
from ${\msr G}_-$ to $\msr A$.  Define a bilinear form
$(\cdot|\cdot)$ on ${\msr G}^{E_7}$ by
$$(h_1|h_2)=(h_1,h_2),\;\; (h|E_{\al})=0,\;\;
 (E_{\al}|E_{\be})=-\dlt_{\al+\be,0}\eqno(12.4.12)$$
for $h_1,h_2\in H$ and $\al,\be\in \Phi_{E_7}$. It can be verified
that $(\cdot|\cdot)$ is a ${\msr G}^{E_7}$-invariant form; that is,
$$([u,v]|w)=-(v|[u,w])\qquad\for\;\;u,v\in{\msr
G}^{E_7}\eqno(12.4.13)$$(cf. (4.4.24) and (4.4.25)). Hence the
linear map determined by $\mfk a_i\mapsto \ptl_{x_i}$ for
$i\in\ol{1,27}$ gives a ${\msr G}^{E_6}$-module isomorphism from
${\msr G}_+$ to $\msr T$. Hence we define the action of ${\msr G}_+$
on ${\msr A}$ by
$$\mfk a_i|_{\msr
A}=\ptl_{x_i}\qquad\for\;\;i\in\ol{1,27}.\eqno(12.4.14)$$

Recall the Witt algebra
$$\mbb W_{27}=\sum_{i=1}^{27}{\msr A}\ptl_{x_i}\eqno(12.4.15)$$
with the Lie bracket (6.7.9). Now we want to find the differential
operators $P_1,P_2,...,P_{27}\in \mbb W_{27}$ such that the
following action matches the structure of ${\msr G}^{E_7}$:
$$\mfk b_i|_{\msr A}=P_i\qquad\for\;\;i\in\ol{1,27}.\eqno(12.4.16)$$
Comparing the weights in Table 11.1.1 and Table 11.2.1, we use
(11.2.1) and (11.2.3)-(11.2.28) to assume
\begin{eqnarray*}\qquad P_1&=&x_1D+c_1\zeta_1\ptl_{x_{14}}+c_2\zeta_2\ptl_{x_{17}}+c_3\zeta_3\ptl_{x_{19}}
+c_4\zeta_4\ptl_{x_{21}}+c_5\zeta_5\ptl_{x_{22}}
\\&&+c_6\zeta_6\ptl_{x_{23}}+c_7\zeta_7\ptl_{x_{24}}+c_9\zeta_9\ptl_{x_{25}}+c_{11}\zeta_{11}\ptl_{x_{26}}
+c_{14}\zeta_{14}\ptl_{x_{27}},\hspace{2.15cm}(12.4.17)\end{eqnarray*}
where $c_i\in\mbb{C}$. Imposing
\begin{eqnarray*}\qquad [\ptl_{x_1},P_1]&=&\mfk r([\mfk a_1,\mfk b_1])=-\mfk r(\al_7)=2x_1\ptl_{x_1}+\sum_{i=2}^{13}x_i\ptl_{x_i}\\&
&+x_{15}\ptl_{x_{15}} +x_{16}\ptl_{x_{16}}+x_{18}\ptl_{x_{18}}
+x_{20}\ptl_{x_{20}}\hspace{4.4cm}(12.4.18)\end{eqnarray*} by
(4.4.25) and (12.4.8), we get
$$c_1=c_2=c_3=-c_4=c_5=c_6=-c_7=c_9=c_{11}=c_{14}=-1\eqno(12.4.19)$$
by (11.2.1), (11.2.3)-(11.2.8), (11.2.10), (11.2.12) and (11.2.15).
Thus
\begin{eqnarray*}\qquad\qquad
P_1&=&x_1D-\zeta_1\ptl_{x_{14}}-\zeta_2\ptl_{x_{17}}-\zeta_3\ptl_{x_{19}}
+\zeta_4\ptl_{x_{21}}-\zeta_5\ptl_{x_{22}}
\\&&-\zeta_6\ptl_{x_{23}}+\zeta_7\ptl_{x_{24}}-\zeta_9\ptl_{x_{25}}-\zeta_{11}\ptl_{x_{26}}
-\zeta_{14}\ptl_{x_{27}}.\hspace{3.4cm}(12.4.20)\end{eqnarray*}

According to (11.2.1), (11.2.3)-(11.2.8), (11.2.10), (11.2.12) and
(11.2.15),
 we find
$$[\ptl_{x_2},P_1]=x_1\ptl_{x_2}-x_{11}\ptl_{x_{14}}-x_{13}\ptl_{x_{17}}
-x_{16}\ptl_{x_{19}}-x_{18}\ptl_{x_{21}}-x_{20}\ptl_{x_{23}},\eqno(12.4.21)$$
$$[\ptl_{x_3},P_1]=x_1\ptl_{x_3}-x_9\ptl_{x_{14}}-x_{12}\ptl_{x_{17}}
-x_{15}\ptl_{x_{19}}+x_{18}\ptl_{x_{22}}+x_{20}\ptl_{x_{24}},\eqno(12.4.22)$$
$$[\ptl_{x_4},P_1]=x_1\ptl_{x_4}-x_7\ptl_{x_{14}}-x_{10}\ptl_{x_{17}}
+x_{15}\ptl_{x_{21}}+x_{16}\ptl_{x_{22}}-x_{20}\ptl_{x_{25}},\eqno(12.4.23)$$
$$[\ptl_{x_5},P_1]=x_1\ptl_{x_5}+x_6\ptl_{x_{14}}-x_{10}\ptl_{x_{19}}
-x_{12}\ptl_{x_{21}}-x_{13}\ptl_{x_{22}}+x_{20}\ptl_{x_{25}},\eqno(12.4.24)$$
$$[\ptl_{x_6},P_1]=x_1\ptl_{x_6}+x_5\ptl_{x_{14}}+x_8\ptl_{x_{17}}
-x_{15}\ptl_{x_{23}}-x_{16}\ptl_{x_{24}}-x_{18}\ptl_{x_{25}},\eqno(12.4.25)$$
$$[\ptl_{x_7},P_1]=x_1\ptl_{x_7}-x_4\ptl_{x_{14}}+x_8\ptl_{x_{19}}
+x_{12}\ptl_{x_{23}}+x_{13}\ptl_{x_{24}}+x_{18}\ptl_{x_{26}},\eqno(12.4.26)$$
$$[\ptl_{x_8},P_1]=x_1\ptl_{x_8}+x_6\ptl_{x_{17}}+x_7\ptl_{x_{19}}
+x_9\ptl_{x_{21}}+x_{11}\ptl_{x_{22}}+x_{20}\ptl_{x_{27}},\eqno(12.4.27)$$
$$[\ptl_{x_9},P_1]=x_1\ptl_{x_9}-x_3\ptl_{x_{14}}+x_8\ptl_{x_{21}}
-x_{10}\ptl_{x_{23}}+x_{13}\ptl_{x_{25}}-x_{16}\ptl_{x_{26}},\eqno(12.4.28)$$
$$[\ptl_{x_{10}},P_1]=x_1\ptl_{x_{10}}-x_4\ptl_{x_{17}}-x_5\ptl_{x_{19}}
-x_9\ptl_{x_{23}}-x_{11}\ptl_{x_{24}}+x_{18}\ptl_{x_{27}},\eqno(12.4.29)$$
$$[\ptl_{x_{11}},P_1]=x_1\ptl_{x_{11}}-x_2\ptl_{x_{14}}+x_8\ptl_{x_{22}}
-x_{10}\ptl_{x_{24}}-x_{12}\ptl_{x_{25}}+x_{15}\ptl_{x_{26}},\eqno(12.4.30)$$
$$[\ptl_{x_{12}},P_1]=x_1\ptl_{x_{12}}-x_3\ptl_{x_{17}}-x_5\ptl_{x_{21}}
+x_7\ptl_{x_{23}}-x_{11}\ptl_{x_{25}}-x_{16}\ptl_{x_{27}},\eqno(12.4.31)$$
$$[\ptl_{x_{13}},P_1]=x_1\ptl_{x_{13}}-x_2\ptl_{x_{17}}-x_5\ptl_{x_{22}}
+x_7\ptl_{x_{24}}+x_9\ptl_{x_{25}}+x_{15}\ptl_{x_{27}},\eqno(12.4.32)$$
$$[\ptl_{x_{15}},P_1]=x_1\ptl_{x_{15}}-x_3\ptl_{x_{19}}+x_4\ptl_{x_{21}}
-x_6\ptl_{x_{23}}+x_{11}\ptl_{x_{26}}+x_{13}\ptl_{x_{27}},\eqno(12.4.33)$$
$$[\ptl_{x_{16}},P_1]=x_1\ptl_{x_{16}}-x_2\ptl_{x_{19}}+x_4\ptl_{x_{22}}
-x_6\ptl_{x_{24}}-x_9\ptl_{x_{26}}-x_{12}\ptl_{x_{27}},\eqno(12.4.34)$$
$$[\ptl_{x_{18}},P_1]=x_1\ptl_{x_{18}}-x_2\ptl_{x_{21}}+x_3\ptl_{x_{22}}
-x_6\ptl_{x_{25}}+x_7\ptl_{x_{26}}+x_{10}\ptl_{x_{27}},\eqno(12.4.35)$$
$$[\ptl_{x_{20}},P_1]=x_1\ptl_{x_{20}}-x_2\ptl_{x_{23}}+x_3\ptl_{x_{24}}
-x_4\ptl_{x_{25}}+x_5\ptl_{x_{26}}+x_8\ptl_{x_{27}},\eqno(12.4.36)$$
$$[x_r,P_1]=0\qquad\for\;\;r=14,17,19,21,22,23,24,25,26,27.\eqno(12.4.37)$$
On the other hand, if $[\mfk a_i,\mfk b_1]\neq 0$ for $i\geq 2$,
then it is equal to the vector obtained by deleting the last
coordinate $1$ (cf. (11.1.4)-(11.1.12) and (12.4.1)). For instance,
$\mfk a_{18}=E_{(1,1,2,3,2,1,1)}$ and $[\mfk a_{18},\mfk
b_1]=E_{(1,1,2,3,2,1)}$. By  (12.4.21)-(12.4.36) and correspondingly
(11.1.23), (11.1.28), (11.1.33), (11.1.38), (11.1.37), (11.1.42),
(11.1.40), (11.1.45), (11.1.44), (11.1.48), (11.1.47), (11.1.50),
(11.1.49),  (11.1.51)-(11.1.53), we have:
$$[\ptl_{x_i},P_1]=[\mfk a_i,\mfk b_1]|_{\msr
A}\qquad\for\;\;i\in\ol{1,27}.\eqno(12.4.38)$$

Expressions (6.7.9), (11.1.18)-(11.1.23) and (11.2.29)-(11.2.34)
imply that $P_1$ is a ${\msr G}^{E_6}$-singular vector in $\mbb
W_{27}$ with weight $\lmd_6$. We set
\begin{eqnarray*}
\qquad P_2&=&-[E_{-\al_6},P_1]=
x_2D-\zeta_1\ptl_{x_{11}}-\zeta_2\ptl_{x_{13}}-\zeta_3\ptl_{x_{16}}
+\zeta_4\ptl_{x_{18}}\\ &
&-\zeta_6\ptl_{x_{20}}-\zeta_8\ptl_{x_{22}}
+\zeta_{10}\ptl_{x_{24}}-\zeta_{12}\ptl_{x_{25}}-\zeta_{15}\ptl_{x_{26}}
+\zeta_{17}\ptl_{x_{27}},\hspace{2.2cm}(12.4.39)\end{eqnarray*}
\begin{eqnarray*}
\qquad P_3&=&-[E_{-\al_5},P_2]=
x_3D-\zeta_1\ptl_{x_9}-\zeta_2\ptl_{x_{12}}-\zeta_3\ptl_{x_{15}}
+\zeta_5\ptl_{x_{18}}\\ &
&-\zeta_7\ptl_{x_{20}}-\zeta_8\ptl_{x_{21}}
+\zeta_{10}\ptl_{x_{23}}-\zeta_{13}\ptl_{x_{25}}-\zeta_{16}\ptl_{x_{26}}
+\zeta_{19}\ptl_{x_{27}},\hspace{2.2cm}(12.4.40)\end{eqnarray*}
\begin{eqnarray*}
\qquad P_4&=&-[E_{-\al_4},P_3]=
x_4D-\zeta_1\ptl_{x_7}-\zeta_2\ptl_{x_{10}}-\zeta_4\ptl_{x_{15}}
+\zeta_5\ptl_{x_{16}}\\ &
&-\zeta_8\ptl_{x_{19}}-\zeta_9\ptl_{x_{20}}
+\zeta_{12}\ptl_{x_{23}}-\zeta_{13}\ptl_{x_{24}}-\zeta_{18}\ptl_{x_{26}}
+\zeta_{21}\ptl_{x_{27}},\hspace{2.3cm}(12.4.41)\end{eqnarray*}
\begin{eqnarray*}
\qquad P_5&=&-[E_{-\al_3},P_4]=
x_5D+\zeta_1\ptl_{x_6}-\zeta_3\ptl_{x_{10}}+\zeta_4\ptl_{x_{12}}
-\zeta_5\ptl_{x_{13}}\\ &
&+\zeta_8\ptl_{x_{17}}+\zeta_{11}\ptl_{x_{20}}
-\zeta_{15}\ptl_{x_{23}}+\zeta_{16}\ptl_{x_{24}}-\zeta_{18}\ptl_{x_{25}}
-\zeta_{22}\ptl_{x_{27}},\hspace{2.2cm}(12.4.42)\end{eqnarray*}
\begin{eqnarray*}
\qquad P_6&=&-[E_{-\al_2},P_4]=
x_6D+\zeta_1\ptl_{x_5}+\zeta_2\ptl_{x_8}-\zeta_6\ptl_{x_{15}}
+\zeta_7\ptl_{x_{16}}\\ &
&-\zeta_9\ptl_{x_{18}}-\zeta_{10}\ptl_{x_{19}}
+\zeta_{12}\ptl_{x_{21}}-\zeta_{13}\ptl_{x_{22}}+\zeta_{20}\ptl_{x_{26}}
-\zeta_{23}\ptl_{x_{27}},\hspace{2.2cm}(12.4.43)\end{eqnarray*}
\begin{eqnarray*}
\qquad P_7&=&-[E_{-\al_3},P_6]=
x_7D-\zeta_1\ptl_{x_4}+\zeta_3\ptl_{x_8}+\zeta_6\ptl_{x_{12}}
-\zeta_7\ptl_{x_{13}}\\ &
&+\zeta_{10}\ptl_{x_{17}}+\zeta_{11}\ptl_{x_{18}}
-\zeta_{15}\ptl_{x_{21}}+\zeta_{16}\ptl_{x_{22}}+\zeta_{20}\ptl_{x_{25}}
+\zeta_{24}\ptl_{x_{27}},\hspace{2.1cm}(12.4.44)\end{eqnarray*}
\begin{eqnarray*}
\qquad P_8&=&-[E_{-\al_2},P_5]=
x_8D+\zeta_2\ptl_{x_6}+\zeta_3\ptl_{x_7}-\zeta_4\ptl_{x_9}
+\zeta_5\ptl_{x_{11}}\\ &
&-\zeta_8\ptl_{x_{14}}+\zeta_{14}\ptl_{x_{20}}
+\zeta_{17}\ptl_{x_{23}}-\zeta_{19}\ptl_{x_{24}}+\zeta_{21}\ptl_{x_{25}}
+\zeta_{22}\ptl_{x_{26}},\hspace{2.1cm}(12.4.45)\end{eqnarray*}
 \begin{eqnarray*}
\qquad P_9&=&[E_{-\al_4},P_7]=
x_9D-\zeta_1\ptl_{x_3}-\zeta_4\ptl_{x_8}-\zeta_6\ptl_{x_{10}}
+\zeta_9\ptl_{x_{13}}\\ &
&-\zeta_{12}\ptl_{x_{17}}-\zeta_{11}\ptl_{x_{16}}
+\zeta_{15}\ptl_{x_{19}}-\zeta_{18}\ptl_{x_{22}}-\zeta_{20}\ptl_{x_{24}}
+\zeta_{25}\ptl_{x_{27}},\hspace{2cm}(12.4.46)\end{eqnarray*}
\begin{eqnarray*}
\qquad P_{10}&=&-[E_{-\al_1},P_7]=
x_{10}D-\zeta_2\ptl_{x_4}-\zeta_3\ptl_{x_5}-\zeta_6\ptl_{x_9}
+\zeta_7\ptl_{x_{11}}\\ &
&-\zeta_{10}\ptl_{x_{14}}+\zeta_{14}\ptl_{x_{18}}
+\zeta_{17}\ptl_{x_{21}}-\zeta_{19}\ptl_{x_{22}}-\zeta_{23}\ptl_{x_{25}}
-\zeta_{24}\ptl_{x_{26}},\hspace{1.9cm}(12.4.47)\end{eqnarray*}\begin{eqnarray*}
\qquad P_{11}&=&[E_{-\al_5},P_9]=
x_{11}D-\zeta_1\ptl_{x_2}+\zeta_5\ptl_{x_8}+\zeta_7\ptl_{x_{10}}
-\zeta_9\ptl_{x_{12}}\\ &
&+\zeta_{11}\ptl_{x_{15}}+\zeta_{13}\ptl_{x_{17}}
-\zeta_{16}\ptl_{x_{19}}+\zeta_{18}\ptl_{x_{21}}+\zeta_{20}\ptl_{x_{23}}
+\zeta_{26}\ptl_{x_{27}},\hspace{1.9cm}(12.4.48)\end{eqnarray*}
 \begin{eqnarray*}
\qquad
P_{12}&=&-[E_{-\al_1},P_9]=x_{12}D-\zeta_2\ptl_{x_3}+\zeta_4\ptl_{x_5}+\zeta_6\ptl_{x_7}
-\zeta_9\ptl_{x_{11}}\\ &
&+\zeta_{12}\ptl_{x_{14}}-\zeta_{14}\ptl_{x_{16}}
-\zeta_{17}\ptl_{x_{19}}+\zeta_{21}\ptl_{x_{22}}+\zeta_{23}\ptl_{x_{24}}
-\zeta_{25}\ptl_{x_{26}},\hspace{2cm}(12.4.49)\end{eqnarray*}
 \begin{eqnarray*}
\qquad
P_{13}&=&[E_{-\al_5},P_{12}]=x_{13}D-\zeta_2\ptl_{x_2}-\zeta_5\ptl_{x_5}-\zeta_7\ptl_{x_7}
+\zeta_9\ptl_{x_9}\\ &
&-\zeta_{13}\ptl_{x_{14}}+\zeta_{14}\ptl_{x_{15}}
+\zeta_{19}\ptl_{x_{19}}-\zeta_{21}\ptl_{x_{21}}-\zeta_{23}\ptl_{x_{23}}
-\zeta_{26}\ptl_{x_{26}},\hspace{2cm}(12.4.50)\end{eqnarray*}
\begin{eqnarray*}
\qquad P_{14}&=&[E_{-\al_6},P_{11}]=
x_{14}D-\zeta_1\ptl_{x_1}-\zeta_8\ptl_{x_8}-\zeta_{10}\ptl_{x_{10}}
+\zeta_{12}\ptl_{x_{12}}\\ &
&-\zeta_{13}\ptl_{x_{13}}-\zeta_{15}\ptl_{x_{15}}
+\zeta_{16}\ptl_{x_{16}}-\zeta_{18}\ptl_{x_{18}}-\zeta_{20}\ptl_{x_{20}}
+\zeta_{27}\ptl_{x_{27}},\hspace{2cm}(12.4.51)\end{eqnarray*}
\begin{eqnarray*}
\qquad
P_{15}&=&-[E_{-\al_3},P_{12}]=x_{15}D-\zeta_3\ptl_{x_3}-\zeta_4\ptl_{x_4}-\zeta_6\ptl_{x_6}
+\zeta_{11}\ptl_{x_{11}}\\ &
&+\zeta_{14}\ptl_{x_{13}}-\zeta_{15}\ptl_{x_{14}}
+\zeta_{17}\ptl_{x_{17}}-\zeta_{22}\ptl_{x_{22}}-\zeta_{24}\ptl_{x_{24}}
-\zeta_{25}\ptl_{x_{25}},\hspace{2cm}(12.4.52)\end{eqnarray*}
\begin{eqnarray*}
\qquad
P_{16}&=&[E_{-\al_5},P_{15}]=x_{16}D-\zeta_3\ptl_{x_2}+\zeta_5\ptl_{x_4}+\zeta_7\ptl_{x_6}
-\zeta_{11}\ptl_{x_9}\\ &
&-\zeta_{14}\ptl_{x_{12}}+\zeta_{16}\ptl_{x_{14}}
-\zeta_{19}\ptl_{x_{17}}+\zeta_{22}\ptl_{x_{21}}+\zeta_{24}\ptl_{x_{23}}
-\zeta_{26}\ptl_{x_{25}},\hspace{2cm}(12.4.53)\end{eqnarray*}
 \begin{eqnarray*}
\qquad P_{17}&=&=-[E_{-\al_1},P_{14}]=
x_{17}D-\zeta_2\ptl_{x_1}+\zeta_8\ptl_{x_5}+\zeta_{10}\ptl_{x_7}
-\zeta_{12}\ptl_{x_9}\\ &
&+\zeta_{13}\ptl_{x_{11}}+\zeta_{17}\ptl_{x_{15}}
-\zeta_{19}\ptl_{x_{16}}+\zeta_{21}\ptl_{x_{18}}+\zeta_{23}\ptl_{x_{20}}
-\zeta_{27}\ptl_{x_{26}},\hspace{2cm}(12.4.54)\end{eqnarray*}
  \begin{eqnarray*}
\qquad
P_{18}&=&[E_{-\al_4},P_{16}]=x_{18}D+\zeta_4\ptl_{x_2}+\zeta_5\ptl_{x_3}-\zeta_9\ptl_{x_6}
+\zeta_{11}\ptl_{x_7}\\ &
&+\zeta_{14}\ptl_{x_{10}}-\zeta_{18}\ptl_{x_{14}}
+\zeta_{21}\ptl_{x_{17}}-\zeta_{22}\ptl_{x_{19}}+\zeta_{25}\ptl_{x_{23}}
+\zeta_{26}\ptl_{x_{24}},\hspace{2cm}(12.4.55)\end{eqnarray*}
 \begin{eqnarray*}
\qquad
P_{19}&=&[E_{-\al_6},P_{16}]=x_{19}D-\zeta_3\ptl_{x_1}-\zeta_8\ptl_{x_4}-\zeta_{10}\ptl_{x_6}
+\zeta_{15}\ptl_{x_9}\\ &
&-\zeta_{17}\ptl_{x_{12}}-\zeta_{16}\ptl_{x_{11}}
+\zeta_{19}\ptl_{x_{13}}-\zeta_{22}\ptl_{x_{18}}-\zeta_{24}\ptl_{x_{20}}
-\zeta_{27}\ptl_{x_{25}},\hspace{1.8cm}(12.4.56)\end{eqnarray*}
  \begin{eqnarray*}
\qquad
P_{20}&=&[E_{-\al_2},P_{18}]=x_{20}D-\zeta_6\ptl_{x_2}-\zeta_7\ptl_{x_3}-\zeta_9\ptl_{x_4}
+\zeta_{11}\ptl_{x_5}\\ &
&+\zeta_{14}\ptl_{x_8}-\zeta_{20}\ptl_{x_{14}}
+\zeta_{23}\ptl_{x_{17}}-\zeta_{24}\ptl_{x_{19}}-\zeta_{25}\ptl_{x_{21}}
-\zeta_{26}\ptl_{x_{22}},\hspace{2cm}(12.4.57)\end{eqnarray*}
 \begin{eqnarray*}
\qquad
P_{21}&=&[E_{-\al_4},P_{19}]=x_{21}D+\zeta_4\ptl_{x_1}-\zeta_8\ptl_{x_3}+\zeta_{12}\ptl_{x_6}
-\zeta_{15}\ptl_{x_7}\\ &
&+\zeta_{17}\ptl_{x_{10}}+\zeta_{18}\ptl_{x_{11}}
-\zeta_{21}\ptl_{x_{13}}+\zeta_{22}\ptl_{x_{16}}-\zeta_{25}\ptl_{x_{20}}
+\zeta_{27}\ptl_{x_{24}},\hspace{1.9cm}(12.4.58)\end{eqnarray*}
  \begin{eqnarray*}
\qquad
P_{22}&=&[E_{-\al_5},P_{21}]=x_{22}D-\zeta_5\ptl_{x_1}-\zeta_8\ptl_{x_2}-\zeta_{13}\ptl_{x_6}
+\zeta_{16}\ptl_{x_7}\\ &
&-\zeta_{19}\ptl_{x_{10}}-\zeta_{18}\ptl_{x_9}
+\zeta_{21}\ptl_{x_{12}}-\zeta_{22}\ptl_{x_{15}}-\zeta_{26}\ptl_{x_{20}}
-\zeta_{27}\ptl_{x_{23}},\hspace{2cm}(12.4.59)\end{eqnarray*}
 \begin{eqnarray*}
\qquad
P_{23}&=&[E_{-\al_2},P_{21}]=x_{23}D-\zeta_6\ptl_{x_1}+\zeta_{10}\ptl_{x_3}+\zeta_{12}\ptl_{x_4}
-\zeta_{15}\ptl_{x_5}\\ &
&+\zeta_{17}\ptl_{x_8}+\zeta_{20}\ptl_{x_{11}}
-\zeta_{23}\ptl_{x_{13}}+\zeta_{24}\ptl_{x_{16}}+\zeta_{25}\ptl_{x_{18}}
-\zeta_{27}\ptl_{x_{22}},\hspace{2cm}(12.4.60)\end{eqnarray*}
 \begin{eqnarray*}
\qquad
P_{24}&=&[E_{-\al_2},P_{22}]=x_{24}D+\zeta_7\ptl_{x_1}+\zeta_{10}\ptl_{x_2}-\zeta_{13}\ptl_{x_4}
+\zeta_{16}\ptl_{x_5}\\ &
&-\zeta_{19}\ptl_{x_8}-\zeta_{20}\ptl_{x_9}
+\zeta_{23}\ptl_{x_{12}}-\zeta_{24}\ptl_{x_{15}}+\zeta_{26}\ptl_{x_{18}}
+\zeta_{27}\ptl_{x_{21}},\hspace{2cm}(12.4.61)\end{eqnarray*}
  \begin{eqnarray*}
\qquad
P_{25}&=&[E_{-\al_4},P_{24}]=x_{25}D-\zeta_9\ptl_{x_1}-\zeta_{12}\ptl_{x_2}-\zeta_{13}\ptl_{x_3}
-\zeta_{18}\ptl_{x_5}\\ &
&+\zeta_{21}\ptl_{x_8}+\zeta_{20}\ptl_{x_7}
-\zeta_{23}\ptl_{x_{10}}-\zeta_{25}\ptl_{x_{15}}-\zeta_{26}\ptl_{x_{16}}
-\zeta_{27}\ptl_{x_{19}},\hspace{2cm}(12.4.62)\end{eqnarray*}
  \begin{eqnarray*}
\qquad
P_{26}&=&[E_{-\al_3},P_{25}]=x_{26}D-\zeta_{11}\ptl_{x_1}-\zeta_{15}\ptl_{x_2}-\zeta_{16}\ptl_{x_3}
-\zeta_{18}\ptl_{x_4}\\ &
&+\zeta_{22}\ptl_{x_8}+\zeta_{20}\ptl_{x_6}
-\zeta_{24}\ptl_{x_{10}}-\zeta_{25}\ptl_{x_{12}}-\zeta_{26}\ptl_{x_{13}}
-\zeta_{27}\ptl_{x_{17}},\hspace{2.1cm}(12.4.63)\end{eqnarray*}
  \begin{eqnarray*}
\qquad
P_{27}&=&-[E_{-\al_1},P_{26}]=x_{27}D-\zeta_{14}\ptl_{x_1}+\zeta_{17}\ptl_{x_2}+\zeta_{19}\ptl_{x_3}
+\zeta_{21}\ptl_{x_4}\\ &
&-\zeta_{22}\ptl_{x_5}-\zeta_{23}\ptl_{x_6}
+\zeta_{24}\ptl_{x_7}+\zeta_{25}\ptl_{x_9}+\zeta_{26}\ptl_{x_{11}}
+\zeta_{27}\ptl_{x_{14}}.\hspace{2.4cm}(12.4.64)\end{eqnarray*}

Write
$${\msr P}=\sum_{i=1}^{27}\mbb FP_i,\qquad\mfk C_0=\mfk r({\msr G}^{E_6})+\mbb FD\eqno(12.4.65)$$
(cf. (11.1.18)-(11.1.54) and (12.4.9)) and
$$\mfk C={\msr P}+\mfk C_0+{\msr T}\eqno(12.4.66)$$ (cf. (12.4.10)). Then
we have:\psp

{\bf Theorem 12.4.1}. {\it The space ${\mfk C}$ forms a Lie
subalgebra of the Witt algebra $\mbb W_{27}$ (cf. (12.4.15)).
Moreover, the linear map $\vt$ determined by
$$\vt({\mfk a}_i)=\ptl_{x_i},\;\;\vt({\mfk b_i})=P_i,\;\;\vt(u)=u|_{\msr
A}\qquad\for\;\;i\in\ol{1,27},\;u\in{\msr G}_0\eqno(12.4.67)$$
 (cf. (11.1.18)-(11.1.54) and (12.4.8))
gives a Lie algebra isomorphism from ${\msr G}^{E_7}$ to $\mfk C$.}

{\it Proof}. Since $\msr T\cong {\msr G}_+$ as ${\msr
G}^{E_6}$-modules, we have
$${\msr G}_0+{\msr G}_+\stl{\vt}{\cong}{\mfk C}_0+{\msr T}\eqno(12.4.68)$$ as Lie algebras.
Recall that $U({\msr G})$ stands for the universal enveloping
algebra of a Lie algebra ${\msr G}$. Note that
$${\msr B}_-={\msr G}_0+{\msr G}_-,\qquad{\msr B}_+={\msr G}_0+{\msr G}_+\eqno(12.4.69)$$
are also  Lie subalgebras of ${\msr G}^{E_7}$ and
$${\msr G}^{E_7}={\msr B}_-\oplus {\msr G}_+={\msr G}_-\oplus {\msr B}_+.\eqno(12.4.70)$$
We define a one-dimensional ${\msr B}_-$-module $\mbb{C}u_0$ by
$$w(u_0)=0\qquad\for\;\;w\in\msr G_-+{\msr G}^{E_6},\;\;\widehat\al(u_0)=54u_0\eqno(12.4.71)$$
(cf. (12.4.9)). Let
$$\Psi=U({\msr G}^{E_7})\otimes_{{\msr B}_-}\mbb{C}u_0\cong U({\msr G}_+)\otimes_\mbb{F} \mbb{F}u_0\eqno(12.4.72)$$
be the induced ${\msr G}^{E_7}$-module.

Recall that $\mbb{N}$ is the set of nonnegative integers. Let
$$\wht{\msr
A}=\mbb{F}[\ptl_{x_1},\ptl_{x_2},...,\ptl_{x_{27}}].\eqno(12.4.73)$$
We define an action of the associative algebra $\mbb{A}$ (cf.
(11.2.35)) on $\wht{\msr A}$ by
$$\ptl_{x_i}(\prod_{j=1}^{27}\ptl_{x_j}^{\be_j})=\ptl_{x_i}^{\be_i+1}\prod_{i\neq
j\in\ol{1,27}}\prod_{j=1}^{27}\ptl_{x_j}^{\be_j}\eqno(12.4.74)$$ and
$$x_i(\prod_{j=1}^{27}\ptl_{x_j}^{\be_j})=-\be_i\ptl_{x_i}^{\be_i-1}\prod_{i\neq
j\in\ol{1,27}}\prod_{j=1}^{27}\ptl_{x_j}^{\be_j}\eqno(12.4.75)$$ for
$i\in\ol{1,27}$. Since
$$[-x_i,\ptl_{x_j}]=[\ptl_{x_i},x_j]=\dlt_{i,j}\qquad\for\;\;i,j\in\ol{1,27},\eqno(12.4.76)$$
the above action gives an associative algebra representation of
$\mbb{A}$. Thus it also gives a Lie algebra representation of
$\mbb{A}$, whose Lie bracket is the commutator.
 It is straightforward to verify that
$$[d|_{\wht{\msr A}},\ptl|_{\wht{\msr A}}]=[d,\ptl]|_{\wht{\msr
A}}\qquad\for\;\;d\in{\mfk C}_0,\;\ptl\in{\msr T}.\eqno(12.4.77)$$

Define linear map $\vs: \Psi\rta {\msr A}'$ by
$$\vs(\prod_{i=1}^{27}\xi_i^{\ell_i}\otimes
u_0)=\prod_{i=1}^{27}\ptl_{x_i}^{\ell_i}\qquad(\ell_1,...,\ell_{27})\in\mbb{N}^{27}.\eqno(12.4.78)$$
According to (11.1.18)-(11.1.54), (12.4.74) and (12.4.75),
$$D(1)=-27,\;\;d(1)=0\qquad\for\;\;d\in \mfk r(\msr G^{E_6}).\eqno(12.4.79)$$
Moreover, (12.4.71), (12.4.72),  (12.4.74), (12.4.75) and (12.4.79)
imply
$$\vs(\xi(v))=\vt(\xi)\vs(v)\qquad\for\;\;\xi\in{\msr
G}_0,\;v\in\Psi.\eqno(12.4.80)$$

Now (12.4.72) and (12.4.74) imply
$$\vs(w(u))=\vt(w)(\vs(u))\qquad\for\;\;w\in{\msr
B}_+,\;\;u\in\Psi.\eqno(12.4.81)$$ Thus  we have
$$\vs w|_{\Psi}\vs^{-1}=\vt(w)|_{\wht{\msr A}}\qquad\for\;\;w\in{\msr
B}_+.\eqno(12.4.82)$$ On the other hand, the linear map
$$\psi(v)=\vs v|_{\Psi}\vs^{-1}\qquad\for\;\; v\in{\msr
G}^{E_7}\eqno(12.4.83)$$ is a Lie algebra monomorphism from ${\msr
G}^{E_7}$ to $\mbb{A}|_{\wht{\msr A}}$. According to (12.4.38) and
(12.4.77),
$$\psi(\mfk b_1)=P_1|_{\wht{\msr A}}.\eqno(12.4.84)$$
By the constructions of $P_2,...,P_{27}$ in (12.4.39)-(12.4.64), we
have
$$\psi(\mfk b_i)=P_i|_{\wht{\msr
A}}\qquad\for\;\;i\in\ol{2,27}.\eqno(12.4.85)$$ Therefore, we have
$$\psi(v)=\vt(v)|_{\wht{\msr A}}\qquad\for\;\;v\in{\msr
G}^{E_7}.\eqno(12.4.86)$$ In particular, ${\mfk C}|_{\wht{\msr
A}}=\vt({\msr G}^{E_7})|_{\wht{\msr A}}=\psi({\msr G}^{E_7})$ forms
a Lie algebra. Since the linear map $d\mapsto d|_{\wht{\msr A}}$ for
$d\in {\mfk C}$ is injective, we have that ${\mfk C}$ forms a Lie
subalgebra of $\mbb{A}$ and $\vt$ is a Lie algebra
isomorphism.$\qquad\Box$ \psp

By the above theorem, a Lie group of type $E_7$ is generated by the
linear transformations $\{e^{b\mfk r(E_\al)}\mid b\in\mbb{R},\;
\al\in\Phi_{E_6}\}$ associated with (11.1.18)-(11.1.54), the real
translations and dilations in $\mbb R^{27}$ with $x_r$ as the $r$th
coordinate function, and the fractional transformations
$\{e^{bP_s}:(x_1,x_2,...,x_{27})\mapsto
(e^{bP_s}(x_1),e^{bP_s}(x_2),...,e^{bP_s}(x_{27}))\mid
b\in\mbb{R},\; s\in\ol{1,27}\}$, where $e^{bP_s}(x_i)$ are of the
 forms as the following forms of the case $s=1$:
$$e^{bP_1}(x_i)=\frac{x_i}{1-bx_1},\;\;i\in\{\ol{1,13},15,16,18,20\},\eqno(12.4.87)$$
$$e^{bP_1}(x_{14})
=x_{14}-\frac{b(x_2x_{11}+x_3x_9+x_4x_7-x_5x_6)}{1-bx_1},\eqno(12.4.88)$$
$$e^{bP_1}(x_{17})
=x_{17}-\frac{b(x_2x_{13}+x_3x_{12}+x_4x_{10}-x_6x_8)}{1-bx_1},\eqno(12.4.89)$$
$$e^{bP_1}(x_{19})
=x_{19}-\frac{b(x_2x_{16}+x_3x_{15}+x_5x_{10}-x_7x_8)}{1-bx_1},\eqno(12.4.90)$$
$$e^{bP_1}(x_{21})
=x_{21}-\frac{b(x_2x_{18}-x_4x_{15}+x_5x_{12}
-x_8x_9)}{1-bx_1},\eqno(12.4.91)$$
$$e^{bP_1}(x_{22})
=x_{22}+\frac{b(x_3x_{18}+x_4x_{16}-x_5x_{13}
+x_8x_{11})}{1-bx_1},\eqno(12.4.92)$$
$$e^{bP_1}(x_{23})
=x_{23}-\frac{b(x_2x_{20}+x_6x_{15}-x_7x_{12}
+x_9x_{10})}{1-bx_1},\eqno(12.4.93)$$
$$e^{bP_1}(x_{24})
=x_{24}+\frac{b(x_3x_{20}-x_6x_{16}+x_7x_{13}
-x_{10}x_{11})}{1-bx_1},\eqno(12.4.94)$$
$$e^{bP_1}(x_{25})
=x_{25}-\frac{b(x_4x_{20}+x_6x_{18}-x_9x_{13}
+x_{11}x_{12})}{1-bx_1},\eqno(12.4.95)$$
$$e^{bP_1}(x_{26})
=x_{26}+\frac{b(x_5x_{20}+x_7x_{18}-x_9x_{16}+x_{11}x_{15})}{1-bx_1},\eqno(12.4.96)$$
$$e^{bP_1}(x_{27})
=x_{27}+\frac{b(x_8x_{20}+x_{10}x_{18}-x_{12}x_{16}+x_{13}x_{15}
)}{1-bx_1}\eqno(12.4.97)$$ by (11.2.1), (11.2.3)-(11.2.8),
(11.2.10), (11.2.12), (11.2.15) and (12.4.20).

\section{Functor from $E_6$-Mod to $E_7$-Mod}

In this section, we construct a new functor from $E_6$-{\bf Mod} to
$E_7$-{\bf Mod}.

Note that
$${\msr G}^{E_6}_{\msr A}=(\bigoplus_{i=1}^6{\msr A}\al_i)\oplus\bigoplus_{\al\in
\Phi_{E_6}}{\msr A}E_{\al}\eqno(12.5.1)$$ forms a Lie algebra with
the Lie bracket:
$$[fu_1,gu_2]=fg[u_1,u_2]\qquad\for\;\;f,g\in{\msr A},\;u_1,u_2\in
{\msr G}^{E_6}.\eqno(12.5.2)$$ Moreover, we define the Lie algebra
$${\msr K}={\msr G}^{E_6}_{\msr A}\oplus {\msr A}\kappa\eqno(12.5.3)$$
with the Lie bracket:
$$[\xi_1+f\kappa,\xi_2+g\kappa]=[\xi_1,\xi_2]\qquad\for\;\;\xi_1,\xi_2\in
{\msr G}^{E_6}_{\msr A},\;f,g\in{\msr A}.\eqno(12.5.4)$$  Recall the
Witt algebra $\mbb W_{27}=\sum_{i=1}^{27}{\msr A}\ptl_{x_i}$, and
Shen's monomorphism $\Im$ from the Lie algebra $\mbb W_{27}$ to the
Lie algebra of semi-product $\mbb W_{27}+gl(27,{\msr A})$ (cf.
(6.7.9)-(6.7.12)) defined by
$$\Im(\sum_{i=1}^{27}f_i\ptl_{x_i})=\sum_{i=1}^{27}f_i\ptl_{x_i}+\Im_1(\sum_{i=1}^{27}f_i\ptl_{x_i}),\;\;
\Im_1(\sum_{i=1}^{27}f_i\ptl_{x_i})=\sum_{i,j=1}^{27}\ptl_{x_i}(f_j)E_{i,j}.
\eqno(12.5.5)$$ According to our construction of $P_1$-$P_{27}$ in
(12.4.20) and (12.4.39)-(12.4.64),
$$\Im_1(P_i)=\sum_{r=1}^{27}x_r\Im_1(\mfk r([\mfk a_r,\mfk b_i]))\qquad\for\;\;i\in\ol{1,27}.\eqno(12.5.6)$$

On the other hand,
$$\wht{\msr K}=\mbb W_{27}\oplus {\msr K}\eqno(12.5.7)$$
becomes a Lie algebra with the Lie bracket
\begin{eqnarray*}&&[d_1+f_1u_1+f_2\kappa,d_2+g_1u_2+g_2\kappa]\\
&=&[d_1,d_2]+f_1g_1[u_1,u_2]+d_1(g_2)u_2
-d_2(g_1)u_1+(d_1(g_2)-d_2(g_1))\kappa\hspace{2.3cm}(12.5.8)\end{eqnarray*}
for $f_1,f_2,g_1,g_2\in{\msr A},\;u_1,u_2\in {\msr G}^{E_6}$ and
$d_1,d_2\in\mbb W_{27}$. Note
$${\msr G}_0={\msr
G}^{E_6}\oplus\mbb F\wht\al\eqno(12.5.9)$$ by (12.4.2) and (12.4.5).
So there exists a Lie algebra monomorphism $\mfk k:{\msr G}_0\rta
{\msr K}$ determined by
$$\mfk k(\wht\al)=3\kappa,\;\;\mfk k(u)=u\qquad\for\;\;u\in{\msr
G}^{E_6}.\eqno(12.5.10)$$ Since $\Im$ is a Lie algebra monomorphism,
our construction of $P_1$-$P_{27}$ in (12.4.20) and
(12.4.39)-(12.4.64) show that we have a Lie algebra monomorphism
$\iota: {\msr G}^{E_7}\rta \wht{\msr K}$ given by
$$\iota(u)=\mfk r(u)+\mfk k(u)\qquad\for\;\;u\in{\msr
G}_0,\eqno(12.5.11)$$
$$\iota(\mfk a_i)=\ptl_{x_i},\;\;\iota(\mfk b_i)=P_i+\sum_{r=1}^{27}x_r\mfk k([\mfk a_r,\mfk b_i])\qquad\for
\;\;i\in\ol{1,27}.\eqno(12.5.12)$$

According to (4.4.18), (4.4.25), (11.1.4)-(11.1.12) and (12.4.1),
\begin{eqnarray*} \iota(\mfk b_1)&=&P_1+x_2E_{\al_6}
+x_3E_{(0,0,0,0,1,1)}+x_4E_{(0,0,0,1,1,1)}
+x_5E_{(0,0,1,1,1,1)}+x_6E_{(0,1,0,1,1,1)}\\
& &-x_1\mfk
k(\al_7)+x_7E_{(0,1,1,1,1,1)}+x_8E_{(1,0,1,1,1,1)}+x_9E_{(0,1,1,2,1,1)}+x_{10}E_{(1,1,1,1,1,1)}
\\& & +x_{11}E_{(0,1,1,2,2,1)}+x_{12}E_{(1,1,1,2,1,1)}+x_{13}E_{(1,1,1,2,2,1)}+x_{15}E_{(1,1,2,2,1,1)}
\\& &
+x_{16}E_{(1,1,2,2,2,1)}+x_{18}E_{(1,1,2,3,2,1)}+x_{20}E_{(1,2,2,3,2,1)}.
\hspace{3.9cm}(12.5.13)\end{eqnarray*} Moreover, (4.4.23), (11.1.23)
and (11.1.56) yield
\begin{eqnarray*}
\iota(\mfk b_2)&=&-\iota([E_{-\al_6},\mfk
b_1])=-[\iota(E_{-\al_6}),\iota(\mfk b_1)]=
-[\mfk r(E_{-\al_6})+E_{-\al_6},\iota(\mfk b_1)]\\
&=&P_2-x_1E_{-\al_6} +x_3E_{\al_5}+x_4E_{(0,0,0,1,1)}
+x_5E_{(0,0,1,1,1)}+x_6E_{(0,1,0,1,1)}\\
&&-x_2\mfk
k(\al_{(0,0,0,0,0,1,1)})+x_7E_{(0,1,1,1,1)}+x_8E_{(1,0,1,1,1)}+x_9E_{(0,1,1,2,1)}+x_{10}E_{(1,1,1,1,1)}
\\& & +x_{12}E_{(1,1,1,2,1)}-x_{14}E_{(0,1,1,2,2,1)}+x_{15}E_{(1,1,2,2,1)}-x_{17}E_{(1,1,1,2,2,1)}
\\& &
-x_{19}E_{(1,1,2,2,2,1)}-x_{21}E_{(1,1,2,3,2,1)}-x_{23}E_{(1,2,2,3,2,1)}.
\hspace{3.9cm}(12.5.14)\end{eqnarray*}  Similarly, we have:
\begin{eqnarray*}
\iota(\mfk b_3) &=&P_3-x_1E'_{(0,0,0,0,1,1)}
-x_2E_{-\al_5}+x_4E_{\al_4}
+x_5E_{(0,0,1,1)}+x_6E_{(0,1,0,1)}\\
&&-x_3\mfk
k(\al_{(0,0,0,0,1,1,1)})+x_7E_{(0,1,1,1)}+x_8E_{(1,0,1,1)}+x_{10}E_{(1,1,1,1)}-x_{11}E_{(0,1,1,2,1)}
\\& & -x_{13}E_{(1,1,1,2,1)}-x_{14}E_{(0,1,1,2,1,1)}-x_{16}E_{(1,1,2,2,1)}-x_{17}E_{(1,1,1,2,1,1)}
\\& &
-x_{19}E_{(1,1,2,2,1,1)}+x_{22}E_{(1,1,2,3,2,1)}+x_{24}E_{(1,2,2,3,2,1)},
\hspace{3.8cm}(12.5.15)\end{eqnarray*}
\begin{eqnarray*}
\iota(\mfk b_4) &=&P_4-x_1E'_{(0,0,0,1,1,1)}
-x_2E'_{(0,0,0,1,1)}-x_3E_{-\al_4}
+x_5E_{\al_3}+x_6E_{\al_2}\\
&&-x_4\mfk
k(\al_{(0,0,0,1,1,1,1)})-x_9E_{(0,1,1,1)}+x_8E_{(1,0,1)}-x_{12}E_{(1,1,1,1)}-x_{11}E_{(0,1,1,1,1)}
\\& & -x_{13}E_{(1,1,1,1,1)}-x_{14}E_{(0,1,1,1,1,1)}-x_{17}E_{(1,1,1,1,1,1)}+x_{18}E_{(1,1,2,2,1)}
\\& &
+x_{21}E_{(1,1,2,2,1,1)}+x_{22}E_{(1,1,2,2,2,1)}-x_{25}E_{(1,2,2,3,2,1)},
\hspace{3.8cm}(12.5.16)\end{eqnarray*}
\begin{eqnarray*}
\iota(\mfk b_5) &=&P_5-x_1E'_{(0,0,1,1,1,1)}
-x_2E'_{(0,0,1,1,1)}-x_3E'_{(0,0,1,1)}
-x_4E_{-\al_3}+x_7E_{\al_2}\\
&&-x_5\mfk
k(\al_{(0,0,1,1,1,1,1)})+x_8E_{\al_1}+x_9E_{(0,1,0,1)}-x_{15}E_{(1,1,1,1)}+x_{11}E_{(0,1,0,1,1)}
\\& & -x_{16}E_{(1,1,1,1,1)}+x_{14}E_{(0,1,0,1,1,1)}-x_{19}E_{(1,1,1,1,1,1)}-x_{18}E_{(1,1,1,2,1)}
\\& &
-x_{21}E_{(1,1,1,2,1,1)}-x_{22}E_{(1,1,1,2,2,1)}+x_{26}E_{(1,2,2,3,2,1)},
\hspace{3.8cm}(12.5.17)\end{eqnarray*}
\begin{eqnarray*}
\iota(\mfk b_6) &=&P_6-x_1E'_{(0,1,0,1,1,1)}
-x_2E'_{(0,1,0,1,1)}-x_3E'_{(0,1,0,1)}
-x_4E_{-\al_2}+x_7E_{\al_3}\\
&&-x_6\mfk
k(\al_{(0,1,0,1,1,1,1)})+x_9E_{(0,0,1,1)}+x_{10}E_{(1,0,1)}+x_{12}E_{(1,0,1,1)}+x_{11}E_{(0,0,1,1,1)}
\\& & +x_{13}E_{(1,0,1,1,1)}+x_{14}E_{(0,0,1,1,1,1)}+x_{17}E_{(1,0,1,1,1,1)}-x_{20}E_{(1,1,2,2,1)}
\\& &
-x_{23}E_{(1,1,2,2,1,1)}-x_{24}E_{(1,1,2,2,2,1)}-x_{25}E_{(1,1,2,3,2,1)},
\hspace{3.9cm}(12.5.18)\end{eqnarray*}
\begin{eqnarray*}
\iota(\mfk b_7) &=&P_7-x_1E'_{(0,1,1,1,1,1)}
-x_2E'_{(0,1,1,1,1)}-x_3E'_{(0,1,1,1)}
-x_5E_{-\al_2}-x_6E_{-\al_3}\\
&&-x_7\mfk
k(\al_{(0,1,1,1,1,1,1)})-x_9E_{\al_4}+x_{10}E_{\al_1}+x_{15}E_{(1,0,1,1)}-x_{11}E_{(0,0,0,1,1)}
\\& & +x_{16}E_{(1,0,1,1,1)}-x_{14}E_{(0,0,0,1,1,1)}+x_{19}E_{(1,0,1,1,1,1)}+x_{20}E_{(1,1,1,2,1)}
\\& &
+x_{23}E_{(1,1,1,2,1,1)}+x_{24}E_{(1,1,1,2,2,1)}+x_{26}E_{(1,1,2,3,2,1)},
\hspace{3.8cm}(12.5.19)\end{eqnarray*}
\begin{eqnarray*}
\iota(\mfk b_8) &=&P_8-x_1E'_{(1,0,1,1,1,1)}
-x_2E'_{(1,0,1,1,1)}-x_3E'_{(1,0,1,1)}
-x_4E'_{(1,0,1)}-x_5E_{-\al_1}\\
&&-x_8\mfk
k(\al_{(1,0,1,1,1,1,1)})+x_{10}E_{\al_2}+x_{12}E_{(0,1,0,1)}+x_{15}E_{(0,1,1,1)}+x_{13}E_{(0,1,0,1,1)}
\\& & +x_{16}E_{(0,1,1,1,1)}+x_{17}E_{(0,1,0,1,1,1)}+x_{19}E_{(0,1,1,1,1,1)}+x_{18}E_{(0,1,1,2,1)}
\\& &
+x_{21}E_{(0,1,1,2,1,1)}+x_{22}E_{(0,1,1,2,2,1)}+x_{27}E_{(1,2,2,3,2,1)},
\hspace{3.8cm}(12.5.20)\end{eqnarray*}
\begin{eqnarray*}
\iota(\mfk b_9) &=&P_9-x_1E'_{(0,1,1,2,1,1)}
-x_2E'_{(0,1,1,2,1)}+x_4E'_{(0,1,1,1)}
-x_5E'_{(0,1,0,1)}-x_6E'_{(0,0,1,1)}\\
&&-x_9\mfk
k(\al_{(0,1,1,2,1,1,1)})+x_7E_{-\al_4}+x_{12}E_{\al_1}-x_{15}E_{(1,0,1)}-x_{11}E_{\al_5}
\\& & +x_{18}E_{(1,0,1,1,1)}-x_{14}E_{(0,0,0,0,1,1)}+x_{21}E_{(1,0,1,1,1,1)}-x_{20}E_{(1,1,1,1,1)}
\\& &
-x_{23}E_{(1,1,1,1,1,1)}+x_{25}E_{(1,1,1,2,2,1)}-x_{26}E_{(1,1,2,2,2,1)},
\hspace{3.9cm}(12.5.21)\end{eqnarray*}
\begin{eqnarray*}
\iota(\mfk b_{10}) &=&P_{10}-x_1E'_{(1,1,1,1,1,1)}
-x_2E'_{(1,1,1,1,1)}-x_3E'_{(1,1,1,1)}
-x_8E_{-\al_2}-x_6E'_{(1,0,1)}\\
&&-x_{10}\mfk
k(\al_{(1,1,1,1,1,1,1)})-x_7E_{-\al_1}-x_{12}E_{\al_4}-x_{15}E_{(0,0,1,1)}-x_{13}E_{(0,0,0,1,1)}
\\& & -x_{16}E_{(0,0,1,1,1)}-x_{17}E_{(0,0,0,1,1,1)}-x_{19}E_{(0,0,1,1,1,1)}-x_{20}E_{(0,1,1,2,1)}
\\& &
-x_{23}E_{(0,1,1,2,1,1)}-x_{24}E_{(0,1,1,2,2,1)}+x_{27}E_{(1,1,2,3,2,1)},
\hspace{3.8cm}(12.5.22)\end{eqnarray*}
\begin{eqnarray*}
\iota(\mfk b_{11}) &=&P_{11}-x_1E'_{(0,1,1,2,2,1)}
+x_3E'_{(0,1,1,2,1)}+x_4E'_{(0,1,1,1,1)}
-x_5E'_{(0,1,0,1,1)}\\
&&-x_{11}\mfk
k(\al_{(0,1,1,2,2,1,1)})-x_6E'_{(0,0,1,1,1)}+x_7E'_{(0,0,0,1,1)}+x_9E_{-\al_5}\\&
&+x_{13}E_{\al_1}-x_{16}E_{(1,0,1)}
 -x_{18}E_{(1,0,1,1)}-x_{14}E_{\al_6}+x_{20}E_{(1,1,1,1)}\\& &
+x_{22}E_{(1,0,1,1,1,1)}-x_{24}E_{(1,1,1,1,1,1)}-x_{25}E_{(1,1,1,2,1,1)}+x_{26}E_{(1,1,2,2,1,1)},
\hspace{0.8cm}(12.5.23)\end{eqnarray*}
\begin{eqnarray*}
\iota(\mfk b_{12}) &=&P_{12}-x_1E'_{(1,1,1,2,1,1)}
-x_2E'_{(1,1,1,2,1)}+x_4E'_{(1,1,1,1)}-x_6E'_{(1,0,1,1)}
\\
&&-x_{12}\mfk
k(\al_{(1,1,1,2,1,1,1)})-x_8E'_{(0,1,0,1)}-x_9E_{-\al_1}+x_{10}E_{-\al_4}\\&
&-x_{13}E_{\al_5}+x_{15}E_{\al_3}
 -x_{17}E_{(0,0,0,0,1,1)}-x_{18}E_{(0,0,1,1,1)}+x_{20}E_{(0,1,1,1,1)}\\& &-x_{21}E_{(0,0,1,1,1,1)}
+x_{23}E_{(0,1,1,1,1,1)}-x_{25}E_{(0,1,1,2,2,1)}-x_{27}E_{(1,1,2,2,2,1)},
\hspace{0.8cm}(12.5.24)\end{eqnarray*}
 \begin{eqnarray*}
\iota(\mfk b_{13}) &=&P_{13}-x_1E'_{(1,1,1,2,2,1)}
+x_3E'_{(1,1,1,2,1)}+x_4E'_{(1,1,1,1,1)}
-x_8E'_{(0,1,0,1,1)}\\
&&-x_{13}\mfk
k(\al_{(1,1,1,2,2,1,1)})-x_6E'_{(1,0,1,1,1)}+x_{10}E'_{(0,0,0,1,1)}-x_{11}E_{-\al_1}\\&
&+x_{12}E_{-\al_5}+x_{16}E_{\al_3}
 -x_{17}E_{\al_6}+x_{18}E_{(0,0,1,1)}-x_{20}E_{(0,1,1,1)}\\& &-x_{22}E_{(0,0,1,1,1,1)}
+x_{24}E_{(0,1,1,1,1,1)}+x_{25}E_{(0,1,1,2,1,1)}+x_{27}E_{(1,1,2,2,1,1)},
\hspace{0.8cm}(12.5.25)\end{eqnarray*}
\begin{eqnarray*}
\iota(\mfk b_{14}) &=&P_{14}+x_2E'_{(0,1,1,2,2,1)}
+x_3E'_{(0,1,1,2,1,1)}+x_4E'_{(0,1,1,1,1,1)}
-x_5E'_{(0,1,0,1,1,1)}\\
&&-x_{14}\mfk
k(\al_{(0,1,1,2,2,2,1)})-x_6E'_{(0,0,1,1,1,1)}+x_7E'_{(0,0,0,1,1,1)}+x_9E'_{(0,0,0,0,1,1)}\\&
&+x_{11}E_{-\al_6}+x_{17}E_{\al_1}
 -x_{19}E_{(1,0,1)}-x_{21}E_{(1,0,1,1)}-x_{22}E_{(1,0,1,1,1)}\\& &+x_{23}E_{(1,1,1,1)}
+x_{24}E_{(1,1,1,1,1)}+x_{25}E_{(1,1,1,2,1)}-x_{26}E_{(1,1,2,2,1)},
\hspace{2cm}(12.5.26)\end{eqnarray*}
\begin{eqnarray*}
\iota(\mfk b_{15}) &=&P_{15}-x_1E'_{(1,1,2,2,1,1)}
-x_2E'_{(1,1,2,2,1)}+x_5E'_{(1,1,1,1)}
-x_7E'_{(1,0,1,1)}\\
&&-x_{15}\mfk
k(\al_{(1,1,2,2,1,1,1)})-x_8E'_{(0,1,1,1)}+x_9E'_{(1,0,1)}+x_{10}E'_{(0,0,1,1)}\\&
&-x_{12}E_{-\al_3}-x_{16}E_{\al_5} +x_{18}E_{(0,0,0,1,1)}
-x_{19}E_{(0,0,0,0,1,1)}-x_{20}E_{(0,1,0,1,1)}\\&
&+x_{21}E_{(0,0,0,1,1,1)}
-x_{23}E_{(0,1,0,1,1,1)}+x_{26}E_{(0,1,1,2,2,1)}+x_{27}E_{(1,1,1,2,2,1)},
\hspace{0.8cm}(12.5.27)\end{eqnarray*}
 \begin{eqnarray*}\iota(\mfk b_{16})
&=&P_{16}-x_1E'_{(1,1,2,2,2,1)}
+x_3E'_{(1,1,2,2,1)}+x_5E'_{(1,1,1,1,1)}
-x_7E'_{(1,0,1,1,1)}\\
&&-x_{16}\mfk
k(\al_{(1,1,2,2,2,1,1)})-x_8E'_{(0,1,1,1,1)}+x_{10}E'_{(0,0,1,1,1)}+x_{11}E'_{(1,0,1)}\\&
&-x_{13}E_{-\al_3}+x_{15}E_{-\al_5}
 -x_{18}E_{\al_4}-x_{19}E_{\al_6}+x_{20}E_{(0,1,0,1)}+x_{22}E_{(0,0,0,1,1,1)}\\& &
-x_{24}E_{(0,1,0,1,1,1)}-x_{26}E_{(0,1,1,2,1,1)}-x_{27}E_{(1,1,1,2,1,1)},
\hspace{3.8cm}(12.5.28)\end{eqnarray*}
\begin{eqnarray*}
\iota(\mfk b_{17}) &=&P_{17}+x_2E'_{(1,1,1,2,2,1)}
+x_3E'_{(1,1,1,2,1,1)}+x_4E'_{(1,1,1,1,1,1)}
-x_6E'_{(1,0,1,1,1,1)}\\
&&-x_{17}\mfk
k(\al_{(1,1,1,2,2,2,1)})-x_8E'_{(0,1,0,1,1,1)}+x_{10}E'_{(0,0,0,1,1,1)}-x_{14}E_{-\al_1}\\&
&+x_{12}E'_{(0,0,0,0,1,1)}+x_{13}E_{-\al_6}+x_{19}E_{\al_3}
 +x_{21}E_{(0,0,1,1)}-x_{23}E_{(0,1,1,1)}\\& &+x_{22}E_{(0,0,1,1,1)}
-x_{24}E_{(0,1,1,1,1)}-x_{25}E_{(0,1,1,2,1)}-x_{27}E_{(1,1,2,2,1)},
\hspace{1.7cm}(12.5.29)\end{eqnarray*}
\begin{eqnarray*}\iota(\mfk b_{18})
&=&P_{18}-x_1E'_{(1,1,2,3,2,1)}
-x_4E'_{(1,1,2,2,1)}+x_5E'_{(1,1,1,2,1)}
-x_8E'_{(0,1,1,2,1)}\\
&&-x_{18}\mfk
k(\al_{(1,1,2,3,2,1,1)})-x_9E'_{(1,0,1,1,1)}+x_{11}E'_{(1,0,1,1)}+x_{12}E'_{(0,0,1,1,1)}\\&
&-x_{13}E'_{(0,0,1,1)}-x_{15}E'_{(0,0,0,1,1)}
 +x_{16}E_{-\al_4}-x_{21}E_{\al_6}+x_{22}E_{(0,0,0,0,1,1)}\\& &
-x_{20}E_{\al_2}-x_{25}E_{(0,1,0,1,1,1)}+x_{26}E_{(0,1,1,1,1,1)}+x_{27}E_{(1,1,1,1,1,1)},
\hspace{2cm}(12.5.30)\end{eqnarray*}
 \begin{eqnarray*}\iota(\mfk b_{19})
&=&P_{19}+x_2E'_{(1,1,2,2,2,1)}
+x_3E'_{(1,1,2,2,1,1)}+x_5E'_{(1,1,1,1,1,1)}
-x_7E'_{(1,0,1,1,1,1)}\\
&&-x_{19}\mfk
k(\al_{(1,1,2,2,2,2,1)})-x_8E'_{(0,1,1,1,1,1)}+x_{10}E'_{(0,0,1,1,1,1)}+x_{14}E'_{(1,0,1)}\\&
&+x_{15}E'_{(0,0,0,0,1,1)}+x_{16}E_{-\al_6}-x_{17}E_{-\al_3}
 -x_{21}E_{\al_4}+x_{23}E_{(0,1,0,1)}\\& &-x_{22}E_{(0,0,0,1,1)}
+x_{24}E_{(0,1,0,1,1)}+x_{26}E_{(0,1,1,2,1)}+x_{27}E_{(1,1,1,2,1)},
\hspace{1.8cm}(12.5.31)\end{eqnarray*}
 \begin{eqnarray*}\iota(\mfk b_{20})
&=&P_{20}-x_1E'_{(1,2,2,3,2,1)}
+x_6E'_{(1,1,2,2,1)}-x_7E'_{(1,1,1,2,1)}
+x_{10}E'_{(0,1,1,2,1)}\\
&&-x_{20}\mfk
k(\al_{(1,2,2,3,2,1,1)})+x_9E'_{(1,1,1,1,1)}-x_{11}E'_{(1,1,1,1)}-x_{12}E'_{(0,1,1,1,1)}\\&
&+x_{13}E'_{(0,1,1,1)}+x_{15}E'_{(0,1,0,1,1)}
 -x_{16}E_{(0,1,0,1)}+x_{18}E_{-\al_2}-x_{23}E_{\al_6}\\& &+x_{24}E_{(0,0,0,0,1,1)}
-x_{25}E_{(0,0,0,1,1,1)}+x_{26}E_{(0,0,1,1,1,1)}+x_{27}E_{(1,0,1,1,1,1)},
\hspace{0.9cm}(12.5.32)\end{eqnarray*}
\begin{eqnarray*}\iota(\mfk b_{21})
&=&P_{21}+x_2E'_{(1,1,2,3,2,1)}
-x_4E'_{(1,1,2,2,1,1)}+x_5E'_{(1,1,1,2,1,1)}
-x_8E'_{(0,1,1,2,1,1)}\\
&&-x_{21}\mfk
k(\al_{(1,1,2,3,2,2,1)})-x_9E'_{(1,0,1,1,1,1)}+x_{12}E'_{(0,0,1,1,1,1)}+x_{14}E'_{(1,0,1,1)}\\&
&-x_{15}E'_{(0,0,0,1,1,1)}-x_{17}E'_{(0,0,1,1)}+x_{18}E_{-\al_6}
 +x_{19}E_{-\al_4}-x_{23}E_{\al_2}\\& &-x_{22}E_{\al_5}
+x_{25}E_{(0,1,0,1,1)}-x_{26}E_{(0,1,1,1,1)}-x_{27}E_{(1,1,1,1,1)},
\hspace{2.8cm}(12.5.33)\end{eqnarray*}
\begin{eqnarray*}\iota(\mfk b_{22})
&=&P_{22}-x_3E'_{(1,1,2,3,2,1)} -x_4E'_{(1,1,2,2,2,1)} -x_{22}\mfk
k(\al_{(1,1,2,3,3,2,1)})+x_5E'_{(1,1,1,2,2,1)}\\
&&-x_8E'_{(0,1,1,2,2,1)}-x_{11}E'_{(1,0,1,1,1,1)}+x_{13}E'_{(0,0,1,1,1,1)}+x_{14}E'_{(1,0,1,1,1)}\\&
&-x_{16}E'_{(0,0,0,1,1,1)}-x_{17}E'_{(0,0,1,1.1)}-x_{18}E'_{(0,0,0,0,1,1)}
 +x_{19}E'_{(0,0,0,1,1)}\\& &+x_{21}E_{-\al_5}-x_{24}E_{\al_2}
-x_{25}E_{(0,1,0,1)}+x_{26}E_{(0,1,1,1)}+x_{27}E_{(1,1,1,1)},
\hspace{1.5cm}(12.5.34)\end{eqnarray*}
 \begin{eqnarray*}\iota(\mfk b_{23})
&=&P_{23}+x_2E'_{(1,2,2,3,2,1)}
+x_6E'_{(1,1,2,2,1,1)}-x_7E'_{(1,1,1,2,1,1)}
+x_9E'_{(1,1,1,1,1,1)}\\
&&-x_{23}\mfk
k(\al_{(1,2,2,3,2,2,1)}+x_{10}E'_{(0,1,1,2,1,1)}-x_{12}E'_{(0,1,1,1,1,1)}-x_{14}E'_{(1,1,1,1)}\\&
&+x_{15}E'_{(0,1,0,1,1,1)}+x_{17}E'_{(0,1,1,1)}+x_{20}E_{-\al_6}
 -x_{19}E_{(0,1,0,1)}+x_{21}E_{-\al_2}\\& &-x_{24}E_{\al_5}
+x_{25}E_{(0,0,0,1,1)}-x_{26}E_{(0,0,1,1,1)}-x_{27}E_{(1,0,1,1,1)},
\hspace{2.8cm}(12.5.35)\end{eqnarray*}
 \begin{eqnarray*}\iota(\mfk b_{24})
&=&P_{24}-x_3E'_{(1,2,2,3,2,1)} +x_6E'_{(1,1,2,2,2,1)} -x_{24}\mfk
k(\al_{(1,2,2,3,3,2,1)})-x_7E'_{(1,1,1,2,2,1)}\\
&&+x_{10}E'_{(0,1,1,2,2,1)}+x_{11}E'_{(1,1,1,1,1,1)}-x_{13}E'_{(0,1,1,1,1,1)}-x_{14}E'_{(1,1,1,1,1)}\\&
&+x_{16}E'_{(0,1,0,1,1,1)}+x_{17}E'_{(0,1,1,1,1)}-x_{19}E'_{(0,1,0,1,1)}-x_{20}E'_{(0,0,0,0,1,1)}
 \\& &+x_{22}E_{-\al_2}+x_{23}E_{-\al_5}
-x_{25}E_{\al_4}+x_{26}E_{(0,0,1,1)}+x_{27}E_{(1,0,1,1)},
\hspace{2cm}(12.5.36)\end{eqnarray*}
\begin{eqnarray*}\iota(\mfk b_{25})
&=&P_{25}+x_4E'_{(1,2,2,3,2,1)} +x_6E'_{(1,1,2,3,2,1)} -x_{25}\mfk
k(\al_{(1,2,3,4,3,2,1)}-x_9E'_{(1,1,1,2,2,1)}\\
&&+x_{12}E'_{(0,1,1,2,2,1)}+x_{11}E'_{(1,1,1,2,1,1)}-x_{13}E'_{(0,1,1,2,1,1)}-x_{14}E'_{(1,1,1,2,1)}\\&
&+x_{17}E'_{(0,1,1,2,1)}+x_{18}E'_{(0,1,0,1,1,1)}+x_{20}E'_{(0,0,0,1,1,1)}-x_{21}E'_{(0,1,0,1,1)}
 \\& &+x_{22}E'_{(0,1,0,1)}-x_{23}E'_{(0,0,0,1,1)}
+x_{24}E_{-\al_4}-x_{26}E_{\al_3}-x_{27}E_{(1,0,1)},
\hspace{1.5cm}(12.5.37)\end{eqnarray*}
\begin{eqnarray*}\iota(\mfk b_{26})
&=&P_{26}-x_5E'_{(1,2,2,3,2,1)} -x_7E'_{(1,1,2,3,2,1)} -x_{26}\mfk
k(\al_{(1,2,3,4,3,2,1)})+x_9E'_{(1,1,2,2,2,1)}\\
&&-x_{15}E'_{(0,1,1,2,2,1)}-x_{11}E'_{(1,1,2,2,1,1)}+x_{16}E'_{(0,1,1,2,1,1)}+x_{14}E'_{(1,1,2,2,1)}\\&
&-x_{19}E'_{(0,1,1,2,1)}-x_{18}E'_{(0,1,1,1,1,1)}-x_{20}E'_{(0,0,1,1,1,1)}+x_{21}E'_{(0,1,1,1,1)}
 \\& &-x_{22}E'_{(0,1,1,1)}+x_{23}E'_{(0,0,1,1,1)}
-x_{24}E'_{(0,0,1,1)}+x_{25}E_{-\al_3}+x_{27}E_{\al_1},
\hspace{1.2cm}(12.5.38)\end{eqnarray*}
\begin{eqnarray*}\iota(\mfk b_{27})
&=&P_{27}-x_8E'_{(1,2,2,3,2,1)} -x_{10}E'_{(1,1,2,3,2,1)}
-x_{27}\mfk
k(\al_{(2,2,3,4,3,2,1)})+x_{12}E'_{(1,1,2,2,2,1)}\\
&&-x_{15}E'_{(1,1,1,2,2,1)}-x_{13}E'_{(1,1,2,2,1,1)}+x_{16}E'_{(1,1,1,2,1,1)}+x_{17}E'_{(1,1,2,2,1)}\\&
&-x_{19}E'_{(1,1,1,2,1)}-x_{18}E'_{(1,1,1,1,1,1)}-x_{20}E'_{(1,0,1,1,1,1)}+x_{21}E'_{(1,1,1,1,1)}
 \\& &-x_{22}E'_{(1,1,1,1)}+x_{23}E'_{(1,0,1,1,1)}
-x_{24}E'_{(1,0,1,1)}+x_{25}E'_{(1,0,1)}-x_{26}E_{-\al_1}.
\hspace{0.7cm}(12.5.39)\end{eqnarray*}
  Note that the coefficient of $x_i$ in the above
 $\iota(\mfk b_i)$ is exactly $\mfk k(\be)$ if $\mfk b_i=E_\be$.

Recall ${\msr A}=\mbb F[x_1,...,x_{27}]$. Let $M$ be a ${\msr
G}^{E_6}$-module and set
$$\wht{M}={\msr A}\otimes_{\mbb{C}}M.\eqno(12.5.40)$$
We identify
$$f\otimes v=fv\qquad\for\;\;f\in{\msr A},\;v\in M.\eqno(12.5.41)$$
Recall the Lie algebra $\wht{\msr K}$ defined via (12.5.1)-(12.5.8).
Fix $c\in\mbb F$. Then $\wht M$ becomes a $\wht{\msr K}$-module with
the action defined by
$$d(fw)=d(f)w,\;\;\kappa(fw)=cfw,\;\;(gu)(fw)=fg
u(w)\eqno(12.5.42)$$ for $d\in \mbb W_{27},\;f,g\in{\msr A},\;w\in
M$ and $u\in{\msr G}^{E_6}$.

Since the linear map $\iota: {\msr G}^{E_7}\rta \wht{\msr K}$
defined in (12.5.10)-(12.5.12) is a  Lie algebra monomorphism,
$\wht{M}$ becomes a ${\msr G}^{E_7}$-module with the action defined
by
$$\xi(w)=\iota(\xi)(w)\qquad\for\;\;\xi\in{\msr
G}^{E_7},\;w\in\wht{M}.\eqno(12.5.43)$$ In fact, we have:\psp

{\bf Theorem 12.5.1}. {\it The map $M\mapsto \wht M$ gives a functor
from ${\msr G}^{E_6}$-{\bf Mod} to ${\msr G}^{E_7}$-{\bf Mod}, where
the morphisms are ${\msr G}^{E_6}$-module homomorphisms and ${\msr
G}^{E_7}$-module homomorphisms, respectively.}\psp

We remark that  the module $\wht M$ is not a generalized module in
general because it may not be equal to $U({\msr G})(M)=U({\msr
G}_-)(M)$. \psp

{\bf Proposition 12.5.2}. {\it If $M$ is an irreducible ${\msr
G}^{E_6}$-module, then $U({\msr G}_-)(M)$ is an irreducible ${\msr
G}^{E_7}$-module.}

{\it Proof.} Note that for any $i\in\ol{1,27}$, $f\in{\msr A}$ and
$v\in M$, (12.5.12), (12.5.42) and (12.5.43) imply
$$\mfk a_i(fv)=\ptl_{x_i}(f)v.\eqno(12.5.44)$$
Let $W$ be any nonzero ${\msr G}^{E_7}$-submodule. The above
expression shows that $W\bigcap M\neq\{0\}$. According to (12.5.42),
$W\bigcap M$ is a ${\msr G}^{E_6}$-submodule. By the irreducibility
of $M$, $M\subset W$. Thus $U({\msr G}_-)(M)\subset W$. So $U({\msr
G}_-)(M)=W$ is irreducible. $\qquad\Box$\psp

By the above proposition, the map $M\mapsto U({\msr G}_-)(M)$ is a
polynomial extension from  irreducible ${\msr G}^{E_6}$-modules to
irreducible ${\msr G}^{E_7}$-modules. We can use it to derive
Gel'fand-Zetlin bases for $E_7$ from those for $E_6$.

\section{Irreducibility of the Functor}

In this section, we want to determine the irreducibility of ${\msr
G}^{E_6}$-modules $\wht M$.

Note that $\wht{M}$ can be viewed as a ${\msr G}^{E_6}$-module.
Indeed, (12.5.11) and (12.5.43) show
$$u(fv)=u(f)v+fu(v)\qquad\for\;\;u\in{\msr
G}^{E_6},\;f\in{\msr A},\;v\in M\eqno(12.6.1)$$ (cf.
(11.1.18)-(11.1.56)). So $\wht{M}={\msr A}\otimes_{\mbb{C}}M$ is a
tensor module of ${\msr G}^{E_7}$. Write
$$\mfk b^\al=\prod_{i=1}^{27}\mfk b_{i}^{\alpha_i},\;\;
 |\al|=\sum_{i=1}^{27}\al_i\;\; \for\;\;
\al=(\al_1,\al_2,...,\al_{27})\in\mbb{N}^{27}\eqno(12.6.2)$$ (cf.
(11.1.4)-(11.1.12)).
 Recall the Lie subalgebras ${\msr G}_\pm$ and ${\msr G}_0$ of ${\msr G}^{E_6}$ defined
 in (12.4.2).  For $k\in\mbb{N}$, we set
$${\msr A}_k=\mbox{Span}_{\mbb F}\{x^\al\mid
\al \in\mbb{N}^{27};|\al|=k\},\;\; \wht M_{k}={\msr
A}_kM\eqno(12.6.3)$$ (cf.  (12.5.41)) and
 $$(U({\msr G}_-)(
M))_{k}=\mbox{Span}_{\mbb F}\{\mfk b^\al( M)\mid \al
\in\mbb{N}^{27}, \ |\al|=k\}. \eqno(12.6.4)$$
 Moreover,
$$(U({\msr G}_-)(
M))_{0}=\wht M_{0}= M.\eqno(12.6.5)$$ Furthermore,
 $$\wht M=\bigoplus\limits_{k=0}^\infty\wht M_{
 k},\qquad
 U({\msr G}_-)(M)=\bigoplus\limits_{k=0}^\infty(U({\msr G}_-)(M))_{ k}.\eqno(12.6.6)$$

Next we define a linear transformation $\vf$ on  $\wht M$ determined
by
$$\vf(x^\al v)=\mfk b^\al(
v)\qquad\for\;\;\al\in\mbb{N}^{27},\;v\in M.\eqno(12.6.7)$$ Note
that ${\msr A}_1=\sum_{i=1}^{27}\mbb Fx_i$ forms the
$27$-dimensional ${\msr G}_0$-module of highest weight $\lmd_6$.
According to (4.4.24) and (4.4.25), ${\msr G}_-$ forms a ${\msr
G}_0$-module with respect to the adjoint representation, and the
linear map from ${\msr A}_1$ to ${\msr G}_-$ determined by
$x_i\mapsto \mfk b_i$ for $i\in\ol{1,27}$ gives a ${\msr
G}_0$-module isomorphism. Thus $\vf$ can also be viewed as a ${\msr
G}_0$-module homomorphism from $\wht M$ to $U({\msr G}_0)(M)$.
Moreover,
$$\vf(\wht M_{ k})=(U({\msr G}_-)(
M))_{k}\qquad\for\;\;k\in\mbb{N}.\eqno(12.6.8)$$

Observe that the fundamental weights of ${\msr G}^{E_6}$ are:
$$\lmd_1=\frac{1}{3}(4\al_1+3\al_2+5\al_3+6\al_4+4\al_5+2\al_6),\eqno(12.6.9)$$
$$\lmd_2=\al_1+2\al_2+2\al_3+3\al_4+2\al_5+\al_6,\eqno(12.6.10)$$
$$\lmd_3=\frac{1}{3}(5\al_1+6\al_2+10\al_3+12\al_4+8\al_5+4\al_6),\eqno(12.6.11)$$
$$\lmd_4=2\al_1+3\al_2+4\al_3+6\al_4+4\al_5+2\al_6,\eqno(12.6.12)$$
$$\lmd_5=\frac{1}{3}(4\al_1+6\al_2+8\al_3+12\al_4+10\al_5+5\al_6),\eqno(12.6.13)$$
$$\lmd_6=\frac{1}{3}(2\al_1+3\al_2+4\al_3+6\al_4+5\al_5+4\al_6).\eqno(12.6.14)$$
Using the above expressions, we can view $\lmd_i\in H_{E_6}$ by
(11.1.2). Then Casimir element of ${\msr G}^{E_6}$ is
$$\omega=\sum_{i=1}^6\lmd_i\al_i-\sum_{\be\in\Phi_{E_6}}E_{\be}E_{-\be}\in
U({\msr G}^{E_6})\eqno(12.6.15)$$ due to (12.4.12).
 The algebra
$U({\msr G}^{E_6})$  is imbedded into the tensor algebra $U({\msr
G}^{E_6})\otimes_\mbb{F} U({\msr G}^{E_6})$ by the associative
algebra homomorphism $\mfk{d}: U({\msr G}^{E_6}) \rightarrow U({\msr
G}^{E_6})\otimes_\mbb{F} U({\msr G}^{E_6}))$ determined  by
$$\mfk{d}(u)=u\otimes 1 +1 \otimes u \qquad \mbox{ for} \ u\in
{\msr G}^{E_6}.\eqno(12.6.16)$$ Set
$$\td\omega=\frac{1}{2}(\mfk{d}(\omega)-\omega\otimes 1-1\otimes
\omega)\in U({\msr G}^{E_6})\otimes_\mbb{F} U({\msr
G}^{E_6}).\eqno(12.6.17)$$ Since $\sum_{i=1}^6\lmd_i\otimes \al_i$
is symmetric with respect to $\{\al_1,...,\al_6\}$ by
(12.6.9)-(12.6.14),
$$\td\omega=\sum_{i=1}^6\lmd_i\otimes\al_i-\sum_{\be\in\Phi_{E_6}}E_{\be}\otimes
E_{-\be}.\eqno(12.6.18)$$ Furthermore, $\td{\omega}$ acts on $\wht
M$ as a ${\msr G}^{E_6}$-module homomorphism via
$$(u_1\otimes u_2)(fw)=u_1(f)u_2(w)\qquad\for\;\;u_1,u_2\in
{\msr G}^{E_6},\;f\in{\msr A},\;w\in M.\eqno(12.6.19)$$

{\bf Lemma 12.6.1}. {\it We have $\vf|_{\wht M_{
1}}=(\td\omega-c)|_{\wht M_{ 1}}$.}

{\it Proof.} By (11.1.54) with Table 11.1.1 and (12.6.9)-(12.6.14)
$$\lmd_r|_{\msr
A}=\sum_{i=1}^{27}c_{i,r}x_i\ptl_{x_i}\qquad\for\;\;r\in\ol{1,6}\eqno(12.6.20)$$
with $c_{i,r}\in\mbb{F}$, for instance,
$$c_{1,1}=\frac{2}{3},\;c_{1,2}=1,\;c_{1,3}=\frac{4}{3},\;c_{1,4}=2,\;c_{1,5}=
\frac{5}{3},\;c_{1,6}=\frac{4}{3}.\eqno(12.6.21)$$ On the other
hand,
$$\al_7=\frac{1}{3}[\wht\al-(2\al_1+3\al_2+4\al_3+6\al_4+5\al_5+4\al_6)]\eqno(12.6.22)$$
by (12.4.5). According to (12.5.10), (12.5.42), (12.6.19) and
(12.6.21),
$$(\sum_{i=1}^6\lmd_i\otimes\al_i)(x_1w)=cx_1w-x_1\al_7(w)\qquad\for\;w\in M.\eqno(12.6.23)$$
Thus (11.1.18)-(11.1.56) yield
\begin{eqnarray*}\td\omega(x_1w)&=&[cx_1+x_2E_{\al_6}+x_3E_{(0,0,0,0,1,1)}
+x_4E_{(0,0,0,1,1,1)}+x_5E_{(0,0,1,1,1,1)} +x_6E_{(0,1,0,1,1,1)}\\ &
&-x_1\al_7(w)+x_7E_{(0,1,1,1,1,1)}+x_8E_{(1,0,1,1,1,1)}+x_9E_{(0,1,1,2,1,1)}+x_{10}E_{(1,1,1,1,1,1)}
\\ & &+x_{11}E_{(0,1,1,2,2,1)}+x_{12}E_{(1,1,1,2,1,1)}+x_{13}E_{(1,1,1,2,2,1)}+x_{15}E_{(1,1,2,2,1,1)}\\
&
&+x_{16}E_{(1,1,2,2,2,1)}+x_{18}E_{(1,1,2,3,2,1)}+x_{20}E_{(1,2,2,3,2,1)}](w).\hspace{2.7cm}(12.6.24)
\end{eqnarray*}
Comparing (12.5.13) and (12.6.24), we get
$\vf(x_1w)=(\td\omega-c)(x_1w)$. According to (11.1.23) and
(11.1.56), $-E_{-\al_6}(x_1)=x_2$. So for $v\in M$,
\begin{eqnarray*}\qquad& &\vf(x_2v)-\vf(x_1E_{-\al_6}(v))=
-(E_{-\al_6})(\vf(x_1v))\\
&=&-(E_{-\al_6})[(\td\omega-c)(x_1v)]=(\td\omega-c)[-E_{-\al_6}(x_1v)]
\\
&=&(\td\omega-c)[x_2v-x_1E_{-\al_6}(v)]\\
&=&(\td\omega-c)(x_2v)-(\td\omega-c)(x_1E_{-\al_6})(v))
\\ &=&(\td\omega-c)(x_2v)-\vf(x_1E_{-\al_6}(v)),\hspace{7.2cm}(12.6.25)\end{eqnarray*}
or equivalently, $\vf(x_2v)=(\td\omega-c)(x_2v)$. By
(11.1.18)-(11.1.23) and (11.1.56), similar arguments as (12.6.25)
and induction on $i$, we can prove
$$\vf(x_iv)=(\td\omega-c)(x_iv)\qquad\for\;\;i\in\ol{1,27},\;v\in M;\eqno(12.6.26)$$
that is, the lemma holds.$\qquad\Box$\psp

Recall that $\mbb{N}$ is the set of nonnegative integers and the set
of dominate integral weights is
$$\Lmd^+=\sum_{r=1}^6\mbb{N}\lmd_r.\eqno(12.6.27)$$
Denote
$$\rho=\sum_{r=1}^6\lmd_r\eqno(12.6.28)$$
(cf. (12.6.9)-(12.6.14). For any $\mu\in \Lmd^+$, we denote by
$V(\mu)$ the finite-dimensional irreducible ${\msr G}^{E_6}$-module
with the highest weight $\mu$ and have
$$\omega|_{V(\mu)}=(\mu+2\rho,\mu)\mbox{Id}_{V(\mu)}\eqno(12.6.29)$$
by (12.6.15). According to (11.1.54) and Table 11.1.1, the weight
set of the ${\msr G}^{E_6}$-module ${\msr A}_1$ is
$$\Pi({\msr A}_1)=\{\sum_{r=1}^6a_{i,r}\lmd_r\mid
i\in\ol{1,27}\}.\eqno(12.6.30)$$
 Fixing
$\lmd\in\Lmd^+$, we define
$$\Upsilon(\lmd)=\{\lmd+\mu\mid\mu\in\Pi({\msr
A}_1),\;\lmd+\mu\in\Lmd^+\}.\eqno(12.6.31)$$ \pse

{\bf Lemma 12.6.2}. {\it We have}:
$${\msr A}_1\otimes V(\lmd)\cong \bigoplus_{\lmd'\in
\Upsilon(\lmd)}V(\lmd').\eqno(12.6.32)$$

{\it Proof}. Note that all the weight subspaces of ${\msr A}_1$ are
one-dimensional. Thus all the irreducible components of ${\msr
A}_1\otimes V(\lmd)$ are of multiplicity one. Since
$$\rho+\lmd+\mu\in\Lmd^+\qquad\for\;\;\mu\in \Pi({\msr
A}_1),\eqno(12.6.33)$$ Theorem 5.4.3 says that $V(\lmd')$ is a
component of ${\msr A}_1\otimes V(\lmd)$ if and only if $\lmd'\in
\Upsilon(\lmd).\qquad\Box$\psp

Recall
$$\mbox{the highest weight of}\;{\msr
A}_1=\lmd_6\eqno(12.6.34)$$ by Table 11.1.1. Thus the eigenvalues of
$\td{\omega}|_{\wht{V(\lmd)}_{1}}$ are
$$\{[(\lmd'+2\rho,\lmd')-(\lmd+2\rho,\lmd)-(\lmd_6+2\rho,\lmd_6)]/2\mid\lmd'\in
\Upsilon(\lmd)\}\eqno(12.6.35)$$ by (12.6.17) and (12.6.19). We
remark that the above fact is equivalent to a special explicit
version of Kostant's characteristic identity (cf. [Kb]). Define
$$\ell_\omega(\lmd)=\min\{[(\lmd'+2\rho,\lmd')-(\lmd+2\rho,\lmd)-(\lmd_6+2\rho,\lmd_6)]/2\mid\lmd'\in
\Upsilon(\lmd)\},\eqno(12.6.36)$$ which will be used to determine
the irreducibility of $\wht{V(\lmd)}$. If
$\lmd'=\lmd+\lmd_6-\al\in\Upsilon(\lmd)$ with $\al\in\Phi_{D_6}^+$,
then
$$(\lmd'+2\rho,\lmd')-(\lmd+2\rho,\lmd)-(\lmd_6+2\rho,\lmd_6)
=2[(\lmd,\lmd_6)+1-(\rho+\lmd+\lmd_6,\al)].\eqno(12.6.37)$$

Recall the differential operators $P_1,...,P_{27}$ given in
(12.4.20) and (12.4.39)-(12.4.64). We also view the elements of
${\msr A}$ as the multiplication operators on ${\msr A}$. Recall
$\zeta_1$ in (11.2.1). It turns out that we need the following lemma
in order to determine the irreducibility of $\wht{V(\lmd)}$. Recall
the Dickson's invariant $\eta$ given in (11.2.42).\psp

{\bf Lemma 12.6.3}. {\it As operators on} ${\msr A}$:
\begin{eqnarray*}& &P_{14}x_1+P_1x_{14}+P_{11}x_2+P_2x_{11}+P_9x_3+P_3x_9+P_7x_4+P_4x_7
-P_6x_5-P_5x_6\\
&=&\zeta_1(D-8)+\eta\ptl_{x_{27}}.\hspace{10.1cm}(12.6.38)\end{eqnarray*}

{\it Proof}. According to  (12.4.20), (12.4.39)-(12.4.44),
(12.4.46), (12.4.48) and (12.4.51), we find that
\begin{eqnarray*}&
&P_{14}x_1+P_1x_{14}+P_{11}x_2+P_2x_{11}+P_9x_3+P_3x_9+P_7x_4+P_4x_7
-P_6x_5-P_5x_6\\
&=&-8\zeta_1+x_1P_{14}+x_{14}P_1+x_2P_{11}+x_{11}P_2+x_3P_9\\
& &+x_9P_3+x_4P_7+x_7P_4
-x_5P_6-x_6P_5\hspace{6.9cm}(12.6.39)\end{eqnarray*} and
\begin{eqnarray*}&
&x_1P_{14}+x_{14}P_1+x_2P_{11}+x_{11}P_2+x_3P_9+x_9P_3+x_4P_7+x_7P_4
-x_5P_6-x_6P_5
\\\!\!\!\!\!\!&=&\!\!\!x_1(x_{14}D-\zeta_1\ptl_{x_1}-\zeta_8\ptl_{x_8}-\zeta_{10}\ptl_{x_{10}}
+\zeta_{12}\ptl_{x_{12}}-\zeta_{13}\ptl_{x_{13}}-\zeta_{15}\ptl_{x_{15}}
+\zeta_{16}\ptl_{x_{16}}-\zeta_{18}\ptl_{x_{18}}\\ &
&-\zeta_{20}\ptl_{x_{20}}
+\zeta_{27}\ptl_{x_{27}})+x_{14}(x_1D-\zeta_1\ptl_{x_{14}}-\zeta_2\ptl_{x_{17}}-\zeta_3\ptl_{x_{19}}
+\zeta_4\ptl_{x_{21}}-\zeta_5\ptl_{x_{22}}
-\zeta_6\ptl_{x_{23}}\\&&+\zeta_7\ptl_{x_{24}}-\zeta_9\ptl_{x_{25}}-\zeta_{11}\ptl_{x_{26}}
-\zeta_{14}\ptl_{x_{27}})+x_2(x_{11}D-\zeta_1\ptl_{x_2}+\zeta_5\ptl_{x_8}+\zeta_7\ptl_{x_{10}}
-\zeta_9\ptl_{x_{12}}\\
& &+\zeta_{11}\ptl_{x_{15}}+\zeta_{13}\ptl_{x_{17}}
-\zeta_{16}\ptl_{x_{19}}+\zeta_{18}\ptl_{x_{21}}+\zeta_{20}\ptl_{x_{23}}
+\zeta_{26}\ptl_{x_{27}})+x_{11}(x_2D-\zeta_1\ptl_{x_{11}}\\
& &-\zeta_2\ptl_{x_{13}}-\zeta_3\ptl_{x_{16}}
+\zeta_4\ptl_{x_{18}}-\zeta_6\ptl_{x_{20}}-\zeta_8\ptl_{x_{22}}
+\zeta_{10}\ptl_{x_{24}}-\zeta_{12}\ptl_{x_{25}}-\zeta_{15}\ptl_{x_{26}}
+\zeta_{17}\ptl_{x_{27}})
\\ &&+x_3(x_9D-\zeta_1\ptl_{x_3}-\zeta_4\ptl_{x_8}-\zeta_6\ptl_{x_{10}}
+\zeta_9\ptl_{x_{13}}-\zeta_{12}\ptl_{x_{17}}-\zeta_{11}\ptl_{x_{16}}
+\zeta_{15}\ptl_{x_{19}}-\zeta_{18}\ptl_{x_{22}}\\ &
&-\zeta_{20}\ptl_{x_{24}}
+\zeta_{25}\ptl_{x_{27}})+x_9(x_3D-\zeta_1\ptl_{x_9}-\zeta_2\ptl_{x_{12}}-\zeta_3\ptl_{x_{15}}
+\zeta_5\ptl_{x_{18}}-\zeta_7\ptl_{x_{20}}-\zeta_8\ptl_{x_{21}}
\\ &&+\zeta_{10}\ptl_{x_{23}}-\zeta_{13}\ptl_{x_{25}}-\zeta_{16}\ptl_{x_{26}}
+\zeta_{19}\ptl_{x_{27}})+x_4(x_7D-\zeta_1\ptl_{x_4}+\zeta_3\ptl_{x_8}+\zeta_6\ptl_{x_{12}}
-\zeta_7\ptl_{x_{13}}\\ &
&+\zeta_{10}\ptl_{x_{17}}+\zeta_{11}\ptl_{x_{18}}
-\zeta_{15}\ptl_{x_{21}}+\zeta_{16}\ptl_{x_{22}}+\zeta_{20}\ptl_{x_{25}}
+\zeta_{24}\ptl_{x_{27}})+x_7(x_4D-\zeta_1\ptl_{x_7}\\
& &-\zeta_2\ptl_{x_{10}}-\zeta_4\ptl_{x_{15}}
+\zeta_5\ptl_{x_{16}}-\zeta_8\ptl_{x_{19}}-\zeta_9\ptl_{x_{20}}
+\zeta_{12}\ptl_{x_{23}}-\zeta_{13}\ptl_{x_{24}}-\zeta_{18}\ptl_{x_{26}}
+\zeta_{21}\ptl_{x_{27}})\\ &
&-x_5(x_6D+\zeta_1\ptl_{x_5}+\zeta_2\ptl_{x_8}-\zeta_6\ptl_{x_{15}}
+\zeta_7\ptl_{x_{16}}-\zeta_9\ptl_{x_{18}}-\zeta_{10}\ptl_{x_{19}}
+\zeta_{12}\ptl_{x_{21}}\\
& &-\zeta_{13}\ptl_{x_{22}}+\zeta_{20}\ptl_{x_{26}}
-\zeta_{23}\ptl_{x_{27}})-x_6(x_5D+\zeta_1\ptl_{x_6}-\zeta_3\ptl_{x_{10}}+\zeta_4\ptl_{x_{12}}
\\ &&-\zeta_5\ptl_{x_{13}}+\zeta_8\ptl_{x_{17}}+\zeta_{11}\ptl_{x_{20}}
-\zeta_{15}\ptl_{x_{23}}+\zeta_{16}\ptl_{x_{24}}-\zeta_{18}\ptl_{x_{25}}
-\zeta_{22}\ptl_{x_{27}})
\\\!\!\!\!\!\!&=&\!\!\!\zeta_1D+\eta\ptl_{x_{27}}.\qquad\Box\hspace{10.1cm}(12.6.40)
\end{eqnarray*}
\psp

{\bf Lemma 12.6.4}. {\it As operators on ${\msr A}$,}
$$\sum_{14\neq
i\in\ol{1,26}}P_i\zeta_{28-i}-P_{14}\zeta_{14}-P_{27}\zeta_1=\eta(24-5D).\eqno(12.6.41)$$

{\it Proof}. We calculate it by (11.2.1), (11.2.3)-(11.2.28),
(12.4.20) and (12.4.39)-(12.4.64). In particular,
$$2(\zeta_1\zeta_{14}-\zeta_2\zeta_{11}-\zeta_3\zeta_9+\zeta_4\zeta_7-\zeta_5\zeta_6)=-2\eta x_1,\eqno(12.6.42)$$
which is the coefficient of $\ptl_{x_1}$ in addition to the term
containing $D$. According to (11.2.41), the operator on the left
hand side of (12.6.41) is a ${\msr G}^{E_6}$-invariant differential
operator. By symmetry,
$$\mbox{the coefficient of}\;\ptl_{x_i}=-2\eta x_i\;\;\for\;\;i\in\ol{1,27}.\qquad\Box\eqno(12.6.43)$$
\pse

We define the multiplication
$$f(gv)=(fg)v\qquad\for\;\;f,g\in{\msr A},\;v\in M.\eqno(12.6.44)$$
 Furthermore, we have:\psp

{\bf Lemma 12.6.5}. {\it As operators on $\wht M$,}
$$[(\sum_{14\neq
i\in\ol{1,26}}\iota(\mfk b_i)\zeta_{28-i}-\iota(\mfk
b_{14})\zeta_{14}-\iota(\mfk b_{27})\zeta_1]|_{\wht
M}=\eta(24-5D+3c).\eqno(12.6.45)$$

{\it Proof}. By (12.5.13)-(12.5.39), the coefficient of $E_{\al_6}$
is
$$x_2\zeta_{27}-x_{14}\zeta_{17}-x_{17}\zeta_{15}-x_{19}\zeta_{12}-x_{21}\zeta_{10}-x_{23}\zeta_8.\eqno(12.6.46)$$
Moreover, we use (11.2.9), (11.2.11), (11.2.13), (11.2.16),
(11.2.18) and (11.2.28) to find that (12.6.46) is equal to zero.
Since the left hand side of (12.6.45) is a ${\msr
G}^{E_6}$-invariant differential operator, it is invariant under the
action of the $E_6$ Weyl group. The transitivity of the Weyl group
on $\Phi_{E_6}$ yields that
$$\mbox{the coefficient of}\;E_\be=0\;\;\mbox{for
any}\;\be\in\Phi_{E_6}.\eqno(12.6.47)$$ By (12.5.10), (12.5.42),
(12.6.22) and Lemma 12.6.4,
$$\sum_{14\neq
i\in\ol{1,26}}\iota(\mfk b_i)\zeta_{28-i}-\iota(\mfk
b_{14})\zeta_{14}-\iota(\mfk b_{27})\zeta_1
=\eta(24-5D+3c)+\sum_{r=1}^6f_r\al_r\eqno(12.6.48)$$ as the elements
of $\wht {\msr K}$ acting on $\wht M$ (cf. (12.5.7)). The ${\msr
G}^{E_6}$-invariancy implies
$$[E_\be|_{\wht M},\sum_{r=1}^6f_r\al_r]=0\;\;\mbox{for
any}\;\be\in\Phi_{E_6}.\eqno(12.6.49)$$ Thus
$\sum_{r=1}^6f_r\al_r=0$; that is, (12.6.45) holds. $\qquad\Box$\psp

Next we calculate
\begin{eqnarray*}\qquad T_1&=&\iota(\mfk b_{14})x_1+\iota(\mfk b_1)x_{14}+\iota(\mfk b_{11})x_2+\iota(\mfk b_2)x_{11}+\iota(\mfk b_9)x_3
\\ &&+\iota(\mfk b_3)x_9+\iota(\mfk b_7)x_4+\iota(\mfk b_4)x_7
-\iota(\mfk b_6)x_5-\iota(\mfk b_5)x_6\\
&=&\zeta_1(D-8)+\eta\ptl_{x_{27}}+\zeta_2E_{\al_1}-\zeta_3E_{(1,0,1)}+\zeta_4E_{(1,0,1,1)}-\zeta_5E_{(1,0,1,1,1)}
\\ & &+\zeta_6E_{(1,1,1,1)}
-\zeta_7E_{(1,1,1,1,1)}+\zeta_8E_{(1,0,1,1,1,1)}+\zeta_9E_{(1,1,1,2,1)}\\
& &+\zeta_{10}E_{(1,1,1,1,1,1)}
-\zeta_{11}E_{(1,1,2,2,1)}-\zeta_{12}E_{(1,1,1,2,1,1)}+\zeta_{13}E_{(1,1,1,2,2,1)}
\\ & &+\zeta_{15}E_{(1,1,2,2,1,1)}-\zeta_{16}E_{(1,1,2,2,2,1)}+\zeta_{18}E_{(1,1,2,3,2,1)}
+\zeta_{20}E_{(1,2,2,3,2,1)}\\& &-2\zeta_1
c+\frac{\zeta_1}{3}(4\al_1+3\al_2+5\al_3+6\al_4+4\al_5+2\al_6)
 \hspace{3.5cm}(12.6.50)\end{eqnarray*}
by (11.2.1), (11.2.3)-(11.2.14), (11.2.16), (11.2.17), (11.2.19),
(11.2.21), (12.5.13)-(12.5.19), (12.5.21), (12.5.23), (12.5.26),
(12.6.20) and Lemma 12.6.3.

We define a ${\msr G}^{E_6}$-module
 structure on the space $\mbox{End}\:\wht M$ of linear
 transformations on $\wht M$ by
 $$\iota(u)(T)=[\iota(u),T]=\iota(u)T-T\iota(u)\qquad\for\;\;u\in
{\msr G}^{E_6},\;T\in \mbox{End}\:\wht M\eqno(12.6.51)$$ (cf.
(12.6.1)). It can be verified that $T_1$ is a ${\msr
G}^{E_6}$-singular vector with weight $\lmd_1$ in $\mbox{End}\:\wht
M$. So it generates the 27-dimensional  module of highest weight
$\lmd_1$.
  We set
\begin{eqnarray*}T_2&=&-[\iota(E_{-\al_1}),T_1]=
\zeta_2(D-8)-\eta\ptl_{x_{26}}-\zeta_1E_{-\al_1}+\zeta_3E_{\al_3}-\zeta_4E_{(0,0,1,1)}\\&
&+\zeta_5E_{(0,0,1,1,1)}-\zeta_6E_{(0,1,1,1)}
+\zeta_7E_{(0,1,1,1,1)}-\zeta_8E_{(0,0,1,1,1,1)}-\zeta_9E_{(0,1,1,2,1)}\\
& &-\zeta_{10}E_{(0,1,1,1,1,1)}-\zeta_{14}E_{(1,1,2,2,1)}
+\zeta_{12}E_{(0,1,1,2,1,1)}-\zeta_{13}E_{(0,1,1,2,2,1)}\\&
&-\zeta_{17}E_{(1,1,2,2,1,1)}+\zeta_{19}E_{(1,1,2,2,2,1)}-\zeta_{21}E_{(1,1,2,3,2,1)}
-\zeta_{23}E_{(1,2,2,3,2,1)}\\& &-2\zeta_2
c+\frac{\zeta_2}{3}(\al_1+3\al_2+5\al_3+6\al_4+4\al_5+2\al_6),
 \hspace{4.4cm}(12.6.52)\end{eqnarray*}
\begin{eqnarray*}T_3&=&-[\iota(E_{-\al_3}),T_2]=
\zeta_3(D-8)-\eta\ptl_{x_{25}}+\zeta_1E'_{(1,0,1)}-\zeta_2E_{-\al_3}+\zeta_4E_{\al_4}\\&
&-\zeta_5E_{(0,0,0,1,1)}+\zeta_6E_{(0,1,0,1)}
-\zeta_7E_{(0,1,0,1,1)}+\zeta_8E_{(0,0,0,1,1,1)}+\zeta_{11}E_{(0,1,1,2,1)}\\
& &+\zeta_{10}E_{(0,1,0,1,1,1)}+\zeta_{14}E_{(1,1,1,2,1)}
-\zeta_{15}E_{(0,1,1,2,1,1)}+\zeta_{16}E_{(0,1,1,2,2,1)}\\&
&+\zeta_{17}E_{(1,1,1,2,1,1)}-\zeta_{19}E_{(1,1,1,2,2,1)}+\zeta_{22}E_{(1,1,2,3,2,1)}
+\zeta_{24}E_{(1,2,2,3,2,1)}\\& &-2\zeta_3
c+\frac{\zeta_3}{3}(\al_1+3\al_2+2\al_3+6\al_4+4\al_5+2\al_6),
 \hspace{4.4cm}(12.6.53)\end{eqnarray*}
\begin{eqnarray*}T_4&=&-[\iota(E_{-\al_4}),T_3]=
\zeta_4(D-8)-\eta\ptl_{x_{24}}-\zeta_1E'_{(1,0,1,1)}+\zeta_2E'_{(0,0,1,1)}-\zeta_3E_{-\al_4}\\&
&+\zeta_5E_{\al_5}+\zeta_6E_{\al_2}
-\zeta_9E_{(0,1,0,1,1)}-\zeta_8E_{(0,0,0,0,1,1)}+\zeta_{11}E_{(0,1,1,1,1)}\\
& &+\zeta_{12}E_{(0,1,0,1,1,1)}+\zeta_{14}E_{(1,1,1,1,1)}
-\zeta_{15}E_{(0,1,1,1,1,1)}+\zeta_{18}E_{(0,1,1,2,2,1)}\\&
&+\zeta_{17}E_{(1,1,1,1,1,1)}-\zeta_{21}E_{(1,1,1,2,2,1)}+\zeta_{22}E_{(1,1,2,2,2,1)}
-\zeta_{25}E_{(1,2,2,3,2,1)}\\& &-2\zeta_4
c+\frac{\zeta_4}{3}(\al_1+3\al_2+2\al_3+3\al_4+4\al_5+2\al_6),
 \hspace{4.4cm}(12.6.54)\end{eqnarray*}
\begin{eqnarray*}T_5&=&-[\iota(E_{-\al_5}),T_4]=
\zeta_5(D-8)-\eta\ptl_{x_{23}}+\zeta_1E'_{(1,0,1,1,1)}-\zeta_2E'_{(0,0,1,1,1)}\\&
&+\zeta_3E'_{(0,0,0,1,1)}-\zeta_4E_{-\al_5}+\zeta_7E_{\al_2}
-\zeta_9E_{(0,1,0,1)}+\zeta_8E_{\al_6}+\zeta_{11}E_{(0,1,1,1)}\\
& &+\zeta_{13}E_{(0,1,0,1,1,1)}+\zeta_{14}E_{(1,1,1,1)}
-\zeta_{16}E_{(0,1,1,1,1,1)}+\zeta_{18}E_{(0,1,1,2,1,1)}\\&
&+\zeta_{19}E_{(1,1,1,1,1,1)}-\zeta_{21}E_{(1,1,1,2,1,1)}+\zeta_{22}E_{(1,1,2,2,1,1)}
+\zeta_{26}E_{(1,2,2,3,2,1)}\\& &-2\zeta_5
c+\frac{\zeta_5}{3}(\al_1+3\al_2+2\al_3+3\al_4+\al_5+2\al_6),
 \hspace{4.5cm}(12.6.55)\end{eqnarray*}
\begin{eqnarray*}T_6&=&-[\iota(E_{-\al_2}),T_4]=
\zeta_6(D-8)-\eta\ptl_{x_{22}}-\zeta_1E'_{(1,1,1,1)}+\zeta_2E'_{(0,1,1,1)}-\zeta_3E'_{(0,1,0,1)}\\&
&+\zeta_7E_{\al_5}-\zeta_4E_{-\al_2}
+\zeta_9E_{(0,0,0,1,1)}-\zeta_{10}E_{(0,0,0,0,1,1)}-\zeta_{11}E_{(0,0,1,1,1)}\\
& &-\zeta_{12}E_{(0,0,0,1,1,1)}-\zeta_{14}E_{(1,0,1,1,1)}
+\zeta_{15}E_{(0,0,1,1,1,1)}-\zeta_{20}E_{(0,1,1,2,2,1)}\\&
&-\zeta_{17}E_{(1,0,1,1,1,1)}+\zeta_{23}E_{(1,1,1,2,2,1)}-\zeta_{24}E_{(1,1,2,2,2,1)}
-\zeta_{25}E_{(1,1,2,3,2,1)}\\& &-2\zeta_6
c+\frac{\zeta_6}{3}(\al_1+2\al_3+3\al_4+4\al_5+2\al_6),
 \hspace{5.5cm}(12.6.56)\end{eqnarray*}
\begin{eqnarray*}T_7&=&-[\iota(E_{-\al_5}),T_6]=
\zeta_7(D-8)-\eta\ptl_{x_{21}}+\zeta_1E'_{(1,1,1,1,1)}-\zeta_2E'_{(0,1,1,1,1)}\\&
&+\zeta_3E'_{(0,1,0,1,1)}-\zeta_6E_{-\al_5}-\zeta_5E_{-\al_2}
+\zeta_9E_{\al_4}+\zeta_{10}E_{\al_6}-\zeta_{11}E_{(0,0,1,1)}\\
& &-\zeta_{13}E_{(0,0,0,1,1,1)}-\zeta_{14}E_{(1,0,1,1)}
+\zeta_{16}E_{(0,0,1,1,1,1)}-\zeta_{20}E_{(0,1,1,2,1,1)}\\&
&-\zeta_{19}E_{(1,0,1,1,1,1)}+\zeta_{23}E_{(1,1,1,2,1,1)}-\zeta_{24}E_{(1,1,2,2,1,1)}
+\zeta_{26}E_{(1,1,2,3,2,1)}\\& &-2\zeta_7
c+\frac{\zeta_7}{3}(\al_1+2\al_3+3\al_4+\al_5+2\al_6),
 \hspace{5.7cm}(12.6.57)\end{eqnarray*}
\begin{eqnarray*}T_8&=&-[\iota(E_{-\al_6}),T_5]=
\zeta_8(D-8)-\eta\ptl_{x_{20}}-\zeta_1E'_{(1,0,1,1,1,1)}+\zeta_2E'_{(0,0,1,1,1,1)}\\&
&-\zeta_3E'_{(0,0,0,1,1,1)}+\zeta_4E'_{(0,0,0,0,1,1)}-\zeta_5E_{-\al_6}+\zeta_{10}E_{\al_2}
-\zeta_{12}E_{(0,1,0,1)}\\
&&+\zeta_{13}E_{(0,1,0,1,1)}+\zeta_{15}E_{(0,1,1,1)}-\zeta_{17}E_{(1,1,1,1)}
-\zeta_{16}E_{(0,1,1,1,1)}+\zeta_{18}E_{(0,1,1,2,1)}\\&
&+\zeta_{19}E_{(1,1,1,1,1)}-\zeta_{21}E_{(1,1,1,2,1)}+\zeta_{22}E_{(1,1,2,2,1)}
-\zeta_{27}E_{(1,2,2,3,2,1)}\\& &-2\zeta_8
c+\frac{\zeta_8}{3}(\al_1+3\al_2+2\al_3+3\al_4+\al_5-\al_6),
 \hspace{4.7cm}(12.6.58)\end{eqnarray*}
\begin{eqnarray*}T_9&=&-[\iota(E_{-\al_4}),T_7]=
\zeta_9(D-8)-\eta\ptl_{x_{19}}-\zeta_1E'_{(1,1,1,2,1)}+\zeta_2E'_{(0,1,1,2,1)}\\&
&+\zeta_4E'_{(0,1,0,1,1)}-\zeta_6E'_{(0,0,0,1,1)}+\zeta_5E'_{(0,1,0,1)}
-\zeta_7E_{-\al_4}-\zeta_{11}E_{\al_3}+\zeta_{12}E_{\al_6}\\
& &+\zeta_{13}E_{(0,0,0,0,1,1)}-\zeta_{14}E_{(1,0,1)}
+\zeta_{18}E_{(0,0,1,1,1,1)}-\zeta_{20}E_{(0,1,1,1,1,1)}\\&
&-\zeta_{21}E_{(1,0,1,1,1,1)}+\zeta_{23}E_{(1,1,1,1,1,1)}+\zeta_{25}E_{(1,1,2,2,1,1)}
+\zeta_{26}E_{(1,1,2,2,2,1)}\\& &-2\zeta_9
c+\frac{\zeta_9}{3}(\al_1+2\al_3+\al_5+2\al_6),
 \hspace{6.8cm}(12.6.59)\end{eqnarray*}
\begin{eqnarray*}T_{10}&=&-[\iota(E_{-\al_2}),T_8]=
\zeta_{10}(D-8)-\eta\ptl_{x_{18}}-\zeta_1E'_{(1,1,1,1,1,1)}+\zeta_2E'_{(0,1,1,1,1,1)}\\&
&-\zeta_3E'_{(0,1,0,1,1,1)}+\zeta_6E'_{(0,0,0,0,1,1)}-\zeta_7E_{-\al_6}-\zeta_8E_{-\al_2}
+\zeta_{12}E_{\al_4}\\
&&-\zeta_{13}E_{(0,0,0,1,1)}-\zeta_{15}E_{(0,0,1,1)}+\zeta_{17}E_{(1,0,1,1)}
+\zeta_{16}E_{(0,0,1,1,1)}-\zeta_{20}E_{(0,1,1,2,1)}\\&
&-\zeta_{19}E_{(1,0,1,1,1)}+\zeta_{23}E_{(1,1,1,2,1)}-\zeta_{24}E_{(1,1,2,2,1)}
-\zeta_{27}E_{(1,1,2,3,2,1)}\\& &-2\zeta_{10}
c+\frac{\zeta_{10}}{3}(\al_1+2\al_3+3\al_4+\al_5-\al_6),
 \hspace{5.4cm}(12.6.60)\end{eqnarray*}
\begin{eqnarray*}T_{11}&=&[\iota(E_{-\al_3}),T_9]=
\zeta_{11}(D-8)-\eta\ptl_{x_{17}}+\zeta_1E'_{(1,1,2,2,1)}-\zeta_3E'_{(0,1,1,2,1)}\\&
&-\zeta_4E'_{(0,1,1,1,1)}+\zeta_6E'_{(0,0,1,1,1)}-\zeta_5E'_{(0,1,1,1)}
+\zeta_7E'_{(0,0,1,1)}+\zeta_9E_{-\al_3}\\
&&+\zeta_{14}E_{\al_1}+\zeta_{15}E_{\al_6}+\zeta_{16}E_{(0,0,0,0,1,1)}
+\zeta_{18}E_{(0,0,0,1,1,1)}-\zeta_{20}E_{(0,1,0,1,1,1)}\\&
&-\zeta_{22}E_{(1,0,1,1,1,1)}+\zeta_{24}E_{(1,1,1,1,1,1)}+\zeta_{25}E_{(1,1,1,2,1,1)}
+\zeta_{26}E_{(1,1,1,2,2,1)}\\& &-2\zeta_{11}
c+\frac{\zeta_{11}}{3}(\al_1-\al_3+\al_5+2\al_6),
 \hspace{6.6cm}(12.6.61)\end{eqnarray*}
\begin{eqnarray*}T_{12}&=&-[\iota(E_{-\al_4}),T_{10}]=
\zeta_{12}(D-8)-\eta\ptl_{x_{16}}+\zeta_1E'_{(1,1,1,2,1,1)}-\zeta_2E'_{(0,1,1,2,1,1)}\\&
&-\zeta_4E'_{(0,1,0,1,1,1)}+\zeta_6E'_{(0,0,0,1,1,1)}+\zeta_8E'_{(0,1,0,1)}-\zeta_9E_{-\al_6}
-\zeta_{10}E_{-\al_4}\\
&&+\zeta_{13}E_{\al_5}-\zeta_{15}E_{\al_3}+\zeta_{17}E_{(1,0,1)}
+\zeta_{18}E_{(0,0,1,1,1)}-\zeta_{20}E_{(0,1,1,1,1)}\\&
&-\zeta_{21}E_{(1,0,1,1,1)}+\zeta_{23}E_{(1,1,1,1,1)}+\zeta_{25}E_{(1,1,2,2,1)}
-\zeta_{27}E_{(1,1,2,2,2,1)}\\& &-2\zeta_{12}
c+\frac{\zeta_{12}}{3}(\al_1+2\al_3+\al_5-\al_6),
 \hspace{6.5cm}(12.6.62)\end{eqnarray*}
\begin{eqnarray*}T_{13}&=&-[\iota(E_{-\al_5}),T_{12}]=
\zeta_{13}(D-8)-\eta\ptl_{x_{15}}-\zeta_1E'_{(1,1,1,2,2,1)}-\zeta_2E'_{(0,1,1,2,2,1)}\\&
&-\zeta_5E'_{(0,1,0,1,1,1)}+\zeta_7E'_{(0,0,0,1,1,1)}-\zeta_8E'_{(0,1,0,1,1)}-\zeta_9E'_{(0,0,0,0,1,1)}
+\zeta_{10}E'_{(0,0,0,1,1)}\\
&&-\zeta_{12}E_{-\al_5}-\zeta_{16}E_{\al_3}+\zeta_{19}E_{(1,0,1)}
+\zeta_{18}E_{(0,0,1,1)}-\zeta_{20}E_{(0,1,1,1)}\\&
&-\zeta_{21}E_{(1,0,1,1)}+\zeta_{23}E_{(1,1,1,1)}-\zeta_{26}E_{(1,1,2,2,1)}
-\zeta_{27}E_{(1,1,2,2,1,1)}\\& &-2\zeta_{13}
c+\frac{\zeta_{13}}{3}(\al_1+2\al_3-2\al_5-\al_6),
 \hspace{6.2cm}(12.6.63)\end{eqnarray*}
\begin{eqnarray*}T_{14}&=&-[\iota(E_{-\al_1}),T_{11}]=
\zeta_{14}(D-8)+\eta\ptl_{x_{14}}+\zeta_2E'_{(1,1,2,2,1)}-\zeta_3E'_{(1,1,1,2,1)}\\&
&-\zeta_4E'_{(1,1,1,1,1)}+\zeta_6E'_{(1,0,1,1,1)}-\zeta_5E'_{(1,1,1,1)}
+\zeta_7E'_{(1,0,1,1)}+\zeta_9E'_{(1,0,1)}\\
&&-\zeta_{11}E_{-\al_1}-\zeta_{17}E_{\al_6}-\zeta_{19}E_{(0,0,0,0,1,1)}
-\zeta_{21}E_{(0,0,0,1,1,1)}+\zeta_{23}E_{(0,1,0,1,1,1)}\\&
&+\zeta_{22}E_{(0,0,1,1,1,1)}-\zeta_{24}E_{(0,1,1,1,1,1)}-\zeta_{25}E_{(0,1,1,2,1,1)}
-\zeta_{26}E_{(0,1,1,2,2,1)}\\& &-2\zeta_{14}
c+\frac{\zeta_{14}}{3}(-2\al_1-\al_3+\al_5+2\al_6),
 \hspace{6.1cm}(12.6.64)\end{eqnarray*}
\begin{eqnarray*}T_{15}&=&[\iota(E_{-\al_3}),T_{12}]=
\zeta_{15}(D-8)-\eta\ptl_{x_{13}}-\zeta_1E'_{(1,1,2,2,1,1)}+\zeta_3E'_{(0,1,1,2,1,1)}\\&
&+\zeta_4E'_{(0,1,1,1,1,1)}-\zeta_6E'_{(0,0,1,1,1,1)}-\zeta_8E'_{(0,1,1,1)}
+\zeta_{10}E'_{(0,0,1,1)}-\zeta_{11}E_{-\al_6}\\
&&+\zeta_{12}E_{-\al_3}+\zeta_{16}E_{\al_5}-\zeta_{17}E_{\al_1}
+\zeta_{18}E_{(0,0,0,1,1)}-\zeta_{20}E_{(0,1,0,1,1)}\\&
&-\zeta_{22}E_{(1,0,1,1,1)}+\zeta_{24}E_{(1,1,1,1,1)}+\zeta_{25}E_{(1,1,1,2,1)}
-\zeta_{27}E_{(1,1,1,2,2,1)}\\& &-2\zeta_{15}
c+\frac{\zeta_{15}}{3}(\al_1-\al_3+\al_5-\al_6),
 \hspace{6.8cm}(12.6.65)\end{eqnarray*}
\begin{eqnarray*}T_{16}&=&[\iota(E_{-\al_3}),T_{13}]=
\zeta_{16}(D-8)-\eta\ptl_{x_{12}}+\zeta_1E'_{(1,1,2,2,2,1)}+\zeta_3E'_{(0,1,1,2,2,1)}\\&
&+\zeta_5E'_{(0,1,1,1,1,1)}-\zeta_7E'_{(0,0,1,1,1,1)}+\zeta_8E'_{(0,1,1,1,1)}
-\zeta_{10}E'_{(0,0,1,1,1)}-\zeta_{11}E'_{(0,0,0,0,1,1)}\\
&&+\zeta_{13}E_{-\al_3}-\zeta_{15}E_{-\al_5}-\zeta_{19}E_{\al_1}
+\zeta_{18}E_{\al_4}-\zeta_{20}E_{(0,1,0,1)}\\&
&-\zeta_{22}E_{(1,0,1,1)}+\zeta_{24}E_{(1,1,1,1)}-\zeta_{26}E_{(1,1,1,2,1)}
-\zeta_{27}E_{(1,1,1,2,1,1)}\\& &-2\zeta_{16}
c+\frac{\zeta_{16}}{3}(\al_1-\al_3-2\al_5-\al_6),
 \hspace{6.5cm}(12.6.66)\end{eqnarray*}
\begin{eqnarray*}T_{17}&=&[\iota(E_{-\al_1}),T_{15}]=
\zeta_{17}(D-8)-\eta\ptl_{x_{11}}+\zeta_2E'_{(1,1,2,2,1,1)}-\zeta_3E'_{(1,1,1,2,1,1)}\\&
&-\zeta_4E'_{(1,1,1,1,1,1)}+\zeta_6E'_{(1,0,1,1,1,1)}+\zeta_8E'_{(1,1,1,1)}
-\zeta_{10}E'_{(1,0,1,1)}-\zeta_{12}E'_{(1,0,1)}\\
&&+\zeta_{14}E_{-\al_6}+\zeta_{15}E_{-\al_1}+\zeta_{19}E_{\al_5}
+\zeta_{21}E_{(0,0,0,1,1)}-\zeta_{23}E_{(0,1,0,1,1)}\\&
&-\zeta_{22}E_{(0,0,1,1,1)}+\zeta_{24}E_{(0,1,1,1,1)}+\zeta_{25}E_{(0,1,1,2,1)}
-\zeta_{27}E_{(0,1,1,2,2,1)}\\& &-2\zeta_{17}
c-\frac{\zeta_{17}}{3}(2\al_1+\al_3-\al_5+\al_6),
 \hspace{6.6cm}(12.6.67)\end{eqnarray*}
\begin{eqnarray*}T_{18}&=&-[\iota(E_{-\al_4}),T_{16}]=
\zeta_{18}(D-8)-\eta\ptl_{x_{10}}-\zeta_1E'_{(1,1,2,3,2,1)}+\zeta_4E'_{(0,1,1,2,2,1)}\\&
&-\zeta_5E'_{(0,1,1,2,1,1)}-\zeta_8E'_{(0,1,1,2,1)}-\zeta_9E'_{(0,0,1,1,1,1)}
-\zeta_{11}E'_{(0,0,0,1,1,1)}-\zeta_{12}E'_{(0,0,1,1,1)}\\
&&-\zeta_{13}E'_{(0,0,1,1)}-\zeta_{15}E'_{(0,0,0,1,1)}-\zeta_{16}E_{-\al_4}
-\zeta_{20}E_{\al_2}-\zeta_{21}E_{\al_1}\\&
&-\zeta_{22}E_{(1,0,1)}-\zeta_{25}E_{(1,1,1,1)}-\zeta_{26}E_{(1,1,1,1,1)}
-\zeta_{27}E_{(1,1,1,1,1,1)}\\& &-2\zeta_{16}
c+\frac{\zeta_{18}}{3}(\al_1-\al_3-3\al_4-2\al_5-\al_6),
 \hspace{5.4cm}(12.6.68)\end{eqnarray*}
\begin{eqnarray*}T_{19}&=&[\iota(E_{-\al_1}),T_{16}]=
\zeta_{19}(D-8)-\eta\ptl_{x_9}-\zeta_2E'_{(1,1,2,2,2,1)}-\zeta_3E'_{(1,1,1,2,2,1)}\\&
&-\zeta_5E'_{(1,1,1,1,1,1)}+\zeta_7E'_{(1,0,1,1,1,1)}-\zeta_8E'_{(1,1,1,1,1)}
+\zeta_{10}E'_{(1,0,1,1,1)}-\zeta_{13}E'_{(1,0,1)}\\
&&+\zeta_{14}E'_{(0,0,0,0,1,1)}+\zeta_{16}E_{-\al_1}-\zeta_{17}E_{-\al_5}
+\zeta_{21}E_{\al_4}-\zeta_{23}E_{(0,1,0,1)}\\&
&-\zeta_{22}E_{(0,0,1,1)}+\zeta_{24}E_{(0,1,1,1)}-\zeta_{26}E_{(0,1,1,2,1)}
-\zeta_{27}E_{(0,1,1,2,1,1)}\\& &-2\zeta_{19}
c-\frac{\zeta_{19}}{3}(2\al_1+\al_3+2\al_5+\al_6),
 \hspace{6.4cm}(12.6.69)\end{eqnarray*}
\begin{eqnarray*}T_{20}&=&[\iota(E_{-\al_2}),T_{18}]=
\zeta_{20}(D-8)-\eta\ptl_{x_8}-\zeta_1E'_{(1,2,2,3,2,1)}-\zeta_6E'_{(0,1,1,2,2,1)}\\&
&+\zeta_7E'_{(0,1,1,2,1,1)}+\zeta_{10}E'_{(0,1,1,2,1)}+\zeta_9E'_{(0,1,1,1,1,1)}
+\zeta_{11}E'_{(0,1,0,1,1,1)}+\zeta_{12}E'_{(0,1,1,1,1)}\\
&&+\zeta_{13}E'_{(0,1,1,1)}+\zeta_{15}E'_{(0,1,0,1,1)}+\zeta_{16}E'_{(0,1,0,1)}
+\zeta_{18}E_{-\al_2}-\zeta_{23}E_{\al_1}\\&
&-\zeta_{24}E_{(1,0,1)}-\zeta_{25}E_{(1,0,1,1)}-\zeta_{26}E_{(1,0,1,1,1)}
-\zeta_{27}E_{(1,0,1,1,1,1)}\\& &-2\zeta_{20}
c+\frac{\zeta_{20}}{3}(\al_1-3\al_2-\al_3-3\al_4-2\al_5-\al_6),
 \hspace{4.3cm}(12.6.70)\end{eqnarray*}
\begin{eqnarray*}T_{21}&=&-[\iota(E_{-\al_4}),T_{19}]=
\zeta_{21}(D-8)-\eta\ptl_{x_7}+\zeta_2E'_{(1,1,2,3,2,1)}-\zeta_4E'_{(1,1,1,2,2,1)}\\&
&+\zeta_5E'_{(1,1,1,2,1,1)}+\zeta_8E'_{(1,1,1,2,1)}+\zeta_9E'_{(1,0,1,1,1,1)}
+\zeta_{12}E'_{(1,0,1,1,1)}+\zeta_{13}E'_{(1,0,1,1)}\\
&&+\zeta_{14}E'_{(0,0,0,1,1,1)}-\zeta_{17}E'_{(0,0,0,1,1)}+\zeta_{18}E_{-\al_1}
-\zeta_{19}E_{-\al_4}-\zeta_{22}E_{\al_3}\\&
&-\zeta_{23}E_{\al_2}-\zeta_{25}E_{(0,1,1,1)}-\zeta_{26}E_{(0,1,1,1,1)}
-\zeta_{27}E_{(0,1,1,1,1,1)}\\& &-2\zeta_{21}
c-\frac{\zeta_{21}}{3}(2\al_1+\al_3+3\al_4+2\al_5+\al_6),
 \hspace{5.3cm}(12.6.71)\end{eqnarray*}
\begin{eqnarray*}T_{22}&=&[\iota(E_{-\al_3}),T_{21}]=
\zeta_{22}(D-8)-\eta\ptl_{x_6}-\zeta_3E'_{(1,1,2,3,2,1)}+\zeta_4E'_{(1,1,2,2,2,1)}\\&
&-\zeta_5E'_{(1,1,2,2,1,1)}-\zeta_8E'_{(1,1,2,2,1)}+\zeta_{11}E'_{(1,0,1,1,1,1)}
+\zeta_{15}E'_{(1,0,1,1,1)}+\zeta_{16}E'_{(1,0,1,1)}\\
&&-\zeta_{14}E'_{(0,0,1,1,1,1)}+\zeta_{17}E'_{(0,0,1,1,1)}+\zeta_{18}E'_{(1,0,1)}
+\zeta_{19}E'_{(0,0,1,1)}+\zeta_{21}E_{-\al_3}\\&
&-\zeta_{24}E_{\al_2}-\zeta_{25}E_{(0,1,0,1)}-\zeta_{26}E_{(0,1,0,1,1)}
-\zeta_{27}E_{(0,1,0,1,1,1)}\\& &-2\zeta_{22}
c-\frac{\zeta_{22}}{3}(2\al_1+4\al_3+3\al_4+2\al_5+\al_6),
 \hspace{5.1cm}(12.6.72)\end{eqnarray*}
\begin{eqnarray*}T_{23}&=&[\iota(E_{-\al_2}),T_{21}]=
\zeta_{23}(D-8)-\eta\ptl_{x_5}+\zeta_2E'_{(1,2,2,3,2,1)}+\zeta_6E'_{(1,1,1,2,2,1)}\\&
&-\zeta_7E'_{(1,1,1,2,1,1)}-\zeta_{10}E'_{(1,1,1,2,1)}-\zeta_9E'_{(1,1,1,1,1,1)}
-\zeta_{12}E'_{(1,1,1,1,1)}-\zeta_{13}E'_{(1,1,1,1)}\\
&&-\zeta_{14}E'_{(0,1,0,1,1,1)}+\zeta_{17}E'_{(0,1,0,1,1)}+\zeta_{19}E'_{(0,1,0,1)}+\zeta_{20}E_{-\al_1}
+\zeta_{21}E_{-\al_2}\\&
&-\zeta_{24}E_{\al_3}-\zeta_{25}E_{(0,0,1,1)}-\zeta_{26}E_{(0,0,1,1,1)}
-\zeta_{27}E_{(0,0,1,1,1,1)}\\& &-2\zeta_{23}
c-\frac{\zeta_{23}}{3}(2\al_1+3\al_2+\al_3+3\al_4+2\al_5+\al_6),
 \hspace{4.1cm}(12.6.73)\end{eqnarray*}
\begin{eqnarray*}T_{24}&=&[\iota(E_{-\al_3}),T_{23}]=
\zeta_{24}(D-8)-\eta\ptl_{x_4}-\zeta_3E'_{(1,2,2,3,2,1)}-\zeta_6E'_{(1,1,2,2,2,1)}\\&
&+\zeta_7E'_{(1,1,2,2,1,1)}+\zeta_{10}E'_{(1,1,2,2,1)}-\zeta_{11}E'_{(1,1,1,1,1,1)}
+\zeta_{14}E'_{(0,1,1,1,1,1)}-\zeta_{15}E'_{(1,1,1,1,1)}\\
&&-\zeta_{16}E'_{(1,1,1,1)}-\zeta_{17}E'_{(0,1,1,1,1)}
-\zeta_{19}E'_{(0,1,1,1)}+\zeta_{20}E'_{(1,0,1)}+\zeta_{22}E_{-\al_2}\\&
&+\zeta_{23}E_{-\al_3}-\zeta_{25}E_{\al_4}-\zeta_{26}E_{(0,0,0,1,1)}
-\zeta_{27}E_{(0,0,0,1,1,1)}\\& &-2\zeta_{24}
c-\frac{\zeta_{24}}{3}(2\al_1+3\al_2+4\al_3+3\al_4+2\al_5+\al_6),
 \hspace{3.9cm}(12.6.74)\end{eqnarray*}
\begin{eqnarray*}T_{25}&=&[\iota(E_{-\al_4}),T_{24}]=
\zeta_{25}(D-8)-\eta\ptl_{x_3}+\zeta_4E'_{(1,2,2,3,2,1)}-\zeta_6E'_{(1,1,2,3,2,1)}\\&
&-\zeta_9E'_{(1,1,2,2,1,1)}-\zeta_{11}E'_{(1,1,1,2,1,1)}-\zeta_{12}E'_{(1,1,2,2,1)}
+\zeta_{14}E'_{(0,1,1,2,1,1)}-\zeta_{15}E'_{(1,1,1,2,1)}\\
&&+\zeta_{18}E'_{(1,1,1,1)}-\zeta_{17}E'_{(0,1,1,2,1)}
+\zeta_{20}E'_{(1,0,1,1)}+\zeta_{21}E'_{(0,1,1,1)}+\zeta_{22}E'_{(0,1,0,1)}\\&
&+\zeta_{23}E'_{(0,0,1,1)}+\zeta_{24}E_{-\al_4}-\zeta_{26}E_{\al_5}
-\zeta_{27}E_{(0,0,0,0,1,1)}\\& &-2\zeta_{25}
c-\frac{\zeta_{25}}{3}(2\al_1+3\al_2+4\al_3+6\al_4+2\al_5+\al_6),
 \hspace{3.9cm}(12.6.75)\end{eqnarray*}
\begin{eqnarray*}T_{26}&=&[\iota(E_{-\al_5}),T_{25}]=
\zeta_{26}(D-8)-\eta\ptl_{x_2}-\zeta_5E'_{(1,2,2,3,2,1)}+\zeta_7E'_{(1,1,2,3,2,1)}\\&
&-\zeta_9E'_{(1,1,2,2,2,1)}-\zeta_{11}E'_{(1,1,1,2,2,1)}+\zeta_{13}E'_{(1,1,2,2,1)}
+\zeta_{14}E'_{(0,1,1,2,2,1)}\\&
&+\zeta_{16}E'_{(1,1,1,2,1)}+\zeta_{18}E'_{(1,1,1,1,1)}+\zeta_{19}E'_{(0,1,1,2,1)}
+\zeta_{20}E'_{(1,0,1,1,1)}\\&
&+\zeta_{21}E'_{(0,1,1,1,1)}+\zeta_{22}E'_{(0,1,0,1,1)}+\zeta_{23}E'_{(0,0,1,1,1)}+\zeta_{24}E'_{(0,0,0,1,1)}+\zeta_{25}E_{-\al_5}
\\& &-\zeta_{27}E_{\al_6}-2\zeta_{26}
c-\frac{\zeta_{26}}{3}(2\al_1+3\al_2+4\al_3+6\al_4+5\al_5+\al_6),
 \hspace{2.2cm}(12.6.76)\end{eqnarray*}
\begin{eqnarray*}T_{27}&=&[\iota(E_{-\al_6}),T_{26}]=
\zeta_{27}(D-8)-\eta\ptl_{x_1}+\zeta_8E'_{(1,2,2,3,2,1)}-\zeta_{10}E'_{(1,1,2,3,2,1)}\\&
&+\zeta_{12}E'_{(1,1,2,2,2,1)}+\zeta_{13}E'_{(1,1,2,2,1,1)}+\zeta_{15}E'_{(1,1,1,2,2,1)}
+\zeta_{17}E'_{(0,1,1,2,2,1)}\\&
&+\zeta_{16}E'_{(1,1,1,2,1,1)}+\zeta_{18}E'_{(1,1,1,1,1,1)}+\zeta_{19}E'_{(0,1,1,2,1,1)}
+\zeta_{20}E'_{(1,0,1,1,1,1)}\\&
&+\zeta_{21}E'_{(0,1,1,1,1,1)}+\zeta_{22}E'_{(0,1,0,1,1,1)}+\zeta_{23}E'_{(0,0,1,1,1,1)}+\zeta_{24}E'_{(0,0,0,1,1,1)}+\zeta_{26}E_{-\al_6}
\\& &+\zeta_{25}E'_{(0,0,0,0,1,1)}-2\zeta_{27}
c-\frac{\zeta_{27}}{3}(2\al_1+3\al_2+4\al_3+6\al_4+5\al_5+4\al_6).
 \hspace{0.8cm}(12.6.77)\end{eqnarray*}
Then ${\mfk T}=\sum_{r=1}^{27}\mbb{F}T_r$ forms the 27-dimensional
${\msr G}^{E_6}$ of highest weight $\lmd_1$. Indeed, the map
$\zeta_r\mapsto T_r$ for $r\in\ol{1,27}$ determines a ${\msr
G}^{E_6}$-module isomorphism from  $U$ to ${\mfk T}$ (cf. (11.2.2)).

Denote
$$T'_r=T_r-\zeta_r(D-2c-8)+\eta\ptl_{x_{28-r}}\qquad\for\;\;r\in\ol{1,27}.\eqno(12.6.78)$$
Easily see that  ${\mfk T}'=\sum_{r=1}^{27}\mbb{F}T_r'$ forms the
27-dimensional ${\msr G}^{E_6}$-module of highest weight $\lmd_1$.
So we have the ${\msr G}^{E_6}$-module isomorphism from
$U=\sum_{r=1}^{27}\mbb{F}\zeta_r$ to ${\mfk T}'$ determined by
$\zeta_r\mapsto T_r'$ for $r\in\ol{1,27}$. The weight set of $U$ is
$$\Pi(U)=\{\sum_{s=1}^6b_{r,s}\lmd_s\mid r\in\ol{1,27}\}\eqno(12.6.79)$$
(cf. (12.6.9)-(12.6.14) and Table 11.2.1). Let $\lmd\in\Lmd^+$.
Denote
$$\Upsilon'(\lmd)=\{\lmd+\mu\mid
\mu\in\Pi(U),\;\;\lmd+\mu\in\Lmd^+\}\eqno(12.6.80)$$ (cf.
(12.6.27)). Take $M=V(\lmd)$, the irreducible ${\msr
G}^{E_6}$-module of highest weight $\lmd$. Theorem 5.4.3 implies
$$UV(\lmd)=U\otimes_{\mbb{F}}V(\lmd)\cong
\bigoplus_{\lmd'\in \Upsilon'(\lmd)}V(\lmd').\eqno(12.6.81)$$

Given $\lmd'\in\Upsilon'(\lmd)$, we pick a singular vector
$$u=\sum_{r=1}^{27}\zeta_ru_r\eqno(12.6.82)$$
of weight $\lmd'$ in $UV(\lmd)$, where $u_r\in V(\lmd)$. Moreover,
any singular vector of  weight $\lmd'$ in $UV(\lmd)$ is a scalar
multiple of $u$. Note that the vector
$$w=\sum_{r=1}^{27}T'_r(u_r)\eqno(12.6.83)$$
is also a  singular vector of  weight $\lmd'$ if it is not zero.
Thus
$$w=\flat_{\lmd'}u,\qquad \flat_{\lmd'}\in\mbb{F}.\eqno(12.6.84)$$
Set
$$\flat(\lmd)=\min\{\flat_{\lmd'}\mid\lmd'\in\Upsilon'(\lmd)\}.\eqno(12.6.85)$$
Recall the notion $\ell_\omega(\lmd)$ in (12.6.36). \psp

{\bf Theorem 12.6.6}. {\it The ${\msr G}^{E_7}$-module
$\wht{V(\lmd)}$ is irreducible if}
$$c\in\mbb{F}\setminus\{-8+5\mbb{N}/3,
(1/2)(\flat(\lmd)+\mbb{N})-4,\ell_\omega(\lmd)+\mbb{N}\}.\eqno(12.6.86)$$

{\it Proof}. Recall that the ${\msr G}^{E_7}$-submodule $U({\msr
G}_-)(V(\lmd))$ is irreducible by Proposition 12.5.2. It is enough
to prove $\wht{V(\lmd)}=U({\msr G}_-)(V(\lmd))$. It is obvious that
$$\wht{V(\lmd)}_{ 0}=V(\lmd)=(U({\msr G}_-)(V(\lmd)))_{
0}\eqno(12.6.87)$$ (cf. (12.6.3) and (12.6.4) with $M=V(\lmd)$).
Moreover, Lemma 12.6.1 with $M=V(\lmd)$,  (12.6.36) and (12.6.86)
imply that $\vf|_{\wht{V(\lmd)}_{ 1}}$ is invertible, or
equivalently,
$$\wht{V(\lmd)}_{ 1}=(U({\msr G}_-)(V(\lmd)))_{
1}.\eqno(12.6.88)$$ Suppose that
$$\wht{V(\lmd)}_{ i}=(U({\msr G}_-)(V(\lmd)))_{
i}\eqno(12.6.89)$$ for $i\in\ol{0,k}$ with $1\leq k\in\mbb{N}$.

For any $v\in V(\lmd)$ and $\al\in\mbb{N}^{16}$ such that
$|\al|=k-1$, we have
$$T_r(x^\al v)=x^\al[(|\al|-2c-8)\zeta_r+T'_r](v)+\eta \ptl_{x_{28-r}}(x^\al)v\in (U({\msr G}_-)(V(\lmd)))_{
k+1}\eqno(12.6.90)$$ for $r\in\ol{1,27}$ by (12.6.78), (12.6.89)
with $i=k-1,k$. If $k=1$, then $\al=0$. So $\eta
\ptl_{x_{28-r}}(x^\al)v=0$. When $k>1$,
$$\ptl_{x_{28-r}}(x^\al)v\in
\wht{V(\lmd)}_{ k-2}=(U({\msr G}_-)(V(\lmd)))_{ k-2}\eqno(12.6.91)$$
and so \begin{eqnarray*}\qquad & &[(\sum_{14\neq
p\in\ol{1,26}}\iota(\mfk b_p)\zeta_{28-p}-\iota(\mfk
b_{14})\zeta_{14}-\iota(\mfk
b_{27})\zeta_1](\ptl_{x_{28-r}}(x^\al)v)
\\&=&\eta(24-5(|\al|-1)+3c)(\ptl_{x_{28-r}}(x^\al)v)\in (U({\msr G}_-)(V(\lmd)))_{
k+1}\hspace{2.4cm}(12.6.92)\end{eqnarray*} by Lemma 12.6.5. Thus
(12.6.86) gives
$$\eta\ptl_{x_{28-r}}(x^\al)v\in (U({\msr G}_-)(V(\lmd)))_{
k+1}\qquad\for\;\;r\in\ol{1,27}.\eqno(12.6.93)$$ Hence in any case,
$$T_r(x^\al v)=x^\al[(|\al|-2c-8)\zeta_r+T'_r](v)\in (U({\msr G}_-)(V(\lmd)))_{
k+1}\qquad\for\;\;r\in\ol{1,27}.\eqno(12.6.94)$$

On the other hand,
$$V'=\mbox{Span}\{[(|\al|-2c-8)\zeta_r+T'_r](v)\mid r\in\ol{1,27},\;v\in
V(\lmd)\}\eqno(12.6.95)$$ forms a ${\msr G}^{E_6}$-submodule of
$UV(\lmd)$ with respect to the action in (12.6.1). Let $u$ be a
${\msr G}^{E_6}$-singular vector in (12.6.82). Then
$$V'\ni
\sum_{r=1}^{27}[(|\al|-2c-8)\zeta_r+T'_r](u_r)=(|\al|-2c-8)u+w=(|\al|-2c-8+\flat_{\lmd'})u\eqno(12.6.96)$$
by (12.6.83) and (12.6.84). Moreover, (12.6.85) and (12.6.86) yield
$u\in V'$. Since $UV(\lmd)$ is a ${\msr G}^{E_6}$-module generated
by all the singular vectors, we have $V'=UV(\lmd)$. So
$$x^\al UV(\lmd)\subset (U({\msr G}_-)(V(\lmd)))_{
k+1}.\eqno(12.6.97)$$ The arbitrariness of $\al$ implies
$$\zeta_r\wht{V(\lmd)}_{ k-1}\subset (U({\msr G}_-)(V(\lmd)))_{
k+1}\qquad\for\;\;r\in\ol{1,27}.\eqno(12.6.98)$$

Given any $f\in{\msr A}_k$ and $v\in V(\lmd)$, we have
$$\zeta_r\ptl_{x_s}(f)v\in \zeta_r\wht{V(\lmd)}_{ k-1}\subset (U({\msr G}_-)(V(\lmd)))_{
k+1}\qquad\for\;\;r,s\in\ol{1,27}.\eqno(12.6.99)$$ Moreover,
\begin{eqnarray*}\qquad\mfk b_s(fv)&=&\iota(\mfk b_s)(fv)=P_s(fv)+f(\td\omega-c)(x_sv)\\
&\equiv&
f(k+\td\omega-c)(x_sv)\;\;(\mbox{mod}\;\sum_{r=1}^{27}\zeta_r\wht{V(\lmd)}_{
k-1})\hspace{3.9cm}(12.6.100)\end{eqnarray*} for $s\in\ol{1,27}$ by
(12.4.20), (12.4.39)-(12.4.64), (12.5.13)-(12.5.39) and Lemma
12.6.1. According to (12.6.36), (12.6.86), (12.6.96) and (12.6.98),
we get
$$x_sfv\in (U({\msr G}_-)(V(\lmd)))_{
k+1}\qquad\for\;\;s\in\ol{1,27}.\eqno(12.6.101)$$ Thus (12.6.89)
holds for $i=k+1$. By induction on $k$, (12.6.89) holds for any
$i\in\mbb{N}$; that is, $\wht{V(\lmd)}=U({\msr
G}_-)(V(\lmd)).\qquad\Box$\psp

When $\lmd=0$, $V(0)$ is the one-dimensional trivial module and
$\ell_\omega(0)=\flat(0)=0$. So we have:\psp

{\bf Corollary 12.6.7}. {\it The ${\msr G}^{E_7}$-module
$\wht{V(0)}$ is irreducible if} $c\in\mbb{F}\setminus\{5\mbb{N}/3-8,
\mbb{N}/2-4\}.$ \psp

Next we consider $\lmd=k\lmd_1$. In this case,
$$\Upsilon(\lmd)=\{k\lmd_1+\lmd_6,(k-1)\lmd_1+\lmd_2,(k-1)\lmd_1\}
\eqno(12.6.102)$$ by (12.6.31) and Tables 11.1.1 and 11.2.1. Thus we
have
$$\ell_\omega(\lmd)=-16-\frac{4k}{3}\eqno(12.6.103)$$ by (12.6.36).
Moreover,
$$\Upsilon'(\lmd)=\{(k+1)\lmd_1,(k-1)\lmd_1+\lmd_3,(k-1)\lmd_1+\lmd_6\}\eqno(12.6.104)$$
by Table 11.2.1 and (12.6.80).

We define a representation of ${\msr G}^{E_6}$ on ${\msr
Z}=\mbb{F}[z_1,...,z_{27}]$ determined via (11.2.36)-(11.2.40) with
$U$ replaced by ${\msr Z}$ and $\zeta_i$ replaced by $z_i$.
 Then the ${\msr G}^{E_6}$-submodule ${\msr N}_k$
generated by $z_1^k$ is isomorphic to $V(k\lmd_1)$.

We calculate:
$$z_1^{k-1}z_2=-\frac{1}{k}E_{-\al_1}(z_1^k)\in{\msr
N}_k,\;\;z_1^{k-1}z_3=-E_{-\al_3}(z_1^{k-1}z_2)\in{\msr
N}_k,\eqno(12.6.105)$$
$$z_1^{k-1}z_4=-E_{-\al_4}(z_1^{k-1}z_3)\in{\msr N}_k,\;\;z_1^{k-1}z_5=-E_{-\al_5}(z_1^{k-1}z_4)\in{\msr
N}_k,\eqno(12.6.106)$$
$$z_1^{k-1}z_6=-E_{-\al_2}(z_1^{k-1}z_4)\in{\msr N}_k,\;\;z_1^{k-1}z_7=-E_{-\al_2}(z_1^{k-1}z_5)\in{\msr
N}_k,\eqno(12.6.107)$$
$$z_1^{k-1}z_9=-E_{-\al_4}(z_1^{k-1}z_7)\in{\msr
N}_k,\;\;z_1^{k-1}z_{11}=E_{-\al_3}(z_1^{k-1}z_9)\in{\msr
N}_k,\eqno(12.6.108)$$
$$-E_{-\al_1}(z_1^{k-1}z_{11})=z_1^{k-1}z_{14}+(k-1)z_1^{k-2}z_2z_{11}\in{\msr
N}_k,\eqno(12.6.109)$$
$$-E_{-\al_3}(z_1^{k-1}z_{14}+(k-1)z_1^{k-2}z_2z_{11})=(k-1)z_3z_{11}\in{\msr
N}_k,\eqno(12.6.110)$$
$$E_{\al_3}((k-1)z_3z_{11})=(k-1)(z_2z_{11}-z_3z_9)\in{\msr
N}_k,\eqno(12.6.111)$$
$$E_{\al_4}E_{-\al_4}[(k-1)(z_2z_{11}-z_3z_9)]=(k-1)(z_3z_9+z_4z_7)\in{\msr
N}_k,\eqno(12.6.112)$$
$$-E_{\al_5}E_{-\al_5}[(k-1)(z_3z_9+z_4z_7)]=(k-1)(z_4z_7+z_5z_6)\in{\msr
N}_k.\eqno(12.6.113)$$

Now we take $M={\msr N}_k$ in our earlier settings. First
$\zeta_1z_1^k$ is a singular vector of weight $(k+1)\lmd_1$ in
$U{\msr N}_k$. By (12.6.50), (11.2.40) and Table 11.2.1,
$$T_1'(z_1^k)=\frac{\zeta_1}{3}(4\al_1+3\al_2+5\al_3+6\al_4+4\al_5+2\al_6)(z_1^k)=\frac{4k}{3}\zeta_1z_1^k.\eqno(12.6.114)$$
So $\flat_{(k+1)\lmd_1}=4k/3$. Next
$\zeta_1z_1^{k-1}z_2-\zeta_2z_1^k$ is a singular vector of weight
$(k-1)\lmd_1+\lmd_3$ in $U{\msr N}_k$. According to (12.6.50),
(12.6.52), (11.2.40) and Table 11.2.1:
\begin{eqnarray*}& &T'_1(z_1^{k-1}z_2)-T'_2(z_1^k)\\
&=&\zeta_2z_1^k+\frac{4k-3}{3}\zeta_1z_1^{k-1}z_2-k\zeta_1z_1^{k-1}z_2-\frac{k}{3}\zeta_2z_1^k
=\frac{k-3}{3}(\zeta_1z_1^{k-1}z_2-\zeta_2z_1^k),
\hspace{0.7cm}(12.6.115)\end{eqnarray*} which gives
$\flat_{(k-1)\lmd_1+\lmd_3}=k/3-1$.

Expressions (11.2.29)-(11.2.34) and (12.6.105)-(12.6.113) show that
\begin{eqnarray*}\qquad
u&=&\zeta_1[4z_1^{k-1}z_{14}+(k-1)z_1^{k-2}(z_2z_{11}+z_3z_9-z_4z_7+z_5z_6)]-(k+3)z_1^{k-1}\\&&\times[\zeta_2z_{11}+\zeta_3z_9
-\zeta_4z_7+\zeta_5z_6+\zeta_6z_5-\zeta_7z_4+\zeta_9z_3+\zeta_{11}z_2-\zeta_{14}z_1]\hspace{0.7cm}(12.6.116)\end{eqnarray*}
is a singular vector of weight $(k-1)\lmd_1+\lmd_6$ in $U{\msr
N}_k$. We find
\begin{eqnarray*}&
&T_1'[4z_1^{k-1}z_{14}+(k-1)z_1^{k-2}(z_2z_{11}+z_3z_9-z_4z_7+z_5z_6)]\\
&
&-(k+3)[T_2'(z_1^{k-1}z_{11})+T_3'(z_1^{k-1}z_9)-T_4'(z_1^{k-1}z_7)
+T_5'(z_1^{k-1}z_6)\\
&&+T_6'(z_1^{k-1}z_5)-T_7'(z_1^{k-1}z_4)+T_9'(z_1^{k-1}z_3)+T_{11}'(z_1^{k-1}z_2)-T_{14}'(z_1^k)]\\&
=&(k+3)z_1^{k-1}(\zeta_2z_{11}+\zeta_3z_9-\zeta_4z_7+\zeta_5z_6+\zeta_6z_5-\zeta_7z_4+\zeta_9z_3+\zeta_{11}z_2)
\\ & &+(4k/3-2)\zeta_1[4z_1^{k-1}z_{14}+(k-1)z_1^{k-2}(z_2z_{11}+z_3z_9-z_4z_7+z_5z_6)]
-(k+3)z_1^{k-2}\\ &
&\times[2(4\zeta_1z_1z_{14}+(k-1)\zeta_1(z_2z_{11}+z_3z_9-z_4z_7+z_5z_6))
+(2k/3+8)\zeta_{14}z_1^2 \\ & &-(2k/3+7)z_1(\zeta_2z_{11}+\zeta_3z_9
-\zeta_4z_7+\zeta_5z_6+\zeta_6z_5-\zeta_7z_4+\zeta_9z_3+\zeta_{11}z_2)]
\\ &=&-(2k/3+8)u
 \hspace{10.7cm}(12.6.117)\end{eqnarray*}
by  (11.2.29)-(11.2.34), (11.2.40) with Table 11.2.1, (12.6.50),
(12.6.52)-(12.6.64), (12.6.78). Thus
$\flat_{(k-1)\lmd_1+\lmd_6}=-(2k/3+8)$. Therefore,
$$\flat(k\lmd_1)=-(2k/3+8).\eqno(12.6.118)$$
\pse

{\bf Corollary 12.6.8}. {\it The ${\msr G}^{E_7}$-module
$\wht{V(k\lmd_1)}$ is irreducible if}
$$c\in\mbb{F}\setminus\{-8+5\mbb{N}/3,
\mbb{N}/2-8-k/3,\mbb{N}-16-4k/3\}.\eqno(12.6.119)$$ \pse

By (12.4.5), (12.5.10) and (12.5.42), the above irreducible ${\msr
G}^{E_7}$-module $\wht{V(k\lmd_1)}$ is a highest-weight module with
highest weight $k\lmd_1+(c-2k/3)\lmd_7$.

\section{Combinatorics of the  Representation of
$E_6$}

In this section, we want to study the symmetry of $\ol{1,27}$ with
respect to the representation $\mfk r$ in (11.1.18)-(11.1.56).

Set
$$\mfk W=\{\mfk r(E_\al)\mid\al\in\Phi_{E_6}\}\eqno(12.7.1)$$ (cf. (4.4.42)-(4.4.44), (11.1.18)-(11.1.53) and (11.1.56)).
 Write
$$\mfk Z=\{\zeta_i\mid i\in\ol{1,27}\}\eqno(12.7.2)$$
(cf. (11.2.1) and (11.2.3)-(11.2.28)). First we have:\psp

{\bf Lemma 12.7.1}. {\it The set $\ol{1,27}$ is symmetric with
respect to the sets $\mfk W$ and $\mfk Z$.} \psp

{\it Proof}.  Let $\msr G_1$ be the Lie subalgebra of $\msr G^{E_6}$
generated by
$\{E_{\pm\al_1},E_{\pm\al_3},E_{\pm\al_4},E_{\pm\al_5},E_{\pm\al_6}\}$.
According to the Dynkin diagram, $\msr G_1\cong sl(6,\mbb C)$. Table
11.1.1 tells us that  $x_1$ and $x_{20}$ are  $\msr G_1$-singular
vectors of weight $\lmd_5$ and $x_6$ is a $\msr G_1$-singular vector
of weight
 $\lmd_2$. By (11.1.18), (11.1.20)-(11.1.23) and (11.1.56),
the $\msr G_1$-submodule generated by $x_1$ is $\sum_{i=1}^5\mbb
Fx_i+\mbb Fx_8$ and the $\msr G_1$-submodule generated by $x_{20}$
is $\mbb Fx_{20}+\sum_{r=23}^{27}\mbb Fx_r$. Thus the following
 subsets of indices:
$$\ol{1,5}\bigcup
\{8\},\;\;\ol{9,19}\bigcup\{6,7,21,22\},\;\;\{20\}\bigcup\ol{23,27}.\eqno(12.7.3)$$
 are symmetric with respect to the restricted
representation of $\msr G_1$

Let $\msr G_2$ be the Lie subalgebra of $\msr G^{E_6}$ generated by
$\{E_{\pm\al_1},E_{\pm\al_2},E_{\pm\al_3},E_{\pm\al_4},E_{\pm\al_5}\}$.
According to the Dynkin diagram, $\msr G_2\cong o(10,\mbb C)$. Table
11.1.1 says that $x_1$ is a $\msr G_2$-singular vector of weight 0,
$x_2$ is a $\msr G_2$-singular vector of weight $\lmd_5$ and
$x_{14}$ is a $\msr G_2$-singular vector of weight $\lmd_1$. So
$\mbb Fx_1$ is a trivial $\msr G_2$-module. The $\msr G_2$-submodule
generated by $x_{14}$ is $\mbb Fx_{14}+\mbb Fx_{17}+\mbb
Fx_{19}+\sum_{i=21}^{27}\mbb Fx_i$ by (11.1.18)-(11.1.22) and
(11.1.56). Now $\sum_{r=2}^{13}\mbb Fx_r+\sum_{s=15,16,18,20}\mbb
Fx_s$ form a spin module and $x_{20}$ is its lowest-weight vector.
Thus the subsets
$$\{1\},\;\;\{2\},\;\;\ol{3,10}\bigcup\{12,15\},\;\;\{11,13,16,18,20\},\;\;\{14,17,19\}
\bigcup\ol{21,27}\eqno(12.7.4)$$
are symmetric with respect to the restricted representation of $\msr
G_2$. Since $\msr G^{E_6}$ is generated by $\msr G_1$ and $\msr
G_2$, we have that the following subsets of indices:
$$\{1\},\;\;\{2\},\;\;\{3,4,5,8\},\;\;\{6,7,9,10,12,15\},\eqno(12.7.5)$$
$$\{11,13,16,18\},\;\;\{14,17,19,21,22\},\;\;\{20\},\;\;\ol{23,27}\eqno(12.7.6)$$
are symmetric with respect to the representation $\mfk r$ of $\msr
G^{E_6}$.

Denote
$$\es_i=(0,...,0,\stl{i}{1},0,...,0)\in\mbb Z^{27}.\eqno(12.7.7)$$
If $d=\sum_{r=1}^6a_rx_{i_r}\ptl_{x_{j_r}}\in\mfk W$ with
$a_r\in\mbb F$, we write
$$d^\clubsuit=\sum_{r=1}^6(\es_{i_r}-\es_{j_r})\in\mbb
Z^{27}.\eqno(12.7.8)$$ It turns out that the map $d\mapsto
d^\clubsuit$ is an injective map from $\mfk W$ to $\mbb F_3^{27}$.
Moreover, if $[d_1,d_2]=bd_3$ for some $d_1,d_2,d_3\in\mfk W$ and
$0\neq b\in\mbb F$, then
$$d_3^\clubsuit=d_1^\clubsuit+d_2^\clubsuit.\eqno(12.7.9)$$
According to (11.1.56),
$$(\mfk r(E_{-\al}))^\clubsuit=-(\mfk
r(E_{\al}))^\clubsuit\qquad\for\;\;\al\in\Phi_{E_6}^+.\eqno(12.7.10)$$

Let $\sgm$ be a permutation on $\ol{1,27}$ such that
$$\sgm(i)=i,\;\sgm(r)=28-r\qquad\for\;\;i\in\ol{13,15},\;i\in\ol{1,12}\bigcup\ol{16,27}.\eqno(12.7.11)$$
We define a linear transformation $\wht\sgm$ on $\mbb F_3^{27}$ by
$$\wht\sgm(\es_i)=\es_{\sgm(i)}\qquad\for\;\;i\in\ol{1,27}.\eqno(12.7.12)$$
We find that
$$\wht\sgm[(\mfk r(E_{\al_1}))^\clubsuit]=(\mfk
r(E_{-\al_6}))^\clubsuit,\;\;\wht\sgm[(\mfk
r(E_{\al_2}))^\clubsuit]=(\mfk
r(E_{-\al_2}))^\clubsuit,\eqno(12.7.13)$$
$$\wht\sgm[(\mfk
r(E_{\al_3}))^\clubsuit]=(\mfk
r(E_{-\al_5}))^\clubsuit,\;\;\wht\sgm[(\mfk
r(E_{\al_4}))^\clubsuit]=(\mfk
r(E_{-\al_4}))^\clubsuit\eqno(12.7.14)$$
$$\wht\sgm[(\mfk
r(E_{\al_5}))^\clubsuit]=(\mfk
r(E_{-\al_3}))^\clubsuit,\;\;\wht\sgm[(\mfk
r(E_{\al_6}))^\clubsuit]=(\mfk
r(E_{-\al_1}))^\clubsuit.\eqno(12.7.15)$$ According to (12.7.9),
(12.7.10) and (12.7.13)-(12.7.15), we find
$$\wht\sgm(\mfk W^\clubsuit)=\mfk W^\clubsuit.\eqno(12.7.16)$$
This shows that $\mfk W$ is symmetric under $\sgm$.

For $\zeta=\sum_{s=1}^5c_sx_{\ell_s}x_{k_s}\in\mfk Z$ and $0\neq
b\in\mbb F$, we define
$$(b\zeta)^\sharp=\sum_{s=1}^5(\es_{\ell_s}+\es_{k_s})\in\mbb
Z^{27}.\eqno(12.7.17)$$ We find that the map $\zeta\mapsto
\zeta^\sharp$ is an injective map from $\mfk Z$ to $\mbb Z^{27}$.
Moreover, if $d(\zeta)\neq 0$ for some $d\in\mfk W$ and
$\zeta\in\mfk Z$, then we have
$$(d(\zeta))^\sharp=d^\clubsuit+\zeta^\sharp.\eqno(12.7.18)$$
Furthermore, we define
$$\wht\sgm(\zeta)=\sum_{s=1}^5c_sx_{\sgm(\ell_s)}x_{\sgm(k_s)}\qquad\for\;\;
\zeta=\sum_{s=1}^5c_sx_{\ell_s}x_{k_s}\in\mfk Z.\eqno(12.7.19)$$ We
calculate
$$(\wht\sgm(\zeta_i))^\sharp=\zeta_i^\sharp,\;\;(\wht\sgm(\zeta_r))^\sharp=\zeta_{28-r}^\sharp
\qquad\for\;\;i\in\ol{13,15},\;\;r\in\ol{1,12}\bigcup\ol{16,27}.\eqno(12.7.20)$$
Therefore, $\mfk W$ and $\mfk Z$ are symmetric under $\sgm$. Thus
the subsets
$$\{1,27\},\;\{2,26\},\;\{3,25\},\;\{6,22\},\;\{11,17\},\;\{12,16\}\eqno(12.7.21)$$
are symmetric with respect to $\mfk W$ and $\mfk Z$. Expressions
(12.7.5), (12.7.6) and (12.7.21) imply that $\ol{1,27}$ is symmetric
with respect to $\mfk W$ and $\mfk Z$. $\qquad\Box$\psp

For any $i\in\ol{1,27}$, we define
$$\Upsilon_i=\{(r,s)\mid
\zeta_r\;\mbox{contains}\;x_ix_s\}.\eqno(12.7.22)$$ We have the
following equivalent combinatorial property:\psp

{\bf Lemma 12.7.2}. {\it Each polynomial $\zeta_r$ contains exactly
ten $x_i$'s. Moreover,
$$|\Upsilon_i|=10\qquad\for\;\;i\in\ol{16}.\eqno(12.7.23)$$}

{\it Proof}. According to (11.2.1) and (11.2.3)-(11.2.28), we
calculate
\begin{eqnarray*}\hspace{2.5cm}\Upsilon_1&=&\{(1,14), (2,17), (3,19), (4,21),
(5,22),\\& &(6,23), (7,24), (9,25), (11,26),
(14,27)\}.\hspace{3.4cm}(12.7.24)\end{eqnarray*} Then (12.7.23)
follows from Lemma 12.7.1. The first statement is obtained by
checking (11.2.1) and (11.2.3)-(11.2.28) one-by-one.$\qquad\Box$\psp

If we represent $\{x_i\mid i\in\ol{1,27}\}$ by 27 vertices and
represent $\{\zeta_r\mid r\in\ol{1,27}\}$ by 27 lines, then we
obtain a graph of $27$ vertices and 27 lines such that each line
contains 10 vertices and each vertex is exactly on 10 lines.

For $i\in\ol{1,27}$, we denote
$$I_i=\{r\in\ol{1,27}\mid\mbox{some}\;d\in\mfk W\;\mbox{contains}\;x_r\ptl_{x_i}\}\eqno(12.7.25)$$ and $$\mfk W_i=
\{d\in\mfk W\mid d\;\mbox{does not
contain}\;x_i\;\mbox{and}\;\ptl_{x_i}\}.\eqno(12.7.26)$$ Write
$$J_i=\{s\in\ol{1,27}\mid (r,s)\in\Upsilon_i\;\mbox{for
some}\;r\in\ol{1,27}\}.\eqno(12.7.27)$$ The following equivalent
combinatorial properties are crucial to our main result.\psp

{\bf Lemma 12.7.3}. {\it For any $i\in\ol{1,27}$, we have
$$|I_i|=16,\;\;|\mfk W_i|=40,\;\;|J_i|=10.\eqno(12.7.28)$$
In fact,
$$I_i\bigcup J_i=\ol{1,27}\setminus\{i\}.\eqno(12.7.29)$$
Moreover, every element in $\mfk W_i$ contains exactly two $x_s$
with $s\in J_i$, and for any $r\in J_i$, $x_r$ is contained in
exactly eight elements in $\mfk W_i$.}

{\it Proof}. Note that the elements in $\mfk W$ containing
$\ptl_{x_1}$ are
$$\mfk r(E_{-\al_6}),\;\;\mfk r(E'_{(0,0,0,0,1,1)}),\;\;\mfk r(E'_{(0,0,0,1,1,1)}),\;\;
\mfk r(E'_{(0,1,0,1,1,1)}),\eqno(12.7.30)$$
$$\mfk r(E'_{(0,0,1,1,1,1)}),\;\;\mfk r(E'_{(1,0,1,1,1,1)}),\;\;\mfk r(E'_{(0,1,1,1,1,1)}),\;\;\mfk r(E'_{(1,1,1,1,1,1)}),
\eqno(12.7.31)$$
$$\mfk r(E'_{(0,1,1,2,1,1)}),\;\;\mfk r(E'_{(1,1,1,2,1,1)}),\;\;\mfk r(E'_{(0,1,1,2,2,1)}),\;\;\mfk r(E'_{(1,1,2,2,1,1)}),
\eqno (12.7.32)$$
$$\mfk r(E'_{(1,1,1,2,2,1)}),\;\;\mfk r(E'_{(1,1,2,2,2,1)}),\;\;\mfk r(E'_{(1,1,2,3,2,1)}),\;\;\mfk r(E'_{(1,2,2,3,2,1)}).
\eqno(12.7.33)$$ Thus
$$I_1=\ol{2,13}\bigcup\{15,16,18,20\}.\eqno(12.7.34)$$
According to (12.7.24),
$$J_1=\{14,17,19\}\bigcup\ol{21,27}.\eqno(12.7.35)$$
Since $\tau(\mfk W)=-\mfk W$ (cf. (11.1.55)), $|\mfk W_1|=40$. So
(12.7.27) and (12.7.28) hold for $i=1$. The elements in $\mfk W_1$
containing $x_{14}$ are
$$\mfk r(E_{\al_1}),\;\;\mfk r(E_{(1,0,1)}),\;\;\mfk r(E_{(1,0,1,1)}),\;\;\mfk r(E_{(1,1,1,1)}),\eqno(12.7.36)$$
$$\mfk r(E_{(1,0,1,1,1)}),\;\;\mfk r(E_{(1,1,1,1,1}),\;\;\mfk r(E_{(1,1,1,2,1)}),\;\;\mfk r(E_{(1,1,2,2,1)}).\eqno(12.7.37)$$
Lemma 12.7.1 implies that  $x_r$ is contained in exactly eight
elements in $\mfk W_1$ for any $r\in J_1$. The fact that every
element in $\mfk W_1$ contains exactly two $x_s$ with $s\in J_1$ is
checked case by case. By Lemma 12.7.1, the lemma holds for any
$i\in\ol{1,27}$. $\qquad\Box$\psp

If we view $J_i$ as the set of vertices and $\mfk W_i$ as the set of
edges, the last statement in Lemma 12.7.3 gives a graph consisting
of 10 vertices and 40 edges such that each edge contain exactly two
vertices and each vertex is contained exactly 8 edges. Finally, we
have the following  duality:\psp

{\bf Lemma 12.7.4}. {\it Let $i\in\ol{1,27}$. For any $r\in J_i$,
there exists a unique $\Im_i(r)\in J_i$ such that an element in
$\mfk W_i$ contains $x_r$ if and only if it contains
$\ptl_{x_{\Im_i(r)}}$. Moreover, $\Im_i^2=1$.}

{\it Proof}. According to Lemma 12.7.3,  $x_r$ is contained in
exactly eight elements in ${\mfk W}_i$ for any $r\in J_i$. Since
$\tau(\mfk W_i)=-\mfk W_i$, we have that $\ptl_{x_r}$ is contained
in exactly eight elements in $\mfk W_i$ for any $r\in J_i$.

By Lemma 12.7.1, we only need to consider $i=1$. Every differential
operator in (12.7.35) and (12.7.36) contains $\ptl_{x_{27}}$. So
$\Im_1(14)=27$ and $\Im_1(27)=14$. By symmetry, there exists a
unique $\Im_1(r)\in J_1$ for any $r\in J_1$ such that an element in
$\mfk W_1$ contains $x_r$ if and only if it contains
$\ptl_{x_{\Im_1(r)}}$. Moreover, $\Im_1^2=1$. So the lemma holds for
$i=1$. $\qquad\Box$

\section{Representations on Exponential-Polynomial Functions}

In this section, we study a family of representations of the simple
Lie algebra of type $E_7$ on a space of exponential-polynomial
functions and prove that their irreducibility is related to an
explicit given algebraic variety.

Recall
$$D=\sum_{i=1}^{27}x_i\ptl_{x_i}.\eqno(12.8.1)$$
We identify $\widehat{V(0)}={\msr A}\otimes v_0$ with ${\msr A}$ by
$$f\otimes v_0\leftrightarrow f\qquad\for\;\;f\in{\msr
B},\eqno(12.8.2)$$ where $V(0)=\mbb Fv_0$. Then we have the
following one-parameter  inhomogeneous first-order differential
operator representation $\mfk r_c$ of $\msr G^{E_7}$:
$$\mfk r_c(u)=\mfk r(u)\qquad\for\;\;u\in \msr
G^{E_6}\eqno(12.8.3)$$(cf. (11.1.18)-(11.1.56)),
$$\mfk r_c(\hat\al)=-2D+3c,\qquad\mfk r_c(\mfk a_i)=\ptl_{x_i}\qquad\for\;\;i\in\ol{1,27},\eqno(12.8.4)$$
\begin{eqnarray*}\qquad\mfk r_c(\mfk b_1)&=&x_1(D-c)-\zeta_1\ptl_{x_{14}}-\zeta_2\ptl_{x_{17}}-\zeta_3\ptl_{x_{19}}
+\zeta_4\ptl_{x_{21}}-\zeta_5\ptl_{x_{22}}
\\&&-\zeta_6\ptl_{x_{23}}+\zeta_7\ptl_{x_{24}}-\zeta_9\ptl_{x_{25}}-\zeta_{11}\ptl_{x_{26}}
-\zeta_{14}\ptl_{x_{27}}.\hspace{3.8cm}(12.8.5)\end{eqnarray*}
\begin{eqnarray*}
\qquad\mfk r_c(\mfk b_2) &=&
x_2(D-c)-\zeta_1\ptl_{x_{11}}-\zeta_2\ptl_{x_{13}}-\zeta_3\ptl_{x_{16}}
+\zeta_4\ptl_{x_{18}}-\zeta_6\ptl_{x_{20}}\\ &
&-\zeta_8\ptl_{x_{22}}
+\zeta_{10}\ptl_{x_{24}}-\zeta_{12}\ptl_{x_{25}}-\zeta_{15}\ptl_{x_{26}}
+\zeta_{17}\ptl_{x_{27}},\hspace{3.5cm}(12.8.6)\end{eqnarray*}
\begin{eqnarray*}
\qquad\mfk r_c(\mfk b_3)&=&
x_3(D-c)-\zeta_1\ptl_{x_9}-\zeta_2\ptl_{x_{12}}-\zeta_3\ptl_{x_{15}}
+\zeta_5\ptl_{x_{18}}-\zeta_7\ptl_{x_{20}}\\ &
&-\zeta_8\ptl_{x_{21}}
+\zeta_{10}\ptl_{x_{23}}-\zeta_{13}\ptl_{x_{25}}-\zeta_{16}\ptl_{x_{26}}
+\zeta_{19}\ptl_{x_{27}},\hspace{3.5cm}(12.8.7)\end{eqnarray*}
\begin{eqnarray*}
\qquad\mfk r_c(\mfk b_4)&=&
x_4(D-c)-\zeta_1\ptl_{x_7}-\zeta_2\ptl_{x_{10}}-\zeta_4\ptl_{x_{15}}
+\zeta_5\ptl_{x_{16}}-\zeta_8\ptl_{x_{19}}\\ &
&-\zeta_9\ptl_{x_{20}}
+\zeta_{12}\ptl_{x_{23}}-\zeta_{13}\ptl_{x_{24}}-\zeta_{18}\ptl_{x_{26}}
+\zeta_{21}\ptl_{x_{27}},\hspace{3.5cm}(12.8.8)\end{eqnarray*}
\begin{eqnarray*}\qquad
\mfk r_c(\mfk b_5)&=&
x_5(D-c)+\zeta_1\ptl_{x_6}-\zeta_3\ptl_{x_{10}}+\zeta_4\ptl_{x_{12}}
-\zeta_5\ptl_{x_{13}}+\zeta_8\ptl_{x_{17}}\\ &
&+\zeta_{11}\ptl_{x_{20}}
-\zeta_{15}\ptl_{x_{23}}+\zeta_{16}\ptl_{x_{24}}-\zeta_{18}\ptl_{x_{25}}
-\zeta_{22}\ptl_{x_{27}},\hspace{3.4cm}(12.8.9)\end{eqnarray*}
\begin{eqnarray*}
\qquad\mfk r_c(\mfk b_6)&=&
x_6(D-c)+\zeta_1\ptl_{x_5}+\zeta_2\ptl_{x_8}-\zeta_6\ptl_{x_{15}}
+\zeta_7\ptl_{x_{16}}-\zeta_9\ptl_{x_{18}}\\ &
&-\zeta_{10}\ptl_{x_{19}}
+\zeta_{12}\ptl_{x_{21}}-\zeta_{13}\ptl_{x_{22}}+\zeta_{20}\ptl_{x_{26}}
-\zeta_{23}\ptl_{x_{27}},\hspace{3.2cm}(12.8.10)\end{eqnarray*}
\begin{eqnarray*}
\qquad\mfk r_c(\mfk b_7)&=&
x_7(D-c)-\zeta_1\ptl_{x_4}+\zeta_3\ptl_{x_8}+\zeta_6\ptl_{x_{12}}
-\zeta_7\ptl_{x_{13}}+\zeta_{10}\ptl_{x_{17}}\\ &
&+\zeta_{11}\ptl_{x_{18}}
-\zeta_{15}\ptl_{x_{21}}+\zeta_{16}\ptl_{x_{22}}+\zeta_{20}\ptl_{x_{25}}
+\zeta_{24}\ptl_{x_{27}},\hspace{3.2cm}(12.8.11)\end{eqnarray*}
 \begin{eqnarray*}
\qquad\mfk r_c(\mfk b_8)&=&
x_8(D-c)+\zeta_2\ptl_{x_6}+\zeta_3\ptl_{x_7}-\zeta_4\ptl_{x_9}
+\zeta_5\ptl_{x_{11}}-\zeta_8\ptl_{x_{14}}\\ &
&+\zeta_{14}\ptl_{x_{20}}
+\zeta_{17}\ptl_{x_{23}}-\zeta_{19}\ptl_{x_{24}}+\zeta_{21}\ptl_{x_{25}}
+\zeta_{22}\ptl_{x_{26}},\hspace{3.2cm}(12.8.12)\end{eqnarray*}
 \begin{eqnarray*}
\qquad\mfk r_c(\mfk b_9)&=&
x_9(D-c)-\zeta_1\ptl_{x_3}-\zeta_4\ptl_{x_8}-\zeta_6\ptl_{x_{10}}
+\zeta_9\ptl_{x_{13}}-\zeta_{12}\ptl_{x_{17}}\\ &
&-\zeta_{11}\ptl_{x_{16}}
+\zeta_{15}\ptl_{x_{19}}-\zeta_{18}\ptl_{x_{22}}-\zeta_{20}\ptl_{x_{24}}
+\zeta_{25}\ptl_{x_{27}},\hspace{3.2cm}(12.8.13)\end{eqnarray*}
\begin{eqnarray*}
\qquad\mfk r_c(\mfk b_{10})&=&
x_{10}(D-c)-\zeta_2\ptl_{x_4}-\zeta_3\ptl_{x_5}-\zeta_6\ptl_{x_9}
+\zeta_7\ptl_{x_{11}}-\zeta_{10}\ptl_{x_{14}}\\ &
&+\zeta_{14}\ptl_{x_{18}}
+\zeta_{17}\ptl_{x_{21}}-\zeta_{19}\ptl_{x_{22}}-\zeta_{23}\ptl_{x_{25}}
-\zeta_{24}\ptl_{x_{26}},\hspace{3.1cm}(12.8.14)\end{eqnarray*}
\begin{eqnarray*}
\qquad\mfk r_c(\mfk b_{11})&=&
x_{11}(D-c)-\zeta_1\ptl_{x_2}+\zeta_5\ptl_{x_8}+\zeta_7\ptl_{x_{10}}
-\zeta_9\ptl_{x_{12}}+\zeta_{11}\ptl_{x_{15}}\\ &
&+\zeta_{13}\ptl_{x_{17}}
-\zeta_{16}\ptl_{x_{19}}+\zeta_{18}\ptl_{x_{21}}+\zeta_{20}\ptl_{x_{23}}
+\zeta_{26}\ptl_{x_{27}},\hspace{3.1cm}(12.8.15)\end{eqnarray*}
 \begin{eqnarray*}
\qquad\mfk r_c(\mfk
b_{12})&=&x_{12}(D-c)-\zeta_2\ptl_{x_3}+\zeta_4\ptl_{x_5}+\zeta_6\ptl_{x_7}
-\zeta_9\ptl_{x_{11}}+\zeta_{12}\ptl_{x_{14}}\\ &
&-\zeta_{14}\ptl_{x_{16}}
-\zeta_{17}\ptl_{x_{19}}+\zeta_{21}\ptl_{x_{22}}+\zeta_{23}\ptl_{x_{24}}
-\zeta_{25}\ptl_{x_{26}},\hspace{3.1cm}(12.8.16)\end{eqnarray*}
 \begin{eqnarray*}
\qquad\mfk r_c(\mfk
b_{13})&=&x_{13}(D-c)-\zeta_2\ptl_{x_2}-\zeta_5\ptl_{x_5}-\zeta_7\ptl_{x_7}
+\zeta_9\ptl_{x_9}-\zeta_{13}\ptl_{x_{14}}\\ &
&+\zeta_{14}\ptl_{x_{15}}
+\zeta_{19}\ptl_{x_{19}}-\zeta_{21}\ptl_{x_{21}}-\zeta_{23}\ptl_{x_{23}}
-\zeta_{26}\ptl_{x_{26}},\hspace{3.1cm}(12.8.17)\end{eqnarray*}
\begin{eqnarray*}\qquad\mfk r_c(\mfk b_{14})&=&
x_{14}(D-c)-\zeta_1\ptl_{x_1}-\zeta_8\ptl_{x_8}-\zeta_{10}\ptl_{x_{10}}
+\zeta_{12}\ptl_{x_{12}}-\zeta_{13}\ptl_{x_{13}}\\ &
&-\zeta_{15}\ptl_{x_{15}}
+\zeta_{16}\ptl_{x_{16}}-\zeta_{18}\ptl_{x_{18}}-\zeta_{20}\ptl_{x_{20}}
+\zeta_{27}\ptl_{x_{27}},\hspace{3.1cm}(12.8.18)\end{eqnarray*}
\begin{eqnarray*}
\qquad\mfk r_c(\mfk
b_{15})&=&x_{15}(D-c)-\zeta_3\ptl_{x_3}-\zeta_4\ptl_{x_4}-\zeta_6\ptl_{x_6}
+\zeta_{11}\ptl_{x_{11}}+\zeta_{14}\ptl_{x_{13}}\\ &
&-\zeta_{15}\ptl_{x_{14}}
+\zeta_{17}\ptl_{x_{17}}-\zeta_{22}\ptl_{x_{22}}-\zeta_{24}\ptl_{x_{24}}
-\zeta_{25}\ptl_{x_{25}},\hspace{3.1cm}(12.8.19)\end{eqnarray*}
\begin{eqnarray*}
\qquad\mfk r_c(\mfk
b_{16})&=&x_{16}(D-c)-\zeta_3\ptl_{x_2}+\zeta_5\ptl_{x_4}+\zeta_7\ptl_{x_6}
-\zeta_{11}\ptl_{x_9}-\zeta_{14}\ptl_{x_{12}}\\ &
&+\zeta_{16}\ptl_{x_{14}}
-\zeta_{19}\ptl_{x_{17}}+\zeta_{22}\ptl_{x_{21}}+\zeta_{24}\ptl_{x_{23}}
-\zeta_{26}\ptl_{x_{25}},\hspace{3.1cm}(12.8.20)\end{eqnarray*}
\begin{eqnarray*}\qquad\mfk r_c(\mfk b_{17})&=&
x_{17}(D-c)-\zeta_2\ptl_{x_1}+\zeta_8\ptl_{x_5}+\zeta_{10}\ptl_{x_7}
-\zeta_{12}\ptl_{x_9}+\zeta_{13}\ptl_{x_{11}}\\ &
&+\zeta_{17}\ptl_{x_{15}}
-\zeta_{19}\ptl_{x_{16}}+\zeta_{21}\ptl_{x_{18}}+\zeta_{23}\ptl_{x_{20}}
-\zeta_{27}\ptl_{x_{26}},\hspace{3.1cm}(12.8.21)\end{eqnarray*}
  \begin{eqnarray*}
\qquad\mfk r_c(\mfk
b_{18})&=&x_{18}(D-c)+\zeta_4\ptl_{x_2}+\zeta_5\ptl_{x_3}-\zeta_9\ptl_{x_6}
+\zeta_{11}\ptl_{x_7}+\zeta_{14}\ptl_{x_{10}}\\ &
&-\zeta_{18}\ptl_{x_{14}}
+\zeta_{21}\ptl_{x_{17}}-\zeta_{22}\ptl_{x_{19}}+\zeta_{25}\ptl_{x_{23}}
+\zeta_{26}\ptl_{x_{24}},\hspace{3.1cm}(12.8.22)\end{eqnarray*}
 \begin{eqnarray*}
\qquad\mfk r_c(\mfk
b_{19})&=&x_{19}(D-c)-\zeta_3\ptl_{x_1}-\zeta_8\ptl_{x_4}-\zeta_{10}\ptl_{x_6}
+\zeta_{15}\ptl_{x_9}-\zeta_{17}\ptl_{x_{12}}\\ &
&-\zeta_{16}\ptl_{x_{11}}
+\zeta_{19}\ptl_{x_{13}}-\zeta_{22}\ptl_{x_{18}}-\zeta_{24}\ptl_{x_{20}}
-\zeta_{27}\ptl_{x_{25}},\hspace{2.9cm}(12.8.23)\end{eqnarray*}
\begin{eqnarray*}
\qquad\mfk r_c(\mfk
b_{20})&=&x_{20}(D-c)-\zeta_6\ptl_{x_2}-\zeta_7\ptl_{x_3}-\zeta_9\ptl_{x_4}
+\zeta_{11}\ptl_{x_5}+\zeta_{14}\ptl_{x_8}\\ &
&-\zeta_{20}\ptl_{x_{14}}
+\zeta_{23}\ptl_{x_{17}}-\zeta_{24}\ptl_{x_{19}}-\zeta_{25}\ptl_{x_{21}}
-\zeta_{26}\ptl_{x_{22}},\hspace{2.9cm}(12.8.24)\end{eqnarray*}
 \begin{eqnarray*}
\qquad\mfk r_c(\mfk
b_{21})&=&x_{21}(D-c)+\zeta_4\ptl_{x_1}-\zeta_8\ptl_{x_3}+\zeta_{12}\ptl_{x_6}
-\zeta_{15}\ptl_{x_7}+\zeta_{17}\ptl_{x_{10}}\\ &
&+\zeta_{18}\ptl_{x_{11}}
-\zeta_{21}\ptl_{x_{13}}+\zeta_{22}\ptl_{x_{16}}-\zeta_{25}\ptl_{x_{20}}
+\zeta_{27}\ptl_{x_{24}},\hspace{2.9cm}(12.8.25)\end{eqnarray*}
  \begin{eqnarray*}\qquad
\mfk r_c(\mfk
b_{22})&=&x_{22}(D-c)-\zeta_5\ptl_{x_1}-\zeta_8\ptl_{x_2}-\zeta_{13}\ptl_{x_6}
+\zeta_{16}\ptl_{x_7}-\zeta_{19}\ptl_{x_{10}}\\ &
&-\zeta_{18}\ptl_{x_9}
+\zeta_{21}\ptl_{x_{12}}-\zeta_{22}\ptl_{x_{15}}-\zeta_{26}\ptl_{x_{20}}
-\zeta_{27}\ptl_{x_{23}},\hspace{3cm}(12.8.26)\end{eqnarray*}
 \begin{eqnarray*}\qquad
\mfk r_c(\mfk
b_{23})&=&x_{23}(D-c)-\zeta_6\ptl_{x_1}+\zeta_{10}\ptl_{x_3}+\zeta_{12}\ptl_{x_4}
-\zeta_{15}\ptl_{x_5}+\zeta_{17}\ptl_{x_8}\\ &
&+\zeta_{20}\ptl_{x_{11}}
-\zeta_{23}\ptl_{x_{13}}+\zeta_{24}\ptl_{x_{16}}+\zeta_{25}\ptl_{x_{18}}
-\zeta_{27}\ptl_{x_{22}},\hspace{2.9cm}(12.8.27)\end{eqnarray*}
 \begin{eqnarray*}\qquad
\mfk r_c(\mfk
b_{24})&=&x_{24}(D-c)+\zeta_7\ptl_{x_1}+\zeta_{10}\ptl_{x_2}-\zeta_{13}\ptl_{x_4}
+\zeta_{16}\ptl_{x_5}-\zeta_{19}\ptl_{x_8}\\ &
&-\zeta_{20}\ptl_{x_9}
+\zeta_{23}\ptl_{x_{12}}-\zeta_{24}\ptl_{x_{15}}+\zeta_{26}\ptl_{x_{18}}
+\zeta_{27}\ptl_{x_{21}},\hspace{3cm}(12.8.28)\end{eqnarray*}
  \begin{eqnarray*}\qquad\mfk r_c(\mfk b_{25})&=&x_{25}(D-c)-\zeta_9\ptl_{x_1}-\zeta_{12}\ptl_{x_2}-\zeta_{13}\ptl_{x_3}
-\zeta_{18}\ptl_{x_5}+\zeta_{21}\ptl_{x_8}\\ &
&+\zeta_{20}\ptl_{x_7}
-\zeta_{23}\ptl_{x_{10}}-\zeta_{25}\ptl_{x_{15}}-\zeta_{26}\ptl_{x_{16}}
-\zeta_{27}\ptl_{x_{19}},\hspace{3cm}(12.8.29)\end{eqnarray*}
  \begin{eqnarray*}\qquad
\mfk r_c(\mfk
b_{26})&=&x_{26}(D-c)-\zeta_{11}\ptl_{x_1}-\zeta_{15}\ptl_{x_2}-\zeta_{16}\ptl_{x_3}
-\zeta_{18}\ptl_{x_4}+\zeta_{22}\ptl_{x_8}\\ &
&+\zeta_{20}\ptl_{x_6}
-\zeta_{24}\ptl_{x_{10}}-\zeta_{25}\ptl_{x_{12}}-\zeta_{26}\ptl_{x_{13}}
-\zeta_{27}\ptl_{x_{17}},\hspace{3cm}(12.8.30)\end{eqnarray*}
  \begin{eqnarray*}
\qquad\mfk r_c(\mfk
b_{27})&=&x_{27}(D-c)-\zeta_{14}\ptl_{x_1}+\zeta_{17}\ptl_{x_2}+\zeta_{19}\ptl_{x_3}
+\zeta_{21}\ptl_{x_4}-\zeta_{22}\ptl_{x_5}\\ &
&-\zeta_{23}\ptl_{x_6}
+\zeta_{24}\ptl_{x_7}+\zeta_{25}\ptl_{x_9}+\zeta_{26}\ptl_{x_{11}}
+\zeta_{27}\ptl_{x_{14}}.\hspace{3.4cm}(12.8.31)\end{eqnarray*}

Let $\vec a=(a_1,...,a_{27})\in\mbb F^{27}\setminus\{\vec 0\}$.
Define
$$\vec a\cdot \vec x=\sum_{i=1}^{27}a_ix_i.\eqno(12.8.32)$$
Recall $\msr A$ in (11.1.16) and set
$${\msr A}_{\vec a}=\{fe^{\vec a\cdot\vec
x}\mid f\in{\msr A}\}.\eqno(12.8.33)$$ For $i\in\ol{1,27}$, we write
$$\Upsilon_i=\{(\iota(r),r)\mid r\in J_i\}.\eqno(12.8.34)$$

Recall $\eta$ in (11.2.42) and the map $\Im_i$ in Lemma 12.7.4.
Denote
$$\msr F=\{x_i\zeta_{\iota(r)}\zeta_{\iota(s)}\eta\mid i\in\ol{1,27};\: r,s\in
J_i,\;r\neq s,\Im_i(s)\}.\eqno(12.8.35)$$ For any function
$f(x_1,...,x_{27})$, we define
$$f(\vec b)=f(b_1,...,b_{27})\qquad\for\;\;\vec b=(b_1,...,b_{27})\in\mbb
F^{27}.\eqno(12.8.36)$$ Now we define
$$\msr V=\{\vec b\in\mbb F^{27}\setminus\{\vec 0\}\mid f(\vec b)=0\;\for\;f\in\msr
F\},\eqno(12.8.37)$$which gives rise to a projective algebraic
variety. The following is our main theorem in this section.\psp

{\bf Theorem 12.8.1}. {\it With respect to the representation
$\pi_c$, ${\msr A}_{\vec a}$ forms an irreducible $\msr
G^{E_7}$-module if $\vec a\not\in\msr V$.}

{\it Proof}. Let $\msr{A}_k$ be the subspace of homogeneous
polynomials with degree $k$. Set
$$\msr{A}_{\vec a,k}=\msr{A}_ke^{\vec a\cdot\vec
x}\qquad\for\;k\in\mbb{N}.\eqno(12.8.38)$$ Let  ${M}$ be a nonzero
$\msr G^{E_7}$-submodule of ${\msr A}_{\vec a}$. Take any $0\neq
fe^{\vec a\cdot\vec x}\in  M$ with $f\in \msr{A}$. By the second
equation in (12.8.4),
$$(\xi_i-a_i)(fe^{\vec a\cdot\vec x})=\ptl_{x_i}(f)e^{\vec a\cdot\vec
x}\in M\qquad\for\;\;i\in\ol{1,27}.\eqno(12.8.39)$$ Repeatedly
applying (12.8.39) if necessarily, we obtain $e^{\vec a\cdot\vec
x}\in M$; that is, $\msr{A}_{\vec a,0}\subset M$.

Suppose $\msr{A}_{\vec a,\ell}\subset M$ for some $\ell\in\mbb{N}$.
Let $ge^{\vec a\cdot\vec x}$ be any element in $\msr{A}_{\vec
a,\ell}$. First we assume $a_1\neq 0$. Applying (12.7.30)-(12.7.33)
to $ge^{\vec a\cdot\vec x}$, we get by (11.1.18)-(11.1.53) and
(11.1.56) that
$$(a_1x_2-a_{11}x_{14}-a_{13}x_{17}
-a_{16}x_{19}-a_{18}x_{21}-a_{20}x_{23})ge^{\vec a\cdot\vec x}\equiv
0\;\;(\mbox{mod}\;M),\eqno(12.8.40)$$
$$(a_1x_3-a_9x_{14}-a_{12}x_{17}
-a_{15}x_{19}+a_{18}x_{22}+a_{20}x_{24})ge^{\vec a\cdot\vec x}\equiv
0\;\;(\mbox{mod}\;M),\eqno(12.8.41)$$
$$(a_1x_4-a_7x_{14}-a_{10}x_{17}+a_{15}x_{21}+a_{16}x_{22}-a_{20}x_{25})ge^{\vec a\cdot\vec x}\equiv
0\;\;(\mbox{mod}\;M),\eqno(12.8.42)$$
$$(a_1x_6+a_5x_{14}+a_8x_{17}-a_{15}x_{23}-a_{16}x_{24}-a_{18}x_{25})ge^{\vec a\cdot\vec x}\equiv
0\;\;(\mbox{mod}\;M),\eqno(12.8.43)$$
$$(a_1x_5+a_6x_{14}-a_{10}x_{19}-a_{12}x_{21}-a_{13}x_{22}
+a_{20}x_{26})ge^{\vec a\cdot\vec x}\equiv 0\;\;(\mbox{mod}\;
M),\eqno(12.8.44)$$
$$(a_1x_8+a_6x_{17}+a_7x_{19}+a_9x_{21}+a_{11}x_{22}
+a_{20}x_{27})ge^{\vec a\cdot\vec x}\equiv 0\;\;(\mbox{mod}\;
M),\eqno(12.8.45)$$
$$(a_1x_7-a_4x_{14}+a_8x_{19}+a_{12}x_{23}+a_{13}x_{24}
+a_{18}x_{26})ge^{\vec a\cdot\vec x}\equiv 0\;\;(\mbox{mod}\;
M),\eqno(12.8.46)$$
$$(a_1x_{10}
-a_4x_{17}-a_5x_{19}-a_9x_{23}-a_{11}x_{24}+a_{18}x_{27})ge^{\vec
a\cdot\vec x}\equiv 0\;\;(\mbox{mod}\;M),\eqno(12.8.47)$$
$$(a_1x_9-a_3x_{14}+a_8x_{21}-a_{10}x_{23}
+a_{13}x_{25}-a_{16}x_{26})ge^{\vec a\cdot\vec x}\equiv
0\;\;(\mbox{mod}\;M),\eqno(12.8.48)$$
$$(a_1x_{12}
-a_3x_{17}-a_5x_{21}+a_7x_{23}-a_{11}x_{25}-a_{16}x_{27})ge^{\vec
a\cdot\vec x}\equiv 0\;\;(\mbox{mod}\;M),\eqno(12.8.49)$$
$$(a_1x_{11}
-a_2x_{14}+a_8x_{22}-a_{10}x_{24}
-a_{12}x_{25}+a_{15}x_{26})ge^{\vec a\cdot\vec x}\equiv
0\;\;(\mbox{mod}\;M),\eqno(12.8.50)$$
$$(a_1x_{15}-a_3x_{19}+a_4x_{21}-a_6x_{23}+a_{11}x_{26}+a_{13}x_{27})ge^{\vec a\cdot\vec x}\equiv
0\;\;(\mbox{mod}\;M),\eqno(12.8.51)$$
$$(a_1x_{13}
-a_2x_{17}-a_5x_{22}+a_7x_{24} +a_9x_{25}+a_{15}x_{27})ge^{\vec
a\cdot\vec x}\equiv 0\;\;(\mbox{mod}\;M),\eqno(12.8.52)$$
$$(a_1x_{16}
-a_2x_{19}+a_4x_{22}-a_6x_{24} -a_9x_{26}-a_{12}x_{27})ge^{\vec
a\cdot\vec x}\equiv 0\;\;(\mbox{mod}\;M),\eqno(12.8.53)$$
$$(a_1x_{18}
-a_2x_{21}+a_3x_{22}-a_6x_{25}+a_7x_{26}+a_{10}x_{27})ge^{\vec
a\cdot\vec x}\equiv 0\;\;(\mbox{mod}\;M),\eqno(12.8.54)$$
$$(a_1x_{20}-a_2x_{23}+a_3x_{24}-a_4x_{25}+a_5x_{26}+a_8x_{27})ge^{\vec a\cdot\vec x}\equiv
0\;\;(\mbox{mod}\;M).\eqno(12.8.55)$$

Next we multiply $a_1$ to  (12.7.36) and (12.7.37), and apply them
to $ge^{\vec a\cdot\vec x}$:
$$a_1(a_8x_5+a_{10}x_7+a_{12}x_9+a_{13}x_{11}
+a_{17}x_{14}+a_{27}x_{26})ge^{\vec a\cdot\vec x}\equiv
0\;\;(\mbox{mod}\;M),\eqno(12.8.56)$$
$$a_1(a_8x_4+a_{10}x_6-a_{15}x_9-a_{16}x_{11}-a_{19}x_{14}
-a_{27}x_{25})ge^{\vec a\cdot\vec x}\equiv 0\;\;(\mbox{mod}\;
M),\eqno(12.8.57)$$
$$a_1(a_8x_3+a_{12}x_6+a_{15}x_7-a_{18}x_{11}
-a_{21}x_{14}+a_{27}x_{24})ge^{\vec a\cdot\vec x}\equiv
0\;\;(\mbox{mod}\;M),\eqno(12.8.58)$$
$$a_1(a_{10}x_3-a_{12}x_4-a_{15}x_5+a_{20}x_{11}+a_{23}x_{14}
+a_{27}x_{22})ge^{\vec a\cdot\vec x}\equiv 0\;\;(\mbox{mod}\;
M),\eqno(12.8.59)$$
$$a_1(a_8x_2+a_{13}x_6+a_{16}x_7+a_{18}x_9-a_{22}x_{14}-a_{27}x_{23})ge^{\vec a\cdot\vec x}\equiv 0\;\;(\mbox{mod}\;
M),\eqno(12.8.60)$$
$$a_1(a_{10}x_2-a_{13}x_4-a_{16}x_5-a_{20}x_9+a_{24}x_{14}-a_{27}x_{21})ge^{\vec a\cdot\vec x}\equiv 0\;\;(\mbox{mod}\;
M),\eqno(12.8.61)$$
$$a_1(a_{12}x_2-a_{13}x_3
-a_{18}x_5+a_{20}x_7+a_{25}x_{14}+a_{27}x_{19})ge^{\vec a\cdot\vec
x}\equiv 0\;\;(\mbox{mod}\;M),\eqno(12.8.62)$$
$$a_1(a_{15}x_2-a_{16}x_3
+a_{18}x_4-a_{20}x_6-a_{26}x_{14}-a_{27}x_{17})ge^{\vec a\cdot\vec
x}\equiv 0\;\;(\mbox{mod}\;M).\eqno(12.8.63)$$

Note that (12.8.44), (12.8.46), (12.8.48) and (12.8.50) yield
$$a_1x_5ge^{\vec a\cdot\vec x}\equiv
-(a_6x_{14}-a_{10}x_{19}-a_{12}x_{21}-a_{13}x_{22}
+a_{20}x_{26})ge^{\vec a\cdot\vec x}\;\;(\mbox{mod}\;
M),\eqno(12.8.64)$$
$$a_1x_7ge^{\vec a\cdot\vec x}\equiv(a_4x_{14}-a_8x_{19}-a_{12}x_{23}-a_{13}x_{24} -a_{18}x_{26})ge^{\vec
a\cdot\vec x}\;\;(\mbox{mod}\;M),\eqno(12.8.65)$$
$$a_1x_9ge^{\vec a\cdot\vec x}\equiv(a_3x_{14}-a_8x_{21}+a_{10}x_{23}
-a_{13}x_{25}+a_{16}x_{26})ge^{\vec a\cdot\vec
x}\;\;(\mbox{mod}\;M),\eqno(12.8.66)$$
$$a_1x_{11}ge^{\vec a\cdot\vec x}\equiv(a_2x_{14}-a_8x_{22}+a_{10}x_{24}
+a_{12}x_{25}-a_{15}x_{26})ge^{\vec a\cdot\vec
x}\;\;(\mbox{mod}\;M).\eqno(12.8.67)$$ Substituting
(12.8.33)-(12.8.37) into (12.8.26), we get
\begin{eqnarray*}\qquad\qquad& &[a_1a_{17}x_{14}-a_8(a_6x_{14}-a_{10}x_{19}-a_{12}x_{21}-a_{13}x_{22}
+a_{20}x_{26})\\ & &+a_1a_{27}x_{26}+
+a_{10}(a_4x_{14}-a_8x_{19}-a_{12}x_{23}-a_{13}x_{24} -a_{18}x_{26})
\\& &+a_{12}(a_3x_{14}-a_8x_{21}+a_{10}x_{23}
-a_{13}x_{25}+a_{16}x_{26})\\
& &+a_{13}(a_2x_{14}-a_8x_{22}+a_{10}x_{24}
+a_{12}x_{25}-a_{15}x_{26})]ge^{\vec a\cdot\vec x}
\\ &=&[(a_1a_{17}+a_2a_{13}+a_3a_{12}+a_4a_{10}-a_6a_8)x_{14}\\ & &
+(a_1a_{27}-a_8a_{20}-a_{10}a_{18}+a_{12}a_{16}-a_{13}a_{15})x_{26}a_{12}a_{16}
\\&=& (\zeta_2(\vec a)x_{14}+\zeta_{14}(\vec a)x_{26})ge^{\vec a\cdot\vec
x}\equiv 0\;\;(\mbox{mod}\; M)\hspace{4.3cm}(12.8.68)\end{eqnarray*}
by(11.2.3) and (11.2.15).

Substituting (12.8.40)-(12.8.55) into (12.8.57)-(12.8.63), we can
similarly obtain
$$(\zeta_3(\vec a)x_{14}+\zeta_{14}(\vec a)x_{25})ge^{\vec a\cdot\vec
x}\equiv 0\;\;(\mbox{mod}\;M),\eqno(12.8.69)$$
$$(\zeta_4(\vec a)x_{14}+\zeta_{14}(\vec a)x_{24})ge^{\vec a\cdot\vec
x}\equiv 0\;\;(\mbox{mod}\;M),\eqno(12.8.70)$$
$$(\zeta_6(\vec a)x_{14}+\zeta_{14}(\vec a)x_{22})ge^{\vec a\cdot\vec
x}\equiv 0\;\;(\mbox{mod}\;M),\eqno(12.8.71)$$
$$(\zeta_5(\vec a)x_{14}+\zeta_{14}(\vec a)x_{23})ge^{\vec a\cdot\vec
x}\equiv 0\;\;(\mbox{mod}\;M),\eqno(12.8.72)$$
$$(\zeta_7(\vec a)x_{14}+\zeta_{14}(\vec a)x_{21})ge^{\vec a\cdot\vec
x}\equiv 0\;\;(\mbox{mod}\;M),\eqno(12.8.73)$$
$$(\zeta_9(\vec a)x_{14}+\zeta_{14}(\vec a)x_{19})ge^{\vec a\cdot\vec
x}\equiv 0\;\;(\mbox{mod}\;M),\eqno(12.8.74)$$
$$(\zeta_{11}(\vec a)x_{14}+\zeta_{14}(\vec a)x_{17})ge^{\vec a\cdot\vec
x}\equiv 0\;\;(\mbox{mod}\;M).\eqno(12.8.75)$$

By Lemma 12.7.4,  we dually have
$$(\zeta_1(\vec a)x_{17}+\zeta_{11}(\vec a)x_{27})ge^{\vec a\cdot\vec
x}\equiv 0\;\;(\mbox{mod}\;M),\eqno(12.8.76)$$
$$(\zeta_1(\vec a)x_{19}+\zeta_9(\vec a)x_{27})ge^{\vec a\cdot\vec
x}\equiv 0\;\;(\mbox{mod}\;M),\eqno(12.8.77)$$
$$(\zeta_1(\vec a)x_{21}+\zeta_7(\vec a)x_{27})ge^{\vec a\cdot\vec
x}\equiv 0\;\;(\mbox{mod}\;M),\eqno(12.8.78)$$
$$(\zeta_1(\vec a)x_{23}+\zeta_5(\vec a)x_{27})ge^{\vec a\cdot\vec
x}\equiv 0\;\;(\mbox{mod}\;M),\eqno(12.8.79)$$
$$(\zeta_1(\vec a)x_{22}+\zeta_6(\vec a)x_{27})ge^{\vec a\cdot\vec
x}\equiv 0\;\;(\mbox{mod}\;M),\eqno(12.8.80)$$
$$(\zeta_1(\vec a)x_{24}+\zeta_4(\vec a)x_{27})ge^{\vec a\cdot\vec
x}\equiv 0\;\;(\mbox{mod}\;M),\eqno(12.8.81)$$
$$(\zeta_1(\vec a)x_{25}+\zeta_3(\vec a)x_{27})ge^{\vec a\cdot\vec
x}\equiv 0\;\;(\mbox{mod}\;M),\eqno(12.8.82)$$
$$(\zeta_1(\vec a)x_{26}+\zeta_2(\vec a)x_{27})ge^{\vec a\cdot\vec
x}\equiv 0\;\;(\mbox{mod}\;M).\eqno(12.8.83)$$

Now  if $\zeta_2(\vec a)\neq 0$, then (12.8.68) and (12.8.83) gives
$$\zeta_1(\vec a)x_{14}\equiv\zeta_{14}(\vec a)x_{27}\;\;(\mbox{mod}\;M).\eqno(12.8.84)$$
Moreover, (12.8.84) can be also derived from the following triples:
$$\{(12.8.69),(12.8.82),\zeta_3(\vec a)\neq
0\},\;\;\{(12.8.70),(12.8.81),\zeta_4(\vec a)\neq 0\},
\eqno(12.8.85)$$ $$\{(12.8.71),(12.8.80),\zeta_6(\vec a)\neq
0\},\;\;\{(12.8.42),(12.8.49),\zeta_5(\vec a)\neq
0\},\eqno(12.8.86)$$ $$\{(12.8.43),(12.8.48),\zeta_7(\vec a)\neq
0\},\;\;\{(12.8.74),(12.8.77),\zeta_9(\vec a)\neq
0\},\eqno(12.8.87)$$
$$\{(12.8.75),(12.8.76),\zeta_{11}(\vec a)\neq
0\}.\eqno(12.8.88)$$ So
$$(12.8.83)\;\;\mbox{holds if}\;\;\zeta_r(\vec a)\neq 0\;\;\mbox{for
some}\;\;r\in\{2,3,4,5,6,7,9,11\}.\eqno(12.8.89)$$

Note $\Im_1(14)=27$. Moreover,
$$\iota(J_1)=\{1,2,3,4,5,6,7,9,11,14\}\eqno(12.8.90)$$
by (12.7.24), where $\iota_1(14)=1$
 and $\iota_1(27)=14$.
Suppose
$$\{(x_1\zeta_r\zeta_{14}\vt)(\vec a)\mid
r\in\{2,3,4,5,6,7,9,11\}\}\neq\{0\}.\eqno(12.8.91)$$ By the above
paragraph, (12.8.84) holds. Substituting (12.8.68)-(12.8.75) into
$\zeta_{14}(\vec a)\times$(12.8.40), we get
\begin{eqnarray*} a_1\zeta_{14}(\vec a)x_2ge^{\vec
a\cdot\vec x}&\equiv &[a_{11}\zeta_{14}(\vec a)
-a_{13}\zeta_{11}(\vec a)-a_{16}\zeta_9(\vec a)\\
& &-a_{18}\zeta_7(\vec a)-a_{20}\zeta_5(\vec a)]x_{14}ge^{\vec
a\cdot\vec x}
\\&\equiv&[a_{11}(a_1a_{27}-a_8a_{20}-a_{10}a_{18}+a_{12}a_{16}
-a_{13}a_{15})\\& &-a_{13}(a_1a_{26}-a_5a_{20}-a_7a_{18}+a_9a_{16}
-a_{11}a_{15})\\ & &-a_{16}(a_1a_{25}+a_4a_{20}+a_6a_{18}-a_9a_{13}
+a_{11}a_{12})\\ & &-a_{18}(-a_1a_{24}+a_3a_{20}-a_6a_{16}+a_7a_{13}
-a_{10}a_{11})\\ &&-a_{20}(a_1a_{22}-a_3a_{18}-a_4a_{16}+a_5a_{13}
-a_8a_{11})]x_{14}ge^{\vec a\cdot\vec x}
\\&\equiv& a_1(a_{11}a_{27}-a_{13}a_{26}-a_{16}a_{25}+a_{18}a_{24}-a_{20}a_{22})x_{14}ge^{\vec a\cdot\vec x}
\\&\equiv &-a_1\zeta_{26}(\vec a)x_{14}ge^{\vec a\cdot\vec x}\;\;
(\mbox{mod}\;M)\hspace{5.6cm}(12.8.92)\end{eqnarray*} by (11.2.6),
(11.2.8), (11.2.10), (11.2.12), (11.2.15) and (11.2.27). Since
$a_1\neq 0$, we have
$$\zeta_{14}(\vec a)x_2ge^{\vec a\cdot\vec x}\equiv-\zeta_{26}(\vec a)x_{14}ge^{\vec a\cdot\vec x}\;\;(\mbox{mod}\;
M).\eqno(12.8.93)$$ Moreover, we substitute (12.8.68)-(12.8.75) and
(12.8.84) into $\zeta_{14}(\vec a)\times$[(12.8.41)-(12.8.55)]  and
can similarly derive
$$\zeta_{14}(\vec a)x_3ge^{\vec a\cdot\vec x}\equiv-\zeta_{25}(\vec a)x_{14}ge^{\vec a\cdot\vec x}\;\;(\mbox{mod}\;
M),\eqno(12.8.94)$$
$$\zeta_{14}(\vec a)x_4ge^{\vec a\cdot\vec x}\equiv-\zeta_{24}(\vec a)x_{14}ge^{\vec a\cdot\vec x}\;\;(\mbox{mod}\;
M),\eqno(12.8.95)$$
$$\zeta_{14}(\vec a)x_6ge^{\vec a\cdot\vec x}\equiv-\zeta_{22}(\vec a)x_{14}ge^{\vec a\cdot\vec x}\;\;(\mbox{mod}\;
M),\eqno(12.8.96)$$
$$\zeta_{14}(\vec a)x_5ge^{\vec a\cdot\vec x}\equiv-\zeta_{23}(\vec a)x_{14}ge^{\vec a\cdot\vec x}\;\;(\mbox{mod}\;
M),\eqno(12.8.97)$$
$$\zeta_{14}(\vec a)x_8ge^{\vec a\cdot\vec x}\equiv-\zeta_{20}(\vec a)x_{14}ge^{\vec a\cdot\vec x}\;\;(\mbox{mod}\;
M),\eqno(12.8.98)$$
$$\zeta_{14}(\vec a)x_7ge^{\vec a\cdot\vec x}\equiv-\zeta_{21}(\vec a)x_{14}ge^{\vec a\cdot\vec x}\;\;(\mbox{mod}\;
M),\eqno(12.8.99)$$
$$\zeta_{14}(\vec a)x_{10}ge^{\vec a\cdot\vec x}\equiv-\zeta_{18}(\vec a)x_{14}ge^{\vec a\cdot\vec x}\;\;(\mbox{mod}\;
M),\eqno(12.8.100)$$
$$\zeta_{14}(\vec a)x_9ge^{\vec a\cdot\vec x}\equiv-\zeta_{19}(\vec a)x_{14}ge^{\vec a\cdot\vec x}\;\;(\mbox{mod}\;
M),\eqno(12.8.101)$$
$$\zeta_{14}(\vec a)x_{12}ge^{\vec a\cdot\vec x}\equiv-\zeta_{16}(\vec a)x_{14}ge^{\vec a\cdot\vec x}\;\;(\mbox{mod}\;
M),\eqno(12.8.102)$$
$$\zeta_{14}(\vec a)x_{11}ge^{\vec a\cdot\vec x}\equiv-\zeta_{17}(\vec a)x_{14}ge^{\vec a\cdot\vec x}\;\;(\mbox{mod}\;
M),\eqno(12.8.103)$$
$$\zeta_{14}(\vec a)x_{15}ge^{\vec a\cdot\vec x}\equiv-\zeta_{13}(\vec a)x_{14}ge^{\vec a\cdot\vec x}\;\;(\mbox{mod}\;
M),\eqno(12.8.104)$$
$$\zeta_{14}(\vec a)x_{13}ge^{\vec a\cdot\vec x}\equiv-\zeta_{15}(\vec a)x_{14}ge^{\vec a\cdot\vec x}\;\;(\mbox{mod}\;
M),\eqno(12.8.105)$$
$$\zeta_{14}(\vec a)x_{16}ge^{\vec a\cdot\vec x}\equiv-\zeta_{12}(\vec a)x_{14}ge^{\vec a\cdot\vec x}\;\;(\mbox{mod}\;
M),\eqno(12.8.106)$$
$$\zeta_{14}(\vec a)x_{18}ge^{\vec a\cdot\vec x}\equiv-\zeta_{10}(\vec a)x_{14}ge^{\vec a\cdot\vec x}\;\;(\mbox{mod}\;
M),\eqno(12.8.107)$$
$$\zeta_{14}(\vec a)x_{20}ge^{\vec a\cdot\vec x}\equiv-\zeta_8(\vec a)x_{14}ge^{\vec a\cdot\vec x}\;\;(\mbox{mod}\;
M).\eqno(12.8.108)$$

Applying $\zeta_{14}(\vec a)\al_6$  to $ge^{\vec a\cdot\vec x}$, we
have
\begin{eqnarray*}\qquad& &\zeta_{14}(\vec
a)(a_1x_1-a_2x_2+a_{11}x_{11}
+a_{13}x_{13}-a_{14}x_{14}+a_{16}x_{16}-a_{17}x_{17}\\
&
&+a_{18}x_{18}-a_{19}x_{19}+a_{20}x_{20}-a_{21}x_{21}-a_{23}x_{23})ge^{\vec
a\cdot\vec x}\equiv 0\;\;(\mbox{mod}\;
M)\hspace{1.2cm}(12.8.109)\end{eqnarray*} by (11.1.54) and Table
11.1.1. According to (12.8.93), (12.8.103), (12.8.105), (12.8.106),
(12.8.75), (12.8.107), (12.8.74), (12.8.108), (12.8.73) and
(12.8.72), the above equation is equivalent to
\begin{eqnarray*}\hspace{2cm}& &a_1\zeta_{14}(\vec
a)x_1ge^{\vec a\cdot\vec x}+[a_2\zeta_{26}(\vec
a)-a_{11}\zeta_{17}(\vec a) -a_{13}\zeta_{15}(\vec
a)-a_{14}\zeta_{14}(\vec a)\\ & &-a_{16}\zeta_{12}(\vec
a)+a_{17}\zeta_{11}(\vec a)-a_{18}\zeta_{10}(\vec
a)+a_{19}\zeta_9(\vec a)-a_{20}\zeta_8(\vec a)\\
& &+a_{21}\zeta_7(\vec a)+a_{23}\zeta_5(\vec a)]x_{14}ge^{\vec
a\cdot\vec x}\equiv 0\;\;(\mbox{mod}\;
M)\hspace{3.6cm}(12.8.110)\end{eqnarray*} Moreover, (11.2.1) and
(11.2.3)-(11.2.28) imply
\begin{eqnarray*}\hspace{2cm}& &-a_{14}\zeta_{14}(\vec a)+a_{17}\zeta_{11}(\vec a)+a_{19}\zeta_9(\vec a)+a_{21}\zeta_7(\vec a)+a_{23}\zeta_5(\vec
a)\\ &=&-a_{14}(a_1a_{27}-a_8a_{20}-a_{10}a_{18}+a_{12}a_{16}
-a_{13}a_{15})\\ & &+a_{17}(a_1a_{26}-a_5a_{20}-a_7a_{18}+a_9a_{16}
-a_{11}a_{15})\\ & &+a_{19}(a_1a_{25}+a_4a_{20}+a_6a_{18}-a_9a_{13}
+a_{11}a_{12})\\& &+a_{21}(-a_1a_{24}+a_3a_{20}-a_6a_{16}+a_7a_{13}
-a_{10}a_{11})\\ & &+a_{23}(a_1a_{22}-a_3a_{18}-a_4a_{16}+a_5a_{13}
-a_8a_{11})\\&=&a_1\zeta_{27}(\vec a)-a_{14}(-a_8a_{20}-a_{10}a_{18}+a_{12}a_{16} -a_{13}a_{15})\\ &&+a_{17}(-a_5a_{20}-a_7a_{18}+a_9a_{16} -a_{11}a_{15})\\
& &+a_{19}(a_4a_{20}+a_6a_{18}-a_9a_{13}
+a_{11}a_{12})\\ & &+a_{21}(a_3a_{20}-a_6a_{16}+a_7a_{13} -a_{10}a_{11})\\
& &+a_{23}(-a_3a_{18}-a_4a_{16}+a_5a_{13}
-a_8a_{11}),\hspace{4.4cm}(12.8.111)
\end{eqnarray*}\qquad
\begin{eqnarray*}& &a_2\zeta_{26}(\vec
a)-a_{11}\zeta_{17}(\vec a) -a_{13}\zeta_{15}(\vec a)
-a_{16}\zeta_{12}(\vec a)-a_{18}\zeta_{10}(\vec
a)-a_{20}\zeta_8(\vec a)\\
&=&a_2(-a_{11}a_{27}+a_{13}a_{26}+a_{16}a_{25}-a_{18}a_{24}+a_{20}a_{22})
\\ &
&-a_{11}(-a_2a_{27}-a_8a_{23}-a_{10}a_{21}+a_{12}a_{19}-a_{15}a_{17})
\\& &-a_{13}(a_2a_{26}+a_5a_{23}+a_7a_{21}-a_9a_{19}+a_{14}a_{15})\\
& &-a_{16}(a_2a_{25}-a_4a_{23}-a_6a_{21}+a_9a_{17} -a_{12}a_{14})\\
&&-a_{18}(-a_2a_{24}-a_3a_{23}+a_6a_{19}-a_7a_{17} +a_{10}a_{14})\\
& &-a_{20}(a_2a_{22}+a_3a_{21}+a_4a_{19}-a_5a_{17} +a_8a_{14})
\\&=& -a_{11}(-a_8a_{23}-a_{10}a_{21}+a_{12}a_{19}-a_{15}a_{17})
\\& &-a_{13}(+a_5a_{23}+a_7a_{21}-a_9a_{19}+a_{14}a_{15})\\
& &-a_{16}(-a_4a_{23}-a_6a_{21}+a_9a_{17} -a_{12}a_{14})\\
&&-a_{18}(-a_3a_{23}+a_6a_{19}-a_7a_{17} +a_{10}a_{14})\\
& &-a_{20}(+a_3a_{21}+a_4a_{19}-a_5a_{17}
+a_8a_{14}).\hspace{6.3cm}(12.8.112)
\end{eqnarray*}
Thus we have
$$\zeta_{14}(\vec
a)x_1ge^{\vec a\cdot\vec x}\equiv-\zeta_{27}(\vec a)x_{14}ge^{\vec
a\cdot\vec x}\;\;(\mbox{mod}\;M).\eqno(12.8.113)$$

Applying $\zeta_{14}(\vec a)\times$ [the first equation in (12.8.4)]
to $ge^{\vec a\cdot\vec x}$, we get
$$\sum_{i=1}^{27}a_i\zeta_{14}(\vec a)x_ige^{\vec a\cdot\vec x}\equiv 0\;\;
(\mbox{mod}\;M).\eqno(12.8.114)$$ By (12.8.68)-(12.8.75), (12.8.84),
(12.8.93)-(12.8.108) and (12.8.113), the above equation is
equivalent to
$$[a_{14}\zeta_{14}(\vec a)+a_{27}\zeta_1(\vec a)-\sum_{14\neq
i\in\ol{1,26}}a_i\zeta_{28-i}(\vec a)x_{14}ge^{\vec a\cdot\vec
x}\equiv 0\;\;(\mbox{mod}\;M).\eqno(12.8.115)$$ Equivalently,
$$\eta(\vec a)x_{14}ge^{\vec a\cdot\vec
x}\equiv 0\;\;(\mbox{mod}\;M)\eqno(12.8.116)$$ by (11.2.41).
Assumption (12.8.91) implies $\vt(\vec a)\neq 0$ and
$\zeta_{14}(\vec a)\neq 0$. Thus $x_{14}ge^{\vec a\cdot\vec x}\in
M$. Substituting it to (12.8.68)-(12.8.75), (12.8.84),
(12.8.93)-(12.8.108) and (12.8.113), we get $x_ige^{\vec a\cdot\vec
x}\in M$ for any $i\in\ol{1,27}$. Hence $\msr{A}_{\vec
a,\ell+1}\subset M$. By induction,
$$\msr{A}_{\vec
a,k}\subset M\qquad\for\;\;k\in\mbb N.\eqno(12.8.117)$$ Therefore,
$M=\sum_{k=0}^\infty\msr{A}_{\vec a,k}=\msr A_{\vec a}$. So $\msr
A_{\vec a}$ forms an irreducible $\msr G^{E_7}$-module.

Now the theorem follows from Lemmas 12.7.1-11.7.4.$\qquad\Box$

\part{Related Topics}

\chapter{Oscillator Representations of $gl(n|m)$ and $osp(n|2m)$ }

Classical harmonic analysis says that the spaces of homogeneous
harmonic polynomials (solutions of Laplace equation) are irreducible
modules of the corresponding orthogonal Lie group (algebra) and the
whole polynomial algebra is a free module over the invariant
polynomials generated by harmonic polynomials.  In this chapter, we
first establish two-parameter $\mbb{Z}^2$-graded supersymmetric
oscillator generalizations of the above theorem for the Lie
superalgebra $gl(n|m)$. Then we extend the result to two-parameter
$\mbb{Z}$-graded supersymmetric oscillator generalizations of the
above theorem for the Lie superalgebras $osp(2n|2m)$ and
$osp(2n+1|2m)$. This is a reformulation of Luo and the author's work
[LX4].

In Section 13.1, we study the canonical supersymmetric orthogonal
oscillator representation of general linear Lie superalgebras.
Section 13.2 is a $\mbb Z^2$-graded oscillator generalizations of
the above representation obtained by partially swapping differential
operators and multiplication operators. Using the results in the
above two sections, we determine in Section 13.3 the structure of
the corresponding representations for even ortho-symplectic Lie
superalgebras. Finally in Section 13.4, we investigate the
oscillator representations of odd ortho-symplectic Lie superalgebras
based on the above three sections.

\section{Canonical Oscillator Representation of $gl(n|m)$}

In this section, we study the canonical supersymmetric orthogonal
oscillator representations of general linear Lie superalgebras.

Fix two positive integers $m$ and $n$.  Set
$$gl(n|m)_0=\sum_{i,j=1}^n\mbb{F}E_{i,j}+\sum_{r,s=1}^m\mbb{F}E_{n+r,n+s}\eqno(13.1.1)$$
and
$$gl(n|m)_1=\sum_{i=1}^n\sum_{r=1}^m(\mbb{F}E_{i,n+r}+\mbb{F}E_{n+r,i}).\eqno(13.1.2)$$
The general linear Lie superalgebra $gl(n|m)=gl(n|m)_0+gl(n|m)_1$
with the algebraic operation $[\cdot,\cdot]$ defined by
$$[A,B]=AB-(-1)^{i_1i_2}BA\qquad\for\;\;A\in gl(n|m)_{i_1},\;B\in
gl(n|m)_{i_2}.\eqno(13.1.3)$$

Let ${\msr A}$ be the polynomial algebra in bosonic variables
$\{x_i\mid i\in\ol{1,2n}\}$ and fermionic variables $\{\sta_j\mid
j\in\ol{1,2m}\}$; i.e.,
$$x_rx_s=x_sx_r,\;\sta_p\sta_q=-\sta_q\sta_p,\;\;x_r\sta_p=\sta_px_r\eqno(13.1.4)
$$ for $r,s\in\ol{1,2n}$ and $p,q\in\ol{1,2m}$. Set
$$\Theta=\sum_{p=1}^{2m}\mbb{F}\sta_p.\eqno(13.1.5)$$  Write
$${\msr A}_{(0)}=\sum_{q=0}^m\mbb{F}[x_1,...,x_{2n}]\Theta^{2q},\qquad{\msr
A}_{(1)}=\sum_{q=1}^{m-1}\mbb{F}[x_1,...,x_{2n}]\Theta^{2q+1}.\eqno(13.1.6)$$
Then ${\msr A}={\msr A}_{(0)}\oplus {\msr A}_{(1)}$ is a
$\mbb{Z}_2$-graded algebra.

For $r\in\ol{1,n}$, the usual differential operator $\ptl_{x_r}$
acts on ${\msr A}$ as a derivation  such that
$$\ptl_{x_r}(x_s)=\dlt_{r,s},\;\;\ptl_{x_r}(\sta_p)=0\qquad\for\;
s\in\ol{1,2n},\;p\in\ol{1,2m}.\eqno(13.1.7)$$  Moreover, for
$p\in\ol{1,2m}$, we define $\ptl_{\sta_p}$ as a linear operator on
${\msr A}$ with
$$\ptl_{\sta_p}(x_r)=0,\;\;\ptl_{\sta_p}(\sta_q)=\dlt_{p,q}\qquad\for\;r\in\ol{1,n},\;q\in\ol{1,2m}\eqno(13.1.8)$$ such that
$$\ptl_{\sta_p}(fg)=\ptl_{\sta_p}(f)g+(-1)^\iota f
\ptl_{\sta_p}(g)\qquad\for\;\;f\in{\msr A}_{(\iota)},\;g\in{\msr
A}.\eqno(13.1.9)$$ For later notational convenience, we redenote
$$y_i=x_{n+i},\qquad\vt_j=\sta_{m+j}\qquad\for\;\;i\in\ol{1,n},\;j\in\ol{1,m}.\eqno(13.1.10)$$
 Define a representation of $gl(n|m)$
on ${\msr A}$ determined by
$$E_{i,j}|_{\msr
A}=x_i\ptl_{x_j}-y_j\ptl_{y_i},\qquad E_{i,n+r}|_{\msr
A}=x_i\ptl_{\sta_r}-\vt_r\ptl_{y_i},\eqno(13.1.11)$$
$$E_{n+r,i}|_{\msr A}=\sta_r\ptl_{x_i}+y_i\ptl_{\vt_r},\qquad
E_{n+r,n+s}|_{\msr
A}=\sta_r\ptl_{\sta_s}-\vt_s\ptl_{\vt_r}\eqno(13.1.12)$$ for
$i,j\in\ol{1,n}$ and $r,s\in\ol{1,m}$.

\subsection{Representations of the Even Subalgebra}

Set
$$\bar{\msr A}=\mbb{F}[x_1,...,x_n,y_1,...,y_n],\qquad \check{\msr
A}=\sum_{i=0}^{2m}\Theta^i.\eqno(13.1.13)$$ Then $\bar{\msr A}$ and
$\check{\msr A}$ are subalgebras of ${\msr A}$, and
$${\msr
A}=\bar{\msr A}\check{\msr A}.\eqno(13.1.14)$$ Write
$$\bar{\msr G}=\sum_{i,j=1}^n\mbb{F}E_{i,j},\qquad\check{\msr
G}=\sum_{r,s=1}^m\mbb{F}E_{n+r,n+s}.\eqno(13.1.15)$$ Then they are
Lie subalgebras of $gl(n|m)$. In fact, the even subalgebra
$gl(n|m)_0=\bar{\msr G}+\check{\msr G}$. Let
$$\bar H=\sum_{i=1}^n\mbb{F}E_{i,i},\qquad\check H=\sum_{r=1}^m\mbb{F}E_{n+r,n+r},\eqno(13.1.16)$$
$$\bar{\msr G}_+=\sum_{1\leq i<j\leq n}\mbb{F}E_{i,j},\qquad\check{\msr
G}_+=\sum_{1\leq r<s\leq m}\mbb{F}E_{n+r,n+s}.\eqno(13.1.17)$$ We
take $\bar H$ as a Cartan subalgebra of $\bar{\msr G}$ and
$\bar{\msr G}_+$ as the subalgebra spanned by positive root vectors
in $\bar{\msr G}$. Similarly, we take $\check H$ as a Cartan
subalgebra of $\check {\msr G}$ and $\check {\msr G}_+$ as the
subalgebra spanned by positive root vectors in $\check {\msr G}$.

Let
$$\bar\Dlt=\sum_{i=1}^n\ptl_{x_i}\ptl_{y_i},\qquad
\bar\eta=\sum_{i=1}^nx_iy_i.\eqno(13.1.18)$$
 Recall
$$x^\al=x_1^{\al_1}\cdots x_n^{\al_n},\;\;y^\be=y_1^{\be_1}\cdots
y_n^{\be_n}\;\;\for\;\;\al=(\al_1,...,\al_n),\be=(\be_1,...,\be_n)\in\mbb{N}^{\:n}.\eqno(13.1.19)$$
For $\ell_1,\ell_2\in\mbb{N}$, we denote
$$\bar{\msr A}_{\ell_1,\ell_2}=\mbox{Span}\{x^\al
y^\be\mid\al,\be\in\mbb{N}^{\:n};\sum_{i=1}^n\al_i=\ell_1;\sum_{i=1}^n\be_i=\ell_2\}.\eqno(13.1.20)$$
Define
$$\bar{\msr H}_{\ell_1,\ell_2}=\{f\in\bar{\msr
A}_{\ell_1,\ell_2}\mid\bar\Dlt(f)=0\}.\eqno(13.1.21)$$ For
$i\in\ol{1,n}$, we define $\bar\ves_i\in\bar H^\ast$ by
$$\bar\ves_i(E_{j,j})=\dlt_{i,j}\qquad\for\;\;j\in\ol{1,n}.\eqno(13.1.22)$$
According to Theorem 6.2.4, we have:\psp

{\bf Lemma 13.1.1}. {\it Suppose $n>1$. For any $\ell_1,\ell_2$,
$\bar{\msr H}_{\ell_1,\ell_2}$ is a finite-dimensional irreducible
$\bar {\msr G}$-module with  highest-weight vector
$x_1^{\ell_1}y_n^{\ell_2}$ of weight $\ell_1\ves_1+\ell_2\ves_n$.
Moreover,
$$\bar{\msr A}=\bigoplus_{\ell_1,\ell_2,\ell_3=0}^\infty\bar
\eta^{\ell_1}\bar{\msr H}_{\ell_2,\ell_3}\eqno(13.1.23)$$ is a
decomposition of irreducible $\bar {\msr G}$-submodules.}\psp

When $n=1$, we have
$$\bar{\msr H}_{\ell,0}=\mbb{F}x_1^\ell,\qquad\bar{\msr
H}_{0,\ell}=\mbb{F}y_1^\ell\qquad\mbox{for}\;\;\ell\in\mbb{N}.\eqno(13.1.24)$$
Moreover,
$$\bar{\msr A}=\bigoplus_{\ell_1,\ell_2=0}^\infty(\bar\eta^{\ell_1}\bar{\msr
H}_{\ell_2,0}\oplus\bar\eta^{\ell_1}\bar{\msr
H}_{0,\ell_2+1}).\eqno(13.1.25)$$ \psp

Denote
$$\Theta_1=\sum_{i=1}^m\mbb{F}\sta_i,\qquad
\Theta_2=\sum_{i=1}^m\mbb{F}\vt_i.\eqno(13.1.26)$$ For
$\ell_1,\ell_2\in\ol{1,m}$, we define
$$\check{\msr
A}_{\ell_1,\ell_2}=\Theta_1^{\ell_1}\Theta_2^{\ell_2}.\eqno(13.1.27)$$
Then $\check{\msr A}_{\ell_1,\ell_2}$ is a finite-dimensional
$\check{\msr G}$-module and
$$\check{\msr
A}=\bigoplus_{\ell_1,\ell_2=0}^m\check{\msr
A}_{\ell_1,\ell_2}.\eqno(13.1.28)$$

Write
$$\check\Dlt=\sum_{r=1}^m\ptl_{\sta_r}\ptl_{\vt_r},\qquad
\check\eta=\sum_{r=1}^m\sta_r\vt_r.\eqno(13.1.29)$$ Recall the
notions
$$\vec\sta_r=\sta_1\cdots\sta_r,\qquad
\vec\vt_r=\vt_m\cdots\vt_r\qquad\for\;\;r\in\ol{1,m}\eqno(13.1.30)$$
and
$$\vec\sta_0=1=\vec\vt_{m+1}.\eqno(13.1.31)$$
For $\ell_1,\ell_2\in\ol{0,m}$, we define
$$\check{\msr H}_{\ell_1,\ell_2}=\{f\in\check{\msr
A}_{\ell_1,\ell_2}\mid\check\Dlt(f)=0\}.\eqno(13.1.32)$$
 For
$r\in\ol{1,m}$, we define $\ves_r'\in\check H^\ast$ by
$$\ves_r'(E_{n+s,n+s})=\dlt_{r,s}\qquad\for\;\;s\in\ol{1,m}.\eqno(13.1.33)$$
Moreover, we treat $\ves'_0=\ves'_{m+1}=0$. By Theorem 6.2.5. we
have:\psp

{\bf Lemma 13.1.2}. {\it For $0\leq r<s\leq m+1$, $\check{\msr
H}_{r,m+1-s}$ is a finite-dimensional irreducible $\check{\msr
G}$-module with the highest-weight vector $\vec\sta_r\vec\vt_s$ of
weight $\sum_{p=0}^r\ves_p'-\sum_{q=s}^{m+1}\ves'_q$. Moreover,
$$\check{\msr A}=\bigoplus_{0\leq r<s\leq
m+1}\bigoplus_{\ell=0}^{s-r-1}\check\eta^\ell \check{\msr
H}_{r,m+1-s}.\eqno(13.1.34)$$}\pse

According to (6.2.38) and (6.2.57), we have:
$$\bar\Dlt\bar\eta=\bar\eta\bar\Dlt+n+\sum_{i=1}^n(x_i\ptl_{x_i}
+y_i\ptl_{y_i}),\;\;\check\Dlt\check\eta=\check\eta\check\Dlt-m+\sum_{r=1}^m(\sta_r\ptl_{\sta_r}
+\vt_r\ptl_{\vt_r}).\eqno(13.1.35)$$ Set
$$\Dlt=\bar\Dlt+\check\Dlt=\sum_{i=1}^n\ptl_{x_i}\ptl_{y_i}+\sum_{r=1}^m\ptl_{\sta_r}\ptl_{\vt_r},\;\;
\eta=\bar\eta+\check\eta=\sum_{i=1}^nx_iy_i+\sum_{r=1}^m\sta_r\vt_r.\eqno(13.1.36)$$
For $\ell_1,\ell_2,\ell_3\in\mbb{N}$, $0\leq r<s\leq m+1$ and
$\ell\in\ol{0,s-r-1}$, $f\in\bar{\msr H}_{\ell_1,\ell_2}$ and
$g\in\check{\msr H}_{r,m+1-s}$, we have
$$\Dlt(\bar\eta^{\ell_3}\check\eta^\ell
fg)=\ell_3(n+\ell_1+\ell_2+\ell_3-1)
\bar\eta^{\ell_3-1}\check\eta^\ell
fg+\ell(\ell+r-s)\bar\eta^{\ell_3}\check\eta^{\ell-1}
fg.\eqno(13.1.37)$$ Suppose $r+1<s$ and $\ell\in\ol{1,s-r-1}$. If
$$\Dlt(\sum_{p=0}^\ell a_p\bar\eta^p\check\eta^{\ell-p}
fg)=0,\eqno(13.1.38)$$ then
$$(\ell-p)(p+s-r-\ell)a_p=(p+1)(n+\ell_1+\ell_2+p)a_{p+1}
\qquad\for\;\;p\in\ol{0,\ell-1}.\eqno(13.1.39)$$ Thus we can take
$a_0=(\ell+1)![\prod_{\iota_2=1}^{\ell+1}
(\iota_2+n+\ell_1+\ell_2-1)]$ and
$$a_{p+1}=[\prod_{\iota_1=0}^p(\ell-\iota_1)(\iota_1+s-r-\ell)][\prod_{\iota_2=p+2}^{\ell+1}
\iota_2(\iota_2+n+\ell_1+\ell_2-1)]\eqno(13.1.40)$$ for
$p\in\ol{0,\ell-1}$. Denote
\begin{eqnarray*} & &\Im(\ell_1,\ell_2;r,s,\ell)\\ &=&(\ell+1)![\prod_{\iota_2=1}^{\ell+1}
(\iota_2+n+\ell_1+\ell_2-1)]\check\eta^\ell+\sum_{p=0}^{\ell-1}
[\prod_{\iota_1=0}^p(\ell-\iota_1)(\iota_1+s-r-\ell)]\\
& &\times[\prod_{\iota_2=p+2}^{\ell+1}
\iota_2(\iota_2+n+\ell_1+\ell_2-1)]\bar\eta^{p+1}\check\eta^{\ell-p-1}.\hspace{5.6cm}(13.1.41)\end{eqnarray*}
 For
convenience, we treat
$$\Im(\ell_1,\ell_2;r,s,0)=n+\ell_1+\ell_2.\eqno(13.1.42)$$

 Write
$${\msr G}=\bar{\msr G}+\check{\msr G}\cong gl(n,\mbb{F})\oplus
gl(m,\mbb{F}). \eqno(13.1.43)$$
 Then
$${\msr A}=\bigoplus_{\ell_1,\ell_2,\ell_3=0}^\infty\;\bigoplus_{0\leq r<s\leq
m+1}\bigoplus_{\ell=0}^{s-r-1}(\bar \eta^{\ell_1}\bar{\msr
H}_{\ell_2,\ell_3})(\check\eta^\ell \check{\msr
H}_{r,m+1-s})\eqno(13.1.44)$$ is a direct sum of irreducible ${\msr
G}$-submodules. Moreover,
$${\msr H}=\{f\in{\msr A}\mid\Dlt(f)=0\}\eqno(13.1.45)$$
it is a ${\msr G}$-submodule. According to (13.1.37)-(13.1.42),
$${\msr H}= \bigoplus_{\ell_2,\ell_3=0}^\infty\;\bigoplus_{0\leq
r<s\leq m+1}\bigoplus_{\ell=0}^{s-r-1}
\Im(\ell_2,\ell_3;r,s,\ell)(\bar{\msr H}_{\ell_2,\ell_3}\check{\msr
H}_{r,m+1-s})\eqno(13.1.46)$$ is a direct sum of irreducible ${\msr
G}$-submodules.

For $\ell,\ell'\in\mbb{N}$, we let
$${\msr
A}_{\ell,\ell'}=\sum_{\ell_1,\ell_2\in\mbb{N},\;\ell_3,\ell_4\in\ol{0,m};\;\ell_1+\ell_3=\ell,\;
\ell_2+\ell_4=\ell'}\bar{\msr A}_{\ell_1,\ell_2}\check{\msr
A}_{\ell_3,\ell_4}\eqno(13.1.47)$$ and
$${\msr H}_{\ell,\ell'}={\msr A}_{\ell,\ell'}\bigcap{\msr
H}.\eqno(13.1.48)$$ It is straightforward to verify
$$E_{i,j}\Dlt=\Dlt E_{i,j},\qquad E_{i,j}\eta=\eta
E_{i,j}\qquad\for\;\;i,j\in\ol{1,m+n}\eqno(13.1.49)$$ by (13.1.11),
(13.1.12) and (13.1.36). Thus ${\msr A}_{\ell,\ell'}$ and ${\msr
H}_{\ell,\ell'}$ are $gl(n|m)$-submodules. Moreover, $${\msr
H}_{\ell,\ell'}=
\bigoplus_{\ell_1+r+\ell_3=\ell,\;\ell_2+\ell_3+m+1-s=\ell'}
\Im(\ell_1,\ell_2;r,s,\ell_3)\bar{\msr H}_{\ell_1,\ell_2}\check{\msr
H}_{r,m+1-s}\eqno(13.1.50)$$ is a direct sum of irreducible ${\msr
G}$-submodules. Take the Cartan subalgebra $H=\bar H+\check H$ of
${\msr G}$ (cf. (13.1.16)) and the subspace
 ${\msr G}_+=\bar{\msr G}_++\check{\msr G}_+$ (cf. (13.1.17)) spanned by positive
 root vectors in ${\msr G}$.  A ${\msr G}$-{\it singular
vector} $v$ is a nonzero weight vector of ${\msr G}$ such that
${\msr G}_+(v)=0$. We count singular vector up to a nonzero scalar
multiple. Hence we have:\psp

{\bf Lemma 13.1.3}. {\it The set
\begin{eqnarray*}\qquad&&\{\Im(\ell_1,\ell_2;r,s,\ell_3)(x_1^{\ell_1}y_n^{\ell_2}\vec\sta_r\vec\vt_s)
\mid\ell_1,\ell_2\in\mbb{N};0\leq r<s\leq m+1;\\
& &\ell_3\in\ol{0,s-r-1};\ell_1+r+\ell_3=\ell,\;
\ell_2+\ell_3+m+1-s=\ell'\}\hspace{2.4cm}(13.1.51)\end{eqnarray*} is
the set of  all the ${\msr G}$-singular vectors in ${\msr
H}_{\ell,\ell'}$, where $\ell_1\ell_2=\ell_1'\ell_2'=0$ if $n=1$.}

\subsection{Singular Vectors of $gl(n|m)$}

Take $H=\bar H+\check H$ as a Cartan subalgebra of the Lie
superalgebra $gl(n|m)$ and
$$gl(n|m)_+={\msr
G}_++\sum_{r=1}^n\sum_{s=1}^m\mbb{F}E_{r,n+s}\eqno(13.1.52)$$ as the
subalgebra generated by positive root vectors.  A $gl(n|m)$-{\it
singular vector} $v$ is a nonzero weight vector of $gl(n|m)$ such
that $gl(n|m)_+(v)=0$. We count singular vector up to a nonzero
scalar multiple. Assume that
$x_1^{\ell_1'}y_n^{\ell_2'}\vec\sta_{r'}\vec\vt_{s'}$ is a
$gl(n|m)$-singular vector, where $\ell_1'\ell_2'=0$ when $n=1$. If
$r'\neq 0$, then
\begin{eqnarray*}\qquad & &E_{1,n+r'}(x_1^{\ell_1'}y_n^{\ell_2'}\vec\sta_{r'}\vec\vt_{s'})
\\ &=&
(x_1\ptl_{\sta_{r'}}-\vt_{r'}\ptl_{y_1})(x_1^{\ell_1'}y_n^{\ell_2'}\vec\sta_{r'}\vec\vt_{s'})
\\&=&(-1)^{r'-1}x_1^{\ell_1'+1}y_n^{\ell_2'}\vec\sta_{r'-1}\vec\vt_{s'}-\dlt_{1,n}\ell_2'
x_1^{\ell_1'}y_n^{\ell_2'-1}\vt_{r'}\vec\sta_{r'}\vec\vt_{s'} \neq
0\hspace{3.2cm}(13.1.53)\end{eqnarray*}
 by the second equation in
(13.1.11), which contradicts the definition of singular vector. So
$r'=0$. Suppose $\ell_2'>0$ and $s'>1$. Again the second equation in
(13.1.11) gives
\begin{eqnarray*}\qquad\quad
E_{n,n+s'-1}(x_1^{\ell_1'}y_n^{\ell_2'}\vec\vt_{s'})&=&
(x_n\ptl_{\sta_{s'-1}}-\vt_{s'-1}\ptl_{y_n})(x_1^{\ell_1'}y_n^{\ell_2'}\vec\vt_{s'})
\\ &=&-\ell_2'
x_1^{\ell_1'}y_n^{\ell_2'-1}\vec\sta_{r'-1}\vec\vt_{s'-1}\neq
0,\hspace{4cm}(13.1.54)\end{eqnarray*} which is absurd. Thus
 $x_1^{\ell_1'}y_n^{\ell_2'}\vec\sta_{r'}\vec\vt_{s'}$ is a
$gl(n|m)$-singular vector if and only if $r'=0$ and
$\ell_2'(s'-1)=0$.

For $\ell_1,\ell_2\in\mbb{N}$, $1\leq r<s-1\leq m$ and
$\ell\in\ol{1,s-r-1}$,
\begin{eqnarray*}& &E_{1,n+s-1}[\Im(\ell_1,\ell_2;r,s,\ell)]\\ &=&(x_1\ptl_{\sta_{s-1}}-\vt_{s-1}\ptl_{y_1})[\Im(\ell_1,\ell_2;r,s,\ell)]
\\ &=&x_1\vt_{s-1}\{
(\ell+1)![\prod_{\iota_2=1}^{\ell+1}
(\iota_2+n+\ell_1+\ell_2-1)]\ell \check\eta^{\ell-1}
\\ & &+\sum_{p=0}^{\ell-2}
[\prod_{\iota_1=0}^p(\ell-\iota_1)(\iota_1+s-r-\ell)]\\
& &\times(\ell-p-1)[\prod_{\iota_2=p+2}^{\ell+1}
\iota_2(\iota_2+n+\ell_1+\ell_2-1)]\bar\eta^{p+1}\check\eta^{\ell-p-2}
\\ & &-\sum_{p=0}^{\ell-1}
(p+1)[\prod_{\iota_1=0}^p(\ell-\iota_1)(\iota_1+s-r-\ell)]
\\ & &\times [\prod_{\iota_2=p+2}^{\ell+1}
\iota_2(\iota_2+n+\ell_1+\ell_2-1)]\bar\eta^p\check\eta^{\ell-p-1}\}\\
&=&\ell(\ell+1)(n+\ell+\ell_1+\ell_2+r-s)\Im(\ell_1+1,\ell_2;r,s-1,\ell-1)x_1\vt_{s-1}.\hspace{1.5cm}(13.1.55)
\end{eqnarray*}
Moreover, (13.1.41) yields
$$\Im(\ell_1,\ell_2;r,s,\ell)=(\ell+1)![\prod_{\iota_1=0}^\ell(\iota_1+s-r-\ell)]\eta^\ell\;\;\mbox{if}\;\;
n+\ell+\ell_1+\ell_2+r-s=0.\eqno(13.1.56)$$

Suppose that
$\Im(\ell_1,\ell_2;r,s,\ell_3)(x_1^{\ell_1}y_n^{\ell_2}\vec\sta_r\vec\vt_s)$
is a $gl(n|m)$-singular vector, then
$$n+\ell_3+\ell_1+\ell_2+r-s=0\eqno(13.1.57)$$
and
$$\Im(\ell_1,\ell_2;r,s,\ell_3)(x_1^{\ell_1}y_n^{\ell_2}\vec\sta_r\vec\vt_s)
=c\eta^{\ell_3}
x_1^{\ell_1}y_n^{\ell_2}\vec\sta_r\vec\vt_s,\eqno(13.1.58)$$ where
$c=(\ell_3+1)![\prod_{\iota_1=0}^\ell(\iota_1+s-r-\ell_3)]$ by
(13.1.54). Since
$$E_{i,n+r}(c\eta^{\ell_3}
x_1^{\ell_1}y_n^{\ell_2}\vec\sta_r\vec\vt_s)=c\eta^{\ell_3}
E_{i,n+r}(x_1^{\ell_1}y_n^{\ell_2}\vec\sta_r\vec\vt_s)\eqno(13.1.59)$$
by (13.1.11), the arguments in (13.1.53) and (13.1.54) show
$r=\ell_2(s-1)=0$. On the other hand, $1\leq 1+r<s$. So $\ell_2=0$.
According to (13.1.57),
$$\ell_3=s-\ell_1-n.\eqno(13.1.60)$$
Thus we only get the singular vector
$$\eta^{s-\ell_1-n}
x_1^{\ell_1}\vec\vt_s\in{\msr
H}_{s-n,m+1-n-\ell_1}\;\mbox{with}\;\;s>\ell_1+n.\eqno(13.1.61)$$

Note $n\leq m$ by $\ell_3>0$. Moreover, $s\leq m+1$ implies
$$s-n,m+1-n-\ell_1\leq
m+1-n.\eqno(13.1.62)$$ Furthermore,
$$(s-n)+(m+1-n-\ell_1)=m+1-n+\ell_3>m+1-n.\eqno(13.1.63)$$
This shows that \begin{eqnarray*}\qquad& &{\msr
H}_{\ell,\ell'}\;\mbox{has a unique $gl(n|m)$-singular vector if}\\
& &\ell>m+1-n\;\mbox{or}\;\ell'>m+1-n\;\mbox{or}\;\ell+\ell'\leq
m+1-n.\hspace{3cm}(13.1.64)\end{eqnarray*}

Suppose $n\leq m$ and $\ell,\ell'\in\ol{0,m+1-n}$ such that
$\ell+\ell'>m+1-n$. We take
$$s=n+\ell,\qquad
\ell_1=m+1-n-\ell',\qquad\ell_3=\ell+\ell'+n-m-1.\eqno(13.1.65)$$
Then
$$\eta^{\ell+\ell'+n-m-1}(x_1^{m+1-n-\ell'}\vec\vt_{n+\ell})\in
{\msr H}_{\ell,\ell'}.\eqno(13.1.66)$$ Hence $${\msr
H}_{\ell,\ell'}\;\mbox{has exactly two $gl(n|m)$-singular vectors
if}\; \ell+\ell'> m+1-n.\eqno(13.1.67)$$ In summary, we have: \psp

{\bf Lemma 13.1.4}. {\it Let $\ell,\ell'\in\mbb{N}$. If
$\ell>m+1-n\;\mbox{or}\;\ell'>m+1-n\;\mbox{or}\;\ell+\ell'\leq
m+1-n,$ the $gl(n|m)$-module ${\msr H}_{\ell,\ell'}$ has a unique
$gl(n|m)$-singular vector. When $\ell,\ell'\leq m+1-n$ and
$\ell+\ell'>m+1-n$, the $gl(n|m)$-module ${\msr H}_{\ell,\ell'}$ has
exactly two $gl(n|m)$-singular vectors.}\psp

\subsection{Structure of the Representation for $gl(n|m)$}

Fix $\ell,\ell'\in\mbb{N}$. Let
$$v_{\ell,\ell'}=x_1^\ell y_n^{\ell_1}\vec\vt_s,\qquad
\ell_1+m+1-s=\ell',\;\;\ell_1(s-1)=0,\;\;\dlt_{n,1}\ell\ell_1=0.\eqno(13.1.68)$$

{\bf Lemma 13.1.5}. {\it The $gl(n|m)$-module ${\msr
H}_{\ell,\ell'}$ is generated by $v_{\ell,\ell'}$.}

{\it Proof.} Let $M$ be the $gl(n|m)$-submodule of ${\msr
H}_{\ell,\ell'}$ generated by $v_{\ell,\ell'}$.  First consider
$n>1$ or $\ell=0$. For any $s'\in\ol{s+1,m+1}$, we have
$$E_{n+s'-1,n}E_{n+s'-2,n}\cdots E_{n+s,n}(v_{\ell,\ell'})=
x_1^\ell y_n^{\ell_1+s'-s}\vec\vt_{s'}\in M\eqno(13.1.69)$$ by
(13.1.12). In other words, \begin{eqnarray*}\qquad\quad & &x_1^\ell
y_n^{\ell_2}\vec\vt_{s'}\in M\;\;\mbox{for
any}\;\ell_2\in\mbb{N}\;\mbox{and}\;s'\in\ol{1,m+1}\\ & &\mbox{such
that}\;\;\ell_2+m+1-s'=\ell'\;\mbox{and}\;\dlt_{n,1}\ell\ell_2=0.
\hspace{4.1cm}(13.1.70)\end{eqnarray*} If $\ell >0$, then for any
$1\leq r\leq \mbox{min}\{\ell,s'-1\}$, we have
$$E_{n+1,1}E_{n+2,1}\cdots E_{n+r,1}(x_1^\ell
y_n^{\ell_2}\vec\vt_{s'})=[\prod_{p=0}^{r-1}(\ell-p)]x_1^{\ell-r}y_n^{\ell_2}\vec\sta_r\vec\vt_{s'}\in
M\eqno(13.1.71)$$ by (13.1.12) again. Thus we have showed that
$$x_1^{\ell_3}y_n^{\ell_2}\vec\sta_r\vec\vt_{s'}\in
M\;\;\mbox{whenever}\;\;r+\ell_3=\ell,\;\ell_2+m+1-s'=\ell',\;\;\dlt_{n,1}\ell_2\ell_3=0\eqno(13.1.72)$$
for $\ell_2,\ell_3\in\mbb{N}$ and $0\leq r<s'\leq m+1$. Recall the
Lie algebra ${\msr G}$ defined in (13.1.42). As ${\msr G}$-modules,
$$\sum_{r+\ell_3=\ell,\;\ell_2+m+1-s'=\ell'}\check{\msr
H}_{\ell_3,\ell_2}\bar{\msr H}_{r,m+1-s'}\subset M.\eqno(13.1.73)$$

For $i\in\mbb{N}+1$, we define
$${\msr H}^{(i)}_{\ell,\ell'}=\{f\in{\msr H}_{\ell,\ell'}\mid
\bar\Dlt^i(f)=0.\}\eqno(13.1.74)$$ First
$${\msr H}^{(1)}_{\ell,\ell'}=\sum_{r+\ell_3=\ell,\;\ell_2+m+1-s'=\ell'}\check{\msr
H}_{\ell_3,\ell_2}\bar{\msr H}_{r,m+1-s'}\subset M.\eqno(13.1.75)$$
Denote
$$k_{\ell,\ell'}=\min\{\ell,\ell'\}.\eqno(13.1.76)$$
We have
$${\msr H}_{\ell,\ell'}^{(k_{\ell,\ell'}+1)}={\msr
H}_{\ell,\ell'}\eqno(13.1.77)$$ by (13.1.48). Let
$\ell_1,\ell_2,\in\mbb{N}$, $0\leq r'+1<s'\leq m+1$ and
$\ell_3\in\ol{1,s'-r'-1}$ such that $\ell_1+\ell_3+r'=\ell$ and
$\ell_2+\ell_3+m+1-s'=\ell'$. Then
$$\Im(\ell_1,\ell_2;r',s',\ell_3)(x_1^{\ell_1}y_n^{\ell_2}\vec\sta_{r'}\vec\vt_{s'})
\in {\msr H}_{\ell,\ell'}^{(\ell_3+1)}\eqno(13.1.78)$$ and
$$\bar\Dlt^{\ell_3}\Im(\ell_1,\ell_2;r',s',\ell_3)
(x_1^{\ell_1}y_n^{\ell_2}\vec\sta_{r'}\vec\vt_{s'})=
c(\ell_1,\ell_2;r',s',\ell_3)x_1^{\ell_1}y_n^{\ell_2}\vec\sta_{r'}\vec\vt_{s'}\eqno(13.1.79)$$
with
$$c(\ell_1,\ell_2;r',s',\ell_3)=
\ell_3(n+\ell_1+\ell_3+\ell_3-1)(n+\ell_1+\ell_2+\ell_3)(\ell_3+1)!
[\prod_{p=1}^{\ell_3}(s'-r'-p)]\eqno(13.1.80)$$ by (13.1.35) and
(13.1.41). On the other hand,
$$x_1^{\ell_1+\ell_3}y_n^{\ell_2}\vec\sta_{r'}\vec\vt_{s'-\ell_3}\in
M\;\;\mbox{with}\;\;\dlt_{n,1}\ell_2=0\eqno(13.1.81)$$ by (13.1.73)
and
\begin{eqnarray*}\qquad & &M\ni f=E_{n+s'-1,1}E_{n+s'-2,1}\cdots
E_{n+s'-\ell_3,1}(x_1^{\ell_1+\ell_3}y_n^{\ell_2}\vec\sta_{r'}\vec\vt_{s'-\ell_3})
\\ &=&(-1)^{r'\ell_3}x_1^{\ell_1+\ell_3}y_1^{\ell_3}
y_n^{\ell_2}\vec\sta_{r'}\vec\vt_{s'}+y_n^{\ell_2}\sum_{i=0}^{\ell_3-1}\zeta_iy_1^i
\hspace{5.8cm}(13.1.83)\end{eqnarray*} with
$\zeta_0,...,\zeta_{\ell_3-1}\in\mbb{F}[x_1]\check{\msr A}$ (cf.
(13.1.15)). Moreover,
\begin{eqnarray*}\qquad\qquad & &\bar\Dlt^{\ell_3}[(\ell_3!)^2{\ell_1+\ell_3\choose\ell_3}\Im(\ell_1,\ell_2;r',s',\ell_3)
(x_1^{\ell_1}y_n^{\ell_2}\vec\sta_{r'}\vec\vt_{s'})\\
& &-(-1)^{r'\ell_3}
c(\ell_1,\ell_2;r',s',\ell_3)f]=0.\hspace{6.3cm}(13.1.84)\end{eqnarray*}
Hence
$$\Im(\ell_1,\ell_2;r',s',\ell_3)(\bar{\msr
H}_{\ell_1,\ell_2}\check{\msr H}_{r',m+1-s'})\subset {\msr
H}^{(\ell_3)}_{\ell,\ell'}+M.\eqno(13.1.85)$$ By (13.1.50) and
induction on $i$, we have
$${\msr H}^{(i)}_{\ell,\ell'}\subset M\qquad\mbox{for
any}\;i\in\mbb{N}+1.\eqno(13.1.86)$$ According to (13.1.77), ${\msr
H}_{\ell,\ell'}=M$. So ${\msr H}_{\ell,\ell'}$ is generated by
$v_{\ell,\ell'}.\qquad\Box$\psp

{\bf Theorem 13.1.6}. {\it Let $\ell,\ell'\in\mbb{N}$. The space
${\msr H}_{\ell,\ell'}$ is an irreducible $gl(n|m)$-module  if and
only if
$\ell>m+1-n\;\mbox{or}\;\ell'>m+1-n\;\mbox{or}\;\ell+\ell'\leq
m+1-n.$  When $|\ell-\ell'|>m+1-n$ or $\ell+\ell'\leq m+1-n$, ${\msr
A}_{\ell,\ell'}=\bigoplus_{i=0}^{k_{\ell,\ell'}}\eta^i{\msr
H}_{\ell-i,\ell'-i}$ is a decomposition of irreducible
$gl(n|m)$-submodules}.

 {\it Proof.} According to (13.1.67), a necessary condition
for ${\msr H}_{\ell,\ell'}$ to be an irreducible $gl(n|m)$-module is
$\ell>m+1-n\;\mbox{or}\;\ell'>m+1-n\;\mbox{or}\;\ell+\ell'\leq
m+1-n.$ To prove the sufficiency, we suppose that
$\ell>m+1-n\;\mbox{or}\;\ell'>m+1-n\;\mbox{or}\;\ell+\ell'\leq
m+1-n.$ Let $V$ be a nonzero $gl(n|m)$-submodule of ${\msr
H}_{\ell,\ell'}$. According to Lemma 13.1.4, ${\msr H}_{\ell,\ell'}$
has a unique singular vector $v_{\ell,\ell'}$ (cf. (13.1.68)). Since
$V$ is finite-dimensional, it contains a singular vector. So
$v_{\ell,\ell'}\in V$. Lemma 13.1.5 says $V={\msr H}_{\ell,\ell'}$.
Hence ${\msr H}_{\ell,\ell'}$ is an irreducible $gl(n|m)$-module.

Let $\ell_1,\ell_2\in\mbb{N}$ and $0\leq r<s\leq m+1$ such that
$n+\ell_1+\ell_2+r-s\geq 0$ and $\dlt_{n,1}\ell_1\ell_2=0$.
 For any
$\ell_4\in\mbb{N}+1$ and $\ell_3\in\ol{0,s-r-1}$, $$\Dlt^{\ell_4+1}
(\eta^{\ell_4}\Im(\ell_1,\ell_2;r,s,\ell_3)x_1^{\ell_1}y_n^{\ell_2}\vec\sta_r\vec\vt_s)=0\eqno(13.1.87)$$
and
\begin{eqnarray*} & &\Dlt^{\ell_4}
(\eta^{\ell_4}\Im(\ell_1,\ell_2;r,s,\ell_3)x_1^{\ell_1}y_n^{\ell_2}\vec\sta_r\vec\vt_s)
\\ &=&\ell_4![\prod_{i=1}^{\ell_4}(n+i+\ell_1+2\ell_3+r-s)]
\Im(\ell_1,\ell_2;r,s,\ell_3)(x_1^{\ell_1}y_n^{\ell_2}\vec\sta_r\vec\vt_s)\neq
0\hspace{2cm}(13.1.88)\end{eqnarray*} by (13.1.35). Thus the set
$$\{\eta^{\ell_4}\Im(\ell_1,\ell_2;r,s,\ell_3)x_1^{\ell_1}y_n^{\ell_2}\vec\sta_r\vec\vt_s\mid\ell_4\in\mbb{N},
\;\ell_3\in\ol{0,s-r-1}\}\eqno(13.1.89)$$ is linearly independent.

Note that
$$\check\eta^{s-r}x_1^{\ell_1}y_n^{\ell_2}\vec\sta_r\vec\vt_s=0\eqno(13.1.90)$$
by (13.1.29) and (13.1.30). So for any $k\in\mbb{N}$,
\begin{eqnarray*}& &\mbox{Span}\{\eta^{\ell_4}\Im(\ell_1,\ell_2;r,s,\ell_3)x_1^{\ell_1}y_n^{\ell_2}
\vec\sta_r\vec\vt_s\mid\ell_4\in\mbb{N},\;\ell_3\in\ol{0,s-r-1};\ell_3+\ell_4=k\}
\\ & \subset& \mbox{Span}\{
\bar\eta^{\ell_5}\check\eta^{\ell_6}x_1^{\ell_1}y_n^{\ell_2}
\vec\sta_r\vec\vt_s\mid\ell_5\in\mbb{N};\ell_6\in\ol{0,s-r-1};\ell_5+\ell_6=k\}.
\hspace{2.4cm}(13.1.91)\end{eqnarray*} But the linear independency
of (13.1.89) implies that the above subspaces have the same
dimension. Thus
\begin{eqnarray*}& &\mbox{Span}\{\eta^{\ell_4}\Im(\ell_1,\ell_2;r,s,\ell_3)x_1^{\ell_1}y_n^{\ell_2}
\vec\sta_r\vec\vt_s\mid\ell_4\in\mbb{N},\;\ell_3\in\ol{0,s-r-1}\}
\\ & =& \mbox{Span}\{
\bar\eta^{\ell_5}\check\eta^{\ell_6}x_1^{\ell_1}y_n^{\ell_2}
\vec\sta_r\vec\vt_s\mid\ell_5\in\mbb{N};\ell_6\in\ol{0,s-r-1}\}.
\hspace{4.5cm}(13.1.92)\end{eqnarray*} Therefore, as ${\msr
G}$-modules,
$$\bigoplus_{\ell_3=0}^{s-r-1}\sum_{\ell_4=0}^\infty
\eta^{\ell_4}\Im(\ell_1,\ell_2;r,s,\ell_3)\bar{\msr
H}_{\ell_1,\ell_2}\check{\msr H}_{r,m+1-s}
=\bigoplus_{\ell_5=0}^{s-r-1}\sum_{\ell_6=0}^\infty
\bar\eta^{\ell_5}\check\eta^{\ell_6}\bar{\msr
H}_{\ell_1,\ell_2}\check{\msr H}_{r,m+1-s}.\eqno(13.1.93)$$

Assume that $|\ell-\ell'|>m+1-n$ or $\ell+\ell'\leq m+1-n$.
According to (13.1.44) and (13.1.92), the ${\msr G}$-module
\begin{eqnarray*}\qquad &&{\msr A}_{\ell,\ell'}=
\mbox{Span}\{\bar\eta^{\ell_5}\check\eta^{\ell_6}\bar{\msr
H}_{\ell_1,\ell_2}\check{\msr
H}_{r,m+1-s}\mid\ell_1,\ell_2,\ell_6\in\mbb{N};\\ & &0\leq r<s\leq
m+1; \ell_5\in\ol{0,s-r-1}; \dlt_{n,1}\ell_1\ell_2=0;\\
& &\ell_1+\ell_5+\ell_6+r=\ell,\;\ell_2+\ell_5+\ell_6+m+1-s=\ell'\}
\\&=&\mbox{Span}\{\eta^{\ell_4}\Im(\ell_1,\ell_2;r,s,\ell_3)\bar{\msr
H}_{\ell_1,\ell_2}\check{\msr
H}_{r,m+1-s}\mid\ell_1,\ell_2,\ell_6\in\mbb{N};\\
& &\dlt_{n,1}\ell_1\ell_2=0; 0\leq r<s\leq m+1;
\ell_5\in\ol{0,s-r-1};\\
& &\ell_1+\ell_3+\ell_4+r=\ell,\;\ell_2+\ell_3+\ell_4+m+1-s=\ell'\}.
\hspace{3.3cm}(13.1.94)\end{eqnarray*}

According to  (13.1.50), (13.1.93) and (13.1.94),
$${\msr A}_{\ell,\ell'}=\bigoplus_{i=0}^{k_{\ell,\ell'}}\eta^i{\msr
H}_{\ell-i,\ell'-i}\eqno(13.1.95)$$ (cf. (13.1.76)). Let
$i\in\ol{0,k_{\ell,\ell'}}$. If $|\ell-\ell'|>m+1-n$, then
$$\ell-i\geq |\ell-\ell'|>m+1-n\;\;\mbox{or}\;\;\ell'-i\geq
|\ell-\ell'|>m+1-n.\eqno(13.1.96)$$ When $\ell+\ell'\leq m+1-n$,
$$(\ell-i)+(\ell'-i)=\ell+\ell'-2i
\leq m+1-n.\eqno(13.1.97)$$ Thus all ${\msr H}_{\ell-i,\ell'-i}$ are
irreducible $gl(n|m)$-submodules. Hence (13.1.94) is a direct sum of
irreducible $gl(n|m)$-submodules. $\qquad\Box$ \psp

Denote $\G_0=\emptyset$ and
$$\G_\ell=\{\vec j=(j_1,j_2,..,j_\ell)\mid 1\leq j_1<j_2<\cdots<
j_\ell\leq m\}\qquad\for\;\;\ell\in\ol{1,m}.\eqno(13.1.98)$$
Moreover, we set
$$\sta_\emptyset=\vt_\emptyset=1,\;\;\sta_{\vec
j}=\sta_{j_1}\sta_{j_2}\cdots\sta_{j_\ell},\;\;\vt_{\vec
j}=\vt_{j_1}\vt_{j_2}\cdots\vt_{j_\ell}.\eqno(13.1.99)$$ Then the
set
\begin{eqnarray*}& &\{\sum_{i=0}^\infty\frac{(-1)^ix_1^iy_1^i}{\prod_{r=1}^i(\al_1+i)(\be_1+i)}
(\sum_{s=2}^n\ptl_{x_2}\ptl_{y_2}+\sum_{r=1}^m\ptl_{\sta_r}\ptl_{\vt_r})^i(x^\al
y^\be\sta_{\vec j}\vt_{\vec k})\mid\al,\be\in\mbb{N}^n;\\ & &\vec
j\in\G_{\ell_1};\vec
k\in\G_{\ell_2};\al_1\be_1=0;\ell_1,\ell_2\in\ol{0,m};|\al|+\ell_1=\ell;|\be|+\ell_2=\ell'\}\hspace{1.9cm}(13.1.100)\end{eqnarray*}
forms a basis of ${\msr H}_{\ell,\ell'}$ by Lemma 6.1.1 with
$T_1=\ptl_{x_1}\ptl_{y_1},\;T_2=\Dlt-T_1$ and
$T_1^-=\int_{(x_1)}\int_{(y_1)}$.\psp

{\bf Remark 13.1.7}. If  $\ell,\ell'\leq m+1-n$ and
$\ell+\ell'>m+1-n$, then ${\msr H}_{\ell,\ell'}$ is an
indecomposable $gl(n|m)$-module. In fact, ${\msr
H}_{\ell,\ell'}\bigcap\eta{\msr A}_{\ell-1,\ell'-1}\neq\{0\}$. This
also shows that ${\msr A}_{\ell,\ell'}$ is not completely reducible
when $|\ell-\ell'|\leq m+1-n$ and $\ell+\ell'>m+1-n$. It can be
verified that the space $\mbb{F}\Dlt+\mbb{F}[\Dlt,\eta]+\mbb{F}\eta$
also forms an operator Lie algebra isomorphic to $sl(2,\mbb{F})$.
The above theorem establishes a supersymmetric
$(sl(2,\mbb{F}),gl(n|m))$ Howe duality on the homogeneous subspaces
${\msr A}_{\ell,\ell'}$ with $\ell,\ell'\in\mbb{N}$ such that
$|\ell-\ell'|>m+1-n$ or $\ell+\ell'\leq m+1-n$.

\section{Noncanonical oscillator Representations of $gl(n|m)$}

In this section, we study $\mbb Z^2$-graded oscillator
generalizations of the canonical oscillator representations of
$gl(n|m))$ presented in the last section.

 Fix the integers $1<n_1+1<n_2\leq n$.
  Recall the symmetry
$$[\ptl_{x_r},x_r]=1=[-x_r,\ptl_{x_r}],\;\;[\ptl_{y_s},y_s]=1=[-y_s,\ptl_{y_s}].\eqno(13.2.1)$$
Changing operators $\ptl_{x_r}\mapsto -x_r,\;
 x_r\mapsto
\ptl_{x_r}$  for $r\in\ol{1,n_1}$ and $\ptl_{y_s}\mapsto -y_s,\;
 y_s\mapsto\ptl_{y_s}$  for $s\in\ol{n_2+1,n}$ in (13.1.11) and
 (13.1.12), we get a new representation of $gl(n|m)$ on
  ${\msr A}$ determined by
$$E_{i,j}|_{\msr
A}=E_{i,j}^x-E_{j,i}^y,\qquad E_{n+r,n+s}|_{\msr
A}=\sta_r\ptl_{\sta_s}-\vt_s\ptl_{\vt_r}\eqno(13.2.2)$$ with
$$E_{i,j}^x|_{\msr A}=\left\{\begin{array}{ll}-x_j\ptl_{x_i}-\delta_{i,j}&\mbox{if}\;
i,j\in\ol{1,n_1};\\ \ptl_{x_i}\ptl_{x_j}&\mbox{if}\;i\in\ol{1,n_1},\;j\in\ol{n_1+1,n};\\
-x_ix_j &\mbox{if}\;i\in\ol{n_1+1,n},\;j\in\ol{1,n_1};\\
x_i\partial_{x_j}&\mbox{if}\;i,j\in\ol{n_1+1,n}
\end{array}\right.\eqno(13.2.3)$$
and
$$E_{i,j}^y|_{\msr A}=\left\{\begin{array}{ll}y_i\ptl_{y_j}&\mbox{if}\;
i,j\in\ol{1,n_2};\\ -y_iy_j&\mbox{if}\;i\in\ol{1,n_2},\;j\in\ol{n_2+1,n};\\
\ptl_{y_i}\ptl_{y_j} &\mbox{if}\;i\in\ol{n_2+1,n},\;j\in\ol{1,n_2};\\
-y_j\partial_{y_i}-\delta_{i,j}&\mbox{if}\;i,j\in\ol{n_2+1,n};
\end{array}\right.\eqno(13.2.4)$$
$$ E_{i,n+r}|_{\msr
A}=\left\{\begin{array}{ll}
\ptl_{x_i}\ptl_{\sta_r}-\vt_r\ptl_{y_i}&\mbox{if}\;
i\in\ol{1,n_1};\\ x_i\ptl_{\sta_r}-\vt_r\ptl_{y_i}&\mbox{if}\;
i\in\ol{n_1+1,n_2};\\ x_i\ptl_{\sta_r}+y_i\vt_r&\mbox{if}\;
i\in\ol{n_2+1,n};\end{array}\right.\eqno(13.2.5)$$
$$ E_{n+r,i}|_{\msr
A}=\left\{\begin{array}{ll}-x_i\sta_r+y_i\ptl_{\vt_r} &\mbox{if}\;
i\in\ol{1,n_1};\\ \sta_r\ptl_{x_i}+y_i\ptl_{\vt_r}&\mbox{if}\; i\in\ol{n_1+1,n_2};\\
\sta_r\ptl_{x_i}+\ptl_{y_i}\ptl_{\vt_r}&\mbox{if}\;
i\in\ol{n_2+1,n}\end{array}\right.\eqno(13.2.6)$$ for
$i,j\in\ol{1,n}$ and $r,s\in\ol{1,m}$.

The Laplace operator in (13.1.18) changes to
$$\bar\Dlt_{n_1,n_2}=-\sum_{i=1}^{n_1}x_i\ptl_{y_i}+\sum_{r=n_1+1}^{n_2}\ptl_{x_r}\ptl_{y_r}-\sum_{s=n_2+1}^n
y_s\ptl_{x_s}\eqno(13.2.7)$$ and its dual changes to
$$\bar\eta=\sum_{i=1}^{n_1}y_i\ptl_{x_i}+\sum_{r=n_1+1}^{n_2}x_ry_r+\sum_{s=n_2+1}^n
x_s\ptl_{y_s}.\eqno(13.2.8)$$ We take (13.1.29) and then
supersymmetric Laplace operator
$$\Dlt_{n_1,n_2}=\bar\Dlt_{n_1,n_2}+\check\Dlt=
-\sum_{i=1}^{n_1}x_i\ptl_{y_i}+\sum_{r=n_1+1}^{n_2}\ptl_{x_r}\ptl_{y_r}-\sum_{s=n_2+1}^n
y_s\ptl_{x_s}+\sum_{r=1}^m\ptl_{\sta_r}\ptl_{\vt_r}\eqno(13.2.9)$$
and its dual
$$\eta=\bar\eta+\check \eta=\sum_{i=1}^{n_1}y_i\ptl_{x_i}+\sum_{r=n_1+1}^{n_2}x_ry_r+\sum_{s=n_2+1}^n
x_s\ptl_{y_s}+\sum_{r=1}^m\sta_r\vt_r.\eqno(13.2.10)$$
 Then with respect to the
representation in (13.2.2)-(13.2.6), we have
$$E_{i,j}\Dlt_{n_1,n_2}=\Dlt_{n_1,n_2} E_{i,j},\qquad E_{i,j}\eta=\eta
E_{i,j}\qquad\for\;\;i,j\in\ol{1,m+n}.\eqno(13.2.11)$$ Moreover, we
take the settings in (13.1.13), (13.1.15)-(13.1.7) and  (13.1.19).

Denote \begin{eqnarray*}\qquad\qquad& & \bar{\msr A}_{\la
\ell_1,\ell_2\ra}=\mbox{Span}\{x^\al
y^\be\mid\al,\be\in\mbb{N}^n;\\
& &\sum_{r=n_1+1}^n\al_r-\sum_{i=1}^{n_1}\al_i=\ell_1;
\sum_{i=1}^{n_2}\be_i-\sum_{r=n_2+1}^n\be_r=\ell_2\}\hspace{3.6cm}(13.2.12)\end{eqnarray*}
for $\ell_1,\ell_2\in\mbb{Z}$. Then $\bar{\msr
A}_{\la\ell_1,\ell_2\ra}$ forms a $\bar{\msr G}$-submodule.
Moreover, for $\ell,\ell'\in\mbb{Z}$, we let
$${\msr
A}_{\la\ell,\ell'\ra}=\sum_{\ell_1,\ell_2\in\mbb{Z},\;\ell_3,\ell_4\in\ol{0,m};\;\ell_1+\ell_3=\ell,\;
\ell_2+\ell_4=\ell'}\bar{\msr A}_{\la\ell_1,\ell_2\ra}\check{\msr
A}_{\ell_3,\ell_4}.\eqno(13.2.13)$$ It can be verified that ${\msr
A}_{\la\ell,\ell'\ra}$ forms a $gl(n|m)$-submodule. Define
$${\msr H}=\{f\in {\msr A}\mid \Dlt_{n_1,n_2}(f)=0\},\;\;{\msr H}_{\la\ell_1,\ell_2\ra}={\msr
H}\bigcap{\msr A}_{\la\ell_1,\ell_2\ra}.\eqno(13.2.14)$$ By
(13.2.11), ${\msr H}_{\la\ell,\ell'\ra}$ forms a $gl(n|m)$-submodule
of ${\msr A}_{\la\ell,\ell'\ra}$. According to (6.6.16) and
(6.6.17), we have:\psp

{\bf Lemma 13.2.1}. {\it The  nonzero vectors in
$$\{\mbb{F}[\bar\eta](x_i^{m_1}y_j^{m_2})\mid
m_1,m_2\in\mbb{N};i=n_1,n_1+1;j=n_2,n_2+1-\dlt_{n_2,n}\}\eqno(13.2.15)$$
are all the singular vectors of $\bar{\msr G}$ (cf.
(13.1.15)-(13.1.17)) in $\bar{\msr A}$ (cf. (13.1.13))}.\psp

Recall ${\msr G}=\bar{\msr G}+\check{\msr G}$ (cf. (13.1.15)). Take
the Cartan subalgebra $H=\bar H+\check H$ of ${\msr G}$ (cf.
(13.1.16)) and the subspace
 ${\msr G}_+=\bar{\msr G}_++\check{\msr G}_+$ (cf. (13.1.17)) spanned by positive
 root vectors in ${\msr G}$. Then
the  nonzero vectors in
\begin{eqnarray*}\hspace{2cm}& &\{\mbb{F}[\bar\eta,\check\eta](x_i^{m_1}y_j^{m_2}\vec\sta_r\vec\vt_s)\mid
m_1,m_2\in\mbb{N};i=n_1,n_1+1;\\ & &j=n_2,n_2+1-\dlt_{n_2,n};0\leq
r<s\leq m+1\}\hspace{4cm}(13.2.16)\end{eqnarray*} are all the ${\msr
G}$-singular vectors   in ${\msr A}$.
 Choose $H$ as a Cartan
subalgebra of the Lie superalgebra $gl(n|m)$ and $gl(n|m)_+={\msr
G}_++\sum_{r=1}^n\sum_{s=1}^m\mbb{F}E_{r,n+s}$ as the subalgebra
generated by positive root vectors.

Fix $\ell,\ell'\in\mbb{Z}$. Then a $gl(n|m)$-singular vector in
${\msr A}_{\ell,\ell'}$ must be of the form
$$f=\sum_{p=0}^{\min\{s-r-1,\ell_1\}}b_p\bar\eta^{\ell_1-p}\check\eta^p
(x_i^{m_1}y_j^{m_2}\vec\sta_r\vec\vt_s),\eqno(13.2.17)$$ where
$\ell_1,m_1,m_2\in\mbb{N},\;\;0\leq r<s\leq m+1,\;b_p\in\mbb{F}$ and
$$(i,j)\in\{(n_1,n_2),
(n_1,n_2+1-\dlt_{n_2,n}),(n_1+1,n_2),(n_1+1,n_2+1-\dlt_{n_2,n})\}.\eqno(13.2.18)$$

Suppose that $\ell_1=0$ or $s-r-1=0$. Then we can assume
$f=x_i^{m_1}y_j^{m_2}\vec\sta_r\vec\vt_s$. If $r\neq 0$, then
(13.2.5) implies
$$E_{n_1+1,n+r}(f)=(-1)^{r-1}x_{n_1+1}x_i^{m_1}y_j^{m_2}\vec\sta_{r-1}\vec\vt_s\neq 0,\eqno(13.2.19)$$
which is absurd. So $r=0$. Assume that $m_2>0,\;s>1$ and $j=n_2$.
According to (13.2.5),
$$E_{n_2,n+1}(f)=-m_2x_i^{m_1}y_j^{m_2-1}\vt_1\vec\vt_s\neq
0,\eqno(13.2.20)$$ which leads a contradiction. Hence $m_2(s-1)=0$
if $j=n_2$. If $n_2<n$, then we have
$$E_{n_2+1,n+1}(f)=x_i^{m_1}y_j^{m_2}y_{n_2+1}\vt_1\vec\vt_s\neq
0\eqno(13.2.21)$$ by (13.2.5), which is absurd. Thus $s=1$. In
summary,
$$f=x_i^{m_1}y_j^{m_2}\vec\vt_s\qquad\mbox{with}\;s=1\;\mbox{or}\;n_2=n\;\mbox{and}\;m_2=0.\eqno(13.2.22)$$

Consider the case $\ell_1>0$ and $s-r-1>0$. By (13.2.5) and the fact
$\check\eta^{s-r-1}\vt_{r+1}\vec\sta_r\vec\vt_s=0$, we have
\begin{eqnarray*}0&=&E_{n_1+1,n+r+1}(f)=(x_{n_1+1}\ptl_{\sta_{r+1}}-\vt_{r+1}\ptl_{y_{n_1+1}})(f)
\\
&=&[\sum_{p=1}^{\min\{s-r-1,\ell_1\}}pb_p\bar\eta^{\ell_1-p}\check\eta^{p-1}-
\sum_{p=0}^{\min\{s-r-1,\ell_1\}-1}(\ell_1-p)b_p\bar\eta^{\ell_1-p-1}\check\eta^p
(x_i^{m_1}y_j^{m_2}\vec\sta_r\vec\vt_s)]\\& &\times
x_{n_1+1}\vt_{r+1}x_i^{m_1}y_j^{m_2}\vec\sta_r\vec\vt_s,
\hspace{8.8cm}(13.2.23)\end{eqnarray*} which implies
$$f=b_0(\bar\eta+\check\eta)^{\ell_1}(x_i^{m_1}y_j^{m_2}\vec\sta_r\vec\vt_s)=b_0\eta^{\ell_1}(x_i^{m_1}y_j^{m_2}\vec\sta_r\vec\vt_s).
\eqno(13.2.24)$$ The arguments in the previous paragraph and
(13.2.11) give:\psp

{\bf Lemma 13.2.2}. {\it Any $gl(n|m)$-singular vector in ${\msr
A}_{\la\ell,\ell'\ra}$ must be of the form:
$$\eta^{\ell_1}(x_i^{m_1}y_j^{m_2}\vec\vt_s)\;\;\mbox{with}\;\ell_1,m_1,m_2\in\mbb{N},\;s\in\ol{1,m+1}\;\mbox{and (13.2.18)}
\eqno(13.2.25)$$ such that $s=1$ or $n_2=n$ and $m_2=0$.}\psp

We define
$$\td D=\sum_{r=n_1+1}^nx_r\ptl_{x_r}-\sum_{i=1}^{n_1}x_i\ptl_{x_i},\;\;
\td
D'=\sum_{i=1}^{n_2}y_i\ptl_{y_i}-\sum_{r=n_2+1}^ny_r\ptl_{y_r}.\eqno(13.2.26)$$
Then
$$\bar{\msr A}_{\la\ell,\ell'\ra}=\{f\in\bar{\msr A}\mid
\td D(f)=\ell f;\td D'(f)=\ell' f\}.\eqno(13.2.27)$$
 We calculate
$$\bar\Dlt_{n_1,n_2}\bar\eta=\bar\eta\bar\Dlt_{n_1,n_2}+n_2-n_1+\td D+\td D'.\eqno(13.2.28)$$
For $\ell_1\in\mbb{N}+1$ and $f\in{\msr H}_{\la\ell,\ell'\ra}$ (cf.
(13.2.19)), (13.1.35) and (13.2.28) imply
$$\Dlt_{n_1,n_2}\eta^{\ell_1}(f)=\ell_1(n_2-n_1-m+\ell+\ell'+\ell_1-1)\eta^{\ell_1-1}(f).\eqno(13.2.29)$$
Thus
$$\Dlt_{n_1,n_2}\eta^{\ell_1}(f)=0\dar \ell+\ell'\leq
n_1+m-n_2\;\;\mbox{and}\;\;\ell_1=n_1+m-n_2-\ell-\ell'+1.\eqno(13.2.30)$$
If the condition holds, then
$$\eta^{\ell_1}(f)\in{\msr H}_{\la
n_1+m-n_2-\ell'+1,n_1+m-n_2-\ell+1\ra}.\eqno(13.2.31)$$ Moreover,
$$(n_1+m-n_2-\ell'+1)+(n_1+m-n_2-\ell+1)\geq
n_1+m-n_2+2.\eqno(13.2.32)$$

Observe that
$$x_{n_1}^{m_1}y_{n_2}^{m_2}\vec\vt_s\in{\msr
H}_{\la-m_1,m+1+m_2-s\ra},\qquad
x_{n_1}^{m_1}y_{n_2+1}^{m_2}\vec\vt_1\in{\msr
H}_{\la-m_1,m-m_2\ra},\eqno(13.2.33)$$
$$x_{n_1+1}^{m_1}y_{n_2}^{m_2}\vec\vt_s\in{\msr
H}_{\la m_1,m+1+m_2-s\ra},\qquad
x_{n_1+1}^{m_1}y_{n_2+1}^{m_2}\vec\vt_1\in{\msr H}_{ \la
m_1,m-m_2\ra}.\eqno(13.2.34)$$ By Lemma 13.2.2,
$$\mbox{any nonzero}\;{\msr H}_{\la\ell,\ell'\ra}\;\;\mbox{contains a singular
vector of the form}\;x_i^{m_1}y_j^{m_2}\vec\vt_s,\eqno(13.2.35)$$
where $s=1$ or $n_2=n$ and $m_2=0$.

Now we consider $f=x_i^{m_1}y_j^{m_2}\vec\vt_s$ with
$m_1,m_2\in\mbb{N},\;s\in\ol{1,m+1}$ and (13.2.18) such that $s=1$
or $n_2=n$ and $m_2=0$. Assume $\Dlt_{n_1,n_2}\eta^{\ell_1}(f)=0$
for some $\ell_1\in\mbb{N}+1$.\psp

{\it Case 1}. $(i,j)=(n_1,n_2)$.\psp

In this subcase,
 $\ell=-m_1$ and
$\ell'=m_2+m+1-s$ by (13.2.33). Thus $m_2-m_1+1-s\leq n_1-n_2$ and
$\ell_1= n_1+m_1+s-n_2-m_2$ by (13.2.30). So
$$\eta^{\ell_1}(f)\in{\msr H}_{\la
n_1+s-n_2-m_2,n_1+m_1-n_2+m+1\ra}.\eqno(13.2.36)$$

{\it Case 2}. $(i,j)=(n_1,n_2+1)$.\psp

In this subcase $s=1$, $\ell=-m_1$ and $\ell'=m-m_2$ by (13.2.33).
Thus $m_1+m_2\geq n_2-n_1$ and $ \ell_1=n_1+m_1+m_2-n_2+1$ by
(13.2.30). Hence
$$\eta^{\ell_1}(f)\in{\msr H}_{\la
n_1+m_2-n_2+1,n_1+m_1-n_2+m+1\ra}.\eqno(13.2.37)$$

{\it Case 3}. $(i,j)=(n_1+1,n_2)$.\psp

In this subcase,
 $\ell=m_1$ and
$\ell'=m_2+m+1-s$ by (13.2.34). Thus $m_2+m_1+1-s\leq n_1-n_2$ and
$\ell_1= n_1-m_1+s-n_2-m_2$ by (13.2.30). So
$$\eta^{\ell_1}(f)\in{\msr H}_{\la
n_1+s-n_2-m_2,n_1-m_1-n_2+m+1\ra}.\eqno(13.2.38)$$

{\it Case 4}. $(i,j)=(n_1+1,n_2+1)$.\psp

In this subcase $s=1$, $\ell=m_1$ and $\ell'=m-m_2$ by (13.2.34).
Thus $m_1-m_2\leq n_1-n_2$ and $ \ell_1=n_1-m_1+m_2-n_2+1$ by
(13.2.30). Hence
$$\eta^{\ell_1}(f)\in{\msr H}_{\la
n_1+m_2-n_2+1,n_1-m_1-n_2+m+1\ra}.\eqno(13.2.39)$$

Thus we obtain:\psp

{\bf Lemma 13.2.3}. {\it A nonzero $gl(n|m)$-module ${\msr
H}_{\la\ell,\ell'\ra}$ has a unique singular vector if and only if
$\ell+\ell'\leq n_1+m+1-n_2$ or $\ell\not\in\ol{n_1+1-n,n_1+m+1-n}$
and $n_2=n$. If the condition holds, the unique singular vector is
of the form $x_i^{m_1}y_j^{m_2}\vec\vt_s$ with (13.2.11), where
$s=1$ or $n_2=n$ and $m_2=0$.}\psp

Fix ${\msr H}_{\la\ell,\ell'\ra}\neq \{0\}$. Assume
$$v_{\ell,\ell'}=x_i^{m_1}y_j^{m_2}\vec\vt_s\in {\msr
H}_{\la\ell,\ell'\ra}\eqno(13.2.40)$$ for some $(i,j)$ in (13.2.18),
$m_1,m_2\in\mbb{N}$ and $s\in\ol{1,m+1}$ such that $s=1$ or $n_2=n$
and $m_2=0$.\psp

{\bf Lemma 13.2.4}. {\it As a $gl(n|m)$-module, ${\msr
H}_{\la\ell,\ell'\ra}$ is generated by $v_{\ell,\ell'}$.}

{\it Proof}. Denote
$$\td\G=\{\td\al=(\al_{n_1+1},...,\al_{n_2})\in\mbb{N}^{n_2-n_1}\},\;\;|\td\al|=\sum_{i=1}^{n_2-n_1}\al_{n_1+i}.
\eqno(13.2.41)$$ Set
$$\td{\msr A}=\check{\msr
A}[x_{n_1+1},...,x_{n_2},y_{n_1+1},...,y_{n_2}],\eqno(13.2.42)$$
$$\td{\msr H}_{\la\ell_1,\ell_2\ra}=\td{\msr A}\bigcap{\msr
H}_{\la\ell_1,\ell_2\ra},\eqno(13.2.43)$$
$$\td\Dlt_{n_1,n_2}=-\sum_{i=1}^{n_1}x_i\ptl_{y_i}+\sum_{r=n_1+2}^{n_2}\ptl_{x_r}\ptl_{y_r}-\sum_{s=n_2+1}^n
y_s\ptl_{x_s}+\check\Dlt.\eqno(13.2.44)$$ Then
$\Dlt_{n_1,n_2}=\ptl_{x_{n_1+1}}\ptl_{y_{n_1+1}}+\td\Dlt_{n_1,n_2}$.
For any $k_1,k_2\in\mbb{N}$, we define the operator
$$T_{k_1,k_2}=\sum_{i=0}^\infty
\frac{(-1)^ix_{n_1+1}^{k_1+i}y_{n_1+1}^{k_2+i}}{\prod_{r=1}^i(k_1+r)(k_2+r)}\td\Dlt_{n_1,n_2}^i.\eqno(13.2.45)$$
Then Lemma 6.1.1 yields
\begin{eqnarray*}{\msr
H}_{\la\ell,\ell'\ra}&=&\mbox{Span}\{T_{\al_{n_1+1},\be_{n_1+1}}([\prod_{i\neq
n_1+1}x_i^{\al_i}y_i^{\be_i}])\sta_{\vec j}\vt_{\vec
k})\mid\al,\be\in\mbb{N}^n;\\ & &\qquad\vec j\in\G_{\ell_1};\vec
k\in\G_{\ell_2};\ell_1,\ell_2\in\ol{0,m};\al_{n_1+1}\be_{n_1+1}=0;\\
& &\qquad\sum_{s=n_1+1}^n\al_s-\sum_{r=1}^{n_1}\al_r+\ell_1=\ell;
\sum_{r=1}^{n_2}\be_r-\sum_{s=n_2+1}^n\be_s+\ell_2=\ell'\},\hspace{1.3cm}(13.2.46)\end{eqnarray*}
\begin{eqnarray*}\td{\msr
H}_{\la\ell,\ell'\ra}\!\!\!&=&\!\!\!\mbox{Span}\{T_{\al_{n_1+1},\be_{n_1+1}}([\prod_{i=n_1+2}^{n_2}
x_i^{\al_i}y_i^{\be_i}])\sta_{\vec j}\vt_{\vec
k})\mid\td\al,\td\be\in\td\G;\vec j\in\G_{\ell_1};\vec
k\in\G_{\ell_2};\\ & &\qquad
\al_{n_1+1}\be_{n_1+1}=0;\ell_1,\ell_2\in\ol{0,m};|\td\al|+\ell_1=\ell;
|\td\be|+\ell_2=\ell'\}.\hspace{1.7cm}(13.2.47)\end{eqnarray*}

Write
$${\msr G}_1=\sum_{i,j=1}^{n_1}\mbb{F}E_{i,j},\qquad {\msr
G}_2=\sum_{r,s=n_2+1}^n\mbb{F}E_{r,s},\eqno(13.2.48)$$
$${\msr G}_3=\sum_{n_1+1\neq i,j\in\ol{1,n_2}}\mbb{F}E_{i,j},\qquad {\msr
G}_4=\sum_{r,s=n_1+2}^n\mbb{F}E_{r,s}.\eqno(13.2.49)$$ Then
$$\xi\td\Dlt_{n_1,n_2}=\td\Dlt_{n_1,n_2}\xi\qquad\for\;\;\xi\in {\msr
G}_i,\;i\in\ol{1,4}.\eqno(13.2.50)$$
 Denote by $V$
the $gl(n|m)$-submodule of ${\msr H}_{\la\ell,\ell'\ra}$ generated
by $v_{\ell,\ell'}$. By (13.2.2)-(13.2.4)
$$-E_{n_1+1,n_1}|_{\msr
A}=x_{n_1}x_{n_1+1}+y_{n_1}\ptl_{y_{n_1+1}},\;\;E_{n_2+1,n_2}|_{\msr
A}=x_{n_2+1}\ptl_{x_{n_2}}+y_{n_2}y_{n_2+1}.\eqno(13.2.51)$$
According to (13.2.5) and (13.2.6),
$$E_{n_2,n+r}=-x_{n_2}\ptl_{\sta_r}+\vt_r\ptl_{y_{n_2}},\;\;E_{n+r,n_2}=\sta_r\ptl_{x_{n_2}}+y_{n_2}\ptl_{\vt_r}
\qquad\for\;\;r\in\ol{1,m}.\eqno(13.2.52)$$ Repeatedly applying the
operators in (13.2.51) and (13.2.52) to $v_{\ell,\ell'}$, we obtain
$$x_{n_1}^{p_1}x_{n_1+1}^{p_2}y_{n_2}^{p_3}y_{n_2+1}^{p_4}\vec\vt_{s'}\in
V\eqno(13.2.53)$$ for $p_i\in\mbb{N}$ and $s'\in\ol{1,m+1}$ such
that
$$p_2-p_1=\ell,\;p_3-p_4+m+1-s'=\ell',\;p_3(s'-1)=0.\eqno(13.2.54)$$
Lemma 13.1.5, (13.1.69),  and (13.2.43)  tell us that
$$x_{n_1}^{p_1}y_{n_2+1}^{p_4}\td{\msr H}_{
p_2,p_3+m+1-s'}\subset V.\eqno(13.2.55)$$

Recall that $U({\msr G}_i)$ stands for the universal enveloping of
the Lie algebra ${\msr G}_i$. Applying $U({\msr G}_1)$ and $U({\msr
G}_2)$ to (13.2.55), we get
$$[\prod_{\iota_1=1}^{n_1}x_{\iota_1}^{\al_{\iota_1}}][\prod_{\iota_2=n_2+1}^n y_{\iota_2}^{\be_{\iota_2}}]
\td{\msr H}_{\la p_2,p_3+m+1-s'\ra}\subset V\eqno(13.2.56)$$ for any
$(\al_1,...,\al_{n_1})\in\mbb{N}^{n_1}$ and
$(\be_{n_2+1},...,\be_n)\in\mbb{N}^{n-n_2}$ such that
$\sum_{\iota_1=1}^{n_1}\al_{\iota_1}=p_1$ and
$\sum_{\iota_2=n_2+1}^n \be_{\iota_2}=p_4$. Applying $U({\msr G}_3)$
to (13.2.56), we obtain
$$T_{\al_{n_1+1},\be_{n_1+1}}([\prod_{i\neq
n_1+1}x_i^{\al_i}y_i^{\be_i}])\sta_{\vec j}\vt_{\vec k})\in
V\eqno(13.2.57)$$ by (13.2.47) and (13.2.50), where $\al,\be,\vec
j,\vec k$ are as those in (13.2.46) and
$\al_{n_2+1}=\cdots=\al_n=0$. Finally, we get (13.2.57) with any
$\al,\be,\vec j,\vec k$ in (13.2.46) by (13.2.50) and applying
$U({\msr G}_4)$. According to (13.2.46), $V={\msr
H}_{\la\ell,\ell'\ra}.\qquad\Box$\psp

{\bf Theorem 13.2.5}. {\it Let $\ell,\ell'\in\mbb{Z}$ such that
$\ell'\geq 0$ if $n_2=n$.  The $gl(n|m)$-module ${\msr
H}_{\la\ell,\ell'\ra}$ is irreducible if and only if $\ell+\ell'\leq
n_1+m+1-n_2$ or $\ell\not\in\ol{n_1+1-n,n_1+m+1-n}$ and $n_2=n$.
When $\ell+\ell'\leq n_1+m+1-n_2$, ${\msr
A}_{\la\ell,\ell'\ra}=\bigoplus_{i=0}^\infty\eta^i({\msr
H}_{\la\ell-i,\ell'-i\ra})$ is the decomposition of irreducible
$gl(n|m)$-submodules if  $n_2<n$, and ${\msr
A}_{\la\ell,\ell'\ra}=\bigoplus_{i=0}^{\ell'}\eta^i({\msr
H}_{\la\ell-i,\ell'-i\ra})$ is the decomposition of irreducible
$gl(n|m)$-submodules if $n_2=n$.}

{\it Proof}. Suppose $\ell+\ell'\leq n_1+m+1-n_2$ or
$\ell\not\in\ol{n_1+1-n,n_1+m+1-n}$ and $n_2=n$. Let $V$ be a
nonzero submodule of ${\msr H}_{\la\ell,\ell'\ra}$. According to
Lemma 13.2.3, the vector $v_{\ell,\ell'}$ in (13.2.40) is the unique
singular vector of ${\msr H}_{\la\ell,\ell'\ra}$. Since $gl(n|m)_+$
in (13.1.52) is locally nilpotent by (13.2.2)-(13.2.6), $V$ contains
a singular vector. So $v_{\ell,\ell'}\in V$. By Lemma 13.2.4,
$V={\msr H}_{\la\ell,\ell'\ra}$; that is, ${\msr
H}_{\la\ell,\ell'\ra}$ is irreducible. The necessity also follows
from Lemma 12.2.3.

Assume $\ell+\ell'\leq n_1+m+1-n_2$. Since $\Dlt_{n_1,n_2}$ is
locally nilpotent by (13.2.9) and (13.2.10), for any $0\neq u\in
{\msr A}_{\la\ell,\ell'\ra}$, there exists an element
$\kappa(u)\in\mbb{N}$ such that
$$\Dlt_{n_1,n_2}^{\kappa(u)}(u)\neq
0\;\;\mbox{and}\;\;\Dlt_{n_1,n_2}^{\kappa(u)+1}(u)=0.\eqno(13.2.58)$$
Set
$$\Psi=\left\{\begin{array}{ll}\sum_{i=0}^\infty\eta^i({\msr
H}_{\la\ell-i,\ell'-i\ra})&\mbox{if}\;n_2<n,\\
\sum_{i=0}^{\ell'}\eta^i({\msr
H}_{\ell-i,\ell'-i})&\mbox{if}\;n_2=n.\end{array}\right.\eqno(13.2.59)$$
Given $0\neq u\in {\msr A}_{\la\ell,\ell'\ra}$, $\kappa(u)=1$
implies $u\in {\msr H}_{\la\ell,\ell'\ra}\subset\Psi$. Suppose that
$u\in \Psi$ whenever $\kappa(u)<r$ for some positive integer $r$.
Assume $\kappa(u)=r$. First
$$v=\Dlt_{n_1,n_2}^r(u)\in {\msr
H}_{\la\ell-r,\ell'-r\ra}\subset\Psi.\eqno(13.2.60)$$ Note
$$\Dlt_{n_1,n_2}^r[\eta^r(v)]=r![\prod_{i=1}^r(n_2-n_1-m+\ell+\ell'-r-i)]v\eqno(13.2.61)$$
by (13.2.29). Thus we  have either
$$u=\frac{1}{r![\prod_{i=1}^r(n_2-n_1-m+\ell+\ell'-r-i)]}\eta^r(v)\in\Psi\eqno(13.2.62)$$
or
$$\kappa\left(u-\frac{1}{r![\prod_{i=1}^r(n_2-n_1-m+\ell+\ell'-r-i)]}\eta^r(v)\right)<r.\eqno(13.2.63)$$
By induction,
$$u-\frac{1}{r![\prod_{i=1}^r(n_2-n_1-m+\ell+\ell'-r-i)]}\eta^r(v)\in\Psi,\eqno(13.2.64)$$
which implies $u\in\Psi$. Therefore, we have $\Psi={\msr
A}_{\la\ell,\ell'\ra}$. Since all $\eta^i({\msr
H}_{\la\ell-i,\ell'-i\ra})$ have distinct highest weights, the sums
in (13.2.58) are direct sums. $\qquad\Box$\psp

{\bf Remark 13.2.6}. If $\ell+\ell'> n_1+m+1-n_2$  and
$\ell\in\ol{n_1+1-n,n_1+m+1-n}$ when $n_2=n$, the $gl(n|m)$-module
${\msr H}_{\la\ell,\ell'\ra}$ is indecomposable. When $n_2<n$ and
$\ell+\ell'> n_1+m+1-n_2$, ${\msr A}_{\la\ell,\ell'\ra}$ is not
completely reducible. With the settings (13.2.9) and (13.2.10),  the
space
$\mbb{F}\Dlt_{n_1,n_2}+\mbb{F}[\Dlt_{n_1,n_2},\eta]+\mbb{F}\eta$
again forms an operator Lie algebra isomorphic to $sl(2,\mbb{F})$.
The above theorem establishes a supersymmetric
$(sl(2,\mbb{F}),gl(n|m))$ Howe duality on the homogeneous subspaces
${\msr A}_{\la\ell,\ell'\ra}$ with $\ell,\ell'\in\mbb{Z}$ such that
 $\ell+\ell'\leq n_1+m+1-n_2$ and $\ell'\geq 0$ if $n_2=n$.\psp

Recall the notations in (13.1.99) and (13.1.100). The set
\begin{eqnarray*}\hspace{1cm}& &\{\sum_{i=0}^\infty\frac{(-1)^ix_{n_1+1}^iy_{n_1+1}^i}
{\prod_{r=1}^i(\al_{n_1+1}+r)(\be_{n_1+1}+r)}
(\Dlt_{n_1,n_2}-\ptl_{x_{n_1+1}}\ptl_{y_{n_1+1}})^i(x^\al
y^\be\sta_{\vec j}\vt_{\vec k})\\ & &\qquad
\mid\al,\be\in\mbb{N}^n;\vec j\in\G_{\ell_1};\vec
k\in\G_{\ell_2};\al_{n_1+1}\be_{n_1+1}=0;\ell_1,\ell_2\in\ol{0,m};\\
& &\qquad \sum_{r=n_1+1}^n\al_r-\sum_{i=1}^{n_1}\al_i+\ell_1=\ell;
\sum_{i=1}^{n_2}\be_i-\sum_{r=n_2+1}^n\be_r+\ell_2=\ell'\}\hspace{2cm}(13.2.65)\end{eqnarray*}
 forms a basis of ${\msr H}_{\la\ell,\ell'\ra}$ by Lemma 6.1.1.

\section{Oscillator Representations of  $osp(2n|2m)$}

In this section, we use the results in last two sections to study
the corresponding the even ortho-symplectic Lie superalgebras.

 Set
$${\msr K}_0=\sum_{i,j=1}^n\mbb{F}(E_{i,j}-E_{n+j,n+i})+\sum_{r,s=1}^m\mbb{F}
(E_{2n+r,2n+s}-E_{2n+m+s,2n+m+r}),\eqno(13.3.1)$$
$${\msr
K}_1=\sum_{i=1}^n\sum_{r=1}^m[\mbb{F}(E_{i,2n+r}-E_{2n+m+r,n+i})+\mbb{F}
(E_{2n+r,i}+E_{n+i,2n+m+r})].\eqno(13.3.2)$$ Then ${\msr K}={\msr
K}_0+{\msr K}_1$ forms a Lie sub-superalgebra of $gl(2n|2m)$
isomorphic to $gl(n|m)$. Let \begin{eqnarray*} \qquad
osp(2n|2m)_0&=&{\msr K}_0+\sum_{1\leq i<j\leq
n}[\mbb{F}(E_{i,n+j}-E_{j,n+i})+\mbb{F}(E_{n+i,j}-E_{n+j,i})]
\\ & &+\sum_{1\leq r\leq s\leq m}[\mbb{F} (E_{2n+r,2n+m+s}+E_{2n+s,2n+m+r})
\\ & &\quad+\mbb{F}
(E_{2n+m+r,2n+s}+E_{2n+m+s,2n+r})],\hspace{3.8cm}(13.3.3)\end{eqnarray*}
$$osp(2n|2m)_1={\msr
K}_1+\sum_{i=1}^n\sum_{r=1}^m[\mbb{F}(E_{i,2n+m+r}+E_{2n+r,n+i})+\mbb{F}
(E_{n+i,2n+r}-E_{2n+m+r,i})].\eqno(13.3.4)$$ The space
$osp(2n|2m)=osp(2n|2m)_0+osp(2n|2m)_1$ forms a simple Lie
sub-superalgebra of  $gl(2n|2m)$, which is called an {\it even
ortho-symplectic Lie superalgebra}. \index{even ortho-symplectic Lie
superalgebra} Moreover, its Lie subalgebra
$$osp(2n|2m)_0\cong o(2n,\mbb{F})\oplus sp(2n,\mbb{F}).\eqno(13.3.5)$$

Take settings in (13.1.4)-(13.1.10).  Define a representation of
$osp(2n|2m)$ on ${\msr A}$ determined by
$$(E_{i,j}-E_{n+j,n+i})|_{\msr
A}=x_i\ptl_{x_j}-y_j\ptl_{y_i},\;\;(E_{i,2n+r}-E_{2n+m+r,n+i})|_{\msr
A}=x_i\ptl_{\sta_r}-\vt_r\ptl_{y_i},\eqno(13.3.6)$$
$$(E_{i,n+j}-E_{j,n+i})|_{\msr A}=x_i\ptl_{y_j}-x_j\ptl_{y_i},\;\;(E_{n+i,j}-E_{n+j,i})|_{\msr
A}=y_i\ptl_{x_j}-y_j\ptl_{x_i},\eqno(13.3.7)$$
$$(E_{2n+r,i}+E_{n+i,2n+m+r})|_{\msr
A}=\sta_r\ptl_{x_i}+y_i\ptl_{\vt_r},\eqno(13.3.8)$$
$$(E_{i,2n+m+r}+E_{2n+r,n+i})|_{\msr
A}=x_i\ptl_{\vt_r}+\sta_r\ptl_{y_i},\eqno(13.3.9)$$
$$(E_{2n+r,2n+s}-E_{2n+m+s,2n+m+r})|_{\msr
A}=\sta_r\ptl_{\sta_s}-\vt_s\ptl_{\vt_r},\eqno(13.3.10)$$
$$(E_{2n+m+r,2n+s}+E_{2n+m+s,2n+r})]|_{\msr
A}=\vt_r\ptl_{\sta_s}+\vt_s\ptl_{\sta_r},\eqno(13.3.11)$$
$$\;\;(E_{2n+r,2n+m+s}+E_{2n+s,2n+m+r})|_{\msr
A}=\sta_r\ptl_{\vt_s}+\sta_s\ptl_{\vt_r},\eqno(13.3.12)$$
$$
(E_{n+i,2n+r}-E_{2n+m+r,i})|_{\msr
A}=y_i\ptl_{\sta_r}-\vt_r\ptl_{x_i}\eqno(13.3.13)$$
 for $i,j\in\ol{1,n}$ and $r,s\in\ol{1,m}$.

Recall that we write $\Theta_1=\sum_{r=1}^m\mbb{F}\sta_r$ and
$\Theta_2=\sum_{s=1}^m\mbb{F}\vt_s$. For $k\in\mbb{N}$, we denote
$${\msr A}_k=\mbox{Span}\{x^\al
y^\al\Theta_1^{\ell'_1}\Theta_2^{\ell'_2}\mid\al,\be\in\mbb{N}^n;\ell_1',\ell_2'\in\mbb{N};|\al|+\ell_1'+
|\be|+\ell_2'=k\}.\eqno(13.3.14)$$ Again we take
$$\Dlt=\sum_{i=1}^n\ptl_{x_i}\ptl_{y_i}+\sum_{r=1}^m\ptl_{\sta_r}\ptl_{\vt_r},\qquad
\eta=\sum_{i=1}^nx_iy_i+\sum_{r=1}^m\sta_r\vt_r.\eqno(13.3.15)$$ Set
$${\mbb
W}_0=[\sum_{i=1}^n(\mbb{F}x_i+\mbb{F}y_i)+\sum_{r=1}^m(\mbb{F}\sta_r+\mbb{F}\vt_r)]
[\sum_{j=1}^n(\mbb{F}\ptl_{x_{_j}}+\mbb{F}\ptl_{y_{_j}})+\sum_{s=1}^m(\mbb{F}\ptl_{\sta_s}+\mbb{F}\ptl_{\vt_s})].
\eqno(13.3.16)$$ Then
$$osp(2n|2m)|_{\msr A}=\{T\in{\mbb W}_0\mid T(\eta)=0\}.\eqno(13.3.17)$$
Moreover,
$$\xi\Dlt=\Dlt \xi,\qquad \xi\eta=\eta
\xi\qquad\for\;\;\xi\in osp(2n|2m)\eqno(13.3.18)$$ as operators on
${\msr A}$. For $k\in\mbb{N}$, the subspace
$${\msr H}_k=\{f\in{\msr A}_k\mid\Dlt_{n_1,n_2}(f)=0\}\eqno(13.3.19)$$ forms an
$osp(2n|2m)$-submodule.  We have the first main theorem:\psp

{\bf Theorem 13.3.1}. {\it Suppose $n>1$. For $k\in\mbb{N}$, ${\msr
H}_k$ is an irreducible $osp(2n|2m)$-module if and only if $k\leq
m+1-n$ or $k>2(m+1-n)$. When $k\leq m+1-n$, ${\msr
A}_k=\bigoplus_{i=0}^{\llbracket k/2 \rrbracket}\eta^i{\msr
H}_{k-2i}$ is a decomposition of irreducible
$osp(2n|2m)$-submodules}.

{\it Proof}. We take the subspace of diagonal matrices in
$osp(2n|2m)$ as a Carten subalgebras and the subspace
\begin{eqnarray*}osp(2n|2m)_+&=&\sum_{1\leq i<j\leq
n}[\mbb{F}(E_{i,j}-E_{n+j,n+i})+\mbb{F}(E_{i,n+j}-E_{j,n+i})] \\ &
&+\sum_{i=1}^n\sum_{r=1}^m[
\mbb{F}(E_{i,2n+r}-E_{2n+m+r,n+i})+\mbb{F}(E_{i,2n+m+r}+E_{2n+r,n+i})]\\
& &
+\sum_{1\leq r<s\leq m}\mbb{F} (E_{2n+r,2n+s}-E_{2n+m+s,2n+m+r})\\
& &+\sum_{1\leq r\leq s\leq m}\mbb{F}
(E_{2n+r,2n+m+s}+E_{2n+s,2n+m+r})
\hspace{3.5cm}(13.3.20)\end{eqnarray*}
 as the space
spanned by positive root vectors. An $osp(2n|2m)$-{\it singular
vector} $v$ is a nonzero weight vector of $osp(2n|2m)$ such that
$osp(2n|2m)_+(v)=0$. We count singular vector up to a nonzero scalar
multiple.

 Observe ${\msr K}|_{\msr
A}=gl(n|m)|_{\msr A}$. According to the arguments in
(13.1.52)-(13.1.61), the homogeneous singular vectors of ${\msr K}$
are:
$$\{\eta^{\ell_3}x_1^{\ell_1}y_n^{\ell_2}\vec\vt_{s}\mid\ell_i\in\mbb{N};s\in\ol{1,m+1};
\ell_2(s-1)=0\}.\eqno(13.3.21)$$ By (13.3.9),
$$(E_{n,2n+m+s}+E_{2n+s,2n})|_{\msr
A}=x_n\ptl_{\vt_s}+\sta_s\ptl_{y_n}.\eqno(13.3.22)$$ Thus the
homogeneous singular vectors of $osp(2n|2m)$ are
$$\{\eta^{\ell_2}x_1^{\ell_1}\mid\ell_1,\ell_2\in\mbb{N}\}.\eqno(13.3.23)$$
Moreover, (13.1.35) imply
$$\Dlt(\eta^{\ell_2}x_1^{\ell_1})=\ell_2(n-m+\ell_1+\ell_2-1)
=0\lra \ell_1+n\leq
m\;\mbox{and}\;\ell_2=m+1-n-\ell_1.\eqno(13.3.24)$$ In this case,
$$\eta^{m+1-n-\ell_1}x_1^{\ell_1}\in{\msr H}_{2(m+1-n)-\ell_1}.\eqno(13.3.25)$$ Thus
\begin{eqnarray*}\qquad\qquad& &{\msr H}_k\;\mbox{has a unique singular vector if and
only if}\\ & &k\leq m+1-n\;\mbox{or}\;
k>2(m+1-n),\hspace{5.5cm}(13.3.26)\end{eqnarray*} and
$${\msr H}_k\;\mbox{has two singular vectors when}\;
 m+1-n< k\leq 2(m+1-n).\eqno(13.3.27)$$

Note $x_1^k\in{\msr H}_k$. Let $U$ be the $osp(2n|2m)$-submodule
generated by $x_1^k$. Repeatedly applying
$$(E_{n+1,2n+r}-E_{2n+m+r,1})|_{\msr
A}=y_1\ptl_{\sta_r}-\vt_r\ptl_{x_1}\eqno(13.3.28)$$  by (13.3.13),
we get
$$x_1^\ell\vec\vt_s\in U\qquad\for\;\;\ell+m+1-s=k.\eqno(13.3.29)$$
According to (13.3.7),
$$(E_{n+1,n}-E_{2n,1})|_{\msr
A}=y_1\ptl_{x_n}-y_n\ptl_{x_1}.\eqno(13.3.30)$$
 Thus
$$x_1^{\ell_1}y_n^{\ell_2}\vec\vt_1\in
U\qquad\for\;\;\ell_1+\ell_2+m=k\eqno(13.3.31)$$ when $k\geq m$.
Since
$${\msr H}_k=\sum_{i=0}^k{\msr H}_{i,k-i},\eqno(13.3.32)$$
Lemma 13.1.5, (13.3.29) and (13.3.31) imply $U={\msr H}_k$. So
${\msr H}_k$ is an $osp(2n|2m)$-module generated by $x_1^k$.

Suppose $k\leq m+1-n$ or $k>2(m+1-n)$. Let $M$ be a nonzero
$osp(2n|2m)$-submodule of ${\msr H}_k$. By (13.3.26), ${\msr H}_k$
contains a unique singular vector $x_1^k$. Thus $x_1^k\in M$. By the
above paragraph, ${\msr H}_k\subset M$. Hence ${\msr H}_k$ is
irreducible.

If ${\msr H}_k$ is irreducible, (13.3.27) implies $k\leq m+1-n$ or
$k>2(m+1-n)$.

Assume $k\leq m+1-n$. Note
$${\msr A}_j=\sum_{\ell=0}^j{\msr A}_{\ell,j-\ell}\qquad\for\;\;j\in\mbb{N}.\eqno(13.3.33)$$ By Theorem
13.1.6,
$${\msr A}_k=\bigoplus_{\ell=0}^k{\msr A}_{\ell,k-\ell}=
\bigoplus_{\ell=0}^k\bigoplus_{i=0}^{\llbracket k/2
\rrbracket}\eta^i{\msr
H}_{\ell-i,k-\ell-i}=\bigoplus_{i=0}^{\llbracket k/2
\rrbracket}\eta^i{\msr H}_{k-2i}.\qquad\Box\eqno(13.3.34)$$
 \pse

The above theorem establishes a supersymmetric
$(sl(2,\mbb{F}),osp(2n|2m))$ Howe duality on the homogeneous
subspaces ${\msr A}_k$ with $k\in\mbb{N}$ such that $k\leq m+1-n$.

We remark that if $m+1-n<k\leq 2(m+1-n)$, then ${\msr H}_k$ is an
indecomposable $gl(n|m)$-module. In fact, ${\msr
H}_k\bigcap\eta{\msr A}_{k-2}\neq\{0\}$. This also shows that ${\msr
A}_k$ is not completely reducible when $k>m+1-n$.  The conclusion
with $n>m+1$ was first obtained by Zhang [Z].

 Fix $n_1,n_2\in\ol{1,n}$ with $n_1+1<n_2$.
  Changing operators $\ptl_{x_r}\mapsto -x_r,\;
 x_r\mapsto
\ptl_{x_r}$  for $r\in\ol{1,n_1}$ and $\ptl_{y_s}\mapsto -y_s,\;
 y_s\mapsto\ptl_{y_s}$  for $s\in\ol{n_2+1,n}$ in (13.3.6)-(13.3.13), we get a new representation of $osp(2n|2m)$ on
  ${\msr A}$ determined by
$$E_{2n+r,2n+s}|_{\msr A}=\sta_r\ptl_{\sta_s},\qquad
E_{2n+m+r,2n+m+s}|_{\msr A}=\vt_r\ptl_{\vt_s},\eqno(13.3.35)$$
$$E_{2n+r,2n+m+s}|_{\msr A}=\sta_r\ptl_{\vt_s},\qquad
E_{2n+m+r,2n+s}|_{\msr A}=\vt_r\ptl_{\sta_s},\eqno(13.3.36)$$
$$E_{i,j}|_{\msr A}=\left\{\begin{array}{ll}-x_j\ptl_{x_i}-\delta_{i,j}&\mbox{if}\;
i,j\in\ol{1,n_1};\\ \ptl_{x_i}\ptl_{x_j}&\mbox{if}\;i\in\ol{1,n_1},\;j\in\ol{n_1+1,n};\\
-x_ix_j &\mbox{if}\;i\in\ol{n_1+1,n},\;j\in\ol{1,n_1};\\
x_i\partial_{x_j}&\mbox{if}\;i,j\in\ol{n_1+1,n};
\end{array}\right.\eqno(13.3.37)$$
$$E_{n+i,n+j}|_{\msr A}=\left\{\begin{array}{ll}y_i\ptl_{y_j}&\mbox{if}\;
i,j\in\ol{1,n_2};\\ -y_iy_j&\mbox{if}\;i\in\ol{1,n_2},\;j\in\ol{n_2+1,n};\\
\ptl_{y_i}\ptl_{y_j} &\mbox{if}\;i\in\ol{n_2+1,n},\;j\in\ol{1,n_2};\\
-y_j\partial_{y_i}-\delta_{i,j}&\mbox{if}\;i,j\in\ol{n_2+1,n};
\end{array}\right.\eqno(13.3.38)$$
$$E_{i,n+j}|_{\msr
A}=\left\{\begin{array}{ll}
\ptl_{x_i}\ptl_{y_j}&\mbox{if}\;i\in\ol{1,n_1},\;j\in\ol{1,n_2};\\
-y_j\ptl_{x_i}&\mbox{if}\;i\in\ol{1,n_1},\;j\in\ol{n_2+1,n};\\
x_i\ptl_{y_j}&\mbox{if}\;i\in\ol{n_1+1,n},\;j\in\ol{1,n_2};\\
-x_iy_j&\mbox{if}\;i\in\ol{n_1+1,n},\;j\in\ol{n_2+1,n};\end{array}\right.\eqno(13.3.39)$$
$$E_{n+i,j}|_{\msr A}=\left\{\begin{array}{ll}
-x_jy_i&\mbox{if}\;j\in\ol{1,n_1},\;i\in\ol{1,n_2};\\
-x_j\ptl_{y_i}&\mbox{if}\;j\in\ol{1,n_1},\;i\in\ol{n_2+1,n};\\
y_i\ptl_{x_j}&\mbox{if}\;j\in\ol{n_1+1,n},\;i\in\ol{1,n_2};\\
\ptl_{x_j}\ptl_{y_i}&\mbox{if}\;j\in\ol{n_1+1,n},\;i\in\ol{n_2+1,n};
\end{array}\right.\eqno(13.3.40)$$
$$E_{i,2n+r}|_{\msr
A}=\left\{\begin{array}{ll}\ptl_{x_i}\ptl_{\sta_r}&\mbox{if}\;i\in\ol{1,n_1};\\
x_i\ptl_{\sta_r}&\mbox{if}\;i\in\ol{n_1+1,n};\end{array}\right.\eqno(13.3.41)$$
$$E_{i,2n+m+r}|_{\msr
A}=\left\{\begin{array}{ll}\ptl_{x_i}\ptl_{\vt_r}&\mbox{if}\;i\in\ol{1,n_1};\\
x_i\ptl_{\vt_r}&\mbox{if}\;i\in\ol{n_1+1,n};\end{array}\right.\eqno(13.3.42)$$
$$E_{2n+r,i}|_{\msr
A}=\left\{\begin{array}{ll}-x_i\sta_r&\mbox{if}\;i\in\ol{1,n_1};\\\sta_r\ptl_{x_i}&\mbox{if}\;i\in\ol{n_1+1,n};
\end{array}\right.\eqno(13.3.43)$$
$$E_{2n+m+r,i}|_{\msr
A}=\left\{\begin{array}{ll}-x_i\vt_r&\mbox{if}\;i\in\ol{1,n_1};\\
\vt_r\ptl_{x_i}&\mbox{if}\;i\in\ol{n_1+1,n};\end{array}\right.\eqno(13.3.44)$$
$$E_{n+i,2n+r}|_{\msr
A}=\left\{\begin{array}{ll}y_i\ptl_{\sta_r}&\mbox{if}\;i\in\ol{1,n_2};\\
\ptl_{y_i}\ptl_{\sta_r}&\mbox{if}\;i\in\ol{n_2+1,n};\end{array}\right.\eqno(13.3.45)$$
$$E_{n+i,2n+m+r}|_{\msr
A}=\left\{\begin{array}{ll}y_i\ptl_{\vt_r}&\mbox{if}\;i\in\ol{1,n_2};\\
\ptl_{y_i}\ptl_{\vt_r}&\mbox{if}\;i\in\ol{n_2+1,n};\end{array}\right.\eqno(13.3.46)$$
$$E_{2n+r,n+i}|_{\msr
A}=\left\{\begin{array}{ll}\sta_r\ptl_{y_i}&\mbox{if}\;i\in\ol{1,n_2};\\
-y_i\sta_r&\mbox{if}\;i\in\ol{n_2+1,n};\end{array}\right.
\eqno(13.3.47)$$
$$E_{2n+m+r,n+i}|_{\msr
A}=\left\{\begin{array}{ll}\vt_r\ptl_{y_i}&\mbox{if}\;i\in\ol{1,n_2};\\
-y_i\vt_r&\mbox{if}\;i\in\ol{n_2+1,n};\end{array}\right.\eqno(13.3.48)$$
for $i,j\in\ol{1,n}$ and $r,s\in\ol{1,m}$.

The related Laplace operator becomes
$$\Dlt_{n_1,n_2}=-\sum_{i=1}^{n_1}x_i\ptl_{y_i}+\sum_{r=n_1+1}^{n_2}\ptl_{x_r}\ptl_{y_r}-\sum_{s=n_2+1}^n
y_s\ptl_{x_s}+\sum_{r=1}^m\ptl_{\sta_r}\ptl_{\vt_r}\eqno(13.3.49)$$
and its dual
$$\eta=\sum_{i=1}^{n_1}y_i\ptl_{x_i}+\sum_{r=n_1+1}^{n_2}x_ry_r+\sum_{s=n_2+1}^n
x_s\ptl_{y_s}+\sum_{r=1}^m\sta_r\vt_r.\eqno(13.3.50)$$ It can be
verified that (13.3.18) holds again. Denote
\begin{eqnarray*}{\msr A}_{\la
k\ra}&=&\mbox{Span}\{x^\al
y^\be\Theta_1^{\ell_1'}\Theta_2^{\ell_2'}\mid\al,\be\in\mbb{N}^n;\ell_1',\ell_2'\in\mbb{N};\\
& &\sum_{r=n_1+1}^n\al_r-\sum_{i=1}^{n_1}\al_i+
\sum_{i=1}^{n_2}\be_i-\sum_{r=n_2+1}^n\be_r+\ell_1'+\ell_2'=\ell_2\}\hspace{3.1cm}(13.3.51)\end{eqnarray*}
for $k\in\mbb{Z}$.

Again we set
$${\msr H}_{\la k\ra}=\{f\in {\msr A}_{\la k\ra}\mid
\Dlt_{n_1,n_2}(f)=0\}. \eqno(13.3.52)$$  We have the second main
theorem:\psp

{\bf Theorem 13.3.2}. {\it Let $k\in\mbb{Z}$.  The
$osp(2n|2m)$-module ${\msr H}_{\la k\ra}$ is irreducible if and only
if $k\leq n_1+m+1-n_2$. When $k\leq n_1+m+1-n_2$, ${\msr A}_{\la
k\ra}=\bigoplus_{i=0}^\infty\eta^i({\msr H}_{\la k-2i\ra})$ is the
decomposition of irreducible $osp(2n|2m)$-submodules.}

{\it Proof}. Observe ${\msr K}|_{\msr A}=gl(n|m)|_{\msr A}$ in terms
of the representation of $gl(n|m)$ given in (13.2.2)-(13.2.6). Lemma
13.2.2 says that the homogeneous singular vectors of ${\msr K}$ are
of the form:
$$\eta^{\ell_1}(x_i^{m_1}y_j^{m_2}\vec\vt_s)\;\;
\mbox{with}\;\ell_1,m_1,m_2\in\mbb{N},\;s\in\ol{1,m+1}\eqno(13.3.53)$$
 and
$$(i,j)\in\{(n_1,n_2),
(n_1,n_2+1-\dlt_{n_2,n}),(n_1+1,n_2),(n_1+1,n_2+1-\dlt_{n_2,n})\}.\eqno(13.3.54)$$

{\it Claim 1}. For $k\in\mbb{N}$, ${\msr H}_{\la k\ra}$ is an
$osp(2n|2m)$-module generated by $x_{n_1+1}^k$ and ${\msr H}_{\la
-k\ra}$ is an $osp(2n|2m)$-module generated by $x_{n_1}^k$.\psp

Let $V$ be the $osp(2n|2m)$-module generated by $x_{n_1+1}^k\in
{\msr H}_{\la k\ra}$. By (13.3.38),
$$(E_{n_2+1,n+n_1+1}-E_{n_1+1,n+n_2+1})|_{\msr
A}=x_{n_2+1}\ptl_{y_{n_1+1}}+x_{n_1+1}y_{n_2+1}.\eqno(13.3.55)$$
Thus
$$(E_{n_2+1,n+n_1+1}-E_{n_1+1,n+n_2+1})^{k_1}(x_{n_1+1}^k)=x_{n_1+1}^{k+k_1}y_{n_2+1}^{k_1}\in
V.\eqno(13.3.56)$$ According to (13.3.40),
$$(E_{n+n_2,n_1+1}-E_{n+n_1+1,n_2})|_{\msr
A}=y_{n_2}\ptl_{x_{n_1+1}}-y_{n_1+1}\ptl_{x_{n_2}}.\eqno(13.3.57)$$
Repeatedly applying (13.3.57) to $x_{n_1+1}^k$, we obtain
$$x_{n_1+1}^{k_1}y_{n_2}^{k_2}\in
V\qquad\for\;k_1,k_2\in\mbb{N}\;\mbox{such
that}\;k_1+k_2=k.\eqno(13.3.58)$$ Note that (13.3.44) and (13.3.45)
imply
$$(E_{n+n_1+1,2n+r}-E_{2n+m+r,n_1+1})|_{\msr
A}=y_{n_1+1}\ptl_{\sta_r}-\vt_r\ptl_{x_{n_1+1}}\eqno(13.3.59)$$
Applying (13.3.59) with various $r$ to (13.3.56) and (13.3.58), we
obtain
$$x_{n_1+1}^{\ell_1}y_j^{\ell_2}\vec\vt_s\in {\msr H}_{\la
k\ra}\;\;\mbox{with}\;j\in\{n_2,n_2+1\}\lra
x_{n_1+1}^{\ell_1}y_j^{\ell_2}\vec\vt_s\in V.\eqno(13.3.60)$$ In
terms of (13.2.40),
$$v_{\ell,\ell'}\in V,\qquad \ell+\ell'=k.\eqno(13.3.61)$$
By Lemma 13.2.4,
$${\msr H}_{\la \ell,\ell'\ra}\subset V,\qquad \ell+\ell'=k.\eqno(13.3.62)$$
Therefore
$${\msr H}_{\la k\ra}=\bigoplus_{\ell,\ell'\in\mbb{Z};\;\ell+\ell'=k}
{\msr H}_{\la \ell,\ell'\ra}\subset V.\eqno(13.3.63)$$

Suppose $k>0$. Let $U$ be  the $osp(2n|2m)$-module generated by
$x_{n_1}^k\in{\msr H}_{\la -k\ra}$. Observe
$$(E_{n_2+1,n+n_1}-E_{n_1,n+n_2+1})|_{\msr
A}=x_{n_2+1}\ptl_{y_{n_1}}+y_{n_2+1}\ptl_{x_{n_1}}\eqno(13.3.64)$$
by (13.3.39), and
$$(E_{n+n_1,n_2}-E_{n+n_2,n_1}|_{\msr
A}=y_{n_1}\ptl_{x_{n_1}}+x_{n_1}y_{n_1}\eqno(13.3.65)$$ by
(13.3.40). Repeatedly applying the above two equations to
$x_{n_1}^k$, we have
$$x_{n_1}^{k_1}y_{n_2+1}^{k_2},x_{n_1}^{k+k_3}y_{n_2}^{k_3}\in
U\qquad\for\;\;k_1,k_2,k_3\in\mbb{N}\;\mbox{such
that}\;k_1+k_2=k.\eqno(13.3.66)$$ According to (13.3.44) and
(13.3.45),
$$(E_{n+n_1,2n+r}-E_{2n+m+r,n_1})|_{\msr
A}=y_{n_1}\ptl_{\sta_r}+x_{n_1}\vt_r.\eqno(13.3.67)$$ Applying
(13.3.67) with various $r$ to (13.3.66), we find
$$x_{n_1}^{\ell_1}y_j^{\ell_2}\vec\vt_s\in {\msr H}_{\la-k\ra}\;\;\mbox{with}\;j\in\{n_2,n_2+1\}\lra
x_{n_1}^{\ell_1}y_j^{\ell_2}\vec\vt_s\in V.\eqno(13.3.68)$$ In terms
of (13.2.40),
$$v_{\ell,\ell'}\in V,\qquad \ell+\ell'=-k.\eqno(13.3.69)$$
Lemma 13.2.4 gives
$${\msr H}_{\la\ell,\ell'\ra}\subset V,\qquad \ell+\ell'=-k.\eqno(13.3.70)$$
Thus
$${\msr H}_{\la-k\ra}=\bigoplus_{\ell,\ell'\in\mbb{Z};\;\ell+\ell'=-k}
{\msr H}_{\la\ell,\ell'\ra}\subset V.\eqno(13.3.71)$$ This prove
Claim 1.\psp

{\it Claim 2}. For $n_1+m+1-n_2\geq k\in\mbb{Z}$,  any nonzero
$osp(2n|2m)$-submodule of ${\msr H}_{\la k\ra}$ contains
$x_{n_1+1}^k$ if $k\geq 0$ or $x_{n_1}^{-k}$ when $k<0$.\psp

Note
$$x_{n_1}^{m_1}y_{n_2}^{m_2}\vec\vt_s\in{\msr H}_{\la
m+m_2+1-s-m_1\ra},\qquad
x_{n_1}^{m_1}y_{n_2+1}^{m_2}\vec\vt_s\in{\msr H}_{\la
m+1-s-m_1-m_2\ra},\eqno(13.3.72)$$
$$x_{n_1+1}^{m_1}y_{n_2}^{m_2}\vec\vt_s\in{\msr H}_{\la
m+m_1+m_2+1-s\ra},\qquad
x_{n_1+1}^{m_1}y_{n_2+1}^{m_2}\vec\vt_s\in{\msr H}_{\la
m+m_1+1-s-m_2\ra}.\eqno(13.3.73)$$ For $\ell_1\in\mbb{N}+1$ and
$f\in{\msr H}_{\la k'\ra}$ with $k'\in\mbb{Z}$, (13.2.28) and the
second equation in (13.1.35)  imply
$$\Dlt_{n_1,n_2}\eta^{\ell_1}(f)=\ell_1(n_2-n_1-m+k'+\ell_1-1)\eta^{\ell_1-1}(f).\eqno(13.3.74)$$
Thus
$$\Dlt_{n_1,n_2}\eta^{\ell_1}(f)=0\dar k'\leq
n_1+m-n_2\;\;\mbox{and}\;\;\ell_1=n_1+m-n_2-k'+1.\eqno(13.3.75)$$ If
the condition holds, then
$$\eta^{\ell_1}(f)\in{\msr H}_{\la
2(n_1+m+1-n_2)-k'\ra}.\eqno(13.3.76)$$ Moreover,
$$2(n_1+m+1-n_2)-k'=n_1+m+2-n_2+(n_1+m-n_2-k')\geq n_1+m+2-n_2.\eqno(13.3.77)$$
This shows that
$${\msr H}_{\la k_1\ra}\bigcap(\bigcup_{i=1}^\infty \eta^i({\msr
H}))\neq \{0\}\dar k_1\geq  n_1+m+2-n_2.\eqno(13.3.78)$$

Suppose $k\leq n_1+m+1-n_2$. Then the singular vectors of ${\msr K}$
in ${\msr H}_{\la k\ra}$ are of the form
$$x_i^{m_1}y_j^{m_2}\vec\vt_s\;\;
\mbox{with}\;m_1,m_2\in\mbb{N},\;s\in\ol{1,m+1}\eqno(13.3.79)$$ with
(13.3.54). Observe that
$$(E_{n_1+1,n+n_2}-E_{n_2,n+n_1+1})|_{\msr
A}=x_{n_1+1}\ptl_{y_{n_2}}-x_{n_2}\ptl_{y_{n_1+1}}\eqno(13.3.80)$$
by (13.3.39), and
$$(E_{n+n_2+1,n_1}-E_{n+n_1,n_2+1})|_{\msr
A}=-x_{n_1}\ptl_{y_{n_2+1}}-y_{n_1}\ptl_{x_{n_2+1}}\eqno(13.3.81)$$
by (13.3.40). According to (13.3.42) and (13.3.47),
$$(E_{n_1+1,2n+m+r}+E_{2n+r,n+n_1+1})|_{\msr
A}=x_{n_1+1}\ptl_{\vt_r}+\sta_r\ptl_{y_{n_1+1}}.\eqno(13.3.82)$$ Let
$M$ be any nonzero $osp(2n|2m)$-submodule of ${\msr H}_{\la k\ra}$.
Then $M$ contains at least one of the ${\msr K}$-singular vectors in
(13.3.79). Applying (13.3.70)-(13.3.82), we get
$$x_{n_1}^{k_1}x_{n_2}^{k_2}\in M\;\;\mbox{for
some}\;k_1,k_2\in\mbb{N}\;\mbox{such
that}\;k_2-k_1=k.\eqno(13.3.83)$$ By (13.3.37) and (13.3.38),
$$(E_{n_1,n_1+1}-E_{n+n_1+1,n+n_1})|_{\msr
A}=\ptl_{x_{n_1}}\ptl_{x_{n_1+1}}-y_{n_1+1}\ptl_{y_{n_1}}.\eqno(13.3.84)$$
Repeatedly applying (13.3.84) to (13.3.83), we obtain
$x_{n_1+1}^k\in M$ if $k\geq 0$ or $x_{n_1}^{-k}\in M$ when
$k<0$.\psp

The above claims show that ${\msr H}_{\la k\ra}$ is an irreducible
$osp(2n|2m)$-module if $k\leq  n_1+m+1-n_2$.

Note
$${\msr A}_{\la j\ra}=\sum_{\ell\in\mbb{Z}}{\msr A}_{\la\ell,j-\ell\ra},\qquad
{\msr H}_{\la j\ra}=\sum_{\ell\in\mbb{Z}}{\msr
H}_{\la\ell,j-\ell\ra}\qquad\for\;\;j\in\mbb{Z}.\eqno(13.3.85)$$ By
(13.3.78), $k\leq  n_1+m+1-n_2$ if ${\msr H}_{\la k\ra}$ is
irreducible. When $k\leq  n_1+m+1-n_2$,  Theorem 13.2.5 implies
$${\msr A}_{\la k\ra}=\bigoplus_{\ell\in\mbb{Z}}{\msr A}_{\la\ell,k-\ell\ra}=
\bigoplus_{\ell\in\mbb{Z}}\bigoplus_{i=0}^\infty\eta^i{\msr
H}_{\la\ell-i,k-\ell-i\ra}=\bigoplus_{i=0}^\infty\eta^i{\msr H}_{\la
k-2i\ra}.\qquad\Box\eqno(13.3.86)$$
 \psp

The above theorem establishes a supersymmetric
$(sl(2,\mbb{F}),osp(2n|2m))$ Howe duality on the homogeneous
subspaces ${\msr A}_{\la k\ra}$ with $k\in\mbb{Z}$ such that $k\leq
n_1+m+1-n_2$.

We remark that ${\msr H}_{\la k\ra}$ is an indecomposable
$osp(2n|2m)$-module if $k \geq n_1+m+2-n_2$ by Claim 1 and
(13.3.78). This also implies that  ${\msr A}_{\la k\ra}$ is not
completely reducible $osp(2n|2m)$-module when $k\geq n_1+m+2-n_2$.

\section{Oscillator Representations of $osp(2n+1|2m)$}
In this section, we use the results in last three sections to study
the corresponding the odd ortho-symplectic Lie superalgebras.

Write
$$osp(2n+1|2m)_0=osp(2n|2m)_0+\sum_{i=1}^n[\mbb{F}(E_{0,i}-E_{n+i,0})
+\mbb{F}(E_{0,n+i}-E_{i,0})],\eqno(13.4.1)$$
\begin{eqnarray*}\qquad \qquad osp(2n+1|2m)_1&=&osp(2n|2m)_1+\sum_{r=1}^m[\mbb{F}(E_{0,2n+r}-E_{2n+m+r,0})
\\ & &+\mbb{F}(E_{0,2n+m+r}+E_{2n+r,0})].\hspace{4.3cm}(13.4.2)\end{eqnarray*} The {\it odd ortho-symplectic
 Lie superalgebra} $osp(2n+1|2m)=osp(2n+1|2m)_0+osp(2n+1|2m)_1$ \index{ odd ortho-symplectic
 Lie superalgebra} is a Lie
sub-superalgebra of $gl(2n+1|2m)$. Let $x_0$ be a bosonic
(commuting) variable. Set
$${\msr B}={\msr A}[x_0]=\bigoplus_{k=0}^\infty{\msr B}_k,\qquad{\msr
B}_k=\sum_{i=0}^k{\msr A}_{k-i}x_0^i.\eqno(13.4.3)$$ Moreover, the
corresponding supersymmetric Laplace operator and invariant become
$$\Dlt'=\ptl_{x_0}^2+2\sum_{i=1}^n\ptl_{x_i}\ptl_{y_i}+2\sum_{r=1}^m\ptl_{\sta_r}\ptl_{\vt_r},\qquad
\eta'=x_0^2+2\sum_{i=1}^nx_iy_i+2\sum_{r=1}^m\sta_r\vt_r.\eqno(13.4.4)$$
 Now $osp(2n|2m)$ acts on ${\msr
B}$ by the differential operators in (13.3.6)-(13.3.13); namely, we
change the subindex $|_{\msr A}$ to $|_{\msr B}$. Extend the
representation of $osp(2n|2m)$ on ${\msr B}$ to a representation of
$osp(2n+1|2m)$ on ${\msr B}$ by:
$$(E_{0,i}-E_{n+i,0})|_{\msr B}=x_0\ptl_{x_i}-y_i\ptl_{x_0},
\;\;(E_{0,2n+m+r}+E_{2n+r,0})|_{\msr
B}=x_0\ptl_{\vt_r}+\sta_r\ptl_{x_0} \eqno(13.4.5)$$
$$(E_{0,n+i}-E_{i,0})|_{\msr
B}=x_0\ptl_{y_i}-x_i\ptl_{x_0},,\;\;(E_{0,2n+r}-E_{2n+m+r,0})|_{\msr
B}=x_0\ptl_{\sta_r}-\vt_r\ptl_{x_0},\eqno(13.4.6)$$ for
$i\in\ol{1,n}$ and $r\in\ol{1,m}$.

Set
$${\mbb W}_0'=[\sum_{i=0}^n(\mbb{F}x_i+\mbb{F}y_i)+\sum_{r=1}^m(\mbb{F}\sta_r+\mbb{F}\vt_r)]
[\sum_{j=0}^n(\mbb{F}\ptl_{x_{_j}}+\mbb{F}\ptl_{y_{_j}})+\sum_{s=1}^m(\mbb{F}\ptl_{\sta_s}+\mbb{F}\ptl_{\vt_s})].
\eqno(13.4.7)$$ Then
$$osp(2n+1|2m)|_{\msr B}=\{T\in {\mbb W}_0'\mid
T(\eta')=0\}.\eqno(13.4.8)$$ Define
$${\msr H}'=\{f\in{\msr B}'\mid \Dlt'(f)=0\},\qquad {\msr
H}'_k={\msr H}'\bigcap {\msr B}_k.\eqno(13.4.9)$$

Again we take the subspace of diagonal matrices in $osp(2n+1|2m)$ as
a Cartan subalgebra and take the space spanned by positive roots:
$$osp(2n+1|2m)_+=\mbb{F}(E_{0,n+i}-E_{i,0})+\mbb{F}(E_{0,2n++m+r}+E_{2n+r,0})+
osp(2n|2m)_+.\eqno(13.4.10)$$ An $osp(2n+1|2m)$-{\it singular
vector} $v$ is a nonzero weight vector of $osp(2n+1|2m)$ such that
$osp(2n+1|2m)_+(v)=0$. We count singular vector up to a nonzero
scalar multiple. According to the proof of Theorem 13.2.1, any
singular vector of $osp(2n+1|2m)$ must be in
$\mbb{F}[x_0,x_1,\eta]$, where
$$\eta=\sum_{i=1}^nx_iy_i+\sum_{r=1}^m\sta_r\vt_r.\eqno(13.4.11)$$
Note that $\eta'=x_0^2+2\eta$. By (13.4.10), $x_1^k$ is a singular
vector of $osp(2n+1|2m)$ for any $k\in\mbb{N}$. Thus a homogeneous
singular vector of $osp(2n+1|2m)$ must be of the form
$$f=\sum_{i=0}^\ell
b_ix_0^{2i+\iota}\eta^{\ell-i}x_1^k,\eqno(13.4.12)$$ where
$b_i\in\mbb{F}$, $\ell,k\in\mbb{N}$ and $\iota=0,1.$ Note
$$(E_{0,n+i}-E_{i,0})(f)=(x_0\ptl_{y_i}-x_i\ptl_{x_0})(f)=0\dar
f=b_0{\eta'}^\ell x_1^k.\eqno(13.4.13)$$ Thus $\{{\eta'}^\ell
x_1^k\mid \ell,k\in\mbb{N}\}$ are all the homogeneous singular
vectors of $osp(2n+1|2m)$ in ${\msr B}$.

Observe that
$$[\Dlt',\eta']=2+4(n-m)+4[x_0\ptl_{x_0}+\sum_{i=1}^n(x_i\ptl_{x_i}
+y_i\ptl_{y_i})+\sum_{r=1}^m(\sta_r\ptl_{\sta_r}+\vt_r\ptl_{\vt_r})]\eqno(13.4.14)$$
 by (13.1.35). So
$$\Dlt'({\eta'}^\ell g)=2\ell(1+2(n-m+k+\ell-1)){\eta'}^{\ell-1}
g\qquad\for\;\;g\in{\msr H}_k\eqno(13.4.15)$$ Thus
$${\msr H}'_k\;\mbox{has a unique singular vector}\;x_1^k\;\mbox{for
any}\;k\in\mbb{N}.\eqno(13.4.16)$$ Indeed we have the first main
theorem in this section:\psp

{\bf Theorem 13.4.1}. {\it For any $k\in\mbb{N}$, ${\msr H}_k'$ is
an irreducible $osp(2n+1|2m)$-module. Moreover,
 $${\msr
B}=\bigoplus_{\ell,k=0}^\infty {\eta'}^\ell {\msr
H}_k\eqno(13.4.17)$$ is a direct sum of irreducible
$osp(2n+1|2m)$-submodules.}

{\it Proof}. By the arguments in (13.2.58)-(13.2.64), we only need
to prove that ${\msr H}'_k$ is an $osp(2n+1|2m)$-module generated by
$x_1^k$. For $\iota=0,1$, we define
$$T_\iota=\sum_{i=0}^\infty\frac{(-2)^ix_0^{2i+\iota}}{(2i+\iota)!}\Dlt^i,\qquad\Dlt=\sum_{i=1}^n\ptl_{x_i}\ptl_{y_i}
+\sum_{r=1}^m\ptl_{\sta_r}\ptl_{\vt_r} .\eqno(13.4.18)$$ We take the
notations in (13.1.98) and (13.1.99). By Lemma 6.1.1,
\begin{eqnarray*}\qquad{\msr H}_k'&=&\mbox{Span}\{T_\iota(x^\al y^\be
\sta_{\vec j}\vt_{\vec j'})\mid\al,\be\in\mbb{N}^n;\vec j\in
\G_{k_1},\;\vec j'\in
\G_{k_2};\\
&&\qquad\iota\in\{0,1\};|\al|+|\be|+k_1+k_2+\iota=k\}.\hspace{4.2cm}(13.4.19)\end{eqnarray*}
Let $U$ be the $osp(2n+1|2m)$-module generated by $x_1^k\in{\msr
H}_k$. Since $o(2n+1,\mbb{F})$ is a subalgebra of $osp(2n+1|2m)$,
Theorem 8.1.1 yields
$$T_\iota(x^\al y^\be)\in
U\qquad\for\;\;\al,\be\in\mbb{N}^n;|\al|+|\be|=k.\eqno(13.4.20)$$
Repeatedly applying
$$(E_{2n+r,i}+E_{n+i,2n+m+r})|_{\msr
A}=\sta_r\ptl_{x_i}+y_i\ptl_{\vt_r}\eqno(13.4.21)$$
 to (13.4.20) with $i\in\ol{1,n}$
and $r\in\ol{1,m}$ (cf. (13.3.8)), we obtain
$$T_\iota(x^\al y^\be\sta_{\vec j})\in
U\qquad\for\;\;\al,\be\in\mbb{N}^n,\;\vec
j\in\G_{k_1};|\al|+|\be|+k_1=k.\eqno(13.4.22)$$ Finally, we get
$U={\msr H}_k'$ by repeatedly applying
$(E_{i,2n+r}-E_{2n+m+r,n+i})|_{\msr
B}=x_i\ptl_{\sta_r}-\vt_r\ptl_{y_i}$ to (13.4.18) with
$i\in\ol{1,n}$ and $r\in\ol{1,m}$ (cf. (13.3.6)).$\qquad\Box$\psp

It can be verified by (13.4.14) that the space
$\mbb{F}\Dlt'+\mbb{F}[\Dlt',\eta']+\mbb{F}\eta'$ also forms an
operator Lie algebra isomorphic to $sl(2,\mbb{F})$. The above
theorem establishes a supersymmetric $(sl(2,\mbb{F}),osp(2n+1|2m))$
 Howe duality. We remark that the conclusion in the above theorem with $n>m$ was first obtained
  by Zhang [Z].

 Next
$osp(2n|2m)$ acts on ${\msr B}$ via the differential operators in
(13.3.35)-(13.3.48); namely, we change the subindex $|_{\msr A}$ to
$|_{\msr B}$. Moreover, we extend the representation of $osp(2n|2m)$
on ${\msr B}$ to a representation of $osp(2n+1|2m)$ on ${\msr B}$
by:
$$(E_{0,i}-E_{n+i,0})|_{\msr
B}=\left\{\begin{array}{ll}-x_0x_i-y_i\ptl_{x_0}&\mbox{if}\;i\in\ol{1,n_1},\\
x_0\ptl_{x_i}-y_i\ptl_{x_0}&\mbox{if}\;i\in\ol{n_1+1,n_2},\\
x_0\ptl_{x_i}-\ptl_{x_0}\ptl_{y_i}&\mbox{if}\;i\in\ol{n_2+1,n};\end{array}\right.
\eqno(13.4.23)$$
$$(E_{0,n+i}-E_{i,0})|_{\msr
B}=\left\{\begin{array}{ll}x_0\ptl_{y_i}-\ptl_{x_i}\ptl_{x_0}&\mbox{if}\;i\in\ol{1,n_1},\\
x_0\ptl_{y_i}-x_i\ptl_{x_0}&\mbox{if}\;i\in\ol{n_1+1,n_2},\\
-x_0y_i-x_i\ptl_{x_0}&\mbox{if}\;i\in\ol{n_2+1,n};\end{array}\right.
\eqno(13.4.24)$$
$$(E_{0,2n+r}-E_{2n+m+r,0})|_{\msr B}=x_0\ptl_{\sta_r}-\vt_r\ptl_{x_0},
\eqno(13.4.25)$$ $$ (E_{0,2n++m+r}+E_{2n+r,0})|_{\msr
B}=x_0\ptl_{\vt_r}+\sta_r\ptl_{x_0}\eqno(13.4.26)$$ for
$i\in\ol{1,n}$ and $r\in\ol{1,m}$.

Now the corresponding Laplace operator becomes
$$\Dlt_{n_1,n_2}'=\ptl_{x_0}^2+2\Dlt_{n_1,n_2},\;\;\eta'=x_0^2+2\eta\eqno(13.4.27)$$
with
$$\Dlt_{n_1,n_2}=-\sum_{i=1}^{n_1}x_i\ptl_{y_i}+\sum_{r=n_1+1}^{n_2}
\ptl_{x_r}\ptl_{y_r}-\sum_{s=n_2+1}^n
y_s\ptl_{x_s}+\sum_{r=1}^m\ptl_{\sta_r}\ptl_{\vt_r}\eqno(13.4.28)$$
and
$$\eta=\sum_{i=1}^{n_1}y_i\ptl_{x_i}+\sum_{r=n_1+1}^{n_2}x_ry_r+\sum_{s=n_2+1}^n
x_s\ptl_{y_s}+\sum_{r=1}^m\sta_r\vt_r.\eqno(13.4.29)$$ We take the
notation in (13.3.51) and set
$${\msr B}_{\la k\ra}=\sum_{i=0}^\infty {\msr A}_{\la
k-i\ra}x_0^i,\qquad {\msr H}'_{\la k\ra}=\{f\in {\msr B}_{ k}\mid
\Dlt_{n_1,n_2}'(f)=0\}.\eqno(13.4.30)$$ Then we have the second main
theorem in this section:\psp

{\bf Theorem 13.4.2}.  {\it For any $k\in\mbb{Z}$, ${\msr H}_{\la
k\ra}'$ is an irreducible $osp(2n+1|2m)$-module. Moreover,
$${\msr
B}=\bigoplus_{\ell,k=0}^\infty {\eta'}^\ell ({\msr H}_{\la
k\ra})\eqno(13.4.31)$$ is a direct sum of irreducible
$osp(2n+1|2m)$-submodules.}

{\it Proof}. We define $T_\iota$ as in (13.4.18) with $\Dlt$
replaced by $\Dlt_{n_1,n_2}$ in (13.4.28). By Lemma 6.1.1,
$${\msr H}'_{\la k\ra}=T_0({\msr A}_{\la k\ra})+T_1({\msr A}_{\la
k-1\ra})\qquad\for\;\;k\in\mbb{Z}.\eqno(13.4.32)$$ Since
$\Dlt_{n_1,n_2}\xi=\xi\Dlt_{n_1,n_2}$ for $\xi\in osp(2n|2m)$, we
have
$$\xi(T_\iota(f))=T_\iota(\xi(f))\qquad\for\;\;\xi\in osp(2n|2m),\;
f\in{\msr A}.\eqno(13.4.33)$$

First we consider ${\msr H}'_{\la k\ra}$ with $k\in\mbb{N}$. Let $V$
be any nonzero $osp(2n+1|2m)$-submodule of ${\msr H}'_{\la k\ra}$.
According to the arguments in paragraph of (13.3.79)-(13.3.84), $V$
contains some $T_\iota(\eta^\ell(x_{n_1+1}^{k-\iota-2\ell})).$ By
(13.4.24),
$$(E_{n_1+1,0}-E_{0,n+n_1+1})|_{\msr B}=x_{n_1+1}\ptl_{x_0}-x_0\ptl_{y_{n_1+1}}.\eqno(13.4.34)$$
Moreover, as operators on ${\msr B}$,
\begin{eqnarray*}\hspace{1cm}& &[E_{n_1+1,0}-E_{0,n+n_1+1},T_0]\\&=&
[x_{n_1+1}\ptl_{x_0}-x_0\ptl_{y_{n_1+1}},\sum_{i=0}^\infty\frac{(-2)^ix_0^{2i}}{(2i)!}\Dlt_{n_1,n_2}^i]\\
&=&[x_{n_1+1},\sum_{i=1}^\infty\frac{(-2)^ix_0^{2i}}{(2i)!}\Dlt_{n_1,n_2}^i]\ptl_{x_0}+x_{n_1+1}
\sum_{i=1}^\infty\frac{(-2)^ix_0^{2i-1}}{(2i-1)!}\Dlt_{n_1,n_2}^i\\&=&
-[\sum_{i=1}^\infty\frac{i(-2)^ix_0^{2i}}{(2i)!}\Dlt_{n_1,n_2}^{i-1}\ptl_{x_0}
+\sum_{i=1}^\infty\frac{i(-2)^ix_0^{2i-1}}{(2i-1)!}\Dlt_{n_1,n_2}^{i-1}]\ptl_{y_{n_1+1}}
\\ & &-2(T_1\Dlt_{n_1,n_2}) x_{n_1+1},\hspace{8.8cm}(13.4.35)
\end{eqnarray*}
\begin{eqnarray*}& &[E_{n_1+1,0}-E_{0,n+n_1+1},T_1]\\&=&
[x_{n_1+1}\ptl_{x_0}-x_0\ptl_{y_{n_1+1}},\sum_{i=0}^\infty\frac{(-2)^ix_0^{2i+1}}{(2i+1)!}\Dlt_{n_1,n_2}^i]\\
&=&[x_{n_1+1},\sum_{i=1}^\infty\frac{(-2)^ix_0^{2i+1}}{(2i+1)!}\Dlt_{n_1,n_2}^i]\ptl_{x_0}+x_{n_1+1}
\sum_{i=0}^\infty\frac{(-2)^ix_0^{2i}}{(2i)!}\Dlt_{n_1,n_2}^i\\&=&
-[\sum_{i=1}^\infty\frac{i(-2)^ix_0^{2i+1}}{(2i+1)!}\Dlt_{n_1,n_2}^{i-1}\ptl_{x_0}
+\sum_{i=1}^\infty\frac{i(-2)^ix_0^{2i}}{(2i)!}\Dlt_{n_1,n_2}^{i-1}]\ptl_{y_{n_1+1}}
+T_0 x_{n_1+1}.\hspace{1.4cm}(13.4.36)
\end{eqnarray*}
If $T_0(\eta^\ell(x_{n_1+1}^{k-2\ell}))\in V$ for some
$\ell\in\mbb{N}+1$, we have
\begin{eqnarray*}& &(E_{n_1+1,0}-E_{0,n+n_1+1})T_0(\eta^\ell(x_{n_1+1}^{k-2\ell}))\\
&=&
([E_{n_1+1,0}-E_{0,n+n_1+1},T_0]+T_0(E_{n_1+1,0}-E_{0,n+n_1+1}))(\eta^\ell(x_{n_1+1}^{k-2\ell}))
\\&=&-[\sum_{i=1}^\infty\frac{i(-2)^ix_0^{2i-1}}{(2i-1)!}\Dlt_{n_1,n_2}^{i-1}\ptl_{y_{n_1+1}}+x_0T_0\ptl_{y_{n_1+1}}
+2(T_1\Dlt_{n_1,n_2})x_{n_1+1}](\eta^\ell(x_{n_1+1}^{k-2\ell}))
\\&=& [T_1\ptl_{y_{n_1+1}}-2(T_1\Dlt_{n_1,n_2})x_{n_1+1}](\eta^\ell(x_{n_1+1}^{k-2\ell}))
\\ &=&\ell[1-2(m+n_1-n_2+\ell-1)](\eta^{\ell-1}(x_{n_1+1}^{k-2(\ell-1)-1}))
\\ &=&\ell[3-2(m+n_1-n_2+\ell)]T_1(\eta^{\ell-1}(x_{n_1+1}^{k-2(\ell-1)-1}))\in
V\hspace{4.1cm}(13.4.37)\end{eqnarray*} by (13.3.74), (13.4.34) and
(13.4.35). So $T_1(\eta^{\ell-1}(x_{n_1+1}^{k-2(\ell-1)-1}))\in V$.
When $T_1(\eta^\ell(x_{n_1+1}^{k-2\ell-1}))\in V$ for some
$\ell\in\mbb{N}$, (13.4.34) and (13.4.36) yield
$$(E_{n_1+1,0}-E_{0,n+n_1+1})T_1(\eta^\ell(x_{n_1+1}^{k-2\ell-1}))=T_0(\eta^\ell(x_{n_1+1}^{k-2\ell}))\in
V.\eqno(13.4.38)$$ By induction on $\ell$, we have
$x_{n_1+1}^k=T_0(x_{n_1+1}^k)\in V$.

Note
$$(E_{n+i,n+n_1+1}-E_{n_1+1,i})|_{\msr
B}=y_i\ptl_{y_{n_1+1}}+x_ix_{n_1+1}\qquad\for\;\;i\in\ol{1,n_1}\eqno(13.4.39)$$
and
$$(E_{n_2+r,n_2}-E_{n+n_2,n+m_2+r})|_{\msr
A}=x_{n_2+r}\ptl_{x_{n_2}}+y_{n_2}y_{n_2+r}\;\;\for\;\;r\in\ol{1,n-n_2}\eqno(13.4.40)$$
if $n_2<n$ by (13.3.37) and (13.3.38). Repeatedly applying (13.4.39)
and (13.4.40) to (13.4.38) with various $i\in\ol{1,n_1}$ and
$r\in\ol{1,n-n_2}$ if $n_2<n$, we have
$$[\prod_{i=1}^{n_1+1}x_i^{\al_i}][\prod_{j=n_2}^ny_j^{\be_j}]\in
V\qquad\for\;\;\al_i,\be_j\in\mbb{N};\al_{n_1+1}+\be_{n_2}-\sum_{i=1}^{n_1}\al_i
-\sum_{r=n_2+1}^n\be_r=k.\eqno(13.4.41)$$

Denote
$$I=\{0,\ol{n_1+1,n_2},\ol{n+n_1+1,n+n_2},\ol{2n+1,2n+2m}\}.\eqno(13.4.42)$$
Then the Lie subalgebra
$${\msr G}=osp(2n+1|2m)\bigcap(\sum_{i,j\in I}\mbb{F}E_{i,j})\cong
osp(2(n_2-n_1)+1|2m).\eqno(13.4.43)$$ Applying Theorem 13.4.1 to
${\msr G}$ and
$\mbb{F}[x_0,x_{n_1+1},...,x_{n_2},y_{n_1+1},...,y_{n_2},\sta_1,...,\sta_m,\vt_1,...,\vt_m]$,
we get
$$T_\iota(x^\al y^\be\sta_{\vec j}\vt_{\vec j'})\in V\eqno(13.4.44)$$
for $\al,\be\in\mbb{N}^n$, $\vec j\in\G_{k_1}$ and $\vec
j'\in\G_{k_2}$ such that $\be_i=0$ if $i\leq n_1$ and $\al_j=0$ if
$j>n_2$, and
$$\iota+k_1+k_2+\sum_{r=n_1+1}^{n_2}(\al_r+\be_r)
-\sum_{i=1}^{n_1}\al_i-\sum_{j=n_2+1}^n\be_j=k.\eqno(13.4.45)$$

Repeatedly applying (13.4.39) to (13.4.45) under above conditions
with various $i\in\ol{1,n_1}$, we obtain (13.4.45) for
$\al,\be\in\mbb{N}^n$, $\vec j\in\G_{k_1}$ and $\vec j'\in\G_{k_2}$
such that $\al_i=0$ if $i>n_2$, and
$$\iota+k_1+k_2+\sum_{r=n_1+1}^{n_2}\al_r+\sum_{s=1}^{n_2}\be_s-\sum_{i=1}^{n_1}\al_i-\sum_{j=n_2+1}^n\be_j=k.
\eqno(13.4.46)$$
 Observe
$$(E_{n_2+r,n_1+s}-E_{n+n_1+s,n+n_2+r})|_{\msr
B}=y_{n_1+s}y_{n_2+r}+x_{n_2+r}\ptl_{x_{n_1+s}}\eqno(13.4.47)$$ for
$r\in\ol{1,n-n_2}$ and $s\in\ol{1,n_2-n_1}$ by (13.3.37) and
(13.3.38). Repeatedly applying (13.4.47) to (13.4.45) with $\al_i=0$
if $i>n_2$, we obtain ${\msr H}_{\la k\ra}\subset V$ by (13.4.32).
So ${\msr H}_{\la k\ra}$ is an irreducible $osp(2n+1|2m)$-module.

Next we consider ${\msr H}_{\la -k\ra}$ with $k\in\mbb{N}+1$.  Let
$U$ be any nonzero $osp(2n+1|2m)$-submodule of ${\msr
H}'_{\la-k\ra}$. According to the arguments in paragraph of
(13.3.79)-(13.3.84), $U$ contains some
$T_\iota(\eta^\ell(x_{n_1}^{k+\iota+2\ell})).$ Observe Note
$$(E_{n_1,0}-E_{0,n+n_1})=\ptl_{x_0}\ptl_{x_{n_1}}-x_0\ptl_{y_{n_1}}\eqno(13.4.48)$$
by (13.3.37) and (13.3.38). As operators on ${\msr B}$,
\begin{eqnarray*}\hspace{2cm}& &[E_{n_1,0}-E_{0,n+n_1},T_0]\\&=&
[\ptl_{x_0}\ptl_{x_{n_1}}-x_0\ptl_{y_{n_1}}
,\sum_{i=0}^\infty\frac{(-2)^ix_0^{2i}}{(2i)!}\Dlt_{n_1,n_2}^i]\\
&=&\sum_{i=1}^\infty\frac{(-2)^ix_0^{2i-1}}{(2i-1)!}\Dlt_{n_1,n_2}^i\ptl_{x_{n_1}}
-\ptl_{x_0}\sum_{i=1}^\infty\frac{i(-2)^ix_0^{2i}}{(2i)!}\Dlt_{n_1,n_2}^{i-1}\ptl_{y_{n_1}}
\\&=&-2T_1\Dlt_{n_1,n_2}\ptl_{x_{n_1}}-
\sum_{i=1}^\infty\frac{i(-2)^ix_0^{2i-1}}{(2i-1)!}\Dlt_{n_1,n_2}^{i-1}\ptl_{y_{n_1}}\\
& &-
\sum_{i=1}^\infty\frac{i(-2)^ix_0^{2i}}{(2i)!}\Dlt_{n_1,n_2}^{i-1}\ptl_{y_{n_1}}\ptl_{x_0},
\hspace{6.1cm}(13.4.49)
\end{eqnarray*}
\begin{eqnarray*}& &[E_{n_1,0}-E_{0,n+n_1},T_1]\\&=&
[\ptl_{x_0}\ptl_{x_{n_1}}-x_0\ptl_{y_{n_1}},\sum_{i=0}^\infty\frac{(-2)^ix_0^{2i+1}}{(2i+1)!}\Dlt_{n_1,n_2}^i]\\
&=&\sum_{i=0}^\infty\frac{(-2)^ix_0^{2i}}{(2i)!}\Dlt_{n_1,n_2}^i\ptl_{x_{n_1}}-
\ptl_{x_0}\sum_{i=1}^\infty\frac{i(-2)^ix_0^{2i+1}}{(2i+1)!}\Dlt_{n_1,n_2}^{i-1}\ptl_{y_{n_1}}
\\ &=& T_0\ptl_{x_{n_1}}-\sum_{i=1}^\infty\frac{i(-2)^ix_0^{2i}}{(2i)!}\Dlt_{n_1,n_2}^{i-1}\ptl_{y_{n_1}}
-\sum_{i=1}^\infty\frac{i(-2)^ix_0^{2i+1}}{(2i+1)!}\Dlt_{n_1,n_2}^{i-1}\ptl_{y_{n_1}}\ptl_{x_0}
.\hspace{1.7cm}(13.4.50)
\end{eqnarray*}

If $T_0(\eta^\ell(x_{n_1}^{k+2\ell}))\in U$ for some
$\ell\in\mbb{N}+1$, we have
\begin{eqnarray*}& &(E_{n_1,0}-E_{0,n+n_1})T_0(\eta^\ell(x_{n_1+1}^{k-2\ell}))
\\&=& [T_1\ptl_{y_{n_1}}-2T_1\Dlt_{n_1,n_2}\ptl_{x_{n_1}}](\eta^\ell(x_{n_1+1}^{k+2\ell}))
\\ &=&(k+2\ell)\ell[1+2(n_1+k+\ell-n_2)]T_1(\eta^{\ell-1}(x_{n_1+1}^{k+2\ell-1}))
\in V\hspace{3.6cm}(13.4.51)\end{eqnarray*} by (13.3.74), (13.4.48)
and (13.4.49). So $T_1(\eta^{\ell-1}(x_{n_1}^{k+2(\ell-1)+1}))\in
U$. When $T_1(\eta^\ell(x_{n_1}^{k+2\ell+1}))\in V$ for some
$\ell\in\mbb{N}$, (13.4.48) and (13.4.50) yield
$$(E_{n_1,0}-E_{0,n+n_1})T_1(\eta^\ell(x_{n_1}^{k+2\ell+1}))=(k+2\ell+1)T_0(\eta^\ell(x_{n_1}^{k+2\ell}))\in
U.\eqno(13.4.52)$$ By induction on $\ell$, we have
$x_{n_1}^k=T_0(x_{n_1}^k)\in U$.

According to (13.4.39) with $i=n_1$,
$$x_{n_1}^{k+k'}x_{n_1+1}^{k'}\in
U\qquad\for\;\;k'\in\mbb{N}.\eqno(13.4.53)$$
 Moreover,
$$(E_{i,n_1}-E_{n+n_1,n+i})|_{\msr
B}=x_i\ptl_{x_{n_1}}-y_{n_1}\ptl_{y_i}\qquad\for
i\in\ol{1,n_1-1}\eqno(13.4.54)$$ by (13.3.37) and (13.3.38).
Repeatedly applying (13.4.54) to (13.4.53) with various
$i\in\ol{1,n_1-1}$, we have
$$\prod_{i=1}^{n_1+1}x_i^{\al_i}\in
U\qquad\for\;\;\al_i\in\mbb{N};\al_{n_1+1}-\sum_{i=1}^{n_1}\al_i=-k.\eqno(13.4.55)$$
Observe
$$(E_{n_2+r,n+1}-E_{1,n+n_2+r})|_{\msr
A}=x_{n_2+r}\ptl_{y_1}+y_{n_2+r}\ptl_{x_1}\;\;\for\;\;r\in\in\ol{1,n-n_2}\;\;
\mbox{if}\;\;n_2<n\eqno(13.4.56)$$ by (13.3.39). Repeatedly applying
(13.4.56) to (13.4.55) with various $r\in\ol{1,n-n_2}$ if $n_2<n$,
we find
$$[\prod_{i=1}^{n_1+1}x_i^{\al_i}][\prod_{j=n_2+1}^ny_j^{\be_j}]\in
U\qquad\for\;\;\al_i,\be_j\in\mbb{N};\al_{n_1+1}-\sum_{i=1}^{n_1}\al_i-\sum_{j=n_2+1}^n=-k.
\eqno(13.4.57)$$ By the same arguments from (13.4.43) to the end of
the paragraph below (13.4.47) with $k$ replaced by $-k$, we prove
that ${\msr H}_{\la-k\ra}$ is an irreducible $osp(2n+1|2m)$-module.

We calculate
$$[\Dlt_{n_1,n_2}',\eta']=2+4(n_2-n_1-m)+4[x_0\ptl_{x_0}+\sum_{i=1}^n(x_i\ptl_{x_i}
+y_i\ptl_{y_i})+\sum_{r=1}^m(\sta_r\ptl_{\sta_r}+\vt_r\ptl_{\vt_r})].\eqno(13.4.58)$$
By the arguments in (13.2.58)-(13.2.64),  ${\msr
B}=\bigoplus_{\ell,k=0}^\infty {\eta'}^\ell {\msr H}_{\la k\ra}$ is
a direct sum of irreducible $osp(2n+1|2m)$-submodules for any
$k\in\mbb{Z}.\qquad\Box$\psp

It can be verified  by (13.4.58) that the space
$\mbb{F}\Dlt_{n_1,n_2}'+\mbb{F}[\Dlt_{n_1,n_2}',\eta']+\mbb{F}\eta'$
in terms of (13.4.28) and (13.4.29) also forms an operator Lie
algebra isomorphic to $sl(2,\mbb{F})$. The above theorem establishes
a supersymmetric $(sl(2,\mbb{F}),osp(2n+1|2m))$
 Howe duality.

\chapter{Lie Theoretic Codes}

Linear codes with large minimal distances are important error
correcting codes in information theory. Orthogonal codes have more
applications in the other fields of mathematics. In this chapter, we
study the binary and ternary orthogonal codes generated by the
weight matrices on finite-dimensional  modules of
 the simple Lie algebras. The Weyl groups of the Lie algebras act on these codes isometrically.
  It turns out
 that certain weight matrices of $sl(n,\mbb{C})$ and $o(2n,\mbb{C})$ generate doubly-even  binary orthogonal codes
 and ternary orthogonal codes with large minimal distances. Moreover, we prove that the
 weight matrices of $F_4$, $E_6$, $E_7$
 and $E_8$ on their minimal irreducible modules and adjoint modules all generate ternary orthogonal codes
 with large minimal distances.
 In determining the minimal distances, we have used the Weyl groups and branch rules of
 the irreducible representations of the related simple Lie algebras.

\section{Basics Definitions}

Let $m$ be a positive integer and denote
$\mbb{Z}_m=\mbb{Z}/m\mbb{Z}$, which is a {\it ring}\index{ring} with
the algebraic operations induced from $\mbb Z$. A {\it code ${\msr
C}$ of length $n$} is a subset of $(\mbb{Z}_m)^n$ for some $m$,
where the ring structure of $\mbb{Z}_m$ may not be used. The
elements of ${\msr C}$ are called {\it codewords}.\index{codeword}
The {\it (Hamming) distance}\index{Hamming distance} between two
codewords is the number of different coordinates. The {\it minimal
distance} of a code is the minimal number among the distances of all
its pairs of codewords in the code. A code with minimal distance $d$
can be used to correct $\llbracket(d-1)/2\rrbracket$ errors in
signal transmissions.

A {\it linear code}\index{linear code} ${\msr C}$ over the ring
$\mbb{Z}_m$ is a $\mbb{Z}_m$-submodule of $(\mbb{Z}_m)^n$. The {\it
(Hamming) weight} of a codeword in a linear code ${\msr C}$ is the
number of its nonzero coordinates. In this case, the minimal
distance of ${\msr C}$ is exactly the minimal weight of the nonzero
codewords in ${\msr C}$. The inner product in $(\mbb{Z}_m)^n$ is
defined by
$$(a_1,...,a_n)\cdot(b_1,...,b_n)=\sum_{i=1}^na_ib_i.\eqno(14.1.1)$$
Moreover, ${\msr C}$ is called {\it orthogonal}\index{orthogonal
code} if
$${\msr C}\subseteq\{\vec a\in (\mbb{Z}_m)^n\mid\vec\al\cdot\vec b=0\;\for\;\vec b\in{\msr C}\}.
\eqno(14.1.2)$$ When the equality holds, we call ${\msr C}$ a {\it
self-dual} code.\index{self-dual code} Orthogonal linear codes
(especially, self-dual codes) have important applications to the
other mathematical fields such as sphere packing, integral linear
lattices, finite group theory, etc.  A code is called {\it
binary}\index{binary code} if $m=2$ and {\it ternary} when $m=3$.
\index{ternary code} A binary linear code is called  {\it even
(doubly-even)}\index{even (doubly-even) code} if the weights of all
its codewords are divisible by 2 (by 4).

Let ${\msr G}$ be a finite-dimensional simple Lie algebras over
$\mbb{C}$, the field of complex numbers. Take a Cartan subalgebra
$H$ and simple positive roots $\{\al_1,\al_2,...,\al_n\}$. Moreover,
we denote by $\{h_1,h_2,...,h_n\}$ the elements of $H$ such that the
matrix
$$(\al_i(h_j))_{n\times n}\;\;\mbox{is the Cartan matrix of}\;\;{\msr G}\eqno(14.1.3)$$
((2.5.13), (2.5.18), (3.1.8) and (3.3.2) ). For a finite-dimensional
${\msr G}$-module $V$,  $V$ has a weight-subspace decomposition:
$$V=\bigoplus_{\mu\in H^\ast}V_\mu,\qquad V_\mu=\{v\in V\mid h(v)=\mu(h)v\;\for\;h\in H\}.\eqno(14.1.4)$$
Take a maximal linearly independent set $\{u_1,u_2,...,u_k\}$ of
weight vectors with nonzero weights in $V$ such
 that the order is compatible with the partial order $\lhd$ of weights in (3.4.49). Write
 $$h_i(u_j)=c_{i,j}u_j,\qquad C(V)=(c_{i,j})_{n\times k}.\eqno(14.1.5)$$
 By the representation theory of simple Lie algebras, all $c_{i,j}$ are integers.
  We call $C(V)$ the {\it weight matrix
 of ${\msr G}$ on}\index{weight matrix} $V$. Identify  integers with their images in
 $\mbb{Z}_m$ when the context is clear. Denote by ${\msr C}_m(V)$ the linear code
 over $\mbb{Z}_m$ generated by
$C(V)$.

Suppose that the weight of $u_i$ is $\mu_i$. Set
 $$\bar H_m=\sum_{i=1}^n\mbb{Z}_mh_i.\eqno(14.1.6)$$
We define a map $\jmath: \bar H_m\rta (\mbb{Z}_m)^k$ by
$$\jmath(\sum_{i=1}^nl_ih_i)=(\sum_{i=1}^nl_i\mu_1(h_i),\sum_{i=1}^nl_i\mu_2(h_i),...,\sum_{i=1}^nl_i\mu_k(h_i)).\eqno(14.1.7)$$
Then
$${\msr C}_m(V)=\jmath(\bar H_m).\eqno(14.1.8)$$
Let ${\msr W}$ be the Weyl group of the simple Lie algebra ${\msr
G}$. For any $\sgm\in {\msr W}$, there exists a linear automorphism
$\hat\sgm$ of $V$ such that
$$\hat\sgm (V_{\mu})=V_{\sgm(\mu)},\qquad\sgm(\mu)(\sgm(h))=\mu(h)\qquad\for\;\;h\in H\eqno(14.1.9)$$
(cf. (5.2.41)), where $\sgm(t_\al)=t_{\sgm(\al)}$ in terms of
(2.5.13). Moreover, we define an action of ${\msr W}$ on $\bar H_m$
by
$$\sgm(\sum_{i=1}^nl_ih_i)=\sum_{i=1}^nl_i\sgm(h_i)\qquad\for\;\;\sgm\in {\msr W}.\eqno(14.1.10)$$
According to (14.1.9),
$$\wt\jmath(\sgm(h))=\wt\jmath(h)\qquad\for\;\;\sgm \in {\msr W}({\msr G}),\;h\in \bar H_m.
\eqno(14.1.11)$$ So the number of the distinct weights of codewords
in ${\msr C}_m(V)$ is less than or equal to the number of ${\msr
W}$-orbits in $\bar H_m$. Expression (14.1.11) will be used later in
determining minimal distances.

Let $\Lmd(V)$ be the set of nonzero weights of $V$. The module $V$
is called {\it self-dual}\index{self-dual module} if
$\Lmd(V)=-\Lmd(V)$. In this chapter, we are only interested in the
binary and ternary codes. We call ${\msr C}_2(V)$ the {\it binary
weight code of} ${\msr G}$ {\it on} $V$.\index{binary weight code}
If $V$ is self-dual, then the weight matrix $C(V)=(-A,A)$ and ${\msr
C}_3(V)$ is orthogonal if and only if $A$ generates a ternary
orthogonal code (e.g., cf. [P4]). For this reason, we call the
ternary code generated by $A$ the {\it ternary weight code of}
${\msr G}$ {\it on} $V$ if $V$ is self-dual. When $V$ is not
self-dual, then ${\msr C}_3(V)$ is the {\it ternary weight code of}
${\msr G}$ {\it on} $V$. \index{ternary weight code}

Denote by $V_X(\lmd)$ the finite-dimensional irreducible module of a
simple Lie algebra of type $X$ with the highest weight $\lmd$. Let
$p$ be a prime number. Then $\mbb{Z}_p$ is a finite field, which is
traditionally denoted as $\mbb{F}_p$. A linear code ${\msr C}$ of
length $n$ over $\mbb{F}_p$ is a vector subspace of $\mbb{F}_p^n$
over $\mbb{F}_p$. If $\dim {\msr C}=k$, we say that ${\msr C}$ is of
{\it type} $[n,k]$. When $d$ is the minimal distance of ${\msr C}$,
 we call ${\msr C}$ an {\it $[n,k,d]$-code.}

\section{Codes Related to Representations of $sl(n,\mbb C)$}

In this section, we study the binary and ternary codes related to
representations of $sl(n,\mbb C)$, where $n>1$ is an integer. Again
we denote
$$\ves_i=(0,...,\stl{i}{1},0,...,0)\in\mbb{R}^n.\eqno(14.2.1)$$
So
$$\mbb{R}^n=\sum_{i=1}^n\mbb{R}\ves_i.\eqno(14.2.2)$$
Then inner product ``$(\cdot,\cdot)$" is Euclidian; that is,
$$(\sum_{i=1}^nk_i\ves_i,\sum_{j=1}^nl_j\ves_j)=\sum_{i=1}^nk_il_i.\eqno(14.2.3)$$
Recall the special linear Lie algebra $sl(n,\mbb C)$ and related
settings in (6.2.1)-(6.2.4). Its root system
$$\Phi_{A_{n-1}}=\{\ves_i-\ves_j\mid i,j\in\ol{1,n},\;i\neq
j\}.\eqno(14.2.4)$$ Take the simple positive roots
$$\al_i=\ves_i-\ves_{i+1}\qquad\for\;\;i\in\ol{1,n-1}.\eqno(14.2.5)$$
 The Weyl group ${\msr W}_{A_{n-1}}$ of $sl(n,\mbb C)$ is
exactly the full permutation group $S_n$ on $\ol{1,n}$, which acts
on $H_{A_{n-1}}$ and $\mbb{R}^n$ by permuting sub-indices of
$E_{i,i}$ and $\ves_i$, respectively.

Let $\Psi$ be the exterior algebra generated by
$\{\sta_1,\sta_2,...,\sta_n\}$ (cf. (6.2.15)) and take the related
settings (6.2.16)-(6.2.19). Proposition 6.2.2 says that $\Psi_r$
forms an irreducible $sl(n,\mbb C)$-submodule of highest weight
$\lmd_r$ for $r\in\ol{1,n-1}$. The Weyl group ${\msr W}_{A_{n-1}}$
acts on ${\msr A}$ by permuting sub-indices of $\sta_i$.

Two $k_1\times k_2$ matrices $A_1$ and $A_2$ with entries in
$\mbb{Z}_m$ are called {\it equivalent} in the sense of coding
theory if there exist an invertible $k_1\times k_1$ matrix $K_1$ and
an invertible $k_2\times k_2$ monomial matrix $K_2$ such that
$A_1=K_1A_2K_2.$ Equivalent matrices generate isomorphic codes.
 Take any order of the basis
$$\{x_{r,1},x_{r,2},...,x_{r,{n\choose r}}\}=\{
\sta_{i_1}\sta_{i_2}\cdots\sta_{i_r}\mid 1\leq i_1<i_2<\cdots<
i_r\leq n\}\eqno(14.2.6)$$ of $\Psi_r$. Then we have
$$h_i(x_{r,j})=a_{i,j}(r)x_{r,j},\qquad
a_{i,j}(r)\in\mbb{Z}.\eqno(14.2.7)$$ Modulo equivalence, the weight
matrix
$$C(\Psi_r)=[a_{i,j}(r)]_{(n-1)\times {n\choose
r}}.\eqno(14.2.8)$$

{\bf Theorem 14.2.1}. {\it When $n=2m\geq 4$ is even, ${\msr
C}_2(\Psi_2)$ is a doubly-even binary orthogonal
$[m(2m-1),2(m-1),4(m-1)]$-code}.

{\it Proof}. Denote by $\xi_i$ the $i$th row $C_2(\Psi_2)$. Then
$$\wt\xi_i=2(n-2)\qquad\for\;\;i\in\ol{1,n-1}.\eqno(14.2.9)$$
Moreover,
$$\sum_{i=0}^{m-1}\xi_{2i+1}=0\qquad\mbox{in}\;\;{\msr C}_2(\Psi_2). \eqno(14.2.10)$$ Furthermore,
$$\xi_i\cdot\xi_j=4\equiv 0\qquad \mbox{if}\;\;i+1<j\eqno(14.2.11)$$
and
$$\xi_i\cdot\xi_{i+1}=2(m-1)\equiv 0.\eqno(14.2.12)$$

Write
$$E_{i,i}(x_{r,j})=b_{i,j}(r)x_{r,j},\qquad
B_r=[b_{i,j}(r)]_{n\times {n\choose r}}.\eqno(14.2.13)$$ Denote by
$\zeta_i$ the $i$th row of $B_2$. By symmetry (cf.
(14.1.9)-(14.1.11)), any nonzero codeword in ${\msr C}_2(\Psi_2)$
has the same weight as the codeword
$$u=\sum_{s=1}^{2t}\zeta_s\in\mbb{F}_2^{n(n-1)/2}\qquad\mbox{for
some}\;\;t\in\ol{1,m-1}.\eqno(14.2.14)$$ We calculate
$$\wt u=4t(m-t)=-4t^2+4mt.\eqno(14.2.15)$$ Since the function $-4t^2+t(4m-1)$ attains
maximal at $t=m/2$, $\wt u$ is minimal at $t=1$ or $m-1$. Note
$$\wt u=4(m-1)\qquad\mbox{if}\;\;t=1\;\mbox{or}\;m-1.\eqno(14.2.16)$$
Thus the code ${\msr C}_2(\Psi_2)$ has the minimal distance
$4(m-1).\qquad\Box$\psp

When $m=2$, ${\msr C}_2(\Psi_2)$ is a doubly-even binary orthogonal
$[6,2,4]$-code. If $m=3$,
 ${\msr C}_2(\Psi_2)$ becomes a doubly-even binary orthogonal $[15,4,8]$-code.
These two code are optimal linear
 codes (e.g., cf. [PHB]).
In  the case of $m=4$,  ${\msr C}_2(\Psi_2)$ is a doubly-even binary
orthogonal $[28,6,12]$-code.\psp

{\bf Theorem 2.2}. {\it The code ${\msr C}_2(\Psi_3)$ is a
doubly-even binary orthogonal $[{n\choose 3}, n-1,(n-2)(n-3)]$-code
if $n>9$ and $n\equiv 2,3\;(\mbox{mod}\;4)$.}

{\it Proof}. Denote by $\xi_i$ the $i$th row the weight matrix
$C(\Psi_3)$. Then
$$\wt\xi_i=(n-2)(n-3)\qquad\for\;\;i\in\ol{1,n-1}.\eqno(14.2.17)$$
Moreover,
$$\xi_i\cdot\xi_j=4(n-4)\qquad \mbox{if}\;\;i+1<j\eqno(14.2.18)$$
and
$$\xi_i\cdot\xi_{i+1}=n-3+{n-3\choose 2}=\frac{(n-2)(n-3)}{2}.\eqno(14.2.19)$$
So ${\msr C}_2(\Psi_3)$ is a doubly-even binary orthogonal code
under the assumption.

Denote by $\zeta_i$ the $i$th row of $B_3$ (cf. (14.2.13)). By
symmetry (cf. (14.1.9)-(14.1.11)), any nonzero codeword in ${\msr
C}_2(\Psi_3)$ has the same weight as the codeword
$$u(t)=\sum_{s=1}^{2t}\zeta_s\in\mbb{F}_2^n\qquad\mbox{for
some}\;\;t\in\ol{1,\llbracket n/2\rrbracket}.\eqno(14.2.20)$$ We
calculate
$$f(t)=3\wt u(t)=3{2t\choose 3}+6t{n-2t\choose
2}=t[16t^2-12nt+3n(n-1)+2].\eqno(14.2.21)$$ Moreover,
$$f'(t)=48t^2-24nt+3n(n-1)+2=48\left(t-\frac{n}{4}\right)^2-3n+2.\eqno(14.2.22)$$
Thus
$$f'(t_0)=0\lra
t_0=\frac{n}{4}\pm\frac{1}{4}\sqrt{n-\frac{2}{3}}.\eqno(14.2.23)$$
Since $f'(0)=3n(n-1)+2>0$, $f(t)$ attains local maximum at
$$t=\frac{n}{4}-\frac{1}{4}\sqrt{n-\frac{2}{3}}\eqno(14.2.24)$$
and local minimum at
$$t=\frac{n}{4}+\frac{1}{4}\sqrt{n-\frac{2}{3}}.\eqno(14.2.25)$$
According to (14.2.17) and (14.2.21), $f(1)=3(n-2)(n-3)$.
Furthermore,
\begin{eqnarray*}\qquad\qquad& &f\left(\frac{n}{4}+\frac{1}{4}\sqrt{n-\frac{2}{3}}\right)
\\ &=&\left(\frac{n}{4}+\frac{1}{4}\sqrt{n-\frac{2}{3}}\right)
\big[16\left(\frac{n}{4}+\frac{1}{4}\sqrt{n-\frac{2}{3}}\right)^2
\\& &-12n\left(\frac{n}{4}+\frac{1}{4}\sqrt{n-\frac{2}{3}}\right)
+3n(n-1)+2\big]\\&=&\frac{1}{4}\left(n+\sqrt{n-\frac{2}{3}}\right)
\big[\left(n+\sqrt{n-\frac{2}{3}}\right)^2\\ & &-3n\left(n+
\sqrt{n-\frac{2}{3}}\right)
+3n(n-1)+2\big]\\&=&\frac{1}{4}\left(n+\sqrt{n-\frac{2}{3}}\right)
\left[n\left(n-\sqrt{n-\frac{2}{3}}\right)-2n+\frac{4}{3}\right]\\
&=&\frac{1}{4}\left[n^3-3n^2+2n+\left(\frac{4}{3}-2n\right)
\sqrt{n-\frac{2}{3}}\right]\\
&>&\frac{1}{4}(n^3-5n^2+2n).\hspace{8.2cm}(14.2.26)\end{eqnarray*}
Thus
\begin{eqnarray*}& &
f\left(\frac{n}{4}+\frac{1}{4}\sqrt{n-\frac{2}{3}}\right)-f(1)\\
&>&\frac{1}{4}(n^3-5n^2+2n)-3(n-2)(n-3)=\frac{1}{4}(n^3-17n^2+62n-72)
\\ &>&\frac{n^2(n-17)}{4}.\hspace{11.2cm}(14.2.27)\end{eqnarray*}
If $n\geq 17$, we have
$$f\left(\frac{n}{4}+\frac{1}{4}\sqrt{n-\frac{2}{3}}\right)>f(1)
\eqno(14.2.28)$$ and \begin{eqnarray*}\hspace{2cm}&
&f(n/2)-f(1)\\&=& \frac{n}{2}[4n^2-6n^2+3n(n-1)+2]-3(n-2)(n-3)
\\ &=&\frac{n(n-1)(n-2)}{2}-3(n-2)(n-3)\\&=&\frac{(n-2)(n^2-7n+9)}{2}>0\;\;\mbox{if}\;\;n\geq 6.
\hspace{5.1cm}(14.2.29)\end{eqnarray*} Thus the minimal weight is
$f(1)/3=(n-2)(n-3)$ when $n\geq 17$.

When $n=10$, we calculate
\begin{center}{\bf Table 14.2.1}\end{center}
\begin{center}\begin{tabular}{|c|c|c|c|c|c|}\hline
$t$&1&2&3&4&5
\\\hline\wt u(t)&56&64&56&64&120\\\hline\end{tabular}\end{center}
If $n=11$, we find
\begin{center}{\bf Table 14.2.2}\end{center}
\begin{center}\begin{tabular}{|c|c|c|c|c|c|}\hline
$t$&1&2&3&4&5
\\\hline\wt u(t)&72&88&80&80&120\\\hline\end{tabular}\end{center}
When $n=14$, we obtain
\begin{center}{\bf Table 14.2.3}\end{center}
\begin{center}\begin{tabular}{|c|c|c|c|c|c|c|c|}\hline
$t$&1&2&3&4&5&6&7
\\\hline\wt u(t)&132&184&188&176&180&232&364\\\hline\end{tabular}\end{center}
If $n=15$, we get
\begin{center}{\bf Table 14.2.4}\end{center}
\begin{center}\begin{tabular}{|c|c|c|c|c|c|c|c|}\hline
$t$&1&2&3&4&5&6&7
\\\hline\wt u(t)&156&224&216&224&220&256&364\\\hline\end{tabular}\end{center}
This prove the conclusion in the theorem.$\qquad\Box$\psp

Note that when $n=6$, we find \begin{center}{\bf Table
14.2.5}\end{center}
\begin{center}\begin{tabular}{|c|c|c|c|}\hline
$t$&1&2&3
\\\hline\wt u(t)&12&8&20\\\hline\end{tabular}\end{center}
So ${\msr C}_2(\Psi_3)$ is a doubly-even binary orthogonal
$[20,5,8]$-code. Moreover, if $n=7$, we find
\begin{center}{\bf Table 14.2.6}\end{center}
\begin{center}\begin{tabular}{|c|c|c|c|}\hline
$t$&1&2&3
\\\hline\wt u(t)&20&16&20\\\hline\end{tabular}\end{center}
 Hence ${\msr C}_2(\Psi_3)$ a doubly-even binary orthogonal
$[35,6,16]$-code. In both cases, the above theorem fails and both
codes are the best even codes among the binary codes with the same
length and dimension (e.g., cf. [PHB]).

  According
to the above theorem, ${\msr C}_2(\Psi_3)$ is a doubly-even binary
orthogonal $[120,9,56]$-code when $n=10$, $[165,10,72]$-code if
$n=11$, $[364,13,132]$-code when $n=14$ and
 $[455,14,156]$-code if $n=15$.

Next let us consider the ternary codes.  Again by symmetry, any
nonzero codeword in ${\msr C}_3(\Psi_r)$ has the same weight as the
codeword
$$u(s,t)=\sum_{r=1}^s\zeta_r-\sum_{i=1}^t\zeta_{s+i}\in\mbb{F}_3^{{n\choose r}}\eqno(14.2.30)$$
for some nonnegative integers $s,t$, where $\zeta_\iota$ is the
$\iota$th row of the matrix $B_r$ in (14.2.13).
 Moreover,
$$\wt u(s,t)=\wt u(t,s).\eqno(14.2.31)$$
Furthermore, we have
$$\wt u(s,t)=(s+t)(n-s-t)+{s\choose 2}+{t\choose 2}\qquad\mbox{in}\;\;
{\msr C}_3(\Psi_2)\eqno(14.2.32)$$ and
$$\wt u(s,t)=(s+t){n-s-t\choose 2}+(n-s){s\choose 2}+(n-t){t\choose 2}\qquad\mbox{in}\;\;
{\msr C}_3(\Psi_3).\eqno(14.2.33)$$

For convenience, we denote \begin{eqnarray*}
\hspace{1cm}f(s,t)&=&2\wt u(s,t)=2(s+t)(n-s-t)+s(s-1)+t(t-1)\\
&=&(2n-1)(s+t)-s^2-t^2-4st\hspace{5.9cm}(14.2.34)\end{eqnarray*} in
${\msr C}_3(\Psi_2)$ and
 \begin{eqnarray*}\qquad\qquad g(s,t)&=&2\wt u(s,t)\\ &=&
(s+t)(n-s-t)(n-s-t-1)\\ & &+(n-s)s(s-1)+(n-t)t(t-1)
\\
&=&(s+t)^3-(2n-1)(s+t)^2+n(n-1)(s+t)\\ & &-s^3-t^3
+(n+1)(s^2+t^2)-n(s+t)
\\&=&3st^2+3s^2t+(2-n)(s^2+t^2)\\& &-2(2n-1)st+n(n-2)(s+t)\hspace{5.3cm}(14.2.35)\end{eqnarray*}
in ${\msr C}_3(\Psi_3)$.

Note
$$f(3,0)=3(2n-1)-9=6(n-2),\;\;f(n,0)=n(2n-1)-n^2=n(n-1),\eqno(14.2.36)$$
$$f(1,1)=2(2n-1)-6=4(n-2),\;\;f(1,n-1)=(n-1)(n-2).\eqno(14.2.37)$$
Since geometrically $f(s,t)$ has only local minimum, it attains the
absolute minimum at boundary points. Thus
$$\min \{f(s,t)\mid s\equiv t\;(\mbox{mod}\;3)\}=4(n-2)\qquad\mbox{if}\;\;n\geq 5.\eqno(14.2.38)$$

Now
$$g_s(s,t)=3t^2+6st+2(2-n)s-2(2n-1)t+n(n-2),\eqno(14.2.39)$$
$$g_t(s,t)=3s^2+6st+2(2-n)t-2(2n-1)s+n(n-2). \eqno(14.2.40)$$ Suppose that
$g_s(s_0,t_0)=g_t(s_0,t_0)=0$ for $s_0,t_0\geq 0$; that is,
$$3t_0^2+6s_0t_0+2(2-n)s_0-2(2n-1)t_0+n(n-2)=0,\eqno(14.2.41)$$
$$3s_0^2+6s_0t_0+2(2-n)t_0-2(2n-1)s_0+n(n-2)=0.\eqno(14.2.42)$$
By $(14.2.41)-(14.2.42)$, we get
$$(t_0-s_0)(3t_0+3s_0-2(n+1))=0\lra t_0=s_0\;\;\mbox{or}\;\;
3t_0+3s_0=2(n+1).\eqno(14.2.43)$$

If $s_0=t_0$, then
 we find
$$9s_0^2-2(n-1)s_0+n(n-2)=0\sim
8s_0^2+(s_0-n+1)^2-1=0,\eqno(14.2.44)$$ which leads to a
contradiction because $n>1$. Thus $3t_0+3s_0=2(n-1)$. Denote
$s_1=3t_0$ and $t_1=3t_0$. Then $s_1+t_1=2(n+1)$ and (14.2.41)
becomes
\begin{eqnarray*}\qquad\quad& &t_1^2+2(2(n+1)-t_1)t_1+2(2-n)(2(n+1)-t_1)\\ & &-2(2n-1)t_1+3n(n-2)=0,
\hspace{6.9cm}(14.2.45)\end{eqnarray*}or equivalently,
\begin{eqnarray*}\qquad\quad& &t_1^2-2(n+1)t_1+(n-2)(n+4)=0\\ & &\sim (t_1-n-1)^2-9=0 \lra
t_1=n+4,\;n-2.\hspace{4.1cm}(14.2.46)\end{eqnarray*} Therefore,
$$s_0=\frac{n+4}{3},\;\;t_0=\frac{n-2}{3}\qquad\mbox{or}\qquad
t_0=\frac{n+4}{3},\;\;s_0=\frac{n-2}{3}.\eqno(14.2.47)$$ We
calculate
$$g(s_0,t_0)=\frac{2(n-2)(n^2-n-3)}{9},\eqno(14.2.48)$$
$$g(1,0)=g(n-1,0)=(n-1)(n-2),\eqno(14.2.49)$$
$$g(3,0)=3(n-2)(n-3),\;\;g(n,0)=0,\eqno(14.2.50)$$
$$g(1,1)=g(n-2,1)=2(n-2)(n-3),\;\;g(n-2,0)=2(n-2)^2.\eqno(14.2.51)$$
Moreover,
$$g(s_0,t_0)\geq g(1,0),\;g(1,1)\qquad\mbox{if}\;\;n\geq 6.\eqno(14.2.52)$$
When $n=5$, we calculate
$$g(1,0)=g(1,1)=g(2,1)=g(2,2)=g(3,1)=g(4,0)=g(4,1)=12,\eqno(14.2.53)$$
$$g(2,0)=g(3,0)=g(3,2)=18.\eqno(14.2.54)$$
In summary, we have:\psp

{\bf Theorem 14.2.3}. {\it Let $n\geq 5$. The  matrix $B_3$ (cf.
(14.2.13)) generates a ternary $\left[{n\choose 3},n-1,{n-1\choose
2}\right]$-code, which is equal to ${\msr C}_3(\Psi_3)$ if
$n\not\equiv 0\;(\mbox{\it mod}\;3)$. If $n=3m+2$ for some positive
integer $m$, ${\msr C}_3(\Psi_2)$ is a ternary orthogonal
$[{3m+2\choose 2},3m+1,6m]$-code and ${\msr C}_3(\Psi_3)$ is a
ternary orthogonal $[{3m+2\choose 3},3m+1,3m(3m+1)/2]$-code.  The
code ${\msr C}_3(\Psi_3)$ is a ternary orthogonal $\left[{n\choose
3},n-2,(n-2)(n-3)\right]$-code when $n\equiv 0\;(\mbox{\it
mod}\;3)$. }

{\it Proof}. The part of minimal distances has been proved by the
above arguments. We only need to prove orthogonality.

Suppose $n=3m+2$. In ${\msr C}_3(\Psi_2)$, $\xi_r$ stands for the
$r$th row of  the weight matrix $C(\Psi_2)$ and
$$\xi_i\cdot \xi_j=2-2=0\qquad\for\;\;1\leq i<j-1\leq
n-2,\eqno(14.2.55)$$
$$\xi_r\cdot
\xi_{r+1}=-(n-2)=-3m,\;\;\xi_s\cdot\xi_s=2(n-2)=6m\eqno(14.2.56)$$
for $r\in\ol{1,n-2}$ and $s\in\ol{1,n-1}$. So ${\msr C}_3(\Psi_2)$
is orthogonal. Now $\zeta_r$ stands for the $r$th row of $B_3$ (cf.
(14.2.18)). Observe
$$\sum_{i=1}^n\zeta_i=0\in \mbb C_3^{{n\choose 3}}\eqno(14.2.57)$$
by (6.2.16) and (6.2.19). Moreover,
$$\zeta_i\cdot\zeta_j=n-2=3m,\;\;\zeta_i\cdot\zeta_i={n-1\choose
2}=\frac{3m(3m+1)}{2},\qquad i\neq j.\eqno(14.2.58)$$ Thus $B_3$
generate a ternary  orthogonal code.

Assume that $n=3m$ for some nonnegative integer $m$. In ${\msr
C}_3(\Psi_3)$, we also use $\xi_r$ for the $r$th row of the weight
code $C(\Psi_3)$ and
$$\xi_i\cdot\xi_j=2(n-4)-2(n-4)=0\qquad\for\;\;1\leq i<j-1\leq
n-2,\eqno(14.2.59)$$
$$\xi_s\cdot\xi_s=2\xi_r\cdot\xi_{r+1}=(n-2)(n-3)=3(3m-2)(m-1)\equiv
0\eqno(14.2.60)$$ for $r\in\ol{1,n-2}$ and $s\in\ol{1,n-1}$. So
${\msr C}_3(\Psi_3)$ is orthogonal.$\qquad\Box$ \psp

According to the above theorem, ${\msr C}_3(\Psi_2)$ is a ternary
orthogonal $[10,4,6]$-code when $n=5$ (which is optimal (e.g., cf.
[PHB])), $[28,7,12]$-code when $n=8$, and $[55,10,18]$-code when
$n=11$. Moreover, ${\msr C}_3(\Psi_3)$ is a ternary orthogonal
$[10,4,6]$-code when $n=5$, $[15,4,12]$-code if $n=6$,
$[56,7,21]$-code when $n=8$, $[84,7,42]$-code if $n=9$,
$[165,10,45]$-code when $n=11$ and $[220,10,90]$-code when $n=12$.

Finally, we consider the adjoint representation of $sl(n,\mbb C)$.
Note that $\{E_{i,j}\mid 1\leq i< j\leq n\}$ are positive root
vectors. Given an order
$$\{y_1,...,y_{{n\choose 2}}\}=\{E_{i,j}\mid 1\leq i< j\leq
n\},\eqno(14.2.61)$$ we write
$$[h_i,y_j]=k_{i,j}y_j,\qquad
[E_{r,r},y_j]=l_{r,j}y_j.\eqno(14.2.62)$$  Denote
$$K=(k_{i,j})_{(n-1)\times {n\choose r}},
\qquad L=(l_{i,j})_{n\times {n\choose r}}.\eqno(14.2.63)$$ Let
${\msr C}_K$ be the ternary code generated by $K$ and let ${\msr
C}_L$ be the ternary code generated by $L$. Moreover, $\vec k_i$
stands for the $i$th row of $K$ and $\vec l_r$ stands for the $r$th
row of $L$. Set
$$u(s,t)=\sum_{i=1}^s\vec l_i-\sum_{j=1}^t\vec
l_{s+j}.\eqno(14.2.64)$$ For any nonzero codeword $v\in{\msr C}_L$,
using negative root vectors, we can prove
$$\wt (v,-v)=\wt (u(s,t),-u(s,t))\eqno(14.2.65)$$
for some $s$ and $t$ by symmetry (cf. (14.1.9)-(14.1.11)). Thus
$$\wt v=\wt u(s,t)=(s+t)(n-s-t)+st=\phi(s,t).\eqno(14.2.66)$$
Note
$$\phi(s,t)=n^2-\frac{1}{2}[(s-n)^2+(t-n)^2+(s-t)^2].\eqno(14.2.67)$$
So $\phi(s,t)$ has only local maximum. Thus it attains the absolute
minimum at the boundary points. We calculate
$$\phi(1,0)=\phi(n-1,0)=n-1,\;\;\phi(n-3,0)=3(n-3),\eqno(14.2.68)$$
$$\phi(1,1)=2n-3,\qquad\phi(n-2,1)=2(n-1).\eqno(14.2.69)$$
Since
$$\sum_{i=1}^n\vec l_i=0,\eqno(14.2.70)$$
$${\msr C}_K={\msr C}_L\qquad\mbox{if}\;\;n\neq
0\;(\mbox{mod}\;3).\eqno(14.2.71)$$
$$\vec k_i\cdot\vec k_j=2-2=0\qquad1\leq i<j-1\leq n,\eqno(14.2.72)$$
$$\vec k_r\cdot\vec k_{r+1}=6-n,\;\;\vec k_s\cdot\vec
k_s=2n-3.\eqno(14.2.73)$$ In summary, we have:\psp

{\bf Theorem 14.2.4}. {\it The code ${\msr C}_L$ is a ternary
$[{n\choose 2},n-1,n-1]$-code if $n\geq 4$, which is also the
ternary weight code on the adjoint module $sl(n,\mbb C)$ when $n\neq
0\;(\mbox{mod}\;3)$. If $n=3m$ for some integer $m>1$, then the
ternary weight code ${\msr C}_K$ on $sl(3m,\mbb C)$ is an orthogonal
$[{3m\choose 2}, 3m-2, 3(m-1)]$-code.}

\section{Codes Related to Representations of $o(2m,\mbb C)$}

In this section, we only study ternary codes related to certain
representations of $so(2m,\mbb C)$, some of which will be used to
investigate the codes related to exceptional  simple Lie algebras.

 Let $n=2m$ be a positive even integer. Take the settings in (14.2.1)-(14.2.3) (with $n\rta m$).  The
orthogonal Lie algebra
\begin{eqnarray*}o(2m,\mbb C)&=&\sum_{1\leq i<j\leq
m}[\mbb C(E_{i,j}-E_{m+j,m+i})+\mbb C(E_{j,i}-E_{m+i,m+j}) +\mbb
C(E_{i,m+j}-E_{j,m+i})\\ & &+\mbb C(E_{m+i,j}-E_{m+j,i})]
+\sum_{r=1}^m\mbb Ch_r,\hspace{6.1cm}(14.3.1)\end{eqnarray*} where
$$h_s=E_{s,s}-E_{s+1,s+1}-E_{m+s,m+s}+E_{m+s+1,m+s+1}\qquad\for\;\;s
\in\ol{1,m-1}\eqno(14.3.2)$$ and
$$h_m=E_{m-1,m-1}+E_{m,m}-E_{2m-1,2m-1}-E_{2m,2m}.\eqno(14.3.3)$$
Indeed, we take the Cartan subalgebra
$$ H_{D_m}=\sum_{i=1}^m\mbb Ch_i\eqno(14.3.4)$$
of $o(2m,\mbb C)$.  The root system
$$\Phi_{D_m}=\{\pm\ves_i\pm\ves_j\mid i,j\in\ol{1,m},\;i\neq
j\}\eqno(14.3.5)$$ and simple positive roots are:
$$\al_i=\ves_i-\ves_{i+1},\;\;\al_m=\ves_{m-1}+\ves_m,\qquad
i\in\ol{1,m-1}.\eqno(14.3.6)$$ The Weyl group is
$S_m\ltimes\mbb{Z}_2^{m-1}$, which acts $H_{D_m}$ and $\mbb{R}^m$ by
permuting sub-indices of $\ves_i$ and $E_{i,i}-E_{m+i,m+i}$, and
changing sign on even number of their coefficients.

Take the settings in (6.2.15)-(6.2.18) and (14.2.6)-(14.2.8).
Moreover, the representation of
 $o(2m,\mbb C)$ on $\Psi$ determined by (6.2.19).
 For any $\vec\iota=(\iota_1,...,\iota_m)$ with $\iota_i\in\{0,1\}$ and $\tau\in S_m$, we have an associative algebra
automorphism $\sgm_{\tau,\vec\iota}$ of $\Psi$ determined by
$$\sgm_{\tau,\vec\iota}(\sta_i)=
\sta_{m\dlt_{\iota_i,1}+\tau(i)},\;\;
\sgm_{\tau,\vec\iota}(\sta_{m+i})=\sta_{m\dlt_{\iota_i,0}+\tau(i)}
\qquad\for\;\;i\in\ol{1,m}.\eqno(14.3.7)$$ Moreover, we define a
linear map $\sgm_{\tau,\vec\iota}$ on ${\msr H}$ by
$$\sgm_{\tau,\vec\iota}(E_{i,i}-E_{m+i,m+i})=(-1)^{\iota_i}(E_{\tau(i),\tau(i)}-E_{m+\tau(i),m+\tau(i)})\qquad\for
\;\;i\in\ol{1,m}.\eqno(14.3.8)$$ Then
$$\sgm_{\tau,\vec\iota}(h(w))=\sgm_{\tau,\vec\iota}(h)[\sgm_{\tau,\vec\iota}(w)]\qquad\for\;\;h\in{\msr H},\;
w\in\Psi.\eqno(14.3.9)$$

Note that all $\Psi_r\cong V_{D_m}(\lmd_r)$ are self-dual $o(2m,\mbb
C)$-submodules for $r\in\ol{1,m-2}$ by Theorem 7.1.2. In particular,
the ternary weight code ${\msr C}_2$ of $o(2m,\mbb C)$ on $\Psi_2$
is given by the weight matrix on its subspace
 $$\Psi_{2,1}=\sum_{1\leq i<j\leq m}(\mbb C\sta_i\sta_j+\mbb C\sta_i\sta_{m+j}).\eqno(14.3.10)$$
We take any order
$$\{x_1,x_2,\cdots,x_{m(m-1)}\}=\{\sta_i\sta_j,\sta_i\sta_{m+j}\mid 1\leq i<j\leq m\}\eqno(14.3.11)$$
and write
$$(E_{i,i}-E_{m+i,m+i})(x_j)=c_{i,j}(2)x_j,\qquad C_2=(c_{i,j}(2))_{m\times m(m-1)}.\eqno(14.3.12)$$
Moreover,
$$\mbox{the weight matrix on}\; \Psi_2\;\mbox{is equivalent to}\;(C_2,-C_2). \eqno(14.3.13)$$
Since
$$\sum_{i=1}^m\mbb C_3h_i=\sum_{i=1}^m\mbb C_3(E_{i,i}-E_{m+i,m+i}),\eqno(14.3.14)$$
$C_2$ is a generator matrix of the ternary code ${\msr C}_2$. Denote
by $\zeta_i$ the $i$th row of $C_2$.  By (14.3.9) and (14.3.13), any
nonzero codeword in ${\msr C}_2$ has the same weight as the codeword
$$u(t)=\sum_{i=1}^t\zeta_t\qquad\mbox{for some}\;\;t\in\ol{1,m}.
\eqno(14.3.15)$$ Moreover,
$$f(t)=\wt u(t)={t\choose 2}+2t(m-t)=\frac{(4m-1)t-3t^2}{2}\eqno(14.3.16)$$ So $f(t)$ has only local maximum and it
attains the absolute  minimum at the boundary points. Note that
$$f(1)=2(m-1),\qquad f(m)=\frac{m(m-1)}{2}.\eqno(14.3.17)$$
Hence
$$\mbox{the minimal distance of}\;\;{\msr
C}_2\;\;\mbox{is}\;\;2(m-1)\;\;\mbox{if}\;m\geq 4.\eqno(14.3.18)$$
\pse

{\bf Theorem 14.3.1}. {\it When $m=3m_1+1$ for some positive integer
$m_1$, the ternary weight code ${\msr C}_2$ of $o(2m,\mbb C)$ on
$\Psi_2$ is an orthogonal $[m(m-1),m,2(m-1)]$-code.}

{\it Proof}. Note that for $i,j\in\ol{1,m}$ with $i\neq j$,
$$\zeta_i\cdot\zeta_i=f(1)=6m_1,\eqno(14.3.19)$$
$$(\zeta_i+\zeta_j)\cdot(\zeta_i+\zeta_j)=f(2)=1+4(m-2)=4m-7=12(m_1-1).\eqno(14.3.20)$$
Thus
$$\zeta_i\cdot \zeta_j=\frac{f(2)-2f(1)}{2}=-6.\eqno(14.3.21)$$
Hence  ${\msr C}_2$ is an orthogonal ternary code. $\qquad\Box$ \psp

In particular, ${\msr C}_2$ is an orthogonal ternary $[12,4,6]$-code
when $m_1=1$, $[42,7,12]$-code when $m_1=2$ and $[90,10,18]$-code
when $m_1=3$. It can be proved that ${\msr C}_2$ is also the weight
code on the adjoint module of $o(2m,\mbb C)$.

The ternary weight code ${\msr C}_3$ of  $o(2m,\mbb C)$ on $\Psi_3$
is given by the weight matrix on its subspace
 $$\Psi_{3,1}=\sum_{1\leq i<j<l\leq m}\mbb C\sta_i\sta_j\sta_l+\sum_{1\leq i<j\leq m}\;\sum_{l=1}^m
 \mbb C\sta_i\sta_j\sta_{m+l}.\eqno(14.3.22)$$
We take any order
\begin{eqnarray*} & &\{y_1,y_2,\cdots,y_{{m\choose 3}+m{m\choose 2}}\}\\ &=&\{\sta_i\sta_j\sta_l,
\sta_r\sta_s\sta_{m+q}\mid 1\leq i<j<l\leq m;\;1\leq r<s\leq
m;\;q\in\ol{1,m}\}\hspace{1.9cm}(14.3.23)\end{eqnarray*} and write
$$(E_{i,i}-E_{m+i,m+i})(y_j)=c_{i,j}(3)y_j,\qquad C_3=(c_{i,j}(3))_{m\times\left({m\choose 3}+m{m\choose 2}\right)}.
\eqno(14.3.24)$$ Moreover,
$$\mbox{the weight matrix on}\; \Psi_3\;\mbox{is equivalent to}\;(C_3,-C_3).
\eqno(14.3.25)$$

Denote by $\mfk c_i$ the $i$th row of $C_3$.  By (14.3.9) and
(14.3.25), any nonzero codeword in ${\msr C}_3$ has the same weight
as the codeword
$$u(t)=\sum_{i=1}^t\mfk c_i\qquad\mbox{for some}\;\;t\in\ol{1,m}.
\eqno(14.3.26)$$ Moreover,
\begin{eqnarray*}\hspace{2cm}g(t)&=&\wt u(t)=(2m-t){t\choose 2}+2t{m-t\choose 2}+t(m-t)^2
\\ &=&\frac{t(t-1)(2m-t)+2t(m-t)(2m-2t-1)}{2}\\ &=&\frac{t}{2}[3t^2+3(1-2m)t+4(m^2-m)].
\hspace{4.6cm}(14.3.27)\end{eqnarray*} Observe that
$$g'(t)=\frac{1}{2}[9t^2+6(1-2m)t+4(m^2-m)]=\frac{1}{2}[(3t+1-2m)^2-1].\eqno(14.3.28)$$
Thus
$$g'(t_0)=0\lra t_0=\frac{2(m-1)}{3},\;\frac{2m}{3}.\eqno(14.3.29)$$
Since $g'(0)=(m^2-m)/2\geq 0$, $t=2(m-1)/3$ is a point of local
maximum and $t=2m/3$ is a point of local minimum. We calculate
$$g(1)=(m-1)(2m-3),\qquad g(m)=\frac{m^2(m-1)}{2},\qquad
g(2m/3)=\frac{2}{9}m^2(2m-3).\eqno(14.3.30)$$ Note that $g(m)\geq
g(1)$ and $g(2m/3)\geq g(1)$ if $m\geq 3$. \psp

{\bf Theorem 14.3.2}. {\it Let $m\geq 3$. The ternary weight code
${\msr C}_3$ of  $o(2m,\mbb C)$ on $\Psi_3$ is of type
$[m(m-1)(2m-1)/3,m,(m-1)(2m-3)]$. Moreover, it is orthogonal if
 $m\not\equiv -1\;(\mbox{mod}\;3)$.}

{\it Proof}. Note
$$\mfk c_i\cdot\mfk c_i=g(1)=(m-1)(2m-3)\eqno(14.3.31)$$
and
$$(\mfk c_i+\mfk c_j)\cdot(\mfk c_i+\mfk c_j)=g(2)=2(2(m-2)^2+1)\eqno(14.3.32)$$
for $i,j\in\ol{1,m}$ such that $i\neq j$. Thus
$$\mfk c_i\cdot\mfk c_j=\frac{g(2)-2g(1)}{2}=3(2-m).\eqno(14.3.33)$$
So  ${\msr C}_3$ is orthogonal if $m\not\equiv
-1\;(\mbox{mod}\;3).\qquad\Box$\psp

Remark that ${\msr C}_3$ is an orthogonal $[10,3,6]$-code when
$m=3$, $[28,4,15]$-code when $m=4$, $[110,6,45]$-code when $m=6$ and
$[182,7,66]$-code when $m=7$.

Let ${\msr B}$ be the subalgebra of $\Psi$ generated by $\{1_{\msr
A},\sta_i\mid i\in\ol{1,m}\}$ and
$${\msr B}_r=\Psi_r\bigcap {\msr B}\qquad\for\;\;r\in\ol{0,m}.\eqno(14.3.34)$$
The spin representation of $o(2m,\mbb C)$ is given by
$$E_{i,j}-E_{m+j,m+i}=\sta_i\ptl_{\sta_j}-\frac{\dlt_{i,j}}{2}\qquad
\for\;\;i,j\in\ol{1,m},\eqno(14.3.35)$$
$$E_{m+s,r}-E_{m+r,s}=\ptl_{\sta_s}\ptl_{\sta_r},\qquad
 E_{r,m+s}-E_{s,m+r}=\sta_r\sta_s\eqno(14.3.36)$$
 for $1\leq r<s\leq m$. Then the subspace
$${\cal V}=\sum_{i=1}^{\llbracket m/2\rrbracket}{\msr B}_{m-i}\eqno(14.3.37)$$
is the irreducible module with highest weight $\lmd_m$; that is,
${\cal V}\cong V_{D_m}(\lmd_m)$.

If $m=2m_1+1$ is odd, then
$$\{\sta_{i_1}\cdots\sta_{i_{m-2r}}\mid r\in
\ol{0,m_1};\;1\leq i_1<\cdots< i_{m-2r}\leq m\}\eqno(14.3.38)$$
forms a weight-vector basis of ${\cal V}$. When $m=2m_1$ is even,
$$\{1,\sta_{i_1}\cdots\sta_{i_{m-2r}}\mid r\in
\ol{0,m_1-1};\;1\leq i_1<\cdots< i_{m-2r}\leq m\}\eqno(14.3.39)$$ is
a weight-vector basis of ${\msr V}$. Take any order
$\{z_1,z_2,...,z_{2^{m-1}}\}$ of the above base vectors. Denote
$$(E_{r,r}-E_{m+r,m+r})(z_i)=q_{r,i}z_i,\qquad
C({\cal V})=(q_{r,i})_{m\times 2^{m-1}}.\eqno(14.3.40)$$ Note that
$$\frac{1}{2}\equiv -1\qquad\mbox{in}\;\;\mbb F_3.\eqno(14.3.41)$$

Denote by $\xi_r$ the $r$th row of the weight matrix $C({\cal V})$.
Set
$$\bar u=\sum_{r=1}^{m-1}\xi_r-\xi_m,\;\;
u(t)=\sum_{i=1}^t\xi_i\qquad\for\;\;t\in\ol{1,m}.\eqno(14.3.42)$$
Then any  nonzero codeword in ${\msr C}_3({\cal V})$ is conjugated
to some $u(t)$ or $\bar u$ under the action of the Weyl group of
$o(2m,\mbb C)$ (cf. (14.1.10) and (14.1.11)). It has the same weight
as $u(t)$ or $\bar u$. We calculate
$$\wt u(1)=2^{m-1},\qquad\wt u(2)=2^{m-2}.\eqno(14.3.43)$$
Moreover, we have the following more general estimates. For any
positive integer $k>2$, we always have
$${k\choose l-1}+{k\choose l+1}>{k\choose
l}\qquad\for\;\;l\in\ol{0,k},\eqno(14.3.44)$$ where we treat
${k\choose -1}={k\choose k+1}=0$. If $t=3t_1$ for some positive
integer $t_1$, we have
\begin{eqnarray*} \wt
u(t)&=&2^{m-3t_1-1}\sum_{i=0}^{t_1}\left[{3t_1\choose
6i+1}+{3t_1\choose 6i+2}+{3t_1\choose 6i+4}+{3t_1\choose 6i+5}\right]\\
&>& 2^{m-3t_1-1}\sum_{i=0}^{t_1}\left[{3t_1\choose
6i+1}+{3t_1\choose 6i+3}+{3t_1\choose 6i+5}\right]=
2^{m-2}.\hspace{1.8cm}(14.3.45)\end{eqnarray*} When $t=3t_1+1$ for
some positive integer $t_1$, we obtain
\begin{eqnarray*} \wt
u(t)&=&2^{m-3t_1-2}\sum_{i=0}^{t_1}\left[{3t_1+1\choose
6i}+{3t_1+1\choose 6i+1}+{3t_1+1\choose 6i+3}+{3t_1+1\choose 6i+4}\right]\\
&>& 2^{m-3t_1-2}\sum_{i=0}^{t_1}\left[{3t_1+1\choose
6i}+{3t_1+1\choose 6i+2}+{3t_1+1\choose 6i+4}\right]=
2^{m-2}.\hspace{1.2cm}(14.3.46)\end{eqnarray*} If $t=3t_1+2$ for
some positive integer $t_1$, we get
\begin{eqnarray*} \wt
u(t)&=&2^{m-3t_1-3}\sum_{i=0}^{t_1}\left[{3t_1+2\choose
6i}+{3t_1+2\choose 6i+2}+{3t_1+2\choose 6i+3}+{3t_1+2\choose 6i+5}\right]\\
&>& 2^{m-3t_1-3}\sum_{i=0}^{t_1}\left[{3t_1+2\choose
6i}+{3t_1+2\choose 6i+2}+{3t_1+2\choose 6i+4}\right]=
2^{m-2}.\hspace{1.2cm}(14.3.47)\end{eqnarray*}

Let $k$ be  a positive integer. We have
$${2k\choose i}+{2k\choose i+4}>{2k\choose i+1}\eqno(14.3.48)$$
if $i\leq k-3$ or $i\geq k$. Moreover,
$${2k\choose k-2}+{2k\choose k+2}-{2k\choose k-1}=\frac{k-4}{k-1}{2k\choose
k-2},\eqno(14.3.49)$$
$${2k\choose k-1}+{2k\choose k+3}-{2k\choose
k}=\frac{k^3-4k^2-3k-6}{k(k-1)(k-2)}{2k\choose k-3}.\eqno(14.3.50)$$
Thus (14.3.48) always holds if $k\geq 5$. Furthermore,
$${2k+1\choose i}+{2k+1\choose i+4}>{2k+1\choose i+1}\eqno(14.3.51)$$
if $i\neq  k-1$. Observe that
$${2k+1\choose k-1}+{2k+1\choose i+3}-{2k+1\choose
k}=\frac{k^2-3k-6}{k(k-1)}{2k+1\choose k-2}.\eqno(14.3.52)$$ So
(14.3.51) holds whenever $k\geq 5$. Therefore,
$${k\choose i}+{k\choose i+4}>{k\choose i+1}\qquad\mbox{if}\;\;k\geq 10.\eqno(14.3.53)$$

If $m=3m_1$ for some positive integer $m_1$,
\begin{eqnarray*} \wt\bar u&=&\sum_{i=0}^m\left[{m\choose 6i}+{m-1\choose
6i+1}+{m-1\choose 6i+4}\right]\\ &=&\sum_{i=0}^m\left[{m-1\choose
6i}+{m-1\choose 6i+5}+{m-1\choose 6i+1}+{m-1\choose
6i+4}\right],\hspace{2.9cm}(14.3.54)\end{eqnarray*} which is $>
2^{m-2}$ if $m_1\geq 4$ by (14.3.53). When $m=3m_1+1$ for some
positive integer $m_1$,
\begin{eqnarray*} \wt\bar u&=&\sum_{i=0}^m\left[{m-1\choose
6i}+{m\choose 6i+2}+{m-1\choose 6i+3}\right]\\
&=&\sum_{i=0}^m\left[{m-1\choose 6i}+{m-1\choose 6i+1}+{m-1\choose
6i+2}+{m-1\choose 6i+3}\right]\\
&=&1+\sum_{i=0}^m\left[{m-1\choose 6i+1}+{m-1\choose
6i+3}+{m-1\choose 6i+2}+{m-1\choose
6i+6}\right],\hspace{2.2cm}(14.3.55)\end{eqnarray*} which is again
$>2^{m-2}$ if $m_1\geq 4$ by (14.3.53). Assuming  $m=3m_1+2$ for
some positive integer $m_1$, we have
\begin{eqnarray*} \wt\bar u&=&\sum_{i=0}^m\left[{m-1\choose
6i+2}+{m\choose 6i+4}+{m-1\choose 6i+5}\right]\\ &=&{m-1\choose
3}+\sum_{i=0}^m\left[{m-1\choose 6i+2}+{m-1\choose 6i+4}+{m-1\choose
6i+5}+{m-1\choose 6i+9}\right]
,\hspace{0.7cm}(14.3.56)\end{eqnarray*} which is $> 2^{m-2}$ if
$m_1\geq 3$ by (14.3.53).
 Moreover, we have the
following table:
\begin{center}{\bf Table 14.3.1}\end{center}
\begin{center}\begin{tabular}{|c|c|c|c|c|c|c|c|}\hline
m&4&5&6&7&8&9&10
\\\hline $\wt\bar u$&8&11&12&43&112&171&260\\\hline\end{tabular}\end{center}
In summary, we have:\psp

{\bf Theorem 14.3.3}. {\it Let $m>3$ be an integer. The ternary code
${\msr C}_3({\cal V})$ is of type $[2^{m-1},m,2^{m-2}]$ if $m\neq 6$
and of type $[32,6,12]$ when $n=6$}.\psp

We remark that the spin module ${\cal V}$ is self-dual if and only
if $m$ is even. \psp

{\bf Corollary 14.3.4}. {\it When $m=6m_1+2$ for some positive
integer $m_1$, the ternary weight code of $o(2m,\mbb C)$ on
$o(2m,\mbb C)+{\cal V}$ is an orthogonal ternary
$[m(m-1)+2^{m-2},m,4m-7+2^{m-3}]$-code. If
 $m=6m_1+3$ for some positive integer $m_1$, the ternary weight code of
 $o(2m,\mbb C)$ on
$o(2m,\mbb C)+{\cal V}$ is an orthogonal ternary
$[2m(m-1)+2^{m-1},m,8m-14+2^{m-2}]$-code. In the case $m=6m_1+5$ and
$m=6m_1+12$ for some nonnegative integer $m_1$, the code ${\msr
C}_2\oplus {\msr C}_3({\cal V})$ is an orthogonal ternary
$[m(m-1)+2^{m-1},m, 4m-7+2^{m-2}]$-code. When $m=6$, the code ${\msr
C}_2\oplus {\msr C}_3({\cal V})$ is an orthogonal ternary $[62,6,
27]$-code.}

{\it Proof}. Suppose $m=6m_1+2$ for some positive integer $m_1$.
Then the weight matrix of $o(2m,\mbb C)$ on $o(2m,\mbb C)+{\cal V}$
is equivalent to $(A,-A)$, where $A$ generates the weight code
${\msr C}$ of $o(2m,\mbb C)+{\cal V}$. Moreover, ${\msr C}$ is
orthogonal if and only if the matrix $(A,-A)$ generates an
orthogonal code. But
$$(A,-A)\sim (C_2,C_2,C({\cal V})).\eqno(14.3.57)$$
Note that
$$\wt (\zeta_i,\zeta_i,\xi_i)=2f(1)+2^{m-1}=4(m-1)+2^{m-1}\equiv 1+(-1)^{6m_1+1}\equiv 0
\;\;(\mbox{mod}\;3),\eqno(14.3.58)$$
\begin{eqnarray*} & &\wt (\zeta_i+\zeta_j,\zeta_i+\zeta_j,\xi_i+\xi_j)
\\ &=&2f(2)+2^{m-2}=8m-14+2^{m-2}\equiv 2+(-1)^{6m_1}\equiv 0
\;\;(\mbox{mod}\;3)\hspace{2.4cm}(14.3.59)\end{eqnarray*} for
$i,j\in\ol{1,m}$ with $i\neq j$ by (14.3.17) and (14.3.43). Thus
$$(\zeta_i,\zeta_i,\xi_i)\cdot(\zeta_i,\zeta_i,\xi_i)\equiv \wt (\zeta_i,\zeta_i,\xi_i)\equiv 0,\eqno(14.3.60)$$
\begin{eqnarray*} & &(\zeta_i,\zeta_i,\xi_i)\cdot(\zeta_j,\zeta_j,\xi_j)
\\ &\equiv& -[\wt (\zeta_i+\zeta_j,\zeta_i+\zeta_j,\xi_i+\xi_j)-\wt (\zeta_i,\zeta_i,\xi_i)-
\wt (\zeta_j,\zeta_j,\xi_j)]\equiv
0\hspace{2cm}(14.3.61)\end{eqnarray*} by (14.3.41). Thus ${\msr C}$
is orthogonal. Note
$$f(2)=4m-7\leq \frac{m(m-1)}{2}=f(m)\qquad\mbox{if}\;\;m\geq 7.\eqno(14.3.62)$$
Thus
$$f(2)\leq f(t)\qquad\for\;\;t\in\ol{2,m}.\eqno(14.3.63)$$
By (14.3.9),
$$\wt (\sum_{i=1}^{m-1}\zeta_i-\zeta_m)=f(m)\geq f(2).\eqno(14.3.64)$$
Thus the minimum distance of ${\msr C}$ is
$$\min\{f(1)+2^{m-2},f(2)+2^{m-3}\}=4m-7+2^{m-3}\qquad\mbox{if}\;\;m\geq 6.\eqno(14.3.65)$$
This proves the first conclusion. The other conclusions for $m\geq
7$ can be proved similarly.

In the case $m=5$, we have
\begin{center}{\bf Table 14.3.2}\end{center}
\begin{center}\begin{tabular}{|c|c|c|c|c|c|}\hline
t&1&2&3&4&5
\\\hline f(t)&8&13&15&14&10\\\hline\end{tabular}\end{center}
and on the ${\cal V}$,
\begin{center}{\bf Table 14.3.3}\end{center}
\begin{center}\begin{tabular}{|c|c|c|c|c|c|}\hline
t&1&2&3&4&5
\\\hline \wt u(t)&16&8&12&10&11\\\hline\end{tabular}\end{center}
By Tables 14.3.1-14.3.3 and the fact $\wt
(\sum_{i=1}^4\zeta_i-\zeta_5)=f(5)$ in ${\msr C}_3(\Psi_2)$, the
third conclusion holds for $m=5$.

If $m=6$,
\begin{center}{\bf Table 14.3.4}\end{center}
\begin{center}\begin{tabular}{|c|c|c|c|c|c|c|}\hline
t&1&2&3&4&5&6
\\\hline f(t)&10&17&21&22&20&15\\\hline\end{tabular}\end{center}
and on the ${\cal V}$,
\begin{center}{\bf Table 14.3.5}\end{center}
\begin{center}\begin{tabular}{|c|c|c|c|c|c|c|}\hline
t&1&2&3&4&5&6
\\\hline \wt u(t)&32&16&24&20&22&21\\\hline\end{tabular}\end{center}
By Tables 14.3.1, 14.3.4, and 14.3.5,  and the fact $\wt
(\sum_{i=1}^5\zeta_i-\zeta_6)=f(6)$ in ${\msr C}_3(\Psi_2)$, the
last conclusion holds. $\qquad\Box$ \psp

When $m=8$, the ternary weight code of  $o(16,\mbb C)$ on $o(16,\mbb
C)+{\cal V}$ is a ternary orthogonal  $[120,8,57]$-code, which will
later be proved also to be the ternary weight code of $E_8$ on its
adjoint module. If $m=9$, the ternary weight code of $o(18,\mbb C)$
on $o(18,\mbb C)+{\cal V}$ is a ternary orthogonal
$[400,8,186]$-code. When $m=5$, the code ${\msr C}_2\oplus {\msr
C}_3({\cal V})$ is a  ternary orthogonal $[36,5,21]$-code, which
will later be proved also to be the ternary weight code of $E_6$ on
its adjoint module. In the case $m=11$, the code ${\msr C}_2\oplus
{\msr C}_3({\cal V})$ is a ternary orthogonal $[1134,8,549]$-code.

\section{Exceptional Lie Algebras and Ternary Codes}

In this section, we study the ternary weight codes of
$F_4,E_6,E_7,E_8$ on its minimal irreducible module and adjoint
module.

We go back to the settings in (14.2.1)-(14.2.3) with $n=4$.  The
root system of $F_4$ is
$$\Phi_{F_4}=\left\{\pm \ves_i,\pm \ves_i\pm \ves_j,\frac{1}{2}(\pm \ves_1\pm \ves_2\pm
\ves_3\pm \ves_4)\mid i\neq j\right\}\eqno(14.4.1)$$  and the
positive simple roots are
$$\bar\al_1=\ves_2-\ves_3,\bar\al_2=\ves_3-\ves_4,\bar\al_3=\ves_4,\bar\al_4=\frac{1}{2}(\ves_1-\ves_2-\ves_3-\ves_4).
\eqno(14.4.2)$$ The Weyl group ${\msr W}_{F_4}$ of $F_4$ contains
the permutation group $S_4$ on the sub-indices of $\ves_i$ and all
reflections with respect to the coordinate hyperplanes. Moreover,
there is an identification:
$$h_1\leftrightarrow\al_1,\;h_2\leftrightarrow\al_2,\;h_3\leftrightarrow 2\al_3,\;h_4\leftrightarrow 2\al_4.\eqno(14.4.3)$$
Thus
$${\msr H}_3=\sum_{i=1}^4\mbb{F}_3h_i=\sum_{i=1}^4\mbb{F}_3\ves_i.\eqno(14.4.4)$$
Moreover,
$${\msr H}_3=\{{\msr W}_{F_4}(h_1),{\msr W}_{F_4}(h_1+h_3),{\msr W}_{F_4}(h_3),{\msr W}_{F_4}(h_4)\}.\eqno(14.4.5)$$

Recall that the basic (minimal) irreducible module $V_{F_4}$ of the
52-dimensional Lie algebra ${\msr G}^{F_4}$ has a basis
$\{x_i\mid\ol{1,26}\}$ and with the representation determined by
(10.2.44)-(10.2.74). The module $V_{F_4}$ is self-dual. The weight
matrix of $V_{F_4}$ is $(A_{F_4},-A_{F_4})$ with
$$A_{F_4}=\left[\begin{array}{rrrrrrrrrrrr}0&0&0&1&1&-1&1&-1&-1&0&0&0\\
0&0&1&-1&0&0&0&1&1&-1&-1&0\\
0&1&-1&1&-1&1&0&-1&0&1&2&-1\\
1&-1&0&0&1&0&-1&1&-1&1&-1&2\end{array}\right]\eqno(14.4.6)$$ by
(10.2.68)-(10.2.71).\psp

{\bf Theorem 14.4.1}. {\it The ternary weight code ${\msr
C}_{F_4,1}$ (generated by $A_{F_4}$) of $F_4$ on $V_{F_4}$ is an
orthogonal [12,4,6]-code}.

 {\it Proof}. Denote by $\xi_i$ the $i$th row of the matrix $A_{F_4}$. Then
$$\wt\xi_1=6,\qquad \wt (\xi_1+\xi_3)=\wt\xi_3=\wt\xi_4=9.\eqno(14.4.7)$$
According to (14.4.5), any nonzero codeword in ${\msr C}_{F_4,1}$
has weight 6 or 9. By an argument as (14.3.31)-(14.3.33), ${\msr
C}_{F_4,1}$ is orthogonal.$\qquad\Box$\psp

Next we consider the adjoint representation of $F_4$. Its weight
code ${\msr C}_{F_4,2}$ is determined by the set $\Phi_{F_4}^+$ of
positive roots, the corresponding root vectors are constructed in
(10.2.29)-(10.2.42). The weight matrix $B_{F_4}$ on
$\sum_{\al\in\Phi_{F_4}^+}\mbb CE_\al$ is given by
$${\tiny \left[\begin{array}{rrrrrrrrrrrrrrrrrrrrrrrr}2&$-1$&0&0&1&$-1$&0&1&$-1$&$-1$&1&$-1$&1&0&1&$-1$&0&1&0&0&0&0&$-1$&1\\
$-1$&2&$-1$&0&1&1&$-1$&0&1&0&$-1$&0&0&1&$-1$&0&1&$-1$&1&0&0&$-1$&1&0\\
0&$-2$&2&$-1$&$-2$&0&1&0&$-1$&2&2&1&$-1$&0&1&0&$-1$&0&$-2$&1&0&2&0&0\\
0&0&$-1$&2&0&$-1$&1&$-1$&1&$-2$&$-2$&0&1&$-2$&0&2&0&2&2&$-1$&1&0&0&0\end{array}\right].}\eqno(14.4.8)$$
\pse

{\bf Theorem 14.4.2}. {\it The ternary weight code ${\msr
C}_{F_4,2}$ (generated by $B_{F_4}$) of $F_4$ on its adjoint module
is an orthogonal $[24,4,15]$-code.}\psp

{\it Proof}. Denote by $\eta_i$ the $i$th row of the above matrix.
Then
$$\wt \eta_i=15,\qquad\wt(\eta_1+\eta_3)=18.\eqno(14.4.9)$$
According to (14.4.5), any nonzero codeword in ${\msr C}_{F_4,2}$
has weight 15 or 18. By an argument as (14.3.31)-(14.3.33) , ${\msr
C}_{F_4,2}$ is orthogonal.$\qquad\Box$\psp

Now we take the settings in (14.2.1)-(14.2.3) with $n=7$ and have
the following root system of $E_6$:
\begin{eqnarray*}\qquad\qquad \Phi_{E_6}&=&\big\{\ves_i-\ves_j,\frac{1}{2}(\sum_{s=1}^6\iota_s\ves_s\pm
\sqrt{2}\ves_7),\pm \sqrt{2}\ves_7\\ & &\mid i,j\in\ol{1,6},i\neq
j;\iota_s=\pm
1;\sum_{i=1}^6\iota_i=0\big\}\hspace{4.3cm}(14.4.10)\end{eqnarray*}
and the simple positive roots are
$$\al_1=\ves_1-\ves_2,\;\;\al_2=\frac{1}{2}(\sum_{j=1}^3(\ves_{3+j}-\ves_j)+\sqrt{2}\ves_7),\;\;
\al_i=\ves_{i-1}-\ves_i,\qquad i\in\ol{3,6}.\eqno(14.4.11)$$

Note that
$${\msr H}_{E_6,3}=\sum_{i=1}^6\mbb{F}_3h_i=\{\sum_{i=1}^6\iota_i\ves_i+\iota_7\sqrt{2}\ves_7\mid
\iota_r\in\mbb{F}_3,\;\sum_{i=1}^6\iota_i=0\}.\eqno(14.4.12)$$
Moreover, the Weyl group ${\msr W}_{E_6}$ contains the permutation
group $S_6$ on the first six sub-indices  of $\ves_i$ and the
reflection
$$\sum_{i=1}^6\iota_i\ves_i+\iota_7\sqrt{2}\ves_7\mapsto \sum_{i=1}^6\iota_i\ves_i-\iota_7\sqrt{2}\ves_7.\eqno(14.4.13)$$
So
$${\msr H}_{E_6,3}={\msr W}_{E_6}(\{\sum_{i=1}^s\ves_i-\sum_{j=1}^t\ves_{s+j}+\iota\sqrt{2}\ves_7,\sqrt{2}\ves_7\mid
\iota=0,1;s-t\equiv 0\;(\mbox{mod}\;3)\}).\eqno(14.4.14)$$

Recall that the 27-dimensional basic irreducible module $V_{E_6}$ of
weight $\lmd_1$ for $E_6$ has a basis $\{x_i\mid i\in\ol{1,27}\}$
with the representation formulas determined by (11.1.31)-(11.1.70).
Moreover,
$$E_{\al_r}(x_i)\neq 0\Leftrightarrow a_{r,i}<0,\;\;
E_{-\al_r}(x_i)\neq 0\Leftrightarrow a_{r,i}>0.\eqno(14.4.15)$$
 \pse

{\bf Theorem 14.4.3}. {\it The ternary weight code ${\msr
C}_{E_6,1}$ of $E_6$ on $V_{E_6}$ is an orthogonal
$[27,6,12]$-code.}\psp

{\it Proof}. According to Table 11.1.1 and (11.1.72), we write
$$A_{E_6}=(a_{r,i})_{6\times 27}.\eqno(14.4.16)$$
Denote by $\xi_r$ the $r$th row of the matrix $A_{E_6}$. Then
$$\wt\xi_r=12\qquad\for\;\;r\in\ol{1,6}.\eqno(14.4.17)$$
Moreover,
$$\wt(\xi_1+\xi_3)=\wt(\xi_2+\xi_4)=12,\;\;\wt (\xi_1+\xi_4)=18, \eqno(14.4.18)$$
$$\wt (\xi_1+\xi_2)=\wt(\xi_2+\xi_3)=\wt(\xi_2+\xi_5)=\wt(\xi_2+\xi_6)=18.\eqno(14.4.19)$$
By an argument as (14.3.31)-(14.3.33) and symmetry, we have
$$\xi_i\cdot\xi_j\equiv 0\;(\mbox{mod}\;3)\qquad\for\;\;i,j\in\ol{1,6};\eqno(14.4.20)$$
that is ${\msr C}_{E_6,1}$ is orthogonal.

 Note that
the Lie subalgebra ${\msr G}^{E_6}_{A,1}$ generated by
$\{E_{\pm\al_i}\mid 2\neq i\in\ol{1,6}\}$ is isomorphic to
$sl(6,\mbb C)$. Recall that a singular vector in a module of simple
Lie algebra is a weight vector annihilated by its positive root
vectors. By Table 11.1.1 and (14.4.15), the ${\msr
G}^{E_6}_{A,1}$-singular vectors are $x_1$ of weight $\lmd_1$, $x_6$
of  weight $\lmd_4$ and $x_{20}$ of weight $\lmd_1$. So the $({\msr
G}^{E_6},{\msr G}^{E_6}_{A,1})$-branch rule on $V_{E_6}$ is
$$V_{E_6}\cong V_{A_5}(\lmd_1)\oplus V_{A_5}(\lmd_4)\oplus V_{A_5}(\lmd_1).\eqno(14.4.21)$$
Denote by ${\msr G}^{E_6}_{A,2}$ the Lie subalgebra of ${\msr
G}^{E_6}$ generated by $\{E_{\pm\al_r},\;E_{\pm(\al_2+\al_4)}\mid
2,4\neq r\in\ol{1,6}\}$. The algebra ${\msr G}^{E_6}_{A,2}$ is also
isomorphic to $sl(6,\mbb C)$. According to Table 11.1.1 and
(14.4.15), the ${\msr G}^{E_6}_{A,2}$-singular vectors are $x_1$ of
weight $\lmd_1$, $x_4$ of  weight $\lmd_4$ and $x_{18}$ of  weight
$\lmd_1$. Hence (14.4.21) is also the $({\msr G}^{E_6},{\msr
G}^{E_6}_{A,2})$-branch rule. Since the module $V_{A_5}(\lmd_2)$ is
contragredient to $V_{A_5}(\lmd_4)$, they have the same ternary
weight code. By (14.2.34) and (14.2.38) with $n=6$, the minimal
distances of the subcodes $\sum_{2\neq i\in\ol{1,6}}\mbb{F}_3\xi_i$
and $\mbb{F}_3(\xi_2+\xi_4)+\sum_{2,4\neq
i\in\ol{1,6}}\mbb{F}_3\xi_i$ are $\wt\xi_1=12$.

Recall $\frac{1}{2}=-1$ in $\mbb{F}_3$. Moreover,
$$-(\al_2+\al_4)=-\ves_1-\ves_2+\ves_3-\ves_4+\ves_5+\ves_6+\sqrt{2}\ves_7\;\;\mbox{in}\;\;{\msr H}_{E_6,3}.
\eqno(14.4.22)$$ Thus in ${\msr H}_{E_6,3}$,
$$\al_1-(\al_2+\al_4)=\ves_2+\ves_3-\ves_4+\ves_5+\ves_6+\sqrt{2}\ves_7,\eqno(14.4.23)$$
$$\al_1-\al_2-(\al_2+\al_4)=-\ves_3-\ves_4+\ves_5+\ves_6+\sqrt{2}\ves_7,\eqno(14.4.24)$$
$$\al_1-\al_2-(\al_2+\al_4)+\al_6=-\ves_3-\ves_4-\ves_5+\sqrt{2}\ves_7,\eqno(14.4.25)$$
$$\al_1-\al_2-(\al_2+\al_4)-\al_5+\al_6=-\ves_3+\ves_4+\sqrt{2}\ves_7.\eqno(14.4.26)$$
Note that
$$\wt (\xi_1-(\xi_2+\xi_4)),\;\wt(\xi_1-\xi_2-(\xi_2+\xi_4))\geq 12,\eqno(14.4.27)$$
$$\wt(\xi_1-\xi_2-(\xi_2+\xi_4)+\xi_6),\;\;\wt(\xi_1-\xi_2-(\xi_2+\xi_4)-\xi_5+\xi_6)\geq 12\eqno(14.4.28)$$
because the minimal distance of
$\mbb{F}_3(\xi_2+\xi_4)+\sum_{2,4\neq i\in\ol{1,6}}\mbb{F}_3\xi_i$
is 12. Furthermore,
$$-\sum_{i=1}^6\ves_i+\sqrt{2}\ves_7=\al_1-\al_2-\al_3\qquad\mbox{in}\;\;{\msr H}_{E_6,3}.\eqno(14.4.29)$$
We calculate
$$\wt(\xi_1-\xi_2-\xi_3)=21.\eqno(14.4.30)$$
By (14.4.13), the minimal distance of the ternary code ${\msr
C}_{E_6,1}$ is 12. $\qquad\Box$\psp

Next we consider the ternary weight code ${\msr C}_{E_6,2}$ of $E_6$
on its adjoint module. Take any order
$$\{y_1,...,y_{36}\}=\{E_\al\mid \al\in \Phi^+_{E_6}\}.\eqno(14.4.31)$$
Write
$$[\al_i,y_j]=b_{i,j}y_j,\qquad B_{E_6}=(b_{i,j})_{6\times 36}.\eqno(14.4.32)$$
\pse

{\bf Theorem 14.4.4}. {\it The ternary weight code ${\msr
C}_{E_6,2}$ (generated $B_{E_6}$) of $E_6$ on its adjoint module is
an orthogonal $[36,5,21]$-code.}

{\it Proof}. Denote by $\zeta_i$ the $i$th row of $B_{E_6}$. Note
$$(\al_1-\al_3+\al_5-\al_6,\al_i)\equiv
0\;\;(\mbox{mod}\;3)\qquad\for\;\;i\in\ol{1,6}.\eqno(14.4.33)$$ Thus
$$\zeta_1-\zeta_3+\zeta_5-\zeta_6\equiv 0\qquad\mbox{in}\;\;\mbb{F}_3^{36}.\eqno(14.4.34)$$
Hence
$${\msr C}_{E_6,2}=\sum_{i=2}^6\mbb{F}_3\zeta_i.\eqno(14.4.35)$$
Denote by ${\msr G}^{E_6}_D$ the Lie subalgebra  of ${\msr G}^{E_6}$
generated by $\{E_{\pm \al_r}\mid r\in\ol{2,6}\}$. According to the
Dynkin diagram of $E_6$,
$${\msr G}^{E_6}_D\cong o(10,\mbb C).\eqno(14.4.36)$$
Let ${\msr G}^{E_6}_+=\sum_{i=1}^{36}\mbb Cy_i$ and denote by ${\msr
G}^{E_6}_{D,+}$ the subspace spanned by the  root vectors
$E_\al\in{\msr G}^{E_6}_D$ with $\al\in\Phi_{E_6}^+$. Then $[{\msr
G}^{E_6}_{D,+},{\msr G}^{E_6}_+]\subset {\msr G}^{E_6}_+$. Moreover,
 the space
${\msr G}^{E_6}_+$ contains  ${\msr G}^{E_6}_D$-singular vectors
$E_{\al_4+\al_5+\sum_{i=2}^6\al_i}$ of weight $\lmd_2$ (the highest
root) and $E_{\al_2+\al_4+\sum_{r=3}^5\al_r+\sum_{i=1}^6\al_i}$ of
weight $\lmd_5$. Hence, we have the partial $({\msr G}_{E_6},{\msr
G}^{E_6}_D)$-branch rule on ${\msr G}_{E_6}$:
$${\msr G}_{E_6}^+\cong{\msr G}^{E_6}_{D+}\oplus V_{D_5}(\lmd_5).
\eqno(14.4.37)$$ Thus the ternary weight code ${\msr C}_{E_6,2}$ of
$E_6$ on its adjoint module  is exactly the code ${\msr C}_2\oplus
{\msr C}_3({\cal V})$ with $n=5$ in Corollary 14.3.4, which is a
ternary orthogonal $[36,5,21]$-code.$\qquad\Box$ \psp

Take the settings in (14.2.1)-(14.2.3) with $n=8$. We have the root
system of $E_7$:
 $$\Phi_{E_7}=\left\{\ves_i-\ves_j,\frac{1}{2}\sum_{s=1}^8\iota_s\ves_s
\mid i,j\in\ol{1,8},\;i\neq j;\;\iota_s=\pm
1,\;\sum_{s=1}^8\iota_s=0\right\}\eqno(14.4.38)$$ and the simple
positive roots are:
$$\al_1=\ves_2-\ves_3,\;\al_2=\frac{1}{2}\sum_{j=1}^4(\ves_{4+j}-\ves_j),\;\;\al_i=\ves_i-\ves_{i+1},\qquad i\in\ol{3,7}.
\eqno(14.4.39)$$
 The minimal module $V_{E_7}$ of $E_7$ is of  56-dimensional and
has a basis $\{x_i\mid i\in\ol{1,56}\}$ with the representation
given by (12.1.30)-(12.1.96). Again we have
$$E_{\al_r}(x_i)\neq 0\Leftrightarrow a_{r,i}<0,\;\;
E_{-\al_r}(x_i)\neq 0\Leftrightarrow a_{r,i}>0.\eqno(14.4.40)$$
According to (12.1.96) and Table 12.1.1, we define
$$A_{E_7}=(a_{r,i})_{7\times 28}.\eqno(14.4.41)$$
\pse

{\bf Theorem 14.4.5}. {\it The ternary weight code ${\msr
C}_{E_7,1}$ of $E_7$ on $V_{E_7}$ is an orthogonal
$[28,7,12]$-code}.

{\it Proof}. Note that the root system of $A_7$:
$$\Phi_{A_7}=\{\ves_i-\ves_j\mid i,j\in\ol{1,8},\;i\neq j\}\subset \Phi_{E_7}.\eqno(14.4.42)$$
Thus we have the Lie subalgebra of ${\msr G}^{E_7}$ (cf.
(12.1.1)-(12.1.3)):
$${\msr G}^{E_7}_A=\sum_{i=1}^7\mbb C\al_i+\sum_{\al\in\Phi_{A_7}}\mbb CE_\al\cong sl(8,\mbb C)
.\eqno(14.4.43)$$ Moreover,
$$\al_1'=\ves_1-\ves_2=-2\al_2-2\al_1-3\al_3-4\al_4-3\al_5-2\al_6-\al_7.\eqno(14.4.44)$$
Note that $x_{23}$ is a ${\msr G}^{E_7}_A$-singular vector of weight
$\lmd_6$ and $x_{49}$ is a ${\msr G}^{E_7}_A$-singular vector of
weight $\lmd_2$ by (12.1.96), (14.4.40) and Table 12.1.1. Thus the
$({\msr G}^{E_7},{\msr G}^{E_7}_A)$-branch rule on $V_{E_7}$ is
$$V_{E_7}\cong V_{A_7}(\lmd_2)\oplus V_{A_7}(\lmd_6).\eqno(14.4.45)$$
Since  $V_{A_7}(\lmd_6)$ is contragredient to $V_{A_7}(\lmd_2)$,
they have the same ternary weight code of ${\msr G}^{E_7}_A$, which
is the ${\msr C}_3(\Psi_2)$ with $m=2$ in Theorem 14.2.3. Hence the
weight matrix of ${\msr G}^{E_7}_A$ on $V_{E_7}$ generates a ternary
orthogonal $[56,7,24]$-code.

On the other hand,
$$\sum_{i=1}^7\mbb{F}_3\al_i=\mbb{F}_3\al_1'+\sum_{2\neq i\in\ol{1,7}}\mbb{F}_3\al_i\eqno(14.4.46)$$
by (14.4.44) and  the fact $1/2\equiv -1$ in $\mbb{F}_3$. Thus the
weight matrix $(A_{E_7},-A_{E_7})$ of $E_7$ on $V_{E_7}$ generates
the same ternary code as the weight matrix of ${\msr G}^{E_7}_A$ on
$V_{E_7}$. So $(A_{E_7},-A_{E_7})$ generates a ternary orthogonal
$[56,7,24]$-code. Hence the ternary code ${\msr C}_{E_7,1}$
generated by $A_{E_7}$ is  an orthogonal
$[28,7,12]$-code.$\qquad\Box$\psp

Next we consider the ternary weight code of $E_7$ on its adjoint
module.  By (12.2.88)-(12.2.90) and Table 12.2.1,  it can be proved
that the $({\msr G}^{E_7},{\msr G}^{E_7}_A)$-branch rule on ${\msr
G}^{E_7}$ is
$${\msr G}^{E_7}\cong {\msr G}^{E_7}_A\oplus V_{A_7}(\lmd_4).\eqno(14.4.47)$$
The module $V_{A_7}(\lmd_4)$ of $sl(8,\mbb C)\;(\cong {\msr
G}^{E_7}_A)$ is exactly $\Psi_4$ in (6.2.16) with $n=8$, which is
self-dual. For convenience, we study the ternary code generated by
the weight matrix of $sl(8,\mbb C)$ on $\Psi_4$. Taking any order of
its basis
$$\{z_1,...,z_{70}\}=\{\sta_{i_1}\sta_{i_2}\sta_{i_3}\sta_{i_4}\mid 1\leq i_1<i_2<i_3<i_4\leq 8\},\eqno(14.4.48)$$
we write
$$[E_{r,r},z_i]=b_{r,i}z_i,\qquad B_{E_7}=(b_{r,i})_{7\times 70}.\eqno(14.4.49)$$
Denote by $\eta_r$ the $r$th row of $B_{E_7}$ and by ${\msr C}'$ the
ternary code generated by $B_{E_7}$. Set
$$v(s,t)=\sum_{i=1}^s\eta_i-\sum_{j=1}^t\eta_{s+j}\qquad\in {\msr C}'.\eqno(14.4.50)$$
Moreover, we only calculate the related weights:
\begin{center}{\bf Table 14.4.1}\end{center}
\begin{center}\begin{tabular}{|c|c|c|c|c|c|c|c|c|}\hline
(s,t)&(1,1)&(2,2)&(3,3)&(4,4)&(3,0)&(6,0)&(4,1)&(5,2)
\\\hline \wt v(s,t)&40&44&48&34&60&30&46&50\\\hline\end{tabular}\end{center}
Recall (14.2.65)-(14.2.70). We have
\begin{center}{\bf Table 14.4.2}\end{center}
\begin{center}\begin{tabular}{|c|c|c|c|c|c|c|c|c|}\hline
(s,t)&(1,1)&(2,2)&(3,3)&(4,4)&(3,0)&(6,0)&(4,1)&(5,2)
\\\hline 2\wt u(s,t)&26&40&42&32&30&24&38&34\\\hline\end{tabular}\end{center}
According to (14.4.1), the Weyl group ${\msr W}_{E_7}$ contains the
permutation group $S_8$ on the sub-indices of $\ves_i$. By (14.1.9),
(14.1.11) and the values of $\wt v(s,t)+2\wt u(s,t)$ from the above
tables, 54, 66, 84 and 90 are the only weights of the nonzero
codewords in ${\msr C}_3({\msr G}^{E_7})$, the ternary code
generated by the weight matrix of ${\msr G}^{E_7}_A$ on ${\msr
G}^{E_7}$. By (14.4.46) and an argument as (14.3.31)-(14.3.33), we
have: \psp

{\bf Theorem 14.4.6}. {\it The ternary weight code of $E_7$ on its
adjoint module is an orthogonal $[63,7,27]$-code}.\psp

The minimal representation of $E_8$ is its adjoint module given
(12.1.1)-(12.1.3). Recall the settings in (14.2.1)-(14.2.3). We have
the $E_8$ root system
$$\Phi_{E_8}=\left\{\pm \ves_i\pm \ves_j,\frac{1}{2}\sum_{i=1}^8\iota_i\ves_i
\mid i,j\in\ol{1,8},\;i\neq j;\;\iota_i=\pm
1,\;\sum_{i=1}^8\iota_i\in 2\mbb{Z}\right\}\eqno(14.4.51)$$ and
positive simple roots:
$$\al_1=\frac{1}{2}(\sum_{j=2}^7\ves_j-\ves_1-\ves_8),\;\al_2=-\ves_1-\ves_2,\;\al_r=\ves_{r-2}-\ves_{r-1},\qquad
r\in\ol{3,8}.\eqno(14.4.52)$$

 Observe that the root system of
$o(16,\mbb C)$:
$$\Phi_{D_8}=\{\pm \ves_i\pm \ves_j
\mid i,j\in\ol{1,8},\;i\neq j\}\subset \Phi_{E_8}.\eqno(14.4.53)$$
So the Lie subalgebra
$${\msr G}^{E_8}_D=H_{E_8}+\sum_{\al\in\Phi_{D_8}}\mbb CE_\al\eqno(14.4.54)$$
of ${\msr G}^{E_8}$ is exactly isomorphic to $o(16,\mbb C)$.
Moreover, it can be proved that the $({\msr G}^{E_8},{\msr
G}^{E_8}_D)$-branch rule on ${\msr G}^{E_8}$ is
$${\msr G}^{E_8}\cong {\msr G}^{E_8}_D\oplus V_{D_8}(\lmd_8).\eqno(14.4.55)$$
In fact, $V_{D_8}(\lmd_8)$ is exactly the spin module ${\cal V}$ in
(14.3.37). Since
$$\sum_{i=1}^8\mbb{F}_3\al_i=\sum_{\al\in\Phi_{D_8}}\mbb{F}_3\al,\eqno(14.4.56)$$
the ternary weight code of $E_8$ on ${\msr G}^{E_8}$ is the same as
that of ${\msr G}^{E_8}_D$ on ${\msr G}^{E_8}$. By Corollary 14.3.4
with $m=8$, we have: \psp

{\bf Theorem 14.4.7}. {\it The ternary weight code of $E_8$ on its
adjoint module is an orthogonal $[120,8,57]$-code }.

\chapter{Root-Related Integrable Systems}

The Calogero-Sutherland model is an exactly solvable quantum
many-body system in one-dimension (cf. [C], [S]), which is related
to a root system of A. Olshanesky and Perelomov [OP] generalized it
the integrable systems related the root syustems of types B, C and
$D$. In this chapter, we want to show their connections with
representation theory.

 In Section 15.1, we prove that certain variations of Weyl denominator
functions in (5.3.13) of classical types are solutions of  the
corresponding Olshanesky-Perelomov models. In Section 15.2, we show
that the Etingof's trance functions of types A and C give rise to
the solutions of the Calogero-Sutherland model and the corresponding
Olshanesky-Perelomov model, respectively. In Section 15.3, we
present a connection between the Calogero-Sutherland model of two
particles and the classical Gauss hypergeometric function. A new
family of multi-variable hypergeometric functions are found from the
Calogero-Sutherland model of more particles via the Etingof's trance
functions of types A. In Section 15.4, we determine another family
of multi-variable hypergeometric functions from the
Olshanesky-Perelomov model  via the Etingof's trance functions of
types C. In Section 15.5, we find the differential properties, Euler
integral representations and differential equations  for the
hypergeometric functions of type A,  and the differential properties
and differential equations  for the hypergeometric functions of type
C. Moreover, we define our multi-variable hypergeometric functions
of type B and D analogously as those of type C and find the
corresponding differential properties and differential equations.

The results in Sections 15.1, 15.3, 15.4 and 15.5 were due to our
work [X]. The results in Section 15.2 are taken from Etingof's work
[Ep] and Etingof-Styrkas' work [ES].

\section{Integrable Systems and Weyl Functions }

The Calogero-Sutherland model is an exactly solvable quantum
many-body system in one-dimension
 (cf. [C], [S]), whose Hamiltonian is given by
$$H_{CS}=\sum_{i=1}^n\ptl_{x_i}^2+K\sum_{1\leq i<j\leq n}\frac{1}{\sinh^2(x_i-x_j)},\eqno(15.1.1)$$
where $K$ is a constant. The model was used to study long-range
interactions of $n$ particles. It has become an important
mathematical object partly because its potential is related to the
root system of the special linear algebra $sl(n,\mbb C)$. Olshanesky
and Perelomov [OP] generalized the
  Calogero-Sutherland model to the system whose Hamiltonian is given by
\begin{eqnarray*}\qquad H_{OP}&=&\sum_{i=1}^n\ptl_{x_i}^2+\sum_{1\leq i<j\leq n}\left(\frac{K_1}
{\sinh^2(x_i-x_j)}+\frac{K_2}{\sinh^2(x_i+x_j)}\right)\\ &
&+\sum_{i=1}^n\left(\frac{K_3}
{\sinh^2x_i}+\frac{K_4}{\sinh^22x_i}\right),\hspace{6.5cm}(15.1.2)\end{eqnarray*}
where $K_1,K_2,K_3,K_4$ are constants. The potential of the
Olshanesky-Perelomov model is related to the root systems of all
four families of finite-dimensional classical simple Lie algebras.
Solving a system is to find eigenfunctions of its Hamiltonians.

 Changing
variables
$$z_i=e^{2x_i}\qquad\for\;\;i\in\ol{1,n},\eqno(15.1.3)$$
we have the equation for the Calogero-Sutherland model (15.1.1):
$$\sum_{i=1}^n(z_i\ptl_{z_i})^2(\Psi)+K\left(\sum_{1\leq i<j\leq n}\frac{z_iz_j}{(z_i-z_j)^2}
\right)\Psi=\nu \Psi\eqno(15.1.4)$$ and  the equation for the
Olshanesky-Perelomov model (15.1.2):
\begin{eqnarray*}\qquad& & \sum_{i=1}^n\ptl_{z_i}^2(\Psi)+[\sum_{1\leq i<j\leq n}
\left(K_1\frac{z_iz_j}{(z_i-z_j)^2}+K_2\frac{z_iz_j}{(z_iz_j-1)^2}\right)\\
& &+\sum_{i=1}^n
\left(K_3\frac{z_i}{(z_i-1)^2}+K_4\frac{z_i^2}{(z_i^2-1)^2}\right)]\Psi=\nu\Psi.\hspace{4.9cm}
(15.1.5)\end{eqnarray*}

Let $f_{i,j}(z)\mid i,j\in\ol{1,n}\}$ be a set of one-variable
differentiable functions and
 let $d_i$ be a one-variable differential operator in $z_i$ for $i\in\ol{1,n}$. It is easy to
 verified the following lemma:
\psp

{\bf Lemma 15.1.1}. {\it We have the following equation on
differentiation of determinants}:
\begin{eqnarray*}\qquad & &(\sum_{i=1}^nd_i)\left(\left|\begin{array}{cccc}f_{1,1}(z_1)&
f_{1,2}(z_2)&\cdots & f_{1,n}(z_n)\\ f_{2,1}(z_1)& f_{2,2}(z_2)&\cdots & f_{2,n}(z_n)\\
 \vdots&\vdots&\vdots&\vdots\\ f_{n,1}(z_1)& f_{n,2}(z_2)&\cdots &f_{n,n}(z_n)\end{array}
 \right|\right)\\&=&\sum_{i=1}^n\left|\begin{array}{cccc}f_{1,1}(z_1)& f_{1,2}(z_2)&
\cdots & f_{1,n}(z_n)\\ \vdots&\vdots&\vdots&\vdots\\
f_{i-1,1}(z_1)& f_{i-1,2}(z_2)&\cdots & f_{i-1,n}(z_n)\\
d_1(f_{i,1}(z_1))& d_2(f_{i,2}(z_2))&\cdots & d_n(f_{i,n}(z_n))
\\ f_{i+1,1}(z_1)& f_{i+1,2}(z_2)&\cdots & f_{i+1,n}(z_n)\\ \vdots&\vdots&\vdots&\vdots\\
 f_{n,1}(z_1)& f_{n,2}(z_2)&\cdots & f_{n,n}(z_n)\end{array}\right|.\hspace{3.8cm}(15.1.6)
 \end{eqnarray*}\pse

Suppose that $R^+_X\subset \mbb{R}^n$ is the set of positive roots
of the finite-dimensional
 simple Lie algebra of type X. We denote
$$z^\al=\prod_{i=1}^nz_i^{\al_i}\qquad\for\;\;\al\in\mbb{R}^n.\eqno(15.1.7)$$
The Weyl function
$${\cal W}_X=\prod_{\al\in R^+_X}(z^{\al/2}-z^{-\al/2}).\eqno(15.1.8)$$

First  the Weyl function of $sl(n,\mbb C)$ is
$${\cal W}_{A_{n-1}}=\prod_{1\leq i<j\leq n}(z_i^{1/2}z_j^{-1/2}-z_i^{-1/2}z_j^{1/2})=
(z_1z_2\cdots z_n)^{(1-n)/2}\prod_{1\leq i<j\leq
n}(z_i-z_j),\eqno(15.1.9)$$ Denote the Vandermonde determinant
$$W(z_1,z_2,...,z_n)=\left|\begin{array}{cccc}1&1&\cdots &1\\ z_1&z_2&\cdots &z_n\\ z_1^2&z_2^2
&\cdots &z_n^2 \\ \vdots&\vdots&\vdots&\vdots\\
z_1^{n-1}&z_2^{n-1}&\cdots &z_n^{n-1}
\end{array}\right|=\prod_{1\leq i<j\leq n}(z_i-z_j).\eqno(15.1.10)$$
Then $W(z_1,z_2,...,z_n)$ is the fundamental part of ${\cal
W}_{A_{n-1}}$. By Lemma 15.1.1,
 we have
\begin{eqnarray*} (\sum_{i=1}^n(z_i\ptl_{z_i})^2)(W(z_1,z_2,...,z_n))&=&(\sum_{i=1}^{n-1}i^2)
W(z_1,z_2,...,z_n)\\
&=&\frac{(n-1)n(2n-1)}{6}W(z_1,z_2,...,z_n).\hspace{2cm}(15.1.11)
\end{eqnarray*}
On the other hand,
\begin{eqnarray*} & &(\sum_{i=1}^n(z_i\ptl_{z_i})^2)(W(z_1,z_2,...,z_n))\\ &=&(\sum_{r=1}^n
(z_r\ptl_{z_r})^2)(\prod_{1\leq i<j\leq n}(z_i-z_j))\\
&=&\sum_{r=1}^nz_r\left[-\sum_{s=1}^{r-1} \left(\prod_{1\leq i<j\leq
n;\:(i,j)\neq (s,r)}(z_i-z_j)\right)+\sum_{s=r+1}^n
\left(\prod_{1\leq i<j\leq n;\:(i,j)\neq
(r,s)}(z_i-z_j)\right)\right]\\ & &+2\sum_{r=1}^n z_r^2(\sum_{1\leq
s_1<s_2<r}\left[\prod_{1\leq i<j\leq n;\:(i,j)\neq
(s_1,r),(s_2,r)}(z_i-z_j) \right]\\ & &+ \sum_{r<s_1<s_2\leq
n}\left[\prod_{1\leq i<j\leq n;\:(i,j)\neq
(r,s_1),(r,s_2)}(z_i-z_j)\right]
\\ & &-\sum_{1\leq s_1<r<s_2\leq n}\left[\prod_{1\leq i<j\leq n;\:(i,j)\neq (s_1,r),(r,s_2)}
(z_i-z_j)\right])\\ &=&\sum_{1\leq s<r\leq n}(z_s-z_r)\prod_{1\leq
i<j\leq n;\:(i,j)\neq (s,r)} (z_i-z_j)\\ &
&+2W(z_1,z_2,...,z_n)\sum_{r=1}^n\sum_{1\leq s_1<s_2\leq
n;\:s_1,s_2\neq r} \frac{z_r^2}{(z_{s_1}-z_r)(z_{s_2}-z_r)}\\
&=&\left(\frac{n(n-1)}{2}+2 \sum_{r=1}^n\sum_{1\leq s_1<s_2\leq
n;\:s_1,s_2\neq r}\frac{z_r^2}{(z_{s_1}-z_r)(z_{s_2}-z_r)} \right)\\
& &\times W(z_1,z_2,...,z_n).\hspace{9.9cm}(15.1.12)\end{eqnarray*}
Thus (15.1.11) and (15.1.12) yield
\begin{eqnarray*}\qquad\qquad& &\sum_{r=1}^n\sum_{1\leq s_1<s_2\leq n;\:s_1,s_2\neq r}
\frac{z_r^2}{(z_{s_1}-z_r)(z_{s_2}-z_r)}\\
&=&\frac{1}{2}\left[\frac{(n-1)n(2n-1)}{6}
-\frac{n(n-1)}{2}\right]\\
&=&\frac{(n-1)n(n-2)}{6}={n\choose
3}.\hspace{7.2cm}(15.1.13)\end{eqnarray*}

Let
$$\phi^A_{\mu_1,\mu_2}=(z_1z_2\cdots z_n)^{\mu_1}W^{\mu_2}(z_1,z_2,...,z_n)\qquad\for\;\;
\mu_1,\mu_2\in\mbb{C}.\eqno(15.1.14)$$ Then
$$z_r\ptl_{z_r}(\phi^A_{\mu_1,\mu_2})=\left(\mu_1-\mu_2\sum_{s=1}^{r-1}\frac{z_r}{z_s-z_r}
+\mu_2\sum_{s=r+1}^n\frac{z_r}{z_r-z_s}\right)\phi^A_{\mu_1,\mu_2}\eqno(15.1.15)$$
for $r\in\ol{1,n}$. Thus
\begin{eqnarray*}& &\sum_{r=1}^n(z_r\ptl_{z_r})^2(\phi^A_{\mu_1,\mu_2})\\ &=&\sum_{r=1}^n
[\mu_1^2-2\mu_1\mu_2\sum_{s=1}^{r-1}\frac{z_r}{z_s-z_r}+2\mu_1\mu_2\sum_{s=r+1}^n\frac{z_r}
{z_r-z_s}-\mu_2\sum_{r\neq s\in\ol{1,n}}\frac{z_sz_r}{(z_s-z_r)^2}\\
&&+\mu_2^2\sum_{r\neq
s\in\ol{1,n}}\frac{z_r^2}{(z_s-z_r)^2}+2\mu_2^2 \sum_{1\leq
s_1<s_2\leq n;\:s_1,s_2\neq
r}\frac{z_r^2}{(z_{s_1}-z_r)(z_{s_2}-z_r)}] \phi^A_{\mu_1,\mu_2}\\
&=&[n\mu_1^2+n(n-1)\mu_1\mu_2+\mu_2^2 {n\choose
3}-2\mu_2\sum_{1\leq r<s\leq n} \frac{z_sz_r}{(z_s-z_r)^2}\\
& &+2\mu_2^2\sum_{1\leq r<s\leq n}\frac{z_r^2+z_s^2}
{(z_s-z_r)^2}]\phi^A_{\mu_1,\mu_2}\\&=&[n\mu_1^2+n(n-1)\mu_1\mu_2+2\mu_2^2
{n\choose 3}-2\mu_2\sum_{1\leq r<s\leq n}
\frac{z_sz_r}{(z_s-z_r)^2}\\ & &+\mu_2^2\sum_{1\leq r<s\leq
n}\frac{z_r^2+z_s^2-2z_rz_s+2z_rz_s}{(z_s-z_r)^2}]\phi^A_{\mu_1,\mu_2}\\
&=&[n\mu_1^2+n(n-1)(\mu_1+\mu_2/2)\mu_2+2\mu_2^2{n\choose 3}\\ &
&+2\mu_2(\mu_2-1)\sum_{1\leq r<s\leq
n}\frac{z_sz_r}{(z_s-z_r)^2}]\phi^A_{\mu_1,\mu_2}\hspace{6.7cm}(15.1.16)\end{eqnarray*}
by (15.1.13) and (15.1.15). Therefore, we have: \psp

{\bf Theorem 15.1.2}. {\it The function $\phi^A_{\mu_1,\mu_2}$
satisfies:
\begin{eqnarray*}\qquad& &\sum_{r=1}^n(z_r\ptl_{z_r})^2(\phi^A_{\mu_1,\mu_2})+2\mu_2(1-\mu_2)
\left(\sum_{1\leq i<j\leq n}\frac{z_iz_j}{(z_i-z_j)^2}\right)\phi^A_{\mu_1,\mu_2}\\
&=&\left[n\mu_1^2+n(n-1)(\mu_1+\mu_2/2)\mu_2+2{n\choose
3}\mu_2^2\right]\phi^A_{\mu_1,\mu_2},\hspace{3.7cm}(15.1.17)\end{eqnarray*}
which is a Calogero-Sutherland equation.} \psp

We remark that above result was  known when $\mu_1=\mu_2$ or
$\mu_1=0$ before our work [X14]. Moreover, ${\cal
W}_{A_{n-1}}=\phi^A_{(1-n)/2,1}$.

 According to simplicity, we secondly consider the Weyl function of  $o(2n,\mbb{C})$. By the root
 structure,
\begin{eqnarray*}\qquad{\cal W}_{D_n}&=&\prod_{1\leq i<j\leq n}(z_i^{1/2}z_j^{-1/2}-z_i^{-1/2}
z_j^{1/2})(z_i^{1/2}z_j^{1/2}-z_i^{-1/2}z_j^{-1/2})\\ &=&
\prod_{1\leq i<j\leq n}(z_i+z_i^{-1}-(z_j+z_j^{-1}))\\
&=&W(z_1+z_1^{-1},z_2+z_2^{-1},...,z_n
+z_n^{-1})\hspace{5.5cm}(15.1.18)\end{eqnarray*}(cf. (15.1.10)).
Note that
\begin{eqnarray*}(z_i\ptl_{z_i})^2[(z_i+z_i^{-1})^k] &=&k(z_i+z_i^{-1})^k+k(k-1)(z_i-z_i^{-1})^2
(z_i+z_i^{-1})^{k-2}\\ &=& k(z_i+z_i^{-1})^k+k(k-1)(z_i^2-2+z_i^{-2})(z_i+z_i^{-1})^{k-2}\\
 &=& k^2(z_i+z_i^{-1})^k-4k(k-1)(z_i+z_i^{-1})^{k-2}.\hspace{3cm}(15.1.19)\end{eqnarray*}
By Lemma 15.1.1 and a similar calculation as (15.1.10),
$$(\sum_{i=1}^n(z_i\ptl_{z_i})^2)({\cal W}_{D_n})=\frac{(n-1)n(2n-1)}{6}{\cal W}_{D_n}.
\eqno(15.1.20)$$ On the other hand, by (15.1.12) and (15.1.19), we
have
\begin{eqnarray*} & &\sum_{r=1}^n\sum_{1\leq s_1<s_2\leq n;\:s_1,s_2\neq r}
\frac{(z_r-z_r^{-1})^2}{(z_{s_1}+z_{s_1}^{-1}-z_r-z_r^{-1})(z_{s_2}+z_{s_2}^{-1}-z_r-z_r^{-1})}
\\ &=&\frac{(n-1)n(n-2)}{6}={n\choose 3}.
\hspace{8.9cm}(15.1.21)\end{eqnarray*} Moreover,
\begin{eqnarray*} & &\frac{(z_r-z_r^{-1})^2+(z_s-z_s^{-1})^2}{(z_s+z_s^{-1}-z_r-z_r^{-1})^2}\\
 &=& \frac{z_r^2+z_r^{-2}+z_s^2+z_s^{-2}-4}{(z_rz_s)^{-2}(z_s-z_r)^2(z_rz_s-1)^2}\\ &=&
 \frac{z_r^4z_s^2+z_s^2+z_r^2z_s^4+z_s^2-4z_r^2z_s^2}{(z_s-z_r)^2(z_rz_s-1)^2}\\ &=&
\frac{1}{(z_s-z_r)^2(z_rz_s-1)^2}
 [z_r^4z_s^2+z_s^2+z_r^2z_s^4+z_s^2-2z_rz_s(z_r^2+z_s^2+(z_rz_s)^2\\ & &-2z_rz_s+1)
 +2z_rz_s(z_r^2+z_s^2+(z_rz_s)^2-4z_rz_s+1)]\\ &=&\frac{(z_s-z_r)^2(z_rz_s-1)^2
 +2z_rz_s[(z_r-z_s)^2+(z_rz_s-1)^2]}{(z_s-z_r)^2(z_rz_s-1)^2}\\ &=&1 +2\frac{z_rz_s}
 {(z_r-z_s)^2}+2\frac{z_rz_s}{(z_rz_s-1)^2}\hspace{7.9cm}(15.1.22)\end{eqnarray*}
and
\begin{eqnarray*}\qquad& &z_r\ptl_{z_r}\left(\frac{z_r-z_r^{-1}}{z_r+z_r^{-1}-z_s-z_s^{-1}}
\right)\\
&=&\frac{(z_r+z_r^{-1})(z_r+z_r^{-1}-z_s-z_s^{-1})-(z_r-z_r^{-1})^2}{(z_r+z_r^{-1}
-z_s-z_s^{-1})^2}\\ &=&
\frac{(z_r+z_r^{-1})^2-(z_r+z_r^{-1})(z_s+z_s^{-1})-(z_r-z_r^{-1})^2}
{(z_r+z_r^{-1}-z_s-z_s^{-1})^2}\\ &=&
 \frac{4-z_rz_s-z_rz_s^{-1}-z_r^{-1}z_s-z_r^{-1}z_s^{-1}}{(z_rz_s)^{-2}(z_s-z_r)^2(z_rz_s-1)^2}
 \\ &=& -z_rz_s\frac{(z_rz_s)^2+z_r^2+z_s^2+1-4z_rz_s}{(z_s-z_r)^2(z_rz_s-1)^2}\\ &=&
 -\frac{z_rz_s}{(z_s-z_r)^2}-\frac{z_rz_s}{(z_rz_s-1)^2}\hspace{7.9cm}(15.1.23)\end{eqnarray*}
for $r,s\in\ol{1,n}$ such that $r\neq s$.

Set
$$\phi_\mu^D=({\cal W}_{D_n})^\mu=\prod_{1\leq i<j\leq n}(z_i+z_i^{-1}-z_j-z_j^{-1})^{\mu}
\qquad\for\;\;\mu\in\mbb{C}.\eqno(15.1.24)$$ Then
\begin{eqnarray*} & &(\sum_{r=1}^n(z_r\ptl_{z_r})^2)(\phi_\mu^D)\\ &=&\sum_{r=1}^n
[\sum_{r\neq s\in\ol{1,n}}\left[\mu
z_r\ptl_{z_r}\left(\frac{z_r-z_r^{-1}}
{z_r+z_r^{-1}-z_s-z_s^{-1}}\right)+\mu^2\frac{(z_r-z_r^{-1})^2}{(z_s+z_s^{-1}-z_r-z_r^{-1})^2}
\right]\\ & &+2\mu^2\sum_{1\leq s_1<s_2\leq n;\:s_1,s_2\neq
r}\frac{(z_r-z_r^{-1})^2}{(z_{s_1}
+z_{s_1}^{-1}-z_r-z_r^{-1})(z_{s_2}+z_{s_2}^{-1}-z_r-z_r^{-1})}]\phi_\mu^D\\
&=& [\sum_{1\leq r<s\leq
n}\left[-2\mu\left(\frac{z_rz_s}{(z_s-z_r)^2}+\frac{z_rz_s}{(z_rz_s-1)^2}
\right)+\mu^2\frac{(z_r-z_r^{-1})^2+(z_s-z_s^{-1})^2}{(z_s+z_s^{-1}-z_r-z_r^{-1})^2}\right]
\phi_\mu^D\\&&+2\mu^2\sum_{r=1}^n\;\sum_{1\leq s_1<s_2\leq
n;\:s_1,s_2\neq r}
\frac{(z_r-z_r^{-1})^2}{(z_{s_1}+z_{s_1}^{-1}-z_r-z_r^{-1})(z_{s_2}+z_{s_2}^{-1}-z_r-z_r^{-1})}]
\\ &=& \left[2\mu(\mu-1)\sum_{1\leq r<s\leq n}\left(\frac{z_rz_s}{(z_s-z_r)^2}+\frac{z_rz_s}
{(z_rz_s-1)^2}\right)+\mu^2\left(\frac{n(n-1)}{2}+2\left(\!\!\begin{array}{c}n\\
3\end{array} \!\!\right)\right)\right]\phi_\mu^D\\
&=&\big[2\mu(\mu-1)\sum_{1\leq r<s\leq
n}\left(\frac{z_rz_s}{(z_s-z_r)^2}+\frac{z_rz_s}{(z_rz_s-1)^2}\right)
\\ & &+\frac{n(n-1)(2n-1)}{6}\mu^2\big]\phi_\mu^D.\hspace{8.8cm}(15.1.25)\end{eqnarray*}
Thus we obtain: \psp

{\bf Theorem 15.1.3}. {\it The function $\phi_\mu^D$ the equation:
\begin{eqnarray*}& &\sum_{r=1}^n(z_r\ptl_{z_r})^2(\phi^D_\mu)+2\mu(1-\mu)
\left(\sum_{1\leq i<j\leq n}\frac{z_iz_j}{(z_i-z_j)^2}+\sum_{1\leq
i<j\leq n}\frac{z_iz_j} {(z_iz_j-1)^2}\right)\phi^D_\mu\\
&=&\frac{n(n-1)(2n-1)\mu^2}{6}\phi_\mu^D,
\hspace{9.2cm}(15.1.26)\end{eqnarray*} which is an
Olshanesky-Perelomov equation}. \psp

We remark that the above result with $n=4$ was given in [NFP].

The Weyl function of $sp(2n,\mbb{C})$ is
$${\cal W}_{C_n}=\left[\prod_{r=1}^n(z_r-z_r^{-1})\right]\prod_{1\leq i<j\leq n}
(z_i+z_i^{-1}-(z_j+z_j^{-1})),\eqno(15.1.27)$$ In terms of
determinant,
$${\cal W}_{C_n}=\left|\begin{array}{ccc}z_1-z_1^{-1}&\cdots& z_n-z_n^{-1}\\ (z_1-z_1^{-1})
(z_1+z_1^{-1})&\cdots& (z_n-z_n^{-1})(z_n+z_n^{-1})\\
(z_1-z_1^{-1})(z_1+z_1^{-1})^2& \cdots&
(z_n-z_n^{-1})(z_n+z_n^{-1})^2\\ \vdots&\vdots&\vdots\\
(z_1-z_1^{-1}) (z_1+z_1^{-1})^{n-1}&\cdots&
(z_n-z_n^{-1})(z_n+z_n^{-1})^{n-1}\end{array}\right|.
\eqno(15.1.28)$$ Moreover,
\begin{eqnarray*} (z_r\ptl_{z_r})^2[(z_r-z_r^{-1})(z_r+z_r^{-1})^k]&=&(k+1)^2(z_r-z_r^{-1})
(z_r+z_r^{-1})^k\\ &
&+4k(k-1)(z_r-z_r^{-1})(z_r+z_r^{-1})^{k-2}\hspace{1.7cm}(15.1.29)
\end{eqnarray*}
by (15.1.19). Thus
$$(\sum_{i=1}^n(z_i\ptl_{z_i})^2)({\cal W}_{C_n})=\frac{n(n+1)(2n+1)}{6}{\cal W}_{C_n}
\eqno(15.1.30)$$ by Lemma 15.1.1 and (15.1.29).

Set
$$\psi^C_\mu=\prod_{r=1}^n(z_r-z_r^{-1})^{\mu}\qquad\for\;\;\mu\in\mbb{C}.\eqno(15.1.31)$$
Then
\begin{eqnarray*}\sum_{r=1}^n(z_r\ptl_{z_r})^2(\psi^C_\mu)&=&\mu\psi^C_\mu\sum_{r=1}^n
\left[z_r\ptl_{z_r}\left(\frac{z_r+z_r^{-1}}{z_r-z_r^{-1}}\right)+\mu\left(\frac{z_r+z_r^{-1}}
{z_r-z_r^{-1}}\right)^2\right]\\ &=&
\mu\psi^C_\mu\sum_{r=1}^n\left[\frac{(z_r-z_r^{-1})^2-(z_r+z_r^{-1})^2}{(z_r-z_r^{-1})^2}
+\mu\frac{(z_r+z_r^{-1})^2}{(z_r-z_r^{-1})^2}\right]\\
&=&\mu\sum_{r=1}^n\left[\mu+4(\mu-1)
\frac{1}{(z_r-z_r^{-1})^2}\right]\psi^C_\mu\\
&=&\left[n\mu^2+4\mu(\mu-1)\sum_{r=1}^n
\frac{z_r^2}{(z_r^2-1)^2}\right]\psi^C_\mu.\hspace{4.1cm}(15.1.32)\end{eqnarray*}
Moreover,
\begin{eqnarray*}\qquad& &\sum_{r=1}^n(z_r\ptl_{z_r})(\psi^C_{\mu_1})(z_r\ptl_{z_r})(\phi_{\mu_2}
^D)\\
&=&\mu_1\mu_2\sum_{r=1}^n\frac{z_r+z_r^{-1}}{z_r-z_r^{-1}}\psi^C_{\mu_1}
\sum_{r\neq
s\in\ol{1,n}}\frac{z_r-z_r^{-1}}{z_r+z_r^{-1}-z_s-z_s^{-1}}\phi_{\mu_2}^D
\\&=&\mu_1\mu_2\psi^C_{\mu_1}\phi_{\mu_2}^D\sum_{r=1}^n\sum_{r\neq
s\in\ol{1,n}} \frac{z_r+z_r^{-1}}{z_r+z_r^{-1}-z_s-z_s^{-1}}\\
&=&\frac{n(n-1)}{2}\mu_1\mu_2\psi^C_{\mu_1}
\phi_{\mu_2}^D.\hspace{9cm}(15.1.33)\end{eqnarray*} Set
$$\phi_{\mu_1,\mu_2}^C=\psi^C_{\mu_1}\phi_{\mu_2}^D=\left[\prod_{r=1}^n(z_r-z_r^{-1})^{\mu_1}
\right]\prod_{1\leq i<j\leq
n}(z_i+z_i^{-1}-z_j-z_j^{-1})^{\mu_2}\eqno(15.1.34)$$ for
$\mu_1,\mu_2\in\mbb{C}$. By Theorem 15.1..3 and (15.1.32)-(15.1.34),
we have: \psp

{\bf Theorem 15.1.4}. {\it The function $\phi^C_{\mu_1,\mu_2}$
satisfies:
\begin{eqnarray*}& &\sum_{r=1}^n(z_r\ptl_{z_r})^2(\phi^C_{\mu_1,\mu_2})+[4\mu_1(1-\mu_1)
\sum_{r=1}^n\frac{z_r^2}{(z_r^2-1)^2}\\ &
&+2\mu_2(1-\mu_2)\left(\sum_{1\leq i<j\leq n}
\frac{z_iz_j}{(z_i-z_j)^2}+\sum_{1\leq i<j\leq
n}\frac{z_iz_j}{(z_iz_j-1)^2}\right)] \phi^C_{\mu_1,\mu_2}\\
&=&\left(n\mu^2_1+n(n-1)\mu_1\mu_2+\frac{n(n-1)(2n-1)\mu_2^2}{6}\right)
\phi^C_{\mu_1,\mu_2},\hspace{3.9cm}(15.1.35)\end{eqnarray*} which is
an Olshanesky-Perelomov equation}. \psp

Observe that ${\cal W}_{C_n}=\phi^C_{1,1}$. The Weyl function of
$o(2n+1,\mbb C)$ is
$${\cal W}_{B_n}=\left[\prod_{r=1}^n(z_r^{1/2}-z_r^{-1/2})\right]\prod_{1\leq i<j\leq n}
(z_i+z_i^{-1}-(z_j+z_j^{-1})),\eqno(15.1.36)$$ In terms of
determinant,
\begin{eqnarray*}& &{\cal W}_{B_n}=\\
& &{\small \left|\begin{array}{ccc}z_1^{1/2}-z_1^{-1/2}&\cdots&
z_n^{1/2}-z_n^{-1/2}\\ (z_1^{1/2}
-z_1^{-1/2})(z_1^{1/2}+z_1^{-1/2})^2&\cdots& (z_n^{1/2}-z_n^{-1/2})(z_n^{1/2}+z_n^{-1/2})^2\\
(z_1^{1/2}-z_1^{-1/2})(z_1^{1/2}+z_1^{-1/2})^4&\cdots&
(z_n^{1/2}-z_n^{-1/2})(z_n^{1/2} +z_n^{-1/2})^4\\
\vdots&\vdots&\vdots\\
(z_1^{1/2}-z_1^{-1/2})(z_1^{1/2}+z_1^{-1/2})^{2(n-1)} &\cdots&
(z_n^{1/2}-z_n^{-1/2})(z_n^{1/2}+z_n^{-1/2})^{2(n-1)}\end{array}\right|}.
\hspace{0.6cm}(15.1.37)\end{eqnarray*} Moreover,
\begin{eqnarray*}\qquad & &(z_r\ptl_{z_r})^2[(z_r^{1/2}-z_r^{-1/2})(z_r^{1/2}+z_r^{-1/2})^{2k}]
\\ &=&\left(k^2+k+\frac{1}{4}\right)(z_r^{1/2}-z_r^{-1/2})(z_r^{1/2}+z_r^{-1/2})^{2k}\\
&
&+2k(2k-1)(z_r^{1/2}-z_r^{-1/2})(z_r^{1/2}+z_r^{-1/2})^{2(k-1)}.\hspace{4.4cm}(15.1.38)
\end{eqnarray*}
 Thus
$$(\sum_{i=1}^n(z_i\ptl_{z_i})^2)({\cal W}_{B_n})=\frac{n(4n^2-1)}{12}{\cal W}_{B_n}\eqno(15.1.39)$$
by Lemma 15.1.1 and (15.1.38).

Set
$$\psi^B_\mu=\prod_{r=1}^n(z_r^{1/2}-z_r^{-1/2})^{\mu}\qquad\for\;\;\mu\in\mbb{C}.\eqno(15.1.40)$$
Then
\begin{eqnarray*}& &\sum_{r=1}^n(z_r\ptl_{z_r})^2(\psi^B_\mu)\\&=&\frac{\mu}{2}\psi^B_\mu
\sum_{r=1}^n\left[z_r\ptl_{z_r}\left(\frac{z_r^{1/2}+z_r^{-1/2}}{z_r^{1/2}-z_r^{-1/2}}\right)
+\frac{\mu}{2}\left(\frac{z_r^{1/2}+z_r^{-1/2}}{z_r^{1/2}-z_r^{-1/2}}\right)^2\right]\\
&=&
\frac{\mu}{4}\psi^B_\mu\sum_{r=1}^n\left[\frac{(z_r^{1/2}-z_r^{-1/2})^2-(z_r^{1/2}+
z_r^{-1/2})^2}{(z_r^{1/2}-z_r^{-1/2})^2}+\mu\frac{(z_r^{1/2}+z_r^{-1/2})^2}{(z_r^{1/2}
-z_r^{-1/2})^2}\right]\\
&=&\left[\frac{n\mu^2}{4}+\mu(\mu-1)\sum_{r=1}^n\frac{z_r}
{(z_r-1)^2}\right]\psi^B_\mu.\hspace{6.9cm}(15.1.41)\end{eqnarray*}
Moreover,
\begin{eqnarray*}\qquad& &\sum_{r=1}^n(z_r\ptl_{z_r})(\psi^B_{\mu_1})(z_r\ptl_{z_r})
(\phi_{\mu_2}^D)\\
&=&\frac{\mu_1\mu_2}{2}\sum_{r=1}^n\frac{z_r^{1/2}+z_r^{-1/2}}
{z_r^{1/2}-z_r^{-1/2}}\psi^B_{\mu_1}\sum_{r\neq
s\in\ol{1,n}}\frac{z_r-z_r^{-1}}
{z_r+z_r^{-1}-z_s-z_s^{-1}}\phi_{\mu_2}^D\\
&=&\frac{\mu_1\mu_2}{2}\psi^B_{\mu_1}
\phi_{\mu_2}^D\sum_{r=1}^n\sum_{r\neq
s\in\ol{1,n}}\frac{(z_r^{1/2}+z_r^{-1/2})^2}
{z_r+z_r^{-1}-z_s-z_s^{-1}}\\
&=&\frac{n(n-1)}{4}\mu_1\mu_2\psi^B_{\mu_1}\phi_{\mu_2}^D.
\hspace{8.7cm}(15.1.42)\end{eqnarray*} Set
$$\phi_{\mu_1,\mu_2}^B=\psi^B_{\mu_1}\phi_{\mu_2}^D=\left[\prod_{r=1}^n(z_r^{1/2}
-z_r^{-1/2})^{\mu_1}\right]\prod_{1\leq i<j\leq
n}(z_i+z_i^{-1}-z_j-z_j^{-1})^{\mu_2} \eqno(15.1.43)$$ for
$\mu_1,\mu_2\in\mbb{C}$. By Theorem 15.1.3 and (15.1.41)-(15.1.43),
we obtain: \psp

{\bf Theorem 15.1.5}. {\it The function $\phi^B_{\mu_1,\mu_2}$
satisfies:
\begin{eqnarray*}& &\sum_{r=1}^n(z_r\ptl_{z_r})^2(\phi^B_{\mu_1,\mu_2})+
[\mu_1(1-\mu_1)\sum_{r=1}^n\frac{z_r}{(z_r-1)^2}\\ &
&+2\mu_2(1-\mu_2) \left(\sum_{1\leq i<j\leq
n}\frac{z_iz_j}{(z_i-z_j)^2}+\sum_{1\leq i<j\leq n}\frac{z_iz_j}
{(z_iz_j-1)^2}\right)]\phi^C_{\mu_1,\mu_2}\\
&=&\left(\frac{n\mu^2_1}{4}+\frac{n(n-1)}{2}
\mu_1\mu_2+\frac{n(n-1)(2n-1)\mu_2^2}{6}\right)\phi^B_{\mu_1,\mu_2},\hspace{3.7cm}(15.1.44)
\end{eqnarray*}
which is an Olshanesky-Perelomov equation}.

\section{Etingof Traces}

In this section, we present the trace functions of the intertwining
operators among certain modules  introduced by Etingof [Ep].

Recall that the general Lie algebra
$$gl(n,\mbb R)=\sum_{i,j=1}^n\mbb RE_{i,j}.\eqno(15.2.1)$$
The subspace
$$H_A=\sum_{i=1}^n\mbb RE_{i,i}\eqno(15.2.2)$$
is a toral Cartan subalgebra of the Lie algebra $gl(n,\mbb R)$. For
a weight $\lmd\in H_A^\ast$, we denote
$$\lmd(E_{i,i})=\lmd_i\qquad \for\;\;i\in\ol{1,n}.\eqno(15.2.3)$$
Note that
$$\omega_{_A}=\sum_{i,j=1}^nE_{i,j}E_{j,i}\eqno(15.2.4)$$
is the Casimier element in $U(gl(n,\mbb R)$. Let
$$M=\bigoplus_{\lmd\in H_A^\ast}M_\lmd\eqno(15.2.5)$$
be any weight $gl(n,\mbb R)$-module with the representation
$\pi_{_M}$ such that
$$\pi_{_M}(\omega_{_A})=b\:\mbox{Id}_M\qquad\mbox{for some}\;\;b\in \mbb{C},\eqno(15.2.6)$$
where
$$M^\lmd=\{w\in M\mid h(w)=\lmd (h)w\;\for\;h\in H\}.\eqno(15.2.7)$$
In particular, we can take $M$ to be any highest weight module; say
the Verma module $M(\lmd)$ given in (5.2.8)-(5.2.10). The module $M$
is not necessarily irreducible.

Let $\{x_1,x_2,...,x_n\}$ be $n$ indeterminates. Set
$$x^{\vec i}=x_1^{i_1}x_2^{i_2}\cdots x_n^{i_n}\qquad\vec
i=(i_1,i_2,...,i_n)\in\mbb{Z}^n. \eqno(15.2.8)$$ Fix $\kappa\in\mbb
R$ and denote
$$\bar x^\kappa=(x_1x_2\cdots x_n)^\kappa.\eqno(15.2.9)$$
Define the space of formal Laurent series in $\{x_1,x_2,...,x_n\}$
with the coefficients in
 $M$:
$$\td{M}=\{\sum_{\vec i\in\mbb{Z}^n}w_{\vec i}x^{\vec i}\bar x^\kappa\mid w_{\vec i}\in M\}.\eqno(15.2.10)$$
We extend the representation $\pi_{_M}$ to $\td M$ by
$$\pi_{_M}(\xi)(\sum_{\vec i\in\mbb{Z}^n}w_{\vec i}x^{\vec i}\bar
x^\kappa)=\sum_{\vec i\in\mbb{Z}^n}\pi_{_M}(\xi)(w_{\vec i})x^{\vec
i}\bar x^\kappa\;\;\for\;\xi\in gl(n,\mbb R).\eqno(15.2.11)$$
 Moreover, we define another representation $\pi'$ of
$gl(n,\mbb R)$ on $\td{M}$ by
$$\pi'(E_{r,s})(\sum_{\vec i\in\mbb{Z}^n}w_{\vec i}x^{\vec i}\bar x^\kappa)=\sum_{\vec i\in\mbb{Z}^n}
[ w_{\vec i}x_r\ptl_{x_s}(x^{\vec i}\bar x^\kappa)-\dlt_{r,s}\kappa
w_{\vec i}x^{\vec i} \bar x^\kappa]\eqno(15.2.12)$$ for
$r,s\in\ol{1,n}$. Define the representation of $gl(n,\mbb R)$ on
$\td{M}$ by
$$\pi=\pi_{_M}+\pi'.\eqno(15.2.13)$$

 Suppose that $\Phi:M\rta \td{M}$ is a linear map
such that
$$\Phi(\pi_{_M}(\xi)(w))=\Dlt(\xi)\Phi(w)\eqno(15.2.14)$$
for $\xi\in gl(n,\mbb R)$ and $w\in M$. View $\Phi$ as a function in
$\{x_1,x_2,...,x_n\}$ taking value
 in the spaces of  linear transformations on $M$.
 The {\it Etingof trace function}
$$E_A(z_1,z_2,...,z_n)=\mbox{tr}_{_{M}}\:\Phi z_1^{E_{1,1}}
z_2^{E_{2,2}}\cdots z_n^{E_{n,n}}\eqno(15.2.15)$$ (cf. [Ep]). By
considering weights, we can prove
$$(\bar x^\kappa)^{-1}E_A(z_1,z_2,...,z_n)\;\mbox{is
independent of}\;x_1,x_2,...,x_n.\eqno(15.2.16)$$ Etingof [Ep]
proved: \psp

{\bf Theorem 15.2.1}. {\it The function
$$\Psi_A(z_1,z_2,...,z_n)={\cal W}_{A_{n-1}}(\bar x^\kappa)^{-1}E_A(z_1,z_2,...,z_n)\eqno(15.2.17)$$ is a solution of
(15.1.4) for suitable constants $K$ and $\nu$.}

{\it Proof}. Recall ${\cal W}_{A_{n-1}}=\phi^A_{(1-n)/2,1}$. By
(15.1.15), $$z_r\ptl_{z_r}({\cal W}_{A_{n-1}})=\left(\frac{1-n}{2}
+\sum_{s\neq
r}\frac{z_r}{z_r-z_s}\right){\cal W}_{A_{n-1}}\\
=\left(\frac{n-1}{2} +\sum_{s\neq r}\frac{z_s}{z_r-z_s}\right){\cal
W}_{A_{n-1}}\eqno(15.2.18)$$ for $r\in\ol{1,n}$. Set
$$E_A^{i_1,i_2}=\mbox{tr}_M\Phi \pi_{_M}( E_{i_1,i_2}E_{i_2,i_1})
z_1^{\pi_{_M}(E_{1,1})}z_2^{\pi_{_M}(E_{2,2})}\cdots
z_n^{\pi_{_M}(E_{n,n})}\eqno(15.2.19)$$ for $i_1,i_2\in\ol{1,n}$.
Observe that
\begin{eqnarray*}\qquad\qquad& &\mbox{tr}_M\Phi \pi_{_M}( E_{i,i}))
z_1^{\pi_{_M}(E_{1,1})}z_2^{\pi_{_M}(E_{2,2})}\cdots
z_n^{\pi_{_M}(E_{n,n})}\\&=&\mbox{tr}_M\Phi
z_1^{\pi_{_M}(E_{1,1})}z_2^{\pi_{_M}(E_{2,2})}\cdots
z_n^{\pi_{_M}(E_{n,n})}\pi_{_M}( E_{i,i}))\\
&=&z_i\ptl_{z_i}[\mbox{tr}_M\Phi
z_1^{\pi_{_M}(E_{1,1})}z_2^{\pi_{_M}(E_{2,2})}\cdots
z_n^{\pi_{_M}(E_{n,n})}]\hspace{4.5cm}(15.2.20)\end{eqnarray*} for
$i\in\ol{1,n}$. In particular,
$$E_A^{i,i}=(z_i\ptl_{z_i})^2 (E_A)\eqno(15.2.21)$$
and for $i_1\neq i_2$,
\begin{eqnarray*}  E_A^{i_1,i_2}&=&\mbox{tr}_M\pi(E_{i_1,i_2})\Phi \pi_{_M}
(E_{i_2,i_1})z_1^{\pi_{_M}(E_{1,1})}z_2^{\pi_{_M}(E_{2,2})}\cdots
z_n^{\pi_{_M}(E_{n,n})}
\\ &=& \mbox{tr}_M\pi'(E_{i_1,i_2})\Phi \pi_{_M}(E_{i_2,i_1})z_1^{\pi_{_M}(E_{1,1})}
z_2^{\pi_{_M}(E_{2,2})}\cdots  z_n^{\pi_{_M}(E_{n,n})}\\ &
&+\mbox{tr}_M\pi_{_M}(E_{i_1,i_2}) \Phi
\pi_{_M}(E_{i_2,i_1})z_1^{\pi_{_M}(E_{1,1})}z_2^{\pi_{_M}(E_{2,2})}\cdots
z_n^{\pi_{_M} (E_{n,n})}\\ &=&\mbox{tr}_M\pi'(E_{i_1,i_2})\Phi
\pi_{_M}(E_{i_2,i_1})z_1^{\pi_{_M}(E_{1,1})}
z_2^{\pi_{_M}(E_{2,2})}\cdots  z_n^{\pi_{_M}(E_{n,n})}\\ &
&+\mbox{tr}_M \Phi
\pi_{_M}(E_{i_2,i_1})z_1^{\pi_{_M}(E_{1,1})}z_2^{\pi_{_M}(E_{2,2})}\cdots
z_n^{\pi_{_M}(E_{n,n})}\pi_{_M}(E_{i_1,i_2})\\
&=&\mbox{tr}_M\pi'(E_{i_1,i_2}) \Phi
\pi_{_M}(E_{i_2,i_1})z_1^{\pi_{_M}(E_{1,1})}z_2^{\pi_{_M}(E_{2,2})}\cdots
z_n^{\pi_{_M}(E_{n,n})}\\ & &+\frac{z_{i_1}}{z_{i_2}}\mbox{tr}_M\Phi
\pi_{_M}(E_{i_2,i_1}
E_{i_1,i_2})z_1^{\pi_{_M}(E_{1,1})}z_2^{\pi_{_M}(E_{2,2})}\cdots
z_n^{\pi_{_M}(E_{n,n})}
 \\&=& \mbox{tr}_M\pi'(E_{i_1,i_2})\Phi \pi_{_M}(E_{i_2,i_1})z_1^{\pi_{_M}(E_{1,1})}
 z_2^{\pi_{_M}(E_{2,2})}\cdots  z_n^{\pi_{_M}(E_{n,n})}\\ & &+\frac{z_{i_1}}{z_{i_2}}
 \mbox{tr}_M\Phi \pi_{_M}(E_{i_1,i_2}E_{i_2,i_1}+[E_{i_2,i_1},E_{i_1,i_2}])z_1^{\pi_{_M}
 (E_{1,1})}z_2^{\pi_{_M}(E_{2,2})}\cdots
 z_n^{\pi_{_M}(E_{n,n})}\\&=& \mbox{tr}_M\pi'
 (E_{i_1,i_2})\Phi \pi_{_M}(E_{i_2,i_1})z_1^{\pi_{_M}(E_{1,1})}z_2^{\pi_{_M}(E_{2,2})}\cdots
  z_n^{\pi_{_M}(E_{n,n})}\\ & &+\frac{z_{i_1}}{z_{i_2}}\mbox{tr}_M\Phi \pi_{_M}(E_{i_1,i_2}
  E_{i_2,i_1}+E_{i_2,i_2}-E_{i_1,i_1})z_1^{\pi_{_M}(E_{1,1})}z_2^{\pi_{_M}(E_{2,2})}\cdots
   z_n^{\pi_{_M}(E_{n,n})}\\ &=&
\mbox{tr}_M\pi'(E_{i_1,i_2})\Phi
\pi_{_M}(E_{i_2,i_1})z_1^{\pi_{_M}(E_{1,1})}z_2^{\pi_{_M}
(E_{2,2})}\cdots  z_n^{\pi_{_M}(E_{n,n})}+\frac{z_{i_1}}{z_{i_2}} E_A^{i_1,i_2}\\
 & &+\frac{z_{i_1}}{z_{i_2}}(z_{i_2}\ptl_{z_{i_2}}-z_{i_1}\ptl_{z_{i_1}})(E_A).
 \hspace{7.5cm}(15.2.22)\end{eqnarray*}
Thus
\begin{eqnarray*}E_A^{i_1,i_2}&=&\frac{z_{i_2}}{z_{i_2}-z_{i_1}}\mbox{tr}_M
\pi'(E_{i_1,i_2})\Phi
\pi_{_M}(E_{i_2,i_1})z_1^{\pi_{_M}(E_{1,1})}z_2^{\pi_{_M}(E_{2,2})}
\cdots  z_n^{\pi_{_M}(E_{n,n})}\\ &
&+\frac{z_{i_1}}{z_{i_2}-z_{i_1}}(z_{i_2}
\ptl_{z_{i_2}}-z_{i_1}\ptl_{z_{i_1}})
(E_A).\hspace{6.6cm}(15.2.23)\end{eqnarray*} Furthermore,
\begin{eqnarray*}& &\mbox{tr}_M\pi'(E_{i_1,i_2})\Phi \pi_{_M}(E_{i_2,i_1})
z_1^{\pi_{_M}(E_{1,1})}z_2^{\pi_{_M}(E_{2,2})}\cdots
z_n^{\pi_{_M}(E_{n,n})}\\ &=&
\mbox{tr}_M\pi'(E_{i_1,i_2})\pi(E_{i_2,i_1})\Phi
z_1^{\pi_{_M}(E_{1,1})}z_2^{\pi_{_M} (E_{2,2})}\cdots
z_n^{\pi_{_M}(E_{n,n})}\\ &=&
 \mbox{tr}_M\pi'(E_{i_1,i_2})\pi'(E_{i_2,i_1})\Phi z_1^{\pi_{_M}(E_{1,1})}z_2^{\pi_{_M}
 (E_{2,2})}\cdots  z_n^{\pi_{_M}(E_{n,n})}\\ & &+ \mbox{tr}_M\pi'(E_{i_1,i_2})\pi_{_M}
 (E_{i_2,i_1})\Phi z_1^{\pi_{_M}(E_{1,1})}z_2^{\pi_{_M}(E_{2,2})}\cdots  z_n^{\pi_{_M}
 (E_{n,n})}\\ &=&x_{i_1}\ptl_{x_{i_2}}x_{i_2}\ptl_{x_{i_1}} (E_A)+
 \mbox{tr}_M\pi'(E_{i_1,i_2})\Phi z_1^{\pi_{_M}(E_{1,1})}z_2^{\pi_{_M}(E_{2,2})}\cdots
 z_n^{\pi_{_M}(E_{n,n})}\pi_{_M}(E_{i_2,i_1})\\ &=& \frac{z_{i_2}}
 {z_{i_1}}\mbox{tr}_M\pi'(E_{i_1,i_2})\Phi \pi_{_M}(E_{i_2,i_1}) z_1^{\pi_{_M}(E_{1,1})}
 z_2^{\pi_{_M}(E_{2,2})}\cdots  z_n^{\pi_{_M}(E_{n,n})}\\ & &+\kappa(\kappa+1)E_A\hspace{11cm}(15.2.24)\end{eqnarray*}
by (15.2.18). Hence
$$\mbox{tr}_M\pi'(E_{i_1,i_2})\Phi \pi_{_M}(E_{i_2,i_1})
z_1^{\pi_{_M}(E_{1,1})}z_2^{\pi_{_M}(E_{2,2})}\cdots
z_n^{\pi_{_M}(E_{n,n})}=
\frac{\kappa(\kappa+1)z_{i_1}}{z_{i_1}-z_{i_2}}E_A.\eqno(15.2.25)$$
Substituting it to (15.2.23), we obtain
$$E_A^{i_1,i_2}=\frac{z_{i_1}}{z_{i_2}-z_{i_1}}(z_{i_2}\ptl_{z_{i_2}}-z_{i_1}
\ptl_{z_{i_1}})(E_A)-\frac{\kappa(\kappa+1)z_{i_1}z_{i_2}}{(z_{i_1}-z_{i_2})^2}E_A\eqno(15.2.26)$$
for $i_1,i_2\in\ol{1,n}$ such that $i_1\neq i_2.$

According to (15.1.17) with $\mu_1=(1-n)/2$ and $\mu_2=1$, we have
$$\sum_{r=1}^n(z_r\ptl_{z_r})^2({\cal W}_{A_{n-1}})
=\left[\frac{1}{2}+2{n\choose 3}\right]{\cal
W}_{A_{n-1}}\eqno(15.2.27)$$ Now (15.2.6) yields
\begin{eqnarray*} b\Psi_A&=&{\cal W}_{A_{n-1}}(x_1x_2\cdots x_n)^{-\kappa}
\mbox{tr}_{_{M}}\:\Phi\pi_{_M}(\omega_{_A})z_1^{E_{1,1}}
z_2^{E_{2,2}}\cdots z_n^{E_{n,n}}\\&=&\sum_{i_1,i_2=1}^n {\cal
W}_{A_{n-1}}(x_1x_2\cdots x_n)^{-\kappa}E_A^{i_1,i_2}
\\&=&\sum_{i=1}^n{\cal
W}_{A_{n-1}}(x_1x_2\cdots
x_n)^{-\kappa}(z_i\ptl_{z_i})^2(E_A)-\sum_{i_1\neq
i_2}\frac{\kappa(\kappa+1)z_{i_1}z_{i_2}}{(z_{i_1}-z_{i_2})^2}\Psi_A\\&
&+\sum_{i_1\neq i_2}{\cal W}_{A_{n-1}}(x_1x_2\cdots
x_n)^{-\kappa}\frac{z_{i_1}}{z_{i_2}-z_{i_1}}(z_{i_2}\ptl_{z_{i_2}}-z_{i_1}
\ptl_{z_{i_1}})(E_A)
\\&=&
\sum_{i=1}^n{\cal W}_{A_{n-1}}(x_1x_2\cdots
x_n)^{-\kappa}(z_i\ptl_{z_i})^2(E_A)-\sum_{i_1\neq
i_2}\frac{\kappa(\kappa+1)z_{i_1}z_{i_2}}{(z_{i_1}-z_{i_2})^2}\Psi_A\\&
&+\sum_{i_2=1}^n((1-n)/2+z_{i_2}\ptl_{z_{i_2}})({\cal
W}_{A_{n-1}})(x_1x_2\cdots
x_n)^{-\kappa}z_{i_2}\ptl_{z_{i_2}}(E_A)\\ & &+
\sum_{i_1=1}^n((n-1)/2+z_{i_1}\ptl_{z_{i_1}})({\cal
W}_{A_{n-1}})(x_1x_2\cdots x_n)^{-\kappa}z_{i_1}\ptl_{z_{i_1}}(E_A)
\\&=&
\sum_{i=1}^n{\cal W}_{A_{n-1}}(x_1x_2\cdots
x_n)^{-\kappa}(z_i\ptl_{z_i})^2(E_A)-\sum_{i_1\neq
i_2}\frac{\kappa(\kappa+1)z_{i_1}z_{i_2}}{(z_{i_1}-z_{i_2})^2}\Psi_A\\&
&+\sum_{i=1}^n2z_i\ptl_{z_i}({\cal W}_{A_{n-1}})(x_1x_2\cdots
x_n)^{-\kappa}z_i\ptl_{z_i}(E_A)
\\&=&
\sum_{i=1}^n(z_i\ptl_{z_i})^2(\Psi_A)-\left[\frac{1}{2}+2{n\choose
3}\right]\Psi_A-\sum_{i_1\neq
i_2}\frac{\kappa(\kappa+1)z_{i_1}z_{i_2}}{(z_{i_1}-z_{i_2})^2}\Psi_A\hspace{2cm}(15.2.28)
\end{eqnarray*}
by (15.2.4), (15.2.6), (15.2.18), (15.2.19), (15.2.21), (15.2.26)
and (15.2.27). Thus
$$\sum_{i=1}^n(z_i\ptl_{z_i})^2(\Psi_A)-2\kappa(\kappa+1)\left(\sum_{1\leq i<j\leq n}\frac{z_iz_j}{(z_i-z_j)^2}
\right)\Psi_A= \left[b+\frac{1}{2}+2{n\choose
3}\right]\Psi_A,\eqno(15.2.29)$$ which shows that $\Psi_A$ is a
solution of (15.1.4) for
$$K=-2\kappa(\kappa+1),\;\;\nu=b+\frac{1}{2}+2{n\choose
3}.\qquad\Box\eqno(15.2.30)$$\pse

For convenience, we denote
$$C_{i,j}=E_{i,j}-E_{n+j,n+i}\qquad\for\;\;i,j\in\ol{1.n},\eqno(15.2.31)$$
$$C_{j,n+k}=E_{j,n+k}+E_{k,n+j},\qquad C_{n+j,k}=E_{n+j,k}+E_{n+k,j}\eqno(15.2.32)$$
for $j,k\in\ol{1,n}$ such that $j\neq k$, and
$$C_{r,n+r}=E_{r,n+r},\qquad C_{n+r,r}=E_{n+r,r}\qquad\for\;\;r\in\ol{1,n}.\eqno(15.2.33)$$
Then the sympletic Lie algebra
$$sp(2n,\mbb R)=\sum_{i,j=1}^n\mbb RC_{i,j}+\sum_{1\leq r\leq s\leq n}(\mbb RC_{r,n+s}+\mbb RC_{n+r,s}).\eqno(15.2.34)$$
Moreover,
\begin{eqnarray*}\qquad\qquad \omega_{_C}&=&\sum_{i,j=1}^nC_{i,j}C_{j,i}+\sum_{1\leq s<r\leq n}
(C_{n+r,s}C_{s,n+r}+C_{s,n+r}C_{n+r,s})\\ &
&+2\sum_{p=1}^n(C_{n+p,p}C_{p,n+p}+C_{p,n+p}
C_{n+p,p})\hspace{4.7cm}(15.2.35)\end{eqnarray*} is a Casimier
element in $U(sp(2n,\mbb R))$. Now
$$H_C=\sum_{i=1}^n\mbb RC_{i,i}\eqno(15.2.36)$$
is a toral Cartan subalgebra of $sp(2n,\mbb R)$.  Let
$$M=\bigoplus_{\lmd\in H_C^\ast}M_\lmd\eqno(15.2.37)$$
be a weight $sp(2n,\mbb R)$-module with the representation
$\pi_{_M}$ such that
$$\pi_{_M}(\omega_{_C})=b \:\mbox{Id}_M\qquad\mbox{for some}\;\;b\in \mbb R,\eqno(15.2.38)$$
In particular, we can take $M$ to be any highest weight module; say
the Verma module $M(\lmd)$ given in (5.2.8)-(5.2.10). The module $M$
is not necessarily irreducible.

Denote
$$x^\ast=x_1^{-1/2}x_2^{-1/2}\cdots x_n^{-1/2}.\eqno(15.2.39)$$
Denote
 the space of formal Laurent series in $\{x_1,x_2,...,x_n\}$ with the coefficients in $M$ by
$$\td{M}=\{\sum_{\vec i\in\mbb{Z}^n}w_{\vec i}x^{\vec i}x^\ast\mid w_n\in M\}.\eqno(15.2.40)$$
We extend the $\pi_{_M}$ to $\td{M}$ by
$$\pi_{_M}(u)(\sum_{\vec i\in\mbb{Z}^n}w_{\vec i}x^{\vec i}x^\ast)=\sum_{\vec i\in\mbb{Z}^n}
\pi_{_M}(\xi)(w_{\vec i})x^{\vec i}x^\ast\qquad\for\;\;\xi\in
sp(2n,\mbb R).\eqno(15.2.41)$$ Moreover, we define another
representation $\pi'$ of $sp(2n,\mbb R)$ on $\td{M}$ by
$$\pi'(C_{p,q})(\sum_{\vec i\in\mbb{Z}^n}w_{\vec i}x^{\vec i}x^\ast)=\sum_{\vec i\in\mbb{Z}^n}
\left(C_{p,q}(w_{\vec i}x_p\ptl_{x_q}(x^{\vec
i}x^\ast)+\frac{\dlt_{p,q}}{2}w_{\vec i}x^{\vec i}
x^\ast\right),\eqno(15.2.42)$$
$$\pi'(C_{p,n+q})(\sum_{\vec i\in\mbb{Z}^n}w_{\vec i}x^{\vec i}x^\ast)=-\frac{1}{1+\dlt_{p,q}}
\sum_{\vec i\in\mbb{Z}^n}w_{\vec i}x_px_qx^{\vec
i}x^\ast,\eqno(15.2.43)$$
$$\pi'(C_{n+p,q})(\sum_{\vec i\in\mbb{Z}^n}w_{\vec i}x^{\vec i}x^\ast)=
\frac{1}{1+\dlt_{p,q}}\sum_{\vec i\in\mbb{Z}^n}w_{\vec
i}\ptl_{x_p}\ptl_{x_q} (x^{\vec i}x^\ast)\eqno(15.2.44)$$ for
$p,q\in\ol{1 n}$. Set
$$\pi=\pi_{_M}+\pi'.\eqno(15.2.45)$$

Suppose that $\Phi:M\rta \td{M}$ is a linear map such that
$$\Phi(\pi_{_M}(u)(w))=\pi(u)\Phi(w)\eqno(15.2.46)$$
for $u\in sp(2n,\mbb R)$ and $w\in M$. View $\Phi$ as a function in
$\{x_1,x_2,...,x_n\}$ taking value
 in the spaces of  linear transformations on $M_\lmd$. Define the  trace function:
$$E_C(z_1,z_2,...,z_n)=\mbox{tr}_M\Phi z_1^{\pi_{_M}(C_{1,1})} z_2^{\pi_{_M}(C_{2,2})}
\cdots  z_n^{\pi_{_M}(C_{n,n})}.\eqno(15.2.47)$$ By considering
weights, we can prove
$$(x^\ast)^{-1}E_C(z_1,z_2,...,z_n)\;\mbox{is
independent of}\;x_1,x_2,...,x_n.\eqno(15.2.48)$$ The following
theorem was essentially proved by Etingof and Styrkas [ES].\psp

{\bf Theorem 15.2.2}. {\it The function
$$\Psi_C={\cal W}_{C_n}(x^\ast)^{-1}E_C(z_1,z_2,...,z_n)\eqno(15.2.49)$$
satisfies the Olshanesky-Perelomov equation (15.1.5) with $K_3=0$
and proper constants $K_1,K_2,K_4$.}

{\it Proof}. Note that (15.1.27) implies
\begin{eqnarray*}z_i\ptl_{z_i}({\cal W}_{C_n})&=&\sum_{i=1}^n\left[\frac{z_i+z_i^{-1}}{z_i
-z_i^{-1}}+\sum_{i\neq
r\in\ol{1,n}}\frac{z_i-z_i^{-1}}{z_i+z_i^{-1}-z_r-z_r^{-1}}\right]
{\cal W}_{C_n}\\
&=&\sum_{i=1}^n\left[\frac{z_i^2+1}{z_i^2-1}+\sum_{i\neq
r\in\ol{1,n}}
\left(\frac{z_i}{z_i-z_r}+\frac{1}{z_iz_r-1}\right)\right]{\cal
W}_{C_n}.\hspace{2.5cm} (15.2.50)\end{eqnarray*}

Set
$$E_C^{i_1,i_2}=(x^\ast)^{-1}\mbox{tr}_M\Phi \pi_{_M}( C_{i_1,i_2}C_{i_2,i_1})
z_1^{\pi_{_M}(C_{1,1})}z_2^{\pi_{_M}(C_{2,2})}\cdots
z_n^{\pi_{_M}(C_{n,n})}\eqno(15.2.51)$$ for $i_1,i_2\in\ol{1,n}$.
Then
$$E_C^{i,i}=(z_i\ptl_{z_i})^2 (E_C)\eqno(15.2.52)$$
by (15.2.22), and
$$E_C^{i_1,i_2}=\frac{z_{i_1}}{z_{i_2}-z_{i_1}}(z_{i_2}\ptl_{z_{i_2}}-z_{i_1}
\ptl_{z_{i_1}})(E_C)+\frac{z_{i_1}z_{i_2}}{4(z_{i_1}-z_{i_2})^2}E_C\eqno(15.2.53)$$
for $i_1,i_2\in\ol{1,n}$ such that $i_1\neq i_2$ by
(15.2.24)-(15.2.28).

Let
$$E_C^{s,n+r}=(x^\ast)^{-1}\mbox{tr}_M\Phi \pi_{_M}( C_{s,n+r}C_{n+r,s})z_1^{\pi_{_M}
(C_{1,1})}z_2^{\pi_{_M}(C_{2,2})}\cdots
z_n^{\pi_{_M}(C_{n,n})}\eqno(15.2.54)$$ and
$$E_C^{n+r,s}=(x^\ast)^{-1}\mbox{tr}_M\Phi \pi_{_M}(C_{n+r,s}C_{s,n+r})z_1^{\pi_{_M}
(C_{1,1})}z_2^{\pi_{_M}(C_{2,2})}\cdots
z_n^{\pi_{_M}(C_{n,n})}\eqno(15.2.55)$$ for $1\leq s\leq r\leq n$.
By similar calculations as (15.2.22)-(15.2.26), we have
$$E_C^{s,n+r}=\frac{z_rz_s}{z_rz_s-1}(z_r\ptl_{z_r}+z_s\ptl_{z_s})E_C+\frac{z_rz_s}
{4(z_rz_s-1)^2}E_C.\eqno(15.2.56)$$
$$E_C^{n+r,s}=\frac{1}{z_rz_s-1}(z_r\ptl_{z_r}+z_s\ptl_{z_s})E_C+\frac{z_rz_s}
{4(z_rz_s-1)^2}E_C\eqno(15.2.57)$$ for $1\leq s<r\leq n$ and
$$E_C^{p,n+p}=\frac{z_p^2}{z_p^2-1}z_p\ptl_{z_p}E_C+\frac{3z_p^2}{16(z_p^2-1)^2}
E_C.\eqno(15.2.58)$$
$$E_C^{n+p,p}=\frac{1}{z_p^2-1}z_p\ptl_{z_p}E_C+\frac{3z_p^2}{16(z_p^2-1)^2}
E_C\eqno(15.2.59)$$ for $p\in\ol{1,n}$. Thus
$$(\sum_{i=1}^n(z_i\ptl_{z_i})^2({\cal
W}_{C_n})=\frac{n(n+1)(2n+1)}{6}{\cal W}_{C_n} \eqno(15.2.60)$$ by
(15.1.30).

 Now (15.2.38) gives
\begin{eqnarray*}b \Psi&=&{\cal W}_{C_n}(x^\ast)^{-1}\mbox{tr}_M\Phi \pi_{_M}(\omega_{_C}) z_1^{\pi_{_M}
(C_{1,1})} z_2^{\pi_{_M}(C_{2,2})}\cdots  z_n^{\pi_{_M}{C_{n,n}}}
\\ &=&{\cal W}_{C_n}x^\ast)^{-1}[\sum_{i_1,i_2=1}^nE_C^{i_1,i_2}+\sum_{1\leq s\leq r\leq n}
(\dlt_{r,s}+1)(E_C^{n+r,s}+E_C^{s,n+r})]\\ &=& {\cal
W}_{C_n}\big[\sum_{i_1,i_2\in\ol{1,n}\:i_1\neq i_2}
\left(\frac{z_{i_1}}{z_{i_2}-z_{i_1}}(z_{i_2}\ptl_{z_{i_2}}-z_{i_1}\ptl_{z_{i_1}})
(E_C) +\frac{z_{i_1}z_{i_2}}{4(z_{i_1}-z_{i_2})^2}E_C\right)\\ &
&+\sum_{1\leq s<r\leq n}
\left(\frac{z_rz_s+1}{z_rz_s-1}(z_r\ptl_{z_r}+z_s\ptl_{z_s})(E_C)+\frac{z_rz_s}{2(z_rz_s-1)^2}E_C\right)\\&
&+2\sum_{p=1}^n\left(\frac{z_p^2+1}{z_p^2-1}z_p\ptl_{z_p}{\cal
E}+\frac{3z_p^2}{8(z_p^2-1)^2}E_C\right)+\sum_{i=1}^n(z_i\ptl_{z_i})^2(E_C)\big]\\
&=&\left[\sum_{1\leq i_1<i_2\leq
n}\frac{1}{2}\left(\frac{z_{i_1}z_{i_2}}{(z_{i_1}-z_{i_2})^2}+
\frac{z_{i_1}z_{i_2}}{(z_{i_1}z_{i_2}-1)^2}\right)+\frac{3}{4}\sum_{p=1}^n
\frac{z_p^2}{(z_p^2-1)^2}\right]\Psi_C\\
& &+ 2{\cal W}_{C_n}(x^\ast)^{-1}\sum_{r,s\in\ol{1,n};\;r\neq
s}\left(\frac{z_r}{z_r-z_s}+\frac{1}{z_rz_s-1}\right)z_r\ptl_{z_r}(E_C)\\
& &+2{\cal
W}_{C_n}(x^\ast)^{-1}\sum_{p=1}^n\frac{z_p^2+1}{z_p^2-1}z_p\ptl_{z_p}(E_C)+{\cal
W}_{C_n}(x^\ast)^{-1}\sum_{i=1}^n(z_i\ptl_{z_i})^2(E_C)\\&=&
\left[\sum_{1<i_1<i_2\leq
n}\frac{1}{2}\left(\frac{z_{i_1}z_{i_2}}{(z_{i_1}-z_{i_2})^2}+\frac{z_{i_1}z_{i_2}}
{(z_{i_1}z_{i_2}-1)^2}\right)+\frac{3}{4}\sum_{p=1}^n\frac{z_p^2}{(z_p^2-1)^2}\right]\Psi_C\\
& &+\sum_{i=1}^n[{\cal
W}_{C_n}(x^\ast)^{-1}(z_i\ptl_{z_i})^2(E_C)+2z_i\ptl_{z_i}({\cal
W}_{C_n})(x^\ast)^{-1}(z_i\ptl_{z_i})(E_C)]\\ &=&
\left[\sum_{1<i_1<i_2\leq
n}\frac{1}{2}\left(\frac{z_{i_1}z_{i_2}}{(z_{i_1}-z_{i_2})^2}+\frac{z_{i_1}z_{i_2}}
{(z_{i_1}z_{i_2}-1)^2}\right)+\frac{3}{4}\sum_{p=1}^n\frac{z_p^2}{(z_p^2-1)^2}\right]\Psi_C\\
& &+\sum_{i=1}^n{\cal
W}_{C_n}(x^\ast)^{-1}(z_i\ptl_{z_i})^2(\Psi_C)-\sum_{i=1}^n(z_i\ptl_{z_i})^2({\cal
W}_{C_n})(x^\ast)^{-1}E_C\\ &=&\left[\sum_{1<i_1<i_2\leq
n}\frac{1}{2}\left(\frac{z_{i_1}z_{i_2}}{(z_{i_1}-z_{i_2})^2}+\frac{z_{i_1}z_{i_2}}
{(z_{i_1}z_{i_2}-1)^2}\right)+\frac{3}{4}\sum_{p=1}^n\frac{z_p^2}{(z_p^2-1)^2}\right]\Psi_C\\
&
&+\sum_{i=1}^n(z_i\ptl_{z_i})^2(\Psi_C)-\frac{n(n+1)(2n+1)}{6}\Psi_C\hspace{5.5cm}(15.2.61)
\end{eqnarray*}
by (15.2.50), (15.2.52), (15.2.53) and (15.2.56)-(15.2.60).
Therefore, we have:
\begin{eqnarray*}\qquad& &\left[\sum_{1<i_1<i_2\leq n}
\frac{1}{2}\left(\frac{z_{i_1}z_{i_2}}{(z_{i_1}-z_{i_2})^2}+\frac{z_{i_1}z_{i_2}}
{(z_{i_1}z_{i_2}-1)^2}\right)+\frac{3}{4}\sum_{p=1}^n\frac{z_p^2}{(z_p^2-1)^2}\right]\Psi_C\\
 &=&\left(b+\frac{n(n+1)(2n+1)}{6}\right)\Psi_C-\sum_{i=1}^n(z_i\ptl_{z_i})^2(\Psi_C).\hspace{4.1cm}(15.2.62)\end{eqnarray*}
So $\Psi_C$ satisfies the Olshanesky-Perelomov equation (15.1.5)
with $K_3=0$ and
$$K_1=K_2=\frac{1}{2},\;\;K_4=\frac{3}{4},\;\;\nu=b+\frac{n(n+1)(2n+1)}{6}.\qquad\Box\eqno(15.2.63)$$

\section{Path Hypergeometric Functions of Type A}

In this section, we present our path hypergeometric functions of
type A found in [X14] via the Etingof trace functions of type A.

First we calculate the Etingof trace functions for $gl(2,\mbb R)$ in
terms of classical Gauss hypergeometric functions. Let $\lmd$ be a
weight of $gl(2,\mbb R)$ such that
$$\sgm=\lmd_1-\lmd_2\not\in\mbb{N}\eqno(15.3.1)$$
(cf. (15.2.3)). Choose $E_{1,2}$ as  the positive root vector.  The
Verma $gl(2,\mbb R)$-module with the highest-weight vector $v_\lmd$
of weight $\lmd$ is given by
$$M(\lmd) =\mbox{Span}\{E_{2,1}^iv_\lmd\mid i\in\mbb{N}\}\eqno(15.3.2)$$
with the action
$$E_{1,,2}(E_{2,1}^iv_\lmd)=i(\sgm+1-i)E_{2,1}^{i-1}v_\lmd,\qquad E_{2,1}(E_{2,1}^iv_\lmd)
=E_{2,1}^{i+1}v_\lmd,\eqno(15.3.3)$$
$$E_{1,1}(E_{2,1}^iv_\lmd)=(\lmd_1-i)E_{2,1}^iv_\lmd,\qquad E_{2,2}(E_{2,1}^iv_\lmd)
=(\lmd_2+i)E_{2,1}^iv_\lmd.\eqno(15.3.4)$$ Under the assumption
(15.4.1), $M(\lmd)$ is an irreducible $gl(2,\mbb R)$-module.

 Let $\kappa\in\mbb{C}$ be a fixed constant and recall
$$\td{M}=\{\sum_{i_1,i_2\in\mbb{Z}}v_{i_1,i_2}x_1^{i_1+\kappa}x_2^{i_2+\kappa}\mid v_{i_1,i_2}\in
M(\lmd)\},\eqno(15.3.5)$$ a space of Laurent series with
coefficients in $M(\lmd)$. Moreover, the representation $\pi$ of
$gl(2,\mbb R)$ on $\td{M}$ is given by
\begin{eqnarray*}E_{j_1,j_2}(\sum_{i_1,i_2\in\mbb{Z}}v_{i_1,i_2}x_1^{i_1+\kappa}x_2^{i_2+\kappa})
&=&\sum_{i_1,i_2\in\mbb{Z}}[E_{j_1,j_2}(v_{i_1,i_2})x_1^{i_1+\kappa}x_2^{i_2+\kappa}
\\ & &+v_{i_1,i_2}(x_{j_1}\ptl_{x_{j_2}}-\kappa\dlt_{j_1,j_2})(x_1^{i_1+\kappa}x_2^{i_2+\kappa})]
\hspace{2.2cm}(15.3.6)\end{eqnarray*} for $j_1,j_2\in\{1,2\}$ by
(15.2.11)-(15.2.13). Recall that a singular vector in a Verma module
is a nonzero weight vector annihilated by positive root vectors.
 Any singular vector of weight $\lmd$ in $\td{M}$ is of the form:
$$u_\lmd=(x_1x_2)^\kappa\sum_{i=0}^{\infty}a_iE_{2,1}^iv_\lmd x_1^ix_2^{-i}\eqno(15.3.7)$$
Observe that
\begin{eqnarray*}0&=&E_{1,2}(u_\lmd)\\ &=&(x_1x_2)^\kappa\sum_{i=0}^{\infty}a_i[i(\sgm+1-i)
E_{2,1}^{i-1}v_\lmd x_1^ix_2^{-i}+(\kappa-i)E_{2,1}^iv_\lmd
x_1^{i+1}x_2^{-i-1}]
\\ &=&(x_1x_2)^\kappa\sum_{i=0}^{\infty}((\kappa-i)a_i+(i+1)(\sgm-i)a_{i+1})v_\lmd
x_1^{i+1}x_2^{-i-1}.\hspace{3.6cm}(15.3.8)\end{eqnarray*} Thus
$$(\kappa-i)a_i+(i+1)(\sgm-i)a_{i+1}=0\qquad\for\;\;i\in\mbb{N},\eqno(15.3.9)$$
or equivalently,
$$a_{i+1}=-\frac{(\kappa-i)a_i}{(i+1)(\sgm-i)}\qquad\for\;\;i\in\mbb{N},\eqno(15.3.10)$$
For any constant $c$ and nonnegative integer $i$, we denote
$$ \la c\ra_i=c(c-1)\cdots (c-(i-1)).\eqno(15.3.11)$$
We normalize $u_\lmd$ by taking
$$a_0=1.\eqno(15.3.12)$$
By induction,
$$a_i=\frac{(-1)^i \la\kappa\ra_i}{i! \la\sgm\ra_i}.\eqno(15.3.13)$$
Thus
$$u_\lmd =\sum_{i=0}^{\infty}\frac{(-1)^i \la\kappa\ra_i}{i! \la\sgm\ra_i}E_{2,1}^iv_\lmd
x_1^{i+\kappa}x_2^{-i+\kappa}.\eqno(15.3.14)$$

Let $\Phi: M_\lmd\rta \td{M}$ be the Lie algebra module homomorphism
such that $\Phi(v_\lmd)=u_\lmd$. Note that
$$E_{2,1}^j(E_{2,1}^iv_\lmd
x_1^{i+\kappa}x_2^{-i+\kappa})=\sum_{k=0}^j{j\choose
 k}\la\kappa+i\ra_kE_{2,1}^{i+j-k}v_\lmd x_1^{i-k+\kappa}x_2^{-i+k+\kappa}.
 \eqno(15.3.15)$$
Moreover, for any $\iota\in\mbb{N}$, we have
\begin{eqnarray*}\qquad\qquad\sum_{k=\iota}^{\infty}{k\choose
\iota}z^k&=& z^{\iota}\sum_{k=\iota}^{\infty}{k\choose
\iota}z^{k-\iota} =\frac{z^{\iota}}{\iota
!}\frac{d^\iota}{dz^\iota}(\sum_{k=0}^{\infty}z^k)
\\ &=&\frac{z^{\iota}}{\iota
!}\frac{d^\iota}{dz^\iota}\left(\frac{1}{1-z}\right)=\frac{z^{\iota}}
{(1-z)^{\iota+1}}.\hspace{4.5cm}(15.3.16)\end{eqnarray*} View $\Phi$
as a ``function taking values in $\mbox{End}\:M_\lmd$.'' By
(15.3.15) and (15.3.16), the Etingof trace function
\begin{eqnarray*}\qquad& &E_A(z_1,z_2)=\mbox{tr}_{_{M_\lmd}}\:\Phi
z_1^{E_{1,1}}z_2^{E_{2,2}}\\
&=&(x_1x_2)^\kappa
z_1^{\lmd_1}z_2^{\lmd_2}\sum_{j=0}^{\infty}\sum_{i=0}^\infty
{j\choose i}\la\kappa+i\ra_i\frac{(-1)^i \la\kappa\ra_i}{i!
\la\sgm\ra_i}\left(\frac{z_2}{z_1}\right)^j\hspace{8cm}
\end{eqnarray*}
\begin{eqnarray*} \qquad&=&(x_1x_2)^\kappa
z_1^{\lmd_1}z_2^{\lmd_2}\sum_{i=0}^\infty\la\kappa+i\ra_i\frac{(-1)^i
\la\kappa\ra_i}{i! \la\sgm\ra_i} \sum_{j=0}^{\infty}{j\choose
i}\left(\frac{z_2}{z_1} \right)^j\\&=&(x_1x_2)^\kappa
z_1^{\lmd_1}z_2^{\lmd_2}\sum_{i=0}^\infty\la\kappa+i\ra_i\frac{(-1)^i
\la\kappa\ra_i}{i!
\la\sgm\ra_i}\left(\frac{z_2}{z_1}\right)^i\left(\frac{1}{1-z_2/z_1}\right)^{i+1}\\
&=&(x_1x_2)^\kappa\frac{z_1^{\lmd_1+1}z_2^{\lmd_2}}{z_1-z_2}\sum_{i=0}^\infty\frac{\la\kappa+i\ra_i
\la\kappa\ra_i}{i!
\la\sgm\ra_i}\left(\frac{z_2}{z_2-z_1}\right)^i.\hspace{4.4cm}(15.3.17)\end{eqnarray*}

Denote
$$(c)_i=c(c+1)(c+2)\cdots (c+i-1)\qquad \for\;\;c\in\mbb{C},\;i\in\mbb{N}.\eqno(15.3.18)$$
The  classical Gauss hypergeometric function
$$_2F_1(a,b;c; z)=\sum_{m=0}^{\infty}\frac{(a)_m(b)_m}{m!(c)_m}z^m,\qquad\mbox{where}\;\;a,b,c
\in\mbb{C},\;-c\not\in\mbb{N}.\eqno(15.3.19)$$
 Note
$$ \la c+n\ra_n=(c+1)_n,\;\; \la c\ra_n=(-1)^n(-c)_n\qquad\for\;\;c\in\mbb{C},\;n\in\mbb{N}.
\eqno(15.3.20)$$ By  (15.3.1) and (15.3.17)-(15.4.20), we have: \psp

{\bf Theorem 15.3.1}. {\it The Etingof trace function for $gl(2,\mbb
R)$ is
$$E_A(z_1,z_2)=(x_1x_2)^\kappa\frac{z_1^{\lmd_1+1}z_2^{\lmd_2}}{z_1-z_2}\:_2F_1\left(\kappa+1,-\kappa;\lmd_2-\lmd_1:
\frac{z_2}{z_2-z_1}\right).\eqno(15.3.21)$$ The function in
(15.2.17)
$$\Psi_A(z_1,z_2)=z_1^{\lmd_1+1/2}z_2^{\lmd_2-1/2}\:_2F_1\left(\kappa+1,-\kappa;\lmd_2-\lmd_1:
\frac{z_2}{z_2-z_1}\right)\eqno(15.3.22)$$}\pse

Next we want to find the Etingof trace functions for the general
case of $gl(n,\mbb R)$ with $n>2$ under a certain condition and
their master hypergeometric functions.

 As the case of $sl(n,\mbb R)$, we choose
$$\{E_{i,j}\mid 1\leq i<j\leq n\}\;\;\mbox{as positive root vectors}.\eqno(15.3.23)$$
In particular, we have
$$\{E_{i,i+1}\mid i=1,2,...,n-1\}\;\;\mbox{as  positive simple root vectors}.\eqno(15.3.24)$$
Accordingly,
$$\{E_{i,j}\mid 1\leq j<i\leq n\}\;\;\mbox{are negative root vectors}\eqno(15.3.25)$$
and
$$\{E_{i+1,i}\mid i=1,2,...,n-1\}\;\;\mbox{are negative simple root vectors}.\eqno(15.3.26)$$
Let
$$\G_A=\sum_{1\leq j<i\leq n}\mbb{N}\ves_{i,j}\eqno(15.3.27)$$
be the torsion-free additive semigroup  of  rank $n(n-1)/2$ with
$\ves_{i,j}$ as base elements.
 For
$$\al =\sum_{1\leq j<i\leq n}\al_{ij}\ves_{i,j}\in \G_A,\eqno(15.3.28)$$
we denote
$$E^\al=E_{2,1}^{\al_{2,1}}E_{3,1}^{\al_{3,1}}E_{3,2}^{\al_{3,2}}E_{4,1}^{\al_{4,1}}\cdots
E_{n,1}^{\al_{n,1}}\cdots E_{n,n-1}^{\al_{n,n-1}}\in
U(gl(n)).\eqno(15.3.29)$$ Denote by  ${\msr G}_-$ the Lie subalgebra
spanned by (15.3.25). Then
$$\{E^\al\mid \al\in\G_A\}\;\;\mbox{forms a PBW basis of}\;\;U({\msr G}_-).\eqno(15.3.30)$$

Let $\lmd$ be a weight of $gl(n,\mbb R)$ such that
$$\lmd_1-\lmd_2=\cdots =\lmd_{n-2}-\lmd_{n-1}=\kappa\;\;\mbox{and}\;\;\lmd_{n-1}-\lmd_n
\not\in\mbb{N}, \eqno(15.3.31)$$ for some constant $\kappa$ (cf.
(15.2.3)). Set
$$\sgm_i=\lmd_i-\lmd_{i+1}\qquad\for\;\;i\in\ol{1,n-1}.\eqno(15.3.32)$$
 The Verma $gl(n,\mbb R)$-module with the highest-weight
vector $v_\lmd$ of weight $\lmd$ is given by
$$M(\lmd) =\mbox{Span}\{E^\al v_\lmd\mid\al\in \G_A\},\eqno(15.3.33)$$
with the action determined by
\begin{eqnarray*}E_{i,i+1}(E^\al v_\lmd)&=&(\sum_{p=1}^{i-1}\al_{i+1,p}E^{\al+\ves_{i,p}
-\ves_{i+1,p}}-\sum_{p=i+2}^n\al_{p,i}E^{\al+\ves_{p,i+1}-\ves_{p,i}}\\
& &+\al_{i+1,i}
(\sgm_i+1-\sum_{p=i+1}^n\al_{p,i}+\sum_{p=i+2}^n\al_{p,i+1})E^{\al-\ves_{i+1,i}})v_\lmd,
\hspace{1.7cm}(15.3.34)\end{eqnarray*}
$$E_{j,i}(E^\al v_\lmd)=(E^{\al+\ves_{j,i}}+\sum_{p=1}^{i-1}\al_{i,p}E^{\al+\ves_{j,p}
-\ves_{i,p}}) v_\lmd,\eqno(15.3.35)$$
$$E_{k,k}(E^\al v_\lmd)=(\lmd_k+\sum_{p=1}^{k-1}\al_{k,p}-\sum_{q=k+1}^n\al_{q,k})
E^\al v_\lmd \eqno(15.3.36)$$ for $1\leq i<j\leq n$ and
$k\in\ol{1,n}$.

Recall the notions in (15.2.7)-(15.2.9) and for $M=M(\lmd)$,
$$\td{M}=\{\sum_{\vec i\in\mbb{Z}^n}v_{\vec i}x^{\vec i}\bar x^\kappa\mid v_{\vec i}\in M(\lmd)\}.\eqno(15.3.37)$$
 Moreover, the representation $\pi$ of $gl(n,\mbb R)$ on $\td{M}$ is given by
$$E_{j_1,j_2}(\sum_{\vec i\in\mbb{Z}^n}(v_{\vec i})x^{\vec i}\bar
x^\kappa) =\sum_{\vec i\in\mbb{Z}^n}[E_{j_1,j_2}(v_{\vec i})x^{\vec
i}\bar x^\kappa+v_{\vec i}
(\ptl_{x_{j_2}}-\kappa\dlt_{j_1,j_2})(x^{\vec i}\bar x^\kappa)]
\eqno(15.3.38)$$ for $j_1,j_2\in\ol{1,n}$.

For
$$\vec i=(i_1,i_2,...,i_{n-1})\in\mbb{N}^{n-1},\eqno(15.3.39)$$
we denote
$$ x^{(\vec i)}=x_1^{i_1}x_2^{i_2-i_1}\cdots x_{n-1}^{i_{n-1}-i_{n-2}}x_n^{-i_{n-1}},
\qquad \es(\vec i)=\sum_{p=1}^{n-1}i_p\ves_{p+1,p}.\eqno(15.3.40)$$
Then we have: \psp

{\bf Lemma 15.3.2}. {\it The vector
$$u_\lmd=\bar x^\kappa\sum_{i_1,...,i_{n-1}=0}^{\infty}\frac{(-1)^{i_1+i_2+\cdots
+i_{n-1}} \la\kappa\ra_{i_{n-1}}}{i_1!i_2!\cdots i_{n-1}!
\la\sgm\ra_{i_{n-1}}}E^{\es(\vec i)} v_\lmd x^{(\vec
i)}\eqno(15.3.41)$$ is a singular vector of weight $\lmd$, where
$\sgm=\sgm_{n-1}=\lmd_{n-1}-\lmd_n$.}

{\it Proof}. For $p\in\{1,2,...,n-2\}$, we have
\begin{eqnarray*} & &E_{p,p+1}(u_\lmd)\\ &=&
\bar
x^\kappa\sum_{i_1,...,i_{n-1}=0}^{\infty}\frac{(-1)^{i_1+i_2+\cdots+i_{n-1}}
 \la\kappa\ra_{i_{n-1}}}{i_1!i_2!\cdots i_{n-1}!\la
\sgm\ra_{i_{n-1}}}[(\kappa+i_{p+1}-i_p) E^{\es(\vec i)}v_\lmd
x_px_{p+1}^{-1}x^{(\vec i)}\\ & &+(\kappa+1-i_p+i_{p+1})i_p
E^{\es(\vec i)-\ves_p}v_\lmd x^{(\vec i)}]\\ &=& \bar
x^\kappa[\sum_{i_1,...,i_{n-1}=0}^{\infty}\frac{(-1)^{i_1+i_2+\cdots+i_{n-1}}
 \la\kappa\ra_{i_{n-1}}}{i_1!i_2!\cdots i_{n-1}!
\la\sgm\ra_{i_{n-1}}}(\kappa+i_{p+1}-i_p) E^{\es(\vec i)}v_\lmd
x_px_{p+1}^{-1}x^{(\vec i)}\\ &&+
\sum_{i_1,...,i_{n-1}=0}^{\infty}\frac{(-1)^{i_1+i_2+\cdots+i_{n-1}}
\la\kappa\ra_{i_{n-1}}} {i_1!i_2!\cdots i_{n-1}!
\la\sgm\ra_{i_{n-1}}}(\kappa+1-i_p+i_{p+1})i_pE^{\es(\vec i)-\ves_p}
v_\lmd x^{(\vec i)}]\\ &=&\bar
x^\kappa[\sum_{i_1,...,i_{n-1}=0}^{\infty}
\frac{(-1)^{i_1+i_2+\cdots+i_{n-1}}
\la\kappa\ra_{i_{n-1}}(\kappa+i_{p+1}-i_p)}{i_1!i_2!\cdots i_{n-1}!
\la\sgm\ra_{i_{n-1}}}E^{\es(\vec i)}v_\lmd x_px_{p+1}^{-1}x^{(\vec
i)}
\\ &&-\sum_{i_1,...,i_{n-1}=0}^{\infty}\frac{(-1)^{i_1+i_2+\cdots+i_{n-1}}
 \la\kappa\ra_{i_{n-1}}(\kappa+i_{p+1}-i_p)} {i_1!i_2!\cdots
i_{n-1}!\la \sgm\ra_{i_{n-1}}}E^{\es(\vec i)}v_\lmd x_p
x_{p+1}^{-1}x^{(\vec i)}]=0.\hspace{0.5cm}(15.3.42)\end{eqnarray*}

Moreover,
\begin{eqnarray*}& & E_{n-1,n}(u_\lmd)\\ &=&
\bar
x^\kappa\sum_{i_1,...,i_{n-1}=0}^{\infty}\frac{(-1)^{i_1+i_2+\cdots+i_{n-1}}
\la\kappa\ra_{i_{n-1}}}{i_1!i_2!\cdots i_{n-1}!
\la\sgm\ra_{i_{n-1}}}[(\kappa-i_{n-1})E^{\es(\vec i)}v_\lmd
x_{n-1}x_n^{-1}x^{(\vec i)}\\ & &+(\sgm+1-i_{n-1})i_{n-1}E^{\es(\vec
i)-\ves_{n-1}}v_\lmd x^{(\vec i)}]\\
&=&\bar
x^\kappa[\sum_{i_1,...,i_{n-1}=0}^{\infty}\frac{(-1)^{i_1+i_2+\cdots+i_{n-1}}
\la\kappa\ra_{i_{n-1}}}{i_1!i_2!\cdots i_{n-1}!
\la\sgm\ra_{i_{n-1}}}(\kappa-i_{n-1})E^{\es(\vec i)}v_\lmd
x_{n-1}x_n^{-1}x^{(\vec i)}\\ &&+
\sum_{i_1,...,i_{n-1}=0}^{\infty}\frac{(-1)^{i_1+i_2+\cdots+i_{n-1}}
\la\kappa\ra_{i_{n-1}}}{i_1!i_2!\cdots i_{n-1}!
\la\sgm\ra_{i_{n-1}}}(\sgm+1-i_{n-1})i_{n-1}E^{\es(\vec
i)-\ves_{n-1}}v_\lmd x^{(\vec i)}]\\
&=&\bar
x^\kappa[\sum_{i_1,...,i_{n-1}=0}^{\infty}\frac{(-1)^{i_1+i_2+\cdots+i_{n-1}}
\la\kappa\ra_{i_{n-1}}}{i_1!i_2!\cdots i_{n-1}!
\la\sgm\ra_{i_{n-1}}}(\kappa-i_{n-1})E^{\es(\vec i)}v_\lmd
x_{n-1}x_n^{-1}x^{(\vec i)}\\
&&-\sum_{i_1,...,i_{n-1}=0}^{\infty}\frac{(-1)^{i_1+i_2+\cdots+i_{n-1}}
\la\kappa\ra_{i_{n-1}+1}}{i_1!i_2!\cdots i_{n-1}!
\la\sgm\ra_{i_{n-1}}}E^{\es(\vec i)}v_\lmd x_px_{p+1}^{-1}x^{(\vec
i)}]=0.\hspace{2.6cm}(15.3.43)\end{eqnarray*} Hence $u_\lmd$ is a
singular vector. $\qquad\Box$ \psp

 Let $\Phi: M(\lmd)\rta\td{M}$ be the Lie algebra module homomorphism
such that $\Phi(v_\lmd)=u_\lmd$. Let $\vec i$ as in (15.2.8) and fix
$\be\in\G_A$. By (15.3.35),  we have
\begin{eqnarray*}& &E_{n,n-1}^{\be_{n,n-1}}(E^{\es(\vec i)}v_\lmd)=
\sum_{p_{n,n-1}=0}^{\be_{n,n-1}} {\be_{n,n-1}\choose p_{n,n-1}}\la
i_{n-2}\ra_{p_{n,n-1}}
\\ & &E^{\es(\vec i)-p_{n,n-1}\ves_{n-1,n-2}+p_{n,n-1}\ves_{n,n-2}+(\be_{n,n-1}-p_{n,n-1})
\ves_{n,n-1}}v_\lmd,\hspace{4.5cm}(15.3.44)\end{eqnarray*}
\begin{eqnarray*}& &E_{n,n-2}^{\be_{n,n-2}}E_{n,n-1}^{\be_{n,n-1}}(E^{\es(\vec i)}v_\lmd)
=\sum_{p_{n,n-2},p_{n,n-1}=0}^{\infty} {\be_{n,n-2}\choose
p_{n,n-2}}{\be_{n,n-1}\choose p_{n,n-1}}\la i_{n-3}\ra_{p_{n,n-2}}
\la i_{n-2} \ra_{p_{n,n-1}}\\ & & E^{\es(\vec
i)+\sum_{r=n-2,n-1}[-p_{n,r}\ves_{r,r-1}+p_{n,r}\ves_{n,r-1}+(\be_{n,r}-p_{n,r}
\ves_{n,r}]}v_\lmd,\hspace{4.7cm}(15.3.45)\end{eqnarray*}
\begin{eqnarray*}\qquad& &\left[\prod_{r=1}^{n-1}E_{n,r}^{\be_{n,r}}\right]
(E^{\es(\vec
i)}v_\lmd)=\sum_{p_{n,2},...,,p_{n,n-1}=0}^{\infty}\left[\prod_{r=2}^{n-1}
{\be_{n,r}\choose p_{n,r}} \la i_{r-1}\ra_{p_{n,r}}\right]\\
& & E^{\es(\vec
i)+\be_{n,1}\ves_{n,1}+\sum_{s=2}^{n-1}(-p_{n,s}\ves_{s,s-1}+p_{n,s}\ves_{n,s-1}
+(\be_{n,s}-p_{n,s})\ves_{n,s})}v_\lmd,\hspace{3.4cm}(15.3.46)\end{eqnarray*}
\begin{eqnarray*}&&\left[\prod_{r=1}^{n-2}E_{n-1,r}^{\be_{n-1,r}}\right]\left[\prod_{r=1}^{n-1}
E_{n,r}^{\be_{n,r}}\right](E^{\es(\vec
i)}v_\lmd)=\sum_{i=n-1,n}\;\sum_{j=1}^{i-1}
\sum_{p_{i,j}=0}^{\infty}\\ &
&\left[\prod_{r=1}^{n-2}{\be_{n-1,r}\choose p_{n-1,r}}\right]
\left[\prod_{r=1}^{n-1}{\be_{n,r}\choose p_{n,r}} \right]
\left[\prod_{r=1}^{n-3} \la i_r\ra_{p_{n,r+1}+p_{n-1,r+1}}\right]
\la i_{n-2}\ra_{p_{n,n-1}}
\\ & &E^{\es(\vec i)+\sum_{r=n-1,n}\be_{r,1}\ves_{r,1}+\sum_{r=n-1,n}\sum_{s=2}^{r-1}(-p_{r,s}
\ves_{s,s-1}+p_{r,s}\ves_{r,s-1}+(\be_{r,s}-p_{r,s})\ves_{r,s})}v_\lmd,\hspace{1.8cm}(15.3.47)
\end{eqnarray*}
\begin{eqnarray*}& &E^\be(E^{\es(\vec i)}v_\lmd)=\sum_{2\leq k<j\leq n}\;
\sum_{p_{j,k}=0}^{\infty} \left[\prod_{2\leq m<l\leq
n}{\be_{l,m}\choose p_{l,m}} \right]\left[\prod_{s=1}^{n-2} \la
i_s\ra_{_{\sum_{r=s+2}^np_{r,s+1}}}\right]\\ & & E^{\es(\vec
i)+\sum_{1\leq m<l\leq
n}(\be_{l,m}+(1-\dlt_{l,m+1})p_{l,m+1}-p_{l,m}-\dlt_{l,m+1}
\sum_{s=l+1}^np_{s,l})\ves_{l,m}}v_\lmd.\hspace{3.2cm}(15.3.48)\end{eqnarray*}
where we treat $p_{2,1}=p_{3,1}=\cdots=p_{n,1}=0$. Hence by
(15.3.37) and (15.3.47), we obtain
\begin{eqnarray*}& &E^\al(E^{\es(\vec i)}v_\lmd  x^{(\vec i)}\bar x^\kappa)\\
&=&\sum_{p_{2,1}^0,...,p_{n,1}^0=0}^{\infty} \;\sum_{2\leq k<j\leq
n}\;\sum_{p_{j,k}^0,p_{j,k}=0}^{\infty}\left[\prod_{m=2}^n{\al_{m,1}\choose
p_{m,1}^0}\right]\left[\prod_{2\leq m<l\leq n}{\al_{l,m}\choose
p_{l,m}^0,p_{l,m}}\right]\\
&&\left[\prod_{s=1}^{n-2}\la
i_s\ra_{_{\sum_{r=s+2}^np_{r,s+1}}}\right]\left[\prod_{s=1}^{n-1}
\la\kappa+i_s-i_{s-1}\ra_{_{\sum_{r=s+1}^np_{r,s}^0}}\right]\\ & &
E^{\es(\vec i)+\sum_{1\leq m<l\leq
n}(\al_{l,m}+p_{l,m+1}-p_{l,m}^0-p_{l,m}
-\dlt_{l,m+1}\sum_{s=l+1}^np_{s,l})\ves_{l,m}}v_\lmd\\ & & \bar
x^\kappa\prod_{r=1}^nx_r^{i_r-i_{r-1}+\sum_{s=1}^{r-1}p^0_{r,s}
-\sum_{q=r+1}^np_{q,r}^0}\hspace{7.5cm}(15.3.49)\end{eqnarray*} for
$\al\in\G_A$,  where we treat
$$i_{-1}=i_n=p_{2,1}=p_{3,1}=\cdots =p_{n,1}=p_{2,2}=p_{3,3}=\cdots p_{n,n}=0.\eqno(15.3.50)$$

View $\Phi$  as a function in $\{x_1,x_2,...,x_n\}$ taking values in
$\mbox{End}\:M(\lmd)$, the space of linear transformations on
$M(\lmd)$.
 In order to
calculate the Etingof trace function $E_A(z_1,z_2,...,z_n)$ in
(15.2.15), we need to take
$$p_{l,m+1}=p_{l,m}^0+p_{l,m},\;\; i_k=p_{k+1,k}^0+p_{k+1,k}+\sum_{s=k+2}^np_{s,k+1}\eqno(15.3.51)$$
in (15.3.49) for $1\leq m<l-1\leq n-1$ and $k\in\ol{1,n-1}$. By
(15.3.50), (15.3.51) and induction,
$$p_{l,m}=\sum_{r=1}^{m-1}p_{l,r}^0\qquad\for\;\;1\leq m<l\leq n.\eqno(15.3.52)$$
Moreover,
$$i_k=\sum_{r=k+1}^n(p_{r,k}^0+p_{r,k})=\sum_{r=k+1}^n\sum_{s=1}^kp^0_{r,s}\qquad\for\;\;
k\in\ol{1,n-1}.\eqno(15.3.53)$$
$$i_k-i_{k-1}=\sum_{r=k+1}^n\sum_{s=1}^kp^0_{r,s}-\sum_{r=k}^n\sum_{s=1}^{k-1}p^0_{r,s}=
\sum_{r=k+1}^np^0_{r,k}-\sum_{s=1}^{k-1}p^0_{k,s}\eqno(15.3.54)$$
for $k\in\ol{1,k}$. Furthermore,
$${\al_{j,k}\choose p_{j,k}^0,p_{j,k}}={p^0_{j,k}+p_{j,k}\choose
p_{j,k}}{ \al_{j,k}\choose p^0_{j,k}+p_{j,k}} \eqno(15.3.55)$$ and
$$ {p^0_{j,k}+p_{j,k}\choose p_{j,k}}=
\frac{(p^0_{j,k}+p_{j,k})!}{p^0_{j,k}!p_{j,k}!}=\frac{(p^0_{j,k}+p_{j,k})!}{p^0_{j,k}!
(p^0_{j,k-1}+p_{j,k-1})!}\eqno(15.3.56)$$ by (15.3.51).

For convenience, we denote
$$\al_{_{1\ast}}=\al_n^\ast=0,\;\;\al_{_{i\ast}}=\sum_{r=1}^{i-1}\al_{i,r},\;\;
\al_{_{j}}^\ast=\sum_{s=j+1}^n\al_{s,j}\eqno(15.3.57)$$ for
$\al\in\G_A,\;i\in \ol{2,n}$ and $j\in\ol{1,n-1}$. By (15.3.40),
(15.3.41) and (15.3.49)-(15.3.57), the Etingof's trace function
\begin{eqnarray*}& &E_A(z_1,z_2,...,z_n)
\\ &=&\sum_{\al\in\G_A}\;\sum_{p_{r_1,r_2}^0\in\mbb{N},
\;1\leq r_2<r_1\leq n}(-1)^{\sum_{s=1}^{n-1}i_s}\left[\prod_{2\leq
k<j\leq n}{p^0_{j,k}+p_{j,k}\choose p_{j,k}}\right]\\
&&\times
\frac{\prod_{s=1}^{n-1}\la\kappa+i_s-i_{s-1}\ra_{_{\sum_{r=s+1}^np^0_{r,s}}}}
{\la\sgm\ra_{_{i_{n-1}}}\prod_{s=1}^{n-1}(p_{s+1,s}^0+p_{s+1,s})!}
\left[\prod_{1\leq k<j\leq n}{\al_{j,k}\choose
p^0_{j,k}+p_{j,k}}\right] \la\kappa\ra_{i_{n-1}}\\& &\times \bar
x^\kappa
\left(\prod_{r=1}^nz_r^{\lmd_r}\right)\prod_{1\leq k<j\leq n}\left(\frac{z_j}{z_k}\right)^{\al_{j,k}}\\
&=&\sum_{p_{r_1,r_2}^0\in\mbb{N},\;1\leq r_2<r_1\leq
n}(-1)^{\sum_{s=1}^{n-1}i_s}\frac{\prod_{s=1}^{n-1}\la
\kappa+i_s-i_{s-1}\ra_{_{\sum_{r=s+1}^n
p^0_{r,s}}}\la\kappa\ra_{i_{n-1}}}{\la
\sgm\ra_{_{i_{n-1}}}\prod_{1\leq s_1<s_2\leq n}p_{s_2,s_1}^0!}\\ &
&\times \left(\prod_{r=1}^nz_r^{\lmd_r}\right)\left[\prod_{1\leq
r_1<r_2\leq n}\left(\frac{z_{r_2}}{z_{r_1}}\right)^{p_{r_2,r_1}^0
+p_{r_2,r_1}}\right]\\ & &\times\left[\prod_{1\leq s_1<s_2\leq
n}\left(\frac{1}{1-z_{s_2}/z_{s_1}}
\right)^{p_{s_2,s_1}^0+p_{s_2,s_1}+1}\right]\bar x^\kappa \\
&=&\frac{\bar x^\kappa\prod_{r=1}^nz_r^{\lmd_r+n-r}}{\prod_{1\leq
k<j\leq
n}(z_k-z_j)}\sum_{\be\in\G_A}\frac{\left[\prod_{s=1}^{n-1}\la\kappa+\be_s^\ast-\be_{s\ast}
\ra_{_{\be_s^\ast}}\right]\la\kappa\ra_{_{\be_{n\ast}}}}{\la\sgm\ra_{_{\be_{n\ast}}}
\prod_{1\leq s_1<s_2\leq n}\be_{s_2,s_1}!}\\ & &\times
\left[\prod_{1\leq k<j\leq
n}\left(\frac{z_j}{z_j-z_k}\right)^{\sum_{r=1}^k\be_{j,r}}\right].\hspace{7.6cm}(15.3.58)
\end{eqnarray*}
Set
$$\xi_{r_2,r_1}^A=\prod_{s=r_1}^{r_2-1}\frac{z_{r_2}}{z_{r_2}-z_s}
\qquad \for\;\;1\leq r_1<r_2\leq n.\eqno(15.3.59)$$

Motivated by (15.3.17)-(15.3.19) and (15.3.58), we define our
$(n(n-1)/2)$-variable {\it path hypergeometric function of type
A}\index{path hypergeometric function! of type A} by
$${\msr X}_A(\tau_1,..,\tau_n;\vt)\{z_{j,k}\}=\sum_{\be\in\G_A}\frac{\left[\prod_{s=1}^{n-1}
(\tau_s-\be_{s\ast})_{_{\be_s^\ast}}\right](\tau_n)_{_{\be_{n\ast}}}}
{\be!(\vt)_{_{\be_{n\ast}}}}z^\be\eqno(15.3.60)$$ (cf. (15.3.57)),
where
$$\be!=\prod_{1\leq k<j\leq n}\be_{j,k}!,\qquad z^\be=\prod_{1\leq k<j\leq n}z_{j,k}^{\be_{j,k}}.
\eqno(15.3.61)$$
 According to (15.3.20) and (15.3.58)-(15.3.60), we have:
\psp

{\bf Theorem 15.3.3}. {\it Under the condition (15.3.31), the
Etingof trace function   is
$$E_A(z_1,z_2,...,z_n)=\left[\frac{\bar x^\kappa\prod_{r=1}^nz_r^{\lmd_r+n-r}}{\prod_{1\leq k<j\leq n}
(z_k-z_j)}\right]{\msr
X}_A(\kappa+1,..,\kappa+1,-\kappa;-\sgm)\{\xi_{r_2,r_1}^A\}.\eqno(15.3.62)$$
Moreover, (15.1.9) implies that the function in (15.2.17):
$$\Psi_A(z_1,z_2,\cdots,z_n)=\prod_{r=1}^nz_r^{\lmd_r+(n+1)/2-r}
{\msr
X}_A(\kappa+1,..,\kappa+1,-\kappa;-\sgm)\{\xi_{r_2,r_1}^A\}.\eqno(15.3.63)$$}

\section{Path Hypergeometric Functions of Type
C}

In this section, we want to introduce and calculate  the Etingof
trace function of $sp(2n,\mbb R)$.

Take the settings in (15.2.31)-(15.2.36).  Set
$$h_n=C_{n,n},\;\;h_i=C_{i,i}-C_{i+1,i+1}\eqno(15.4.1)$$
for $i=1,2,...,n-1.$ We choose positive root vectors
$$\{C_{i,j},\;C_{i,n+j},\;C_{k,n+k} \mid 1\leq i<j\leq n,\; k\ol{1,n}\}.
\eqno(15.4.2)$$ In particular, we take
$$\{C_{i,i+1},\;C_{n,2n}\mid i=1,2,...,n-1\}\;\;\mbox{as positive simple root
vectors}.\eqno(15.4.3)$$ Accordingly,
$$\{C_{i,j},\;C_{i,n+j},\;C_{n+k,k} \mid 1\leq j<i\leq n,\;k\in\ol{1,n}\}
\eqno(15.4.4)$$ are negative root vectors and
$$\{C_{i+1,i},\;C_{2n,n}\mid i\in\ol{1,n-1}\}\;\;\mbox{are negative simple root
vectors}.\eqno(15.4.5)$$

Let
$$\G_C=\sum_{1\leq j<i\leq n}\mbb{N}\ves_{i,j}+\sum_{1\leq j\leq i\leq n}\mbb{N}\ves_{n+i,j}
\eqno(15.4.6)$$ be the torsion-free additive semigroup  of  rank
$n^2$ with $\ves_{p,q}$ as base elements, and let $U(sp(2n,\mbb R))$
be the universal enveloping algebra of $sp(2n,\mbb R)$. For any
$$\al=\sum_{1\leq j<i\leq n}\al_{i,j}\ves_{i,j}+\sum_{1\leq j\leq i\leq n}\al_{n+i,j}\ves_{n+i,j}
\in \G_C,\eqno(15.4.7)$$ we denote
\begin{eqnarray*}\qquad C^\al&=&C_{2,1}^{\al_{2,1}}C_{3,1}^{\al_{3,1}}C_{3,2}^{\al_{3,2}}
C_{4,1}^{\al_{4,1}}\cdots C_{n,1}^{\al_{n,1}}\cdots
C_{n,n-1}^{\al_{n,n-1}}\\ & &\times
 C_{n+1,1}^{\al_{n+1,1}}C_{n+2,1}^{\al_{n+2,1}}C_{n+2,2}^{\al_{n+2,2}}C_{n+3,1}^{\al_{n+3,1}}
 \cdots C_{2n,1}^{\al_{2n,1}}\cdots C_{2n,n}^{\al_{2n,n}}\hspace{3.5cm}(15.4.8)\end{eqnarray*}
Denote by  ${\msr G}_-$  the Lie subalgebra spanned by (15.4.4).
Then
$$\{C^\al\mid \al\in\G_C\}\;\;\mbox{forms a PBW basis of}\;\;U({\msr G}_-)).\eqno(15.4.9)$$

For convenience, we treat
$$\ves_{n+i,j}=\ves_{n+j,i}\qquad\for\;\;i,j=1,2,...,n.\eqno(15.4.10)$$
Note that
$$[C_{n,2n},C_{2n,n}]=h_n,\;\; [C_{i,i+1},C_{i+1,i}]=h_i,\eqno(15.4.11)$$
$$[h_i,C_{n+i+1,i}]=0,\;\; [h_i,C_{n+i,i}]=-2C_{n+i,i},\;\; [h_i,C_{n+i+1,i+1}]=2C_{n+i+1,i+1},\eqno(15.4.12)$$
$$[h_i,C_{n+k,i}]=-C_{n+k,i},\;\;[h_i,C_{n+k,i+1}]=C_{n+k,i+1},\eqno(15.4.13)$$
$$[C_{i,i+1},C_{n+p,i}]=-C_{n+p,i+1},\;\; [C_{i,i+1},C_{n+i+1,i}]=-2C_{n+i+1,i+1},\eqno(15.4.14)$$
$$[C_{n,2n},C_{2n,i}]=C_{n,i},\;\; [C_{n+j,r}, C_{r,q}]=(1+\dlt_{j,q})C_{n+j,q},
\;\;[C_{n+r,r},C_{r,q}]=C_{n+r,q}\eqno(15.4.15)$$ for
$i\in\ol{1,n-1}$ and $ j,k,p,q,r\in\ol{1,n}$ such that  $k\neq
i,i+1,\;p\neq i+1$ and $q< r$.

 Let $\lmd$ be a weight, which is a linear function on $H_C$, such that
$$\lmd(h_i)=-\frac{1}{2}\qquad\for\;\;i=1,2,...,n-1;\qquad \lmd(h_n)=\lmd_n\in\mbb{C}. \eqno(15.4.16)$$
Recall that $sp(2n,\mbb R)$ is generated by
$$\{C_{i,i+1},C_{i+1,i},C_{n,2n},C_{2n,n},\mid i=1,2,...,n-1\}\eqno(15.4.17)$$
as a Lie algebra. The Verma  $sp(2n,\mbb R)$-module with the
highest-weight vector $v_\lmd$ of weight $\lmd$ is given by
$$M(\lmd) =\mbox{Span}\{E^\al v_\lmd\mid\al\in \G_C\},\eqno(15.4.18)$$
with the action determined by
\begin{eqnarray*}& &C_{i,i+1}(C^\al v_\lmd)\\ &=&[\sum_{j=1}^{i-1}\al_{i+1,j}C^{\al+\ves_{i,j}
-\ves_{i+1,j}}-\sum_{j=i+2}^n\al_{j,i}C^{\al+\ves_{j,i+1}-\ves_{j,i}}-\sum_{k\neq
i+1}\al_{n+k,i} C^{\al-\ves_{n+k,i}+\ves_{n+k,i+1}}\\ &
&-2\al_{n+i+1,i}C^{\al-\ves_{n+i+1,i}+\ves_{n+i+1,i+1}}+
\al_{i+1,i}(1/2-\sum_{j=i+1}^n\al_{j,i}+\sum_{j=i+2}^n\al_{j,i+1}\\
& &+\sum_{k\neq i,i+1}
(\al_{n+k,i+1}-\al_{n+k,i})-2\al_{n+i,i}+2\al_{n+i+1,i+1})C^{\al-\ves_{i+1,i}}]v_\lmd
\hspace{2.4cm}(15.4.19)\end{eqnarray*} for $i\in \ol{1,n-1}$,
\begin{eqnarray*}& &C_{n,2n}(C^\al v_\lmd)=[\sum_{i=1}^{n-1}\al_{2n,i}(C^{\al-\ves_{2n,i}
+\ves_{n,i}}+(\al_{2n,i}-1)C^{\al-2\ves_{2n,i}+\ves_{n+i,i}}\\ &
&+\sum_{j=1}^{i-1}\al_{2n,j}
C^{\al-\ves_{2n,i}-\ves_{2n,j}+\ves_{n+i,j}})+\al_{2n,n}(\lmd_n+1-\al_{2n,n})C^{\al-\ves_{2n,n}}]
v_\lmd\hspace{2cm}(15.4.20)\end{eqnarray*} and
$$C_{j,i}(E^\al v_\lmd)=(C^{\al+\ves_{j,i}}+\sum_{p=1}^{i-1}\al_{i,p}C^{\al+\ves_{j,p}
-\ves_{i,p}})v_\lmd,\eqno(15.4.21)$$
\begin{eqnarray*}C_{n+j,i}(E^\al v_\lmd)&=&(C^{\al+\ves_{n+j,i}}+\sum_{q=1}^{j-1}\al_{n+j,q}
C^{\al+\ves_{n+q,i}-\ves_{j,q}}\\ & &+
\sum_{p=1}^{i-1}\al_{i,p}(C^{\al+\ves_{n+j,p}-\ves_{i,p}}+\sum_{q=1}^{j-1}\al_{n+j,q}C^{\al+
\ves_{n+q,p}-\ves_{i,p}-\ves_{j,q}})\hspace{1.6cm}(15.4.22)\end{eqnarray*}
for $1\leq i<j\leq n$,
\begin{eqnarray*}\hspace{1cm}C_{n+k,k}(C^\al v_\lmd)&=&[C^{\al+\ves_{n+k,k}}+\sum_{r=1}^{k-1}
\al_{k,r}(C^{\al-\ves_{k,r}+\ves_{n+k,r}}+(\al_{k,r}-1)C^{\al-2\ves_{k,r}+\ves_{n+r,r}}\\
&&
+\sum_{s=r+1}^{k-1}\al_{k,s}C^{\al-\ves_{k,r}-\ves_{k,s}+\ves_{n+s,r}})]v_\lmd\hspace{4.3cm}
(15.4.23)\end{eqnarray*}
 by (15.4.11)-(15.4.15).

Recall the notions in (15.2.8) and (15.2.39). Moreover, we take
$M=M(\lmd)$ and
$$\td{M}=\{\sum_{\vec i\in\mbb{Z}^n}w_{\vec i}x^{\vec i}x^\ast\mid w_{\vec i}\in
M(\lmd)\}. \eqno(15.4.24)$$ The representation $\pi$ in
(15.2.41)-(15.2.45) of $sp(2n,\mbb R)$ on $\td{M}$ is given by
$$C_{p,q}(\sum_{\vec i\in\mbb{Z}^n}w_{\vec i}x^{\vec i}x^\ast)=\sum_{\vec i\in\mbb{Z}^n}
\left(C_{p,q}(w_{\vec i})x^{\vec i}x^\ast+w_{\vec
i}x_p\ptl_{x_q}(x^{\vec i}x^\ast) +\frac{\dlt_{p,q}}{2}w_{\vec
i}x^{\vec i}x^\ast\right),\eqno(15.4.25)$$
$$C_{p,n+q}(\sum_{\vec i\in\mbb{Z}^n}w_{\vec i}x^{\vec i}x^\ast)=\sum_{\vec i\in\mbb{Z}^n}
\left(C_{p,n+q}(w_{\vec i})x^{\vec
i}x^\ast-\frac{1}{1+\dlt_{p,q}}w_{\vec i}x_px_qx^{\vec i}
x^\ast\right),\eqno(15.4.26)$$
$$C_{n+p,q}(\sum_{\vec i\in\mbb{Z}^n}w_{\vec i}x^{\vec i}x^\ast)=\sum_{\vec i\in\mbb{Z}^n}
(C_{n+p,q}(w_{\vec i})x^{\vec i}x^\ast+\frac{1}{1+\dlt_{p,q}}w_{\vec
i}\ptl_{x_p}\ptl_{x_q} (x^{\vec i}x^\ast))\eqno(15.4.27)$$ for
$1\leq p,q\leq n$.

 For $\vec i\in \mbb{N}^n$, we set
$$x^{(\vec i)}=x_1^{i_1}x_2^{i_2-i_1}x_3^{i_3-i_2}\cdots x_{n-1}^{i_{n-1}-i_{n-2}}x_n^{2i_n
-i_{n-1}},\eqno(15.4.28)$$
$$\es(\vec i)=\sum_{p=1}^{n-1}i_p\ves_{p+1,p}+i_n\ves_{2n,n}.\eqno(15.4.29)$$
Let
$$u_\lmd=\sum_{\vec i\in\mbb{N}^n}\frac{(-1)^{i_1+i_2+\cdots+i_{n-1}}}{i_1!i_2!\cdots i_n!
2^{i_n} \la\lmd_n\ra_{i_n}}C^{\es(\vec i)}v_\lmd x^{(\vec
i)}x^\ast.\eqno(15.4.30)$$ It is verified that $u_\lmd$ is a
singular vector of weight $\lmd$ in $\td{M}$; that is,
$$C_{i,j}(u_\lmd)=C_{p,n+q}(u_\lmd)=0\qquad\for\;\;1\leq i<j\leq n;\;1\leq p\leq q\leq n.
\eqno(15.4.31)$$ Hence there exists $sp(2n,\mbb R)$-module
homomorphism $\Phi: M(\lmd)\rta \td{M}$ such that
$$\Phi(v_\lmd)=u_\lmd.\eqno(15.4.32)$$
View $\Phi$ as a function in $\{x_1,x_2,...,x_n\}$ taking value in
the spaces of  linear transformations on $M(\lmd)$. In the rest of
this section, we want to calculate $E_C(z_1,z_2,...,z_n)$ in
(15.2.47).

For $\vec i\in\mbb{N}^n$ and $\be \in \G_C$ (cf. (15.4.6)), we have
\begin{eqnarray*}\qquad\;\;& &C^\be(C^{\es(\vec i)}v_\lmd)\\ &=&\sum_{2\leq j_1<j_2\leq n}\;\sum_{k_1=1}^n
\;\sum_{k_2=2}^n\;\sum_{k_3=k_2}^n(\sum_{p_{j_2,j_1},q_{k_1,k_2},q_{k_2,k_2}^1=0}^{\infty}\\
& &\left[\prod_{r=2}^n{\be_{n+r,r}\choose
q_{r,r},q_{r,r}^1}{\be_{m+r,1}\choose q_{1,r}}\right]\la
i_{n-1}\ra_{_{\sum_{r=1}^nq_{r,n}+\sum_{r=2}^nq^1_{n,r}+q_{n,n}^1}}
\\& &\times
\left[\prod_{2\leq l_1<l_2\leq n}{\be_{l_2,l_1}\choose
p_{l_2,l_1}}{\be_{n+l_2,l_1}\choose q_{l_1,l_2},q_{l_2,l_1},
q^1_{l_2,l_1}}\right] \\
& &\times \left[\prod_{s=1}^{n-2}\la
i_s\ra_{_{\sum_{r=s+2}^np_{r,s+1}+\sum_{r=1}^nq_{r,s+1}+\sum_{r=2}^{s+1}q^1_{s+1,r}+
\sum_{r=s+1}^nq^1_{r,s+1} }}\right]
\\ & &\times
E^{-\dlt_{l_2,l_1+1}(\sum_{r=l_2+1}^n
p_{r,l_2}+\sum_{s=1}^nq_{s,l_2}+\sum_{s=2}^{l_2}q^1_{l_2,s}+\sum_{s=l_2}^nq^1_{s,l_2})
\ves_{l_2,l_1}}
\\ & &\times E^{\sum_{1\leq s_1<s_2\leq n}(q_{s_2,s_1+1}+q_{s_1,s_2+1}+q_{s_2+1,s_1+1}^1
-q_{s_1,s_2}-q_{s_2,s_1}-q_{s_2,s_1}^1)\ves_{n+s_2,s_1}}\\ & &\times
E^{\es(\vec i)+\be
+\sum_{r=1}^n(q_{r,r+1}+q_{r+1,r+1}^1-q_{r,r}-q_{r,r}^1)\ves_{n+r,r}}v_\lmd)\\
& &\times E^{\sum_{1\leq l_1<l_2\leq
n}(p_{l_2,l_1+1}-p_{l_2,l_1})\ves_{l_2,l_1}}\hspace{7cm}
(15.4.33)\end{eqnarray*}
 by the calculations in (15.3.44)-(15.3.49)
and (15.4.19)-(15.4.23), where we treat
$$p_{r,r}=q_{1,1}=q^1_{n+1,n+1}=q^1_{r,1}=q^1_{n+1,r}=q_{s,n+1}=q_{n+1,s}=0\eqno(15.4.34)$$
for $2\leq r\leq n$ and $1\leq s\leq n.$ Hence by (15.4.25),
(15.4.27) and (15.4.33), we obtain
\begin{eqnarray*}& &C^\al(C^{\es(\vec i)}v_\lmd x^{(\vec i)}x^\ast)\\ &=&
\sum_{2\leq j_1<j_2\leq
n}\;\sum_{k_1=1}^n\;\sum_{k_2=2}^n\;\sum_{k_3=k_1}^n\;\sum_{k_4=k_2}^n\{\sum_{p^0_{k_2,1},
p_{j_2,j_1}^0,p_{j_2,j_1},q^0_{k_3,k_1},q^1_{k_3,k_1},q_{k_1,k_2}=0}^{\infty}\\
& &\times \frac{\la
i_{n-1}\ra_{_{\sum_{r=1}^nq_{r,n}+\sum_{r=2}^nq^1_{n,r}+q_{n,n}^1}}}
{2^{\sum_{r=1}^nq^0_{r,r}}}\la
2i_n-i_{n-1}-1/2\ra_{_{2q^0_{n,n}+\sum_{r=1}^{n-1}q^0_{n,r}}}
\\& &\times \left[\prod_{2\leq l_1<l_2\leq
n}{\al_{l_2,l_1}\choose  p_{l_2,l_1}^0,
p_{l_2,l_1}}{\al_{n+l_2,l_1}\choose q_{l_2,l_1}^0,q_{l_2,l_1}^1,
q_{l_1,l_2},q_{l_2,l_1}}\right]\\ &
&\times\left[\prod_{m=2}^n{\al_{m,1}\choose
p_{m,1}^0}\right]{\al_{n+1,1}\choose
q_{1,1}^0}\left[\prod_{r=2}^n{\al_{n+r,r}\choose q_{r,r}^0,
q_{r,r}^1,q_{r,r}}{\al_{n+r,1}\choose q_{r,1}^0,q_{1,r}}\right]\\
& &\times \left[\prod_{s=1}^{n-2} \la
i_s\ra_{_{\sum_{r=s+2}^np_{r,s+1}+\sum_{r=1}^nq_{r,s+1}+\sum_{r=2}^{s+1}q^1_{s+1,r}+
\sum_{r=s+1}^nq^1_{r,s+1} }}\right]\\
& &\times \left[\prod_{s=1}^{n-1} \la
i_s-i_{s-1}-1/2\ra_{_{\sum_{r=s+1}^np_{r,s}^0+\sum_{r=1}^sq_{s,r}^0+\sum_{r=s}^nq_{r,s}^0}}
\right]
\\ & &\times
E^{-\sum_{1\leq l_1<l_2\leq n}\dlt_{l_2,l_1+1}
(\sum_{r=l_2+1}^np_{r,l_2}+\sum_{s=1}^nq_{s,l_2}+\sum_{s=2}^{l_2}q^1_{l_2,s}+\sum_{s=l_2}^n
q^1_{s,l_2})\ves_{l_2,l_1}}
\\ & &\times E^{\sum_{1\leq s_1<s_2\leq n}(q_{s_2,s_1+1}+q_{s_1,s_2+1}+q_{s_2+1,s_1+1}^1
-q_{s_1,s_2}-q_{s_2,s_1}-q_{s_2,s_1}^0-q_{s_2,s_1}^1)\ves_{n+s_2,s_1}}\\
& &\times E^{\al+\es(\vec
i)+\sum_{r=1}^n(q_{r,r+1}+q_{r+1,r+1}^1-q_{r,r}-q_{r,r}^0-q_{r,r}^1)
\ves_{n+r,r}}\\ & & \times E^{\sum_{1\leq l_1<l_2\leq
n}(p_{l_2,l_1+1}-p_{l_2,l_1}-p_{l_2,l_1}^0)\ves_{l_2,l_1}} v_\lmd\\
& &
\left[\prod_{s=1}^{n-1}x_s^{i_s-i_{s-1}-1/2+\sum_{r=1}^{s-1}p^0_{s,r}-\sum_{r=s+1}^np_{r,s}^0
-\sum_{r=1}^sq_{s,r}^0-\sum_{r=s}^nq_{r,s}^0}\right]\\ & &\times
x_n^{2i_n-i_{n-1}-1/2+
\sum_{r=1}^{n-1}p^0_{n,r}-2q^0_{n,n}-\sum_{r=1}^{n-1}q^0_{n,r}}\},\hspace{6.2cm}(15.4.35)
\end{eqnarray*}
where we treat $i_{-1}=0$.

 In order to  calculate  the Etingof trace (15.2.47), we have to let in (15.4.35):
$$ p_{r_2,r_1}=p_{r_2,r_1-1}^0+p_{r_2,r_1-1}\qquad\for\;\;2\leq r_1<r_2\leq n,\eqno(15.4.36)$$
$$q_{r,r+1}+q_{r+1,r+1}^1=q^0_{r,r}+q_{r,r}++q_{r,r}^1\qquad\for\;\;r\in\ol{1,n-1},\eqno(15.4.37)$$
$$q_{s_2,s_1+1}+q_{s_1,s_2+1}+q_{s_2+1,s_1+1}^1=q_{s_1,s_2}+q_{s_2,s_1}+q_{s_2,s_1}^0
+q_{s_2,s_1}^1\eqno(15.4.38)$$ for $1\leq s_1<s_2\leq n$ and
$$i_r=p^0_{r+1,r}+p_{r+1,r}+\sum_{s=r+2}^np_{s,r+1}+\sum_{s=1}^nq_{s,r+1}
+\sum_{s=2}^{r+1}q_{r+1,s}^1+\sum_{s=r+1}^nq_{s,r+1}^1\eqno(15.4.39)$$
for $r\in\ol{1,n-1}$,
$$i_n=q_{n,n}^0+q_{n,n}+q^1_{n,n}.\eqno(15.4.40)$$
In particular, (15.4.36) implies
$$p_{r_2,r_1}=\sum_{s=1}^{r_1-1}p^0_{r_2,s}\qquad\for\;\;2\leq r_1<r_2\leq n.\eqno(15.4.41)$$

Set
$$\Im_r=\sum_{s=1}^nq_{s,r+1}+\sum_{s=2}^{r+1}q_{r+1,s}^1+\sum_{s=r+1}^nq_{s,r+1}^1\eqno(15.4.42)$$
for $r\in\ol{1,n-1}$. By (15.4.34) and (15.4.37), we have
$$q_{1,2}+q^1_{2,2}=q_{1,1}^0.\eqno(15.4.43)$$
Moreover, (15.4.38) implies
$$q_{i,2}+q_{1,i+1}+q_{i+1,2}^1=q_{1.i}+q_{i,1}^0\qquad\for\;\;i\in\ol{2,n}.\eqno(15.4.44)$$
$$2(q_{1,2}+q^1_{2,2})+\sum_{i=2}^nq_{i,2}+\sum_{r=3}^n(q_{1,r}+q_{r,2}^1)=2q^0_{1,1}
+\sum_{i=2}^n(q_{1,i}+q^0_{i,1}),\eqno(15.4.45)$$ or equivalently,
$$\sum_{i=1}^nq_{i,2}+q^1_{2,2}+\sum_{r=2}^nq_{r,2}^1=q^0_{1,1}+\sum_{i=1}^nq^0_{i,1}.\eqno(15.4.46)$$
Thus
$$\Im_1=q^0_{1,1}+\sum_{i=1}^nq^0_{i,1}.\eqno(15.4.47)$$

Let $r\in\ol{2,n-1}$. By (15.4.38), we have
$$q_{i,r+1}+q_{r,i+1}+q_{i+1,r+1}^1=q_{r,i}+q_{i,r}+q_{i,r}^0+q_{i,r}^1\eqno(15.4.48)$$
for $i\in\ol{r+1,n}$ and
$$q_{r,j+1}+q_{j,r+1}+q_{r+1,j+1}^1=q_{j,r}+q_{r,j}+q_{r,j}^0+q_{r,j}^1\eqno(15.4.49)$$
for $j\in\ol{1,r-1}$. Expressions (15.4.37), (15.4.48) and (15.4.49)
imply
\begin{eqnarray*} & &2(q_{r,r+1}+q_{r+1,r+1}^1)+\sum_{i=r+1}^n(q_{i,r+1}+q_{r,i+1}
+q_{i+1,r+1}^1)+\sum_{j=1}^{r-1}(q_{r,j+1}+q_{j,r+1}+q_{r+1,j+1}^1)\\
&=& 2(q^0_{r,r}+q_{r,r}++q_{r,r}^1)
+q_{i,r}^1)+\sum_{j=1}^{r-1}(q_{j,r}+q_{r,j}+q_{r,j}^0+q_{r,j}^1)\\&
&
+\sum_{s=r+1}^n(q_{r,s}+q_{s,r}+q_{s,r}^0+q_{s,r}^1),\hspace{7.9cm}(15.4.50)\end{eqnarray*}
or equivalently,
\begin{eqnarray*}\qquad & &\sum_{i=1}^nq_{i,r+1}+\sum_{j=2}^{r+1}q_{r+1,j}^1+\sum_{i=r+1}^n q_{i,r+1}
\\ &=&\sum_{j=1}^rq_{r,j}^0+\sum_{s=r}^nq_{s,r}^0+\sum_{s=1}^nq_{s,r}+\sum_{j=2}^rq_{r,j}^1
+\sum_{s=r}^nq_{s,r}^1\hspace{4.4cm}(15.4.51)\end{eqnarray*} by
(15.4.34). So
$$\Im_r-\Im_{r-1}=\sum_{j=1}^rq_{r,j}^0+\sum_{i=r}^nq_{i,r}^0.\eqno(15.4.52)$$
Furthermore, by (15.4.39), (15.4.40), (15.4.42), (15.4.52) and
induction, we have
$$
i_r=\sum_{s=r+1}^n\sum_{j=1}^rp_{s,j}^0+\sum_{s=1}^r(\sum_{j=1}^sq_{s,j}^0+\sum_{j=s}^nq_{j,s}^0)
\eqno(15.4.53)$$ for $r\in\ol{1,n-1}$. In particular,
$$
i_r-i_{r-1}=\sum_{s=r+1}^np^0_{s,r}-\sum_{s=1}^{r-1}p^0_{r,s}+\sum_{j=1}^rq_{r,j}^0
+\sum_{s=r}^nq_{s,r}^0\eqno(15.4.54)$$ for $r\in\ol{1,n-1}$.

Next we want to calculate $i_n$. First we have
$$q_{n,n}+q^1_{n,n}=q_{n,n-1}+q_{n-1,n}+q_{n,n}^1+q_{n,n-1}^1+q_{n,n-1}^0\eqno(15.4.55)$$
by (15.4.38). Moreover, for $r\in\ol{1,n-2}$, we have
\begin{eqnarray*}& &\sum_{i=0}^{2r}q_{n-i,n-2r+i}+\sum_{i=0}^{r-1}q_{n-i,n-2r+1+i}
\\ &=&q_{n,n-2r}+\sum_{i=1}^r(q_{n-i,n-2r+i}+q_{n-2r-1+i,n-i+1}+q_{n-i+1,n-2r+i}^1)\\ &=&
\sum_{i=1}^r(q_{n-i,n-2r-1+i}+q_{n-2r-1+i,n-i}+q_{n-i,n-2r-1+i}^0+q_{n-i,n-2r-1+i}^1)
\\ && +q_{n,n-2r-1}+q_{n-2r-1,n}+q_{n,n-2r-1}^0+q_{n,n-2r-1}^1
\\ &=& \sum_{i=0}^r(q_{n-i,n-2r-1+i}+q_{n-2r-1+i,n-i})+\sum_{i=0}^rq_{n-i,n-2r-1+i}^1+
 \sum_{i=0}^rq_{n-i,n-2r-1+i}^0\\ &=&\sum_{i=0}^{2r+1}q_{n-i,n-2r-1+i}+ \sum_{i=0}^r
 q_{n-i,n-2r-1+i}^0+ \sum_{i=0}^rq_{n-i,n-2r-1+i}^1.\hspace{3.2cm}(15.4.56)\end{eqnarray*}
by (15.4.38). Furthermore, for $r\in\ol{0,n-2}$, we get
\begin{eqnarray*}& &\sum_{i=0}^{2r+1}q_{n-i,n-2r-1+i}+\sum_{i=0}^rq_{n-i,n-2r+i}\\
&=&q_{n,n-2r-1}+q_{n-r-1,n-r}+q_{n-r,n-r}^1+\sum_{i=1}^r(q_{n-i,n-2r-1+i}\\
& & +q_{n-2r-2+i,n+1-i} +q_{n-i+1,n-2r-1+i}^1)
\\ &=& q_{n,n-2r-2}+q_{n-2r-2,n}+q^0_{n,n-2r-2}+q^1_{n,n-2r-2}+q_{n-r-1,n-r-1}
+q_{n-r-1,n-r-1}^0\\ & &+q_{n-r-1,n-r-1}^1
+\sum_{i=1}^r(q_{n-i,n-2r-2+i}+q_{n-2r-2+i,n-i}+q_{n-i,n-2r-2+i}^0+q_{n-i,n-2r-2+i}^1)
\\ &=&
q_{n-r-1,n-r-1}+\sum_{i=0}^r(q_{n-i,n-2r-2+i}+q_{n-2r-2+i,n-i})
\\ & &+\sum_{i=0}^{r+1}q_{n-i,n-2r-2+i}^0+\sum_{i=0}^{r+1}q_{n-i,n-2r-2+i}^1\\
&=& \sum_{i=0}^{2r+2}
q_{n-i,n-2r-2+i}+\sum_{i=0}^{r+1}q_{n-i,n-2r-2+i}^0+\sum_{i=0}^{r+1}q_{n-i,n-2r-2+i}^1
\hspace{3.4cm}(15.4.57)\end{eqnarray*} based on (15.4.37) and
(15.4.38).

By (15.4.55)-(15.4.57) and induction, we obtain
\begin{eqnarray*}\qquad \qquad q_{n,n}+q_{n,n}^1&=&\sum_{r=0}^{n-1}\sum_{s=0}^rq_{n-s,n-2r-1+s}^0+
\sum_{r=0}^{n-2}\sum_{s=0}^rq_{n-s,n-2r-2+s}^0\\ &=&\sum_{1\leq
r_1\leq r_2\leq n}q^0_{r_2,r_1}
-q_{n,n}^0.\hspace{6.1cm}(15.4.58)\end{eqnarray*} Therefore,
$$i_n=\sum_{1\leq r_1\leq r_2\leq n}q_{s_2,s_1}^0.\eqno(15.4.59)$$
Now
\begin{eqnarray*}\qquad \qquad 2i_n-i_{n-1}&=&\sum_{1\leq r_1\leq r_2\leq
n}2q_{s_2,s_1}^0-\sum_{r=1}^{n-1}(p^0_{n,r}
+\sum_{p=1}^rq_{r,i}^0+\sum_{s=r}^nq_{s,r}^0)\\
&=&q_{n,n}^0+\sum_{r=1}^nq_{n,r}^0
-\sum_{s=1}^{n-1}p^0_{n,s}.\hspace{5.6cm}(15.4.60)\end{eqnarray*}

By (15.4.37) and induction on $r$, we get
$$q_{r,r}^1=\sum_{s=1}^{r-1}(q_{s,s}^0+q_{s,s}-q_{s,s+1})\qquad\for\;\;r\in\ol{2,n}.\eqno(15.4.61)$$
Moreover, (15.4.38) and induction imply
$$q_{n,r}=\sum_{s=1}^{r-1}(q^0_{n,s}+q_{s,n}+q_{n,s}^1)\qquad\for\;\;r\in\ol{2,n}\eqno(15.4.62)$$
and
\begin{eqnarray*}\qquad \qquad q_{r_2,r_1}^1&=&\sum_{s=1}^{r_2-1}(q_{r_2-s,r_1-s}^0+q_{r_2-s,r_1-s}+q_{r_1-s,r_2-s})
\\ & &-\sum_{s=1}^{r_1-1}(q_{r_2-s,r_1-s+1}+q_{r_1-s,r_2-s+1})\hspace{5cm}(15.4.63)\end{eqnarray*}
for $2\leq r_1<r_2\leq n$. Furthermore, (15.4.62) and (15.4.63) show
\begin{eqnarray*}\qquad q_{n,r}&=&\sum_{s=1}^{r-1}(q_{n,i}^0+q_{i,n})+\sum_{s=2}^{r-1}\sum_{j=1}^{s-1}(q_{n-j,s-j}^0+
q_{n-j,s-j}+q_{s-j,n-j})\\ &
&-\sum_{s=2}^{r-1}\sum_{j=1}^{s-1}(q_{n-j,s-j+1}+q_{s-j,n-j+1})\hspace{6cm}
(15.4.64)\end{eqnarray*} for $r\in\ol{2,n}$. Set
$$\vec{p^0}!=\prod_{1\leq r_1<r_2\leq n}p_{r_2,r_1}^0!,\qquad \vec{q^0}!=
\prod_{1\leq r_1\leq r_2\leq n}q_{r_2,r_1}^0!.\eqno(15.4.65)$$

By (15.2.47), (15.4.35), (15.4.53), (15.4.59), (15.4.61), (15.4.63)
and (15.4.64), we have
\begin{eqnarray*}& &E_C(z_1,z_2,...,z_n)\\ &=& x^\ast\sum_{\al\in
\G_C}\; \sum_{1\leq j_1<j_2\leq
n}\;\sum_{k_1=1}^n\;\sum_{k_2=2}^n\;\sum_{k_3=k_1}^n\;\sum_{k_4=1}^{n-1}\;
\sum_{p_{j_2,j_1}^0, q_{k_3,k_1}^0,q_{k_4,k_2}=0}^{\infty}\\
& &\frac{(-1)^{i_1+i_2+\cdots+i_{n-1}}}{i_1!i_2!\cdots i_n!
2^{i_n+\sum_{r=1}^nq^0_{r,r}} \la\lmd_n\ra_{i_n}} \la
i_{n-1}\ra_{_{\sum_{r=1}^nq_{r,n}+\sum_{r=2}^nq_{n,r}^1+q^1_{n,n}}}
\\ & &\times \left[\prod_{2\leq l_1\leq l_2\leq
n}{\al_{l_2,l_1}\choose p^0_{l_2,l_1},
p_{l_2,l_1}}{\al_{n+l_2,l_1}\choose
q^0_{l_2,l_1},q_{l_2,l_1}^1,q_{l_1,l_2},q_{l_2,l_1}}\right]\\
& &\times\left[\prod_{m=2}^n{\al_{m,1}\choose
p^0_{m,1}}\right]{\al_{n+1,1}\choose
q^0_{1,1}}\left[\prod_{r=2}^n{\al_{n+r,r}\choose
q^0_{r,r},q_{r,r}^1,q_{r,r}}{\al_{n+r,1}\choose
q^0_{r,1},q_{1,r}}\right]\\ & &\times \left[\prod_{s=1}^{n-2} \la
i_s\ra_{_{\sum_{r_1=s+2}^np_{r_1,s+1}+\sum_{r=1}^nq_{r,s+1}+\sum_{r=2}^{s+1}q_{s+1,r}^1
+\sum_{r=s+1}^nq_{r,s+1}^1}}\right] \\
& &\times \left[\prod_{s=1}^{n-1}\la i_s
-i_{s-1}-1/2\ra_{_{\sum_{r=s+1}^np_{r,s}^0+\sum_{r=1}^sq_{s,r}^0+\sum_{r=s}^nq_{r,s}^0}}\right]\\
& &\times\la
2i_n-i_{n-1}-1/2\ra_{_{2q^0_{n,n}+\sum_{r=1}^{n-1}q^0_{n,r}}}(z_1z_2\cdots
z_{n-1})^{-1/2}z_n^{\lmd_n}\\
& &\times\left[\prod_{1\leq r_1<r_2\leq
n}\left(\frac{z_{r_2}}{z_{r_1}}\right)^{\al_{r_2,r_1}}\frac{1}{(z_{r_1}z_{r_2})
^{\al_{n+r_2,r_1}}}\right]\left[\prod_{r=1}^n\frac{1}{z_r^{\al_{n+r,r}}}\right]\\
&=& x^\ast\sum_{1\leq j_1<j_2\leq
n}\;\sum_{k_1=1}^n\;\sum_{k_2=2}^n\;\sum_{k_3=k_1}^n\;\sum_{k_4=1}^{n-1}
\;\sum_{p_{j_2,j_1}^0,
q_{k_3,k_1}^0,q_{k_4,k_2}=0}^{\infty}\frac{(-1)^{i_1+i_2+\cdots+i_{n-1}}}{\vec
p^0! 2^{i_n+\sum_{r=1}^nq^0_{r,r}}\la \lmd_n\ra_{i_n}}\\ &
&\times\frac{\la
2i_n-i_{n-1}-1/2\ra_{_{2q^0_{n,n}+\sum_{r=1}^{n-1}q^0_{n,r}}}}{q_{n,n}^0!q_{n,n}!q_{n,n}^1!}\left[\prod_{r=2}^{n-1}
{q^0_{r,r}+q_{r,r}^1+q_{r,r}\choose
q_{r,r}^1,q_{r,r}}\right]\hspace{9cm}
\end{eqnarray*}
\begin{eqnarray*} &
&\times \left[\prod_{1\leq l_1\leq l_2\leq n}
{q^0_{l_2,l_1}+q_{l_2,l_1}^1+q_{l_1,l_2}+q_{l_2,l_1} \choose
q_{l_2,l_1}^1,q_{l_1,l_2},q_{l_2,l_1}}\right]\\ &
&\times\left[\prod_{s=1}^{n-1} \la
i_s-i_{s-1}-1/2\ra_{_{\sum_{r=s+1}^np_{r,s}^0+\sum_{r=1}^sq_{s,r}^0+\sum_{r=s}^nq_{r,s}^0}}
\right] \\& &\times \left[\prod_{1\leq r_1<r_2\leq
n}\frac{z_{r_1}^2z_{r_2}^{1+\sum_{s=1}^{r_1}p_{r_2,s}^0}}{(z_{r_1}-z_{r_2})
^{1+\sum_{s=1}^{r_1}p_{r_2,s}^0}(z_{r_1}z_{r_2}-1)^{1+q^0_{r_2,r_1}+q_{r_2,r_1}^1
+q_{r_1,r_2}+q_{r_2,r_1}}}\right]\\ & &\times
\left[\prod_{r=1}^n\frac{z_r^2}{(z_r^2-1)^{1+q^0_{r,r}+q^1_{r,r}+q_{r,r}}}\right]
(z_1z_2\cdots
z_{n-1})^{-1/2}z_n^{\lmd_n}.\hspace{4.5cm}(15.4.66)\end{eqnarray*}
Note that
\begin{eqnarray*}& &\frac{\prod_{r=2}^{n-1}
{q^0_{r,r}+q_{r,r}^1+q_{r,r}\choose
q_{r,r}^1,q_{r,r}}}{q_{n,n}^0!q_{n,n}!q_{n,n}^1!} \left[\prod_{1\leq
l_1\leq l_2\leq n}
{q^0_{l_2,l_1}+q_{l_2,l_1}^1+q_{l_1,l_2}+q_{l_2,l_1}
\choose q_{l_2,l_1}^1,q_{l_1,l_2},q_{l_2,l_1}}\right]\\
&=&\frac{\prod_{r=2}^{n-1} {q_{r+1,r+1}^1+q_{r,r+1}\choose
q_{r,r}^1,q_{r,r}}}{q_{n,n}^0!q_{n,n}!q_{n,n}^1!} \left[\prod_{1\leq
l_1\leq l_2\leq n} {q_{l_2+1,l_1+1}^1+q_{l_1,l_2+1}+q_{l_2,l_1+1}
\choose
q_{l_2,l_1}^1,q_{l_1,l_2},q_{l_2,l_1}}\right]\\
&=&\frac{\prod_{r=1}^{n-1} {q_{r+1,r+1}^1+q_{r,r+1}\choose
q_{r,r+1}}}{\vec q^0!} \left[\prod_{1\leq l_1\leq l_2\leq n-1}
{q_{l_2+1,l_1+1}^1+q_{l_1,l_2+1}+q_{l_2,l_1+1} \choose
q_{l_1,l_2+1},q_{l_2,l_1+1}}\right]\\
&=&\frac{\prod_{r=1}^{n-1} {q^0_{r,r}+q_{r,r}^1+q_{r,r}\choose
q_{r,r+1}}}{\vec q^0!} \left[\prod_{1\leq l_1\leq l_2\leq n-1}
{q^0_{l_2,l_1}+q_{l_2,l_1}^1+q_{l_1,l_2}+q_{l_2,l_1} \choose
q_{l_1,l_2+1},q_{l_2,l_1+1}}\right]\hspace{2.2cm}(15.4.67)\end{eqnarray*}
according to (15.4.37) and (15.4.38). Furthermore,
$$q_{r,r}^0+q^1_{r,r}+q_{r,r}=\sum_{s=1}^r(q_{s,s}^0+q_{s,s})-\sum_{s=1}^{r-1}q_{s,s+1}
\eqno(15.4.68)$$ by (15.4.61) and
\begin{eqnarray*}& & q_{n,r}^0+q_{n,r}^1+q_{r,n}+q_{n,r}\\
&=&\sum_{s=1}^r(q^0_{n,s}+q_{s,n})+\sum_{s=0}^{r-2}\sum_{j=1}^{r-s-1}(q^0_{n-j,r-s-j}
+q_{n-j,r-s-j}+q_{r-s-j,n-j})\\ &
&-\sum_{s=0}^{r-2}\sum_{j=1}^{r-s-1}(q_{n-j,r-s-j+1}+q_{r-s-j,n-j+1})\\
&=&\sum_{s=0}^{r-1}\sum_{j=0}^{r-s-1}q^0_{n-j,r-s-j}+\sum_{j=0}^{r-1}q_{r-j,n-j}-\sum_{j=1}^{r-1}
q_{r-j,n-j+1}\hspace{4.5cm}(15.4.69)\end{eqnarray*} by (15.4.62) and
(15.4.64) for $r\in\ol{1,n-1}$, and
\begin{eqnarray*}
q_{r_2,r_1}^0+q_{r_2,r_1}^1+q_{r_1,r_2}+q_{r_2,r_1}
&=&\sum_{s=0}^{r_1-1}(q_{r_2-s,r_1-s}^0+q_{r_2-s,r_1-s}+q_{r_1-s,r_2-s})
\\ & &-\sum_{s=1}^{r_1-1}(q_{r_2-s,r_1-s+1}+q_{r_1-s,r_2-s+1})
\hspace{2.6cm}(15.4.70)\end{eqnarray*} by (15.4.63) for $1\leq r_1<
r_2\leq n-1.$

Let
$$ y_{r_2,r_1}=\frac{1}{z_{r_1}z_{r_2}-1}\qquad\for\;\;1\leq
r_1\leq r_2\leq n.\eqno(15.4.71)$$ Set
$$w_{r_2,r_1}=\frac{y_{r_2,r_1}y_{r_2+1,r_1+1}\cdots
y_{n-1,n+r_1-r_2-1}}{y_{r_2+1,r_1}y_{r_2+2,r_1+1}\cdots
y_{n,n+r_1-r_2-1}}\qquad\for\;\;2\leq r_1\leq r_2\leq
n-1,\eqno(15.4.72)$$
$$ w_{k,n}=y_{n,k}\qquad\for\;\;k\in\ol{1,n-1}\eqno(15.4.73)$$
and
$$w_{r_1,r_2}=\frac{y_{r_2,r_1}y_{r_2+1,r_1+1}\cdots
y_{n,n+r_1-r_2}}{y_{r_2,r_1+1}y_{r_2+1,r_1+2}\cdots
y_{n-1,n+r_1-r_2}}\qquad\for\;\;1\leq r_1<r_2\leq
n-1.\eqno(15.4.74)$$ Based on (15.4.68)-(15.4.70), we have
\begin{eqnarray*}& &\prod_{1\leq r_1\leq r_2\leq
n}y_{r_2,r_1}^{q^0_{r_2,r_1}+q^1_{r_2,r_1}+q_{r_2,r_1}+q_{r_1,r_2}}=\left[\prod_{r_1=1}^{n-1}
\prod_{r_2=2}^nw_{r_1,r_2}^{q_{r_1,r_2}}\right]\\ &&\times
\left[\prod_{r=1}^n(y_{r,r}y_{r+1,r+1}\cdots
y_{n,n})^{q_{r,r}^0}\right]\left[\prod_{r=1}^{n-1}(y_{n,r}y_{n,r+1}\cdots
y_{n,n})^{q_{n,r}^0}\right]\prod_{1\leq r_1<r_2\leq n-1}\\ &
&(y_{r_2,r_1}y_{r_2+1,r_1+1}\cdots
y_{n-1,n+r_1-r_2-1}y_{n,n+r_1-r_1}y_{n,n+r_1-r_1+1}\cdots
y_{n,n})^{q^0_{r_2,r_1}}.\hspace{1.9cm}(15.4.75)\end{eqnarray*}

Next we have
\begin{eqnarray*}\hspace{2cm}& &\sum_{q_{n-1,n}=0}^{\infty}
\left(\!\!\begin{array}{c}q^0_{n-1,n-1}+q_{n-1,n-1}^1+q_{n-1,n-1}\\
q_{n-1,n}\end{array}\!\!\right)w_{n-1,n}^{q_{n-1,n}}\\
&=&(1+w_{n-1,n})^{q^0_{n-1,n-1}+q_{n-1,n-1}^1+q_{n-1,n-1}}\\
&=&(1+w_{n-1,n})^{\sum_{s=1}^{n-1}(q_{s,s}^0+q_{s,s})-\sum_{s=1}^{n-2}q_{s,s+1}}
\hspace{5cm}(15.4.76)\end{eqnarray*} by (15.4.68). Set
$$ w^{(2n-1)}_{r,r+1}=\frac{w_{r,r+1}}{1+w_{n-1,n}},\qquad
w_{s,s}^{(2n-1)}=(1+w_{n-1,n})w_{s,s}\eqno(15.4.77)$$ for
$r\in\ol{1,n-2},\;s\in\ol{1,n-1}$ and
$$
w^{(2n-1)}_{r_1,r_2}=w_{r_1,r_2}\qquad\for\;\;r_1\in\ol{1,n-1},\;r_1,r_1+1\neq
r_2\in\ol{2,n}.\eqno(15.4.78)$$ Suppose that we have defined
$\{w_{r_1,r_2}^{(2k+1)}\mid r_1+r_2<2k+1\}$. Set
$$\ell=\max\{1,2k-n\}.\eqno(15.4.79)$$
Note
\begin{eqnarray*}&&
\sum_{q_{\ell,2k-\ell},q_{\ell+1,2k-\ell-1},\cdots,q_{2k-\ell-1,\ell+1}=0}^{\infty}\;
\prod_{r=\ell}^{k-1}{q^0_{2k-1-r,r}+q^1_{2k-1-r,r}+q_{2k-1-r,r}+q_{r,2k-1-r}\choose
q_{r,2k-r},q_{2k-r-1,r+1}}\\ & &\times
w_{r,2k-r}^{(2k+1)q_{r,2k-r}}w_{2k-r-1,r+1}^{(2k+1)q_{2k-r-1,r+1}}\\
&=&\prod_{r=\ell}^{k-1}(1+w_{r,2k-r}^{(2k+1)}+w_{2k-r-1,r+1}^{(2k+1)})
^{q^0_{2k-1-r,r}+q^1_{2k-1-r,r}+q_{2k-1-r,r}+q_{r,2k-1-r}}\\
&=&\prod_{r=\ell}^{k-1}(1+w_{r,2k-r}^{(2k+1)}+w_{2k-r-1,r+1}^{(2k+1)})^{\sum_{s=1}^r
(q^0_{2k-1-2r+s,s}+q_{2k-1-2r+s,s}+q_{s,2k-1-2r+s})}\\ &
&\times(1+w_{r,2k-r}^{(2k+1)}+w_{2k-r-1,r+1}^{(2k+1)})^{-\sum_{s=1}^r
(q_{2k-1-2r+s,s+1}+q_{s,2k-2r+s})}\hspace{3.2cm}(15.4.80)\end{eqnarray*}
by (15.4.70). Set
$$w_{2k-1-2r+s,s}^{(2k)}=w_{2k-1-2r+s,s}^{(2k+1)}(1+w_{r,2k-r}^{(2k+1)}+w_{2k-r-1,r+1}^{(2k+1)}),
\eqno(15.4.81)$$
$$w_{s,2k-1-2r+s}^{(2k)}=w_{s,2k-1-2r+s}^{(2k+1)}(1+w_{r,2k-r}^{(2k+1)}+w_{2k-r-1,r+1}^{(2k+1)})
\eqno(15.4.82)$$ for $1\leq s\leq r\leq k-1$;
$$w_{2k-1-2r+s,s+1}^{(2k)}=\frac{w_{2k-1-2r+s,s+1}^{(2k+1)}}{1+w_{r,2k-r}^{(2k+1)}
+w_{2k-r-1,r+1}^{(2k+1)}}, \eqno(15.4.83)$$
$$w_{s,2k-2r+s}^{(2k)}=\frac{w_{s,2k-2r+s}^{(2k+1)}}{1+w_{r,2k-r}^{(2k+1)}+
w_{2k-r-1,r+1}^{(2k+1)}} \eqno(15.4.84)$$ for $1\leq s\leq r\leq
k-1$ and
$$w^{(2k)}_{r_1,r_2}=w^{(2k+1)}_{r_1,r_2}\;\;\mbox{for the other
pairs}\;\;(r_1,r_2)\;\;\mbox{such that}\;\;r_1+r_2\leq
2k.\eqno(15.4.85)$$

Let
$$\iota=\max\{1,2k-1-n\}.\eqno(15.4.86)$$
Moreover,
\begin{eqnarray*}&&\sum_{q_{\iota,2k-\iota-1},q_{\iota+1,2k-\iota-2},\cdots,
q_{2k-\iota-2,\iota+1}=0}
^{\infty}[\prod_{r=\iota}^{k-2}{q^0_{2k-2-r,r}+q^1_{2k-2-r,r}+q_{2k-2-r,r}+q_{r,2k-2-r}
\choose q_{r,2k-r-1},q_{2k-r-2,r+1}}\\ &&\times
w_{r,2k-r-1}^{(2k)q_{r,2k-r-1}}w_{2k-r-2,r+1}^{(2k)q_{2k-r-2,r+1}}]{
q^0_{k-1,k-1}+q^1_{k-1,k-1}
+q_{k-1,k-1}\choose q_{k-1,k}}w_{k-1,k}^{(2k)q_{k-1,k}}\\
&=&\left[\prod_{r=\iota}^{k-2}(1+w_{r,2k-r-1}^{(2k)}+w^{(2k)}_{2k-r-2,r+1})
^{q^0_{2k-2-r,r}+q^1_{2k-2-r,r}+q_{2k-2-r,r}+q_{r,2k-2-r}}\right]\\
&&\times (1+w_{k-1,k}^{(2k)})^{q^0_{k-1,k-1}+q^1_{k-1,k-1}
+q_{k-1,k-1}}\\
&=&[\prod_{r=\iota}^{k-2}(1+w_{r,2k-r-1}^{(2k)}+w^{(2k)}_{2k-r-2,r+1})
^{\sum_{s=1}^r(q^0_{2k-2-2r+s,s}+q_{2k-2-2r+s,s}+q_{s,2k-2-2r+s})}\\
& &\times(1+w_{r,2k-r-1}^{(2k)}+w^{(2k)}_{2k-r-2,r+1})
^{-\sum_{s=1}^{r-1}(q_{2k-2-2r+s,s+1}+q_{s,2k-1-2r+s})}\\ &&\times
(1+w_{k-1,k}^{(2k)})^{\sum_{r=1}^{k-1}(q^0_{r,r}+q_{r,r})-\sum_{r=1}^{k-1}q_{r,r+1}}
\hspace{6.7cm}(15.4.87)\end{eqnarray*} by (15.4.68) and (15.4.70).
Set
$$w_{2k-2-2r+s,s}^{(2k-1)}=w_{2k-2-2r+s,s}^{(2k)}(1+w_{r,2k-r-1}^{(2k)}+w_{2k-r-2,r+1}^{(2k)}),
\eqno(15.4.88)$$
$$w_{s,2k-2-2r+s}^{(2k-1)}=w_{s,2k-2-2r+s}^{(2k)}(1+w_{r,2k-r-1}^{(2k)}+w_{2k-r-2,r+1}^{(2k)})
\eqno(15.4.89)$$ for $1\leq s\leq r\leq k-2$;
$$w_{2k-2-2r+s,s+1}^{(2k-1)}=\frac{w_{2k-2-2r+s,s+1}^{(2k)}}{1+w_{r,2k-r-1}^{(2k)}
+w_{2k-r-2,r+1}^{(2k)}}, \eqno(15.4.90)$$
$$w_{s,2k-1-2r+s}^{(2k-1)}=\frac{w_{s,2k-1-2r+s}^{(2k)}}{1+w_{r,2k-r-1}^{(2k)}+
w_{2k-r-2,r+1}^{(2k)}} \eqno(15.4.91)$$ for $1\leq s\leq r\leq k-2$;
$$w_{r,r}^{(2k-1)}=w_{r,r}^{(2k)}(1+w_{k-1,k}^{(2k)}),\qquad
w_{s,s+1}^{(2k-1)}=\frac{w_{s,s+1}^{(2k)}}{1+w_{k-1,k}^{(2k)}}\eqno(15.4.92)$$
for $r\in\ol{1,k-1},\;s\in\ol{1,k-2}$; and
$$w^{(2k-1)}_{r_1,r_2}=w^{(2k)}_{r_1,r_2}\;\;\mbox{for the other
pairs}\;\;(r_1,r_2)\;\;\mbox{such that}\;\;r_1+r_2\leq
2k-1.\eqno(15.4.93)$$ By induction, we have defined
$$\{w^{(k)}_{r_1,r_2}\mid
k\in\ol{1,2n-1};\;r_1\in\ol{1,n-1},\;r_2\in\ol{1,n},\;r_1+r_2\leq
k\}.\eqno(15.4.94)$$ For convenience, we denote
$$
w_{r_1,r_2}^{(2n)}=w_{r_1,r_2}\qquad\for\;\;r_1\in\ol{1,n-1},\;r_2\in\ol{2,n}.\eqno(15.4.95)$$
Based on (15.4.76), (15.4.80) and (15.4.87), we have
\begin{eqnarray*}&&\sum_{r_1=1}^{n-1}\;\sum_{r_2=2}^n\;\sum_{q_{r_1,r_2}=0}^{\infty}
[\prod_{1\leq l_1\leq l_2\leq n-1}
\left(\!\!\begin{array}{c}q^0_{l_2,l_1}+q_{l_2,l_1}^1+q_{l_1,l_2}+q_{l_2,l_1}
\\ q_{l_1,l_2+1},q_{l_2,l_1+1}\end{array}\!\!\right)\\ & &\times
w_{l_1,l_2+1}^{q_{l_1,l_2+1}}w_{l_2,l_1+1}^{q_{l_2,l_1+1}}]\left[\prod_{r=1}^{n-1}
\left(\!\!\begin{array}{c}q^0_{r,r}+q_{r,r}^1+q_{r,r}\\
q_{r,r+1}\end{array}\!\!\right)w_{r,r+1}^{q_{r,r+1}}\right]\\
&=&\left[\prod_{k=2}^n(1+w_{k-1,k}^{(2k)})^{\sum_{r=1}^{k-1}q_{r,r}^0}\right]\left[\prod_{k=2}
^{n-1}\prod_{r=\ell}^{k-1}(1+w_{r,2k-r}^{(2k+1)}+w_{2k-r-1,r+1}^{(2k+1)})^{\sum_{s=1}^r
q_{2k-1-2r+s,s}^0}\right]\\ &&\times
\left[\prod_{k=3}^{n-1}\prod_{r=\iota}^{k-2}(1+w_{r,2k-r-1}^{(2k)}+w_{2k-r-2,r+1}^{(2k)})
^{\sum_{s=1}^rq_{2k-2-2r+s,s}^0}\right]\\ &=&[\prod_{1\leq
r_1<r_2\leq n-1;r_2-r_1\;\mbox{\small is
odd}}\;(\prod_{k=(r_2-r_1+3)/2}^{n-1}(1+w^{(2k+1)}_{(2k+r_1-r_2-1)/2,(2k-r_1+r_2+1)/2}\\
&&+w^{(2k+1)}_{(2k-r_1+r_2-1)/2,(2k+r_1-r_2+1)/2})^{q^0_{r_2,r_1}})]
[\prod_{1\leq r_1<r_2\leq n-1;r_2-r_1\;\mbox{\small is even}}\\& &
(\prod_{k=(r_2-r_1+2)/2}^{n-1}(1+w^{(2k+2)}_{(2k+r_1-r_2)/2,(2k-r_1+r_2+2)/2}+w^{(2k+2)}
_{(2k-r_1+r_2)/2,(2k+r_1-r_2+2)/2})^{q^0_{r_2,r_1}})]\\ & &\times
\left[\prod_{r=1}^{n-1}\left(\prod_{k=r}^{n-1}(1+w_{k,k+1}^{(2k)})\right)^{q_{r,r}^0}\right].
\hspace{7.7cm}(15.4.96)\end{eqnarray*}

According to (15.4.66), (15.4.75) and (15.4.96), we define
$$ \xi_{n,i}^C=(-1)^{n+i+1}\frac{y_{n,i}y_{n,i+1}\cdots
y_{n,n}}{2^{1+\dlt_{n,i}}}\qquad\for\;\;1\leq i\leq
n,\eqno(15.4.97)$$
$$\xi_{i,i}^C=-\frac{y_{i,i}\cdots
y_{n,n}\prod_{k=i}^{n-1}(1+w_{k,k+1}^{(2k)})}{4}\eqno(15.4.98)$$ for
$1\leq i\leq n-1$,
\begin{eqnarray*}\xi_{r_2,r_1}^C&=&\frac{(-1)^{r_1+r_2+1}}{2}y_{r_2,r_1}y_{r_2+1,r_1+1}\cdots
y_{n-1,n+r_1-r_2-1} y_{n,n+r_1-r_2}\\ & &\times
y_{n,n+r_1-r_2+1}\cdots
y_{n,n}\prod_{k=(r_2-r_1+3)/2}^{n-1}(1+w^{(2k+1)}_{(2k+r_1-r_2-1)/2,(2k-r_1+r_2+1)/2}
\\ & &+w^{(2k+1)}_{(2k-r_1+r_2-1)/2,(2k+r_1-r_2+1)/2})\hspace{7.2cm}(15.4.99)\end{eqnarray*}
for $1\leq r_1<r_2\leq n-1$ such that $r_2-r_1$ is odd, and
\begin{eqnarray*}\xi_{r_2,r_1}^C&=&\frac{(-1)^{r_1+r_2+1}}{2}y_{r_2,r_1}y_{r_2+1,r_1+1}\cdots
y_{n-1,n+r_1-r_2-1}y_{n,n+r_1-r_2}\\ & &\times
y_{n,n+r_1-r_2+1}\cdots
y_{n,n}\prod_{k=(r_2-r_1+2)/2}^{n-1}(1+w^{(2k+2)}_{(2k+r_1-r_2)/2,(2k-r_1+r_2+2)/2}\\
& &+w^{(2k+2)}
_{(2k-r_1+r_2)/2,(2k+r_1-r_2+2)/2})\hspace{7.3cm}(15.4.100)\end{eqnarray*}
for $1\leq r_1<r_2\leq n-1$ such that $r_2-r_1$ is even.

 For
$\al\in\G_C$, we set
$$\al^c_i=\sum_{r=1}^i\al_{n+r,i}+\sum_{s=i}^n\al_{n+s,i}\qquad\for\;\;i\in\ol{1,n}\eqno(15.4.101)$$
and
$$\al^c=\sum_{1\leq r_1\leq r_2\leq n}\al_{n+r_2,r_1}.\eqno(15.4.102)$$
By (15.3.57), (15.3.62), (15.4.66), (15.4.67), (15.4.71)-(15.4.75)
and (15.4.96)-(15.4.102), we have
\begin{eqnarray*}& &E_C(z_1,z_2,...,z_n)\\ &=&\frac{x^\ast z_n^{n+\lmd_n+1}\prod_{i=1}^{n-1}
z_i^{2n-i+1/2}}{\left[\prod_{1\leq r_1<r_2\leq
n}(z_{r_1}-z_{r_2})\right]\left[\prod_{1\leq s_1 \leq s_2\leq
n}(z_{s_1}z_{s_2}-1)\right]}\sum_{\gm\in\G_C}\frac{1}{\gm!\la\lmd_n\ra_{_{\gm^c}}}
\\& &\times\prod_{r=1}^n
\la\gm_r^\ast-\gm_{r\ast}+\gm_r^c-1/2\ra_{_{\gm_r^\ast+\gm_r^c}}
\xi^\gm,\hspace{7cm}(15.4.103)\end{eqnarray*} where
$$\gm !=\left[\prod_{1\leq r_1<r_2\leq n}\gm_{r_2,r_1}!\right]\left[\prod_{1\leq s_1\leq s_2
\leq n}\gm_{n+s_2,s_1}!\right]\eqno(15.4.104)$$ and
$$\xi^\gm=
\left[\prod_{1\leq r_1<r_2\leq
n}(\xi^A_{r_2,r_1})^{\gm_{r_2,r_1}}\right] \left[\prod_{1\leq
s_1\leq s_2\leq
n}(\xi^C_{n+s_2,s_1})^{\gm_{n+s_2,s_1}}\right].\eqno(15.4.105)$$
Recall the notations in (15.3.18) and (15.3.57). We define our {\it
path hypergeometric function of type C}\index{path hypergeometric
function! of type C} by
$${\msr X}_C(\tau_1,...,\tau_n;\vt)\{z_{r_2,r_1},z_{n+s_2,s_1}\}=
\sum_{\al\in\G_C}\frac{\prod_{r=1}^n(\tau_r-\al_{r\ast})_{_{\al_r^\ast+\al^c_r}}}
{\al!(\vt)_{\al^c}}z^\al.\eqno(15.4.106)$$
 By (15.1.18), (15.1.27),(15.3.20), (15.4.102) and (15.4.106), we obtain the following main theorem in this section:
\psp

{\bf Theorem 15.4.1}. {\it The Etingof trace function in (15.2.47):
\begin{eqnarray*}E_C(z_1,z_2,...,z_n)&=&
\frac{x^\ast z_n^{n+\lmd_n+1}\prod_{i=1}^{n-1}z_i^{2n-i+1/2}}
{\left[\prod_{1\leq r_1<r_2\leq
n}(z_{r_1}-z_{r_2})\right]\left[\prod_{1\leq s_1\leq s_2\leq
n}(z_{s_1}z_{s_2}-1)\right]}\\ & &\times {\msr
X}_C(1/2,...,1/2;-\lmd_n)\{\xi^A_{r_2,r_1},\xi^C_{n+s_2,s_1}\}.
\hspace{3.2cm}(15.4.107)\end{eqnarray*} Moreover, the function in
(15.2.49):
$$\Psi_C=z_n^{\lmd_n+2-n}\big[\prod_{i=1}^{n-1}z_i^{-i-1/2}\big]{\msr
X}_C(1/2,...,1/2;-\lmd_n)\{\xi^A_{r_2,r_1},\xi^C_{n+s_2,s_1}\}.\eqno(15.4.108)$$}

\section{Properties of Path Hypergeometric Functions}

In this section, we find the differential properties, integral
representations and differential equations  for the hypergeometric
functions of type A in (15.3.60),  and the differential properties
and differential equations  for the hypergeometric functions of type
C in (15.4.106). Moreover, we define our hypergeometric functions of
type B and D analogously as those of type C.

Recall the differentiation property of the classical Gauss
hypergeometric function:
$$\frac{d}{dz}\:
_2F_1(\al,\be;\gm;z)=\frac{\al\be}{\gm}\:
_2F_1(\al+1,\be+1;\gm+1;z)\eqno(15.5.1)$$(e.g., cf. [X21]). For two
positive integers $k_1$ and $k_2$ such that $k_1<k_2$, a {\it path}
from $k_1$ to
 $k_2$ is a sequence $(m_0,m_1....,m_r)$ of positive integers such that
$$k_1=m_0<m_1<m_2<\cdots <m_{r-1}<m_r=k_2.\eqno(15.5.2)$$
One can imagine a path from $k_1$ to $k_2$ is a way of a super man
going from $k_1$th floor to $k_2$th floor through a stairway. Let
$${\msr P}_{k_1}^{k_2}=\mbox{the set of all paths from}\;k_1\;\mbox{to}\;k_2.\eqno(15.5.3)$$
The {\it path polynomial} from $k_1$ to $k_2$ is defined as
$$P_{[k_1,k_2]}=\sum_{(m_0,m_1,...,m_r)\in {\msr P}_{k_1}^{k_2}}(-1)^rz_{m_1,m_0}z_{m_2,m_1}
\cdots z_{m_{r-1},m_{r-2}}z_{m_r,m_{r-1}}.\eqno(15.5.4)$$ Moreover,
we set
$$P_{[k,k]}=1\qquad\for\;\;0<k\in\mbb{N}.\eqno(15.5.5)$$
For convenience, we simply denote
$${\msr X}_A={\msr X}_A(\tau_1,..,\tau_n;\vt)\{z_{j,k}\},\eqno(15.5.6)$$
$${\msr X}_A[i,j]={\msr X}_A(\tau_1,...,\tau_i+1,...,\tau_j-1,...\tau_n;\vt)\{z_{r_2,r_1}\}
\eqno(15.5.7)$$ obtained from ${\msr X}_A$ by changing $\tau_i$ to
$\tau_i+1$ and $\tau_j$ to $\tau_j-1$ for $1\leq i<j\leq n-1$ and
$${\msr X}_A[k,n]={\msr X}_A(\tau_1,..,\tau_k+1,...,\tau_n+1;\vt+1)\{z_{r_2,r_1}\}\eqno(15.5.8)$$
obtained from ${\msr X}_A$ by changing $\tau_i$ to $\tau_i+1$,
$\tau_n$ to $\tau_n+1$ and $\vt$ to $\vt+1$ for $k\in\ol{1,n-1}$. We
have the following natural generalization of (15.5.1): \psp

{\bf Theorem 15.5.1}. {\it For $1\leq r_1<r_2\leq n-1$ and $r\in
\ol{1,n-1}$, we have
$$\ptl_{z_{r_2,r_1}}({\msr X}_A)=\sum_{s=1}^{r_1}\tau_sP_{[s,r_1]}{\msr X}_A[s,r_2],\eqno(15.5.9)$$
$$\ptl_{z_{n,r}}({\msr X}_A)=\frac{\tau_n}{\vt}\sum_{s=1}^r\tau_sP_{[s,r]}{\msr X}_A[s,n].
\eqno(15.5.10)$$}

{\it Proof}. Note that
\begin{eqnarray*}\qquad\qquad\;\;& & \ptl_{z_{r_2,r_1}}({\msr X}_A)\\ &=& \sum_{\be\in\G_A}
\frac{\left[\prod_{s=1}^{n-1}(\tau_s-\be_{s\ast})_{_{\be_s^\ast}}\right]
(\tau_n)_{_{\be_{n\ast}}}}{\be!(\vt)_{_{\be_{n\ast}}}}\be_{r_2,r_1}z^{\be-\ves_{r_2,r_1}}
\\ &=&\sum_{\be\in\G_A}\frac{(\tau_{r_1}-\be_{r_1\ast})(\tau_{r_1}+1
 -\be_{r_1\ast})_{\be_{r_1}^\ast}(\tau_{r_2}-1-\be_{r_2\ast})_{\be_{r_2}^\ast}
 } {\be!(\vt)_{_{\be_{n\ast}}}}\\ & &\times \big[\prod_{s\neq
 r_1,r_2}(\tau_s-\be_{s\ast})_{_{\be_s^\ast}}\big](\tau_n)_{_{\be_{n\ast}}}z^\be\\ &=&(\tau_{r_1}-\sum_{s=1}^{r_1-1}z_{r_1,s}\ptl_{z_{r_1,s}})
 ({\msr X}_A[r_1,r_2]).\hspace{5.7cm}(15.5.11)\end{eqnarray*}
In particular,
$$ \ptl_{z_{r,1}}({\msr X}_A)=\tau_1{\msr X}_A[1,r]\qquad\for\;\;r\in\ol{1,n-1}.\eqno(15.5.12)$$
By (15.5.7), (15.5.11), (15.5.12) and  induction on $r_1$, we get
(15.5.9). Moreover,
\begin{eqnarray*}& & \ptl_{z_{n,r}}({\msr X}_A)\\&=& \sum_{\be\in\G_A}\frac{\left[
\prod_{s=1}^{n-1}(\tau_s-\be_{s\ast})_{_{\be_s^\ast}}\right](\tau_n)_{_{\be_{n\ast}}}}
{\be!(\vt)_{_{\be_{n\ast}}}}\be_{n,r}z^{\be-\ves_{n,r}}\\ &=&
\sum_{\be\in\G_A}\frac{\left[(\tau_r-\be_{r\ast})(\tau_r+1-\be_{r\ast})_{\be_r^\ast}
\prod_{s\neq
r}(\tau_s-\be_{s\ast})_{_{\be_s^\ast}}\right]\tau_n(\tau_n+1)_{_{\be_{n\ast}}}}
{\be!\vt(\vt+1)_{_{\be_{n\ast}}}}z^\be\\
&=&\frac{\tau_n}{\vt}(\tau_r-\sum_{s=1}^rz_{r,s}
\ptl_{z_{r,s}})({\msr
X}_A[r,n]),\hspace{7.8cm}(15.5.13)\end{eqnarray*} which leads to
(15.5.10) by (15.5.7), (15.5.8) and (15.5.9).$\qquad\Box$ \psp

For any $z\in \mbb{C}\setminus (-\mbb{N})$, the gamma function
$$\G(z)=\left[z e^{cz}\prod_{m=1}^{\infty}\left\{\left(1+\frac{z}{m}\right)e^{-z/m}\right\}
\right]^{-1},\eqno(15.5.14)$$ where $c$ is Euler's constant given by
$$c=\lim_{m\rta\infty}\left(\sum_{k=1}^m\frac{1}{k}-\ln m\right).\eqno(15.5.15)$$
When $\mbox{Re}\:z>0$, we have
$$\G(z)=\int_0^{\infty}t^{z-1}e^{-t}dt.\eqno(15.5.16)$$
Recall the integral representation of the classical Gauss
hypergeometric function:
$$
_2F_1(\al,\be;\gm;z)=\frac{\G(\gm)}{\G(\be)\G(\gm-\be)}\int_0^1t^{\be
-1}(1-t)^{\gm-\be-1}(1-zt)^{-\al}dt\eqno(15.5.17)$$ Now we have the
following analogous integral representation for our path
hypergeometric function of type A: \psp

{\bf Theorem 15.5.2}. {\it Suppose $\mbox{\it Re}\:\tau_n>0$ and
$\mbox{\it Re}\:(\vt-\tau_n)>0$.
 We have}
$${\msr X}_A=\frac{\G(\vt)}{\G(\vt-\tau_n)\G(\tau_n)}\int_0^1\left[\prod_{r=1}^{n-1}
(\sum_{s=r}^{n-1}P_{[r,s]}+tP_{[r,n]})^{-\tau_r}\right]t^{\tau_n-1}(1-t)^{\vt-\tau_n-1}dt.
\eqno(15.5.18)$$

{\it Proof}.  For $\kappa_1,\kappa_2\in\mbb{C}$ with
$\mbox{Re}\:\kappa_1,\mbox{Re}\:\kappa_1>0$, we have the following
Euler integral
$$\int_0^1t^{\kappa_1-1}(1-t)^{\kappa_2-1}dt=\frac{\G(\kappa_1)\G(\kappa_2)}{\G(\kappa_1+\kappa_2)}.\eqno(15.5.19)$$
Moreover,
$$(\kappa)_m=\frac{\G(\kappa+m)}{\G(\kappa)}\qquad\for\;\;m\in\mbb{N},\;\kappa\in\mbb{C}.\eqno(15.5.20)$$
Recall
 $${\msr X}_A=\sum_{\be\in\G_A}\frac{\left[\prod_{s=1}^{n-1}(\tau_s-\be_{s\ast})_{_{\be_s^\ast}}
 \right](\tau_n)_{_{\be_{n\ast}}}}{\be!(\vt)_{_{\be_{n\ast}}}}z^\be.\eqno(15.5.21)$$
We have
\begin{eqnarray*}\qquad & &\frac{(\tau_n)_{_{\be_{n\ast}}}}{(\vt)_{_{\be_{n\ast}}}}=
\frac{\G(\vt)}{\G(\vt-\tau_n)\G(\tau_n)}\frac{\G(\tau_n+\be_{n\ast})\G(\vt-\tau_n)}{\G(\vt
+\be_{n\ast})}\\
&=&\frac{\G(\vt)}{\G(\vt-\tau_n)\G(\tau_n)}\int_0^1t^{\tau_n+\be_{n\ast}-1}
(1-t)^{\vt-\tau_n-1}dt\hspace{4.8cm}(15.5.22)\end{eqnarray*} by
(15.5.19) and (15.5.20).

Denote
$$\G_{(s)}=\sum_{1\leq r_1<r_2\leq s}\mbb{N}\ves_{r_2,r_1},\;\;\be_r^s=\sum_{p=r+1}^s
\be_{p,r}\;\;\for\;\;\be\in\G_A,\;1\leq r<s\leq n-1.\eqno(15.5.23)$$
Note that
\begin{eqnarray*}& &\sum_{\be\in\G_A}\frac{\prod_{s=1}^{n-1}(\tau_s-\be_{s\ast})_{_{\be_s^\ast}}}
{\be!}z^\be
\\&=&\sum_{\be\in\G_{(n-1)}}\frac{\prod_{s=1}^{n-2}(\tau_s-\be_{s\ast})
_{_{\be^{n-1}_s}}}{\be!}z^\be\sum_{\be_{n,1},...,\be_{n,n-1}=0}^n\frac{(\tau_{n-1}
-\be_{(n-1)\ast})_{\be_{n,n-1}}z_{n,n-1}^{\be_{n,n-1}}}{\be_{n,n-1}!}\\
& &\times\prod_{r=1}^{n-2}
\frac{(\tau_r-\be_{r\ast}+\be_r^{n-1})_{\be_{n,r}}z_{n,r}^{\be_{n,r}}}{\be_{n,r}!}
\hspace{12cm}\end{eqnarray*}
 \begin{eqnarray*} &=&
\sum_{\be\in\G_{(n-1)}}\frac{\prod_{s=1}^{n-2}(\tau_s-\be_{s\ast})_{_{\be_s^{n-1}}}}
{\be!}z^\be(1-z_{n,n-1})^{\be_{(n-1)\ast}-\tau_{n-1}}\\&
&\times\prod_{r=1}^{n-2}(1-z_{n,r})^{\be_{r\ast}
-\be_r^{n-1}-\tau_r}\\
&=& \left[\prod_{s=1}^{n-1}(1-z_{n,s})^{-\tau_s}\right]
\sum_{\be\in\G_{(n-1)}}\frac{\prod_{s=1}^{n-2}(\tau_s-\be_{s\ast})_{_{\be_s^{n-1}}}}
{\be!}\\
& &\times\prod_{1\leq r_1<r_2\leq
n-1}\left[\frac{(1-z_{n.r_2})z_{r_2,r_1}}{1-z_{n,r_1}}\right]^{\be_{r_2,r_1}}.
\hspace{6.7cm}(15.5.24)\end{eqnarray*} Observe that
$$1-\frac{(1-z_{n,n-1})z_{n-1,r}}{1-z_{n,r}}=\frac{1-z_{n-1,r}-z_{n,r}+z_{n,n-1}z_{n-1,r}}{1
-z_{n,r}}\eqno(15.5.25)$$for $r\in\ol{1,n-1}$. By (15.5.24) and
(15.5.25), we have
\begin{eqnarray*}& &\sum_{\be\in\G_A}\frac{\prod_{s=1}^{n-1}(\tau_s-\be_{s\ast})_{_{\be_s^\ast}}}
{\be!}z^\be=(1-z_{n,n-1})^{-\tau_{n-1}}(1-z_{n-1,n-2}+P_{[n-2,n]})^{-\tau_{n-2}}\\
& &\times
\left[\prod_{s=1}^{n-2}(1-z_{n-1,s}-z_{n,s}+z_{n,n-1}z_{n-1,s})^{-\tau_s}\right]
\sum_{\be\in\G_{(n-2)}}\frac{\prod_{s=1}^{n-3}(\tau_s-\be_{s\ast})_{_{\be^{n-2}_s}}}
{\be !}\\ & &\times \prod_{1\leq r_1<r_2\leq
n-2}\left[\frac{(1-z_{n-1,r_2}-z_{n,r_2}
+z_{n,n-1}z_{n-1,r_2})z_{r_2,r_1}}{1-z_{n-1,r_1}-z_{n,r_1}+z_{n,n-1}z_{n-1,r_1}}
\right]^{\be_{r_2,r_1}}.\hspace{2.7cm}(15.5.26)\end{eqnarray*} Now
\begin{eqnarray*} & &1-\frac{(1-z_{n-1,n-2}-z_{n,r-2}+z_{n,n-1}z_{n-1,n-2})z_{n-2,r}}
{1-z_{n-1,r}-z_{n,r}+z_{n,n-1}z_{n-1,r}}\\
&=&\frac{1}{1-z_{n-1,r}-z_{n,r}+z_{n,n-1}z_{n-1,r}}
[1-z_{n-2,r}-z_{n-1,r}+z_{n-1,n-2}z_{n-2,r}\\ &
&-z_{n,r}+z_{n,n-1}z_{n-1,r}+z_{n,r-2}z_{n-2,r}
-z_{n,n-1}z_{n-1,n-2}z_{n-2,r}].\hspace{3.2cm}(15.5.27)\end{eqnarray*}
By induction, we can prove that
$$\sum_{\be\in\G_A}\frac{\prod_{s=1}^{n-1}(\tau_s-\be_{s\ast})_{_{\be_s^\ast}}}{\be!}z^\be
=\prod_{r=1}^{n-1}(\sum_{s=r}^nP_{[r,s]})^{-\tau_r}.\eqno(15.5.28)$$
Thus
$$\sum_{\be\in\G_A}\frac{\prod_{s=1}^{n-1}(\tau_s-\be_{s\ast})_{_{\be_s^\ast}}}{\be!}
z^\be
t^{\be_{n\ast}}=\prod_{r=1}^{n-1}(\sum_{s=r}^{n-1}P_{[r,s]}+tP_{[r,n]})^{-\tau_r}.
\eqno(15.5.29)$$ Hence we obtain (15.5.18) by (15.5.22) and
(15.5.29).$\qquad\Box$ \psp

{\bf Remark 15.5.3}. According to [M],
$$\left(\begin{array}{ccccc}1&0&0&\cdots&0\\ z_{2,1}&
1&0&\cdots&0\\ z_{3,1}&z_{3,2}&1&\ddots&\vdots\\
\vdots&\vdots&\ddots&\ddots&0\\ z_{n,1}&z_{n,2}&\cdots&z_{n,n-1}&
1\end{array}\right)^{-1}=\left(\begin{array}{ccccc}1&0&0&\cdots&0\\
P_{[1,2]}&
1&0&\cdots&0\\ P_{[1,3]}&P_{[2,3]}&1&\ddots&\vdots\\
\vdots&\vdots&\ddots&\ddots&0\\
P_{[1,n]}&P_{[2,n]}&\cdots&P_{[n-1,n]}&
1\end{array}\right).\eqno(15.5.30)$$ \pse

 Recall the classical Gauss
hypergeometric equation
$$z(1-z){y'}'+[\gm-(\al+\be+1)z]y'-\al\be y=0.\eqno(15.5.31)$$
Note
$$D=z\frac{d}{dz}\lra
D^2=z^2\frac{d^2}{dz^2}+z\frac{d}{dz}.\eqno(15.5.32)$$ Then
(15.5.31) can be rewritten as
$$(\gm+D)\frac{d}{dz}(y)=(\al+D)(\be+D)(y).\eqno(15.5.33)$$
Denote
$${\msr D}_{i\ast}^A=\sum_{r=1}^{i-1}z_{i,r}\ptl_{z_{i,r}},\qquad {\msr D}_i^{A\ast}=
\sum_{s=i+1}^nz_{s,i}\ptl_{z_{s,i}}\qquad\for\;\;i\in\ol{1,n}.\eqno(15.5.34)$$
Then we have the following analogue of (15.5.33): \psp

{\bf Theorem 15.5.4}. {\it We have:
$$(\tau_{r_2}-1-{\msr D}_{r_2\ast}^A+{\msr D}_{r_2}^{A\ast})\ptl_{z_{r_2,r_1}}({\msr X}_A)
=(\tau_{r_2}-1-{\msr D}_{r_2\ast}^A)(\tau_{r_1}-{\msr
D}_{r_1\ast}^A+{\msr D}_{r_1}^{A\ast}) ({\msr X}_A)\eqno(15.5.35)$$
for $1\leq r_1<r_2\leq n-1$ and
$$(\vt+{\msr D}_{n\ast}^A)\ptl_{z_{n,r}}({\msr X}_A)=(\tau_n+{\msr D}_{n\ast}^A)
(\tau_r-{\msr D}_{r\ast}^A+{\msr D}_r^{A\ast})({\msr
X}_A)\eqno(15.5.36)$$ for $r\in\ol{1,n-1}$.} \pse

{\it Proof}. Let
$$a_\be=\frac{\left[\prod_{s=1}^{n-1}(\tau_s-\be_{s\ast})_{_{\be_s^\ast}}\right]
(\tau_n)_{_{\be_{n\ast}}}}{\be!(\vt)_{_{\be_{n\ast}}}}\qquad\for\;\;\be\in\G_A.\eqno(15.5.37)$$
Then
$$\frac{a_{\be+\ves_{r_2,r_1}}}{a_\be}=\frac{(\tau_{r_2}-\be_{r_2\ast}-1)(\tau_{r_1}-
\be_{r_1\ast}+\be_{r_1}^\ast)}{(\tau_{r_2}-1-\be_{r_2\ast}+\be_{r_2}^\ast)(\be_{r_2,r_1}+1)},
\eqno(15.5.38)$$  or equivalently,
\begin{eqnarray*}\qquad\qquad &&(\tau_{r_2}-1-\be_{r_2\ast}+\be_{r_2}^\ast)(\be_{r_2,r_1}+1)a_{\be+\ves_{r_2,r_1}}
\\&=&(\tau_{r_2}
-\be_{r_2\ast}-1)(\tau_{r_1}-\be_{r_1\ast}+\be_{r_1}^\ast)a_\be.\hspace{5.5cm}(15.5.39)\end{eqnarray*}
Thus
\begin{eqnarray*}\qquad\qquad &&\sum_{\be\in\G_A}(\tau_{r_2}-1-\be_{r_2\ast}+\be_{r_2}^\ast)(\be_{r_2,r_1}+1)
a_{\be+\ves_{r_2,r_1}}z^\be\\&=&\sum_{\be\in\G_A}(\tau_{r_2}-\be_{r_2\ast}-1)(\tau_{r_1}
-\be_{r_1\ast}+\be_{r_1}^\ast)a_\be
z^\be,\hspace{4.2cm}(15.5.40)\end{eqnarray*} which is equivalent to
(15.5.35). Similarly, (15.5.36) follows from
$$\frac{a_{\be+\ves_{n,r}}}{a_\be}=\frac{(\tau_n+\be_{n\ast})(\tau_r-\be_{r\ast}+\be_r^\ast)}
{(\vt+\be_{n\ast})(\be_{n,r}+1)}.\qquad\Box\eqno(15.5.41)$$ \psp

Next we study our hypergeometric functions of type C. For
convenience, we simply denote
$${\msr X}_C={\msr X}_C(\tau_1,...,\tau_n;\vt)\{z_{r_2,r_1},z_{n+s_2,s_1}\}\eqno(15.5.42)$$
(cf. (15.4.106)). First we have
\begin{eqnarray*}& &\ptl_{z_{r_2,r_1}}({\msr X}_C)\\ &=&\sum_{\al\in\G_C}\frac{\prod_{r=1}^n
(\tau_r-\al_{r\ast})_{_{\al_r^\ast+\al^c_r}}}{\al!(\vt)_{\al^c}}\al_{r_2,r_1}z^{\al
-\ves_{r_2,r_1}}\\
 &=& \sum_{\al\in\G_C}\frac{\prod_{r\neq
r_1,r_2}(\tau_r-\al_{r\ast})_{_{\al_r^\ast+\al^c_r}}}{\al!(\vt)_{\al^c}}(\tau_{r_1}
-\al_{r_1\ast})(\tau_{r_1}+1-\al_{r_1\ast})_{_{\al_{r_1}^\ast+\al^c_{r_1}}}\\
&
&\times(\tau_{r_2}-1-\al_{r_2\ast})_{_{\al_{r_2}^\ast+\al^c_{r_2}}}z^\al.
\hspace{8.7cm}(15.5.43)\end{eqnarray*} or $1\leq r_1<r_2\leq n$.
Moreover,
\begin{eqnarray*}& &\ptl_{z_{n+r_2,r_1}}({\msr X}_C)\\ &=&\sum_{\al\in\G_C}
\frac{\prod_{r=1}^{n-1}(\tau_r-\al_{r\ast})_{_{\al_r^\ast+\al^c_r}}}{\al!(\vt)_{\al^c}}
\al_{n+r_2,r_1}z^{\al-\ves_{n+r_2,r_1}}\\ &=&
\frac{1}{\vt}\sum_{\al\in\G_C} \frac{\prod_{r\neq
r_1,r_2}(\tau_r-\al_{r\ast})_{_{\al_r^\ast+\al^c_r}}}
{\al!(\vt+1)_{\al^c}}(\tau_{r_1}-\al_{r_1\ast})(\tau_{r_2}-\al_{r_2\ast})\\
& &\times(\tau_{r_1}
+1-\al_{r_1\ast})_{_{\al_{r_1}^\ast+\al^c_{r_1}}} (\tau_{r_2}+1-
\al_{r_2\ast})_{_{\al_{r_2}^\ast+\al^c_{r_2}}}z^{\al}\hspace{5cm}(15.5.44)\end{eqnarray*}
for $1\leq r_1<r_2\leq n$ and
\begin{eqnarray*}\ptl_{z_{n+s,s}}({\msr X}_C)&=&\sum_{\al\in\G_C}\frac{\prod_{r=1}^n(\tau_r
-\al_{r\ast})_{_{\al_r^\ast+\al^c_r}}}{\al!(\vt)_{\al^c}}\al_{n+s,s}z^{\al-\ves_{n+s,s}}\\
&=&\frac{1}{\vt}\sum_{\al\in\G_C}\frac{\prod_{r\neq
s}(\tau_r-\al_{r\ast})_{_{\al_r^\ast
+\al^c_r}}}{\al!(\vt+1)_{\al^c}}(\tau_{s}-\al_{s\ast})(\tau_{s}+1-\al_{s\ast})\\
& &\times
(\tau_s+2-\al_{s\ast})_{_{\al_s^\ast+\al^c_{s}}}z^{\al}\hspace{7.5cm}(15.5.45)\end{eqnarray*}
for $s\in\ol{1,n}$.

Expressions (15.5.43)-(15.5.45) motivate us to define
$${\msr X}_C[i_1,i_2]={\msr X}_C(\tau_1,...,\tau_{i_1}+1,...,\tau_{i_2}-1,...\tau_n;\vt)
\{z_{r_2,r_1},z_{n+s_2,s_1}\}\eqno(15.5.46)$$ obtained from ${\msr
X}_C$  by changing $\tau_{i_1}$ to $\tau_{i_1}+1$ and $\tau_{i_2}$
to $\tau_{i_2}-1$ for $1\leq i_1<i_2\leq n$,
$${\msr X}_C[j_1,j_2;1]={\msr X}_C(\tau_1,...,\tau_{j_1}+1,...,\tau_{j_2}+1,...\tau_n;\vt+1)
\{z_{r_2,r_1},z_{n+s_2,s_1}\}\eqno(15.5.47)$$ obtained from ${\msr
X}_C$  by changing $\tau_{j_1}$ to $\tau_{j_1}+1$, $\tau_{j_2}$ to
$\tau_{j_2}+1$ and $\vt$ to $\vt+1$ for $1\leq j_1<j_2\leq n$ and
$${\msr X}_C[k(2)]={\msr X}_C(\tau_1,...,\tau_k+2,...,\tau_n;\vt+1)\{z_{r_2,r_1},z_{n+s_2,s_1}\}
\eqno(15.5.48)$$ obtained from ${\msr X}_C$  by changing $\tau_k$ to
$\tau_k+2$ and $\vt$ to $\vt+1$ for $k\in\ol{1,n}$. Now (15.5.43)
can be written as
$$\ptl_{z_{r_2,r_1}}({\msr X}_C)=(\tau_{r_1}-\sum_{s=1}^{r_1-1}z_{r_1,s}\ptl_{z_{r_1,s}})
({\msr X}_C[i_1,i_2])\eqno(15.5.49)$$ and (15.5.44) becomes
$$\ptl_{z_{n+r_2,r_1}}({\msr X}_C)=\frac{1}{\vt}(\tau_{r_1}-\sum_{s_1=1}^{r_1-1}z_{r_1,s_1}
\ptl_{z_{r_1,s_1}})(\tau_{r_2}-\sum_{s_2=1}^{r_2-1}z_{r_2,s_2}\ptl_{z_{r_2,s_2}})
({\msr X}_C[i_1,i_2;1])\eqno(15.5.50)$$ for $1\leq r_1<r_2\leq n.$
Moreover, (15.5.45) is equivalent to
$$\ptl_{z_{n+s,s}}({\msr X}_C)=\frac{1}{\vt}(\tau_s-\sum_{r=1}^{s-1}z_{s,r}\ptl_{z_{s,r}})
(\tau_s+1-\sum_{r=1}^{s-1}z_{s,r}\ptl_{z_{s,r}})({\msr
X}_C[s(2)])\eqno(15.5.51)$$ for $s\in\ol{1,n}$. By
(15.5.49)-(15.5.51) and induction, we obtain: \psp

{\bf Theorem 15.5.5}. {\it The following equations hold for ${\msr
X}_C$:
$$\ptl_{z_{r_2,r_1}}({\msr X}_C)=\sum_{s=1}^{r_1}\tau_sP_{[s,r_1]}{\msr X}_C[s,r_2],
\eqno(15.5.52)$$
\begin{eqnarray*}& &\ptl_{z_{n+r_2,r_1}}({\msr X}_C)\\ &=&\frac{1}{\vt}[\sum_{i=1}^{r_1}\tau_i^2
P_{[i,r_1]}P_{[i,r_2]}{\msr
X}_C[i(2)]+\sum_{s_1=1}^{r_1}\sum_{s_2=r_1+1}^{r_2}\tau_{s_1}
\tau_{s_2}P_{[s_1,r_1]}P_{[s_2,r_2]}{\msr X}_C[s_1,s_2;1]\\ &
&+\sum_{1\leq s_1<s_2\leq r_1}
\tau_{s_1}\tau_{s_2}(P_{[s_1,r_1]}P_{[s_2,r_2]}+P_{[s_2,r_1]}P_{[s_1,r_2]})
{\msr X}_C[s_1,s_2;1]]\hspace{2.8cm}(15.5.53)\end{eqnarray*} for
$1\leq r_1<r_2\leq n$ and
\begin{eqnarray*}\ptl_{z_{n+s,s}}({\msr X}_C)&=&\frac{1}{\vt}[\sum_{i=1}^s\tau_i^2P_{[i,s]}^2
{\msr X}_C[i(2)]+\tau_s{\msr X}_C[s(2)]+
\sum_{i=1}^{s-1}\tau_iP_{[i,s]}{\msr X}_C[i,s;1]\\ & & +2\sum_{1\leq
r_1<r_2\leq s}\tau_{r_1}\tau_{r_2}P_{[r_1,s]}P_{[r_2,s]}{\msr
X}_C[r_1,r_2;1]] \hspace{4cm}(15.5.54)\end{eqnarray*} for
$s\in\ol{1,n}$.} \psp

Up to this stage, we have not found a nice integral representation
for ${\msr X}_C$. In fact, there is no simple integral formula of
Euler type with an elementary integrand for Lauricella third
multiple hypergeometric function (e.g., cf. [Eh]). It might also be
the case for our hypergeometric function ${\msr X}_C$.

By (15.4.106), we set
$$c_\al=\frac{\prod_{r=1}^n(\tau_r-\al_{r\ast})_{_{\al_r^\ast+\al^c_r}}}{\al!(\vt)_{\al^c}}
z^\al\qquad\for\;\;\al\in\G_C.\eqno(15.5.55)$$ Note that
$$\frac{c_{\al+\ves_{r_2,r_1}}}{c_\al}=\frac{(\tau_{r_1}-\al_{_{r_1\ast}}+\al_{_{r_1}^\ast}
+\al^c_{r_1})(\tau_{r_2}-1-\al_{_{r_2\ast}})}{(\al_{r_2,r_1}+1)(\tau_{r_2}-1-\al_{_{r_2\ast}}
+\al_{_{r_2}^\ast}+\al^c_{r_2})},\eqno(15.5.56)$$
$$\frac{c_{\al+\ves_{n+r_2,r_1}}}{c_\al}=\frac{(\tau_{r_1}-\al_{_{r_1\ast}}+\al_{_{r_1}^\ast}
+\al^c_{r_1})(\tau_{r_2}-\al_{_{r_2\ast}}+\al_{_{r_2}^\ast}+\al^c_{r_2})}{(\al_{n+r_2,r_1}+1)
(\vt+\al^c)}\eqno(15.5.57)$$ for $1\leq r_1<r_2\leq n$ and
$$\frac{c_{\al+\ves_{n+s,s}}}{c_\al}=\frac{(\tau_s-\al_{_{s\ast}}+\al_{_s^\ast}+\al^c_{s})
(\tau_s+1-\al_{_{s\ast}}+\al_{_s^\ast}+\al^c_{s})}{(\al_{n+s,s}+1)(\vt+\al^c)}\eqno(15.5.58)$$
for $s\in\ol{1,n}$. Let
$${\msr D}^C_r=\sum_{i=1}^rz_{n+r,i}\ptl_{z_{n+r,i}}+\sum_{s=r}^nz_{n+s,r}\ptl_{z_{n+s,r}}
\qquad\for\;\;r\in\ol{1,n}\eqno(15.5.59)$$ and
$${\msr D}^C=\sum_{1\leq r_1\leq r_2\leq n}z_{n+r_2,r_1}\ptl_{z_{n+r_2,r_1}}.\eqno(15.5.60)$$
By the proof of Theorem 15.5.4, we have: \psp

{\bf Theorem 15.5.6}. {\it The function ${\msr X}_C$ satisfies:
\begin{eqnarray*}& &(\tau_{r_2}-1-{\msr D}_{r_2\ast}^A+{\msr D}_{r_2}^{A\ast}+{\msr
D}^C_{r_2})\ptl_{z_{r_2,r_1}} ({\msr X}_C)\\ &=&(\tau_{r_2}-1-{\msr
D}_{r_2\ast}^A)(\tau_{r_1}-{\msr D}_{r_1\ast}^A +{\msr
D}_{r_1}^{A\ast}+{\msr D}_{r_1}^C)({\msr
X}_C),\hspace{5.1cm}(15.5.61)\end{eqnarray*}
\begin{eqnarray*}& &(\vt+{\msr D}^C)\ptl_{z_{n+r_2,r_1}}({\msr X}_C)\\ &=&(\tau_{r_2}-{\msr
D}_{r_2\ast}^A +{\msr D}_{r_2}^{A\ast}+{\msr
D}^C_{r_2})(\tau_{r_1}-{\msr D}_{r_1\ast}^A+{\msr D}_{r_1}^{A\ast}
+{\msr D}_{r_1}^C)({\msr X}_C)\hspace{3.6cm}(15.5.62)\end{eqnarray*}
for $1\leq r_1<r_2\leq n$ and
\begin{eqnarray*}& &(\vt+{\msr D}^C)\ptl_{z_{n+s,s}}({\msr
X}_C)\\ &=&(\tau_s-{\msr D}_{s\ast}^A+{\msr D}_s^{A\ast} +{\msr
D}^C_s)(\tau_s+1-{\msr D}_{s\ast}^A+{\msr D}_s^{A\ast}+{\msr
D}_{s}^C)({\msr X}_C)\hspace{3.6cm}(15.5.63)\end{eqnarray*}for
$s\in\ol{1,n}$.} \psp

Set
$$\G_D=\sum_{1\leq r_1<r_2\leq n}(\mbb{N}\ves_{r_2,r_1}+\mbb{N}\ves_{n+r_2,r_1})\subset \G_C.
\eqno(15.5.64)$$ Moreover, we let
$$\be_r^D=\sum_{i=1}^{r-1}\be_{n+r,i}+\sum_{s=r+1}^n\be_{n+s,r},\;\;\be^D=\sum_{1\leq r_1<r_2
\leq n}\be_{n+r_2,r_1}\qquad\for\;\;\be\in\G_D\eqno(15.5.65)$$ and
$$\be_r^B=\sum_{i=1}^r\be_{n+r,i}+\sum_{s=r+1}^n\be_{n+s,r},
\;\;\be^B=\be^c=\sum_{1\leq r_1\leq r_2\leq
n}\be_{n+r_2,r_1}\qquad\for\;\;\be\in\G_C. \eqno(15.5.66)$$ We
define the following hypergeometric functions:
$${\msr X}_D(\tau_1,...,\tau_n;\vt)\{z_{r_2,r_1},z_{n+r_2,r_1}\}=\sum_{\al\in\G_D}
\frac{\prod_{r=1}^n(\tau_r-\al_{r\ast})_{_{\al_r^\ast+\al^D_r}}}{\al!(\vt)_{\al^D}}z^\al.
\eqno(15.5.67)$$
$${\msr X}_B(\tau_1,...,\tau_n;\vt)\{z_{r_2,r_1},z_{n+s_2,s_1}\}=\sum_{\al\in\G_C}
\frac{\prod_{r=1}^n(\tau_r-\al_{r\ast})_{_{\al_r^\ast+\al^B_r}}}{\al!(\vt)_{\al^B}}z^\al.
\eqno(15.5.68)$$ Furthermore, we define ${\msr X}_B,{\msr
X}_B[i_1,i_2],{\msr X}_B[j_1,j_2;1],{\msr X}_B[k(2)]$ and ${\msr
X}_D,{\msr X}_D[i_1,i_2],{\msr X}_D[j_1,j_2;1],$ ${\msr X}_D[k(2)]$
as those of type C in (15.5.42) and (15.5.46)-(15.5.48). We let
$${\msr X}_B[k]={\msr X}_N(\tau_1,...,\tau_k+1,...,\tau_n;\vt+1)\{z_{r_2,r_1},z_{n+s_2,s_1}\}
\eqno(15.5.69)$$ obtained from ${\msr X}_B$  by changing $\tau_k$ to
$\tau_k+1$ and $\vt$ to $\vt+1$ for $k\in\ol{1,n}$. \psp

{\bf Theorem 15.5.7}. {\it The following equations hold for ${\msr
X}_B$ and ${\msr X}_D$:
$$\ptl_{z_{r_2,r_1}}({\msr X}_B)=\sum_{s=1}^{r_1}\tau_sP_{[s,r_1]}{\msr X}_B[s,r_2],
\eqno(15.5.70)$$
\begin{eqnarray*}& &\ptl_{z_{n+r_2,r_1}}({\msr X}_B)\\&=&\frac{1}{\vt}[\sum_{i=1}^{r_1}
\tau_i^2P_{[i,r_1]}P_{[i,r_2]}{\msr
X}_B[i(2)]+\sum_{s_1=1}^{r_1}\sum_{s_2=r_1+1}^{r_2}
\tau_{s_1}\tau_{s_2}P_{[s_1,r_1]}P_{[s_2,r_2]}{\msr
X}_B[s_1,s_2;1]\\ & &+\sum_{1\leq s_1<s_2 \leq
r_1}\tau_{s_1}\tau_{s_2}(P_{[s_1,r_1]}P_{[s_2,r_2]}+P_{[s_2,r_1]}P_{[s_1,r_2]})
{\msr X}_B[s_1,s_2;1]],\hspace{2.6cm}(15.5.71)\end{eqnarray*}
$$\ptl_{z_{r_2,r_1}}({\msr X}_D)=\sum_{s=1}^{r_1}\tau_sP_{[s,r_1]}{\msr X}_D[s,r_2],\eqno(15.5.72)$$
\begin{eqnarray*}&& \ptl_{z_{n+r_2,r_1}}({\msr X}_D)\\&=&\frac{1}{\vt}[\sum_{i=1}^{r_1}
\tau_i^2P_{[i,r_1]}P_{[i,r_2]}{\msr
X}_D[i(2)]+\sum_{s_1=1}^{r_1}\sum_{s_2=r_1+1}^{r_2}
\tau_{s_1}\tau_{s_2}P_{[s_1,r_1]}P_{[s_2,r_2]}{\msr
X}_D[s_1,s_2;1]\\ & &+ \sum_{1\leq s_1<s_2\leq
r_1}\tau_{s_1}\tau_{s_2}(P_{[s_1,r_1]}P_{[s_2,r_2]}+
P_{[s_2,r_1]}P_{[s_1,r_2]}){\msr
X}_D[s_1,s_2;1]],\hspace{2.4cm}(15.5.73)\end{eqnarray*} for $1\leq
r_1<r_2\leq n$ and
$$\ptl_{z_{n+s,s}}({\msr X}_B)=\sum_{r=1}^s\tau_rP_{[s,r]}\msr X_B[s]\eqno(15.5.74)$$
for $s\in\ol{1,n}$.}\psp

 Let
$${\msr D}^B_r=\sum_{i=1}^rz_{n+r,i}\ptl_{z_{n+r,i}}+\sum_{s=r+1}^nz_{n+s,r}\ptl_{z_{n+s,r}},
\eqno(15.5.75)$$
$${\msr D}^D_r=\sum_{i=1}^{r-1}z_{n+r,i}\ptl_{z_{n+r,i}}+\sum_{s=r+1}^nz_{n+s,r}
\ptl_{z_{n+s,r}}\eqno(15.5.76)$$ for $r\in\ol{1,n}$ and
$${\msr D}^B={\msr D}^C,\;\;{\msr D}^D=\sum_{1\leq r_1< r_2\leq n}z_{n+r_2,r_1}
\ptl_{z_{n+r_2,r_1}}.\eqno(15.5.77)$$ As Theorem 15.5.4, we have:
\psp

{\bf Theorem 15.5.8}. {\it The functions ${\msr X}_B$ and ${\msr
X}_D$ satisfy: \begin{eqnarray*}\qquad\qquad&&(\tau_{r_2}-1-{\msr
D}_{r_2\ast}^A+{\msr D}_{r_2}^{A\ast}+{\msr
D}^B_{r_2})\ptl_{z_{r_2,r_1}} ({\msr X}_B)\\&=&(\tau_{r_2}-1-{\msr
D}_{r_2\ast}^A)(\tau_{r_1}-{\msr D}_{r_1\ast}^A+ {\msr
D}_{r_1}^{A\ast}+{\msr D}_{r_1}^B)({\msr
X}_B),\hspace{3.1cm}(15.5.78)\end{eqnarray*}
\begin{eqnarray*}\qquad\qquad&&(\tau_{r_2}-1-{\msr D}_{r_2\ast}^A+{\msr D}_{r_2}^{A\ast}+{\msr
D}^D_{r_2})\ptl_{z_{r_2,r_1}} ({\msr X}_D)\\&=&(\tau_{r_2}-1-{\msr
D}_{r_2\ast}^A)(\tau_{r_1}-{\msr D}_{r_1\ast}^A+ {\msr
D}_{r_1}^{A\ast}+{\msr D}_{r_1}^D)({\msr
X}_D),\hspace{3.1cm}(15.5.79)\end{eqnarray*}
\begin{eqnarray*}\qquad&&(\vt+{\msr D}^B)\ptl_{z_{n+r_2,r_1}}({\msr X}_B)\\ &=&(\tau_{r_2}-{\msr
D}_{r_2\ast}^A +{\msr D}_{r_2}^{A\ast}+{\msr
D}^B_{r_2})(\tau_{r_1}-{\msr D}_{r_1\ast}^A+{\msr D}_{r_1}^{A\ast}
+{\msr D}_{r_1}^B)({\msr
X}_B),\hspace{2.3cm}(15.5.80)\end{eqnarray*}
\begin{eqnarray*}\qquad&&(\vt+{\msr D}^D)\ptl_{z_{n+r_2,r_1}}({\msr
X}_D)\\ &=&(\tau_{r_2}-{\msr D}_{r_2\ast}^A +{\msr
D}_{r_2}^{A\ast}+{\msr D}^D_{r_2})(\tau_{r_1}-{\msr
D}_{r_1\ast}^A+{\msr D}_{r_1}^{A\ast} +{\msr D}_{r_1}^D)({\msr
X}_D)\hspace{2.5cm}(15.5.81)\end{eqnarray*} for $1\leq r_1<r_2\leq
n$ and
$$(\vt+{\msr D}^B)\ptl_{z_{n+s,s}}({\msr X}_B)=(\tau_s-{\msr D}_{s\ast}^A+{\msr D}_s^{A\ast}
+{\msr D}^B_s)({\msr X}_B)\eqno(15.5.82)$$ for $s\in\ol{1,n}$.}\psp

Heckman and Opdam [Hg1-Hg3, HO, O1-O6,BO] introduced hypergeometric
equations related to root systems and analogous to (15.1.4). They
proved the existence of
 solutions (hypergeometric functions) of their equations.
In our case, the interesting functions are the path hypergeometric
functions like ${\cal X}_A$
 in (15.3.60) and the functions like $\Psi_A$ in (15.3.63) are not interesting from pure function
  point of view. Gel'fand and Graev studied analogues of classical hypergeometric functions
  (so called GG-functions) by generalizing the differential property of the classical
  hypergeometric functions (e.g. cf. [GG]).

\addcontentsline{toc}{chapter}{\numberline{}Index}

\printindex


\begin{thebibliography}{}
\addcontentsline{toc}{chapter}{\numberline{}Bibliography}

\bibitem[A]{} J. Adams, {\it Lectures on Exceptional Lie Groups}, The
University of Chicago Press Ltd., London, 1996.


\bibitem[AAR]{} G. Andrews, R. Askey and R. Roy, {\it Special
Functions}, Cambridge University Press, 1999.

\bibitem[B1]{}Z. Bai, Gelfand-Kirillov dimensionas of the $\mbb
Z^2$-graded oscillator representations of $sl(n)$, {\it Acta Math.
Sinica (Engl. Ser.)} {\bf 31} (2015), no. 6, 921--937.

\bibitem[B2]{}Z. Bai, Gelfand-Kirillov dimensionas of the $\mbb
Z$-graded oscillator representations of $o(n,\mbb C)$ and
$sp(2n,\mbb C)$, {\it arXiv: 1501,04781 v1 [Math.RT].}



\bibitem[BDK]{}Barakat A., De Sole A., Kac V. G., Poisson vertex algebras in the
theory of Hamiltonian equations, {\it Jpn. J. Math.} {\bf 4} (2009),
no. 2, 141--252.





\bibitem[BO]{} R. Beerends and E. Opdam, Certain hypergeometric series related to the root
 system $BC$, {\it Trans. Amer. Math. Soc. } {\bf 119} (1993), 581--609.


\bibitem[BGG]{} I. N. Bernstein, I. M. Gel'fand and S. I. Gel'fand,
Structure of representations generated by vectors of highest weight,
(Russian) {\it Funktsional. Anal. i Prilozhen.} {\bf 8} (1971), no.
1, 1--9.

\bibitem[BH]{} F. Beukers and G. Heckman, Monodromy for a hypergeometric function $_nF_{n-1}$,
{\it Invent. Math.} {\bf 95} (1989), no. 2, 325--354.


\bibitem[B-N]{} J. Bion-Nadal, Subfactor of the hyperfinite $\Pi_1$
factor with Coxeter graph $E_6$ as invariant, {\it J. Operator
Theory} {\bf 28} (1992), 27--50.



\bibitem[BBL]{} G. Benkart, D. Britten and F. W. Lemire, Modules with
bounded multiplicities for simple Lie algebras, {\it Math. Z.} {\bf
225} (1997), 333--353.

\bibitem[Br]{} R. Block, On torsion-free abelian groups and Lie
algebras, {\it Proc. Amer. Soc.} {\bf 9} (1958), 613--620.

\bibitem[{[Bo]}]{} R. E. Borcherds, Vertex algebras, Kac-Moody
algebras, and the Monster, {\it Proc. Natl. Acad. Sci. USA} {\bf 83}
(1986), 3068--3071.




\bibitem[B]{} M. Brion, Invariants d'un sous-groupe unipotent maxaimal d'un groupe
semi-simple, {\it Ann. Inst. Fourier (Grenoble)} {\bf 33} (1983),
no. 1, 1--27.

\bibitem[BKR]{} L. Brink, S. Kim and P. Ramond, $E_{7(7)}$ on the light
cone, {\it J. High Energy Phys.} (2008), no. 6, 034, 18pp.

\bibitem[BFL]{} D. Britten, V. Futorny and F. W. Lemire, Simple
$A_2$-modules with a finite-dimensional weight space, {\it Commun.
Algebra} {\bf 23} (1995), 467--510.

\bibitem[BHL]{} D. Britten , J. Hooper and F. W. Lemire, Simple
$C_n$-modules with multiplicities 1 and applications, {\it Canad. J.
Phys.} {\bf 72} (1994), 326--335.

\bibitem[BL1]{} D. Britten and F. W. Lemire, A classification of
simple Lie modules having 1-dimensional weight space, {\it Trans.
Amer. Math. Soc. }{\bf 299} (1987), 683--697.

\bibitem[BL2]{} D. Britten and F. W. Lemire, On modules of bounded
multiplicities for symplectic algebras, {\it Trans. Amer. Math. Soc.
}{\bf 351} (1999), 3413--3431.

\bibitem[B]{} R. Brown, A minimal representation for the Lie algebra of
type $E_7$, {\it Illiois J. Math.} {\bf 12} (1968), 190--200.




\bibitem[Cf]{} F. Calogero, Solution of the one-dimensional $n$-body problem with quadratic and
 /or inversely quadratic pair potentials, {\it J. Math. Phys.} {\bf 12} (1971), 419--432.

\bibitem[Cb]{} B. Cao, Solutions of Navier Equations and Their
Representation Structure, {\it Advances in Applied Mathematics,}
\textbf{43} (2009), 331--374.

\bibitem[Cl]{} L. Chen, Twisted Hamiltonian Lie algebras and their
multiplicity-free representations, {\it Acta Math. Din. (Engl.
Ser.)} {\bf 27} (2011), no 1. 45--72.



\bibitem[DES]{} M. Davidson, T.
Enright, and R. Stanke, {\it Differential Operators and Highest
Weight Representations}, Memoirs of American Mathematical Society
{\bf 94}, no. 455, 1991.


\bibitem[D]{} L. Dickson, A class of groups in an arbitrary realm
connected with the configuration of the 27 lines on a cubic surface,
{\it J. Math.} {\bf 33} (1901), 145--123.

\bibitem[Ep]{} P. Etingof, Quantum integrable systems and representations of
Lie algebras, {\it J. Math. Phys.} {\bf 36} (1995), no.
6,2637--2651.

\bibitem[ES]{} P. Etingof and K. Styrkas, Algebraic integrability
of Schr\"{o}dinger operators and representations of Lie algebras,
{\it Compositio Math}. {\bf 98} (1995), no. 1, 91--112.

\bibitem[Eh]{} H. Exton, Multiple hypergeometric functions and applications, Halsted Pree
(John Wiley \& Sons, Inc,), New York, 1976.

\bibitem[FC]{} F. M. Fern\'{a}ndez and E. A. Castro, {\it Algebraic
Methods in Quantum Chemistry and Physics}, CRC Press, Inc., 1996.



\bibitem[Fs]{} S. L. Fernando, Lie algebra modules with
finite-dimensional weight spaces, I, {\it Trans. Amer. Math. Soc.
}{\bf 322} (1990), 757--781.





\bibitem[FSS]{} L.Frappat, A.Sciarrino and P.sorba, Dictionary on Lie
Algebras and Superalgebras, Academic Press,2000

\bibitem[Fv]{} V. Futorny, The weight representations of semisimple
finite-dimensional Lie algebras, {\it Ph.D. Thesis, Kiev University,
1987.}



\bibitem[GG]{} I. M. Gel'fand and M. I. Graev, GG-functions and
their relation to general hypergeometric functions, {\it Russian
Math. Surveys} {\bf 52}(1997), no. 4, 639--684.

\bibitem[GDi1]{} I. M. Gel'fand and L. A. Dikii, Asymptotic
behaviour of the resolvent of Sturm-Liouville equations and the
algebra of the Korteweg-de Vries equations, {\it Russian Math.
Surveys} {\bf 30:5} (1975), 77--113.

\bibitem[GDi2]{} I. M. Gel'fand and L. A. Dikii, A Lie algebra
structure in a formal variational Calculation, {\it Func. Anal.
Appl.}  {\bf 10} (1976), 16--22.

\bibitem[GDo1]{} I. M. Gel'fand and I. Ya. Dorfman, Hamiltonian
operators and algebraic structures related to them, {\it Funkts.
Anal. Prilozhen} {\bf 13} (1979), 13--30.

\bibitem[GDo2]{} I. M. Gel'fand and I. Ya. Dorfman, Schouten
brackets and Hamiltonian operators, {\it Funkts. Anal. Prilozhen}
{\bf 14} (1980),  71--74.

\bibitem[GDo3]{} I. M. Gel'fand and I. Ya. Dorfman,
 Hamiltonian operators and infinite-dimensional Lie algebras, {\it Funkts. Anal. Prilozhen}  {\bf 14} (1981),  23--40.


\bibitem[GG]{} I. M. Gel'fand and M. I. Graev, GG-functions and their relation to general
hypergeometric functions, {\it Russian Math. Surveys} {\bf
52}(1997), no. 4, 639--684.

\bibitem[GT1]{} I. Gel'fand and M. Tsetlin, Finite-dimensional
representations of the group of unimodular matrices, {\it Doklady
Akad. Nauk SSR (N.S.)} {\bf 71} (1950), 825--828.

\bibitem[GT2]{} I. Gel'fand and M. Tsetlin, Finite-dimensional
representations of the group of orthogonal matrices, {\it Doklady
Akad. Nauk SSR (N.S.)} {\bf 71} (1950), 1017--1020.




\bibitem[G]{} H. Georgi, {\it Lie Algebras in Particle Physics},
Second Edition, Perseus Books Group, 1999.




\bibitem[Hg1]{} G. Heckman, Root systems and hypergeometric functions II, {\it Compositio. Math.}
 {\bf 64} (1987), 353--373.

\bibitem[Hg2]{} G. Heckman, Heck algebras and hypergeometric functions, {\it Invent. Math.}
{\bf 100} (1990), 403--417.

\bibitem[Hg3]{} G. Heckman, An elementary approach to the hypergeometric shift operators of
 Opdam, {\it Invent. Math.} {\bf 103} (1990), 341--350.

\bibitem[HO]{} G. Heckman and E. Opdam, Root systems and hypergeometric functions I,
{\it Compositio. Math.} {\bf 64} (1987), 329--352.





\bibitem[HH]{} J. E. M. Homos and Y. M. M. Homos, Algebraic model for
the evolution of the generic code, {\it Phys. Rev. Lett.} {\bf 71}
(1991), 4401--4404 (A reference for preface).


\bibitem[Hr1]{} R. Howe, Dual pairs in physics:harmonic oscillators,
photons, electrons, and singletons, {\it Applications of group
theory in physics and mathematical physics (Chicago, 1982)},
179--207, Lectures in Appl. Mathe. {\bf 21}, {\it Amer. Math. Soc.},
Providence, RI. 1985.

\bibitem[Hr2]{} R. Howe, Remarks on classical invariant theory, {\it
Trans. Amer. Math. Soc.} {\bf 313} (1989), 539--570.

\bibitem[Hr3]{}R. Howe, Transcending classical invariant theory, {\it
J. Amer. Math. Soc.} {\bf 2} (1989), 535--552.

\bibitem[Hr4]{} R. Howe, Perspectives on invariant theory: Schur
duality, multiplicity-free actions and beyond, {\it The Schur
lectures} (1992) ({\it Tel Aviv}), 1--182, {\it Israel Math. Conf.
Proc.,} 8, {\it Bar-Ilan Univ., Ramat Gan,} 1995.



\bibitem[H]{} J. E. Humphreys, {\it Introduction to Lie Algebras and Representation Theory},
 Springer-Verlag New York Inc., 1972.

\bibitem[I1]{} N. H. Ibragimov, {\it Transformation groups applied to
mathematical physics}, Nauka, 1983.

\bibitem[I2]{} N. H. Ibragimov, {\it Lie Group Analysis of
Differential Equations}, Volume 2, CRC Handbook, CRC Press, 1995.

\bibitem[J]{} N. Jacobson, {\it Lie algebras}, Interscience
Publishers, New York/London, 1962.


\bibitem[Kv1]{} V. G. Kac, {\it Infinite-Dimensional Lie algebras}, 3rd Edition, Cambridge
University Press, 1990.

\bibitem[Kv2]{} V. G. Kac,{\it Vertex algebras for beginners}, University lecture series, Vol {\bf 10}, AMS, Providendence RI, 1996.


\bibitem[Kb]{} B. Kostant, On the tensor product of a finite and an infinite
dimensional representation, {\it J. Fund. Anal.} {\bf 20} (1975),
257--285.


\bibitem[KK]{} I.  Krasil'shchik and P. Kersten. {\it Symmetries and
Recursion Operators for Classical Supersymmetric Differential
Equations}, Mathematics and Its Applications, Kluwer Academic
Publishers, 2000.


\bibitem[Lj]{} J. Lepowsky, Generalized Verma modules, the Cartan-Helgason theorem, and the Harish-Chandra homomorphism, {\it J. Algebra} {\bf 49} (1977), no. 2, 470--495.


\bibitem[Lie]{} S. Lie, Begr\"{u}dung einer Invariantentheorie der
Ber\"{u}hrungstransformationen, {\it Math. Ann.} {\bf 8} (1874),
215--288.

\bibitem[LF]{} W. Ludwig and C. Falter, {\it
Symmetries in Physics}, Second Edition, Springer-Verlag,
Berlin/Heidelberg, 1996.

\bibitem[LX1]{}C. Luo and X. Xu, $\mbb{Z}^2$-Graded oscillator generalizations of $sl(n)$,
{\it Commun.  Algebra} {\bf 41} (2013), 3147--3173.

\bibitem[LX2]{}C. Luo and X. Xu, Z-Graded oscillator generalizations of the
classical theorem on harmonic polynomials, {\it J. Lie Theory} {\bf
23} (2013), 979--1003.

\bibitem[LX3]{} C. Luo and X. Xu, $\mbb{Z}$-Graded  oscillator
 representations  of Symplectic Lie algebras, {\it J. Algebra} {\bf 403}
 (2014), 401--425.

\bibitem[LX4]{} C. Luo and X. Xu, Supersymmetric analogues of the classical
theorem on harmonic polynomials, {\it J. Algebra Appl.} {\bf 13}
(2014) no. 6, 1450011, 42pp.

\bibitem[Ll1]{} L. Luo, Cohomology of oriented tree diagram Lie
algebras, {\it Comm. Algebra} {\bf 37} (2009), no. 3, 965--984.

\bibitem[Ll2]{} L. Luo, Oriented tree diagram Lie
algebras and their abelian ideals, {\it Acta Math. Sin. (Engl.
Ser.)} {\bf 20} (2010), no. 11, 2041--2058.


\bibitem[M]{} O. Mathieu, Classification of irreducible weight
modules, {\it Ann. Inst. Fourier (Grenoble)} {\bf 50} (2000),
537--592.



\bibitem[M]{} M. Marvan, Reducibility of zero curvature
representations with application to recursion operators, {\it Acta
Appl. Math.} {\bf 83} (2004), 39--68.



\bibitem[Ma]{} A. Molev, A basis for representations of symplectic
Lie algebras, {\it Comm. Math. Phys.} {\bf 201} (1999), 591--618.


\bibitem[OP]{} M. Olslanetsky and A. Perelomov, Completely integrable Hamiltonian systems
connected with semisimple Lie algebras, {\it Invent. Math.} {\bf 37}
(1976), 93--108.

\bibitem[Op]{} P. J. Olver, {\it Equivalence, Invariants and
Symmetry}, Cambridge University Press, Cambridge, 1995.


\bibitem[O1]{}  E. Opdam, Root systems and hypergeometric functions III, {\it Compositio. Math.}
 {\bf 67} (1988), 21--49.

\bibitem[O2]{}  E. Opdam, Root systems and hypergeometric functions IVI, {\it Compositio. Math.}
 {\bf 67} (1988), 191--209.

\bibitem[O3]{} E. Opdam, Some applications of hypergeometric shift operators,
{\it Invent.  Math.} {\bf 98} (1989), 1--18.

\bibitem[O4]{}  E. Opdam, An analogue of the Gauss summation formula for hypergeometric
 functions related to root systems, {\it Math. Z.} {\bf 212} (1993), 313--336.

\bibitem[O5]{}  E. Opdam, Harmonic analysis for certain representations of graded Hecke algebras,
{\it Acta Math.} {\bf 175} (1995), 75--121.

\bibitem[O6]{}  E. Opdam, Cuspital hypergeometric functions, {\it Methods Appl. Anal.} {\bf 6}
(1999), 67--80.



\bibitem[P4]{} V. Pless, {\it Introduction to the Theory of Error-Correcting Codes},
 3rd Edition, John Wiley \& Sons, Inc., 1998.

\bibitem[PHB]{} V. Pless, W. C. Huffman and R. A. Brualdi, An introduction to algebraic codes, in: {\it Handbook
of Coding Theory} ed. by V. Pless and W. Huffman, Chap.1, Elsevier
Science B.V., 1998.


\bibitem[Sg]{} G. Shen, Graded modules of graded Lie algebras of
Cartan type (I)---mixed product of modules, {\it Science in China A}
{\bf  29} (1986), 570--581.


\bibitem[Sb]{} B. Sutherland, Exact results for a quantum many-body problem in
one-dimension, {\it Phys. Rev. A} {\bf 5} (1972), 1372--1376.

\bibitem[Vj]{}  J. H. van Lint, {\it Introduction to Coding Theory}, 3rd Edition, Springer-Verlag, Berlin Heidelberg, 1999.

\bibitem[V1]{} D,-N. Verma, Structure of certain induced representations of complex semisimple
 Lie algebras, {\it thesis, Yale University,} 1966.

\bibitem[V2]{} D,-N. Verma, Structure of certain induced representations of complex semisimple
 Lie algebras, {\it Bull. Amer. Math. Soc.} {\bf 74}(1968), 160--166.



\bibitem[WG]{} Z. Wang and D. Guo, {\it Special functions}, World Scientific,
Singapore, 1998.


\bibitem[X1]{} X. Xu, A gluing technique for constructing relatively self-dual codes,
 {\it J. Comb. Theory, Ser. A} {\bf 66} (1994), 137--159.

\bibitem[X2]{} X. Xu, Self-dual lattices of type $A$,  {\it Acta
Math.} {\bf 175} (1995), 123--150.

\bibitem[X3]{} X. Xu, Hamiltonian superoperators, {\it J. Phys. A} {\bf 28} (1995),
1681-1698.

\bibitem[X4]{} X. Xu, Hamiltonian operators and associative algebras with a
derivation, {\it Lett. Math. Phys.} {\bf 33} (1995), 1-6.

\bibitem[X5]{}X. Xu, Differential invariants of classical groups, {\it Duke Math. J.} {\bf 94} (1998), 543-572.

\bibitem[X6]{} X. Xu,  {\it Introduction to Vertex Operator
Superalgebras and Their Modules},
 Kluwer Academic Publishers, Dordrecht/Boston/London, 1998.

\bibitem[X7]{}X. Xu, Generalizations of the Block algebras, {\it Manuscripta
Mathematica} {\bf 100} (1999), 489--518.

\bibitem[X8]{}X. Xu, New generalized simple Lie algebras of Cartan type over a
field with characteristic 0, {\it J. Algebra} {\bf 224} (2000),
23--58.

\bibitem[X9]{}X. Xu, Variational calculus of supervariables and related algebraic structures,
{\it J. Algebra} {\bf 223} (2000), 386-437.

\bibitem[X10]{}X. Xu, Equivalence of conformal superalgebras to Hamiltonian superoperators, {\it Algebra Colloq.}
{\bf 9} (2001), 63-92.

\bibitem[X11]{} X. Xu, Poisson and Hamiltonian superpairs over polarized associative
algebras, {\it J. Phys. A} {\bf 34} (2001), 4241-4265.

\bibitem[X12]{}X. Xu, Differential equations for singular vectors of sp(2n),
 {\it Commun. Algebra} {\bf 33} (2005), 4177--4196.

\bibitem[X13]{}X. Xu, Tree diagram Lie algebras of  differential operators and
evolution partial differential equations, {\it J. Lie Theory} {\bf
16} (2006),  691--718.

\bibitem[X14]{}X. Xu, Path hypergeometric functions,  {\it J. Algebra
Appl.} {\bf 6} (2007),  595--653.

\bibitem[X15]{}X. Xu, {\it Kac-Moody Algebras and Their
Representations}, China Science Press, 2007.

\bibitem[X16]{} X.Xu, Flag partial differential equations and
representations of Lie algebras, {\it Acta Appl. Math.} {\bf 102}
(2008), 249--280.

\bibitem[X17]{} X. Xu, A cubic $E_6$-generalization of the classical theorem on harmonic polynomials,
{\it J. Lie Theory} 21 (2011), 145--164.

\bibitem[X18]{} X. Xu, Partial differential equation approach to $F_4$,
{\it Front. Math. China} {\bf 6} (2011),  759--774.

\bibitem[X19]{} X.Xu, Differential-operator representation of $S_n$ and
singular vectors in Verma modules, {\it Algebr. Represent. Theor.}
{\bf 15} (2012), 211--231.

\bibitem[X20]{}X. Xu, Representations of Lie algebras and coding theory,  {\it
J. Lie Theory}  {\bf 22} (2012), no. 3, 647--682.


\bibitem[X21]{}X. Xu, {\it Algebraic Approaches to Partial
Differential Equations}, Springer, Heidelberg/New York/Dordrecht/
London, 2013.

\bibitem[X22]{}X. Xu, Partial Differential equation approach to
$E_7$, {\it Acta Math Sin. (Engl. Ser.)} {\bf 31} (2015), 177--200.

\bibitem[X23]{}X. Xu, A new functor from $D_5$-{\bf Mod} to $E_6$-{\bf
Mod}, {\it arXiv:1112.3792v1[math.RT]}.


\bibitem[X24]{}X. Xu,  A new functor from $E_6$-{\bf Mod} to $E_7$-{\bf
Mod}, {\it Commun. Algebra} {\bf 43} (2015), no. 9, 3589-3636.

\bibitem[X25]{}X. Xu, Projective oscillator representations of $sl(n+1)$ and
$sp(2m+2)$, {\it J. Lie Theory} {\bf 26} (2016),  96--114.

\bibitem[X26]{}X. Xu,  Conformal oscillator representations of orthogonal Lie
algebras, {\it Sci. China: Math. (Engl. Ser.)} {\bf 59}, No.1,
37-48.

\bibitem[X27]{}X. Xu, Representations of $E_6$, combinatorics  and algebraic
Varieties, {\it preprint}.

\bibitem[X28]{}X. Xu, Representations of $E_7$, combinatorics  and algebraic
Varieties, {\it preprint}.


\bibitem[XZ]{} X. Xu and Y. Zhao, Extensions of the conformal
representations for orthogonal Lie algebras, {\it J. Algebra} {\bf
377} (2013), 97--124.

\bibitem[XW1]{}W. Xiao, Differential equations and singular vectors
in Verma modules over $sl(n,\mbb C)$, {\it Acta Math. Sinica, (Engl.
Ser.)} {\bf 31} (2015), no. 7, 1057--1066.

\bibitem[XW2]{}W. Xiao, Differential-operator representations of Weyl
groups and singular vectors in Verma modules, {\it Preprint}.




\bibitem[Z]{} R. B. Zhang, Orthosymplectic Lie superalgebras in
superspace analogues of quantum Kepler problems, {\it Commun. Math.
Phys.} {\bf 280} (2008), 545--562.

\bibitem[ZX]{} Y. Zhao and X. Xu, Generalized  projective  representations
 for sl(n+1), {\it J. Algebra} {\bf 328} (2011), 132--154.




\end{thebibliography}
\end{document}